\documentclass[10pt]{book}
\title{Compressible Flow and Euler's Equations\footnote{Work supported by ERC Advanced Grant 246574 \emph{Partial Differential Equations of Classical Physics}. }}
\usepackage{amsmath,bbold}
\usepackage{geometry}
\usepackage{framed}
\usepackage{amsfonts,amssymb,amsmath,verbatim}
\usepackage{mathrsfs}
\usepackage{amssymb}
\usepackage{amsmath}
\usepackage{amsfonts}
\usepackage{mathrsfs}
\usepackage{amsfonts, amssymb, amsmath}
\usepackage{a4,graphicx}
\usepackage{graphicx}
\usepackage{slashed}
\usepackage{amsthm}
\usepackage{authblk}
\newcommand{\leftexp}[2]{{\vphantom{#2}}^{#1}{#2}}
\author[1]{Demetrios Christodoulou}
\author[1,2]{Shuang Miao}
\affil[1]{Department of Mathematics, ETH Zurich}
\affil[2]{Academy of Mathematics and Systems Sciences,
CAS}
\begin{document}
\maketitle
\chapter*{Introduction}
The equations describing the motion of a perfect fluid were first formulated by Euler
in 1752 (see \cite{Eu1}, \cite{Eu2}), based, in part, on the earlier work of D.Bernoulli \cite{Be}. These equations were among the first partial differential equations to be written down, preceded, it seems, only by D'Alembert's 1749 formulation \cite{DA} of the one-dimensional wave equation describing the motion of a vibrating string in the linear approximation. In contrast to D'Alembert's equation however, we are still, after the lapse of two and a half centuries, far from having achieved an adequate understanding 
of the observed phenomena which are supposed to lie within the domain of validity of Euler's equations.

The phenomena displayed in the interior of a fluid fall into two broad classes, the phenomena of sound, the linear theory of which is acoustics, and the phenomena of vortex motion. The sound phenomena depend on the compressibility of a fluid, while the vortex phenomena occur even in a regime where the fluid may be considered to be incompressible. The formation of shocks, the subject of the present monograph, belongs to the class of sound phenomena, but lies in nonlinear regime, beyond the range covered by the linear acoustics.

Let us make a short review of the history of the study of the sound phenomena in fluids, in particular the phenomena of the formation of shocks in the nonlinear regime. At the time when the equations of the fluid motion were formulated, thermodynamics was in its infancy, however, it was already clear that the local state of a fluid as seen by a comoving observer is determined by two thermodynamic variables, say pressure and temperature. Of these, only pressure entered the equations of motion, while the equations involve also the density of the fluid. Density was already known to be a function of pressure and temperature for a given type of fluid. However, in the absence of an additional equation, the system of equations at the time of Euler, which consisted of the momentum equations and the equation of continuity, was underdetermined, except in the incompressible limit. The additional equation was supplied by Laplace in 1816 \cite{La} in the form of what was later to be called adiabatic condition, and allowed him to make the first correct calculation of the sound speed.

The first work on the formation of shocks was done by Riemann in 1858 \cite{Ri}. Riemann considered the case of isentropic flow with plane symmetry, where the equations of fluid mechanics reduces to a system of conservation laws for two unknowns and with two independent variables, a single space coordinate and time. He introduced for such systems the so-called Riemann invariants, and with the help of these showed that solutions which arise from smooth initial conditions develop infinite gradients in finite time.

In 1865 the concept of entropy was introduced into theoretical physics by Clausius \cite{Cl2}, and the adiabatic condition was understood to be the requirement that the entropy of each fluid element remains constant during its evolution.

The first general result on the formation of singularity in three dimensional fluids was obtained by Sideris in 1985 \cite{S}. Sideris considered the compressible Euler equations in the case of a classical ideal gas with adiabatic index $\gamma>1$ and with initial data which coincide with those of a constant state outside a ball. The assumptions of his theorem on the initial data were that there is an annular region bounded by the sphere outside which the constant state holds, and a concentric sphere in its interior, such that a certain integral in this annular region of $\rho-\rho_{0}$, the departure of the density $\rho$ from its value $\rho_{0}$ in the constant state, is positive, while another integral in the same region of $\rho v^{r}$, the radial momentum density, is non-negative. These integrals involve kernels which are functions of the distance from the center. It is also assumed that everywhere in the annular region the specific entropy $s$ is not less than its value $s_{0}$ in the constant state. The conclusion of the theorem is that the maximal time interval of existence of a smooth solution is finite. The chief drawback of this theorem is that it tells us nothing about the nature of the breakdown. Also the method relies the strict convexity of the pressure as a function of density displayed by the equation of state of an ideal gas, and does not extend to more general equation of state.

The most recent and complete results on the formation of shocks in three dimensional fluids were obtained by Christodoulou in 2007 \cite{Chr}. Christodoulou considered the relativistic Euler equations in three space dimensions for a perfect fluid with an arbitrary equation of state. He considered the regular initial data on a space-like hyperplane $\Sigma_{0}$ in Minkowski spacetime which outside a sphere coincide with the data corresponding to a constant state. He considered the restriction of the initial data to the exterior of a concentric sphere in $\Sigma_{0}$ and the maximal classical development of this data. Under suitable restriction on the size of the departure of the initial data from those of the constant state, he proved certain theorems which give a complete description of the 
maximal classical development. In particular, theorems give a detailed description of the geometry of the boundary of the domain of the maximal classical solution and a detailed analysis of the behavior of the solution at this boundary.

The aim of the present monograph is to derive analogous results for the classical, non-relativistic, compressible Euler's equations taking the data to be irrotational and isentropic, and to give a proof of these results which is considerably simpler and completely self-contained.
The present monograph in fact not only gives simpler proofs but also sharpens some of the results. In addition the present monograph explains in depth the ideas on which the approach is based. Finally certain geometric aspects which pertain only to the non-relativistic theory are discussed.

We shall presently explain the basis of the approach of the present monograph (also of the previous one dealing with the relativistic case). This basic idea can be thought as an extension of the method of Riemann invariants combined with the method of the partial hodograph transformation, to the case of more than one space dimension. We first recall some basic facts about Riemann invariants. In the case of one space dimension, the isentropic Euler system reads:
\begin{align*}
\partial_{t}\rho+\partial_{x}(\rho v)=0\\
\partial_{t}v+v\partial_{x}v=-\frac{1}{\rho}\partial_{x}p
\end{align*}
and it can be written as a single equation of the velocity potential $\phi$:
\begin{align*}
(g^{-1})^{\mu\nu}\partial_{\mu}\partial_{\nu}\phi=0
\end{align*}
or, setting $\psi_{\mu}=\partial_{\mu}\phi$,
\begin{align*}
(g^{-1})^{\mu\nu}\partial_{\mu}\psi_{\nu}=0
\end{align*}
where $g$ is the following Lorentzian metric on the spacetime manifold $\mathcal{M}$:
\begin{align*}
g=-\eta^{2}dt^{2}+(dx-vdt)^{2},
\end{align*}
Here 
\begin{align*}
v=-\psi_{x}
\end{align*}
is the fluid velocity and 
\begin{align*}
h=\psi_{t}-\frac{1}{2}\psi_{x}^{2}
\end{align*}
is the enthalpy. The pressure $p$ is a given function of $h$ and 
\begin{align*}
\rho=dp/dh
\end{align*}
while the sound speed $\eta$ is given by
\begin{align*}
\eta^{2}=dp/d\rho
\end{align*}

Riemann invariants are the functions $R_{1}, R_{2}$ defined on the cotangent space $T^{*}_{p}\mathcal{M}$, which are the two functionally independent solutions of the following eikonal equation:
\begin{align*}
g_{\mu\nu}\frac{\partial R}{\partial\psi_{\mu}}\frac{\partial R}{\partial \psi_{\nu}}=0
\end{align*}
i.e.
\begin{align*}
-\eta^{2}(\frac{\partial R}{\partial\psi_{t}})^{2}+(\frac{\partial R}{\partial\psi_{x}}+\psi_{x}\frac{\partial R}{\partial\psi_{t}})^{2}=0
\end{align*}
We define the vectorfields $N_{1}$ $N_{2}$ on $\mathcal{M}$ by:
\begin{align*}
N_{1}=\frac{\partial R_{1}}{\partial\psi_{\mu}}\frac{\partial}{\partial x^{\mu}},\quad N_{2}=\frac{\partial R_{2}}{\partial\psi_{\mu}}\frac{\partial}{\partial x^{\mu}}
\end{align*}
These are null vectorfields with respect to $g$. We choose $R_{1}$ and $R_{2}$ so that the integral curves of $N_{1}$ and $N_{2}$ are the incoming and outgoing null curves respectively. Then $R_{1}$, $R_{2}$ as functions on $\mathcal{M}$ satisfy:
\begin{align}
N_{2}R_{1}=0,\quad N_{1}R_{2}=0
\end{align}
Let us introduce the acoustical coordinates $(t,u)$ so that $u$ is constant along the outgoing null curves. Then the vectorfields:
\begin{align*}
L=\frac{\partial}{\partial t},\quad \underline{L}=\eta^{-1}\kappa L+2T
\end{align*}
where 
\begin{align}
T=\frac{\partial}{\partial u}\quad\textrm{and}\quad \kappa=-\frac{\partial x}{\partial u}
\end{align}
are null vectorfields with respect to $g$, and the integral curves of $L$ and $\underline{L}$ are outgoing and incoming null curves respectively. Therefore $L$ and $\underline{L}$ are colinear to $N_{2}$ and $N_{1}$ respectively. Therefore equations (1) are equivalent to:  
\begin{align}
LR_{1}=0,\quad \underline{L}R_{2}=0
\end{align}
To write down explicit expressions for $R_{1}, R_{2}$ we use, instead of $(\psi_{t},\psi_{x})$, the following variables in $T^{*}_{p}\mathcal{M}$:
\begin{align*}
h=\psi_{t}-\frac{1}{2}\psi_{x}^{2},\quad v=-\psi_{x}
\end{align*}
Then for any function $f=f(h,v)$ defined on $T^{*}_{p}\mathcal{M}$, we have:
\begin{align*}
\frac{\partial f}{\partial\psi_{x}}+\psi_{x}\frac{\partial f}{\partial\psi_{t}}=-\frac{\partial f}{\partial v}
\end{align*}
Let us introduce a function $r=r(h)$ by:
\begin{align*}
\frac{dr}{dh}=\frac{1}{\eta},\quad r(0)=0
\end{align*}
Then 
\begin{align*}
R_{1}=r+v,\quad R_{2}=r-v
\end{align*}
are the two functionally independent solutions of the eikonal equation, therefore the two Riemann invariants. From the first equation of (2), we know that:
\begin{align*}
R_{1}=R_{1}(u)
\end{align*}
is determined by initial data, while to obtain $R_{2}$, we consider the second 
of equations (3), namely the equation:
\begin{align*}
\eta^{-1}\kappa\frac{\partial R_{2}}{\partial t}+2\frac{\partial R_{2}}{\partial u}=0
\end{align*}
Here $\kappa$ enters, which is defined by (2). To obtain an equation for $\kappa$ we consider the equation:
\begin{align*}
\frac{\partial x}{\partial t}=c_{+},\quad c_{+}=v+\eta
\end{align*}
We shall derive an equation for $\kappa$. From (3) we have:
\begin{align*}
\frac{\partial\kappa}{\partial t}=-\frac{\partial c_{+}}{\partial u}
=-\frac{1}{2}\frac{\partial R_{1}}{\partial u}+\frac{1}{2}\frac{\partial R_{2}}{\partial u}-\frac{d\eta}{dh}\frac{\partial h}{\partial u}
\end{align*}
While
\begin{align*}
\frac{\partial h}{\partial u}=\frac{dh}{dr}\frac{\partial r}{\partial u}=
\frac{1}{2}\eta(\frac{\partial R_{1}}{\partial u}
+\frac{\partial R_{2}}{\partial u})
\end{align*}
Substituting the above we obtain:
\begin{align*}
\frac{\partial\kappa}{\partial t}=\frac{1}{2}(-1-\eta\frac{d\eta}{dh})\frac{\partial R_{1}}{\partial u}+\frac{1}{2}(1-\eta\frac{d\eta}{dh})\frac{\partial R_{2}}{\partial u}
\end{align*}
Since 
\begin{align*}
\frac{\partial R_{1}}{\partial u}=\frac{2}{\eta}\frac{\partial h}{\partial u}-\frac{\partial R_{2}}{\partial u}
\end{align*}
we have:
\begin{align*}
\frac{\partial\kappa}{\partial t}=\frac{1}{2\eta}(-2-2\eta\frac{d\eta}{dh})\frac{\partial h}{\partial u}+\frac{\partial R_{2}}{\partial u}
\end{align*}
Let us define:
\begin{align*}
H=-2h-\eta^{2}
\end{align*}
Noting that by the second of (3)
\begin{align*}
\frac{\partial R_{2}}{\partial u}=-\frac{\kappa}{2\eta}\frac{\partial R_{2}}{\partial t}
\end{align*}
We conclude that $\kappa$ satisfies the following equation:
\begin{align*}
\frac{\partial\kappa}{\partial t}=m^{\prime}+\kappa e^{\prime}
\end{align*}
with
\begin{align*}
m^{\prime}=\frac{1}{2\eta}\frac{dH}{dh}\frac{\partial h}{\partial u},\quad e^{\prime}=-\frac{1}{2\eta}\frac{\partial R_{2}}{\partial t}
\end{align*}

The main idea in the one dimensional case is that $R_{1}, R_{2}$ as well as the rectangular coordinate $x$ are smooth functions of $(t,u)$. The partial hodograph transformation is the transformation:
\begin{align*}
(t,u)\longmapsto(t,x),
\end{align*}
from acoustical to rectangular coordinates. The Jacobian is:
\begin{align*}
\frac{\partial(t,x)}{\partial(t,u)}=\begin{vmatrix}
1&0\\v+\eta&-\kappa
\end{vmatrix}
=-\kappa,
\end{align*}
and vanishes when $\kappa$ vanishes. This means $R_{1}, R_{2}$ are not smooth in $(t,x)$ when shocks form.

In the case of more than one space dimension, in particular, the case of three space dimensions, we do not have Riemann invariants. We work instead with the first order variations, which are defined through the \emph{variation fields}:
\begin{align*}
\partial_{\mu},\quad \mathring{R}_{i}=\epsilon_{ijk}x^{j}\frac{\partial}{\partial x^{k}},\quad x^{\mu}\partial_{\mu}-I,\quad \mu=0,1,2,3;\quad 1\leq i<j\leq 3
\end{align*} 
Here $I$ is the multiplication operator by $1$.
These fields are the generators the of the subgroup of the scale-extended Galilean group, the invariance group of the compressible Euler system, which leaves the constant state invariant. These first order variations satisfy the linear wave equation:
\begin{align}
\Box_{\tilde{g}}\psi=0
\end{align}
where $\tilde{g}$ is the conformal acoustical metric:
\begin{align*}
\tilde{g}=\Omega g,\quad \Omega=\frac{\rho}{\eta},\quad g=-\eta^{2}dt^{2}
+\sum_{i}(dx^{i}-v^{i}dt)^{2}
\end{align*}
Like the equations satisfied by $R_{1}, R_{2}$ in the case of one space dimension, equation (4) does not depend on the Galilean structure, but depends only on the properties of $(\mathcal{M},g)$ as a Lorentzian manifold. Actually, it depends more sensitively on the conformal class of $(\mathcal{M},g)$. This is similar to the fact that in the case of one space dimension, null curves depend only on the conformal class of the acoustical metric. 
 Also like in the case of one space dimension, we shall work in the \emph{acoustical coordinates} $(t,u,\vartheta)$. Here $u$ is the \emph{acoustical function} in $\mathcal{M}$, whose level sets $C_{u}$ are outgoing null hypersurfaces. The level curves of $\vartheta\in\mathbb{S}^{2}$ on each $C_{u}$ are the generators of $C_{u}$. Like in the case of one space dimension, we denote by $L$ the tangent vectorfield of the generator of $C_{u}$ parametrized by $t$, and $\underline{L}$ the incoming null normals of $S_{t,u}:=C_{u}\bigcap\Sigma_{t}$, whose definitions are formally the same as one space dimensional case. We also denote by $T$ the tangent vectorfield of the inward normal curves to the $S_{t,u}$ in $\Sigma_{t}$ parametrized by $u$. Here $\Sigma_{t}$ is the level set of the function $t$, which is isometric to the Euclidean space.

To obtain a fundamental energy estimate for this linear equation in $(t,u,\vartheta)$, we need two \emph{multiplier vectorfields} $K_{0}, K_{1}$. These are non-spacelike and future-directed with respect to $g$ (a requirement which actually depends only on the conformal class of $g$). They are linear combinations of $L$ and $\underline{L}$, with coefficients which are smooth in $(t,u,\vartheta)$.
The concept of multiplier vectorfields originates from 
Noether's theorem \cite{N} on conserved currents. A modern more general treatment of \emph{compatible currents} is found in \cite{Chr1}.
In order to obtain higher order energy estimates, we consider the $n$th order variations by applying a string of \emph{commutation vectorfields} of length $n-1$ to the first order variations:
\begin{align*}
T, \quad Q:=(1+t)L,\quad R_{i}=\Pi \mathring{R}_{i}
\end{align*}
Here $\Pi$ is the orthogonal projection from $T_{p}\Sigma_{t}$ to $T_{p}S_{t,u}$. Then we obtain an inhomogeneous wave equation for the $n$th order variation $\psi_{n}$:
\begin{align}
\Box_{\tilde{g}}\psi_{n}=\rho_{n}
\end{align}
where $\rho_{n}$ is determined by the deformation tensors of the commutation vectorfields. Actually $\rho_{n}$ depends on up to the $n-1$th order derivatives of the deformation tensors, and $\rho_{1}=0$. The use of commutation fields originates in \cite{Kl}. Multiplier fields and commutation fields on general curved spacetimes have first been used in \cite{CK}.

After we solve (5) in the acoustical coordinates, we need to go back to the original rectangular coordinates. Again we must consider the inverse of the transformation:
\begin{align*}
(t,u,\vartheta)\longmapsto(t,x)
\end{align*}
where $\vartheta\in\mathbb{S}^{2}$ and $x\in\mathbb{R}^{3}$. This is what replaces the partial hodograph transformation in higher dimensions. The Jacobian of this transformation is:
\begin{align}
\kappa\sqrt{\det\slashed{g}}
\end{align}
Here $\slashed{g}$ is the induced acoustical metric on $S_{t,u}$. So we consider the system satisfied by the rectangular coordinates on each $C_{u}$:
\begin{align*}
\frac{\partial x^{i}}{\partial t}=L^{i}=-\eta\hat{T}^{i}-\psi_{i},
\end{align*}
 which is a fully nonlinear system for $x^{i}$. Here $\hat{T}$ is the inward unit normal of $S_{t,u}$ in $\Sigma_{t}$, whose expression is the ratio of a homogeneous quadratic polynomial in $\frac{\partial x^{i}}{\partial\vartheta^{A}}$ to the square root of a homogeneous quartic polynomial in $\frac{\partial x^{i}}{\partial\vartheta^{A}}$. The estimates of the derivatives of $x^{i}$ reduce to the estimates of the derivatives of $\chi$ and $\mu$. These are defined as follows:
 \begin{align*}
 2\chi=\slashed{\mathcal{L}}_{L}\slashed{g},\quad \mu=\eta\kappa
\end{align*}  
where $\kappa$ is the magnitude of $T$. Thus $\chi$ is the 2nd fundamental form of $S_{t,u}$ in $C_{u}$. Finally the way we estimate $\chi, \mu$ is to study the geometric structure equations of the foliation of $\mathcal{M}$ by surfaces $S_{t,u}$.

To summarize, in the case of $n$ space dimensions, the role of Riemann invariants is played by the first order variations $\psi$, 
which shall be proved to be smooth functions of $(t,u,\vartheta)$.
Moreover, we shall show that the $x^{i}$ are also smooth functions of $(t,u,\vartheta)$. This shall be done through estimates on $\chi$ and $\mu$ based on the geometric structure equations. One of these equations is the same as in the one space dimensional case:
\begin{align*}
L\mu=m+\mu e
\end{align*}
where
\begin{align*}
m=\frac{1}{2}\frac{dH}{dh}Th,\quad H=-2h-\eta^{2}
\end{align*}
and
\begin{align*}
e=\frac{1}{2\eta^{2}}(\frac{\rho}{\rho^{\prime}})^{\prime}Lh+\frac{1}{\eta}\hat{T}^{i}(L\psi_{i})
\end{align*}

As a consequence of these facts, the boundary of the maximal classical development contains a singular part $\mathcal{H}$, where the Jacobian (6) vanishes. Since $\sqrt{\det\slashed{g}}$ is, by virtue of the estimates for $\chi$, bounded from below by a positive constant, $\kappa$, equivalently $\mu$, vanishes on $\mathcal{H}$. Thus, the inverse of the transformation:
\begin{align*}
(t,u,\vartheta)\longmapsto(t,x)
\end{align*}
is not differentiable at $\mathcal{H}$. Therefore the $\psi$ are not differentiable with respect to the rectangular coordinates at $\mathcal{H}$. Nevertheless, $\mathcal{H}$, the zero-level set of $\mu$, a smooth function of $(t,u,\vartheta)$, is a smooth hypersurface in $\mathcal{M}$ relative to the differential structure induced by the acoustical coordinates. This is because we can show that $L\mu$ is bounded from above by a strictly negative function at $\mathcal{H}$, therefore $\mathcal{H}$ is a non-critical level set of $\mu$.

The first main result in the present monograph can be thought as an existence theorem \\($\textbf{Theorem 17.1}$) for the nonlinear wave equation of the velocity potential:
\begin{align*}
(g^{-1})^{\mu\nu}\partial_{\mu}\partial_{\nu}\phi=0
\end{align*}
with small initial data: 

\emph{The solution of the above nonlinear wave equation can be extended smoothly to the boundary of the maximal solution in acoustical coordinates $(t,u,\vartheta)$, and the solution is also smooth in rectangular coordinates before $\mu$ becomes $0$, since the differential structures induced by acoustical coordinates and rectangular coordinates are equivalent to each other when $\mu>0$.}

 Theorem 17.1 also gives a lower bound for the time when the shock forms (i.e. $\mu=0$), and the energy estimates for solution as well as various geometric quantities associated to the acoustical spacetime $(\mathcal{M}, g)$.  

Based on this existence theorem, we find some conditions on initial data which guarantee the formation of shocks in finite time (See $\textbf{Theorem 18.1}$). The conditions are imposed on a $\Sigma_{t}$-integral of the following function:
\begin{align}
(1-u+t)\underline{L}\psi_{0}-\psi_{0}
\end{align}
where
\begin{align*}
\psi_{0}:=\partial_{t}\phi
\end{align*}
is one of the first order variations.
The principal part of function (7), namely, $(1-u+t)\underline{L}\psi_{0}$ determines the properties of function $m$, which, in turn, determines the formation of shocks. Moreover, the spherical mean of function (7) on $S_{t,u}$ satisfies an ordinary differential inequality in the parameter $t$. Then we can connect the properties of $m$ near the point where the shocks form and the properties of $m$ on the initial hypersurface $\Sigma_{0}$ by using this ordinary differential inequality. Then the necessary properties of $m$ then follow from its properties on the initial hypersurface $\Sigma_0$.

Also based on the existence theorem, we can give a geometric description of the boundary of the maximal classical solution in acoustical differential structure ($\textbf{Proposition 19.1}$):

\emph{The boundary contains a regular part $\underline{C}$, which is an incoming null hypersurface in $(\mathcal{M},g)$, and a singular part $\mathcal{H}$, on which the function $\mu$ vanishes. $\mathcal{H}$ is a spacelike hypersurface in $(\mathcal{M},g)$, and it has the common past boundary with $\underline{C}$, denoted by $\partial_{-}\mathcal{H}$, which is a 2-dimensional spacelike surface in $(\mathcal{M},g)$. However, the singular boundary, from the intrinsic point of view, is a null hypersurface in $(\mathcal{M},g)$, on which the acoustical metric $g$ degenerates in acoustical coordinates.} 

The corresponding 
description of the singular boundary in the standard differential structure (that is, in rectangular coordinates) is given in $\textbf{Proposition 19.3}$. 

Moreover, we establish a trichotomy theorem ($\textbf{Theorem 19.1}$) describing the behavior of the past-directed null geodesics initiated at the singular boundary:

 \emph{For each point $q$ of the singular boundary, the intersection of the past null geodesic conoid of $q$ with any $\Sigma_{t}$ in the past of $q$ splits into three parts, the parts corresponding to the outgoing and to the incoming sets of null geodesics ending at $q$ being embedded discs with a common boundary, an embedded circle, which corresponds to the set of the remaining null geodesics ending at $q$. All outgoing null geodesics ending at $q$ have the same tangent vector at $q$.}

Finally, considering the transformation from one acoustical function to another, we show that the foliations corresponding to different families of outgoing null hypersurfaces have equivalent geometric properties and degenerate in precisely the same way on the same singular boundary (See $\textbf{Proposition 19.2}$).  

Let us now give an outline of the present monograph. The first four chapters concern the geometric set up. Then in Chapter 5, we obtain energy estimates for the linear wave equation associated to the conformal acoustical metric. Chapter 6 deals with the preliminary estimates for the deformation tensor of the commutation vectorfields, the precise estimates of which are given in Chapter 10 and 11. We also introduce the basic bootstrap assumptions on variations as well as on $\chi$ and $\mu$ in Chapter 5 and Chapter 6. Chapter 8 and Chapter 9 are crucial in the whole work, because it is here that estimates for $\chi$ and $\mu$ are
derived which do not lose derivatives, thus allowing us to close the
bootstrap argument. Chapter 8 concerns the estimates for the top order spatial derivatives of $\chi$. In fact only the top order angular derivatives are involved. 
While Chapter 9 concerns the estimates for the top order spatial derivatives of $\mu$. In Chapter 12, based on the bootstrap assumptions on variations, we recover the bootstrap assumptions on $\chi$ and $\mu$, except $\textbf{C1}$, $\textbf{C2}$ and $\textbf{C3}$, which are recovered in Chapter 13. In Chapter 14, based on a crucial lemma ($\textbf{Lemma 8.11}$) established in Chapter 8, we estimate the \emph{borderline} contribution from the top order spatial derivatives of $\chi$ and $\mu$. Then in Chapter 15, we obtain the energy estimates for the top order variations. These are allowed to blow up as shocks begin to form. We then revisit the lower order energy estimates and show that the estimates of each preceding order blow up successively more slowly
until we finally reach energy estimates of a certain order which do not
blow up at all. These bounded energy estimates allow us to close the bootstrap argument (See Chapter 16-17). 

In regard to the notational conventions, Latin indices take the values $1,2,3$, while Greek indices take the values $0,1,2,3$. Repeated indices are meant to be summed, unless otherwise specified.

\tableofcontents

\chapter{Compressible Fluids and Non-linear Wave Equations}



\section{The Euler's Equations}

The mathematical description of the state of a moving fluid is determined by the distribution of the fluid velocity 
$\textbf{v}=\textbf{v}(\textbf{x},t)$ and of any two thermodynamic quantities pertaining to the fluid, for instance the pressure $p=p(\textbf{x},\textsl{t})$ 
and the mass density $\rho=\rho(\textbf{x},\textsl{t})$. All the thermodynamics are determined by the value of any two of them, together with the equation of state.
The equation of motion of a perfect fluid can be derived as follows.

We begin with the equation which expresses the conservation of mass. We consider the domain $\Omega(t)$ in $\mathbb{R}^3$  occupied by some fluid at time $t$. 
The mass of fluid in $\Omega(t)$ is $\int_{\Omega(t)}\rho dV$. This integral should not depend on $t$ due to the mass conservation. So we get

\begin{equation}
    \frac{d}{dt}\int_{\Omega(t)}\rho dV=0
\end{equation}

This means

\begin{equation}
    \int_{\Omega(t)}\frac{\partial\rho}{\partial t}dV+\int_{\partial\Omega(t)}
\rho\textbf{v}\cdot\textbf{n}dS=0
\end{equation}
where $\textbf{n}$ is the out normal of $\partial\Omega({t})$. By Green's formula, we have
\begin{equation}
 \int_{\Omega({t})}[\frac{\partial\rho}{\partial{t}}+\textrm{div}(\rho\textbf{v})]{dV} =0
\end{equation}
Since $\Omega({t})$ is arbitrary, we have

\begin{equation}
\frac{\partial\rho}{\partial{t}}+\textrm{div}(\rho\textbf{v})=0
\end{equation}
This is the equation of continuity.

Then we consider the conservation of momentum. By Newton's law,

\begin{equation}
\frac{d\textbf{P}}{dt}=\textbf{F}
\end{equation}

where $\textbf{P}$ is the momentum, $\textbf{F}$ is the force. In our case, equation (1.5) becomes

\begin{equation}
\frac{d}{dt}\int_{\Omega({t})}\rho\textbf{v} {dV}=-\int_{\partial\Omega({t})}{p}\textbf{n}{dS}
\end{equation}
 i.e.

\begin{equation}
\frac{d}{dt}\int_{\Omega({t})}\rho\textbf{v} {dV}=-\int_{\Omega({t})}\nabla{p}{dV}
\end{equation}

A component of this equation is

 \begin{equation}
    \frac{d}{dt}\int_{\Omega({t})}\rho{v}^i {dV}=-\int_{\Omega({t})}\frac{\partial{p}}{\partial {x}^i}{dV}
 \end{equation}

 Then we get

 \begin{equation}
    \int_{\Omega({t})}[\frac{\partial(\rho{v}^i)} {\partial{t}}+\textrm{div}(\rho{v}^i \textbf{v})]{dV}=
-\int_{\Omega({t})}\frac{\partial{p}}{\partial {x}^i}{dV}
 \end{equation}

 Since $\Omega({t})$ is arbitrary, as well as the equation of continuity, we get

 \begin{equation}
    \frac{\partial{v}^i}{\partial{t}}+\textbf{v}\cdot\nabla{v}^i=-\frac{1}{\rho}\frac{\partial {p}}{\partial{x}^i}
 \end{equation}

 In a perfect fluid, heat exchange between different parts of the fluid is absent. So the motion is adiabatic throughout the fluid.

In the adiabatic motion the entropy of any particle of fluid remains constant as that particle moves about in space. 
Denoting by ${s}$ the entropy per unit mass. We can express the condition for adiabatic motion as

\begin{equation}
    \frac{ds}{dt}=0
\end{equation}

which can be written as

\begin{equation}
    \frac{\partial{s}}{\partial{t}}+\textbf{v}\cdot\nabla{s}=0
\end{equation}

\section{Irrotational Flow and the Nonlinear Wave Equation}
The adiabatic condition may take a much simpler form. If the entropy is constant throughout the fluid at some initial instant, 
it remains everywhere the same constant value at all times. In this case, we can write the adiabatic condition simply as

\begin{equation}
    {s}=constant
\end{equation}

Such a motion is said to be isentropic.
We may in this case put (1.10) into a different form. To do this, we employ the thermodynamic relation

\begin{equation}
    {dh}={V}{dp}+\theta{ds}
\end{equation}

where ${h}$ is the enthalpy, defined in terms of the energy per unit mass:
\begin{align*}
 h=e+pV
\end{align*}

${V}=\frac{1}{\rho}$ is the specific volume, and $\theta$ is the temperature. Since $s$ is constant, we have:

\begin{equation}
    {dh}={V}{dp}=\frac{1}{\rho}{dp}
\end{equation}

So (1.10) becomes

\begin{equation}
    \frac{\partial{v}^i}{\partial {t}}+\textbf{v}\cdot\nabla{v}^i=
    -\frac{\partial{h}}{\partial{x}^i}
\end{equation}

This implies

\begin{equation}
    \frac{\partial\omega}{\partial{t}}+\textbf{v}\cdot\nabla
    \omega=\omega\cdot\nabla\textbf{v}-(\textrm{div}\textbf{v})\omega
\end{equation}

where $\omega=\nabla\times\textbf{v}$.

If $\omega\mid_{{t}=0}=0$, then by (1.17) $\omega\equiv0$. 
A flow for which $\omega=\nabla\times\textbf{v}\equiv0$ in all space and time is called an irrotational flow. In this case, since $\mathbb{R}^3$ is simply-connected, 
there exists a function $\phi$, such that $\textbf{v}=-\nabla\phi$. The equation (1.16) becomes

\begin{equation}
    \frac{\partial^2 \phi}{\partial{x}^i\partial{t}}-\sum_{j}\frac{\partial \phi}{\partial{x}^j}\frac{\partial^2 \phi}
{\partial {x}^j \partial {x}^i}=\frac{\partial {h}}{\partial {x}^i}
\end{equation}
i.e.
\begin{equation}
    \frac{\partial}{\partial {x}^i}[\frac{\partial \phi}{\partial {t}}-\frac{1}{2}|\nabla\phi|^{2}-{h}]=0
\end{equation}

Since $\phi$ is defined up to a constant which may depend on ${t}$, without loss of generality, we may set

 \begin{equation}
    {h}=\frac{\partial\phi}{\partial{t}}-\frac{1}{2}|\nabla\phi|^{2}
\end{equation}

Then only the continuity equation remains.
While $\frac{dh}{dp}={V}>0$, then by inverse function theorem, ${p}$ can be viewed as a function of ${h}$. 
This is because ${h}$ can be viewed as a function of ${p}$ and ${s}$, and in the isentropic case, 
${h}$ is a function of ${p}$. Then $\rho=\frac{dp}{dh}$, $\rho$ is also a function of ${h}$.

The continuity equation is

\begin{equation}
    \frac{\partial\rho}{\partial{t}}-\textrm{div}(\rho\nabla\phi)=0
\end{equation}

i.e.

\begin{equation}
    \frac{{d}\rho}{dh}\frac{\partial}{\partial{t}}(\frac{\partial \phi}{\partial{t}}
-\frac{1}{2}\sum_{i}(\frac{\partial\phi}{\partial{x}^i})^2)-
\sum_{j}\frac{{d}\rho}{dh}\frac{\partial}{\partial {x}^j}(\frac{\partial \phi}{\partial{t}}-
\frac{1}{2}\sum_{i}(\frac{\partial\phi}{\partial{x}^i})^2)\frac{\partial \phi}{\partial{x}^j}-\rho\Delta\phi=0
\end{equation}

\begin{equation}
    \frac{\partial^2\phi}{\partial{t}^2}-
2\sum_{i}\frac{\partial\phi}{\partial {x}^i}\frac{\partial^2\phi}{\partial{t}\partial {x}^i }+
\sum_{ij}\frac{\partial\phi}{\partial {x}^i}\frac{\partial\phi}{\partial {x}^j}\frac{\partial^2\phi}{\partial {x}^i\partial {x}^j}-\eta^2\Delta\phi=0
\end{equation}

where we have used the fact that $\frac{dp}{dh}=\rho$ and the definition of $\eta$, the sound speed, which is the following:

\begin{equation}
    \eta^2=\eta^2(h):=(\frac{\partial{p}}{\partial\rho})_{s}
\end{equation}

We define the function
\begin{equation}
    {H}={H}({h}):=-\eta^2-2{h}
\end{equation}

From (1.24) we know that
\begin{equation}
    \eta^{-2}=(\frac{{d}\rho}{dp})_{s}=
    -\frac{1}{{V}^2}(\frac{dV}{dp})_{s}
\end{equation}

Since ${V}=(\frac{dh}{dp})_{s}$, by direct calculation we obtain

\begin{equation}
    \frac{dH}{dh}=-\frac{\eta^4}{{V}^3}(\frac{{d}^2{V}}{{dp}^2})_{{s}}
\end{equation}

So $\frac{dH}{dh}$ vanishes if and only if ${V}$ is linear in ${p}$. This fact will be used later.
The meaning of the function $H$ will become apparent in the sequel.

Equation (1.21) is the Euler-Lagrange equation of the Lagrangian

\begin{equation}
    {L}={p}({h})
\end{equation}

In general, ${L}={L}({t},\textbf{x};\phi,\frac{\partial\phi}{\partial {t}},\nabla\phi)$, and the Euler-lagrange equation is

\begin{equation}
    \sum_{\alpha}\frac{\partial}{\partial{x}^\alpha}(\frac{\partial{L}({t},\textbf{x};\phi,\frac{\partial\phi}{\partial{t}},\nabla\phi) }{\partial {u}_{\alpha}})
=\frac{\partial{L}}{\partial{q}}(t,\textbf{x};\phi,\frac{\partial\phi}{\partial t},\nabla\phi)
\end{equation}

where ${q}=\phi$, ${x}^0={t}$, ${u}_{0}=\frac{\partial\phi}{\partial{t}}$, ${u}_{i}=\frac{\partial\phi}{\partial{x}^i}$  $i=1,2,3$

In the case of ${L}={p}({h})$,

\begin{equation}
    \frac{\partial{L}}{\partial{u}_{\alpha}}=\frac{{dp}}{{dh}}\frac{\partial{h} }{\partial{u}_{\alpha}}=\rho\frac{\partial{h}}{\partial{u}_{\alpha}}
\end{equation}

where ${h}={u}_{0}-\frac{1}{2}\sum_{i}({u}_{i})^2$

So we have 
\begin{align*}
 \frac{\partial{h}}{\partial{u}_{0}}=1, \frac{\partial{h}}{\partial{u}_{i}}=-{u}_{i}, \frac{\partial{h}}{\partial{q}}=0
\end{align*}

Then the Euler-lagrange equation is
\begin{equation}
    \frac{\partial\rho}{\partial{t}}-\sum_{i}\frac{\partial}{\partial{x}^i}(\rho{u}_{i})=0
\end{equation}

This is the equation (1.21).

\section{The Equation of Variations and the Acoustical Metric}
Next, we shall discuss the linearized equation corresponding to (1.21).

Let $\phi$ be a given solution of (1.21) and let ${\phi_{\tau}:\tau\in I}$, $I$ an open interval of the real line containing $0$, 
be a differentiable 1-parameter family of solution such that $\phi_{0}=\phi$.

Then

\begin{equation}
    \psi=(\frac{{d}\phi_{\tau}}{{d}\tau})_{\tau=0}
\end{equation}

is a variation of $\phi$ through solutions.

Consider the Lagrangian of the unknown $\psi$:
\begin{equation}
    {L}_{\tau}[\psi]:={L}[\phi+\tau\psi]
\end{equation}

Then the linearized Lagrangian of ${L}[\phi]$ is
\begin{equation}
    \dot{{L}}[\psi]:=\frac{1}{2}\frac{{d}^2{L}_{\tau}[\psi]}{{d}\tau^2}\mid_{\tau=0}
\end{equation}

In this case, ${L}_{\tau}[\psi]={p}({h}_{\tau})$, where 
\begin{align*}
 {h}_{\tau}=\frac{\partial}{\partial{t}}(\phi+\tau\psi)-\frac{1}{2}\sum_{i}[\frac{\partial}{\partial{x}^i}(\phi+\tau\psi)]^2
\end{align*}

By direct calculation, we have

\begin{equation}
    \frac{{d}^2{L}_{\tau}}{{d}\tau^2}[\psi]\mid_{\tau=0}=\frac{{d}^2{p}}{{dh}^2}
({h})[\frac{\partial\psi}{\partial{t}}-\sum_{i}\frac{\partial\phi}{\partial{x}^i}
\frac{\partial\psi}{\partial{x}^i}]^2+\frac{{dp}}{{dh}}({h})(-\sum_{i}(\frac{\partial\psi}{\partial{x}^i})^2)
\end{equation}
i.e.
\begin{equation}
    \frac{{d}^2{L}_{\tau}}{{d}\tau^2}[\psi]\mid_{\tau=0}=
\frac{{d}\rho}{{dh}}({h})(\frac{\partial\psi}{\partial{t}}
-\sum_{i}\frac{\partial\phi}{\partial{x}^i}\frac{\partial\psi}{\partial{x}^i})^2-\rho(\sum_{i}(\frac{\partial\psi}{\partial{x}^i})^2)
\end{equation}
i.e.
\begin{equation}
  \frac{{d}^2{L}_{\tau}}{{d}\tau^2}[\psi]\mid_{\tau=0}=
\rho[\eta^{-2}(\partial_{t}\psi-\sum_{i}\partial_{i}\phi\partial_{i}\psi)^{2}-\sum_{i}(\partial_{i}\psi)^2]
\end{equation}
i.e.
\begin{equation}
    \frac{{d}^2{L}_{\tau}}{{d}\tau^2}[\psi]\mid_{\tau=0}=-\rho({g}^{-1})^{\mu\nu}\partial_{\mu}\psi\partial_{\nu}\psi
\end{equation}

where the metric $g$ is:
\begin{equation}
    {g}=-\eta^2{dt}^2+\sum_{i}({dx}^{i}-{v}^{i}{dt})^2
\end{equation}
and
\begin{equation}
 {g}^{-1}=-\eta^{-2}(\frac{\partial}{\partial t}+v^{i}\frac{\partial}{\partial x^{i}})\otimes(\frac{\partial}{\partial t}+v^{j}\frac{\partial}{\partial x^{j}})
+\sum_{i}\frac{\partial}{\partial x^{i}}\otimes \frac{\partial}{\partial x^{i}}
\end{equation}

Here,${v}^i=-\partial_{i}\phi$.

We have used the fact: 
\begin{align*}
 \frac{\rho'}{\rho}=\frac{{d}\rho}{{dh}}\frac{1}{\rho}=\frac{{d}\rho}{{dp}}\frac{{dp}}{{dh}}\frac{1}{\rho}=\eta^{-2}
\end{align*}

Consider the conformal metric

\begin{equation}
    \tilde{g}_{\mu\nu}=\Omega{g}_{\mu\nu}
\end{equation}

which satisfy the following conditions:

\begin{equation}
    \rho({g}^{-1})^{\mu\nu}\partial_{\mu}\psi\partial_{\nu}\psi=(\tilde{g}^{-1})^{\mu\nu}
\partial_{\mu}\psi\partial_{\nu}\psi\sqrt{-{\det}\tilde{g}}
\end{equation}

Since $\sqrt{-{\det g}}=\eta$, we get $\sqrt{-{\det}\tilde{{g}}}=\Omega^2\eta$, we have $\Omega=\frac{\rho}{\eta}$

Since the linearized equation of (1.21) is the Euler-Lagrange equation of the linearized Lagrangian $\dot{{L}}[\psi]$, the linearized equation of (1.21) is
\begin{equation}
    \Box_{\tilde{{g}}}\psi=0
\end{equation}

Since $\rho_{0},\eta_{0}$, which are the mass density and the sound speed corresponding to the constant state, are constants, and the Euler-Lagrangian equations
 are not affected if the Lagrangian is multiplied by a constant, we may choose
\begin{equation}
\Omega=\frac{\rho/\rho_{0}}{\eta/\eta_{0}}
\end{equation}
so that $\Omega=1$ in the constant state.

We may in fact choose the unit of time in relation to the unit of length so that 
\begin{align*}
 \eta_{0}=1
\end{align*}
and we may further choose the unit of mass in relation to the unit of length so that
\begin{align*}
 \rho_{0}=1
\end{align*}
These choices shall be understood from the next chapter to the end of the book.

\section{The Fundamental Variations}
Now if $\phi$ is a solution of (1.21) and ${f}_{\tau}$ is a 1-parameter subgroup of the isometry group of $\mathbb{E}^3$, that is $\mathbb{R}^{3}$ with the standard
Euclidean metric
\begin{align*}
 \sum_{i}(dx^{i})^{2}
\end{align*}
 ,or the translation on ${t}$-direction, 
then for each $\tau$,

\begin{equation}
    \phi_{\tau}=\phi\circ{f}_{\tau}
\end{equation}

is also a solution. It follows that

\begin{equation}
    \psi=(\frac{{d}\phi\circ{f}_{\tau}}{{d}\tau})_{\tau=0}={X}\phi
\end{equation}

satisfies the linear equation (1.43)

${X}$ can be the following:
\begin{equation}
    \Omega_{ij}:={x}^{i}\partial_{j}-{x}^{j}\partial_{i}
\end{equation}
$1\leq i < j \leq 3$

and

\begin{equation}
    {T}_{\mu}:=\frac{\partial}{\partial{x}^\mu}
\end{equation}
$\mu=0,1,2,3$, which are the generators of ${f}_{\tau}$.

Consider now the scaling group for equation (1.21). It is easy to see that the physical dimension of $\phi$ is $\frac{length^{2}}{time}$, so we consider
the transformation, with constants $a,b>0$:
\begin{align}
 \tilde{\phi}(t,\textbf{x}):=\frac{b^{2}}{a}\phi(\frac{t}{a},\frac{\textbf{x}}{b})
\end{align}
A direct calculation implies
\begin{align}
 \tilde{h}(t,\textbf{x}):=\frac{\partial\tilde{\phi}}{\partial t}-\frac{1}{2}|\nabla\tilde{\phi}|^{2}
=\frac{b^{2}}{a^{2}}h(\frac{t}{a},\frac{\textbf{x}}{b})
\end{align}
In general, $p$ is an arbitrary function of $h$:
\begin{align}
 p=p(h)
\end{align}
so 
\begin{align}
 \frac{dp}{dh}=\rho>0\quad\textrm{i.e}\quad p^{\prime}>0\quad\textrm{and also}\quad p>0
\end{align}
and
\begin{align}
 \frac{dp}{d\rho}=\eta^{2}=\frac{dp}{dh}\frac{dh}{d\rho}=\rho/\frac{d^{2}p}{dh^{2}}>0
\end{align}

Under the transformation (1.49), we have:
\begin{align}
 \tilde{p}(t,\textbf{x})=p(\tilde{h}(t,\textbf{x}))=p(\frac{b^{2}}{a^{2}}h(\frac{t}{a},\frac{\textbf{x}}{b}))
\end{align}
and then
\begin{align}
 \tilde{\rho}(t,\textbf{x})=\frac{d\tilde{p}}{d\tilde{h}}=p^{\prime}(\frac{b^{2}}{a^{2}}h(\frac{t}{a},\frac{\textbf{x}}{b}))=\rho(\frac{b^{2}}{a^{2}}h(\frac{t}{a},
\frac{\textbf{x}}{b}))
\end{align}
Now we can transform the Euler-Lagrange equation (1.21) under the transformation (1.49).

We can see that in general, only in the case $a=b$ the Euler-Lagrange equation is invariant under the transformation (1.49). But still, there are some 
special cases that we have 2-dimensional scaling group.

     First we consider the case:
\begin{align}
 p(h)=\frac{1}{2}h^{2}
\end{align}
Note that this function $p(h)$ satisfies (1.51)-(1.53)

Now we may take $a\slashed{=}b$:
\begin{align*}
 \tilde{h}(t,\textbf{x})=\frac{b^{2}}{a^{2}}h(\frac{t}{a},\frac{\textbf{x}}{b})
\end{align*}
Then
\begin{align}
 p(\tilde{h}(t,\textbf{x}))=\frac{1}{2}\tilde{h}^{2}(t,\textbf{x})=\frac{b^{4}}{a^{4}}\frac{1}{2}h^{2}(\frac{t}{a},\frac{\textbf{x}}{b})
=\frac{b^{4}}{a^{4}}p(h(\frac{t}{a},\frac{\textbf{x}}{b}))
\end{align}
That is 
\begin{align}
 \tilde{p}(t,\textbf{x})=\frac{b^{4}}{a^{4}}p(\frac{t}{a},\frac{\textbf{x}}{b})
\end{align}
and also:
\begin{align}
 \tilde{\rho}(t,\textbf{x})=\tilde{h}(t,\textbf{x})=\frac{b^{2}}{a^{2}}h(\frac{t}{a},\frac{\textbf{x}}{b})
=\frac{b^{2}}{a^{2}}\rho(\frac{t}{a},\frac{\textbf{x}}{b})
\end{align}
So now we have:
\begin{align}
 \frac{\partial\tilde{\rho}}{\partial t}-\textrm{div}(\tilde{\rho}\nabla\tilde{\phi})\\\notag
=\frac{b^{2}}{a^{3}}\frac{\partial\rho}{\partial t}-\frac{b^{2}}{a^{2}}\frac{1}{a}\rho\Delta\phi-\frac{b}{a^{2}}\nabla\rho\cdot\frac{b}{a}\nabla\phi
=\frac{b^{2}}{a^{3}}(\frac{\partial\rho}{\partial t}-\textrm{div}(\rho\nabla\phi))
\end{align}
This means the equation (1.21) is invariant under the transformation (1.49). So in this case, the scaling group is 2-dimensional.

     More generally, we can consider the case 
\begin{align}
 p(h)=Ch^{\alpha}, \alpha>1, C>0, h>0
\end{align}
Obviously in this case the function $p(h)$ satisfies (1.51)-(1.53) and we have:
\begin{align}
 p=C(\frac{\rho}{\alpha C})^{\frac{\alpha}{\alpha-1}}
\end{align}
If we set:
\begin{align*}
 \gamma=\frac{\alpha}{\alpha-1},\quad k=\frac{C}{(\alpha C)^{\frac{\alpha}{\alpha-1}}}
\end{align*}
then
\begin{align}
 p=k\rho^{\gamma}
\end{align}
with $\gamma>1$ and $k$ depending on $\gamma$.

This is the polytropic case. We can calculate as above to see that in this case, we also have the 2-dimensional scaling group. Also, one can easily see 
that we have the 2-dimensional scaling group only when the function $p$ is a homogeneous function of $h$.

     In present book, we shall consider the general case $p=p(h)$, so the scaling group is 1-dimensional. In this case, we consider the 1-parameter group of 
dilations of $\mathbb{R}\times\mathbb{E}^{3}$, given by:
\begin{align}
(t,\textbf{x})\mapsto (e^{\tau}t, e^{\tau}\textbf{x})
\end{align}
where $\tau\in\mathbb{R}$.

Let $\phi$ be a solution of (1.21) and define for each $\tau\in\mathbb{R}$ the function
\begin{align}
 \phi_{\tau}:=e^{-\tau}\phi(e^{\tau}t, e^{\tau}\textbf{x})
\end{align}
As we have seen before, $\phi_{\tau}$ is also a solution of (1.21). Then
\begin{align}
 \psi=(\frac{d\phi_{\tau}}{d\tau})|_{\tau=0}:=S\phi
\end{align}
satisfies the linear equation (1.43). Here $S$ is the differential operator
\begin{align}
 S=D-I,\quad D=x^{\mu}\frac{\partial}{\partial x^{\mu}}
\end{align}

\chapter{The Basic Geometric Construction}



\section{Null Foliation Associated with the Acoustical Metric}

\subsection{Galilean Spacetime}
Since the Galilean spacetime is the frame of classical fluid mechanics, we shall mention some basic facts about it.

 A $Galilean$ $spacetime$ is a quadruple $(\mathcal{G},\mathcal{A},\pi,e)$ where:

(1) $\mathcal{G}$ is a 4-dimensional affine space;

(2) $\mathcal{A}$ is a 1-dimensional affine space with unit displacement, which represents time;

(3) The map $\pi$: $\mathcal{G}\rightarrow \mathcal{A}$ is a surjective affine map;

So $\Sigma_{t}:=\pi^{-1}(t)$ is an affine subspace of $\mathcal{G}$ for each $t\in \mathcal{A}$, and we can introduce a Euclidean metric $e_{t}$ on each $\Sigma_{t}$.
The $\Sigma_{t}$ is the hypersurfaces of absolute simultaneity.

(4) A $Galilean$ $frame$ is an affine foliation of $\mathcal{G}$ by a family of parallel lines transversal to $\Sigma_{t}$. A Galilean frame represents a family of
inertial observers at rest relative to each other. The metric $e$ satisfies the condition that parallel transport along the lines of a Galilean frame induces an isometry
between $(\Sigma_{t_{1}},e_{t_{1}})$ and $(\Sigma_{t_{2}},e_{t_{2}})$ for any $t_{1}, t_{2}\in\mathcal{A}$. 

Two Galilean spacetimes, $(\mathcal{G},\mathcal{A},\pi,e)$ and $(\mathcal{G}^{\prime},\mathcal{A}^{\prime},\pi^{\prime},e^{\prime})$ are equivalent if and 
only if there is a pair $(I,J)$, where $I$ is an affine spaces isomorphism of $\mathcal{G}$ onto $\mathcal{G}^{\prime}$, and $J$ is an affine spaces isomorphism
of $\mathcal{A}$ onto $\mathcal{A}^{\prime}$ such that

$(1)\quad \pi^{\prime}\circ I=J\circ\pi$ ;

(2) $J$ takes the unit displacement of $\mathcal{A}$ into the unit displacement of $\mathcal{A}^{\prime}$ ;

(3) $I|_{\Sigma_{t}}$ is an isometry of $(\Sigma_{t},e_{t})$ onto $(\Sigma^{\prime}_{t^{\prime}}, e^{\prime}_{t^{\prime}})$, for each $t\in\mathcal{A}$. Here 
$\Sigma_{t}=\pi^{-1}(t)$, $\Sigma^{\prime}_{t^{\prime}}=\pi^{\prime -1}(t^{\prime})$, and $t^{\prime}=J(t)$.

 Classical mechanics including classical continuum mechanics are invariant under such Galilean transformations. We have used in fact this Galilean invariance to define the 
first order variations in Chapter 1.

However, in this book, we shall adopt a different point of view in analyzing the equation of variations. We shall retain from the Galilean spacetime $\mathcal{G}$ only
its differentiale manifold structure $\mathbb{R}^{4}$, a manifold endowed with the acoustical metric $g$. We call the pair $(\mathbb{R}^{4},g)$ the acoustical spacetime,
which is a Lorentzian manifold. Two such Lorentzian manifolds $(\mathcal{M},g)$ and $(\mathcal{M}^{\prime},g^{\prime})$ are equivalent if and only if there is 
a diffeomorphism $f$ of $\mathcal{M}^{\prime}$ onto $\mathcal{M}$ such that $g^{\prime}=f^{*}g$, i.e. $(\mathcal{M}^{\prime},g^{\prime})$ is isometric to
$(\mathcal{M},g)$, through  an arbitrary  diffeomorphism $f$.

\subsection{Null Foliation and Acoustical Coordinates}
The initial data for the equation of motion
\begin{equation}
\frac{\partial\rho}{\partial t}+\textrm{div}(\rho\textbf{v})=0
\end{equation} 
\begin{equation}
\frac{\partial\textbf{v}}{\partial t}+\textbf{v}\cdot\nabla\textbf{v}=-\frac{1}{\rho}\nabla p
\end{equation}
\begin{equation} 
\frac{\partial s}{\partial t}+\textbf{v}\cdot\nabla s=0
\end{equation} 
is to be given on the hyperplane $\Sigma_{0}$ and consists in specification of the triplet $(p,s,\textbf{v})$. 
We assume that there is a sphere $S_{0,0}$ outside which the initial data coincide with those of a constant state, that is we have 
\begin{equation}
 p=p_{0},\quad s=s_{0},\quad \textbf{v}=0
\end{equation}
By choosing the unit of length equal to the radius of $S_{0,0}$, we then take $S_{0,0}$ to be the unit sphere centered at the origin in $\Sigma_{0}$.
We consider an annular interior neighborhood of $S_{0,0}$ in $\Sigma_{0}$:
\begin{equation}
 \Sigma^{\epsilon_{0}}_{0}=\{\textbf{x} \in \Sigma_{0}: 1-\epsilon_{0}\leqslant r(\textbf{x})\leqslant1\}
\end{equation}
where $r$ is the Euclidean distance function from the origin in $\Sigma_{0}$ and ${\epsilon_{0}}$ is a positive constant, satisfying the condition:
\begin{equation}
 \epsilon_{0}\leqslant\frac{1}{2}
\end{equation}
We define in $\Sigma_{0}$ the function
\begin{equation}
 u=1-r
\end{equation}
which on $\Sigma_{0}$ minus the origin is a smooth function without critical points, vanishing on $S_{0,0}$ and increasing inward. 
For each value of $u$ in the closed interval $[0,\epsilon_{0}]$, the corresponding level set $S_{0,u}$ of $u$ 
is a sphere of radius $1-u$ in the interval $[1-\epsilon_{0},1]$, and we have:
\begin{equation}
 \Sigma^{\epsilon_{0}}_{0}=\bigcup_{u \in [0,\epsilon_{0}]}  S_{0,u}
\end{equation}
In the case that the initial data is irrotational and isentropic outside $S_{0,\epsilon_{0}}$, the data in the exterior of $S_{0,\epsilon_{0}}$ in
$\Sigma_{0}$ specify the pair $(\phi,\partial_{t}\phi)$ in this region. This is initial data for the nonlinear wave equation (1.23) in the region in question 
and we have:
\begin{equation}
 \phi=0,\quad \partial_{0}\phi=h_{0}
\end{equation}
$h_{0}$ is the enthalpy of constant state. in the exterior of $S_{0,0}$ in $\Sigma_{0}$.
As we have seen in Chapter 1, the enthalpy $h$ is defined up to an additive constant, we may use this freedom to set:
\begin{align*}
 h_{0}=0
\end{align*}
To any given initial data set as above there corresponds a unique $maximal$ $solution$ of the equation (2.1)-(2.3), of the nonlinear wave equation (1.23) 
in the irrotational isentropic case. The notion of maximal solution or maximal development of an initial data set is the following. Given an initial data set, 
the local existence theorem asserts the existence of a $development$ of this set, namely of a domain $\mathcal{D}$ in Galilean spacetime, whose past boundary is 
the domain of the initial data, and of a solution defined in $\mathcal{D}$ and taking the given data at the past boundary, such that if we consider any point 
$p \in \mathcal{D}$ and any curve issuing at $p$ with the property that its tangent vector at any point $q$ belongs to $\bar{I}^{-}_{q}$, the closure of the past 
component of the double cone defined by $g_{q}$, the acoustical metric at $q$, then the curve terminates in the past at a point of the domain of the initial data.
The local uniqueness theorem asserts that if $(\mathcal{D}_{1},(p_{1}, s_{1}, \textbf{v}_{1}))$ and $(\mathcal{D}_{2},(p_{2}, s_{2}, \textbf{v}_{2}))$
are two developments of the same initial data set ($(\mathcal{D}_{1},\phi_{1})$ and $(\mathcal{D}_{2},\phi_{2})$ in the irrotational isentropic case), then $(p_{1},
s_{1}, \textbf{v}_{1})$ coincides with $(p_{2}, s_{2}, \textbf{v}_{2})$ in $\mathcal{D}_{1} \bigcap \mathcal{D}_{2}$ ($\phi_{1}$ coincides with 
$\phi_{2}$ in $\mathcal{D}_{1} \bigcap \mathcal{D}_{2}$ in the irrotational isentropic case). It follows that the union of all developments of a given initial data
set is itself a development, the unique $maximal$ $development$ of the initial data set.
We consider, in the domain of maximal solution, the family $\{C_{u}: u \in [0,\epsilon_{0}]\}$ of outgoing characteristic hypersurfaces corresponding to the family
$\{S_{0,u}: u \in [0,\epsilon_{0}]\}$:
\begin{equation}
 C_{u}\bigcap\Sigma_{0}=S_{0,u}: \forall u \in [0,\epsilon_{0}]
\end{equation}
Each bicharacteristic generator of each $C_{u}$ is to extend in the domain of maximal solution as long as it remains on the boundary of the domain of dependence of the
exterior of the surface $S_{0,u}$ in $\Sigma_{0}$. if we denote by $W_{\epsilon_{0}}$ the spacetime domain :
\begin{equation}
 W_{\epsilon_{0}}=\bigcup_{u \in [0,\epsilon_{0}]} C_{u}
\end{equation}
then the domain $M_{\epsilon_{0}}$ of the maximal solution corresponding to the given data set is the union of $W_{\epsilon_{0}}$ with the domain in Galilean spacetime
bounded by the exterior of the unit sphere $S_{0,0}$ in $\Sigma_{0}$ and by the outgoing characteristic hypersurface $C_{0}$ corresponding to
 $S_{0,0}$. By the domain of dependence theorem the solution coincide in $M_{\epsilon_{0}} \setminus W_{\epsilon_{0}}$ with the constant state. 
This implies that $C_{0}$ is a complete cone, each of its bicharacteristic generators extending to infinity in the parameter $t$.
We extend the function $u$ to $W_{\epsilon_{0}}$ by requiring that its level sets are precisely the outgoing characteristic hypersurfaces $C_{u}$. The
function $u$ is then the solution of the equation:
\begin{equation}
 (g^{-1})^{\mu\nu}\partial_{\mu}u\partial_{\nu}u=0
\end{equation}
We shall call such a function an $acoustical$ $function$. The vector field $\hat{L}$ given by:
\begin{equation}
 \hat{L}=-(g^{-1})^{\mu\nu}\partial_{\nu}u
\end{equation}
is then a future-directed null geodesic vectorfield with respect to the Lorenzian metric $g$. Its integral curves are the bicharacteristic generators of each 
$C_{u}$. Now the parametrization of these given by $\hat{L}$ is affine; however, for reasons which shall become apparent in the following we wish 
the generators
to be parametrized by the function $t$ instead. For this reason we choose the colinear vectorfield 
\begin{equation}
 L=\mu\hat{L}
\end{equation}
where the proportionality factor $\mu$, a positive function, is chosen so that
\begin{equation}
 Lt=1
\end{equation}
Thus
\begin{equation}
 \frac{1}{\mu}=-(g^{-1})^{\mu\nu}\partial_{\mu}t\partial_{\nu}u
\end{equation}
For each $u \in [0,\epsilon_{0}]$ there is a greatest lower bound $t_{\ast}(u)$ of the extent of the generators of $C_{u}$,
in the parameter $t$, in the domain of the maximal solution. This $t_{\ast}(u)$ is either a positive real number or $\infty$. Acoording to the 
above we have $t_{\ast}(0)=\infty$. Let us denote
\begin{equation}
 t_{\ast \epsilon_{0}}=\inf_{u \in [0,\epsilon_{0}]} t_{\ast}(u)
\end{equation}
In the following we shall confine attention to the spacetime domain 
\begin{equation}
 W^\ast_{\epsilon_{0}}=\bigcup_{(t,u)\in [0,t_{\ast \epsilon_{0}})\times[0,\epsilon_{0}]}S_{t,u}
\end{equation}
where
\begin{equation}
 S_{t,u}=C_{u}\bigcap\Sigma_{t}
\end{equation}
Define in $W^{\ast}_{\epsilon_{0}}$ the vectorfield $T$ such that it is tangential to $\Sigma_{t}$, orthogonal to the surfaces $\{S_{t,u}:
u \in [0,\epsilon_{0}]\}$ in metric $g$, and it verifies:
\begin{equation}
 Tu=1
\end{equation}
Consider the commutator:
\begin{equation}
\Lambda=[L,T]
\end{equation}
From (2.12)-(2.14); we know
\begin{equation}
Lu=0
\end{equation}
and by the fact that $T$ is tangential to $\Sigma_{t}$, we have:
\begin{equation}
 Tt=0
\end{equation}
then from (2.15),(2.20), we get:
\begin{equation}
 \Lambda u=\Lambda t=0
\end{equation}
Therefore $\Lambda$ is tangential to $S_{t,u}$.
From (2.14) and (2.22), we know that
\begin{equation}
 g(L,L)=0
\end{equation}
i.e. $L$ is a null vector with respect to the acoustical metric $g$. Also, from (2.13),(2.14) and (2.20),
\begin{equation}
 g(L,T)=\mu g(\hat{L},T)=-\mu Tu=-\mu
\end{equation}
Since the restriction of $g$ to $\Sigma_{t}$ is the Euclidean metric, $T$ is a spacelike vector with respect to $g$ and we can write 
\begin{equation}
 g(T,T)=\kappa^2
\end{equation}
where $\kappa$ is a positive function, the Euclidean magnitude of $T$.
Let $N$ be any vector tangent to one of $C_{u}$, then we have 
\begin{equation}
 Nu=0
\end{equation}
i.e.
\begin{equation}
 g(L,N)=0
\end{equation}
So for any $X$ tangent to one of $S_{t,u}$, 
\begin{equation}
 g(L,X)=0
\end{equation}
Since by definition, $X$ is also $g$ orthogonal to $T$, it follows that
\begin{equation}
 T_{p}W^{\ast}_{\epsilon_{0}}=\Pi_{p} \oplus T_{p}S_{t,u}
\end{equation}
 where $\Pi_{p}$ is spanned by $L$ and $T$.
The metric $g$ in $\Pi_{p}$ is given by (2.25)-(2.27) in terms of $\mu$ and $\kappa$. By (2.31) $g$ is then completely specified 
once we give the metric $\slashed{g}$ induced on $S_{t,u}$:
\begin{equation}
\slashed{g}(X,Y)=g(X,Y)\quad \forall X, Y \in T_{p}S_{t,u}
\end{equation}

For each $u \in [0,\epsilon_{0}]$ the generators of $C_{u}$ defines a smooth mapping of $S_{0,u}$ to 
$S_{t,u}$ for each $t \in [0,t_{\ast\epsilon_{0}})$. While each $S_{0,u}$ is diffeomorphic to $S^{2}$. 
Therefore, we get a diffeomorphism from $S_{t,u}$ to $S^{2}$:
\begin{equation}
 p\mapsto\vartheta=\vartheta(p) 
\end{equation}
where $\vartheta\in S^{2}$, $p \in S_{t,u}$. If local coordinates $(\vartheta^{1},\vartheta^{2})$ 
are chosen on $S^{2}$, this diffeomorphsim defines local coordinates on $S_{t,u}$ for every $(t,u) 
\in [0,t_{\ast\epsilon_{0}})\times[0,\epsilon_{0}]$
Since the diffeomorphism from $S_{0,u}$ to $S^{2}$ is arbitrary, the diffeomorphsim (2.33) is arbitrary as it may be composed on the left by a 
transformation of the form:
\begin{equation}
 \vartheta \mapsto \tilde{\vartheta}=\tilde{\vartheta}(u,\vartheta)
\end{equation}
The local coordinates $(\vartheta^{1},\vartheta^{2})$, together with the functions $(t,u)$ define a complete system  of local coordinates 
$(t,u,\vartheta^{1},\vartheta^{2})$ for $W^\ast_{\epsilon_{0}}$. We shall call these $acoustical$ $coordinates$ and we shall derive 
an expression for $g$ in $W^{\ast}_{\epsilon_{0}}$ in these coordinates.

First, the integral curves of $L$ are the lines of constant $\vartheta$ and $u$, parametrized by $t$. Therefore
\begin{equation}
 L=\frac{\partial}{\partial t}
\end{equation}
Next, by (2.20),(2.23), we have
\begin{equation}
 T=\frac{\partial}{\partial u}-\Xi
\end{equation}
where $\Xi$ is a vectorfield tangential to $S_{t,u}$. Thus $\Xi$ can be expanded in terms of the frame field 
($\frac{\partial}{\partial\vartheta^{A}}$: $\textsl{A}=1,2$).and we have
\begin{equation}
 \Xi=\Xi^{A}\frac{\partial}{\partial\vartheta^{A}}
\end{equation}
By (2.21),(2.35),(2.36) we have
\begin{equation}
 [L,\Xi]=-\Lambda
\end{equation}
or, in terms of components,
\begin{equation}
 \frac{\partial\Xi^{A}}{\partial t}=-\Lambda^{A}
\end{equation}
We can set $\Xi=0$ along any one of the hypersurfaces $\Sigma_{t}$. However the non-vanishing of $\Lambda$ forbids setting $\Xi=0$ everywhere.
With the components
\begin{equation}
 \slashed{g}_{AB}=\slashed{g}(\frac{\partial}{\partial\vartheta^{A}},\frac{\partial}{\partial\vartheta^{B}})
\end{equation}
where $A,B=1,2$, we have
\begin{equation}
g=-2\mu dudt+\kappa^{2}du^{2}+\slashed{g}_{AB}(d\vartheta^{A}+\Xi^{A}du)(d\vartheta^{B}+\Xi^{B}du)
\end{equation}
We define in $W^{*}_{\epsilon_{0}}$ the vectorfield $B$ by the conditions that it be orthogonal to $\Sigma_{t}$ with respect to $g$ 
and that it verifies :
\begin{equation}
 Bt=1
\end{equation}
Since $\Sigma_{t}$ is spacelike relative to $g$, $B$ is future-directed timelike relative to $g$. Thus, there is a positive function 
$\alpha$ such that
\begin{equation}
 g(B,B)=-\alpha^{2}
\end{equation}
In fact, we have
\begin{equation}
 B^{\mu}=-\alpha^{2}(g^{-1})^{\mu\nu}\partial_{\nu}t
\end{equation}
   and:
\begin{equation}
 \alpha^{-2}=-(g^{-1})^{\mu\nu}\partial_{\mu}t\partial_{\nu}t=\eta^{-2}
\end{equation}
Since $\alpha$ and $\eta$ are both positive, $\alpha=\eta$.
Then we have
\begin{equation}
 B^{0}=-\alpha^{2}(g^{-1})^{00}=1, \quad B^{i}=-\alpha^{2}(g^{-1})^{i0}=-\alpha^{2}(-\eta^{-2})v^{i}=v^{i}
\end{equation}
So 
\begin{align}
 B=\frac{\partial}{\partial t}+v^{i}\frac{\partial}{\partial x^{i}}
\end{align}
and in irrotational isentropic case,
\begin{equation*}
 B=(1,-\nabla\phi)
\end{equation*}

Now at any point $p$, $B \in \Pi_{p}$, a timelike plane relative to $g$. Therefore, $B$ is a linear combination of $L$ 
and $T$. From (2.36),(2.35),(2.42), and the fact that $B$ is future-directed timelike relative to $g$, we have
\begin{equation}
 B=L+fT
\end{equation}
where $f$ is a positive function. Taking $g$-inner product of (2.48) with $T$, using(2.26) and the fact that $g(B,T)=0$ we 
have
\begin{equation}
 \mu=f\kappa^{2}
\end{equation}
On the other hand, $L$ is null relative to $g$, we have
\begin{equation}
 0=g(B,B)+f^{2}g(T,T)=-\alpha^2+f^{2}\kappa^{2}
\end{equation}
Since $f$ is positive, we have
\begin{equation}
 f=\frac{\alpha}{\kappa}
\end{equation}
Then we get
\begin{equation}
 \mu=\alpha\kappa
\end{equation}
Let us introduce the vectorfield $\hat{T}:=\kappa^{-1}T$. Then from (2.48),(2.51), we have
\begin{equation}
 L=B-\alpha\hat{T}
\end{equation}
By (2.46), (2.47),we have
\begin{equation}
 L=\frac{\partial}{\partial t}-(\alpha\hat{T}^{i}-v^{i})\frac{\partial}{\partial x^{i}}
\end{equation}
in the irrotational isentropic case:
\begin{equation}
 L=\frac{\partial}{\partial t}-(\alpha\hat{T}^{i}+\partial_{i}\phi)\frac{\partial}{\partial x^{i}}
\end{equation}
We finally calculate the Jacobian of the mapping:
\begin{align}
(t,u,\vartheta^{1},\vartheta^{2})\mapsto(x^{0},x^{1},x^{2},x^{3})
\end{align}
Since $x^{0}=t$, we have
\begin{equation}
 \frac{\partial x^{0}}{\partial t}=1, \quad\frac{\partial x^{0}}{\partial u}=\frac{\partial x^{0}}{\partial\vartheta^A}=0
\end{equation}
where $A=1,2$.
Also, since
\begin{equation}
 \frac{\partial x^{\mu}}{\partial t}=L^{\mu}
\end{equation}
we have
\begin{equation}
 \frac{\partial x^{i}}{\partial t}=L^{i}
\end{equation}
where $i=1,2,3$
Similarly, we have
\begin{equation}
 \frac{\partial x^{i}}{\partial u}=T^{i}+\xi^{i}
\end{equation}
where
\begin{equation}
 \xi^{i}=\Xi^{A}X^{i}_{A}, \quad  X^{i}_{A}=\frac{\partial x^{i}}{\partial\vartheta^{A}}
\end{equation}
where $i=1,2,3; A=1,2$.

Then we get the Jacobian determinant of (2.56):
\begin{equation}
 \Delta=\begin{vmatrix} 
         1&0&0&0\\ L^{1}& T^{1}+\xi^{1}& X^{1}_{1}& X^{1}_{2}\\ L^{2}& T^{2}+\xi^{2}& X^{2}_{1}& X^{2}_{2}\\
L^{3}& T^{3}+\xi^{3}& X^{3}_{1}& X^{3}_{2}
        \end{vmatrix}
\end{equation}
 We have,
\begin{equation}
 \Delta=\begin{vmatrix}
         T^{1}+\xi^{1}& X^{1}_{1}& X^{1}_{2}\\ T^{2}+\xi^{2}& X^{2}_{1}& X^{2}_{2}\\
T^{3}+\xi^{3}& X^{3}_{1}& X^{3}_{2}
        \end{vmatrix}=\dot{\Delta}+\ddot{\Delta}
\end{equation}
where
\begin{equation}
 \dot{\Delta}=\begin{vmatrix}
         T^{1}&X^{1}_{1}& X^{1}_{2}\\ T^{2}& X^{2}_{1}& X^{2}_{2}\\
T^{3}& X^{3}_{1}& X^{3}_{2}
        \end{vmatrix}
\end{equation}
and
\begin{equation}
 \ddot{\Delta}=\begin{vmatrix}
         \xi^{1}& X^{1}_{1}& X^{1}_{2}\\\xi^{2}& X^{2}_{1}& X^{2}_{2}\\
\xi^{3}& X^{3}_{1}& X^{3}_{2}
        \end{vmatrix}
\end{equation}
Now from (2.61), we have
\begin{equation}
 \ddot{\Delta}=\sum^{2}_{A=1}\Xi^{A}\begin{vmatrix}
         X^{1}_{A}&X^{1}_{1}&X^{1}_{2}\\ X^{2}_{A}& X^{2}_{1}& X^{2}_{2}\\
X^{3}_{A}&X^{3}_{1}&X^{3}_{2}
        \end{vmatrix}=0
\end{equation}
then we can simply write
\begin{equation}
 \Delta=(T, X_{1}, X_{2})
\end{equation}
i.e.
\begin{equation}
 \Delta=|T||X_{1}\wedge X_{2}|
\end{equation}
where the magnitude is with respect to $\bar{g}$, the restriction of $g$ on $\Sigma_{t}$, which is the Euclidean metirc.
then we know that 
\begin{equation}
 |X_{1}\wedge X_{2}|=\sqrt{|X_{1}|^{2}|X_{2}|^{2}-(X_{1},X_{2})^{2}}=\sqrt{\det\slashed{g}}
\end{equation}
then we get
\begin{equation}
 \Delta=\kappa\sqrt{\det\slashed{g}}
\end{equation}

\section{A Geometric Interpretation for Function $H$}
Recall the definition of $H$ in Chapter 1:
\begin{align}
 H=-2h-\eta^{2}
\end{align}
Here we just consider the case $H=const$. Actually, we just need to consider the case $H=0$ due to the following transformation:
\begin{align}
 \tilde{\phi}:=\phi+\frac{k}{2}t,\quad \tilde{h}:=h+\frac{k}{2},\quad \tilde{H}:=H-k\quad k\quad\textrm{is}\quad\textrm{the}\quad ``const''
\end{align}
which leaves the non-linear wave equation (1.23) invariant.

     Since 
\begin{align*}
 H=0
\end{align*}
implies
\begin{align*}
 h=-\frac{1}{2}\eta^{2}
\end{align*}
in this case the acoustical metric is:
\begin{align*}
 g=2\frac{\partial\phi}{\partial t}(dt)^{2}+2\sum_{i}\frac{\partial\phi}{\partial x^{i}}dx^{i}dt+\sum_{i}(dx^{i})^{2}
\end{align*}
i.e.
\begin{align}
 g=2d\phi dt+\sum_{i}(dx^{i})^{2}
\end{align}
Introducing a new coordinate $s$, we consider the Lorentzian manifold $(\mathbb{R}^{1+4}, \tilde{g})$ where:
\begin{align}
 \tilde{g}=2dsdt+\sum_{i}(dx^{i})^{2}
\end{align}
So the manifold $(\mathbb{R}\times\mathbb{E}^{3}, g)$ can be realized as a submanifold of $(\mathbb{R}^{1+4}, \tilde{g})$ with the induced metric, in fact 
as the graph over the null hyperplane:
\begin{align}
 s=\phi(\textbf{x},t)
\end{align}
Also, by direct calculation, we find:
\begin{align}
 \det g=2h=-\eta^{2}<0
\end{align}
which implies that this submanifold is a time-like hypersurface in $(\mathbb{R}^{1+4}, \tilde{g})$.

Moreover, the manifold $(\mathbb{R}\times\mathbb{E}^{3},g)$ is a minimal submanifold of $(\mathbb{R}^{1+4},\tilde{g})$. Since in general
\begin{align*}
 \frac{d\rho}{dp}=\eta^{-2}\quad\textrm{and}\quad\frac{d\rho}{dp}=\frac{1}{\rho}\frac{d\rho}{dh}
\end{align*}
In the present case we have by (2.71):
\begin{align*}
 \frac{1}{\rho}\frac{d\rho}{dh}=-\frac{1}{2h}
\end{align*}
Hence we can solve this ODE to obtain:
\begin{align}
 \rho=\frac{C}{\sqrt{-2h}}=\frac{C}{\eta}
\end{align}
where $C$ is a positive constant. The equation:
\begin{align*}
 \frac{dp}{dh}=\rho
\end{align*}
gives
\begin{align*}
 p=p_{0}-C\sqrt{-2h}\quad\textrm{i.e}\quad
\end{align*}
By (2.77) we can write:
\begin{align}
 p=p_{0}-C^{2}V\quad\textrm{or}\quad p=p_{0}-C\eta
\end{align}
Thus, the action in a domain $\Omega$ in $\mathbb{R}\times\mathbb{E}^{3}$ is:
\begin{align}
 \int_{\Omega}pdtd^{3}x=p_{0}\int_{\Omega}dtd^{3}x-C\int_{\Omega}\eta dtd^{3}x
\end{align}
 While
\begin{align}
 \int_{\Omega}\eta dtd^{3}x=\int_{\Omega}\sqrt{-\det g}dtd^{3}x
\end{align}
is the volume of $\Omega$ with respect to $g$.

It follows that the Euler-Lagrange equation coincides with that corresponding to the Lagrangian $\sqrt{-\det g}$. But the Euler-Lagrange equation
corresponding to the volume of a domain as the action in the domain is the minimal surface equation.

\chapter{The Acoustical Structure Equations}



\section{The Acoustical Structure Equations}
Define $k$ by 
\begin{equation}
 2\alpha k=\bar{\mathcal{L}}_{B}g
\end{equation}
where $\bar{\mathcal{L}}$ is the restriction of $\mathcal{L}$ to $\Sigma_{t}$.

If $X$, $Y$ are two vectors tangent to $\Sigma_{t}$ at a point, we then have 
\begin{equation}
 \alpha k(X,Y)=g(D_{X}B, Y)=g(D_{Y}B, X)
\end{equation}
where $k$ is the 2nd fundamental form of $\Sigma_{t}$ relative to $g$ and $D$ is the covariant derivative associated to $g$.
Recall from Chapter 2, we have
\begin{equation}
 B=\partial_{0}+v^{i}\partial_{i}
\end{equation}
Then
\begin{equation}
 2\alpha k_{ij}=(\bar{\mathcal{L}}_{B}g)_{ij}=(B\bar{g}_{ij})+
\bar{g}_{im}\partial_{j}v^{m}+\bar{g}_{jm}\partial_{i}v^{m}
\end{equation}
Since $\bar{g}_{ij}=\delta_{ij}$, we have
\begin{equation}
 2\alpha k_{ij}=\partial_{j}v^{i}+\partial_{i}v^{j}=-2\partial_{i}\psi_{j}
\end{equation}
The last equality above is in the irrotational isentropic case.

Any vector $X$ tangent to $C_{u}$ at a point $p$ can be uniquely decomposed into a vector colinear to $L$ and a vector tangent
to $S_{t,u}$.
\begin{equation}
 X=c_{X}L+\Pi X
\end{equation}
If $X$, $Y$ are tangent to $C_{u}$, then $g(X,Y)=\slashed{g}(\Pi X,\Pi Y)$.

Let $X_{A}=\frac{\partial}{\partial\vartheta^{A}}$, given any $Z\in T_{p}C_{u}$, we can expand $\Pi Z=Z^{A}X_{A}$.
Taking inner product with $X_{B}$:
\begin{equation}
 \slashed{g}(X_{B},\Pi Z)=\slashed{g}_{AB}Z^{A}
\end{equation}
Define $\chi$ by
\begin{equation}
 2\chi=\slashed{\mathcal{L}}_{L}g
\end{equation}
Here, $\slashed{\mathcal{L}}$ is the restriction of $\mathcal{L}$ to $S_{t,u}$.

If $X$, $Y$ are tangent to $S_{t,u}$, then we have
\begin{equation}
 \chi(X,Y)=g(D_{X}L,Y)=g(D_{Y}L,X)
\end{equation}
Since $\hat{L}=\mu^{-1}L$ is an affinely parametrized geodesic field, we have:
\begin{equation}
 D_{L}L=\mu D_{\hat{L}}(\mu\hat{L})=L(\mu)\hat{L}=\mu^{-1}(L\mu)L
\end{equation}
Denoting $\chi_{AB}=\chi(X_{A},X_{B})$.
We shall derive a propagation equation for $\chi_{AB}$ along $C_{u}$. First we have
\begin{equation}
 L\chi_{AB}=g(D_{L}D_{{X}_{A}}L,X_{B})+g(D_{{X}_{A}}L,D_{L}X_{B})
\end{equation}
Now $D_{L}X_{A}-D_{X_{A}}L=[L,X_{A}]=0$, if we denote by $R$ the curvature of $g$, we have: 
\begin{align*}
D_{L}D_{{X}_{A}}L-D_{X_{A}}D_{L}L=R(L,X_{A})L
\end{align*}
where the curvature transformation is defined by:
\begin{align*}
 R(X,Y)Z=D_{X}D_{Y}Z-D_{Y}D_{X}Z-D_{[X,Y]}Z
\end{align*}
Then
\begin{equation}
 g(D_{L}D_{{X}_{A}}L,X_{B})=\mu^{-1}(L\mu)\chi_{AB}-\alpha_{AB}
\end{equation}
where 
$\alpha_{AB}=R(X_{A},L,X_{B},L)$.
Recall that the curvature tensor is defined through the curvature transformation as follows:
\begin{equation}
 R(W,Z,X,Y):=g(W,R(X,Y)Z)
\end{equation}
The second term on the right hand side of (3.11) is equal to $g(D_{X_{A}}L,D_{X_{B}}L)$.

Since $L$ is null with respect to $g$, then for any vectorfield $X$, $D_{X}L$ is $g$-orthogonal to 
$L$, therefore tangential to $C_{u}$.
Thus, setting $W_{A}=D_{X_{A}}L$: $A=1,2$, we have
\begin{equation}
 g(W_{A},W_{B})=\slashed{g}(\Pi\textsl{W}_{A},\Pi\textsl{W}_{B})
\end{equation}
We now expand $\Pi W_{A}$ in the basis $X_{B},B=1,2$.The coefficient of $X_{B}$ in this expansion is 
\begin{equation}
 (\slashed{g}^{-1})^{BC}\slashed{g}(\Pi W_{A},X_{C})=(\slashed{g}^{-1})^{BC}g(D_{X_{A}}L,X_{C})
=(\slashed{g}^{-1})^{BC}\chi_{AC}
\end{equation}
hence, 
\begin{align*}
\slashed{g}(\Pi W_{A},\Pi W_{B})=\slashed{g}(\chi_{A}^{C}X_{C},\chi_{B}^{D}X_{D})=\chi_{A}^{C}\chi_{BC}
\end{align*}
Here the capital indices are raised and lowered with respect to $\slashed{g}_{AB}$. Thus,
\begin{align*}
 \chi_{A}^{B}=(\slashed{g}^{-1})^{BC}\chi_{AC}
\end{align*}

So we get the propagation equation:
\begin{equation}
 L\chi_{AB}=\mu^{-1}(L\mu)\chi_{AB}+\chi_{A}^{C}\chi_{BC}-\alpha_{AB}
\end{equation}
From this, we can get a propagation equation for $\textrm{tr}\chi$.
We denote by $S_{\mu\nu}$ the Ricci tensor:
\begin{equation}
 S_{\mu\nu}:=(g^{-1})^{\kappa\lambda}R_{\mu\kappa\nu\lambda}
\end{equation}
We can express $(g^{-1})^{\kappa\lambda}$ in terms of the frame $L,T,X_{1},X_{2}$.
We have
\begin{equation}
 (g^{-1})^{\mu\nu}=-\alpha^{-2}L^{\mu}L^{\nu}-\mu^{-1}(L^{\mu}T^{\nu}+L^{\nu}T^{\mu})
+(\slashed{g}^{-1})^{AB}X^{\mu}_{A}X^{\nu}_{B}
\end{equation}

From (3.18), we have:
\begin{equation}
 \textrm{tr}\alpha=(\slashed{g}^{-1})^{AB}R_{\mu\kappa\nu\lambda}X^{\mu}_{A}L^{\kappa}X^{\nu}_{B}L^{\lambda}=
(g^{-1})^{\mu\nu}R_{\mu\kappa\nu\lambda}L^{\kappa}L^{\lambda}=S_{\kappa\lambda}L^{\kappa}L^{\lambda}
\end{equation}
On the other hand,
\begin{equation}
 L(\textrm{tr}\chi)=L[(\slashed{g}^{-1})^{AB}\chi_{AB}]=(\slashed{g}^{-1})^{AB}L(\chi_{AB})
+\chi_{AB}L[(\slashed{g}^{-1})^{AB}]
\end{equation}
This equals 
\begin{equation}
 (\slashed{g}^{-1})^{AB}L(\chi_{AB})-\chi_{AB}(\slashed{g}^{-1})^{AC}(\slashed{g}^{-1})^{BD}L(\slashed{g}_{CD})
\end{equation}
i.e.
\begin{equation}
 (\slashed{g}^{-1})^{AB}L(\chi_{AB})-2\chi_{AB}(\slashed{g}^{-1})^{AC}(\slashed{g}^{-1})^{BD}\chi_{CD}
\end{equation}
where we have used the following fact:
\begin{align*}
L[(\slashed{g}^{-1})^{AB}]=-(\slashed{g}^{-1})^{AC}(\slashed{g}^{-1})^{BD}L(\slashed{g}_{CD})
\end{align*}
This is simply the derivative of the reciprocal of a non-degenerate matrix: 
\begin{align*}
\frac{dM^{-1}}{dt}=-M^{-1}\frac{dM}{dt}M^{-1}
\end{align*}
Since we have
\begin{equation}
 2\chi_{AB}=L(\slashed{g}_{AB})
\end{equation}
then taking trace with (3.16), we get
\begin{equation}
 L(\textrm{tr}\chi)=\mu^{-1}(L\mu)\textrm{tr}\chi-|\chi|^{2}_{\slashed{g}}-S(L,L)
\end{equation}
We shall write down the Gauss and Codazzi equations of the embedding of $S_{t,u}$ in the acoustical spacetime.
First, we consider $S_{t,u}$ as a submanifold of $\Sigma_{t}$. In this case, the normal vectorfield is $T$.

We define $\theta$, the second fundamental form of $S_{t,u}$ relative to $\Sigma_{t}$, by:
\begin{equation}
 2\kappa\theta:=\slashed{\mathcal{L}}_{T}\bar{g}
\end{equation}
For any two vectors tangent to $S_{t,u}$ at a point, we have 
\begin{equation}
 \kappa\theta(X,Y)=\bar{g}(\bar{D}_{X}T,Y)=\bar{g}(\bar{D}_{Y}T,X)
\end{equation}
Denote by $\slashed{k}$, the restriction of $k$, the second fundamental form of $\Sigma_{t}$ in spacetime manifold, to $TS_{t,u}$,
then by (2.53), we have
\begin{equation}
 \chi=\alpha(\slashed{k}-\theta)
\end{equation}
The Gauss equation of $S_{t,u}$ in $\Sigma_{t}$ is expressed in terms of $(X_{1},X_{2})$ by
\begin{equation}
 \slashed{R}_{ABCD}-\theta_{AC}\theta_{BD}+\theta_{AD}\theta_{BC}=0
\end{equation}
$\slashed{R}$ is the curvature tensor of $\slashed{g}$. Since $\textrm{dim}S_{t,u}=2$, we have
\begin{equation}
 \slashed{R}_{ABCD}=K(\slashed{g}_{AC}\slashed{g}_{BD}-\slashed{g}_{AD}\slashed{g}_{BC})
\end{equation}
where $K$ is the Gauss curvature of $S_{t,u}$.
Then we have
\begin{equation}
 \theta_{AC}\theta_{BD}-\theta_{AD}\theta_{BC}=K(\slashed{g}_{AC}\slashed{g}_{BD}-\slashed{g}_{AD}\slashed{g}_{BC})
\end{equation}
contracting this with $\frac{1}{2}(\slashed{g}^{-1})^{AC}(\slashed{g}^{-1})^{BD}$:
\begin{equation}
 \frac{1}{2}(\textrm{tr}_{\slashed{g}}\theta)^{2}-\frac{1}{2}|\theta|^{2}_{\slashed{g}}=K
\end{equation}
Also from the Gauss equations of the embedding of $\Sigma_{t}$ in the whole manifold, we have
\begin{equation}
 \slashed{k}_{AC}\slashed{k}_{BD}-\slashed{k}_{BC}\slashed{k}_{AD}=R_{ABCD}
\end{equation}
In view of the symmetries of the curvature tensor we can express:
\begin{equation}
 R_{ABCD}=\rho\epsilon_{AB}\epsilon_{CD}
\end{equation}
where $\epsilon_{AB}$ are the components of the area form of $(S_{t,u},\slashed{g})$:
\begin{equation}
 \epsilon_{AB}=\sqrt{\det \slashed{g}}[AB]
\end{equation}
Contracting (3.33) with $\frac{1}{2}(\slashed{g}^{-1})^{AC}(\slashed{g}^{-1})^{BD}$ and taking into account the fact that
\begin{equation}
 (\slashed{g}^{-1})^{BD}\epsilon_{AB}\epsilon_{CD}=\slashed{g}_{AC}
\end{equation}
we obtain
\begin{equation}
 \frac{1}{2}(\textrm{tr}_{\slashed{g}}\slashed{k})^{2}-\frac{1}{2}|\slashed{k}|^{2}_{\slashed{g}}=\rho
\end{equation}

This is the Gauss equation of $S_{t,u}$ in acoustical spacetime.

We then derive the Codazzi equation:
Let $X,Y,Z$ be vectorfields tangent to $S_{t,u}$, we can extend them to $C_{u}$ by pushing forward by the flow of $L$. They are then
tangential to each of the sections $S_{t,u}$, $t\in [0,t_{*\varepsilon_{0}})$. 
 We have:
\begin{equation}
 X(\chi(Y,Z))=X(g(D_{Y}L,Z))=g(D_{X}D_{Y}L,Z)+g(D_{Y}L,D_{X}Z)
\end{equation}
On the other hand,
\begin{equation}
 (\slashed{D}_{X}\chi)(Y,Z)=X(\chi(Y,Z))-\chi(\slashed{D}_{X}Y,Z)-\chi(Y,\slashed{D}_{X}Z)
\end{equation}
Here $\slashed{D}$ is the covariant derivative operator on the $S_{t,u}$ associated to the induced metric $\slashed{g}$.
Hence,
\begin{equation}
 (\slashed{D}_{X}\chi)(Y,Z)=g(D_{X}D_{Y}L,Z)+g(D_{Y}L,D_{X}Z)-\chi(\slashed{D}_{X}Y,Z)-\chi(Y,\slashed{D}_{X}Z)
\end{equation}
Exchanging $X,Y$, and noting that 
\begin{equation}
 D_{X}D_{Y}L-D_{Y}D_{X}L=R(X,Y)L+D_{[X,Y]}L
\end{equation}
while $g(R(X,Y)L,Z)=R(Z,L,X,Y)$ and $\slashed{D}_{X}Y-\slashed{D}_{Y}X=[X,Y]$,we obtain
\begin{align*}
(\slashed{D}_{X}\chi)(Y,Z)-(\slashed{D}_{Y}\chi)(X,Z)=
R(Z,L,X,Y)+g(D_{[X,Y]}L,Z)+g(D_{Y}L,D_{X}Z)-g(D_{X}L,D_{Y}Z)\\
-\chi(\slashed{D}_{X}Y,Z)+\chi(\slashed{D}_{Y}X,Z)-\chi(Y,\slashed{D}_{X}Z)+\chi(X,\slashed{D}_{Y}Z)\\
=R(Z,L,X,Y)+g(D_{Y}L,D_{X}Z)-g(D_{X}L,D_{Y}Z)-\chi(Y,\slashed{D}_{X}Z)+\chi(X,\slashed{D}_{Y}Z)
\end{align*}

Define the 1-form $\zeta$ on $S_{t,u}$, by
\begin{equation}
 \zeta(X)=g(D_{X}L,T)
\end{equation}
for any vector $X$ tangent to $S_{t,u}$ at some point. We shall derive an expression for $\zeta$.
For this, we set $D_{X}L=xL+yT+z^{A}X_{A}$.
Then
\begin{equation}
 g(D_{X}L,T)=-{\mu}x+\kappa^{2}y, \quad g(D_{X}L,L)=-{\mu}y,\quad 
g(D_{X}L,X_{A})=z^{B}\slashed{g}_{AB}
\end{equation}
Solving for the coefficients $x,y,z^{A}$, we obtain:
\begin{equation}
 D_{X}L=-\mu^{-1}\zeta(X)L+\chi\cdot X
\end{equation}
where $\slashed{g}(\chi\cdot X,Y)=\chi(X,Y)$.
Then, 
\begin{equation}
 g(D_{X}L,D_{Y}Z)=\mu^{-1}\zeta(X)\chi(Y,Z)+g(\chi\cdot X,D_{Y}Z)=\mu^{-1}\zeta(X)\chi(Y,Z)+\chi(X,\slashed{D}_{Y}Z)
\end{equation}
Substituting this into the expression above, we get 
\begin{equation}
 (\slashed{D}_{X}\chi)(Y,Z)-(\slashed{D}_{Y}\chi)(X,Z)=R(Z,L,X,Y)-\mu^{-1}(\zeta(X)\chi(Y,Z)-\zeta(Y)\chi(X,Z))
\end{equation}
In this equation, we set $X=X_{A}, Y=X_{B}, Z=X_{C}$. Defining the 1-form $\beta$ on the $S_{t,u}$ by:
\begin{equation}
 R(X_{C},L,X_{A},X_{B})=\beta_{C}\epsilon_{AB}
\end{equation}
we obtain the Codazzi equation:
\begin{equation}
 \slashed{D}_{A}\chi_{BC}-\slashed{D}_{B}\chi_{AC}=\beta_{C}\epsilon_{AB}-\mu^{-1}(\zeta_{A}\chi_{BC}-\zeta_{B}\chi_{AC})
\end{equation}
Denoting by $\slashed{\textrm{curl}}\chi$ the 1-form on $S_{t,u}$ with components:
\begin{equation}
 \slashed{\textrm{curl}}\chi_{C}:=\frac{1}{2}\epsilon^{AB}(\slashed{D}_{A}\chi_{BC}-\slashed{D}_{B}\chi_{AC})
\end{equation}
where $\epsilon^{AB}=(\slashed{g}^{-1})^{AC}(\slashed{g}^{-1})^{BD}\epsilon_{CD}=(\sqrt{\det\slashed{g}})^{-1}[AB]$,
 we can write Codazzi equation as follows:
\begin{equation}
 \slashed{\textrm{curl}}\chi=\beta-\mu^{-1}\zeta\wedge\chi
\end{equation}
Here $(\zeta\wedge\chi)_{C}=\frac{1}{2}\epsilon^{AB}(\zeta_{A}\chi_{BC}-\zeta_{B}\chi_{AC})$, and this equation has the whole content of equations (3.47). 
Contracting (3.47) with $(\slashed{g}^{-1})^{AC}$, we obtain the contracted form of the Codazzi equation.
\begin{equation}
 \slashed{\textrm{div}}\chi-\slashed{d}\textrm{tr}\chi=\beta^{*}-\mu^{-1}(\zeta\cdot\chi-\zeta\textrm{tr}\chi)
\end{equation}
which is equivalent to the equation (3.49). Here $(\zeta\cdot\chi)_{B}=(\slashed{g}^{-1})^{AC}\zeta_{A}\chi_{BC}, 
\beta^{*}_{B}=(\slashed{g}^{-1})^{AC}\beta_{C}\epsilon_{AB}$

We shall now derive an expression for $zeta$, from the expression for $L$,
\begin{equation}
 \zeta(X)=g(D_{X}B,T)-\alpha\kappa^{-1}g(D_{X}T,T)-g(T,T)X(\alpha\kappa^{-1})
\end{equation}
While $g(T,T)=\kappa^{2}$, we have $g(D_{X}T,T)=\kappa X\kappa$, so this reduces to
\begin{equation}
 \zeta(X)=g(D_{X}B,T)-\kappa X\alpha
\end{equation}
Define 1-form $\varepsilon$ on $S_{t,u}$ by: $\kappa\varepsilon(X)=k(X,T)$. Then
\begin{equation}
 \zeta=\kappa(\alpha\varepsilon-\slashed{d}\alpha)
\end{equation}

We now define $S_{t,u}$ 1-form:
\begin{equation}
 \eta:=\zeta+\slashed{d}\mu
\end{equation}
Obviously, we have 
\begin{equation}
 \eta(X)=-g(D_{X}T,L)
\end{equation}
We can express the commutator $\Lambda$ of $L$ and $T$, a vectorfield tangential to $S_{t,u}$, in terms of $\zeta$ and $\eta$. We have
\begin{equation}
 \Lambda=\Lambda^{A}X_{A}, \Lambda^{A}\slashed{g}_{AB}=g(X_{B},D_{L}T)-g(X_{B},D_{T}L)
\end{equation}
and
\begin{equation}
 g(X_{B},D_{L}T)=-g(D_{L}X_{B},T)=-g(D_{X_{B}}L,T)=-\zeta_{B}
\end{equation}
To compute $g(X_{B},D_{T}L)$, we express $L$ in terms of the geodesic vectorfield $\hat{L}$.
We have 
\begin{equation}
 g(X_{B},D_{T}L)=g(X_{B},D_{T}(\mu\hat{L}))=\mu g(X_{B},D_{T}\hat{L})=\mu g(T,D_{X_{B}}\hat{L})=\mu g(T,D_{X_{B}}(\mu^{-1}L))
\end{equation}
Here we have used the fact that $\hat{L}^{\mu}=-(g^{-1})^{\mu\nu}\partial_{\nu}u$, that is, $\hat{L}$ is gradient of a function, so we can exchange
$T$ and $X_{B}$. This equals 
\begin{equation}
 g(T,D_{X_{B}}L)-\mu^{-1}(X_{B}\mu)g(T,L)=\zeta_{B}+X_{B}\mu=\eta_{B}
\end{equation}
So we get
\begin{equation}
 \Lambda^{A}=-(\slashed{g}^{-1})^{AB}(\zeta_{B}+\eta_{B})
\end{equation}
We shall express $D_{T}L$ and $D_{L}T$ in the frames $(L,T,X_{1},X_{2})$
First, by (3.59) and (3.60),
\begin{equation}
 g(D_{T}L,X_{B})=\eta_{B}
\end{equation}
and
\begin{equation}
 g(D_{T}L,T)=g(D_{L}T,T)=L(\frac{1}{2}\kappa^{2}), \quad g(D_{T}L,L)=0
\end{equation}
It follows that
\begin{equation}
 D_{T}L=(\slashed{g}^{-1})^{AB}\eta_{B}X_{A}-\alpha^{-1}(L\kappa)L
\end{equation}
Also,
\begin{equation}
 g(D_{L}T,X_{B})=-g(D_{L}X_{B},T)=-g(D_{X_{B}}L,T)=-\zeta_{B}
\end{equation}
and
\begin{equation}
 g(D_{L}T,L)=g(D_{T}L,L)=0,\quad g(D_{L}T,T)=\frac{1}{2}L(\kappa^{2})
\end{equation}
where we have used the fact that $[L,T]$ is tangential to $S_{t,u}$.
We thus obtain
\begin{equation}
 D_{L}T=-\slashed{g}^{AB}\zeta_{B}X_{A}-\alpha^{-1}(L\kappa)L
\end{equation}

Next, we compute $D_{X_{A}}T$. Set
\begin{equation}
 D_{X_{A}}T=a_{A}L+b_{A}T+c^{B}_{A}X_{B}
\end{equation}
Taking inner product with $L$:
\begin{equation}
 -\mu{b_{A}}=g(D_{X_{A}}T,L)=-\eta_{A} \Rightarrow b_{A}=\mu^{-1}\eta_{A}
\end{equation}
Taking inner product with $T$:
\begin{equation}
 -\mu{a_{A}}+\kappa^{2}b_{A}=g(D_{X_{A}}T,T)=\frac{1}{2}X_{A}(\kappa^{2})
\end{equation}
Thus we obtain:
\begin{equation}
 a_{A}=\alpha^{-2}\eta_{A}-\alpha^{-1}X_{A}(\kappa)=\alpha^{-2}(\eta_{A}-X_{A}(\mu)+\kappa X_{A}(\alpha))
=\alpha^{-2}(\zeta_{A}+\kappa X_{A}(\alpha))=\alpha^{-1}\kappa\varepsilon_{A}
\end{equation}
Finally, taking inner product with $X_{C}$:
\begin{equation}
 c^{B}_{A}\slashed{g}_{BC}=g(D_{X_{A}}T,X_{C})=\kappa\theta_{AC}
\end{equation}
We conclude that
\begin{equation}
 D_{X_{A}}T=\alpha^{-1}\kappa\varepsilon_{A}L+\mu^{-1}\eta_{A}T+\kappa\theta_{AB}(\slashed{g}^{-1})^{BC}X_{C}
\end{equation}
In the following, we will compute $D_{T}T$. We may decompose:
\begin{equation}
 D_{T}T=\bar{D}_{T}T+aB
\end{equation}
where $\bar{D}$ is the covariant derivative operator on the $\Sigma_{t}$, associated to the induced metric $\bar{g}$, which, as we have seen, coincides 
with the Euclidean metric.
Taking inner product with $B$:
\begin{equation}
 -\alpha^{2}a=g(D_{T}T,B)=g(D_{T}T,L)+\alpha\kappa^{-1}g(D_{T}T,T)=-T\mu-g(T,D_{T}L)+\alpha\kappa^{-1}T(\frac{1}{2}\kappa^{2})
\end{equation}
where we have used the fact that $g(D_{T}T,L)=T(g(T,L))-g(T,D_{T}L)$.
The righthand side of (3.74) equals:
\begin{equation}
 -T\mu-\kappa L\kappa+\alpha T\kappa=-\kappa(T\alpha+L\kappa)
\end{equation}
Using the fact that $[L,T]$ is tangential to $S_{t,u}$, so
$g(T,D_{T}L)=g(T,D_{L}T)=\frac{1}{2}L(\kappa^{2})$

We may use the rectangular coordinates on $\Sigma_{t}$, the induced metric $\bar{g}$ being Euclidean. For any pair of vectorfields $X,Y$ tangential to $\Sigma_{t}$
we then have, in terms of components in rectangular coordinates,
\begin{align}
 (\bar{D}_{X}Y)^{i}=X^{j}\partial_{j}Y^{i}
\end{align}
 So we have:
\begin{equation}
 (\bar{D}_{T}T)^{i}=T^{j}\bar{D}_{j}T^{i}=T^{j}\partial_{j}T^{i}
\end{equation}
and 
\begin{equation}
 \bar{g}_{ij}T^{j}=\kappa^{2}\partial_{i}u, \quad \bar{g}_{ij}=\delta_{ij} \Rightarrow T^{i}=\kappa^{2}\partial_{i}u
\end{equation}
So,
\begin{equation}
 (\bar{D}_{T}T)^{i}=\kappa^{2}\partial_{j}u\partial_{j}(\kappa^{2}\partial_{i}u)=\kappa^{2}\partial_{j}u\partial_{j}\kappa^{2}\partial_{i}u
+\kappa^{4}\partial_{j}u\partial_{j}\partial_{i}u
\end{equation}
By (3.78), $\sum_{i}(\partial_{i}u)^{2}=\kappa^{-2}$ then we have
\begin{equation}
 (\bar{D}_{T}T)^{i}=\kappa^{-2}T^{i}T^{j}\partial_{j}(\kappa^{2})+\frac{1}{2}\kappa^{4}\partial_{i}(\kappa^{-2})
\end{equation}
 This equals to
\begin{equation}
 \kappa^{-2}T^{i}T^{j}\partial_{j}(\kappa^{2})-\frac{1}{2}\partial_{i}(\kappa^{2})
\end{equation}
\begin{equation}
 \Rightarrow\bar{D}_{T}T=\frac{1}{2}\kappa^{-2}T(\kappa^{2})T-\frac{1}{2}\slashed{g}^{AB}X_{B}(\kappa^{2})X_{A}
\end{equation}
In view of (3.74), (3.75) and (3.82), we conclude that:
\begin{align*}
 D_{T}T=\kappa\alpha^{-2}(T\alpha+L\kappa)L+[\alpha^{-1}(T\alpha+L\kappa)+\kappa^{-1}T\kappa]T-\frac{1}{2}(\slashed{g}^{-1})^{AB}X_{B}(\kappa^{2})X_{A}
\end{align*}

In regard to (3.74),
\begin{equation}
 g(D_{T}T,B)=-g(T,D_{T}B)=-\alpha k(T,T)
\end{equation}
Comparing with (3.75), we see that:
\begin{equation}
 L(\frac{1}{2}\kappa^{2})=-\kappa T\alpha+\alpha k(T,T)
\end{equation}
Using this equation, we shall derive a propagation equation along the generator of $C_{u}$ for $\mu$:
\begin{equation}
 L\mu=\kappa L\alpha+\alpha L\kappa=\kappa L\alpha-\alpha T\alpha+\alpha\mu k(\hat{T},\hat{T})
\end{equation}
First, 
\begin{equation}
 \alpha\mu k(\hat{T},\hat{T})=-\mu\hat{T}^{i}\hat{T}^{j}\partial_{i}\psi_{j}=-\alpha\hat{T}^{i}(T\psi_{i})
\end{equation}
We have used the fact that $\alpha k_{ij}=-\partial_{i}\psi_{j}$ (See (3.5)).

To compute the first two terms on the right hand side of (3.85), we shall use: $\alpha d\alpha=\frac{1}{2}d\alpha^{2}=\frac{1}{2}d\eta^{2}$. 
While
\begin{equation}
 \eta^{2}=\frac{\rho}{{\rho}^{\prime}} \Rightarrow d\eta^{2}={(\frac{\rho}{{\rho}^{\prime}})}^{\prime}dh
\end{equation}
It follows that 
\begin{equation}
 \kappa L\alpha-\alpha T\alpha=\frac{1}{2}\frac{\mu}{{\alpha}^{2}}{(\frac{\rho}{{\rho}^{\prime}})}^{\prime}Lh-\frac{1}{2}{(\frac{\rho}{{\rho}^{\prime}})}^{\prime}Th
\end{equation}
Substituting this and (3.86) into (3.85), we get
\begin{equation}
 L\mu=\frac{1}{2}\frac{\mu}{{\alpha}^{2}}{(\frac{\rho}{{\rho}^{\prime}})}^{\prime}Lh-\frac{1}{2}{(\frac{\rho}{{\rho}^{\prime}})}^{\prime}Th-\alpha\hat{T}^{i}(T\psi_{i})
\end{equation}
Recalling that 
\begin{equation}
 L=\partial_{t}-(\alpha\hat{T}^{i}+\psi_{i})\partial_{i} \Rightarrow -\alpha\hat{T}^{i}=L^{i}+\psi_{i}
\end{equation}
we express
\begin{equation}
 -\alpha\hat{T}^{i}(T\psi_{i})=L^{i}(T\psi_{i})+\sum_{i}\psi_{i}T\psi_{i}=T^{i}(L\psi_{i})-Th=\alpha^{-1}\mu{\hat{T}}^{i}(L\psi_{i})-Th
\end{equation}
We thus obtain the following propagation equation for $\mu$:
\begin{equation}
 L\mu=m+\mu{e}
\end{equation}
where 
\begin{equation}
 m=\frac{1}{2}\frac{dH}{dh}Th, \quad H=-2h-\eta^{2}
\end{equation}
and 
\begin{equation}
 e=\frac{1}{2{\alpha}^{2}}{(\frac{\rho}{{\rho}^{\prime}})}^{\prime}Lh+\alpha^{-1}{\hat{T}}^{i}(L\psi_{i})
\end{equation}

To complete the set of connection coefficients of the frame field $(L,T,X_{1},X_{2})$, we should compute $D_{{X}_{A}}{{X}_{B}}$.
We decompose:
\begin{equation}
 D_{{X}_{A}}{{X}_{B}}=\slashed{D}_{{X}_{A}}{{X}_{B}}+a_{AB}L+b_{AB}T
\end{equation}
then
\begin{equation}
 -\mu{b_{AB}}=g(D_{{X}_{A}}{{X}_{B}},L)=-\chi_{AB}
\end{equation}
and
\begin{equation}
 -\mu{a_{AB}}+\kappa^{2}b_{AB}=g(T, D_{{X}_{A}}{{X}_{B}})=-\kappa\theta_{AB}
\end{equation}
By $\chi=\alpha(\slashed{k}-\theta)$, we get
\begin{equation}
 b_{AB}=\mu^{-1}\chi_{AB}, a_{AB}=\alpha^{-1}\slashed{k}_{AB}
\end{equation}
So we get the following table:
\begin{align}
D_{L}L=\mu^{-1}(L\mu)L\\
D_{T}L=(\slashed{g}^{-1})^{AB}\eta_{B}X_{A}-\alpha^{-1}(L\kappa)L\\
D_{L}T=-(\slashed{g}^{-1})^{AB}\zeta_{B}X_{A}-\alpha^{-1}(L\kappa)L\\
D_{X_{A}}L=-\mu^{-1}\zeta_{A}L+\chi_{A}^{B}X_{B}\\
D_{T}T=\alpha^{-2}\kappa(T\alpha+L\kappa)L+(\alpha^{-1}L\kappa+\mu^{-1}T\mu)T-
\frac{1}{2}(\slashed{g}^{-1})^{AB}X_{B}(\kappa^{2})X_{A}\\
D_{X_{A}}T=\alpha^{-1}\kappa\varepsilon_{A}L+\mu^{-1}\eta_{A}T+\kappa\theta_{AB}(\slashed{g}^{-1})^{BC}X_{C}\\
D_{L}X_{A}=D_{X_{A}}L\\
D_{{X}_{A}}{{X}_{B}}=\slashed{D}_{{X}_{A}}{{X}_{B}}+\alpha^{-1}\slashed{k}_{AB}L+\mu^{-1}\chi_{AB}T
\end{align}
We now investigate the connection between the Lie derivative of $\chi$ with respect to $T$ and the derivative of $\eta$ tangential to $S_{t,u}$.

First, we extend $\chi$ from $TS_{t,u}$ for each $t,u$ to $TW^{*}_{\epsilon_{0}}$ by the conditions:
\begin{equation}
 \chi(X,L)=\chi(X,T)=0
\end{equation}
 We define $\slashed{\mathcal{L}}_{T}\chi$ to be the restriction of $\mathcal{L}_{T}\chi$ to $S_{t,u}$:
\begin{equation}
 (\slashed{\mathcal{L}}_{T}\chi)(X_{A},X_{B})=(D_{T}\chi)(X_{A},X_{B})+\chi(X_{A},D_{X_{B}}T)+\chi(X_{B},D_{X_{A}}T)
\end{equation}
Due to the extension of $\chi$ and using (3.99)-(3.106), this equals to
\begin{equation}
 (D_{T}\chi)(X_{A},X_{B})+\kappa\theta_{BC}\chi_{A}^{C}+\kappa\theta_{AC}\chi_{B}^{C}
\end{equation}
To compute $(D_{T}\chi)(X_{A},X_{B})$:
\begin{equation}
 (D_{T}\chi)(X_{A},X_{B})=T(\chi_{AB})-\chi(\Pi D_{T}X_{A},X_{B})-\chi(\Pi D_{T}X_{B},X_{A})
\end{equation}
where $\Pi$ is the $g$-projection to $S_{t,u}$.
Now 
\begin{equation}
 T(\chi_{AB})=T(g(D_{X_{A}}L,X_{B}))=g(D_{T}D_{X_{A}}L,X_{B})+g(D_{X_{A}}L,D_{T}X_{B})
\end{equation}
While
\begin{equation}
 g(D_{T}D_{X_{A}}L,X_{B})=g(D_{X_{A}}D_{T}L,X_{B})+g(D_{[T,X_{A}]}L,X_{B})+R(X_{B},L,T,X_{A})
\end{equation}
From (3.100),
\begin{equation}
 g(D_{X_{A}}D_{T}L,X_{B})=g(D_{X_{A}}(\eta^{C}X_{C}-\alpha^{-1}(L\kappa)L),X_{B})
\end{equation}
This equals to
\begin{equation}
 -\alpha^{-1}(L\kappa)\chi_{AB}+X_{A}(\eta_{B})-\eta(\slashed{D}_{X_{A}}X_{B})=-\alpha^{-1}(L\kappa)\chi_{AB}+\slashed{D}_{A}\eta_{B}
\end{equation}
Since $[T,X_{A}]$ is tangential to $S_{t,u}$, we have
\begin{align}
 g(D_{[T,X_{A}]}L,X_{B})=\chi([T,X_{A}],X_{B})=\chi(\Pi D_{T}X_{A},X_{B})-\chi(\Pi D_{X_{A}}T,X_{B})\\
=\chi(\Pi D_{T}X_{A},X_{B})-\kappa\theta_{AC}\chi_{B}^{C}
\end{align}
Also,
\begin{equation}
 g(D_{X_{A}}L,D_{T}X_{B})=-\mu^{-1}\zeta_{A}g(L,D_{T}X_{B})+g(\chi_{AC}X^{C},\Pi D_{T}X_{B})
\end{equation}
which equals 
\begin{equation}
 \mu^{-1}\zeta_{A}\eta_{B}+\chi(X_{A},\Pi D_{T}X_{B})
\end{equation}
So we obtain:
\begin{align}
 (D_{T}\chi)(X_{A},X_{B})=\slashed{D}_{A}\eta_{B}+\mu^{-1}\zeta_{A}\eta_{B}-\alpha^{-1}(L\kappa)\chi_{AB}-\kappa\theta_{AC}\chi_{B}^{C}-R(X_{A},T,X_{B},L)
\end{align}
The symmetric part of this equation is
\begin{align}
(D_{T}\chi)(X_{A},X_{B})=\frac{1}{2}(\slashed{D}_{A}\eta_{B}+\slashed{D}_{B}\eta_{A})
+\frac{1}{2}\mu^{-1}(\zeta_{A}\eta_{B}+\zeta_{B}\eta_{A})\\-\alpha^{-1}(L\kappa)\chi_{AB}-\frac{1}{2}\kappa(\theta_{AC}\chi_{B}^{C}+\theta_{BC}\chi_{A}^{C})\\
-\frac{1}{2}(R(X_{A},T,X_{B},L)+R(X_{B},T,X_{A},L))
\end{align}
The antisymmetric part of this equation contains no new information.

Denoting
\begin{equation}
 (\slashed{\mathcal{L}}_{T}\chi)(X_{A},X_{B})=\slashed{\mathcal{L}}_{T}\chi_{AB}
\end{equation}
our conclusion is
\begin{align}
 \slashed{\mathcal{L}}_{T}\chi_{AB}=\frac{1}{2}(\slashed{D}_{A}\eta_{B}+\slashed{D}_{B}\eta_{A})
+\frac{1}{2}\mu^{-1}(\zeta_{A}\eta_{B}+\zeta_{B}\eta_{A})\\-\alpha^{-1}(L\kappa)\chi_{AB}+\frac{1}{2}\kappa(\theta_{AC}\chi_{B}^{C}+
\theta_{BC}\chi_{A}^{C})-\gamma_{AB}
\end{align}
where 
\begin{equation}
 \gamma_{AB}=\frac{1}{2}(R(X_{A},T,X_{B},L)+R(X_{B},T,X_{A},L))
\end{equation}

We shall also use the null frame $(L,\underline{L},X_{1},X_{2})$, where 
\begin{equation}
 \underline{L}=\alpha^{-1}\kappa L+2T
\end{equation}
One can check that $\underline{L}$ is a incoming future directed null vector satisfying the following condition:
\begin{equation}
 g(L,\underline{L})=-2\mu
\end{equation}
 The inverse of the acoustical metric is expressed in this frame by
\begin{equation}
 (g^{-1})^{\mu\nu}=-\frac{1}{2\mu}(L^{\mu}\underline{L}^{\nu}+L^{\nu}\underline{L}^{\mu})
+(\slashed{g}^{-1})^{AB}X^{\mu}_{A}X^{\nu}_{B}
\end{equation}
From (3.99)-(3.106), we have the following table:
\begin{align}
D_{L}L=\mu^{-1}(L\mu)L\\
D_{\underline{L}}L=-L(\alpha^{-1}\kappa)L+2\eta^{A}X_{A}\\
D_{X_{A}}L=-\mu^{-1}\zeta_{A}L+\chi_{A}^{B}X_{B}\\
D_{L}\underline{L}=-2\zeta^{A}X_{A}\\
D_{\underline{L}}\underline{L}=(\mu^{-1}\underline{L}\mu+L(\alpha^{-1}\kappa))\underline{L}-2\mu(\slashed{g}^{-1})^{AB}X_{B}(\alpha^{-1}\kappa)X_{A}\\
D_{X_{A}}\underline{L}=\mu^{-1}\eta_{A}\underline{L}+\underline{\chi}_{A}^{B}X_{B}\\
D_{L}X_{A}=D_{X_{A}}L\\
D_{{X}_{A}}{{X}_{B}}=\slashed{D}_{{X}_{A}}{{X}_{B}}+\frac{1}{2}\mu^{-1}\underline{\chi}_{AB}L+\frac{1}{2}\mu^{-1}\chi_{AB}\underline{L}
\end{align}
Similarly, we have 
\begin{equation}
 \underline{\chi}=\kappa(\slashed{k}+\theta)
\end{equation}
This completes the exposition of the acoustical structure equations.

We shall presently derive expressions for the operators $\Box_{g}$ and $\Box_{\tilde{g}}$. we have :
\begin{equation}
 \Box_{g}f=\textrm{tr}(D^{2}f)
\end{equation}
where $\textrm{tr}$ denotes the trace with respect to $g$ and $D^{2}f$ denotes the Hessian of $f$ with respect to $g$.
In view of the expression (3.129) for $g^{-1}$ we have
\begin{equation}
 \Box_{g}f=(g^{-1})^{\mu\nu}(D^{2}f)_{\mu\nu}=-\mu^{-1}(D^{2}f)_{\underline{L}L}+(\slashed{g}^{-1})^{AB}(D^{2}f)_{AB}
\end{equation}
Also we can consider the Hessian of the restriction of $f$ to $S_{t,u}$ with respect to the induced metric $\slashed{g}$:
\begin{equation}
 \slashed{D}^{2}f
\end{equation}
then the operator $\slashed{\Delta}=\Delta_{\slashed{g}}$ is given by
\begin{equation}
 \slashed{\Delta}=\textrm{tr}(\slashed{D}^{2}f)
\end{equation}
Here $\textrm{tr}$ is the trace with respect to $\slashed{g}$. In terms of the frame, we have 
\begin{equation}
 (D^{2}f)_{AB}=X_{A}(X_{B}f)-(D_{X_{A}}X_{B})f
\end{equation}
and 
\begin{equation}
 (\slashed{D}^{2}f)_{AB}=X_{A}(X_{B}f)-(\slashed{D}_{X_{A}}X_{B})f
\end{equation}
By (3.137), we have 
\begin{equation}
 (D^{2}f)_{AB}=(\slashed{D}^{2}f)_{AB}-\frac{1}{2}\mu^{-1}\underline{\chi}_{AB}(Lf)-\frac{1}{2}\mu^{-1}\chi_{AB}(\underline{L}f)
\end{equation}
Hence:
\begin{equation}
 (\slashed{g}^{-1})^{AB}(D^{2}f)_{AB}=\slashed{\Delta}f-\frac{1}{2}\mu^{-1}\textrm{tr}\underline{\chi}(Lf)
-\frac{1}{2}\mu^{-1}\textrm{tr}\chi(\underline{L}f)
\end{equation}
Also, by (3.133),
\begin{equation}
(D^{2}f)_{\underline{L}L}=L(\underline{L}f)-(D_{L}\underline{L})f=L(\underline{L}f)+2\zeta\cdot\slashed{d}f
\end{equation}
So we get the following expressions for $\Box_{g}$:
\begin{equation}
 \Box_{g}f=\slashed{\Delta}f-\frac{1}{2}\mu^{-1}\textrm{tr}\underline{\chi}(Lf)
-\frac{1}{2}\mu^{-1}\textrm{tr}\chi(\underline{L}f)-\mu^{-1}L(\underline{L}f)-2\mu^{-1}\zeta\cdot\slashed{d}f
\end{equation}
Since the conformal acoustical metric $\tilde{g}$ is given by:
\begin{equation}
 \tilde{g}_{\mu\nu}=\Omega g_{\mu\nu}
\end{equation}
 we then get
\begin{equation}
 \Box_{\tilde{g}}=\Omega^{-1}\Box_{g}f+\Omega^{-2}(g^{-1})^{\mu\nu}\partial_{\mu}\Omega\partial_{\nu}f
=\Omega^{-1}\Box_{g}f+\Omega^{-2}\frac{d\Omega}{dh}(g^{-1})^{\mu\nu}\partial_{\mu}h\partial_{\nu}f
\end{equation}
In view of (3.129), we have
\begin{equation}
 (g^{-1})^{\mu\nu}\partial_{\mu}h\partial_{\nu}f=-\frac{1}{2}\mu^{-1}(\underline{L}h)(Lf)-\frac{1}{2}\mu^{-1}(\underline{L}f)(Lh)+\slashed{d}h\cdot\slashed{d}f
\end{equation}
Thus, setting
\begin{equation}
 \nu=\frac{1}{2}(\textrm{tr}\chi+\Omega^{-1}\frac{d\Omega}{dh}(Lh))
\end{equation}
\begin{equation}
 \underline{\nu}=\frac{1}{2}(\textrm{tr}\underline{\chi}+\Omega^{-1}\frac{d\Omega}{dh}(\underline{L}h))
\end{equation}
we obtain the formula
\begin{equation}
 \Omega\Box_{\tilde{g}}f=\slashed{\Delta}f-\mu^{-1}L(\underline{L}f)-\mu^{-1}(\nu\underline{L}f+\underline{\nu}Lf)
-2\mu^{-1}\zeta\cdot\slashed{d}f+\Omega^{-1}\frac{d\Omega}{dh}\slashed{d}h\cdot\slashed{d}f
\end{equation}

\section{The Derivatives of the Rectangular Components\\
 of $L$ and $\hat{T}$}
We shall now draw the conclusion that the derivatives of $L^{i}$, $\hat{T}^{i}$ , the rectangular components of $L$ and $\hat{T}$, with respect to $L$, $T$ and $X_{A}$
are regular as $\mu\rightarrow0$.
In the following we shall denote by $\nabla$ the covariant differentiation with respect to the natural connection associated to the affine space $\mathcal{G}$, the Galilean
spacetime. The $t$ and $x^{i} : i=1,2,3$ being linear coordinates on $\mathcal{G}$, $\nabla$ amounts simply to partial differentiation in these coordinates. We have
\begin{align}
 D_{\mu}W^{\nu}=\nabla_{\mu}W^{\nu}+\Gamma_{\mu\lambda}^{\nu}W^{\lambda}\\
\Gamma_{\mu\lambda}^{\nu}=(g^{-1})^{\nu\kappa}\Gamma_{\mu\lambda\kappa}\\
\Gamma_{\mu\lambda\kappa}=\frac{1}{2}(\partial_{\mu}g_{\lambda\kappa}+\partial_{\lambda}g_{\mu\kappa}-\partial_{\kappa}g_{\mu\lambda})
\end{align}
In the irrotational and isentropic case, recall that $v^{i}=-\partial_{i}\phi$, then the acoustical metric can be written as follows
\begin{equation}
 g=-\eta^{2}dt^{2}+\sum_{i}({dx}^{i}+\partial_{i}\phi{dt})^{2}
\end{equation}
By direct calculation, we have the Christoffel symbols:
\begin{align}
 \Gamma_{000}=\frac{1}{2}\partial_{t}(-\eta^{2}+|\textbf{v}|^{2})\\
\Gamma_{0i0}=\frac{1}{2}\partial_{i}(-\eta^{2}+|\textbf{v}|^{2})\\
\Gamma_{ij0}=\partial_{i}\partial_{j}\phi\\
\Gamma_{00k}=\partial_{0}\partial_{k}\phi-\frac{1}{2}\partial_{k}(-\eta^{2}+|\textbf{v}|^{2})\\
\Gamma_{i0k}=\Gamma_{ijk}=0
\end{align}
Consider now the vectorfield 
\begin{equation}
 L(L^{\mu})\partial_{\mu}
\end{equation}
Since $L(L^{0})=0$, it can be expanded as 
\begin{equation}
 a_{L}\hat{T}+b_{L}
\end{equation}
where $b_{L}\in TS_{t,u}$.
On the other hand, (3.165) is also $\nabla_{L}L$.
That is $\nabla_{L}L= a_{L}\hat{T}+b_{L}$.

Taking inner product with $\hat{T}$: 
\begin{equation}
 a_{L}=g(\nabla_{L}L, \hat{T})
\end{equation}
Writing 
\begin{equation}
 b_{L}=b_{L}^{A}X_{A}
\end{equation}
and taking inner product with $X_{B}$:
\begin{equation}
 g(\nabla_{L}L,X_{B})=b_{L}^{A}\slashed{g}_{AB}
\end{equation}
While
\begin{equation}
 g(\nabla_{L}L,\hat{T})=g(D_{L}L,\hat{T})-\Gamma_{\alpha\beta\nu}L^{\alpha}L^{\beta}\hat{T}^{\nu}
\end{equation}
From (3.92) and (3.99):
\begin{equation}
 g(D_{L}L,\hat{T})=-\kappa^{-1}(L\mu)=-\kappa^{-1}m-\alpha{e}
\end{equation}
Also, 
\begin{equation}
 \Gamma_{\alpha\beta\nu}L^{\alpha}L^{\beta}\hat{T}^{\nu}=\Gamma_{00k}\hat{T}^{k}
=\hat{T}(\psi_{0})-\frac{1}{2}\hat{T}(-\eta^{2}+|\textbf{v}|^{2})=-\frac{1}{2}\frac{dH}{dh}\hat{T}(h)
\end{equation}
\begin{equation}
 \Rightarrow g(\nabla_{L}L,\hat{T})=-\alpha{e}
\end{equation}
The $\kappa^{-1}$ term cancels!
From (3.94), we have
\begin{equation}
 a_{L}=-\alpha e=-\frac{1}{2\alpha}L(\eta^{2})-\hat{T}^{i}(L\psi_{i})
\end{equation}
Since $g(D_{L}L,X_{B})=0$, we have
\begin{equation}
 g(\nabla_{L}L,X_{B})=-\Gamma_{\alpha\beta\nu}L^{\alpha}L^{\beta}X_{B}^{\nu}=-\Gamma_{00k}X_{B}^{k}=
-X_{B}(\psi_{0})+\frac{1}{2}X_{B}(-\eta^{2}+|\textbf{v}|^{2})=\frac{1}{2}\frac{dH}{dh}X_{B}(h)
\end{equation}
So 
\begin{equation}
 b_{L}^{A}=(\slashed{g}^{-1})^{AB}g(\nabla_{L}L,X_{B})=\frac{1}{2}\frac{dH}{dh}X^{A}(h)
\end{equation}
Consider next the vectorfield $L(\hat{T}^{\mu})\partial_{\mu}$.
Also it can be expanded as
\begin{equation}
 \nabla_{L}\hat{T}=p_{L}\hat{T}+q_{L}
\end{equation}
Obviously,
\begin{equation}
 p_{L}=g(\nabla_{L}\hat{T},\hat{T})
\end{equation}
Writing 
\begin{equation}
 q_{L}=q_{L}^{A}X_{A}
\end{equation}
and taking inner product with $X_{B}$:
\begin{equation}
 q_{L}^{A}\slashed{g}_{AB}=g(\nabla_{L}\hat{T},X_{B})
\end{equation}
Since $g(D_{L}\hat{T},\hat{T})=\frac{1}{2}L(g(\hat{T},\hat{T}))=0$, we have
\begin{equation}
 g(\nabla_{L}\hat{T},\hat{T})=-\Gamma_{\alpha\beta\nu}L^{\alpha}\hat{T}^{\beta}\hat{T}^{\nu}=0 \Rightarrow p_{L}=0
\end{equation}
Also we have:
\begin{equation}
 g(\nabla_{L}\hat{T},X_{B})=g(D_{L}\hat{T},X_{B})-\Gamma_{\alpha\beta\nu}L^{\alpha}\hat{T}^{\beta}X^{\nu}_{B}=
g(D_{L}\hat{T},X_{B})=\kappa^{-1}g(D_{L}T,X_{B})=-\kappa^{-1}\zeta^{A}\slashed{g}_{AB}
\end{equation}
\begin{equation}
 \Rightarrow q_{L}^{A}=-\kappa^{-1}\zeta^{A}=-\alpha\varepsilon^{A}+X^{A}(\alpha)=-\alpha k_{ij}X^{Ai}\hat{T}^{j}+X^{A}(\alpha)=-\hat{T}^{j}(X^{A}\psi_{j})+X^{A}(\alpha)
\end{equation}
This is regular as $\mu\rightarrow0$.
Consider the vectorfield 
\begin{equation}
 T(L^{\mu})\partial_{\mu}
\end{equation}
Also, it can be expanded as 
\begin{equation}
 \nabla_{T}L=a_{T}\hat{T}+b_{T}
\end{equation}
\begin{equation}
 \Rightarrow  a_{T}=g(\nabla_{T}L,\hat{T}), \quad b_{T}=b_{T}^{A}X_{A},\quad
 b_{T}^{A}\slashed{g}_{AB}=g(\nabla_{T}L,X_{B})
\end{equation}

We have: $g(\nabla_{T}L,\hat{T})=g(D_{T}L,\hat{T})-\Gamma_{\alpha\beta\nu}T^{\alpha}L^{\beta}\hat{T}^{\nu}$.
From (3.100), we have 
\begin{equation}
 g(D_{T}L,\hat{T})=g(-\alpha^{-1}(L\kappa)L,\hat{T})=L\kappa
\end{equation}
While $-\Gamma_{\alpha\beta\nu}T^{\alpha}L^{\beta}\hat{T}^{\nu}=0$, we get
\begin{equation}
 a_{T}=L\kappa
\end{equation}
Also
\begin{equation}
 g(\nabla_{T}L,X_{B})=g(D_{T}L,X_{B})-\Gamma_{\alpha\beta\nu}T^{\alpha}L^{\beta}X_{B}^{\nu}=g(D_{T}L,X_{B})=\eta_{B}
\end{equation}
\begin{equation}
 \Rightarrow b_{T}^{A}=\eta^{A}
\end{equation}
Consider finally the vectorfield
\begin{equation}
 T(\hat{T}^{\mu})\partial_{\mu}
\end{equation}
It can be expanded as 
\begin{equation}
 \nabla_{T}\hat{T}=p_{T}\hat{T}+q_{T}
\end{equation}
\begin{equation}
 \Rightarrow p_{T}=g(\nabla_{T}\hat{T},\hat{T})=g(D_{T}\hat{T},\hat{T})-\Gamma_{\alpha\beta\nu}T^{\alpha}\hat{T}^{\beta}\hat{T}^{\nu}=0
\end{equation}
Writing $q_{T}=q_{T}^{A}X_{A}$, we have
\begin{equation}
 q_{T}^{A}\slashed{g}_{AB}=g(\nabla_{T}\hat{T},X_{B})=g(D_{T}\hat{T},X_{B})-\Gamma_{\alpha\beta\nu}T^{\alpha}\hat{T}^{\beta}X^{\nu}_{B}
\end{equation}
This equals to
\begin{equation}
 g(D_{T}\hat{T},X_{B})=\kappa^{-1}g(D_{T}T,X_{B})=-\frac{1}{2\kappa}(\slashed{g}^{-1})^{AC}X_{C}(\kappa^{2})\slashed{g}_{AB}=-X_{B}(\kappa)
\end{equation}
We proceed to derive expressions for 
\begin{equation}
 \nabla_{X_{A}}L=X_{A}(L^{\mu})\partial_{\mu}
\end{equation}
and
\begin{equation}
 \nabla_{X_{A}}\hat{T}=X_{A}(\hat{T}^{\mu})\partial_{\mu}
\end{equation}
First, we expand
\begin{equation}
 \nabla_{X_{A}}L=\slashed{a}_{A}\hat{T}+\slashed{b}_{A}
\end{equation}
then we have
\begin{equation}
 \slashed{a}_{A}=g(\nabla_{X_{A}}L,\hat{T})=g(D_{X_{A}}L,\hat{T})-\Gamma_{\alpha\beta\nu}X_{A}^{\alpha}L^{\beta}\hat{T}^{\nu}=g(D_{X_{A}}L,\hat{T})=\kappa^{-1}\zeta_{A}
\end{equation}
This is regular as $\mu\rightarrow0$ and we have used (3.102).

On the other hand
\begin{equation}
\slashed{g}_{BC}\slashed{b}^{B}_{A}=g(\nabla_{X_{A}}L,X_{C})=g(D_{X_{A}}L,X_{C})-\Gamma_{\alpha\beta\nu}X^{\alpha}_{A}L^{\beta}X^{\nu}_{C}=g(D_{X_{A}}L,X_{C})=\chi_{AC}
\end{equation}
Finally
\begin{equation}
 \nabla_{X_{A}}\hat{T}=\slashed{p}_{A}\hat{T}+\slashed{q}_{A}
\end{equation}
\begin{equation}
 \Rightarrow \slashed{p}_{A}=g(\nabla_{X_{A}}\hat{T},\hat{T})=g(D_{X_{A}}\hat{T},\hat{T})-\Gamma_{\alpha\beta\nu}X^{\alpha}_{A}\hat{T}^{\beta}\hat{T}^{\nu}=0
\end{equation}
\begin{equation}
 \slashed{q}^{B}_{A}\slashed{g}_{BC}=
g(\nabla_{X_{A}}\hat{T},X_{C})=g(D_{X_{A}}\hat{T},X_{C})-\Gamma_{\alpha\beta\nu}X^{\alpha}_{A}\hat{T}^{\beta}X^{\nu}_{C}=\kappa^{-1}g(D_{X_{A}}T,X_{C})=\theta_{AC}
\end{equation}

\chapter{The Acoustical Curvature}




We presently derive expressions for the curvature tensor of the acoustical metric $g$.

In this chapter, we need the following expressions of Christoffel symbols derived in Chapter 3:
\begin{align}
 \Gamma_{000}=\frac{1}{2}\partial_{0}(-\eta^{2}+|\textbf{v}|^{2})\\\notag
 \Gamma_{0i0}=\frac{1}{2}\partial_{i}(-\eta^{2}+|\textbf{v}|^{2})\\\notag
 \Gamma_{ij0}=\partial_{i}\partial_{j}\phi\\\notag
 \Gamma_{00k}=\frac{1}{2}(2\partial_{0}\partial_{k}\phi-\partial_{k}(-\eta^{2}+|\textbf{v}|^{2}))\\\notag
 \Gamma_{i0k}=\Gamma_{ijk}=0
\end{align}
\section{Expressions for Curvature Tensor}
We have the expressions for the curvature tensor:
\begin{equation}
 R_{\mu\nu\alpha\beta}=R^{(2)}_{\mu\nu\alpha\beta}+R^{(1)}_{\mu\nu\alpha\beta}
\end{equation}
where 
\begin{equation}
 R^{(2)}_{\mu\nu\alpha\beta}=\frac{1}{2}(\partial_{\alpha}\partial_{\nu}g_{\beta\mu}
+\partial_{\beta}\partial_{\mu}g_{\alpha\nu}-\partial_{\alpha}\partial_{\mu}g_{\beta\nu}
-\partial_{\beta}\partial_{\nu}g_{\alpha\mu})
\end{equation}
\begin{equation}
 R^{(1)}_{\mu\nu\alpha\beta}=-(g^{-1})^{\kappa\lambda}(\Gamma_{\alpha\mu\kappa}\Gamma_{\beta\nu\lambda}-\Gamma_{\beta\mu\kappa}
\Gamma_{\alpha\nu\lambda})
\end{equation}
From these expressions and those for the Christoffel symbols, we have
\begin{equation}
 R_{ijkl}=R_{ijkl}^{(1)}=\eta^{-2}(\partial_{i}\partial_{k}\phi\partial_{j}\partial_{l}\phi-\partial_{i}\partial_{l}\phi\partial_{j}\partial_{k}\phi)
\end{equation}
Also, we have:
\begin{equation}
 R_{0i0j}=R_{0i0j}^{(2)}+R_{0i0j}^{(1)}
\end{equation}
where:
\begin{equation}
 R^{(2)}_{0i0j}=\frac{1}{2}(\partial_{0}\partial_{i}g_{0j}+\partial_{0}\partial_{j}g_{0i}-\partial^{2}_{0}g_{ij}-\partial_{i}\partial_{j}g_{00})=
\partial^{2}_{ij}(\partial_{0}\phi-\frac{1}{2}\sum_{k}(\partial_{k}\phi)^{2}+\frac{1}{2}\eta^{2})
\end{equation}
i.e.
\begin{equation}
 R^{(2)}_{0i0j}=\partial^{2}_{ij}h+\frac{1}{2}\partial^{2}_{ij}\eta^{2}=-\frac{1}{2}\partial^{2}_{ij}H=-\frac{1}{2}\frac{dH}{dh}\partial^{2}_{ij}h
-\frac{1}{2}\frac{d^{2}H}{dh^{2}}\partial_{i}h\partial_{j}h
\end{equation}
So the real principal part of $R_{0i0j}$ is $-\frac{1}{2}\frac{dH}{dh}\partial^{2}_{ij}h$. We can write this in terms of $D^{2}_{ij}h$ we have:
\begin{equation}
 R_{0i0j}=R^{[P]}_{0i0j}+R^{[N]}_{0i0j}
\end{equation}
where
\begin{equation}
 R^{[P]}_{0i0j}=-\frac{1}{2}\frac{dH}{dh}D^{2}_{ij}h
\end{equation}
Since
\begin{align*}
 R_{0i0j}^{(2)}=-\frac{1}{2}\frac{dH}{dh}D^{2}_{ij}h-\frac{1}{2}\frac{dH}{dh}\Gamma_{ij}^{\alpha}\partial_{\alpha}h
-\frac{1}{2}\frac{d^{2}H}{dh^{2}}\partial_{i}h\partial_{j}h
\end{align*}
we have:
\begin{equation}
 R^{[N]}_{0i0j}=-\frac{1}{2}\frac{d^{2}H}{dh^{2}}\partial_{i}h\partial_{j}h
-\frac{1}{2}\frac{dH}{dh}\Gamma^{\alpha}_{ij}\partial_{\alpha}h+R^{(1)}_{0i0j}
\end{equation}
where
\begin{align}
 R^{(1)}_{0i0j}=\eta^{-2}[\frac{1}{2}\partial_{0}(-\eta^2+|\textbf{v}|^2)\partial_{i}\partial_{j}\phi-
\frac{1}{4}\partial_{i}(-\eta^2+|\textbf{v}|^2)\partial_{j}(-\eta^2+|\textbf{v}|^2)]\\\notag
-\eta^{-2}\partial_{k}\phi(\partial_{i}\partial_{j}\phi)\frac{1}{2}(2\partial_{0}\partial_{k}\phi-\partial_{k}(-\eta^2+|\textbf{v}|^2))
\end{align}

Finally, we have the component 
\begin{equation}
 R_{0ijk}=R^{(2)}_{0ijk}+R^{(1)}_{0ijk}
\end{equation}
and
\begin{equation}
 R^{(2)}_{0ijk}=\frac{1}{2}(\partial^{2}_{ij}g_{0k}+\partial^{2}_{0k}g_{ij}-\partial^{2}_{0j}g_{ki}-\partial^{2}_{ki}g_{0j})=0
\end{equation}
while
\begin{equation}
 R^{(1)}_{0ijk}=-(g^{-1})^{\sigma\lambda}(\Gamma_{j0\sigma}\Gamma_{ik\lambda}-\Gamma_{k0\sigma}\Gamma_{ij\lambda})
\end{equation}
Hence
\begin{equation}
 R_{0ijk}=R^{(1)}_{0ijk}=\eta^{-2}[\frac{1}{2}\partial_{j}(-\eta^{2}+|\textbf{v}|^{2})\partial^{2}_{ik}\phi
-\frac{1}{2}\partial_{k}(-\eta^{2}+|\textbf{v}|^{2})\partial^{2}_{ij}\phi]
\end{equation}
Now we have calculated all the components of the curvature tensor, and we shall express all the components in a more systematic way. To do this, we 
need to introduce the following notations:
\begin{align*}
 \tau_{\mu}=\partial_{\mu}h\\
v_{\mu\nu}=D_{\mu}\tau_{\nu}=(D^{2}h)_{\mu\nu}\\
w_{\mu\nu}=\partial_{\mu}\psi_{\nu}=w_{\nu\mu}
\end{align*}
First, since only the component $R_{0i0j}$ has the principal part, we have:
\begin{align}
 R^{[P]}_{\mu\nu\alpha\beta}=\frac{1}{2}\frac{dH}{dh}(\delta^{0}_{\mu}\delta^{0}_{\beta}v_{\alpha\nu}+\delta^{0}_{\alpha}\delta^{0}_{\nu}v_{\mu\beta}
-\delta^{0}_{\nu}\delta^{0}_{\beta}v_{\alpha\mu}-\delta^{0}_{\alpha}\delta^{0}_{\mu}v_{\nu\beta})
\end{align}
The lower order terms are much more complicated. First from (4.5):
\begin{align}
 R^{[N]}_{ijkl}=\eta^{-2}(w_{ik}w_{jl}-w_{il}w_{jk})
\end{align}
From (4.11) and (4.12):
\begin{align*}
 R^{[N]}_{0i0j}=-\frac{1}{2}\frac{d^{2}H}{dh^{2}}\tau_{i}\tau_{j}-\frac{1}{2}\frac{dH}{dh}(\Gamma^{0}_{ij}\tau_{0}+\Gamma^{k}_{ij}\tau_{k})\\
+\eta^{-2}[\frac{1}{2}\partial_{0}(-\eta^2+|\textbf{v}|^2)w_{ij}-
\frac{1}{4}\partial_{i}(-\eta^2+|\textbf{v}|^2)\partial_{j}(-\eta^2+|\textbf{v}|^2)]\\\notag
-\eta^{-2}\tau_{k}(w_{ij})\frac{1}{2}(2w_{0k}-\partial_{k}(-\eta^2+|\textbf{v}|^2))
\end{align*}
From (4.1):
\begin{align*}
 \Gamma_{ij}^{0}=(g^{-1})^{00}\Gamma_{ij0}=-\eta^{-2}\omega_{ij}\\
\Gamma_{ij}^{k}=(g^{-1})^{k0}\Gamma_{ij0}=-\eta^{-2}v^{k}\omega_{ij}
\end{align*}
Also, since 
\begin{align*}
 h=\psi_{0}-\frac{1}{2}|\textbf{v}|^{2}\quad\textrm{and}\quad H=-2h-\eta^{2}
\end{align*}
we have:
\begin{align}
 -\eta^{2}+|\textbf{v}|^{2}=H+2\psi_{0}
\end{align}
Substituting in the expression for $R^{[N]}_{0i0j}$, we obtain, in view of the fact that $\psi_{k}=-v^{k}$,
\begin{align}
 R^{[N]}_{0i0j}=-\frac{1}{2}[\frac{d^{2}H}{dh^{2}}+\frac{1}{2}\eta^{-2}(\frac{dH}{dh})^{2}]\tau_{i}\tau_{j}\\\notag
+\frac{1}{2}\eta^{-2}\frac{dH}{dh}(2\tau_{0}w_{ij}-\tau_{i}w_{j0}-\tau_{j}w_{i0})\\\notag
+\eta^{-2}(w_{00}w_{ij}-w_{i0}w_{j0})
\end{align}
(4.16) reads:
\begin{align*}
 R^{[N]}_{0ijk}=\frac{1}{2}\eta^{-2}\{\partial_{j}(-\eta^{2}+|\textbf{v}|^{2})w_{ik}-\partial_{k}(-\eta^{2}+|\textbf{v}|^{2})w_{ij}\}
\end{align*}
Substituting gives:
\begin{align}
 R^{[N]}_{0ijk}=\frac{1}{2}\eta^{-2}\frac{dH}{dh}(\tau_{j}w_{ik}-\tau_{k}w_{ij})+\eta^{-2}(w_{j0}w_{ik}-w_{k0}w_{ij})
\end{align}
Finally, the results for $R^{[N]}_{ijkl}, R^{[N]}_{0ijk}$ and $R^{[N]}_{0i0j}$, combine to produce the formula for $R^{[N]}_{\mu\nu\alpha\beta}$:
\begin{align}
 R^{[N]}_{\mu\nu\alpha\beta}=\eta^{-2}A_{\mu\nu\alpha\beta}+\frac{1}{2}\eta^{-2}\frac{dH}{dh}C_{\mu\nu\alpha\beta}-\frac{1}{2}H_{2}B_{\mu\nu\alpha\beta}
\end{align}
where:
\begin{align*}
 A_{\mu\nu\alpha\beta}=w_{\mu\alpha}w_{\nu\beta}=w_{\mu\beta}w_{\nu\alpha}\\
B_{\mu\nu\alpha\beta}=(\delta^{0}_{\mu}\tau_{\nu}-\delta^{0}_{\nu}\tau_{\mu})(\delta^{0}_{\alpha}\tau_{\beta}-\delta^{0}_{\beta}\tau_{\alpha})
\end{align*}
and, with:
\begin{align*}
 \xi_{\mu\alpha\beta}=\delta^{0}_{\alpha}w_{\beta\mu}-\delta^{0}_{\beta}w_{\alpha\mu}
\end{align*}
we have:
\begin{align*}
 C_{\mu\nu\alpha\beta}=\tau_{\mu}\xi_{\nu\alpha\beta}-\tau_{\nu}\xi_{\mu\alpha\beta}+\tau_{\alpha}\xi_{\beta\mu\nu}
-\tau_{\beta}\xi_{\alpha\mu\nu}
\end{align*}
Also, the function $H_{2}$ is defined by:
\begin{align*}
 H_{2}=\frac{d^{2}H}{dh^{2}}+\frac{1}{2}\eta^{-2}(\frac{dH}{dh})^{2}
\end{align*}
The expressions $A_{\mu\nu\alpha\beta}, B_{\mu\nu\alpha\beta}, C_{\mu\nu\alpha\beta}$ all possess the algebraic properties of the curvature 
tensor, namely the antisymmetriy in the first as well as the second pair of indices, and the cyclic property, which imply the symmetry under exchange of the first
with the second pair of indices. In particular the cyclic property of $C_{\mu\nu\alpha\beta}$ follows from that of $\xi_{\mu\alpha\beta}$:
\begin{align*}
 \xi_{\mu\alpha\beta}+\xi_{\alpha\beta\mu}+\xi_{\beta\mu\alpha}=0
\end{align*}

\section{Regularity for the Acoustical Structure Equations as $\mu\rightarrow0$}
First we compute $\alpha_{AB}$, which appears in the propagation equation of $\chi_{AB}$.
By definition, 
\begin{align*}
 \alpha_{AB}=R_{\mu\nu\alpha\beta}X^{\mu}_{A}L^{\nu}X^{\alpha}_{B}L^{\beta}=
R_{i0j0}X^{i}_{A}X^{j}_{B}+R_{ijkl}X^{i}_{A}L^{j}X^{k}_{B}L^{l}
+R_{i0jk}X^{i}_{A}X^{j}_{B}L^{k}+R_{i0jk}X^{j}_{A}X^{i}_{B}L^{k}
\end{align*}

From (3.106), the principal part of $\alpha_{AB}$ is 
\begin{align}
 \alpha_{AB}^{[P]}=R^{[P]}_{0i0j}X^{i}_{A}X^{j}_{B}=
-\frac{1}{2}\frac{dH}{dh}D^{2}h(X_{A},X_{B})=
-\frac{1}{2}\frac{dH}{dh}[X_{A}(X_{B}h)-(D_{{X}_{A}}X_{B})h]\\\notag
=-\frac{1}{2}\frac{dH}{dh}[\slashed{D}^{2}h(X_{A},X_{B})-\alpha^{-1}\slashed{k}_{AB}Lh-\mu^{-1}\chi_{AB}Th]
\end{align}
So the singular term as $\mu\rightarrow0$ in $\alpha_{AB}^{[P]}$ is 
\begin{equation}
 \frac{1}{2\mu}\frac{dH}{dh}\chi_{AB}T(h)
\end{equation}

For the lower order term, we have:
\begin{align}
 \alpha^{[N]}_{AB}=\eta^{-2}\alpha^{[A]}_{AB}+\frac{1}{2}\eta^{-2}\frac{dH}{dh}\alpha^{[C]}_{AB}-\frac{1}{2}H_{2}\alpha^{[B]}_{AB}
\end{align}
and:
\begin{align}
 \alpha^{[A]}_{AB}=(w_{\mu\alpha}w_{\nu\beta}-w_{\mu\beta}w_{\nu\alpha})X^{\mu}_{A}L^{\nu}X_{B}^{\alpha}L^{\beta}
=\slashed{w}_{AB}w_{LL}-(\slashed{w}_{L})_{A}(\slashed{w}_{L})_{B}\\\notag
\alpha^{[B]}_{AB}=(\delta^{0}_{\mu}\tau_{\nu}-\delta^{0}_{\nu}\tau_{\mu})(\delta^{0}_{\alpha}\tau_{\beta}-\delta^{0}_{\beta}\tau_{\alpha})
X^{\mu}_{A}L^{\nu}X^{\alpha}_{B}L^{\beta}=\slashed{\tau}_{A}\slashed{\tau}_{B}\\\notag
\alpha^{[C]}_{AB}=(\tau_{\mu}\xi_{\nu\alpha\beta}-\tau_{\nu}\xi_{\mu\alpha\beta}+\tau_{\alpha}\xi_{\beta\mu\nu}-\tau_{\beta}\xi_{\alpha\mu\nu})
X^{\mu}_{A}L^{\nu}X_{B}^{\alpha}L^{\beta}
=2\tau_{L}\slashed{w}_{AB}-\slashed{\tau}_{A}(\slashed{w}_{L})_{B}-\slashed{\tau}_{B}(\slashed{w}_{L})_{A}
\end{align}
So $\alpha^{[C]}_{AB}$ is regular as $\mu\rightarrow0$. The only singular term in $\alpha_{AB}$ is (4.24).

Recall the propagation equation for $\chi$ in chapter 3:
\begin{equation}
 L\chi_{AB}=\mu^{-1}(L\mu)\chi_{AB}+\chi_{A}^{C}\chi_{BC}-\alpha_{AB}=e\chi_{AB}+\mu^{-1}m\chi_{AB}+\chi_{A}^{C}\chi_{BC}-\alpha_{AB}
\end{equation}
We have used the propagation equation for $\mu$
\begin{equation}
 L\mu=m+\mu e
\end{equation}
where 
\begin{equation}
 m=\frac{1}{2}\frac{dH}{dh}T(h),\quad 
e=\frac{1}{2\eta}\frac{d(\eta^{2})}{dh}L(h)+\frac{1}{\eta}\hat{T}^{i}(L\psi_{i})
\end{equation}
Then we know that the right hand side of (4.27) is regular as $\mu\longrightarrow0$.

Next, we will analyze the regularity of the right hand side of Codazzi equation and Gauss equation.

Recall the Codazzi equation (3.51):
\begin{equation}
 (\slashed{\textrm{div}}\chi)_{A}-\slashed{d}_{A}\textrm{tr}\chi=\beta^{*}_{A}-\mu^{-1}(\zeta_{B}\chi_{A}^{B}-\zeta_{A}\textrm{tr}\chi)
\end{equation}
Here $\beta^{*}_{B}=(\slashed{g}^{-1})^{AC}\beta_{C}\epsilon_{AB}$

Recalling (3.54) and the definition of $k_{ij}$, we know that $\mu^{-1}\zeta_{A}$ is regular as $\mu\longrightarrow0$. Therefore the regularity
of the right hand side of (4.35) reduces to the regularity of
\begin{equation}
 R(X_{C},L,X_{A},X_{B})=\beta_{C}\epsilon_{AB}
\end{equation}
We have:
\begin{align}
 \epsilon_{AB}\beta^{[P]}_{C}=R^{[P]}_{\mu\nu\alpha\beta}X^{\mu}_{C}L^{\nu}X^{\alpha}_{A}X^{\beta}_{B}=0
\end{align}
and:
\begin{align}
 \beta^{[N]}_{A}=\eta^{-2}\beta^{[A]}_{A}+\frac{1}{2}\eta^{-2}\frac{dH}{dh}\beta^{[C]}_{A}-\frac{1}{2}H_{2}\beta^{[B]}_{A}
\end{align}
where:
\begin{align*}
 \epsilon_{AB}\beta^{[A]}_{C}=\slashed{w}_{CA}(\slashed{w}_{L})_{B}-\slashed{w}_{CB}(\slashed{w}_{L})_{A}\\
\epsilon_{AB}\beta^{[B]}_{C}=0\\
\epsilon_{AB}\beta^{[C]}_{C}=-\slashed{\tau}_{A}\slashed{w}_{BC}+\slashed{\tau}_{B}\slashed{w}_{AC}
\end{align*}
So $\beta_{A}$ is regular as $\mu\rightarrow0$.

We turn to Gauss equation (3.36). Since
\begin{equation}
  R_{ABCD}=\rho\epsilon_{AB}\epsilon_{CD}
\end{equation}
it suffice to prove that $R_{ABCD}$ is regular as $\mu\longrightarrow0$. We have:
\begin{align}
 R^{[P]}_{ABCD}=R^{[P]}_{\mu\nu\alpha\beta}X^{\mu}_{A}X^{\nu}_{B}X^{\alpha}_{C}X^{\beta}_{D}=0
\end{align}
and
\begin{align}
 R^{[N]}_{ABCD}=-\eta^{-2}R^{[A]}_{ABCD}+\frac{1}{2}\eta^{-2}\frac{dH}{dh}R^{[C]}_{ABCD}-\frac{1}{2}H_{2}R^{[B]}_{ABCD}
\end{align}
where
\begin{align*}
 R^{[A]}_{ABCD}=\slashed{w}_{AC}\slashed{w}_{BD}-\slashed{w}_{AD}\slashed{w}_{BC}\\
R^{[B]}_{ABCD}=R^{[C]}_{ABCD}=0
\end{align*}
$\rho$ is regular as $\mu\rightarrow0$.

Finally, we shall discuss about the equation (3.125)-(3.126). As we have seen, $\mu^{-1}\zeta$ is regular as $\mu\longrightarrow0$. We shall thus discuss about 
$\gamma_{AB}$. We have:
\begin{align}
 \gamma^{[P]}_{AB}=\frac{1}{2}(R^{[P]}(X_{A},T,X_{B},L)+R^{[P]}(X_{B},T,X_{A},L))
\end{align}
and
\begin{align}
 \gamma^{[N]}_{AB}=\eta^{-2}\gamma^{[A]}_{AB}+\frac{1}{2}\eta^{-2}\frac{dH}{dh}\gamma^{[C]}_{AB}-\frac{1}{2}H_{2}\gamma^{[B]}_{AB}
\end{align}
where:
\begin{align*}
 \gamma^{[A]}_{AB}=\frac{\kappa}{2}(2\slashed{w}_{AB}w_{L\hat{T}}-(\slashed{w}_{L})_{A}(\slashed{w}_{\hat{T}})_{B}
-(\slashed{w}_{L})_{B}(\slashed{w}_{\hat{T}})_{A})\\
\gamma^{[B]}_{AB}=0\\
\gamma^{[C]}_{AB}=\frac{1}{2}(2\tau_{T}\slashed{w}_{AB}-\kappa\slashed{\tau}_{A}(\slashed{w}_{\hat{T}})_{B}-\kappa\slashed{\tau}_{B}(\slashed{w}_{\hat{T}})_{A})
\end{align*}
$\gamma_{AB}$ is regular as $\mu\rightarrow0$.

Note that all the components of $\tau$ in the $(L,T,X_{1},X_{2})$ frame i.e. $\tau_{L}, \tau_{T}, \slashed{\tau}_{A}$ are regular as $\mu\rightarrow0$.
Also all the components of $w$ in the $(L,\hat{T},X_{1},X_{2})$ frame i.e. $w_{LL}, w_{L\hat{T}}, (\slashed{w}_{\hat{T}})_{A}, \slashed{w}_{AB}$, except 
$w_{\hat{T}\hat{T}}$ are regular as $\mu\rightarrow0$. The last is $\kappa^{-1}w_{T\hat{T}}$ and $w_{T\hat{T}}$ is regular as $\mu\rightarrow0$.

\section{A Remark}
     We end this chapter with a remark. By (4.15), we can see that if $H=const$, then the principal part of the acoustical curvature vanishes. 
This fact can be also obtained by the geometric interpretation of $H$ in Chapter 2.
By the Gauss equation:
\begin{align}
 R_{\alpha\beta\mu\nu}-\tilde{k}_{\alpha\mu}\tilde{k}_{\beta\nu}+\tilde{k}_{\beta\mu}\tilde{k}_{\alpha\nu}=\tilde{R}_{\alpha\beta\mu\nu}
\end{align}
 where $R$ and $\tilde{R}$ are the curvature tensor of $(\mathbb{R}\times\mathbb{E}^{3},g)$ and $(\mathbb{R}^{1+4},\tilde{g})$ respectively. And 
$\tilde{k}$ is the second fundamental form of $(\mathbb{R}\times\mathbb{E}^{3},g)$ in $(\mathbb{R}^{1+4},\tilde{g})$, which contains only the first derivatives of the 
acoustical metric. Now since $(\mathbb{R}^{1+4},\tilde{g})$ is a flat spacetime, $R$ contains only the first derivative of the acoustical metric. I.e. The 
principal part of $R$ vanishes.

\chapter{The Fundamental Energy Estimate}



\section{Bootstrap Assumptions. Statement of the Theorem}

In this chapter we consider the wave equation
\begin{equation}
 \Box_{\tilde{g}}\psi=0
\end{equation}
 in a spacetime whose metric is the conformal acoustical metric $\tilde{g}$. As we have seen in Chapter 1, this equation is satisfied by a variation of the wave
function $\phi$ through solutions. $\tilde{g}$ is related to the acoustical metric $g$ by
\begin{equation}
 \tilde{g}_{\mu\nu}=\Omega{g_{\mu\nu}}
\end{equation}
where
\begin{equation}
 \Omega=\frac{\rho/\rho_{0}}{\eta/\eta_{0}}=\frac{\rho}{\eta}
\end{equation}
where $\rho_{0}$ and $\eta_{0}$ are the mass density and sound speed corresponding to the constant state which we can both set equal to unity
by appropriate choice of units. More generally, we consider the inhomogeneous wave equation
\begin{equation}
 \Box_{\tilde{g}}\psi=\rho
\end{equation}
where the metric $\tilde{g}$, $\psi$ and $\rho$ are all defined in the domain $M_{\epsilon_{0}}$ of the maximal solution.

In application to first order estimates we shall take each $\psi$ to be one of:
\begin{equation}
 \psi_{1}=\partial_{0}\phi-h_{0}=\psi_{0}; \quad \partial_{i}\phi: i=1,2,3; \quad \Omega_{ij}\phi: i<j=1,2,3; \quad S\phi.
\end{equation}
Remember in Chapter 2 we have set:
\begin{align*}
 h_{0}=0
\end{align*}

We know that each $\psi_{1}$ vanishes in the exterior of the outgoing characteristic hypersurface $C_{0}$. In application to higher order estimates we shall 
take each $\psi$ to be one of:
\begin{equation}
 \psi_{n}=Y_{i_{1}}...Y_{i_{n-1}}\psi_{1}
\end{equation}
where $\psi_{1}$ is any one of the above first order variations and the indices $i_{1}...i_{n-1}$ take values in the set ${1,2,3,4,5}$, with $\{Y_{i}:
i=1,2,3,4,5\}$ the set of vectorfields to be discussed in the next chapter. Each of the $\psi_{n}$ also vanishes in the exterior of $C_{0}$.

We shall therefore assume in general that $\psi$  vanishes in the exterior of $C_{0}$. We may then confine attention to the spacetime domain $W_{\epsilon_{0}}$.
For any $s\in (0,t_{*\epsilon_{0}})$, we set:
\begin{equation}
 W^{s}_{\epsilon_{0}}=\bigcup_{(t,u)\in[0,s]\times[0,\epsilon_{0}]}S_{t,u}
\end{equation}
In particular we have:
\begin{equation}
 \bigcup_{s\in(0,t_{*\epsilon_{0}})}W^{s}_{\epsilon_{0}}=W^{*}_{\epsilon_{0}}
\end{equation}
Everything which follows depends on the geometric construction in Chapter 2. We note that the characteristic hypersurfaces $C_{u}$ depend only on the conformal
 class of the metric $g$. Thus, the 2-parameter foliation of $W^{*}_{\epsilon_{0}}$ given by the surfaces $S_{t,u}$ likewise depends only on the conformal class of 
$g$. The geometric properties of this foliation may be referred either to $g$ or to the conformal metric $\tilde{g}$. Since we are studying the wave equation in 
the metric $\tilde{g}$, it would appear more natural if we refer these properties to $\tilde{g}$. However, our bootstrap argument is constructed relative to $g$, 
and we introduce additional assumptions referring to the conformal factor $\Omega$.

The bootstrap assumptions are the followings:

There is a positive constant $C$ independent of $s$ such that in $W^{*}_{\epsilon_{0}}$,
\begin{align*}
 \textbf{A1}: C^{-1}\leq\Omega\leq{C}\\
\textbf{A2}: C^{-1}\leq\eta\leq{C}\\
\textbf{A3}: \mu\leq{C[1+\log(1+t)]}
\end{align*}
 and
\begin{align*}
 \textbf{B1}: C^{-1}(1+t)^{-1}\leq \nu \leq C(1+t)^{-1}\\
\textbf{B2}: |\underline{\nu}| \leq C(1+t)^{-1}[1+\log(1+t)]^{4}\\
\textbf{B3}: |\hat{\chi}|\leq C(1+t)^{-1}[1+\log(1+t)]^{-2}\\
\textbf{B4}: |\hat{\underline{\chi}}|\leq C(1+t)[1+\log(1+t)]^{-6}
\end{align*}
Here,
\begin{align}
 \nu=\frac{1}{2}\tilde{\textrm{tr}}\chi=\frac{1}{2}(\textrm{tr}\chi+L\log\Omega)\\
\underline{\nu}=\frac{1}{2}\tilde{\textrm{tr}}\underline{\chi}=\frac{1}{2}(\textrm{tr}\underline{\chi}+\underline{L}\log\Omega)
\end{align}
$\hat{\chi}$ and $\hat{\underline{\chi}}$ are the trace-free part of $\chi$ and $\underline{\chi}$ respectively, and the pointwise norms of the tensors on 
$S_{t,u}$ are with respect to the induced metric $\slashed{g}$.

Next, we have the assumptions:
\begin{align*}
 \textbf{B5}: |L\log\Omega|\leq C(1+t)^{-1}[1+\log(1+t)]^{-2}\\
\textbf{B6}: |\underline{L}\log\Omega|\leq C(1+t)[1+\log(1+t)]^{-6}
\end{align*}
and
\begin{align*}
 \textbf{B7}: |\zeta+\eta|\leq C(1+t)^{-1}[1+\log(1+t)]\\
\textbf{B8}: |\slashed{d}(\eta^{-1}\kappa)|\leq C(1+t)^{-1}[1+\log(1+t)]\\
\textbf{B9}: |L\eta|\leq C(1+t)^{-1}[1+\log(1+t)]^{2}\\
\textbf{B10}:|L(\eta^{-1}\kappa)|\leq C(1+t)^{-1}[1+\log(1+t)]^{3}\\
\textbf{B11}:|\underline{L}(\eta^{-1}\kappa)|\leq C(1+t)[1+\log(1+t)]^{-2}
\end{align*}
as well as
\begin{align*}
 \textbf{B12}:|L\nu+\nu^{2}|\leq C(1+t)^{-2}[1+\log(1+t)]^{-2}\\
\textbf{B13}:|\underline{L}\nu|\leq C(1+t)^{-2}[1+\log(1+t)]^{3}\\
\textbf{B14}:|\slashed{d}\nu|\leq C(1+t)^{-2}[1+\log(1+t)]^{\frac{1}{2}}
\end{align*}
The next set of bootsrap assumptions concern the behavior of the function $\mu$. In the following we denote by $f_{+}$ and $f_{-}$ respectively the positive
and negative parts of the function $f$:
\begin{equation}
 f_{+}=\max{(f(x),0)}, f_{-}=\min{(f(x),0)}
\end{equation}
\begin{equation*}
 \textbf{C1}: \mu^{-1}(L\mu)_{+}\leq (1+t)^{-1}[1+\log(1+t)]^{-1}+A(t)
\end{equation*}
where $A(t)$ is a nonnegative function such that:
\begin{equation}
 \int^{s}_{0}A(t)dt \leq C
\end{equation}
($C$ is independent of $s$.)
\begin{equation*}
 \textbf{C2}:\mu^{-1}(L\mu+\underline{L}\mu)_{+}\leq B(t)
\end{equation*}
where $B(t)$ is a nonnegative function such that
\begin{equation}
 \int^{s}_{0}(1+t)^{-2}[1+\log(1+t)]^{4}B(t)dt \leq C
\end{equation}
($C$ is also independent of $s$).

Moreover, denoting by $\mathcal{U}$ the region
\begin{equation}
 \mathcal{U}=\{x\in W^{*}_{\epsilon_{0}}: \mu < 1/4\}
\end{equation}
we have:
\begin{equation*}
 \textbf{C3}: L\mu \leq -C^{-1}(1+t)^{-1}[1+\log(1+t)]^{-1}
\end{equation*}
in $\mathcal{U}\bigcap W^{s}_{\epsilon_{0}}$

The final set of bootstrap assumptions concerns the existence of a function $\omega$ verifying the following conditions:
\begin{align*}
 \textbf{D1}:C^{-1}(1+t) \leq \omega \leq C(1+t)\\
\textbf{D2}: |L\omega-\nu\omega| \leq C[1+\log(1+t)]^{-2}\\
\textbf{D3}: |\underline{L}\omega| \leq C[1+\log(1+t)]^{3}\\
\textbf{D4}: |\slashed{d}\omega| \leq C[1+\log(1+t)]^{\frac{1}{2}}\\
\textbf{D5}: \int^{s}_{0}\{\int_{0}^{\epsilon_{0}}\sup_{S_{t,u}}(\mu|\Box_{\tilde{g}}\omega|)du\}dt\leq C[1+\log(1+s)]^{4}
\end{align*}
Here we note that by (5.3) $\textbf{A1}$ is equivalent to the following modulo $\textbf{A2}$:
\begin{align*}
 \textbf{A1}^{\prime}: C^{-1}\leq\rho\leq C
\end{align*}

The main theorem in this chapter is 

$\textbf{Theorem 5.1}$ Let $\psi$ be a solution of the wave equation associated to the metric $\tilde{g}$, defined in $M_{\epsilon_{0}}$ and vanishing in the 
exterior of $C_{0}$. Suppose that assumptions $\textbf{A1}-\textbf{A3}$, $\textbf{B1}-\textbf{B14}$, $\textbf{C1}-\textbf{C3}$, $\textbf{D1}-\textbf{D5}$ hold 
in $W^{s}_{\epsilon_{0}}$, for some $s\in (0,t_{*\epsilon_{0}}]$. Let us denote:
\begin{equation}
 D_{0}=\int^{\epsilon_{0}}_{0}\{\int_{S_{0,u}}[(L\psi)^{2}+(\underline{L}\psi)^{2}+|\slashed{d}\psi|^{2}]d\mu_{\slashed{g}}\}du.
\end{equation}
Then there exist constants $C$ independent of $s$ such that:
\begin{align*}
 (\textbf{i}) \sup_{t\in[0,s]}\int_{0}^{\epsilon_{0}}\{\int_{S_{t,u}}[\mu(1+\mu)((L\psi)^{2}+|\slashed{d}\psi|^{2})
+(\underline{L}\psi)^{2}]d\mu_{\slashed{g}}\}du \leq CD_{0}\\
(\textbf{ii}) \int_{S_{t,u}}\psi^{2}d\mu_{\slashed{g}} \leq C\epsilon_{0}D_{0}\\
(\textbf{iii}) \sup_{u\in[0,\epsilon_{0}]}\int^{s}_{0}\{\int_{S_{t,u}}[(1+\mu)(L\psi)^{2}+\mu|\slashed{d}\psi|^{2}]d\mu_{\slashed{g}}
\}dt \leq CD_{0}\\
(\textbf{iv}) \sup_{t\in[0,s]}[1+\log(1+t)]^{-4}(1+t)^{2}\int^{\epsilon_{0}}_{0}\{\int_{S_{t,u}}\mu[(L\psi+\nu\psi)^{2}+|\slashed{d}\psi|^{2}]
d\mu_{\slashed{g}}\}du \leq CD_{0}\\
(\textbf{v}) \sup_{u\in[0,\epsilon_{0}]}[1+\log(1+t)]^{-4}\int^{s}_{0}(1+t)^{2}\{\int_{S_{t,u}}(L\psi+\nu\psi)^{2}d\mu_{\slashed{g}}\}dt \leq CD_{0}\\
(\textbf{vi}) \int_{\mathcal{U}\bigcup W^{s}_{\epsilon_{0}}}(1+t)[1+\log(1+t)]^{-1}|\slashed{d}\psi|^{2}d\mu_{\slashed{g}}dudt \leq CD_{0}[1+\log(1+s)]^{4}
\end{align*}

\section{The Multiplier Fields $K_{0}$ and $K_{1}$. The Associated Energy-Momentum Density Vectorfields}
 We begin with the energy-momentum-stress tensor $\tilde{T}_{\mu\nu}$ associated to the function $\psi$ through the metric $\tilde{g}$:
\begin{align}
 \tilde{T}_{\mu\nu}:=\partial_{\mu}\psi\partial_{\nu}\psi-\frac{1}{2}\tilde{g}_{\mu\nu}(\tilde{g}^{-1})^{\kappa\lambda}\partial_{\kappa}\psi
\partial_{\lambda}\psi\\\notag
=\partial_{\mu}\psi\partial_{\nu}\psi-\frac{1}{2}g_{\mu\nu}(g^{-1})^{\kappa\lambda}\partial_{\kappa}\psi\partial_{\lambda}\psi=T_{\mu\nu}
\end{align}
We have:
\begin{equation*}
 \tilde{D}^{\mu}\tilde{T}_{\mu\nu}:=(\tilde{g}^{-1})^{\mu\lambda}\tilde{D}_{\mu}\tilde{T}_{\lambda\nu}=\partial_{\nu}\psi\Box_{\tilde{g}}\psi
\end{equation*}
where $\tilde{D}$ is the covariant derivative operator associated to the metric $\tilde{g}_{\mu\nu}$. Thus, for a solution of (5.1), $\tilde{T}_{\mu\nu}$
is divergence-free with respect to the metric $\tilde{g}$, while for a solution of the inhomogenneous wave equation (5.4),
\begin{equation}
 \tilde{D}^{\mu}\tilde{T}_{\mu\nu}=\rho\partial_{\nu}\psi
\end{equation}
We consider the future-directed, time-like vectorfield $K_{0}$:
\begin{equation}
 K_{0}=(1+\alpha^{-1}\kappa)L+\underline{L}
\end{equation}
Also, given a function $\omega$ satisfying $\textbf{D1}$-$\textbf{D5}$ we consider the future-directed, null vectorfield $K_{1}$:
\begin{equation}
 K_{1}=(\omega/\nu)L
\end{equation}
At this point, we discuss why we choose these two multiplier fields. A multiplier should be non-spacelike and future-directed so that it gives positive energies. Also the multiplier should be smooth with respect to the basis $(\frac{\partial}{\partial t},\frac{\partial}{\partial u},\frac{\partial}{\partial\vartheta^{A}})$ in acoustical coordinates $(t,u,\vartheta^{A})$, because only the smooth vectorfields give error terms which can be estimated, due to their expressions in terms of deformation tensors. So in the flat spacetime, $\partial_{t}$ is a natural 
choice. While in the present case, $L+\underline{L}$ fits the above requirements. Indeed $L+\underline{L}$ is the analogue of time translation, i.e. $\partial_{t}$ in flat spacetime. However, in contrast to this, $K_{0}$ becomes null when $\mu\longrightarrow0$  
On the other hand, before we set $\eta_{0}=1$, $L$ and $\underline{L}$ 
have different physical dimensions, $L$ being $\frac{1}{time}$ and $\underline{L}$ being $\frac{1}{length}$. However, after we set $\eta_{0}=1$, the unit
of time is the time required for sound to travel the unit length in the constant state. Therefore time is the length and $L$, $\underline{L}$ are 
both $\frac{1}{length}$. To deal with the growth of $\mu$, we need to plus the term 
$\alpha^{-1}\kappa L$. This is how we choose $K_{0}$.  

     Concerning the choice of $K_{1}$, we remind you the generator of the inverted time translation in Minkowski spacetime:
\begin{align*}
 K=(u^{2}\underline{L}+\underline{u}^{2}L)
\end{align*}
where 
\begin{align*}
 L=\partial_{t}+\partial_{r}\quad \underline{L}=\partial_{t}-\partial_{r}
\end{align*}
and
\begin{align*}
 \underline{u}=t+r\quad u=t-r
\end{align*}
Now at the present case, we want to choose a similar vectorfield. The analogue of $t-r$ is our acoustical function $u$, which is bounded in $W^{*}_{\epsilon_{0}}$.
While when $t$ is large, $\underline{u}\sim t$.
Since the first term in the expression for $K$ is contained in $K_{0}$, so by $\textbf{B1}$ and $\textbf{D1}$, it is reasonale to choose:
\begin{align*}
 K_{1}=(\omega/\nu)L
\end{align*}
We will see later more clearly why we exactly choose this form.

We denote by $\tilde{\pi}_{0}$ and $\tilde{\pi}_{1}$ the Lie derivatives of the conformal metric $\tilde{g}$ with respect to the vectorfields $K_{0}$ and
$K_{1}$:
\begin{equation}
 \tilde{\pi}_{0}=\mathcal{L}_{K_{0}}\tilde{g}, \tilde{\pi}_{1}=\mathcal{L}_{K_{1}}\tilde{g}
\end{equation}
To $K_{0}$ we associate the vectorfield 
\begin{equation}
 \tilde{P}^{\mu}_{0}=-\tilde{T}^{\mu}_{\nu}K^{\nu}_{0}
\end{equation}
where
\begin{align*}
 \tilde{T}^{\mu}_{\nu}=(\tilde{g}^{-1})^{\mu\kappa}T_{\kappa\nu}
\end{align*}

To $K_{1}$ we associate the vectorfield 
\begin{equation}
 \tilde{P}^{\mu}_{1}=-\tilde{T}^{\mu}_{\nu}K^{\nu}_{1}-(\tilde{g}^{-1})^{\mu\nu}(\omega\psi\partial_{\nu}\psi-\frac{1}{2}\psi^{2}\partial_{\nu}\omega)
\end{equation}
Let
\begin{equation}
 \tilde{T}^{\mu\nu}=(\tilde{g}^{-1})^{\nu\lambda}\tilde{T}^{\mu}_{\lambda}=(\tilde{g}^{-1})^{\mu\kappa}(\tilde{g}^{-1})^{\nu\lambda}\tilde{T}_{\kappa\lambda}
\end{equation}
For any vectorfield $X$, we have , by the virtue of the symmetry of $\tilde{T}^{\mu\nu}$,
\begin{align}
\tilde{D}_{\mu}(\tilde{T}^{\mu}_{\nu}X^{\nu})=(\tilde{D}_{\mu}\tilde{T}^{\mu}_{\nu})X^{\nu}
+\tilde{T}^{\mu\lambda}(\tilde{g}_{\lambda\nu}\tilde{D}_{\mu}X^{\nu})\\\notag
=\rho X\psi+\frac{1}{2} \tilde{T}^{\mu\lambda}(\tilde{g}_{\lambda\nu}\tilde{D}_{\mu}X^{\nu}+\tilde{g}_{\mu\nu}\tilde{D}_{\lambda}X^{\nu})=
\rho X\psi+\frac{1}{2} \tilde{T}^{\mu\lambda}\mathcal{L}_{X}\tilde{g}_{\mu\lambda}
\end{align}
Thus, we have
\begin{equation}
 \tilde{D}_{\mu}\tilde{P}^{\mu}_{0}=-\rho K_{0}\psi-\frac{1}{2}\tilde{T}^{\mu\nu}\tilde{\pi}_{0,\mu\nu}:=\tilde{Q}_{0}
\end{equation}
Also, we have
\begin{equation*}
 \tilde{D}_{\mu}\tilde{P}^{\mu}_{1}=-\rho(K_{1}\psi+\omega\psi)-\frac{1}{2}\tilde{T}^{\mu\nu}\tilde{\pi}_{1,\mu\nu}
-\omega(\tilde{g}^{-1})^{\mu\nu}\partial_{\mu}\psi\partial_{\nu}\psi+\frac{1}{2}\psi^{2}\Box_{\tilde{g}}\omega
\end{equation*}
Taking into account the fact that:
\begin{equation*}
 \textrm{tr}\tilde{T}:=\tilde{g}_{\mu\nu}\tilde{T}^{\mu\nu}=-(\tilde{g}^{-1})^{\mu\nu}\partial_{\mu}\psi\partial_{\nu}\psi
\end{equation*}
and introducing:
\begin{equation}
 \tilde{\pi}^{\prime}_{1}=\tilde{\pi}_{1}-2\omega\tilde{g}
\end{equation}
we can write the above in the form
\begin{equation}
 \tilde{D}_{\mu}\tilde{P}^{\mu}_{1}=-\rho(K_{1}\psi+\omega\psi)-\frac{1}{2}\tilde{T}^{\mu\nu}\tilde{\pi}^{\prime}_{1,\mu\nu}
+\frac{1}{2}\psi^{2}\Box_{\tilde{g}}\omega:=\tilde{Q}_{1}
\end{equation}
Consider now the equation
\begin{equation}
 \tilde{D}_{\mu}\tilde{P}^{\mu}=\tilde{Q}
\end{equation}
for an arbitrary vectorfield $\tilde{P}$ and function $\tilde{Q}$. In arbitrary local coordinates we have
\begin{equation*}
 \tilde{D}_{\mu}\tilde{P}^{\mu}=\frac{1}{\sqrt{-\det\tilde{g}}}\partial_{\mu}(\sqrt{-\det\tilde{g}}\tilde{P}^{\mu})=
\frac{1}{\Omega^{2}\sqrt{-\det g}}\partial_{\mu}(\Omega^{2}\sqrt{-\det{g}}\tilde{P}^{\mu})=\Omega^{-2}D_{\mu}P^{\mu}
\end{equation*}
where
\begin{equation}
 P^{\mu}=\Omega^{2}\tilde{P}^{\mu}
\end{equation}
Hence, with
\begin{equation}
 Q=\Omega^{2}\tilde{Q}
\end{equation}
equation (5.28) is equivalent to the equation:
\begin{equation}
 D_{\mu}P^{\mu}=Q
\end{equation}
We may express (5.31) in acoustical coordinates $(t,u,\vartheta^{1},\vartheta^{2})$. Expanding $P$ in the associated coordinate frame field,
\begin{align*}
 (\frac{\partial}{\partial t}, \frac{\partial}{\partial u}, \frac{\partial}{\partial \vartheta^{1}}, \frac{\partial}{\partial \vartheta^{2}})\\
P=P^{t}\frac{\partial}{\partial t}+P^{u}\frac{\partial}{\partial u}+\sum_{A}(P^{\vartheta})^{A}\frac{\partial}{\partial \vartheta^{A}}
\end{align*}
and noting that, by (2.41), we have
\begin{equation}
 \sqrt{-\det g}=\mu\sqrt{\det\slashed{g}}
\end{equation}
equation (5.31) takes the form:
\begin{equation}
 \frac{1}{\sqrt{\det\slashed{g}}}\{\frac{\partial}{\partial t}(\mu\sqrt{\det\slashed{g}}P^{t})+\frac{\partial}{\partial u}(\mu\sqrt{\det\slashed{g}}P^{u})\}
+\slashed{\textrm{div}}M=\mu Q
\end{equation}
where $M$ is the vectorfield on $S^{2}$:
\begin{equation*}
 M=\mu\sum_{A}(P^{\vartheta})^{A}\frac{\partial}{\partial \vartheta^{A}}
\end{equation*}
This follows from the fact that
\begin{equation*}
 \slashed{\textrm{div}}M=\frac{1}{\sqrt{\det\slashed{g}}}\frac{\partial}{\partial\vartheta^{A}}(\sqrt{\det\slashed{g}}M^{A})
\end{equation*}
We integrate (5.33) on $S^{2}$ with respect to the measure
\begin{equation}
 d\mu_{\slashed{g}}=\sqrt{\det\slashed{g}}d\vartheta^{1}d\vartheta^{2}
\end{equation}
to obtain the equation
\begin{equation}
 \frac{\partial}{\partial t}(\int_{S_{t,u}}\mu P^{t}d\mu_{\slashed{g}})
+\frac{\partial}{\partial u}(\int_{S_{t,u}}\mu P^{u}d\mu_{\slashed{g}})=\int_{S_{t,u}}\mu Q d\mu_{\slashed{g}}
\end{equation}
Replacing $(t,u)$ by $(t^{\prime},u^{\prime})$ and integrating with respect to $(t^{\prime},u^{\prime})$ on $[0,t]\times[0,u]$, we obtain, under the hypothesis 
that $P$ vanishes in the closure of the exterior of $C_{0}$.
\begin{equation}
 \mathcal{E}^{u}(t)-\mathcal{E}^{u}(0)+\mathcal{F}^{t}(u)=\int_{W^{t}_{u}}Q d\mu_{g}
\end{equation}
Here $\mathcal{E}^{u}$ is the ``energy'':
\begin{equation}
 \mathcal{E}^{u}(t)=\int_{\Sigma^{u}_{t}}\mu P^{t}d\mu_{\slashed{g}}du^{\prime}=\int^{u}_{0}(\int_{S_{t,u^{\prime}}}\mu P^{t}d\mu_{\slashed{g}})du^{\prime}
\end{equation}
and $\mathcal{F}^{t}$ is the ``flux'':
\begin{equation}
 \mathcal{F}^{t}(u)=\int_{C_{u}^{t}}\mu P^{u}d\mu_{\slashed{g}}dt^{\prime}=\int_{0}^{t}(\int_{S_{t^{\prime},u}}\mu P^{u}d\mu_{\slashed{g}})dt^{\prime}
\end{equation}
In obtaining (5.36) we have used the fact that
\begin{equation}
 \int_{W^{t}_{u}}Qd\mu_{g}=\int\int_{[0,t]\times[0,u]}(\int_{S_{t^{\prime},u^{\prime}}}\mu Qd\mu_{\slashed{g}})du^{\prime}dt^{\prime}
\end{equation}
In the above, $\Sigma^{u}_{t}$ is the annular region:
\begin{equation}
 \Sigma^{u}_{t}=\bigcup_{u^{\prime}\in [0,u]}S_{t,u^{\prime}}
\end{equation}
in the space-like hyperplane $\Sigma_{t}$, $C^{t}_{u}$ is the characteristic hypersurface $C_{u}$ truncated from above by $\Sigma_{t}$:
\begin{equation}
 C^{t}_{u}=\bigcup_{t^{\prime}\in [0,t]}S_{t^{\prime},u}
\end{equation}
and $W^{t}_{u}$ is the spacetime domain:
\begin{equation}
 W^{t}_{u}=\bigcup_{(t^{\prime},u^{\prime})\in [0,t]\times[0,u]}S_{t^{\prime},u^{\prime}}
\end{equation}
bounded by the characteristic hypersurfaces $C_{u}$ and $C_{0}$ and the space-like hypersurfaces $\Sigma_{t}$ and $\Sigma_{0}$.

Now we expand $P$ in the frame field $(L,\underline{L},X_{1},X_{2})$:
\begin{equation*}
 P=P^{L}L+P^{\underline{L}}\underline{L}+\sum_{A}P^{A}X_{A}
\end{equation*}
 From the conclusion of Chapter 2 and Chapter 3:
\begin{align}
L=\frac{\partial}{\partial t}\\
\underline{L}=\eta^{-1}\kappa\frac{\partial}{\partial t}+2(\frac{\partial}{\partial u}-\Xi^{A}\frac{\partial}{\partial\vartheta^{A}})\\
X_{A}=\frac{\partial}{\partial \vartheta^{A}}
\end{align}
By direct calculations, we get
\begin{align}
P^{t}=P^{L}+\eta^{-1}\kappa P^{\underline{L}}\\
P^{u}=2P^{\underline{L}}\\
(P^{\vartheta})^{A}=P^{A}-2P^{\underline{L}}\Xi^{A}
\end{align}

We shall write down the expressions for the energy and flux integrals corresponding to $P_{0}$ and $P_{1}$.

We begin with the components of the energy-momentum-stress tensor in the null frame 

$(L,\underline{L},X_{1},X_{2})$:
\begin{align}
 T_{LL}=(L\psi)^{2},\quad T_{\underline{L}\underline{L}}=(\underline{L}\psi)^{2},\quad
T_{L\underline{L}}=\mu|\slashed{d}\psi|^{2}\\\notag
T_{LA}=(L\psi)(\slashed{d}_{A}\psi),\quad T_{\underline{L}A}=(\underline{L}\psi)(\slashed{d}_{A}\psi)\\\notag
T_{AB}=(\slashed{d}_{A}\psi)(\slashed{d}_{B}\psi)-\frac{1}{2}\slashed{g}_{AB}(-\mu^{-1}(L\psi)(\underline{L}\psi)+|\slashed{d}\psi|^{2})
\end{align}
Here, $\slashed{d}_{A}\psi=X_{A}\psi$, and the pointwise norms of tensors on $S_{t,u}$ are with respect to the induced metric $\slashed{g}$.

We also give the components of the tensor $T^{\mu\nu}$. From (3.130), we have:
\begin{align}
 T^{LL}=\frac{(\underline{L}\psi)^{2}}{4\mu^{2}},\quad T^{\underline{L}\underline{L}}=\frac{(L\psi)^{2}}{4\mu^{2}},\quad
T^{L\underline{L}}=\frac{|\slashed{d}\psi|^{2}}{4\mu}\\\notag
T^{LA}=-\frac{1}{2\mu}(\underline{L}\psi)(\slashed{d}^{A}\psi), \quad T^{\underline{L}A}=-\frac{1}{2\mu}(L\psi)(\slashed{d}^{A}\psi)\\\notag
T^{AB}=\frac{1}{2\mu}(\slashed{g}^{-1})^{AB}(L\psi)(\underline{L}\psi)+\{(\slashed{d}^{A}\psi)(\slashed{d}^{B}\psi)-\frac{1}{2}(\slashed{g}^{-1})
^{AB}|\slashed{d}\psi|^{2}\}
\end{align}
where $\slashed{d}^{A}\psi=(\slashed{g}^{-1})^{AB}\slashed{d}_{B}\psi$

We consider $P_{0}$. From (5.66) and (5.54) we have
\begin{align}
 P^{\mu}_{0}=\Omega^{2}\tilde{P}^{\mu}_{0}=-\Omega^{2}\tilde{T}^{\mu}_{\nu}K^{\nu}_{0}=-\Omega T^{\mu}_{\nu}K^{\nu}_{0}
\end{align}
We get
\begin{align*}
 P^{\underline{L}}_{0}=-\Omega T^{\underline{L}}_{L}(1+\eta^{-1}\kappa)-\Omega T^{\underline{L}}_{\underline{L}}
=\frac{\Omega}{2\mu}((1+\eta^{-1}\kappa)T_{LL}+T_{L\underline{L}})\\
 P^{L}_{0}=-\Omega T^{L}_{L}(1+\eta^{-1}\kappa)-\Omega T^{L}_{\underline{L}}
=\frac{\Omega}{2\mu}((1+\eta^{-1}\kappa)T_{L\underline{L}}+T_{\underline{L}\underline{L}})
\end{align*}
By (3.91), we have
\begin{align}
 P^{\underline{L}}_{0}=\frac{\Omega}{2\mu}((1+\eta^{-1}\kappa)(L\psi)^{2}+\mu|\slashed{d}\psi|^{2})\\\notag
P^{L}_{0}=\frac{\Omega}{2\mu}((1+\eta^{-1}\kappa)\mu|\slashed{d}\psi|^{2}+(\underline{L}\psi)^{2})
\end{align}
So we get
\begin{align}
 \mathcal{E}^{u}_{0}(t)=\int_{\Sigma^{u}_{t}}\frac{\Omega}{2}\{\eta^{-1}\kappa(1+\eta^{-1}\kappa)(L\psi)^{2}+(\underline{L}\psi)^{2}
+(1+2\eta^{-1}\kappa)\mu|\slashed{d}\psi|^{2}\}d\mu_{\slashed{g}}du^{\prime}\\
\mathcal{F}^{t}_{0}(u)=\int_{C^{t}_{u}}\Omega\{(1+\eta^{-1}\kappa)(L\psi)^{2}+\mu|\slashed{d}\psi|^{2}\}d\mu_{\slashed{g}}dt^{\prime}
\end{align}
We now consider $P_{1}$. By some similar calculation we have:
\begin{equation}
 P^{\mu}_{1}=-\Omega\{T^{\mu}_{\nu}K^{\nu}_{1}+(g^{-1})^{\mu\nu}(\omega\psi\partial_{\nu}\psi-\frac{1}{2}\psi^{2}\partial_{\nu}\omega)\}
\end{equation}
Then we get
\begin{align*}
 P^{\underline{L}}_{1}=\frac{\Omega}{2\mu}\{\omega\nu^{-1}T_{LL}+\omega\psi(L\psi)-\frac{1}{2}\psi^{2}(L\omega)\}\\
P^{L}_{1}=\frac{\Omega}{2\mu}\{\omega\nu^{-1}T_{L\underline{L}}+\omega\psi(\underline{L}\psi)-
\frac{1}{2}\psi^{2}(\underline{L}\omega)\}
\end{align*}
Then from (5.49) we get 
\begin{align}
  P^{\underline{L}}_{1}=\frac{\Omega}{2\mu}\{\omega\nu^{-1}(L\psi)^{2}+\omega\psi(L\psi)-\frac{1}{2}\psi^{2}(L\omega)\}\\\notag
P^{L}_{1}=\frac{\Omega}{2\mu}\{\omega\nu^{-1}\mu|\slashed{d}\psi|^{2}+\omega\psi(\underline{L}\psi)-
\frac{1}{2}\psi^{2}(\underline{L}\omega)\}
\end{align}
and
\begin{align}
 \mathcal{E}^{u}_{1}(t)=\int_{\Sigma^{u}_{t}}\frac{\Omega}{2}\{\omega\nu^{-1}[\eta^{-1}\kappa(L\psi)^{2}+\mu|\slashed{d}\psi|^{2}]+
\omega\psi[\eta^{-1}\kappa(L\psi)+(\underline{L}\psi)]\\\notag
-\frac{1}{2}\psi^{2}[\eta^{-1}\kappa(L\omega)+(\underline{L}\omega)]\}d\mu_{\slashed{g}}du^{\prime}\\
\mathcal{F}^{t}_{1}(u)=\int_{C^{t}_{u}}\Omega\{\omega\nu^{-1}(L\psi)^{2}+\omega\psi(L\psi)-\frac{1}{2}\psi^{2}(L\omega)\}d\mu_{\slashed{g}}dt^{\prime}
\end{align}
We actually define the flux integral associated to the vectorfield $K_{1}$ to be:
\begin{equation}
 \mathcal{F}^{\prime t}_{1}(u)=\int_{C^{t}_{u}}\Omega\omega\nu^{-1}(L\psi+\nu\psi)^{2}d\mu_{\slashed{g}}dt^{\prime}
\end{equation}
By a direct calculation, we deduce:
\begin{equation}
 \mathcal{F}^{t}_{1}(u)-\mathcal{F}^{\prime t}_{1}(u)=-\int_{C^{t}_{u}}\frac{1}{2}\Omega[L(\omega\psi^{2})+2\nu\omega\psi^{2}]d\mu_{\slashed{g}}dt^{\prime}
\end{equation}
Now we have
\begin{equation}
 \slashed{\mathcal{L}}_{L}d\mu_{\slashed{g}}=\textrm{tr}\chi d\mu_{\slashed{g}}
\end{equation}
In view of (5.9), we have for an arbitrary function $f$ defined in $W^{*}_{\epsilon_{0}}$:
\begin{align}
 \frac{\partial}{\partial t}(\int_{S_{t,u}}\Omega fd\mu_{\slashed{g}})=\int_{S_{t,u}}\Omega(Lf+2\nu f)d\mu_{\slashed{g}}
\end{align}
Setting $f=\frac{1}{2}\omega\psi^{2}$, we have from (5.60),
\begin{align}
 \mathcal{F}^{t}_{1}(u)-\mathcal{F}^{\prime t}_{1}(u)=-\int^{t}_{0}\frac{\partial}{\partial t^{\prime}}(\int_{S_{t^{\prime},u}}\frac{1}{2}\Omega\omega\psi^{2}
d\mu_{\slashed{g}})dt^{\prime}\\\notag
=-\int_{S_{t,u}}\frac{1}{2}\Omega\omega\psi^{2}d\mu_{\slashed{g}}+\int_{S_{0,u}}\frac{1}{2}\Omega\omega\psi^{2}d\mu_{\slashed{g}}
\end{align}
Consider next the energy integral (5.57). We actually define the energy integral associated to the vectorfield $K_{1}$ to be:
\begin{align}
 \mathcal{E}^{\prime u}_{1}(t)
=\int_{\Sigma^{u}_{t}}\frac{\Omega}{2}\omega\nu^{-1}\{\eta^{-1}\kappa(L\psi+\nu\psi)^{2}+\mu|\slashed{d}\psi|^{2}\}d\mu_{\slashed{g}}du^{\prime}
\end{align}
Taking into account the fact that,
\begin{equation}
 \underline{L}-\eta^{-1}\kappa L=2T
\end{equation}
we find:
\begin{align}
 \mathcal{E}^{u}_{1}(t)-\mathcal{E}^{\prime u}_{1}(t)=\int_{\Sigma^{u}_{t}}\frac{\Omega}{2}\{2\omega\psi(T\psi)-\frac{1}{2}[\eta^{-1}\kappa L\omega
+\underline{L}\omega+2\eta^{-1}\kappa\nu\omega]\psi^{2}\}d\mu_{\slashed{g}}du^{\prime}
\end{align}
From the definition of $\theta$, we have
\begin{equation}
 \slashed{\mathcal{L}}_{T}d\mu_{\slashed{g}}=\kappa\textrm{tr}\theta d\mu_{\slashed{g}}
\end{equation}
From (5.67), we know that for an arbitrary function $f$ defined in $W^{*}_{\epsilon_{0}}$:
\begin{align}
\frac{\partial}{\partial u}(\int_{S_{t,u}}\Omega fd\mu_{\slashed{g}})=\int_{S_{t,u}}\slashed{\mathcal{L}}_{T}(\Omega fd\mu_{\slashed{g}})
=\int_{S_{t,u}}\Omega\{Tf+[\kappa\textrm{tr}\theta+T(\log\Omega)]f\}d\mu_{\slashed{g}}
\end{align}
 Setting $f=\frac{1}{2}\omega\psi^{2}$ yields:
\begin{align*}
 \int_{\Sigma^{u}_{t}}\frac{\Omega}{2}\{2\omega\psi(T\psi)+(T\omega)\psi^{2}+[\kappa\textrm{tr}\theta+T(\log\Omega)]\omega\psi^{2}\}
d\mu_{\slashed{g}}du^{\prime}\\
=\int^{u}_{0}\frac{\partial}{\partial u^{\prime}}(\int_{S_{t,u^{\prime}}}\frac{\Omega}{2}\omega\psi^{2}d\mu_{\slashed{g}})du^{\prime}
=\int_{S_{t,u}}\frac{1}{2}\Omega\omega\psi^{2}d\mu_{\slashed{g}}
\end{align*}
Here we have used the fact that $\psi$ vanishes on $S_{t,0}\subset C_{0}$. Comparing this with (5.66) and using (5.65), we get:
\begin{equation}
 \mathcal{E}^{u}_{1}(t)-\mathcal{E}^{\prime u}_{1}(t)=\int_{S_{t,u}}\frac{1}{2}\Omega\omega\psi^{2}d\mu_{\slashed{g}}-I^{u}(t)
\end{equation}
where
\begin{equation*}
 I^{u}(t)=\int_{\Sigma^{u}_{t}}\frac{1}{2}\Omega\{2T\omega+\eta^{-1}\kappa[L\omega+(\nu+\eta\textrm{tr}\theta)\omega]+\omega T(\log\Omega)\}
\psi^{2}d\mu_{\slashed{g}}du^{\prime}
\end{equation*}
Taking into account the fact that:
\begin{equation}
 \textrm{tr}\chi=\eta(\textrm{tr}\slashed{k}-\textrm{tr}\theta),\quad \textrm{tr}\underline{\chi}=\kappa(\textrm{tr}\slashed{k}+\textrm{tr}\theta)
\end{equation}
we have
\begin{equation*}
 \textrm{tr}\underline{\chi}=\eta^{-1}
 \kappa[\textrm{tr}\chi+2\eta\textrm{tr}\theta]
\end{equation*}
Then we get
\begin{equation}
 I^{u}(t)=\int_{\Sigma^{u}_{t}}\frac{1}{2}\Omega(\underline{L}\omega+\underline{\nu}\omega)
 \psi^{2}d\mu_{\slashed{g}}du^{\prime}
\end{equation}
So (5.69) takes the form 
\begin{equation}
 \mathcal{E}^{u}_{1}(t)-\mathcal{E}^{\prime u}_{1}(t)=\int_{S_{t,u}}\frac{1}{2}\Omega\omega\psi^{2}d\mu_{\slashed{g}}
-I^{u}(t)
\end{equation}
Substituting this and (5.63) into (5.36), the surface integrals cancel and we get
\begin{equation}
 \mathcal{E}^{\prime u}_{1}(t)+\mathcal{F}^{\prime t}_{1}(u)=I^{u}(t)-I^{u}(0)+\mathcal{E}^{\prime u}_{1}(0)+\int_{W^{t}_{u}}Q_{1}d\mu_{g}
\end{equation}
Also for $K_{0}$ we have 
\begin{equation}
 \mathcal{E}^{u}_{0}(t)+\mathcal{F}^{t}_{0}(u)=\mathcal{E}^{u}_{0}(0)+\int_{W^{t}_{u}}Q_{0}d\mu_{g}
\end{equation}
Obviously, by $\textbf{A1}$ and $\textbf{A2}$ there is a positive constant $C$ such that
\begin{equation}
 C^{-1}\mathcal{E}^{u}_{0}(t)\leq \int_{0}^{u}\{\int_{S_{t,u^{\prime}}}[\mu(1+\mu)((L\psi)^{2}+|\slashed{d}
\psi|^{2})+(\underline{L}\psi)^{2}]d\mu_{\slashed{g}}\}du^{\prime}\leq C\mathcal{E}^{u}_{0}(t)
\end{equation}
 and
\begin{equation}
 C^{-1}\mathcal{F}^{t}_{0}(u)\leq \int_{0}^{t}\{\int_{S_{t^{\prime},u}}[(1+\mu)(L\psi)^{2}+\mu|\slashed{d}\psi|^{2}]d\mu_{\slashed{g}}
\}dt'\leq C\mathcal{F}^{t}_{0}(u)
\end{equation}
By $\textbf{A1}, \textbf{A2}$ and $\textbf{B1}, \textbf{D1}$ there is a positive constant $C$ such that:
\begin{equation}
 C^{-1}\mathcal{E}^{\prime u}_{1}(t)\leq (1+t)^{2}\int_{0}^{u}\{\int_{S_{t,u^{\prime}}}\mu[(L\psi+\nu\psi)^{2}+|\slashed{d}\psi|^{2}]
d\mu_{\slashed{g}}\}du^{\prime}\leq C\mathcal{E}^{\prime u}_{1}(t)
\end{equation}
 and
\begin{equation}
 C^{-1}\mathcal{F}^{\prime t}_{1}(u)\leq \int_{0}^{t}(1+t^{\prime})^{2}\{\int_{S_{t^{\prime},u}}(L\psi+\nu\psi)^{2}d\mu_{\slashed{g}}\}dt^{\prime}
\leq C\mathcal{F}^{\prime t}_{1}(u)
\end{equation}

$\textbf{Lemma 5.1}$ There is a numerical constant $C$ such that for all $u\in[0,\epsilon_{0}]$:
\begin{equation*}
 \int_{S_{t,u}}\psi^{2}d\mu_{\slashed{g}}\leq \epsilon_{0}C\mathcal{E}^{u}_{0}(t)
\end{equation*}
$Proof$. We shall use the acoustical coordinates $(t,u,\vartheta^{1},\vartheta^{2})$. Now on a given hyperplane $\Sigma_{t}$ we may set $\Xi=0$, then on this 
$\Sigma_{t}$ we have
\begin{equation*}
 T=\frac{\partial}{\partial u}
\end{equation*}
hence, in view of the fact that $\psi$ vanishes on $C_{0}$,
\begin{equation}
 \psi(t,u,\vartheta)=\int^{u}_{0}(T\psi)(t,u^{\prime},\vartheta)du^{\prime}
\end{equation}
It follows that:
\begin{align}
 \int_{S_{t,u}}\psi^{2}d\mu_{\tilde{\slashed{g}}}=\int_{S^{2}}\psi^{2}(t,u,\vartheta)d\mu_{\tilde{\slashed{g}}}(t,u,\vartheta)\\\notag
=\int_{S^{2}}\{\int_{0}^{u}(T\psi)(t,u^{\prime},\vartheta)du^{\prime}\}^{2}d\mu_{\tilde{\slashed{g}}}(t,u,\vartheta)\\\notag
\leq \epsilon_{0}\int_{S^{2}}\{\int_{0}^{u}(T\psi)^{2}(t,u^{\prime},\vartheta)du^{\prime}\}d\mu_{\tilde{\slashed{g}}}(t,u,\vartheta)
\end{align}
In view of the fact $\tilde{\slashed{g}}=\Omega\slashed{g}$,
\begin{align}
 \slashed{\mathcal{L}}_{T}d\mu_{\tilde{\slashed{g}}}=\Omega\slashed{\mathcal{L}}_{T}d\mu_{\slashed{g}}
+T(\Omega)d\mu_{\slashed{g}}=(\kappa\textrm{tr}\theta+T\log\Omega)d\mu_{\tilde{\slashed{g}}}
\end{align}
From (5.9) and (5.10) we have
\begin{equation}
 \kappa\textrm{tr}\theta+T\log\Omega=-\eta^{-1}\kappa\nu+\underline{\nu}
\end{equation}
Thus
\begin{equation}
 \slashed{\mathcal{L}}_{T}d\mu_{\tilde{\slashed{g}}}=(-\eta^{-1}\kappa\nu+\underline{\nu})d\mu_{\tilde{\slashed{g}}}
\end{equation}
From $\textbf{B1}, \textbf{B2}$ and $\textbf{A3}, \textbf{A2}$ we have
\begin{equation}
 |-\eta^{-1}\kappa\nu+\underline{\nu}|\leq C(1+t)^{-1}[1+\log(1+t)]^{4} \leq C^{\prime}
\end{equation}
Integrate (5.83) along $T$ on $\Sigma_{t}$ and using the fact $\epsilon_{0} \leq \frac{1}{2}$ we get
\begin{equation}
 C^{-1} \leq \frac{d\mu_{\tilde{\slashed{g}}(t,u,\vartheta)}}{d\mu_{\tilde{\slashed{g}}(t,0,\vartheta)}}\leq C
\end{equation}
for all $(u,\vartheta)\in [0,\epsilon_{0}]\times S^{2}$.
Then (5.80) is bounded by 
\begin{equation}
 C\epsilon_{0}\int_{0}^{u}\{\int_{S^{2}}(T\psi)^{2}(t,u^{\prime},\vartheta)d\mu_{\tilde{\slashed{g}}}(t,u^{\prime},\vartheta)\}du^{\prime}\leq
C\epsilon_{0}\int_{0}^{u}\{\int_{S_{t,u^{\prime}}}(T\psi)^{2}d\mu_{\slashed{g}}\}du^{\prime}
\end{equation}
Obviously
\begin{equation}
 (T\psi)^{2}\leq C[(\underline{L}\psi)^{2}+\eta^{-2}\kappa^{2}(L\psi)^{2}]
\end{equation}
The result follows. $\qed$

We now introduce the following quantities which are non-decreasing functions of $t$ at each $u$:
\begin{align}
 \bar{\mathcal{E}}^{u}_{0}(t)=\sup_{t'\in [0,t]}\mathcal{E}^{u}_{0}(t')\\
\mathcal{F}^{t}_{0}(u)\\
\bar{\mathcal{E}}^{\prime u}_{1}(t)=\sup_{t'\in[0,t]}[1+\log(1+t^{\prime})]^{-4}\mathcal{E}^{\prime u}_{1}(t^{\prime})\\
\bar{\mathcal{F}}^{\prime t}_{1}(u)=\sup_{t'\in[0,t]}[1+\log(1+t^{\prime})]^{-4}\mathcal{F}^{\prime t^{\prime}}_{1}(u)
\end{align}
Obviously, $\bar{\mathcal{E}}^{u}_{0}(t)$, $\bar{\mathcal{E}}^{\prime u}_{1}(t)$ are also non-decreasing functions of $u$ at each $t$. The statements $(\textbf{i})
-(\textbf{v})$ of $\textbf{Theorem 5.1}$ is equivalent to 
\begin{equation}
 \bar{\mathcal{E}}^{\epsilon_{0}}_{0}(s), \sup_{u\in [0,\epsilon_{0}]}\mathcal{F}^{s}_{0}(u), 
\bar{\mathcal{E}}^{\prime\epsilon_{0}}_{1}(s), \sup_{u\in [0,\epsilon_{0}]}\bar{\mathcal{F}}^{\prime s}_{1}(u) \leq CD_{0}
\end{equation}
 In view of (5.73) and (5.74), to prove the theorem we must estimate properly the spacetime integrals
\begin{equation}
 \int_{W^{t}_{u}}Q_{0}d\mu_{g}, \quad \int_{W^{t}_{u}}Q_{1}d\mu_{g}
\end{equation}
in terms of $\bar{\mathcal{E}}^{u}_{0}(t),\mathcal{F}^{t}_{0}(u),\bar{\mathcal{E}}^{\prime u}_{1}(t),\bar{\mathcal{F}}^{\prime t}_{1}(u)$. When
we estimate the error spacetime integral we shall see the reason why we allow the growth like $[1+\log(1+t)]^{4}$ for $\mathcal{E}^{\prime u}_{1}(t)$
and $\mathcal{F}^{\prime t}_{1}(u) $

\section{The Error Integrals}
From (5.25) and (5.30), we have:
\begin{align}
 Q_{0}=-\Omega^{2}\rho K_{0}\psi-\frac{1}{2}T^{\mu\nu}\tilde{\pi}_{0,\mu\nu}\\\notag
=Q_{0,0}+Q_{0,1}+Q_{0,2}+Q_{0,3}+Q_{0,4}+Q_{0,5}+Q_{0,6}+Q_{0,7}
\end{align}
where 
\begin{align}
 Q_{0,0}=-\Omega^{2}\rho K_{0}\psi\\
Q_{0,1}=-\frac{1}{2}T^{LL}\tilde{\pi}_{0,LL}
=-\frac{1}{8}\mu^{-2}(\underline{L}\psi)^{2}\tilde{\pi}_{0,LL}\\
Q_{0,2}=-\frac{1}{2}T^{\underline{L}\underline{L}}\tilde{\pi}_{0,\underline{L}\underline{L}}
=-\frac{1}{8}\mu^{-2}(L\psi)^{2}\tilde{\pi}_{0,\underline{L}\underline{L}}\\
Q_{0,3}=-T^{L\underline{L}}\tilde{\pi}_{0,L\underline{L}}
=-\frac{1}{4}\mu^{-1}|\slashed{d}\psi|^{2}\tilde{\pi}_{0,L\underline{L}}\\
Q_{0,4}=-T^{LA}\tilde{\pi}_{0,LA}
=\frac{1}{2}\mu^{-1}(\underline{L}\psi)(\slashed{d}^{A}\psi)\tilde{\pi}_{0,LA}\\
Q_{0,5}=-T^{\underline{L}A}\tilde{\pi}_{0,\underline{L}A}
=\frac{1}{2}\mu^{-1}(L\psi)(\slashed{d}^{A}\psi)\tilde{\pi}_{0,\underline{L}A}
\end{align}
and,
\begin{equation*}
 -\frac{1}{2}T^{AB}\tilde{\pi}_{0,AB}=Q_{0,6}+Q_{0,7}
\end{equation*}
where
\begin{align}
 Q_{0,6}=-\frac{1}{2}\{\slashed{d}^{A}\psi\slashed{d}^{B}\psi-\frac{1}{2}(\slashed{g}^{-1})^{AB}|\slashed{d}\psi|^{2}\}\hat{\tilde{\slashed{\pi}}}_{0,AB}\\
Q_{0,7}=-\frac{1}{4}\mu^{-1}(L\psi)(\underline{L}\psi)\textrm{tr}\tilde{\slashed{\pi}}_{0}
\end{align}
Here $\tilde{\slashed{\pi}}_{0}$ denotes the restriction to $S_{t,u}$ of $\tilde{\pi}_{0}$, and $\hat{\tilde{\slashed{\pi}}}_{0}$ denotes the trace-free
part of $\tilde{\slashed{\pi}}_{0}$ and we have made use of the trace-free nature of 
\begin{equation*}
 \slashed{d}^{A}\psi\slashed{d}^{B}\psi-\frac{1}{2}(\slashed{g}^{-1})^{AB}|\slashed{d}\psi|^{2}
\end{equation*}
We have
\begin{equation*}
 \tilde{\pi}_{0}=\Omega\pi_{0}+(K_{0}\Omega)g
\end{equation*}
Using (3.131)-(3.138) we have
\begin{align}
 \tilde{\pi}_{0,LL}=0\\
\tilde{\pi}_{0,\underline{L}\underline{L}}=-4\Omega\mu\{\underline{L}(\eta^{-1}\kappa)
-(\eta_{0}^{-1}+\alpha^{-1}\kappa)L(\eta^{-1}\kappa)\}\\
\tilde{\pi}_{0,L\underline{L}}=-2\Omega\mu\{\mu^{-1}(K_{0}\mu)+(K_{0}\log\Omega)+2L(\eta^{-1}\kappa)\}\\
\tilde{\pi}_{0,LA}=-2\Omega(\zeta_{A}+\eta_{A})\\
\tilde{\pi}_{0,\underline{L}A}=2\Omega\{(\eta_{0}^{-1}+\alpha^{-1}\kappa)(\zeta_{A}+\eta_{A})-\mu\slashed{d}_{A}(\eta^{-1}\kappa)\}\\
\hat{\tilde{\slashed{\pi}}}_{0,AB}=2\Omega\{(\eta^{-1}_{0}+\eta^{-1}\kappa)\hat{\chi}_{AB}+\hat{\underline{\chi}}_{AB}\}\\
\textrm{tr}\tilde{\slashed{\pi}}_{0}=4\Omega\{(\eta_{0}^{-1}+\eta^{-1}\kappa)\nu+\underline{\nu}\}
\end{align}
From (5.27) and (5.30) we have:
\begin{align}
 Q_{1}=-\Omega^{2}\rho K_{1}\psi-\Omega^{2}\rho\omega\psi-\frac{1}{2}T^{\mu\nu}\tilde{\pi}'_{1,\mu\nu}+\frac{1}{2}\Omega^{2}\psi^{2}\Box_{\tilde{g}}\omega\\\notag
=Q_{1,0}+Q_{1,1}+Q_{1,2}+Q_{1,3}+Q_{1,4}+Q_{1,5}+Q_{1,6}+Q_{1,7}+Q_{1,8}
\end{align}
where
\begin{align}
 Q_{1,0}=-\Omega^{2}\rho K_{1}\psi-\Omega^{2}\rho\omega\psi\\
Q_{1,1}=-\frac{1}{2}T^{LL}\tilde{\pi}^{\prime}_{1,LL}=
-\frac{1}{8}\mu^{-2}(\underline{L}\psi)^{2}\tilde{\pi}^{\prime}_{1,LL}\\
Q_{1,2}=-\frac{1}{2}T^{\underline{L}\underline{L}}\tilde{\pi}^{\prime}_{1,\underline{L}\underline{L}}=
-\frac{1}{8}\mu^{-2}(L\psi)^{2}\tilde{\pi}^{\prime}_{1,\underline{L}\underline{L}}\\
Q_{1,3}=-T^{L\underline{L}}\tilde{\pi}^{\prime}_{1,L\underline{L}}
=-\frac{1}{4}\mu^{-1}|\slashed{d}\psi|^{2}\tilde{\pi}^{\prime}_{1,L\underline{L}}\\
Q_{1,4}=-T^{LA}\tilde{\pi}^{\prime}_{1,LA}
=\frac{1}{2}\mu^{-1}(\underline{L}\psi)(\slashed{d}^{A}\psi)\tilde{\pi}^{\prime}_{1,LA}\\
Q_{1,5}=-T^{\underline{L}A}\tilde{\pi}^{\prime}_{1,\underline{L}A}
=\frac{1}{2}\mu^{-1}(L\psi)(\slashed{d}^{A}\psi)\tilde{\pi}^{\prime}_{1,\underline{L}A}
\end{align}
and
\begin{equation*}
 -\frac{1}{2}T^{AB}\tilde{\pi}^{\prime}_{1,AB}=Q_{1,6}+Q_{1,7}
\end{equation*}
where
\begin{align}
 Q_{1,6}=-\frac{1}{2}\{\slashed{d}^{A}\psi\slashed{d}^{B}\psi-\frac{1}{2}(\slashed{g}^{-1})^{AB}|\slashed{d}\psi|^{2}\}\hat{\tilde{\slashed{\pi}}}^{\prime}_{1,AB}\\
Q_{1,7}=-\frac{1}{4}\mu^{-1}(L\psi)(\underline{L}\psi)\textrm{tr}\tilde{\slashed{\pi}}^{\prime}_{1}
\end{align}
Also,
\begin{equation}
 Q_{1,8}=\frac{1}{2}\Omega^{2}\psi^{2}\Box_{\tilde{g}}\omega
\end{equation}
The symbols are similar as (5.101)-(5.102).

Since $\tilde{\pi}_{1}=\Omega\pi_{1}+(K_{1}\Omega)g$, we have $\tilde{\pi}^{\prime}_{1}=\Omega\{\pi_{1}+(K_{1}\log\Omega-2\omega)g\}$.
So we get the following:
\begin{align}
 \tilde{\pi}^{\prime}_{1,LL}=0\\
\tilde{\pi}^{\prime}_{1,\underline{L}\underline{L}}=
-4\Omega\mu\{\underline{L}(\nu^{-1}\omega)-\nu^{-1}\omega L(\eta^{-1}\kappa)\}\\
\tilde{\pi}^{\prime}_{1,L\underline{L}}=-2\Omega\mu\{\mu^{-1}K_{1}\mu+L(\nu^{-1}\omega)-2\omega+K_{1}\log\Omega\}\\
\tilde{\pi}^{\prime}_{1,LA}=0\\
\tilde{\pi}^{\prime}_{1,\underline{L}A}=2\Omega\{\nu^{-1}\omega(\zeta_{A}+\eta_{A})-\mu\slashed{d}_{A}(\nu^{-1}\omega)\}\\
\hat{\tilde{\slashed{\pi}}}^{\prime}_{1,AB}=2\Omega\nu^{-1}\omega\hat{\chi}_{AB}\\
\textrm{tr}\tilde{\slashed{\pi}}^{\prime}_{1}=0
\end{align}
At this point, we can see more clearly why we choose 
\begin{align*}
 K_{1}=(\omega/\nu)L
\end{align*}
To see this, let us first set:
\begin{align*}
 K_{1}=fL
\end{align*}
 where $f$ is a function to be determined. 

Recall (5.26):
\begin{align*}
 \tilde{\pi}^{\prime}_{1,AB}=\tilde{\pi}_{1,AB}-2\omega\tilde{\slashed{g}}_{AB}
\end{align*}
while
\begin{align*}
 \tilde{\pi}_{1,AB}=f\{\tilde{g}(D_{A}L,e_{B})+\tilde{g}(D_{B}L,e_{A})\}=2f\tilde{\chi}_{AB}\\
=2f\hat{\tilde{\chi}}_{AB}+f\textrm{tr}\tilde{\chi}\tilde{\slashed{g}}_{AB}
\end{align*}
Since
\begin{align*}
 \textrm{tr}\tilde{\chi}=2\nu
\end{align*}
we have:
\begin{align*}
 \textrm{tr}\tilde{\slashed{\pi}}^{\prime}_{1}=4(\nu f-\omega)
\end{align*}
So if we choose
\begin{align*}
 f=\omega/\nu
\end{align*}
then 
\begin{align*}
 \textrm{tr}\tilde{\slashed{\pi}}^{\prime}_{1}=0
\end{align*}
Hence $Q_{1,7}$ vanishes. This fact is very important, otherwise, the growth of $\omega$ would prevent us from appropriately estimating this term.

\section{The Estimates for the Error Integrals}
We now consider the spacetime integral of $Q_{0}$.

In view of (5.32), we have for an arbitrary function $f$,
\begin{equation}
 \int_{W^{t}_{u}}\mu^{-1}fd\mu_{g}=\int_{0}^{t}\{\int_{\Sigma_{t}^{\prime u}}f\}dt'=\int_{0}^{u}\{\int_{C^{t}_{u'}}f\}du'
\end{equation}
where
\begin{equation}
 \int_{\Sigma_{t}^{\prime u}}f=\int_{0}^{u}[\int_{S_{t',u'}}fd\mu_{\slashed{g}}]du',\quad \int_{C^{t}_{u'}}f
=\int_{0}^{t}[\int_{S_{t',u'}}fd\mu_{\slashed{g}}]dt'
\end{equation}
In view of the fact that $\rho$ vanishes and (5.103), 
\begin{equation}
 Q_{0,0}=Q_{0,1}=0
\end{equation}
By $\textbf{B10},\textbf{B11}$ and $\textbf{A1}$,
\begin{equation}
 \mu^{-1}|\tilde{\pi}_{0,\underline{L}\underline{L}}|\leq C(1+t)[1+\log(1+t)]^{-2}
\end{equation}
hence, from (5.97), we have
\begin{equation}
 \mu|Q_{0,2}|\leq C(1+t)[1+\log(1+t)]^{-2}(L\psi)^{2}
\end{equation}
Writing
\begin{equation}
 (L\psi)^{2}\leq 2(L\psi+\nu\psi)^{2}+2(\nu\psi)^{2}
\end{equation}
then
\begin{equation}
 \int_{W_{u}^{t}}|Q_{0,2}|d\mu_{g} \leq C(J_{0}+J_{1})
\end{equation}
where 
\begin{equation}
 J_{0}=\int_{0}^{t}(1+t')[1+\log(1+t')]^{-2}\{\int_{\Sigma_{t'}^{u}}(\nu\psi)^{2}\}dt'
\end{equation}
\begin{equation}
 J_{1}=\int_{0}^{t}(1+t')[1+\log(1+t')]^{-2}\{\int_{\Sigma_{t'}^{u}}(L\psi+\nu\psi)^{2}\}dt'
\end{equation}
By $\textbf{B1}$,
\begin{equation}
 J_{0}\leq C\int_{0}^{t}(1+t')^{-1}[1+\log(1+t')]^{-2}\{\int_{\Sigma_{t'}^{u}}\psi^{2}\}dt'
\end{equation}
Using Lemma 5.1 we obtain:
\begin{equation}
 J_{0}\leq C\int_{0}^{t}(1+t')^{-1}[1+\log(1+t')]^{-2}\mathcal{E}^{u}_{0}(t')dt'
\end{equation}
On the other hand,
\begin{equation}
 J_{1}=\int_{0}^{t}f_{0}(t')\frac{dg(t')}{dt'}dt'
\end{equation}
where
\begin{equation}
 g(t)=\int_{0}^{t}(1+t')^{2}\{\int_{\Sigma^{u}_{t'}}(L\psi+\nu\psi)^{2}\}dt'
\end{equation}
and
\begin{equation}
 f_{0}(t)=(1+t)^{-1}[1+\log(1+t)]^{-2}
\end{equation}
Integrating by parts:
\begin{equation}
 J_{1}=f_{0}(t)g(t)-\int_{0}^{t}g(t')\frac{df_{0}(t')}{dt'}dt'
\end{equation}
By (5.76)
\begin{equation}
 \int_{0}^{t}\int_{S_{t',u}}(\textsl{L}\psi)^{2}d\mu_{\slashed{g}}dt' \leq C\mathcal{F}^{t}_{0}(u)
\end{equation}
Hence
\begin{equation}
 \int_{0}^{t}\{\int_{\Sigma_{t'}^{u}}(L\psi)^{2}\}dt'\leq C\int_{0}^{u}\mathcal{F}_{0}^{t}(u')du'
\end{equation}
Also by Lemma 5.1 and $\textbf{B1}$,
\begin{equation}
 \int_{0}^{t}\int_{\Sigma_{t'}^{u}}(\nu\psi)^{2}d\mu_{\slashed{g}}dt'\leq C\epsilon_{0}\int_{0}^{t}(1+t')^{-2}
\sup_{u'\in [0,u]}(\int_{S_{t',u'}}\psi^{2}d\mu_{\slashed{g}})dt'\leq C\epsilon_{0}^{2}
\int_{0}^{t}(1+t')^{-2}\mathcal{E}^{u}_{0}(t')dt'
\end{equation}
Set:
\begin{equation}
 F(t,u)=\{\int_{\Sigma_{t}^{u}}(L\psi+\nu\psi)^{2}\}^{\frac{1}{2}}
\end{equation}
We then have:
\begin{equation}
 g(t)=\int_{0}^{t}(1+t')^{2}F^{2}(t',u)dt'
\end{equation}
So 
\begin{align}
f_{0}(t)g(t)\leq [1+\log(1+t)]^{-2}\int_{0}^{t}(1+t')F^{2}(t',u)dt'\\\notag
\leq \{\int_{0}^{t}F^{2}(t',u)dt'\}^{\frac{1}{2}}\{[1+\log(1+t)]^{-4}\int_{0}^{t}(1+t')^{2}F^{2}(t',u)dt'\}^{\frac{1}{2}} 
\end{align}
By (5.143) and (5.144),
\begin{equation}
 \int_{0}^{t}F^{2}(t',u)dt' \leq C\int_{0}^{u}\mathcal{F}_{0}^{t}(u')du'+C\epsilon_{0}^{2}\int_{0}^{t}(1+t')^{-2}\mathcal{E}^{u}_{0}(t')dt'
\end{equation}
On the other hand, by (5.78) we have:
\begin{equation}
 \int_{0}^{t}(1+t')^{2}F^{2}(t',u)dt'\leq C\int_{0}^{u}\mathcal{F}'^{t}_{1}(u')du'
\end{equation}
So we have:
\begin{equation}
 f_{0}(t)g(t)\leq C\{\int_{0}^{u}\mathcal{F}_{0}^{t}(u')du'+C\epsilon_{0}^{2}\int_{0}^{t}(1+t')^{-2}\mathcal{E}^{u}_{0}(t')dt'\}^{\frac{1}{2}}
\{[1+\log(1+t)]^{-4}\int_{0}^{u}\mathcal{F}'^{t}_{1}(u')du'\}^{\frac{1}{2}}
\end{equation}
Moreover, since
\begin{equation}
 |\frac{df_{0}}{dt}|\leq C(1+t)^{-2}[1+\log(1+t)]^{-2}
\end{equation}
we have by (5.146) and (5.149),
\begin{equation}
 \int_{0}^{t}g(t')|\frac{df_{0}}{dt}|dt'\leq C\int_{0}^{t}(1+t')^{-2}[1+\log(1+t']^{-2}\{\int_{0}^{u}\mathcal{F}'^{t'}_{1}(u')du'\}dt'
\end{equation}
We thus obtain
\begin{align}
 J_{1}\leq C\{\int_{0}^{u}\mathcal{F}_{0}^{t}(u')du'+C\epsilon_{0}^{2}\int_{0}^{t}(1+t')^{-2}\mathcal{E}^{u}_{0}(t')dt'\}^{\frac{1}{2}}
\{[1+\log(1+t)]^{-4}\int_{0}^{u}\mathcal{F}'^{t}_{1}(u')du'\}^{\frac{1}{2}}\\\notag
+ C\int_{0}^{t}(1+t')^{-2}[1+\log(1+t')]^{-2}\{\int_{0}^{u}\mathcal{F}'^{t'}_{1}(u')du'\}dt'
\end{align}
which together with (5.137) yields
\begin{align}
 \int_{W^{t}_{u}} |Q_{0,2}|d\mu_{g}\leq C\int_{0}^{t}(1+t')^{-1}[1+(\log(1+t')]^{-2}\mathcal{E}^{u}_{0}(t')dt'\\\notag
+C\{\int_{0}^{u}\mathcal{F}_{0}^{t}(u')du'+C\epsilon_{0}^{2}\int_{0}^{t}(1+t')^{-2}\mathcal{E}^{u}_{0}(t')dt'\}^{\frac{1}{2}}\
\{[1+\log(1+t)]^{-4}\int_{0}^{u}\mathcal{F}'^{t}_{1}(u')du'\}^{\frac{1}{2}}\\\notag
+ C\int_{0}^{t}(1+t')^{-2}[1+\log(1+t')]^{-2}\{\int_{0}^{u}\mathcal{F}'^{t'}_{1}(u')du'\}dt'
\end{align}
We now turn to $Q_{0,6}$ and $Q_{0,7}$. The terms $Q_{0,3}, Q_{0,4}, Q_{0,5}$ shall be treated later.
 
From $\textbf{B3},\textbf{B4}$ and $\textbf{A1},\textbf{A3}$, we have
\begin{equation}
 |\hat{\tilde{\slashed{\pi}}}_{0}| \leq C(1+t)[1+\log(1+t)]^{-6}
\end{equation}
In view of (5.77) and (5.90) we can then estimate
\begin{align}
 \int_{W^{t}_{u}}|Q_{0,6}|d\mu_{g} \leq C\int_{0}^{t}(1+t')^{-1}[1+\log(1+t')]^{-6}\mathcal{E}'^{u}_{1}(t')dt'\\\notag
\leq C\int_{0}^{t}(1+t')^{-1}[1+\log(1+t')]^{-2}\bar{\mathcal{E}}^{\prime u}_{1}(t')dt'
\end{align}
Also, by $\textbf{B1},\textbf{B2}, \textbf{A1},\textbf{A3}$, we have
\begin{equation}
 |\textrm{tr}\tilde{\slashed{\pi}}_{0}|\leq C(1+t)^{-1}[1+\log(1+t)]^{4}
\end{equation}
Writing $|L\psi|\leq |L\psi+\nu\psi|+|\nu\psi|$, we can then estimate: 
\begin{equation}
 \int_{W^{t}_{u}}|Q_{0,7}|d\mu_{g}\leq C(J_{2}+J_{3})
\end{equation}
where:
\begin{align}
 J_{2}=\int_{0}^{t}f_{1}(t')\int_{\Sigma_{t'}^{u}}|\nu\psi||\underline{L}\psi|\\
J_{3}=\int_{0}^{t}f_{1}(t')\int_{\Sigma_{t'}^{u}}|L\psi+\nu\psi||\underline{L}\psi|
\end{align}
and 
\begin{equation}
 f_{1}(t)=(1+t)^{-1}[1+\log(1+t)]^{4}
\end{equation}
We have:
\begin{equation}
 J_{2}\leq I_{0}^{\frac{1}{2}}I_{1}^{\frac{1}{2}}
\end{equation}
where:
\begin{align}
 I_{0}=\int_{0}^{t}(1+t')f_{1}(t')\{\int_{\Sigma_{t'}^{u}}(\nu\psi)^{2}\}dt'\\
I_{1}=\int_{0}^{t}(1+t')^{-1}f_{1}(t')\{\int_{\Sigma_{t'}^{u}}(\underline{L}\psi)^{2}\}dt'
\end{align}
By $\textbf{B1}$ and Lemma 5.1, we have ,with
\begin{equation}
 f_{2}(t)=(1+t)^{-1}f_{1}(t)=(1+t)^{-2}[1+\log(1+t)]^{4}
\end{equation}
\begin{equation}
 I_{0}\leq C\epsilon_{0}\int_{0}^{t}f_{2}(t')\mathcal{E}^{u}_{0}(t')dt'
\end{equation}
Also, by (5.75),
\begin{equation}
 I_{1}\leq C\int_{0}^{t}f_{2}(t')\mathcal{E}^{u}_{0}(t')dt'
\end{equation}
Thus we obtain:
\begin{equation}
 J_{2}\leq C\epsilon^{\frac{1}{2}}_{0}\int_{0}^{t}f_{2}(t')\mathcal{E}^{u}_{0}(t')dt'
\end{equation}
We have:
\begin{equation}
 J_{3}\leq \frac{1}{2}(I_{1}+I_{2})
\end{equation}
where:
\begin{equation}
 I_{2}=\int_{0}^{t}(1+t')f_{1}(t')\{\int_{\Sigma^{u}_{t'}}(L\psi+\nu\psi)^{2}\}dt'
\end{equation}
Recalling (5.139), we have:
\begin{equation}
 I_{2}=\int_{0}^{t}f_{2}(t')\frac{dg(t')}{dt'}dt'
\end{equation}
We can see that $I_{2}$ is similar to $J_{1}$, moreover, $f_{2}$ decays faster than $f_{0}$, thus $I_{2}$ has the same bound as $J_{1}$:
\begin{align}
 I_{2}\leq C\{\int_{0}^{u}\mathcal{F}_{0}^{t}(u')du'+C\epsilon_{0}^{2}\int_{0}^{t}(1+t')^{-2}\mathcal{E}^{u}_{0}(t')dt'\}^{\frac{1}{2}}
\{[1+\log(1+t)]^{-4}\int_{0}^{u}\mathcal{F}'^{t}_{1}(u')du'\}^{\frac{1}{2}}\\\notag
+ C\int_{0}^{t}(1+t')^{-2}[1+\log(1+t')]^{-2}\{\int_{0}^{u}\mathcal{F}^{\prime t'}_{1}(u')du'\}dt'
\end{align}
We conclude that
\begin{align}
\int_{W^{t}_{t}}|Q_{0,7}|d\mu_{g} \leq C\int_{0}^{t}f_{2}(t')\mathcal{E}^{u}_{0}(t')dt'\\\notag
+C\{\int_{0}^{u}\mathcal{F}_{0}^{t}(u')du'+C\epsilon_{0}^{2}\int_{0}^{t}(1+t')^{-2}\mathcal{E}^{u}_{0}(t')dt'\}^{\frac{1}{2}}
\{[1+\log(1+t)]^{-4}\int_{0}^{u}\mathcal{F}^{\prime t}_{1}(u')du'\}^{\frac{1}{2}}\\\notag
+ C\int_{0}^{t}(1+t')^{-2}[1+\log(1+t')]^{-2}\{\int_{0}^{u}\mathcal{F}^{\prime t'}_{1}(u')du'\}dt'
\end{align}
We turn to $Q_{1}$. Since we are considering the case that $\rho$ vanishes, we have:
\begin{equation}
 Q_{1,0}=0
\end{equation}
Also, from (5.120),(5.123) and (5.126) we have
\begin{equation}
 Q_{1,1}=Q_{1,4}=Q_{1,7}=0
\end{equation}
So we first consider $Q_{1,2}$. By $\textbf{B1},\textbf{B10},\textbf{B13},\textbf{D1},\textbf{D3},\textbf{A1}$, we have
\begin{equation}
 \mu^{-1}|\tilde{\pi}'_{1,\underline{L}\underline{L}}|\leq C(1+t)[1+\log(1+t)]^{3}
\end{equation}
Hence,
\begin{equation}
 \mu|Q_{1,2}|\leq C(1+t)[1+\log(1+t)]^{3}(L\psi)^{2}
\end{equation}
Writing again: $(L\psi)^{2} \leq 2(L\psi+\nu\psi)^{2}+2(\nu\psi)^{2}$, it then follows:
\begin{equation}
 \int_{W^{t}_{u}}|Q_{1,2}|d\mu_{g}\leq C(J_{4}+J_{5})
\end{equation}
where
\begin{align}
 J_{4}=\int_{0}^{t}(1+t')[1+\log(1+t')]^{3}\{\int_{\Sigma_{t'}^{u}}(\nu\psi)^{2}\}dt'\\
J_{5}=\int_{0}^{t}(1+t')[1+\log(1+t')]^{3}\{\int_{\Sigma_{t'}^{u}}(L\psi+\nu\psi)^{2}\}dt'
\end{align}
By $\textbf{B1}$ and Lemma 5.1 we have:
\begin{equation}
 J_{4}\leq C\epsilon_{0}^{2}\int_{0}^{t}(1+t')^{-1}[1+\log(1+t')]^{3}\mathcal{E}^{u}_{0}(t')dt'
\leq C\epsilon^{2}_{0}\bar{\mathcal{E}}^{u}_{0}(t)[1+\log(1+t)]^{4}
\end{equation}
From (5.139) we have
\begin{equation}
 J_{5}=\int_{0}^{t}f_{3}(t')\frac{dg(t')}{dt'}dt'
\end{equation}
where
\begin{equation}
 f_{3}(t)=(1+t)^{-1}[1+\log(1+t)]^{3}
\end{equation}
Then we have
\begin{equation}
 J_{5}=f_{3}(t)g(t)-\int_{0}^{t}g(t')\frac{df_{3}(t')}{dt'}dt'
\leq g(t)\{f_{3}(t)+\int_{0}^{t}|\frac{df_{3}(t')}{dt'}|dt'\}\leq C\int_{0}^{u}\mathcal{F}'^{t}_{1}(u')du'
\end{equation}
We have used the fact that $g(t)$ is non-decreasing. In conclusion,
\begin{equation}
 \int_{W^{t}_{u}}|Q_{1,2}|d\mu_{g}\leq C\int_{0}^{u}\mathcal{F}'^{t}_{1}(u')du'+C\epsilon^{2}_{0}\bar{\mathcal{E}}^{u}_{0}(t)[1+\log(1+t)]^{4}
\end{equation}
Consider next $Q_{1,6}$. By $\textbf{A1},\textbf{B1},\textbf{B3},\textbf{D1}$, we have
\begin{equation}
 |\hat{\tilde{\slashed{\pi}}}'_{1}|\leq C(1+t)[1+\log(1+t)]^{-2}
\end{equation}
Hence, by (5.77),
\begin{align}
 \int_{W^{t}_{u}}|Q_{1,6}|d\mu_{g}\leq C\int_{0}^{t}(1+t')[1+\log(1+t')]^{-2}\{\int_{\Sigma_{t'}^{u}}\mu|\slashed{d}\psi|^{2}\}dt'\\\notag
\leq C\int_{0}^{t}(1+t')^{-1}[1+\log(1+t')]^{-2}\mathcal{E}'^{u}_{1}(t')dt'
\end{align}
To estimate $Q_{1,8}$, we shall use $\textbf{D5}$ and Lemma 5.1:
\begin{align}
 \int_{W^{t}_{u}}|Q_{1,8}|d\mu_{g}\leq \int_{0}^{t}\{\int_{0}^{u}\sup_{S_{t',u'}}(\mu|\Box_{\tilde{g}}\omega|)
\cdot[\int_{S_{t',u'}}\psi^{2}d\mu_{\slashed{g}}]du'dt'\\\notag
\leq C\epsilon_{0}\int_{0}^{t}\{\int_{0}^{u}\sup_{S_{t',u'}}(\mu|\Box_{\tilde{g}}\omega|)du'\}\mathcal{E}^{u}_{0}(t')dt'
\end{align}
 So from $\textbf{D5}$ we have 
\begin{equation}
 \int_{W^{t}_{u}}|Q_{1,8}|d\mu_{g}\leq C\epsilon_{0}[1+\log(1+t)]^{4}\bar{\mathcal{E}}^{u}_{0}(t)
\end{equation}
We now turn to the crucial terms, $Q_{0,3}$ and $Q_{1,3}$.

Consider first $Q_{1,3}$. We have:
\begin{equation}
 \tilde{\pi}^{\prime}_{1,L\underline{L}}=-2\Omega\mu(\omega\nu^{-1})(\mu^{-1}L\mu+r_{1})
\end{equation}
where 
\begin{equation}
 r_{1}=\omega^{-1}(L\omega-\nu\omega)-\nu^{-1}(L\nu+\nu^{2})+L\log\Omega
\end{equation}
Thus,
\begin{equation}
 \int_{W^{t}_{u}}Q_{1,3}d\mu_{g}=\int_{0}^{t}\{\int_{\Sigma_{t'}^{u}}\frac{\Omega}{2}(\omega\nu^{-1})(\mu^{-1}L\mu+r_{1})\mu|\slashed{d}\psi|^{2}\}dt'
\end{equation}
Decomposing $\mu^{-1}L\mu$ into its positive and negative parts, we write:
\begin{equation}
 \int_{W^{t}_{u}}Q_{1,3}d\mu_{g}=\int_{0}^{t}\{\int_{\Sigma_{t'}^{u}}\frac{\Omega}{2}(\omega\nu^{-1})[\mu^{-1}(L\mu)_{+}+\mu^{-1}(L\mu)_{-}+r_{1}]
\mu|\slashed{d}\psi|^{2}\}dt'
\end{equation}
Then by $\textbf{C1}$ and (5.64) we have
\begin{equation}
 \int_{W^{t}_{u}}Q_{1,3}d\mu_{g}\leq \int_{0}^{t}\{(1+t')^{-1}[1+\log(1+t')]^{-1}+A(t')
+\sup_{\Sigma^{u}_{t'}}|r_{1}|\}\mathcal{E}^{\prime u}_{1}(t')dt'-K(t,u)
\end{equation}
where $K(t,u)$ is the non-negative spacetime integral:
\begin{equation}
 K(t,u)=-\int_{W^{t}_{u}}\frac{\Omega}{2}\omega\nu^{-1}\mu^{-1}(L\mu)_{-}|\slashed{d}\psi|^{2}d\mu_{g}
\end{equation}
Moreover, by $\textbf{B1},\textbf{D1},\textbf{B5},\textbf{B12},\textbf{D2}$ we have:
\begin{equation}
 \sup_{\Sigma^{\epsilon_{0}}_{t}}|r_{1}|\leq C(1+t)^{-1}[1+\log(1+t)]^{-2}
\end{equation}
Consider next $Q_{0,3}$. From (5.105) and (5.18) we have:
\begin{equation}
 \tilde{\pi}_{0,L\underline{L}}=-2\Omega\mu(\mu^{-1}(L\mu+\underline{L}\mu)+r_{0})
\end{equation}
where
\begin{align}
 r_{0}=\alpha^{-2}L\mu+2L(\alpha^{-1}\kappa)+K_{0}\log\Omega\\\notag
=3L(\alpha^{-1}\kappa)+2\alpha^{-2}\kappa L\alpha+(1+\alpha^{-1}\kappa)L\log\Omega+\underline{L}\log\Omega
\end{align}
Thus,
\begin{equation}
 \int_{W^{t}_{u}}Q_{0,3}d\mu_{g}=\int_{0}^{t}\{\int_{\Sigma_{t'}^{u}}\frac{\Omega}{2}(\mu^{-1}(L\mu+\underline{L}\mu)
+r_{0})\mu|\slashed{d}\psi|^{2}\}dt'
\end{equation}
Decomposing $\mu^{-1}(L\mu+\underline{L}\mu)$ into its positive and negative parts, and then dropping the negative part, we obtain:
\begin{equation}
 \int_{W^{t}_{u}}Q_{0,3}d\mu_{g}\leq \int_{0}^{t}\{\int_{\Sigma_{t'}^{u}}\frac{\Omega}{2}(\mu^{-1}(L\mu+\underline{L}\mu)_{+}
+|r_{0}|)\mu|\slashed{d}\psi|^{2}\}dt'
\end{equation}
Then by $\textbf{C2}$ and (5.64) we have:
\begin{equation}
 \int_{W^{t}_{u}}Q_{0,3}d\mu_{g}\leq \int_{0}^{t}(1+t')^{-2}[B(t')+\sup_{\Sigma_{t'}^{u}}|r_{0}|]\mathcal{E}'^{u}_{1}(t')dt'
\end{equation}
Moreover, by $\textbf{B5},\textbf{B6},\textbf{B9},\textbf{B10}$ we have:
\begin{equation}
 \sup_{\Sigma_{t}^{\epsilon_{0}}}|r_{0}|\leq C(1+t)[1+\log(1+t)]^{-6}
\end{equation}

Finally, we consider $Q_{0,4}$,$Q_{0,5}$ and $Q_{1,5}$.

First, we have:
\begin{equation}
 \int_{W^{t}_{u}}|Q_{1,5}|d\mu_{g}\leq M_{1}+R_{1}
\end{equation}
where
\begin{align}
 M_{1}=\int_{W^{t}_{u}}\Omega(\omega\nu^{-1})|L\psi||\slashed{d}\psi||\zeta+\eta|\mu^{-1}d\mu_{g}\\
R_{1}=\int_{W^{t}_{u}}\Omega|L\psi||\slashed{d}\psi||\slashed{d}(\omega\nu^{-1})|d\mu_{g}
\end{align}
First, we estimate $M_{1}$. We decompose:
\begin{equation}
 M_{1}=M'_{1}+M''_{1}
\end{equation}
where
\begin{align}
 M_{1}^{\prime}=\int_{\mathcal{U}\bigcap W_{u}^{t}}\Omega(\omega\nu^{-1})|\textsl{L}\psi||\slashed{d}\psi||\zeta+\eta|\mu^{-1}d\mu_{g}\\
M_{1}^{\prime\prime}=\int_{\mathcal{U}^{c}\bigcap W_{u}^{t}}\Omega(\omega\nu^{-1})|\textsl{L}\psi||\slashed{d}\psi||\zeta+\eta|\mu^{-1}d\mu_{g}
\end{align}
The region $\mathcal{U}$ is defined by (5.14).

According to $\textbf{C3}$ we have:
\begin{equation}
 -(L\mu)_{-}\geq C^{-1}(1+t')^{-1}[1+\log(1+t')]^{-1}
\end{equation}
in $\mathcal{U}\bigcap W_{u}^{t}$.

Comparing with (5.195) we have
\begin{align}
 K(t,u)\geq -\int_{\mathcal{U}\bigcap W_{u}^{t}}\frac{\Omega}{2}\omega\nu^{-1}\mu^{-1}(L\mu)_{-}|\slashed{d}\psi|^{2}d\mu_{g}\\\notag
\geq \frac{1}{2C}\int_{\mathcal{U}\bigcap W_{u}^{t}}\Omega\omega\nu^{-1}(1+t')^{-1}[1+\log(1+t')]^{-1}\mu^{-1}|\slashed{d}\psi|^{2}d\mu_{g}
\end{align}
Thus we can estimate:
\begin{equation}
 M^{\prime}_{1}\leq CK^{\frac{1}{2}}N_{1}^{\frac{1}{2}}
\end{equation}
where: 
\begin{equation}
 N_{1}=\int_{W_{u}^{t}}\Omega\omega\nu^{-1}(1+t')[1+\log(1+t')]|\zeta+\eta|^{2}|L\psi|^{2}\mu^{-1}d\mu_{g}
\end{equation}
By virtue of $\textbf{B7}$,
\begin{equation}
 N_{1}\leq \int_{W_{u}^{t}}\Omega\omega\nu^{-1}(1+t')^{-1}[1+\log(1+t')]^{3}|L\psi|^{2}\mu^{-1}d\mu_{g}
\end{equation}
Thus 
\begin{equation}
 N_{1}\leq C(N_{1,0}+N_{1,1})
\end{equation}
where
\begin{align}
 N_{1,0}=\int_{\mathcal{U}\bigcap W_{u}^{t}}\Omega\omega\nu^{-1}(1+t')^{-1}[1+\log(1+t')]^{3}|\nu\psi|^{2}\mu^{-1}d\mu_{g}\\
N_{1,1}=\int_{\mathcal{U}\bigcap W_{u}^{t}}\Omega\omega\nu^{-1}(1+t')^{-1}[1+\log(1+t')]^{3}|L\psi+\nu\psi|^{2}\mu^{-1}d\mu_{g}
\end{align}
Now by $\textbf{B1},\textbf{D1}$ and $\textbf{A1}$,
\begin{equation}
 N_{1,0}\leq C\int_{0}^{t}(1+t')^{-1}[1+\log(1+t')]^{3}\{\int_{0}^{u}[\int_{S_{t',u'}}\psi^{2}d\mu_{\slashed{g}}]du'\}dt'
\end{equation}
hence by Lemma 5.1,
\begin{align}
 N_{1,0}\leq C\epsilon_{0}^{2}\int_{0}^{t}(1+t')^{-1}[1+\log(1+t')]^{3}\mathcal{E}^{u}_{0}(t')dt'\\\notag
\leq C\epsilon_{0}^{2}\bar{\mathcal{E}}^{u}_{0}(t)[1+\log(1+t)]^{4}
\end{align}
Also, from (5.59),
\begin{equation}
 N_{1,1}\leq C\int_{0}^{u}\{\int_{C_{u'}^{t}}\Omega\omega\nu^{-1}(L\psi+\nu\psi)^{2}\}du'=C\int_{0}^{u}\mathcal{F}'^{t}_{1}(u')du'
\end{equation}
We thus obtain: 
\begin{equation}
 N_{1}\leq C\epsilon_{0}^{2}\bar{\mathcal{E}}^{u}_{0}(t)[1+\log(1+t)]^{4}+C\int_{0}^{u}\mathcal{F}^{\prime t}_{1}(u')du'
\end{equation}
We conclude that
\begin{equation}
 M'_{1}\leq CK(t,u)^{\frac{1}{2}}\{\epsilon_{0}^{2}\bar{\mathcal{E}}^{u}_{0}(t)[1+\log(1+t)]^{4}+C\int_{0}^{u}\mathcal{F}^{\prime t}_{1}(u')du'\}^{\frac{1}{2}}
\end{equation}
On the other hand, since in $\mathcal{U}^{c}\bigcap W_{u}^{t}$ we have $\mu\geq 1/4$,
\begin{align}
 M''_{1}\leq 2\int_{\mathcal{U}^{c}\bigcap W_{u}^{t}}\Omega\omega\nu^{-1}|\textsl{L}\psi||\slashed{d}\psi||\zeta+\eta|\mu^{-\frac{1}{2}}d\mu_{g}\\\notag
\leq 2\int_{W_{u}^{t}}\Omega\omega\nu^{-1}|\textsl{L}\psi||\slashed{d}\psi||\zeta+\eta|\mu^{-\frac{1}{2}}d\mu_{g}\\\notag
\leq 2\{\int_{W_{u}^{t}}(1+t')^{-1}[1+\log(1+t')]^{-1}\Omega\omega\nu^{-1}|\slashed{d}\psi|^{2}d\mu_{g}\}^{\frac{1}{2}}N_{1}
^{\frac{1}{2}}
\end{align}
where $N_{1}$ is given by (5.212). In view of (5.64),
\begin{equation}
 M''_{1}\leq C\{\int_{0}^{t}(1+t')^{-1}[1+\log(1+t')]^{-1}\mathcal{E}'^{u}_{1}(t')dt'\}^{\frac{1}{2}}N_{1}
^{\frac{1}{2}}
\end{equation}
So by (5.220) we conclude:
\begin{equation}
 M''_{1}\leq C\{\int_{0}^{t}(1+t')^{-1}[1+\log(1+t')]^{-1}\mathcal{E}'^{u}_{1}(t')dt'\}^{\frac{1}{2}}
\{\epsilon_{0}^{2}\bar{\mathcal{E}}^{u}_{0}(t)[1+\log(1+t)]^{4}+C\int_{0}^{u}\mathcal{F}'^{t}_{1}(u')du'\}^{\frac{1}{2}}
\end{equation}
In the following we estimate $R_{1}$ given by (5.205). Using $\textbf{B14},\textbf{D4},\textbf{B1},\textbf{D1}$, we have
\begin{equation}
 \omega^{-1}\nu|\slashed{d}(\omega\nu^{-1})|\leq C(1+t)^{-1}[1+\log(1+t)]^{\frac{1}{2}}
\end{equation}
so $\mu^{\frac{1}{2}}\omega^{-1}\nu|\slashed{d}(\omega\nu^{-1})|$ enjoys the same bound as $|\zeta+\eta|$, namely by $C(1+t)^{-1}[1+\log(1+t)]$. It follows 
that $R_{1}$ is bounded in the same way as $M''_{1}$.
So we conclude:
\begin{align}
 \int_{W_{u}^{t}}|Q_{1,5}|d\mu_{g}\leq C\{\epsilon_{0}^{2}\bar{\mathcal{E}}^{u}_{0}(t)[1+\log(1+t)]^{4}+C\int_{0}^{u}\mathcal{F}'^{t}_{1}(u')du'\}
^{\frac{1}{2}}\\\notag
\{K(t,u)+\int_{0}^{t}(1+t')^{-1}[1+\log(1+t')]^{-1}\mathcal{E}'^{u}_{1}(t')dt'\}^{\frac{1}{2}}
\end{align}

Next we consider $Q_{0,4}$. Obviously, we have
\begin{equation}
 \int_{W_{u}^{t}}|Q_{0,4}|d\mu_{g}\leq M_{0}
\end{equation}
where
\begin{equation}
 M_{0}=\int_{W_{u}^{t}}\Omega|\underline{L}\psi||\slashed{d}\psi||\zeta+\eta|\mu^{-1}d\mu_{g}
\end{equation}
Again, we decompose:
\begin{equation}
 M_{0}=M'_{0}+M''_{0}
\end{equation}
where
\begin{align}
 M'_{0}=\int_{\mathcal{U}\bigcap W_{u}^{t}}\Omega|\underline{L}\psi||\slashed{d}\psi||\zeta+\eta|\mu^{-1}d\mu_{g}\\
M''_{0}=\int_{\mathcal{U}^{c}\bigcap W_{u}^{t}}\Omega|\underline{L}\psi||\slashed{d}\psi||\zeta+\eta|\mu^{-1}d\mu_{g}
\end{align}
By $\textbf{B7}$, we estimate:
\begin{align}
 \int_{\mathcal{U}\bigcap W^{t}_{u}}\Omega|\underline{L}\psi||\slashed{d}\psi||\zeta+\eta|\mu^{-1}d\mu_{g}
\leq C\int_{\mathcal{U}\bigcap W^{t}_{u}}(1+t^{\prime})^{-1}[1+\log(1+t^{\prime})]|\underline{L}\psi||\slashed{d}\psi|d\mu_{\slashed{g}}du^{\prime}dt^{\prime}
\end{align}
By Holder Inequality, the right hand side of the above is bounded by:
\begin{align}
 C\int_{0}^{t}(1+t^{\prime})^{-1}[1+\log(1+t^{\prime})]\sqrt{\mathcal{E}^{u}_{0}(t^{\prime})}\sqrt{\int_{\mathcal{U}\bigcap\Sigma_{t^{\prime}}^{u}}
|\slashed{d}\psi|^{2}}dt^{\prime}\\\notag
\leq (\int_{0}^{t}(1+t^{\prime})^{-2}[1+\log(1+t^{\prime})]^{2}\mathcal{E}^{u}_{0}(t^{\prime})dt^{\prime})^{1/2}\cdot
(\int_{\mathcal{U}\bigcap W^{t}_{u}}|\slashed{d}\psi|^{2}d\mu_{\slashed{g}}du^{\prime}dt^{\prime})^{1/2}
\end{align}
Let us define:
\begin{align}
 F(t,u):=\int_{0}^{t}(1+t^{\prime})[1+\log(1+t^{\prime})]^{-1}(\int_{\mathcal{U}\bigcap\Sigma_{t^{\prime}}^{u}}|\slashed{d}\psi|^{2})dt^{\prime}
\end{align}
By (5.210) and $\textbf{B1},\textbf{D1}$,
\begin{equation}
 K(t,u)\geq \frac{1}{C}\int_{\mathcal{U}\bigcap W_{u}^{t}}\Omega(1+t')[1+\log(1+t')]^{-1}\mu^{-1}|\slashed{d}\psi|^{2}d\mu_{g}
\end{equation}

Let us define
\begin{align*}
 \bar{K}(t,u)=\sup_{t'\in[0,t]}[1+\log(1+t')]^{-4}K(t',u)
\end{align*}
Therefore
\begin{align}
 F(t,u)\leq CK(t,u)\leq C\bar{K}(t,u)[1+\log(1+t)]^{4}
\end{align}
and
\begin{align}
 \frac{dF}{dt}(t,u)=(1+t)[1+\log(1+t)]^{-1}\int_{\mathcal{U}\bigcap \Sigma_{t}^{u}}|\slashed{d}\psi|^{2}
\end{align}
What we need to estimate is the following integral:
\begin{align}
 I(t):=\int_{0}^{t}(\int_{\mathcal{U}\bigcap\Sigma_{t^{\prime}}^{u}}|\slashed{d}\psi|^{2})dt^{\prime}
\end{align}
We have:
\begin{align}
 I(t)=\int_{0}^{t}(1+t^{\prime})^{-1}[1+\log(1+t^{\prime})]\frac{dF}{dt^{\prime}}(t^{\prime})dt^{\prime}\\\notag
=(1+t)^{-1}[1+\log(1+t)]F(t)+\int_{0}^{t}F(t^{\prime})\{(1+t^{\prime})^{-2}[1+\log(1+t^{\prime})]-(1+t^{\prime})^{-2}\}dt^{\prime}
\end{align}
Substitute (5.236) we obtain:
\begin{align}
 I(t)\leq C^{\prime}\bar{K}(t,u)
\end{align}
Therefore
\begin{align}
 M^{\prime}_{0}\leq C\bar{K}(t,u)^{1/2}\{\int_{0}^{t}(1+t^{\prime})^{-2}[1+\log(1+t^{\prime})]^{2}\bar{\mathcal{E}}^{u}_{0}(t^{\prime})dt^{\prime}\}^{1/2}
\end{align}

To estimate $M''_{0}$ we note that since $\mu\geq \eta_{0}/4$ in $\mathcal{U}^{c}$:
\begin{align}
 M''_{0}\leq 2\int_{\mathcal{U}^{c}\bigcap W_{u}^{t}}\Omega|\underline{L}\psi||\slashed{d}\psi||\zeta+\eta|\mu^{-1/2}d\mu_{g}\\\notag
\leq 2\int_{W_{u}^{t}}\Omega|\underline{L}\psi||\slashed{d}\psi||\zeta+\eta|\mu^{-1/2}d\mu_{g}\\\notag
=2\int_{0}^{t}\{\int_{\Sigma_{t'}^{u}}\Omega|\underline{L}\psi||\slashed{d}\psi||\zeta+\eta|\mu^{1/2}\}dt'
\end{align}
Hence 
\begin{align}
 M''_{0}\leq 2\int_{0}^{t}\{\int_{\Sigma_{t'}^{u}}\Omega\mu|\slashed{d}\psi|^{2}\}^{1/2}
\cdot\{\int_{\Sigma_{t'}^{u}}\Omega|\underline{L}\psi|^{2}|\zeta+\eta|^{2}\}^{1/2}dt'\\\notag
\leq C\int_{0}^{t}\frac{\mathcal{E}'^{u}_{1}(t')^{1/2}}{(1+t')}\cdot\frac{[1+\log(1+t')]\mathcal{E}^{u}_{0}(t')^{1/2}}{(1+t')}dt'\\\notag
\leq C\{\int_{0}^{t}(1+t')^{-2}[1+\log(1+t')]^{3}\bar{\mathcal{E}}^{\prime u}_{1}(t')dt'\}^{1/2}
\cdot\{\int_{0}^{t}(1+t')^{-2}[1+\log(1+t')]^{3}\bar{\mathcal{E}}^{u}_{0}(t')dt'\}^{1/2}
\end{align}
From (5.227),(5.229),(5.241) and (5.243) we conclude that:
\begin{align}
 \int_{W_{u}^{t}}|Q_{0,4}|d\mu_{g}\leq C\bar{K}(t,u)^{\frac{1}{2}}(\int_{0}^{t}(1+t')^{-2}[1+\log(1+t')]^{2}\bar{\mathcal{E}}^{u}_{0}(t')dt')^{\frac{1}{2}}\\\notag
+C\{\int_{0}^{t}(1+t')^{-2}[1+\log(1+t')]^{3}\bar{\mathcal{E}}^{\prime u}_{1}(t')dt'\}^{1/2}
\cdot\{\int_{0}^{t}(1+t')^{-2}[1+\log(1+t')]^{3}\bar{\mathcal{E}}^{u}_{0}(t')dt'\}^{1/2}
\end{align}
We are left with $Q_{0,5}$. We have:
\begin{equation}
 \int_{W_{u}^{t}}|Q_{0,5}|d\mu_{g}\leq \tilde{M_{0}}+R_{0}
\end{equation}
where
\begin{align}
 \tilde{M_{0}}=\int_{W_{u}^{t}}\Omega|L\psi||\slashed{d}\psi|(1+\eta^{-1}\kappa)|\zeta+\eta|\mu^{-1}d\mu_{g}\\\notag
R_{0}=\int_{W_{u}^{t}}\Omega|L\psi||\slashed{d}\psi||\slashed{d}(\eta^{-1}\kappa)|d\mu_{g}
\end{align}
Decomposing:
\begin{align}
 \tilde{M}_{0}=\tilde{M}'_{0}+\tilde{M}''_{0}\\\notag
\tilde{M}'_{0}=\int_{\mathcal{U}\bigcap W^{t}_{u}}\Omega|L\psi||\slashed{d}\psi|(1+\eta^{-1}\kappa)|\zeta+\eta|\mu^{-1}d\mu_{g}\\\notag
\tilde{M}''_{0}=\int_{\mathcal{U}^{c}\bigcap W^{t}_{u}}\Omega|L\psi||\slashed{d}\psi|(1+\eta^{-1}\kappa)|\zeta+\eta|\mu^{-1}d\mu_{g}
\end{align}
Similarly we have:
\begin{align}
 \int_{\mathcal{U}\bigcap W^{t}_{u}}\Omega|L\psi||\slashed{d}\psi|(1+\eta^{-1}\kappa)|\zeta+\eta|d\mu_{\slashed{g}}du^{\prime}dt^{\prime}\\\notag
\leq \int_{0}^{t}(1+t^{\prime})^{-1}[1+\log(1+t^{\prime})]^{3/2}\sqrt{\mathcal{E}^{u}_{0}(t^{\prime})}\sqrt{\int_{\mathcal{U}\bigcap \Sigma_{t^{\prime}}^{u}}
|\slashed{d}\psi|^{2}}dt^{\prime}
\end{align}
This is similar to (5.233).
We then proceed as we estimate $M^{\prime}_{0}$ to obtain:
\begin{equation}
 \tilde{M}'_{0}\leq C\bar{K}(t,u)^{1/2}(\int_{0}^{u}\mathcal{F}^{t}_{0}(u')du')^{1/2}
\end{equation}
Next we estimate $\tilde{M}''_{0}$. Since $\mu\geq 1/4$ in $\mathcal{U}^{c}\bigcap W_{u}^{t}$, we have, by$\textbf{A2},\textbf{A3}$ and $\textbf{B7}$,
\begin{align}
 \tilde{M}''_{0}\leq \frac{2}{\sqrt{\eta_{0}}}\int_{\mathcal{U}\bigcap W_{u}^{t}}\Omega|L\psi||\slashed{d}\psi||\zeta+\eta|
(1+\eta^{-1}\kappa)\mu^{-1/2}d\mu_{g}\\\notag
\leq C\int_{W_{u}^{t}}(1+t')^{-1}[1+\log(1+t')]^{3/2}\Omega|L\psi||\slashed{d}\psi|(1+\eta^{-1}\kappa)^{1/2}\mu^{-1/2}d\mu_{g}
\end{align}
Since the factor $(1+t')^{-1}[1+\log(1+t')]^{3/2}$ is bounded by a numerical constant we obtain:
\begin{equation}
 \tilde{M}''_{0}\leq C\int_{W_{u}^{t}}\Omega\{(1+\eta^{-1}\kappa)(L\psi)^{2}+\mu|\slashed{d}\psi|^{2}\}\mu^{-1}d\mu_{g}
\end{equation}
hence in view of (5.54) we get
\begin{equation}
 \tilde{M}''_{0}\leq C\int_{0}^{u}\mathcal{F}^{t}_{0}(u')du'
\end{equation}
Finally we estimate $R_{0}$.
By $\textbf{B8}$ we have:
\begin{equation}
 R_{0}\leq C\int_{W_{u}^{t}}(1+t')^{-1}[1+\log(1+t')]\Omega|L\psi||\slashed{d}\psi|(1+\eta^{-1}\kappa)^{1/2}\mu^{-1/2}d\mu_{g}
\end{equation}
We have used the fact that $(1+\eta^{-1}\kappa)^{1/2}\mu^{-1/2}=((1+\eta^{-2}\mu)/\mu)^{1/2}\geq \eta^{-1}\geq C^{-1}$.

Since the factor $(1+t')^{-1}[1+\log(1+t')]$ is bounded by a numerical constant we again obtain:
\begin{equation}
 R_{0}\leq C\int_{W_{u}^{t}}\Omega\{(1+\eta^{-1}\kappa)(L\psi)^{2}+\mu|\slashed{d}\psi|^{2}\}\mu^{-1}d\mu_{g}
\leq C\int_{0}^{u}\mathcal{F}_{0}^{t}(u')du'
\end{equation}
So we conclude that 
\begin{equation}
 \int_{W_{u}^{t}}|Q_{0,5}|d\mu_{g}\leq C\bar{K}(t,u)^{1/2}(\int_{0}^{u}\mathcal{F}^{t}_{0}(u')du')^{1/2}+C\int_{0}^{u}\mathcal{F}_{0}^{t}(u')du'
\end{equation}

\section{Treatment of the Integral Inequalities Depending on $t$ and $u$. Completion of the Proof}
We now focus on the identity (5.73). By (5.174),(5.175),(5.185),(5.187),(5.189),(5.194),(5.196) and (5.226), the spacetime integral on right hand side of 
(5.73) is bounded from above by:
\begin{align}
 \int_{W_{u}^{t}}Q_{1}d\mu_{g}\leq CM(t,u)+L(t,u)+\int_{0}^{t}\tilde{A}(t')\mathcal{E}'^{u}_{1}(t')dt'\\\notag
-K(t,u)+C(K(t,u)^{1/2}+L(t,u)^{1/2})M(t,u)^{1/2}
\end{align}
Here,
\begin{align}
 M(t,u)=\bar{\mathcal{E}}^{u}_{0}(t)[1+\log(1+t)]^{4}+\int_{0}^{u}\mathcal{F}'^{t}_{1}(u')du'\\\notag
L(t,u)=\int_{0}^{t}(1+t')^{-1}[1+\log(1+t')]^{-1}\mathcal{E}'^{u}_{1}(t')dt'\\\notag
\tilde{A}(t)=A(t)+C(1+t)^{-1}[1+\log(1+t)]^{-2}
\end{align}
Now from $\textbf{C1}$ we have:
\begin{align}
 \int_{0}^{t}\tilde{A}(t')dt'\leq C 
\end{align}
where $C$ is independent of $t$.

Using the inequalities:
\begin{align}
 -K+CK^{1/2}M^{1/2}\leq -\frac{1}{2}K+CM,\quad CL^{1/2}M^{1/2}\leq \frac{1}{2}L+CM
\end{align}
we get the following:
\begin{align}
 \int_{W_{u}^{t}}Q_{1}d\mu_{g}\leq -\frac{1}{2}K(t,u)+CM(t,u)+\frac{3}{2}L(t,u)+\int_{0}^{t}\tilde{A}(t')\mathcal{E}'^{u}_{1}(t')dt'
\end{align}
We have:
\begin{align}
 L(t,u)\leq \int_{0}^{t}(1+t')^{-1}[1+\log(1+t')]^{3}\bar{\mathcal{E}}'^{u}_{1}(t')dt'\leq \frac{1}{4}
[1+\log(1+t)]^{4}\bar{\mathcal{E}}'^{u}_{1}(t)
\end{align}
Also, since 
\begin{equation}
 \mathcal{F}^{\prime t}_{1}(u)\leq [1+\log(1+t)]^{4}\bar{\mathcal{F}}'^{t}_{1}(u)
\end{equation}
defining:
\begin{equation}
 V'_{1}(t,u)=\int_{0}^{u}\bar{\mathcal{F}}'^{t}_{1}(u')du'
\end{equation}
we have:
\begin{equation}
 \int_{0}^{u}\mathcal{F}^{\prime t}_{1}(u')du'\leq [1+\log(1+t)]^{4}V'_{1}(t,u)
\end{equation}
Note that $V'_{1}(t,u)$ is a non-decreasing function of $t$ at each $u$ as well as a non-decreasing function of $u$ at each $t$.
Hence:
\begin{equation}
 M(t,u)\leq [1+\log(1+t)]^{4}(\bar{\mathcal{E}}^{u}_{0}(t)+V'_{1}(t,u))
\end{equation}
and also,
\begin{equation}
\int_{0}^{t}\tilde{A}(t')\mathcal{E}'^{u}_{1}(t')dt'\leq [1+\log(1+t)]^{4}\int_{0}^{t}\tilde{A}(t')\bar{\mathcal{E}}'^{u}_{1}(t')dt'      
\end{equation}
From $\textbf{A1},\textbf{B2},\textbf{D1},\textbf{D3}$ and Lemma 5.1, the space-like hypersurface integrals on the right hand side of (5.73) are bounded by:
\begin{equation}
 |\int_{\Sigma_{t}^{u}}\frac{1}{2}\Omega(\underline{L}\omega+\underline{\nu}\omega)\psi^{2}|\leq C\bar{\mathcal{E}}^{u}_{0}(t)[1+\log(1+t)]^{4},\quad
|\int_{\Sigma_{0}^{u}}\frac{1}{2}\Omega(\underline{L}\omega+\underline{\nu}\omega)\psi^{2}|\leq CD_{0}
\end{equation}
Also, by (5.75) and (5.77) at $t=0$, the remaining term on the right hand side of (5.73) is bounded by:
\begin{equation}
 \mathcal{E}'^{u}_{1}(0)\leq CD_{0}
\end{equation}
 In view of (5.260)-(5.268), the identity (5.73) implies
\begin{align}
 \mathcal{E}^{\prime u}_{1}(t)+\mathcal{F}^{\prime t}_{1}(u)+\frac{1}{2}K(t,u)\\\notag
\leq [1+\log(1+t)]^{4}\{\frac{3}{8}\bar{\mathcal{E}}^{\prime u}_{1}(t)+C(\bar{\mathcal{E}}^{u}_{0}(t)+V'_{1}(t,u))
+\int_{0}^{t}\tilde{A}(t')\bar{\mathcal{E}}^{\prime u}_{1}(t')dt'\}+CD_{0}
\end{align}
Keeping only the first term on the left we have:
\begin{align}
 [1+\log(1+t)]^{-4}\mathcal{E}^{\prime u}_{1}(t)\leq \frac{3}{8}\bar{\mathcal{E}}^{\prime u}_{1}(t)+C(\bar{\mathcal{E}}^{u}_{0}(t)+V'_{1}(t,u))
+\int_{0}^{t}\tilde{A}(t')\bar{\mathcal{E}}^{\prime u}_{1}(t')dt'+CD_{0}
\end{align}
The same holds with $t$ replaced by $t'\in [0,t]$. Since the right hand side is a non-decreasing function of $t$ at each $u$, we deduce:
\begin{equation}
 \bar{\mathcal{E}}^{\prime u}_{1}(t)\leq C(\bar{\mathcal{E}}^{u}_{0}(t)+V'_{1}(t,u))
+\int_{0}^{t}\tilde{A}(t')\bar{\mathcal{E}}^{\prime u}_{1}(t')dt'+CD_{0}
\end{equation}
Since $\bar{\mathcal{E}}^{u}_{0}(t)+V'_{1}(t,u)$ is a non-decreasing function of $t$ at each $u$ 
we can use Gronwall's Inequality to obtain in view of (5.259),
\begin{equation}
 \bar{\mathcal{E}}'^{u}_{1}(t)\leq C(\bar{\mathcal{E}}^{u}_{0}(t)+V'_{1}(t,u))+CD_{0}
\end{equation}
hence also:
\begin{equation}
 \int_{0}^{t}\tilde{A}(t')\bar{\mathcal{E}}^{\prime u}_{1}(t')dt'\leq \bar{\mathcal{E}}^{\prime u}_{1}(t)\int_{0}^{t}\tilde{A}(t')dt'
\leq C(\bar{\mathcal{E}}^{u}_{0}(t)+V'_{1}(t,u))
\end{equation}
Substituting (5.272)-(5.273) to (5.269) and keeping only the second term on the left, we get
\begin{equation}
 [1+\log(1+t)]^{-4}\mathcal{F}^{\prime t}_{1}(u)\leq C(\bar{\mathcal{E}}^{u}_{0}(t)+V'_{1}(t,u))+CD_{0}
\end{equation}
The same holds with $t$ replaced by $t^{\prime}\in[0,t]$. Since the right hand side is a non-decreasing function of $t$ at each $u$, we deduce:
\begin{equation}
 \bar{\mathcal{F}}^{\prime t}_{1}(u)\leq C(\bar{\mathcal{E}}^{u}_{0}(t)+V'_{1}(t,u))+CD_{0}
\end{equation}
In view of the definition of $V'_{1}(t,u)$, this is
\begin{equation}
 \bar{\mathcal{F}}'^{t}_{1}(u)\leq C\bar{\mathcal{E}}^{u}_{0}(t)+C\int_{0}^{u}\bar{\mathcal{F}}'^{t}_{1}(u')du'+CD_{0}
\end{equation}
In view of the fact that $\bar{\mathcal{E}}^{u}_{0}(t)$ is a non-decreasing function of $u$ at each $t$ and $[0,\epsilon_{0}]$ is a bounded interval,
we can use Gronwall's inequality to get
\begin{equation}
 \bar{\mathcal{F}}'^{t}_{1}(u)\leq C\bar{\mathcal{E}}^{u}_{0}(t)+CD_{0}
\end{equation}
hence also:
\begin{equation}
 V'_{1}(t,u)\leq C\epsilon_{0}\bar{\mathcal{E}}^{u}_{0}(t)+CD_{0}
\end{equation}
Substituting (5.278) into (5.272) we obtain:
\begin{equation}
 \bar{\mathcal{E}}^{\prime u}_{1}(t)\leq C\bar{\mathcal{E}}^{u}_{0}(t)+CD_{0}
\end{equation}
Substituting (5.273), (5.278) and (5.279) into (5.269) and keeping only the third term on the left we obtain:
\begin{equation}
 [1+\log(1+t)]^{-4}K(t,u)\leq C\bar{\mathcal{E}}^{u}_{0}(t)+CD_{0}
\end{equation}
which implies:
\begin{equation}
 \bar{K}(t,u)\leq C\bar{\mathcal{E}}^{u}_{0}(t)+CD_{0}
\end{equation}

We now turn to the identity (5.74). By (5.129), (5.154), (5.156), (5.173), (5.201)-(5.202), (5.244) and (5.255), the spacetime integral 
on the right hand side of (5.74) is bounded from above by, in view of (5.264): 
\begin{align}
 \int_{W_{u}^{t}}Q_{0}d\mu_{g}\leq \int_{0}^{t}(1+t')^{-2}[1+\log(1+t')]^{4}B(t')\bar{\mathcal{E}}^{\prime u}_{1}(t')dt'\\\notag
+C\int_{0}^{t}(1+t')^{-1}[1+\log(1+t')]^{-2}(\bar{\mathcal{E}}^{\prime u}_{1}(t')+\bar{\mathcal{E}}^{u}_{0}(t'))dt'+CV_{0}(t,u)\\\notag
+C\{V_{0}(t,u)+\int_{0}^{t}(1+t')^{-1}[1+\log(1+t')]^{-2}\bar{\mathcal{E}}^{\prime u}_{0}(t')dt'\}^{1/2}(V'_{1}(t,u))^{1/2}\\\notag
+C\int_{0}^{t}(1+t')^{-2}[1+\log(1+t')]^{2}V'_{1}(t',u)dt'\\\notag
+C\bar{K}(t,u)^{1/2}(\int_{0}^{t}(1+t')^{-1}[1+\log(1+t')]^{-2}\bar{\mathcal{E}}^{u}_{0}(t')dt')^{1/2}\\\notag
+C\bar{K}(t,u)^{1/2}(V_{0}(t,u))^{1/2}
\end{align}
Here we have defined:
\begin{equation}
 V_{0}(t,u)=\int_{0}^{u}\mathcal{F}^{t}_{0}(u')du'
\end{equation}

We now substitute (5.278), (5.279) and (5.281) to the above. In doing this we estimate the fourth term on the right by:
\begin{align}
 C\bar{\mathcal{E}}^{u}_{0}(t)^{1/2}\{V_{0}(t,u)+\int_{0}^{t}(1+t')^{-1}[1+\log(1+t')]^{-2}\bar{\mathcal{E}}'^{u}_{0}(t')dt'\}^{1/2}\\\notag
\leq \frac{\delta}{2}\bar{\mathcal{E}}^{u}_{0}(t)+\frac{C^{2}}{2\delta}\{V_{0}(t,u)+\int_{0}^{t}(1+t')^{-1}[1+\log(1+t')]^{-2}
\bar{\mathcal{E}}^{\prime u}_{0}(t')dt'\}
\end{align}
and the sixth term by:
\begin{align}
 C\bar{\mathcal{E}}^{u}_{0}(t)^{1/2}(\int_{0}^{t}(1+t')^{-1}[1+\log(1+t')]^{-2}\bar{\mathcal{E}}^{u}_{0}(t')dt')^{1/2}\\\notag
\leq \frac{\delta}{2}\bar{\mathcal{E}}^{u}_{0}(t)+\frac{C^{2}}{2\delta}(\int_{0}^{t}(1+t')^{-1}[1+\log(1+t')]^{-2}\bar{\mathcal{E}}^{u}_{0}(t')dt')
\end{align}
and the seventh term by:
\begin{align}
 C\bar{\mathcal{E}}^{u}_{0}(t)^{1/2}(V_{0}(t,u))^{1/2}\leq \frac{\delta}{2}\bar{\mathcal{E}}^{u}_{0}(t)+\frac{C^{2}}{2\delta}V_{0}(t,u)
\end{align}
After these substitutions, we obtain:
\begin{align}
 \int_{W_{u}^{t}}Q_{0}d\mu_{g}\leq \frac{3\delta}{2}\bar{\mathcal{E}}^{u}_{0}(t)+C(1+\frac{1}{\delta})V_{0}(t,u)+C(1+\frac{1}{\delta})\int_{0}^{t}\tilde{B}(t')
\bar{\mathcal{E}}^{u}_{0}(t')dt'+CD_{0}
\end{align}
The constants $C$ are independent of $\delta$. Here:
\begin{align}
 \tilde{B}(t)=(1+t)^{-2}[1+\log(1+t)]^{4}B(t)+C(1+t)^{-1}[1+\log(1+t)]^{-2}
\end{align}
By $\textbf{C2}$ we have:
\begin{equation}
 \int_{0}^{t}\tilde{B}(t')dt'\leq C
\end{equation}
In view of (5.287) the identity (5.74) implies:
\begin{align}
 \mathcal{E}^{u}_{0}(t)+\mathcal{F}^{t}_{0}(u)\leq \mathcal{E}^{u}_{0}(0)+\frac{3\delta}{2}\bar{\mathcal{E}}^{u}_{0}(t)+C(1+\frac{1}{\delta})
\{V_{0}(t,u)+\int_{0}^{t}\tilde{B}(t')\bar{\mathcal{E}}^{u}_{0}(t')dt'\}
\end{align}
We set: $\delta=\frac{1}{3}$. Replacing $t$ by $t'\in[0,t]$ in the above, noting that the right hand side is a non-decreasing function of $t$ at each $u$, 
and keeping only the first term on the left, we obtain:
\begin{align}
 \bar{\mathcal{E}}^{u}_{0}(t)\leq \mathcal{E}^{u}_{0}(0)+CV_{0}(t,u)+C\int_{0}^{t}\tilde{B}(t')\bar{\mathcal{E}}^{u}_{0}(t')dt'
\end{align}
Since $V_{0}(t,u)$ is a non-decreasing function of $t$ at each $u$, we can apply Gronwall's inequality to obtain:
\begin{align}
 \bar{\mathcal{E}}^{u}_{0}(t)\leq C(\mathcal{E}^{u}_{0}(0)+V_{0}(t,u))
\end{align}
hence also:
\begin{align}
 \int_{0}^{t}\tilde{B}(t')\bar{\mathcal{E}}^{u}_{0}(t')dt'\leq \bar{\mathcal{E}}^{u}_{0}(t)\int_{0}^{t}\tilde{B}(t')dt'
\leq C(\mathcal{E}^{u}_{0}(0)+V_{0}(t,u))
\end{align}
Keeping only the second term on the left in (5.290) we obtain:
\begin{align}
 \mathcal{F}^{t}_{0}(u)\leq C(\mathcal{E}^{u}_{0}(0)+V_{0}(t,u))=C(\mathcal{E}^{u}_{0}(0)+\int_{0}^{u}\mathcal{F}^{t}_{0}(u')du')
\end{align}
In view that $\mathcal{E}^{u}_{0}(0)$ is a non-decreasing function of $u$ and $[0,\epsilon_{0}]$ is a bounded interval, we can use Gronwall's inequality to obtain:
\begin{equation}
 \mathcal{F}^{t}_{0}(u)\leq C\mathcal{E}^{u}_{0}(0)
\end{equation}
hence also:
\begin{equation}
 V_{0}(t,u)\leq C\epsilon_{0}\mathcal{E}^{u}_{0}(0)
\end{equation}
Substituting this into (5.292) we obtain:
\begin{equation}
 \bar{\mathcal{E}}^{u}_{0}(t)\leq C\mathcal{E}^{u}_{0}(0)
\end{equation}
Substituting finally (5.297) to (5.277), (5.279), (5.281) we conclude that: 
\begin{equation}
 \bar{\mathcal{E}}^{\prime u}_{1}(t),\bar{\mathcal{F}}^{\prime t}_{1}(u),\bar{K}(t,u)\leq C\mathcal{E}^{u}_{0}(0)
\end{equation}
The proof of $\textbf{Theorem 5.1}$ is now complete. $\qed$

We conclude this chapter by noting certain key features of the above estimates.

(i) When we estimate $Q_{1,3}$, we do not just bound the contribution from 
$L\mu$ simply in absolute value. Instead, we split $L\mu$ into positive and negative 
parts. Because by $\textbf{C3}$ and (5.195), we can use the negative part of $L\mu$ to estimate the tangential derivative of variations on the region 
$\mathcal{U}$, where $\mu$ is very small. In estimating $Q_{0,3}$, we also split the term $L\mu+\underline{L}\mu$ into positive part and negative 
part so that we can use the bootstrap assumption $\textbf{C1}$ and $\textbf{C2}$.

(ii) Let us revisit the estimate (5.218) for $N_{1,0}$, which is original from $Q_{1,5}$. There is a growth of $[1+\log(1+t)]^{4}$ which can not be
improved on the right hand side, this is why we introduce the weight for $\mathcal{E}^{\prime u}_{1}(t)$ and $\mathcal{F}^{\prime t}_{1}(u)$, and also
for $K(t,u)$ see (5.90), (5.91) and (5.235). Then we can complete the error estimates associated to $K_{1}$, since also this weight is enough to bound the
 contribution of $L(t,u)$ (see (5.257)), whose growth in time is like $\log[\log(1+t)]$.


(iii) Finally, the treatment in section 5.5 is not standard, because we have two variables $t, u$ in the integral inequalities with $u$ in a bounded interval. 
However, we can use Gronwall inequality with respect to one variable when the other is fixed to complete the estimates.

\chapter{Construction of Commutation Vectorfields}



\section{Commutation Vectorfields and Their Deformation Tensors}

In this chapter, we shall construct the vectorfields $Y_{i}: i=1,2,3,4,5$ used to define the higher order variations of the wave function $\phi$. That is,
\begin{equation}
 \psi_{n}=Y_{i_{1}}...Y_{i_{n-1}}\psi_{1}
\end{equation}
where $\psi_{1}$ is a first order variation, namely one of the functions (5.5). Here $i_{1},...,i_{n-1} \in {1,2,3,4,5}$. Since $\psi_{1}$ satisfies:
\begin{equation}
 \Box_{\tilde{g}}\psi_{1}=0
\end{equation}
$\psi_{n}$ satisfies
\begin{equation}
 \Box_{\tilde{g}}\psi_{n}=\rho_{n}
\end{equation}
where $\rho_{n}$ is obtained by commuting the vectorfields $Y_{i_{1}}...Y_{i_{n-1}}$ with $\Box_{\tilde{g}}$. So we call $Y_{i}: i=1,2,3,4,5$ the 
$commutation$  $vectorfields$.

We require that, at each point, these commutation vectors span the tangent space to the spacetime manifold at the corresponding point. Thus, the set
of all $\psi_{n}$ for fixed $n$ corresponding to a given $\psi_{1}$ contains all derivatives of $\psi_{1}$ of order $n-1$.

We take
\begin{equation}
 Y_{1}=T
\end{equation}
Since this is transversal to $C_{u}$, we require $Y_{i}:i=2,3,4,5$ to be tangential to $C_{u}$. Moreover, for each $u\in[0,\epsilon_{0}]$ and each point 
$x\in C_{u}$, the set $\{Y_{i}(x):i=2,3,4,5\}$ spans the tangent space to $C_{u}$ at $x$.

Next, we take
\begin{equation}
 Y_{2}=Q:=(1+t)L
\end{equation}
In view of the fact that the range of $u$ on $W^{*}_{\epsilon_{0}}$ is the bounded interval $[0,\epsilon_{0}]$, $Q$ is an analogue in the acoustical
spacetime of the dilation field $D$ of flat spacetime introduced in Chapter 1.
Since $Y_{2}$ is transversal to the surfaces $S_{t,u}$, we require that $Y_{i}: i=3,4,5$ to be tangential to $S_{t,u}$. Moreover,
for each $t,u$ and each point $x\in S_{t,u}$, the set $\{Y_{i}(x):i=3,4,5\}$ spans the tangent plane to $S_{t,u}$

We set:
\begin{equation}
 Y_{i+2}=R_{i}:=\Pi \mathring{R}_{i}: i=1,2,3
\end{equation}
Here $\mathring{R}_{i}$ are the generators of rotations about the three rectangular coordinate axes:
\begin{equation}
 \mathring{R}_{i}=\epsilon_{ijk}x^{j}\frac{\partial}{\partial x^{k}}=\frac{1}{2}\epsilon_{ijk}\Omega_{jk}
\end{equation}
where $\epsilon_{ijk}$ is the fully antisymmetric 3-dimensional symbol. In (6.6) $\Pi$ is the projection to the tangent plane to the surfaces $S_{t,u}$
with respect to the induced acoustical metric $\bar{g}$ on $\Sigma_{t}$, which, as we have seen, coincides with the Euclidean metric.

When we consider the commutator $[Y,\Box_{\tilde{g}}]$, the deformation tensor $\leftexp{(Y)}{\tilde{\pi}}$ appears. Since:
\begin{equation}
 \leftexp{(Y)}{\tilde{\pi}}=\Omega\leftexp{(Y)}{\pi}+(Y\Omega)g
\end{equation}
we must consider the expression for $\leftexp{(Y)}{\pi}$:
\begin{equation}
 \leftexp{(Y)}{\pi}(Z_{1},Z_{2})=g(D_{Z_{1}}Y,Z_{2})+g(D_{Z_{2}}Y,Z_{1})
\end{equation}
The components of the deformation tensors $\leftexp{(T)}{\pi}, \leftexp{(Q)}{\pi}$ in the null frame 
$(L, \underline{L},X_{1},X_{2})$ can be directly computed from (3.99)-(3.106) and (3.131)-(3.138).
\begin{align}
 \leftexp{(T)}{\tilde{\pi}_{LL}}=0\\
\leftexp{(T)}{\tilde{\pi}}_{\underline{L}\underline{L}}=4\Omega\mu T(\eta^{-1}\kappa)\\
\leftexp{(T)}{\tilde{\pi}}_{\underline{L}L}=-2\Omega(T\mu+\mu T\log\Omega)\\
\leftexp{(T)}{\tilde{\pi}}_{LA}=-\Omega(\zeta_{A}+\eta_{A})\\
\leftexp{(T)}{\tilde{\pi}}_{\underline{L}A}=-\Omega\alpha^{-1}\kappa(\zeta_{A}+\eta_{A})\\
\leftexp{(T)}{\hat{\tilde{\slashed{\pi}}}}_{AB}=\Omega(\underline{\hat{\chi}}_{AB}-\eta^{-1}\kappa\hat{\chi}_{AB})\\
\textrm{tr}\leftexp{(T)}{\tilde{\slashed{\pi}}}=2\Omega(\underline{\nu}-\eta^{-1}\kappa\nu)
\end{align}
and:
\begin{align}
 \leftexp{(Q)}{\tilde{\pi}}_{LL}=0\\
\leftexp{(Q)}{\tilde{\pi}}_{\underline{L}\underline{L}}=4\Omega\mu\{Q(\eta^{-1}\kappa)-\eta^{-1}\kappa\}\\
\leftexp{(Q)}{\tilde{\pi}}_{L\underline{L}}=-2\Omega\{Q\mu+\mu Q\log\Omega+\mu\}\\
\leftexp{(Q)}{\tilde{\pi}}_{LA}=0\\
\leftexp{(Q)}{\tilde{\pi}}_{\underline{L}A}=2\Omega(1+t)(\zeta_{A}+\eta_{A})\\
\leftexp{(Q)}{\hat{\tilde{\slashed{\pi}}}}_{AB}=2\Omega(1+t)\hat{\chi}_{AB}\\
\textrm{tr}\leftexp{(Q)}{\tilde{\slashed{\pi}}}=4\Omega(1+t)\nu
\end{align}
Here we denote by $\leftexp{(Y)}{\tilde{\slashed{\pi}}}$ the restriction of $\leftexp{(Y)}{\tilde{\pi}}$ to $S_{t,u}$, and by
$\leftexp{(Y)}{\hat{\tilde{\slashed{\pi}}}}$ the trace-free part of $\leftexp{(Y)}{\tilde{\slashed{\pi}}}$.

Noting that the definition of $R_{i}: i=1,2,3$ is intrinsic to each space-like hypersurface $\Sigma_{t}$, we shall derive expressions for the components of 
$\leftexp{(R_{i})}{\pi}$ in the frame $(L,T,X_{1},X_{2})$, taking advantage of the fact that $(T,X_{1},
X_{2})$ is a frame field for each $\Sigma_{t}$. From (3.99)-(3.106), we have:
\begin{align}
 \leftexp{(R_{i})}{\pi}_{LL}=2g(D_{L}R_{i},L)=-2g(R_{i},D_{L}L)=0
\end{align}
This is because $R_{i}$ is tangential to $S_{t,u}$, while $L$ is orthogonal to $S_{t,u}$. For the same reason, 
$T$ is also orthogonal to $R_{i}$, hence:
\begin{equation}
 \leftexp{(R_{i})}{\pi}_{TT}=2g(D_{T}R_{i},T)=-2g(R_{i},D_{T}T)
\end{equation}
Since $R_{i}$ is tangential to $S_{t,u}$, we can expand
\begin{equation}
 R_{i}=R^{A}_{i}X_{A}
\end{equation}
we obtain
\begin{equation}
 \leftexp{(R_{i})}{\pi}_{TT}=2\kappa R_{i}\kappa
\end{equation}
Next,
\begin{align}
 \leftexp{(R_{i})}{\pi}_{LT}=g(D_{L}R_{i},T)+g(D_{T}R_{i},L)\\\notag
=-g(R_{i},D_{L}T)-g(R_{i},D_{T}L)=-\eta_{A}R^{A}_{i}+\zeta_{A}R^{A}_{i}=-R_{i}\mu
\end{align}
where we have used (3.55) in chapter 3.

Next we have:
\begin{equation}
 \leftexp{(R_{i})}{\pi}_{LA}=g(D_{L}R_{i},X_{A})+g(D_{X_{A}}R_{i},L)
\end{equation}
Now, 
\begin{equation}
 g(D_{X_{A}}R_{i},L)=-g(R_{i},D_{X_{A}}L)=-\chi_{AB}R^{B}_{i}
\end{equation}
On the other hand, from (6.6) we have
\begin{equation}
 D_{L}R_{i}=(D_{L}\Pi)\mathring{R}_{i}+\Pi(D_{L}\mathring{R}_{i})
\end{equation}
Noting that for any vectorfield $Z$
\begin{equation}
 g(\Pi Z,X_{A})=g(Z,X_{A})
\end{equation}
we obtain:
\begin{equation}
 g(D_{L}R_{i},X_{A})=g((D_{L}\Pi)\mathring{R}_{i},X_{A})+g(D_{L}\mathring{R}_{i},X_{A})
\end{equation}
Thus,
\begin{equation}
 \leftexp{(R_{i})}{\pi}_{LA}=g((D_{L}\Pi)\mathring{R}_{i},X_{A})+g(D_{L}\mathring{R}_{i},X_{A})-\chi_{AB}R^{B}_{i}
\end{equation}
Similarly, we have
\begin{equation}
 \leftexp{(R_{i})}{\pi}_{TA}=g((D_{T}\Pi)\mathring{R}_{i},X_{A})+g(D_{T}\mathring{R}_{i},X_{A})-\kappa\theta_{AB}R^{B}_{i}
\end{equation}
and
\begin{align}
 \leftexp{(R_{i})}{\pi}_{AB}=g((D_{X_{A}}\Pi)\mathring{R}_{i},X_{B})+g((D_{X_{B}}\Pi)\mathring{R}_{i},X_{A})\\\notag
+g(D_{X_{A}}\mathring{R}_{i},X_{B})+g(D_{X_{B}}\mathring{R}_{i},X_{A})
\end{align}

Now we must compute $(D\Pi)\mathring{R}_{i}$ and $D\mathring{R}_{i}$. Since $\mathring{R}_{i}$ are tangential to $\Sigma_{t}$ we can expand:
\begin{equation}
 \mathring{R}_{i}=R^{A}_{i}X_{A}+\lambda_{i}\hat{T}, \quad \hat{T}=\kappa^{-1}T
\end{equation}
 for some functions $\lambda_{i}$. Since $\Pi T=0$, for any vectorfield $Z$ we have:
\begin{equation}
 (D_{Z}\Pi)T=D_{Z}(\Pi T)-\Pi(D_{Z}T)=-\Pi(D_{Z}T)
\end{equation}
Setting $Z$ equal to $L,T,X_{A}: A=1,2$, we obtain, using (3.99)-(3.106),
\begin{align}
 (D_{L}\Pi)T=\zeta^{A}X_{A}\\
(D_{T}\Pi)T=\kappa(\slashed{d}^{A}\kappa)X_{A}\\
(D_{X_{A}}\Pi)T=-\kappa\theta_{A}^{B}X_{B}
\end{align}
Since $\Pi X_{A}=X_{A}$, for any vectorfield $Z$ we have:
\begin{align}
 (D_{Z}\Pi)X_{A}=D_{Z}(\Pi X_{A})-\Pi(D_{Z}X_{A})=
D_{Z}X_{A}-\Pi(D_{Z}X_{A})
\end{align}
Obviously, the right hand side of the above is orthogonal to $S_{t,u}$. So we get
\begin{align}
 g((D_{L}\Pi)\mathring{R}_{i},X_{A})=\kappa^{-1}\lambda_{i}\zeta_{A}\\
g((D_{T}\Pi)\mathring{R}_{i},X_{A})=\lambda_{i}\slashed{d}_{A}\kappa\\
g((D_{X_{A}}\Pi)\mathring{R}_{i},X_{B})=-\lambda_{i}\theta_{AB}
\end{align}
 Next, we shall compute $D\mathring{R}_{i}$. In an arbitrary coordinate system we have:
\begin{equation}
 D_{\mu}\mathring{R}^{\nu}_{i}=\frac{\partial \mathring{R}^{\nu}_{i}}{\partial x^{\mu}}+\Gamma^{\nu}_{\mu\lambda}\mathring{R}^{\lambda}_{i}
\end{equation}
where $\Gamma^{\nu}_{\mu\lambda}$ are the Christoffel symbols of the acoustical metric $g$ in the given coordinate system, and we have:
\begin{equation}
 \Gamma^{\nu}_{\mu\lambda}=(g^{-1})^{\nu\kappa}\Gamma_{\mu\lambda\kappa}
\end{equation}
Recall from Chapter 3, that in a Galilean coordinate system we have:
\begin{equation}
 \Gamma_{000}=\frac{1}{2}\partial_{t}(-\eta^{2}+|\textbf{v}|^{2})
\end{equation}
\begin{equation}
 \Gamma_{0i0}=\frac{1}{2}\partial_{i}(-\eta^{2}+|\textbf{v}|^{2})
\end{equation}
\begin{equation}
 \Gamma_{ij0}=\partial_{i}\partial_{j}\phi
\end{equation}
\begin{equation}
 \Gamma_{00k}=\frac{1}{2}(2\partial_{t}\partial_{k}\phi-\partial_{k}(-\eta^{2}+|\textbf{v}|^{2}))
\end{equation}
\begin{equation}
 \Gamma_{i0k}=\Gamma_{ijk}=0
\end{equation}
Now we can compute $g(D_{L}\mathring{R}_{i},X_{A}), g(D_{T}\mathring{R}_{i},X_{A})$ and 
$g(D_{X_{A}}\mathring{R}_{i},X_{B})$.
\begin{align}
 g(D_{L}\mathring{R}_{i},X_{A})=g_{\mu\nu}L^{\lambda}D_{\lambda}\mathring{R}^{\mu}_{i}X^{\nu}_{A}
=g_{\mu\nu}L^{\lambda}(\frac{\partial \mathring{R}^{\mu}_{i}}{\partial x^{\lambda}}+\Gamma_{\lambda\gamma}^{\mu}\mathring{R}^{\gamma}_{i})X^{\nu}_{A}
\end{align}
This is:
\begin{align}
L^{\lambda}(g_{\mu\nu}\frac{\partial \mathring{R}^{\mu}_{i}}{\partial x^{\lambda}}+\Gamma_{\lambda\mu\nu}\mathring{R}^{\mu}_{i})X^{\nu}_{A}
=L^{l}\bar{g}_{mn}\frac{\partial \mathring{R}^{m}_{i}}{\partial x^{l}}X^{n}_{A}=L^{l}\delta_{mn}\epsilon_{ilm}X^{n}_{A}
=L^{l}\epsilon_{ilm}X^{m}_{A}
\end{align}
where we have used the expressions for Christoffel symbols and the fact that $\mathring{R}_{i}^{0}=X^{0}_{A}=0$

Similarly, we have 
\begin{align}
 g(D_{T}\mathring{R}_{i},X_{A})=T^{l}\epsilon_{ilm}X^{m}_{A}
\end{align}
and 
\begin{align}
 g(D_{X_{A}}\mathring{R}_{i},X_{B})=X^{l}_{A}\epsilon_{ilm}X^{m}_{B}
\end{align}
Substituting the above in (6.35)-(6.37) we obtain the following:
\begin{align}
 \leftexp{(R_{i})}{\pi_{LA}}=-\chi_{AB}R^{B}_{i}+L^{l}\epsilon_{ilm}X^{m}_{A}+\kappa^{-1}\lambda_{i}\zeta_{A}\\
\leftexp{(R_{i})}{\pi_{TA}}=-\kappa\theta_{AB}R^{B}_{i}+T^{l}\epsilon_{ilm}X^{m}_{A}+\lambda_{i}\slashed{d}_{A}\kappa
\end{align}
and noting that $\epsilon_{ilm}$ is antisymmetric in $l,m$, while $X^{l}_{A}X^{m}_{B}+X^{m}_{A}X^{l}_{B}$ is symmetric in
$l,m$:
\begin{equation}
 \leftexp{(R_{i})}{\pi_{AB}}=-2\lambda_{i}\theta_{AB}
\end{equation}
Since the coefficient of $\lambda_{i}$ in (6.57) is $\kappa^{-1}\zeta_{A}$, 
and recalling from Chapter 3 that $\kappa^{-1}\zeta_{A}$, $\theta_{AB}$ are regular as $\mu\rightarrow0$, these expressions are regular as $\mu\rightarrow0$.

We introduce the functions $y^{i}$ by setting:
\begin{equation}
 \hat{T}^{i}=-\frac{x^{i}}{1-u+t}+y^{i}
\end{equation}
using this and the fact that:
\begin{equation}
 \epsilon_{ilm}x^{l}X^{m}_{A}=\sum_{m}\mathring{R}^{m}_{i}X^{m}_{A}=\slashed{g}_{AB}R^{B}_{i}
\end{equation}
we rewrite the first two terms on the right hand side in (6.58) as:
\begin{equation}
 -\kappa\theta_{AB}R^{B}_{i}+T^{l}\epsilon_{ilm}X^{m}_{A}=-\kappa(\theta_{AB}+\frac{\slashed{g}_{AB}}{1-u+t})R^{B}_{i}
+\kappa\epsilon_{ilm}y^{l}X^{m}_{A}
\end{equation}
Then (6.58) takes the form:
\begin{equation}
 \leftexp{(R_{i})}{\pi_{TA}}=-\kappa(\theta_{AB}+\frac{\slashed{g}_{AB}}{1-u+t})R^{B}_{i}
+\kappa\epsilon_{ilm}y^{l}X^{m}_{A}+\lambda_{i}\slashed{d}_{A}\kappa
\end{equation}

Also we define the function $z^{i}$ by setting:
\begin{equation}
 L^{i}=\frac{x^{i}}{1-u+t}+z^{i}
\end{equation}
By (6.61) and (6.65) and recalling: $L=\partial_{0}-(\eta\hat{T}^{i}+\psi_{i})\partial_{i}$, we get:
\begin{equation}
 z^{i}=-\eta y^{i}+\frac{(\eta-1)x^{i}}{1-u+t}-\psi_{i}
\end{equation}
Using (6.62) and (6.65), we rewrite the first two terms of (6.57):
\begin{equation}
 -\chi_{AB}R^{B}_{i}+L^{l}\epsilon_{ilm}X^{m}_{A}=-(\chi_{AB}-\frac{\slashed{g}_{AB}}{1-u+t})R^{B}_{i}
+\epsilon_{ilm}z^{l}X^{m}_{A}
\end{equation}
So (6.57) takes the form:
\begin{equation}
 \leftexp{(R_{i})}{\pi_{LA}}=-(\chi_{AB}-\frac{\slashed{g}_{AB}}{1-u+t})R^{B}_{i}
+\epsilon_{ilm}z^{l}X^{m}_{A}+\lambda_{i}(\kappa^{-1}\zeta_{A})
\end{equation}

\section{Preliminary Estimates for the Deformation Tensors}
In the remainder of this chapter, we shall show how the deformation tensors $\leftexp{(Y)}{\pi}$ of the commutation vectorfields are to 
be controlled in terms of $\chi, \mu$ and $\psi_{\mu}$.

In the following, the assumptions $\textbf{A1},\textbf{A2},\textbf{A3}$ are assumed to hold in $W^{s}_{\epsilon_{0}}$.
\begin{align*}
 \textbf{A1}: C^{-1}\leq \Omega \leq C\\
\textbf{A2}: C^{-1}\leq \eta\leq C\\
\textbf{A3}: \mu\leq C[1+\log(1+t)]
\end{align*}
where the constant $C$ is independent of $s$.

Recall that 
\begin{equation}
 \psi_{\mu}=\partial_{\mu}\phi
\end{equation}
and we shall use use the following bootstrap assumptions on $\psi_{\mu}$ and their first derivatives. In the following we denote by $\delta_{0}$ a positive 
constant, which is by definition less or equal to unity. This constant is used to keep track of the relative size of various quantities.

On the other hand, no smallness condition is assumed on $\epsilon_{0}$, other than the condition: $\epsilon_{0}\leq \frac{1}{2}$.

There is a positive constant $C$ independent of $s$ such that in $W^{s}_{\epsilon_{0}}$,
\begin{align*}
 \textbf{E1}:|\psi_{0}-h_{0}|, \quad |\psi_{i}|\leq C\delta_{0}(1+t)^{-1}\\
\textbf{E2}: |T\psi_{\mu}|\leq C\delta_{0}(1+t)^{-1},\quad |L\psi_{\mu}|,\quad |\slashed{d}\psi_{\mu}|\leq C\delta_{0}(1+t)^{-2}
\end{align*}
The second of $\textbf{E2}$ shall follow from $|Q\psi_{\mu}|\leq C\delta_{0}(1+t)^{-1}$, while the third shall follow from 
$|R_{i}\psi_{\mu}|\leq C\delta_{0}(1+t)^{-1}: i=1,2,3$ in the actual estimates of Chapter 10,11 and 12.

We also have the bootstrap assumptions on $\mu$ and $\chi$ in $W^{s}_{\epsilon_{0}}$:
\begin{align*}
 \textbf{F1}: |T\mu|\leq C\delta_{0}[1+\log(1+t)], \quad |\slashed{d}\mu|\leq C\delta_{0}(1+t)^{-1}[1+\log(1+t)]\\
\textbf{F2}:|\chi-\frac{\slashed{g}}{1-u+t}|\leq C\delta_{0}(1+t)^{-2}[1+\log(1+t)]
\end{align*}
Here the positive constant $C$ is also independent of $s$.

Under the assumptions $\textbf{A},\textbf{E},\textbf{F}$, the estimates which we shall derive will hold 
in $W^{s}_{\epsilon_{0}}$ with constants which are independent of $s$.

First, recall that the enthalpy is $h=\psi_{0}-\frac{1}{2}\sum_{i=1}^{3}\psi_{i}^{2}$, so by $\textbf{E}$ we have:
\begin{align}
 |h|\leq C\delta_{0}(1+t)^{-1},\quad |Th|\leq C\delta_{0}(1+t)^{-1},\quad |Lh|,\quad |\slashed{d}h|\leq C\delta_{0}(1+t)^{-2}
\end{align}
Recall that $\hat{T}$ is the unit normal to the surface $S_{t,u}$ in $\Sigma_{t}$ with respect to the induced acoustical metric $\bar{g}$, 
which is the Euclidean metric so:
\begin{equation}
 |\hat{T}|= 1
\end{equation}
 It follows by $\textbf{E1}$ that 
\begin{equation}
 |\psi_{\hat{T}}|\leq C\delta_{0}(1+t)^{-1}
\end{equation}
Recalling that $L=\partial_{0}-(\alpha\hat{T}^{i}+\psi_{i})\partial_{i}$, from $\textbf{E1}$ and (6.70) we have:
\begin{equation}
 |L^{i}|\leq C
\end{equation}
Since $L^{0}=1$, we have :
\begin{equation}
 |\psi_{L}|\leq C\delta_{0}(1+t)^{-1}
\end{equation}
Let us express $\bar{g}^{-1}$ in the frame $\hat{T},X_{1},X_{2}$
\begin{align}
 \delta^{ij}=(\bar{g}^{-1})^{ij}=\hat{T}^{i}\hat{T}^{j}+(\slashed{g}^{-1})^{AB}X_{A}^{i}X_{B}^{j}
\end{align}
It follows from $\textbf{E1}$ that:
\begin{align}
 |\slashed{\psi}|^{2}=(\slashed{g}^{-1})^{AB}\slashed{\psi}_{A}\slashed{\psi}_{B}=(\slashed{g}^{-1})^{AB}X_{A}^{i}X_{B}^{j}\psi_{i}\psi_{j}\\\notag
\leq (\bar{g}^{-1})^{ij}\psi_{i}\psi_{j}=\sum_{i}(\psi_{i})^{2}\leq C\delta_{0}^{2}(1+t)^{-2}
\end{align}
that is:
\begin{equation}
 |\slashed{\psi}|\leq C\delta_{0}(1+t)^{-1}
\end{equation}
For a bilinear form $\omega$ in the tangent space to $S_{t,u}$ at a point, the magnitude $|\omega|$ is given by:
\begin{align}
 |\omega|^{2}=(\slashed{g}^{-1})^{AC}(\slashed{g}^{-1})^{BD}\omega_{AB}\omega_{CD}
\end{align}
Recall that $\slashed{k}_{AB}=X^{i}_{A}X^{j}_{B}k_{ij}=-\eta^{-1}X^{i}_{A}X^{j}_{B}\partial_{i}\psi_{j}$. So
\begin{align}
 |\slashed{k}|^{2}=(\slashed{g}^{-1})^{AC}(\slashed{g}^{-1})^{BD}\slashed{k}_{AB}\slashed{k}_{CD}
=(\slashed{g}^{-1})^{AC}(\slashed{g}^{-1})^{BD}X^{i}_{A}X^{j}_{B}k_{ij}X^{l}_{C}X^{m}_{D}k_{lm}
\\\notag=\eta^{-2}(\slashed{g}^{-1})^{AC}(\slashed{g}^{-1})^{BD}X^{i}_{A}\slashed{d}_{B}\psi_{i}X^{l}_{C}\slashed{d}_{D}\psi_{l}\\\notag
\leq C\eta^{-2}(\bar{g}^{-1})^{ij}(\slashed{g}^{-1})^{BD}(\slashed{d}_{B}\psi_{i})(\slashed{d}_{D}\psi_{j})\\\notag=\eta^{-2}\sum_{i}
(\slashed{g}^{-1})^{BD}(\slashed{d}_{B}\psi_{i})
(\slashed{d}_{D}\psi_{i})=\eta^{-2}\sum_{i}|\slashed{d}\psi_{i}|^{2}\leq C\delta_{0}^{2}(1+t)^{-4}
\end{align}
That is 
\begin{equation}
 |\slashed{k}|\leq C\delta_{0}(1+t)^{-2}
\end{equation}
By assumptions $\textbf{A2}$ and $\textbf{A3}$ this implies:
\begin{equation}
 \kappa|\slashed{k}|\leq C\delta_{0}(1+t)^{-2}[1+\log(1+t)]
\end{equation}
For a linear form $v$ in the tangent space to $S_{t,u}$ at a point, the magnitude $|v|$ is given by:
\begin{equation}
 |v|^{2}=(\slashed{g}^{-1})^{AB}v_{A}v_{B}
\end{equation}
From the definition in Chapter 3, $\epsilon_{A}=\kappa^{-1}k(X_{A},T)$, we have 
\begin{align}
 |\epsilon|^{2}=(\slashed{g}^{-1})^{AB}\epsilon_{A}\epsilon_{B}=(\slashed{g}^{-1})^{AB}k_{A\hat{T}}k_{B\hat{T}}=
(\slashed{g}^{-1})^{AB}X^{i}_{A}X^{j}_{B}\hat{T}^{l}\hat{T}^{m}k_{il}k_{jm}\\\notag
=\eta^{-2}(\slashed{g}^{-1})^{AB}\hat{T}^{l}\hat{T}^{m}\slashed{d}_{A}\psi_{l}\slashed{d}_{B}\psi_{m}
\leq \eta^{-2}\sum_{i}|\slashed{d}\psi_{i}|^{2}\leq C\delta^{2}_{0}(1+t)^{-4}
\end{align}
That is 
\begin{equation}
 |\epsilon|\leq C\delta_{0}(1+t)^{-2}
\end{equation}

Next, we shall estimate $\eta$. Recall that 
\begin{equation}
\eta=(\frac{\rho'(h)}{\rho(h)})^{-\frac{1}{2}}
\end{equation}
From the first of (6.69), it then follows that: 
\begin{equation}
 |\eta-1|\leq C\delta_{0}(1+t)^{-1}
\end{equation}
Differentiating (6.84) tangentially to $S_{t,u}$, we obtain $\slashed{d}\eta=(\frac{\rho'(h)}{\rho(h)})^{-\frac{1}{2}\prime}\slashed{d}h$,
hence by the last of (6.69):  
\begin{equation}
 |\slashed{d}\eta|\leq C\delta_{0}(1+t)^{-2}
\end{equation}
Similarly, we have
\begin{equation}
 |T\eta|\leq C\delta_{0}(1+t)^{-1},\quad |L\eta|\leq C\delta_{0}(1+t)^{-2}
\end{equation}
In view of the fact that $\kappa^{-1}\zeta_{A}=\eta\epsilon_{A}-\slashed{d}_{A}\alpha$, we then obtain:
\begin{equation}
 |\kappa^{-1}\zeta|\leq C\delta_{0}(1+t)^{-2}
\end{equation}
hence
\begin{equation}
  |\zeta|\leq C\delta_{0}(1+t)^{-2}[1+\log(1+t)]
\end{equation}
We can also bound $L\mu$, given by (3.92) in Chapter 3:
\begin{equation}
 L\mu=m+\mu e
\end{equation}
where
\begin{equation}
 m=\frac{1}{2}\frac{dH}{dh}Th,\quad e=\frac{1}{2\alpha^{2}}(\frac{\rho^{\prime}}{\rho})^{\prime}Lh+\alpha^{-1}\hat{T}^{i}(L\psi_{i})
\end{equation}
From (6.69) we have:
\begin{equation}
 |m|\leq C\delta_{0}(1+t)^{-1}
\end{equation}
From (6.69), (6.70) and $\textbf{E2}$ we have
\begin{equation}
 |e|\leq C\delta_{0}(1+t)^{-2}
\end{equation}
Hence
\begin{equation}
 |L\mu|\leq C\delta_{0}(1+t)^{-1}
\end{equation}
In view of $\kappa=\alpha^{-1}\mu$, and the second of (6.87), we have:
\begin{equation}
 |L\kappa|\leq C\delta_{0}(1+t)^{-1}
\end{equation}
Integrating this along the integral curves of $L$ yields
\begin{equation}
 |\kappa-\kappa_{0}|\leq C\delta_{0}\log(1+t)
\end{equation}
Here the value of $\kappa_{0}$ at a point in $W^{*}_{\epsilon_{0}}$ is the value of $\kappa$ 
at the corresponding point on $\Sigma_{0}$ along the same integral curve of $L$. Combining (6.96) and the fact that on $\Sigma_{0}$ we have
\begin{align*}
\kappa_{0}=1
\end{align*} 
we obtain the following estimates in $W^{s}_{\epsilon_{0}}$:
\begin{equation}
 |\kappa-1|\leq C\delta_{0}[1+\log(1+t)]
\end{equation}
We now combine $\textbf{F1}$ with (6.86) and the first of (6.87) to estimate the first derivatives of $\kappa$ on $\Sigma_{t}$:
\begin{equation}
 |T\kappa|\leq C\delta_{0}[1+\log(1+t)], \quad |\slashed{d}\kappa|\leq C\delta_{0}(1+t)^{-1}[1+\log(1+t)]
\end{equation}

In view of the fact that $\eta=\zeta+\slashed{d}\mu$, $\textbf{F1}$ and (6.89) yield:
\begin{equation}
 |\eta|\leq C\delta_{0}(1+t)^{-1}[1+\log(1+t)]
\end{equation}
 Now we use $\textbf{F2}$ to estimate
\begin{align}
 \nu=\frac{1}{2}(\textrm{tr}\chi+\frac{d\log\Omega}{dh}Lh),\quad
\underline{\nu}=\frac{1}{2}(\textrm{tr}\underline{\chi}+\frac{d\log\Omega}{dh}\underline{L}h)
\end{align}
If $\delta_{0}$ is small enough, $\textbf{F2}$ implies:
\begin{equation}
 \textrm{tr}\chi\geq C^{-1}(1+t)^{-1}
\end{equation}
This restriction on the size of $\delta_{0}$ is imposed from now on.

Recall that:
\begin{equation}
 \chi=\eta(\slashed{k}-\theta),\quad \underline{\chi}=\kappa(\slashed{k}+\theta)
\end{equation}
Thus we have:
\begin{equation}
 \textrm{tr}\underline{\chi}=2\kappa\textrm{tr}\slashed{k}-\eta^{-1}\kappa\textrm{tr}\chi,
\hat{\underline{\chi}}=2\kappa\hat{\slashed{k}}-\eta^{-1}\kappa\hat{\chi}
\end{equation}
From (6.80), (6.96) and $\textbf{F2}$ we have:
\begin{align}
 |\textrm{tr}\underline{\chi}|\leq C(1+t)^{-1}[1+\log(1+t)],\quad
|\hat{\underline{\chi}}|\leq C\delta_{0}(1+t)^{-2}[1+\log(1+t)]^{2}
\end{align}
Also, using (6.79), (6.85) and $\textbf{F2}$ we deduce:
\begin{align}
 |\theta+\frac{\slashed{g}}{1-u+t}|\leq C\delta_{0}(1+t)^{-2}[1+\log(1+t)]
\end{align}
Finally, combining $\textbf{F2}$, (6.69) and (6.101), we obtain:
\begin{align}
 C^{-1}(1+t)^{-1}\leq \nu\leq C(1+t)^{-1},\quad
|\nu-(1-u+t)^{-1}|\leq C\delta_{0}(1+t)^{-2}[1+\log(1+t)]
\end{align}
Also, from (6.69) and (6.104) we have:
\begin{equation}
 |\underline{\nu}|\leq C(1+t)^{-1}[1+\log(1+t)]
\end{equation}

We can now give the bounds of the deformation tensor of $T$ and $Q$. In the following, we denote by $\leftexp{(Y)}{\slashed{\pi}_{L}}$ 
and $\leftexp{(Y)}{\slashed{\pi}_{\underline{L}}}$ the 1-forms on each surface $S_{t,u}$ with components:
\begin{equation}
 \leftexp{(Y)}{\slashed{\pi}_{L}}(X_{A})=\leftexp{(Y)}{\slashed{\pi}_{LA}}
, \leftexp{(Y)}{\slashed{\pi}_{\underline{L}}}(X_{A})=\leftexp{(Y)}{\slashed{\pi}_{\underline{L}A}}
\end{equation}
and similarly for $\leftexp{(Y)}{\tilde{\slashed{\pi}}_{L}}$ and $\leftexp{(Y)}{\tilde{\slashed{\pi}}_{\underline{L}}}$.

First, we estimate the deformation tensor of $T$.
By (6.11), (6.97), (6.98) and (6.87) we have:
\begin{equation}
 \mu^{-1}|\leftexp{(T)}{\tilde{\pi}}_{\underline{L}\underline{L}}|\leq C\delta_{0}[1+\log(1+t)]
\end{equation}
By (6.12), $\textbf{F1}, \textbf{A}$, and (6.69) we have:
\begin{equation}
 |\leftexp{(T)}{\tilde{\pi}}_{L\underline{L}}|\leq C\delta_{0}[1+\log(1+t)]
\end{equation}
By (6.13), (6.89) and (6.99) we have 
\begin{equation}
 |\leftexp{(T)}{\tilde{\slashed{\pi}}}_{L}|\leq C\delta_{0}(1+t)^{-1}[1+\log(1+t)]
\end{equation}
By (6.14), (6.89), (6.97) and (6.99) we have:
\begin{equation}
 \mu^{-1}|\leftexp{(T)}{\tilde{\slashed{\pi}}}_{\underline{L}}|\leq C\delta_{0}(1+t)^{-1}[1+\log(1+t)]
\end{equation}
By (6.15), (6.97), (6.104) and $\textbf{F2}$,
we have:
\begin{equation}
 |\leftexp{(T)}{\hat{\tilde{\slashed{\pi}}}}|\leq C\delta_{0}(1+t)^{-2}[1+\log(1+t)]^{2}
\end{equation}
Finally by (6.16), (6.97), (6.106) and (6.107) we have:
\begin{equation}
 |\textrm{tr}\leftexp{(T)}{\tilde{\slashed{\pi}}}|\leq C(1+t)^{-1}[1+\log(1+t)]
\end{equation}

Next, we estimate the deformation tensor of $Q$.
By (6.18), (6.87), (6.95) and (6.97) we have:
\begin{equation}
 \mu^{-1}|\leftexp{(Q)}{\tilde{\pi}_{\underline{L}\underline{L}}}|\leq C[1+\log(1+t)]
\end{equation}
By (6.19), (6.69), (6.94) we have:
\begin{equation}
 |\leftexp{(Q)}{\tilde{\pi}}_{L\underline{L}}|\leq C[1+\log(1+t)]
\end{equation}
By (6.21), (6.89) and (6.99) we have:
\begin{equation}
 |\leftexp{(Q)}{\tilde{\slashed{\pi}}}_{\underline{L}}|\leq C\delta_{0}[1+\log(1+t)]
\end{equation}
By (6.22), $\textbf{F2}$ we have:
\begin{equation}
 |\leftexp{(Q)}{\hat{\tilde{\slashed{\pi}}}}|\leq C\delta_{0}(1+t)^{-1}[1+\log(1+t)]
\end{equation}
By (6.23), (6.106) we have:
\begin{equation}
 |\textrm{tr}\leftexp{(Q)}{\tilde{\slashed{\pi}}}|\leq C
\end{equation}
We can obtain, in fact, a more precise estimate for $\textrm{tr}\leftexp{(Q)}{\tilde{\slashed{\pi}}}$:
\begin{align}
 \textrm{tr}\leftexp{(Q)}{\tilde{\slashed{\pi}}}-4=4(\Omega-1)(1+t)\nu+4(1+t)(\nu-\frac{1}{1+t})
\end{align}

Recalling that $\Omega(0)=1$, i.e. at the constant state, the conformal factor $\Omega=1$, (6.69) implies:
\begin{equation}
 |\Omega-1|\leq C\delta_{0}(1+t)^{-1}
\end{equation}
Then by (6.106), we get
\begin{align}
 |\textrm{tr}\leftexp{(Q)}{\tilde{\slashed{\pi}}}-4|\leq C(1+t)^{-1}[1+\log(1+t)]
\end{align}
where we have used the fact that the behavior of difference between $(1+t)^{-1}$ and $(1-u+t)^{-1}$ is like $(1+t)^{-2}$.

To bound the deformation tensor of $R_{i}$, we must get estimates for $\lambda_{i}$ defined by (6.37), as well as the functions $y_{i}$ and $z_{i}$ defined by 
(6.60) and (6.64).

First we shall get upper and lower bounds for the function:
\begin{equation}
 r=\sqrt{\sum_{i=1}^{3}(x^{i})^{2}}
\end{equation}
Since $r$ achieves its maximum value in $\Sigma_{t}\bigcap W^{*}_{\epsilon_{0}}$ at the outer boundary $\Sigma_{t}\bigcap C_{0}=S_{t,0}$ where it is equal to
$1+t$, we have:
\begin{equation}
 r\leq 1+t 
\end{equation}
in $W^{*}_{\epsilon_{0}}$.

For the lower bound, we have:
\begin{equation}
 Tr=\sum_{i=1}^{3}\frac{x^{i}T^{i}}{r}=\kappa\sum_{i=1}^{3}\frac{x^{i}\hat{T}^{i}}{r}
\end{equation}
In view of (6.70), we get 
\begin{equation}
 |Tr|\leq \kappa 
\end{equation}
 then by (6.97), we have:
\begin{equation}
 |Tr|\leq 1+C\delta_{0}[1+\log(1+t)]
\end{equation}
Integrating (6.127) along the integral curves of $T$ from $\Sigma_{t}\bigcap C_{0}=S_{t,0}$ where $r=1+t$ to $S_{t,u}$, we get
\begin{equation}
 r\geq 1-u+t-C\delta_{0}u[1+\log(1+t)]
\end{equation}
In $W_{\epsilon_{0}}^{*}$, $u\in [0,\epsilon_{0}]$, hence taking $\delta_{0}$ suitably small we get:
\begin{equation}
 r\geq C^{-1}(1+t)
\end{equation}
We proceed to derive the estimates for $\lambda_{i}$, from (6.37) we have:
\begin{equation}
 \lambda_{i}=\bar{g}(\mathring{R}_{i},\hat{T})
\end{equation}
First, we derive the estimates for $\mathring{R}_{i}$ and $R_{i}$
\begin{equation}
 |\mathring{R}_{i}|=\sqrt{r^{2}-(x_{i})^{2}}\leq r\leq 1+t
\end{equation}
From (6.37), we have:
\begin{equation}
 |R_{i}|\leq |\mathring{R}_{i}|\leq 1+t
\end{equation}
Next, we derive an ordinary differential equation for $\lambda_{i}$ along the integral curves of $L$. Note that the surfaces $S_{0,u}$ being spheres centered 
at the origin and we have $\hat{T}=-\frac{\partial}{\partial r}$ on $\Sigma_{0}$, hence 
\begin{align}
 \lambda_{i}=0\quad: \quad\textrm{on}\quad\Sigma_{0}
\end{align}
Now we have:
\begin{equation}
 L\lambda_{i}=\sum_{m}(L\mathring{R}^{m}_{i})\hat{T}^{m}+\sum_{m}\mathring{R}^{m}_{i}L\hat{T}^{m}
\end{equation}
Recall that $L=\partial_{t}-(\eta\hat{T}^{i}+\psi_{i})\partial_{i}$ and $\mathring{R}_{i}^{m}=\epsilon_{ijm}x^{j}$.
\begin{equation}
 \Rightarrow L\mathring{R}_{i}^{m}=-(\eta\hat{T}^{j}+\psi_{j})\epsilon_{ijm}
\end{equation}
While $L\hat{T}^{m}$ is given by (3.177):
\begin{equation}
 L\hat{T}^{m}=p_{L}\hat{T}^{m}+q^{m}_{L}
\end{equation}
Recall from Chapter 3, $p_{L}=0$, and from (3.179) we have:
\begin{align}
 q_{L}=q^{A}_{L}X_{A}=q^{A}_{L}X_{A}^{m}\partial_{m} \Rightarrow q_{L}^{m}=q_{L}^{A}X^{m}_{A}
\end{align}
we thus obtain: 
\begin{align}
 L\lambda_{i}=-\psi_{j}\hat{T}^{m}\epsilon_{ijm}+\epsilon_{ijm}q^{A}_{L}X^{m}_{A}x^{j}
\end{align}
where we have used the fact that $\epsilon_{ijm}\hat{T}^{j}\hat{T}^{m}=0$.

Recalling from (3.183) that
\begin{equation}
 q^{A}_{L}=-\kappa^{-1}\zeta^{A}
\end{equation}
the estimate (6.88) gives
\begin{align}
 |q_{L}|\leq C\delta_{0}(1+t)^{-2}
\end{align}
The second term on the right in (6.138) is $\bar{g}(\mathring{R}_{i},q_{L})$.

By (6.131) and (6.140): 
\begin{equation}
 |\bar{g}(\mathring{R}_{i},q_{L})|\leq |\mathring{R}_{i}||q_{L}|\leq C\delta_{0}(1+t)^{-1}
\end{equation}
using also $\textbf{E1}$ we obtain: 
\begin{equation}
 |L\lambda_{i}|\leq C\delta_{0}(1+t)^{-1}
\end{equation}
Integrating this along the integral curve of $L$, we conclude that:
\begin{equation}
 |\lambda_{i}|\leq C\delta_{0}[1+\log(1+t)]
\end{equation}

Next, we should derive the estimates of the derivatives of $\lambda_{i}$ tangent to the $\Sigma_{t}$. We have:
\begin{equation}
 T\lambda_{i}=\sum_{m}(T\mathring{R}_{i}^{m})\hat{T}^{m}+\sum_{m}\mathring{R}^{m}_{i}T(\hat{T}^{m})
\end{equation}
From $\mathring{R}^{m}_{i}=\epsilon_{ijm}x^{j}$, we have:
\begin{equation}
 \hat{T}\mathring{R}^{m}_{i}=\epsilon_{ijm}\hat{T}^{j}
\end{equation}
\begin{equation*}
 \Rightarrow \sum_{m}(T\mathring{R}_{i}^{m})\hat{T}^{m}=\epsilon_{ijm}\hat{T}^{j}\hat{T}^{m}=0
\end{equation*}
Hence:
\begin{equation}
  T\lambda_{i}=\sum_{m}\mathring{R}^{m}_{i}T(\hat{T}^{m})
\end{equation}
From (3.192), we know that 
\begin{equation}
 T(\hat{T}^{m})=p_{T}\hat{T}^{m}+q^{m}_{T}
\end{equation}
 From (3.193), we have $p_{T}=0$, and from (3.194) and (3.195), we have $q^{A}_{T}=-X^{A}(\kappa)$. So
\begin{equation}
 T\lambda_{i}=\bar{g}(\mathring{R}_{i},q_{T}) \Rightarrow |T\lambda_{i}|\leq |\mathring{R}_{i}||q_{T}|
\end{equation}
and 
\begin{equation}
 |q_{T}|^{2}=\slashed{g}_{AB}q^{A}_{T}q^{B}_{T}=|\slashed{d}\kappa|^{2}\leq C\delta_{0}^{2}(1+t)^{-2}[1+\log(1+t)]^{2}
\end{equation}
by the second of (6.98). In view also of (6.131) we conclude that:
\begin{equation}
 |T\lambda_{i}|\leq C\delta_{0}[1+\log(1+t)]
\end{equation}
Next, we have:
\begin{equation}
 \slashed{d}_{A}\mathring{R}^{m}_{i}=\epsilon_{ijm}\slashed{d}_{A}x^{j}=\epsilon_{ijm}\textsl{X}^{j}_{A}
\end{equation}
\begin{align}
 \Rightarrow \slashed{d}_{A}\lambda_{i}=\sum_{m}(\slashed{d}_{A}\mathring{R}^{m}_{i})\hat{T}^{m}
+\sum_{m}\mathring{R}^{m}_{i}\slashed{d}_{A}\hat{T}^{m}
=\epsilon_{ijm}\textsl{X}^{j}_{A}\hat{T}^{m}+\sum_{m}\mathring{R}^{m}_{i}(\slashed{p}_{A}\hat{T}^{m}+\slashed{q}^{m}_{A})
\end{align}
where we have used (3.201). Also from (3.202) and (3.203) we have:
\begin{equation}
 \slashed{p}_{A}=0, \quad \slashed{q}_{A}^{m}=\theta^{B}_{A}X_{B}^{m}
\end{equation}
Therefore: 
\begin{equation}
 \slashed{d}_{A}\lambda_{i}=\epsilon_{ijm}\textsl{X}^{j}_{A}\hat{T}^{m}+\bar{g}(\mathring{R}_{i},\slashed{q}_{A})
\end{equation}
By (6.60), the first term on the right in (6.154) is:
\begin{align}
 -\frac{\epsilon_{ijm}X^{j}_{A}x^{m}}{1-u+t}+\epsilon_{ijm}X^{j}_{A}y^{m}
\end{align}
On the other hand, the second term on the right in (6.154) is:
\begin{align}
 \mathring{R}^{m}_{i}\theta_{A}^{B}X^{m}_{B}=\epsilon_{ijm}x^{j}\theta_{A}^{B}X_{B}^{m}\\\notag
=-\frac{\epsilon_{ijm}x^{j}X_{A}^{m}}{1-u+t}+\epsilon_{ijm}x^{j}(\theta_{A}^{B}+\frac{\delta_{A}^{B}}{1-u+t})X_{B}^{m}
\end{align}
The first terms in (6.155) and (6.156) cancel each other, hence we obtain:
\begin{align}
 \slashed{d}_{A}\lambda_{i}=\epsilon_{ijm}[X^{j}_{A}y^{m}+x^{j}(\theta_{A}^{B}+\frac{\delta_{A}^{B}}{1-u+t})X^{m}_{B}]
\end{align}
In view of (6.105) and (6.124), the estimate
\begin{align*}
 |y^{i}|\leq C\delta_{0}(1+t)^{-1}[1+\log(1+t)]
\end{align*}
would imply:
\begin{align}
 |\slashed{d}\lambda_{i}|\leq C\delta_{0}(1+t)^{-1}[1+\log(1+t)]
\end{align}
To derive the above estimate for the functions $y^{i}$, and a similar estimate for the functions $x^{i}$, we shall need another bound for $r$.

Note that for any vector $V$ in euclidean space $\mathbb{R}^{3}$, we have:
\begin{equation}
 \sum_{i}\bar{g}(\mathring{R}_{i},V)^{2}=r^{2}|\Sigma V|^{2}
\end{equation}
where $\Sigma$ is the Euclidean projection operator to the Euclidean spheres:
\begin{equation}
 \Sigma^{i}_{j}=\delta_{j}^{i}-r^{-2}x^{i}x^{j}
\end{equation}
We can check this with the identities:
\begin{equation}
 \sum_{i}\epsilon_{ijm}\epsilon_{ikn}=\delta_{jk}\delta_{mn}-\delta_{jn}\delta_{km}, \quad \sum_{i}\Sigma_{m}^{i}\Sigma^{i}_{n}=\Sigma^{m}_{n}
\end{equation}
Taking in (6.159) $V$ to be $\hat{T}$, we obtain in view of (6.130):
\begin{equation}
 |\Sigma\hat{T}|^{2}=r^{-2}\sum_{i}(\lambda_{i})^{2}
\end{equation}
Then by (6.129) and (6.143) we have:
\begin{equation}
 |\Sigma\hat{T}|\leq C\delta_{0}(1+t)^{-1}[1+\log(1+t)]
\end{equation}
Consider $N$, the Euclidean outward unit normal to the Euclidean spheres:
\begin{equation}
 N=\frac{x^{i}}{r}\partial_{i}
\end{equation}
Then (6.125) takes the form:
\begin{equation}
 Tr=\kappa\bar{g}(N,\hat{T})
\end{equation}
We decompose:
\begin{equation}
 \hat{T}=\bar{g}(N,\hat{T})N+\Sigma\hat{T}
\end{equation}
\begin{equation*}
 \Rightarrow |\hat{T}|^{2}=|\bar{g}(N,\hat{T})|^{2}+|\Sigma\hat{T}|^{2}
\end{equation*}
Recalling that $|\hat{T}|=1$, we have:
\begin{equation}
 1-\bar{g}(N,\hat{T})^{2}=|\Sigma\hat{T}|^{2}\leq C\delta_{0}(1+t)^{-2}[1+\log(1+t)]^{2}
\end{equation}
If $\delta_{0}$ is small enough, the right hand side above is no bigger than $\frac{1}{2}$. While on $S_{t,0}$, 
$\bar{g}(N,\hat{T})=-1$, by continuity in $u$ this implies the angle between $\hat{T}$ and $N$ is greater than $\frac{\pi}{2}$. Hence
\begin{align*}
 \bar{g}(\hat{T},N)<0
\end{align*}
Therefore
\begin{equation}
  0\leq 1+\bar{g}(N,\hat{T})\leq C\delta_{0}(1+t)^{-2}[1+\log(1+t)]^{2}
\end{equation}
Going back to (6.165),
\begin{equation}
 Tr+1=\kappa(1+\bar{g}(N,\hat{T}))-(\kappa-1)
\end{equation}
From (6.97) and (6.168) we have:
\begin{equation}
 |Tr+1|\leq C\delta_{0}[1+\log(1+t)]
\end{equation}
Since $r=1+t$ on $S_{t,0}$, we have, on $S_{t,u}$:
\begin{equation}
 r-1-t+u=\int_{0}^{u}(Tr+1)du'
\end{equation}
Using (6.170), we get:
\begin{equation}
 r\leq 1+t-u[1-C\delta_{0}(1+\log(1+t))]
\end{equation}

If $C\delta_{0}[1+\log(1+t)]<1$, (6.172) is stronger than (6.124), otherwise weaker.

From (6.128) and (6.172), we get 
\begin{equation}
 |\frac{1}{r}-\frac{1}{1-u+t}|\leq C\delta_{0}u(1+t)^{-2}[1+\log(1+t)]
\end{equation}
Consider now the vectorfield:
\begin{equation}
 y'=\hat{T}+N=\Sigma\hat{T}+(1+\bar{g}(N,\hat{T}))N
\end{equation}
(6.163) and (6.168) imply:
\begin{equation}
 |y'|\leq C\delta_{0}(1+t)^{-1}[1+\log(1+t)]
\end{equation}
Recall (6.60) and $N^{i}=\frac{x^{i}}{r}$, we get:
\begin{equation}
 y'^{i}=-\frac{x^{i}}{1-u+t}+y^{i}+\frac{x^{i}}{r} \Rightarrow
y^{i}=y'^{i}-x^{i}(\frac{1}{r}-\frac{1}{1-u+t})
\end{equation}
Using (6.124), (6.172) and (6.174) we have:
\begin{equation}
 |y^{i}|\leq C\delta_{0}(1+t)^{-1}[1+\log(1+t)]
\end{equation}
By this and (6.65), $\textbf{E1}$, we have:
\begin{equation}
 |z^{i}|\leq C\delta_{0}(1+t)^{-1}[1+\log(1+t)]
\end{equation}
We also used the fact that $|\alpha-1|\leq C\delta_{0}(1+t)^{-1}$.

We turn to estimate the deformation tensor of $R_{i}$. Since $\underline{L}=\alpha^{-1}\kappa L+2T$, we get from (6.24) and (6.27)-(6.29):
\begin{align}
 \leftexp{(R_{i})}{\pi}_{LL}=0\\
\leftexp{(R_{i})}{\pi}_{L\underline{L}}=-2R_{i}\mu\\
\leftexp{(R_{i})}{\pi}_{\underline{L}\underline{L}}=4\mu R_{i}(\alpha^{-1}\kappa)
\end{align}
By $\textbf{F1}$, (6.85), (6.86), (6.97), (6.98) and (6.132), we have:
\begin{align}
 |\leftexp{(R_{i})}{\pi}_{L\underline{L}}|\leq C\delta_{0}[1+\log(1+t)]\\
\mu^{-1}|\leftexp{(R_{i})}{\pi_{\underline{L}\underline{L}}}|\leq C\delta_{0}[1+\log(1+t)]
\end{align}
 Next, we estimate the $S_{t,u}$ 1-form $\leftexp{(R_{i})}{\slashed{\pi}_{L}}$, whose component $\leftexp{(R_{i})}{\pi_{LA}}$ is 
given by (6.67). In view of $\textbf{F2}$ and (6.132), the first term on the right of (6.67) is bounded by $C\delta_{0}(1+t)^{-1}[1+\log(1+t)]$.
The square magnitude of the second term is $(\slashed{g}^{-1})^{AB}\epsilon_{ilm}z^{l}\textsl{X}^{m}_{A}\epsilon_{irp}z^{r}\textsl{X}^{p}_{B}$.
From (6.74), this is bounded by $\sum_{m}\epsilon_{ilm}z^{l}\epsilon_{irm}z^{r}=\sum_{l\slashed{=}i}(z^{l})^{2}$. We have used the identity:
$\sum_{i}\epsilon_{ijm}\epsilon_{ikn}=\delta_{jk}\delta_{mn}-\delta_{jn}\delta_{km}$. Then by (6.178), the second term in (6.67) is bounded by:
$C\delta_{0}(1+t)^{-1}[1+\log(1+t)]$. From (6.88) and (6.143), the third term is bounded by: $C\delta_{0}^{2}(1+t)^{-2}[1+\log(1+t)]$.
\begin{equation}
 \Rightarrow |\leftexp{(R_{i})}{\slashed{\pi}_{L}}|\leq C\delta_{0}(1+t)^{-1}[1+\log(1+t)]
\end{equation}
Next, we estimate $\leftexp{(R_{i})}{\pi_{TA}}$, which is given by (6.63).

By (6.105), the first term is bounded by $C\delta_{0}(1+t)^{-1}[1+\log(1+t)]^{2}$.

By (6.97) and (6.177) we know that the second term is bounded by $C\delta_{0}(1+t)^{-1}[1+\log(1+t)]^{2}$

Lastly, Using (6.98) and (6.143) the third term $\lambda_{i}\slashed{d}_{A}\kappa$ is bounded by $C\delta_{0}(1+t)^{-1}[1+\log(1+t)]^{2}$.

So we get
\begin{equation}
 |\leftexp{(R_{i})}{\slashed{\pi}}_{T}|\leq C\delta_{0}(1+t)^{-1}[1+\log(1+t)]^{2}
\end{equation}
Since $\underline{L}=\alpha^{-1}\kappa L+2T$, we get:
\begin{equation}
 |\leftexp{(R_{i})}{\slashed{\pi}}_{\underline{L}}|\leq C\delta_{0}(1+t)^{-1}[1+\log(1+t)]^{2}
\end{equation}
Finally, we estimate $\leftexp{(R_{i})}{\slashed{\pi}}_{AB}$, which is given by (6.59). Using $\textbf{F2}$, (6.80) and (6.143), this is bounded by
$C\delta_{0}(1+t)^{-1}[1+\log(1+t)]$.
But we would like to bound its trace and trace-free part separately.

Since $\textrm{tr}\leftexp{(R_{i})}{\slashed{\pi}}=2\lambda_{i}(\alpha^{-1}\textrm{tr}\chi-\textrm{tr}\slashed{k})$, also by 
$\textbf{F2}$, (6.79) and (6.143),
\begin{equation}
 |\textrm{tr}\leftexp{(R_{i})}{\slashed{\pi}}|\leq C\delta_{0}(1+t)^{-1}[1+\log(1+t)]
\end{equation}
as well as 
\begin{equation}
 |\leftexp{(R_{i})}{\hat{\slashed{\pi}}}|\leq C\delta_{0}(1+t)^{-2}[1+\log(1+t)]^{2}
\end{equation}
Then we can substitute all these estimates to the expressions of $\leftexp{(R_{i})}{\tilde{\pi}}$:
\begin{equation}
 \leftexp{(R_{i})}{\tilde{\pi}}=\Omega\leftexp{(R_{i})}{\pi}+\frac{d\Omega}{dh}(R_{i}h)g
\end{equation}
By (6.69) and (6.132) we have: $|R_{i}h|\leq C\delta_{0}(1+t)^{-1}$, then we get:
\begin{align}
 \leftexp{(R_{i})}{\tilde{\pi}}_{LL}=0\\
|\leftexp{(R_{i})}{\tilde{\pi}}_{L\underline{L}}|\leq C\delta_{0}[1+\log(1+t)]\\
\mu^{-1}|\leftexp{(R_{i})}{\tilde{\pi}}_{\underline{L}\underline{L}}|\leq C\delta_{0}[1+\log(1+t)]\\
|\leftexp{(R_{i})}{\tilde{\slashed{\pi}}}_{L}|\leq C\delta_{0}(1+t)^{-1}[1+\log(1+t)]\\
|\leftexp{(R_{i})}{\tilde{\slashed{\pi}}}_{\underline{L}}|\leq C\delta_{0}(1+t)^{-1}[1+\log(1+t)]^{2}\\
|\textrm{tr}\leftexp{(R_{i})}{\tilde{\slashed{\pi}}}|\leq C\delta_{0}(1+t)^{-1}[1+\log(1+t)]\\
|\leftexp{(R_{i})}{\hat{\tilde{\slashed{\pi}}}}|\leq C\delta_{0}(1+t)^{-2}[1+\log(1+t)]^{2}
\end{align}
Finally, we shall estimate the functions $\leftexp{(R_{i})}{\delta}$, defined by:
\begin{equation}
 \leftexp{(R_{i})}{\delta}=\frac{1}{2}\tilde{\textrm{tr}}\leftexp{(Y)}{\tilde{\pi}}-\mu^{-1}Y\mu-2\Omega^{-1}Y\Omega
\end{equation}
which will play a role in Chaptter 7.
Here $\tilde{\textrm{tr}}$ is the trace with respect to $\tilde{g}$. Since $\tilde{g}^{-1}=\Omega^{-1}g^{-1}$, we have:
\begin{equation}
 \tilde{\textrm{tr}}\leftexp{(Y)}{\tilde{\pi}}=\Omega^{-1}\textrm{tr}\leftexp{(Y)}{\tilde{\pi}}
=\Omega^{-1}(-\mu^{-1}\leftexp{(Y)}{\tilde{\pi}_{L\underline{L}}}+\textrm{tr}\leftexp{(Y)}{\tilde{\slashed{\pi}}})
\end{equation}
Hence,
\begin{equation}
 \leftexp{(Y)}{\delta}=-\frac{1}{2}\Omega^{-1}\mu^{-1}\leftexp{(Y)}{\tilde{\pi}_{L\underline{L}}}-\mu^{-1}Y\mu
+\Omega^{-1}(\frac{1}{2}\textrm{tr}\leftexp{(Y)}{\tilde{\slashed{\pi}}}-2\frac{d\Omega}{dh}Yh)
\end{equation}
From (6.12), we find:
\begin{equation}
 \leftexp{(T)}{\delta}=\Omega^{-1}(\frac{1}{2}\textrm{tr}\leftexp{(T)}{\tilde{\slashed{\pi}}}-\frac{d\Omega}{dh}Th)
\end{equation}
From (6.19),
\begin{equation}
 \leftexp{(Q)}{\delta}=\Omega^{-1}(\frac{1}{2}\textrm{tr}\leftexp{(Q)}{\tilde{\slashed{\pi}}}-\frac{d\Omega}{dh}Qh)+1
\end{equation}
From (6.69) and (6.114), we have:
\begin{equation}
 |\leftexp{(T)}{\delta}|\leq C(1+t)^{-1}[1+\log(1+t)]
\end{equation}
from (6.69) and (6.122), we have:
\begin{equation}
 |\leftexp{(Q)}{\delta}|\leq C
\end{equation}
From (6.180) and (6.189), we have:
\begin{equation}
 \leftexp{(R_{i})}{\tilde{\pi}}_{L\underline{L}}=-2\Omega R_{i}\mu-2\mu R_{i}\Omega
\end{equation}
\begin{equation}
 \Rightarrow \leftexp{(R_{i})}{\delta}=\Omega^{-1}(\frac{1}{2}\textrm{tr}\leftexp{(R_{i})}{\tilde{\slashed{\pi}}}
-\frac{d\Omega}{dh}R_{i}h)
\end{equation}
By (6.69) and (6.195), we have:
\begin{equation}
 |\leftexp{(R_{i})}{\delta}|\leq C\delta_{0}(1+t)^{-1}[1+\log(1+t)]
\end{equation}

\chapter{Outline of the Derived Estimates of Each Order}

\section{The Inhomogeneous Wave Equations for the Higher\\Order Variations.
The Recursion Formula for the Source Functions}
$\textbf{Proposition 7.1}$ Let $\psi$ be a solution of the inhomogeneous wave equation
\begin{equation}
 \Box_{\tilde{g}}\psi=\rho
\end{equation}
and let $Y$ be an arbitrary vectorfield. Then
\begin{equation}
 \psi'=Y\psi
\end{equation}
satisfies the inhomogeneous wave equation:
\begin{equation}
 \Box_{\tilde{g}}\psi'=\rho'
\end{equation}
where the source function $\rho'$ is related to the source function $\rho$ by:
\begin{equation}
 \rho'=\tilde{\textrm{div}}\leftexp{(Y)}{\tilde{J}}+Y\rho+\frac{1}{2}\tilde{\textrm{tr}}\leftexp{(Y)}{\tilde{\pi}}\rho
\end{equation}
Here
\begin{equation}
 \leftexp{(Y)}{\tilde{J}^{\mu}}=\frac{1}{2}[(\tilde{g}^{-1})^{\mu\alpha}(\tilde{g}^{-1})^{\nu\beta}
+(\tilde{g}^{-1})^{\nu\alpha}(\tilde{g}^{-1})^{\mu\beta}-(\tilde{g}^{-1})^{\mu\nu}(\tilde{g}^{-1})^{\alpha\beta}]
\leftexp{(Y)}{\tilde{\pi}}_{\alpha\beta}\partial_{\nu}\psi
\end{equation}
is the commutation current associated to $\psi$ and $Y$.

$Proof$. Let $f_{s}$ be the local 1-parameter group generated by $Y$. We denote by $f_{s*}$ the corresponding pullback. We then have:
\begin{equation}
 \Box_{f_{s*}\tilde{g}}(f_{s*}\psi)=f_{s*}(\Box_{\tilde{g}}\psi)=f_{s*}\rho
\end{equation}
Now, in an arbitrary system of local coordinates,
\begin{equation}
 \Box_{\tilde{g}}\psi=\frac{1}{\sqrt{-\det\tilde{g}}}\partial_{\mu}((\tilde{g}^{-1})^{\mu\nu}\sqrt{-\det\tilde{g}}\partial_{\nu}\psi)
\end{equation}
and
\begin{equation}
 \Box_{f_{s*}\tilde{g}}(f_{s*}\psi)=\frac{1}{\sqrt{-\det f_{s*}\tilde{g}}}\partial_{\mu}(((f_{s*}\tilde{g})^{-1})^{\mu\nu}
\sqrt{-\det f_{s*}\tilde{g}}\partial_{\nu}(f_{s*}\psi))
\end{equation}
Let us differentiate (7.8)  with respect to $s$ at $s=0$. In view of the facts that:
\begin{align}
 (\frac{d}{ds}f_{s*}\tilde{g})_{s=0}=\mathcal{L}_{Y}\tilde{g}=\leftexp{(Y)}{\tilde{\pi}}\\
(\frac{d}{ds}\sqrt{-\det f_{s*}\tilde{g}})_{s=0}=\frac{1}{2}\sqrt{-\det\tilde{g}}\tilde{\textrm{tr}}\leftexp{(Y)}{\tilde{\pi}}\\
(\frac{d}{ds}((f_{s*}\tilde{g})^{-1})^{\mu\nu})_{s=0}=-(\tilde{g}^{-1})^{\mu\alpha}(\tilde{g}^{-1})^{\nu\beta}\leftexp{(Y)}{\tilde{\pi}}_{\alpha\beta}
\end{align}
and:
\begin{equation}
 (\frac{d}{ds}f_{s*}\psi)_{s=0}=(\frac{d}{ds}\psi\circ f_{s})_{s=0}=Y\psi
\end{equation}
 we obtain:
\begin{align}
 (\frac{d}{ds}\Box_{f_{s*}\tilde{g}}(f_{s*}\psi))_{s=0}=-\frac{1}{2}\frac{1}{\sqrt{-\det\tilde{g}}}\tilde{\textrm{tr}}
\leftexp{(Y)}{\tilde{\pi}}\partial_{\mu}((\tilde{g}^{-1})^{\mu\nu}\sqrt{-\det\tilde{g}}\partial_{\nu}\psi)\\\notag
+\frac{1}{\sqrt{-\det\tilde{g}}}\partial_{\mu}\{\sqrt{-\det\tilde{g}}(
-(\tilde{g}^{-1})^{\mu\alpha}(\tilde{g}^{-1})^{\nu\beta}\leftexp{(Y)}{\tilde{\pi}}_{\alpha\beta}+
\frac{1}{2}(\tilde{g}^{-1})^{\mu\nu}\tilde{\textrm{tr}}\leftexp{(Y)}{\tilde{\pi}})\partial_{\nu}\psi\\\notag
+(\tilde{g}^{-1})^{\mu\nu}\sqrt{-\det\tilde{g}}\partial_{\nu}(Y\psi)\}
\end{align}
On the other hand, by (7.6), we have:
\begin{equation}
 (\frac{d}{ds}\Box_{f_{s*}\tilde{g}}(f_{s*}\psi))_{s=0}=(\frac{d}{ds}f_{s*}\rho)_{s=0}=(\frac{d}{ds}\rho\circ f_{s})_{s=0}=Y\rho
\end{equation}
Comparing (7.13) and (7.14) and in view of the expression (7.7) with $Y\psi$ in the role of $\psi$, the proposition follows. $\qed$

Let $X$ be an arbitrary vectorfield. In an arbitrary system of local coordinates the divergence of $X$ with respect to the 
acoustical metric $g$ is expressed by:
\begin{equation}
 \textrm{div}X=D_{\mu}X^{\mu}=\frac{1}{\sqrt{-\det g}}\frac{\partial}{\partial x^{\mu}}(\sqrt{-\det g}X^{\mu})
\end{equation}
while its divergence with respect to the conformal acoustical metric $\tilde{g}$ is expressed by:
\begin{equation}
 \tilde{\textrm{div}} X=\tilde{D}_{\mu}X^{\mu}=\frac{1}{\sqrt{-\det \tilde{g}}}\frac{\partial}{\partial x^{\mu}}(\sqrt{-\det \tilde{g}}X^{\mu})
\end{equation}
Since
\begin{equation}
 \tilde{g}_{\mu\nu}=\Omega g_{\mu\nu}, \quad\sqrt{-\det\tilde{g}}=\Omega^{2}\sqrt{-\det g}
\end{equation}
we have:
\begin{equation}
 \tilde{D}_{\mu}X^{\mu}=\Omega^{-2}D_{\mu}(\Omega^{2}X^{\mu})
\end{equation}
Applying this to the vectorfield $\leftexp{(Y)}{\tilde{J}}$, we obtain:
\begin{equation}
 \tilde{\textrm{div}}\leftexp{(Y)}{\tilde{J}}=\Omega^{-2}\textrm{div}\leftexp{(Y)}{J}
\end{equation}
where:
\begin{equation}
 \leftexp{(Y)}{J}=\Omega^{2}\leftexp{(Y)}{\tilde{J}}
\end{equation}
In an arbitrary local coordinates:
\begin{align}
 \leftexp{(Y)}{J}^{\mu}=\frac{1}{2}((g^{-1})^{\mu\alpha}(g^{-1})^{\nu\beta}+(g^{-1})^{\nu\alpha}(g^{-1})^{\mu\beta}
-(g^{-1})^{\mu\nu}(g^{-1})^{\alpha\beta})\leftexp{(Y)}{\tilde{\pi}_{\alpha\beta}}\partial_{\nu}\psi
\end{align}
Setting
\begin{equation}
 \leftexp{(Y)}{\tilde{\pi}}^{\mu\nu}=(g^{-1})^{\mu\alpha}(g^{-1})^{\nu\beta}\leftexp{(Y)}{\tilde{\pi}}_{\alpha\beta}
\end{equation}
we can write:
\begin{equation}
 \leftexp{(Y)}{J}^{\mu}=(\leftexp{(Y)}{\tilde{\pi}}^{\mu\nu}-\frac{1}{2}(g^{-1})^{\mu\nu}\textrm{tr}\leftexp{(Y)}{\tilde{\pi}})\partial_{\nu}\psi
\end{equation}

We now consider the $n$th order variations $\psi_{n}$ of the wave function $\phi$, as defined in the begining of Chapter 6, by applying
to a first order variation $\psi_{1}$ a string of commutation vectorfields $Y_{i}: i=1,2,3,4,5$ of length $n-1$. We shall use proposition 7.1 
to derive a recursion formula for the corresponding source functions $\rho_{n}$:
\begin{equation}
 \Box_{\tilde{g}}\psi_{n}=\rho_{n}
\end{equation}
First, we have:
\begin{equation}
 \rho_{1}=0
\end{equation}
Denoting by $Y$ one of the $Y_{i}: i=1,2,3,4,5$, we have:
\begin{equation}
 \psi_{n}=Y\psi_{n-1}
\end{equation}
Applying Proposition 7.1 to $\psi_{n-1}$ and $\psi_{n}$ we get:
\begin{equation}
 \rho_{n}=\tilde{\textrm{div}}\leftexp{(Y)}{\tilde{J}}_{n-1}+Y\rho_{n-1}+\frac{1}{2}\tilde{\textrm{tr}}\leftexp{(Y)}{\tilde{\pi}}\rho_{n-1}
\end{equation}
Here,
\begin{align}
 \leftexp{(Y)}{\tilde{J}}_{n-1}^{\mu}=\frac{1}{2}[(\tilde{g}^{-1})^{\mu\alpha}(\tilde{g}^{-1})^{\nu\beta}
+(\tilde{g}^{-1})^{\nu\alpha}(\tilde{g}^{-1})^{\mu\beta}-(\tilde{g}^{-1})^{\mu\nu}(\tilde{g}^{-1})^{\alpha\beta}]
\leftexp{(Y)}{\tilde{\pi}}_{\alpha\beta}\partial_{\nu}\psi_{n-1}
\end{align}
is the commutation current associated to $\psi_{n-1}$ and $Y$.
Equation (7.27) is a recursion formula for the source functions $\rho_{n}$, but it is not quite in the form which can be used in our estimats.
To obtain the appropriate form we consider instead the re-scaled sources:
\begin{equation}
 \tilde{\rho}_{n}=\Omega^{2}\mu\rho_{n}
\end{equation}
By direct computation and (7.27), (7.28) as well as (7.19)-(7.23), we have:
\begin{equation}
 \tilde{\rho}_{n}=\leftexp{(Y)}{\sigma}_{n-1}+Y\tilde{\rho}_{n-1}+\leftexp{(Y)}{\delta}\tilde{\rho}_{n-1}
\end{equation}
Here,
\begin{align}
 \leftexp{(Y)}{\sigma}_{n-1}=\mu \textrm{div}\leftexp{(Y)}{J}_{n-1}\\
\leftexp{(Y)}{J}^{\mu}_{n-1}=(\leftexp{(Y)}{\tilde{\pi}}^{\mu\nu}-\frac{1}{2}(g^{-1})^{\mu\nu}\textrm{tr}\leftexp{(Y)}{\tilde{\pi}})\partial_{\nu}\psi_{n-1}
\end{align}
and $\leftexp{(Y)}{\delta}$ are defined in Chapter 6:
\begin{equation}
 \leftexp{(Y)}{\delta}=\frac{1}{2}\tilde{\textrm{tr}}\leftexp{(Y)}{\tilde{\pi}}-\mu^{-1}Y\mu-2\Omega^{-1}Y\Omega
\end{equation}
Moreover, by (7.25):
\begin{equation}
 \tilde{\rho}_{1}=0
\end{equation}

\section{The First Term in $\tilde{\rho}_{n}$}
Since $\psi_{n}$ is the solution of (7.24), we can apply Theorem 5.1 to $\psi_{n}$, provided that we can estimate:
\begin{align}
 \int_{W^{t}_{u}}Q_{0,0,n}d\mu_{g}=-\int_{W^{t}_{u}}\Omega^{2}\rho_{n}K_{0}\psi_{n}d\mu_{g}=
-\int_{W^{t}_{u}}\tilde{\rho}_{n}(K_{0}\psi_{n})dt'du'd\mu_{\slashed{g}}
\end{align}
and
\begin{align}
\int_{W^{t}_{u}}Q_{1,0,n}d\mu_{g}=-\int_{W^{t}_{u}}\Omega^{2}\rho_{n}(K_{1}\psi_{n}+\omega\psi_{n})d\mu_{g}=
-\int_{W^{t}_{u}}\tilde{\rho}_{n}(K_{1}\psi_{n}+\omega\psi_{n})dt'du'd\mu_{\slashed{g}}
\end{align}
We have used here the fact that:
\begin{equation}
 d\mu_{g}=\mu dtdud\mu_{\slashed{g}}
\end{equation}

First, we consider the contribution of the term $\leftexp{(Y)}{\sigma}_{n-1}$ in $\tilde{\rho}_{n}$.

Let $V$ be an arbitrary vectorfield defined in the spacetime domain $W_{\epsilon_{0}}$. We decompose:
\begin{equation}
 V=V^{L}L+V^{\underline{L}}\underline{L}+\slashed{V}
\end{equation}
where $\slashed{V}=\Pi V=V^{A}X_{A}$ is tangent to $S_{t,u}$. We have:
\begin{align}
 V^{L}=-\frac{1}{2\mu}g(V,\underline{L})\\\notag
V^{\underline{L}}=-\frac{1}{2\mu}g(V,L)\\\notag
V^{A}=(\slashed{g}^{-1})^{AB}g(V,X_{B})
\end{align}
The divergence of $V$ is then expressed by:
\begin{equation}
 \textrm{div}V=(D_{L}V)^{L}+(D_{\underline{L}}V)^{\underline{L}}+(D_{X_{A}}V)^{A}
\end{equation}
Replacing $V$ in (7.39) by $D_{L}V, D_{\underline{L}}V, D_{X_{A}}V$, we obtain:
\begin{align}
 (D_{L}V)^{L}=-\frac{1}{2\mu}g(D_{L}V,\underline{L})\\\notag
(D_{\underline{L}}V)^{\underline{L}}=-\frac{1}{2\mu}g(D_{\underline{L}}V,L)\\\notag
(D_{X_{A}}V)^{A}=(\slashed{g}^{-1})^{AB}g(D_{X_{A}}V,X_{B})
\end{align}
Moreover, substituting (7.38) we obtain:
\begin{align}
 g(D_{X_{A}}V,X_{B})=V^{L}g(D_{X_{A}}L,X_{B})+V^{\underline{L}}g(D_{X_{A}}\underline{L},X_{B})+g(D_{X_{A}}\slashed{V},X_{B})\\\notag
=\chi_{AB}V^{L}+\underline{\chi}_{AB}V^{\underline{L}}+\slashed{g}(\slashed{D}_{X_{A}}\slashed{V},X_{B})
\end{align}
hence:
\begin{equation}
 (D_{X_{A}}V)^{A}=\textrm{tr}\chi V^{L}+\textrm{tr}\underline{\chi}V^{\underline{L}}+\slashed{\textrm{div}}\slashed{V}
\end{equation}
Thus we obtain:
\begin{equation}
 \textrm{div}V=-\frac{1}{2\mu}[g(D_{L}V,\underline{L})+g(D_{\underline{L}}V,L)]+
\textrm{tr}\chi V^{L}+\textrm{tr}\underline{\chi}V^{\underline{L}}+\slashed{\textrm{div}}\slashed{V}
\end{equation}
We can express the first term on the right-hand side in terms of:
\begin{equation}
 V_{L}=g(V,L),\quad V_{\underline{L}}=g(V,\underline{L})
\end{equation}
Using (3.131)-(3.138):
\begin{align}
 g(D_{L}V,\underline{L})=L(g(V,\underline{L}))-g(V,D_{L}\underline{L})
=L(V_{\underline{L}})+g(V,2\zeta^{A}X_{A})=L(V_{\underline{L}})+2\zeta_{A}V^{A}\\
g(D_{\underline{L}}V,L)=\underline{L}(g(V,L))-g(V,D_{\underline{L}}L)
=\underline{L}(V_{L})-g(V,-L(\eta^{-1}\kappa)L+2\eta^{A}X_{A})\\\notag
=\underline{L}(V_{L})+L(\eta^{-1}\kappa)V_{L}-2\eta_{A}V^{A}
\end{align}
Substituting the equation:
\begin{equation}
 \eta_{A}-\zeta_{A}=\slashed{d}_{A}\mu
\end{equation}
we have:
\begin{equation}
 -2\mu \textrm{div}V=L(V_{\underline{L}})+\underline{L}(V_{L})-2\slashed{\textrm{div}}(\mu\slashed{V})+
L(\eta^{-1}\kappa)V_{L}+\textrm{tr}\chi V_{\underline{L}}+\textrm{tr}\underline{\chi}V_{L}
\end{equation}
We shall apply the above formula to $\leftexp{(Y)}{J_{n-1}}$, given by (7.32). 
We introduce the vectorfields $\leftexp{(Y)}{\tilde{Z}}$ and $\leftexp{(Y)}{\underline{\tilde{Z}}}$, associated to 
the commutation vectorfield $Y$, and tangent to $S_{t,u}$, by the condition that:
\begin{equation}
 g(\tilde{Z},V)=\leftexp{(Y)}{\tilde{\pi}}(L,\Pi V), \quad g(\underline{\tilde{Z}},V)=\leftexp{(Y)}{\tilde{\pi}}(\underline{L},\Pi V)
\end{equation}
for any vector $V\in TW_{\epsilon_{0}}$. In terms of the null frames we have:
\begin{equation}
 \leftexp{(Y)}{\tilde{Z}}=\leftexp{(Y)}{\tilde{Z}}^{A}X_{A}, 
\leftexp{(Y)}{\tilde{\underline{Z}}}=\leftexp{(Y)}{\tilde{\underline{Z}}}^{A}X_{A}
\end{equation}
where:
\begin{equation}
 \leftexp{(Y)}{\tilde{Z}}^{A}=\leftexp{(Y)}{\tilde{\pi}}_{LB}(\slashed{g}^{-1})^{AB},\quad
\leftexp{(Y)}{\tilde{\underline{Z}}}^{A}=\leftexp{(Y)}{\tilde{\pi}}_{\underline{L}B}(\slashed{g}^{-1})^{AB}
\end{equation}
Note that the $S_{t,u}$-tangential vectorfields $\leftexp{(Y)}{\tilde{Z}}$ and $\leftexp{(Y)}{\tilde{\underline{Z}}}$ correspond, through the induced
metric $\slashed{g}$, to the 1-forms $\leftexp{(Y)}{\slashed{\pi}}_{L}$ and $\leftexp{(Y)}{\slashed{\pi}}_{\underline{L}}$, respectively.
Taking into account these definitions as well as the fact that:
\begin{equation}
 \textrm{tr}\leftexp{(Y)}{\tilde{\pi}}=-\frac{1}{\mu}\leftexp{(Y)}{\tilde{\pi}}_{L\underline{L}}
+\textrm{tr}\leftexp{(Y)}{\tilde{\slashed{\pi}}}
\end{equation}
we deduce the following the components of $\leftexp{(Y)}{J_{n-1}}$:
\begin{align}
 \leftexp{(Y)}{J_{n-1,L}}=-\frac{1}{2}\textrm{tr}\leftexp{(Y)}{\tilde{\slashed{\pi}}}(L\psi_{n-1})
+\leftexp{(Y)}{\tilde{Z}}\cdot\slashed{d}\psi_{n-1}\\
\leftexp{(Y)}{J_{n-1,\underline{L}}}=-\frac{1}{2}\textrm{tr}\leftexp{(Y)}{\tilde{\slashed{\pi}}}(\underline{L}\psi_{n-1})
+\leftexp{(Y)}{\underline{\tilde{Z}}}\cdot\slashed{d}\psi_{n-1}
-\frac{1}{2\mu}\leftexp{(Y)}{\tilde{\pi}}_{\underline{L}\underline{L}}(L\psi_{n-1})\\
\mu\leftexp{(Y)}{\slashed{J}}^{A}_{n-1}=-\frac{1}{2}\leftexp{(Y)}{\tilde{Z}}^{A}(\underline{L}\psi_{n-1})
-\frac{1}{2}\leftexp{(Y)}{\tilde{\underline{Z}}}^{A}(L\psi_{n-1})\\\notag
+\frac{1}{2}
(\leftexp{(Y)}{\tilde{\pi}}_{L\underline{L}}-\mu\textrm{tr}\leftexp{(Y)}{\tilde{\slashed{\pi}}})\slashed{d}^{A}\psi_{n-1}
+\mu\leftexp{(Y)}{\tilde{\slashed{\pi}}}^{A}_{B}\slashed{d}^{B}\psi_{n-1}
\end{align}
We note here that there is no $\frac{1}{\mu}\leftexp{(Y)}{\tilde{\pi}}_{L\underline{L}}$ in the expressions for 
$\leftexp{(Y)}{J_{n-1,L}},\leftexp{(Y)}{J_{n-1,\underline{L}}}$, although such a term appears in (7.53). This is related to the fact that
the operator $\Delta_{g}$ on the 2-dimensional Riemannian manifold $(M,g)$ is conformally covariant. Here in the role of such a 2-dimensional 
manifold we have the 2-dimensional distribution of time-like planes spanned by $L$ and $\underline{L}$ at each point. This distribution is not integrable, 
the obstruction to integrability being
\begin{equation}
 \Pi[L,\underline{L}]=2\Lambda, \quad \Lambda=-(\slashed{g}^{-1})^{AB}(\zeta_{B}+\eta_{B})X_{A}
\end{equation}
However, the conformal covariance is still reflected by the fact that the restriction of $\leftexp{(Y)}{J_{n-1}}$ to the plane spanned by
$L$, $\underline{L}$ depends only on the trace-free, relative to this plane, part of the restriction of $\tilde{\pi}$ to the plane, therefore 
not on $\tilde{\pi}_{L\underline{L}}$. Similarly, one can easily check that $\leftexp{(Y)}{\slashed{J}}^{A}_{n-1}$ does not actually depend on $\textrm{tr}\leftexp{(Y)}
{\tilde{\slashed{\pi}}}$.

Applying (7.49) to $\leftexp{(Y)}{J_{n-1}}$ we obtain:
\begin{align}
 \leftexp{(Y)}{\sigma_{n-1}}=-\frac{1}{2}L(\leftexp{(Y)}{J_{n-1,\underline{L}}})
-\frac{1}{2}\underline{L}(\leftexp{(Y)}{J_{n-1,L}})+\slashed{\textrm{div}}(\mu\leftexp{(Y)}{\slashed{J}_{n-1}})\\\notag
-\frac{1}{2}L(\eta^{-1}\kappa)\leftexp{(Y)}{J_{n-1,L}}-\frac{1}{2}\textrm{tr}\chi\leftexp{(Y)}{J_{n-1,\underline{L}}}
-\frac{1}{2}\textrm{tr}\underline{\chi}\leftexp{(Y)}{J_{n-1,L}}
\end{align}
The right hand side of (7.58) contains the principal terms, the derivatives of products of components of $\leftexp{(Y)}{\tilde{\pi}}$ 
with derivatives of $\psi_{n-1}$, while the second line contains lower order terms of the form of triple products of connection coefficients
of the null frame with components of $\leftexp{(Y)}{\tilde{\pi}}$ and derivatives of $\psi_{n-1}$. We decompose:
\begin{equation}
 \leftexp{(Y)}{\sigma_{n-1}}=\leftexp{(Y)}{\sigma_{1,n-1}}+\leftexp{(Y)}{\sigma_{2,n-1}}+\leftexp{(Y)}{\sigma_{3,n-1}}
\end{equation}
where $\leftexp{(Y)}{\sigma_{1,n-1}}$ contains the terms which are products of components of $\leftexp{(Y)}{\tilde{\pi}}$ with 
2nd derivatives of $\psi_{n-1}$, $\leftexp{(Y)}{\sigma_{2,n-1}}$ contains the terms which are products of 1st derivatives of 
$\leftexp{(Y)}{\tilde{\pi}}$ with 1st derivatives of $\psi_{n-1}$, and $\leftexp{(Y)}{\sigma_{3,n-1}}$ contains the lower order terms.

Now, in view of (7.54)-(7.56), the first two terms on the right of (7.58) contain $L(\leftexp{(Y)}{\underline{\tilde{Z}}}\cdot\slashed{d}\psi_{n-1})$
and $\underline{L}(\leftexp{(Y)}{\tilde{Z}}\cdot\slashed{d}\psi_{n-1})$. In decomposing each of these into three parts as above, we need the following
lemma.

$\textbf{Lemma 7.1}$ Let $V$ be a vectorfield defined on the spacetime domain $W^{*}_{\epsilon_{0}}$ and tangent to $S_{t,u}$.
Also, let $f$ be a function on $W^{*}_{\epsilon_{0}}$. We then have:
\begin{align*}
 L(V\cdot\slashed{d}f)=V\cdot\slashed{d}Lf+\slashed{\mathcal{L}}_{L}V\cdot\slashed{d}f\\
\underline{L}(V\cdot\slashed{d}f)=V\cdot\slashed{d}\underline{L}f+\slashed{\mathcal{L}}_{\underline{L}}V\cdot\slashed{d}f
-(V\cdot\slashed{d}(\eta^{-1}\kappa))Lf
\end{align*}
Here, we denote:
\begin{equation}
 \slashed{\mathcal{L}}_{L}V=\Pi\mathcal{L}_{L}V=\Pi[L,V],\quad
\slashed{\mathcal{L}}_{\underline{L}}V=\Pi\mathcal{L}_{\underline{L}}V=\Pi[\underline{L},V]
\end{equation}
$Proof$. Since $V$ is tangent to $S_{t,u}$ we have:
\begin{equation}
 V\cdot\slashed{d}f=V\cdot df
\end{equation}
Hence,
\begin{align}
 L(V\cdot\slashed{d}f)=L(V\cdot df)=\mathcal{L}_{L}V\cdot df+V\cdot dLf=\mathcal{L}_{L}V\cdot df
+V\cdot\slashed{d}Lf\\\notag
\underline{L}(V\cdot\slashed{d}f)=\underline{L}(V\cdot df)=\mathcal{L}_{\underline{L}}V\cdot df+V\cdot d\underline{L}f
=\mathcal{L}_{\underline{L}}V\cdot df+V\cdot\slashed{d}\underline{L}f
\end{align}
Since
\begin{align*}
 Lt=1\quad Lu=0
\end{align*}
and
\begin{align}
 S_{t,u}=\Sigma_{t}\bigcap C_{u}
\end{align}
and $V$ is $S_{t,u}$-tangential, thus
\begin{align}
 V(u)=V(t)=0
\end{align}
it follows that
\begin{align*}
 (\mathcal{L}_{L}V)(t)=L(Vt)-V(Lt)=0\\
(\mathcal{L}_{L}V)(u)=L(Vu)-V(Lu)=0
\end{align*}
So we have:
\begin{align}
 \mathcal{L}_{L}V=\slashed{\mathcal{L}}_{L}V
\end{align}
Concerning the term,
\begin{align}
 \mathcal{L}_{\underline{L}}V
\end{align}
 we use the fact that:
\begin{align}
 \underline{L}=\eta^{-1}\kappa L+2T
\end{align}
Then we have:
\begin{align}
 \mathcal{L}_{\underline{L}}V=(\eta^{-1}\kappa)\mathcal{L}_{L}V+2\mathcal{L}_{T}V-(V\cdot\slashed{d}(\eta^{-1}\kappa))L
\end{align}
By an argument similar to that leading to (7.65), we can see that
\begin{align}
 \mathcal{L}_{T}V=\slashed{\mathcal{L}}_{T}V
\end{align}
Therefore
\begin{align}
 \slashed{\mathcal{L}}_{\underline{L}}V=(\eta^{-1}\kappa)\slashed{\mathcal{L}}_{L}V+2\slashed{\mathcal{L}}_{T}V
\end{align}
So we obtain:
\begin{align}
 \mathcal{L}_{\underline{L}}V=\slashed{\mathcal{L}}_{\underline{L}}V-(V(\eta^{-1}\kappa))L
\end{align}
In view of (7.62), (7.65) and (7.71), the lemma follows. $\qed$

Applying Lemma 7.1 to the case $V=\leftexp{(Y)}{\underline{\tilde{Z}}}, f=\psi_{n-1}$, we obtain:
\begin{equation}
 L(\leftexp{(Y)}{\underline{\tilde{Z}}}\cdot\slashed{d}\psi_{n-1})=
\leftexp{(Y)}{\underline{\tilde{Z}}}\cdot\slashed{d}L\psi_{n-1}+
\slashed{\mathcal{L}}_{L}\leftexp{(Y)}{\underline{\tilde{Z}}}\cdot\slashed{d}\psi_{n-1}
\end{equation}
while applying Lemma 7.1 to the case $V=\leftexp{(Y)}{\tilde{Z}}, f=\psi_{n-1}$, yields:
\begin{equation}
 \underline{L}(\leftexp{(Y)}{\tilde{Z}}\cdot\slashed{d}\psi_{n-1})=
\leftexp{(Y)}{\tilde{Z}}\cdot\slashed{d}\underline{L}\psi_{n-1}
+\slashed{\mathcal{L}}_{\underline{L}}\leftexp{(Y)}{\tilde{Z}}\cdot\slashed{d}\psi_{n-1}
-\leftexp{(Y)}{\tilde{Z}}\cdot\slashed{d}(\eta^{-1}\kappa)L\psi_{n-1}
\end{equation}
Substituting (7.54)-(7.56), we get the following expressions in regard to the decomposition (7.59):
\begin{align}
\leftexp{(Y)}{\sigma_{1,n-1}}=\frac{1}{2}\textrm{tr}\leftexp{(Y)}{\tilde{\slashed{\pi}}}(L\underline{L}\psi_{n-1}
+\nu\underline{L}\psi_{n-1})\\\notag
+\frac{1}{4}(\mu^{-1}\leftexp{(Y)}{\tilde{\pi}}_{\underline{L}\underline{L}})L^{2}\psi_{n-1}\\\notag
-\leftexp{(Y)}{\tilde{Z}}\cdot\slashed{d}\underline{L}\psi_{n-1}
-\leftexp{(Y)}{\underline{\tilde{Z}}}\cdot\slashed{d}L\psi_{n-1}\\\notag
+\frac{1}{2}\leftexp{(Y)}{\tilde{\pi}}_{L\underline{L}}\slashed{\Delta}\psi_{n-1}+
\mu\leftexp{(Y)}{\hat{\tilde{\slashed{\pi}}}}\cdot\slashed{D}^{2}\psi_{n-1}
\end{align}
\begin{align}
 \leftexp{(Y)}{\sigma_{2,n-1}}=\frac{1}{4}L(\textrm{tr}\leftexp{(Y)}{\tilde{\slashed{\pi}}})\underline{L}\psi_{n-1}
+\frac{1}{4}\underline{L}(\textrm{tr}\leftexp{(Y)}{\tilde{\slashed{\pi}}})L\psi_{n-1}\\\notag
+\frac{1}{4}L(\mu^{-1}\leftexp{(Y)}{\tilde{\pi}}_{\underline{L}\underline{L}})L\psi_{n-1}\\\notag
-\frac{1}{2}\slashed{\mathcal{L}}_{L}\leftexp{(Y)}{\underline{\tilde{Z}}}\cdot\slashed{d}\psi_{n-1}
-\frac{1}{2}\slashed{\mathcal{L}}_{\underline{L}}\leftexp{(Y)}{\tilde{Z}}\cdot\slashed{d}\psi_{n-1}\\\notag
-\frac{1}{2}\slashed{\textrm{div}}\leftexp{(Y)}{\tilde{Z}}\underline{L}\psi_{n-1}-
\frac{1}{2}\slashed{\textrm{div}}\leftexp{(Y)}{\underline{\tilde{Z}}}L\psi_{n-1}\\\notag
+\frac{1}{2}\slashed{d}\leftexp{(Y)}{\tilde{\pi}}_{L\underline{L}}\cdot\slashed{d}\psi_{n-1}+
\slashed{\textrm{div}}(\mu\leftexp{(Y)}{\hat{\tilde{\slashed{\pi}}}})\cdot\slashed{d}\psi_{n-1}
\end{align}
and 
\begin{equation}
 \leftexp{(Y)}{\sigma_{3,n-1}}=\leftexp{(Y)}{\sigma_{3,n-1}^{L}}L\psi_{n-1}+
\leftexp{(Y)}{\sigma_{3,n-1}^{\underline{L}}}\underline{L}\psi_{n-1}+
\leftexp{(Y)}{\slashed{\sigma}_{3,n-1}}\cdot\slashed{d}\psi_{n-1}
\end{equation}
where
\begin{align}
 \leftexp{(Y)}{\sigma_{3,n-1}^{L}}=\frac{1}{4}\textrm{tr}\underline{\chi}\textrm{tr}\leftexp{(Y)}{\tilde{\slashed{\pi}}}
+\frac{1}{4}\textrm{tr}\chi(\mu^{-1}\leftexp{(Y)}{\tilde{\pi}}_{\underline{L}\underline{L}})\\\notag
+\frac{1}{2}\leftexp{(Y)}{\tilde{Z}}\cdot\slashed{d}(\eta^{-1}\kappa)\\
\leftexp{(Y)}{\sigma_{3,n-1}^{\underline{L}}}=-\frac{1}{4}(L\log\Omega)\textrm{tr}\leftexp{(Y)}{\tilde{\slashed{\pi}}}\\
\leftexp{(Y)}{\slashed{\sigma}_{3,n-1}}=-\frac{1}{2}(\textrm{tr}\leftexp{(Y)}{\tilde{\slashed{\pi}}})\Lambda\\\notag
-\frac{1}{2}(\textrm{tr}\underline{\chi}+L(\eta^{-1}\kappa))\leftexp{(Y)}{\tilde{Z}}-\frac{1}{2}\textrm{tr}\chi\leftexp{(Y)}{\underline{\tilde{Z}}}
\end{align}
Note that we have incorporated the term $\frac{1}{2}\textrm{tr}\leftexp{(Y)}{\tilde{\slashed{\pi}}}\nu\underline{L}\psi_{n-1}$ into the first term on the 
right in (7.74). This term comes from the contribution $\frac{1}{4}\textrm{tr}\chi\textrm{tr}\leftexp{(Y)}{\tilde{\slashed{\pi}}}\underline{L}\psi_{n-1}$
of the first term on the right in (7.55) to the second term in the second line in (7.58). This leaves as a remainder the coefficient (7.78) of 
$\underline{L}\psi_{n-1}$ in (7.76).

\section{The Estimates of the Contribution of the First Term in $\tilde{\rho}_{n}$ to the Error Integrals}
We shall begin our estimates of the contribution of $\leftexp{(Y)}{\sigma_{n-1}}$ to the error integrals (7.35), (7.36), with the estimates of 
the partial contribution of $\leftexp{(Y)}{\sigma_{1,n-1}}$. We shall use the following assumptions:

$\textbf{G0}$: There is a positive constant $C$ independent of $s$ such that in $W^{s}_{\epsilon_{0}}$, for all five commutation fields $Y$, we have:
\begin{align}
 \leftexp{(Y)}{\tilde{\pi}}_{LL}=0\\
\mu^{-1}|\leftexp{(Y)}{\tilde{\pi}}_{\underline{L}\underline{L}}|\leq C[1+\log(1+t)]\\
|\leftexp{(Y)}{\tilde{\pi}}_{L\underline{L}}|\leq C[1+\log(1+t)]\\
|\leftexp{(Y)}{\tilde{\slashed{\pi}}}_{L}|\leq C(1+t)^{-1}[1+\log(1+t)]\\
|\leftexp{(Y)}{\tilde{\slashed{\pi}}}_{\underline{L}}|\leq C[1+\log(1+t)]\\
|\leftexp{(Y)}{\hat{\tilde{\slashed{\pi}}}}|\leq C(1+t)^{-1}[1+\log(1+t)]\\
|\leftexp{(Y)}{\textrm{tr}\tilde{\slashed{\pi}}}|\leq C
\end{align}
This assumption comes from the results in Chapter 6.

We shall also use the assumptions concerning the set of rotation fields $\{R_{i}:i=1,2,3\}$.

$\textbf{H0}$: There is a positive constant $C$ independent of $s$ such that for any function $f$ differentiable on each surface $S_{t,u}$ we have:
\begin{equation}
 |\slashed{d}f|^{2}\leq C(1+t)^{-2}\sum_{i}(R_{i}f)^{2}
\end{equation}

$\textbf{H1}$: There is a positive constant $C$ independent of $s$ such that for any 1-form $\xi$ differentiable on each surface $S_{t,u}$ we have:
\begin{equation}
 |\slashed{D}\xi|^{2}\leq C(1+t)^{-2}\sum_{i}|\slashed{\mathcal{L}}_{R_{i}}\xi|^{2}
\end{equation}

In particular, taking $\xi=\slashed{d}f$, where $f$ is any function twice differentiable on each $S_{t,u}$, we have:
\begin{equation}
 |\slashed{D}^{2}f|^{2}\leq C(1+t)^{-2}\sum_{i}|\slashed{d}R_{i}f|^{2}
\end{equation}

$\textbf{H2}$: There is a positive constant $C$ independent of $s$ such that for any differentiable traceless symmetric 2-covariant tensorfield 
$\vartheta$ on each surface $S_{t,u}$ we have:
\begin{equation}
 |\slashed{D}\vartheta|^{2}\leq C(1+t)^{-2}\sum_{i}|\slashed{\mathcal{L}}_{R_{i}}\vartheta|^{2}
\end{equation}

Similarly as in Chapter 5 we define:
\begin{align*}
 \mathcal{E}^{u}_{0,n}(t)=\sum\int_{\Sigma^{u}_{t}}\frac{\Omega}{2}\{\eta^{-1}\kappa(1+\eta^{-1}\kappa)(L\psi_{n})^{2}
+(\underline{L}\psi_{n})^{2}+(1+2\eta^{-1}\kappa)\mu|\slashed{d}\psi_{n}|^{2}\}d\mu_{\slashed{g}}du^{\prime}\\
\mathcal{F}^{t}_{0,n}(u)=\sum\int_{C^{t}_{u}}\Omega\{(1+\eta^{-1}\kappa)(L\psi_{n})^{2}+\mu|\slashed{d}\psi_{n}|^{2}\}d\mu_{\slashed{g}}dt^{\prime}\\
  \mathcal{E}^{\prime u}_{1,n}(t)=\sum\int_{\Sigma^{u}_{t}}\frac{\Omega}{2}\omega\nu^{-1}\{\eta^{-1}\kappa(L\psi_{n}
+\nu\psi_{n})^{2}+\mu|\slashed{d}\psi_{n}|^{2}\}d\mu_{\slashed{g}}du^{\prime}\\
\mathcal{F}^{\prime t}_{1,n}(u)=\sum\int_{C^{t}_{u}}\Omega\omega\nu^{-1}(L\psi_{n}+\nu\psi_{n})^{2}d\mu_{\slashed{g}}dt^{\prime}
\end{align*}
In each of the above the sum is over the set of $\psi_{n}$ of the form (7.26) as $Y$ ranges over the set $\{Y_{i}:i=1,2,3,4,5\}$.
Also as in Chapter 5, we have:
\begin{align}
  C^{-1}\mathcal{E}^{u}_{0,n}(t)\leq \sum\int_{0}^{u}\{\int_{S_{t,u^{\prime}}}[\mu(1
+\mu)((L\psi_{n})^{2}+|\slashed{d}
\psi_{n}|^{2})+(\underline{L}\psi_{n})^{2}]d\mu_{\slashed{g}}\}du^{\prime}\leq C\mathcal{E}^{u}_{0,n}(t)\\
C^{-1}\mathcal{F}^{t}_{0,n}(u)\leq \sum\int_{0}^{t}\{\int_{S_{t^{\prime},u}}[(1+\mu)(L\psi_{n})^{2}
+\mu|\slashed{d}\psi_{n}|^{2}]d\mu_{\slashed{g}}
\}dt^{\prime}\leq C\mathcal{F}^{t}_{0,n}(u)\\
 C^{-1}\mathcal{E}^{\prime u}_{1,n}(t)\leq (1+t)^{2}\sum\int_{0}^{u}\{\int_{S_{t,u^{\prime}}}\mu[(L\psi_{n}
+\nu\psi_{n})^{2}+|\slashed{d}\psi_{n}|^{2}]
d\mu_{\slashed{g}}\}du^{\prime}\leq C\mathcal{E}^{\prime u}_{1,n}(t)\\
C^{-1}\mathcal{F}^{\prime t}_{1,n}(u)\leq \sum\int_{0}^{t}(1+t^{\prime})^{2}\{\int_{S_{t^{\prime},u}}(L\psi_{n}
+\nu\psi_{n})^{2}d\mu_{\slashed{g}}\}dt^{\prime}
\leq C\mathcal{F}^{\prime t}_{1,n}(u)
\end{align}
Also, the sum has the same meaning as the above. And we can define the following quantities which are non-decreasing in $t$ at each $u$
as in Chapter 5:
\begin{align}
 \bar{\mathcal{E}}^{u}_{0,n}(t)=\sup_{t'\in [0,t]}\mathcal{E}^{u}_{0,n}(t')\\
\mathcal{F}^{t}_{0,n}(u)\\
\bar{\mathcal{E}}'^{u}_{1,n}(t)=\sup_{t'\in[0,t]}[1+\log(1+t')]^{-4}\mathcal{E}'^{u}_{1,n}(t')\\
\bar{\mathcal{F}}'^{t}_{1,n}(u)=\sup_{t'\in[0,t]}[1+\log(1+t')]^{-4}\mathcal{F}'^{t'}_{1,n}(u)
\end{align}
Obviously, $\bar{\mathcal{E}}^{u}_{0,n}(t)$, $\bar{\mathcal{E}}'^{u}_{1,n}(t)$ are also non-decreasing functions of $u$ at each $t$. 

The contribution of $\leftexp{(Y)}{\sigma_{n-1}}$ to (7.35) and (7.36) are 
\begin{equation}
 -\int_{W^{t}_{u}}(K_{0}\psi_{n})\leftexp{(Y)}{\sigma_{n-1}}dt'du'd\mu_{\slashed{g}}
\end{equation}
and
\begin{equation}
 -\int_{W^{t}_{u}}(K_{1}\psi_{n}+\omega\psi_{n})\leftexp{(Y)}{\sigma_{n-1}}dt'du'd\mu_{\slashed{g}}
\end{equation}
respectively.

We shall first consider the contribution from $\leftexp{(Y)}{\sigma_{1,n-1}}$ to each of these integrals and we begin with the contribution to the integral (7.99).
We first consider the first term on the right-hand side of (7.74):
\begin{equation}
 \frac{1}{2}\textrm{tr}\leftexp{(Y)}{\tilde{\slashed{\pi}}}(L\underline{L}\psi_{n-1}+\nu\underline{L}\psi_{n-1})
\end{equation}
Now $\underline{L}$ can be expressed in terms of commutation vectorfields by:
\begin{equation}
 \underline{L}=\eta^{-1}\kappa(1+t)^{-1}Q+2T
\end{equation}
Thus we have:
\begin{align}
 \frac{1}{2}(L\underline{L}\psi_{n-1}+\nu\underline{L}\psi_{n-1})=(LT\psi_{n-1}+\nu T\psi_{n-1})\\\notag
+\frac{1}{2}\eta^{-1}\kappa(1+t)^{-1}(LQ\psi_{n-1}+\nu Q\psi_{n-1})\\\notag
+\frac{1}{2}(1+t)^{-1}(L(\eta^{-1}\kappa)-(1+t)^{-1}\eta^{-1}\kappa)Q\psi_{n-1}
\end{align}
Here we estimate the contribution of the first term on the right in (7.103), which is the leading term. The contributions of the other two terms are easier to 
estimate because of the presence of the decay factor $(1+t)^{-1}$.

We thus first estimate  
\begin{equation}
 -\int_{W^{t}_{u}}(K_{0}\psi_{n})\textrm{tr}\leftexp{(Y)}{\tilde{\slashed{\pi}}}(LT\psi_{n-1}+\nu T\psi_{n-1})dt'du'd\mu_{\slashed{g}}
\end{equation}
By the last of $\textbf{G0}$, this is bounded in absolute value by: $C(M^{L}+M^{\underline{L}})$, where
\begin{align}
 M^{L}=\int_{W^{t}_{u}}(1+\eta^{-1}\kappa)|L\psi_{n}||LT\psi_{n-1}+\nu T\psi_{n-1}|dt^{\prime}du^{\prime}d\mu_{\slashed{g}}\\
M^{\underline{L}}=\int_{W^{t}_{u}}|\underline{L}\psi_{n}||LT\psi_{n-1}+\nu T\psi_{n-1}|dt^{\prime}du^{\prime}d\mu_{\slashed{g}}
\end{align}
Since $T$ is one of the commutation fields, we have, by (7.94) and (7.98), for all $u\in [0,\epsilon_{0}]$,
\begin{equation}
 \int_{C^{t}_{u}}(1+t')^{2}(LT\psi_{n-1}+\nu T\psi_{n-1})^{2}d\mu_{\slashed{g}}dt'\leq
C\bar{\mathcal{F}}^{\prime t}_{1,n}(u)[1+\log(1+t)]^{4}
\end{equation}
Therefore
\begin{align}
 \int_{W^{t}_{u}}(1+t^{\prime})^{2}(LT\psi_{n-1}+\nu T\psi_{n-1})^{2}d\mu_{\slashed{g}}dt^{\prime}du^{\prime}
\leq C\int_{0}^{u}\bar{\mathcal{F}}^{\prime t}_{1,n}(u^{\prime})du^{\prime}[1+\log(1+t)]^{4}
\end{align}
We estimate:
\begin{align}
 M^{\underline{L}}\leq\int_{0}^{t}[(1+t^{\prime})^{-3/2}\mathcal{E}^{u}_{0}(t^{\prime})]^{1/2}
[(1+t^{\prime})^{3/2}\int_{\Sigma_{t^{\prime}}^{u}}(LT\psi_{n-1}+\nu T\psi_{n-1})^{2}]^{1/2}dt^{\prime}
\end{align}
Let us define:
\begin{align}
 F(t,u):=\int_{0}^{t}(1+t^{\prime})^{2}(\int_{\Sigma_{t^{\prime}}^{u}}(LT\psi_{n-1}+\nu T\psi_{n-1})^{2})dt^{\prime}
\end{align}
Therefore 
\begin{align}
 F(t,u)\leq C\int_{0}^{u}\bar{\mathcal{F}}^{\prime t}_{1,n}(u^{\prime})du^{\prime}[1+\log(1+t)]^{4}
\end{align}
Then
\begin{align}
 \frac{dF}{dt}(t,u)=(1+t)^{2}\int_{\Sigma_{t}^{u}}(LT\psi_{n-1}+\nu T\psi_{n-1})^{2}
\end{align}
What we need to estimate is 
\begin{align}
 I(t):=\int_{0}^{t}(1+t^{\prime})^{3/2}(\int_{\Sigma_{t}^{u}}(LT\psi_{n-1}+\nu T\psi_{n-1})^{2})dt^{\prime}
\end{align}
We have:
\begin{align}
 I(t)=\int_{0}^{t}(1+t^{\prime})^{-1/2}\frac{dF}{dt^{\prime}}(t^{\prime})dt^{\prime}\\\notag
=(1+t)^{-1/2}F(t)+\frac{1}{2}\int_{0}^{t}(1+t^{\prime})^{-3/2}F(t^{\prime})dt^{\prime}
\end{align}
Substitute (7.111) we obtain
\begin{align}
 I(t)\leq C\int_{0}^{u}\bar{\mathcal{F}}^{\prime t}_{1,n}(u^{\prime})du^{\prime}(1+t)^{-1/2}[1+\log(1+t)]^{4}
\leq C^{\prime}\int_{0}^{u}\bar{\mathcal{F}}^{\prime t}_{1,n}(u^{\prime})du^{\prime}
\end{align}
Therefore
\begin{align}
 M^{\underline{L}}\leq C(\int_{0}^{u}\bar{\mathcal{F}}^{\prime t}_{1,n}(u')du')^{1/2}
(\int_{0}^{t}(1+t')^{-3/2}\mathcal{E}^{u}_{0,n}(t')dt')^{1/2}
\end{align}

We turn to $M^{L}$ given by (7.105):
\begin{align}
 M^{L}\leq (\int_{W^{t}_{u}}\Omega(1+\eta^{-1}\kappa)(L\psi_{n})^{2}dt'du'd\mu_{\slashed{g}})^{1/2}\\\notag
\cdot(\int_{W^{t}_{u}}\Omega^{-1}(1+\eta^{-1}\kappa)(LT\psi_{n-1}+\nu T\psi_{n-1})^{2}dt'du'd\mu_{\slashed{g}})^{1/2}\\\notag
\leq C(\int_{0}^{u}\mathcal{F}^{t}_{0,n}(u')du')^{1/2}\cdot(\int_{W^{t}_{u}}[1+\log(1+t')](LT\psi_{n-1}+\nu T\psi_{n-1})^{2}dt'du'd\mu_{\slashed{g}})^{1/2}
\end{align}
where we have used assumption $\textbf{A}$.
We can estimate the integral in the last factor in a similar way as we estimate $M^{\underline{L}}$:
\begin{equation}
 \int_{W^{t}_{u}}[1+\log(1+t')](LT\psi_{n-1}+\nu T\psi_{n-1})^{2}dt'du'd\mu_{\slashed{g}}\leq C\int_{0}^{u}\bar{\mathcal{F}}'^{t}_{1,n}(u')du'
\end{equation}
Therefore we obtain:
\begin{equation}
 M^{L}\leq C(\int_{0}^{u}\mathcal{F}^{t}_{0,n}(u')du')^{1/2}(\int_{0}^{u}\bar{\mathcal{F'}}^{t}_{1,n}(u')du')^{1/2}
\end{equation}
This completes the estimate of the integral (7.104).

We then turn to estimate the contribution of the second term on the right in (7.74)
\begin{equation}
 (1/4)(\mu^{-1}\leftexp{(Y)}{\tilde{\pi}}_{\underline{L}\underline{L}})L^{2}\psi_{n-1}
\end{equation}
to the integral (7.99).

Writing 
\begin{align}
 L^{2}\psi_{n-1}=L((1+t)^{-1}Q\psi_{n-1})=(1+t)^{-1}LQ\psi_{n-1}-(1+t)^{-2}Q\psi_{n-1}
\end{align}
what we must estimate is 
\begin{align}
 -(1/4)\int_{W^{t}_{u}}(K_{0}\psi_{n})(\mu^{-1}\leftexp{(Y)}{\tilde{\pi}}_{\underline{L}\underline{L}})(1+t')^{-1}(LQ\psi_{n-1})dt'du'd\mu_{\slashed{g}}\\\notag
+(1/4)\int_{W^{t}_{u}}(K_{0}\psi_{n})(\mu^{-1}\leftexp{(Y)}{\tilde{\pi}}_{\underline{L}\underline{L}})(1+t')^{-2}(Q\psi_{n-1})dt'du'd\mu_{\slashed{g}}
\end{align}
By (7.81), the first integral is bounded in absolute value by
\begin{equation}
 C(I_{1}^{\underline{L}}+I_{1}^{L})
\end{equation}
where 
\begin{align}
 I^{\underline{L}}_{1}=\int_{W^{t}_{u}}(1+t')^{-1}|\underline{L}\psi_{n}|[1+\log(1+t')]|LQ\psi_{n-1}|dt'du'd\mu_{\slashed{g}}\\
I_{1}^{L}=\int_{W^{t}_{u}}(1+\eta^{-1}\kappa)|L\psi_{n}|(1+t')^{-1}[1+\log(1+t')]|LQ\psi_{n-1}|dt'du'd\mu_{\slashed{g}}
\end{align}
We estimate:
\begin{align}
 I^{\underline{L}}_{1}\leq C\{\int_{W^{t}_{u}}(\underline{L}\psi_{n})^{2}(1+t')^{-2}[1+\log(1+t')]^{2}dt'du'd\mu_{\slashed{g}}\}^{1/2}
\{\int_{W^{t}_{u}}(LQ\psi_{n-1})^{2}dt'du'd\mu_{\slashed{g}}\}^{1/2}\\\notag
\leq C\{\int_{0}^{t}(1+t')^{-2}[1+\log(1+t')]^{2}\mathcal{E}^{u}_{0,n}(t')dt'\}^{1/2}\{\int_{0}^{u}\mathcal{F}^{t}_{0,n}(u')du'\}^{1/2}
\end{align}
and,
\begin{equation}
 I_{1}^{L}\leq C\int_{0}^{u}\mathcal{F}^{t}_{0,n}(u')du'
\end{equation}
where we have used the fact that for any $t\geq 0, (1+t)^{-1}[1+\log(1+t)]$ is bounded. 

Also, by (7.81), the second integral in (7.122) is bounded by:
\begin{align}
 C\int_{W^{t}_{u}}|\underline{L}\psi_{n}||Q\psi_{n-1}|(1+t')^{-2}[1+\log(1+t')]dt'du'd\mu_{\slashed{g}}\\\notag
+C\int_{W^{t}_{u}}(1+\eta^{-1}\kappa)|L\psi_{n}||Q\psi_{n-1}|(1+t')^{-2}[1+\log(1+t')]dt'du'd\mu_{\slashed{g}}\\\notag
\leq C\{\int_{W^{t}_{u}}(\underline{L}\psi_{n})^{2}(1+t')^{-2}[1+\log(1+t')]dt'du'd\mu_{\slashed{g}}\}^{1/2}\\\notag
\cdot\{\int_{W^{t}_{u}}(Q\psi_{n-1})^{2}(1+t')^{-2}[1+\log(1+t')]dt'du'd\mu_{\slashed{g}}\}^{1/2}\\\notag
+C\{\int_{W^{t}_{u}}(1+\eta^{-1}\kappa)(L\psi_{n})^{2}dt'du'd\mu_{\slashed{g}}\}^{1/2}\\\notag
\cdot\{\int_{W^{t}_{u}}(1+\eta^{-1}\kappa)(Q\psi_{n-1})^{2}(1+t')^{-4}[1+\log(1+t')]^{2}dt'du'd\mu_{\slashed{g}}\}^{1/2}\\\notag
\leq C\{\int_{0}^{t}(1+t')^{-2}[1+\log(1+t')]\mathcal{E}^{u}_{0,n}(t')dt'\}^{1/2}\\\notag
\cdot\{\epsilon^{2}_{0}\int_{0}^{t}(1+t')^{-2}[1+\log(1+t')]\mathcal{E}^{u}_{0,n}(t')dt'\}^{1/2}\\\notag
+C\{\int_{0}^{u}\mathcal{F}^{t}_{0,n}(u')du'\}^{1/2}\{\epsilon_{0}^{2}\int_{0}^{t}(1+t')^{-4}[1+\log(1+t')]^{3}\mathcal{E}^{u}_{0,n}(t')dt'\}^{1/2}
\end{align}
Here we have made use of Lemma 5.1 with $\psi_{n}$ in the role of $\psi$, namely
\begin{align}
 \int_{S_{t,u}}\sum\psi_{n}^{2}d\mu_{\slashed{g}}\leq C\epsilon_{0}\mathcal{E}^{u}_{0,n}(t)
\end{align}
the sum being over the set of $\psi_{n}$ of the form (7.26) as $Y$ ranges over the set $\{Y_{i} : i=1,2,3,4,5\}$.

Next, we shall estimate the contribution of the third and fourth terms on the right in (7.74):
\begin{equation}
 -\leftexp{(Y)}{\tilde{Z}}\cdot\slashed{d}\underline{L}\psi_{n-1}-
\leftexp{(Y)}{\underline{\tilde{Z}}}\cdot\slashed{d}L\psi_{n-1}
\end{equation}
to the integral (7.99).

Expressing $\underline{L}$ and $L$ in terms of commutation fields $T$, $Q$,
\begin{equation}
 \underline{L}=\eta^{-1}\kappa(1+t)^{-1}Q+2T, \quad L=(1+t)^{-1}Q
\end{equation}
(7.130) takes the form:
\begin{align}
 -2\leftexp{(Y)}{\tilde{Z}}\cdot\slashed{d}T\psi_{n-1}-(1+t)^{-1}(\eta^{-1}\kappa\leftexp{(Y)}{\tilde{Z}}
+\leftexp{(Y)}{\underline{\tilde{Z}}})\cdot\slashed{d}Q\psi_{n-1}\\\notag
-(1+t)^{-1}\leftexp{(Y)}{\tilde{Z}}
\cdot\slashed{d}(\eta^{-1}\kappa)Q\psi_{n-1}
\end{align}
The last term is a lower order term and its contribution is easily estimated.

Using (7.83), (7.84) and noting that $T\psi_{n-1}, Q\psi_{n-1}$ are among the $\psi_{n}$, the contribution of (7.132) to the error integral (7.99) 
can be bounded by
\begin{align}
 \int_{W^{t}_{u}}|\underline{L}\psi_{n}||\slashed{d}\psi_{n}|(1+t')^{-1}[1+\log(1+t')]dt'du'd\mu_{\slashed{g}}\\\notag
+\int_{W^{t}_{u}}(1+\eta^{-1}\kappa)|L\psi_{n}||\slashed{d}\psi_{n}|(1+t')^{-1}[1+\log(1+t')]dt'du'd\mu_{\slashed{g}}
\end{align}
Here the first integral is bounded by:
\begin{equation}
 \{\int_{W^{t}_{u}}(1+t')^{-2}[1+\log(1+t')]^{2}(\underline{L}\psi_{n})^{2}dt'du'd\mu_{\slashed{g}}\}^{1/2}
\{\int_{W^{t}_{u}}|\slashed{d}\psi_{n}|^{2}dt'du'd\mu_{\slashed{g}}\}^{1/2}
\end{equation}
The integral in the first factor in (7.134) is bounded by:
\begin{equation}
 C\int_{0}^{t}(1+t')^{-2}[1+\log(1+t')]^{2}\mathcal{E}^{u}_{0,n}(t')dt'
\end{equation}
While the integral in the second factor decomposes into:
\begin{equation}
 \int_{\mathcal{U}\bigcap W^{t}_{u}}|\slashed{d}\psi_{n}|^{2}dt'du'd\mu_{\slashed{g}}+
\int_{\mathcal{U}^{c}\bigcap W^{t}_{u}}|\slashed{d}\psi_{n}|^{2}dt'du'd\mu_{\slashed{g}}
\end{equation}
where $\mathcal{U}$ is defined in Chapter 5:
\begin{equation}
 \mathcal{U}=\{x\in W^{*}_{\epsilon_{0}}:\mu <1/4\}
\end{equation}
Since $4\mu\geq 1$ in $\mathcal{U}^{c}$,
\begin{equation}
 \int_{\mathcal{U}^{c}\bigcap W^{t}_{u}}|\slashed{d}\psi_{n}|^{2}dt'du'd\mu_{\slashed{g}}
\leq 4\int_{W^{t}_{u}}\mu|\slashed{d}\psi_{n}|^{2}dt'du'd\mu_{\slashed{g}}
\leq C\int_{0}^{u}\mathcal{F}^{t}_{0,n}(u')du'
\end{equation}
To estimate the integral over $\mathcal{U}\bigcap W^{t}_{u}$ we use the spacetime integral $K_{n}(t,u)$ defined 
by (5.195), with $\psi_{n}$ in the role of $\psi$:
\begin{equation}
 K_{n}(t,u)=-\int_{W^{t}_{u}}\frac{\Omega}{2}\omega\nu^{-1}\mu^{-1}(L\mu)_{-}|\slashed{d}\psi_{n}|^{2}d\mu_{\slashed{g}}
\end{equation}
According to $\textbf{C3}, \textbf{B1}, \textbf{D1}$, we have:
\begin{equation}
 K_{n}(t,u)\geq \frac{1}{C}\int_{\mathcal{U}\bigcap W^{t}_{u}}(1+t')[1+\log(1+t')]^{-1}|\slashed{d}\psi_{n}|^{2}dt'du'd\mu_{\slashed{g}}
\end{equation}

Recalling the definition:
\begin{equation}
 \bar{K_{n}}(t,u)=\sup_{t'\in[0,t]}[1+\log(1+t')]^{-4}K_{n}(t',u)
\end{equation}
and noting that for $t'\in[t_{m-1},t_{m}]$ it holds that:
\begin{equation}
 (1+t')^{-1}[1+\log(1+t')]\leq 2(1+t_{m})^{-1}[1+\log(1+t_{m})]
\end{equation}
we then obtain:
\begin{equation}
 \int_{\mathcal{U}\bigcap W^{t_{m-1},t_{m}}_{u}}|\slashed{d}\psi_{n}|^{2}dt'du'd\mu_{\slashed{g}}\leq CA'_{m}\bar{K_{n}}(t_{m},u)
\end{equation}
where:
\begin{equation}
 A^{\prime}_{m}=(1+t_{m})^{-1}[1+\log(1+t_{m})]^{5}=(1+2^{m+r})^{-1}[1+\log(1+2^{m+r})]^{5}
\end{equation}
Hence, in view of the fact that $\sum_{m=0}^{\infty}A^{\prime}_{m}$ is convergent and $\bar{K_{n}}(t,u)$ is a non-decreasing function of $t$, we conclude that:
\begin{equation}
 \int_{\mathcal{U}\bigcap W^{t}_{u}}|\slashed{d}\psi_{n}|^{2}dt'du'd\mu_{\slashed{g}}=
\sum_{m=0}^{N}\int_{\mathcal{U}\bigcap W^{t_{m-1},t_{m}}_{u}}|\slashed{d}\psi_{n}|^{2}dt'du'd\mu_{\slashed{g}}
\leq C\bar{K_{n}}(t,u)
\end{equation}
Combining with (7.138) we get
\begin{equation}
 \int_{W^{t}_{u}}|\slashed{d}\psi_{n}|^{2}dt'du'd\mu_{\slashed{g}}\leq \bar{K_{n}}(t,u)+\int_{0}^{u}\mathcal{F}^{t}_{0,n}(u')du'
\end{equation}
Thus, the first integral in (7.133) is bounded in absolute value by
\begin{align}
 C(\int_{0}^{t}(1+t')^{-2}[1+\log(1+t')]^{2}\mathcal{E}^{u}_{0,n}(t')dt')^{1/2}(\bar{K_{n}}(t,u)+\int_{0}^{u}\mathcal{F}^{t}_{0,n}(u')du')^{1/2}
\end{align}
Also, in view of (7.92) and the estimate (7.146), the second integral in (7.133) is bounded in absolute value by:
\begin{align}
 \{\int_{W^{t}_{u}}(1+t')^{-2}[1+\log(1+t')]^{2}(1+\alpha^{-1}\kappa)^{2}(L\psi_{n})^{2}dt'du'd\mu_{\slashed{g}}\}^{1/2}
\cdot\{\int_{W^{t}_{u}}|\slashed{d}\psi_{n}|^{2}dt'du'd\mu_{\slashed{g}}\}^{1/2}\\\notag
\leq C\{\int_{0}^{u}\mathcal{F}^{t}_{0,n}(u')du'\}^{1/2}\{\bar{K_{n}}(t,u)+\int_{0}^{u}\mathcal{F}^{t}_{0,n}(u')du'\}^{1/2}
\end{align}
Finally, we estimate the contribution of the last term in (7.74),
\begin{equation}
 (1/2)\leftexp{(Y)}{\tilde{\pi}}_{L\underline{L}}\slashed{\Delta}\psi_{n-1}+\mu
\leftexp{(Y)}{\hat{\tilde{\slashed{\pi}}}}\cdot\slashed{D}^{2}\psi_{n-1}
\end{equation}
to the error integral (7.99).
From (7.82) and (7.85), this contribution can be bounded by
\begin{equation}
 C\int_{W^{t}_{u}}|K_{0}\psi_{n}||\slashed{D}^{2}\psi_{n-1}|[1+\log(1+t')]dt'du'd\mu_{\slashed{g}}
\end{equation}
By $\textbf{H1}$, 
\begin{equation}
 |\slashed{D}^{2}\psi_{n-1}|\leq C(1+t)^{-1}\sum_{i}|\slashed{d}R_{i}\psi_{n-1}|\leq C(1+t)^{-1}\sum|\slashed{d}\psi_{n}|
\end{equation}
since the $R_{i}\psi_{n-1}$ are among the $\psi_{n}$.
So (7.160) can be bounded by
\begin{equation}
 C\int_{W^{t}_{u}}|K_{0}\psi_{n}||\slashed{d}\psi_{n}|(1+t')^{-1}[1+\log(1+t')]dt'du'd\mu_{\slashed{g}}
\end{equation}
which is the integral (7.133) bounded by (7.148).

We collect the above results in the following lemma:

$\textbf{Lemma 7.2}$  We have:
\begin{align}
 |\int_{W^{t}_{u}}(K_{0}\psi_{n})\leftexp{(Y)}{\sigma}_{1,n-1}dt'du'd\mu_{\slashed{g}}|\\
\leq C\{\int_{0}^{t}(1+t')^{-3/2}\mathcal{E}^{u}_{0,n}(t')dt'+\int_{0}^{u}\mathcal{F}^{t}_{0,n}(u')du'\}\\
+C\{(\int_{0}^{u}\bar{\mathcal{F}}^{\prime t}_{1,n}(u')du')^{1/2}+(\bar{K_{n}}(t.u))^{1/2}\}\\
\times \{\int_{0}^{t}(1+t')^{-3/2}\mathcal{E}^{u}_{0,n}(t')dt'+\int_{0}^{u}\mathcal{F}^{t}_{0,n}(u')du'\}^{1/2}
\end{align}

    We now consider the contribution of $\leftexp{(Y)}{\sigma_{1,n-1}}$ to the integral (7.100). Since
\begin{equation}
 K_{1}=(\omega/\nu)L
\end{equation}
(7.100) equals
\begin{equation}
 -\int_{W^{t}_{u}}(\omega/\nu)(L\psi_{n}+\nu\psi_{n})\leftexp{(Y)}{\sigma_{n-1}}dt'du'd\mu_{\slashed{g}}
\end{equation}
What we consider here is
\begin{equation}
 -\int_{W^{t}_{u}}(\omega/\nu)(L\psi_{n}+\nu\psi_{n})\leftexp{(Y)}{\sigma_{1,n-1}}dt'du'd\mu_{\slashed{g}}
\end{equation}
By $\textbf{B1}$ and $\textbf{D1}$ we have:
\begin{equation}
 C^{-1}(1+t)^{2}\leq \omega/\nu\leq C(1+t)^{2}
\end{equation}
So (7.159) is bounded by:
\begin{equation}
 \int_{W^{t}_{u}}(1+t')^{2}|L\psi_{n}+\nu\psi_{n}||\leftexp{(Y)}{\sigma_{1,n-1}}|dt'du'd\mu_{\slashed{g}}
\end{equation}
We shall estimate this integral.

We begin with the contribution of the first term on the right in (7.74). We expand this term as in (7.103). So the contribution of the leading term in (7.103) is 
\begin{equation}
 \int_{W^{t}_{u}}(1+t')^{2}|L\psi_{n}+\nu\psi_{n}||\textrm{tr}\leftexp{(Y)}{\tilde{\slashed{\pi}}}|
|LT\psi_{n-1}+\nu T\psi_{n-1}|dt'du'd\mu_{\slashed{g}}
\end{equation}
By (7.86) and the fact that $T\psi_{n-1}$ is one of the $\psi_{n}$ this is bounded by
\begin{equation}
 C\int_{0}^{u}\mathcal{F}'^{t}_{1,n}(u')du'
\end{equation}
To estimate the contribution of the second term on the right in (7.74), we write:
\begin{align}
 L^{2}\psi_{n-1}=(1+t)^{-1}LQ\psi_{n-1}-(1+t)^{-2}Q\psi_{n-1}\\\notag
=(1+t)^{-1}(LQ\psi_{n-1}+\nu Q\psi_{n-1})-(1+t)^{-2}(1+(1+t)\nu)Q\psi_{n-1}
\end{align}
So by (7.81) and the fact that $Q\psi_{n-1}$ is one of the $\psi_{n}$, the contribution in question can be bounded by:
\begin{align}
C\int_{W^{t}_{u}}(1+t')[1+\log(1+t')]|L\psi_{n}+\nu\psi_{n}|^{2}dt'du'd\mu_{\slashed{g}}\\\notag
+C\int_{W^{t}_{u}}[1+\log(1+t')]|L\psi_{n}+\nu\psi_{n}||\psi_{n}|dt'du'd\mu_{\slashed{g}}\\\notag
\leq C\int_{0}^{u}\mathcal{F'}^{t}_{1,n}(u')du'+C\{\int_{0}^{u}\mathcal{F}^{\prime t}_{1,n}(u')du'\}^{1/2}
\{\epsilon^{2}_{0}\int_{0}^{t}(1+t')^{-2}[1+\log(1+t')]^{2}\mathcal{E}^{u}_{0,n}(t')dt'\}^{1/2}
\end{align}
To estimate the contribution of the third and fourth term on the right in (7.74), we shall use the expression (7.132). Using (7.83),(7.84) 
the contribution in question is bounded by:
\begin{equation}
 \int_{W^{t}_{u}}(1+t')[1+\log(1+t')]|L\psi_{n}+\nu\psi_{n}||\slashed{d}\psi_{n}|dt'du'd\mu_{\slashed{g}}
\end{equation}
which is in turn bounded by:
\begin{align}
 \{\int_{W^{t}_{u}}(1+t')^{2}|L\psi_{n}+\nu\psi_{n}|^{2}dt'du'd\mu_{\slashed{g}}\}^{1/2}\\\notag
\times\{\int_{W^{t}_{u}}[1+\log(1+t')]^{2}|\slashed{d}\psi_{n}|^{2}dt'du'd\mu_{\slashed{g}}\}^{1/2}
\end{align}
Now, the integral in the first factor in (7.167) is bounded by
\begin{align}
 C\int_{0}^{u}\mathcal{F}^{\prime t}_{1,n}(u')du^{\prime 1/2}
\end{align}
The integral in the second factor in (7.167) decomposes into:
\begin{align}
 \int_{W^{t}_{u}\bigcap\mathcal{U}}[1+\log(1+t')]^{2}|\slashed{d}\psi_{n}|^{2}dt'du'd\mu_{\slashed{g}}
+\int_{W^{t}_{u}\bigcap\mathcal{U}^{c}}[1+\log(1+t')]^{2}|\slashed{d}\psi_{n}|^{2}dt'du'd\mu_{\slashed{g}}
\end{align}
In $\mathcal{U}^{c}$ we have $4\mu\geq 1$; we can thus estimate, in view of (7.93),
\begin{align}
 \int_{W^{t}_{u}\bigcap\mathcal{U}^{c}}[1+\log(1+t')]^{2}|\slashed{d}\psi_{n}|^{2}dt'du'd\mu_{\slashed{g}}
\leq 4\int_{W^{t}_{u}}[1+\log(1+t')]^{2}\mu|\slashed{d}\psi_{n}|^{2}dt'du'd\mu_{\slashed{g}}\\
\leq C\int_{0}^{t}(1+t')^{-2}[1+\log(1+t')]^{2}\mathcal{E}^{\prime u}_{1,n}(t')dt'
\end{align}
on the other hand, by (7.140),
\begin{equation}
 \int_{\mathcal{U}\bigcap W^{t}_{u}}[1+\log(1+t')]^{2}|\slashed{d}\psi_{n}|^{2}dt'du'd\mu_{\slashed{g}}\leq CK_{n}(t,u)
\end{equation}
So we obtain:
\begin{align}
 \int_{W^{t}_{u}}[1+\log(1+t')]^{2}|\slashed{d}\psi_{n}|^{2}dt'du'd\mu_{\slashed{g}}\\\notag
\leq C\{K_{n}(t,u)+\int_{0}^{t}(1+t')^{-2}[1+\log(1+t')]^{2}\mathcal{E}^{\prime u}_{1,n}(t')dt'\}
\end{align}
Thus (7.166) is bounded by
\begin{equation}
 C\{\int_{0}^{u}\mathcal{F}^{\prime t}_{1,n}(u')du'\}^{1/2}\{K_{n}(t,u)+\int_{0}^{t}(1+t')^{-2}[1+\log(1+t')]^{2}\mathcal{E}^{\prime u}_{1,n}(t')dt'\}^{1/2}
\end{equation}

Finally, we estimate the contribution of the last two terms on the right in (7.74). By (7.82), (7.85) and (7.151) this is bounded by (7.166), which has already 
been estimated.

We collect the above results in the following lemma:

$\textbf{Lemma 7.3}$ We have:
\begin{align}
 |\int_{W^{t}_{u}}(K_{1}\psi_{n}+\omega\psi_{n})\leftexp{(Y)}{\sigma_{1,n-1}}dt'du'd\mu_{\slashed{g}}|\\\notag
\leq C\int_{0}^{u}\mathcal{F}^{\prime t}_{1,n}(u')du'\\\notag
+C\{\int_{0}^{u}\mathcal{F}^{\prime t}_{1,n}(u')du'\}^{1/2}\{K_{n}(t,u)+\int_{0}^{t}(1+t')^{-2}[1+\log(1+t')]^{2}(\mathcal{E}^{\prime u}_{1,n}(t')
+\epsilon_{0}^{2}\mathcal{E}^{u}_{0,n}(t'))dt'\}^{1/2}
\end{align}

To estimate the contribution of $\leftexp{(Y)}{\sigma_{2,n-1}}$ to (7.99) and (7.100), we shall use the following assumptions:

$\textbf{G1}$: There is a positive constant $C$ independent of $s$ such that in $W^{s}_{\epsilon_{0}}$, for all five commutation vectorfields
we have:
\begin{align}
 Y(\leftexp{(Y)}{\tilde{\pi}}_{LL})=0\\
|Y(\mu^{-1}\leftexp{(Y)}{\tilde{\pi}}_{\underline{L}\underline{L}})|\leq C[1+\log(1+t)]\\
|Y(\leftexp{(Y)}{\tilde{\pi}}_{L\underline{L}})|\leq C[1+\log(1+t)]\\
|\slashed{\mathcal{L}}_{Y}(\leftexp{(Y)}{\tilde{\slashed{\pi}}}_{L})|\leq C(1+t)^{-1}[1+\log(1+t)]\\
|\slashed{\mathcal{L}}_{Y}(\leftexp{(Y)}{\tilde{\slashed{\pi}}}_{\underline{L}})|\leq C[1+\log(1+t)]\\
|\slashed{\mathcal{L}}_{Y}(\leftexp{(Y)}{\hat{\tilde{\slashed{\pi}}}})|\leq C(1+t)^{-1}[1+\log(1+t)]\\
|Y(\leftexp{(Y)}{\textrm{tr}\tilde{\slashed{\pi}}})|\leq C(1+t)^{-1}[1+\log(1+t)]
\end{align}
Here in each line $Y$ occurs twice. In each of these occurances $Y$ is meant to range $independently$ over the set of the five commutation fields
$\{Y_{i}:i=1,2,3,4,5\}$.

The last of $\textbf{G1}$ seems stronger than what corresponds to the last of $\textbf{G0}$, however, recall that by (6.114), (6.122) and (6.195) of Chapter 6,
we have:
\begin{equation}
 |\leftexp{(Y_{i})}{\textrm{tr}\tilde{\slashed{\pi}}}|\leq C(1+t)^{-1}[1+\log(1+t)]: i\slashed{=} 2 \notag
\end{equation}
while for $Y_{2}=Q$ we have:
\begin{equation}
 |\textrm{tr}\leftexp{(Q)}{\tilde{\slashed{\pi}}}-4|\leq C(1+t)^{-1}[1+\log(1+t)] \notag
\end{equation}
So the last of $\textbf{G1}$ is in fact in accordance with the other assumptions.

In the following we shall use the bounds:
\begin{align}
 \int_{S_{t,u}}(\underline{L}\psi_{n-1})^{2}d\mu_{\slashed{g}}\leq C\epsilon_{0}\mathcal{E}^{u}_{0,n}(t)\\ \notag
\int_{S_{t,u}}(L\psi_{n-1})^{2}d\mu_{\slashed{g}}\leq C\epsilon_{0}(1+t)^{-2}\mathcal{E}^{u}_{0,n}(t)\\ \notag
\int_{S_{t,u}}|\slashed{d}\psi_{n-1}|^{2}d\mu_{\slashed{g}}\leq C\epsilon_{0}(1+t)^{-2}\mathcal{E}^{u}_{0,n}(t)
\end{align}
These follow from (7.129) using the expressions of $L$ and $\underline{L}$ in terms of the commutation fields $T$, $Q$, as well as $\textbf{H0}$.

We consider now the error integral:
\begin{equation}
 -\int_{W^{t}_{u}}(K_{0}\psi_{n})\leftexp{(Y)}{\sigma_{2,n-1}}dt'du'd\mu_{\slashed{g}}
\end{equation}
The first term in (7.75),
\begin{equation}
 (1/4)L(\textrm{tr}\leftexp{(Y)}{\tilde{\slashed{\pi}}})\underline{L}\psi_{n-1} \notag
\end{equation}
involves $L(\textrm{tr}\leftexp{(Y)}{\tilde{\slashed{\pi}}})$. Writing $L=(1+t)^{-1}Q$, we have, from $\textbf{G1}$,
\begin{equation}
 |L(\textrm{tr}\leftexp{(Y)}{\tilde{\slashed{\pi}}})|\leq C(1+t)^{-2}[1+\log(1+t)]
\end{equation}
It follows that the contribution of the first term to the error integral (7.184) is bounded in absolute value by:
\begin{equation}
 C(I^{\underline{L}}_{2}+I^{L}_{2}) \notag
\end{equation}
where
\begin{align}
 I^{\underline{L}}_{2}=\int_{W^{t}_{u}}|\underline{L}\psi_{n}||\underline{L}\psi_{n-1}|
(1+t')^{-2}[1+\log(1+t')]dt'du'd\mu_{\slashed{g}}\\
I^{L}_{2}=\int_{W^{t}_{u}}(1+\eta^{-1}\kappa)|L\psi_{n}||\underline{L}\psi_{n-1}|
(1+t')^{-2}[1+\log(1+t')]dt'du'd\mu_{\slashed{g}}
\end{align}
Using the first of (7.183) we can estimate:
\begin{align}
 I^{\underline{L}}_{2}\leq \int_{0}^{t}\{\int_{0}^{u}\int_{S_{t',u'}}(\underline{L}\psi_{n})^{2}d\mu_{\slashed{g}}du'\}^{1/2}
\{\int_{0}^{u}\int_{S_{t',u'}}(\underline{L}\psi_{n-1})^{2}d\mu_{\slashed{g}}du'\}^{1/2}\\ \notag
\cdot (1+t')^{-2}[1+\log(1+t')]dt'\\ \notag
\leq C\int_{0}^{t}\{\mathcal{E}^{u}_{0,n}(t')\}^{1/2}\{\epsilon_{0}^{2}\mathcal{E}^{u}_{0,n}(t')\}^{1/2}(1+t')^{-2}[1+\log(1+t')]dt'\\ \notag
=C\epsilon_{0}\int_{0}^{t}(1+t')^{-2}[1+\log(1+t')]\mathcal{E}^{u}_{0,n}(t')dt'
\end{align}
and (using $\textbf{A}$):
\begin{align}
 I^{L}_{2}\leq \{\int_{W^{t}_{u}}(1+\eta^{-1}\kappa)(L\psi_{n})^{2}dt'du'd\mu_{\slashed{g}}\}^{1/2}\\ \notag
\cdot \{\int_{W^{t}_{u}}(1+\eta^{-1}\kappa)(\underline{L}\psi_{n-1})^{2}(1+t')^{-4}[1+\log(1+t')]^{2}dt'du'd\mu_{\slashed{g}}\}^{1/2}\\ \notag
\leq C\epsilon_{0}\{\int_{0}^{u}\mathcal{F}^{t}_{0,n}(u')du'\}^{1/2}\{\int_{0}^{t}(1+t')^{-4}[1+\log(1+t')]^{3}\mathcal{E}^{u}_{0,n}(t')dt'\}^{1/2}
\end{align}
The second term in (7.75),
\begin{equation}
 (1/4)\underline{L}(\textrm{tr}\leftexp{(Y)}{\tilde{\slashed{\pi}}})L\psi_{n-1}
\end{equation}
involves $\underline{L}(\textrm{tr}\leftexp{(Y)}{\tilde{\slashed{\pi}}})$. Writing $\underline{L}=2T+\alpha^{-1}\kappa(1+t)^{-1}Q$, we have, by $\textbf{G1}$,
\begin{equation}
 |\underline{L}(\textrm{tr}\leftexp{(Y)}{\tilde{\slashed{\pi}}})|\leq C(1+t)^{-1}[1+\log(1+t)]
\end{equation}
So the contribution of the second term is bounded in absolute value by:
\begin{equation}
 C(I^{\underline{L}}_{3}+I^{L}_{3}) \notag
\end{equation}
where:
\begin{align}
 I^{\underline{L}}_{3}=\int_{W^{t}_{u}}|\underline{L}\psi_{n}||L\psi_{n-1}|(1+t')^{-1}[1+\log(1+t')]dt'du'd\mu_{\slashed{g}}\\
I^{L}_{3}=\int_{W^{t}_{u}}(1+\eta^{-1}\kappa)|L\psi_{n}||L\psi_{n-1}|(1+t')^{-1}[1+\log(1+t')]dt'du'd\mu_{\slashed{g}}
\end{align}
using the second of (7.183), we can estimate:
\begin{align}
 I^{\underline{L}}_{3}\leq \int_{0}^{t}\{\int_{0}^{u}\int_{S_{t',u'}}(\underline{L}\psi_{n})^{2}d\mu_{\slashed{g}}du'\}^{1/2}
\{\int_{0}^{u}\int_{S_{t',u'}}(L\psi_{n-1})^{2}d\mu_{\slashed{g}}du'\}^{1/2}\\ \notag
\cdot (1+t')^{-1}[1+\log(1+t')]dt'\\ \notag
\leq C\int_{0}^{t}\{\mathcal{E}^{u}_{0,n}(t')\}^{1/2}\{\epsilon_{0}^{2}(1+t')^{-2}\mathcal{E}^{u}_{0,n}(t')\}^{1/2}(1+t')^{-1}[1+\log(1+t')]dt'\\ 
\notag
=C\epsilon_{0}\int_{0}^{t}(1+t')^{-2}[1+\log(1+t')]\mathcal{E}^{u}_{0,n}(t')dt'
\end{align}
and (using $\textbf{A}$):
\begin{align}
 I^{L}_{3}\leq \{\int_{W^{t}_{u}}(1+\eta^{-1}\kappa)(L\psi_{n})^{2}dt'du'd\mu_{\slashed{g}}\}^{1/2}\\ \notag
\cdot \{\int_{W^{t}_{u}}(1+\eta^{-1}\kappa)(L\psi_{n-1})^{2}(1+t')^{-2}[1+\log(1+t')]^{2}dt'du'd\mu_{\slashed{g}}\}^{1/2}\\ \notag
\leq C\epsilon_{0}\{\int_{0}^{u}\mathcal{F}^{t}_{0,n}(u')du'\}^{1/2}\{\int_{0}^{t}(1+t')^{-4}[1+\log(1+t')]^{3}\mathcal{E}^{u}_{0,n}(t')dt'\}^{1/2}
\end{align}

The third term in (7.75),
\begin{equation}
 (1/4)L(\mu^{-1}\leftexp{(Y)}{\tilde{\pi}}_{\underline{L}\underline{L}})L\psi_{n-1} \notag
\end{equation}
involves $L(\mu^{-1}\leftexp{(Y)}{\tilde{\pi}}_{\underline{L}\underline{L}})$. Writing $L=(1+t)^{-1}Q$, we have, from $\textbf{G1}$:
\begin{equation}
 |L(\mu^{-1}\leftexp{(Y)}{\tilde{\pi}}_{\underline{L}\underline{L}})|\leq C(1+t)^{-1}[1+\log(1+t)]
\end{equation}
It follows that the contribution in question can be also bounded in absolute value by:
\begin{equation}
 C(I^{\underline{L}}_{3}+I^{L}_{3}) \notag
\end{equation}

The term fourth and fifth terms in (7.75),
\begin{equation}
 -\frac{1}{2}\slashed{\mathcal{L}}_{L}\leftexp{(Y)}{\underline{\tilde{Z}}}\cdot\slashed{d}\psi_{n-1}
-\frac{1}{2}\slashed{\mathcal{L}}_{\underline{L}}\leftexp{(Y)}{\tilde{Z}}\cdot\slashed{d}\psi_{n-1} \notag
\end{equation}
involve $\slashed{\mathcal{L}}_{L}\leftexp{(Y)}{\underline{\tilde{Z}}}$ and 
$\slashed{\mathcal{L}}_{\underline{L}}\leftexp{(Y)}{\tilde{Z}}$. 
Now from (7.50)-(7.52), we have:
\begin{equation}
 \leftexp{(Y)}{\tilde{\slashed{\pi}}}_{L}=\slashed{g}\cdot\leftexp{(Y)}{\tilde{Z}},
 \leftexp{(Y)}{\tilde{\slashed{\pi}}}_{\underline{L}}=\slashed{g}\cdot\leftexp{(Y)}{\tilde{\underline{Z}}} \notag
\end{equation}
hence:
\begin{align}
 \slashed{\mathcal{L}}_{Y}\leftexp{(Y)}{\tilde{\slashed{\pi}}}_{L}=\slashed{g}\cdot\slashed{\mathcal{L}}_{Y}\leftexp{(Y)}{\tilde{Z}}
+\leftexp{(Y)}{\slashed{\pi}}\cdot\leftexp{(Y)}{\tilde{Z}}\\ \notag
\slashed{\mathcal{L}}_{Y}\leftexp{(Y)}{\tilde{\slashed{\pi}}}_{\underline{L}}=\slashed{g}\cdot\slashed{\mathcal{L}}_{Y}\leftexp{(Y)}{\tilde{\underline{Z}}}
+\leftexp{(Y)}{\slashed{\pi}}\cdot\leftexp{(Y)}{\tilde{\underline{Z}}}
\end{align}
Since 
\begin{equation}
 \leftexp{(Y)}{\tilde{\pi}}=\Omega(\leftexp{(Y)}{\pi}+Y(\log\Omega)g) \notag
\end{equation}
and  by $\textbf{E1}$ and $\textbf{E2}$ we have:
\begin{equation}
 |Yh|\leq C(1+t)^{-1}
\end{equation}
which implies
\begin{equation}
 |Y\log\Omega|\leq C(1+t)^{-1}
\end{equation}
for all five commutation vectorfields $Y$, then by $\textbf{G0}$ and $\textbf{G1}$ we have:
\begin{align}
 |\slashed{\mathcal{L}}_{Y}\leftexp{(Y)}{\tilde{Z}}|\leq C(1+t)^{-1}[1+\log(1+t)]\\ \notag
|\slashed{\mathcal{L}}_{Y}\leftexp{(Y)}{\tilde{\underline{Z}}}|\leq C[1+\log(1+t)]
\end{align}
These imply, writing $\underline{L}=2T+\eta^{-1}\kappa(1+t)^{-1}Q$, $L=(1+t)^{-1}Q$, that
\begin{equation}
 |\slashed{\mathcal{L}}_{\underline{L}}\leftexp{(Y)}{\tilde{Z}}|,
|\slashed{\mathcal{L}}_{L}\leftexp{(Y)}{\tilde{\underline{Z}}}|\leq C(1+t)^{-1}[1+\log(1+t)]
\end{equation}
So the contribution of the fourth and fifth terms in (7.75) can be bounded in absolute value by:
\begin{equation}
 C(I^{\underline{L}}_{4}+I^{L}_{4}) \notag
\end{equation}
where:
\begin{align}
 I^{\underline{L}}_{4}=\int_{W^{t}_{u}}|\underline{L}\psi_{n}||\slashed{d}\psi_{n-1}|(1+t')^{-1}[1+\log(1+t')]dt'du'd\mu_{\slashed{g}}\\
I^{L}_{4}=\int_{W^{t}_{u}}(1+\eta^{-1}\kappa)|L\psi_{n}||\slashed{d}\psi_{n-1}|(1+t')^{-1}[1+\log(1+t')]dt'du'd\mu_{\slashed{g}}
\end{align}
Using the last of (7.183), we can estimate:
\begin{align}
 I^{\underline{L}}_{4}\leq \int_{0}^{t}\{\int_{0}^{u}\int_{S_{t',u'}}(\underline{L}\psi_{n})^{2}d\mu_{\slashed{g}}du'\}^{1/2}
\{\int_{0}^{u}\int_{S_{t',u'}}|\slashed{d}\psi_{n-1}|^{2}d\mu_{\slashed{g}}du'\}^{1/2}\\ \notag
\cdot(1+t')^{-1}[1+\log(1+t')]dt'\\ \notag
\leq C\int_{0}^{t}\{\mathcal{E}^{u}_{0,n}(t')\}^{1/2}\{\epsilon_{0}^{2}(1+t')^{-2}\mathcal{E}^{u}_{0,n}(t')\}^{1/2}(1+t')^{-1}[1+\log(1+t')]dt'\\ \notag
=C\epsilon_{0}\int_{0}^{t}(1+t')^{-2}[1+\log(1+t')]\mathcal{E}^{u}_{0,n}(t')dt'
\end{align}
and (using $\textbf{A}$)
\begin{align}
 I^{L}_{4}\leq\{\int_{W^{t}_{u}}(1+\eta^{-1}\kappa)(L\psi_{n})^{2}dt'du'd\mu_{\slashed{g}}\}^{1/2}\\ \notag
\{\int_{W^{t}_{u}}(1+\eta^{-1}\kappa)|\slashed{d}\psi_{n-1}|^{2}(1+t')^{-2}[1+\log(1+t')]^{2}dt'du'd\mu_{\slashed{g}}\}^{1/2}\\ \notag
\leq C\epsilon_{0}\{\int_{0}^{u}\mathcal{F}^{t}_{0,n}(u')du'\}^{1/2}\{\int_{0}^{t}(1+t')^{-4}[1+\log(1+t')]^{3}\mathcal{E}^{u}_{0,n}(t')dt'\}^{1/2}
\end{align}

The sixth and seventh terms in (7.75) involve $\slashed{\textrm{div}}\leftexp{(Y)}{\tilde{Z}}$ and $\slashed{\textrm{div}}\leftexp{(Y)}{\underline{\tilde{Z}}}$. 
Now by $\textbf{H1}$ applied to 
the 1-forms $\tilde{\slashed{\pi}}_{L}$, $\tilde{\slashed{\pi}}_{\underline{L}}$ on each $S_{t,u}$ yields:
\begin{align}
 |\slashed{D}\leftexp{(Y)}{\tilde{Z}}|^{2}=|\slashed{D}\leftexp{(Y)}{\tilde{\slashed{\pi}}}_{L}|^{2}
\leq C(1+t)^{-2}\sum_{i}|\mathcal{L}_{R_{i}}\leftexp{(Y)}{\tilde{\slashed{\pi}}}_{L}|^{2}\\ \notag
|\slashed{D}\leftexp{(Y)}{\tilde{\underline{Z}}}|^{2}=|\slashed{D}\leftexp{(Y)}{\tilde{\slashed{\pi}}}_{\underline{L}}|^{2}
\leq C(1+t)^{-2}\sum_{i}|\mathcal{L}_{R_{i}}\leftexp{(Y)}{\tilde{\slashed{\pi}}}_{\underline{L}}|^{2}
\end{align}
Then from $\textbf{G1}$:
\begin{align}
 |\slashed{D}\leftexp{(Y)}{\tilde{Z}}|\leq C(1+t)^{-2}[1+\log(1+t)]\\ \notag
|\slashed{D}\leftexp{(Y)}{\tilde{\underline{Z}}}|\leq C(1+t)^{-1}[1+\log(1+t)]
\end{align}
Thus, also:
\begin{align}
 |\slashed{\textrm{div}}\leftexp{(Y)}{\tilde{Z}}|\leq C(1+t)^{-2}[1+\log(1+t)]\\ \notag
|\slashed{\textrm{div}}\leftexp{(Y)}{\tilde{\underline{Z}}}|\leq C(1+t)^{-1}[1+\log(1+t)]
\end{align}
It then follows that the contribution of the sixth term in (7.75) 
\begin{equation}
 -(1/2)\slashed{\textrm{div}}\leftexp{(Y)}{\tilde{Z}}\underline{L}\psi_{n-1} \notag
\end{equation}
is bounded in absolute value by:
\begin{equation}
 C(I^{L}_{2}+I^{\underline{L}}_{2}) \notag
\end{equation}
while the contribution of the seventh term in (7.75)
\begin{equation}
 -(1/2)\slashed{\textrm{div}}\leftexp{(Y)}{\underline{\tilde{Z}}}L\psi_{n-1} \notag
\end{equation}
is bounded in absolute value by:
\begin{equation}
 C(I^{L}_{3}+I^{\underline{L}}_{3}) \notag
\end{equation}
The eighth term in (7.75),
\begin{equation}
 (1/2)\slashed{d}\leftexp{(Y)}{\tilde{\pi}}_{L\underline{L}}\cdot\slashed{d}\psi_{n-1} \notag
\end{equation}
involves $\slashed{d}\leftexp{(Y)}{\tilde{\pi}}_{L\underline{L}}$. By $\textbf{H0}$ and $\textbf{G1}$,
\begin{equation}
 |\slashed{d}\leftexp{(Y)}{\tilde{\pi}}_{L\underline{L}}|\leq C(1+t)^{-1}\sqrt{\sum_{i}(R_{i}\leftexp{(Y)}{\tilde{\pi}}_{L\underline{L}})^{2}}
\leq C(1+t)^{-1}[1+\log(1+t)]
\end{equation}
It then follows that the contribution of this term is bounded by:
\begin{equation}
 C(I^{L}_{4}+I^{\underline{L}}_{4}) \notag
\end{equation}

Finally, the last term in (7.75),
\begin{align}
 \slashed{\textrm{div}}(\mu\leftexp{(Y)}{\hat{\tilde{\slashed{\pi}}}})\cdot\slashed{d}\psi_{n-1} \notag
\end{align}
involves
\begin{equation}
 \slashed{\textrm{div}}(\mu\leftexp{(Y)}{\hat{\tilde{\slashed{\pi}}}})=\mu\slashed{\textrm{div}}\leftexp{(Y)}{\hat{\tilde{\slashed{\pi}}}}
+\slashed{d}\mu\cdot\leftexp{(Y)}{\hat{\tilde{\slashed{\pi}}}} \notag
\end{equation}
Now, $\textbf{H2}$ applied to $\leftexp{(Y)}{\hat{\tilde{\slashed{\pi}}}}$ gives:
\begin{equation}
 |\slashed{D}\leftexp{(Y)}{\hat{\tilde{\slashed{\pi}}}}|^{2}\leq C(1+t)^{-2}\sum_{i}|\slashed{\mathcal{L}}_{R_{i}}\leftexp{(Y)}{\hat{\tilde{\slashed{\pi}}}}|^{2}
\end{equation}
Then $\textbf{G1}$ implies:
\begin{equation}
 |\slashed{D}\leftexp{(Y)}{\hat{\tilde{\slashed{\pi}}}}|\leq C(1+t)^{-2}[1+\log(1+t)]
\end{equation}
thus also:
\begin{equation}
 |\slashed{\textrm{div}}\leftexp{(Y)}{\hat{\tilde{\slashed{\pi}}}}|\leq C(1+t)^{-2}[1+\log(1+t)]
\end{equation}
Using then $\textbf{A3},\textbf{F1}$ and $\textbf{G0}$ we obtain:
\begin{equation}
 |\slashed{\textrm{div}}(\mu\leftexp{(Y)}{\hat{\tilde{\slashed{\pi}}}})|\leq C(1+t)^{-2}[1+\log(1+t)]^{2}
\end{equation}
hence the contribution of this term is also bounded in absolute value by:
\begin{equation}
 C(I^{L}_{4}+I^{\underline{L}}_{4}) \notag
\end{equation}

We collect the above results in the following lemma:

$\textbf{Lemma 7.4}$ We have:
\begin{align}
 |\int_{W^{t}_{u}}(K_{0}\psi_{n})\leftexp{(Y)}{\sigma_{2,n-1}}dt'du'd\mu_{\slashed{g}}|\\  \notag
\leq C\epsilon_{0}\{\int_{0}^{t}(1+t')^{-2}[1+\log(1+t')]\mathcal{E}^{u}_{0,n}(t')dt'+\int_{0}^{u}\mathcal{F}^{t}_{0,n}(u')du'\} \notag
\end{align}

We now consider the error integral
\begin{equation}
 -\int_{W^{t}_{u}}(K_{1}\psi_{n}+\omega\psi_{n})\leftexp{(Y)}{\sigma_{2,n-1}}dt'du'd\mu_{\slashed{g}} \notag
\end{equation}
or
\begin{equation}
 -\int_{W^{t}_{u}}(\omega/\nu)(L\psi_{n}+\nu\psi_{n})\leftexp{(Y)}{\sigma_{2,n-1}}dt'du'd\mu_{\slashed{g}}
\end{equation}
Recalling the bounds for $\omega$ and $\nu$, this is bounded by:
\begin{equation}
 C\int_{W^{t}_{u}}(1+t')^{2}|L\psi_{n}+\nu\psi_{n}||\leftexp{(Y)}{\sigma_{2,n-1}}|dt'du'd\mu_{\slashed{g}}
\end{equation}
We shall estimate this integral.

In view of (7.195), the contribution of the first term in (7.75) to the integral (7.161) is bounded in absolute value by:
\begin{align}
 C\int_{0}^{u}\{\int_{C_{u'}}|L\psi_{n}+\nu\psi_{n}||\underline{L}\psi_{n-1}|[1+\log(1+t')]d\mu_{\slashed{g}}dt'\}du'\\ \notag
\leq C\int_{0}^{u}\{\int_{C_{u'}}(1+t')^{2}(L\psi_{n}+\nu\psi_{n})^{2}d\mu_{\slashed{g}}dt'\}^{1/2}\\ \notag
\cdot\{\int_{C_{u'}}(1+t')^{-2}[1+\log(1+t')]^{2}(\underline{L}\psi_{n-1})^{2}d\mu_{\slashed{g}}dt'\}^{1/2}du'\\ \notag
\leq C\int_{0}^{u}\{\mathcal{F}^{\prime t}_{1,n}(u')\}^{1/2}\{\int_{0}^{t}\epsilon_{0}\mathcal{E}^{u}_{0,n}(t')(1+t')^{-2}[1+\log(1+t')]^{2}dt'\}^{1/2}du'\\ \notag
\leq C\epsilon_{0}\{\int_{0}^{u}\mathcal{F}^{\prime t}_{1,n}(u')du'\}^{1/2}\{\int_{0}^{t}(1+t')^{-2}[1+\log(1+t')]^{2}\mathcal{E}^{u}_{0,n}(t')dt'\}^{1/2}
\end{align}
Here we have used the first of (7.183).

In view of (7.191), the contribution of the second term in (7.75) is bounded in absolute value by:
\begin{align}
 C\int_{0}^{u}\{\int_{C_{u'}}|L\psi_{n}+\nu\psi_{n}||L\psi_{n-1}|(1+t')[1+\log(1+t')]d\mu_{\slashed{g}}dt'\}du'\\ \notag
\leq C\int_{0}^{u}\{\int_{C_{u'}}(1+t')^{2}(L\psi_{n}+\nu\psi_{n})^{2}d\mu_{\slashed{g}}dt'\}^{1/2}\\ \notag
\cdot\{\int_{C_{u'}}[1+\log(1+t')]^{2}(L\psi_{n-1})^{2}d\mu_{\slashed{g}}dt'\}^{1/2}du'\\ \notag
\leq C\int_{0}^{u}\{\mathcal{F}^{\prime t}_{1,n}(u')\}^{1/2}\{\int_{0}^{t}\epsilon_{0}\mathcal{E}^{u}_{0,n}(t')(1+t')^{-2}[1+\log(1+t')]^{2}dt'\}^{1/2}du'\\ \notag
\leq C\epsilon_{0}\{\int_{0}^{u}\mathcal{F}^{\prime t}_{1,n}(u')du'\}^{1/2}\{\int_{0}^{t}(1+t')^{-2}[1+\log(1+t')]^{2}\mathcal{E}^{u}_{0,n}(t')dt'\}^{1/2}
\end{align}
Here we have used the second of (7.183).

In view of (7.196), the contribution of the third term in (7.75) is bounded in a same way as (7.218).

In view of (7.201), the contributions of the fourth and fifth terms in (7.75) is bounded in absolute value by:
\begin{align}
 C\int_{0}^{u}\{\int_{C_{u'}}|L\psi_{n}+\nu\psi_{n}||\slashed{d}\psi_{n-1}|(1+t')[1+\log(1+t')]d\mu_{\slashed{g}}dt'\}du'\\ \notag
\leq C\int_{0}^{u}\{\int_{C_{u'}}(1+t')^{2}(L\psi_{n}+\nu\psi_{n})^{2}d\mu_{\slashed{g}}dt'\}^{1/2}\\ \notag
\cdot\{\int_{C_{u'}}[1+\log(1+t')]^{2}|\slashed{d}\psi_{n-1}|^{2}d\mu_{\slashed{g}}dt'\}^{1/2}du'\\ \notag
\leq C\int_{0}^{u}\{\mathcal{F}^{\prime t}_{1,n}(u')\}^{1/2}\{\int_{0}^{t}\epsilon_{0}\mathcal{E}^{u}_{0,n}(t')(1+t')^{-2}[1+\log(1+t')]^{2}dt'\}^{1/2}du'\\ \notag
\leq C\epsilon_{0}\{\int_{0}^{u}\mathcal{F}^{\prime t}_{1,n}(u')du'\}^{1/2}\{\int_{0}^{t}(1+t')^{-2}[1+\log(1+t')]^{2}\mathcal{E}^{u}_{0,n}(t')dt'\}^{1/2}
\end{align}
Here we have used the last of (7.183).

In view of the bounds (7.208), the contribution of the sixth term in (7.75) is bounded in the same way as (7.217), while 
that of the seventh is bounded in the same way as (7.228).

Finally, in view of (7.209) and (7.213), the contributions of the last two terms of (7.75) is bounded in the same way as (7.219). 

We collect the above results in the following lemma:

$\textbf{Lemma 7.5}$ We have:
\begin{align}
 |\int_{W^{t}_{u}}(K_{1}\psi_{n}+\omega\psi_{n})\leftexp{(Y)}{\sigma_{2,n-1}}dt'du'd\mu_{\slashed{g}}|\\ \notag
\leq C\epsilon_{0}\{\int_{0}^{u}\mathcal{F}^{\prime t}_{1,n}(u')du'\}^{1/2}\{\int_{0}^{t}(1+t')^{-2}[1+\log(1+t')]^{2}\mathcal{E}^{u}_{0,n}(t')dt'\}^{1/2}
\end{align}

Finally, we consider the contribution of $\leftexp{(Y)}{\sigma_{3,n-1}}$ given by (7.77)-(7.79) to (7.99)-(7.100). From $\textbf{E1}, \textbf{E2},
\textbf{F1},\textbf{F2}$ and (6.86), (6.87), (6.89), (6.95), (6.98), (6.99) and (6.104), we have:
\begin{align}
 |L\log\Omega|\leq C(1+t)^{-2}\\\notag
|\textrm{tr}\chi|\leq C(1+t)^{-1}, |\textrm{tr}\underline{\chi}|\leq C(1+t)^{-1}[1+\log(1+t)]\\\notag
|\Lambda|\leq C(1+t)^{-1}[1+\log(1+t)]\\\notag
|L(\eta^{-1}\kappa)|\leq C(1+t)^{-1}, |\slashed{d}(\eta^{-1}\kappa)|\leq C(1+t)^{-1}[1+\log(1+t)]
\end{align}
Together with $\textbf{G0}$ these imply:
\begin{align}
 |\leftexp{(Y)}{\sigma^{\underline{L}}_{3,n-1}}|\leq C(1+t)^{-2}\\ \notag
|\leftexp{(Y)}{\sigma^{L}_{3,n-1}}|\leq (1+t)^{-1}[1+\log(1+t)]\\ \notag
|\leftexp{(Y)}{\slashed{\sigma}_{3,n-1}}|\leq C(1+t)^{-1}[1+\log(1+t)]
\end{align}
Recalling (7.185), (7.191) and (7.201), we conclude that the contribution of $\leftexp{(Y)}{\sigma_{3,n-1}}$ to each of the error integrals
(7.99), (7.100) can be bounded in the same way as the corresponding contribution of $\leftexp{(Y)}{\sigma_{2,n-1}}$.

Combining this with the previous conclusions we arrive at the following lemma:

$\textbf{Lemma 7.6}$  We have:
\begin{align*}
 |\int_{W^{t}_{u}}(K_{0}\psi_{n})\leftexp{(Y)}{\sigma_{n-1}}dt'du'd\mu_{\slashed{g}}|\leq 
C\{\int_{0}^{t}(1+t')^{-3/2}\mathcal{E}^{u}_{0,n}(t')dt'+\int_{0}^{u}\mathcal{F}^{t}_{0,n}(u')du'\}\\
+C\{(\int_{0}^{u}\bar{\mathcal{F}}^{\prime t}_{1,n}(u')du')^{1/2}+(\bar{K_{n}}(t,u))^{1/2}\}\\
\times\{\int_{0}^{t}(1+t')^{-3/2}\mathcal{E}^{u}_{0,n}(t')dt'+\int_{0}^{u}\mathcal{F}^{t}_{0,n}(u')du'\}^{1/2}
\end{align*}
and:
\begin{align*}
 |\int_{W^{t}_{u}}(K_{1}\psi_{n}+\omega\psi_{n})\leftexp{(Y)}{\sigma_{n-1}}dt'du'd\mu_{\slashed{g}}|\\
\leq C\int_{0}^{u}\mathcal{F}^{\prime t}_{1,n}(u')du'+C\{\int_{0}^{u}\mathcal{F}^{\prime}_{1,n}(u')du'\}^{1/2}\{K_{n}(t,u)\}^{1/2}\\
+C\{\int_{0}^{u}\mathcal{F}^{\prime t}_{1,n}(u')du'\}^{1/2}\{\int_{0}^{t}(1+t')^{-2}[1+\log(1+t')]^{2}(\mathcal{E}^{\prime u}_{1,n}(t')
+\epsilon_{0}^{2}\mathcal{E}^{u}_{0,n}(t')dt'\}^{1/2}
\end{align*}

\chapter{Regularization of the Propagation Equation for $\slashed{d}\textrm{tr}\chi$. \\
Estimates for the Top Order Angular Derivatives of $\chi$}
\section{Preliminary}
\subsection{Regularization of The Propagation Equation}
We denote, 
\begin{align}
 \tau_{\mu}=\partial_{\mu}h,\quad  \omega_{\mu\nu}=\partial_{\mu}\psi_{\nu}=\partial_{\nu}\psi_{\mu}
\end{align}
$\textbf{Proposition 8.1}$  The function $h$  satisfies the following inhomogeneous wave equation in the acoustical metric $g$:
\begin{equation}
 \Box_{g}h=-\Omega^{-1}\frac{d\Omega}{dh}a-b
\end{equation}
 where $a$ and $b$ are the functions:
\begin{equation}
 a=(g^{-1})^{\mu\nu}\tau_{\mu}\tau_{\nu}, \quad b=\sum_{i}(g^{-1})^{\mu\nu}\omega_{\mu{i}}\omega_{\nu{i}}
\end{equation}
$Proof$. From Chapter 1, we know that:
\begin{equation}
 \Box_{\tilde{g}}\psi_{\alpha}=0
\end{equation}
It follows that $h$ being given by 
\begin{equation}
 h=\partial_{t}\phi-\frac{1}{2}\sum_{i}(\partial_{i}\phi)^{2}
\end{equation}
satisfies the equation:
\begin{equation}
 \Box_{\tilde{g}}h=-\sum_{i}(\tilde{g}^{-1})^{\mu\nu}\partial_{\mu}\psi_{i}\partial_{\nu}\psi_{i}
\end{equation}
Since 
\begin{equation}
 \tilde{g}_{\mu\nu}=\Omega g_{\mu\nu}
\end{equation}
then for an arbitrary function $f$, we have 
\begin{equation}
 \Box_{\tilde{g}}f=\Omega^{-1}\Box_{g}f+\Omega^{-2}(g^{-1})^{\mu\nu}\partial_{\mu}\Omega\partial_{\nu}f
\end{equation}
and we have
\begin{equation}
 \partial_{\mu}\Omega=\frac{d\Omega}{dh}\partial_{\mu}h
\end{equation}
then the proposition follows. $\qed$

From Chapter 3 $\textrm{tr}\chi$ satisfies the propagation equation
\begin{align}
 L\textrm{tr}\chi=\mu^{-1}(L\mu)\textrm{tr}\chi-|\chi|^{2}-\textrm{tr}\alpha
\end{align}
Since 
\begin{align}
 S_{\mu\nu}=(g^{-1})^{\kappa\lambda}R_{\kappa\mu\lambda\nu}
\end{align}
is the Ricci tensor, we have:
\begin{align}
 \textrm{tr}\alpha=S(L,L)
\end{align}
We decompose:
\begin{align}
 S(L,L)=S^{[P]}(L,L)+S^{[N]}(L,L)=\textrm{tr}\alpha^{[P]}+\textrm{tr}\alpha^{[N]}
\end{align}
From Chapter 4,
\begin{align}
 \textrm{tr}\alpha^{[P]}=-\frac{1}{2}\frac{dH}{dh}(\slashed{g}^{-1})^{AB}v_{AB}
\end{align}
Since
\begin{align}
 -\mu^{-1}v_{L\underline{L}}+(\slashed{g}^{-1})^{AB}v_{AB}=\textrm{tr}v=\Box_{g}h
\end{align}
we have:
\begin{align}
 \textrm{tr}\alpha^{[P]}=-\frac{1}{2}\frac{dH}{dh}(\mu^{-1}v_{L\underline{L}}+\Box_{g}h)
\end{align}
Also from Chapter 4,
\begin{align}
 \textrm{tr}\alpha^{[N]}=\eta^{-2}\textrm{tr}\alpha^{[A]}+\frac{1}{2}H_{1}\textrm{tr}\alpha^{[C]}-\frac{1}{2}H_{2}\textrm{tr}\alpha^{[B]}
\end{align}
where:
\begin{align}
 H_{1}=\eta^{-2}\frac{dH}{dh},\quad H_{2}=\frac{d^{2}H}{dh^{2}}+\frac{1}{2}\eta^{-2}(\frac{dH}{dh})^{2}
\end{align}
and:
\begin{align}
 \textrm{tr}\alpha^{[A]}=w_{LL}\textrm{tr}\slashed{w}-|\slashed{w}_{L}|^{2}\\
\textrm{tr}\alpha^{[B]}=|\slashed{\tau}|^{2}\\
\textrm{tr}\alpha^{[C]}=2(\tau_{L}\textrm{tr}\slashed{w}-\slashed{\tau}\cdot\slashed{w}_{L})
\end{align}
Now,
\begin{align}
 v_{L\underline{L}}=D^{2}h(L,\underline{L})=L(\underline{L}h)-\tau(D_{L}\underline{L})
\end{align}
and from Chapter 3:
\begin{align*}
 D_{L}\underline{L}=-2\zeta^{A}X_{A}
\end{align*}
Hence:
\begin{align}
 v_{L\underline{L}}=L(\underline{L}h)+2\zeta\cdot\slashed{\tau}
\end{align}
Let us define, noting that $\underline{L}h=\tau_{\underline{L}}$,
\begin{align}
 f=-\frac{1}{2\mu}\frac{dH}{dh}\tau_{\underline{L}}
\end{align}
We then have:
\begin{align}
 S(L,L)=Lf+g
\end{align}
where:
\begin{align}
 g=\frac{1}{2\mu}\frac{d^{2}H}{dh^{2}}\tau_{L}\tau_{\underline{L}}-\frac{1}{2\mu^{2}}(L\mu)\frac{dH}{dh}\tau_{\underline{L}}
-\frac{1}{\mu}\frac{dH}{dh}\zeta\cdot\slashed{\tau}-\frac{1}{2}\frac{dH}{dh}\Box_{g}h+\textrm{tr}\alpha^{[N]}
\end{align}
We now set:
\begin{align}
 \check{f}=\mu f=-\frac{1}{2}\frac{dH}{dh}\tau_{\underline{L}}
\end{align}
which is regular as $\mu\rightarrow0$. We then have:
\begin{align}
 \mu S(L,L)=L\check{f}+\check{g}
\end{align}
where
\begin{align}
 \check{g}=\mu g-(L\mu)f\\\notag
=\frac{1}{2}\frac{d^{2}H}{dh^{2}}\tau_{L}\tau_{\underline{L}}-\frac{dH}{dh}\zeta\cdot\slashed{\tau}
-\frac{1}{2}\frac{dH}{dh}\mu\Box_{g}h+\mu\textrm{tr}\alpha^{[N]}
\end{align}
By Proposition 8.1:
\begin{align}
 \mu\Box_{g}h=-\Omega^{-1}\frac{d\Omega}{dh}\mu a-\mu b
\end{align}
where:
\begin{align}
 \mu a=\mu(g^{-1})^{\mu\nu}\tau_{\mu}\tau_{\nu}=-\tau_{L}\tau_{\underline{L}}+\mu|\slashed{\tau}|^{2}
\end{align}
and:
\begin{align}
 \mu b=\sum_{i}\mu(g^{-1})^{\mu\nu}\partial_{\mu}\psi_{i}\partial_{\nu}\psi_{i}
=\sum_{i}\{-(L\psi_{i})(\underline{L}\psi_{i})+\mu|\slashed{d}\psi_{i}|^{2}\}
\end{align}
We see that $\mu a$ and $\mu b$ are both regular as $\mu\rightarrow0$, so is $\mu\Box_{g}h$. It follows that $\check{g}$ is regular as $\mu\rightarrow0$.
(8.10) takes the form:
\begin{equation}
 L(\mu\textrm{tr}\chi+\check{f})=2(L\mu)\textrm{tr}\chi-\frac{1}{2}\mu(\textrm{tr}\chi)^{2}-\mu|\hat{\chi}|^{2}-\check{g}
\end{equation}
Introduce the $S_{t,u}$ 1-form:
\begin{equation}
 x_{0}=\mu\slashed{d}\textrm{tr}\chi+\slashed{d}\check{f}
\end{equation}
We shall obtain a propagation equation for $x_{0}$, by differentiaing (8.33) tangentially to $S_{t,u}$.

$\textbf{Lemma 8.1}$ For an arbitrary function $\phi$, we have:
\begin{equation}
 \slashed{\mathcal{L}}_{L}(\slashed{d}\phi)=\slashed{d}(L\phi)
\end{equation}
$Proof$. Since both sides of the equation are $S_{t,u}$ 1-forms, we need only evaluate both sides on $X_{A}$.
\begin{align*}
\slashed{\mathcal{L}}_{L}(\slashed{d}\phi)(X_{A})=\mathcal{L}_{L}(\slashed{d}\phi)(X_{A})=
L(X_{A}(\phi))-(\slashed{d}\phi)(\mathcal{L}_{L}X_{A})\\
=L(X_{A}(\phi))-[L,X_{A}]\phi=X_{A}(L\phi)
\end{align*}
While $\slashed{d}(L\phi)(X_{A})=X_{A}(L\phi)$, so the lemma follows. $\qed$

Taking $\phi=\mu\textrm{tr}\chi+\check{f}$ in Lemma 8.1, using (8.33), (8.34), we get the propagation equation for $x_{0}$:
\begin{equation}
 \slashed{\mathcal{L}}_{L}x_{0}+(\textrm{tr}\chi-2\mu^{-1}(L\mu))x_{0}=(\frac{1}{2}\textrm{tr}\chi-2\mu^{-1}(L\mu))\slashed{d}\check{f}
-\mu\slashed{d}(|\hat{\chi}|^{2})-g_{0}
\end{equation}
where:
\begin{equation}
 g_{0}=\slashed{d}\check{g}-\frac{1}{2}\textrm{tr}\chi\slashed{d}(\check{f}+2L\mu)+(\slashed{d}\mu)(L\textrm{tr}\chi+|\chi|^{2})
\end{equation}
We may substitute (3.24)
\begin{equation}
 L\textrm{tr}\chi+|\chi|^{2}=\mu^{-1}(L\mu)\textrm{tr}\chi-\textrm{tr}\alpha
\end{equation}
into the expression of $g_{0}$. From Chapter 4,
\begin{align*}
 \textrm{tr}\alpha=\textrm{tr}\alpha^{[P]}+\textrm{tr}\alpha^{[N]}
\end{align*}
where $\textrm{tr}\alpha^{[N]}$ is regular as $\mu\rightarrow0$ while
\begin{align}
 \textrm{tr}\alpha^{[P]}=-\frac{1}{2}\frac{dH}{dh}(\slashed{g}^{-1})^{AB}v_{AB}
\end{align}
where
\begin{align*}
 v_{AB}=D^{2}h(X_{A},X_{B})=\slashed{D}^{2}h(X_{A},X_{B})-\eta^{-1}\slashed{k}_{AB}Lh-\mu^{-1}\chi_{AB}Th
\end{align*}
hence 
\begin{align}
 (\slashed{g}^{-1})^{AB}v_{AB}=\slashed{\Delta}h-\eta^{-1}\textrm{tr}\slashed{k}Lh-\mu^{-1}\textrm{tr}\chi Th
\end{align}
and:
\begin{align}
 \textrm{tr}\alpha^{[P]}=\mu^{-1}m\textrm{tr}\chi-\frac{1}{2}\frac{dH}{dh}(\slashed{\Delta}h-\eta^{-1}\textrm{tr}\slashed{k}Lh)
\end{align}
In view of the propagation equation for $\mu$
\begin{equation}
 L\mu=m+\mu{e}
\end{equation}
the right hand side of (8.38) is 
\begin{align}
 e\textrm{tr}\chi+\frac{1}{2}\frac{dH}{dh}(\slashed{\Delta}h-\eta^{-1}\textrm{tr}\slashed{k}Lh)-\textrm{tr}\alpha^{[N]}
\end{align}
which is regular as $\mu\longrightarrow0$.

Thus $g_{0}$ is of the order of the second derivatives of $\psi_{\mu}$ and bounded as $\mu\longrightarrow0$.

\subsection{Propagation Equations for Higher Order Angular Derivatives}
To control the higher order angular derivatives of $\textrm{tr}\chi$, we introduce the $S_{t,u}$ 1-forms:
\begin{equation}
 \leftexp{(i_{1}...i_{l})}{x_{l}}=\mu\slashed{d}(R_{i_{l}}...R_{i_{1}}\textrm{tr}\chi)+\slashed{d}(R_{i_{l}}...R_{i_{1}}\check{f})
\end{equation}
Thus, for a given positive integer ${l}$ we have the mult-indices $(i_{1}...i_{l})$ of length ${l}$, where each $i_{k} \in \{1,2,3\}$, 
$k=1,..,l$. To the mult-indices $(i_{1}...i_{l})$ there corresponds the string $(R_{i_{l}}...R_{i_{1}})$ of rotation vectorfields. In deriving propagation 
equations for $x_{l}$, we shall use the following lemma.

$\textbf{Lemma 8.2}$ Let $Y$ be an arbitrary $S_{t,u}$-tangential vectorfield in spacetime manifold. We have
\begin{equation}
 [L,Y]=\leftexp{(Y)}{Z}
\end{equation}
where $\leftexp{(Y)}{Z}$ is the $S_{t,u}$-tangential vectorfields, associated to $Y$, 
and defined by the condition that for any vector $V \in TW^{*}_{\epsilon_{0}}$:
\begin{equation}
g(\leftexp{(Y)}{Z},V)=\leftexp{(Y)}{\pi}(L,\Pi V)
\end{equation}
 In terms of $(L,T,X_{1},X_{2})$ or $(L,\underline{L},X_{1},X_{2})$,
\begin{equation}
 \leftexp{(Y)}{Z}=\leftexp{(Y)}{Z}^{A}X_{A},\quad \leftexp{(Y)}{Z}^{A}=\leftexp{(Y)}{\pi}_{LB}(\slashed{g}^{-1})^{AB}
\end{equation}
$Proof$. The vectorfield $[L,Y]$ satisfies
\begin{align*}
 [L,Y]t=L(Yt)-Y(Lt)=0\\
[L,Y]u=L(Yu)-Y(Lu)=0
\end{align*}
It follows that $[L,Y]=\leftexp{(Y)}{Z}$ is a vectorfield which is tangential to the surfaces $S_{t,u}$,
consequently we expand
\begin{equation}
 \leftexp{(Y)}{Z}=\leftexp{(Y)}{Z}^{B}X_{B}
\end{equation}
Taking $g$-inner product with $X_{A}$ we obtain
\begin{align}
 \slashed{g}_{AB}\leftexp{(Y)}{Z}^{B}=\leftexp{(Y)}{Z}_{A}=g(\leftexp{(Y)}{Z},X_{A})\\\notag
=g([L,Y], X_{A})=g(D_{L}Y,X_{A})-g(D_{Y}L,X_{A})\\\notag
=\leftexp{(Y)}{\pi_{LA}}-g(D_{X_{A}}Y,L)-g(D_{Y}L,X_{A})
\end{align}
Now since $g(L,X_{A})=0$ we have
\begin{equation}
 -g(D_{Y}L,X_{A})=g(L,D_{Y}X_{A})
\end{equation}
hence, substituting,
\begin{equation}
 \slashed{g}_{AB}Z^{B}=\leftexp{(Y)}{\pi_{LA}}+g(L,[Y,X_{A}])=\leftexp{(Y)}{\pi_{LA}}
\end{equation}
We have used the fact that $[Y,X_{A}]$ is tangential to $S_{t,u}$, hence orthogonal to $L$. The lemma is proved. $\qed$

$\textbf{Lemma 8.3}$ Let $Y$ be an arbitrary $S_{t,u}$-tangential vectorfield on the spacetime domain $W^{*}_{\epsilon_{0}}$ and 
let $\xi$ be an arbitrary $S_{t,u}$ 1-form on $W^{*}_{\epsilon_{0}}$. We have
\begin{equation}
 \slashed{\mathcal{L}}_{L} \slashed{\mathcal{L}}_{Y}\xi-\slashed{\mathcal{L}}_{Y} \slashed{\mathcal{L}}_{L}\xi=\slashed{\mathcal{L}}_{\leftexp{(Y)}{Z}}\xi
\end{equation}
$Proof$. Since both $L,Y$ are tangential to the hypersurfaces $C_{u}$ we can restrict ourselves to a given $C_{u}$. In defining
 $\slashed{\mathcal{L}}_{Y}\xi,\slashed{\mathcal{L}}_{L}\xi$ we are considering the extension of $\xi$ to $TC_{u}$ by the condition
\begin{equation}
 \xi(L)=0
\end{equation}
We have 
\begin{equation}
 (\mathcal{L}_{L}\xi)(L)=L(\xi(L))-\xi([L,L])=0
\end{equation}
therefore
\begin{equation}
 \slashed{\mathcal{L}}_{L}\xi=\mathcal{L}_{L}\xi
\end{equation}
However
\begin{equation}
 (\mathcal{L}_{Y}\xi)(L)=Y(\xi(L))-\xi([Y,L])=\xi([L,Y])=\xi(\leftexp{(Y)}{Z})
\end{equation}
Since $\slashed{\mathcal{L}}_{Y}\xi$ is defined by restricting $\mathcal{L}_{Y}\xi$ to $TS_{t,u}$ and then extending to $TC_{u}$ by the condition
$(\slashed{\mathcal{L}}_{Y}\xi)(L)=0$, it follows that on the manifold $C_{u}$,
\begin{equation}
 \slashed{\mathcal{L}}_{Y}\xi=\mathcal{L}_{Y}\xi-\xi(\leftexp{(Y)}{Z})dt
\end{equation}
in view of the fact that $dt(L)=Lt=1, \slashed{d}t=0$.

We have
\begin{align}
 (\slashed{\mathcal{L}}_{L} \slashed{\mathcal{L}}_{Y}\xi-\slashed{\mathcal{L}}_{Y}\slashed{\mathcal{L}}_{L}\xi)(X_{A})=
(\mathcal{L}_{L} \slashed{\mathcal{L}}_{Y}\xi-\mathcal{L}_{Y} \slashed{\mathcal{L}}_{L}\xi)(X_{A})
\end{align}
Substituting (8.55) and (8.57), the right hand side of the above is
\begin{equation}
 (\mathcal{L}_{L}\mathcal{L}_{Y}\xi-L(\xi(\leftexp{(Y)}{Z}))dt
-\xi(\leftexp{(Y)}{Z})\mathcal{L}_{L}(dt)-\mathcal{L}_{Y}\mathcal{L}_{L}\xi)(X_{A})
\end{equation}
Now, $dt(X_{A})=0$ and $\mathcal{L}_{L}(dt)=d(Lt)=0$, therefore (8.58) becomes
\begin{equation}
 (\mathcal{L}_{L}\mathcal{L}_{Y}\xi-\mathcal{L}_{Y}\mathcal{L}_{L}\xi)(X_{A})
=(\mathcal{L}_{[L,Y]}\xi)(X_{A})
\end{equation}
and by lemma 8.2, this is $(\mathcal{L}_{\leftexp{(Y)}{Z}}\xi)(X_{A})$. The lemma follows. $\qed$

In the sequel, the proofs of several propositions make use of a general elementary proposition on linear recursions, the proof of which is by a simple 
application of the principle of induction.

$\textbf{Proposition 8.2}$ Let $(y_{n}: n=1,2,...)$ be a given sequence in a space $X$ and $(A_{n}: n=1,2,...)$ a given sequence of operators in $X$. Suppose that 
$(x_{n}: n=0,1,2,...)$ is a sequence in $X$ satisfying the recursion:
\begin{equation}
 x_{n}=A_{n}x_{n-1}+y_{n}
\end{equation}
 Then for each $n=1,2,...$ we have:
\begin{equation}
 x_{n}=A_{n}...A_{1}x_{0}+\sum^{n-1}_{m=0}A_{n}...A_{n-m+1}y_{n-m}
\end{equation}

To present the propagation equation for $\leftexp{(i_{1}...i_{l})}{x_{l}}$, we introduce the functions:
\begin{align}
 \leftexp{(i_{1}..i_{l})}{\check{f}_{l}}=R_{i_{l}}...R_{i_{1}}\check{f}\\
\leftexp{(i_{1}...i_{l})}{h_{l}}=R_{i_{l}}...R_{i_{1}}(|\hat{\chi}|^{2})
\end{align}

$\textbf{Proposition 8.3}$ For each nonnegative integer $l$ and each multi-index $(i_{1}...i_{l})$, the $S_{t,u}$ 1-form 
$\leftexp{(i_{1}...i_{l})}{x_{l}}$ satisfies the propagation equation:
\begin{align*}
 \slashed{\mathcal{L}}_{L}\leftexp{(i_{1}...i_{l})}{x_{l}}+(\textrm{tr}\chi-2\mu^{-1}(L\mu))\leftexp{(i_{1}...i_{l})}{x_{l}}
=(\frac{1}{2}\textrm{tr}\chi-2\mu^{-1}(L\mu))\slashed{d}\leftexp{(i_{1}..i_{l})}{\check{f}_{l}}-\mu\slashed{d}\leftexp{(i_{1}...i_{l})}{h_{l}}-
\leftexp{(i_{1}...i_{l})}{g_{l}}
\end{align*}
where $\leftexp{(i_{1}...i_{l})}{g_{l}}$ is the $S_{t,u}$ 1-form given by:
\begin{align*}
 \leftexp{(i_{1}...i_{l})}{g_{l}}=\slashed{\mathcal{L}}_{R_{i_{l}}}...\slashed{\mathcal{L}}_{R_{i_{1}}}g_{0}-
\sum_{k=0}^{l-1}\slashed{\mathcal{L}}_{R_{i_{l}}}...\slashed{\mathcal{L}}_{R_{i_{l-k+1}}}\slashed{\mathcal{L}}_{\leftexp{(R_{i_{l-k}})}{Z}}
\leftexp{(i_{1}...i_{l-k-1})}{x_{l-k-1}}\\
+\sum_{k=0}^{l-1}\slashed{\mathcal{L}}_{R_{i_{l}}}...\slashed{\mathcal{L}}_{R_{i_{l-k+1}}}\leftexp{(i_{1}..i_{l-k})}{y_{l-k}}
\end{align*}
where $g_{0}$ is given by (8.37). Here for each $j=1...l$, $\leftexp{(i_{1}...i_{j})}{y_{j}}$ is the $S_{t,u}$ 1-form:
\begin{align*}
 \leftexp{(i_{1}...i_{j})}{y_{j}}=(R_{i_{j}}\mu)\leftexp{(i_{1}...i_{j-1})}{a_{j-1}}\\
+(\mu{R_{i_{j}}}\textrm{tr}\chi-R_{i_{j}}L\mu+\leftexp{(R_{i_{j}})}{Z\mu})\slashed{d}(R_{i_{j-1}}...R_{i_{1}}\textrm{tr}\chi)\\
+\frac{1}{2}(R_{i_{j}}\textrm{tr}\chi)\slashed{d}\leftexp{(i_{1}...i_{j-1})}{\check{f}_{j-1}}
\end{align*}
where $\leftexp{(i_{1}...i_{j-1})}{a_{j-1}}$ is the $S_{t,u}$ 1-form:
\begin{equation}
 \leftexp{(i_{1}...i_{j-1})}{a_{j-1}}=\slashed{\mathcal{L}}_{L}\slashed{d}(R_{{i}_{j-1}}...R_{i_{1}}\textrm{tr}\chi)
+\textrm{tr}\chi\slashed{d}(R_{{i}_{j-1}}...R_{i_{1}}\textrm{tr}\chi)+\slashed{d}\leftexp{(i_{1}...i_{j-1})}{h_{j-1}} \notag
\end{equation}
$Proof$. When $l=0$, the propagation equation in the proposition is just the equation (8.36). Thus, by induction, assuming that the propagation equation 
holds with $l$ replaced by $l-1$, that is,
\begin{align}
 \slashed{\mathcal{L}}_{L}\leftexp{(i_{1}...i_{l-1})}{x_{l-1}}+(\textrm{tr}\chi-2\mu^{-1}(L\mu))\leftexp{(i_{1}...i_{l-1})}{x_{l-1}}\\\notag
=(\frac{1}{2}\textrm{tr}\chi-2\mu^{-1}(L\mu))\slashed{d}\leftexp{(i_{1}..i_{l-1})}{\check{f}_{l-1}}-\mu\slashed{d}\leftexp{(i_{1}...i_{l-1})}{h_{l-1}}-
\leftexp{(i_{1}...i_{l-1})}{g_{l-1}}
\end{align}
holds for some $S_{t,u}$ 1-form $\leftexp{(i_{1}...i_{l-1})}{g_{l-1}}$, what we shall show is that a propagation equation of the form given by the proposition 
holds true for $l$, where $\leftexp{(i_{1}...i_{l})}{g_{l}}$ is an $S_{t,u}$ 1-form related to $\leftexp{(i_{1}...i_{l-1})}{g_{l-1}}$ by a certain 
recursion relation. This recursion relation shall then determine $\leftexp{(i_{1}...i_{l})}{g_{l}}$, for each $l$, from the $S_{t,u}$ 1-form $g_{0}$,
given by (8.37).

We begin by re-writing 
\begin{align*}
 2\mu^{-1}(L\mu)(\leftexp{(i_{1}...i_{l-1})}{x_{l-1}}-\slashed{d}\leftexp{(i_{1}...i_{l-1})}{\check{f}}_{l-1})
\end{align*}
as 
\begin{align*}
 2(L\mu)\slashed{d}(R_{i_{l-1}}...R_{i_{1}}\textrm{tr}\chi)
\end{align*}
(see (8.44) and (8.63)) obtaining:
\begin{align}
 \slashed{\mathcal{L}}_{L}\leftexp{(i_{1}...i_{l-1})}{x_{l-1}}+\textrm{tr}\chi\leftexp{(i_{1}...i_{l-1})}{x_{l-1}}\\ \notag
=2(L\mu)\slashed{d}(R_{i_{l-1}}...R_{i_{1}}\textrm{tr}\chi)+\frac{1}{2}\textrm{tr}\chi\slashed{d}\leftexp{(i_{1}...i_{l-1})}{\check{f}}_{l-1}
-\mu\slashed{d}\leftexp{(i_{1}...i_{l-1})}{h}_{l-1}-\leftexp{(i_{1}...i_{l-1})}{g}_{l-1}
\end{align}
We now apply $\slashed{\mathcal{L}}_{R_{i_{l}}}$ to this equation. By Lemma 8.1 applied to $R_{i_{l}}$ and to the functions:
\begin{align*}
 R_{i_{l-1}}...R_{i_{1}}\textrm{tr}\chi,\quad \leftexp{(i_{1}...i_{l-1})}{\check{f}}_{l-1}, \quad \leftexp{(i_{1}...i_{l-1})}{h}_{l-1}
\end{align*}
we have:
\begin{align}
 \slashed{\mathcal{L}}_{R_{i_{l}}}\slashed{d}(R_{i_{l-1}}...R_{i_{1}}\textrm{tr}\chi)=\slashed{d}(R_{i_{l}}...R_{i_{1}}\textrm{tr}\chi)\\\notag
\slashed{\mathcal{L}}_{R_{i_{l}}}\slashed{d}(\leftexp{(i_{1}...i_{l-1})}{\check{f}}_{l-1})=\slashed{d}(\leftexp{(i_{1}...i_{l})}{\check{f}}_{l})\\\notag
\slashed{\mathcal{L}}_{R_{i_{l}}}\slashed{d}(\leftexp{(i_{1}...i_{l-1})}{h}_{l-1})=\slashed{d}(\leftexp{(i_{1}...i_{l})}{h}_{l})
\end{align}
(see (8.63), (8.64)). We then obtain:
\begin{align}
 \slashed{\mathcal{L}}_{R_{i_{l}}}\slashed{\mathcal{L}}_{L}\leftexp{(i_{1}...i_{l-1})}{x_{l-1}} 
+\textrm{tr}\chi\slashed{\mathcal{L}}_{R_{i_{l}}}\leftexp{(i_{1}...i_{l-1})}{x_{l-1}}+
(R_{i_{l}}\textrm{tr}\chi)\leftexp{(i_{1}...i_{l-1})}{x_{l-1}}\\ \notag
=2(L\mu)\slashed{d}(R_{i_{l}}...R_{i_{1}}\textrm{tr}\chi)+2(R_{i_{l}}L\mu)\slashed{d}(R_{i_{l-1}}...R_{i_{1}}\textrm{tr}\chi)
+\frac{1}{2}\textrm{tr}\chi\slashed{d}(\leftexp{(i_{1}...i_{l})}{\check{f}}_{l})\\ \notag+\frac{1}{2}(R_{i_{l}}\textrm{tr}\chi)
\slashed{d}(\leftexp{(i_{1}...i_{l-1})}{\check{f}}_{l-1})
-\mu\slashed{d}(\leftexp{(i_{1}...i_{l})}{h}_{l})-(R_{i_{l}}\mu)\slashed{d}(\leftexp{(i_{1}...i_{l-1})}{h}_{l-1})
-\slashed{\mathcal{L}}_{R_{i_{l}}}\leftexp{(i_{1}...i_{l-1})}{g}_{l-1} \notag
\end{align}
Next, we apply Lemma 8.3 to $Y=R_{i_{l}}, \xi=\leftexp{(i_{1}...i_{l-1})}{x_{l-1}}$, to express
\begin{equation}
 \slashed{\mathcal{L}}_{R_{i_{l}}}\slashed{\mathcal{L}}_{L}\leftexp{(i_{1}...i_{l-1})}{x_{l-1}}=
\slashed{\mathcal{L}}_{L}\slashed{\mathcal{L}}_{R_{i_{l}}}\leftexp{(i_{1}...i_{l-1})}{x_{l-1}}-
\slashed{\mathcal{L}}_{\leftexp{(R_{i_{l}})}{Z}}\leftexp{(i_{1}...i_{l-1})}{x_{l-1}}
\end{equation}
Now, from (8.63) we have:
\begin{align*}
 \leftexp{(i_{1}...i_{l-1})}{x_{l-1}}=\mu\slashed{d}(R_{i_{l-1}}...R_{i_{1}}\textrm{tr}\chi)+\slashed{d}(R_{i_{l-1}}...R_{i_{1}}\check{f})
\end{align*}
Applying $\slashed{\mathcal{L}}_{R_{i_{l}}}$ to this we obtain, using (8.63) and Lemma 8.1,
\begin{align}
 \slashed{\mathcal{L}}_{R_{i_{l}}}\leftexp{(i_{1}...i_{l-1})}{x_{l-1}}=\leftexp{(i_{1}...i_{l})}{x_{l}}
+(R_{i_{l}}\mu)\slashed{d}(R_{i_{l-1}}...R_{i_{1}}\textrm{tr}\chi)
\end{align}
Applying $\slashed{\mathcal{L}}_{L}$ to (8.70) and expressing:
\begin{align}
 LR_{i_{l}}\mu=R_{i_{l}}L\mu+\leftexp{(R_{i_{l}})}{Z}\mu
\end{align}
using the fact that by Lemma 8.2
\begin{equation}
 [L,R_{i_{l}}]=\leftexp{(R_{i_{l}})}{Z}
\end{equation}
yields:
\begin{align}
 \slashed{\mathcal{L}}_{L}\slashed{\mathcal{L}}_{R_{i_{l}}}\leftexp{(i_{1}...i_{l-1})}{x_{l-1}}=
\slashed{\mathcal{L}}_{L}\leftexp{(i_{1}...i_{l})}{x_{l}}+(R_{i_{l}}\mu)\slashed{\mathcal{L}}_{L}\slashed{d}
(R_{i_{l-1}}...R_{i_{1}}\textrm{tr}\chi)\\
+(R_{i_{l}}L\mu+\leftexp{(R_{i_{l}})}{Z}\mu)\slashed{d}(R_{i_{l-1}}...R_{i_{1}}\textrm{tr}\chi)\notag
\end{align}
Substituting (8.72) in (8.69) and the result in (8.68), and using (8.70) to re-write the second term on the left in (8.68), a propagation equation for 
$\leftexp{(i_{1}...i_{l})}{x}_{l}$ of the form given by the proposition results with $\leftexp{(i_{1}...i_{l})}{g}_{l}$ expressed in terms of
$\leftexp{(i_{1}...i_{l})}{g}_{l-1}$ by the recursion formula:
\begin{equation}
 \leftexp{(i_{1}...i_{l})}{g}_{l}=\slashed{\mathcal{L}}_{R_{i_{l}}}\leftexp{(i_{1}...i_{l-1})}{g}_{l-1}
-\slashed{\mathcal{L}}_{\leftexp{(R_{i_{l}})}{Z}}\leftexp{(i_{1}...i_{l-1})}{x}_{l-1}+\leftexp{(i_{1}...i_{l})}{y}_{l}
\end{equation}
Applying Proposition 8.2 to this recursion with the space of $S_{t,u}$ 1-forms in the role of the space $X$, and $\leftexp{(i_{1}...i_{l})}{g}_{l},
\leftexp{(i_{1}...i_{l})}{y}_{l}-\slashed{\mathcal{L}}_{\leftexp{(R_{i_{l}})}{Z}}\leftexp{(i_{1}...i_{l-1})}{x}_{l-1}$, in the role of $x_{n}, y_{n}$,
respectively, the proposition follows. $\qed$

We remark that the $S_{t,u}$ 1-form $\leftexp{(i_{1}...i_{j-1})}{a}_{j-1}$ reduces for $j=1$ to
\begin{equation}
 a_{0}=\slashed{\mathcal{L}}_{L}\slashed{d}\textrm{tr}\chi+\textrm{tr}\chi\slashed{d}\textrm{tr}\chi+\slashed{d}h_{0}
\end{equation}
Define
\begin{equation}
 L\textrm{tr}\chi+|\chi|^{2}=f_{0}
\end{equation}
Applying $\slashed{d}$ to the above, we obtain, in view of Lemma 8.1 and the fact that 
\begin{align*}
 |\chi|^{2}=\frac{1}{2}(\textrm{tr}\chi)^{2}+h_{0},
\end{align*}
simply:
\begin{equation}
 a_{0}=\slashed{d}f_{0}
\end{equation}
For $j\geq 2$, we apply $\slashed{\mathcal{L}}_{R_{i_{j-1}}}...\slashed{\mathcal{L}}_{R_{i_{1}}}$ to (8.75) to obtain, in view of (8.76), the following 
expression for $\leftexp{(i_{1}...i_{j-1})}{a}_{j-1}$:
\begin{equation}
 \leftexp{(i_{1}...i_{j-1})}{a}_{j-1}=\slashed{d}(R_{i_{j-1}}...R_{i_{1}}f_{0})+\leftexp{(i_{1}...i_{j-1})}{b}_{j-1}
\end{equation}
where
\begin{align}
 \leftexp{(i_{1}...i_{j-1})}{b}_{j-1}=[\slashed{\mathcal{L}}_{L}, \slashed{\mathcal{L}}_{R_{i_{j-1}}}...\slashed{\mathcal{L}}_{R_{i_{1}}}]
\slashed{d}\textrm{tr}\chi\\ \notag
+\sum_{m=1}^{j-1}\sum_{k_{1}<...<k_{m}=1}^{j-1}(R_{i_{k_{m}}}...R_{i_{k_{1}}}\textrm{tr}\chi)
(\slashed{d}R_{i_{j-1}}\overset{>i_{k_{m}}...i_{k_{1}}<}{...}R_{i_{1}}\textrm{tr}\chi)
\end{align}
Here and in the sequel the symbol $>\quad<$ signifies that the enclosed indices are absent. 
We can get an explicit expression for the first term on the right hand side of (8.79) by applying 
the following lemma, proved inductively, using Lemma 8.3.

$\textbf{Lemma 8.4}$  For any positive integer $l$ and any $S_{t,u}$ 1-form $\xi$ we have:
\begin{align*}
 [\slashed{\mathcal{L}}_{L}, \slashed{\mathcal{L}}_{R_{i_{l}}}...\slashed{\mathcal{L}}_{R_{i_{1}}}]\xi=
\sum_{k=0}^{l-1}\slashed{\mathcal{L}}_{R_{i_{l}}}...\slashed{\mathcal{L}}_{R_{i_{l-k+1}}}\slashed{\mathcal{L}}_{\leftexp{(R_{i_{l-k}})}{Z}}
\slashed{\mathcal{L}}_{R_{i_{l-k+1}}}...\slashed{\mathcal{L}}_{R_{i_{1}}}\xi
\end{align*}
Applying this to $S_{t,u}$ 1-form $\slashed{d}\textrm{tr}\chi$, we get:
\begin{align}
 [\slashed{\mathcal{L}}_{L}, \slashed{\mathcal{L}}_{R_{i_{j-1}}}...\slashed{\mathcal{L}}_{R_{i_{1}}}]\slashed{d}\textrm{tr}\chi
=\sum_{m=0}^{j-2}\slashed{d}(R_{i_{j-1}}...R_{i_{j-m}}\leftexp{(R_{i_{j-m-1}})}{Z}R_{i_{j-m-2}}...R_{i_{1}}\textrm{tr}\chi)
\end{align}

Let us now investigate the order of the various terms in the propagation equation of Proposition 8.3. We agree to the convention that the order of 
$\psi_{\alpha}: \alpha=0,1,2,3$ is $0$, so the order of enthalpy $h$ and the acoustical metric $g, \bar{g}, \slashed{g}$ is also $0$.
Since $\mu$, $\kappa$ are the components of the acoustical metric $g$ in another coordinate system, they are also of order $0$ . So $L^{\mu}$, $T^{i}$ and 
$\alpha$ are also of order $0$, as well as $\psi_{L}$, $\psi_{\hat{T}}$ and $\slashed{\psi}$. The $\slashed{\omega}$, $\slashed{\omega}_{L}$, $\omega_{LL}$ 
and $\omega_{L\hat{T}}$ are of order $1$, so are the $\kappa^{-1}\zeta$ and $\slashed{k}$. The $\chi$, the focus of this chapter is also of order $1$. 
We are to set $l=n-1$, in estimating the $n$th order angular derivatives of $\chi$ in terms of the $n+1$st order derivatives of $\psi_{\alpha}$.
From the definition of $\check{f}$ (see (8.27)), we see that $\check{f}$ is of order $1$. So $\slashed{d}\leftexp{(i_{1}...i_{l})}{\check{f}}_{l}$ is of order $l+2$,
its principal part contains $l+2$nd derivatives of $h$, and it contains $l+1$st order angular derivatives of $\mu$. 

Next, we shall investigate the term $\leftexp{(i_{1}...i_{l})}{g}_{l}$. First, we must investigate the terms $\slashed{d}\check{g}$, 
$L\textrm{tr}\chi+|\chi|^{2}$, together with $\slashed{d}\check{f}$,  which express $g_{0}$.
From (8.38), we see that $L\textrm{tr}\chi+|\chi|^{2}$ is of order $2$. Its principal term is $\frac{1}{2}\frac{dH}{dh}\slashed{\Delta}h$. The function $\mu$
is not involved in $L\textrm{tr}\chi+|\chi|^{2}$.
We then investigate $\slashed{d}\check{g}$. From the definition (8.29) we know $\slashed{d}\check{g}$ contains $1$st angular derivatives of $\mu$, and 
$2$nd derivatives of the $\psi_{\alpha}$.
It follows that 
\begin{align*}
 \slashed{\mathcal{L}}_{R_{i_{l}}}...\slashed{\mathcal{L}}_{R_{i_{1}}}g_{0}
\end{align*}
the first term in the expression for $\leftexp{(i_{1}...i_{l})}{g}_{l}$, is of principal order $l+2$ and its principal terms are $l+2$ order derivatives of the 
$\psi_{\alpha}$. Also it contains $l+1$ order angular derivatives of $\mu$.

We now consider the last term in the expression for $\leftexp{(i_{1}...i_{l})}{g}_{l}$, namely,
\begin{equation}
 \sum_{k=0}^{l-1}\slashed{\mathcal{L}}_{R_{i_{l}}}...\slashed{\mathcal{L}}_{R_{i_{l-k+1}}}\leftexp{(i_{1}...i_{l-k})}{y}_{l-k}
\end{equation}
$\leftexp{(i_{1}...i_{j})}{y}_{j}$ is the sum of 
\begin{align*}
 \leftexp{(i_{1}...i_{j-1})}{a}_{j-1},\quad \slashed{d}(R_{i_{j-1}}...R_{i_{1}}\textrm{tr}\chi),\quad \slashed{d}\leftexp{(i_{1}...i_{j-1})}{\check{f}}_{j-1}
\end{align*}
with coefficients
\begin{align*}
 R_{i_{j}}\mu,\quad \mu R_{i_{j}}\textrm{tr}\chi-R_{i_{j}}L\mu+\leftexp{(R_{i_{j}})}{Z}\mu
=\mu R_{i_{j}}(\textrm{tr}\chi-e)-R_{i_{j}}m+\leftexp{(R_{i_{j}})}{Z}\mu-eR_{i_{j}}\mu,\quad
\frac{1}{2}(R_{i_{j}}\textrm{tr}\chi)
\end{align*}
of order 1,2,2, respectively. The first two of these coefficients involve the $1$st angular derivatives of $\mu$. The order of 
$\slashed{d}(R_{i_{j-1}}...R_{i_{1}}\textrm{tr}\chi)$ and $\slashed{d}\leftexp{(i_{1}...i_{j-1})}{\check{f}}_{j-1}$ is $j+1$.
On the other hand, $\leftexp{(i_{1}...i_{j-1})}{a}_{j-1}$ is given by (8.78). By the above discussion of $f_{0}$, the first term on the right of 
(8.78) is of order $j+2$, its principal terms being angular derivatives of $h$. This is the principal term in 
$\leftexp{(i_{1}...i_{j-1})}{a}_{j-1}$. The second term on the right, $\leftexp{(i_{1}...i_{j-1})}{b}_{j-1}$, given by (8.79) is of lower order,
namely $j+1$. In fact only the first term in (8.79), given by (8.80), is of order $j+1$ while the second term, the double sum, is of order $j$. Also, 
$\leftexp{(i_{1}...i_{j-1})}{a}_{j-1}$ does not involve the function $\mu$. So we conclude that $\leftexp{(i_{1}...i_{j})}{y}_{j}$ 
is of order $j+2$ and its principal terms are $j+2$ order angular derivatives of $h$. Then the order of (8.81) is $l+2$, its principal terms are $l+2$ 
order angular derivatives of $h$. Moreover, (8.81) contains $l$ order angular derivatives of $\mu$.  

To complete the investigation of $\leftexp{(i_{1}...i_{l})}{g}_{l}$, we must consider:
\begin{equation}
 \sum_{k=0}^{l-1}\slashed{\mathcal{L}}_{R_{i_{l}}}...\slashed{\mathcal{L}}_{R_{i_{l-k+1}}}\slashed{\mathcal{L}}_{\leftexp{(R_{i_{l-k}})}{Z}}
\leftexp{(i_{1}...i_{l-k-1})}{x}_{l-k-1}
\end{equation}
which we can write as:
\begin{equation}
  \sum_{k=0}^{l-1}\slashed{\mathcal{L}}_{\leftexp{(R_{i_{l-k}})}{Z}}\slashed{\mathcal{L}}_{R_{i_{l}}}...\slashed{\mathcal{L}}_{R_{i_{l-k+1}}}
\leftexp{(i_{1}...i_{l-k-1})}{x}_{l-k-1}
\end{equation}
minus the commutator term:
\begin{equation}
 \sum_{k=1}^{l-1}[\slashed{\mathcal{L}}_{\leftexp{(R_{i_{l-k}})}{Z}},   \slashed{\mathcal{L}}_{R_{i_{l}}}...\slashed{\mathcal{L}}_{R_{i_{l-k+1}}}]
\leftexp{(i_{1}...i_{l-k-1})}{x}_{l-k-1}
\end{equation}
We shall show this commutator term is of lower order by using the following lemma.

$\textbf{Lemma 8.5}$ For any $S_{t,u}$-tangential vectorfield $X$ and any $S_{t,u}$ 1-form $\xi$ and positive integer $l$ we have:
\begin{align*}
 [\slashed{\mathcal{L}}_{X},\slashed{\mathcal{L}}_{R_{i_{l}}}...\slashed{\mathcal{L}}_{R_{i_{1}}}]\xi=-\sum_{j=1}^{l}
\sum_{k_{1}<...<k_{j}=1}^{l}\slashed{\mathcal{L}}_{\leftexp{(i_{k_{1}}...i_{k_{j}})}{Y}}
\leftexp{(i_{1}\overset{>i_{k_{1}}...i_{k_{j}}<}{...}i_{l})}{\xi}_{l-j}
\end{align*}
Here $\leftexp{(i_{k_{1}}...i_{k_{j}})}{Y}$ is the $S_{t,u}$-tangential vectorfield:
\begin{align*}
 \leftexp{(i_{k_{1}}...i_{k_{j}})}{Y}=\slashed{\mathcal{L}}_{R_{i_{k_{j}}}}...\slashed{\mathcal{L}}_{R_{i_{k_{1}}}}X
\end{align*}
and $\leftexp{(i_{1}...i_{m})}{\xi}_{m}$ is the $S_{t,u}$ 1-form:
\begin{align*}
 \leftexp{(i_{1}...i_{m})}{\xi}_{m}=\slashed{\mathcal{L}}_{R_{i_{m}}}...\slashed{\mathcal{L}}_{R_{i_{1}}}\xi
\end{align*}
The proof is by induction, starting from the case $l=1$, which is a standard fact in differential geometry, in view of the fact that everything restricts to 
$S_{t,u}$. 

Applying Lemma 8.5 we see that (8.84) is minus:
\begin{align}
 \sum_{k=1}^{l-1}\sum_{j=1}^{k}\sum_{m_{1}<...<m_{j}=l-k+1}^{l}\slashed{\mathcal{L}}_{\leftexp{(i_{m_{1}}...i_{m_{j}};i_{l-k})}{Z}}
\slashed{\mathcal{L}}_{R_{i_{l}}}\overset{>i_{m_{j}}...i_{m_{1}}<}{...}\slashed{\mathcal{L}}_{R_{i_{l-k+1}}}
\leftexp{(i_{1}...i_{l-k-1})}{x}_{l-k-1}
\end{align}
Here the indexes $i_{m_{j}},...,i_{m_{1}}$ are absent, and $\leftexp{(i_{m_{1}}...i_{m_{j}};i_{l-k})}{Z}$ is the $S_{t,u}$-tangential vectorfield:
\begin{equation}
 \leftexp{(i_{m_{1}}...i_{m_{j}};i_{l-k})}{Z}=\slashed{\mathcal{L}}_{R_{i_{m_{j}}}}...\slashed{\mathcal{L}}_{R_{i_{m_{1}}}}
\leftexp{(R_{i_{l-k}})}{Z}
\end{equation}
From (6.57), we know that $\leftexp{(R_{i})}{Z}$ is of order 1 and its principal acoustical term involves $\chi$. Hence 
$\leftexp{(i_{m_{1}}...i_{m_{j}};i_{l-k})}{Z}$ is of order $j+1$ and its principal acoustical term is a $j$th order angular derivative of $\chi$.
In (8.85), the 1-form
\begin{equation}
 \slashed{\mathcal{L}}_{R_{i_{l}}}\overset{>i_{m_{j}}...i_{m_{1}}<}{...}\slashed{\mathcal{L}}_{R_{i_{l-k+1}}}
\leftexp{(i_{1}...i_{l-k-1})}{x}_{l-k-1}
\end{equation}
is of order $l-j+1$ as it is an angular derivative of order $k-j$ of $\leftexp{(i_{1}...i_{l-k-1})}x_{l-k-1}$. If $X$ is an $S_{t,u}$-tangential vectorfield of order
$m$ and $\xi$ an $S_{t,u}$ 1-form of order $n$, then the order of $\slashed{\mathcal{L}}_{X}\xi$ is max$\{m,n\}+1$. This is evident from the formula:
\begin{equation}
 (\slashed{\mathcal{L}}_{X}\xi)(Y)=X(\xi(Y))-\xi(\slashed{\mathcal{L}}_{X}Y) \notag
\end{equation}
Thus, the order of the terms in the triple sum in (8.85) is 
\begin{align*}
 \max\{j+1,l-j+1\}+1
\end{align*}
and the two extreme terms corresponding to $j=1$ and $j=l-1$ are of highest order, namely $l+1$.
So (8.85) is of order $l+1$, thus a lower order term.

We now turn to the principal term (8.83). In view of (8.44) and the fact that for any function $\phi$:
\begin{align*}
 \slashed{\mathcal{L}}_{R_{i_{l}}}...\slashed{\mathcal{L}}_{R_{i_{1}}}\slashed{d}\phi
=\slashed{d}(R_{i_{l}}...R_{i_{1}}\phi)
\end{align*}
we obtain:
\begin{align}
 \slashed{\mathcal{L}}_{R_{i_{l}}}...\slashed{\mathcal{L}}_{R_{i_{l-k+1}}}
\leftexp{(i_{1}...i_{l-k-1})}{x}_{l-k-1}\\ \notag=\leftexp{(i_{l}...i_{l-k+1}i_{l-k-1}...i_{1})}{x}_{l-1}\\ \notag
+\sum_{j=1}^{k}\sum_{m_{1}<...<m_{j}=l-k+1}^{l}(R_{i_{m_{j}}}...R_{i_{m_{1}}}\mu)\slashed{d}(R_{i_{l}}\overset{>i_{m_{j}}...i_{m_{1}}i_{l-k}<}{...}R_{i_{1}}\textrm{tr}\chi)
\end{align}
Here the indexes $i_{m_{j}}...i_{m_{1}}i_{l-k}$ are absent. The first term is of order $l+1$, while the double sum is of order $l$. The contribution of the 
first term in (8.88) to (8.83) is then the principal part of (8.82). This contribution is 
\begin{equation}
 \sum_{k=0}^{l-1}\slashed{\mathcal{L}}_{\leftexp{(R_{i_{l-k}})}{Z}}\leftexp{(i_{1}...i_{l-k-1}i_{l-k+1}...i_{l})}{x}_{l-1}
\end{equation}
We conclude that the principal acoustical terms of $\leftexp{(i_{1}...i_{l})}{g}_{l}$ are all the terms of the sum (8.89). Moreover 
(8.82) does not involve the principal order of function $\mu$.

We shall estimate each term in (8.89) pointwise in terms of 
\begin{equation}
 \max_{j}|\leftexp{(i_{1}...i_{l-k-1}i_{l-k+1}...i_{l}j)}{x}_{l}| \notag
\end{equation}

$\textbf{Lemma 8.6}$ Let $\phi$ and $\rho$ be arbitrary functions and $\mu$ a nonnegative function. Denote by $\xi$ the $S_{t,u}$ 1-form:
\begin{align*}
 \xi=\mu\slashed{d}\phi+\slashed{d}\rho
\end{align*}
and by $\eta_{j}$ the $S_{t,u}$ 1-forms
\begin{align*}
 \eta_{j}=\mu\slashed{d}(R_{j}\phi)+\slashed{d}(R_{j}\rho): j=1,2,3
\end{align*}
Then for any $S_{t,u}$-tangential vectorfield $X$ we have the pointwise estimate:
\begin{align*}
 |\slashed{\mathcal{L}}_{X}\xi|\leq C(1+t)^{-1}\{|X|(\max_{j}|\eta_{j}|+2\max_{j}|\slashed{d}(R_{j}\rho)|
+(1+t)|\slashed{d}\mu||\slashed{d}\phi|)\\
+(\max_{j}|\slashed{\mathcal{L}}_{R_{j}}X|)(\mu|\slashed{d}\phi|+|\slashed{d}\rho|)\}
\end{align*}
$Proof$. We have:
\begin{equation}
 \slashed{\mathcal{L}}_{X}\xi=\mu\slashed{d}(X\phi)+(X\mu)\slashed{d}\phi+\slashed{d}(X\rho)
\end{equation}
hence:
\begin{equation}
 |\slashed{\mathcal{L}}_{X}\xi|\leq \mu|\slashed{d}(X\phi)|+|X||\slashed{d}\mu||\slashed{d}\phi|+|\slashed{d}(X\rho)|
\end{equation}
Now, by $\textbf{H0}$,
\begin{equation}
 |\slashed{d}(X\phi)|\leq C(1+t)^{-1}\max_{j}|R_{j}(X\phi)|
\end{equation}
and we have:
\begin{align}
 R_{j}(X\phi)=X(R_{j}\phi)+[R_{j},X]\phi=X\cdot\slashed{d}(R_{j}\phi)+(\slashed{\mathcal{L}}_{R_{j}}X)\cdot\slashed{d}\phi
\end{align}
hence:
\begin{align}
 \max_{j}|R_{j}(X\phi)|\leq |X|\max_{j}|\slashed{d}(R_{j}\phi)|+(\max_{j}|\slashed{\mathcal{L}}_{R_{j}}X|)|\slashed{d}\phi|
\end{align}
Then we get:
\begin{equation}
 |\slashed{d}(X\phi)|\leq C(1+t)^{-1}\{|X|\max_{j}|\slashed{d}(R_{j}\phi)|+(\max_{j}|\slashed{\mathcal{L}}_{R_{j}}X|)|\slashed{d}\phi|\}
\end{equation}
By the definition of $\eta_{j}$,
\begin{align}
 \mu\max_{j}|\slashed{d}(R_{j}\phi)|=\max_{j}|\mu\slashed{d}(R_{j}\phi)|\\
=\max_{j}|\eta_{j}-\slashed{d}(R_{j}\rho)|\leq \max_{j}|\eta_{j}|+\max_{j}|\slashed{d}(R_{j}\rho)| \notag
\end{align}
Moreover, from (8.95) with $\rho$ in the role of $\phi$ we have:
\begin{equation}
 |\slashed{d}(X\rho)|\leq C(1+t)^{-1}\{|X|\max_{j}|\slashed{d}(R_{j}\rho)|+(\max_{j}|\slashed{\mathcal{L}}_{R_{j}}X|)|\slashed{d}\rho|\}
\end{equation}
In view of (8.95)-(8.97), the lemma follows through (8.91). $\qed$

We now take $\phi=R_{i_{l}}...R_{i_{l-k+1}}R_{i_{l-k-1}}...R_{i_{1}}\textrm{tr}\chi$, and
$\rho=R_{i_{l}}...R_{i_{l-k+1}}R_{i_{l-k-1}}...R_{i_{1}}\check{f}\\
=\leftexp{(i_{1}...i_{l-k-1}i_{l-k+1}...i_{l})}{\check{f}}_{l-1}$.
Then $\xi$ is the $S_{t,u}$ 1-form $\leftexp{(i_{1}...i_{l-k-1}i_{l-k+1}...i_{l})}{x}_{l-1}$ and $\eta_{j}$ the $S_{t,u}$ 1-form
$\leftexp{(i_{1}...i_{l-k-1}i_{l-k+1}...i_{l}j)}{x}_{l}$. Moreover, we take $X$ to be the $S_{t,u}$-tangential vectorfield
$\leftexp{(R_{i_{l-k}})}{Z}$. Noting that:
\begin{align*}
 \mu|\slashed{d}\phi|+|\slashed{d}\rho|\leq |\xi|+2|\slashed{d}\rho|,
\end{align*}
we obtain:
\begin{align}
 |\slashed{\mathcal{L}}_{\leftexp{(R_{i_{l-k}})}{Z}}\leftexp{(i_{1}...i_{l-k-1}i_{l-k+1}...i_{l})}{x}_{l-1}|
\leq C(1+t)^{-1}\{|\leftexp{(R_{i_{l-k}})}{Z}|(\max_{j}|\leftexp{(i_{1}...i_{l-k-1}i_{l-k+1}...i_{l}j)}{x}_{l}|\\ \notag
+2\max_{j}|\slashed{d}(\leftexp{(i_{1}...i_{l-k-1}i_{l-k+1}...i_{l}j)}{\check{f}}_{l})|+(1+t)|\slashed{d}\mu|
|\slashed{d}(R_{i_{l}}...R_{i_{l-k+1}}R_{i_{l-k-1}}...R_{i_{1}}\textrm{tr}\chi)|)\\ \notag
+(\max_{j}|\slashed{\mathcal{L}}_{R_{j}}\leftexp{(R_{i_{l-k}})}{Z}|)(|\leftexp{(i_{1}...i_{l-k-1}i_{l-k+1}...i_{l})}{x}_{l-1}|
+2|\slashed{d}(\leftexp{(i_{1}...i_{l-k-1}i_{l-k+1}...i_{l})}{\check{f}}_{l-1})|)\}
\end{align}

\subsection{Elliptic Theory on $S_{t,u}$}
Returning to the propagation equation for $\leftexp{(i_{1}...i_{l})}{x}_{l}$ given by Proposition 8.3, the term $\mu\slashed{d}
(\leftexp{(i_{1}...i_{l})}{h}_{l})$ remains to be considered. Since $\leftexp{(i_{1}...i_{l})}{h}_{l}$ given by (8.64) is $l$th order 
angular derivative of $|\hat{\chi}|^{2}$, so the term to be considered is a principal acoustical term. Moreover, it involves $l+1$st
order angular derivatives of $\hat{\chi}$, not $\textrm{tr}\chi$, which is what the propagation equation for $\leftexp{(i_{1}...i_{l})}{x}_{l}$ 
allows us to control.

The term in question comes from the right-hand side of propagation equation (8.36) for $x_{0}=\mu\slashed{d}\textrm{tr}\chi+\slashed{d}\check{f}$, 
from the term:
\begin{equation}
 \mu\slashed{d}(|\hat{\chi}|^{2})=2\mu\hat{\chi}\cdot\slashed{D}\hat{\chi} \notag
\end{equation}
Now, the propagation equation (8.36) must be considered in conjunction with the Codazzi equation (4.36):
\begin{equation}
 \slashed{\textrm{div}}\hat{\chi}=\frac{1}{2}\slashed{d}\textrm{tr}\chi+i
\end{equation}
where
\begin{equation}
 i=\beta^{*}_{A}-\mu^{-1}(\zeta_{B}\chi_{A}^{B}-\zeta_{A}\textrm{tr}\chi)
\end{equation}
which is regular as $\mu\longrightarrow 0$. Equation (8.99) allows us to estimate, on $S_{t,u}$, $\hat{\chi}$ and $\slashed{D}\hat{\chi}$ in 
terms of $\slashed{d}\textrm{tr}\chi$ and the derivatives of $\psi_{\mu}$ of up to 1st order. 
The propagation equation (8.36) allows us to estimate, along each generator, $\slashed{d}\textrm{tr}\chi$ in terms of $\slashed{d}(|\hat{\chi}|^{2})$, 
that is in terms of $\hat{\chi}$ and $\slashed{D}\hat{\chi}$, and the derivatives of $\psi_{\mu}$ of up to the 2nd order. Thus, considering on $C_{u}$ 
the propagation equation (8.36) in conjunction with the Codazzi equation (8.99) we can estimate $\chi$ and its 1st order angular derivatives in 
terms of $\psi_{\mu}$ and their derivatives of up to the 2nd order, avoiding in this way loss of differentiability. The idea coupling the propagation 
equations for $\chi$ along the generators of $C_{u}$ with the Codazzi equations at each $S_{t,u}$ was also used in \cite{CK}, where the Ricci curvature
tensor vanishes.

To get a higher order version of (8.99), we shall use the following lemma:

$\textbf{Lemma 8.7}$ Let $(M,g)$ be a 2-dimensional Riemannian manifold, let $X$ be an arbitrary vectorfield on $M$ and let $\theta$ 
be a trace-free symmetric 2-covariant tensorfield on $(M,g)$ satisfying the equation:
\begin{equation}
 \textrm{div}_{g}\theta=f
\end{equation}
 for some 1-form $f$ on $M$. Then $\hat{\mathcal{L}}_{X}\theta$, the trace-free part of $\mathcal{L}_{X}\theta$, satisfies the equation:
\begin{equation}
 \textrm{div}_{g}\hat{\mathcal{L}}_{X}\theta=\leftexp{(X)}{\dot{f}}
\end{equation}
where $\leftexp{(X)}{\dot{f}}$ is the 1-form given, in an arbitrary local frame, by:
\begin{align*}
 \leftexp{(X)}{\dot{f}}_{a}=(\mathcal{L}_{X}f)_{a}+\frac{1}{2}\textrm{tr}\leftexp{(X)}{\pi}f_{a}\\
+\frac{1}{2}\leftexp{(X)}{\hat{\pi}}^{bc}(\nabla_{b}\theta_{ac}+\nabla_{c}\theta_{ab}-\nabla_{a}\theta_{bc})+
(\textrm{div}_{g}\leftexp{(X)}{\hat{\pi}})^{b}\theta_{ab}
\end{align*}
Here, $\leftexp{(X)}{\pi}=\mathcal{L}_{X}g$.

$Proof$. Let $\phi_{t}$ be the local 1-parameter group generated by $X$ and let $\phi_{t*}$ be the corresponding pullback. We then have:
\begin{equation}
 \textrm{div}_{\phi_{t*}g}(\phi_{t*}\theta)=\phi_{t*}(\textrm{div}_{g}\theta)=\phi_{t*}f
\end{equation}
Now in an arbitrary local coordinates:
\begin{align}
 (\textrm{div}_{g}\theta)_{a}=(g^{-1})^{bc}\nabla_{c}\theta_{ab}\\ \notag
=(g^{-1})^{bc}(\frac{\partial\theta_{ab}}{\partial x^{c}}-\Gamma_{ca}^{d}\theta_{db}-\Gamma_{cb}^{d}\theta_{ad})
\end{align}
where $\Gamma^{c}_{ab}$ are the Christoffel symbols of the metric $g$. Similarly,
\begin{align}
 (\textrm{div}_{\phi_{t*}g}(\phi_{t*}\theta))_{a}=((\phi_{t*}g)^{-1})^{bc}(\frac{\partial(\phi_{t*}\theta)_{ab}}{\partial x^{c}}
-\overset{\phi_{t*}g}{{\Gamma}_{ca}^{d}}(\phi_{t*}\theta)_{db}-\overset{\phi_{t*}g}{{\Gamma}_{cb}^{d}}(\phi_{t*}\theta)_{ad})
\end{align}
where $\overset{\phi_{t*}g}{{\Gamma}_{ab}^{c}}$ is the Christoffel symbols of the metric $\phi_{t*}g$.
We have:
\begin{equation}
 \Gamma^{c}_{ab}=\frac{1}{2}(g^{-1})^{cd}(\frac{\partial g_{bd}}{\partial x^{a}}+\frac{\partial g_{ad}}{\partial x^{b}}
-\frac{\partial g_{ab}}{\partial x^{d}})
\end{equation}
and
\begin{equation}
\overset{\phi_{t*}g}{{\Gamma}^{c}_{ab}}=\frac{1}{2}((\phi_{t*}g)^{-1})^{cd}(\frac{\partial(\phi_{t*}g)_{bd}}{\partial x^{a}}
+\frac{\partial(\phi_{t*}g)_{ad}}{\partial x^{b}}-\frac{\partial(\phi_{t*}g)_{ab}}{\partial x^{d}})
\end{equation}
Since, by definition,
\begin{equation}
 (\frac{d}{dt}\phi_{t*}g)_{t=0}=\mathcal{L}_{X}g=\leftexp{(X)}{\pi}
\end{equation}
differentiating (8.107) with respect to $t$ at $t=0$, we obtain , in view of (8.106), 
the following formula:
\begin{equation}
 (\frac{d}{dt}\overset{\phi_{t*}g}{{\Gamma}^{c}_{ab}})_{t=0}=\frac{1}{2}(g^{-1})^{cd}
(\nabla_{a}\leftexp{(X)}{\pi}_{bd}+\nabla_{b}\leftexp{(X)}{\pi}_{ad}-\nabla_{d}\leftexp{(X)}{\pi}_{ab})
\end{equation}
Differentiating then (8.105) with respect to $t$ at $t=0$ and using (8.109) and (8.108) and the fact, by definition:
\begin{equation}
 (\frac{d}{dt}\phi_{t*}\theta)_{t=0}=\mathcal{L}_{X}\theta
\end{equation}
we obtain:
\begin{align}
 (\frac{d}{dt}(\textrm{div}_{\phi_{t*}g}(\phi_{t*}\theta))_{a})_{t=0}=(\textrm{div}_{g}(\mathcal{L}_{X}\theta))_{a}-\leftexp{(X)}{\pi}^{bc}\nabla_{c}\theta_{ab}\\
-\frac{1}{2}\{(\nabla_{a}\leftexp{(X)}{\pi}^{bc})\theta_{bc}+(2\nabla_{c}\leftexp{(X)}{\pi}^{bc}-
\nabla^{b}\textrm{tr}\leftexp{(X)}{\pi})\theta_{ab}\} \notag
\end{align}
From (8.103) and the definition:
\begin{equation}
 (\frac{d}{dt}\phi_{t*}f)_{t=0}=\mathcal{L}_{X}f
\end{equation}
we conclude:
\begin{equation}
 \textrm{div}_{g}(\mathcal{L}_{X}\theta)=\leftexp{(X)}{f'}
\end{equation}
where 
\begin{align}
 \leftexp{(X)}{f'}_{a}=(\mathcal{L}_{X}f)_{a}+\leftexp{(X)}{\pi}^{bc}\nabla_{c}\theta_{ab}\\
+\frac{1}{2}\{(\nabla_{a}\leftexp{(X)}{\pi}^{bc})\theta_{bc}+(2\nabla_{c}\leftexp{(X)}{\pi}^{bc}-
\nabla^{b}\textrm{tr}\leftexp{(X)}{\pi})\theta_{ab}\} \notag
\end{align}
In fact the above hold for any n-dimensional manifold $M$, and irrespective of the trace-free nature of $\theta$. In the present case where $M$ is 2-dimensional, 
decomposing $\leftexp{(X)}{\pi}$ into its trace-free part and its trace,
\begin{equation}
 \leftexp{(X)}{\pi}=\leftexp{(X)}{\hat{\pi}}+\frac{1}{2}g\textrm{tr}\leftexp{(X)}{\pi}
\end{equation}
we may write (8.114) as:
\begin{align}
  \leftexp{(X)}{f'}_{a}=(\mathcal{L}_{X}f)_{a}+\frac{1}{2}\textrm{tr}\leftexp{(X)}{\pi}f_{a}+\leftexp{(X)}{\hat{\pi}}^{bc}\nabla_{c}\theta_{ab}\\
+\frac{1}{2}(\nabla_{a}\leftexp{(X)}{\hat{\pi}}^{bc})\theta_{bc}+(\nabla_{c}\leftexp{(X)}{\hat{\pi}}^{bc})\theta_{ab} \notag
\end{align}
Now,
\begin{equation}
 \textrm{tr}(\mathcal{L}_{X}\theta)=(g^{-1})^{ab}(\mathcal{L}_{X}\theta)_{ab}=-(\mathcal{L}_{X}g^{-1})^{ab}\theta_{ab}
=\leftexp{(X)}{\pi}^{ab}\theta_{ab}=\leftexp{(X)}{\hat{\pi}}\cdot\theta
\end{equation}
where we have used the fact that $\theta$ is trace-free. Therefore, $\hat{\mathcal{L}}_{X}\theta$ is given by:
\begin{equation}
 \hat{\mathcal{L}}_{X}\theta=\mathcal{L}_{X}\theta-\frac{1}{2}g(\leftexp{(X)}{\hat{\pi}}\cdot\theta)
\end{equation}
and noting that for any function $\psi$ we have:
\begin{align*}
 \textrm{div}_{g}(g\psi)=d\psi
\end{align*}
we obtain:
\begin{equation}
 \textrm{div}_{g}(\hat{\mathcal{L}}_{X}\theta)=\textrm{div}_{g}(\mathcal{L}_{X}\theta)-\frac{1}{2}d(\leftexp{(X)}{\hat{\pi}}\cdot\theta)
\end{equation}
Hence, defining:
\begin{equation}
 \leftexp{(X)}{\dot{f}}=\leftexp{(X)}{f'}-\frac{1}{2}d(\leftexp{(X)}{\hat{\pi}}\cdot\theta)
\end{equation}
we conclude that $\textrm{div}_{g}\hat{\mathcal{L}}_{X}\theta=\leftexp{(X)}{\dot{f}}$ and $\leftexp{(X)}{\dot{f}}$ is as given in the lemma. $\qed$

$\textbf{Proposition 8.4}$ For each non-negative integer $l$ and each multi-index $(i_{1},...,i_{l})$, the trace-free symmetric 
2-covariant tensorfield
\begin{align*}
 \leftexp{(i_{1}...i_{l})}{\hat{\chi}}_{l}=\hat{\slashed{\mathcal{L}}}_{R_{i_{l}}}...\hat{\slashed{\mathcal{L}}}_{R_{i_{1}}}\hat{\chi}
\end{align*}
on $S_{t,u}$ satisfies the elliptic system:
\begin{align*}
 \slashed{\textrm{div}}(\leftexp{(i_{1}...i_{l})}{\hat{\chi}}_{l})=\frac{1}{2}\slashed{d}(R_{i_{l}}...R_{i_{1}}\textrm{tr}\chi)+
\leftexp{(i_{1}...i_{l})}{i}_{l}
\end{align*}
where the $S_{t,u}$ 1-form $\leftexp{(i_{1}...i_{l})}{i}_{l}$ is given by:
\begin{align*}
 \leftexp{(i_{1}...i_{l})}{i}_{l}=(\slashed{\mathcal{L}}_{R_{i_{l}}}+\frac{1}{2}\textrm{tr}\leftexp{(R_{i_{l}})}{\slashed{\pi}})...
(\slashed{\mathcal{L}}_{R_{i_{1}}}+\frac{1}{2}\textrm{tr}\leftexp{(R_{i_{1}})}{\slashed{\pi}})i\\
+\sum_{k=0}^{l-1}(\slashed{\mathcal{L}}_{R_{i_{l}}}+\frac{1}{2}\textrm{tr}\leftexp{(R_{i_{l}})}{\slashed{\pi}})...
(\slashed{\mathcal{L}}_{R_{i_{l-k+1}}}+\frac{1}{2}\textrm{tr}\leftexp{(R_{i_{l-k+1}})}{\slashed{\pi}})\leftexp{(i_{1}...i_{l-k})}{q}_{l-k} \notag
\end{align*}
Here $i$ is the $S_{t,u}$ 1-form given by (8.100) and for each $j=1,...,l$, $\leftexp{(i_{1}...i_{j})}{q}_{j}$ is the $S_{t,u}$ 1-form given by:
\begin{align*}
 \leftexp{(i_{1}...i_{j})}{q}_{j,A}=\frac{1}{4}\textrm{tr}\leftexp{(R_{i_{j}})}{\slashed{\pi}}\slashed{d}(R_{i_{j-1}}...R_{i_{1}}\textrm{tr}\chi)\\
+\frac{1}{2}\leftexp{(R_{i_{j}})}{\hat{\slashed{\pi}}}^{BC}\{\slashed{D}_{B}(\leftexp{(i_{1}...i_{j-1})}{\hat{\chi}}_{j-1})_{AC}+
\slashed{D}_{C}(\leftexp{(i_{1}...i_{j-1})}{\hat{\chi}}_{j-1})_{AB}\\-\slashed{D}_{A}(\leftexp{(i_{1}...i_{j-1})}{\hat{\chi}}_{j-1})_{BC}\}
+(\slashed{\textrm{div}}\leftexp{(R_{i_{j}})}{\hat{\slashed{\pi}}})^{B}(\leftexp{(i_{1}...i_{j-1})}{\hat{\chi}}_{j-1})_{AB}
\end{align*}
$Proof$. The proposition reduces for $l=0$, since $i=i_{0}$, to (8.99). By induction, we assume the elliptic system holds with $l$ replaced by $l-1$:
\begin{align*}
 \slashed{\textrm{div}}(\leftexp{(i_{1}...i_{l-1})}{\hat{\chi}}_{l-1})=\frac{1}{2}\slashed{d}(R_{i_{l-1}}...R_{i_{1}}\textrm{tr}\chi)+
\leftexp{(i_{1}...i_{l-1})}{i}_{l-1}
\end{align*}
Taking $(M,g)$ to be $(S_{t,u},\slashed{g})$, $\theta$ to be $\leftexp{(i_{1}...i_{l-1})}{\hat{\chi}}_{l-1}$, $X$ to be $R_{i_{l}}$,
and $f$ to be 
\begin{align*}
 \frac{1}{2}\slashed{d}(R_{i_{l-1}}...R_{i_{1}}\textrm{tr}\chi)+
\leftexp{(i_{1}...i_{l-1})}{i}_{l-1}
\end{align*}
in Lemma 8.7, and taking into account the fact that:
\begin{align*}
 \slashed{\mathcal{L}}_{R_{i_{l}}}\slashed{d}(R_{i_{l-1}}...R_{i_{1}}\textrm{tr}\chi)=\slashed{d}(R_{i_{l}}...R_{i_{1}}\textrm{tr}\chi)
\end{align*}
we get a recursion formula for $\leftexp{(i_{1}...i_{l})}{i}_{l}$: 
\begin{align}
\leftexp{(i_{1}...i_{l})}{i}_{l}=\slashed{\mathcal{L}}_{R_{i_{l}}}\leftexp{(i_{1}...i_{l-1})}{i}_{l-1}
+\frac{1}{2}\textrm{tr}\leftexp{(R_{i_{l}})}{\slashed{\pi}}\leftexp{(i_{1}...i_{l-1})}{i}_{l-1}+\leftexp{(i_{1}...i_{l})}{q}_{l} 
\end{align}
Applying Proposition 8.2 to this recursion with the space of $S_{t,u}$ 1-forms in the role of the space $X$, the operators
\begin{align*}
 \slashed{\mathcal{L}}_{R_{i_{l}}}+\frac{1}{2}\textrm{tr}\leftexp{(R_{i_{l}})}{\slashed{\pi}}
\end{align*}
in the role of the operators $A_{n}$
and $\leftexp{(i_{1}...i_{l})}{i}_{l}$, $\leftexp{(i_{1}...i_{l})}{q}_{l}$, in the role of $x_{n}$, $y_{n}$, respectively, 
the proposition is proved. $\qed$

Let us investigate the order of the terms in $\leftexp{(i_{1}...i_{l})}{i}_{l}$. We begin with $i$, given by (8.100). From the expression (4.39)-(4.40), 
we know that there is no 2nd order term in $i$, so the first term in the expression for $\leftexp{(i_{1}...i_{l})}{i}_{l}$ contains no 
principal term. Moreover, the second term is of order $l+1$. This can be easily got from the fact that 
$\leftexp{(i_{1}...i_{j})}{q}_{j}$ is of order $j+1$. Also $\leftexp{(i_{1}...i_{l})}{i}_{l}$ does not involve $\mu$.

$Remark$:
From the propagation equation in Proposition 8.3, we can control $\slashed{d}(R_{i_{l}}...R_{i_{1}}\textrm{tr}\chi)$ in terms of 
$2\mu\hat{\chi}\slashed{D}(\hat{\slashed{\mathcal{L}}}_{R_{i_{l}}}...\hat{\slashed{\mathcal{L}}}_{R_{i_{1}}}\hat{\chi})$ and the derivatives of 
$\psi_{\mu}$ up to $l+2$nd order. While using the elliptic system in Proposition 8.4, we can control 
$\slashed{D}(\hat{\slashed{\mathcal{L}}}_{R_{i_{l}}}...\hat{\slashed{\mathcal{L}}}_{R_{i_{1}}}\hat{\chi})$ in terms of 
$\slashed{d}(R_{i_{l}}...R_{i_{1}}\textrm{tr}\chi)$ and the derivatives of $\psi_{\mu}$ up to $l+1$st order.

Since the elliptic system is to be considered in conjunction with the propagation equation of Proposition 8.3, the right hand side of the elliptic system for 
$\hat{\slashed{\mathcal{L}}}_{R_{i_{l}}}...\hat{\slashed{\mathcal{L}}}_{R_{i_{1}}}\hat{\chi}$ is thus to be expressed in terms of 
$\leftexp{(i_{1}...i_{l})}{x}_{l}$. From (8.40) we have:
\begin{align*}
 \mu\slashed{d}(R_{i_{l}}...R_{i_{1}}\textrm{tr}\chi)=\leftexp{(i_{1}...i_{l})}{x}_{l}-\slashed{d}\leftexp{(i_{1}...i_{l})}{\check{f}}_{l}
\end{align*}
Then we can only control the product of $\mu$ with the right-hand side of the elliptic system in terms of $\leftexp{(i_{1}...i_{l})}{x}_{l}$.
This shall lead to $\mu$-weighted $L^{2}$ estimates on $S_{t,u}$.

$\textbf{Lemma 8.8}$ Let $(M,g)$ be a compact 2-dimensional Riemannian manifold, and let $\theta$ be a trace-free symmetric 2-covariant tensorfield on 
$(M,g)$ satisfying the equation:
\begin{equation}
 \textrm{div}_{g}\theta=f
\end{equation}
for some 1-form $f$ on $M$. Let also $\mu$ be an arbitrary non-negative function on $M$. Then the following estimate holds on $M$:
\begin{equation}
 \int_{M}\mu^{2}\{\frac{1}{2}|\nabla\theta|^{2}+2K|\theta|^{2}\}d\mu_{g}\leq
3\int_{M}\mu^{2}|f|^{2}d\mu_{g}+3\int_{M}|d\mu|^{2}|\theta|^{2}d\mu_{g} \notag
\end{equation}
where $K$ is the Gauss curvature of $(M,g)$.

$Proof$. Consider the 3-covariant tensorfield given in an arbitrary local frame by:
\begin{equation}
 \omega_{abc}=\nabla_{a}\theta_{bc}-\nabla_{b}\theta_{ac}
\end{equation}
For any vectorfield $X$, $\omega_{abc}X^{c}$ defines a 2-form on $M$. Since $\textrm{dim}M=2$, this 2-form must be a function $\phi$, times the 
volume form $\epsilon$ of $(M,g)$. Since $\phi$ depends linearly on $X$, there must be a 1-form $e$ on $M$ such that $e\cdot X=\phi$.
So we can write:
\begin{equation}
 \omega_{abc}=\epsilon_{ab}e_{c}
\end{equation}
Contracting the right hand side of (8.124) with $(g^{-1})^{ac}$ we get:
\begin{align*}
 (g^{-1})^{ac}\epsilon_{ab}e_{c}=e^{a}\epsilon_{ab}=\leftexp{*}{e}_{b}
\end{align*}
(the Hodge dual of $e$). Contracting the left hand side of (8.124) with $(g^{-1})^{ac}$ and taking into account the fact that $(g^{-1})^{ac}\theta_{ac}=0$
we obtain:
\begin{align*}
 (g^{-1})^{ac}\omega_{abc}=(\textrm{div}\theta)_{b}=f_{b}
\end{align*}
We conclude that:
\begin{equation}
 \leftexp{*}{e}=f, e=-\leftexp{*}{f}
\end{equation}
 and (8.124) becomes:
\begin{equation}
 \omega_{abc}=-\epsilon_{ab}\leftexp{*}{f}_{c}
\end{equation}
In view of the fact that
\begin{align*}
 \epsilon_{ac}\epsilon_{b}{}^{c}=g_{ab}
\end{align*}
It follows that:
\begin{equation}
 \frac{1}{2}\omega_{abc}\omega^{abc}=\frac{1}{2}\epsilon_{ab}\leftexp{*}{f}_{c}\epsilon^{ab}\leftexp{*}{f}^{c}=|\leftexp{*}{f}|^{2}=|f|^{2}
\end{equation}

On the other hand, from the definition (8.123) we have:
\begin{align}
 \frac{1}{2}\omega_{abc}\omega^{abc}=\frac{1}{2}(\nabla_{a}\theta_{bc}-\nabla_{b}\theta_{ac})
(\nabla^{a}\theta^{bc}-\nabla^{b}\theta^{ac})\\
=|\nabla\theta|^{2}-\nabla_{a}\theta_{bc}\nabla^{b}\theta^{ac} \notag
\end{align}
We write:
\begin{equation}
 \nabla_{a}\theta_{bc}\nabla^{b}\theta^{ac}=\nabla^{b}(\theta^{ac}\nabla_{a}\theta_{bc})-\theta^{ac}\nabla^{b}\nabla_{a}\theta_{bc}
\end{equation}
Using a well known fact in differential geometry, we have:
\begin{align}
 \nabla^{b}\nabla_{a}\theta_{bc}-\nabla_{a}\nabla^{b}\theta_{bc}=R_{b}{}^{db}{}_{a}\theta_{dc}+R_{c}{}^{db}{}_{a}\theta_{bd}\\
=S^{d}{}_{a}\theta_{dc}+R_{c}{}^{db}{}_{a}\theta_{bd}\notag
\end{align}
where $R_{abcd}$ is the curvature tensor and $S_{ab}=(g^{-1})^{cd}R_{cadb}$ the Ricci tensor of $(M,g)$. Since $M$ is 2-dimensional, we have:
\begin{align*}
 R_{abcd}=K(g_{ac}g_{bd}-g_{ad}g_{bc}), \quad S_{ab}=Kg_{ab}
\end{align*}
Hence, the right hand side of (8.130) is:
\begin{align*}
 K\delta^{d}_{a}\theta_{dc}+K(\delta^{b}_{c}\delta^{d}_{a}-g_{ca}(g^{-1})^{db})\theta_{bd}=2K\theta_{ac}
\end{align*}
in view of the fact that $\textrm{tr}\theta=0$. Thus, (8.130) reduces to:
\begin{equation}
 \nabla^{b}\nabla_{a}\theta_{bc}-\nabla_{a}\nabla^{b}\theta_{bc}=2K\theta_{ac}
\end{equation}
Substituting in (8.129) we obtain:
\begin{equation}
 \nabla_{a}\theta_{bc}\nabla^{b}\theta^{ac}=\nabla^{b}(\theta^{ac}\nabla_{a}\theta_{bc})-\theta^{ac}\nabla_{a}\nabla^{b}\theta_{bc}-2K|\theta|^{2}
\end{equation}
In reference to the second term on the right, we write:
\begin{align*}
 \theta^{ac}\nabla_{a}\nabla^{b}\theta_{bc}=\nabla_{a}(\theta^{ac}\nabla^{b}\theta_{bc})-(\nabla_{a}\theta^{ac})(\nabla^{b}\theta_{bc})
=\nabla_{a}(\theta^{ac}\nabla^{b}\theta_{bc})-|f|^{2}
\end{align*}
Substituting in (8.132) then yields:
\begin{equation}
 \nabla_{a}\theta_{bc}\nabla^{b}\theta^{ac}=\textrm{div}_{g}J-2K|\theta|^{2}+|f|^{2}
\end{equation}
where $J$ is the vectorfield:
\begin{equation}
 J^{a}=\theta^{b}{}_{c}\nabla_{b}\theta^{ac}-\theta^{a}{}_{c}\nabla_{b}\theta^{bc}
=\theta^{b}{}_{c}\nabla_{b}\theta^{ac}-\theta^{a}{}_{c}f^{c}
\end{equation}
Substituting (8.133) in (8.128) we get:
\begin{equation}
 \frac{1}{2}\omega_{abc}\omega^{abc}=|\nabla\theta|^{2}+2K|\theta|^{2}-|f|^{2}-\textrm{div}_{g}J
\end{equation}
Comparing (8.135) with (8.127) we conclude that:
\begin{equation}
 |\nabla\theta|^{2}+2K|\theta|^{2}=2|f|^{2}+\textrm{div}_{g}J
\end{equation}
Multiplying this equation by $\mu^{2}$ and integrating over $M$, we obtain:
\begin{align}
 \int_{M}\mu^{2}\{|\nabla\theta|^{2}+2K|\theta|^{2}\}d\mu_{g}\\
=2\int_{M}\mu^{2}|f|^{2}d\mu_{g}+\int_{M}\mu^{2}\textrm{div}_{g}Jd\mu_{g} \notag
\end{align}
Writing 
\begin{equation}
 \mu^{2}\textrm{div}_{g}J=\textrm{div}_{g}(\mu^{2}J)-2\mu(J\cdot d\mu) \notag
\end{equation}
we obtain, since $M$ is compact,
\begin{equation}
 \int_{M}\mu^{2}\textrm{div}_{g}Jd\mu_{g}=-2\int_{M}\mu(J\cdot d\mu)d\mu_{g}
\end{equation}
From (8.134), we can estimate
\begin{equation}
 |J|\leq |\theta|(|\nabla\theta|+|f|)
\end{equation}
therefore:
\begin{align*}
 -2\int_{M}\mu(J\cdot d\mu)d\mu_{g}\leq 2\int_{M}\mu|J||d\mu|d\mu_{g}\\
\leq 2\int_{M}|d\mu||\theta|(\mu|\nabla\theta|+\mu|f|)d\mu_{g}
\end{align*}
Applying the inequalities:
\begin{align*}
 2|d\mu||\theta|\mu|\nabla\theta|\leq\frac{1}{2}\mu^{2}|\nabla\theta|^{2}+2|d\mu|^{2}|\theta|^{2}\\
2|d\mu||\theta|\mu|f|\leq\mu^{2}|f|^{2}+|d\mu|^{2}|\theta|^{2}
\end{align*}
we then obtain:
\begin{align}
 -2\int_{M}\mu(J\cdot d\mu)d\mu_{g}\\
\leq \frac{1}{2}\int_{M}\mu^{2}|\nabla\theta|^{2}d\mu_{g}+\int_{M}\mu^{2}|f|^{2}d\mu_{g}+3\int_{M}|d\mu|^{2}|\theta|^{2}d\mu_{g} \notag
\end{align}
So the lemma follows. $\qed$

Next, we shall derive an estimate for $K$, the Gauss curvature of $S_{t,u}$. Recall the bootstrap assumptions $\textbf{A}, \textbf{E}, \textbf{F}$ of 
Chapter 6.

$\textbf{Lemma 8.9}$ Under the bootstrap assumptions $\textbf{A},\textbf{E},\textbf{F}$, the following estimate holds:
\begin{align*}
|K-r^{-2}|\leq C\delta_{0}(1+t)^{-3}[1+\log(1+t)] 
\end{align*}
In particular, taking $\delta_{0}$ suitably small, we have:
\begin{align*}
 K\geq C^{-1}(1+t)^{-2}
\end{align*}
$Proof$. From the equation (3.27), (3.31) and (3.36) in Chapter 3, we have the following expression for $K$:
\begin{align*}
 K=\frac{1}{2}\eta^{-2}[(\textrm{tr}\chi)^{2}-|\chi|^{2}]-\eta^{-1}(\textrm{tr}\slashed{k}\textrm{tr}\chi-\chi\cdot\slashed{k})+\rho
\end{align*}
From (3.5), (3.33) and (4.10), we know that $K$ is regular as $\mu\longrightarrow 0$. Also from $\textbf{E}$ we know that
\begin{equation}
 |\rho|\leq C\delta_{0}^{2}(1+t)^{-4}
\end{equation}
From $\textbf{A2},\textbf{E}$ and $\textbf{F2}$ we have:
\begin{equation}
 |\eta^{-1}(\textrm{tr}\slashed{k}\textrm{tr}\chi-\chi\cdot\slashed{k})|\leq C\delta_{0}(1+t)^{-3}
\end{equation}
Finally, we consider the first term in the expression of $K$. Assumption $\textbf{F2}$ implies
\begin{equation}
 |\frac{1}{2}[(\textrm{tr}\chi)^{2}-|\chi|^{2}]-(\frac{1}{r})^{2}|\leq C\delta_{0}\frac{[1+\log(1+t)]}{(1+t)^{3}}
\end{equation}
In fact, we can decompose:
\begin{align*}
 |\chi|^{2}=|\hat{\chi}|^{2}+\frac{1}{2}(\textrm{tr}\chi)^{2}
\end{align*}
and
\begin{align*}
 |\frac{1}{4}(\textrm{tr}\chi)^{2}-(\frac{1}{r})^{2}|\leq |\frac{1}{2}\textrm{tr}\chi-\frac{1}{r}|
|\frac{1}{2}\textrm{tr}\chi+\frac{1}{r}|
\end{align*}
From (6.124), (6.129) and $\textbf{F2}$ we have:
\begin{align*}
 |\frac{1}{4}(\textrm{tr}\chi)^{2}-(\frac{1}{r})^{2}|\leq C\delta_{0}(1+t)^{-3}[1+\log(1+t)]
\end{align*}
Also from $\textbf{F2}$, we have:
\begin{align*}
 |\hat{\chi}|^{2}\leq C\delta_{0}(1+t)^{-4}[1+\log(1+t)]^{2}
\end{align*}
So (8.143) results.

Then from $\textbf{E}$, we have:
\begin{equation}
 |\eta^{-2}-1|\leq C\delta_{0}(1+t)^{-1}
\end{equation}
so from (8.143) and (6.124), (6.129), we have:
\begin{equation}
 |\frac{1}{2}\eta^{-2}[(\textrm{tr}\chi)^{2}-|\chi|^{2}]-\frac{1}{r^{2}}|\leq C\delta_{0}(1+t)^{-3}[1+\log(1+t)]
\end{equation}
The lemma then follows. $\qed$

\subsection{Preliminary Estimates for the Solutions of the Propagation Equations}
We now apply Lemma 8.8 to Proposition 8.4, taking $(M,g)$ to be $(S_{t,u},\slashed{g})$, the trace-free symmetric 2-covariant 
tensorfields $\theta$ to be $\leftexp{(i_{1}...i_{l})}{\hat{\chi}}_{l}$, and the 1-form $f$ to be:
\begin{align*}
 \frac{1}{2}\slashed{d}(R_{i_{l}}...R_{i_{1}}\textrm{tr}\chi)+\leftexp{(i_{1}...i_{l})}{i}_{l}
\end{align*}
In view of the fact that by Lemma 8.9, taking $\delta_{0}$ suitably small, we have $K\geq 0$, we obtain:
\begin{align}
 \|\mu\slashed{D}\leftexp{(i_{1}...i_{l})}{\hat{\chi}}_{l}\|_{L^{2}(S_{t,u})}
\leq C\|\mu\slashed{d}(R_{i_{l}}...R_{i_{1}}\textrm{tr}\chi)\|_{L^{2}(S_{t,u})}+C\|\mu\leftexp{(i_{1}...i_{l})}{i}_{l}\|_{L^{2}(S_{t,u})}\\
+C\|\slashed{d}\mu\|_{L^{\infty}(S_{t,u})}\|\leftexp{(i_{1}...i_{l})}{\hat{\chi}}_{l}\|_{L^{2}(S_{t,u})} \notag
\end{align}
Now, by $\textbf{F1}$ we have:
\begin{equation}
 |\slashed{d}\mu|\leq C\delta_{0}(1+t)^{-1}[1+\log(1+t)]
\end{equation}
Also,
\begin{equation}
 \mu\slashed{d}(R_{i_{l}}...R_{i_{1}}\textrm{tr}\chi)=\leftexp{(i_{1}...i_{l})}{x}_{l}-\slashed{d}\leftexp{(i_{1}...i_{l})}{\check{f}}_{l}
\end{equation}
Substituting (8.147) and (8.148) in (8.146) yields:
\begin{align}
 \|\mu\slashed{D}\leftexp{(i_{1}...i_{l})}{\hat{\chi}}_{l}\|_{L^{2}(S_{t,u})}\leq C\|\leftexp{(i_{1}...i_{l})}{x}_{l}\|_{L^{2}(S_{t,u})}\\ \notag
+C\|\slashed{d}\leftexp{(i_{1}...i_{l})}{\check{f}}_{l}\|_{L^{2}(S_{t,u})}+C\|\mu\leftexp{(i_{1}...i_{l})}{i}_{l}\|_{L^{2}(S_{t,u})}\\ \notag
+C\delta_{0}(1+t)^{-1}[1+\log(1+t)]\|\leftexp{(i_{1}...i_{l})}{\hat{\chi}}_{l}\|_{L^{2}(S_{t,u})} \notag
\end{align}

Let us return to the propagation equation of Proposition 8.3. Setting:
\begin{equation}
 \leftexp{(i_{1}...i_{l})}{\tilde{g}}_{l}=\mu\slashed{d}\leftexp{(i_{1}...i_{l})}{h}_{l}+\leftexp{(i_{1}...i_{l})}{g}_{l}
\end{equation}
This propagation equation takes the form:
\begin{align}
 \slashed{\mathcal{L}}_{L}\leftexp{(i_{1}...i_{l})}{x}_{l}+(\textrm{tr}\chi-2\mu^{-1}(L\mu))\leftexp{(i_{1}...i_{l})}{x}_{l}\\
=(\frac{1}{2}\textrm{tr}\chi-2\mu^{-1}(L\mu))\slashed{d}\leftexp{(i_{1}...i_{l})}{\check{f}}_{l}-\leftexp{(i_{1}...i_{l})}{\tilde{g}}_{l}\notag
\end{align}

Recall that $S_{0,0}$ is the unit sphere $S^{2}$ in $\mathbb{R}^{3}$ ($\Sigma_{0}$ is identified with $\mathbb{R}^{3}$). We define a diffeomorphism $\Phi_{t,u}$ of 
$S^{2}$ onto $S_{t,u}$ as follows. First, at $t=0$, $\Phi_{0,u}$ is the diffeomorphism of $S^{2}$ onto $S_{0,u}$ defined by the flow of $T$ on $\Sigma_{0}$.Then 
$\Phi_{t,u}=\Phi_{t}\circ\Phi_{0,u}$ where $\Phi_{t}$ is the flow of $L$ on $W^{*}_{\epsilon_{0}}$. That is, $\Phi_{t|S_{0,u}}$ is the diffeomorphism of 
$S_{0,u}$ onto $S_{t,u}$ defined by the generators of $C_{u}$. Given an arbitrary $S_{t,u}$ 1-form $\xi$, we may pull it back by $\Phi_{t,u}$ to $S^{2}$
and consider the 1-form $\xi(t,u)=\Phi^{*}_{t,u}\xi$ as a 1-form  on $S^{2}$ depending on the parameters $t$ and $u$. If $\xi=\slashed{\mathcal{L}}_{L}\zeta$, where 
$\zeta$ is an arbitrary $S_{t,u}$ 1-form, then $\xi(t,u)=\frac{\partial\zeta(t,u)}{\partial t}$. If $\xi=\slashed{d}\phi$ where $\phi$ is a function, then 
$\xi(t,u)=\slashed{d}\phi(t,u)$, where $\phi(t,u)=\phi\circ\Phi_{t,u}$ is the corresponding function on $S^{2}$ and $\slashed{d}\phi(t,u)$ denotes 
its differential on $S^{2}$. Also, we may pull the induced metric 
$\slashed{g}$ on $S_{t,u}$ back to $S^{2}$ and consider $\slashed{g}(t,u)=\Phi^{*}_{t,u}\slashed{g}$ as a metric on $S^{2}$ depending on $t$ and $u$.
The foregoing correspond simply to the description in terms of acoustical coordinates $(t,u,\vartheta^{1},\vartheta^{2})$, where $(\vartheta^{1},\vartheta^{2})$
are arbitrary local coordinates on $S^{2}$ extended to $\Sigma_{0}$ so that the lines corresponding to constant values of $(\vartheta^{1},\vartheta^{2})$ are 
orthogonal to the $S_{0,u}$ (hence $\Xi=0$ on $\Sigma_{0}$).

In view of above, we may think of (8.151) as an equation for the 1-form $\leftexp{(i_{1}...i_{l})}{x}_{l}(t,u)$ on $S^{2}$ depending on $(t,u)$:
\begin{align}
 \frac{\partial}{\partial t}\leftexp{(i_{1}...i_{l})}{x}_{l}+(\textrm{tr}\chi-2\mu^{-1}(L\mu))\leftexp{(i_{1}...i_{l})}{x}_{l}\\
=(\frac{1}{2}\textrm{tr}\chi-2\mu^{-1}(L\mu))\slashed{d}\leftexp{(i_{1}...i_{l})}{\check{f}}_{l}-\leftexp{(i_{1}...i_{l})}{\tilde{g}}_{l}\notag
\end{align}
where $\leftexp{(i_{1}...i_{l})}{x}_{l}=\leftexp{(i_{1}...i_{l})}{x}_{l}(t,u)$ is a 1-form on $S^{2}$ depending on $(t,u)$.

In the following we denote by $(\quad,\quad)$ the pointwise inner product of tensorfields on $S^{2}$, depending on $t$ and $u$, with respect to the metric 
$\slashed{g}(t,u)$. For 1-forms $\xi(t,u)$ and $\zeta(t,u)$, we have:
\begin{equation}
 (\xi(t,u),\zeta(t,u))=(\slashed{g}^{-1})^{AB}(t,u)\xi_{A}(t,u)\zeta_{B}(t,u)
\end{equation}
We also denote by $|\quad|$ the pointwise magnitude of tensorfields on $S^{2}$, depending on $t$ and $u$, with respect to $\slashed{g}(t,u)$. Thus,
\begin{equation}
 |\xi(t,u)|=\sqrt{(\xi(t,u),\xi(t,u))}
\end{equation}
Now by the definition of $\chi$ we have:
\begin{align}
 \frac{\partial\slashed{g}_{AB}(t,u)}{\partial t}=2\chi_{AB}(t,u), \quad \frac{\partial(\slashed{g}^{-1})^{AB}(t,u)}{\partial t}=-2\chi^{AB}(t,u);\\
\chi^{AB}(t,u)=(\slashed{g}^{-1})^{AC}(t,u)(\slashed{g}^{-1})^{BD}(t,u)\chi_{CD}(t,u)\notag
\end{align}
It follows from (8.153) and (8.155) that:
\begin{equation}
 \frac{\partial}{\partial t}(\xi,\zeta)=(\frac{\partial\xi}{\partial t},\zeta)+(\xi,\frac{\partial\zeta}{\partial t})-2\xi\cdot\chi\cdot\zeta
\end{equation}
Here, 
\begin{equation}
\xi\cdot\chi\cdot\zeta=\xi_{A}\chi^{AB}\zeta_{B}       
\end{equation}
In particular, taking $\zeta=\xi$ we obtain, in view of (8.154),
\begin{equation}
 |\xi|\frac{\partial}{\partial t}|\xi|=\frac{1}{2}\frac{\partial}{\partial t}(\xi,\xi)=(\frac{\partial\xi}{\partial t},\xi)-\xi\cdot\chi\cdot\xi
\end{equation}
Decomposing:
\begin{equation}
 \chi_{AB}=\hat{\chi}_{AB}+\frac{1}{2}\slashed{g}_{AB}\textrm{tr}\chi,\quad \chi^{AB}=\hat{\chi}^{AB}+\frac{1}{2}(\slashed{g}^{-1})^{AB}\textrm{tr}\chi
\end{equation}
then we have:
\begin{equation}
 \xi\cdot\chi\cdot\zeta=\xi\cdot\hat{\chi}\cdot\zeta+\frac{1}{2}\textrm{tr}\chi(\xi,\zeta)
\end{equation}
hence (8.158) takes the form:
\begin{equation}
 |\xi|\frac{\partial}{\partial t}|\xi|=(\frac{\partial\xi}{\partial t},\xi)-\xi\cdot\hat{\chi}\cdot\xi-\frac{1}{2}\textrm{tr}\chi|\xi|^{2}
\end{equation}
We shall apply (8.161) taking $\xi=\leftexp{(i_{1}...i_{l})}{x}_{l}$. From (8.152) we have:
\begin{align}
 (\frac{\partial \leftexp{(i_{1}...i_{l})}{x}_{l}}{\partial t},\leftexp{(i_{1}...i_{l})}{x}_{l})
-\frac{1}{2}\textrm{tr}\chi|\leftexp{(i_{1}...i_{l})}{x}_{l}|^{2}=\\ \notag
-(-2\mu^{-1}\frac{\partial\mu}{\partial t}+\frac{3}{2}\textrm{tr}\chi)|\leftexp{(i_{1}...i_{l})}{x}_{l}|^{2}
+(-2\mu^{-1}\frac{\partial\mu}{\partial t}+\frac{1}{2}\textrm{tr}\chi)(\leftexp{(i_{1}...i_{l})}{x}_{l},
\slashed{d}\leftexp{(i_{1}...i_{l})}{\check{f}}_{l})\\ \notag-(\leftexp{(i_{1}...i_{l})}{x}_{l},\leftexp{(i_{1}...i_{l})}{\tilde{g}}_{l})
\end{align}

We shall use the following assumption:

In $W^{s}_{\epsilon_{0}}$ we have:
\begin{equation}
 \textbf{AS}: -2\mu^{-1}\frac{\partial\mu}{\partial t}+\frac{1}{2}\textrm{tr}\chi\geq 0\notag
\end{equation}
From this assumption and (8.162) we have:
\begin{align}
 (\frac{\partial \leftexp{(i_{1}...i_{l})}{x}_{l}}{\partial t},\leftexp{(i_{1}...i_{l})}{x}_{l})
-\frac{1}{2}\textrm{tr}\chi|\leftexp{(i_{1}...i_{l})}{x}_{l}|^{2}\\ \notag
\leq -(-2\mu^{-1}\frac{\partial\mu}{\partial t}+\frac{3}{2}\textrm{tr}\chi)|\leftexp{(i_{1}...i_{l})}{x}_{l}|^{2}
+(-2\mu^{-1}\frac{\partial\mu}{\partial t}+\frac{1}{2}\textrm{tr}\chi)|\leftexp{(i_{1}...i_{l})}{x}_{l}|
|\slashed{d}\leftexp{(i_{1}...i_{l})}{\check{f}}_{l}|\\ \notag
+|\leftexp{(i_{1}...i_{l})}{x}_{l}||\leftexp{(i_{1}...i_{l})}{\tilde{g}}_{l}|
\end{align}
In reference to the second term on the right in (8.161) with $\leftexp{(i_{1}...i_{l})}{x}_{l}$ in the role of $\xi$, we have
\begin{equation}
 |\leftexp{(i_{1}...i_{l})}{x}_{l}\cdot\hat{\chi}\cdot\leftexp{(i_{1}...i_{l})}{x}_{l}|\leq |\hat{\chi}||\leftexp{(i_{1}...i_{l})}{x}_{l}|^{2}
\end{equation}
Substituting (8.164) in (8.163) we obtain, in view of (8.161),
\begin{align}
 \frac{\partial|\leftexp{(i_{1}...i_{l})}{x}_{l}|}{\partial t}\leq 
-(-2\mu^{-1}\frac{\partial\mu}{\partial t}+\frac{3}{2}\textrm{tr}\chi-|\hat{\chi}|)|\leftexp{(i_{1}...i_{l})}{x}_{l}|\\ \notag
+(-2\mu^{-1}\frac{\partial\mu}{\partial t}+\frac{1}{2}\textrm{tr}\chi)|\slashed{d}\leftexp{(i_{1}...i_{l})}{\check{f}}_{l}|
+|\leftexp{(i_{1}...i_{l})}{\tilde{g}}_{l}|
\end{align}
The integrating factor here is:
\begin{align}
 \exp\{\int_{0}^{t}(-2\mu^{-1}\frac{\partial\mu}{\partial t}+\frac{3}{2}\textrm{tr}\chi-|\hat{\chi}|)(t',u)dt'\}=
(\frac{\mu(t,u)}{\mu(0,u)})^{-2}(A(t,u))^{3/2}e^{-S(t,u)}
\end{align}
Here
\begin{equation}
 A(t,u)=\exp(\int_{0}^{t}\textrm{tr}\chi(t',u)dt'), \quad S(t,u)=\int_{0}^{t}|\hat{\chi}|(t',u)dt'
\end{equation}
Since by (8.155):
\begin{equation}
 \textrm{tr}\chi=\frac{1}{\sqrt{\det\slashed{g}}}\frac{\partial}{\partial t}\sqrt{\det\slashed{g}}
\end{equation}
we have:
\begin{equation}
 A(t,u)=\frac{\sqrt{\det\slashed{g}(t,u)}}{\sqrt{\det\slashed{g}(0,u)}}
\end{equation}
that is, $A(t,u)$ is the ratio of the area elements of $S_{t,u}$ and $S_{0,u}$ at corresponding points along the same generator of $C_{u}$.
Integrating (8.165) and using (8.166) we get:
\begin{equation}
 |\leftexp{(i_{1}...i_{l})}{x}_{l}(t,u)|\leq \leftexp{(i_{1}...i_{l})}{F}_{l}(t,u)+\leftexp{(i_{1}...i_{l})}{G}_{l}(t,u)
\end{equation}
where:
\begin{align}
 \leftexp{(i_{1}...i_{l})}{F}_{l}(t,u)=e^{S(t,u)}(A(t,u))^{-3/2}(\mu(t,u))^{2}\{
(\mu(0,u))^{-2}|\leftexp{(i_{1}...i_{l})}{x}_{l}(0,u)|\\+\int_{0}^{t}(\mu(t',u))^{-2}(A(t',u))^{3/2}e^{-S(t',u)}
(-2\mu^{-1}\frac{\partial\mu}{\partial t}+\frac{1}{2}\textrm{tr}\chi)(t',u)|
\slashed{d}\leftexp{(i_{1}...i_{l})}{\check{f}}_{l}(t',u)|\}\notag
\end{align}
and:
\begin{align}
 \leftexp{(i_{1}...i_{l})}{G}_{l}(t,u)=e^{S(t,u)}(A(t,u))^{-3/2}(\mu(t,u))^{2}\\
\cdot\int_{0}^{t}(\mu(t',u))^{-2}(A(t',u))^{3/2}e^{-S(t',u)}|\leftexp{(i_{1}...i_{l})}{\tilde{g}}_{l}(t',u)|dt'\notag
\end{align}

We begin our estimates with $A(t,u)$ and $S(t,u)$. By $\textbf{F2}$ and (8.167) we have:
\begin{align*}
 |\log A(t,u)-\int_{0}^{t}\frac{2}{1-u+t'}dt'|\leq C\delta_{0}\int_{0}^{t}(1+t')^{-2}[1+\log(1+t')]dt'\leq C\delta_{0}
\end{align*}
Since
\begin{align*}
 \int_{0}^{t}\frac{2}{1-u+t'}dt'=2\log(\frac{1-u+t}{1-u})
\end{align*}
it follows that:
\begin{align}
 e^{-C\delta_{0}}(\frac{1-u+t}{1-u})^{2}\leq A(t,u)\leq e^{C\delta_{0}}(\frac{1-u+t}{1-u})^{2}
\end{align}
Also by $\textbf{F2}$ and (8.167) we have:
\begin{equation}
 S(t,u)\leq C\delta_{0}
\end{equation}
the integral
\begin{equation}
 \int_{0}^{\infty}(1+t')^{-2}[1+\log(1+t')]dt'\notag
\end{equation}
being convergent.

In view of (8.173) and (8.174) we obtain from (8.171) and (8.172):
\begin{align}
 \leftexp{(i_{1}...i_{l})}{F}_{l}(t,u)\leq e^{C\delta_{0}}(1-u+t)^{-3}\{\leftexp{(i_{1}...i_{l})}{M}^{0}_{l}(t,u)
+\leftexp{(i_{1}...i_{l})}{M}^{1}_{l}(t,u)+\leftexp{(i_{1}...i_{l})}{M}^{2}_{l}(t,u)\} 
\end{align}
where:
\begin{align}
 \leftexp{(i_{1}...i_{l})}{M}^{0}_{l}(t,u)=
(\frac{\mu(t,u)}{\mu(0,u)})^{2}(1-u)^{3}|\leftexp{(i_{1}...i_{l})}{x}_{l}(0,u)|
\end{align}
\begin{align}
 \leftexp{(i_{1}...i_{l})}{M}^{1}_{l}(t,u)=
\int_{0}^{t}(\frac{\mu(t,u)}{\mu(t',u)})^{2}(1-u+t')^{3}[-2\mu^{-1}(\frac{\partial\mu}{\partial t})_{-}(t',u)]\cdot
|\slashed{d}\leftexp{(i_{1}...i_{l})}{\check{f}}_{l}(t',u)|dt'
\end{align}
\begin{align}
 \leftexp{(i_{1}...i_{l})}{M}^{2}_{l}(t,u)=
\int_{0}^{t}(\frac{\mu(t,u)}{\mu(t',u)})^{2}(1-u+t')^{3}(\frac{1}{2}\textrm{tr}\chi(t',u))
|\slashed{d}\leftexp{(i_{1}...i_{l})}{\check{f}}_{l}(t',u)|dt'
\end{align}
Also:
\begin{align}
 \leftexp{(i_{1}...i_{l})}{G}_{l}(t,u)
\leq e^{C\delta_{0}}(1-u+t)^{-3}\cdot\int_{0}^{t}(\frac{\mu(t,u)}{\mu(t',u)})^{2}(1-u+t')^{3}|\leftexp{(i_{1}...i_{l})}
{\tilde{g}}_{l}(t',u)|dt'
\end{align}
To proceed we must analyze the behavior of $\mu$. We need the following assumption:
\begin{equation}
 \slashed{\textbf{E}}\textbf{3}_{0}:|\slashed{\Delta}\psi_{0}|\leq C\delta_{0}(1+t)^{-3}\notag
\end{equation}

\section{Crucial Lemmas Concerning the Behavior of $\mu$}
$\textbf{Lemma 8.10}$ Under the assumptions $\textbf{E1},\textbf{E2}, \slashed{\textbf{E}}\textbf{3}_{0}, 
\textbf{F2}$ and $\textbf{A}$, the following hold. Let us denote:
 \begin{align*}
  P_{s}(u,\vartheta)=(1+s)(\underline{L}\psi_{0})(s,u,\vartheta)
 \end{align*}
We then have, for all $t\in[0,s]$:
\begin{align*}
 (\underline{L}\psi_{0})(t,u,\vartheta)=\frac{P_{s}(u,\vartheta)}{(1+t)}+R_{s}(t,u,\vartheta)
\end{align*}
and:
\begin{align*}
 |R_{s}(t,u,\vartheta)|\leq C\delta_{0}\frac{[1+\log(1+t)]}{(1+t)^{2}}\frac{(s-t)}{(1+s)}
\end{align*}
where $C$ is a constant independent of $s$.

$Proof$. The function $\psi_{0}$ satisfies the equation:
\begin{equation}
 \Box_{\tilde{g}}\psi_{0}=0\notag
\end{equation}
From the expression for $\Box_{\tilde{g}}$ in Chapter 3, this reads:
\begin{equation}
 L(\underline{L}\psi_{0})+\nu(\underline{L}\psi_{0})=\rho_{0}
\end{equation}
where
\begin{equation}
 \rho_{0}=\mu\slashed{\Delta}\psi_{0}-\underline{\nu}L\psi_{0}-2\zeta\cdot\slashed{d}\psi_{0}+\mu\frac{d\log\Omega}{dh}\slashed{d}h
\cdot\slashed{d}\psi_{0}
\end{equation}
Now, by $\slashed{\textbf{E}}\textbf{3}_{0}$, and $\textbf{A3}$, the first term on the right hand side of (8.181) is bounded in absolute value by:
\begin{align*}
 C\delta_{0}(1+t)^{-3}[1+\log(1+t)]
\end{align*}
By (6.122) and $\textbf{E2}$, the second term on the right hand side of (8.181) is bounded in absolute value by:
\begin{align*}
  C\delta_{0}(1+t)^{-3}[1+\log(1+t)]
\end{align*}
By (6.104) and $\textbf{E2}$, the third term on the right hand side of (8.181) is bounded in absolute value by:
\begin{align*}
 C\delta_{0}(1+t)^{-4}[1+\log(1+t)]
\end{align*}
Finally, by $\textbf{A3}$ and $\textbf{E}$, the last term on the right of (8.181) is bounded in absolute value by:
\begin{align*}
 C\delta_{0}(1+t)^{-4}[1+\log(1+t)]
\end{align*}
We conclude:
\begin{align}
 |\rho_{0}|\leq C\delta_{0}(1+t)^{-3}[1+\log(1+t)]
\end{align}

Consider on $S^{2}$ the function:
\begin{equation}
 \tilde{A}(t,u)=A(t,u)(\frac{\Omega(t,u)}{\Omega(0,u)})
\end{equation}
depending on the parameters $t$ and $u$, where $A(t,u)$ is defined by (8.169). We have, using (8.167) and the definition of $\nu$,
\begin{equation}
 \frac{\partial\tilde{A}(t,u)}{\partial t}=2\nu(t,u)\tilde{A}(t,u), \quad \tilde{A}(0,u)=1
\end{equation}
Thus, setting, in terms of the pullbacks by $\Phi_{t,u}$ to $S^{2}$,
\begin{equation}
 \tau_{0}(t,u)=(\tilde{A}(t,u))^{1/2}(\underline{L}\psi_{0})(t,u),\quad
\tilde{\rho_{0}}(t,u)=(\tilde{A}(t,u))^{1/2}\rho_{0}(t,u)
\end{equation}
(8.180) becomes:
\begin{equation}
 \frac{\partial \tau_{0}}{\partial t}=\tilde{\rho}_{0}
\end{equation}
Integrating (8.184) we obtain:
\begin{equation}
 \tilde{A}(t,u)=e^{2N(t,u)}, \quad N(t,u)=\int_{0}^{t}\nu(t',u)dt'
\end{equation}
Writing,
\begin{equation}
 \nu(t,u)=\frac{1}{1-u+t}+\hat{\nu}
\end{equation}
we have,
\begin{equation}
 N(t,u)=\log(\frac{1-u+t}{1-u})+\hat{N}(t,u), \quad\hat{N}(t,u)=\int_{0}^{t}\hat{\nu}(t',u)dt'
\end{equation}
hence:
\begin{equation}
 \tilde{A}(t,u)=(\frac{1-u+t}{1-u})^{2}e^{2\hat{N}(t,u)}
\end{equation}
Now, by $\textbf{E}$ and $\textbf{F2}$:
\begin{equation}
 |\hat{\nu}(t,u)|\leq C\delta_{0}(1+t)^{-2}[1+\log(1+t)]
\end{equation}
It follows that:
\begin{equation}
 |\hat{N}(s,u)-\hat{N}(t,u)|\leq \int_{t}^{s}|\hat{\nu}(t',u)|dt'\leq C\delta_{0}\int_{t}^{s}(1+t')^{-2}[1+\log(1+t')]dt'
\end{equation}
The last integral is:
\begin{align}
 \int_{t}^{s}(1+t')^{-2}[1+\log(1+t')]dt'\\ \notag
=2[\frac{1}{(1+t)}-\frac{1}{(1+s)}]+[\frac{\log(1+t)}{(1+t)}-\frac{\log(1+s)}{(1+s)}]\\ \notag
\leq 2\frac{[1+\log(1+t)]}{(1+t)}\frac{(s-t)}{(1+s)}
\end{align}
Thus, we obtain:
\begin{equation}
 |\hat{N}(s,u)-\hat{N}(t,u)|\leq C\delta_{0}\frac{[1+\log(1+t)]}{(1+t)}\frac{(s-t)}{(1+s)}
\end{equation}
This implies in particular, replacing $s,t$ by $t,0$ respectively and noting that $N(0,u)=0$,
\begin{equation}
 |\hat{N}(t,u)|\leq C\delta_{0}
\end{equation}

We now return to (8.186). Integrating on $[t,s]$ we obtain:
\begin{equation}
 \tau_{0}(t,u)=\tau_{0}(s,u)-\int_{t}^{s}\tilde{\rho}_{0}(t',u)dt'
\end{equation}
By (8.182), (8.185), (8.190) and (8.195) we have:
\begin{equation}
 |\tilde{\rho}_{0}(t,u)|\leq C\delta_{0}(1+t)^{-2}[1+\log(1+t)]
\end{equation}
Hence, in view of (8.193),
\begin{equation}
 |\int_{t}^{s}\tilde{\rho}_{0}(t',u)dt'|\leq C\delta_{0}\frac{[1+\log(1+t)]}{(1+t)}\frac{(s-t)}{(1+s)}
\end{equation}
From (8.185) and (8.190) we have:
\begin{equation}
 (\underline{L}\psi_{0})(t,u)=(\tilde{A}(t,u))^{-1/2}\tau_{0}(t,u)=(\frac{1-u}{1-u+t})e^{-\hat{N}(t,u)}\tau_{0}(t,u)
\end{equation}
In particular, this holds at $t=s$. On the other hand, according to the definition of $P_{s}(u)$ in the statement of Lemma 8.10,
\begin{align*}
 (\underline{L}\psi_{0})(s,u)=\frac{P_{s}(u)}{(1+s)}
\end{align*}
it follows that:
\begin{equation}
 \tau_{0}(s,u)=\frac{P_{s}(u)}{(1-u)}\frac{(1-u+s)}{(1+s)}e^{\hat{N}(s,u)}
\end{equation}
Substituting (8.200) in (8.196) and the result in (8.199) we get:
\begin{align}
 (\underline{L}\psi_{0})(t,u)=\frac{P_{s}(u)}{(1+s)}\frac{(1-u+s)}{(1-u+t)}e^{\hat{N}(s,u)-\hat{N}(t,u)}\\
-\frac{(1-u)}{(1-u+t)}e^{-\hat{N}(t,u)}\int_{t}^{s}\tilde{\rho}_{0}(t',u)dt'\notag
\end{align}
Setting:
\begin{align}
 R_{s}(t,u)=P_{s}(u)B_{s}(t,u)-\frac{(1-u)}{(1-u+t)}e^{-\hat{N}(t,u)}\int_{t}^{s}\tilde{\rho}_{0}(t',u)dt'
\end{align}
where
\begin{align}
 B_{s}(t,u)=\frac{1}{(1+s)}\frac{(1-u+s)}{(1-u+t)}e^{\hat{N}(s,u)-\hat{N}(t,u)}-\frac{1}{(1+t)}
\end{align}
We obtain:
\begin{equation}
 (\underline{L}\psi_{0})(t,u)=\frac{P_{s}(u)}{(1+t)}+R_{s}(t,u)
\end{equation}
Moreover, (8.203) can be written as:
\begin{align}
 B_{s}(t,u)=\frac{1}{(1+t)}\{e^{\hat{N}(s,u)-\hat{N}(t,u)}-1-\frac{-u}{(1-u+t)}\frac{(s-t)}{(1+s)}e^{\hat{N}(s,u)-\hat{N}(t,u)}\}
\end{align}
Since from (8.194),
\begin{align*}
 e^{|\hat{N}(s,u)-\hat{N}(t,u)|}\leq e^{C\delta_{0}},\\
|e^{\hat{N}(s,u)-\hat{N}(t,u)}-1|\leq e^{C\delta_{0}}|\hat{N}(s,u)-\hat{N}(t,u)|\leq C\delta_{0}
\frac{[1+\log(1+t)]}{(1+t)}\frac{(s-t)}{(1+s)}
\end{align*}
it follows that:
\begin{align*}
 |B_{s}(t,u)|\leq C\frac{[1+\log(1+t)]}{(1+t)^{2}}\frac{(s-t)}{(1+s)}
\end{align*}
Since by $\textbf{E2}$, we also have:
\begin{equation}
 |P_{s}(u)|\leq C\delta_{0}
\end{equation}
it follows from (8.202), in view of (8.198), that:
\begin{equation}
 |R_{s}(t,u)|\leq C\delta_{0}\frac{[1+\log(1+t)]}{(1+t)^{2}}\frac{(s-t)}{(1+s)}
\end{equation}
and the lemma is proved. $\qed$

We now introduce the additional bootstrap assumptions:

There is a positive constant $C$ independent of $s$ such that in $W^{s}_{\epsilon_{0}}$,
\begin{align*}
 \textbf{E}_{\textbf{LT}}\textbf{3}: |LT\psi_{\mu}|\leq C\delta_{0}(1+t)^{-2}\\
\textbf{E}_{\textbf{LL}}\textbf{3}: |LL\psi_{\mu}|\leq C\delta_{0}(1+t)^{-3}
\end{align*}
$\textbf{Proposition 8.5}$ Let the assumptions of Lemma 8.10 hold. Let the additional bootstrap assumptions $\textbf{E}_{\textbf{LT}}\textbf{3}, 
\textbf{E}_{\textbf{LL}}\textbf{3}$, hold as well. Denoting:
\begin{align*}
 E_{s}(u,\vartheta)=(1+s)(\frac{\partial\mu}{\partial t})(s,u,\vartheta)
\end{align*}
we then have, for all $t\in[0,s]$:
\begin{align*}
 (\frac{\partial\mu}{\partial t})(t,u,\vartheta)=\frac{E_{s}(u,\vartheta)}{(1+t)}+Q_{1,s}(t,u,\vartheta)
\end{align*}
and
\begin{align*}
 |Q_{1,s}(t,u,\vartheta)|\leq C\delta_{0}\frac{[1+\log(1+t)]}{(1+t)^{2}}\frac{(s-t)}{(1+s)}
\end{align*}
$Proof$. We begin with the equation for $\mu$ (see (3.92)). We have:
\begin{equation}
 m=\frac{1}{2}\frac{dH}{dh}(Th)
\end{equation}
Setting:
\begin{equation}
 m_{0}=\frac{1}{2}\ell(T\psi_{0}),\quad m_{1}=m-m_{0}
\end{equation}
where $\ell$ is the constant:
\begin{equation}
 \ell=\frac{dH}{dh}(0)
\end{equation}
we have:
\begin{equation}
m_{1}=\frac{1}{2}(\frac{dH}{dh}-\ell)(Th)+\frac{1}{2}\ell[(Th)-(T\psi_{0})]
\end{equation}
From $\textbf{E1}$:
\begin{equation}
 |h|\leq C\delta_{0}(1+t)^{-1}
\end{equation}
which implies:
\begin{equation}
 |\frac{dH}{dh}-\ell|\leq C\delta_{0}(1+t)^{-1}
\end{equation}
Moreover, since 
\begin{align*}
 Th=T\psi_{0}-\sum_{i}\psi_{i}T\psi_{i}
\end{align*}
$\textbf{E}$ implies:
\begin{equation}
 |m_{1}|\leq C\delta_{0}(1+t)^{-2}
\end{equation}
Also using $\textbf{E}_{\textbf{LT}}\textbf{3}$ we deduce:
\begin{equation}
 |Lm_{1}|\leq C\delta_{0}(1+t)^{-3}
\end{equation}
Now we set:
\begin{equation}
 m_{0,1}=\frac{1}{4}\ell(\underline{L}\psi_{0})
\end{equation}
Then, since $\underline{L}=2T+\alpha^{-2}\mu L$, we have:
\begin{equation}
 m_{0,1}=m_{0}+\frac{1}{4}\ell\alpha^{-2}\mu(L\psi_{0})
\end{equation}
Thus, setting also:
\begin{equation}
 e_{1}=e-\frac{1}{4}\ell\alpha^{-2}(L\psi_{0}),\quad n_{1}=m_{1}+\mu e_{1}
\end{equation}
we get:
\begin{equation}
 L\mu=m_{0,1}+n_{1}
\end{equation}
From (8.214), $\textbf{A}$ and $\textbf{E}$, we have:
\begin{equation}
 |e_{1}|\leq C\delta_{0}(1+t)^{-2},\quad |n_{1}|\leq C\delta_{0}(1+t)^{-2}[1+\log(1+t)]
\end{equation}
Moreover, by $\textbf{E}_{\textbf{LL}}\textbf{3}$,
\begin{equation}
 |Le_{1}|\leq C\delta_{0}(1+t)^{-3}
\end{equation}
which, together with (8.215) and the fact that:
\begin{equation}
 |L\mu|\leq C\delta_{0}(1+t)^{-1}
\end{equation}
yields:
\begin{equation}
 |Ln_{1}|\leq C\delta_{0}(1+t)^{-3}[1+\log(1+t)]
\end{equation}
By virtue of Lemma 8.10 we have:
\begin{equation}
 m_{0,1}(t,u)=\frac{1}{4}\ell\frac{P_{s}(u)}{(1+t)}+\frac{1}{4}\ell R_{s}(t,u)
\end{equation}
We define the function $E_{s}(u)$ on $S^{2}$, depending on the parameters $u$ and $s$:
\begin{equation}
 E_{s}(u)=\frac{1}{4}\ell P_{s}(u)+(1+s)n_{1}(s,u)
\end{equation}
Then by the second of (8.220),
\begin{equation}
 |E_{s}(u)-\frac{1}{4}\ell P_{s}(u)|\leq C\delta_{0}(1+s)^{-1}[1+\log(1+s)]
\end{equation}
From (8.219), (8.224) and the fact that $R_{s}(s,u)=0$, we have:
\begin{equation}
 (1+s)(\frac{\partial\mu}{\partial t})(s,u)=E_{s}(u)
\end{equation}
Moreover we have:
\begin{align*}
 (\frac{\partial\mu}{\partial t}(t,u)=\frac{E_{s}(u)}{(1+t)}+Q_{1,s}(t,u)
\end{align*}

\begin{equation}
 Q_{1,s}(t,u)=\frac{1}{4}\ell R_{s}(t,u)-n_{1}(s,u)\frac{(s-t)}{(1+t)}-(n_{1}(s,u)-n_{1}(t,u))
\end{equation}
 From Lemma 8.10, the first term on the right of above satisfies the required bound. From (8.220) and the fact that
$f(t)=(1+t)^{-1}[1+\log(1+t)]$ is decreasing in $t$, the second term on the right of above also satisfies the 
required bound. For the last term, we have, by (8.223),
\begin{equation}
 |n_{1}(s,u)-n_{1}(t,u)|\leq \int_{t}^{s}|\frac{\partial n_{1}}{\partial t}(t',u)|dt'\leq C\delta_{0}\int_{t}^{s}
(1+t')^{-3}[1+\log(1+t')]dt'
\end{equation}
and
\begin{align}
 \int_{t}^{s}(1+t')^{-3}[1+\log(1+t')]dt'\\ \notag
=\frac{3}{4}[\frac{1}{(1+t)^{2}}-\frac{1}{(1+s)^{2}}]+\frac{1}{2}[\frac{\log(1+t)}{(1+t)^{2}}-\frac{\log(1+s)}{(1+s)^{2}}]\\ \notag
\leq \frac{3}{4}\frac{[1+\log(1+t)]}{(1+t)^{2}}\frac{(s-t)}{(1+s)}\notag
\end{align}
We conclude that indeed
\begin{equation}
 |Q_{1,s}(t,u)|\leq C\delta_{0}\frac{[1+\log(1+t)]}{(1+t)^{2}}\frac{(s-t)}{(1+s)}
\end{equation}
and the proposition is proved. $\qed$

From Proposition 8.5, we have:
\begin{equation}
 \mu(t,u)=\mu(0,u)+E_{s}(u)\log(1+t)+\int_{0}^{t}Q_{1,s}(t',u)dt'
\end{equation}
Recall that $\mu=\eta\kappa$. Assumption $\textbf{E1}$ restricted to $\Sigma_{0}^{\epsilon_{0}}$ implies:
\begin{equation}
 |\eta(0,u)-1|\leq C\delta_{0}
\end{equation}
In view of the fact that
\begin{equation}
\kappa=1\quad:\quad\textrm{on}\quad\Sigma_{0}
\end{equation}
we obtain:
\begin{equation}
 |\mu(0,u)-1|\leq C\delta_{0}
\end{equation}
According the bound for $Q_{1,s}(t,u)$ in Proposition 8.5,
\begin{equation}
 |Q_{1,s}(t,u)|\leq C\delta_{0}(1+t)^{-2}[1+\log(1+t)]
\end{equation}
Let us define on $S^{2}$ the function:
\begin{equation}
 \mu_{1,s}(u)=\int_{0}^{s}Q_{1,s}(t,u)dt
\end{equation}
Substituting (8.236) we obtain:
\begin{equation}
 |\mu_{1,s}(u)|\leq C\delta_{0}
\end{equation}
We then define on $S^{2}$ the function:
\begin{equation}
 \mu_{[1],s}(u)=\mu(0,u)+\mu_{1,s}(u)
\end{equation}
From (8.235) and (8.238) we have:
\begin{equation}
 |\mu_{[1],s}(u)-1|\leq C\delta_{0}
\end{equation}
Taking $\delta_{0}$ suitably small, this implies:
\begin{equation}
 \inf_{u\in[0,\epsilon_{0}]}\inf_{S^{2}}\mu_{[1],s}(u)\geq \frac{1}{2}
\end{equation}
In view of (8.237) and (8.239), (8.232) can be written as
\begin{equation}
 \mu(t,u)=\mu_{[1],s}(u)+E_{s}(u)\log(1+t)+Q_{0,s}(t,u)
\end{equation}
where
\begin{equation}
 Q_{0,s}(t,u)=-\int_{t}^{s}Q_{1,s}(t',u)dt'
\end{equation}
We have:
\begin{align}
 |Q_{0,s}(t,u)|\leq \int_{t}^{s}|Q_{1,s}(t',u)|dt'
\leq C\delta_{0}\int_{t}^{s}(1+t')^{-2}[1+\log(1+t')]dt'\\
\leq 2C\delta_{0}\frac{[1+\log(1+t)]}{(1+t)}\frac{(s-t)}{(1+s)}
\end{align}
(see (8.193))

Finally, setting
\begin{align}
 \hat{\mu}_{s}(t,u)=\frac{\mu(t,u)}{\mu_{[1],s}(u)}, \quad \hat{E}_{s}(u)=\frac{E_{s}(u)}{\mu_{[1],s}(u)}\\ \notag
\hat{Q}_{0,s}(t,u)=\frac{Q_{0,s}(t,u)}{\mu_{[1],s}(u)}, \quad \hat{Q}_{1,s}(t,u)=\frac{Q_{1,s}(t,u)}{\mu_{[1],s}(u)},
\end{align}
We collect the above results in the following proposition:

$\textbf{Proposition 8.6}$ Let the assumptions of Proposition 8.5 hold. Then, there is a
function $\mu_{[1],s}(u,\vartheta)$ on $S^{2}$ depending on the parameters $u$ and $s$, satisfying:
\begin{equation}
 |\mu_{[1],s}(u,\vartheta)-1|\leq C\delta_{0}, \quad \inf_{u\in[0,\epsilon_{0}]}\inf_{\vartheta\in S^{2}}\mu_{[1],s}(u,\vartheta)\geq 
\frac{1}{2}\notag
\end{equation}
such that setting:
\begin{equation}
 \hat{\mu}_{s}(t,u,\vartheta)=\frac{\mu(t,u,\vartheta)}{\mu_{[1],s}(u,\vartheta)}\notag
\end{equation}
we have:
\begin{align*}
 \hat{\mu}_{s}(t,u,\vartheta)=1+\hat{E}_{s}(u,\vartheta)\log(1+t)+\hat{Q}_{0,s}(t,u,\vartheta)\\
\frac{\partial \hat{\mu}_{s}}{\partial t}(t,u,\vartheta)=\frac{\hat{E}_{s}(u,\vartheta)}{(1+t)}+\hat{Q}_{1,s}(t,u,\vartheta)
\end{align*}
where:
\begin{align*}
 \hat{Q}_{0,s}(t,u,\vartheta)=-\int_{t}^{s}\hat{Q}_{1,s}(t',u,\vartheta)dt'
\end{align*}
and the functions $\hat{Q}_{0,s}, \hat{Q}_{1,s}$ satisfy the bounds:
\begin{align*}
 |\hat{Q}_{0,s}(t,u,\vartheta)|\leq C\delta_{0}\frac{[1+\log(1+t)]}{(1+t)}\frac{(s-t)}{(1+s)}\\
|\hat{Q}_{1,s}(t,u,\vartheta)|\leq C\delta_{0}\frac{[1+\log(1+t)]}{(1+t)^{2}}\frac{(s-t)}{(1+s)}
\end{align*}

Let us define:
\begin{equation}
 \mu_{m}(t)=\min_{\Sigma^{\epsilon_{0}}_{t}}\mu =\min_{(u,\vartheta)\in[0,\epsilon_{0}]\times S^{2}}\mu(t,u,\vartheta)
\end{equation}
We set:
\begin{equation}
 \bar{\mu}_{m}(t)=\min\{\mu_{m}(t),1\}
\end{equation}
Let us also define:
\begin{equation}
 M(t)=\max_{\Sigma^{\epsilon_{0}}_{t}}\{-\mu^{-1}(L\mu)_{-}\}=
\max_{(u,\vartheta)\in[0,\epsilon_{0}]\times S^{2}}\{-\mu^{-1}(\frac{\partial\mu}{\partial t})_{-}(t,u,\vartheta)\}
\end{equation}
The following lemma plays a crucial role in the sequel.

$\textbf{Lemma 8.11}$ Let the assumptions of Proposition 8.6 hold on $W^{s}_{\epsilon_{0}}$. Then for any constant $a\geq 4$, there
is a positive constant $C$ independent of $s$ and $a$ such that for all $t\in[0,s]$ we have:
\begin{equation}
 I_{a}(t):=\int_{0}^{t}\bar{\mu}_{m}^{-a}(t')M(t')dt'\leq Ca^{-1}\bar{\mu}_{m}^{-a}(t)\notag
\end{equation}
provided $\delta_{0}$ is suitably small depending on $a$.

$Proof$. We can express $M(t)$ as:
\begin{equation}
 M(t)=\max_{(u,\vartheta)\in[0,\epsilon_{0}]\times S^{2}}\{-\frac{1}{\hat{\mu}_{s}}(\frac{\partial\hat{\mu}_{s}}{\partial t})_{-}(t,u,\vartheta)\}
\end{equation}
Let us denote:
\begin{equation}
 \hat{E}_{s,m}=\min_{(u,\vartheta)\in[0,\epsilon_{0}]\times S^{2}}\hat{E}_{s}(u,\vartheta)
\end{equation}
Case 1) $\hat{E}_{s,m}\geq 0$. In this case,
\begin{equation}
 \hat{E}_{s}(u,\vartheta)\geq 0: \forall(u,\vartheta)\in[0,\epsilon_{0}]\times S^{2}
\end{equation}
From Proposition 8.6 we then have:
\begin{align*}
 \frac{\partial\hat{\mu}_{s}}{\partial t}(t,u,\vartheta)\geq\hat{Q}_{1,s}(t,u,\vartheta)
\end{align*}
hence:
\begin{equation}
 -(\frac{\partial\hat{\mu}_{s}}{\partial t})_{-}(t,u,\vartheta)\leq -(\hat{Q}_{1,s})_{-}(t,u,\vartheta)\leq C\delta_{0}
\frac{[1+\log(1+t)]}{(1+t)^{2}}
\end{equation}
Also:
\begin{equation}
 \hat{\mu}_{s}(t,u,\vartheta)\geq 1+\hat{Q}_{0,s}(t,u,\vartheta)\geq 1-C\delta_{0}
\end{equation}
Now (8.253) and (8.254) together imply:
\begin{equation}
 M(t)\leq C\delta_{0}\frac{[1+\log(1+t)]}{(1+t)^{2}}
\end{equation}
provided $\delta_{0}$ is suitably small.

On the other hand, we have:
\begin{align*}
 \mu(t,u,\vartheta)=\mu_{[1],s}(u,\vartheta)\hat{\mu}_{s}(t,u,\vartheta)
\end{align*}
From Proposition 8.6,
\begin{equation}
 |\mu_{[1],s}(u,\vartheta)-1|\leq C\delta_{0}
\end{equation}
The lower bound of (8.254) then implies:
\begin{equation}
 \mu_{m}(t)\geq 1-C\delta_{0}
\end{equation}
hence also:
\begin{equation}
 \bar{\mu}_{m}(t)\geq 1-C\delta_{0}
\end{equation}
Substituting (8.255) and (8.258) into $I_{a}(t)$ we obtain:
\begin{equation}
 I_{a}(t)\leq \int_{0}^{t}(1-C\delta_{0})^{-a}C\delta_{0}\frac{[1+\log(1+t')]}{(1+t')^{2}}dt'\leq C'\delta_{0}
\end{equation}
provided that:
\begin{equation}
 \delta_{0}a\leq \frac{1}{C}
\end{equation}
This is because:
\begin{align*}
 (1-C\delta_{0})^{-a}\leq (1-a^{-1})^{-a}\longrightarrow e,\quad \textrm{as}\quad a\longrightarrow\infty
\end{align*}
On the other hand, we have, in any case, by definition,
\begin{align*}
 \bar{\mu}_{m}(t)\leq 1 \Rightarrow \bar{\mu}_{m}^{-a}(t)\geq 1
\end{align*}
Therefore, the lemma holds in Case 1).

Case 2) $\hat{E}_{s,m}<0$. Let us set in this case:
\begin{equation}
 \hat{E}_{s,m}=-\delta_{1}, \quad \delta_{1}>0
\end{equation}
From the definition of $E_{s}(u)$ and (8.222) we have:
\begin{equation}
 \delta_{1}\leq C\delta_{0}
\end{equation}
In the following we denote:
\begin{equation}
 b(t,s)=\frac{[1+\log(1+t)]}{(1+t)}\frac{(s-t)}{(1+s)},\quad c(t)=\frac{[1+\log(1+t)]}{(1+t)^{2}}
\end{equation}
We also denote:
\begin{equation}
 \hat{\mu}_{s,m}(t)=\min_{(u,\vartheta)\in[0,\epsilon_{0}]\times S^{2}}\hat{\mu}_{s}(t,u,\vartheta)
\end{equation}
From Proposition 8.6 we have, in the present case,
\begin{equation}
 \hat{\mu}_{s}(t,u,\vartheta)\geq 1-\delta_{1}\log(1+t)-C\delta_{0}b(t,s)
\end{equation}
for all $(u,\vartheta)\in[0,\epsilon_{0}]\times S^{2}$, hence:
\begin{equation}
 \hat{\mu}_{s,m}(t)\geq 1-\delta_{1}\log(1+t)-C\delta_{0}b(t,s)
\end{equation}
On the other hand, if
\begin{align*}
 \hat{E}_{s}(u_{m},\vartheta_{m})=\hat{E}_{s,m}
\end{align*}
then
\begin{align*}
 \hat{\mu}_{s,m}(t)\leq\hat{\mu}_{s}(t,u_{m},\vartheta_{m})\\
=1-\delta_{1}\log(1+t)+\hat{Q}_{0,s}(t,u_{m},\vartheta_{m})
\end{align*}
hence:
\begin{equation}
 \hat{\mu}_{s,m}(t)\leq 1-\delta_{1}\log(1+t)+C\delta_{0}b(t,s)
\end{equation}

Moreover, from Proposition 8.6 we have:
\begin{align*}
 -(\frac{\partial\hat{\mu}_{s}}{\partial t})_{-}(t,u,\vartheta)\leq -\frac{(\hat{E}_{s})_{-}(u,\vartheta)}{(1+t)}-
(\hat{Q}_{1,s})_{-}(t,u,\vartheta)
\end{align*}
hence:
\begin{equation}
  -(\frac{\partial\hat{\mu}_{s}}{\partial t})_{-}(t,u,\vartheta)\leq\frac{\delta_{1}}{(1+t)}+C\delta_{0}c(t)
\end{equation}

Consider now the integral $I_{a}(t)$. Let us set:
\begin{equation}
 t_{1}=e^{\frac{1}{2a\delta_{1}}}-1
\end{equation}
We have two subcases to consider:

Subcase 2a) $t\leq t_{1}$. Since $t'\leq t$, where $t^{\prime}$ is the variable of the integration in the integral $I_{a}(t)$, in this subcase we have:
\begin{equation}
 1-\delta_{1}\log(1+t')\geq 1-\frac{1}{2a}
\end{equation}
hence, by (8.266):
\begin{equation}
 \hat{\mu}_{s,m}(t')\geq 1-\frac{1}{a}
\end{equation}
provided that:
\begin{equation}
 \delta_{0}a\leq \frac{1}{2C}
\end{equation}
where we have used the fact that $b(t,s)\leq 1$. In view of (8.256), (8.271) implies
\begin{align*}
 \mu_{m}(t')\geq 1-\frac{2}{a}
\end{align*}
hence also:
\begin{equation}
 \bar{\mu}_{m}(t')\geq 1-\frac{2}{a}
\end{equation}
Now by (8.268) and (8.271) we have:
\begin{equation}
 M(t')\leq (1-\frac{1}{a})^{-1}\{\frac{\delta_{1}}{(1+t')}+C\delta_{0}c(t')\}
\end{equation}
Using (8.273) and (8.274) we obtain:
\begin{equation}
 I_{a}(t)\leq (1-\frac{1}{a})^{-1}(1-\frac{2}{a})^{-a}\int_{0}^{t}\{\frac{\delta_{1}}{(1+t')}+C\delta_{0}c(t')\}dt'
\end{equation}
Here,
\begin{equation}
 \int_{0}^{t}\frac{\delta_{1}}{(1+t')}dt'=\delta_{1}\log(1+t)\leq\delta_{1}\log(1+t_{1})=\frac{1}{2a}
\end{equation}
while,
\begin{equation}
 \int_{0}^{t}C\delta_{0}c(t')dt'\leq C'\delta_{0}\leq \frac{1}{2a}
\end{equation}
provided $\delta_{0}$ is suitably small. In view of the fact that the coefficient in (8.275):
\begin{align*}
 (1-\frac{1}{a})^{-1}(1-\frac{2}{a})^{-a}\longrightarrow e^{2}, \textrm{as} \quad a\longrightarrow\infty
\end{align*}
it follows that for $t\leq t_{1}$:
\begin{equation}
 I_{a}(t)\leq\frac{C}{a}
\end{equation}
Since $\bar{\mu}_{m}(t)\leq 1$, the lemma holds in subcase 2a).

Subcase 2b) $t>t_{1}$. We write:
\begin{equation}
 I_{a}(t)=I_{a}(t_{1})+\int_{t_{1}}^{t}\bar{\mu}_{m}^{-a}(t')M(t')dt'
\end{equation}
By (8.278) with $t$ replaced by $t_{1}$,
\begin{equation}
 I_{a}(t)\leq \frac{C}{a}+\int_{t_{1}}^{t}\bar{\mu}_{m}^{-a}(t')M(t')dt'
\end{equation}
Consider now the upper bound (8.267) evaluated at $t=s$. Since $b(s,s)=0$, we have:
\begin{equation}
 \hat{\mu}_{s,m}(s)\leq 1-\delta_{1}\log(1+s)
\end{equation}
On the other hand, $\mu>0$ in $W^{s}_{\epsilon_{0}}$. It follows that:
\begin{equation}
 s\leq t_{*}, \quad\textrm{where}\quad t_{*}=e^{\frac{1}{\delta_{1}}}-1
\end{equation}
Let us introduce the variable:
\begin{equation}
 \tau'=\log(1+t')
\end{equation}
and correspondingly write:
\begin{equation}
 \tau=\log(1+t),\quad \sigma=\log(1+s)
\end{equation}
and:
\begin{equation}
 \tau_{1}=\log(1+t_{1})=\frac{1}{2a\delta_{1}},\quad \tau_{*}=\log(1+t_{*})=\frac{1}{\delta_{1}}
\end{equation}
Then in (8.280) we have:
\begin{equation}
 \tau_{1}\leq\tau'\leq\tau\leq\sigma\leq\tau_{*}
\end{equation}
Also, we express:
\begin{equation}
 b(t,s)=\frac{(1+\tau)}{e^{\tau}}\frac{(e^{\sigma}-e^{\tau})}{e^{\sigma}}
\end{equation}
Since,
\begin{align*}
 0\leq\frac{(e^{\sigma}-e^{\tau})}{e^{\sigma}}=1-e^{-(\sigma-\tau)}\leq\sigma-\tau
\end{align*}
it follows that
\begin{equation}
 0\leq b(t,s)\leq \frac{(1+\tau)}{e^{\tau}}(\sigma-\tau)
\end{equation}
Substituting (8.288) in (8.266) with $t$ replaced by $t'$ we obtain:
\begin{align}
 \hat{\mu}_{s,m}(t')\geq 1-\delta_{1}\tau'-C\delta_{0}\frac{(1+\tau')}{e^{\tau'}}(\sigma-\tau')\\
=1-\delta_{1}\sigma+\delta_{1}\{1-\frac{C\delta_{0}}{\delta_{1}}\frac{(1+\tau')}{e^{\tau'}}\}(\sigma-\tau')\notag
\end{align}
We have:
\begin{equation}
 \frac{(1+\tau')}{e^{\tau'}}\leq\frac{(1+\tau_{1})}{e^{\tau_{1}}}=(1+\frac{1}{2a\delta_{1}})e^{-\frac{1}{2a\delta_{1}}}
\end{equation}
Hence, under the smallness condition
\begin{align*}
 C\delta_{0}\leq\frac{1}{2a}
\end{align*}
on $\delta_{0}$:
\begin{equation}
 \frac{C\delta_{0}}{\delta_{1}}\frac{(1+\tau')}{e^{\tau'}}\leq\frac{1}{2a\delta_{1}}(1+\frac{1}{2a\delta_{1}})e^{-\frac{1}{2a\delta_{1}}}
\end{equation}
  In view of the fact that
\begin{equation}
 \lim_{x\rightarrow0_{+}}\frac{1}{x}(1+\frac{1}{x})e^{-\frac{1}{x}}=0\notag
\end{equation}
there is a positive constant $K(a)$, depending on $a$, such that:
\begin{align*}
 \frac{1}{x}(1+\frac{1}{x})e^{-\frac{1}{x}}\leq\frac{1}{a}: \forall x\leq\frac{1}{K(a)}
\end{align*}
Here we set:
\begin{align*}
 x=2a\delta_{1}\leq 2aC\delta_{0}
\end{align*}
So the smallness condition:
\begin{equation}
 \delta_{0}\leq\frac{1}{2CaK(a)}
\end{equation}
on $\delta_{0}$ implies:
\begin{equation}
 \frac{1}{2a\delta_{1}}(1+\frac{1}{2a\delta_{1}})e^{-\frac{1}{2a\delta_{1}}}\leq \frac{1}{a}
\end{equation}
It then follows from (8.291) that:
\begin{equation}
 1-\frac{C\delta_{0}}{\delta_{1}}\frac{(1+\tau')}{e^{\tau'}}\geq 1-\frac{1}{a}
\end{equation}
In reference to (8.289), taking into account the fact that:
\begin{equation}
 1-\delta_{1}\sigma\geq 1-\delta_{1}\tau_{*}=0
\end{equation}
we obtain:
\begin{align}
 1-\delta_{1}\sigma+\delta_{1}\{1-\frac{C\delta_{0}}{\delta_{1}}\frac{(1+\tau')}{e^{\tau'}}\}(\sigma-\tau')\\ \notag
\geq 1-\delta_{1}\sigma+(1-\frac{1}{a})\delta_{1}(\sigma-\tau')\\ \notag
\geq (1-\frac{1}{a})(1-\delta_{1}\sigma)+(1-\frac{1}{a})\delta_{1}(\sigma-\tau')
=(1-\frac{1}{a})(1-\delta_{1}\tau')\notag
\end{align}
From (8.289) we then conclude that:
\begin{align}
 \hat{\mu}_{s}(t',u,\vartheta)\geq \hat{\mu}_{s,m}(t')\geq (1-\frac{1}{a})(1-\delta_{1}\tau')\\
:\forall (u,\vartheta)\in[0,\epsilon_{0}]\times S^{2}, \forall t'\geq t_{1}\notag
\end{align}

We now consider the bound (8.268). Since
\begin{equation}
 (1+t')c(t')=\frac{(1+\tau')}{e^{\tau'}}
\end{equation}
the bound (8.268) with $t$ replaced by $t'$ reads:
\begin{equation}
 -(1+t')(\frac{\partial\hat{\mu}_{s}}{\partial t})_{-}(t',u,\vartheta)\leq \delta_{1}
\{1+\frac{C\delta_{0}}{\delta_{1}}\frac{(1+\tau')}{e^{\tau'}}\}
\end{equation}
If (8.292) holds, this implies:
\begin{align}
 -(1+t')(\frac{\partial\hat{\mu}_{s}}{\partial t})_{-}(t',u,\vartheta)\leq \delta_{1}(1+\frac{1}{a})\\
:\forall (u,\vartheta)\in[0,\epsilon_{0}]\times S^{2}, \forall t'\geq t_{1} \notag
\end{align}

The bounds (8.297) and (8.300) together yield:
\begin{equation}
 (1+t')M(t')\leq \frac{1+\frac{1}{a}}{1-\frac{1}{a}}\frac{\delta_{1}}{(1-\delta_{1}\tau')}
\end{equation}
By (8.256), taking $C\delta_{0}\leq a^{-1}$, and (8.297),
\begin{equation}
 \mu_{m}(t')\geq (1-\frac{2}{a})(1-\delta_{1}\tau')
\end{equation}
hence also:
\begin{equation}
 \bar{\mu}_{m}(t')\geq (1-\frac{2}{a})(1-\delta_{1}\tau')
\end{equation}
We now consider the integral in (8.280). In terms of the variable $\tau'$ we have:
\begin{align}
 \int_{t_{1}}^{t}\bar{\mu}_{m}^{-a}(t')M(t')dt'=\int_{\tau_{1}}^{\tau}\bar{\mu}_{m}^{-a}(t')(1+t')M(t')d\tau'\notag
\end{align}
thus (8.301) and (8.303) yields:
\begin{equation}
 \int_{t_{1}}^{t}\bar{\mu}_{m}^{-a}(t')M(t')dt'\leq\frac{1+\frac{1}{a}}{1-\frac{1}{a}}\frac{1}{(1-\frac{2}{a})^{a}}
\int_{\tau_{1}}^{\tau}(1-\delta_{1}\tau')^{-a-1}\delta_{1}d\tau'
\end{equation}
Since:
\begin{align*}
 \int_{\tau_{1}}^{\tau}(1-\delta_{1}\tau')^{-a-1}\delta_{1}d\tau'
=\frac{1}{a}[(1-\delta_{1}\tau)^{-a}-(1-\delta_{1}\tau_{1})^{-a}]\leq\frac{1}{a}(1-\delta_{1}\tau)^{-a}
\end{align*}
while the coefficient:
\begin{align*}
 \frac{1+\frac{1}{a}}{1-\frac{1}{a}}\frac{1}{(1-\frac{2}{a})^{a}} \rightarrow e^{2}, \textrm{as} \quad a\rightarrow\infty
\end{align*}
we conclude that:
\begin{equation}
 \int_{t_{1}}^{t}\bar{\mu}_{m}^{-a}(t')M(t')dt'\leq\frac{C}{a}(1-\delta_{1}\tau)^{-a}
\end{equation}

On the other hand, by (8.267) and (8.288),
\begin{align}
 \hat{\mu}_{s,m}(t)\leq 1-\delta_{1}\tau+C\delta_{0}\frac{(1+\tau)}{e^{\tau}}(\sigma-\tau)\\
=1-\delta_{1}\sigma+\delta_{1}\{1+\frac{C\delta_{0}}{\delta_{1}}\frac{(1+\tau)}{e^{\tau}}\}(\sigma-\tau)\notag
\end{align}
Now by (8.291) with $\tau'$ replaced by $\tau$ and (8.293),
\begin{equation}
 \frac{C\delta_{0}}{\delta_{1}}\frac{(1+\tau)}{e^{\tau}}\leq \frac{1}{a}
\end{equation}
provided (8.292) holds. Hence,
\begin{equation}
 1+\frac{C\delta_{0}}{\delta_{1}}\frac{(1+\tau)}{e^{\tau}}\leq 1+\frac{1}{a}
\end{equation}
Then,
\begin{align}
 1-\delta_{1}\sigma+\delta_{1}\{1+\frac{C\delta_{0}}{\delta_{1}}\frac{(1+\tau)}{e^{\tau}}\}(\sigma-\tau)\\ \notag
\leq 1-\delta_{1}\sigma+(1+\frac{1}{a})\delta_{1}(\sigma-\tau)\\ \notag
\leq (1+\frac{1}{a})(1-\delta_{1}\sigma)+(1+\frac{1}{a})\delta_{1}(\sigma-\tau)\\ \notag
=(1+\frac{1}{a})(1-\delta_{1}\tau)\notag
\end{align}
Thus we obtain:
\begin{equation}
 \hat{\mu}_{s,m}(t)\leq (1+\frac{1}{a})(1-\delta_{1}\tau)
\end{equation}
which together with (8.256) and a smallness condition on $\delta_{0}$ of the form (8.260) yields:
\begin{equation}
 \mu_{m}(t)\leq (1+\frac{2}{a})(1-\delta_{1}\tau)
\end{equation}
hence also:
\begin{equation}
 \bar{\mu}_{m}(t)\leq (1+\frac{2}{a})(1-\delta_{1}\tau)
\end{equation}
Therefore:
\begin{equation}
 \bar{\mu}_{m}^{-a}(t)\geq \frac{1}{(1+\frac{2}{a})^{a}}(1-\delta_{1}\tau)^{-a}
\end{equation}
Since
\begin{align*}
 (1+\frac{2}{a})^{a}\rightarrow e^{2} , \textrm{as}\quad a\rightarrow\infty
\end{align*}
we conclude that:
\begin{equation}
 \bar{\mu}^{-a}_{m}(t)\geq\frac{1}{C}(1-\delta_{1}\tau)^{-a}
\end{equation}
Comparing with (8.305), then yields:
\begin{equation}
 \int_{t_{1}}^{t}\bar{\mu}_{m}^{-a}(t')M(t')dt'\leq\frac{C}{a}\bar{\mu}_{m}^{-a}(t)
\end{equation}
In view of above, (8.280) and the fact that $\bar{\mu}_{m}(t)\leq 1$, the lemma follows as well in subcase 2b).

The proof of the lemma is now complete. $\qed$

$\textbf{Corollary 1}$ Under the assumptions of Lemma 8.11, there is a positive constant $C$ independent of $s$ such that the following holds.
If at some $(u,\vartheta)\in[0,\epsilon_{0}]\times S^{2}$ we have $\hat{E}_{s}(u,\vartheta)<0$, then we have:
\begin{equation}
 \frac{\hat{\mu}_{s}(t,u,\vartheta)}{\hat{\mu}_{s}(t',u,\vartheta)}\leq C\notag
\end{equation}
for all $t'\in[0,t]$ and all $t\in[0,s]$.

$Proof$. Suppose at $(u,\vartheta)$, we have
\begin{align*}
 \hat{E}_{s}(u,\vartheta)=-\delta_{1}, \delta_{1}>0, 
\end{align*}
From (8.266) and (8.267),
\begin{equation}
 1-\delta_{1}\log(1+t)-C\delta_{0}b(t,s)\leq\hat{\mu}_{s}(t,u,\vartheta)\leq 1-\delta_{1}\log(1+t)+C\delta_{0}b(t,s)
\end{equation}
Fixing $a=4$, and setting, in accordance with (8.269),
\begin{equation}
 t_{1}=e^{\frac{1}{8\delta_{1}}}-1
\end{equation}
we have to consider the two subcases according as to whether $t\leq t_{1}$ or $t>t_{1}$.
If $t\leq t_{1}$, from (8.271),
\begin{equation}
 \hat{\mu}_{s}(t',u,\vartheta)\geq\frac{1}{2}
\end{equation}
while by the upper bound in (8.316),
\begin{equation}
 \hat{\mu}_{s}(t,u,\vartheta)\leq\frac{5}{4}
\end{equation}
So the corollary holds when $t\leq t_{1}$.

If $t> t_{1}$, from (8.310) we have:
\begin{equation}
 \hat{\mu}_{s}(t,u,\vartheta)\leq\frac{5}{4}(1-\delta_{1}\tau)
\end{equation}
If then $t'\leq t_{1}$, the corollary holds by (8.318). While if $t'>t_{1}$, from (8.297) we have:
\begin{equation}
 \hat{\mu}_{s}(t',u,\vartheta)\geq\frac{3}{4}(1-\delta_{1}\tau')
\end{equation}
Since $(1-\delta_{1}\tau)\leq(1-\delta_{1}\tau')$, the corollary follows. $\qed$

$\textbf{Corollary 2}$ Under the assumptions of Lemma 8.11, there is a positive constant $C$ independent of $s$ and $a$ such that:
\begin{align*}
 \bar{\mu}_{m}^{-a}(t')\leq C\bar{\mu}_{m}^{-a}(t)
\end{align*}
holds for all $t'\in[0,t]$ and all $t\in[0,s]$.

$Proof$. Case 1) $\hat{E}_{s,m}\geq 0$. In this case, from (8.258):
\begin{equation}
 \bar{\mu}_{m}(t')\geq 1-C\delta_{0}
\end{equation}
On the other hand, by definition,
\begin{align*}
 \bar{\mu}_{m}(t)\leq 1
\end{align*}
It follows that:
\begin{equation}
 \frac{\bar{\mu}_{m}^{-a}(t')}{\bar{\mu}^{-a}_{m}(t)}\leq (1-C\delta_{0})^{-a}
\end{equation}
So the corollary holds provided that
\begin{align*}
 C\delta_{0}\leq\frac{1}{a}
\end{align*}
Case 2) $\hat{E}_{s,m}<0$. In this case we define $\delta_{1}$ as in (8.261) and $t_{1}$ as in (8.269) and we consider the subcases
2a) $t'\leq t_{1}$ and 2b) $t'>t_{1}$ separately.

In subcase 2a), (8.273) holds, thus we have:
\begin{equation}
 \frac{\bar{\mu}^{-a}_{m}(t')}{\bar{\mu}^{-a}_{m}(t)}\leq (1-\frac{2}{a})^{-a}
\end{equation}
which is bounded by a constant independent of $a$ (for $a\geq 4$). 

In subcase 2b), (8.303) as well as (8.312) hold, hence:
\begin{equation}
 \frac{\bar{\mu}^{-a}_{m}(t')}{\bar{\mu}^{-a}_{m}(t)}\leq (\frac{1-\frac{2}{a}}{1+\frac{2}{a}})^{-a}
(\frac{1-\delta_{1}\tau}{1-\delta_{1}\tau'})^{a}\leq (\frac{1-\frac{2}{a}}{1+\frac{2}{a}})^{-a}
\end{equation}
(for,$\tau'\leq\tau$) which is bounded by a constant independent of $a$ (for $a\geq 4$), as required. $\qed$

$\textbf{Corollary 3}$ The assumptions of Lemma 8.11 imply the assumption $\textbf{AS}$.

$Proof$. We consider a given $(u,\vartheta)\in[0,\epsilon_{0}]\times S^{2}$ as in the proof of Corollary 1. Again, we 
have the two cases $\hat{E}_{s}(u,\vartheta)\geq 0$ and $\hat{E}_{s}(u,\vartheta)<0$ to consider.

In the first case, recalling (8.254), we have:
\begin{align*}
 \hat{\mu}_{s}(t,u,\vartheta)\geq 1-C\delta_{0}
\end{align*}
On the other hand, from Proposition 8.6 we have:
\begin{align*}
 \frac{\partial \hat{\mu}_{s}}{\partial t}(t,u,\vartheta)=\frac{\hat{E}_{s}(u,\vartheta)}{(1+t)}+
\hat{Q}_{1,s}(t,u,\vartheta)\leq C\delta_{0}(1+t)^{-1}
\end{align*}
where we have used (8.206), (8.226) and (8.231).
Therefore:
\begin{align*}
 -2\mu^{-1}\frac{\partial\mu}{\partial t}=-2\hat{\mu}_{s}^{-1}\frac{\partial\hat{\mu}_{s}}{\partial t}\geq
-C\delta_{0}(1+t)^{-1}
\end{align*}
On the other hand, $\textbf{F2}$ implies
\begin{align*}
 \textrm{tr}\chi\geq C^{-1}(1+t)^{-1}
\end{align*}
So $\textbf{AS}$ follows in Case 1) if $\delta_{0}$ is suitably small.

In the second case we set as in the proof of Corollary 1, $\hat{E}_{s}(u,\vartheta)=-\delta_{1}, \delta_{1}>0, a=4$ and define $t_{1}$ as (8.317).
Then we have to consider the subcases $t\leq t_{1}$ and $t>t_{1}$.

In the first subcase we have (see (8.318))
\begin{align*}
 \hat{\mu}_{s}(t,u,\vartheta)\geq \frac{3}{4}
\end{align*}
from Proposition 8.6
\begin{align*}
 \frac{\partial\hat{\mu}_{s}}{\partial t}(t,u,\vartheta)\leq\hat{Q}_{1,s}(t,u,\vartheta)\leq C\delta_{0}
(1+t)^{-2}[1+\log(1+t)]
\end{align*}
Therefore:
\begin{align*}
 -2\mu^{-1}\frac{\partial\mu}{\partial t}=-2\hat{\mu}_{s}^{-1}\frac{\partial\hat{\mu}_{s}}{\partial t}\geq
-C\delta_{0}(1+t)^{-2}[1+\log(1+t)]
\end{align*}
So $\textbf{AS}$ follows in subcase 2a) if $\delta_{0}$ is suitably small.

In the second subcase we have from Proposition 8.6,
\begin{align*}
 (1+t)\frac{\partial\hat{\mu}_{s}}{\partial t}(t,u,\vartheta)=-\delta_{1}+(1+t)\hat{Q}_{1,s}(t,u,\vartheta)\\
\leq -\delta_{1}\{1-\frac{C\delta_{0}}{\delta_{1}}\frac{(1+\tau)}{e^{\tau}}\}\leq -\frac{3\delta_{1}}{4}
\end{align*}
if $\delta_{0}$ is suitably small, according to (8.292) with $a=4$. It follows that in subcase 2b)
\begin{align*}
 -2\mu^{-1}\frac{\partial\mu}{\partial t}>0
\end{align*}
So $\textbf{AS}$ again holds. $\qed$

\section{The Actual Estimates for the Solutions of the Propagation Equations}
We now return to (8.175). We would like to obtain an estimate for the $L^{2}$ norm of $\leftexp{(i_{1}...i_{l})}{F}_{l}(t)$
on $[0,\epsilon_{0}]\times S^{2}$. In general, if $\phi$ is a function defined on $W^{*}_{\epsilon_{0}}$ and we consider 
$\phi(t,u)$, the corresponding function on $S^{2}$ depending on $t$ and $u$, then the $L^{2}$ norm of $\phi$ on 
$\Sigma^{\epsilon_{0}}_{t}$ is
\begin{align}
 \|\phi\|_{L^{2}(\Sigma^{\epsilon_{0}}_{t})}=\sqrt{\int_{\Sigma^{\epsilon_{0}}_{t}}\phi^{2}
 d\mu_{\slashed{g}}du}\\
=\sqrt{\int_{[0,\epsilon_{0}]\times S^{2}}(\phi(t,u))^{2}d\mu_{\slashed{g}(t,u)}du}\notag
\end{align}
While the $L^{2}$ norm of $\phi(t,.)$ on $[0,\epsilon_{0}]\times S^{2}$ is:
\begin{equation}
 \|\phi(t)\|_{L^{2}([0,\epsilon_{0}]\times S^{2})}=
\sqrt{\int_{[0,\epsilon_{0}]\times S^{2}}(\phi(t,u))^{2}d\mu_{\slashed{g}(0,0)}du}
\end{equation}
$\slashed{g}(0,0)$ being the standard metric on $S^{2}$. We have:
\begin{equation}
 d\mu_{\slashed{g}(t,u)}=\frac{\sqrt{\det\slashed{g}(t,u)}}{\sqrt{\det\slashed{g}(0,0)}}d\mu_{\slashed{g}(0,0)}
=A(t,u)\frac{\sqrt{\det\slashed{g}(0,u)}}{\sqrt{\det\slashed{g}(0,0)}}d\mu_{\slashed{g}(0,0)}
\end{equation}
Now, by the definition of $\theta$, we have, 
\begin{equation}
 \kappa\textrm{tr}\theta=\frac{1}{\sqrt{\det\slashed{g}}}\frac{\partial}{\partial u}\sqrt{\det\slashed{g}}
\end{equation}
on $\Sigma^{\epsilon_{0}}_{0}$, because $\Xi$ vanishes there. Using $\textbf{F2}$, (6.85) and the relation
\begin{align*}
 \chi=\eta(\slashed{k}-\theta)
\end{align*}
we have:
\begin{equation}
 |\kappa\textrm{tr}\theta+\frac{2}{1-u}|\leq C\delta_{0}
\end{equation}
on $\Sigma^{\epsilon_{0}}_{0}$.

Integrating (8.330) we get:
\begin{equation}
 e^{-C\delta_{0}}(1-u)^{2}\leq\frac{\sqrt{\det\slashed{g}(0,u)}}{\sqrt{\det\slashed{g}(0,0)}}
\leq e^{C\delta_{0}}(1-u)^{2}
\end{equation}
Combining with (8.173), we obtain from (8.328):
\begin{equation}
 e^{-C\delta_{0}}(1-u+t)^{2}d\mu_{\slashed{g}(0,0)}\leq d\mu_{\slashed{g}(t,u)}
\leq e^{C\delta_{0}}(1-u+t)^{2}d\mu_{\slashed{g}(0,0)}
\end{equation}
In view of (8.326) and (8.327) it follows that:
\begin{equation}
 C^{-1}(1+t)\|\phi(t)\|_{L^{2}([0,\epsilon_{0}]\times S^{2})}\leq\|\phi\|_{L^{2}(\Sigma^{\epsilon_{0}}_{t})}
\leq C(1+t)\|\phi(t)\|_{L^{2}([0,\epsilon_{0}]\times S^{2})}
\end{equation}
Moreover, if $\xi$ is an $S_{t,u}$ 1-form defined on $W^{*}_{\epsilon_{0}}$ and we consider $\xi(t,u)$, the corresponding 1-form on 
$S^{2}$ depending on $t$ and $u$, then setting $\phi=|\xi|$, $\phi(t,u)=|\xi(t,u)|$, we have:
\begin{equation}
 \|\phi\|_{L^{2}(\Sigma^{\epsilon_{0}}_{t})}=\|\xi\|_{L^{2}(\Sigma^{\epsilon_{0}}_{t})},
\|\phi(t)\|_{L^{2}([0,\epsilon_{0}]\times S^{2})}=\||\xi(t)|\|_{L^{2}([0,\epsilon_{0}]\times S^{2})}
\end{equation}
hence, by (8.333),
\begin{equation}
 C^{-1}(1+t)\||\xi(t)|\|_{L^{2}([0,\epsilon_{0}]\times S^{2})}\leq\|\xi\|_{L^{2}(\Sigma^{\epsilon_{0}}_{t})}
\leq C(1+t)\||\xi(t)|\|_{L^{2}([0,\epsilon_{0}]\times S^{2})}
\end{equation}

We now consider the terms $\leftexp{(i_{1}...i_{l})}{M}^{0}_{l}(t,u),\leftexp{(i_{1}...i_{l})}{M}^{1}_{l}(t,u),
\leftexp{(i_{1}...i_{l})}{M}^{2}_{l}(t,u)$ on the right hand side of (8.175), given by (8.176), (8.177) and (8.178) respectively.

First, by $\textbf{A3}$, (8.335), and (8.235),
\begin{equation}
 \|\leftexp{(i_{1}...i_{l})}{M}^{0}_{l}(t,u)\|_{L^{2}([0,\epsilon_{0}]\times S^{2})}
\leq C[1+\log(1+t)]^{2}\|\leftexp{(i_{1}...i_{l})}{x}_{l}(0)\|_{L^{2}(\Sigma^{\epsilon_{0}}_{0})}
\end{equation}

We turn to $\leftexp{(i_{1}...i_{l})}{M}^{1}_{l}(t,u)$. We partition $[0,\epsilon_{0}]\times S^{2}$, into
\begin{equation}
 \mathcal{V}_{s-}=\{(u,\vartheta)\in[0,\epsilon_{0}]\times S^{2}:\hat{E}_{s}(u,\vartheta)<0\}
\end{equation}
and
\begin{equation}
 \mathcal{V}_{s+}=\{(u,\vartheta)\in[0,\epsilon_{0}]\times S^{2}:\hat{E}_{s}(u,\vartheta)\geq0\}
\end{equation}
We then have:
\begin{align}
 \|\leftexp{(i_{1}...i_{l})}{M}^{1}_{l}(t)\|^{2}_{L^{2}([0,\epsilon_{0}]\times S^{2})}
=\|\leftexp{(i_{1}...i_{l})}{M}^{1}_{l}(t)\|^{2}_{L^{2}(\mathcal{V}_{s-})}+
\|\leftexp{(i_{1}...i_{l})}{M}^{1}_{l}(t)\|^{2}_{L^{2}(\mathcal{V}_{s+})}
\end{align}
By Minkowski's inequality, we have:
\begin{equation}
 \|\leftexp{(i_{1}...i_{l})}{M}^{1}_{l}(t)\|^{2}_{L^{2}(\mathcal{V}_{s-})}\leq\int_{0}^{t}
\|\leftexp{(i_{1}...i_{l})}{N}^{1}_{l}(t,t')\|_{L^{2}(\mathcal{V}_{s-})}dt'
\end{equation}
where:
\begin{align}
 \leftexp{(i_{1}...i_{l})}{N}^{1}_{l}(t,t',u)=(\frac{\mu(t,u)}{\mu(t',u)})^{2}(1-u+t')^{3}
[-2\mu^{-1}(\frac{\partial\mu}{\partial t})_{-}(t',u)]|\slashed{d}\leftexp{(i_{1}...i_{l})}{\check{f}}_{l}(t',u)|
\end{align}
We have:
\begin{align}
 \|\leftexp{(i_{1}...i_{l})}{N}^{1}_{l}(t,t')\|_{L^{2}(\mathcal{V}_{s-})}\leq
(1+t')^{3}[\max_{\mathcal{V}_{s-}}(\frac{\mu(t)}{\mu(t')})]^{2}
\max_{\mathcal{V}_{s-}}[-2\mu^{-1}(\frac{\partial\mu}{\partial t})_{-}(t')]
\||\slashed{d}\leftexp{(i_{1}...i_{l})}{\check{f}}_{l}(t')|\|_{L^{2}(\mathcal{V}_{s-})}
\end{align}
Now, by Corollary 1 to Lemma 8.11 we have:
\begin{equation}
 \max_{\mathcal{V}_{s-}}(\frac{\mu(t)}{\mu(t')})\leq C
\end{equation}
while by (8.249),
\begin{equation}
 \max_{\mathcal{V}_{s-}}[-2\mu^{-1}(\frac{\partial\mu}{\partial t})_{-}(t')]\leq 2M(t')
\end{equation}
Substituting in (8.342),
\begin{align}
 \|\leftexp{(i_{1}...i_{l})}{N}^{1}_{l}(t,t')\|_{L^{2}(\mathcal{V}_{s-})}\leq
C(1+t')^{3}M(t')\||\slashed{d}\leftexp{(i_{1}...i_{l})}{\check{f}}_{l}(t')|\|_{L^{2}([0,\epsilon_{0}]\times S^{2})}
\end{align}
We define:
\begin{equation}
 \leftexp{(i_{1}...i_{l})}{P}_{l}(t)=(1+t)^{2}\||\slashed{d}\leftexp{(i_{1}...i_{l})}{\check{f}}_{l}(t')|\|_{L^{2}([0,\epsilon_{0}]\times S^{2})}
\end{equation}
Suppose that, for non-negative quantities $\leftexp{(i_{1}...i_{l})}{P}_{l}^{(0)},
\leftexp{(i_{1}...i_{l})}{P}_{l}^{(1)}$, we have:
\begin{equation}
 \leftexp{(i_{1}...i_{l})}{P}_{l}(t)\leq \leftexp{(i_{1}...i_{l})}{P}_{l}^{(0)}(t)
+\leftexp{(i_{1}...i_{l})}{P}_{l}^{(1)}(t)
\end{equation}
Defining then the non-decreasing non-negative quantities $\leftexp{(i_{1}...i_{l})}{\bar{P}}_{l,a}^{(0)},
\leftexp{(i_{1}...i_{l})}{\bar{P}}_{l,a}^{(1)}$, by
\begin{equation}
 \leftexp{(i_{1}...i_{l})}{\bar{P}}_{l,a}^{(0)}=
\sup_{t'\in[0,t]}\{\bar{\mu}^{a}_{m}(t')\leftexp{(i_{1}...i_{l})}{P}_{l}^{(0)}(t')\}
\end{equation}
\begin{equation}
 \leftexp{(i_{1}...i_{l})}{\bar{P}}_{l,a}^{(1)}=
\sup_{t'\in[0,t]}\{(1+t')^{1/2}\bar{\mu}^{a}_{m}(t')\leftexp{(i_{1}...i_{l})}{P}_{l}^{(1)}(t')\}
\end{equation}
Then for $t'\in[0,t]$ we have:
\begin{equation}
 \leftexp{(i_{1}...i_{l})}{P}_{l}(t')\leq
\bar{\mu}_{m}^{-a}(t')\{\leftexp{(i_{1}...i_{l})}{\bar{P}}^{(0)}_{l,a}(t)+
(1+t')^{-1/2}\leftexp{(i_{1}...i_{l})}{\bar{P}}^{(1)}_{l,a}(t)\}
\end{equation}
Substituting in (8.345) we obtain:
\begin{align}
 \|\leftexp{(i_{1}...i_{l})}{N}^{1}_{l}(t,t')\|_{L^{2}(\mathcal{V}_{s-})}\leq
C\{(1+t)\leftexp{(i_{1}...i_{l})}{\bar{P}}^{(0)}_{l,a}(t)+(1+t)^{1/2}\leftexp{(i_{1}...i_{l})}{\bar{P}}^{(1)}_{l,a}(t)\}
\bar{\mu}_{m}^{-a}(t')M(t')
\end{align}
for all $t'\in[0,t]$. Substituting this in (8.340) we obtain:
\begin{align}
 \|\leftexp{(i_{1}...i_{l})}{M}^{1}_{l}(t)\|_{L^{2}(\mathcal{V}_{s-})}\leq
C\{(1+t)\leftexp{(i_{1}...i_{l})}{\bar{P}}^{(0)}_{l,a}(t)+(1+t)^{1/2}\leftexp{(i_{1}...i_{l})}{\bar{P}}^{(1)}_{l,a}(t)\}
I_{a}(t)
\end{align}
We then apply Lemma 8.11 to conclude that:
\begin{align}
 \|\leftexp{(i_{1}...i_{l})}{M}^{1}_{l}(t)\|_{L^{2}(\mathcal{V}_{s-})}\leq
Ca^{-1}\{(1+t)\leftexp{(i_{1}...i_{l})}{\bar{P}}^{(0)}_{l,a}(t)+(1+t)^{1/2}\leftexp{(i_{1}...i_{l})}{\bar{P}}^{(1)}_{l,a}(t)\}
\bar{\mu}_{m}^{-a}(t)
\end{align}
In analogy with (8.340) and (8.342), with $\mathcal{V}_{s+}$ in the role of $\mathcal{V}_{s-}$, we have:
\begin{equation}
 \|\leftexp{(i_{1}...i_{l})}{M}^{1}_{l}(t)\|^{2}_{L^{2}(\mathcal{V}_{s+})}\leq\int_{0}^{t}
\|\leftexp{(i_{1}...i_{l})}{N}^{1}_{l}(t,t')\|_{L^{2}(\mathcal{V}_{s+})}dt'
\end{equation}
and:
\begin{equation}
 \|\leftexp{(i_{1}...i_{l})}{N}^{1}_{l}(t,t')\|_{L^{2}(\mathcal{V}_{s+})}\leq
(1+t')^{3}[\max_{\mathcal{V}_{s+}}(\frac{\mu(t)}{\mu(t')})]^{2}
\max_{\mathcal{V}_{s+}}[-2\mu^{-1}(\frac{\partial\mu}{\partial t})_{-}(t')]
\||\slashed{d}\leftexp{(i_{1}...i_{l})}{\check{f}}_{l}(t')|\|_{L^{2}(\mathcal{V}_{s+})}
\end{equation}
From (8.255) we have:
\begin{equation}
 \max_{\mathcal{V}_{s+}}[-2\mu^{-1}(\frac{\partial\mu}{\partial t})_{-}(t')]
\leq C\delta_{0}(1+t')^{-2}[1+\log(1+t')]
\end{equation}
and from (8.254) we have:
\begin{equation}
 \max_{\mathcal{V}_{s+}}(\frac{\mu(t)}{\mu(t')})\leq C[1+\log(1+t)]
\end{equation}
Substituting (8.356) and (8.357) as well as (8.346) in (8.354) and (8.355) we obtain:
\begin{align}
 \|\leftexp{(i_{1}...i_{l})}{M}^{1}_{l}(t)\|_{L^{2}(\mathcal{V}_{s+})}\leq
C\delta_{0}[1+\log(1+t)]^{2}\int_{0}^{t}\frac{[1+\log(1+t')]}{(1+t')}\leftexp{(i_{1}...i_{l})}{P}_{l}(t')dt'
\end{align}
hence, by (8.350),
\begin{align}
 \|\leftexp{(i_{1}...i_{l})}{M}^{1}_{l}(t)\|_{L^{2}(\mathcal{V}_{s+})}\leq
C\delta_{0}[1+\log(1+t)]^{2}\{\leftexp{(i_{1}...i_{l})}{\bar{P}}^{(0)}_{l,a}(t)+\leftexp{(i_{1}...i_{l})}{\bar{P}}^{(1)}_{l,a}(t)\}\cdot\\
\int_{0}^{t}(1+t')^{-1}[1+\log(1+t')]\bar{\mu}_{m}^{-a}(t')dt'\notag
\end{align}
To estimate the integral in (8.359), we follow the proof of Corollary 2 to Lemma 8.11. We have two cases $\hat{E}_{s,m}\geq 0$
and $\hat{E}_{s,m}<0$ to consider. In the first case (8.322) holds, hence if
\begin{align*}
 \delta_{0}a\leq \frac{1}{C}
\end{align*}
we have:
\begin{equation}
 \bar{\mu}_{m}(t')^{-a}\leq(1-\frac{1}{a})^{-a}\leq C
\end{equation}
and then
\begin{equation}
 \int_{0}^{t}\frac{[1+\log(1+t')]}{(1+t')}\bar{\mu}^{-a}_{m}(t')dt'
\leq C[1+\log(1+t)]^{2}
\end{equation}
In the second case we define as before, $\delta_{1}>0$ by $\hat{E}_{s,m}=-\delta_{1}$, and $t_{1}$ according to (8.269).
Then for $t'\leq t_{1}$ (8.273) holds, so an estimate of the form (8.361) holds for all $t\leq t_{1}$.

If $t'>t_{1}$, then (8.303) holds, under the smallness condition (8.292) on $\delta_{0}$. Thus we have:
\begin{align}
 \int_{t_{1}}^{t}\frac{[1+\log(1+t')]}{(1+t')}\bar{\mu}_{m}^{-a}(t')dt'
\leq C\int_{\tau_{1}}^{\tau}(1+\tau')(1-\delta_{1}\tau')^{-a}d\tau'\leq\\ \notag
C(1+\tau)\int_{\tau_{1}}^{\tau}(1-\delta_{1}\tau')^{-a}d\tau'\leq
C(1+\tau)\cdot\frac{(1-\delta_{1}\tau)^{1-a}}{\delta_{1}(a-1)}\\ \notag
=2C(1-\frac{1}{a})^{-1}\tau_{1}(1+\tau)(1-\delta_{1}\tau)^{1-a}\notag
\end{align}
The last equality comes from the first of (8.285). By the upper bound (8.312) we conclude that the last is bounded by a constant multiple of 
\begin{align*}
 [1+\log(1+t)]^{2}\bar{\mu}_{m}^{1-a}(t)
\end{align*}

We conclude that in general,
\begin{equation}
 \int_{0}^{t}\frac{[1+\log(1+t')]}{(1+t')}\bar{\mu}_{m}^{-a}(t')dt'
\leq C[1+\log(1+t)]^{2}\bar{\mu}_{m}^{1-a}(t)
\end{equation}
which, substituted in (8.359), yields:
\begin{align}
 \|\leftexp{(i_{1}...i_{l})}{M}^{1}_{l}(t)\|_{L^{2}(\mathcal{V}_{s+})}
\leq C\delta_{0}[1+\log(1+t)]^{4}\{\leftexp{(i_{1}...i_{l})}{\bar{P}^{(0)}_{l,a}(t)}
+\leftexp{(i_{1}...i_{l})}{\bar{P}^{(1)}_{l,a}(t)}\}\bar{\mu}^{1-a}_{m}(t)
\end{align}
Combining (8.364) and (8.353), and taking account into the fact that $C\delta_{0}\leq \frac{1}{a}$, we obtain:
\begin{align}
 \|\leftexp{(i_{1}...i_{l})}{M}^{1}_{l}(t)\|_{L^{2}([0,\epsilon_{0}]\times S^{2})}
\leq Ca^{-1}\{(1+t)\leftexp{(i_{1}...i_{l})}{\bar{P}}^{(0)}_{l,a}(t)+(1+t)^{1/2}\leftexp{(i_{1}...i_{l})}{\bar{P}}^{(1)}_{l,a}(t)\}
\bar{\mu}_{m}^{-a}(t)
\end{align}

We turn to $\leftexp{(i_{1}...i_{l})}{M}^{2}_{l}(t,u)$. We have:
\begin{equation}
 \|\leftexp{(i_{1}...i_{l})}{M}^{2}_{l}(t)\|^{2}_{L^{2}([0,\epsilon_{0}]\times S^{2})}=
\|\leftexp{(i_{1}...i_{l})}{M}^{2}_{l}(t)\|^{2}_{L^{2}(\mathcal{V}_{s-})}
+\|\leftexp{(i_{1}...i_{l})}{M}^{2}_{l}(t)\|^{2}_{L^{2}(\mathcal{V}_{s+})}
\end{equation}
and also:
\begin{equation}
 \|\leftexp{(i_{1}...i_{l})}{M}^{2}_{l}(t)\|_{L^{2}(\mathcal{V}_{s-})}
\leq \int_{0}^{t}\|\leftexp{(i_{1}...i_{l})}{N}^{2}_{l}(t,t')\|_{L^{2}(\mathcal{V}_{s-})}dt'
\end{equation}
\begin{equation}
 \|\leftexp{(i_{1}...i_{l})}{M}^{2}_{l}(t)\|_{L^{2}(\mathcal{V}_{s+})}
\leq \int_{0}^{t}\|\leftexp{(i_{1}...i_{l})}{N}^{2}_{l}(t,t')\|_{L^{2}(\mathcal{V}_{s+})}dt'
\end{equation}
where
\begin{equation}
 \leftexp{(i_{1}...i_{l})}{N}^{2}_{l}(t,t',u)=(\frac{\mu(t,u)}{\mu(t',u)})^{2}
(1-u+t')^{3}(\frac{1}{2}\textrm{tr}\chi(t',u))|\slashed{d}\leftexp{(i_{1}...i_{l})}{\check{f}}_{l}(t',u)|
\end{equation}
We have:
\begin{align}
 \|\leftexp{(i_{1}...i_{l})}{N}^{2}_{l}(t,t')\|_{L^{2}(\mathcal{V}_{s-})}\leq
(1+t')^{3}[\max_{\mathcal{V}_{s-}}(\frac{\mu(t)}{\mu(t')})]^{2}\max_{\mathcal{V}_{s-}}(\frac{1}{2}\textrm{tr}\chi(t'))
\||\slashed{d}\leftexp{(i_{1}...i_{l})}{\check{f}}_{l}(t')|\|_{L^{2}(\mathcal{V}_{s-})}
\end{align}
Now, by $\textbf{F2}$,
\begin{equation}
 \max_{[0,\epsilon_{0}]\times S^{2}}(\frac{1}{2}\textrm{tr}\chi(t'))\leq C(1+t')^{-1}
\end{equation}
This together with (8.343) imply:
\begin{equation}
 \|\leftexp{(i_{1}...i_{l})}{N}^{2}_{l}(t,t')\|_{L^{2}(\mathcal{V}_{s-})}\leq
C(1+t')^{2}\||\slashed{d}\leftexp{(i_{1}...i_{l})}{\check{f}}_{l}(t')|\|_{L^{2}([0,\epsilon_{0}]\times S^{2})}
\end{equation}
Then from (8.346), (8.350) and Corollary 2 of Lemma 8.11 we have:
\begin{align}
 \|\leftexp{(i_{1}...i_{l})}{N}^{2}_{l}(t,t')\|_{L^{2}(\mathcal{V}_{s-})}\leq C\leftexp{(i_{1}...i_{l})}{P}_{l}(t')\\ \notag
\leq C\{\leftexp{(i_{1}...i_{l})}{\bar{P}}^{(0)}_{l,a}(t)+(1+t')^{-1/2}\leftexp{(i_{1}...i_{l})}{\bar{P}}^{(1)}_{l,a}(t)\}\bar{\mu}^{-a}_{m}(t')\\ \notag
\leq C'\{\leftexp{(i_{1}...i_{l})}{\bar{P}}^{(0)}_{l,a}(t)+(1+t')^{-1/2}\leftexp{(i_{1}...i_{l})}{\bar{P}}^{(1)}_{l,a}(t)\}\bar{\mu}^{-a}_{m}(t)\notag
\end{align}
Substituting (8.373) in (8.367) we get:
\begin{equation}
 \|\leftexp{(i_{1}...i_{l})}{M}^{2}_{l}(t)\|_{L^{2}(\mathcal{V}_{s-})}
\leq C\{(1+t)\leftexp{(i_{1}...i_{l})}{\bar{P}}^{(0)}_{l,a}(t)+(1+t)^{1/2}\leftexp{(i_{1}...i_{l})}{\bar{P}}^{(1)}_{l,a}(t)\}\bar{\mu}^{-a}_{m}(t)
\end{equation}

In analogy with (8.370) we have:
\begin{align}
 \|\leftexp{(i_{1}...i_{l})}{N}^{2}_{l}(t,t')\|_{L^{2}(\mathcal{V}_{s+})}\leq
(1+t')^{3}[\max_{\mathcal{V}_{s+}}(\frac{\mu(t)}{\mu(t')})]^{2}\max_{\mathcal{V}_{s+}}(\frac{1}{2}\textrm{tr}\chi(t'))
\||\slashed{d}\leftexp{(i_{1}...i_{l})}{\check{f}}_{l}(t')|\|_{L^{2}(\mathcal{V}_{s+})}
\end{align}
To estimate
\begin{align*}
 \max_{\mathcal{V}_{s+}}(\frac{\mu(t)}{\mu(t')})=
\max_{(u,\vartheta)\in\mathcal{V}_{s+}}(\frac{\hat{\mu}_{s}(t,u,\vartheta)}{\hat{\mu}_{s}(t',u,\vartheta)})
\end{align*}
we appeal to Proposition 8.6 to obtain:
\begin{equation}
 \max_{\mathcal{V}_{s+}}(\frac{\mu(t)}{\mu(t')})\leq
\max_{(u,\vartheta)\in\mathcal{V}_{s+}}(\frac{1+C\delta_{0}+\hat{E}_{s}(u,\vartheta)\log(1+t)}
{1-C\delta_{0}+\hat{E}_{s}(u,\vartheta)\log(1+t')})
\end{equation}
We have two cases to distinguish according as to whether $t'<\sqrt{t}$ or $t'\geq\sqrt{t}$. In the second case we have $1+t'\geq 1+\sqrt{t}
\geq\sqrt{1+t}$, hence:
\begin{align*}
 \hat{E}_{s}(u,\vartheta)\log(1+t')\geq\frac{1}{2}\hat{E}_{s}(u,\vartheta)\log(1+t) : \forall (u,\vartheta)\in\mathcal{V}_{s+}
\end{align*}
Then from (8.376) we obtain:
\begin{equation}
 \max_{\mathcal{V}_{s+}}(\frac{\mu(t)}{\mu(t')})\leq 2: t'\geq\sqrt{t}
\end{equation}
provided that 
\begin{align*}
 C\delta_{0}\leq\frac{1}{3}
\end{align*}
On the other hand, in the first case, we have, trivially:
\begin{equation}
 \max_{\mathcal{V}_{s+}}(\frac{\mu(t)}{\mu(t')})\leq C[1+\log(1+t)] : t'<\sqrt{t}
\end{equation}
In view of (8.371) and the definition (8.346) we obtain:
\begin{align}
 \|\leftexp{(i_{1}...i_{l})}{N}^{2}_{l}(t,t')\|_{L^{2}(\mathcal{V}_{s+})}\leq 
C[1+\log(1+t)]^{2}\leftexp{(i_{1}...i_{l})}{P}_{l}(t') : t'<\sqrt{t}\\ \notag
 \|\leftexp{(i_{1}...i_{l})}{N}^{2}_{l}(t,t')\|_{L^{2}(\mathcal{V}_{s+})}\leq 
C\leftexp{(i_{1}...i_{l})}{P}_{l}(t') : t'\geq\sqrt{t}
\end{align}
Substituting (8.379) in (8.368), yields:
\begin{align}
\|\leftexp{(i_{1}...i_{l})}{M}^{2}_{l}(t)\|_{L^{2}(\mathcal{V}_{s+})}\leq
\int_{0}^{\sqrt{t}}\|\leftexp{(i_{1}...i_{l})}{N}^{2}_{l}(t,t')\|_{L^{2}(\mathcal{V}_{s+})}dt'
+\int_{\sqrt{t}}^{t}\|\leftexp{(i_{1}...i_{l})}{N}^{2}_{l}(t,t')\|_{L^{2}(\mathcal{V}_{s+})}dt'\\ \notag
\leq C[1+\log(1+t)]^{2}\int_{0}^{\sqrt{t}}\leftexp{(i_{1}...i_{l})}{P}_{l}(t')dt'
+C\int_{\sqrt{t}}^{t}\leftexp{(i_{1}...i_{l})}{P}_{l}(t')dt'\\ \notag
\leq C[1+\log(1+t)]^{2}\int_{0}^{\sqrt{t}}\bar{\mu}_{m}^{-a}(t')\{
\leftexp{(i_{1}...i_{l})}{\bar{P}}^{(0)}_{l,a}(t)+(1+t')^{-1/2}\leftexp{(i_{1}...i_{l})}{\bar{P}}^{(1)}_{l,a}(t)\}dt'\\ \notag
+C\int_{\sqrt{t}}^{t}\bar{\mu}_{m}^{-a}(t')\{
\leftexp{(i_{1}...i_{l})}{\bar{P}}^{(0)}_{l,a}(t)+(1+t')^{-1/2}\leftexp{(i_{1}...i_{l})}{\bar{P}}^{(1)}_{l,a}(t)\}dt'\\ \notag
\leq C[1+\log(1+t)]^{2}\{(1+t)^{1/2}\leftexp{(i_{1}...i_{l})}{\bar{P}}^{(0)}_{l,a}(t)+(1+t)^{1/4}\leftexp{(i_{1}...i_{l})}{\bar{P}}^{(1)}_{l,a}(t)\}
\bar{\mu}_{m}^{-a}(t)\\ \notag
+C\{(1+t)\leftexp{(i_{1}...i_{l})}{\bar{P}}^{(0)}_{l,a}(t)+(1+t)^{1/2}\leftexp{(i_{1}...i_{l})}{\bar{P}}^{(1)}_{l,a}(t)\}
\bar{\mu}_{m}^{-a}(t)\\ \notag
\leq C'\{(1+t)\leftexp{(i_{1}...i_{l})}{\bar{P}}^{(0)}_{l,a}(t)+(1+t)^{1/2}\leftexp{(i_{1}...i_{l})}{\bar{P}}^{(1)}_{l,a}(t)\}
\bar{\mu}_{m}^{-a}(t)\notag
\end{align}
Here we have used (8.350) and Corollary 2 of Lemma 8.11.

Combining (8.374) and (8.380) we obtain:
\begin{equation}
 \|\leftexp{(i_{1}...i_{l})}{M}^{2}_{l}(t)\|_{L^{2}([0,\epsilon_{0}]\times S^{2})}
\leq C\{(1+t)\leftexp{(i_{1}...i_{l})}{\bar{P}}^{(0)}_{l,a}(t)+(1+t)^{1/2}\leftexp{(i_{1}...i_{l})}{\bar{P}}^{(1)}_{l,a}(t)\}
\bar{\mu}_{m}^{-a}(t)
\end{equation}

Finally, (8.336), (8.365) and (8.381) yield, through (8.175), the following estimate:
\begin{align}
 \|\leftexp{(i_{1}...i_{l})}{F}_{l}(t)\|_{L^{2}([0,\epsilon_{0}]\times S^{2})}\leq 
C(1+t)^{-3}[1+\log(1+t)]^{2}\|\leftexp{(i_{1}...i_{l})}{x}_{l}(0)\|_{L^{2}(\Sigma^{\epsilon_{0}}_{0})}\\ \notag
+C(1+t)^{-2}\{\leftexp{(i_{1}...i_{l})}{\bar{P}}^{(0)}_{l,a}(t)+(1+t)^{-1/2}\leftexp{(i_{1}...i_{l})}{\bar{P}}^{(1)}_{l,a}(t)\}
\bar{\mu}_{m}^{-a}(t)
\end{align}

We now consider the estimate (8.179) for $\leftexp{(i_{1}...i_{l})}{G}_{l}(t,u)$. By (8.343) and (8.357) we have, in general,
\begin{equation}
 \max_{[0,\epsilon_{0}]\times S^{2}}(\frac{\mu(t)}{\mu(t')})\leq C[1+\log(1+t)]
\end{equation}
Hence (8.179) implies:
\begin{equation}
 \leftexp{(i_{1}...i_{l})}{G}_{l}(t,u)\leq C(1+t)^{-3}[1+\log(1+t)]^{2}\int_{0}^{t}(1+t')^{3}
|\leftexp{(i_{1}...i_{l})}{\tilde{g}}_{l}(t',u)|dt'
\end{equation}
It follows by Minkowski's inequality that:
\begin{equation}
 \|\leftexp{(i_{1}...i_{l})}{G}_{l}(t)\|_{L^{2}([0,\epsilon_{0}]\times S^{2})}
\leq C(1+t)^{-3}[1+\log(1+t)]^{2}\int_{0}^{t}(1+t')^{3}\||\leftexp{(i_{1}...i_{l})}{\tilde{g}}_{l}(t')|\|_{L^{2}([0,\epsilon_{0}]\times S^{2})}dt'
\end{equation}

Now the $S_{t,u}$ 1-form $\leftexp{(i_{1}...i_{l})}{\tilde{g}}_{l}$ is given by (8.150). Here we must distinguish the principal acoustical terms.
In reference to (8.64), we define the function:
\begin{equation}
 \leftexp{(i_{1}...i_{l})}{\dot{h}}_{l}=\leftexp{(i_{1}...i_{l})}{h}_{l}-2\hat{\chi}\cdot\leftexp{(i_{1}...i_{l})}{\hat{\chi}}_{l}
\end{equation}
then $\slashed{d} \leftexp{(i_{1}...i_{l})}{\dot{h}}_{l}$ does not contain principal acoustical terms. Also, according to the discussion following 
Lemma 8.4, the principal acoustical part of $\leftexp{(i_{1}...i_{l})}{g}_{l}$ consists of the sum (8.89). Therefore the $S_{t,u}$ 1-form:
\begin{equation}
 \leftexp{(i_{1}...i_{l})}{\dot{g}}_{l}=\leftexp{(i_{1}...i_{l})}{g}_{l}+
\sum_{k=0}^{l-1}\slashed{\mathcal{L}}_{\leftexp{(R_{i_{l-k}})}{Z}}\leftexp{(i_{1}...i_{l-k-1}i_{l-k+1}...i_{l})}{x}_{l-1}
\end{equation}
does not contain principal acoustical terms. We conclude that, writing:
\begin{equation}
 \leftexp{(i_{1}...i_{l})}{\tilde{g}}_{l}=2\mu\hat{\chi}\cdot\slashed{D}\leftexp{(i_{1}...i_{l})}{\hat{\chi}}_{l}
-\sum_{k=0}^{l-1}\slashed{\mathcal{L}}_{\leftexp{(R_{i_{l-k}})}{Z}}\leftexp{(i_{1}...i_{l-k-1}i_{l-k+1}...i_{l})}{x}_{l-1}+
\leftexp{(i_{1}...i_{l})}{\ddot{g}}_{l}
\end{equation}
where
\begin{equation}
 \leftexp{(i_{1}...i_{l})}{\ddot{g}}_{l}=\leftexp{(i_{1}...i_{l})}{\dot{g}}_{l}+
\mu(2\slashed{D}\hat{\chi}\cdot\leftexp{(i_{1}...i_{l})}{\hat{\chi}}_{l}+\slashed{d}\leftexp{(i_{1}...i_{l})}{\dot{h}}_{l})
\end{equation}
does not contain principal acoustical terms.

Let us define
\begin{equation}
 X_{l}(t)=\max_{i_{1}...i_{l}}\||\leftexp{(i_{1}...i_{l})}{x}_{l}(t)|\|_{L^{2}([0,\epsilon_{0}]\times S^{2})}
\end{equation}
We shall estimate $\||\leftexp{(i_{1}...i_{l})}{\tilde{g}}_{l}(t)\|_{L^{2}([0,\epsilon_{0}]\times S^{2})}$ in terms of $X_{l}(t)$. Consider first the
second term on the right of (8.388). By (8.98) we have:
\begin{align}
 |\sum_{k=0}^{l-1}\slashed{\mathcal{L}}_{\leftexp{(R_{i_{l-k}})}{Z}}\leftexp{(i_{1}...i_{l-k-1}i_{l-k+1}...i_{l})}{x}_{l-1}|\leq \\ \notag 
C(1+t)^{-1}\{\sum_{k=0}^{l-1}|\leftexp{(R_{i_{l-k}})}{Z}|\sum_{j}|\leftexp{(i_{1}...i_{l-k-1}i_{l-k+1}...i_{l}j)}{x}_{l}|
+\leftexp{(i_{1}...i_{l})}{z}_{l}\}
\end{align}
where
\begin{align}
 \leftexp{(i_{1}...i_{l})}{z}_{l}=\\\notag
\sum_{k=0}^{l-1}|\leftexp{(R_{i_{l-k}})}{Z}|
(\sum_{j}|\slashed{d}(\leftexp{(i_{1}...i_{l-k-1}i_{l-k+1}...i_{l}j)}{\check{f}}_{l})|
+(1+t)|\slashed{d}\mu||\slashed{d}(R_{i_{l}}...R_{i_{l-k+1}}R_{i_{l-k-1}}...R_{i_{1}}\textrm{tr}\chi|)\\ \notag
\sum_{k=0}^{l-1}(\max_{j}|\slashed{\mathcal{L}}_{R_{j}}\leftexp{(R_{i_{l-k}})}{Z}|)
(|\leftexp{(i_{1}...i_{l-k+1}i_{l-k+1}...i_{l})}{x}_{l-1}|+|\slashed{d}(\leftexp{(i_{1}...i_{l-k-1}i_{l-k+1}...i_{l})}{\check{f}}_{l-1})|)\notag
\end{align}
From (6.184), we have:
\begin{equation}
 |\leftexp{(R_{i})}{Z}|\leq C\delta_{0}(1+t)^{-1}[1+\log(1+t)]
\end{equation}
Substituting this in (8.391), we deduce:
\begin{align}
 \||\sum_{k=0}^{l-1}\slashed{\mathcal{L}}_{\leftexp{(R_{i_{l-k}})}{Z}}\leftexp{(i_{1}...i_{l-k-1}i_{l-k+1}...i_{l})}{x}_{l-1}|\|_{L^{2}
([0,\epsilon_{0}]\times S^{2})}\\\notag
\leq Cl\delta_{0}(1+t)^{-2}[1+\log(1+t)]X_{l}(t)+C(1+t)^{-1}\|\leftexp{(i_{1}...i_{l})}{z}_{l}(t)\|_{L^{2}([0,\epsilon_{0}]\times S^{2})}
\end{align}

Consider next the first term on the right in (8.388). By the 2nd of $\textbf{F2}$ we have:
\begin{align}
 \|2\mu|\hat{\chi}\cdot\slashed{D}\leftexp{(i_{1}...i_{l})}{\hat{\chi}}_{l}|(t)\|_{L^{2}([0,\epsilon_{0}]\times S^{2})}\\\notag
\leq C\delta_{0}(1+t)^{-2}[1+\log(1+t)]\|\mu|\slashed{D}\leftexp{(i_{1}...i_{l})}{\hat{\chi}}_{l}|(t)\|_{L^{2}([0,\epsilon_{0}]\times S^{2})}
\end{align}
Using the elliptic estimate (8.149) and the inequality:
\begin{equation}
 C^{-1}(1+t)\||\xi(t,u)|\|_{L^{2}(S^{2})}\leq\|\xi\|_{L^{2}(S_{t,u})}\leq C(1+t)\||\xi(t,u)|\|_{L^{2}(S^{2})}
\end{equation}
which holds for any covariant $S_{t,u}$ tensorfield $\xi$ of any rank, $\xi(t,u)$ being the corresponding tensorfield on $S^{2}$,
we obtain:
\begin{align}
 \|\mu|\slashed{D}\leftexp{(i_{1}...i_{l})}{\hat{\chi}}_{l}|(t,u)\|_{L^{2}(S^{2})}\\ \notag
\leq C\||\leftexp{(i_{1}...i_{l})}{x}_{l}(t,u)|\|_{L^{2}(S^{2})}
+C\||\slashed{d}\leftexp{(i_{1}...i_{l})}{\check{f}}_{l}(t,u)|\|_{L^{2}(S^{2})}
+C\|\mu|\leftexp{(i_{1}...i_{l})}{i}_{l}|(t,u)\|_{L^{2}(S^{2})}\\ \notag
+C\delta_{0}(1+t)^{-1}[1+\log(1+t)]\||\leftexp{(i_{1}...i_{l})}{\hat{\chi}}_{l}(t,u)|\|_{L^{2}(S^{2})}
\end{align}
Taking the $L^{2}$ norms on $[0,\epsilon_{0}]$ then yields:
\begin{align}
 \|\mu|\slashed{D}\leftexp{(i_{1}...i_{l})}{\hat{\chi}}_{l}|(t)\|_{L^{2}([0,\epsilon_{0}]\times S^{2})}\\ \notag
\leq C\||\leftexp{(i_{1}...i_{l})}{x}_{l}(t)|\|_{L^{2}([0,\epsilon_{0}]\times S^{2})}
+C\||\slashed{d}\leftexp{(i_{1}...i_{l})}{\check{f}}_{l}(t)|\|_{L^{2}([0,\epsilon_{0}]\times S^{2})}
+C\|\mu|\leftexp{(i_{1}...i_{l})}{i}_{l}|(t)\|_{L^{2}([0,\epsilon_{0}]\times S^{2})}\\ \notag
+C\delta_{0}(1+t)^{-1}[1+\log(1+t)]\||\leftexp{(i_{1}...i_{l})}{\hat{\chi}}_{l}(t)|\|_{L^{2}([0,\epsilon_{0}]\times S^{2})}
\end{align}

In view of (8.394), (8.395) and (8.398), we obtain from (8.388):
\begin{align}
 \||\leftexp{(i_{1}...i_{l})}{\tilde{g}}_{l}(t)|\|_{L^{2}([0,\epsilon_{0}]\times S^{2})}\\ \notag
\leq C(l+1)\delta_{0}(1+t)^{-2}[1+\log(1+t)]X_{l}(t)+\leftexp{(i_{1}...i_{l})}{Q}_{l}(t)
\end{align}
where:
\begin{align}
 \leftexp{(i_{1}...i_{l})}{Q}_{l}(t)=C\delta_{0}^{2}(1+t)^{-3}[1+\log(1+t)]^{2}\||\leftexp{(i_{1}...i_{l})}{\hat{\chi}}_{l}(t)|\|
_{L^{2}([0,\epsilon_{0}]\times S^{2})}\\\notag
C\delta_{0}(1+t)^{-2}[1+\log(1+t)]\||\slashed{d}\leftexp{(i_{1}...i_{l})}{\check{f}}_{l}(t)|\|_{L^{2}([0,\epsilon_{0}]\times S^{2})}\\ \notag
+C\delta_{0}(1+t)^{-2}[1+\log(1+t)]\|\mu|\leftexp{(i_{1}...i_{l})}{i}_{l}|(t)\|_{L^{2}([0,\epsilon_{0}]\times S^{2})}\\ \notag
+C(1+t)^{-1}\|\leftexp{(i_{1}...i_{l})}{z}_{l}(t)\|_{L^{2}([0,\epsilon_{0}]\times S^{2})}\\ \notag
+\||\leftexp{(i_{1}...i_{l})}{\ddot{g}}_{l}(t)|\|_{L^{2}([0,\epsilon_{0}]\times S^{2})}\notag
\end{align}
Returning to (8.170) and taking $L^{2}$ norms on $[0,\epsilon_{0}]\times S^{2}$ we obtain:
\begin{equation}
 \||\leftexp{(i_{1}...i_{l})}{x}_{l}(t)|\|_{L^{2}([0,\epsilon_{0}]\times S^{2})}
\leq \|\leftexp{(i_{1}...i_{l})}{F}_{l}(t)\|_{L^{2}([0,\epsilon_{0}]\times S^{2})}+
\|\leftexp{(i_{1}...i_{l})}{G}_{l}(t)\|_{L^{2}([0,\epsilon_{0}]\times S^{2})}
\end{equation}
Substituting (8.382), (8.385) and (8.399) then yields:
\begin{align}
 \||\leftexp{(i_{1}...i_{l})}{x}_{l}(t)|\|_{L^{2}([0,\epsilon_{0}]\times S^{2})}\leq \leftexp{(i_{1}...i_{l})}{B}_{l}(t)\\ \notag
+C(l+1)\delta_{0}(1+t)^{-3}[1+\log(1+t)]^{2}\int_{0}^{t}(1+t')[1+\log(1+t')]X_{l}(t')dt'
\end{align}
where:
\begin{align}
 \leftexp{(i_{1}...i_{l})}{B}_{l}(t)=C(1+t)^{-3}[1+\log(1+t)]^{2}\|\leftexp{(i_{1}...i_{l})}{x}_{l}(0)\|_{L^{2}(\Sigma^{\epsilon_{0}}_{0})}\\ \notag
+C(1+t)^{-2}\{\leftexp{(i_{1}...i_{l})}{\bar{P}}^{(0)}_{l,a}(t)+(1+t)^{-1/2}\leftexp{(i_{1}...i_{l})}{\bar{P}}^{(1)}_{l,a}(t)\}\bar{\mu}_{m}^{-a}(t)\\ \notag
+C(1+t)^{-3}[1+\log(1+t)]^{2}\int_{0}^{t}(1+t')^{3}\leftexp{(i_{1}...i_{l})}{Q}_{l}(t')dt'
\end{align}
Taking in (8.402) the maximum over $i_{1},...,i_{l}$, recalling the definition (8.390), we obtain the integral inequality:
\begin{align}
 X_{l}(t)\leq B_{l}(t)+C(l+1)\delta_{0}(1+t)^{-3}[1+\log(1+t)]^{2}\int_{0}^{t}(1+t')[1+\log(1+t')]X_{l}(t')dt'
\end{align}
where:
\begin{equation}
 B_{l}(t)=\max_{i_{1}...i_{l}}\leftexp{(i_{1}...i_{l})}{B}_{l}(t)
\end{equation}

Setting:
\begin{equation}
 Y_{l}(t)=\int_{0}^{t}(1+t')[1+\log(1+t')]X_{l}(t')dt'
\end{equation}
then (8.404) becomes:
\begin{align}
 \frac{dY_{l}(t)}{dt}\leq (1+t)[1+\log(1+t)]B_{l}(t)
+C(l+1)\delta_{0}(1+t)^{-2}[1+\log(1+t)]^{3}Y_{l}(t)
\end{align}
which implies:
\begin{equation}
 Y_{l}(t)\leq e^{C_{l}(t)}\int_{0}^{t}e^{-C_{l}(t')}(1+t')[1+\log(1+t')]B_{l}(t')dt'
\end{equation}
where
\begin{equation}
 C_{l}(t)=C(l+1)\delta_{0}\int_{0}^{t}\frac{[1+\log(1+t')]^{3}}{(1+t')^{2}}dt'
\end{equation}
Obviously,
\begin{equation}
 0\leq C_{l}(t)\leq C'(l+1)\delta_{0}
\end{equation}
Hence, if
\begin{equation}
 \delta_{0}\leq\frac{\log 2}{C'(l+1)}
\end{equation}
(8.408) implies:
\begin{equation}
 Y_{l}(t)\leq 2\int_{0}^{t}(1+t')[1+\log(1+t')]B_{l}(t')dt'
\end{equation}
Substituting this in (8.402), recalling (8.406), we finally conclude that:
\begin{align}
 \||\leftexp{(i_{1}...i_{l})}{x}_{l}(t)|\|_{L^{2}([0,\epsilon_{0}]\times S^{2})}\leq \leftexp{(i_{1}...i_{l})}{B}_{l}(t)\\ \notag
+2C(l+1)\delta_{0}(1+t)^{-3}[1+\log(1+t)]^{2}\int_{0}^{t}(1+t')[1+\log(1+t')]B_{l}(t')dt'
\end{align}

\chapter{Regularization of the Propagation Equation for $\slashed{\Delta}\mu$.\\
 Estimates for the Top Order Spatial Derivatives of $\mu$}

\section{Regularization of the Propagation Equation}
In this chapter, we shall concentrate on the propagation equation (3.92) in Chapter 3:
\begin{equation}
 L\mu=m+\mu e
\end{equation}
where
\begin{equation}
 m=\frac{1}{2}\frac{dH}{dh}Th
\end{equation}
and
\begin{equation}
 e=\frac{1}{2\eta^{2}}(\frac{\rho}{\rho^{\prime}})^{\prime}Lh+\eta^{-1}\hat{T}^{i}(L\psi_{i})
\end{equation}
Since
\begin{align*}
 \frac{\rho}{\rho^{\prime}}=\frac{\rho}{d\rho/dh}=\frac{dp/dh}{d\rho/dh}=\eta^{2}
\end{align*}
therefore (9.3) takes the form:
\begin{align}
 e=\frac{1}{\eta}\frac{d\eta}{dh}Lh+\frac{1}{\eta}\hat{T}^{i}(L\psi_{i})
\end{align}

$\textbf{Lemma 9.1}$ The following commutation formula holds:
\begin{align*}
 [L,\slashed{\Delta}]\phi+\textrm{tr}\chi\slashed{\Delta}\phi=
-2\hat{\chi}\cdot\hat{\slashed{D}}^{2}\phi-2\slashed{\textrm{div}}\hat{\chi}\cdot\slashed{d}\phi
\end{align*}
and
\begin{align*}
 [T,\slashed{\Delta}]\phi+\kappa\textrm{tr}\theta\slashed{\Delta}\phi=-2\kappa\hat{\theta}\cdot\hat{\slashed{D}}^{2}\phi
-2\slashed{\textrm{div}}(\kappa\hat{\theta})\cdot\slashed{d}\phi
\end{align*}

Here, $\phi$ is an arbitrary function defined on $W^{*}_{\epsilon_{0}}$.

$Proof$. Let us work in acoustical coordinates $(t,u,\vartheta^{1},\vartheta^{2})$. We have:
\begin{align}
 \slashed{\Delta}\phi=(\slashed{g}^{-1})^{AB}\slashed{D}_{A}(\slashed{d}_{B}\phi)
=(\slashed{g}^{-1})^{AB}\{\frac{\partial^{2}\phi}{\partial\vartheta^{A}\partial\vartheta^{B}}-
\slashed{\Gamma}^{C}_{AB}\frac{\partial\phi}{\partial\vartheta^{C}}\}
\end{align}
Since in acoustical coordinates $L=\frac{\partial}{\partial t}$, we obtain:
\begin{align}
 L\slashed{\Delta}\phi=(\slashed{g}^{-1})^{AB}\{\frac{\partial^{2}}{\partial\vartheta^{A}\partial\vartheta^{B}}
(\frac{\partial\phi}{\partial t})-\slashed{\Gamma}^{C}_{AB}\frac{\partial}{\partial\vartheta^{C}}(\frac{\partial\phi}{\partial t})\}\\
+\frac{\partial}{\partial t}(\slashed{g}^{-1})^{AB}\slashed{D}_{A}(\slashed{d}_{B}\phi)-(\slashed{g}^{-1})^{AB}
\frac{\partial\slashed{\Gamma}^{C}_{AB}}{\partial t}\frac{\partial\phi}{\partial\vartheta^{C}}\notag
\end{align}
We have:
\begin{align}
 \frac{\partial}{\partial t}(\slashed{g}^{-1})^{AB}=-2\chi^{AB}
\end{align}
and:
\begin{align}
 \frac{\partial\slashed{\Gamma}^{C}_{AB}}{\partial t}=\slashed{D}_{A}\chi_{B}^{C}
+\slashed{D}_{B}\chi_{A}^{C}-\slashed{D}^{C}\chi_{AB}
\end{align}
Here, we used the fact that for a non-degenerate matrix $M$, we have:
\begin{align*}
 \frac{dM^{-1}}{dt}=-M^{-1}\frac{dM}{dt}M^{-1}
\end{align*}
and the way we use to derive (8.109).

Hence:
\begin{equation}
 (\slashed{g}^{-1})^{AB}\frac{\partial\slashed{\Gamma}^{C}_{AB}}{\partial t}=
2\slashed{\textrm{div}}\chi^{C}-\slashed{d}^{C}\textrm{tr}\chi=2\slashed{\textrm{div}}\hat{\chi}^{C}
\end{equation}
Substituting (9.6) and (9.8) in (9.5), we get:
\begin{equation}
 L\slashed{\Delta}\phi=\slashed{\Delta}L\phi-2\chi^{AB}\slashed{D}_{A}(\slashed{d}_{B}\phi)-
2(\slashed{\textrm{div}}\hat{\chi})^{A}\slashed{d}_{A}\phi
\end{equation}
 Writing then:
\begin{align*}
 2\chi^{AB}\slashed{D}_{A}(\slashed{d}_{B}\phi)=\textrm{tr}\chi\slashed{\Delta}\phi+2\hat{\chi}^{AB}\slashed{D}_{A}(\slashed{d}_{B}\phi)
\end{align*}
the first formula follows.
The second formula can be proved in a similar way. $\qed$

Now we apply the second identity by taking $\phi=h$ to obtain:
\begin{align}
 \slashed{\Delta}Th=T\slashed{\Delta}h+c_{T}
\end{align}
 where:
\begin{align}
 c_{T}=\kappa\textrm{tr}\theta\slashed{\Delta}h+2\kappa\hat{\theta}\cdot\hat{\slashed{D}}^{2}h+2\slashed{\textrm{div}}(\kappa\hat{\theta})\cdot\slashed{d}h
\end{align}
$c_{T}$ is of order 2 and regular as $\mu\rightarrow0$. The last term is a principal (i.e.order 2) acoustical term.

We now appeal to the wave equation for $h$ given in Chapter 8:
\begin{align}
 \mu\Box_{g}h=-\Omega^{-1}\frac{d\Omega}{dh}\check{a}-\check{b}
\end{align}
where
\begin{align}
 \check{a}=\mu a,\quad\check{b}=\mu b
\end{align}
are regular as $\mu\rightarrow0$. Since for any function $\phi$ we have:
\begin{align}
 \mu\Box_{g}\phi=\mu\slashed{\Delta}\phi-L\underline{L}\phi-\frac{1}{2}\textrm{tr}\chi\underline{L}\phi-\frac{1}{2}\textrm{tr}\underline{\chi}L\phi
-2\zeta\cdot\slashed{d}\phi
\end{align}
Then (9.13) reads:
\begin{align}
 \mu\slashed{\Delta}h=L\underline{L}h+\frac{1}{2}\textrm{tr}\chi\underline{L}h+\frac{1}{2}\textrm{tr}\underline{\chi}Lh+j
\end{align}
where:
\begin{align}
 j=2\zeta\cdot\slashed{d}h-\Omega^{-1}\frac{d\Omega}{dh}\check{a}-\check{b}
\end{align}
is of order 1 and regular as $\mu\rightarrow0$. Applying $T$ to (9.16) yields, in view of the fact that
\begin{align}
 [L,T]=\Lambda
\end{align}
we have:
\begin{align}
 T(\mu\slashed{\Delta}h)=L(T\underline{L}h)-\Lambda\cdot\slashed{d}\underline{L}h+
\frac{1}{2}T(\textrm{tr}\chi\underline{L}h)+\frac{1}{2}T(\textrm{tr}\underline{\chi}Lh)+Tj
\end{align}
On the other hand from (9.11) we have:
\begin{align}
 \mu\slashed{\Delta}Th=T(\mu\slashed{\Delta}h)-(T\mu)\slashed{\Delta}h+\mu c_{T}
\end{align}
Combining, we obtain:
\begin{align}
 \mu\slashed{\Delta}Th=L(T\underline{L}h)+\frac{1}{2}T(\textrm{tr}\chi\underline{L}h)+\frac{1}{2}T(\textrm{tr}\underline{\chi}Lh)+j^{\prime}
\end{align}
where:
\begin{align}
 j^{\prime}=Tj-\Lambda\cdot\slashed{d}\underline{L}h+\mu c_{T}-(T\mu)\slashed{\Delta}h
\end{align}
$j^{\prime}$ is of order 2 and regular as $\mu\rightarrow0$. The only principal acoustical term in $j^{\prime}$ is contributed by the last term in $c_{T}$
noted as above.

We have:
\begin{align}
 \slashed{\Delta}m=\frac{1}{2}\frac{dH}{dh}\slashed{\Delta}Th+w_{0}
\end{align}
 where
\begin{align}
 w_{0}=\frac{d^{2}H}{dh^{2}}\slashed{d}h\cdot\slashed{d}Th+\frac{1}{2}(\frac{d^{2}H}{dh^{2}}\slashed{\Delta}h+\frac{d^{3}H}{dh^{3}}|\slashed{d}h|^{2})Th
\end{align}
is of order 2 and regular as $\mu\rightarrow0$. Then by (9.21):
\begin{align}
 \mu\slashed{\Delta}m=Lf_{0}^{\prime}-\frac{1}{2}\frac{d^{2}H}{dh^{2}}(Lh)(T\underline{L}h)
+\frac{1}{4}\frac{dH}{dh}\{(Lh)T\textrm{tr}\underline{\chi}+(\underline{L}h)T\textrm{tr}\chi\}\\\notag
+\frac{1}{4}\frac{dH}{dh}\{\textrm{tr}\underline{\chi}(TLh)+\textrm{tr}\chi(T\underline{L}h)\}
+\frac{1}{2}\frac{dH}{dh}j^{\prime}+\mu w_{0}
\end{align}
where:
\begin{align}
 f^{\prime}_{0}=\frac{1}{2}\frac{dH}{dh}T\underline{L}h
\end{align}
here we note that the 3rd term on the right contains principal acoustical terms. In view of the equation for $\slashed{\mathcal{L}}_{T}\chi$ in Chapter 3,
the principal acoustical part of $T\textrm{tr}\chi$ is $\slashed{\Delta}\mu$. Also, since 
\begin{align}
 \underline{\chi}=\eta^{-1}\kappa(-\chi+2\slashed{k})
\end{align}
the principal acoustical part of $T\textrm{tr}\underline{\chi}$ is $-\eta^{-1}\kappa\slashed{\Delta}\mu$. It follows that the principal acoustical part of the 3rd 
term on the right in (9.25) is:
\begin{align}
 \frac{1}{4}\frac{dH}{dh}\{-\eta^{-1}\kappa(Lh)+(\underline{L}h)\}\slashed{\Delta}\mu
=\frac{1}{2}\frac{dH}{dh}(Th)\slashed{\Delta}\mu=m\slashed{\Delta}\mu
\end{align}
We can therefore write:
\begin{align}
 \mu\slashed{\Delta}m=Lf^{\prime}_{0}+\frac{1}{2}\textrm{tr}\chi f^{\prime}_{0}+m\slashed{\Delta}\mu+w^{\prime}_{0}
\end{align}
where:
\begin{align}
 w^{\prime}_{0}=-\frac{1}{2}\frac{d^{2}H}{dh^{2}}(Lh)(T\underline{L}h)
+\frac{1}{4}\frac{dH}{dh}(T\textrm{tr}\chi-\slashed{\Delta}\mu)\underline{L}h\\\notag
+\frac{1}{4}\frac{dH}{dh}(T\textrm{tr}\underline{\chi}+\eta^{-1}\kappa\slashed{\Delta}\mu)Lh
+\frac{1}{4}\frac{dH}{dh}\textrm{tr}\underline{\chi}TLh+\frac{1}{2}\frac{dH}{dh}j^{\prime}+\mu w_{0}
\end{align}
contains principal acoustical terms only through $j^{\prime}$.

Let us turn to $e$. We have:
\begin{align}
 \slashed{\Delta}e=\frac{1}{\eta}\frac{d\eta}{dh}\slashed{\Delta}Lh+\frac{1}{\eta}\hat{T}^{i}\slashed{\Delta}L\psi_{i}+w_{1}
\end{align}
where:
\begin{align}
 w_{1}=2\frac{d}{dh}(\frac{1}{\eta}\frac{d\eta}{dh})\slashed{d}h\cdot\slashed{d}Lh
+2[\frac{d}{dh}(\frac{1}{\eta})\hat{T}^{i}\slashed{d}h\cdot\slashed{d}L\psi_{i}+\frac{1}{\eta}\slashed{d}\hat{T}^{i}\cdot\slashed{d}L\psi_{i}]\\\notag
+[\frac{d}{dh}(\frac{1}{\eta}\frac{d\eta}{dh})\slashed{\Delta}h+\frac{d^{2}}{dh^{2}}(\frac{1}{\eta}\frac{d\eta}{dh})|\slashed{d}h|^{2}]Lh\\\notag
+\{[\frac{d}{dh}(\frac{1}{\eta})\slashed{\Delta}h+\frac{d^{2}}{dh^{2}}(\frac{1}{\eta})|\slashed{d}h|^{2}]\hat{T}^{i}+
2\frac{d}{dh}(\frac{1}{\eta})\slashed{d}h\cdot\slashed{d}\hat{T}^{i}+\frac{1}{\eta}\slashed{\Delta}\hat{T}^{i}\}L\psi_{i}
\end{align}
is of order 2 and regular as $\mu\rightarrow0$.
The only principal acoustical term in $w_{1}$ is the term
\begin{align}
 \frac{1}{\eta}\slashed{\Delta}\hat{T}^{i}L\psi_{i}
\end{align}
contained in the last term on the right in (9.32).

Applying Lemma 9.1 taking $\phi=h, \psi_{i}$ we then obtain:
\begin{align}
 \slashed{\Delta}e=\frac{1}{\eta}\frac{d\eta}{dh}L\slashed{\Delta}h+\frac{1}{\eta}\hat{T}^{i}L\slashed{\Delta}\psi_{i}+c_{L}+w_{1}
\end{align}
where
\begin{align}
 c_{L}=\frac{2}{\eta}\frac{d\eta}{dh}(\chi\cdot\slashed{D}^{2}h+\slashed{\textrm{div}}\hat{\chi}\cdot\slashed{d}h)
+\frac{2}{\eta}\hat{T}^{i}(\chi\cdot\slashed{D}^{2}\psi_{i}
+\slashed{\textrm{div}}\hat{\chi}\cdot\slashed{d}\psi_{i})
\end{align}
contains principal acoustical terms in the terms involving $\slashed{\textrm{div}}\hat{\chi}$. Defining:
\begin{align}
 f^{\prime}_{1}=\frac{1}{\eta}\frac{d\eta}{dh}\slashed{\Delta}h+\frac{1}{\eta}\hat{T}^{i}\slashed{\Delta}\psi_{i}
\end{align}
we then have:
\begin{align}
 \mu^{2}\slashed{\Delta}e=L(\mu^{2}f^{\prime}_{1})+\mu w_{1}^{\prime}
\end{align}
where:
\begin{align}
 w^{\prime}_{1}=-2(L\mu)f^{\prime}_{1}-\mu\{\frac{d}{dh}(\frac{1}{\eta}\frac{d\eta}{dh})(Lh)\slashed{\Delta}h+
[\frac{d}{dh}(\frac{1}{\eta})(Lh)\hat{T}^{i}+\frac{1}{\eta}L\hat{T}^{i}]\slashed{\Delta}\psi_{i}\}\\\notag
+\mu(c_{L}+w_{1})
\end{align}
is of order 2 and regular as $\mu\rightarrow0$.

We commute $\slashed{\Delta}$ with equation (9.1) and multiply it by $\mu$ to obtain:
\begin{align}
 L(\mu\slashed{\Delta}\mu)+(\mu\textrm{tr}\chi-L\mu)\slashed{\Delta}\mu=
-2\mu\hat{\chi}\cdot\hat{\slashed{\Delta}}^{2}\mu\\\notag
-\mu\slashed{d}\mu\cdot(\slashed{d}\textrm{tr}\chi+2i-2\slashed{d}e)
+\mu\slashed{\Delta}m+\mu^{2}\slashed{\Delta}e+e\mu\slashed{\Delta}\mu
\end{align}
Substituting (9.29) in (9.39), the term $m\slashed{\Delta}\mu$ in (9.29) combines with the last term on the right in (9.39) to $(L\mu)\slashed{\Delta}\mu$.
This is brought to the left hand side of (9.39), making the coefficient of $\slashed{\Delta}\mu$ on the left equal to $\mu\textrm{tr}\chi-2L\mu$. Also,
defining
\begin{align}
 \check{f}^{\prime}=f^{\prime}_{0}+\mu^{2}f^{\prime}_{1}
\end{align}
and
\begin{align}
 x^{\prime}=\mu\slashed{\Delta}\mu-\check{f}^{\prime}
\end{align}
the term $L\check{f}^{\prime}$ is brought to the left hand side of (9.39). Then left hand side becomes:
\begin{align}
 Lx^{\prime}+(\textrm{tr}\chi-2\mu^{-1}(L\mu))(x^{\prime}+f^{\prime})
\end{align}
and the resulting right hand side is:
\begin{align}
 \frac{1}{2}\textrm{tr}\chi f^{\prime}_{0}-2\mu\hat{\chi}\cdot\hat{\slashed{D}}^{2}\mu-\mu\slashed{d}\mu\cdot
(\slashed{d}\textrm{tr}\chi+2i-2\slashed{d}e)\\\notag
+w^{\prime}_{0}+\mu w^{\prime}_{1}
\end{align}
Defining then:
\begin{align}
 \check{g}^{\prime}=-\frac{\mu^{2}}{2}\textrm{tr}\chi f^{\prime}_{1}-\mu\slashed{d}\mu\cdot(\slashed{d}\textrm{tr}\chi+2i-2\slashed{d}e)
+w^{\prime}_{0}+\mu w^{\prime}_{1}
\end{align}
the resulting equation takes the form:
\begin{align}
 Lx^{\prime}+(\textrm{tr}\chi-2\mu^{-1}(L\mu))x^{\prime}=\\\notag
-(\frac{1}{2}\textrm{tr}\chi-2\mu^{-1}(L\mu))\check{f}^{\prime}
-2\mu\hat{\chi}\cdot\hat{\slashed{D}}^{2}\mu+\check{g}^{\prime}
\end{align}
The functions $\check{f}^{\prime}$, $\check{g}^{\prime}$ are of order 2 and regular as $\mu\rightarrow0$. The function $\check{f}^{\prime}$ does not contain
principal acoustical terms. From the above discussion, the principal acoustical part of $\check{g}^{\prime}$, is, besides the term $-\mu\slashed{d}\mu\cdot
\slashed{d}\textrm{tr}\chi$, contained in the last two terms in (9.43). In fact, the principal acoustical part of $w^{\prime}_{0}$ is contained in the term
\begin{align}
 \frac{1}{2}\mu\frac{dH}{dh}c_{T}
\end{align}
contributed by the term $\frac{1}{2}\frac{dH}{dh}j^{\prime}$ in (9.30). Since 
\begin{align}
 \theta=\slashed{k}-\eta^{-1}\chi
\end{align}
taking into account the Codazzi equation
\begin{align}
 \slashed{\textrm{div}}\hat{\chi}=\frac{1}{2}\slashed{d}\textrm{tr}\chi+i
\end{align}
this principal acoustical part is
\begin{align}
 -\frac{\mu^{2}}{\eta^{2}}\frac{dH}{dh}\slashed{d}h\cdot\slashed{d}\textrm{tr}\chi
\end{align}
The principal acoustical part of $\mu w^{\prime}_{1}$ is contained in the term
\begin{align}
 \mu^{2}(c_{L}+w_{1})
\end{align}
Taking into account the Codazzi equation the principal acoustical part of $\mu^{2}c_{L}$ is:
\begin{align}
 \frac{\mu^{2}}{\eta}(\frac{d\eta}{dh}\slashed{d}h+\hat{T}^{i}\slashed{d}\psi_{i})\cdot\slashed{d}\textrm{tr}\chi
\end{align}
Finally, the principal acoustical part of $\mu^{2}w_{1}$ is contained in
\begin{align}
 \frac{\mu^{2}}{\eta}\slashed{\Delta}\hat{T}^{i}L\psi_{i}
\end{align}
By the results in Chapter 3, the principal acoustical part of $\slashed{\Delta}\hat{T}^{i}$ is:
\begin{align}
 -\eta^{-1}\slashed{d}x^{i}\cdot\slashed{d}\textrm{tr}\chi
\end{align}
Now,
\begin{align}
 \frac{\partial x^{i}}{\partial\vartheta^{A}}L\psi_{i}=X^{i}_{A}L\psi_{i}=L^{\mu}X_{A}\psi_{\mu}
\end{align}
hence:
\begin{align}
 (\slashed{d}x^{i})L\psi_{i}=L^{\mu}\slashed{d}\psi_{\mu}=\slashed{d}\psi_{t}+L^{i}\slashed{d}\psi_{i}\\\notag
=\slashed{d}\psi_{t}-\psi^{i}\slashed{d}\psi_{i}-\eta\hat{T}^{i}\slashed{d}\psi_{i}
=\slashed{d}h-\eta\hat{T}^{i}\slashed{d}\psi_{i}
\end{align}
Therefore the principal acoustical part of (9.52) is:
\begin{align}
 \frac{\mu^{2}}{\eta}\{-\frac{1}{\eta}\slashed{d}h+\hat{T}^{i}\slashed{d}\psi_{i}\}\cdot\slashed{d}\textrm{tr}\chi
\end{align}
Combining the above results and noting that
\begin{align}
 -\frac{1}{2}\frac{dH}{dh}=1+\eta\frac{d\eta}{dh}
\end{align}
we conclude that the principal acoustical part of $\check{g}^{\prime}$ is:
\begin{align}
 \xi\cdot(\mu\slashed{d}\textrm{tr}\chi)
\end{align}
where
\begin{align}
 \xi=-\slashed{d}\mu+\frac{2\mu}{\eta}\hat{T}^{i}\slashed{d}\psi_{i}+\frac{2\mu}{\eta}\frac{d\eta}{dh}\slashed{d}h
\end{align}

\section{Propagation Equations for the Higher Order Spatial \\
Derivatives}
We now return to the propagation equation (9.45). To control the higher order spatial derivatives of $\mu$, we introduce the function:
\begin{equation}
 \leftexp{(i_{1}...i_{l})}{x}'_{m,l}=\mu R_{i_{l}}...R_{i_{1}}T^{m}\slashed{\Delta}\mu
-R_{i_{l}}...R_{i_{1}}T^{m}\check{f}'
\end{equation}
We have:
\begin{equation}
 x'_{0,0}=x'
\end{equation}
To derive a propagation equation for $x'_{m,l}$, we shall use the following lemma.

$\textbf{Lemma 9.2}$ Let $\alpha$ and $\beta$ be $S_{t,u}$ trace-free symmetric 2-covariant tensorfields defined in the spacetime 
domain $W^{*}_{\epsilon_{0}}$. Then we have:
\begin{align*}
 T(\alpha\cdot\beta)=(\hat{\slashed{\mathcal{L}}}_{T}\alpha)\cdot\beta+\alpha\cdot(\hat{\slashed{\mathcal{L}}}_{T}\beta)
-\textrm{tr}\leftexp{(T)}{\slashed{\pi}}(\alpha\cdot\beta)
\end{align*}
$Proof$. Since $T$ is tangential to $\Sigma_{t}$, we can confine ourselves to $\Sigma_{t}$. To define $\hat{\slashed{\mathcal{L}}}_{T}\alpha$
and $\hat{\slashed{\mathcal{L}}}_{T}\beta$, we extend $\alpha$ and $\beta$ to $T\Sigma_{t}$ as:
\begin{equation}
 \alpha(V,T)=\beta(V,T)=0 :  \forall V\in T\Sigma_{t}
\end{equation}
Since
\begin{equation}
 \bar{g}^{-1}=\hat{T}\otimes\hat{T}+(\slashed{g}^{-1})^{AB}X_{A}\otimes X_{B}
\end{equation}
we have:
\begin{equation}
 \alpha\cdot\beta=(\slashed{g}^{-1})^{AC}(\slashed{g}^{-1})^{BD}\alpha_{AB}\beta_{CD}
=(\bar{g}^{-1})^{ac}(\bar{g}^{-1})^{bd}\alpha_{ab}\beta_{cd}
\end{equation}
Then,
\begin{equation}
 T(\alpha\cdot\beta)=-\leftexp{(T)}{\bar{\pi}}^{ab}\gamma_{ab}+(\mathcal{L}_{T}\alpha)_{ab}\beta^{ab}
+\alpha^{ab}(\mathcal{L}_{T}\beta)_{ab}
\end{equation}
Here,
\begin{equation}
 \leftexp{(T)}{\bar{\pi}}_{ab}=(\mathcal{L}_{T}\bar{g})_{ab},\quad
\leftexp{(T)}{\bar{\pi}}^{ab}=(\bar{g}^{-1})^{ac}(\bar{g}^{-1})^{bd}\leftexp{(T)}{\bar{\pi}}_{cd}
\end{equation}
Here, $\gamma$ is a symmetric 2-covariant tensorfield on $\Sigma_{t}$ given by:
\begin{equation}
 \gamma_{ab}=(\bar{g}^{-1})^{cd}(\alpha_{ac}\beta_{bd}+\beta_{ac}\alpha_{bd})
\end{equation}
We have:
\begin{align*}
 \gamma_{ab}T^{b}=0
\end{align*}
so, we can view $\gamma$ as an $S_{t,u}$ symmetric 2-covariant tensorfield. Now for any two symmetric trace-free 2-dimensional 
matrices $A$ and $B$ we have:
\begin{align*}
 AB+BA-\textrm{tr}(AB)I=0
\end{align*}
In an orthonormal frame relative to $\slashed{g}$ the components of $\alpha$ and $\beta$ form $A$ and $B$. The components of $\gamma$
form $AB+BA$. It follows that:
\begin{equation}
 \gamma_{AB}=(\alpha\cdot\beta)\slashed{g}_{AB}
\end{equation}
Hence, we have:
\begin{equation}
 \leftexp{(T)}{\bar{\pi}}^{AB}\gamma_{AB}=\leftexp{(T)}{\slashed{\pi}}^{AB}\gamma_{AB}=\textrm{tr}\leftexp{(T)}{\slashed{\pi}}
(\alpha\cdot\beta)
\end{equation}
Since $\alpha$ and $\beta$ are trace-free, the lemma follows. $\qed$

Similarly, we have:

$\textbf{Lemma 9.3}$  Let $\alpha$ and $\beta$ be $S_{t,u}$ trace-free symmetric 2-covariant tensorfields defined in the spacetime 
domain $W^{*}_{\epsilon_{0}}$. Then we have:
\begin{align*}
 R_{i}(\alpha\cdot\beta)=(\hat{\slashed{\mathcal{L}}}_{R_{i}}\alpha)\cdot\beta+\alpha\cdot(\hat{\slashed{\mathcal{L}}}_{R_{i}}\beta)
-\textrm{tr}\leftexp{(R_{i})}{\slashed{\pi}}(\alpha\cdot\beta)
\end{align*}

   To simplify the derivation of the propagation equation for $x^{\prime}_{m,l}$, we introduce the functions:
\begin{equation}
 \leftexp{(i_{1}...i_{l})}{\check{f}'}_{m,l}=R_{i_{l}}...R_{i_{1}}T^{m}\check{f}'
\end{equation}
\begin{equation}
 \hat{\slashed{\mu}}_{2}=\hat{\slashed{D}}^{2}\mu
\end{equation}
and 
\begin{equation}
 \leftexp{(i_{1}...i_{l})}{\hat{\slashed{\mu}}}_{2,m,l}=
\hat{\slashed{\mathcal{L}}}_{R_{i_{l}}}...\hat{\slashed{\mathcal{L}}}_{R_{i_{1}}}(\hat{\slashed{\mathcal{L}}}_{T})^{m}\hat{\slashed{\mu}}_{2}
\end{equation}

$\textbf{Proposition 9.1}$ For each non-negative integer $m$ the function $x'_{m,0}$ satisfies the propagation equation:
\begin{align*}
 Lx'_{m,0}+(\textrm{tr}\chi-2\mu^{-1}(L\mu))x'_{m,0}=-(\frac{1}{2}\textrm{tr}\chi-2\mu^{-1}(L\mu))\check{f}'_{m,0}
-2\mu\hat{\chi}\cdot\hat{\slashed{\mu}}_{2,m,0}+g'_{m,0}
\end{align*}
where $g'_{m,0}$ is given by:
\begin{align*}
 g'_{m,0}=T^{m}\check{g}'+\sum_{k=0}^{m-1}T^{k}\Lambda x'_{m-k-1,0}+\sum_{k=0}^{m-1}T^{k}y'_{m-k,0}
\end{align*}
Here, for each $j=1,...,m$, $y'_{j,0}$ is the function:
\begin{align*}
 y'_{j,0}=-(T\mu)a'_{j-1,0}+(TL\mu-\mu T\textrm{tr}\chi-\Lambda\mu)T^{j-1}\slashed{\Delta}\mu\\
+\frac{1}{2}(T\textrm{tr}\chi)\check{f}'_{j-1,0}-2\mu(\hat{\slashed{\mathcal{L}}}_{T}\hat{\chi}-\textrm{tr}\leftexp{(T)}{\slashed{\pi}}\hat{\chi})
\cdot\hat{\slashed{\mu}}_{2,j-1,0}
\end{align*}
where
\begin{align*}
 a'_{j-1,0}=LT^{j-1}\slashed{\Delta}\mu+\textrm{tr}\chi T^{j-1}\slashed{\Delta}\mu+2\hat{\chi}\cdot\hat{\slashed{\mu}}_{2,j-1,0}
\end{align*}

$Proof$. The proposition reduces for $m=0$ to the propagation equation (9.45). 
By induction on $m$, assuming that the propagation equation holds with $m$ replaced by $m-1$:
\begin{align}
  Lx'_{m-1,0}+(\textrm{tr}\chi-2\mu^{-1}(L\mu))x'_{m-1,0}=-(\frac{1}{2}\textrm{tr}\chi-2\mu^{-1}(L\mu))\check{f}'_{m-1,0}
-2\mu\hat{\chi}\cdot\hat{\slashed{\mu}}_{2,m-1,0}+g'_{m-1,0}
\end{align}
 for some $g^{\prime}_{m-1,0}$, we shall show that a propagation equation of the form given by the proposition holds for $m$, where $g^{\prime}_{m,0}$ is related to 
$g^{\prime}_{m-1,0}$ by a certain recursion relation.
 Rewrite the term
\begin{align*}
 2\mu^{-1}(L\mu)(x'_{m-1,0}+\check{f}'_{m-1,0})
\end{align*}
in (9.73) as:
\begin{align*}
 2(L\mu)T^{m-1}\slashed{\Delta}\mu
\end{align*}
(see (9.60) and (9.70)) obtaining the equation:
\begin{align}
 Lx'_{m-1,0}+\textrm{tr}\chi x'_{m-1,0}=2(L\mu)T^{m-1}\slashed{\Delta}\mu-\frac{1}{2}\textrm{tr}\chi\check{f}'_{m-1,0}
-2\mu\hat{\chi}\cdot\hat{\slashed{\mu}}_{2,m-1,0}+g'_{m-1,0}
\end{align}
We apply $T$ to this equation. Since:
\begin{align*}
 TT^{m-1}\slashed{\Delta}\mu=T^{m}\slashed{\Delta}\mu\\
T\check{f}'_{m-1,0}=\check{f}'_{m,0}
\end{align*}
and by Lemma 9.2:
\begin{equation}
 T(\hat{\chi}\cdot\hat{\slashed{\mu}}_{2,m-1,0})=(\hat{\slashed{\mathcal{L}}}_{T}\hat{\chi})\cdot\hat{\slashed{\mu}}_{2,m-1,0}
+\hat{\chi}\cdot\hat{\slashed{\mu}}_{2,m,0}-\textrm{tr}\leftexp{(T)}{\slashed{\pi}}(\hat{\chi}\cdot\hat{\slashed{\mu}}_{2,m-1,0})
\end{equation}
(see definition (9.72)), we obtain:
\begin{align}
 TLx'_{m-1,0}+\textrm{tr}\chi Tx'_{m-1,0}+(T\textrm{tr}\chi)x'_{m-1,0}\\\notag
=2(L\mu)T^{m}\slashed{\Delta}\mu+2(TL\mu)T^{m-1}\slashed{\Delta}\mu\\\notag
-\frac{1}{2}\textrm{tr}\chi\check{f}'_{m,0}-\frac{1}{2}(T\textrm{tr}\chi)\check{f}'_{m-1,0}\\\notag
-2\mu\hat{\chi}\cdot\hat{\slashed{\mu}}_{2,m,0}-2\mu(\hat{\slashed{\mathcal{L}}}_{T}\hat{\chi})\cdot\hat{\slashed{\mu}}_{2,m-1,0}
+2(\mu\textrm{tr}\leftexp{(T)}{\slashed{\pi}}-(T\mu))(\hat{\chi}\cdot\hat{\slashed{\mu}}_{2,m-1,0})\\\notag
+Tg'_{m-1,0}\notag
\end{align}
Using the fact:
\begin{equation}
 [L,T]=\Lambda=-(\slashed{g}^{-1})^{AB}(\zeta_{B}+\eta_{B})X_{A}
\end{equation}
we express:
\begin{equation}
 TLx'_{m-1,0}=LTx'_{m-1,0}-\Lambda x'_{m-1,0}
\end{equation}
From the definition (9.60) with $l=0$ and $m$ replaced by $m-1$:
\begin{align*}
 x'_{m-1,0}=\mu T^{m-1}\slashed{\Delta}\mu-T^{m-1}\check{f}'
\end{align*}
we get:
\begin{equation}
 Tx'_{m-1,0}=x'_{m,0}+(T\mu)T^{m-1}\slashed{\Delta}\mu
\end{equation}
Applying $L$ to (9.79) and expressing:
\begin{align*}
 LT\mu=TL\mu+\Lambda\mu
\end{align*}
yields:
\begin{align}
 LTx'_{m-1,0}=Lx'_{m,0}+(T\mu)LT^{m-1}\slashed{\Delta}\mu+(TL\mu+\Lambda\mu)T^{m-1}\slashed{\Delta}\mu
\end{align}
Substituting (9.80) to (9.78) and the result in (9.76), substituting also (9.79) on the left, we obtain a propagation equation of the required form with:
\begin{equation}
 g'_{m,0}=Tg'_{m-1,0}+\Lambda x'_{m-1,0}+y'_{m,0}
\end{equation}
Applying then Proposition 8.2, the proposition follows. $\qed$

$\textbf{Proposition 9.2}$ For each pair of non-negative integers $(m,l)$ and each multi-index
$(i_{1}...i_{l})$, the function $\leftexp{(i_{1}...i_{l})}{x}'_{m,l}$ satisfies the propagation equation:
\begin{align*}
 L\leftexp{(i_{1}...i_{l})}{x}'_{m,l}+(\textrm{tr}\chi-2\mu^{-1}(L\mu))\leftexp{(i_{1}...i_{l})}{x}'_{m,l}=
-(\frac{1}{2}\textrm{tr}\chi-2\mu^{-1}(L\mu))\leftexp{(i_{1}...i_{l})}{\check{f}'}_{m,l}\\
-2\mu\hat{\chi}\cdot\leftexp{(i_{1}...i_{l})}{\hat{\slashed{\mu}}}_{2,m,l}+\leftexp{(i_{1}...i_{l})}{g}'_{m,l}
\end{align*}
where
\begin{align*}
 \leftexp{(i_{1}...i_{l})}{g}'_{m,l}=R_{i_{l}}...R_{i_{1}}g'_{m,0}+\sum_{k=0}^{l-1}R_{i_{l}}...R_{i_{l-k+1}}
\leftexp{(R_{i_{l-k}})}{Z}\leftexp{(i_{1}...i_{l-k-1})}{x}'_{m,l-k-1}\\
+\sum_{k=0}^{l-1}R_{i_{l}}...R_{i_{l-k+1}}\leftexp{(i_{1}...i_{l-k})}{y}'_{m,l-k}
\end{align*}
where $g'_{m,0}$ is given by Proposition 9.1, and for each $j=1,...,l$, 
\begin{align*}
 \leftexp{(i_{1}...i_{j})}{y}'_{m,j}=-(R_{i_{j}}\mu)\leftexp{(i_{1}...i_{j-1})}{a}'_{m,j-1}\\
+(R_{i_{j}}L\mu-\mu R_{i_{j}}\textrm{tr}\chi-\leftexp{(R_{i_{j}})}{Z}\mu)(R_{i_{j-1}}...R_{i_{1}}T^{m}\slashed{\Delta}\mu)\\
+\frac{1}{2}(R_{i_{j}}\textrm{tr}\chi)\leftexp{(i_{1}...i_{j-1})}{\check{f}'}_{m,j-1}
-2\mu(\hat{\slashed{\mathcal{L}}}_{R_{i_{j}}}\hat{\chi}-\textrm{tr}\leftexp{(R_{i_{j}})}{\slashed{\pi}}\hat{\chi})
\cdot\leftexp{(i_{1}...i_{j-1})}{\hat{\slashed{\mu}}}_{2,m,j-1}
\end{align*}
where:
\begin{align*}
 \leftexp{(i_{1}...i_{j-1})}{a}'_{m,j-1}=L(R_{i_{j-1}}...R_{i_{1}}T^{m}\slashed{\Delta}\mu)\\
+\textrm{tr}\chi(R_{i_{j-1}}...R_{i_{1}}T^{m}\slashed{\Delta}\mu)+2\hat{\chi}\cdot\leftexp{(i_{1}...i_{j-1})}{\hat{\slashed{\mu}}}_{2,m,j-1}
\end{align*}

$Proof$. The proposition reduces for $l=0$ to Proposition 9.1. By induction on $l$, assuming the proposition holds with $l$ replaced by
$l-1$:
\begin{align}
 L\leftexp{(i_{1}...i_{l})}{x}'_{m,l-1}+(\textrm{tr}\chi-2\mu^{-1}(L\mu))\leftexp{(i_{1}...i_{l-1})}{x}'_{m,l-1}=
-(\frac{1}{2}\textrm{tr}\chi-2\mu^{-1}(L\mu))\leftexp{(i_{1}...i_{l-1})}{\check{f}'}_{m,l-1}\\\notag
-2\mu\hat{\chi}\cdot\leftexp{(i_{1}...i_{l-1})}{\hat{\slashed{\mu}}}_{2,m,l-1}+\leftexp{(i_{1}...i_{l-1})}{g}'_{m,l-1}
\end{align}
 for some $\leftexp{(i_{1}...i_{l-1})}{g}^{\prime}_{m,l-1}$, we shall show that a propagation equation of the form given by the proposition holds for $l$, 
where $\leftexp{(i_{1}...i_{l})}{g}^{\prime}_{m,l}$ is related to $\leftexp{(i_{1}...i_{l-1})}{g}^{\prime}_{m,l-1}$ by a certain recursion relation.
   Rewrite the term
\begin{align*}
 2\mu^{-1}(L\mu)(\leftexp{(i_{1}...i_{l-1})}{x}'_{m,l-1}+\leftexp{(i_{1}...i_{l-1})}{\check{f}'}_{m,l-1})
\end{align*}
in (9.82) as:
\begin{align*}
 2(L\mu)R_{i_{l-1}}...R_{i_{1}}T^{m}\slashed{\Delta}\mu
\end{align*}
to obtain the equation:
\begin{align}
 L\leftexp{(i_{1}...i_{l})}{x}'_{m,l-1}+\textrm{tr}\chi\leftexp{(i_{1}...i_{l-1})}{x}'_{m,l-1}\\\notag
=2(L\mu)R_{i_{l-1}}...R_{i_{l}}T^{m}\slashed{\Delta}\mu-\frac{1}{2}\textrm{tr}\chi\leftexp{(i_{1}...i_{l-1})}{\check{f}'}_{m,l-1}\\\notag
-2\mu\hat{\chi}\cdot\leftexp{(i_{1}...i_{l-1})}{\hat{\slashed{\mu}}}_{2,m,l-1}+\leftexp{(i_{1}...i_{l-1})}{g}'_{m,l-1}\notag
\end{align}
We now apply $R_{i_{l}}$ to this equation. Since
\begin{align*}
 R_{i_{l}}\leftexp{(i_{1}...i_{l-1})}{\check{f}'}_{m,l-1}=\leftexp{(i_{1}...i_{l})}{\check{f}'}_{m,l}
\end{align*}
and by Lemma 9.3
\begin{align}
 R_{i_{l}}(\hat{\chi}\cdot\hat{\slashed{\mu}}_{2,m,l-1})=(\hat{\slashed{\mathcal{L}}}_{R_{i_{l}}}\hat{\chi})\cdot\leftexp{(i_{1}...i_{l-1})}
{\hat{\slashed{\mu}}}_{2,m,l-1}+\hat{\chi}\leftexp{(i_{1}...i_{l})}{\hat{\slashed{\mu}}}_{2,m,l}-\textrm{tr}\leftexp{(R_{i_{l}})}{\slashed{\pi}}
(\hat{\chi}\cdot\leftexp{(i_{1}...i_{l-1})}{\hat{\slashed{\mu}}}_{2,m,l-1})
\end{align}
we obtain:
\begin{align}
 R_{i_{l}}L\leftexp{(i_{1}...i_{l-1})}{x}'_{m,l-1}+\textrm{tr}\chi R_{i_{l}}\leftexp{(i_{1}...i_{l-1})}{x}'_{m,l-1}
+(R_{i_{l}}\textrm{tr}\chi)\leftexp{(i_{1}...i_{l-1})}{x}'_{m,l-1}\\\notag
=2(L\mu)R_{i_{l}}...R_{i_{1}}T^{m}\slashed{\Delta}\mu+2(R_{i_{l}}L\mu)R_{i_{l-1}}...R_{i_{1}}T^{m}\slashed{\Delta}\mu\\\notag
-\frac{1}{2}\textrm{tr}\chi\leftexp{(i_{1}...i_{l})}{\check{f}'}_{m,l}-\frac{1}{2}(R_{i_{l}}\textrm{tr}\chi)\leftexp{(i_{1}...i_{l-1})}{\check{f}'}_{m,l-1}\\\notag
-2\mu\hat{\chi}\cdot\leftexp{(i_{1}...i_{l})}{\hat{\slashed{\mu}}}_{2,m,l}-2\mu(\hat{\slashed{\mathcal{L}}}_{R_{i_{l}}}\hat{\chi})\cdot
\leftexp{(i_{1}...i_{l-1})}{\hat{\slashed{\mu}}}_{2,m,l-1}\\\notag
+2(\mu\textrm{tr}\leftexp{(R_{i_{l}})}{\slashed{\pi}}-(R_{i_{l}}\mu))(\hat{\chi}\cdot\leftexp{(i_{1}...i_{l-1})}{\hat{\slashed{\mu}}}_{2,m,l-1})\\\notag
+R_{i_{l}}\leftexp{(i_{1}...i_{l-1})}{g}'_{m,l-1}\notag
\end{align}
By Lemma 8.2 we have:
\begin{equation}
 [L,R_{i_{l}}]=\leftexp{(R_{i_{l}})}{Z}
\end{equation}
hence:
\begin{equation}
 R_{i_{l}}L\leftexp{(i_{1}...i_{l-1})}{x}'_{m,l-1}=LR_{i_{l}}\leftexp{(i_{1}...i_{l-1})}{x}'_{m,l-1}
-\leftexp{(R_{i_{l}})}{Z}\leftexp{(i_{1}...i_{l-1})}{x}'_{m,l-1}
\end{equation}
From the definition of $\leftexp{(i_{1}...i_{l-1})}{x}'_{m,l-1}$
we obtain:
\begin{equation}
 R_{i_{l}}\leftexp{(i_{1}...i_{l-1})}{x}'_{m,l-1}=\leftexp{(i_{1}...i_{l})}{x}'_{m,l}+(R_{i_{l}}\mu)(R_{i_{l-1}}...R_{i_{1}}T^{m}\slashed{\Delta}\mu)
\end{equation}
Applying $L$ to (9.88) and expressing:
\begin{align*}
 LR_{i_{l}}\mu=R_{i_{l}}L\mu+\leftexp{(R_{i_{l}})}{Z}
\end{align*}
yields:
\begin{align}
 LR_{i_{l}}\leftexp{(i_{1}...i_{l-1})}{x}'_{m,l-1}=L\leftexp{(i_{1}...i_{l})}{x}'_{m,l}+
(R_{i_{l}}\mu)L(R_{i_{l-1}}...R_{i_{1}}T^{m}\slashed{\Delta}\mu)\\\notag
+(R_{i_{l}}L\mu+\leftexp{(R_{i_{l}})}{Z}\mu)(R_{i_{l-1}}...R_{i_{1}}T^{m}\slashed{\Delta}\mu)
\end{align}
We substitute (9.89) in (9.87) and the result in (9.85), substituting also (9.88) on the left, a propagation equation of the required form results, with:
\begin{align}
 \leftexp{(i_{1}...i_{l})}{g}'_{m,l}=R_{i_{l}}\leftexp{(i_{1}...i_{l-1})}{g}'_{m,l-1}+
\leftexp{(R_{i_{l}})}{Z}\leftexp{(i_{1}...i_{l-1})}{x}'_{m,l-1}+\leftexp{(i_{1}...i_{l})}{y}'_{m,l}
\end{align}
Applying then Proposition 8.2, the proposition follows. $\qed$

   The function $\leftexp{(i_{1}...i_{j-1})}{a}'_{m,j-1}$ defined in Proposition 9.2 reduces for $m=0, j=1$, to the function:
\begin{equation}
 a'_{0,0}=L\slashed{\Delta}\mu+\textrm{tr}\chi\slashed{\Delta}\mu+2\hat{\chi}\cdot\hat{\slashed{\mu}}_{2}
\end{equation}
According to equation:
\begin{align*}
 L(\slashed{\Delta}\mu)+(\textrm{tr}\chi-e)\slashed{\mu}=-2\hat{\chi}\cdot\hat{\slashed{D}}^{2}\mu-\slashed{d}\mu\cdot
(\slashed{d}\textrm{tr}\chi+2i-2\slashed{d}e)+\slashed{\Delta}m+\mu\slashed{\Delta}e
\end{align*}
which is derived by Lemma 9.1, we have:
\begin{equation}
 a'_{0,0}=e\slashed{\Delta}\mu-2\slashed{d}\mu\cdot\slashed{\textrm{div}}\hat{\chi}+\slashed{\Delta}m+\mu\slashed{\Delta}e+2\slashed{d}\mu\cdot\slashed{d}e
\end{equation}

$\textbf{Lemma 9.4}$ For each non-negative integer $m$ we have:
\begin{equation}
 a'_{m,0}=T^{m}a'_{0,0}+b'_{m,0}\notag
\end{equation}
where:
\begin{align*}
 b'_{m,0}=\sum_{k=0}^{m-1}T^{k}\Lambda T^{m-k-1}\slashed{\Delta}\mu+\sum_{k=0}^{m-1}T^{k}c'_{m-k}
\end{align*}
and for each $j=1,...,m$,
\begin{align*}
 c'_{j}=-(T\textrm{tr}\chi)T^{j-1}\slashed{\Delta}\mu-2(\hat{\slashed{\mathcal{L}}}_{T}\hat{\chi}-\textrm{tr}\leftexp{(T)}{\slashed{\pi}}\hat{\chi})
\cdot\hat{\slashed{\mu}}_{2,j-1,0}
\end{align*}
Moreover, for each pair of non-negative integers $(m,l)$ and each multi-index $(i_{1}...i_{l})$ we have:
\begin{align*}
 \leftexp{(i_{1}...i_{l})}{a}'_{m,l}=R_{i_{l}}...R_{i_{1}}T^{m}a'_{0,0}+\leftexp{(i_{1}...i_{l})}{b}'_{m,l}
\end{align*}
where:
\begin{align*}
 \leftexp{(i_{1}...i_{l})}{b}'_{m,l}=R_{i_{l}}...R_{i_{1}}b'_{m,0}+\sum_{k=0}^{l-1}R_{i_{l}}...R_{i_{l-k+1}}\leftexp{(R_{i_{l-k}})}{Z}
R_{i_{l-k-1}}...R_{i_{1}}T^{m}\slashed{\Delta}\mu\\
+\sum_{k=0}^{l-1}R_{i_{l}}...R_{i_{l-k+1}}\leftexp{(i_{1}...i_{l-k})}{c}''_{m,l-k}
\end{align*}
and for each $j=1,...,l$ and each multi-index $(i_{1}...i_{j})$,
\begin{align*}
 \leftexp{(i_{1}...i_{j})}{c}''_{m,j}=-(R_{i_{j}}\textrm{tr}\chi)(R_{i_{j-1}}...R_{i_{1}}T^{m}\slashed{\Delta}\mu)
-2(\hat{\slashed{\mathcal{L}}}_{R_{i_{j}}}\hat{\chi}-\textrm{tr}\leftexp{(R_{i_{j}})}{\slashed{\pi}}\hat{\chi})
\cdot\leftexp{(i_{1}...i_{j-1})}{\hat{\slashed{\mu}}}_{2,m,j-1}
\end{align*}
$Proof$. To prove the first part, we note that it holds trivially for $m=0$ with:
\begin{equation}
 b'_{0,0}=0
\end{equation}
 Thus, by induction on $m$, assume that:
\begin{equation}
 a'_{m-1,0}=T^{m-1}a'_{0,0}+b'_{m-1,0}
\end{equation}
holds for some $b'_{m-1,0}$. Applying $T$ to this we obtain:
\begin{equation}
 Ta'_{m-1,0}=T^{m}a'_{0,0}+Tb'_{m-1,0}
\end{equation}
On the other hand, from the definition:
\begin{equation}
 a'_{m-1,0}=LT^{m-1}\slashed{\Delta}\mu+\textrm{tr}\chi T^{m-1}\slashed{\Delta}\mu+2\hat{\chi}\cdot\hat{\slashed{\mu}}_{2,m-1,0}
\end{equation}
Applying $T$ to this and using Lemma 9.2 yields:
\begin{align}
 Ta'_{m-1,0}=TLT^{m-1}\slashed{\Delta}\mu+\textrm{tr}\chi T^{m}\slashed{\Delta}\mu+2\hat{\chi}\cdot\hat{\slashed{\mu}}_{2,m,0}\\\notag
+(T\textrm{tr}\chi)T^{m-1}\slashed{\Delta}\mu+2(\hat{\slashed{\mathcal{L}}}_{T}\hat{\chi}-\textrm{tr}\leftexp{(T)}{\slashed{\pi}}\hat{\chi})
\cdot\hat{\slashed{\mu}}_{2,m-1,0}
\end{align}
Writing
\begin{align*}
 TLT^{m-1}\slashed{\Delta}\mu=LT^{m}\slashed{\Delta}\mu-\Lambda T^{m-1}\slashed{\Delta}\mu,
\end{align*}
 (9.97), in view of the definition of $a'_{m,0}$, reads:
\begin{equation}
 Ta'_{m-1,0}=a'_{m,0}-\Lambda T^{m-1}\slashed{\Delta}\mu-c'_{m}
\end{equation}
where $c_{m}^{\prime}$ is as in the statement of the lemma.
Equating the two expressions for $Ta^{\prime}_{m-1,0}$, that is, (9.95) with (9.98), we obtain the formula for $a^{\prime}_{m,0}$, if we set:
\begin{equation}
 b'_{m,0}=Tb'_{m-1,0}+\Lambda T^{m-1}\slashed{\Delta}\mu+c'_{m}
\end{equation}
 Using then Proposition 8.2, the result follows.

   To prove the second part, we apply induction on $l$. Assume:
\begin{equation}
 a'_{m,l-1}=R_{i_{l-1}}...R_{i_{1}}T^{m}a'_{0,0}+\leftexp{(i_{1}...i_{l-1})}{b}'_{m,l-1}
\end{equation}
Applying $R_{i_{l}}$ to this we obtain:
\begin{equation}
 R_{i_{l}}a'_{m,l-1}=R_{i_{l}}...R_{i_{1}}T^{m}a'_{0,0}+R_{i_{l}}\leftexp{(i_{1}...i_{l-1})}{b}'_{m,l-1}
\end{equation}
On the other hand, from the definition in Proposition 9.2,
\begin{equation}
 \leftexp{(i_{1}...i_{l-1})}{a}'_{m,l-1}=L(R_{i_{l-1}}...R_{i_{1}}T^{m}\slashed{\Delta}\mu)
+\textrm{tr}\chi(R_{i_{l-1}}...R_{i_{1}}T^{m}\slashed{\Delta}\mu)+2\hat{\chi}\cdot\leftexp{(i_{1}...i_{l-1})}{\hat{\slashed{\mu}}}_{2,m,l-1}
\end{equation}
Applying $R_{i_{l}}$ to this and using Lemma 9.3:
\begin{align}
 R_{i_{l}}\leftexp{(i_{1}...i_{l-1})}{a}'_{m,l-1}=R_{i_{l}}L(R_{i_{l-1}}...R_{i_{1}}T^{m}\slashed{\Delta}\mu)
+\textrm{tr}\chi(R_{i_{l}}...R_{i_{1}}T^{m}\slashed{\Delta}\mu)\\\notag
+2\hat{\chi}\cdot\leftexp{(i_{1}...i_{l})}{\hat{\slashed{\mu}}}_{2,m,l}
+(R_{i_{l}}\textrm{tr}\chi)(R_{i_{l-1}}...R_{i_{1}}T^{m}\slashed{\Delta}\mu)\\\notag
+2(\hat{\slashed{\mathcal{L}}}_{R_{i_{l}}}\hat{\chi}-\textrm{tr}\leftexp{(R_{i_{l}})}{\slashed{\pi}}\hat{\chi})\cdot
\leftexp{(i_{1}...i_{l-1})}{\hat{\slashed{\mu}}}_{2,m,l-1}\notag
\end{align}
Writing:
\begin{align*}
 R_{i_{l}}L(R_{i_{l-1}}...R_{i_{1}}T^{m}\slashed{\Delta}\mu)=
L(R_{i_{l}}...R_{i_{1}}T^{m}\slashed{\Delta}\mu)-\leftexp{(R_{i_{l}})}{Z}R_{i_{l-1}}...R_{i_{1}}T^{m}\slashed{\Delta}\mu
\end{align*}
(9.103) becomes, in view of the definition of $a^{\prime}_{m,l}$,:
\begin{equation}
 R_{i_{l}}\leftexp{(i_{1}...i_{l-1})}{a}'_{m,l-1}=\leftexp{(i_{1}...i_{l})}{a}'_{m,l}-\leftexp{(R_{i_{l}})}{Z}
R_{i_{l-1}}...R_{i_{1}}T^{m}\slashed{\Delta}\mu-c''_{m,l}
\end{equation}
Equating the two expressions for $R_{i_{l}}\leftexp{(i_{1}...i_{l-1})}{a}^{\prime}_{m,l-1}$, that is (9.101) with (9.104), we obtain the formula for 
$\leftexp{(i_{1}...i_{l})}{a}^{\prime}_{m,l}$, if we set:
\begin{equation}
 \leftexp{(i_{1}...i_{l})}{b}'_{m,l}=R_{i_{l}}\leftexp{(i_{1}...i_{l-1})}{b}'_{m,l-1}
+\leftexp{(R_{i_{l}})}{Z}R_{i_{l-1}}...R_{i_{1}}T^{m}\slashed{\Delta}\mu+c''_{m,l}
\end{equation}
Using again Proposition 8.2, the result follows. $\qed$

   Next, we shall investigate the order of various terms in the propagation equation of Proposition 9.2. In obtaining estimates for the $n+1$st 
order spatial derivatives of $\mu$, of which at least two are angular, we are to set $m+l=n-1$ in Proposition 9.2. 
The principal terms shall be of order $m+l+2=n+1$, and the principal acoustical terms shall be the principal terms in the spatial derivatives of 
$\chi$ and $\mu$.

   From (9.26), (9.36) and (9.40), we know that $\check{f}'$ is of order 2, and contains no acoustical part of order 2. It follows that
$\leftexp{(i_{1}...i_{l})}{\check{f}}_{m,l}$ is of order $m+l+2$, but contains no principal acoustical part.

   We turn to $\leftexp{(i_{1}...i_{l})}{g}^{\prime}_{m,l}$. From the expression of $\leftexp{(i_{1}...i_{l})}{g}^{\prime}_{m,l}$, we must investigate 
$g'_{m,0}$. From (9.44), (9.36), (9.38), (9.35), (9.32), (9.30) and (9.24) we know that $\check{g}'$ is of order 2, and its principal acoustical part 
is given by (9.58), 
consisting of 1st angular derivatives of $\textrm{tr}\chi$ multiplied by $\mu$. It follows that
\begin{equation}
 T^{m}\check{g}'
\end{equation}
the first term in $g'_{m,0}$ is of order $m+2$ and its principal acoustical part consists of 
\begin{equation}
 \xi\cdot(\mu\slashed{d}T^{m}\textrm{tr}\chi)
\end{equation}
or, in view of equation (3.125)-(3.126),
\begin{align}
 \xi\cdot(\mu\slashed{d}\textrm{tr}\chi) : m=0\\\notag
\xi\cdot(\mu\slashed{d}T^{m-1}\slashed{\Delta}\mu) : m\geq 1
\end{align}
to principal acoustical terms. In view of (9.60) and (8.34), (9.108) can be expressed to principal acoustical parts in terms of:
\begin{align}
 \xi\cdot x_{0} : m=0 \\\notag
\xi\cdot\slashed{d}x'_{m-1,0}  :  m\geq 1
\end{align}

   The second term in $g'_{m,0}$ can be written as:
\begin{equation}
 \sum_{k=0}^{m-1}\Lambda T^{k}x'_{m-k-1,0}-\sum_{k=0}^{m-1}[\Lambda,T^{k}]x'_{m-k-1,0}
\end{equation}
The commutator sum is of lower order $m+1$, while
\begin{align}
 T^{k}x'_{m-k-1,0}=T^{k}(\mu T^{m-k-1}\slashed{\Delta}\mu-T^{m-k-1}\check{f}')\\\notag
=x'_{m-1,0}+\sum_{j=1}^{k}C_{k}^{j}(T^{j}\mu)(T^{m-j-1}\slashed{\Delta}\mu)
\end{align}
The contribution of the second term in (9.111) to the first term of (9.110) is of lower order $m+1$, 
while the contribution of the first term is:
\begin{equation}
 m\Lambda x'_{m-1,0}
\end{equation}
This is the principal part of the second term in $g'_{m,0}$.

   The third term in $g'_{m,0}$ is 
\begin{equation}
 \sum_{k=0}^{m-1}T^{k}y'_{m-k,0}
\end{equation}
Consider the expression for $y'_{j,0}$. Here, all terms except the first term:
\begin{equation}
 -(T\mu)a'_{j-1,0}=-(T\mu)(T^{j-1}a'_{0,0}+b'_{j-1,0})
\end{equation}
are of lower order $j+1$. The principal term in $a'_{0,0}$ is 
\begin{align*}
 \slashed{\Delta}m+\mu\slashed{\Delta}e
\end{align*}
It follows that $T^{j-1}a'_{0,0}$ is of principal order $j+2$ but contains no principal acoustical part, its principal term being:
\begin{align*}
 T^{j-1}\slashed{\Delta}m+\mu T^{j-1}\slashed{\Delta}e
\end{align*}
From Lemma 9.4,
\begin{equation}
 b'_{j-1,0}=\sum_{k=0}^{j-2}T^{k}\Lambda T^{j-k-2}\slashed{\Delta}\mu+\sum_{k=0}^{j-2}T^{k}c'_{j-k-1}
\end{equation}
The first sum in (9.115) is of order $j+1$, while from the expression of $c'_{j}$, the second sum in (9.115) is of order $j$.
We conclude that the principal part of the third term in $g'_{m,0}$ is 
\begin{equation}
 -m(T\mu)T^{m-1}a'_{0,0}
\end{equation}
This is of principal order $m+2$, but contains no principal acoustical part. This completes the investigation of $g'_{m,0}$.

  We conclude from above that the contributions of (9.109), (9.112) and (9.116) to 
\begin{equation}
 R_{i_{l}}...R_{i_{1}}g'_{m,0}
\end{equation}
are all of principal order $l+m+2$. The contribution of (9.109) to (9.117) is
\begin{align}
 \xi\cdot\leftexp{(i_{1}...i_{l})}{x}_{l} : m=0\\\notag
\xi\cdot\slashed{d}\leftexp{(i_{1}...i_{l})}{x}'_{m-1,l} : m\geq 1
\end{align}
where $\xi$ is given by (9.59). By $\textbf{H0}$,
\begin{align}
 |\xi\cdot\leftexp{(i_{1}...i_{l})}{x}_{l}|\leq |\xi||\leftexp{(i_{1}...i_{l})}{x}_{l}|\\\notag
|\xi\cdot\slashed{d}\leftexp{(i_{1}...i_{l})}{x}'_{m-1,l}|\leq C(1+t)^{-1}|\xi|\max_{j}
|\leftexp{(i_{1}...i_{l}j)}{x}'_{m-1,l+1}|
\end{align}
The contribution of (9.112) is:
\begin{equation}
 m\Lambda\leftexp{(i_{1}...i_{l})}{x}'_{m-1,l}
\end{equation}
From $\textbf{H0}$, this is bounded by:
\begin{equation}
 Cm(1+t)^{-1}|\Lambda|\max_{j}|\leftexp{(i_{1}...i_{l}j)}{x}'_{m-1,l+1}|
\end{equation}
Finally, the contribution of (9.116) is:
\begin{equation}
 -m(T\mu)R_{i_{l}}...R_{i_{1}}T^{m-1}a'_{0,0}
\end{equation}
This does not contain principal acoustical part. This completes the investigation of the first term of $\leftexp{(i_{1}...i_{l})}{g}'_{m,l}$.

   We turn to the second term of $\leftexp{(i_{1}...i_{l})}{g}'_{m,l}$. This is:
\begin{equation}
 \sum_{k=0}^{l-1}R_{i_{l}}...R_{i_{l-k+1}}\leftexp{(R_{i_{l-k}})}{Z}\leftexp{(i_{1}...i_{l-k-1})}{x}'_{m,l-k-1}
\end{equation}
which we write as:
\begin{equation}
 \sum_{k=0}^{l-1}\leftexp{(R_{i_{l-k}})}{Z}R_{i_{l}}...R_{i_{l-k+1}}\leftexp{(i_{1}...i_{l-k-1})}{x}'_{m,l-k-1}
\end{equation}
minus:
\begin{equation}
 \sum_{k=0}^{l-1}[\leftexp{(R_{i_{l-k}})}{Z},R_{i_{l}}...R_{i_{l-k+1}}]\leftexp{(i_{1}...i_{l-k-1})}{x}'_{m,l-k-1}
\end{equation}
Obviously, the commutator term is of lower order $l+m+1$, while the principal part of (9.135) is:
\begin{equation}
 \sum_{k=0}^{l-1}\leftexp{(R_{i_{l-k}})}{Z}\leftexp{(i_{1}...i_{l-k-1}i_{l-k+1}...i_{l})}{x}'_{m,l-1}
\end{equation}
This is a principal acoustical term. By $\textbf{H0}$, it can be bounded by:
\begin{equation}
 C(1+t)^{-1}(\sum_{k=1}^{l}|\leftexp{(R_{i_{k}})}{Z}|)\max_{j}|\leftexp{(i_{1}...i_{k-1}i_{k+1}...i_{l}j)}{x}'_{m,l}|
\end{equation}

   Finally, we consider the third term in $\leftexp{(i_{1}...i_{l})}{g}'_{m,l}$. This is 
\begin{equation}
 \sum_{k=0}^{l-1}R_{i_{l}}...R_{i_{l-k+1}}\leftexp{(i_{1}...i_{l-k})}{y}'_{m,l-k}
\end{equation}
Note that only the first term in $\leftexp{(i_{1}...i_{j})}{y}'_{m,j}$ is of principal order $m+l+2$:
\begin{equation}
 -(R_{i_{j}}\mu)\leftexp{(i_{1}...i_{j-1})}{a}'_{m,j-1}=
-(R_{i_{j}}\mu)(R_{i_{j-1}}...R_{i_{1}}T^{m}a'_{0,0}+\leftexp{(i_{1}...i_{j-1})}{b}'_{m,j-1})
\end{equation}
Since $a'_{0,0}$ is of order 3, but contains no acoustical part of order 3, it follows that $R_{i_{j-1}}...R_{i_{1}}T^{m}a'_{0,0}$ is 
of principal order $m+j+2$ but contains no principal acoustical part. $\leftexp{(i_{1}...i_{j-1})}{b}'_{m,j-1}$ is given by Lemma 9.4:
\begin{align}
 \leftexp{(i_{1}...i_{j-1})}{b}'_{m,j-1}=R_{i_{j-1}}...R_{i_{1}}b'_{m,0}\\\notag
+\sum_{k=0}^{j-2}R_{i_{j-1}}...R_{i_{j-k}}\leftexp{(R_{i_{j-k-1}})}{Z}R_{i_{j-k-2}}...R_{i_{1}}T^{m}\slashed{\Delta}\mu\\\notag
+\sum_{k=0}^{j-2}R_{i_{j-1}}...R_{i_{j-k}}\leftexp{(i_{1}...i_{j-k-1})}{c}''_{m,j-k-1}\notag
\end{align}
By the discussions after (9.115), we know that the first term on the right of (9.130) is of order $m+j+1$. Obviously, the second term is also of 
order $m+j+1$. While $\leftexp{(i_{1}...i_{j})}{c}''_{m,j}$, given in Lemma 9.4, is of order $m+j+1$, hence, 
the third term on the right of (9.130) is of order $m+j$. Thus $\leftexp{(i_{1}...i_{j-1})}{b}'_{m,j-1}$ is of lower order $m+j+1$.
We conclude that the principal part of (9.128) is:
\begin{equation}
 -\sum_{k=0}^{l-1}(R_{i_{l-k}}\mu)R_{i_{l}}...R_{i_{l-k+1}}R_{i_{l-k-1}}...R_{i_{1}}T^{m}a'_{0,0}
\end{equation}
This is of principal order $m+l+2$, but contains no principal acoustical part. This completes the investigation of $\leftexp{(i_{1}...i_{l})}{g}'_{m,l}$.

\section{Elliptic Theory on $S_{t,u}$}
Returning to the propagation equation for $\leftexp{(i_{1}...i_{l})}{x}^{\prime}_{m,l}$ of Proposition 9.2, the term 
$2\mu\hat{\chi}\cdot\leftexp{(i_{1}...i_{l})}{\hat{\slashed{\mu}}}_{2,m,l}$ remains to be considered. Of course, 
$\leftexp{(i_{1}...i_{l})}{\hat{\slashed{\mu}}}_{2,m,l}$ is a principal acoustical term, it involves the $m+l$th order
spatial derivatives of $\hat{\slashed{\mu}}_{2}$, not $\slashed{\Delta}\mu$, which is  what the propagation equation  for 
$\leftexp{(i_{1}...i_{l})}{x}^{\prime}_{m,l}$ allows us to control. The term in question comes from (9.45), namely, from the term
\begin{align*}
 2\mu\hat{\chi}\cdot\hat{\slashed{{D}}}^{2}\mu=2\mu\hat{\chi}\cdot\hat{\slashed{\mu}}_{2}
\end{align*}
Now, the propagation equation must be considered in conjunction with the definition (9.41):
\begin{equation}
 \mu\slashed{\Delta}\mu=x'+\check{f}'
\end{equation}
which is an elliptic equation for $\mu$ on each surface $S_{t,u}$. Then we can use this equation to estimate on $S_{t,u}$, 
$\slashed{d}\mu$ as well as $\hat{\slashed{\mu}}_{2}$  in terms of $x'$ and the 2nd derivatives of $\psi_{\mu}$. On the 
other hand, (9.45) allows us to estimate $x'$ along each generator of $C_{u}$ in terms of $\hat{\slashed{\mu}}_{2}$ and 2nd 
derivatives of $\psi_{\mu}$. Thus we can estimate $\slashed{d}\mu$ and $\slashed{D}^{2}\mu$ in terms of 2nd derivatives of 
$\psi_{\mu}$.
 
   Similarly, the propagation equation for $\leftexp{(i_{1}...i_{l})}{x}'_{m,l}$ must be considered in conjunction with
an elliptic equation for 
\begin{equation}
 \leftexp{(i_{1}...i_{l})}{\mu}_{m,l}=R_{i_{l}}...R_{i_{1}}T^{m}\mu
\end{equation}
then we can estimate $\slashed{D}^{2}\leftexp{(i_{1}...i_{l})}{\mu}_{m,l}$ in terms of $m+l+2$nd derivatives of $\psi_{\mu}$. 
The required elliptic equation will be derived from (9.60):
\begin{align}
 \mu R_{i_{l}}...R_{i_{1}}T^{m}\slashed{\Delta}\mu=\leftexp{(i_{1}...i_{l})}{x}'_{m,l}+\leftexp{(i_{1}...i_{l})}{\check{f}'}_{m,l}
\end{align}
with the help of the following lemmas.

$\textbf{Lemma 9.5}$ Let $(M,g)$ be a 2-dimensional Riemannian manifold, let $X$ be an arbitrary vectorfield and $f$ an arbitrary function
on $M$. Then the following commutation formulas hold:
\begin{align*}
 X(\Delta_{g}f)-\Delta_{g}(Xf)=-\leftexp{(X)}{\pi}^{ab}(\nabla^{2}f)_{ab}-\textrm{tr}\leftexp{(X)}{\pi}^{b}_{1}d_{b}f
\end{align*}
and:
\begin{align*}
 (\hat{\mathcal{L}}_{X}\hat{\nabla}^{2}f)_{ab}-(\hat{\nabla}^{2}(Xf))_{ab}=-\frac{1}{2}\leftexp{(X)}{\hat{\pi}}_{ab}\Delta_{g}f-
\leftexp{(X)}{\hat{\pi}}^{c}_{1,ab}d_{c}f
\end{align*}
Here,
\begin{align*}
 \textrm{tr}\leftexp{(X)}{\pi}^{b}_{1}=\nabla^{a}\leftexp{(X)}{\pi}_{a}^{b}-\frac{1}{2}\nabla^{b}\textrm{tr}\leftexp{(X)}{\pi}
\end{align*}
\begin{align*}
 \leftexp{(X)}{\hat{\pi}}^{c}_{1,ab}=\frac{1}{2}(\nabla_{a}\leftexp{(X)}{\hat{\pi}}_{b}^{c}+\nabla_{b}\leftexp{(X)}{\hat{\pi}}_{a}^{c}
-\nabla^{c}\leftexp{(X)}{\hat{\pi}}_{ab}-g_{ab}\nabla^{d}\leftexp{(X)}{\hat{\pi}}_{d}^{c})\\
+\frac{1}{4}(\delta_{a}^{c}d_{b}\textrm{tr}\leftexp{(X)}{\pi}+\delta_{b}^{c}d_{a}\textrm{tr}\leftexp{(X)}{\pi}-g_{ab}d^{c}\textrm{tr}
\leftexp{(X)}{\pi})
\end{align*}
where $\mathcal{L}_{X}g=\leftexp{(X)}{\pi}$.

$Proof$. Let $\phi_{t}$ be the local 1-parameter group generated by $X$ and let $\phi_{t*}$ be the corresponding pull back.
Consider an arbitrary 1-form $\alpha$ on $M$. We then have
\begin{equation}
 \phi_{t*}(\overset{g}{\nabla}\alpha)=\overset{\phi_{t*}g}{\nabla}(\phi_{t*}\alpha)
\end{equation}
Now, in an arbitrary coordinates we have:
\begin{equation}
 (\overset{g}{\nabla}\alpha)_{ab}=\frac{\partial\alpha_{b}}{\partial x^{a}}-\overset{g}{\Gamma^{c}_{ab}}\alpha_{c}
\end{equation}
where $\overset{g}{\Gamma^{c}_{ab}}$ is the Christoffel symbol of metric $g$ in these coordinates.
Similarly, in the same coordinates, 
\begin{equation}
 (\overset{\phi_{t*}g}{\nabla}(\phi_{t*}\alpha))_{ab}=\frac{\partial(\phi_{t*}\alpha)_{b}}{\partial x^{a}}-\overset{\phi_{t*}g}{\Gamma^{c}_{ab}}
(\phi_{t*}\alpha)_{c}
\end{equation}
Differentiating (9.137) with respect to $t$ at $t=0$ we obtain:
\begin{align}
 (\frac{d}{dt}(\overset{\phi_{t*}g}{\nabla}(\phi_{t*}\alpha))_{ab})_{t=0}=
\frac{\partial}{\partial x^{a}}(\frac{d}{dt}(\phi_{t*}\alpha)_{b})_{t=0}-\overset{g}{\Gamma^{c}_{ab}}
(\frac{d}{dt}(\phi_{t*}\alpha)_{c})_{t=0}-(\frac{d}{dt}\overset{\phi_{t*}g}{\Gamma^{c}_{ab}})_{t=0}\alpha_{c}
\end{align}
Now, from the definition,
\begin{align*}
 (\frac{d}{dt}(\phi_{t*}\alpha))_{t=0}=\mathcal{L}_{X}\alpha
\end{align*}
By (9.135), the left hand-side of (9.138) is the $ab$-component of:
\begin{align*}
 (\frac{d}{dt}(\phi_{t*}(\overset{g}{\nabla}\alpha)))_{t=0}=\mathcal{L}_{X}(\nabla\alpha)
\end{align*}
Recalling (8.109), we have:
\begin{equation}
 (\frac{d}{dt}\overset{\phi_{t*}g}{\Gamma^{c}_{ab}})_{t=0}=\leftexp{(X)}{\pi}^{c}_{1,ab}
\end{equation}
where:
\begin{equation}
 \leftexp{(X)}{\pi}^{c}_{1,ab}=\frac{1}{2}(\nabla_{a}\leftexp{(X)}{\pi_{b}^{c}}+\nabla_{b}\leftexp{(X)}{\pi}^{c}_{a}
-\nabla^{c}\leftexp{(X)}{\pi}_{ab})
\end{equation}
we conclude that:
\begin{equation}
 (\mathcal{L}_{X}(\nabla\alpha))_{ab}=(\nabla(\mathcal{L}_{X}\alpha))_{ab}-\leftexp{(X)}{\pi}^{c}_{1,ab}\alpha_{c}
\end{equation}
In particular, in the case $\alpha=df$ for some function $f$, since
\begin{align*}
 \mathcal{L}_{X}(df)=d(Xf)
\end{align*}
we obtain:
\begin{equation}
 (\mathcal{L}_{X}(\nabla^{2}f))_{ab}=(\nabla^{2}(Xf))_{ab}-\leftexp{(X)}{\pi}^{c}_{1,ab}d_{c}f
\end{equation}
for any function $f$ on $M$. Since:
\begin{align*}
 \Delta_{g}f=(g^{-1})^{ab}(\nabla^{2}f)_{ab}
\end{align*}
it follows that:
\begin{equation}
 X(\Delta_{g}f)=\mathcal{L}_{X}(\Delta_{g}f)=-\leftexp{(X)}{\pi}^{ab}(\nabla^{2}f)_{ab}
+(g^{-1})^{ab}(\mathcal{L}_{X}(\nabla^{2}f))_{ab}
\end{equation}
Since
\begin{align*}
 (g^{-1})^{ab}(\nabla^{2}(Xf))_{ab}=\Delta_{g}(Xf)
\end{align*}
then by (9.142) and (9.143), the first part of the lemma follows. In fact, the first part holds for an arbitrary 
$n$-dimensional manifold $M$.

   For the second part, we consider the case of 2-dimensional manifold $M$. The trace-free part of $\nabla^{2}f$ is:
\begin{align*}
 \hat{\nabla}^{2}f=\nabla^{2}f-\frac{1}{2}g\Delta_{g}f
\end{align*}
hence:
\begin{equation}
 \mathcal{L}_{X}(\hat{\nabla}^{2}f)=\mathcal{L}_{X}(\nabla^{2}f)-\frac{1}{2}\leftexp{(X)}{\pi}\Delta_{g}f-
\frac{1}{2}gX(\Delta_{g}f)
\end{equation}
Substituting from (9.142) and the first part of the lemma we then obtain:
\begin{align}
 (\mathcal{L}_{X}(\hat{\nabla}^{2}f))_{ab}=(\nabla^{2}(Xf))_{ab}-\leftexp{(X)}{\pi}^{c}_{1,ab}d_{c}f-\frac{1}{2}\leftexp{(X)}{\pi}_{ab}\Delta_{g}f\\\notag
-\frac{1}{2}g_{ab}\{\Delta_{g}(Xf)-\leftexp{(X)}{\pi}^{cd}(\nabla^{2}f)_{cd}-\textrm{tr}\leftexp{(X)}{\pi}^{c}_{1}d_{c}f\}
\end{align}
Taking the trace-free part on both sides of above, then we get the second part of the lemma. $\qed$

$\textbf{Lemma 9.6}$ Let $f$ be an arbitrary function defined on a given hypersurface $\Sigma_{t}$. Then the following commutation formulas hold:
\begin{align*}
 T(\slashed{\Delta}f)-\slashed{\Delta}(Tf)=-\leftexp{(T)}{\slashed{\pi}}^{AB}(\slashed{D}^{2}f)_{AB}-
\textrm{tr}\leftexp{(T)}{\slashed{\pi}}^{B}_{1}\slashed{d}_{B}f
\end{align*}
and:
\begin{align*}
 (\hat{\slashed{\mathcal{L}}}_{T}\hat{\slashed{D}}^{2}f)_{AB}-(\hat{\slashed{D}}^{2}(Tf))_{AB}
=-\frac{1}{2}\leftexp{(T)}{\hat{\slashed{\pi}}}_{AB}\slashed{\Delta}f-\leftexp{(T)}{\hat{\slashed{\pi}}}^{C}_{1,AB}\slashed{d}_{C}f
\end{align*}
Here,
\begin{align*}
 \leftexp{(T)}{\slashed{\pi}}^{C}_{1,AB}=\frac{1}{2}\{\slashed{D}_{A}\leftexp{(T)}{\slashed{\pi}}_{B}^{C}
+\slashed{D}_{B}\leftexp{(T)}{\slashed{\pi}}_{A}^{C}-\slashed{D}^{C}\leftexp{(T)}{\slashed{\pi}}_{AB}\}
\end{align*}
is the Lie derivative with respect to $T$ of the induced connection on $S_{t,u}$,
and the trace in the lemma is with respect to the induced metric on $S_{t,u}$.

$Proof$. Since $T$ is tangential to $\Sigma_{t}$, we can confine attention to $\Sigma_{t}$. Also, we can choose acoustical coordinates 
$(\vartheta_{1},\vartheta_{2})$ so that $\Xi=0$ on the given $\Sigma_{t}$. So on this $\Sigma_{t}$:
\begin{equation}
 T=\frac{\partial}{\partial u} 
\end{equation}
Consider an arbitrary $S_{t,u}$ 1-form $\alpha$ defined on $\Sigma_{t}$. In terms of $(u,\vartheta_{1},\vartheta_{2})$:
\begin{equation}
 (\slashed{D}\alpha)_{AB}=\frac{\partial\alpha_{B}}{\partial\vartheta^{A}}-\slashed{\Gamma}^{C}_{AB}\alpha_{C}
\end{equation}
Differentiating (9.147) with respect to $u$ we obtain:
\begin{equation}
 \frac{\partial}{\partial u}(\slashed{D}\alpha)_{AB}=\frac{\partial}{\partial\vartheta^{A}}
(\frac{\partial\alpha_{B}}{\partial u})-\slashed{\Gamma}^{C}_{AB}\frac{\partial\alpha_{C}}{\partial u}
-\frac{\partial\slashed{\Gamma}^{C}_{AB}}{\partial u}\alpha_{C}
\end{equation}
In view of (9.146),
\begin{equation}
 \frac{\partial\alpha_{A}}{\partial u}=(\slashed{\mathcal{L}}_{T}\alpha)_{A}, \quad
\frac{\partial}{\partial u}(\slashed{D}\alpha)_{AB}=(\slashed{\mathcal{L}}_{T}(\slashed{D}\alpha))_{AB}
\end{equation}
Moreover, from Chapter 3 we have, in view of (9.146),
\begin{equation}
 \frac{\partial\slashed{g}_{AB}}{\partial u}=\leftexp{(T)}{\slashed{\pi}}_{AB}=2\kappa\theta_{AB},\quad
\frac{\partial}{\partial u}(\slashed{g}^{-1})^{AB}=-\leftexp{(T)}{\slashed{\pi}}^{AB}=-2\kappa\theta^{AB}
\end{equation}
hence:
\begin{align*}
 \frac{\partial\slashed{\Gamma}^{C}_{AB}}{\partial u}=\frac{1}{2}\frac{\partial}{\partial u}
\{(\slashed{g}^{-1})^{CD}(\frac{\partial\slashed{g}_{BD}}{\partial\vartheta_{A}}+\frac{\partial\slashed{g}_{AD}}{\partial\vartheta_{B}}
-\frac{\partial\slashed{g}_{AB}}{\partial\vartheta^{D}})\}\\
=-\leftexp{(T)}{\slashed{\pi}}^{CD}\slashed{\Gamma}^{E}_{AB}\slashed{g}_{DE}+\frac{1}{2}(\slashed{g}^{-1})^{CD}
(\frac{\partial\leftexp{(T)}{\slashed{\pi}}_{BD}}{\partial\vartheta^{A}}+\frac{\partial\leftexp{(T)}{\slashed{\pi}}_{AD}}{\partial\vartheta^{B}}
-\frac{\partial\leftexp{(T)}{\slashed{\pi}}_{AB}}{\partial\vartheta^{D}})\\
=\frac{1}{2}(\slashed{g}^{-1})^{CD}(\slashed{D}_{A}\leftexp{(T)}{\slashed{\pi}}_{BD}+\slashed{D}_{B}\leftexp{(T)}{\slashed{\pi}}_{AD}
-\slashed{D}_{D}\leftexp{(T)}{\slashed{\pi}}_{AB})
\end{align*}
that is:
\begin{equation}
 \frac{\partial\slashed{\Gamma}^{C}_{AB}}{\partial u}=\leftexp{(T)}{\slashed{\pi}}^{C}_{1,AB}
\end{equation}
where:
\begin{equation}
 \leftexp{(T)}{\slashed{\pi}}^{C}_{1,AB}=\frac{1}{2}\{\slashed{D}_{A}\leftexp{(T)}{\slashed{\pi}}_{B}^{C}
+\slashed{D}_{B}\leftexp{(T)}{\slashed{\pi}}_{A}^{C}-\slashed{D}^{C}\leftexp{(T)}{\slashed{\pi}}_{AB}\}
\end{equation}
Substituting (9.151) in (9.148), and in view of (9.149), we then obtain:
\begin{equation}
 (\slashed{\mathcal{L}}_{T}(\slashed{D}\alpha))_{AB}=(\slashed{D}(\slashed{\mathcal{L}}_{T}\alpha))_{AB}-
\leftexp{(T)}{\slashed{\pi}^{C}_{1,AB}}\alpha_{C}
\end{equation}
This holds for any $S_{t,u}$ 1-form $\alpha$ defined on $\Sigma_{t}$. In particular, it holds in the case $\alpha=\slashed{d}f$ for some 
function $f$ defined on $\Sigma_{t}$. In this case, we have,
\begin{align*}
 \alpha_{A}=\frac{\partial f}{\partial\vartheta^{A}}, \quad \frac{\partial\alpha_{A}}{\partial u}
=\frac{\partial}{\partial\vartheta^{A}}(\frac{\partial f}{\partial u})
\end{align*}
i.e.
\begin{equation}
 \slashed{\mathcal{L}}_{T}(\slashed{d}f)=\slashed{d}(Tf)
\end{equation}
So by (9.153) we have:
\begin{equation}
 (\slashed{\mathcal{L}}_{T}(\slashed{D}^{2}f))_{AB}=(\slashed{D}^{2}(Tf))_{AB}-\leftexp{(T)}{\slashed{\pi}}^{C}_{1,AB}\slashed{d}_{C}f
\end{equation}
for any function $f$ on $\Sigma_{t}$. Since:
\begin{align*}
 \slashed{\Delta}f=(\slashed{g}^{-1})^{AB}(\slashed{D}^{2}f)_{AB}
\end{align*}
it follows by (9.150) that:
\begin{equation}
 T(\slashed{\Delta}f)=\slashed{\mathcal{L}}_{T}(\slashed{\Delta}f)=-\leftexp{(T)}{\slashed{\pi}}^{AB}
(\slashed{D}^{2}f)_{AB}+(\slashed{g}^{-1})^{AB}(\slashed{\mathcal{L}}_{T}(\slashed{D}^{2}f))_{AB}
\end{equation}
Since
\begin{align*}
 (\slashed{g}^{-1})^{AB}(\slashed{D}^{2}(Tf))_{AB}=\slashed{\Delta}(Tf)
\end{align*}
the first part of the lemma follows.

  To prove the second part, we note that the trace-free part of $\slashed{D}^{2}f$ is:
\begin{align*}
 \hat{\slashed{D}}^{2}f=\slashed{D}^{2}f-\frac{1}{2}\slashed{g}\slashed{\Delta}f
\end{align*}
hence:
\begin{equation}
 \slashed{\mathcal{L}}_{T}(\hat{\slashed{D}}^{2}f)=\slashed{\mathcal{L}}_{T}(\slashed{D}^{2}f)-\frac{1}{2}\leftexp{(T)}{\slashed{\pi}}
\slashed{\Delta}f-\frac{1}{2}\slashed{g}T(\slashed{\Delta}f)
\end{equation}
So using (9.155) and the first part of the lemma we obtain:
\begin{align}
 (\slashed{\mathcal{L}}_{T}(\hat{\slashed{D}}^{2}f))_{AB}=(\slashed{D}^{2}(Tf))_{AB}-
\leftexp{(T)}{\slashed{\pi}}^{C}_{1,AB}\slashed{d}_{C}f-\frac{1}{2}\leftexp{(T)}{\slashed{\pi}}_{AB}\slashed{\Delta}f\\\notag
-\frac{1}{2}\slashed{g}_{AB}\{\slashed{\Delta}(Tf)-\leftexp{(T)}{\slashed{\pi}}^{CD}(\slashed{D}^{2}f)_{CD}-\textrm{tr}
\leftexp{(T)}{\slashed{\pi}}^{C}_{1}\slashed{d}_{C}f\}
\end{align}
Then taking the trace-free part on the above, we get the second part of the lemma. $\qed$

$\textbf{Proposition 9.3}$ For each pair of non-negative integers $(m,l)$ and each multi-index $(i_{1}...i_{l})$ we have:
\begin{align*}
 \slashed{\Delta}\leftexp{(i_{1}...i_{l})}{\mu}_{m,l}-R_{i_{l}}...R_{i_{1}}T^{m}\slashed{\Delta}\mu=\leftexp{(i_{1}...i_{l})}{d}_{m,l}
\end{align*}
where:
\begin{align*}
 \leftexp{(i_{1}...i_{l})}{d}_{m,l}=R_{i_{l}}...R_{i_{1}}d_{m,0}+\sum_{k=0}^{l-1}R_{i_{l}}...R_{i_{l-k+1}}
(\leftexp{(R_{i_{l-k}})}{\slashed{\pi}}\cdot\slashed{D}^{2}\leftexp{(i_{1}...i_{l-k-1})}{\mu}_{m,l-k-1})\\
+\sum_{k=0}^{l-1}R_{i_{l}}...R_{i_{l-k+1}}(\textrm{tr}\leftexp{(R_{i_{l-k}})}{\slashed{\pi}}_{1}\cdot\slashed{d}
\leftexp{(i_{1}...i_{l-k-1})}{\mu}_{m,l-k-1})
\end{align*}
and:
\begin{align*}
 d_{m,0}=\sum_{k=0}^{m-1}T^{k}(\leftexp{(T)}{\slashed{\pi}}\cdot\slashed{D}^{2}\mu_{m-k-1,0})
+\sum_{k=0}^{m-1}T^{k}(\textrm{tr}\leftexp{(T)}{\slashed{\pi}}_{1}\cdot\slashed{d}\mu_{m-k-1,0})
\end{align*}
$Proof$. We first consider the case $l=0$. If also $m=0$ we have trivially:
\begin{equation}
 d_{0,0}=0
\end{equation}
We shall derive a recursion formula for $d_{m,0}$. Consider then:
\begin{align*}
 d_{m-1,0}=\slashed{\Delta}\mu_{m-1,0}-T^{m-1}\slashed{\Delta}\mu
\end{align*}
Applying $T$ to this we obtain:
\begin{equation}
 Td_{m-1,0}=T(\slashed{\Delta}\mu_{m-1,0})-T^{m}\slashed{\Delta}\mu
\end{equation}
From the first part of Lemma 9.6,
\begin{equation}
 T(\slashed{\Delta}\mu_{m-1,0})=\slashed{\Delta}\mu_{m,0}-\leftexp{(T)}{\slashed{\pi}}\cdot\slashed{D}^{2}\mu_{m-1,0}-
\textrm{tr}\leftexp{(T)}{\slashed{\pi}}_{1}\cdot\slashed{d}\mu_{m-1,0}
\end{equation}
Substituting (9.161) in (9.160) and noting that:
\begin{align*}
 \slashed{\Delta}\mu_{m,0}-T^{m}\slashed{\Delta}\mu=d_{m,0}
\end{align*}
we obtain the recursion relation:
\begin{equation}
 d_{m,0}=Td_{m-1,0}+\leftexp{(T)}{\slashed{\pi}}\cdot\slashed{D}^{2}\mu_{m-1,0}
+\textrm{tr}\leftexp{(T)}{\slashed{\pi}}_{1}\cdot\slashed{d}\mu_{m-1,0}
\end{equation}
Applying Proposition 8.2, we get the expression for $d_{m,0}$.

   Next, we consider the case $l\geq 1$. We have:
\begin{align*}
 \leftexp{(i_{1}...i_{l-1})}{d}_{m,l-1}=\slashed{\Delta}\leftexp{(i_{1}...i_{l-1})}{\mu}_{m,l-1}-R_{i_{l-1}}...R_{i_{1}}T^{m}\slashed{\Delta}\mu
\end{align*}
Applying $R_{i_{l}}$ to this we obtain:
\begin{equation}
 R_{i_{l}}\leftexp{(i_{1}...i_{l-1})}{d}_{m,l-1}=R_{i_{l}}\slashed{\Delta}\leftexp{(i_{1}...i_{l-1})}{\mu}_{m,l-1}
-R_{i_{l}}...R_{i_{1}}T^{m}\slashed{\Delta}\mu
\end{equation}
From the first part of Lemma 9.5 we have:
\begin{align}
 R_{i_{l}}\slashed{\Delta}\leftexp{(i_{1}...i_{l-1})}{\mu}_{m,l-1}=\\\notag
\slashed{\Delta}\leftexp{(i_{1}...i_{l})}{\mu}_{m,l}-\leftexp{(R_{i_{l}})}{\slashed{\pi}}\cdot\slashed{D}^{2}
\leftexp{(i_{1}...i_{l-1})}{\mu}_{m,l-1}-\textrm{tr}\leftexp{(R_{i_{l}})}{\slashed{\pi}}_{1}\cdot\slashed{d}
\leftexp{(i_{1}...i_{l-1})}{\mu}_{m,l-1}
\end{align}
Substituting (9.164) in (9.163) and noting that:
\begin{align*}
 \slashed{\Delta}\leftexp{(i_{1}...i_{l})}{\mu}_{m,l}-R_{i_{l}}...R_{i_{1}}T^{m}\slashed{\Delta}\mu=\leftexp{(i_{1}...i_{l})}{d}_{m,l}
\end{align*}
we obtain the recursion relation:
\begin{align}
 \leftexp{(i_{1}...i_{l})}{d}_{m,l}=\\\notag
R_{i_{l}}\leftexp{(i_{1}...i_{l-1})}{d}_{m,l-1}+\leftexp{(R_{i_{l}})}{\slashed{\pi}}\cdot
\slashed{D}^{2}\leftexp{(i_{1}...i_{l-1})}{\mu}_{m,l-1}+\textrm{tr}\leftexp{(R_{i_{l}})}{\slashed{\pi}}_{1}\cdot
\slashed{d}\leftexp{(i_{1}...i_{l-1})}{\mu}_{m,l-1}
\end{align}
Applying again Proposition 8.2, the result follows. $\qed$

$\textbf{Proposition 9.4}$ For each pair of non-negative integers $(m,l)$ and each multi-index $(i_{1}...i_{l})$ we have:
\begin{align*}
 \hat{\slashed{D}}^{2}\leftexp{(i_{1}...i_{l})}{\mu}_{m,l}-\leftexp{(i_{1}...i_{l})}{\hat{\slashed{\mu}}}_{2,m,l}=
\leftexp{(i_{1}...i_{l})}{e}_{m,l}
\end{align*}
where:
\begin{align*}
 \leftexp{(i_{1}...i_{l})}{e}_{m,l}=\hat{\slashed{\mathcal{L}}}_{R_{i_{l}}}...\hat{\slashed{\mathcal{L}}}_{R_{i_{1}}}e_{m,0}
+\frac{1}{2}\sum_{k=0}^{l-1}\hat{\slashed{\mathcal{L}}}_{R_{i_{l}}}...\hat{\slashed{\mathcal{L}}}_{R_{i_{l-k+1}}}
(\leftexp{(R_{i_{l-k}})}{\hat{\slashed{\pi}}}\slashed{\Delta}\leftexp{(i_{1}...i_{l-k-1})}{\mu}_{m,l-k-1})\\
+\sum_{k=0}^{l-1}\hat{\slashed{\mathcal{L}}}_{R_{i_{l}}}...\hat{\slashed{\mathcal{L}}}_{R_{i_{l-k+1}}}
(\leftexp{(R_{i_{l-k}})}{\hat{\slashed{\pi}}}_{1}\cdot\slashed{d}\leftexp{(i_{1}...i_{l-k-1})}{\mu}_{m,l-k-1})
\end{align*}
and:
\begin{align*}
 e_{m,0}=\frac{1}{2}\sum_{k=0}^{m-1}\hat{\slashed{\mathcal{L}}}^{k}_{T}(\leftexp{(T)}{\hat{\slashed{\pi}}}\slashed{\Delta}
\mu_{m-k-1,0})+\sum_{k=0}^{m-1}\hat{\slashed{\mathcal{L}}}_{T}^{k}(\leftexp{(T)}{\hat{\slashed{\pi}}}_{1}\cdot\slashed{d}\mu_{m-k-1,0})
\end{align*}
$Proof$. We first consider the case $l=0$. If also $m=0$, from the definition, we have:
\begin{equation}
 e_{0,0}=0
\end{equation}
Consider then:
\begin{align*}
 e_{m-1,0}=\hat{\slashed{{D}}}^{2}\mu_{m-1,0}-\hat{\slashed{\mu}}_{2,m-1,0}
\end{align*}
Apply $\hat{\slashed{\mathcal{L}}}_{T}$ to this to obatin:
\begin{equation}
 \hat{\slashed{\mathcal{L}}}_{T}e_{m-1,0}=\hat{\slashed{\mathcal{L}}}_{T}\hat{\slashed{D}}^{2}\mu_{m-1,0}
-\hat{\slashed{\mu}}_{2,m,0}
\end{equation}
From the second part of Lemma 9.6:
\begin{equation}
 \hat{\slashed{\mathcal{L}}}_{T}\hat{\slashed{D}}^{2}\mu_{m-1,0}=\hat{\slashed{D}}^{2}\mu_{m,0}-\frac{1}{2}\leftexp{(T)}{\hat{\slashed{\pi}}}
\slashed{\Delta}\mu_{m-1,0}-\leftexp{(T)}{\hat{\slashed{\pi}}}_{1}\cdot\slashed{d}\mu_{m-1,0}
\end{equation}
Substituting (9.168) in (9.167) and noting that:
\begin{align*}
 \hat{\slashed{D}}^{2}\mu_{m,0}-\hat{\slashed{\mu}}_{2,m,0}=e_{m,0}
\end{align*}
we obtain:
\begin{equation}
 e_{m,0}=\hat{\slashed{\mathcal{L}}}_{T}e_{m-1,0}+\frac{1}{2}\leftexp{(T)}{\hat{\slashed{\pi}}}\slashed{\Delta}\mu_{m-1,0}
+\leftexp{(T)}{\hat{\slashed{\pi}}}_{1}\cdot\slashed{d}\mu_{m-1,0}
\end{equation}
Applying Proposition 8.2, the expression for $e_{m,0}$ follows. 

   Next, we consider the case $l\geq1$. Consider:
\begin{align*}
 \leftexp{(i_{1}...i_{l-1})}{e}_{m,l-1}=\hat{\slashed{D}}^{2}\leftexp{(i_{1}...i_{l-1})}{\mu}_{m,l-1}
-\leftexp{(i_{1}...i_{l-1})}{\hat{\slashed{\mu}}}_{2,m,l-1}
\end{align*}
Applying $\hat{\slashed{\mathcal{L}}}_{R_{i_{l}}}$ to this we obtain:
\begin{equation}
 \hat{\slashed{\mathcal{L}}}_{R_{i_{l}}}\leftexp{(i_{1}...i_{l-1})}{e}_{m,l-1}=
\hat{\slashed{\mathcal{L}}}_{R_{i_{l}}}\hat{\slashed{D}}^{2}\leftexp{(i_{1}...i_{l-1})}{\mu}_{m,l-1}
-\leftexp{(i_{1}...i_{l})}{\hat{\slashed{\mu}}}_{2,m,l}
\end{equation}
From the second part of Lemma 9.5:
\begin{align}
 \hat{\slashed{\mathcal{L}}}_{R_{i_{l}}}\hat{\slashed{D}}^{2}\leftexp{(i_{1}...i_{l-1})}{\mu}_{m,l-1}=
\hat{\slashed{D}}^{2}\leftexp{(i_{1}...i_{l})}{\mu}_{m,l}\\\notag
-\frac{1}{2}\leftexp{(R_{i_{l}})}{\hat{\slashed{\pi}}}\slashed{\Delta}\leftexp{(i_{1}...i_{l-1})}{\mu}_{m,l-1}
-\leftexp{(R_{i_{l}})}{\hat{\slashed{\pi}}}_{1}\cdot\slashed{d}\leftexp{(i_{1}...i_{l-1})}{\mu}_{m,l-1}
\end{align}
Substituting (9.171) in (9.170) and noting that:
\begin{align*}
 \hat{\slashed{D}}^{2}\leftexp{(i_{1}...i_{l})}{\mu}_{m,l}-\leftexp{(i_{1}...i_{l})}{\hat{\slashed{\mu}}}_{2,m,l}
=\leftexp{(i_{1}...i_{l})}{e}_{m,l}
\end{align*}
we obtain:
\begin{align}
 \leftexp{(i_{1}...i_{l})}{e}_{m,l}=\hat{\slashed{\mathcal{L}}}_{R_{i_{l}}}\leftexp{(i_{1}...i_{l-1})}{e}_{m,l-1}\\\notag
+\frac{1}{2}\leftexp{(R_{i_{l}})}{\hat{\slashed{\pi}}}\slashed{\Delta}\leftexp{(i_{1}...i_{l-1})}{\mu}_{m,l-1}
+\leftexp{(R_{i_{l}})}{\hat{\slashed{\pi}}}_{1}\cdot\slashed{d}\leftexp{(i_{1}...i_{l-1})}{\mu}_{m,l-1}
\end{align}
Thus, applying again Proposition 8.2, the result follows. $\qed$

   Let us investigate the order of $\leftexp{(i_{1}...i_{l})}{d}_{m,l}$. First, we should investigate $d_{m,0}$. Here, the leading terms are of order $m+1$ and are contributed
by $\slashed{\mathcal{L}}^{k}_{T}(\slashed{D}^{2}\mu_{m-k-1})$ in the first sum, a 2nd angular derivative of $T^{m-1}\mu$ to principal terms, 
and by $\slashed{\mathcal{L}}^{m-1}_{T}\textrm{tr}\leftexp{(T)}{\slashed{\pi}}_{1}$ in the second sum, a 1st angular derivative of 
$\slashed{\mathcal{L}}_{T}^{m-1}\theta$ to principal terms.The principal acoustical part of the latter is a first angular derivative of $\chi$
if $m=1$, a 3rd angular derivative of $T^{m-2}\mu$ if $m\geq2$ (by (3.125)-(3.126)). It follows that the leading terms
in the first term of $d_{m,l}$:
\begin{align*}
 R_{i_{l}}...R_{i_{1}}d_{m,0}
\end{align*}
are of order $m+l+1$, the principal acoustical terms being $l+1$st angular derivatives of $\chi$ and $l+2$nd angular derivatives of $\mu$ if $m=1$,
$l+3$rd angular derivatives of $T^{m-2}\mu$ and $l+2$nd angular derivatives of $T^{m-1}\mu$ if $m\geq2$. 

The leading terms in the second term of $\leftexp{(i_{1}...i_{l})}{d}_{m,l}$:
\begin{align*}
 \sum_{k=0}^{l-1}R_{i_{l}}...R_{i_{l-k+1}}(\leftexp{(R_{i_{l-k}})}{\slashed{\pi}}\cdot\slashed{D}^{2}\mu_{m,l-k-1})
\end{align*}
are also of order $m+l+1$, being $l+1$st angular derivatives of $T^{m}\mu, l\geq1$. Finally, the leading term in the third term 
of $\leftexp{(i_{1}...i_{l})}{d}_{m,l}$:
\begin{align*}
 \sum_{k=0}^{l-1}R_{i_{l}}...R_{i_{l-k+1}}(\textrm{tr}\leftexp{(R_{i_{l-k}})}{\slashed{\pi}}_{1}\cdot
\slashed{d}\leftexp{(i_{1}...i_{l-k-1})}{\mu}_{m,l-k-1})
\end{align*}
is:
\begin{align*}
 (\slashed{\mathcal{L}}_{R_{i_{l}}}...\slashed{\mathcal{L}}_{R_{i_{2}}}\textrm{tr}\leftexp{(R_{i_{1}})}{\slashed{\pi}}_{1})
\cdot\slashed{d}\mu_{m,0}
\end{align*}
From the expression for $\leftexp{(R_{i})}{\pi}_{AB}$, this term involves the $l$th angular derivatives of $\chi$ and is thus of 
order $l+1$, which coincides with the order of the leading terms in the first two terms of the expression for $\leftexp{(i_{1}...i_{l})}{d}_{m,l}$
if $m=0$.
We conclude that $\leftexp{(i_{1}...i_{l})}{d}_{m,l}$ is of order $m+l+1$-one less than principal-and its leading acoustical terms are $l$th angular 
derivatives of $\chi$ and (for $l\geq1$) $l+1$st angular derivatives of $\mu$ if $m=0$, $l+1$st angular derivatives of $\chi$ and $l+2$nd angular 
derivatives of $\mu$ and (for $l\geq1$) $l+1$st angular derivatives of $T\mu$ if $m=1$, $l+3$rd angular derivatives of $T^{m-2}\mu$ and $l+2$nd
angular derivatives of $T^{m-1}\mu$ and (for $l\geq1$) $l+1$st angular derivatives of $T^{m}\mu$ if $m\geq2$. 

   Let us also investigate the order of $\leftexp{(i_{1}...i_{l})}{e}_{m,l}$. First, we should investigate $e_{m,0}$. Here, the leading terms are 
of order $m+1$ and are contributed by $T^{k}(\slashed{\Delta}\mu_{m-k-1,0})$ in the first sum, a 2nd angular derivative of $T^{m-1}\mu$ to principal terms
and by $\hat{\slashed{\mathcal{L}}}^{m-1}_{T}\leftexp{(T)}{\hat{\slashed{\pi}}}_{1}$ in the second sum, a 1st angular derivative of 
$\slashed{\mathcal{L}}^{m-1}_{T}\theta$ to principal terms. The principal acoustical part of the latter is a 1st angular derivative of $\chi$ if $m=1$,
a 3rd angular derivative of $T^{m-2}\mu$ if $m\geq2$ (by (3.125)-(3.126)). It follows that the leading terms in the first term of 
$\leftexp{(i_{1}...i_{l})}{e}_{m,l}$: 
\begin{align*}
 \hat{\slashed{\mathcal{L}}}_{R_{i_{l}}}...\hat{\slashed{\mathcal{L}}}_{R_{i_{1}}}e_{m,0}
\end{align*}
are of order $m+l+1$, the leading acoustical terms being $l+1$st angular derivatives of $\chi$ and $l+2$nd
angular derivatives of $\mu$ if $m=1$, $l+3$rd angular derivatives of $T^{m-2}\mu$ and $l+2$nd angular derivatives of $T^{m-1}\mu$ if $m\geq2$.
The leading terms in the second term of $\leftexp{(i_{1}...i_{l})}{e}_{m,l}$:
\begin{align*}
 \frac{1}{2}\sum_{k=0}^{l-1}\hat{\slashed{\mathcal{L}}}_{R_{i_{l}}}...\hat{\slashed{\mathcal{L}}}_{R_{i_{l-k+1}}}
(\leftexp{(R_{i_{l-k}})}{\hat{\slashed{\pi}}}\slashed{\Delta}\leftexp{(i_{1}...i_{l-k-1})}{\mu}_{m,l-k-1})
\end{align*}
are also of order $m+l+1$, being $l+1$st angular derivatives of $T^{m}\mu, l\geq1$. Finally, the leading term in the third term of 
$\leftexp{(i_{1}...i_{l})}{e}_{m,l}$:
\begin{align*}
 \sum_{k=0}^{l-1}\hat{\slashed{\mathcal{L}}}_{R_{i_{l}}}...\hat{\slashed{\mathcal{L}}}_{R_{i_{l-k+1}}}
(\leftexp{(R_{i_{l-k}})}{\hat{\slashed{\pi}}}_{1}\cdot\slashed{d}\leftexp{(i_{1}...i_{l-k-1})}{\mu}_{m,l-k-1})
\end{align*}
is 
\begin{align*}
 (\hat{\slashed{\mathcal{L}}}_{R_{i_{l}}}...\hat{\slashed{\mathcal{L}}}_{R_{i_{2}}}
\leftexp{(R_{i_{1}})}{\hat{\slashed{\pi}_{1}}})\cdot\slashed{d}\mu_{m,0}
\end{align*}
This term involves the $l$th order angular derivatives of $\chi$ and is thus of order $l+1$, which coincides with the order of the leading terms in 
the first two terms of $\leftexp{(i_{1}...i_{l})}{e}_{m,l}$ if $m=0$. We conclude that $\leftexp{(i_{1}...i_{l})}{e}_{m,l}$ is of order $m+l+1$-one
less than principal-and its leading acoustical terms are $l$th angular derivatives of $\chi$ and (for $l\geq1$) $l+1$st angular derivatives of $\mu$
if $m=0$, $l+1$st angular derivatives of $\chi$ and $l+2$nd angular derivative of $\mu$ and (for $l\geq1$) $l+1$st angular derivatives of $T\mu$ if 
$m=1$, $l+3$rd angular derivatives of $T^{m-2}\mu$ and $l+2$nd angular derivatives of $T^{m-1}\mu$ and (for $l\geq1$) $l+1$st angular derivatives
of $T^{m}\mu$ if $m\geq2$. 

    By equation (9.134) and Proposition 9.3, $\leftexp{(i_{1}...i_{l})}{\mu}_{m,l}$ satisfies on each $S_{t,u}$ the elliptic
equation:
\begin{equation}
 \slashed{\Delta}\leftexp{(i_{1}...i_{l})}{\mu}_{m,l}=\mu^{-1}
(\leftexp{(i_{1}...i_{l})}{x'}_{m,l}+\leftexp{(i_{1}...i_{l})}{\check{f}'}_{m,l})+\leftexp{(i_{1}...i_{l})}{d}_{m,l}
\end{equation}
We shall need the following $\mu$-weighted $L^{2}$ estimate:

$\textbf{Lemma 9.7}$ Let $(M,g)$ be a compact 2-dimensional Riemannian manifold, and let $\phi$ be a function on $(M,g)$ satisfying 
the equation:
\begin{align*}
 \Delta_{g}\phi=\rho
\end{align*}
for some function $\rho$ on $M$. Let also $\mu$ be an arbitrary non-negative function on $M$. Then the following estimate holds on $M$:
\begin{align*}
 \int_{M}\mu^{2}\{\frac{1}{2}|\nabla^{2}\phi|^{2}+K|d\phi|^{2}\}d\mu_{g}\leq2\int_{M}\mu^{2}\rho^{2}d\mu_{g}+3\int_{M}|d\mu|^{2}|d\phi|^{2}d\mu_{g}
\end{align*}
where $K$ is the Gauss curvature of $(M,g)$.

$Proof$. Consider the 1-form:
\begin{equation}
 \psi=d\phi
\end{equation}
We have:
\begin{equation}
 \textrm{div}_{g}\psi=\Delta_{g}\phi=\rho
\end{equation}
So we get:
\begin{equation}
 \mu d(\textrm{div}_{g}\psi)=d(\mu\rho)-\rho d\mu
\end{equation}
Now, in arbitrary local coordinate $d(\textrm{div}_{g}\psi)$ is:
\begin{align*}
\nabla_{a}(\nabla^{b}\psi_{b})
\end{align*}
and we have:
\begin{equation}
 \nabla_{a}(\nabla^{b}\psi_{b})=\nabla^{b}(\nabla_{a}\psi_{b})-S_{a}{}^{b}\psi_{b}
\end{equation}
where $S_{a}{}^{b}=S_{ac}(g^{-1})^{bc}$ and $S_{ac}$ are the components of the Ricci curvature of $(M,g)$.
Since $\textrm{dim}M=2$, we have:
\begin{align*}
 S_{a}{}^{b}=K\delta^{b}_{a}
\end{align*}
Noting also that:
\begin{align*}
 \nabla_{a}\psi_{b}=\nabla_{b}\psi_{a}
\end{align*}
(9.177) reduces to
\begin{equation}
 \nabla_{a}(\nabla^{b}\psi_{b})=\nabla^{b}(\nabla_{b}\psi_{a})-K\psi_{a}
\end{equation}
Substituting in (9.176) we obtain:
\begin{equation}
 \mu\nabla^{b}(\nabla_{b}\psi_{a})-\mu K\psi_{a}=d_{a}(\mu\rho)-\rho d_{a}\mu
\end{equation}
We multiply this equation by $-\mu\psi^{a}$ and integrate over $M$. We have, integrating by parts,
\begin{equation}
 -\int_{M}\mu^{2}\psi^{a}\nabla^{b}(\nabla_{b}\psi_{a})d\mu_{g}=\int_{M}\mu^{2}|\nabla\psi|^{2}d\mu_{g}
+\int_{M}2\mu I^{a}d_{a}\mu d\mu_{g}
\end{equation}
where
\begin{align}
 I^{a}=\psi^{b}\nabla_{b}\psi^{a}
\end{align}
Moreover, we have, also integrating by parts,
\begin{equation}
 -\int_{M}\mu\psi^{a}d_{a}(\mu\rho)d\mu_{g}=\int_{M}\mu\rho \textrm{div}_{g}(\mu\psi)d\mu_{g}=
\int_{M}\{\mu^{2}\rho^{2}+\mu\rho\psi^{a}d_{a}\mu\}d\mu_{g}
\end{equation}
In view of (9.180) and (9.182), we obtain:
\begin{equation}
 \int_{M}\mu^{2}\{|\nabla\psi|^{2}+K|\psi|^{2}\}d\mu_{g}+\int_{M}2\mu I^{a}d_{a}\mu d\mu_{g}
=\int_{M}\{\mu^{2}\rho^{2}+2\mu\rho\psi^{a}d_{a}\mu\}d\mu_{g}
\end{equation}
Now, we have:
\begin{align*}
 2\mu|I^{a}d_{a}\mu|\leq2\mu|I||d\mu|
\end{align*}
and
\begin{align*}
 |I|\leq |\nabla\psi||\psi|
\end{align*}
hence we can estimate:
\begin{align}
 2\mu|I^{a}d_{a}\mu|\leq2\mu|d\mu||\psi||\nabla\psi|\leq
\frac{1}{2}\mu^{2}|\nabla\psi|^{2}+2|d\mu|^{2}|\psi|^{2}
\end{align}
Also we can estimate:
\begin{equation}
 2\mu|\rho\psi^{a}d_{a}\mu|\leq\mu^{2}\rho^{2}+|d\mu|^{2}|\psi|^{2}
\end{equation}
In view of (9.184)-(9.185), the identity (9.183) implies:
\begin{equation}
 \int_{M}\mu^{2}\{\frac{1}{2}|\nabla\psi|^{2}+K|\psi|^{2}\}d\mu_{g}\leq2\int_{M}\mu^{2}\rho^{2}d\mu_{g}
+3\int_{M}|d\mu|^{2}|\psi|^{2}d\mu_{g}
\end{equation}
The lemma thus follows. $\qed$

    We now apply Lemma 9.7 to (9.173), taking $(M,g)$ to be $(S_{t,u}, \slashed{g})$, $\phi$ to be $\leftexp{(i_{1}...i_{l})}{\mu}_{m,l}$,
and the function $\rho$ to be:
\begin{align*}
 \mu^{-1}(\leftexp{(i_{1}...i_{l})}{x}^{\prime}_{m,l}+\leftexp{(i_{1}...i_{l})}{\check{f}'}_{m,l})+\leftexp{(i_{1}...i_{l})}{d}_{m,l}
\end{align*}
In view of Lemma 8.9, taking $\delta_{0}$ suitably small, we have $K\geq0$, thus we obtain:
\begin{align}
 \|\mu\slashed{D}^{2}\leftexp{(i_{1}...i_{l})}{\mu}_{m,l}\|_{L^{2}(S_{t,u})}\\\notag
\leq C\|\leftexp{(i_{1}...i_{l})}{x}^{\prime}_{m,l}\|_{L^{2}(S_{t,u})}+C\|\leftexp{(i_{1}...i_{l})}{\check{f}'}_{m,l}\|_{L^{2}(S_{t,u})}\\\notag
+C\|\mu\leftexp{(i_{1}...i_{l})}{d}_{m,l}\|_{L^{2}(S_{t,u})}+C\|\slashed{d}\mu\|_{L^{\infty}(S_{t,u})}
\|\slashed{d}\leftexp{(i_{1}...i_{l})}{\mu}_{m,l}\|_{L^{2}(S_{t,u})}
\end{align}
By $\textbf{F1}$,
\begin{equation}
 |\slashed{d}\mu|\leq C\delta_{0}(1+t)^{-1}[1+\log(1+t)]
\end{equation}
Substituting (9.188) in (9.187) yields:
\begin{align}
  \|\mu\slashed{D}^{2}\leftexp{(i_{1}...i_{l})}{\mu}_{m,l}\|_{L^{2}(S_{t,u})}\\\notag
\leq C\|\leftexp{(i_{1}...i_{l})}{x}^{\prime}_{m,l}\|_{L^{2}(S_{t,u})}+C\|\leftexp{(i_{1}...i_{l})}{\check{f}'}_{m,l}\|_{L^{2}(S_{t,u})}
+C\|\mu\leftexp{(i_{1}...i_{l})}{d}_{m,l}\|_{L^{2}(S_{t,u})}\\\notag
+C\delta_{0}(1+t)^{-1}[1+\log(1+t)]\|\slashed{d}\leftexp{(i_{1}...i_{l})}{\mu}_{m,l}\|_{L^{2}(S_{t,u})}
\end{align}
In view of Proposition 9.4, we then obtain:
\begin{align}
 \|\mu\leftexp{(i_{1}...i_{l})}{\hat{\slashed{\mu}}}_{2,m,l}\|_{L^{2}(S_{t,u})}\leq
C\|\leftexp{(i_{1}...i_{l})}{x}^{\prime}_{m,l}\|_{L^{2}(S_{t,u})}+C\|\leftexp{(i_{1}...i_{l})}{\check{f}'}_{m,l}\|_{L^{2}(S_{t,u})}\\\notag
+C\|\mu\leftexp{(i_{1}...i_{l})}{d}_{m,l}\|_{L^{2}(S_{t,u})}+C\|\mu\leftexp{(i_{1}...i_{l})}{e}_{m,l}\|_{L^{2}(S_{t,u})}\\\notag
+C\delta_{0}(1+t)^{-1}[1+\log(1+t)]\|\slashed{d}\leftexp{(i_{1}...i_{l})}{\mu}_{m,l}\|_{L^{2}(S_{t,u})}
\end{align}

\section{The Estimates for the Solutions of the Propagation Equations}
    We now return to the propagation equation of Proposition 9.2. Setting:
\begin{equation}
 \leftexp{(i_{1}...i_{l})}{\tilde{g}'}_{m,l}=-2\mu\hat{\chi}\cdot\leftexp{(i_{1}...i_{l})}{\hat{\slashed{\mu}}}_{2,m,l}
+\leftexp{(i_{1}...i_{l})}{g}^{\prime}_{m,l}
\end{equation}
the propagation equation takes the form:
\begin{align}
 L\leftexp{(i_{1}...i_{l})}{x}^{\prime}_{m,l}+(\textrm{tr}\chi-2\mu^{-1}(L\mu))\leftexp{(i_{1}...i_{l})}{x}^{\prime}_{m,l}=
-(\frac{1}{2}\textrm{tr}\chi-2\mu^{-1}(L\mu))\leftexp{(i_{1}...i_{l})}{\check{f}'}_{m,l}+\leftexp{(i_{1}...i_{l})}{\tilde{g}'}_{m,l}
\end{align}
Defining as in Chapter 8 the diffeomorphisms $\Phi_{t,u}$ of $S^{2}$ onto $S_{t,u}$, if $\omega$ is any $S_{t,u}$ $r$-covariant tensorfield
defined in $W^{*}_{\epsilon_{0}}$, we consider the pullback $\omega(t,u)=\Phi^{*}_{t,u}\omega$, a $r$-covariant tensorfield on $S^{2}$ 
depending on the parameters $t$ and $u$. If in place of $\omega$ we have $\slashed{\mathcal{L}}_{L}\omega$, the corresponding tensorfield on $S^{2}$
is $\frac{\partial\omega(t,u)}{\partial t}$. The propagation equation (9.192) can then be viewed as an equation for the function 
$\leftexp{(i_{1}...i_{l})}{x'}(t,u)$ on $S^{2}$ depending on $t$ and $u$:
\begin{align}
 \frac{\partial}{\partial t}\leftexp{(i_{1}...i_{l})}{x}^{\prime}_{m,l}+(\textrm{tr}\chi-2\mu^{-1}\frac{\partial\mu}{\partial t})
\leftexp{(i_{1}...i_{l})}{x}^{\prime}_{m,l}=-(\frac{1}{2}\textrm{tr}\chi-2\mu^{-1}\frac{\partial\mu}{\partial t})\leftexp{(i_{1}...i_{l})}{\check{f}'}_{m,l}
+\leftexp{(i_{1}...i_{l})}{\tilde{g}'}_{m,l}
\end{align}

    Consider a function $\phi(t,u)$ on $S^{2}$ depending on the parameters $t$ and $u$. We have:
\begin{equation}
 |\phi|\frac{\partial}{\partial t}|\phi|=\frac{1}{2}\frac{\partial}{\partial t}\phi^{2}=\phi\frac{\partial\phi}{\partial t}
\end{equation}
Then from the propagation equation (9.193) we have:
\begin{align}
 \leftexp{(i_{1}...i_{l})}{x}^{\prime}_{m,l}\frac{\partial}{\partial t}\leftexp{(i_{1}...i_{l})}{x}^{\prime}_{m,l}
+(\textrm{tr}\chi-2\mu^{-1}\frac{\partial\mu}{\partial t})(\leftexp{(i_{1}...i_{l})}{x}^{\prime}_{m,l})^{2}\\\notag
=-(\frac{1}{2}\textrm{tr}\chi-2\mu^{-1}\frac{\partial\mu}{\partial t})\leftexp{(i_{1}...i_{l})}{x}^{\prime}_{m,l}\leftexp{(i_{1}...i_{l})}{\check{f}'}_{m,l}
+\leftexp{(i_{1}...i_{l})}{x'}_{m,l}\leftexp{(i_{1}...i_{l})}{\tilde{g}'}_{m,l}
\end{align}
Since
\begin{align*}
 |\leftexp{(i_{1}...i_{l})}{x}^{\prime}_{m,l}\leftexp{(i_{1}...i_{l})}{\check{f}'}_{m,l}|
\leq |\leftexp{(i_{1}...i_{l})}{x}^{\prime}_{m,l}||\leftexp{(i_{1}...i_{l})}{\check{f}'}_{m,l}|
\end{align*}
and
\begin{align*}
  |\leftexp{(i_{1}...i_{l})}{x}^{\prime}_{m,l}\leftexp{(i_{1}...i_{l})}{\tilde{g}'}_{m,l}|
\leq |\leftexp{(i_{1}...i_{l})}{x}^{\prime}_{m,l}||\leftexp{(i_{1}...i_{l})}{\tilde{g}'}_{m,l}|
\end{align*}
the assumption $\textbf{AS}$ implies through (9.195)
\begin{align}
 \leftexp{(i_{1}...i_{l})}{x}^{\prime}_{m,l}\frac{\partial}{\partial t}\leftexp{(i_{1}...i_{l})}{x}^{\prime}_{m,l}
+(\textrm{tr}\chi-2\mu^{-1}\frac{\partial\mu}{\partial t})|\leftexp{(i_{1}...i_{l})}{x}^{\prime}_{m,l}|^{2}\\\notag
\leq(\frac{1}{2}\textrm{tr}\chi-2\mu^{-1}\frac{\partial\mu}{\partial t})|\leftexp{(i_{1}...i_{l})}{x}^{\prime}_{m,l}|
|\leftexp{(i_{1}...i_{l})}{\check{f}'}_{m,l}|+|\leftexp{(i_{1}...i_{l})}{x'}_{m,l}||\leftexp{(i_{1}...i_{l})}{\tilde{g}'}_{m,l}|
\end{align}
Applying (9.194) taking $\phi=\leftexp{(i_{1}...i_{l})}{x}^{\prime}_{m,l}$, we get:
\begin{align}
 \frac{\partial}{\partial t}|\leftexp{(i_{1}...i_{l})}{x}^{\prime}_{m,l}|+(\textrm{tr}\chi-2\mu^{-1}\frac{\partial\mu}{\partial t})
|\leftexp{(i_{1}...i_{l})}{x}^{\prime}_{m,l}|\\\notag
\leq(\frac{1}{2}\textrm{tr}\chi-2\mu^{-1}\frac{\partial\mu}{\partial t})|\leftexp{(i_{1}...i_{l})}{\check{f}'}_{m,l}|
+|\leftexp{(i_{1}...i_{l})}{\tilde{g}'}_{m,l}|
\end{align}
The integrating factor here is 
\begin{equation}
 \exp\{\int_{0}^{t}(\textrm{tr}\chi-2\mu^{-1}\frac{\partial\mu}{\partial t})(t',u)dt'\}=
(\frac{\mu(t,u)}{\mu(0,u)})^{-2}A(t,u)
\end{equation}
where $A(t,u)$ is defined by (8.167) and is given by (8.169). From (9.197) we deduce:
\begin{equation}
 |\leftexp{(i_{1}...i_{l})}{x}^{\prime}_{m,l}|\leq\leftexp{(i_{1}...i_{l})}{F}^{\prime}_{m,l}(t,u)+\leftexp{(i_{1}...i_{l})}{G}^{\prime}_{m,l}(t,u)
\end{equation}
where:
\begin{align}
 \leftexp{(i_{1}...i_{l})}{F}^{\prime}_{m,l}(t,u)=(A(t,u))^{-1}(\mu(t,u))^{2}\{(\mu(0,u))^{-2}|\leftexp{(i_{1}...i_{l})}{x}^{\prime}_{m,l}(0,u)|\\\notag
+\int_{0}^{t}(\mu(t',u))^{-2}A(t',u)(\frac{1}{2}\textrm{tr}\chi-2\mu^{-1}\frac{\partial\mu}{\partial t})(t',u)|\leftexp{(i_{1}...i_{l})}
{\check{f}'}_{m,l}(t',u)|dt'\}
\end{align}
and:
\begin{align}
 \leftexp{(i_{1}...i_{l})}{G}^{\prime}_{m,l}(t,u)=(A(t,u))^{-1}(\mu(t,u))^{2}\cdot\int_{0}^{t}
(\mu(t',u))^{-2}A(t',u)|\leftexp{(i_{1}...i_{l})}{\tilde{g}'}_{m,l}(t',u)|dt'
\end{align}
Using the bound (8.173), we obtain from (9.200) and (9.201):
\begin{align}
 \leftexp{(i_{1}...i_{l})}{F}^{\prime}_{m,l}(t,u)\leq e^{C\delta_{0}}(1-u+t)^{-2}\{\leftexp{(i_{1}...i_{l})}{M}^{\prime0}_{m,l}(t,u)\\\notag
+\leftexp{(i_{1}...i_{l})}{M}^{\prime1}_{m,l}(t,u)+\leftexp{(i_{1}...i_{l})}{M}^{\prime2}_{m,l}(t,u)\}
\end{align}
where:
\begin{equation}
 \leftexp{(i_{1}...i_{l})}{M}^{\prime0}_{m,l}(t,u)=(\frac{\mu(t,u)}{\mu(0,u)})^{2}(1-u)^{2}
|\leftexp{(i_{1}...i_{l})}{x}^{\prime}_{m,l}(0,u)|
\end{equation}
\begin{align}
 \leftexp{(i_{1}...i_{l})}{M}^{\prime1}_{m,l}(t,u)=\int_{0}^{t}(\frac{\mu(t,u)}{\mu(t',u)})^{2}(1-u+t')^{2}
\times[-2\mu^{-1}(\frac{\partial\mu}{\partial t})_{-}(t',u)]\cdot|\leftexp{(i_{1}...i_{l})}{\check{f}'}_{m,l}(t',u)|dt'
\end{align}
\begin{align}
 \leftexp{(i_{1}...i_{l})}{M}^{\prime2}_{m,l}(t,u)=\frac{1}{2}\int_{0}^{t}(\frac{\mu(t,u)}{\mu(t',u)})^{2}
(1-u+t')^{2}\textrm{tr}\chi(t',u)|\leftexp{(i_{1}...i_{l})}{\check{f}'}_{m,l}(t',u)|dt'
\end{align}
Also:
\begin{align}
 \leftexp{(i_{1}...i_{l})}{G}^{\prime}_{m,l}(t,u)\leq e^{C\delta_{0}}(1-u+t)^{-2}
\cdot\int_{0}^{t}(\frac{\mu(t,u)}{\mu(t',u)})^{2}(1-u+t')^{2}|\leftexp{(i_{1}...i_{l})}{\tilde{g}}^{\prime}_{m,l}(t',u)|dt'
\end{align}
We shall estimate the $L^{2}$ norm of $\leftexp{(i_{1}...i_{l})}{F}^{\prime}_{m,l}(t)$ on $[0,\epsilon_{0}]\times S^{2}$.

First, by $\textbf{A3}$ and (8.333) at $t=0$,
\begin{equation}
 \|\leftexp{(i_{1}...i_{l})}{M}^{\prime0}_{m,l}(t)\|_{L^{2}([0,\epsilon_{0}]\times S^{2})}
\leq C[1+\log(1+t)]^{2}\|\leftexp{(i_{1}...i_{l})}{x}^{\prime}_{m,l}(0)\|_{L^{2}(\Sigma_{0}^{\epsilon_{0}})}
\end{equation}

    We turn to $\leftexp{(i_{1}...i_{l})}{M}^{\prime1}_{m,l}$. Partitioning $[0,\epsilon_{0}]\times S^{2}$ into the 
 sets $\mathcal{V}_{s-}$, $\mathcal{V}_{s+}$ defined by (8.337) and (8.338), respectively, we have:
\begin{equation}
 \|\leftexp{(i_{1}...i_{l})}{M}^{\prime1}_{m,l}(t)\|^{2}_{L^{2}([0,\epsilon_{0}]\times S^{2})}=
\|\leftexp{(i_{1}...i_{l})}{M}^{\prime1}_{m,l}(t)\|^{2}_{L^{2}(\mathcal{V}_{s-})}+
\|\leftexp{(i_{1}...i_{l})}{M}^{\prime1}_{m,l}(t)\|^{2}_{L^{2}(\mathcal{V}_{s+})}
\end{equation}
By Minkowski's inequality:
\begin{equation}
 \|\leftexp{(i_{1}...i_{l})}{M}^{\prime1}_{m,l}(t)\|_{L^{2}(\mathcal{V}_{s-})}
\leq\int_{0}^{t}\|\leftexp{(i_{1}...i_{l})}{N}^{\prime1}_{m,l}(t,t')\|_{L^{2}(\mathcal{V}_{s-})}dt'
\end{equation}
where:
\begin{align}
 \leftexp{(i_{1}...i_{l})}{N}^{\prime1}_{m,l}(t,t',u)=(\frac{\mu(t,u)}{\mu(t',u)})^{2}(1-u+t')^{2}[-2\mu^{-1}
(\frac{\partial\mu}{\partial t})_{-}(t',u)]|\leftexp{(i_{1}...i_{l})}{\check{f}}^{\prime}_{m,l}(t',u)|
\end{align}
We have:
\begin{align}
 \|\leftexp{(i_{1}...i_{l})}{N}^{\prime1}_{m,l}(t,t')\|_{L^{2}(\mathcal{V}_{s-})}
\leq(1+t')^{2}[\max_{\mathcal{V}_{s-}}(\frac{\mu(t)}{\mu(t')})]^{2}\max_{\mathcal{V}_{s-}}
[-2\mu^{-1}(\frac{\partial\mu}{\partial t})_{-}(t')]\|\leftexp{(i_{1}...i_{l})}{\check{f}'}_{m,l}(t')\|_{L^{2}(\mathcal{V}_{s-})}
\end{align}
In view of (8.343) and (8.344) (see definition (8.249)):
\begin{align}
 \|\leftexp{(i_{1}...i_{l})}{N}^{\prime1}_{m,l}(t,t')\|_{L^{2}(\mathcal{V}_{s-})}
\leq C(1+t')^{2}M(t')\|\leftexp{(i_{1}...i_{l})}{\check{f}'}_{m,l}(t')\|_{L^{2}([0,\epsilon_{0}]\times S^{2})}
\end{align}
Let us define:
\begin{equation}
 \leftexp{(i_{1}...i_{l})}{P}^{\prime}_{m,l}(t)=(1+t)\|\leftexp{(i_{1}...i_{l})}{\check{f}'}_{m,l}(t)\|_{L^{2}([0,\epsilon_{0}]\times S^{2})}
\end{equation}
Suppose that, for non-negative quantities $\leftexp{(i_{1}...i_{l})}{P}^{\prime(0)}_{m,l}$, 
$\leftexp{(i_{1}...i_{l})}{P}^{\prime(1)}_{m,l}$, we have:
\begin{equation}
 \leftexp{(i_{1}...i_{l})}{P}^{\prime}_{m,l}(t)\leq 
\leftexp{(i_{1}...i_{l})}{P}^{\prime(0)}_{m,l}(t)+\leftexp{(i_{1}...i_{l})}{P}^{\prime(1)}_{m,l}(t)
\end{equation}
We define the non-decreasing non-negative quantities $\leftexp{(i_{1}...i_{l})}{\bar{P}}^{\prime(0)}_{m,l,a}$,
$\leftexp{(i_{1}...i_{l})}{\bar{P}}^{\prime(1)}_{m,l,a}$ by:
\begin{equation}
 \leftexp{(i_{1}...i_{l})}{\bar{P}}^{\prime(0)}_{m,l,a}(t)=
\sup_{t'\in[0,t]}\{\bar{\mu}^{a}_{m}(t')\leftexp{(i_{1}...i_{l})}{P}^{\prime(0)}_{m,l}(t')\}
\end{equation}
\begin{equation}
  \leftexp{(i_{1}...i_{l})}{\bar{P}}^{\prime(1)}_{m,l,a}(t)=
\sup_{t'\in[0,t]}\{(1+t')^{1/2}\bar{\mu}^{a}_{m}(t')\leftexp{(i_{1}...i_{l})}{P}^{\prime(1)}_{m,l}(t')\}
\end{equation}
Then for $t'\in[0,t]$ we have:
\begin{equation}
 \leftexp{(i_{1}...i_{l})}{P}^{\prime}_{m,l}(t')\leq
\bar{\mu}^{-a}_{m}(t')\{\leftexp{(i_{1}...i_{l})}{\bar{P}}^{\prime(0)}_{m,l,a}(t)+
(1+t')^{-1/2}\leftexp{(i_{1}...i_{l})}{\bar{P}}^{\prime(1)}_{m,l,a}(t)\}
\end{equation}
Substituting in (9.212) and in view of (9.213) we have:
\begin{equation}
 \|\leftexp{(i_{1}...i_{l})}{N}^{\prime1}_{m,l}(t,t')\|_{L^{2}(\mathcal{V}_{s-})}
\leq C\{(1+t)\leftexp{(i_{1}...i_{l})}{\bar{P}}^{\prime(0)}_{m,l,a}(t)+(1+t)^{1/2}
\leftexp{(i_{1}...i_{l})}{\bar{P}}^{\prime(1)}_{m,l,a}(t)\}\bar{\mu}^{-a}_{m}(t')M(t')
\end{equation}
for all $t'\in[0,t]$. Substituting this in (9.209) we obtain:
\begin{equation}
 \|\leftexp{(i_{1}...i_{l})}{M}^{\prime1}_{m,l}(t)\|_{L^{2}(\mathcal{V}_{s-})}
\leq C\{(1+t)\leftexp{(i_{1}...i_{l})}{\bar{P}}^{\prime(0)}_{m,l,a}(t)+(1+t)^{1/2}
\leftexp{(i_{1}...i_{l})}{\bar{P}}^{\prime(1)}_{m,l,a}(t)\}I_{a}(t)
\end{equation}
Applying Lemma 8.11 we conclude that:
\begin{equation}
 \|\leftexp{(i_{1}...i_{l})}{M}^{\prime1}_{m,l}(t)\|_{L^{2}(\mathcal{V}_{s-})}
\leq Ca^{-1}\{(1+t)\leftexp{(i_{1}...i_{l})}{\bar{P}}^{\prime(0)}_{m,l,a}(t)+(1+t)^{1/2}
\leftexp{(i_{1}...i_{l})}{\bar{P}}^{\prime(1)}_{m,l,a}(t)\}\bar{\mu}^{-a}_{m}(t)
\end{equation}
  In analogy with (9.209) and (9.211), with $\mathcal{V}_{s+}$ in the role of $\mathcal{V}_{s-}$, we have:
\begin{equation}
 \|\leftexp{(i_{1}...i_{l})}{M}^{\prime1}_{m,l}(t)\|_{L^{2}(\mathcal{V}_{s+})}
\leq\int_{0}^{t}\|\leftexp{(i_{1}...i_{l})}{N}^{\prime1}_{m,l}(t,t')\|_{L^{2}(\mathcal{V}_{s+})}dt'
\end{equation}
and:
\begin{align}
 \|\leftexp{(i_{1}...i_{l})}{N}^{\prime1}_{m,l}(t,t')\|_{L^{2}(\mathcal{V}_{s+})}\leq
(1+t')^{2}[\max_{\mathcal{V}_{s+}}(\frac{\mu(t)}{\mu(t')})]^{2}\max_{\mathcal{V}_{s+}}
[-2\mu^{-1}(\frac{\partial\mu}{\partial t})_{-}(t')]\|\leftexp{(i_{1}...i_{l})}{\check{f}'}_{m,l}(t')\|_{L^{2}(\mathcal{V}_{s+})}
\end{align}
Substituting (8.356), (8.357) and (9.213) in (9.222), and the result in (9.221), we obtain:
\begin{equation}
 \|\leftexp{(i_{1}...i_{l})}{\overset{1}{M'}}_{m,l}(t)\|_{L^{2}(\mathcal{V}_{s+})}
\leq C\delta_{0}[1+\log(1+t)]^{2}\int_{0}^{t}\frac{[1+\log(1+t')]}{(1+t')}\leftexp{(i_{1}...i_{l})}{P}^{\prime}_{m,l}(t')dt'
\end{equation}
hence, by (9.217) we have:
\begin{align}
 \|\leftexp{(i_{1}...i_{l})}{M}^{\prime1}_{m,l}(t)\|_{L^{2}(\mathcal{V}_{s+})}
\leq C\delta_{0}[1+\log(1+t)]^{2}\{\leftexp{(i_{1}...i_{l})}{\bar{P}}^{\prime(0)}_{m,l,a}(t)
+\leftexp{(i_{1}...i_{l})}{\bar{P}}^{\prime(1)}_{m,l,a}(t)\}\\\notag\cdot
\int_{0}^{t}\frac{[1+\log(1+t')]}{(1+t')}\bar{\mu}^{-a}_{m}(t')dt'
\end{align}
In view of (8.363), we conclude:
\begin{align}
 \|\leftexp{(i_{1}...i_{l})}{M}^{\prime1}_{m,l}(t)\|_{L^{2}(\mathcal{V}_{s+})}
\leq C\delta_{0}[1+\log(1+t)]^{4}\{\leftexp{(i_{1}...i_{l})}{\bar{P}}^{\prime(0)}_{m,l,a}(t)
+\leftexp{(i_{1}...i_{l})}{\bar{P}}^{\prime(1)}_{m,l,a}(t)\}\bar{\mu}^{1-a}_{m}(t)
\end{align}
Combining finally (9.220) and (9.225) and taking into account the assumption $a\delta_{0}\leq C^{-1}$, we obtain:
\begin{align}
 \|\leftexp{(i_{1}...i_{l})}{M}^{\prime1}_{m,l}(t)\|_{L^{2}([0,\epsilon_{0}]\times S^{2})}
\leq Ca^{-1}\{(1+t)\leftexp{(i_{1}...i_{l})}{\bar{P}}^{\prime(0)}_{m,l,a}(t)+(1+t)^{1/2}
\leftexp{(i_{1}...i_{l})}{\bar{P}}^{\prime(1)}_{m,l,a}(t)\}\bar{\mu}^{-a}_{m}(t)
\end{align}
Here, we also have used the definition (8.248) of Chapter 8.

     We turn to $\leftexp{(i_{1}...i_{l})}{M}^{\prime2}_{m,l}(t,u)$. We have:
\begin{equation}
 \|\leftexp{(i_{1}...i_{l})}{M}^{\prime2}_{m,l}(t)\|^{2}_{L^{2}([0,\epsilon_{0}]\times S^{2})}
= \|\leftexp{(i_{1}...i_{l})}{M}^{\prime2}_{m,l}(t)\|^{2}_{L^{2}(\mathcal{V}_{s-})}+
\|\leftexp{(i_{1}...i_{l})}{M}^{\prime2}_{m,l}(t)\|^{2}_{L^{2}(\mathcal{V}_{s+})}
\end{equation}
and:
\begin{equation}
 \|\leftexp{(i_{1}...i_{l})}{M}^{\prime2}_{m,l}(t)\|_{L^{2}(\mathcal{V}_{s-})}
\leq\int_{0}^{t}\|\leftexp{(i_{1}...i_{l})}{N}^{\prime2}_{m,l}(t,t')\|_{L^{2}(\mathcal{V}_{s-})}dt'
\end{equation}
\begin{equation}
  \|\leftexp{(i_{1}...i_{l})}{M}^{\prime2}_{m,l}(t)\|_{L^{2}(\mathcal{V}_{s+})}
\leq\int_{0}^{t}\|\leftexp{(i_{1}...i_{l})}{N}^{\prime2}_{m,l}(t,t')\|_{L^{2}(\mathcal{V}_{s+})}dt'
\end{equation}
where
\begin{align}
 \leftexp{(i_{1}...i_{l})}{N}^{\prime2}_{m,l}(t,t',u)=\frac{1}{2}(\frac{\mu(t,u)}{\mu(t',u)})^{2}
(1-u+t')^{2}\textrm{tr}\chi(t',u)|\leftexp{(i_{1}...i_{l})}{\check{f}'}_{m,l}(t',u)|
\end{align}
We have:
\begin{align}
 \|\leftexp{(i_{1}...i_{l})}{N}^{\prime2}_{m,l}(t,t')\|_{L^{2}(\mathcal{V}_{s-})}\leq 
\frac{1}{2}(1+t')^{2}[\max_{\mathcal{V}_{s-}}(\frac{\mu(t)}{\mu(t')})]^{2}\max_{\mathcal{V}_{s-}}(\textrm{tr}\chi(t'))
\|\leftexp{(i_{1}...i_{l})}{\check{f}'}_{m,l}(t')\|_{L^{2}(\mathcal{V}_{s-})}
\end{align}
By (8.343) and (8.371):
\begin{equation}
 \|\leftexp{(i_{1}...i_{l})}{N}^{\prime2}_{m,l}(t,t')\|_{L^{2}(\mathcal{V}_{s-})}\leq 
C(1+t')\|\leftexp{(i_{1}...i_{l})}{\check{f}'}_{m,l}(t')\|_{L^{2}([0,\epsilon_{0}]\times S^{2})}
\end{equation}
Then by (9.213) and (9.217):
\begin{align}
 \|\leftexp{(i_{1}...i_{l})}{N}^{\prime2}_{m,l}(t,t')\|_{L^{2}(\mathcal{V}_{s-})}\leq 
C\leftexp{(i_{1}...i_{l})}{P}^{\prime}_{m,l}(t')\\\notag
\leq C\{\leftexp{(i_{1}...i_{l})}{\bar{P}}^{\prime(0)}_{m,l,a}(t)+(1+t')^{-1/2}
\leftexp{(i_{1}...i_{l})}{\bar{P}}^{\prime(1)}_{m,l,a}(t)\}\bar{\mu}_{m}^{-a}(t')\\\notag
\leq C'\{\leftexp{(i_{1}...i_{l})}{\bar{P}}^{\prime(0)}_{m,l,a}(t)+(1+t')^{-1/2}
\leftexp{(i_{1}...i_{l})}{\bar{P}}^{\prime(1)}_{m,l,a}(t)\}\bar{\mu}^{-a}_{m}(t)
\end{align}
where we have used Corollary 2 to Lemma 8.11. Substituting (9.233) in (9.228) we get:
\begin{equation}
 \|\leftexp{(i_{1}...i_{l})}{M}^{\prime2}_{m,l}(t)\|_{L^{2}(\mathcal{V}_{s-})}
\leq C\{(1+t)\leftexp{(i_{1}...i_{l})}{\bar{P}}^{\prime(0)}_{m,l,a}(t)+
(1+t)^{1/2}\leftexp{(i_{1}...i_{l})}{\bar{P}}^{\prime(1)}_{m,l,a}(t)\}\bar{\mu}^{-a}_{m}(t)
\end{equation}

    In analogy with (9.231) we have:
\begin{align}
 \|\leftexp{(i_{1}...i_{l})}{N}^{\prime2}_{m,l}(t,t')\|_{L^{2}(\mathcal{V}_{s+})}
\leq\frac{1}{2}(1+t')^{2}[\max_{\mathcal{V}_{s+}}(\frac{\mu(t)}{\mu(t')})]^{2}\max_{\mathcal{V}_{s+}}
(\textrm{tr}\chi(t'))\|\leftexp{(i_{1}...i_{l})}{\check{f'}}_{m,l}(t')\|_{L^{2}(\mathcal{V}_{s+})}
\end{align}
Here, we have two cases to distinguish according as to whether $t'$ is $<\sqrt{t}$ or $\geq\sqrt{t}$. In the first case,
(8.378) holds, and in the second case (8.377) holds. In view of these, recalling (8.371), we obtain:
\begin{align}
 \|\leftexp{(i_{1}...i_{l})}{N}^{\prime2}_{m,l}(t,t')\|_{L^{2}(\mathcal{V}_{s+})}
\leq C[1+\log(1+t)]^{2}\leftexp{(i_{1}...i_{l})}{P}^{\prime}_{m,l}(t'): t'\leq\sqrt{t}\\\notag
\|\leftexp{(i_{1}...i_{l})}{N}^{\prime2}_{m,l}(t,t')\|_{L^{2}(\mathcal{V}_{s+})}
\leq C\leftexp{(i_{1}...i_{l})}{P}^{\prime}_{m,l}(t'): t'\geq\sqrt{t}
\end{align}
Substituting (9.236) in (9.229), we then have:
\begin{align}
 \|\leftexp{(i_{1}...i_{l})}{M}^{\prime2}_{m,l}(t)\|_{L^{2}(\mathcal{V}_{s+})}\leq 
\int_{0}^{\sqrt{t}}\|\leftexp{(i_{1}...i_{l})}{N}^{\prime2}_{m,l}(t,t')\|_{L^{2}(\mathcal{V}_{s+})}dt'
+\int_{\sqrt{t}}^{t}\|\leftexp{(i_{1}...i_{l})}{N}^{\prime2}_{m,l}(t,t')\|_{L^{2}(\mathcal{V}_{s+})}dt'\\\notag
\leq C[1+\log(1+t)]^{2}\int_{0}^{\sqrt{t}}\leftexp{(i_{1}...i_{l})}{P}^{\prime}_{m,l}(t')dt'
+C\int_{\sqrt{t}}^{t}\leftexp{(i_{1}...i_{l})}{P}^{\prime}_{m,l}(t')dt'\\\notag
\leq C[1+\log(1+t)]^{2}\int_{0}^{\sqrt{t}}\bar{\mu}^{-a}_{m}(t')\{\leftexp{(i_{1}...i_{l})}{\bar{P}}^{\prime(0)}_{m,l,a}(t)
+(1+t')^{-1/2}\leftexp{(i_{1}...i_{l})}{\bar{P}}^{\prime(1)}_{m,l,a}(t)\}dt'\\\notag
+C\int_{\sqrt{t}}^{t}\bar{\mu}^{-a}_{m}(t')\{\leftexp{(i_{1}...i_{l})}{\bar{P}}^{\prime(0)}_{m,l,a}(t)
+(1+t')^{-1/2}\leftexp{(i_{1}...i_{l})}{\bar{P}}^{\prime(1)}_{m,l,a}(t)\}dt'\\\notag
\leq C[1+\log(1+t)]^{2}\{(1+t)^{1/2}\leftexp{(i_{1}...i_{l})}{\bar{P}}^{\prime(0)}_{m,l,a}(t)+
(1+t)^{1/4}\leftexp{(i_{1}...i_{l})}{\bar{P}}^{\prime(1)}_{m,l,a}(t)\}\bar{\mu}^{-a}_{m}(t)\\\notag
+C\{(1+t)\leftexp{(i_{1}...i_{l})}{\bar{P}}^{\prime(0)}_{m,l,a}(t)+
(1+t)^{1/2}\leftexp{(i_{1}...i_{l})}{\overset{(1)}{\bar{P'}}}_{m,l,a}(t)\}\bar{\mu}^{-a}_{m}(t)\\\notag
\leq C'\{(1+t)\leftexp{(i_{1}...i_{l})}{\bar{P}}^{\prime(0)}_{m,l,a}(t)+
(1+t)^{1/2}\leftexp{(i_{1}...i_{l})}{\bar{P}}^{\prime(1)}_{m,l,a}(t)\}\bar{\mu}^{-a}_{m}(t)\notag
\end{align}
where we have used Corollary 2 to Lemma 8.11 and (9.217).
 
   Combining (9.234) and (9.237) we obtain:
\begin{equation}
 \|\leftexp{(i_{1}...i_{l})}{M}^{\prime2}_{m,l}(t)\|_{L^{2}([0,\epsilon_{0}]\times S^{2})}\leq 
C\{(1+t)\leftexp{(i_{1}...i_{l})}{\bar{P}}^{\prime(0)}_{m,l,a}(t)+
(1+t)^{1/2}\leftexp{(i_{1}...i_{l})}{\bar{P}}^{\prime(1)}_{m,l,a}(t)\}\bar{\mu}^{-a}_{m}(t)
\end{equation}
The estimates (9.202), (9.207), (9.226) and (9.238) yield:
\begin{align}
 \|\leftexp{(i_{1}...i_{l})}{F}^{\prime}_{m,l}(t)\|_{L^{2}([0,\epsilon_{0}]\times S^{2})}
\leq C(1+t)^{-2}[1+\log(1+t)]^{2}\|\leftexp{(i_{1}...i_{l})}{x}^{\prime}_{m,l}(0)\|_{L^{2}(\Sigma^{\epsilon_{0}}_{0})}\\\notag
+C(1+t)^{-1}\{\leftexp{(i_{1}...i_{l})}{\bar{P}}^{\prime(0)}_{m,l,a}(t)+(1+t)^{-1/2}
\leftexp{(i_{1}...i_{l})}{\bar{P}}^{\prime(1)}_{m,l,a}(t)\}\bar{\mu}^{-a}_{m}(t)
\end{align}

     We now consider the estimate (9.206). From (8.383):
\begin{equation}
 \leftexp{(i_{1}...i_{l})}{G}^{\prime}_{m,l}(t,u)\leq C(1+t)^{-2}[1+\log(1+t)]^{2}\int_{0}^{t}
(1+t')^{2}|\leftexp{(i_{1}...i_{l})}{\tilde{g}}^{\prime}_{m,l}(t',u)|dt'
\end{equation}
It follows that:
\begin{align}
 \|\leftexp{(i_{1}...i_{l})}{G}^{\prime}_{m,l}(t)\|_{L^{2}([0,\epsilon_{0}]\times S^{2})}\leq 
C(1+t)^{-2}[1+\log(1+t)]^{2}\\\notag
\cdot\int_{0}^{t}(1+t')^{2}\|\leftexp{(i_{1}...i_{l})}{\tilde{g}}^{\prime}_{m,l}(t')\|_{L^{2}([0,\epsilon_{0}]\times S^{2})}dt'
\end{align}

    Now $\leftexp{(i_{1}...i_{l})}{\tilde{g}}^{\prime}_{m,l}$ is given by (9.191). According to the discussion following Lemma 9.4, the principal 
part of $\leftexp{(i_{1}...i_{l})}{g}^{\prime}_{m,l}$ consists of (9.118), (9.120) and (9.126). Define the function 
$\leftexp{(i_{1}...i_{l})}{\dot{g}'}_{m,l}$ by:
\begin{align}
 \leftexp{(i_{1}...i_{l})}{\dot{g}}^{\prime}_{0,l}=\leftexp{(i_{1}...i_{l})}{g}^{\prime}_{0,l}-\xi\cdot\leftexp{(i_{1}...i_{l})}{x}_{l}
-\sum_{k=0}^{l-1}\leftexp{(R_{i_{l-k}})}{Z}\leftexp{(i_{1}...i_{l-k-1}i_{l-k+1}...i_{l})}{x}^{\prime}_{0,l-1} : m=0
\end{align}
and:
\begin{align}
 \leftexp{(i_{1}...i_{l})}{\dot{g}}^{\prime}_{m,l}=\leftexp{(i_{1}...i_{l})}{g}^{\prime}_{m,l}-\xi\cdot\slashed{d}\leftexp{(i_{1}...i_{l})}{x}^{\prime}_{m-1,l}\\\notag
-m\Lambda\leftexp{(i_{1}...i_{l})}{x}^{\prime}_{m-1,l}-\sum_{k=0}^{l-1}\leftexp{(R_{i_{l-k}})}{Z}\leftexp{(i_{1}...i_{l-k-1}i_{l-k+1}...i_{l})}{x}^{\prime}_{m,l-1}
: m\geq 1 
\end{align}
Then $\leftexp{(i_{1}...i_{l})}{\dot{g}}^{\prime}_{m,l}$ does not contain principal acoustical terms and we have:
\begin{align}
 \leftexp{(i_{1}...i_{l})}{\tilde{g}}^{\prime}_{0,l}=-2\mu\hat{\chi}\cdot\leftexp{(i_{1}...i_{l})}{\hat{\slashed{\mu}}}_{2,0,l}
+\xi\cdot\leftexp{(i_{1}...i_{l})}{x}_{l}\\\notag
+\sum_{k=0}^{l-1}\leftexp{(R_{i_{l-k}})}{Z}\leftexp{(i_{1}...i_{l-k-1}i_{l-k+1}...i_{l})}{x}^{\prime}_{0,l-1}
+\leftexp{(i_{1}...i_{l})}{\dot{g}}^{\prime}_{0,l}: m=0
\end{align}
and:
\begin{align}
 \leftexp{(i_{1}...i_{l})}{\tilde{g}'}_{m,l}=-2\mu\hat{\chi}\cdot\leftexp{(i_{1}...i_{l})}{\hat{\slashed{\mu}}}_{2,m,l}
+\xi\cdot\slashed{d}\leftexp{(i_{1}...i_{l})}{x}^{\prime}_{m-1,l}\\\notag
+m\Lambda\leftexp{(i_{1}...i_{l})}{x}^{\prime}_{m-1,l}+\sum_{k=0}^{l-1}\leftexp{(R_{i_{l-k}})}{Z}\leftexp{(i_{1}...i_{l-k-1}i_{l-k+1}...i_{l})}{x}^{\prime}_{m,l-1}
+\leftexp{(i_{1}...i_{l})}{\dot{g}'}_{m,l} : m\geq 1
\end{align}
Let us define:
\begin{equation}
 X'_{m,l}=\max_{i_{1}...i_{l}}\|\leftexp{(i_{1}...i_{l})}{x}^{\prime}_{m,l}(t)\|_{L^{2}([0,\epsilon_{0}]\times S^{2})}
\end{equation}
We shall estimate $\|\tilde{g}^{\prime}_{m,l}(t)\|_{L^{2}([0,\epsilon_{0}]\times S^{2})}$ in terms of $X'_{0,l}(t)$ and $X_{l}(t)$, defined by (8.390), 
for $m=0$, and in terms of $X'_{m,l}(t)$ and $X'_{m-1,l+1}(t)$, for $m\geq 1$. Consider first the second term on the right of (9.244) and (9.245). The 
$S_{t,u}$ 1-form $\xi$ is given by (9.59). By the bootstrap assumption $\textbf{A}, \textbf{E}$ and $\textbf{F}$ we have:
\begin{equation}
 |\xi|\leq C\delta_{0}(1+t)^{-1}[1+\log(1+t)]
\end{equation}
Thus, by the pointwise estimates (9.119):
\begin{align}
 \|\xi\cdot\leftexp{(i_{1}...i_{l})}{x}_{l}\|_{L^{2}([0,\epsilon_{0}]\times S^{2})}
\leq C\delta_{0}(1+t)^{-1}[1+\log(1+t)]X_{l}(t)\\\notag
\|\xi\cdot\slashed{d}\leftexp{(i_{1}...i_{l})}{x'}_{m-1,l}\|_{L^{2}([0,\epsilon_{0}]\times S^{2})}
\leq C\delta_{0}(1+t)^{-2}[1+\log(1+t)]X'_{m-1,l+1}(t)
\end{align}
Consider next the third term on the right of (9.245). From (6.89) and (6.99) we have:
\begin{equation}
 |\Lambda|\leq C\delta_{0}(1+t)^{-1}[1+\log(1+t)]
\end{equation}
Thus, by the pointwise estimate (9.121):
\begin{equation}
 \|m\Lambda\leftexp{(i_{1}...i_{l})}{x}^{\prime}_{m,l}\|_{L^{2}([0,\epsilon_{0}]\times S^{2})}
\leq Cm\delta_{0}(1+t)^{-2}[1+\log(1+t)]X'_{m-1,l+1}(t)
\end{equation}

    Consider next the term next to the last term on the right in (9.244) and (9.245). By the pointwise estimates (9.127) and (8.393):
\begin{equation}
 \|\sum_{k=0}^{l-1}\leftexp{(R_{i_{l-k}})}{Z}\leftexp{(i_{1}...i_{l-k-1}i_{l-k+1}...i_{l})}{x'}_{m,l-1}\|_{L^{2}([0,\epsilon_{0}]\times S^{2})}
\leq Cl\delta_{0}(1+t)^{-2}[1+\log(1+t)]X'_{m,l}(t)
\end{equation}
Finally, we consider the first term on the right in (9.244) and (9.245). By $\textbf{F2}$ we have:
\begin{equation}
 \|2\mu\hat{\chi}\cdot\leftexp{(i_{1}...i_{l})}{\hat{\slashed{\mu}}}_{2,m,l}\|_{L^{2}([0,\epsilon_{0}]\times S^{2})}
\leq C\delta_{0}(1+t)^{-2}[1+\log(1+t)]\|\mu|\leftexp{(i_{1}...i_{l})}{\hat{\slashed{\mu}}}_{2,m,l}|(t)\|_{L^{2}([0,\epsilon_{0}]\times S^{2})}
\end{equation}
We now appeal to (9.190). In view of (8.396), (9.190) is equivalent to:
\begin{align}
 \|\mu|\leftexp{(i_{1}...i_{l})}{\hat{\slashed{\mu}}}_{2,m,l}|(t,u)\|_{L^{2}(S^{2})}
\leq C\|\leftexp{(i_{1}...i_{l})}{x}^{\prime}_{m,l}(t,u)\|_{L^{2}(S^{2})}+C\|\leftexp{(i_{1}...i_{l})}{\check{f}'}_{m,l}(t,u)\|_{L^{2}(S^{2})}\\\notag
+C\|\mu\leftexp{(i_{1}...i_{l})}{d}_{m,l}(t,u)\|_{L^{2}(S^{2})}+C\|\mu|\leftexp{(i_{1}...i_{l})}{e}_{m,l}|(t,u)\|_{L^{2}(S^{2})}\\\notag
+C\delta_{0}(1+t)^{-1}[1+\log(1+t)]\||\slashed{d}\leftexp{(i_{1}...i_{l})}{\mu}_{m,l}|(t,u)\|_{L^{2}(S^{2})}
\end{align}
Taking $L^{2}$ norms on $[0,\epsilon_{0}]$ then yields:
\begin{align}
 \|\mu|\leftexp{(i_{1}...i_{l})}{\hat{\slashed{\mu}}}_{2,m,l}|(t)\|_{L^{2}([0,\epsilon_{0}]\times S^{2})}
\leq C\|\leftexp{(i_{1}...i_{l})}{x}^{\prime}_{m,l}(t)\|_{L^{2}([0,\epsilon_{0}]\times S^{2})}+
C\|\leftexp{(i_{1}...i_{l})}{\check{f}'}_{m,l}(t)\|_{L^{2}([0,\epsilon_{0}]\times S^{2})}\\\notag
+C\|\mu\leftexp{(i_{1}...i_{l})}{d}_{m,l}(t)\|_{L^{2}([0,\epsilon_{0}]\times S^{2})}+
C\|\mu|\leftexp{(i_{1}...i_{l})}{e}_{m,l}|(t)\|_{L^{2}([0,\epsilon_{0}]\times S^{2})}\\\notag
+C\delta_{0}(1+t)^{-1}[1+\log(1+t)]\||\slashed{d}\leftexp{(i_{1}...i_{l})}{\mu}_{m,l}|(t)\|_{L^{2}([0,\epsilon_{0}]\times S^{2})}
\end{align}
    
     In view of (9.248), (9.250), (9.251), (9.252) and (9.254), we obtain from (9.244) and (9.245):
\begin{align}
 \|\leftexp{(i_{1}...i_{l})}{\tilde{g}}^{\prime}_{0,l}(t)\|_{L^{2}([0,\epsilon_{0}]\times S^{2})}\leq 
C(l+1)\delta_{0}(1+t)^{-2}[1+\log(1+t)]X'_{0,l}(t)\\\notag
+C\delta_{0}(1+t)^{-1}[1+\log(1+t)]X_{l}(t)+\leftexp{(i_{1}...i_{l})}{Q'}_{0,l}(t) :  m=0\\\notag
\|\leftexp{(i_{1}...i_{l})}{\tilde{g}}^{\prime}_{m,l}(t)\|_{L^{2}([0,\epsilon_{0}]\times S^{2})}\leq 
C(l+1)\delta_{0}(1+t)^{-2}[1+\log(1+t)]X'_{m,l}(t)\\\notag
+C(m+1)\delta_{0}(1+t)^{-2}[1+\log(1+t)]X'_{m-1,l+1}(t)\\\notag
+\leftexp{(i_{1}...i_{l})}{Q}^{\prime}_{m,l}(t)  :  m\geq 1
\end{align}
where:
\begin{align}
 \leftexp{(i_{1}...i_{l})}{Q}^{\prime}_{m,l}(t)=C\delta_{0}(1+t)^{-3}[1+\log(1+t)]^{2}\||\slashed{d}\leftexp{(i_{1}...i_{l})}{\mu}_{m,l}|(t)\|
_{L^{2}([0,\epsilon_{0}]\times S^{2})}\\\notag+C\delta_{0}(1+t)^{-2}[1+\log(1+t)]\|\leftexp{(i_{1}...i_{l})}{\check{f'}}_{m,l}(t)\|
_{L^{2}([0,\epsilon_{0}]\times S^{2})}\\\notag+C\delta_{0}(1+t)^{-2}[1+\log(1+t)]\|\mu\leftexp{(i_{1}...i_{l})}{d}_{m,l}(t)\|
_{L^{2}([0,\epsilon_{0}]\times S^{2})}\\\notag+C\delta_{0}(1+t)^{-2}[1+\log(1+t)]\|\mu|\leftexp{(i_{1}...i_{l})}{e}_{m,l}|(t)\|
_{L^{2}([0,\epsilon_{0}]\times S^{2})}\\\notag+\|\leftexp{(i_{1}...i_{l})}{\dot{g}'}_{m,l}(t)\|_{L^{2}([0,\epsilon_{0}]\times S^{2})}
\end{align}

    We now return to (9.193). Taking $L^{2}$ norms on $[0,\epsilon_{0}]\times S^{2}$ yields:
\begin{equation}
 \|\leftexp{(i_{1}...i_{l})}{x}^{\prime}_{m,l}(t)\|_{L^{2}([0,\epsilon_{0}]\times S^{2})}
\leq\|\leftexp{(i_{1}...i_{l})}{F}^{\prime}_{m,l}(t)\|_{L^{2}([0,\epsilon_{0}]\times S^{2})}
+\|\leftexp{(i_{1}...i_{l})}{G}^{\prime}_{m,l}(t)\|_{L^{2}([0,\epsilon_{0}]\times S^{2})}
\end{equation}
Substituting (9.239), (9.241) and (9.255) then yields:
\begin{align}
 \|\leftexp{(i_{1}...i_{l})}{x}^{\prime}_{0,l}(t)\|_{L^{2}([0,\epsilon_{0}]\times S^{2})}\leq\leftexp{(i_{1}...i_{l})}{B}^{\prime}_{0,l}(t)\\\notag
+C\delta_{0}(1+t)^{-2}[1+\log(1+t)]^{2}\int_{0}^{t}(1+t')[1+\log(1+t')]X_{l}(t')dt'\\\notag
+C(l+1)\delta_{0}(1+t)^{-2}[1+\log(1+t)]^{2}\int_{0}^{t}[1+\log(1+t')]X'_{0,l}(t')dt' :  m=0\\\notag
\|\leftexp{(i_{1}...i_{l})}{x}^{\prime}_{m,l}(t)\|_{L^{2}([0,\epsilon_{0}]\times S^{2})}\leq \leftexp{(i_{1}...i_{l})}{B}^{\prime}_{m,l}(t)\\\notag
+C(m+1)\delta_{0}(1+t)^{-2}[1+\log(1+t)]^{2}\int_{0}^{t}[1+\log(1+t')]X'_{m-1,l+1}(t')dt'\\\notag
+C(l+1)\delta_{0}(1+t)^{-2}[1+\log(1+t)]^{2}\int_{0}^{t}[1+\log(1+t')]X'_{m,l}(t')dt' :  m\geq 1
\end{align}
where:
\begin{align}
 \leftexp{(i_{1}...i_{l})}{B}^{\prime}_{m,l}(t)=C(1+t)^{-2}[1+\log(1+t)]^{2}\|\leftexp{(i_{1}...i_{l})}{x}^{\prime}_{m,l}(0)\|_{L^{2}(\Sigma^{\epsilon_{0}}_{0})}\\\notag
+C(1+t)^{-1}\{\leftexp{(i_{1}...i_{l})}{\bar{P}}^{\prime(0)}_{m,l,a}(t)+(1+t)^{-1/2}\leftexp{(i_{1}...i_{l})}{\bar{P}}^{\prime(1)}_{m,l,a}(t)\}
\bar{\mu}^{-a}_{m}(t)\\\notag
+C(1+t)^{-2}[1+\log(1+t)]^{2}\int_{0}^{t}(1+t')^{2}\leftexp{(i_{1}...i_{l})}{Q}^{\prime}_{m,l}(t')dt'
\end{align}
Taking in (9.258) the maximum over $i_{1}...i_{l}$ and recalling the definition (9.246) we obtain:
\begin{align}
 X'_{0,l}\leq B'_{0,l}(t)+C\delta_{0}(1+t)^{-2}[1+\log(1+t)]^{2}\int_{0}^{t}(1+t')[1+\log(1+t')]X_{l}(t')dt'\\\notag
+C(l+1)\delta_{0}(1+t)^{-2}[1+\log(1+t)]^{2}\int_{0}^{t}[1+\log(1+t')]X'_{0,l}(t')dt' :  m=0\\\notag
X'_{m,l}\leq B'_{m,l}(t)+C(m+1)\delta_{0}(1+t)^{-2}[1+\log(1+t)]^{2}\int_{0}^{t}[1+\log(1+t')]X'_{m-1,l+1}(t')dt'\\\notag
+C(l+1)\delta_{0}(1+t)^{-2}[1+\log(1+t)]^{2}\int_{0}^{t}[1+\log(1+t')]X'_{m,l}(t')dt' :  m\geq 1
\end{align}
where:
\begin{equation}
 B'_{m,l}(t)=\max_{i_{1}...i_{l}}\leftexp{(i_{1}...i_{l})}{B}^{\prime}_{m,l}(t)
\end{equation}
Setting:
\begin{align*}
 l=n-1-m
\end{align*}
(9.260) is a set of ordinary integral inequalities for $X'_{m,n-1-m}, m=0,...,n-1$. For $m=0$, we have an integral inequality for $X'_{0,n-1}$
containing on the right hand side the quantity $X_{n-1}$ which has been estimated in Chapter 8. For $m\geq1$, we have an integral inequality for 
$X'_{m,n-1-m}$ containing on the right-hand side the quantity $X'_{m-1,n-1-(m-1)}$. Thus, the integral inequalities, considered successively 
in the order of increasing $m$, depend only on quantities already estimated. 
Setting:
\begin{equation}
 Y'_{l,m}(t)=\int_{0}^{t}[1+\log(1+t')]X'_{m,l}(t')dt'
\end{equation}
and recalling (8.406), since
\begin{align*}
 \frac{dY'_{m,l}(t)}{dt}=[1+\log(1+t)]X'_{m,l}(t)
\end{align*}
the $Y'_{m,l}$ satisfy:
\begin{align}
 \frac{dY'_{0,l}(t)}{dt}\leq[1+\log(1+t)]B'_{0,l}(t)+C\delta_{0}(1+t)^{-2}[1+\log(1+t)]^{3}Y_{l}(t)\\\notag
+C(l+1)\delta_{0}(1+t)^{-2}[1+\log(1+t)]^{3}Y'_{0,l}(t)  : m=0\\\notag
\frac{dY'_{m,l}(t)}{dt}\leq [1+\log(1+t)]B'_{m,l}(t)+C(m+1)\delta_{0}(1+t)^{-2}[1+\log(1+t)]^{3}Y'_{m-1,l+1}(t)\\\notag
+C(l+1)\delta_{0}(1+t)^{-2}[1+\log(1+t)]^{3}Y'_{m,l}(t) : m\geq 1
\end{align}
The integrating factor here is:
\begin{align*}
 e^{-C_{l}(t)}
\end{align*}
where $C_{l}(t)$ is of the form (8.409). We thus obtain:
\begin{align}
 Y'_{0,l}(t)\leq e^{C_{l}(t)}\int_{0}^{t}e^{-C_{l}(t')}[1+\log(1+t')]B'_{0,l}(t')dt'\\\notag
C\delta_{0}\int_{0}^{t}(1+t')^{-2}[1+\log(1+t')]^{3}Y_{l}(t')dt'  :  m=0\\\notag
Y'_{m,l}(t)\leq e^{C_{l}(t)}\int_{0}^{t}e^{-C_{l}(t')}[1+\log(1+t')]B'_{m,l}(t')dt'\\\notag
+C(m+1)\delta_{0}\int_{0}^{t}(1+t')^{-2}[1+\log(1+t')]^{3}Y'_{m-1,l+1}(t')dt' : m\geq 1
\end{align}
Since the integral
\begin{align*}
 \int_{0}^{\infty}\frac{[1+\log(1+t')]^{3}}{(1+t')^{2}}dt'
\end{align*}
is convergent, if $\delta_{0}$ satisfies a smallness condition of the form (8.411), taking into account the fact that $Y_{l}, Y'_{m-1,l+1}$ are 
non-decreasing functions of $t$, 
(9.264) implies:
\begin{align}
 Y'_{0,l}(t)\leq 2\int_{0}^{t}[1+\log(1+t)]B'_{0,l}(t')dt'+C\delta_{0}Y_{l}(t) : m=0\\\notag
Y'_{m,l}(t)\leq 2\int_{0}^{t}[1+\log(1+t')]B'_{m,l}(t')dt'+C(m+1)\delta_{0}Y'_{m-1,l+1}(t) : m\geq 1
\end{align}
Setting $l=n-1-m, m=0,1,...,n-1$, we then have:
\begin{equation}
 Y'_{0,n-1}(t)\leq 2\int_{0}^{t}[1+\log(1+t')]B'_{0,n-1}(t')dt'+C\delta_{0}Y_{n-1}(t)
\end{equation}
and, for $m=1,...,n-1$, using a modified version of Proposition 8.2 appropriate to inequalities, we deduce:
\begin{align}
 Y'_{m,n-1-m}(t)\leq (m+1)!(C\delta_{0})^{m}Y'_{0,n-1}(t)\\\notag
+2\sum_{k=0}^{m-1}\frac{(m+1)!}{(m+1-k)!}(C\delta_{0})^{k}\int_{0}^{t}[1+\log(1+t')]B'_{m-k,n-1-(m-k)}(t')dt'
\end{align}
Substituting the bound (8.412), these imply that:
\begin{align}
 Y'_{m,n-1-m}(t)\leq2\int_{0}^{t}(1+t')[1+\log(1+t')]B_{n-1}(t')dt'\\\notag
+2\sum_{k=0}^{m}\int_{0}^{t}[1+\log(1+t')]B'_{k,n-1-k}(t')dt' : m=0,...,n-1
\end{align}
provided $\delta_{0}$ is suitably small depending on $n$. The bounds (8.412) and (9.268) in turn imply:
\begin{align}
 Y_{n-1}(t)+nY'_{0,n-1}(t)\leq 2(n+1)\int_{0}^{t}(1+t')[1+\log(1+t')]B_{n-1}(t')dt'\\\notag
+2n\int_{0}^{t}[1+\log(1+t')]B'_{0,n-1}(t')dt'  :  m=0\\\notag
(m+1)Y'_{m-1,n-1-(m-1)}(t)+(n-m)Y'_{m,n-1-m}(t)\\\notag
\leq 2(n+1)\int_{0}^{t}(1+t')[1+\log(1+t')]B_{n-1}(t')dt'\\\notag
+2(n+1)\sum_{k=0}^{m-1}\int_{0}^{t}[1+\log(1+t')]B'_{k,n-1-k}(t')dt'\\\notag
+2(n-m)\int_{0}^{t}[1+\log(1+t')]B'_{m,n-1-m}(t')dt'  : m=1,...,n-1
\end{align}
Recalling then (9.258) and setting $l=n-1-m$, we conclude that:
\begin{align}
 \|\leftexp{(i_{1}...i_{l})}{x}^{\prime}_{0,n-1}(t)\|_{L^{2}([0,\epsilon_{0}]\times S^{2})}\leq 
\leftexp{(i_{1}...i_{n-1})}{B}^{\prime}_{0,n-1}(t)\\\notag
+C\delta_{0}(1+t)^{-2}[1+\log(1+t)]^{2}[Y_{n-1}(t)+nY'_{0,n-1}(t)]  :   \\\notag m=0\\\notag
 \|\leftexp{(i_{1}...i_{n-1-m})}{x}^{\prime}_{m,n-1-m}(t)\|_{L^{2}([0,\epsilon_{0}]\times S^{2})}\leq 
\leftexp{(i_{1}...i_{n-1-m})}{B}^{\prime}_{m,n-1-m}(t)\\\notag
C\delta_{0}(1+t)^{-2}[1+\log(1+t)]^{2}\cdot[(m+1)Y'_{m-1,n-1-(m-1)}(t)+(n-m)Y'_{m,n-1-m}(t)] : m\geq 1
\end{align}

So finally we obtain the following estimates by substituting (9.269) in (9.270):
\begin{align}
 \|\leftexp{(i_{1}...i_{l})}{x}^{\prime}_{0,n-1}(t)\|_{L^{2}([0,\epsilon_{0}]\times S^{2})}\leq
\leftexp{(i_{1}...i_{n-1})}{B}^{\prime}_{0,n-1}(t)\\\notag
+C\delta_{0}(1+t)^{-2}[1+\log(1+t)]^{2}\{2(n+1)\int_{0}^{t}(1+t^{\prime})[1+\log(1+t^{\prime})]B_{n-1}(t^{\prime})dt^{\prime}\\\notag
+2n\int_{0}^{t}[1+\log(1+t^{\prime})]B^{\prime}_{0,n-1}(t^{\prime})dt^{\prime}\} : m=0
\end{align}
and
\begin{align}
 \|\leftexp{(i_{1}...i_{n-1-m})}{x}^{\prime}_{m,n-1-m}(t)\|_{L^{2}([0,\epsilon_{0}]\times S^{2})}\leq
\leftexp{(i_{1}...i_{n-1-m})}{B}^{\prime}_{m,n-1-m}(t)\\\notag
+C\delta_{0}(1+t)^{-2}[1+\log(1+t)]^{2}\{2(n+1)\int_{0}^{t}(1+t^{\prime})[1+\log(1+t^{\prime})]B_{n-1}(t^{\prime})dt^{\prime}\\\notag
+2(n+1)\sum_{k=0}^{m-1}\int_{0}^{t}[1+\log(1+t^{\prime}]B^{\prime}_{k,n-1-k}(t^{\prime})dt^{\prime}\\\notag
+2(n-m)\int_{0}^{t}[1+\log(1+t^{\prime})]B^{\prime}_{m,n-1-m}(t^{\prime})dt^{\prime}\} : m\geq 1
\end{align}

\chapter{Control of the Angular Derivatives of the First Derivatives of the $x^{i}$.\\
Assumptions and Estimates in Regard to $\chi$}

\section{Preliminary}
One of the aims of this chapter is concerned with derivation of estimates for angular derivatives of the spatial rectangular coordinates $x^{i}$ and $\hat{T}x^{i}: i=1,2,3$  
with respect to $R_{j}: j=1,2,3$. The estimates are based on bootstrap assumptions in regard to $\chi$. These estimates are then used to obtain estimates for angular 
derivatives of the deformation tensors of the commutation fields.

    Before we proceed, we must discuss  the general definition of the operator $\slashed{\mathcal{L}}_{X}$, with $X$ an arbitrary $S_{t,u}$-tangential vectorfield, as well 
as the operator $\slashed{\mathcal{L}}_{L}$, $\slashed{\mathcal{L}}_{T}$, as applied to an arbitrary $T^{q}_{p}$ $S_{t,u}$ tensorfield. Consider first the case of an 
$S_{t,u}$ 1-form $\xi$. Let $X$ be an $S_{t,u}$-tangential vectorfield. Then $\slashed{\mathcal{L}}_{X}\xi$ is a notion intrinsic to $S_{t,u}$. To define $\slashed{\mathcal{L}}_{L}\xi$ 
however, we must consider the extension of $\xi$ to a given $C_{u}$ as an $S_{t,u}$ 1-form, that is, by requiring $\xi(L)=0$. We then define $\slashed{\mathcal{L}}_{L}\xi$
to be the restriction to $TS_{t,u}$ of the usual Lie derivative $\mathcal{L}_{L}\xi$, a notion intrinsic to $C_{u}$. Note that in any case $(\mathcal{L}_{L}\xi)(L)=0$. 
We can define $\slashed{\mathcal{L}}_{T}\xi$ in a similar way by replacing $C_{u}$ by $\Sigma_{t}$.

    The case of any $p$-covariant $S_{t,u}$ tensorfield in the role of $\xi$, is formally identical to the case of an $S_{t,u}$ 1-form. Consider next the case of an 
$S_{t,u}$-tangential vectorfield $Y$. $X$ be another $S_{t,u}$ tangential vectorfield. Then $\slashed{\mathcal{L}}_{X}Y$ is defined to be simply $\mathcal{L}_{X}Y=[X,Y]$.
From Lemma 8.2, $\mathcal{L}_{L}Y=[L,Y]=\leftexp{(Y)}{Z}$ is an $S_{t,u}$-tangential vectorfield, so $\slashed{\mathcal{L}}_{L}Y$ is defined to be $\mathcal{L}_{L}Y$.
Similarly (see Lemma 10.18), $\slashed{\mathcal{L}}_{T}Y$ is defined to be $\mathcal{L}_{T}Y$.

    The case of any $q$-contravariant $S_{t,u}$ tensorfield is formally identical. And also, the general case of an arbitrary $T^{q}_{p}$ $S_{t,u}$ tensorfield $\vartheta$
is similar. Note that the above definitions of $\slashed{\mathcal{L}}_{X}\vartheta$, $\slashed{\mathcal{L}}_{L}\vartheta$, $\slashed{\mathcal{L}}_{T}\vartheta$, do not
depend on the acoustical metric $g$. So we have:
\begin{equation}
 \slashed{\mathcal{L}}_{X}(\vartheta\otimes\varphi)=(\slashed{\mathcal{L}}_{X}\vartheta)\otimes\varphi+\vartheta\otimes(\slashed{\mathcal{L}}_{X}\varphi)
\end{equation}
with $X$ an arbitrary $S_{t,u}$-tangential vectorfield. And also:
\begin{align}
 \slashed{\mathcal{L}}_{L}(\vartheta\otimes\varphi)=(\slashed{\mathcal{L}}_{L}\vartheta)\otimes\varphi+\vartheta\otimes(\slashed{\mathcal{L}}_{L}\varphi)
\end{align}
\begin{align}
 \slashed{\mathcal{L}}_{T}(\vartheta\otimes\varphi)=(\slashed{\mathcal{L}}_{T}\vartheta)\otimes\varphi+\vartheta\otimes(\slashed{\mathcal{L}}_{T}\varphi)
\end{align}
We also note that the same results hold with $Q$ in the role of $L$.

    From the above discussion, we have, for $X,Y\in TS_{t,u}$:
\begin{align*}
 (\mathcal{L}_{L}g)(X,Y)=L(g(X,Y))-g([L,X],Y)-g(X,[L,Y])\\
=L(\slashed{g}(X,Y))-\slashed{g}([L,X],Y)-\slashed{g}(X,[L,Y])=(\slashed{\mathcal{L}}_{L}\slashed{g})(X,Y)
\end{align*}
That is:
\begin{align}
 \leftexp{(L)}{\slashed{\pi}}=\slashed{\mathcal{L}}_{L}\slashed{g}=2\chi
\end{align}
Similarly, we have:
\begin{align}
 \leftexp{(T)}{\slashed{\pi}}=\slashed{\mathcal{L}}_{T}\slashed{g}
\end{align}
and 
\begin{align*}
 (\mathcal{L}_{T}\bar{g})(X,Y)=T(\bar{g}(X,Y))-\bar{g}([T,X],Y)-\bar{g}(X,[T,Y])\\
=T(\slashed{g}(X,Y))-\slashed{g}([T,X],Y)-\slashed{g}(X,[T,Y])=(\slashed{\mathcal{L}}_{T}\slashed{g})(X,Y)
\end{align*}
thus,
\begin{align}
 \slashed{\mathcal{L}}_{T}\slashed{g}=2\kappa\theta
\end{align}

     Given now a positive integer $l$ let us denote by $\slashed{\textbf{E}}_{l}$ the bootstrap assumption that there is a constant $C$ independent of $s$ such that for all
$t\in[0,s]$:
\begin{align*}
 \slashed{\textbf{E}}_{l}\quad:\quad\max_{i_{1}...i_{l}}\max_{\alpha}\|R_{i_{l}}...R_{i_{1}}\psi_{\alpha}\|_{L^{\infty}(\Sigma_{t}^{\epsilon_{0}})}
\leq C\delta_{0}(1+t)^{-1}
\end{align*}
Also, for convenience of notation, let us denote by $\slashed{\textbf{E}}_{0}$ the basic bootstrap assumption $\textbf{E}\textbf{1}$ of Chapter 6:
\begin{align*}
 \slashed{\textbf{E}}_{0}\quad:\quad\max_{\alpha}|\psi_{\alpha}|\leq C\delta_{0}(1+t)^{-1}
\end{align*}
Given a positive integer $l$ we then denote by $\slashed{\textbf{E}}_{[l]}$ the bootstrap assumption:
\begin{align*}
 \slashed{\textbf{E}}_{[l]}\quad:\quad \slashed{\textbf{E}}_{0}\quad\textrm{and}...\textrm{and}\quad\slashed{\textbf{E}}_{l}
\end{align*}
Given a positive integer $n$ we denote:
\begin{align}
 n_{*}=
\begin{cases}
 \frac{n}{2}&:\textrm{if}\quad n \quad\textrm{is even}\\
\frac{n-1}{2}&:\textrm{if}\quad n \quad\textrm{is odd}
\end{cases}
\end{align}

$\textbf{Lemma 10.1}$ Let $G$ be a smooth function of the variations $(\psi_{\alpha}\quad:\quad \alpha=0,1,2,3)$. Suppose that the bootstrap assumption 
$\slashed{\textbf{E}}_{[l_{*}]}$ holds for some positive integer $l$. Then there are constants $C$, $C_{l}$ independent of $s$ such that the following estimate holds:
\begin{align*}
 \|R_{i_{l}}...R_{i_{1}}G\|_{L^{2}(\Sigma_{t}^{\epsilon_{0}})}\\
\leq C\sum_{\alpha}\|R_{i_{l}}...R_{i_{1}}\psi_{\alpha}\|_{L^{2}(\Sigma_{t}^{\epsilon_{0}})}+
C_{l}\delta_{0}(1+t)^{-1}\sum_{k=1}^{l-1}\max_{i_{1}...i_{k}}\max_{\alpha}\|R_{i_{k}}...R_{i_{1}}\psi_{\alpha}\|_{L^{2}(\Sigma_{t}^{\epsilon_{0}})}
\end{align*}
$Proof$: It is a direct calculation. $\qed$ \vspace{7mm}

     Consider now the functions $\hat{T}^{j}$. By (3.201)-(3.203), we have:
\begin{align}
 R_{i}\hat{T}^{j}=R^{A}_{i}\slashed{q}_{A}=\slashed{q}_{i}\cdot\slashed{d}x^{j}
\end{align}
where
\begin{align}
 \slashed{q}_{i}=m\cdot R_{i}
\end{align}
Here, $m$ is given by:
\begin{align}
 m=m_{b}\cdot\slashed{g}^{-1}
\end{align}
where
\begin{align}
 m_{b}=\theta=-\alpha^{-1}\chi+\slashed{k}
\end{align}
    
     Recall the definition of functions $y^{i}$ in Chapter 6:
\begin{align}
 \hat{T}^{j}=-\frac{x^{j}}{1-u+t}+y^{j}
\end{align}
Since $R_{i}x^{j}=R_{i}\cdot\slashed{d}x^{j}$, we have:
\begin{align}
 R_{i}y^{j}=\slashed{q}^{\prime}_{i}\cdot\slashed{d}x^{j}
\end{align}
where
\begin{align}
 \slashed{q}_{i}^{\prime}=\slashed{q}_{i}+\frac{R_{i}}{1-u+t}
\end{align}
By (10.9)-(10.11) we have:
\begin{align}
 \slashed{q}_{i}^{\prime}=m^{\prime}\cdot R_{i}
\end{align}
where
\begin{align}
 m^{\prime}=m^{\prime}_{b}\cdot\slashed{g}^{-1}
\end{align}
and
\begin{align}
 m^{\prime}_{b}=m_{b}+\frac{\slashed{g}}{1-u+t}
\end{align}
is the 2-covariant $S_{t,u}$ tensorfield:
\begin{align}
 m^{\prime}_{b}=-\frac{(\alpha^{-1}-1)\slashed{g}}{1-u+t}-\alpha^{-1}(\chi-\frac{\slashed{g}}{1-u+t})+\slashed{k}
\end{align}

     We consider next the angular derivatives of the rectangular coordinate functions $x^{j}$. According to (6.6) we have:
\begin{align}
 (R_{i})^{j}=\Pi_{k}^{j}(\mathring{R}_{i})^{k}
\end{align}
where $\Pi^{j}_{k}$ are the components of $\Pi$, the $\bar{g}$-orthogonal projection to $S_{t,u}$ on $\Sigma_{t}$:
\begin{align}
 \Pi^{j}_{k}=\delta^{j}_{k}-\bar{g}_{kl}\hat{T}^{l}\hat{T}^{j}=\delta^{j}_{k}
 -\hat{T}^{k}\hat{T}^{j}
\end{align}
and
\begin{align}
 (\mathring{R}_{i})^{j}=\epsilon_{imj}x^{m}
\end{align}
Recall (6.130):
\begin{align}
 \bar{g}(\mathring{R}_{i},\hat{T})=\lambda_{i}=\epsilon_{imj}x^{m}y^{j}
\end{align}
we then obtain:
\begin{align}
 R_{i}=\mathring{R}_{i}-\lambda_{i}\hat{T}^{j}\partial_{j}
\end{align}

     Define, for each non-negative integer $k$, the functions:
\begin{align}
 \leftexp{(k)}{x}^{j}_{i_{1}...i_{k}}=\mathring{R}_{i_{k}}...\mathring{R}_{i_{1}}x^{j}
\end{align}
These are linear functions of the rectangular coordinates, therefore the
\begin{align}
 \leftexp{(k)}{c}^{j}_{l;i_{1}...i_{k}}=\frac{\partial\leftexp{(k)}{x}^{j}_{i_{1}...i_{k}}}{\partial x^{l}}
\end{align}
are constants. Define for each non-negative integer $k$ the functions $\leftexp{(k)}{\delta}^{j}_{i_{1}...i_{j}}$ by the equation:
\begin{align}
 R_{i_{k}}...R_{i_{1}}x^{j}=\leftexp{(k)}{x}^{j}_{i_{1}...i_{k}}-\leftexp{(k)}{\delta}^{j}_{i_{1}...i_{k}}
\end{align}
In particular,
\begin{align}
 \leftexp{(0)}{\delta}^{j}=0
\end{align}
and by (10.23),
\begin{align}
 \leftexp{(1)}{\delta}^{j}_{i}=\lambda_{i}\hat{T}^{j}
\end{align}
Replacing $k$ by $k-1$ in (10.26) and applying $R_{i_{k}}$ we obtain:
\begin{align}
 R_{i_{k}}R_{i_{k-1}}...R_{i_{1}}x^{j}=R_{i_{k}}\leftexp{(k-1)}{x}^{j}_{i_{1}...i_{k-1}}-R_{i_{k}}\leftexp{(k-1)}{\delta}^{j}_{i_{1}...i_{k-1}}
\end{align}
By (10.23)-(10.25),
\begin{align}
 R_{i_{k}}\leftexp{(k-1)}{x}^{j}_{i_{1}...i_{k-1}}=\leftexp{(k)}{x}^{j}_{i_{1}...i_{k}}-\leftexp{(k-1)}{c}^{j}_{l;i_{1}...i_{k-1}}\lambda_{i_{k}}\hat{T}^{l}
\end{align}
Substituting this in (10.29) and comparing with (10.26), we get:
\begin{align}
 \leftexp{(k)}{\delta}^{j}_{i_{1}...i_{k}}=R_{i_{k}}\leftexp{(k-1)}{\delta}^{j}_{i_{1}...i_{k-1}}+\leftexp{(k-1)}{c}^{j}_{l;i_{1}...i_{k-1}}\lambda_{i_{k}}\hat{T}^{l}
\end{align}
Now we can apply Proposition 8.2 with $A_{n}$ in the role of $R_{i_{k}}$, $x_{n}$ in the role of $\leftexp{(k)}{\delta}^{j}_{i_{1}...i_{k}}$ and $y_{n}$ in the role of
\begin{align*}
 \leftexp{(k-1)}{c}^{j}_{l;i_{1}...i_{k-1}}\lambda_{i_{k}}\hat{T}^{l}
\end{align*}
Since by (10.27) $x_{0}=0$ we obtain:
\begin{align}
 \leftexp{(k)}{\delta}^{j}_{i_{1}...i_{k}}=\sum_{m=1}^{k}\leftexp{(m-1)}{c}^{j}_{l;i_{1}...i_{m-1}}R_{i_{k}}...R_{i_{k-m+2}}(\lambda_{i_{m}}\hat{T}^{l})
\end{align}
     Since the $\leftexp{(k)}{x}$ are linear functions of the rectangular coordinates, we have, from (10.25),
\begin{align}
 \leftexp{(k)}{x}^{j}_{i_{1}...i_{k}}=\leftexp{(k)}{c}^{j}_{l;i_{1}...i_{k}}x^{l}
\end{align}
From (10.21) and (10.24) we obtain:
\begin{align}
 \leftexp{(0)}{x}^{j}=x^{j},\quad \leftexp{(1)}{x}^{j}_{i}=\epsilon_{ilj}x^{l}
\end{align}
Moreover,
\begin{align}
 \leftexp{(2)}{x}^{j}_{ik}=\delta_{jk}x^{i}-\delta_{ik}x^{j}
\end{align}
For $k\geq 2$ we have from (10.24) and (10.35):
\begin{align}
 \leftexp{(k)}{x}^{j}_{i_{1}...i_{k}}=\mathring{R}_{i_{k}}...\mathring{R}_{i_{3}}\leftexp{(2)}{x}^{j}_{i_{1}i_{2}}\\\notag
=\mathring{R}_{i_{k}}...\mathring{R}_{i_{3}}(\delta_{ji_{2}}x^{i_{1}}-\delta_{i_{1}i_{2}}x^{j})\\\notag
=\delta_{ji_{2}}\leftexp{(k-2)}{x}^{i_{1}}_{i_{3}...i_{k}}-\delta_{i_{1}i_{2}}\leftexp{(k-2)}{x}^{j}_{i_{3}...i_{k}}
\end{align}
It follows from (10.36) that for $k\geq 2$:
\begin{align}
 \max_{j;i_{1}...i_{k}}|\leftexp{(k)}{x}^{j}_{i_{1}...i_{k}}|=\max_{j;i_{1}...i_{k-2}}|\leftexp{(k-2)}{x}^{j}_{i_{1}...i_{k-2}}|
\end{align}
while from (10.34):
\begin{align}
 \max_{j;i}|\leftexp{(1)}{x}^{j}_{i}|=\max_{j}|\leftexp{(0)}{x}^{j}|=\max_{j}|x^{j}|
\end{align}
Then we have:
\begin{align}
 \max_{j;i_{1}...i_{k}}|\leftexp{(k)}{x}^{j}_{i_{1}...i_{k}}|=\max_{j}|x^{j}|
\end{align}
which by (6.131) implies:
\begin{align}
 \max_{j;i_{1}...i_{k}}\|\leftexp{(k)}{x}^{j}_{i_{1}...i_{k}}\|_{L^{\infty}(\Sigma_{t}^{\epsilon_{0}})}\leq 1+t
\end{align}

$\textbf{Lemma 10.2}$ Let
\begin{align*}
 \max_{j;i_{1}...i_{k}}\|R_{i_{k}}...R_{i_{1}}y^{j}\|_{L^{\infty}(\Sigma_{t}^{\epsilon_{0}})}\leq C_{l}\delta_{0}(1+t)^{-1}[1+\log(1+t)]
\end{align*}
hold for $k=0,...,l$. Then if $\delta_{0}$ is suitably small (perhaps depending on $l$), we have:
\begin{align*}
 \max_{i;i_{1}...i_{k}}\|R_{i_{k}}...R_{i_{1}}\lambda_{i}\|_{L^{\infty}(\Sigma_{t}^{\epsilon_{0}})}\leq C_{l}\delta_{0}[1+\log(1+t)]
\end{align*}
and:
\begin{align*}
 \max_{j;ii_{1}...i_{k}}\|\leftexp{(k+1)}{\delta}^{j}_{ii_{1}...i_{k}}\|_{L^{\infty}(\Sigma_{t}^{\epsilon_{0}})}\leq C_{l}\delta_{0}
[1+\log(1+t)]
\end{align*}
for all $k=0,...,l$.

$Proof$. Noting (6.160) and (10.27), the lemma holds when $l=0$. Then by (10.12), (10.22), (10.26), (10.40) and (10.32), we can use an induction argument. $\qed$ \vspace{7mm}

$\textbf{Lemma 10.3}$ Let 
\begin{align*}
 \max_{j;i_{1}...i_{k}}\|R_{i_{k}}...R_{i_{1}}y^{j}\|_{L^{\infty}(\Sigma_{t}^{\epsilon_{0}})}\leq C_{l}\delta_{0}(1+t)^{-1}[1+\log(1+t)]
\end{align*}
 hold for $k=0,...,l_{*}$. Then if $\delta_{0}$ is suitably small (perhaps depending on $l$) we have:
\begin{align*}
 \max_{i;i_{1}...i_{k}}\|R_{i_{k}}...R_{i_{1}}\lambda_{i}\|_{L^{2}(\Sigma_{t}^{\epsilon_{0}})}\leq C_{l}(1+t)\mathcal{Y}_{[l]}
\end{align*}
and:
\begin{align*}
 \max_{j;ii_{1}...i_{k}}\|\leftexp{(k+1)}{\delta}^{j}_{ii_{1}...i_{k}}\|_{L^{2}(\Sigma_{t}^{\epsilon_{0}})}\leq C_{l}(1+t)\mathcal{Y}_{[l]}
\end{align*}
for all $k=0,...,l$. Here:
\begin{align*}
 \mathcal{Y}_{k}=\max_{j;i_{1}...i_{k}}\|R_{i_{k}}...R_{i_{1}}y^{j}\|_{L^{2}(\Sigma_{t}^{\epsilon_{0}})},\quad 
\mathcal{Y}_{[l]}=\sum_{k=0}^{l}\mathcal{Y}_{k}
\end{align*}
We also define:
\begin{align*}
 \mathcal{W}_{0}=\max_{\alpha}\|\psi_{\alpha}\|_{L^{2}(\Sigma_{t}^{\epsilon_{0}})}
\end{align*}
and for $k\slashed{=}0$:
\begin{align*}
 \mathcal{W}_{k}=\max_{\alpha;i_{1}...i_{k}}\|R_{i_{k}}...R_{i_{1}}\psi_{\alpha}\|_{L^{2}(\Sigma_{t}^{\epsilon_{0}})}
\end{align*}
and:
\begin{align*}
 \mathcal{W}_{[l]}=\sum_{k=0}^{l}\mathcal{W}_{k}
\end{align*}
$Proof$. By (10.22) (10.32) and (10.40), the lemma holds for $l=0$. We then use an induction argument. $\qed$

$\textbf{Lemma 10.4}$ For every non-negative integer $k$ we have:
\begin{align*}
 R_{i_{k}}...R_{i_{1}}R_{i}y^{j}=\leftexp{(k)}{\slashed{q}}^{\prime}_{ii_{1}...i_{k}}\cdot\slashed{d}x^{j}+\leftexp{(k)}{\slashed{r}}^{\prime}_{ii_{1}...i_{k}}
\end{align*}
Here,
\begin{align*}
 \leftexp{(k)}{\slashed{q}}^{\prime}_{ii_{1}...i_{k}}=\slashed{\mathcal{L}}_{R_{i_{k}}}...\slashed{\mathcal{L}}_{R_{i_{1}}}\slashed{q}^{\prime}_{i}
\end{align*}
and
\begin{align*}
 \leftexp{(k)}{\slashed{r}}^{j}_{ii_{1}...i_{k}}=\sum_{m=0}^{k-1}R_{i_{k}}...R_{i_{k-m+1}}\{\leftexp{(k-m-1)}{\slashed{q}}^{\prime}_{ii_{1}...i_{k-m-1}}
\cdot\slashed{d}(R_{i_{k-m}}x^{j})\}
\end{align*}
$Proof$. Since 
\begin{align*}
 \leftexp{(0)}{\slashed{q}}^{\prime}_{i}=\slashed{q}^{\prime}_{i},\quad \leftexp{(0)}{\slashed{r}}^{j}_{i}=0
\end{align*}
the lemma reduces for $k=0$ to equation (10.13). Then arguing by induction, and using Proposition 8.2, the lemma follows. $\qed$
\vspace{7mm}

     Now it is evident from (10.15) and the formulas for the coefficients $\leftexp{(k)}{\slashed{q}}^{\prime}_{ii_{1}...i_{k}}$ and 
$\leftexp{(k)}{\slashed{r}}^{j}_{ii_{1}...i_{k}}$ given in the statement of Lemma 10.4, that $\leftexp{(k)}{\slashed{q}}^{\prime}_{ii_{1}...i_{k}}$ contains the 
angular derivatives of the $R_{i}$ of order up to $k$ and $\leftexp{(k)}{\slashed{r}}^{j}_{ii_{1}...i_{k}}$ contains the angular derivatives of the $R_{i}$ 
of order up to $k-1$, where by ``the angular derivatives of the $R_{i}$ of order up to $l$'' we mean the $S_{t,u}$-tangential vectorfields $\slashed{\mathcal{L}}_{R_{i_{m}}}
...\slashed{\mathcal{L}}_{R_{i_{1}}}R_{i}$ for $m=0,...,l$. We shall thus derive the expressions below for these vectorfields. The basic formula is the formula for
\begin{align*}
 \slashed{\mathcal{L}}_{R_{i}}R_{j}=[R_{i},R_{j}]
\end{align*}
given by the following lemma. Let $w_{i}$ be
\begin{align}
 w_{i}=(w_{i})_{b}\cdot\slashed{g}^{-1}
\end{align}
where
\begin{align}
 (w_{i})_{b}=\epsilon_{ijm}y^{j}\bar{g}_{mn}\slashed{d}x^{n}
 =\epsilon_{ijm}y^{j}\slashed{d}x^{m}
\end{align}
 Let us then define:
\begin{align}
 \tilde{\slashed{q}^{\prime}}_{i}=\slashed{q}^{\prime}_{i}-w_{i}
\end{align}
$\textbf{Lemma 10.5}$ We have:
\begin{align*}
 [R_{i},R_{j}]=-\epsilon_{ijk}R_{k}+\lambda_{i}\tilde{\slashed{q}}^{\prime}_{j}
 -\lambda_{j}\tilde{\slashed{q}}^{\prime}_{i}
\end{align*}
$Proof$. We express:
\begin{align*}
 R_{i}=(R_{i})^{m}\partial_{m}
\end{align*}
According to (10.23):
\begin{align*}
 (R_{i})^{m}=(\mathring{R}_{i})^{m}-\lambda_{i}\hat{T}^{m}
\end{align*}
Thus, we have:
\begin{align}
 [R_{i},R_{j}]=\{R_{i}((R_{j})^{m})-R_{j}((R_{i})^{m})\}\partial_{m}\\\notag
=\{R_{i}((\mathring{R}_{j})^{m}-\lambda_{j}\hat{T}^{m})-R_{j}((\mathring{R}_{i})^{m}-\lambda_{i}\hat{T}^{m})\}\partial_{m}
\end{align}
Now,
\begin{align*}
 R_{i}((\mathring{R}_{j})^{m})=R_{i}(\epsilon_{jkm}x^{k})=\epsilon_{jkm}(R_{i})^{k}\\\notag
=-\epsilon_{jmk}((\mathring{R}_{i})^{k}-\lambda_{i}\hat{T}^{k})=-\epsilon_{jmk}\epsilon_{ilk}x^{l}+\epsilon_{jmk}\lambda_{i}\hat{T}^{k}
\end{align*}
Since,
\begin{align*}
 \epsilon_{jmk}\epsilon_{ilk}+\epsilon_{mik}\epsilon_{jlk}+\epsilon_{ijk}\epsilon_{mlk}=0
\end{align*}
we then obtain:
\begin{align}
 R_{i}((\mathring{R}_{j})^{m})-R_{j}((\mathring{R}_{i})^{m})=-\epsilon_{ijk}(\mathring{R}_{k})^{m}+(\epsilon_{ikm}\lambda_{j}-\epsilon_{jkm}\lambda_{i})\hat{T}^{k}
\end{align}
Substituting (10.45) in (10.44) yields:
\begin{align}
 [R_{i},R_{j}]=-\epsilon_{ijk}\mathring{R}_{k}-\lambda_{i}v_{j}+\lambda_{j}v_{i}\\\notag
-(R_{i}(\lambda_{j})-R_{j}(\lambda_{i}))\hat{T}-(\lambda_{j}R_{i}(\hat{T}^{m})-\lambda_{i}R_{j}(\hat{T}^{m}))\partial_{m}
\end{align}
where
\begin{align}
 v_{i}=\epsilon_{ikm}\hat{T}^{k}\partial_{m}
\end{align}
By (10.12), (10.13),
\begin{align}
 R_{i}(\hat{T}^{m})\partial_{m}=-(1-u+t)^{-1}R_{i}+\slashed{q}^{\prime}_{i}
\end{align}
Substituting (10.48) in (10.46), the latter becomes:
\begin{align}
 [R_{i},R_{j}]=-\epsilon_{ijk}\mathring{R}_{k}\\\notag
-(R_{i}(\lambda_{j})-R_{j}(\lambda_{i}))\hat{T}\\\notag
+\lambda_{i}(\slashed{q}^{\prime}_{j}-v^{\prime}_{j})-\lambda_{j}(\slashed{q}^{\prime}_{i}-v^{\prime}_{i})
\end{align}
where:
\begin{align}
 v^{\prime}_{i}=v_{i}+(1-u+t)^{-1}R_{i}
\end{align}
Since $[R_{i},R_{j}]$ is $S_{t,u}$-tangential, we may apply the projection operator $\Pi$ to the right hand side of (10.49). This annihilates the second term,
and since $\Pi\mathring{R}_{i}=R_{i}$, we obtain:
\begin{align}
 [R_{i},R_{j}]=-\epsilon_{ijk}R_{k}+\lambda_{i}(\slashed{q}^{\prime}_{j}-\Pi v^{\prime}_{j})-\lambda_{j}(\slashed{q}^{\prime}_{i}-\Pi v^{\prime}_{i})
\end{align}
In view of (10.43), we must show that:
\begin{align}
 \Pi v^{\prime}_{i}=w_{i}
\end{align}
We have, in view of (10.50),
\begin{align}
 \Pi v^{\prime}_{i}=\Pi\tilde{v}_{i}
\end{align}
where
\begin{align}
 \tilde{v}_{i}=v_{i}+(1-u+t)^{-1}\mathring{R}_{i}
\end{align}
In view of (10.47), (10.12) and (10.21), we have:
\begin{align}
 \tilde{v}_{i}=\epsilon_{ikm}y^{k}\partial_{m}
\end{align}
We then have, in terms of the frame $(X_{A}: A=1,2)$,
\begin{align}
 \Pi\tilde{v}_{i}=\bar{g}(\tilde{v}_{i},X_{B})(\slashed{g}^{-1})^{BA}X_{A}
\end{align}
and:
\begin{align}
 \bar{g}(\tilde{v}_{i},X_{B})=\bar{g}_{mn}(\tilde{v}_{i})^{m}X^{n}_{B}=\epsilon_{ikm}y^{k}\bar{g}_{mn}\slashed{d}_{B}x^{n}
 =\epsilon_{ikm}y^{k}\slashed{d}_{B}x^{m}
\end{align}
Comparing (10.57) with (10.42) we see that:
\begin{align}
 \bar{g}(\tilde{v}_{i},X_{B})=(w_{i})_{b}(X_{B})
\end{align}
Substituting (10.58) in (10.56) we conclude that:
\begin{align}
 \Pi\tilde{v}_{i}=w_{i}
\end{align}
Therefore by (10.53), the lemma follows. $\qed$ \vspace{7mm}

     Since by (10.15), (10.43),
\begin{align}
 \tilde{\slashed{q}}^{\prime}_{i}=m^{\prime}\cdot R_{i}-w_{i}
\end{align}
defining the $S_{t,u}$ tensorfields $\mu_{i}$ by:
\begin{align}
 \mu_{i}=\lambda_{i}m^{\prime}
\end{align}
and the $S_{t,u}$-tangential vectorfields $\nu_{ij}$:
\begin{align}
 \nu_{ij}=-\lambda_{i}w_{j}+\lambda_{j}w_{i}
\end{align}
the formula of Lemma 10.5 takes the form:
\begin{align}
 [R_{i},R_{j}]=-\epsilon_{ijk}R_{k}+\mu_{i}\cdot R_{j}-\mu_{j}\cdot R_{i}+\nu_{ij}
\end{align}

$\textbf{Lemma 10.6}$ For every non-negative integer $k$ we have:
\begin{align*}
 \slashed{\mathcal{L}}_{R_{i_{k}}}...\slashed{\mathcal{L}}_{R_{i_{1}}}R_{j}=\leftexp{(k)}{\alpha}^{m}_{j;i_{1}...i_{k}}R_{m}
+\leftexp{(k)}{\beta}_{j;i_{1}...i_{k}}\cdot R_{m}+\leftexp{(k)}{\gamma}_{j;i_{1}...i_{k}}
\end{align*}
The coefficients $\leftexp{(k)}{\alpha}^{m}_{j;i_{1}...i_{k}}$ are constants, given by:
\begin{align*}
 \leftexp{(k)}{\alpha}^{m}_{j;i_{1}...i_{k}}=(-1)^{k}\epsilon_{i_{k}n_{k}m}\epsilon_{i_{k-1}n_{k-1}n_{k}}...\epsilon_{i_{1}n_{1}n_{2}}\delta^{n_{1}}_{j}
\end{align*}
The coefficients $\leftexp{(k)}{\beta}_{j;i_{1}...i_{k}}=(\leftexp{(k)}{\beta}^{m}_{i_{1}...i_{k}}\quad:\quad m=1,2,3)$ are triplets of $T^{1}_{1}$-type $S_{t,u}$ 
tensorfields given by:
\begin{align*}
 \leftexp{(k)}{\beta}_{j;i_{1}...i_{k}}=\sum_{n=0}^{k-1}A_{i_{k}}...A_{i_{k-n+1}}\leftexp{(k-n)}{\rho}_{j;i_{1}...i_{k}}
\end{align*}
Here $A_{i}$ are the operators acting on triplets $z=(z^{m}\quad:\quad m=1,2,3)$ of $T^{1}_{1}$-type $S_{t,u}$ tensorfields given by:
\begin{align*}
 (A_{i}z)^{m}=\slashed{\mathcal{L}}_{R_{i}}z^{m}-\sum_{n}\epsilon_{inm}z^{n}+z^{m}\cdot\mu_{i}-(\sum_{n}z^{n}\cdot\mu_{n})\delta^{m}_{i}
\end{align*}
and the $\leftexp{(k)}{\rho}_{j;i_{1}...i_{k}}=(\leftexp{(k)}{\rho}^{m}_{j;i_{1}...i_{k}}\quad:\quad m=1,2,3)$ are the triplets of $T^{1}_{1}$-type $S_{t,u}$ tensorfields,
given by:
\begin{align*}
 \leftexp{(k)}{\rho}^{m}_{j;i_{1}...i_{k}}=\leftexp{(k-1)}{\alpha}^{m}_{j;i_{1}...i_{k}}\mu_{i_{k}}-(\sum_{n}\leftexp{(k-1)}{\alpha}^{n}_{j;i_{1}...i_{k}}\mu_{n})
\delta^{m}_{i_{k}}
\end{align*}
Finally, the coefficients $\leftexp{(k)}{\gamma}_{j;i_{1}...i_{k}}$ are $S_{t,u}$-tangential vectorfields, given by:
\begin{align*}
 \leftexp{(k)}{\gamma}_{j;i_{1}...i_{k}}=\sum_{n=0}^{k-1}\leftexp{(k-1-n)}{\alpha}^{m}_{j;i_{1}...i_{k-n-1}}\slashed{\mathcal{L}}_{R_{i_{k}}}...
\slashed{\mathcal{L}}_{R_{i_{k-n+1}}}\nu_{i_{k-n}m}\\\notag
+\sum_{n=0}^{k-1}\slashed{\mathcal{L}}_{R_{i_{k}}}...\slashed{\mathcal{L}}_{R_{i_{k-n+1}}}(\leftexp{(k-n-1)}{\beta}^{m}_{j;i_{1}...i_{k-n-1}}\cdot\nu_{i_{k-n}m})
\end{align*}
$Proof$. With: 
\begin{align}
 \leftexp{(0)}{\alpha}^{m}_{j}=\delta^{m}_{j},\quad \leftexp{(0)}{\beta}^{m}_{j}=0,\quad \leftexp{(0)}{\gamma}_{j}=0
\end{align}
the lemma is trivial for $k=0$. Moreover, one can easily check that the lemma reduces for $k=1$ to (10.63). Then arguing by induction, using Proposition 8.2 and 
noting (10.64), the lemma follows. $\qed$ \vspace{7mm}

     To analyze the operators $A_{i}$, we define the multiplication operators $M_{i}$. These are 3-dimensional matrices, the entries $(M_{i})^{m}_{n}:m,n=1,2,3$ of 
which are $T^{1}_{1}$-type of $S_{t,u}$ tensorfields given by:
\begin{align}
 (M_{i})^{m}_{n}=-\epsilon_{inm}I+\delta^{m}_{n}\mu_{i}-\delta_{i}^{m}\mu_{n}
\end{align}
Here $I$ is the $T^{1}_{1}$-type of $S_{t,u}$ tensorfield which is the identity as a transformation in the tangent space of $S_{t,u}$ at each point. These matrices act on
triplets $z=(z^{m}: m=1,2,3)$ of $T^{1}_{1}$-type of $S_{t,u}$ tensorfields on the right:
\begin{align}
 (z\cdot M_{i})^{m}=\sum_{n}z^{n}\cdot(M_{i})^{m}_{n}
\end{align}
so that with matrix multiplication defined according to:
\begin{align}
 (M_{i}\cdot M_{j})^{m}_{n}=\sum_{k}(M_{i})^{k}_{n}\cdot(M_{j})^{m}_{k}
\end{align}
we have:
\begin{align}
 (z\cdot M_{i})\cdot M_{j}=z\cdot(M_{i}\cdot M_{j})
\end{align}
The operators $A_{i}$ are defined in terms of the matrices $M_{i}$ by:
\begin{align}
 A_{i}z=\slashed{\mathcal{L}}_{R_{i}}z+z\cdot M_{i}
\end{align}
We may also define the action of the operators $A_{i}$ on matrices $K$ by:
\begin{align}
 A_{i}K=\slashed{\mathcal{L}}_{R_{i}}K+K\cdot M_{i}
\end{align}
Note that if we also denote by $I$ the identity matrix, then according to (10.70) we have:
\begin{align}
 A_{i}I=M_{i}
\end{align}
From (10.69) and (10.70), we readily deduce the following formula by induction on $k$:
\begin{align}
 A_{i_{k}}...A_{i_{1}}z=\sum_{partitions}(\slashed{\mathcal{L}}_{R})^{s_{1}}z\cdot(M)^{s_{2}}
\end{align}
Here the sum is over all ordered partitions $\{s_{1},s_{2}\}$ of the set $\{1,...,k\}$ into two ordered subsets $s_{1},s_{2}$. Also,
we denote:
\begin{align}
 (M)^{s_{2}}=A_{i_{n_{p}}}...A_{i_{n_{1}}}I\quad\textrm{:if}\quad s_{2}=\{n_{1},...,n_{p}\}, p=|s_{2}|
\end{align}
Introducing the matrices $I_{i}$ and $N_{i}$ by:
\begin{align}
 (I_{i})^{m}_{n}=-\epsilon_{inm}I\\
(N_{i})^{m}_{n}=\delta^{m}_{n}\mu_{i}-\delta^{m}_{i}\mu_{n}
\end{align}
we may write:
\begin{align}
 M_{i}=I_{i}+N_{i}
\end{align}
Let us introduce for any positive integer $k$ and multi-index $(i_{1},...,i_{k})$ the matrices $\leftexp{(k)}{I}_{i_{1}...i_{k}}$ and $\leftexp{(k)}{B}_{i_{1}...i_{k}}$
by:
\begin{align}
 \leftexp{(k)}{I}_{i_{1}...i_{k}}=I_{i_{1}}...I_{i_{k}}
\end{align}
and:
\begin{align}
 \leftexp{(k)}{B}_{i_{1}...i_{k}}=A_{i_{k}}...A_{i_{1}}I-\leftexp{(k)}{I}_{i_{1}...i_{k}}
\end{align}
We also set $\leftexp{(0)}{I}=I$, $\leftexp{(0)}{B}=0$. It is then readily seen that:
\begin{align}
 (\leftexp{(k)}{I}_{i_{1}...i_{k}})^{m}_{n}=\leftexp{(k)}{\alpha}^{m}_{n;i_{1}...i_{k}}I
\end{align}
Obviously, the $\leftexp{(k)}{B}_{i_{1}...i_{k}}$ contain terms of degree $1$ up to $k$ in the $N_{i}$, and the terms of degree $q$ contain up to $k-q$ angular derivatives.
Thus the $\leftexp{(k)}{B}_{i_{1}...i_{k}}$ contain the angular derivatives of the $N_{i}$ of order up to $k-1$. Also, (10.73) can be written as:
\begin{align}
 (M)^{s_{2}}=\leftexp{(p)}{I}_{i_{n_{1}}...i_{n_{p}}}+\leftexp{(p)}{B}_{i_{n_{1}}...i_{n_{p}}}
\end{align}
Finally, by the way in which we derive (10.39), we have:
\begin{align}
 \max_{m,j;i_{1}...i_{k}}|\leftexp{(k)}{\alpha}^{m}_{j;i_{1}...i_{k}}|=\max_{m,j}|\leftexp{(0)}{\alpha}^{m}_{j}|=1
\end{align}
which implies:
\begin{align}
 \max_{n,m;i_{1}...i_{k}}|(\leftexp{(k)}{I}_{i_{1}...i_{k}})^{m}_{n}|=1
\end{align}

\section{Estimates for $y^{i}$}

Given a non-negative integer $l$ let us denote by $\slashed{\textbf{E}}^{Q}_{l}$ the bootstrap assumption that there is a constant $C$ independent of $s$ such that for all
$t\in[0,s]$:
\begin{align*}
 \slashed{\textbf{E}}^{Q}_{l}\quad:\quad \max_{i_{1}...i_{l}}\max_{\alpha}\|R_{i_{l}}...R_{i_{1}}Q\psi_{\alpha}\|_{L^{\infty}(\Sigma_{t}^{\epsilon_{0}})}
\leq C\delta_{0}(1+t)^{-1}
\end{align*}
We then denote by $\slashed{\textbf{E}}^{Q}_{[l]}$ the bootstrap assumption:
\begin{align*}
 \slashed{\textbf{E}}^{Q}_{[l]}\quad:\quad \slashed{\textbf{E}}^{Q}_{0} \quad \textrm{and...and}\quad  \slashed{\textbf{E}}^{Q}_{l}
\end{align*}
Given a positive integer $l$, denote by $\slashed{\textbf{X}}_{l}$ the bootstrap assumption that there is a constant $C$ independent of $s$ such that for all $t\in[0,s]$:
\begin{align*}
 \slashed{\textbf{X}}_{l}\quad:\quad \max_{i_{1}...i_{l}}\|\slashed{\mathcal{L}}_{R_{i_{l}}}...\slashed{\mathcal{L}}_{R_{i_{1}}}\chi\|
_{L^{\infty}(\Sigma_{t}^{\epsilon_{0}})}\leq C\delta_{0}(1+t)^{-2}[1+\log(1+t)]
\end{align*}
Let us also denote by $\slashed{\textbf{X}}_{0}$ the bootstrap assumption:
\begin{align*}
 \slashed{\textbf{X}}_{0}\quad:\quad |\chi-\frac{\slashed{g}}{1-u+t}|\leq C\delta_{0}(1+t)^{-2}[1+\log(1+t)]
\end{align*}
This coincides with $\textbf{F2}$ in chapter 6. We then denote:
\begin{align*}
 \slashed{\textbf{X}}_{[l]}\quad:\quad \slashed{\textbf{X}}_{0} \quad\textrm{and....and}\quad \slashed{\textbf{X}}_{l}
\end{align*}

\subsection{$L^{\infty}$ Estimates for $R_{i_{k}}....R_{i_{1}}y^{j}$}

$\textbf{Proposition 10.1}$ Let hypothesis $\textbf{H0}$ and the estimate (6.177) hold. Let also the bootstrap assumptions $\slashed{\textbf{E}}_{[l]}$, 
$\slashed{\textbf{E}}^{Q}_{[l-1]}$ and $\slashed{\textbf{X}}_{[l-1]}$ hold, for some positive integer $l$. Then if $\delta_{0}$ is suitably small (depending on
$l$) we have:
\begin{align*}
 \max_{j;i_{1}...i_{k}}\|R_{i_{k}}...R_{i_{1}}y^{j}\|_{L^{\infty}(\Sigma_{t}^{\epsilon_{0}})}\leq C_{l}\delta_{0}(1+t)^{-1}[1+\log(1+t)]
\end{align*}
for all $k=0,...,l$.

$Proof$. We shall argue by induction. The proposition reduces for $l=0$ to (6.177), which is assumed here. Let then the proposition hold with $l$ replaced by $l-1$.
We then have:
\begin{align}
 \max_{j;i_{1}...i_{k}}\|R_{i_{k}}...R_{i_{1}}y^{j}\|_{L^{\infty}(\Sigma_{t}^{\epsilon_{0}})}\leq C_{l-1}\delta_{0}(1+t)^{-1}[1+\log(1+t)]
\end{align}
for all $k=0,...,l-1$.

Here, we shall use the expressions in Lemma 10.4. First, note that by Lemma 10.2, 
\begin{align}
 \max_{j;i_{1}...i_{k}}\|\leftexp{(k)}{\delta}^{j}_{i_{1}...i_{k}}\|_{L^{\infty}(\Sigma^{\epsilon_{0}}_{t})}\leq C_{l-1}\delta_{0}[1+\log(1+t)]
\end{align}
 for all $k=0,...,l$.

Then by (10.26) and (10.40) we have:
\begin{align}
 \max_{j;i_{1}...i_{k}}\|R_{i_{k}}...R_{i_{1}}x^{j}\|_{L^{\infty}(\Sigma_{t}^{\epsilon_{0}})}\leq C(1+t)
\end{align}
for all $k=0,...,l$.

So what we need to do is to estimate $\slashed{\mathcal{L}}_{R_{i_{k}}}...\slashed{\mathcal{L}}_{R_{i_{1}}}\slashed{q}^{\prime}_{i}$ for $0\leq k\leq l-1$.
First, by $\slashed{\textbf{E}}_{[l]}$ and $\textbf{H0}$, we have:
\begin{align}
 \max_{i_{1}...i_{k}}\|R_{i_{k}}...R_{i_{1}}(\alpha^{-1})\|_{L^{\infty}(\Sigma_{t}^{\epsilon_{0}})}\leq C\delta_{0}(1+t)^{-1}\\
\max_{i_{1}...i_{k}}\|\slashed{\mathcal{L}}_{R_{i_{k}}}...\slashed{\mathcal{L}}_{R_{i_{1}}}\slashed{k}\|_{L^{\infty}(\Sigma_{t}^{\epsilon_{0}})}\leq C\delta_{0}(1+t)^{-2}
\end{align}
for $1\leq k\leq l-1$. 

Then by (6.62) and $\slashed{\textbf{X}}_{[l-1]}$, we have:
\begin{align}
 \max_{j;i_{1}...i_{k}}\|\slashed{\mathcal{L}}_{R_{i_{k}}}...\slashed{\mathcal{L}}_{R_{i_{1}}}\leftexp{(R_{j})}{\slashed{\pi}}\|_{L^{\infty}(\Sigma_{t}^{\epsilon_{0}})}
\leq C_{l}\delta_{0}(1+t)^{-1}[1+\log(1+t)]
\end{align}
for $0\leq k\leq l-1$.
Here we also used the induction assumption so that we have:
\begin{align}
 \max_{i;i_{1}...i_{k}}\|R_{i_{k}}...R_{i_{1}}\lambda_{i}\|_{L^{\infty}(\Sigma_{t}^{\epsilon_{0}})}\leq C_{l-1}\delta_{0}[1+\log(1+t)]
\end{align}
So by $\slashed{\textbf{X}}_{[l-1]}$, and (10.16), (10.18), we have:
\begin{align}
 \max_{i_{1}...i_{k}}\|\slashed{\mathcal{L}}_{R_{i_{k}}}...\slashed{\mathcal{L}}_{R_{i_{1}}}m^{\prime}\|_{L^{\infty}(\Sigma_{t}^{\epsilon_{0}})}\leq C\delta_{0}
(1+t)^{-2}[1+\log(1+t)]
\end{align}
for $0\leq k\leq l-1$.

Second, by Lemma 10.6, we easily obtain:
\begin{align}
 \|\slashed{\mathcal{L}}_{R_{i_{k}}}...\slashed{\mathcal{L}}_{R_{i_{1}}}R_{j}\|_{L^{\infty}(\Sigma_{t}^{\epsilon_{0}})}\leq C(1+t)
\end{align}
for $0\leq k\leq l$.

In fact, from (10.89) and (10.90) we have:
\begin{align}
 \max_{i;i_{1}...i_{k}}\|\slashed{\mathcal{L}}_{R_{i_{k}}}...\slashed{\mathcal{L}}_{R_{i_{1}}}\mu_{i}\|_{L^{\infty}(\Sigma_{t}^{\epsilon_{0}})}
\leq C_{l}\delta_{0}^{2}(1+t)^{-2}[1+\log(1+t)]^{2}
\end{align}
for all $k=0,...,l-1$.

Then by (10.75):
\begin{align}
 \max_{i;i_{1}...i_{k}}\|\slashed{\mathcal{L}}_{R_{i_{k}}}...\slashed{\mathcal{L}}_{R_{i_{1}}}N_{i}\|_{L^{\infty}(\Sigma_{t}^{\epsilon_{0}})}
\leq C_{l}\delta_{0}^{2}(1+t)^{-2}[1+\log(1+t)]^{2}
\end{align}
for all $k=0,...,l-1$.

By the discussion preceding (10.80),
\begin{align}
 \max_{i_{1}...i_{k}}\|\leftexp{(k)}{B}_{i_{1}...i_{k}}\|_{L^{\infty}(\Sigma_{t}^{\epsilon_{0}})}\leq C_{l}\delta_{0}^{2}(1+t)^{-2}[1+\log(1+t)]^{2}
\end{align}
By (10.80), (10.82) and (10.94) we have:
\begin{align}
 \|(M)^{s_{2}}\|_{L^{\infty}(\Sigma_{t}^{\epsilon_{0}})}\leq C\quad: \textrm{for} |s_{2}|\leq l
\end{align}
(provided that $\delta_{0}$ is suitably small depending on $l$). It then follows from (10.72):
\begin{align}
 \max_{i_{1}...i_{k}}\|A_{i_{k}}...A_{i_{1}}z\|_{L^{\infty}(\Sigma_{t}^{\epsilon_{0}})}\leq C_{k}\sum_{m=0}^{k}\max_{i_{1}...i_{m}}
\|\slashed{\mathcal{L}}_{R_{i_{m}}}...\slashed{\mathcal{L}}_{R_{i_{1}}}z\|_{L^{\infty}(\Sigma_{t}^{\epsilon_{0}})}
\end{align}
for all $k=0,...,l$.

Then by Lemma 10.6 and (10.92):
\begin{align}
 \max_{i_{1}...i_{m+n}}\|\slashed{\mathcal{L}}_{R_{i_{m+n}}}...\slashed{\mathcal{L}}_{R_{i_{m+1}}}\leftexp{(m)}{\rho}_{j;i_{1}...i_{m}}\|
_{L^{\infty}(\Sigma_{t}^{\epsilon_{0}})}\leq C_{l}\delta_{0}^{2}(1+t)^{-2}[1+\log(1+t)]^{2}
\end{align}
for all $m$ and all $n\leq l-1$.

Then by (10.96) we obtain:
\begin{align}
 \|A_{i_{m+n}}...A_{i_{m+1}}\leftexp{(m)}{\rho}_{j;i_{1}...i_{m}}\|_{L^{\infty}(\Sigma_{t}^{\epsilon_{0}})}\leq C_{l}\delta^{2}_{0}(1+t)^{-2}[1+\log(1+t)]^{2}
\end{align}
for all $m$ and all $n\leq l-1$.

So finally we have:
\begin{align}
 \max_{j;i_{1}...i_{m+n}}\|\slashed{\mathcal{L}}_{R_{i_{m+n}}}...\slashed{\mathcal{L}}_{R_{i_{m+1}}}\leftexp{(m)}{\beta}_{j;i_{1}...i_{m}}\|
_{L^{\infty}(\Sigma_{t}^{\epsilon_{0}})}\leq C_{l}\delta_{0}^{2}(1+t)^{-2}[1+\log(1+t)]^{2}
\end{align}
for all $m+n\leq l$.

     We still have to consider $\leftexp{(k)}{\gamma}_{j;i_{1}...i_{k}}$. By the expression in Lemma 10.6 and the induction assumption, we easily deduce:
\begin{align}
 \max_{j;i_{1}...i_{k}}\|\leftexp{(k)}{\gamma}_{j;i_{1}...i_{k}}\|_{L^{\infty}(\Sigma_{t}^{\epsilon_{0}})}\leq C_{l}\delta_{0}^{2}(1+t)^{-1}[1+\log(1+t)]^{2}
\end{align}
for all $k=0,...,l$.

We thus arrive at (10.91). Then from (10.91), (10.90) and (10.15) we have:
\begin{align}
 \max_{i;i_{1}...i_{k}}\|\slashed{\mathcal{L}}_{R_{i_{k}}}...\slashed{\mathcal{L}}_{R_{i_{1}}}\slashed{q}^{\prime}_{i}\|_{L^{\infty}(\Sigma_{t}^{\epsilon_{0}})}
\leq C_{l}\delta_{0}(1+t)^{-1}[1+\log(1+t)]
\end{align}
for all $k=0,...,l-1$.

The proposition then follows. $\qed$ \vspace{7mm}

     There are some immediate corollaries of this proposition. First by Lemma 10.2, we have:

$\textbf{Corollary 10.1.a}$ Under the assumptions of Proposition 10.1 we have:
\begin{align*}
 \max_{i;i_{1}...i_{k}}\|R_{i_{k}}...R_{i_{1}}\lambda_{i}\|_{L^{\infty}(\Sigma_{t}^{\epsilon_{0}})}\leq C_{l}\delta_{0}[1+\log(1+t)]
\end{align*}
and
\begin{align*}
 \max_{j;ii_{1}...i_{k}}\|\leftexp{(k+1)}{\delta}_{ii_{1}...i_{k}}\|_{L^{\infty}(\Sigma_{t}^{\epsilon_{0}})}\leq C_{l}\delta_{0}[1+\log(1+t)]
\end{align*}
for all $k=0,...,l$. 

The second of above implies 
\begin{align*}
 \max_{j;i_{1}...i_{k}}\|R_{i_{k}}...R_{i_{1}}x^{j}\|_{L^{\infty}(\Sigma_{t}^{\epsilon_{0}})}\leq C(1+t)
\end{align*}
for all $k=0,...,l+1$, if $\delta_{0}$ is suitably small (depending on $l$).

     Next, from the proof of Proposition 10.1 we have:

$\textbf{Corollary 10.1.b}$ Under the assumptions of Proposition 10.1 we have:
\begin{align*}
 \max_{j;i_{1}...i_{k}}\|R_{i_{k}}...R_{i_{1}}\hat{T}^{j}\|_{L^{\infty}(\Sigma_{t}^{\epsilon_{0}})}\leq C\\
\textrm{for all} \quad k=0,...,l
\end{align*}
and
\begin{align*}
 \max_{i_{1}...i_{k}}\|R_{i_{k}}...R_{i_{1}}\psi_{\hat{T}}\|_{L^{\infty}(\Sigma_{t}^{\epsilon_{0}})}\leq C_{l}\delta_{0}(1+t)^{-1}\\
\textrm{for all}\quad k=0,...,l
\end{align*}
Moreover,
\begin{align*}
 \max_{i_{1}...i_{k}}\|\slashed{\mathcal{L}}_{R_{i_{k}}}...\slashed{\mathcal{L}}_{R_{i_{1}}}\slashed{\psi}\|_{L^{\infty}(\Sigma_{t}^{\epsilon_{0}})}
\leq C_{l}\delta_{0}(1+t)^{-1}\\
\textrm{for all} \quad k=0,...,l
\end{align*}
and:
\begin{align*}
\max_{i_{1}...i_{k}}\|\slashed{\mathcal{L}}_{R_{i_{k}}}...\slashed{\mathcal{L}}_{R_{i_{1}}}\slashed{w}\|_{L^{\infty}(\Sigma_{t}^{\epsilon_{0}}}\leq C_{l}\delta_{0}(1+t)^{-1}\\
\textrm{for all} \quad k=0,...,l-1
\end{align*}
$\textbf{Corollary 10.1.c}$ Under the assumptions of Proposition 10.1 we have:
\begin{align*}
 \max_{i_{1}...i_{k}}\|\slashed{\mathcal{L}}_{R_{i_{k}}}...\slashed{\mathcal{L}}_{R_{i_{1}}}\slashed{k}\|_{L^{\infty}(\Sigma_{t}^{\epsilon_{0}})}
\leq C_{l}\delta_{0}(1+t)^{-2}\\
\textrm{for all} \quad k=0,...,l-1
\end{align*}
$\textbf{Corollary 10.1.d}$ Under the assumptions of Proposition 10.1 we have:
\begin{align*}
 \max_{i_{1}...i_{k}}\|\slashed{\mathcal{L}}_{R_{i_{k}}}...\slashed{\mathcal{L}}_{R_{i_{1}}}\leftexp{(R_{i})}{\slashed{\pi}}\|_{L^{\infty}(\Sigma_{t}^{\epsilon_{0}})}
\leq C_{l}\delta_{0}(1+t)^{-1}[1+\log(1+t)]\\
\textrm{for all} \quad k=0,...,l-1
\end{align*}
$\textbf{Corollary 10.1.e}$ Under the assumptions of Proposition 10.1 we have, if $\delta_{0}$ is suitably small (depending on $l$),
\begin{align*}
 \max_{j;i_{1}...i_{k}}\|\slashed{\mathcal{L}}_{R_{i_{k}}}...\slashed{\mathcal{L}}_{R_{i_{1}}}R_{j}\|_{L^{\infty}(\Sigma_{t}^{\epsilon_{0}})}
\leq C_{l}\delta_{0}(1+t)\\
\textrm{for all} \quad k=0,...,l
\end{align*}
More precisely, we have:
\begin{align*}
 \max_{m,j;i_{1}...i_{k}}|\leftexp{(k)}{\alpha}^{m}_{j;i_{1}...i_{k}}|=1\quad:\textrm{for all}\quad k
\end{align*}
\begin{align*}
 \max_{m,j;i_{1}...i_{k}}\|\leftexp{(k)}{\beta}^{m}_{j;i_{1}...i_{k}}\|_{L^{\infty}(\Sigma_{t}^{\epsilon_{0}})}\leq C_{l}\delta_{0}^{2}
(1+t)^{-2}[1+\log(1+t)]^{2}\\
\textrm{for all} \quad k=0,...,l\\
(\leftexp{(0)}{\beta}^{m}_{j}=0)
\end{align*}
\begin{align*}
 \max_{j;i_{1}...i_{k}}\|\leftexp{(k)}{\gamma}_{j;i_{1}...i_{k}}\|_{L^{\infty}(\Sigma_{t}^{\epsilon_{0}})}\leq C_{l}\delta_{0}^{2}
(1+t)^{-1}[1+\log(1+t)]^{2}\\
\textrm{for all} \quad k=0,...,l\\
(\leftexp{(0)}{\gamma}_{j}=0)
\end{align*}
Finally, we have:

$\textbf{Corollary 10.1.f}$ Under the assumptions of Proposition 10.1, we have:
\begin{align*}
 \max_{i;i_{1}...i_{k}}\|\leftexp{(k)}{\slashed{q}}^{\prime}_{ii_{1}...i_{k}}\|_{L^{\infty}(\Sigma_{t}^{\epsilon_{0}})}\leq C_{l}
(1+t)^{-1}[1+\log(1+t)]\\
\textrm{for all} \quad k=0,...,l-1
\end{align*}
and
\begin{align*}
 \max_{i;i_{1}...i_{k}}\|\leftexp{(k)}{\slashed{r}}^{j}_{ii_{1}...i_{k}}\|_{L^{\infty}(\Sigma_{t}^{\epsilon_{0}})}\leq C_{l}
(1+t)^{-1}[1+\log(1+t)]\\
\textrm{for all} \quad k=0,...,l-1
\end{align*}

\subsection{$L^{2}$ Estimates for $R_{i_{k}}...R_{i_{1}}y^{j}$}
     We define:
\begin{align*}
 \mathcal{A}_{0}=\|\chi-\frac{\slashed{g}}{1-u+t}\|_{L^{2}(\Sigma_{t}^{\epsilon_{0}})}
\end{align*}
and for $k\slashed{=}0$:
\begin{align*}
 \mathcal{A}_{k}=\max_{i_{1}...i_{k}}\|\slashed{\mathcal{L}}_{R_{i_{k}}}...\slashed{\mathcal{L}}_{R_{i_{1}}}\chi\|_{L^{2}(\Sigma_{t}^{\epsilon_{0}})}
\end{align*}
and:
\begin{align*}
 \mathcal{A}_{[l]}=\sum_{k=0}^{l}\mathcal{A}_{k}
\end{align*}
Also, we define:
\begin{align*}
 \mathcal{W}^{Q}_{k}=\max_{\alpha;i_{1}...i_{k}}\|R_{i_{k}}...R_{i_{1}}Q\psi_{\alpha}\|_{L^{2}(\Sigma_{t}^{\epsilon_{0}})},\quad 
\mathcal{W}^{Q}_{[l]}=\sum_{k=0}^{l}\mathcal{W}^{Q}_{k}
\end{align*}

$\textbf{Proposition 10.2}$ Let the hypothesis $\textbf{H0}$ and the estimate (6.177) hold. Let $l$ be a positive integer and let the bootstrap assumptions 
$\slashed{\textbf{E}}_{[l_{*}]}$, $\slashed{\textbf{E}}^{Q}_{[l_{*}-1]}$ and $\slashed{\textbf{X}}_{[l_{*}-1]}$ hold. Then if $\delta_{0}$ is suitably small 
(depending on $l$) we have:
\begin{align*}
 \mathcal{Y}_{k}\leq C_{l}\{\mathcal{Y}_{0}+(1+t)\mathcal{A}_{[l-1]}+\mathcal{W}_{[l]}\}
\end{align*}
for all $k=0,...,l$.

$Proof$. The proof is again by induction. The proposition extends trivially to the case $l=0$. Let then the proposition hold with $l$ replaced by $1,...,l-1$ and 
consider the case $l$. The assumptions of Proposition 10.2 coincide with those of Proposition 10.1 with $l$ replaced by $l_{*}$. Therefore the conclusion of 
Proposition 10.1 holds with $l$ replaced by $l_{*}$. Therefore, we have the $L^{\infty}$-bounds for $\slashed{\mathcal{L}}_{R_{i_{k}}}...\slashed{\mathcal{L}}_{R_{i_{1}}}
\slashed{q}^{\prime}_{i}$ where $k=0,...,l_{*}-1$ and for $\slashed{d}((R)^{s_{2}}x^{j})$ where $|s_{2}|\leq l_{*}$. Also, all the corollaries of Proposition 10.1 hold 
with $l$ replaced by $l_{*}$. As before, we shall use the expression in Lemma 10.4 and we must get a $L^{2}$ estimate for $\slashed{\mathcal{L}}_{R_{i_{k}}}...
\slashed{\mathcal{L}}_{R_{i_{1}}}\slashed{q}^{\prime}_{i}$ where $k=0,...,l-1$ and for $\slashed{d}((R)^{s_{2}}x^{j})$ where $|s_{2}|\leq l-1$. By (10.18), we have:
\begin{align}
 \|(\slashed{\mathcal{L}}_{R})^{s}m^{\prime}_{b}\|_{L^{2}(\Sigma_{t}^{\epsilon_{0}})}\leq C(\delta_{0}(1+t)^{-1}\mathcal{W}_{[l]}+\mathcal{A}_{[l-1]})+C_{l}(1+t)^{-1}
\mathcal{Y}_{[l-1]}
\end{align}
for all $|s|\leq l-1$. Here we have used Lemma 10.3, which yields:
\begin{align}
 \max_{i;i_{1}...i_{k}}\|R_{i_{1}}...R_{i_{k}}\lambda_{i}\|_{L^{2}(\Sigma_{t}^{\epsilon_{0}})}\leq C_{l-1}(1+t)\mathcal{Y}_{[l-1]}
\end{align}
for all $k=0,...,l-1$. 

This is because the hypothesis of Lemma 10.3 with $l$ replaced by $l-1$ are a fortiori satisfied.

Also we have:
\begin{align}
 \|(\slashed{\mathcal{L}}_{R})^{s}\leftexp{(R)}{\slashed{\pi}}\|_{L^{2}(\Sigma_{t}^{\epsilon_{0}})}\leq C_{l}\mathcal{Y}_{[l-2]}+
C_{l}\delta_{0}(1+t)^{-1}[1+\log(1+t)]\mathcal{W}_{[l-1]}+C_{l}\delta_{0}[1+\log(1+t)]\mathcal{A}_{[l-2]}
\end{align}
 for all $s=0,...,l-2$. Here we also used the induction hypothesis for $\lambda_{i}$.

Finally, we have to deal with the case that more than half of the derivatives hit $R_{i}$ in the expression (10.15). Here we must use Lemma 10.6. The first term on the 
right of this expression is obviously bounded in $L^{\infty}$, as well as the second factor in the second term. So first we consider the term $\leftexp{(k)}{\beta}^{m}
_{j;i_{1}...i_{k}}$ for all $k=0,...,l-1$. Let us now consider expression (10.72) for an arbitrary triplet $z$ of $T^{1}_{1}$-type $S_{t,u}$ tensorfields and $k=0,...,l-1$.
We have two cases to consider in regard to the terms in the sum in (10.72).

     Case 1:$|s_{2}|\leq l_{*}$  and:   Case 2: $|s_{1}|\leq l_{*}-1$

     In Case 1 by (10.95) with $l_{*}$ in the role of $l$, we have:
\begin{align}
 \|(\slashed{\mathcal{L}}_{R})^{s_{1}}z\cdot(M)^{s_{2}}\|_{L^{2}(\Sigma_{t}^{\epsilon_{0}})}\leq \|(\slashed{\mathcal{L}}_{R})^{s_{1}}z\|
_{L^{2}(\Sigma_{t}^{\epsilon_{0}})}\|(M)^{s_{2}}\|_{L^{\infty}(\Sigma_{t}^{\epsilon_{0}})}\\\notag
\leq C\sum_{m=0}^{l-1}\max_{i_{1}...i_{m}}\|\slashed{\mathcal{L}}_{R_{i_{m}}}...\slashed{\mathcal{L}}_{R_{i_{1}}}z\|_{L^{2}(\Sigma_{t}^{\epsilon_{0}})}
\end{align}
In Case 2 we express $(M)^{s_{2}}$ as in (10.80). We estimate the contribution of the first term in (10.80) by placing this term in $L^{\infty}$ using (10.82):
\begin{align}
 \|(\slashed{\mathcal{L}}_{R})^{s_{1}}z\cdot\leftexp{(p)}{I}_{i_{n_{1}}...i_{n_{p}}}\|_{L^{2}(\Sigma_{t}^{\epsilon_{0}})}
\leq \|(\slashed{\mathcal{L}}_{R})^{s_{1}}z\|_{L^{2}(\Sigma_{t}^{\epsilon_{0}})}\|\leftexp{(p)}{I}_{i_{n_{1}}...i_{n_{p}}}\|_{L^{\infty}(\Sigma_{t}^{\epsilon_{0}})}\\\notag
\leq C\sum_{m=0}^{l-1}\max_{i_{1}...i_{m}}\|\slashed{\mathcal{L}}_{R_{i_{m}}}...\slashed{\mathcal{L}}_{R_{i_{1}}}z\|_{L^{2}(\Sigma_{t}^{\epsilon_{0}})}
\end{align}
Concerning the contribution of the second term in (10.80), we need a $L^{2}$ estimate for this term. By the discussion above (10.80), we need an estimate for
\begin{align*}
 \max_{i;i_{1}...i_{k}}\|\slashed{\mathcal{L}}_{R_{i_{k}}}...\slashed{\mathcal{L}}_{R_{i_{1}}}N_{i}\|_{L^{2}(\Sigma_{t}^{\epsilon_{0}})}
\end{align*}
for all $k=0,...,l-2$.

Then by (10.75), we need an estimate for
\begin{align*}
 \max_{i;i_{1}...i_{k}}\|\slashed{\mathcal{L}}_{R_{i_{k}}}...\slashed{\mathcal{L}}_{R_{i_{1}}}\mu_{i}\|_{L^{2}(\Sigma_{t}^{\epsilon_{0}})}
\end{align*}
for all $k=0,...,l-2$.

To obtain this estimate, we use the definition (10.61). By (10.103) (with $l$ replaced by $l-1$), we have a $L^{2}$ estimate for $\lambda_{i}$. For a $L^{2}$ estimate 
of $m^{\prime}$, we use (10.104) and (10.102) as well as (10.88) and (10.90), to obtain:
\begin{align}
 \max_{i_{1}...i_{k}}\|\slashed{\mathcal{L}}_{R_{i_{k}}}...\slashed{\mathcal{L}}_{R_{i_{1}}}m^{\prime}\|_{L^{2}(\Sigma_{t}^{\epsilon_{0}})}
\leq C(\delta_{0}(1+t)^{-1}\mathcal{W}_{[l-1]}+\mathcal{A}_{[l-2]})+C_{l}(1+t)^{-1}\mathcal{Y}_{[l-2]}
\end{align}
for $k=0,...,l-2$.

Combining this and (10.103) with (10.90), (10.89) (with $l$ replaced by $l_{*}$), we get:
\begin{align}
 \max_{i;i_{1}...i_{k}}\|\slashed{\mathcal{L}}_{R_{i_{k}}}...\slashed{\mathcal{L}}_{R_{i_{1}}}\mu_{i}\|_{L^{2}(\Sigma_{t}^{\epsilon_{0}})}\\\notag
\leq C_{l}\delta_{0}[1+\log(1+t)](\delta_{0}(1+t)^{-1}\mathcal{W}_{[l-1]}+\mathcal{A}_{[l-2]})+C_{l}(1+t)^{-1}[1+\log(1+t)]\mathcal{Y}_{[l-2]}
\end{align}
for $k=0,...,l-2$.

Therefore we have:
\begin{align}
 \max_{i;i_{1}...i_{k}}\|\slashed{\mathcal{L}}_{R_{i_{k}}}...\slashed{\mathcal{L}}_{R_{i_{1}}}N_{i}\|_{L^{2}(\Sigma_{t}^{\epsilon_{0}})}\\\notag
\leq C_{l}\delta_{0}[1+\log(1+t)](\delta_{0}(1+t)^{-1}\mathcal{W}_{[l-1]}+\mathcal{A}_{[l-2]})+C_{l}(1+t)^{-1}[1+\log(1+t)]\mathcal{Y}_{[l-2]}
\end{align}
So the contribution of the second term in (10.80) in Case 2 is bounded by
\begin{align}
 \|(\slashed{\mathcal{L}}_{R})^{s_{1}}z\cdot\leftexp{(p)}{B}_{i_{n_{1}}...i_{n_{p}}}\|_{L^{2}(\Sigma_{t}^{\epsilon_{0}})}\leq 
\|(\slashed{\mathcal{L}}_{R})^{s_{1}}z\|_{L^{\infty}(\Sigma_{t}^{\epsilon_{0}})}\|\leftexp{(p)}{B}_{i_{n_{1}}...i_{n_{p}}}\|_{L^{2}(\Sigma_{t}^{\epsilon_{0}})}\\\notag
\leq C_{l}\delta_{0}(1+t)^{-1}[1+\log(1+t)]\max_{m\leq l_{*}-1}\max_{i_{1}...i_{m}}\|\slashed{\mathcal{L}}_{R_{i_{m}}}...\slashed{\mathcal{L}}_{R_{i_{1}}}z\|
_{L^{\infty}(\Sigma_{t}^{\epsilon_{0}})}\cdot\{\mathcal{Y}_{[l-2]}+(1+t)\mathcal{A}_{[l-2]}
+\mathcal{W}_{[l-1]}\}
\end{align}
We thus conclude that for any triplet $z$ of $T^{1}_{1}$-type of $S_{t,u}$ tensorfields:
\begin{align}
 \max_{i_{1}...i_{k}}\|A_{i_{k}}...A_{i_{1}}z\|_{L^{2}(\Sigma_{t}^{\epsilon_{0}})}
\leq C_{l}\sum_{m=0}^{l-1}\max_{i_{1}...i_{m}}\|\slashed{\mathcal{L}}_{R_{i_{m}}}...\slashed{\mathcal{L}}_{R_{i_{1}}}z\|_{L^{2}(\Sigma_{t}^{\epsilon_{0}})}\\\notag
+C_{l}\delta_{0}(1+t)^{-1}[1+\log(1+t)]\max_{m\leq l_{*}-1}\max_{i_{1}...i_{m}}\|\slashed{\mathcal{L}}_{R_{i_{m}}}...\slashed{\mathcal{L}}_{R_{i_{1}}}z\|
_{L^{\infty}(\Sigma_{t}^{\epsilon_{0}})}\cdot\{\mathcal{Y}_{[l-2]}+(1+t)\mathcal{A}_{[l-2]}
+\mathcal{W}_{[l-1]}\}
\end{align}
for all $k=0,...,l-1$.

From the expression for the triplets $\leftexp{(k)}{\rho}_{j;i_{1}...i_{k}}$ of Lemma 10.6 and the bounds (10.81) and (10.108) with $l-1$ replaced by $l$, we obtain:
\begin{align}
 \max_{i_{1}...i_{m+n}}\|\slashed{\mathcal{L}}_{R_{i_{m+n}}}...\slashed{\mathcal{L}}_{R_{i_{m+1}}}\leftexp{(m)}{\rho}_{j;i_{1}...i_{m}}\|_
{L^{2}(\Sigma_{t}^{\epsilon_{0}})}\\\notag
\leq C_{l}\delta_{0}(1+t)^{-1}[1+\log(1+t)]\{\mathcal{Y}_{[l-1]}+(1+t)\mathcal{A}_{[l-1]}+\mathcal{W}_{[l]}\}
\end{align}
for all $m$ and all $n\leq l-1$.

We then apply (10.111) to the triplet $\leftexp{(m)}{\rho}_{j;i_{1}...i_{m}}$ using $L^{2}$ estimate (10.112) as well as $L^{\infty}$ bound (10.97) with $l_{*}$ in 
the role of $l$. This yields:
\begin{align}
 \|A_{i_{m+n}}...A_{i_{m+1}}\leftexp{(m)}{\rho}_{j;i_{1}...i_{m}}\|_{L^{2}(\Sigma_{t}^{\epsilon_{0}})}\\\notag
\leq C_{l}\delta_{0}(1+t)^{-1}[1+\log(1+t)]\{\mathcal{Y}_{[l-1]}+(1+t)\mathcal{A}_{[l-1]}+\mathcal{W}_{[l]}\}
\end{align}
 for all $m$ and all $n\leq l-1$.

It then follows that:
\begin{align}
 \max_{j;i_{1}...i_{k}}\|\leftexp{(k)}{\beta}_{j;i_{1}...i_{k}}\|_{L^{2}(\Sigma_{t}^{\epsilon_{0}})}\\\notag
\leq C_{l}\delta_{0}(1+t)^{-1}[1+\log(1+t)]\{\mathcal{Y}_{[l-1]}+(1+t)\mathcal{A}_{[l-1]}+\mathcal{W}_{[l]}\}
\end{align}
for all $k=0,...,l$.

More generally, we have:
\begin{align}
 \max_{j;i_{1}...i_{m+n}}\|\slashed{\mathcal{L}}_{R_{i_{m+n}}}...\slashed{\mathcal{L}}_{R_{i_{m+1}}}\leftexp{(m)}{\beta}_{j;i_{1}...i_{m}}\|
_{L^{2}(\Sigma_{t}^{\epsilon_{0}})}\\\notag
\leq C_{l}\delta_{0}(1+t)^{-1}[1+\log(1+t)]\{\mathcal{Y}_{[l-1]}+(1+t)\mathcal{A}_{[l-1]}+\mathcal{W}_{[l]}\}
\end{align}
for all $m$ and all $n$ such that $m+n\leq l$.

We turn to consider the expression for the coefficients $\leftexp{(k)}{\gamma}_{j;i_{1}...i_{k}}$ given in the statement of Lemma 10.6. Since we already have (10.115), 
what we need is an estimate for $\|\slashed{\mathcal{L}}_{R_{i_{k}}}...\slashed{\mathcal{L}}_{R_{i_{1}}}\nu_{ij}\|_{L^{2}(\Sigma_{t}^{\epsilon_{0}})}$ for all 
$k=0,...,l-2$. In view of (10.62), what is needed is an estimate for $\|\slashed{\mathcal{L}}_{R_{i_{k}}}...\slashed{\mathcal{L}}_{R_{i_{1}}}w_{i}\|_{L^{2}(\Sigma_{t}^{\epsilon_{0}})}$,
therefore an estimate for $\|\slashed{\mathcal{L}}_{R_{i_{k}}}...\slashed{\mathcal{L}}_{R_{i_{1}}}(w_{i})_{b}\|_{L^{2}(\Sigma_{t}^{\epsilon_{0}})}$，for all $k=0,...,l-2$.
Concerning (10.42), the $L^{2}$ estimate for derivatives of $y^{i}$ is just the induction hypothesis. While since we have:
\begin{align}
 R_{q}(R)^{s_{2}}x^{j}=\leftexp{(p+1)}{x}^{j}_{i_{n_{1}}...i_{n_{p}}q}-\leftexp{(p+1)}{\delta}^{j}_{i_{n_{1}}...i_{n_{p}}q}
\end{align}
where $s_{2}=\{n_{1},...,n_{p}\},\quad p=|s_{2}|$. So by $\textbf{H0}$,
\begin{align}
 |\slashed{d}(R)^{s_{2}}x^{j}|\leq C(1+t)^{-1}\sum_{q=1}^{3}\{|\leftexp{(p+1)}{x}^{j}_{i_{n_{1}}...i_{n_{p}}q}|+|\leftexp{(p+1)}{\delta}^{j}_{i_{n_{1}}...i_{n_{p}}q}|\}
\end{align}
By (10.40), we can estimate the first term on the right of (10.117) in $L^{\infty}$ norm. By virtue of Proposition 10.1 with $l$ replaced by $l_{*}$, the hypothesis of 
Lemma 10.3 is satisfied, and of course, the hypothesis of Lemma 10.3 with $l$ replaced by $l-1$ is a fortiori satisfied. Therefore we have:
\begin{align}
 \max_{j;i_{1}...i_{k}}\|\leftexp{(k+1)}{\delta}^{j}_{ii_{1}...i_{k}}\|_{L^{2}(\Sigma_{t}^{\epsilon_{0}})}\leq C_{l-1}(1+t)\mathcal{Y}_{[l-1]}
\end{align}
for all $k=0,...,l-1$. So we can estimate the second term in $L^{2}$ norm. Combining (10.118), (10.85) with induction hypothesis and the conclusion of Proposition 10.1, we obtain:
\begin{align}
 \max_{i;i_{1}...i_{k}}\|\slashed{\mathcal{L}}_{R_{i_{k}}}...\slashed{\mathcal{L}}_{R_{i_{1}}}(w_{i})_{b}\|_{L^{2}(\Sigma_{t}^{\epsilon_{0}})}
\leq C_{l}\mathcal{Y}_{[l-1]}
\end{align}
for all $k=0,...,l-1$. 
Then combining this with (10.104), (10.88) and $L^{\infty}$ bounds for $\slashed{\mathcal{L}}_{R_{i_{k}}}...
\slashed{\mathcal{L}}_{R_{i_{1}}}(w_{i})_{b}$, which is the same as the conclusion of Proposition 10.1, we obtain:
\begin{align}
 \max_{i;i_{1}...i_{k}}\|\slashed{\mathcal{L}}_{R_{i_{k}}}...\slashed{\mathcal{L}}_{R_{i_{1}}}w_{i}\|_{L^{2}(\Sigma_{t}^{\epsilon_{0}})}
\leq C_{l}\{\mathcal{Y}_{[l-1]}+\delta_{0}(1+t)^{-2}[1+\log(1+t)]^{2}[(1+t)\mathcal{A}_{[l-2]}+\mathcal{W}_{[l-1]}]\}
\end{align}
for all $k=0,...,l-1$. Then we can easily obtain by (10.62):
\begin{align}
 \max_{i,j;i_{1}...i_{k}}\|\slashed{\mathcal{L}}_{R_{i_{k}}}...\slashed{\mathcal{L}}_{R_{i_{1}}}\nu_{ij}\|_{L^{2}(\Sigma_{t}^{\epsilon_{0}})}\\\notag
\leq C_{l}\delta_{0}[1+\log(1+t)]\{\mathcal{Y}_{[l-1]}+\delta_{0}(1+t)^{-2}[1+\log(1+t)]^{2}[(1+t)\mathcal{A}_{[l-2]}+\mathcal{W}_{[l-1]}]\}
\end{align}
Combining this with (10.114), we obtain that the last two terms on the right in the expression of Lemma 10.6 are bounded in $L^{2}$ by:
\begin{align}
 C_{l}\delta_{0}[1+\log(1+t)]\{\mathcal{Y}_{[l-1]}+(1+t)\mathcal{A}_{[l-2]}+\mathcal{W}_{[l-1]}\}
\end{align}
Recall that the first term is estimated in $L^{\infty}$ by $1+t$. Combining (10.122) with (10.107), Corollary 10.1.e and (10.90), we obtain:
\begin{align}
 \|\slashed{\mathcal{L}}_{R_{i_{k}}}...\slashed{\mathcal{L}}_{R_{i_{1}}}\slashed{q}^{\prime}_{i}\|_{L^{2}(\Sigma_{t}^{\epsilon_{0}})}
\leq C_{l}\mathcal{Y}_{[l-1]}+C_{l}\delta_{0}[(1+t)\mathcal{A}_{[l-1]}+\mathcal{W}_{[l]}]
\end{align}
for all $k=0,...,l-1$.

Next, we must consider the $L^{2}$ estimate for $\slashed{d}((R)^{s_{2}}x^{j})$ for $|s_{2}|\leq l-1$. Recalling the discussion from (10.116) to (10.118),
the $L^{2}$ contribution is bounded by (10.118). Combining this with Corollary 10.1.f and (10.85) as well as $\textbf{H0}$, we obtain:
\begin{align}
 \|\leftexp{(k)}{\slashed{r}}^{j}_{ii_{1}...i_{k}}\|_{L^{2}(\Sigma_{t}^{\epsilon_{0}})}\leq C_{l}\mathcal{Y}_{[l-1]}+C_{l}\delta_{0}[(1+t)\mathcal{A}_{[l-2]}
+\mathcal{W}_{[l-1]}]
\end{align}
for all $k=0,...,l-1$.

Thus we finally obtain:
\begin{align}
\|R_{i_{k}}...R_{i_{1}}y^{j}\|_{L^{2}(\Sigma_{t}^{\epsilon_{0}})}\leq C_{l}\mathcal{Y}_{[l-1]}+C_{l}\delta_{0}[(1+t)\mathcal{A}_{[l-1]}+\mathcal{W}_{[l]}]
\end{align}
for all $k=0,...,l$.

Recall the induction hypothesis:
\begin{align}
 \|R_{i_{k}}...R_{i_{1}}y^{j}\|_{L^{2}(\Sigma_{t}^{\epsilon_{0}})}\leq C_{l}\{\mathcal{Y}_{0}+(1+t)\mathcal{A}_{[l-2]}+\mathcal{W}_{[l-1]}\}
\end{align}
 for all $k=0,...,l-1$.

Substituting this in (10.125), the proposition is proved. $\qed$ \vspace{7mm}

As before, we have a lot of Corollaries. First, in conjunction with Lemma 10.3 we have:

$\textbf{Corollary 10.2.a}$ Under the assumptions of Proposition 10.2 we have:
\begin{align}
 \max_{i;i_{1}...i_{k}}\|R_{i_{k}}...R_{i_{1}}\lambda_{i}\|_{L^{2}(\Sigma_{t}^{\epsilon_{0}})}
\leq C_{l}(1+t)\{\mathcal{Y}_{0}+(1+t)\mathcal{A}_{[l-1]}+\mathcal{W}_{[l-1]}\}
\end{align}
for all $k=0,...,l$.

and:
\begin{align}
 \max_{j;ii_{1}...i_{k}}\|\leftexp{(k+1)}{\delta}^{j}_{ii_{1}...i_{k}}\|_{L^{2}(\Sigma_{t}^{\epsilon_{0}})}
\leq C_{l}(1+t)\{\mathcal{Y}_{0}+(1+t)\mathcal{A}_{[l-1]}+\mathcal{W}_{[l-1]}\}
\end{align}
for all $k=0,...,l$.

Next, it is not difficult to obtain from the proof of Proposition 10.2 the following:

$\textbf{Corollary 10.2.b}$ Under the assumptions of Proposition 10.2 we have:
\begin{align*}
 \max_{i_{1}...i_{k}}\|R_{i_{k}}...R_{i_{1}}\psi_{\hat{T}}\|_{L^{2}(\Sigma_{t}^{\epsilon_{0}})}
\leq C_{l}\{\mathcal{W}_{[l]}+\delta_{0}(1+t)^{-1}[\mathcal{Y}_{0}+(1+t)\mathcal{A}_{[l-1]}]\}
\end{align*}
for all $k=0,...,l$.

\begin{align*}
 \max_{i_{1}...i_{k}}\|\slashed{\mathcal{L}}_{R_{i_{k}}}...\slashed{\mathcal{L}}_{R_{i_{1}}}\slashed{\psi}\|_{L^{2}(\Sigma_{t}^{\epsilon_{0}})}
\leq C_{l}\{\mathcal{W}_{[l]}+\delta_{0}(1+t)^{-1}[\mathcal{Y}_{0}+(1+t)\mathcal{A}_{[l-1]}]\}
\end{align*}
for all $k=0,...,l$.

and:
\begin{align*}
 \max_{i_{1}...i_{k}}\|\slashed{\mathcal{L}}_{R_{i_{k}}}...\slashed{\mathcal{L}}_{R_{i_{1}}}\slashed{\omega}\|_{L^{2}(\Sigma_{t}^{\epsilon_{0}})}
\leq C_{l}(1+t)^{-1}\{\mathcal{W}_{[l]}+\delta_{0}(1+t)^{-1}[\mathcal{Y}_{0}+(1+t)\mathcal{A}_{[l-2]}]\}
\end{align*}
for all $k=0,...,l-1$.
 
$\textbf{Corollary 10.2.c}$ Under the assumptions of Proposition 10.2 we have:
\begin{align*}
 \max_{i_{1}...i_{k}}\|\slashed{\mathcal{L}}_{R_{i_{k}}}...\slashed{\mathcal{L}}_{R_{i_{1}}}\slashed{k}\|_{L^{2}(\Sigma_{t}^{\epsilon_{0}})}
\leq C_{l}(1+t)^{-1}\{\mathcal{W}_{[l]}+\delta_{0}(1+t)^{-1}[\mathcal{Y}_{0}+(1+t)\mathcal{A}_{[l-2]}]\}
\end{align*}
for all $k=0,...,l-1$.

$\textbf{Corollary 10.2.d}$ Under the assumption of Proposition 10.2 we have:
\begin{align*}
 \max_{i;i_{1}...i_{k}}\|\slashed{\mathcal{L}}_{R_{i_{k}}}...\slashed{\mathcal{L}}_{R_{i_{1}}}\leftexp{(R_{i})}{\slashed{\pi}}\|_{L^{2}(\Sigma_{t}^{\epsilon_{0}})}\\
\leq C_{l}\{\delta_{0}(1+t)^{-1}[1+\log(1+t)][(1+t)\mathcal{A}_{[l-1]}+\mathcal{W}_{[l]}]+\mathcal{Y}_{0}+(1+t)\mathcal{A}_{[l-2]}+\mathcal{W}_{[l-1]}\}
\end{align*}
for all $k=0,...,l-1$.

$\textbf{Corollary 10.2.e}$ Under the assumptions of Proposition 10.2, the coefficients of the expression for $\slashed{\mathcal{L}}_{R_{i_{k}}}...
\slashed{\mathcal{L}}_{R_{i_{1}}}R_{j}$ of Lemma 10.6 satisfy:
\begin{align*}
 \max_{j;i_{1}...i_{k}}\|\leftexp{(k)}{\beta}_{j;i_{1}...i_{k}}\|_{L^{2}(\Sigma_{t}^{\epsilon_{0}})}
\leq C_{l}\delta_{0}(1+t)^{-1}[1+\log(1+t)]\{\mathcal{Y}_{0}+(1+t)\mathcal{A}_{[l-1]}+\mathcal{W}_{[l]}\}
\end{align*}
for all $k=1,...,l$.

\begin{align*}
 \max_{j;i_{1}...i_{k}}\|\leftexp{(k)}{\gamma}_{j;i_{1}...i_{k}}\|_{L^{2}(\Sigma_{t}^{\epsilon_{0}})}
\leq C_{l}\delta_{0}[1+\log(1+t)]\{\mathcal{Y}_{0}+(1+t)\mathcal{A}_{[l-2]}+\mathcal{W}_{[l-1]}\}
\end{align*}
for all $k=1,...,l$.

Finally we have:

$\textbf{Corollary 10.2.f}$ Under the assumptions of Proposition 10.2, the coefficients of the expression of $R_{i_{k}}...R_{i_{1}}R_{j}y^{j}$ of Lemma 10.4 satisfy:
\begin{align*}
 \max_{ii_{1}...i_{k}}\|\leftexp{(k)}{\slashed{q}}^{\prime}_{ii_{1}...i_{k}}\|_{L^{2}(\Sigma_{t}^{\epsilon_{0}})}
\leq C_{l}\{\mathcal{Y}_{0}+(1+t)\mathcal{A}_{[l-1]}+\mathcal{W}_{[l]}\}
\end{align*}
for all $k=0,...,l-1$.

and:
\begin{align*}
 \max_{ii_{1}...i_{k}}\|\leftexp{(k)}{\slashed{r}}^{j}_{ii_{1}...i_{k}}\|_{L^{2}(\Sigma_{t}^{\epsilon_{0}})}
\leq C_{l}\{\mathcal{Y}_{0}+(1+t)\mathcal{A}_{[l-2]}+\mathcal{W}_{[l-1]}\}
\end{align*}
for all $k=0,...,l-1$. \vspace{7mm}

      We shall now estimate the angular derivatives of the functions $L^{i}$. After this we shall proceed with the estimates for the $S_{t,u}$ 1-form $\kappa^{-1}\zeta$
 and the $S_{t,u}$-tangential vectorfields $\leftexp{(R_{i})}{Z}$. According to the equations (6.64) and (6.65):
\begin{align}
 L^{j}=\frac{x^{j}}{1-u+t}+z^{j}
\end{align}
where:
\begin{align}
 z^{j}=-\alpha y^{j}+\frac{(\alpha-1)x^{j}}{1-u+t}-\psi^{j}
\end{align}
Let us define:
\begin{align}
 \slashed{\omega}_{L}=L^{\mu}\slashed{d}\psi_{\mu}
\end{align}
and
\begin{align}
 \omega_{L\hat{T}}=\hat{T}^{i}(L\psi_{i})
\end{align}
As a corollary of Proposition 10.1 we have, in connection with the first two statements of Corollary 10.1.b, the following.

$\textbf{Corollary 10.1.g}$ Under the assumptions of Proposition 10.1 we have:
\begin{align*}
 \max_{j;i_{1}...i_{k}}\|R_{i_{k}}...R_{i_{1}}z^{j}\|_{L^{\infty}(\Sigma_{t}^{\epsilon_{0}})}
\leq C\delta_{0}(1+t)^{-1}[1+\log(1+t)]
\end{align*}
for all $k=0,...,l$.
\begin{align*}
 \max_{j;i_{1}...i_{k}}\|R_{i_{k}}...R_{i_{1}}L^{j}\|_{L^{\infty}(\Sigma_{t}^{\epsilon_{0}})}\leq C
\end{align*}
for all $k=0,...,l$.
Moreover,
\begin{align*}
 \max_{i_{1}...i_{k}}\|R_{i_{k}}...R_{i_{1}}\psi_{L}\|_{L^{\infty}(\Sigma_{t}^{\epsilon_{0}})}\leq C_{l}\delta_{0}(1+t)^{-1}
\end{align*}
for all $k=1,...,l$.
And:
\begin{align*}
 \max_{i_{1}...i_{k}}\|\slashed{\mathcal{L}}_{R_{i_{k}}}...\slashed{\mathcal{L}}_{R_{i_{1}}}\slashed{\omega}_{L}\|_{L^{\infty}(\Sigma_{t}^{\epsilon_{0}})}
\leq C_{l}\delta_{0}(1+t)^{-2}
\end{align*}
for all $k=0,...,l-1$.
\begin{align*}
 \max_{i_{1}...i_{k}}\|\slashed{\mathcal{L}}_{R_{i_{k}}}...\slashed{\mathcal{L}}_{R_{i_{1}}}\omega_{L\hat{T}}\|_{L^{\infty}(\Sigma_{t}^{\epsilon_{0}})}
\leq C_{l}\delta_{0}(1+t)^{-2}
\end{align*}
for all $k=0,...,l-1$.

Also, as a corollary of Proposition 10.2 we have, in connection with the first statement of Corollary 10.2.b:

$\textbf{Corollary 10.2.g}$ Under the assumptions of Proposition 10.2 we have:
\begin{align*}
 \max_{j;i_{1}...i_{k}}\|R_{i_{k}}...R_{i_{1}}z^{j}\|_{L^{2}(\Sigma_{t}^{\epsilon_{0}})}
\leq C_{l}\{\mathcal{Y}_{0}+(1+t)\mathcal{A}_{[l-1]}+\mathcal{W}_{[l]}\}
\end{align*}
for all $k=0,...,l$.

Moreover,
\begin{align*}
 \max_{i_{1}...i_{k}}\|R_{i_{k}}...R_{i_{1}}\psi_{L}\|_{L^{2}(\Sigma_{t}^{\epsilon_{0}})}
\leq C_{l}\{\mathcal{W}_{[l]}+\delta_{0}(1+t)^{-1}[\mathcal{Y}_{0}+(1+t)\mathcal{A}_{[l-1]}]\}
\end{align*}
for all $k=1,..,l$.

And:
\begin{align*}
 \max_{i_{1}...i_{k}}\|\slashed{\mathcal{L}}_{R_{i_{k}}}...\slashed{\mathcal{L}}_{R_{i_{1}}}\slashed{\omega}_{L}\|_{L^{2}(\Sigma_{t}^{\epsilon_{0}})}
\leq C_{l}(1+t)^{-1}\{\mathcal{W}_{[l]}+\delta_{0}(1+t)^{-1}[\mathcal{Y}_{0}+(1+t)\mathcal{A}_{[l-2]}]\}
\end{align*}
for all $k=0,...,l-1$.
\begin{align*}
 \max_{i_{1}...i_{l}}\|\slashed{\mathcal{L}}_{R_{i_{k}}}...\slashed{\mathcal{L}}_{R_{i_{1}}}\omega_{L\hat{T}}\|_{L^{2}(\Sigma_{t}^{\epsilon_{0}})}
\leq C_{l}(1+t)^{-1}\{\mathcal{W}^{Q}_{[l-1]}+\delta_{0}(1+t)^{-1}[\mathcal{Y}_{0}+(1+t)\mathcal{A}_{l-2]}+\mathcal{W}_{[l-1]}]\}
\end{align*}
for all $k=0,...,l-1$.\vspace{7mm}

    We turn to the $S_{t,u}$ 1-form $\kappa^{-1}\zeta$. This is given by (3.54):
\begin{align}
 \kappa^{-1}\zeta=\alpha\varepsilon-\slashed{d}\alpha
\end{align}
We first consider, under the assumptions of Proposition 10.1, $L^{\infty}$ estimates on $\Sigma_{t}^{\epsilon_{0}}$ of the angular derivatives of $\kappa^{-1}\zeta$ 
up to order $l-1$. Using Corollary 10.1.b, we readily deduce:

$\textbf{Corollary 10.1.h}$ Under the assumptions of Proposition 10.1 we have:
\begin{align*}
 \max_{i_{1}...i_{k}}\|\slashed{\mathcal{L}}_{R_{i_{k}}}...\slashed{\mathcal{L}}_{R_{i_{1}}}(\kappa^{-1}\zeta)\|_{L^{\infty}(\Sigma_{t}^{\epsilon_{0}})}
\leq C_{l}\delta_{0}(1+t)^{-2}
\end{align*}
for all $k=0,...,l-1$.

    We consider next, under the assumptions of Proposition 10.2, $L^{2}$ estimates on $\Sigma_{t}^{\epsilon_{0}}$ of the angular derivatives of $\kappa^{-1}\zeta$ up to 
order $l-1$. By Corollary 10.2.g and (10.133), it is easy to see that

$\textbf{Corollary 10.2.h}$ Under the assumptions of Proposition 10.2 we have:
\begin{align*}
 \max_{i_{1}...i_{k}}\|\slashed{\mathcal{L}}_{R_{i_{k}}}...\slashed{\mathcal{L}}_{R_{i_{1}}}(\kappa^{-1}\zeta)\|_{L^{2}(\Sigma_{t}^{\epsilon_{0}})}
\leq C_{l}(1+t)^{-1}\{\mathcal{W}_{[l]}+\delta_{0}(1+t)^{-1}[\mathcal{Y}_{0}+(1+t)\mathcal{A}_{[l-2]}]\}
\end{align*}
for all $k=0,...,l-1$.

    We turn to the $S_{t,u}$-tangential vectorfields $\leftexp{(R_{i})}{Z}$. These are related to the $S_{t,u}$ 1-form $\leftexp{(R_{i})}{\slashed{\pi}}_{L}$ by:
\begin{align}
 \leftexp{(R_{i})}{Z}=\leftexp{(R_{i})}{\slashed{\pi}}_{L}\cdot\slashed{g}^{-1}
\end{align}
According to (6.67):
\begin{align}
 \leftexp{(R_{i})}{\slashed{\pi}}_{L}=-(\chi-\frac{\slashed{g}}{1-u+t})\cdot R_{i}+\epsilon_{ijm}z^{j}\slashed{d}x^{m}+\lambda_{i}(\kappa^{-1}\zeta)
\end{align}
Using the corollaries of Proposition 10.1 we readily deduce:

$\textbf{Corollary 10.1.i}$ Under the assumptions of Proposition 10.1 we have:
\begin{align*}
 \max_{i;i_{1}...i_{k}}\|\slashed{\mathcal{L}}_{R_{i_{k}}}...\slashed{\mathcal{L}}_{R_{i_{1}}}\leftexp{(R_{i})}{\slashed{\pi}}_{L}\|
_{L^{\infty}(\Sigma_{t}^{\epsilon_{0}})}\leq C_{l}\delta_{0}(1+t)^{-1}[1+\log(1+t)]
\end{align*}
for all $k=0,...,l-1$, and:
\begin{align*}
 \max_{i;i_{1}...i_{k}}\|\slashed{\mathcal{L}}_{R_{i_{k}}}...\slashed{\mathcal{L}}_{R_{i_{1}}}\leftexp{(R_{i})}{Z}\|_{L^{\infty}(\Sigma_{t}^{\epsilon_{0}})}
\leq C_{l}\delta_{0}(1+t)^{-1}[1+\log(1+t)]
\end{align*}
for all $k=0,...,l-1$.

    Also, using the corollaries of Proposition 10.1 with $l_{*}$ in the role of $l$, as well as the corollaries of Proposition 10.2 we can easily deduce:

$\textbf{Corollary 10.2.i}$ Under the assumptions of Proposition 10.2 we have:
\begin{align*}
 \max_{i;i_{1}...i_{k}}\|\slashed{\mathcal{L}}_{R_{i_{k}}}...\slashed{\mathcal{L}}_{R_{i_{1}}}\leftexp{(R_{i})}{\slashed{\pi}}_{L}\|_{L^{2}(\Sigma_{t}^{\epsilon_{0}})}\\
\leq C_{l}\{\mathcal{Y}_{0}+(1+t)\mathcal{A}_{[l-1]}+\delta_{0}(1+t)^{-1}[1+\log(1+t)]\mathcal{W}_{[l]}+\mathcal{W}_{[l-1]}\}
\end{align*}
for all $k=0,...,l-1$, and:
\begin{align*}
 \max_{i;i_{1}...i_{k}}\|\slashed{\mathcal{L}}_{R_{i_{k}}}...\slashed{\mathcal{L}}_{R_{i_{1}}}\leftexp{(R_{i})}{Z}\|_{L^{2}(\Sigma_{t}^{\epsilon_{0}})}\\
\leq C_{l}\{\mathcal{Y}_{0}+(1+t)\mathcal{A}_{[l-1]}+\delta_{0}(1+t)^{-1}[1+\log(1+t)]\mathcal{W}_{[l]}+\mathcal{W}_{[l-1]}\}
\end{align*}
for all $k=0,...,l-1$.

\section{Bounds for the quantities $Q_{l}$ and $P_{l}$}
The object of this section is to obtain appropriate bounds for the quantities $\leftexp{(i_{1}...i_{l})}{Q}_{l}$ and $\leftexp{(i_{1}...i_{l})}{P}_{l}$, which,
through $\leftexp{(i_{1}...i_{l})}{B}_{l}$, enter the final estimates of Chapter 8 for the 1-forms $\leftexp{(i_{1}...i_{l})}{x}_{l}$.
\subsection{Estimates for ${Q}_{l}$}
     We begin with the $S_{t,u}$ 1-form $i$ which is given by (8.100):
\begin{align}
 i=\beta^{*}_{A}-\alpha^{-1}[(\chi-\frac{\slashed{g}}{1-u+t})\cdot\slashed{g}^{-1}-(\textrm{tr}\chi-\frac{2}{1-u+t})I]\cdot(\kappa^{-1}\zeta)\\\notag
:=i_{1}+i_{2}
\end{align}
By (4.32)-(4.33), the first term on the right of above is given by:
\begin{align}
 i_{1}=\beta^{*}_{A}=\eta^{-2}(\slashed{g}^{-1})^{BC}[\frac{1}{2}X_{A}(-\eta^{2}+|\textbf{v}|^{2})\slashed{\omega}(X_{B},X_{C})-
\frac{1}{2}X_{B}(-\eta^{2}+|\textbf{v}|^{2})\slashed{\omega}(X_{A},X_{C})]\\\notag
-\eta^{-2}(\slashed{g}^{-1})^{BC}[\slashed{\omega}(X_{A},X_{C})L^{j}X_{B}(\psi_{j})-L^{j}X_{A}(\psi_{j})\slashed{\omega}(X_{C},X_{B})]
\end{align}
 It is obvious that, under the assumptions of Proposition 10.1, with $l$ replaced by $l+1$, we have: 
\begin{align}
 \max_{i_{1}...i_{k}}\|\slashed{\mathcal{L}}_{R_{i_{k}}}...\slashed{\mathcal{L}}_{R_{i_{1}}}i_{1}\|_{L^{\infty}(\Sigma_{t}^{\epsilon_{0}})}
\leq C_{l}\delta_{0}^{2}(1+t)^{-4}
\end{align}
for all $k=0,...,l$.

Also, it is easy to see that
\begin{align}
 \max_{i_{1}...i_{k}}\|\slashed{\mathcal{L}}_{R_{i_{k}}}...\slashed{\mathcal{L}}_{R_{i_{1}}}i_{2}\|_{L^{\infty}(\Sigma_{t}^{\epsilon_{0}})}
\leq C_{l}\delta_{0}^{2}(1+t)^{-4}[1+\log(1+t)]
\end{align}
for all $k=0,...,l$. We thus obtain:

$\textbf{Lemma 10.7}$ Let $\textbf{H0}$ and (6.177) hold. Let also the bootstrap assumptions $\slashed{\textbf{E}}_{[l+1]}$, $\slashed{\textbf{E}}^{Q}_{[l]}$
and $\slashed{\textbf{X}}_{[l]}$ hold, for some positive integer $l$. Then if $\delta_{0}$ is suitably small (depending on $l$) we have:
\begin{align*}
 \max_{i_{1}...i_{k}}\|\slashed{\mathcal{L}}_{R_{i_{k}}}...\slashed{\mathcal{L}}_{R_{i_{1}}}i\|_{L^{\infty}(\Sigma_{t}^{\epsilon_{0}})}
\leq C_{l}\delta_{0}^{2}(1+t)^{-4}[1+\log(1+t)]
\end{align*}
for all $k=0,...,l$. \vspace{7mm}

By direct calculation and the corollaries of Proposition 10.2, with $l$ replaced by $l+1$, it is not difficult to see that
\begin{align}
 \|\slashed{\mathcal{L}}_{R_{i_{k}}}...\slashed{\mathcal{L}}_{R_{i_{1}}}i_{1}\|_{L^{2}(\Sigma_{t}^{\epsilon_{0}})}
\leq C_{l}\delta_{0}(1+t)^{-3}\{\mathcal{W}_{[l+1]}+\delta_{0}(1+t)^{-1}[\mathcal{Y}_{0}+(1+t)\mathcal{A}_{[l-1]}]\}
\end{align}
for all $k=0,...,l$.

Similarly, we have:
\begin{align}
 \|\slashed{\mathcal{L}}_{R_{i_{k}}}...\slashed{\mathcal{L}}_{R_{i_{1}}}i_{2}\|_{L^{2}(\Sigma_{t}^{\epsilon_{0}})}
\leq C_{l}\delta_{0}(1+t)^{-2}\mathcal{A}_{[l]}+\delta_{0}(1+t)^{-3}[1+\log(1+t)]\{\mathcal{W}_{[l+1]}+\delta_{0}(1+t)^{-1}\mathcal{Y}_{0}\}
\end{align}
for all $k=0,...,l$. We thus obtain:

$\textbf{Lemma 10.8}$ Let $\textbf{H0}$ and (6.177) hold. Let $l$ be a positive integer and let the bootstrap assumptions $\slashed{\textbf{E}}_{l_{*}+1}$,
$\slashed{\textbf{E}}^{Q}_{[l_{*}]}$ and $\slashed{\textbf{X}}_{[l_{*}]}$ hold. Then we have:
\begin{align}
 \|\slashed{\mathcal{L}}_{R_{i_{k}}}...\slashed{\mathcal{L}}_{R_{i_{1}}}i\|_{L^{2}(\Sigma_{t}^{\epsilon_{0}})}
\leq C_{l}\delta_{0}(1+t)^{-2}\mathcal{A}_{[l]}+\delta_{0}(1+t)^{-3}[1+\log(1+t)]\{\mathcal{W}_{[l+1]}+\delta_{0}(1+t)^{-1}\mathcal{Y}_{0}\}
\end{align}
for all $k=0,...,l$. \vspace{7mm}

     We proceed to estimate, in $L^{2}(\Sigma_{t}^{\epsilon_{0}})$, the $S_{t,u}$ 1-form $\leftexp{(i_{1}...i_{l})}{i}_{l}$ given by Proposition 8.4,
\begin{align}
 \leftexp{(i_{1}...i_{l})}{i}_{l}=(\slashed{\mathcal{L}}_{R_{i_{l}}}+\frac{1}{2}\textrm{tr}\leftexp{(R_{i_{l}})}{\slashed{\pi}})...
(\slashed{\mathcal{L}}_{R_{i_{1}}}+\frac{1}{2}\textrm{tr}\leftexp{(R_{i_{1}})}{\slashed{\pi}})i\\\notag
+\sum_{k=0}^{l-1}(\slashed{\mathcal{L}}_{R_{i_{l}}}+\frac{1}{2}\textrm{tr}\leftexp{(R_{i_{l}})}{\slashed{\pi}})...
(\slashed{\mathcal{L}}_{R_{i_{l-k+1}}}+\frac{1}{2}\textrm{tr}\leftexp{(R_{i_{l-k+1}})}{\slashed{\pi}})\leftexp{(i_{1}...i_{l-k})}{q}_{l-k}
\end{align}
under the assumptions of Lemma 10.8.

     We begin with the first term on the right of (10.143). We have:
\begin{align}
 (\slashed{\mathcal{L}}_{R_{i_{l}}}+\frac{1}{2}\textrm{tr}\leftexp{(R_{i_{l}})}{\slashed{\pi}})...(\slashed{\mathcal{L}}_{R_{i_{1}}}+
\frac{1}{2}\textrm{tr}\leftexp{(R_{i_{l}})}{\slashed{\pi}})i\\\notag
=\slashed{\mathcal{L}}_{R_{i_{l}}}...\slashed{\mathcal{L}}_{R_{i_{1}}}i+\sum_{j=1}^{l}\sum_{k_{1}<...<k_{j}=1}^{l}
\slashed{\mathcal{L}}_{R_{i_{l}}}...\overset{>\slashed{\mathcal{L}}_{R_{{i_{k}}_{j}}}<}{(\frac{1}{2}\textrm{tr}\leftexp{(R_{{i_{k}}_{j}})}{\slashed{\pi}})}...
\overset{>\slashed{\mathcal{L}}_{R_{{i_{k}}_{1}}}<}{(\frac{1}{2}\textrm{tr}\leftexp{(R_{{i_{k}}_{1}})}{\slashed{\pi}})}...\slashed{\mathcal{L}}_{R_{i_{1}}}i
\end{align}
We have two cases to consider:

Case 1: At most $l_{*}$ angular derivatives fall on one of the factors
\begin{align*}
 \frac{1}{2}\textrm{tr}\leftexp{(R_{{i_{k}}_{1}})}{\slashed{\pi}},...,\frac{1}{2}\textrm{tr}\leftexp{(R_{{i_{k}}_{j}})}{\slashed{\pi}}\quad \textrm{and}:
\end{align*}
Case 2: one of these factors receives more than $l_{*}$ angular derivatives. Then all the other factors receive at most $l_{*}$ angular derivatives.

    In Case 1 we just apply Corollary 10.1.d with $l_{*}+1$ in the role of $l$, then we know that the second term on the right in (10.144) is bounded by:
\begin{align*}
 C_{l}(\delta_{0}(1+t)^{-1}[1+\log(1+t)])^{j}\{\delta_{0}(1+t)^{-2}\mathcal{A}_{[l-j]}\\\notag
+\delta_{0}(1+t)^{-3}[1+\log(1+t)][\mathcal{W}_{[l-j+1]}+\delta_{0}(1+t)^{-1}\mathcal{Y}_{0}]\}
\end{align*}

    In Case 2 the number of angular derivatives falling on $i$ is at most $l_{*}$. So by Lemma 10.7 and Corollary 10.2.d, we have an $L^{2}(\Sigma_{t}^{\epsilon_{0}})$ bound for Case 2, by:
\begin{align*}
 C_{l}(\delta_{0}(1+t)^{-1}[1+\log(1+t)])^{j-1}\cdot\delta_{0}(1+t)^{-4}[1+\log(1+t)]\cdot\\\notag
\{\mathcal{Y}_{0}+(1+t)\mathcal{A}_{[l-j]}+\mathcal{W}_{[l-j+1]}\}
\end{align*}
Here, we have used the following facts that under the assumptions of Proposition 10.2, we have, by Corollary 10.2.d,
\begin{align}
 \max_{i;i_{1}...i_{k}}\|\slashed{\mathcal{L}}_{R_{i_{k}}}...\slashed{\mathcal{L}}_{R_{i_{1}}}\leftexp{(R_{i})}{\slashed{\pi}}\|_{L^{2}(\Sigma_{t}^{\epsilon_{0}})}
\leq C_{l}\{\mathcal{Y}_{0}+(1+t)\mathcal{A}_{[l-1]}+\mathcal{W}_{[l]}\}
\end{align}

Combining the results, we obtain a bound for the second term on the right in (10.144) by:
\begin{align}
 C_{l}(\delta_{0}(1+t)^{-1}[1+\log(1+t)])^{j-1}\{\delta_{0}(1+t)^{-2}\mathcal{A}_{[l-j]}\\\notag
+\delta_{0}(1+t)^{-3}[1+\log(1+t)][\mathcal{W}_{[l-j+1]}+\delta_{0}(1+t)^{-1}\mathcal{Y}_{0}]\}
\end{align}
    Next we consider the second term on the right in (10.143). First, we define the operator $\check{D}$:
\begin{align}
 (\check{\slashed{D}}\vartheta)_{ABC}=\frac{1}{2}(\slashed{D}_{A}\vartheta_{BC}+\slashed{D}_{B}\vartheta_{AC}-\slashed{D}_{C}\vartheta_{AB})
\end{align}
where $\vartheta$ is a symmetric 2-covariant $S_{t,u}$ tensorfield.

So the second term in the expression for $\leftexp{(i_{1}...i_{j})}{q}_{j}$ of Proposition 8.4 is simply
\begin{align*}
 (\check{\slashed{D}}\leftexp{(i_{1}...i_{j-1})}{\hat{\chi}}_{j-1})\cdot(\slashed{g}^{-1}\cdot\leftexp{(R_{j})}{\hat{\slashed{\pi}}}\cdot\slashed{g}^{-1})
\end{align*}
Also, the third term in the same expression is:
\begin{align*}
 \leftexp{(i_{1}...i_{j-1})}{\hat{\chi}}_{j-1}\cdot(\slashed{\textrm{div}}\leftexp{(R_{j})}{\hat{\slashed{\pi}}}\cdot\slashed{g}^{-1})
\end{align*}
While by a direct calculation, we have:
\begin{align}
 \slashed{\textrm{div}}\leftexp{(X)}{\hat{\slashed{\pi}}}\cdot\slashed{g}^{-1}=\textrm{tr}\leftexp{(X)}{\slashed{\pi}}_{1}
\end{align}
Here the definition of $\leftexp{(X)}{\slashed{\pi}}_{1}$ is (9.140):
\begin{align*}
\leftexp{(X)}{\slashed{\pi}}_{1,AB}^{C}=\frac{1}{2}(\slashed{D}_{A}\leftexp{(X)}{\slashed{\pi}}^{C}_{B}+\slashed{D}_{B}\leftexp{(X)}{\slashed{\pi}}_{A}^{C}-\slashed{D}^{C}\leftexp{(X)}{\slashed{\pi}}_{AB})
\end{align*}
We conclude that
\begin{align}
 \leftexp{(i_{1}...i_{j})}{q}_{j}=\frac{1}{4}\textrm{tr}\leftexp{(R_{j})}{\slashed{\pi}}\slashed{d}(R_{i_{j-1}}...R_{i_{1}}\textrm{tr}\chi)
+(\check{\slashed{D}}\leftexp{(i_{1}...i_{j-1})}{\hat{\chi}}_{j-1})\cdot(\slashed{g}^{-1}\cdot\leftexp{(R_{j})}{\hat{\slashed{\pi}}}\cdot\slashed{g}^{-1})\\\notag
+\leftexp{(i_{1}...i_{j-1})}{\hat{\chi}}_{j-1}\cdot\textrm{tr}\leftexp{(R_{i_{j}})}{\slashed{\pi}}_{1}
\end{align}
To estimate the contribution of the above terms, we need some lemmas, which can be proved in a straightforward manner.

$\textbf{Lemma 10.9}$ Let $h$ be an arbitrary smooth function on $S_{t,u}$ and $\vartheta$ be an arbitrary symmetric 2-covariant $S_{t,u}$ tensorfield. 
The following commutation formulas hold for any non-negative
integer $k$:
\begin{align*}
 \leftexp{(i_{1}...i_{k})}{c}_{k}:=\slashed{\mathcal{L}}_{R_{i_{k}}}...\slashed{\mathcal{L}}_{R_{i_{1}}}\slashed{D}^{2}h
-\slashed{D}^{2}(R_{i_{k}}...R_{i_{1}}h)\\
=-\sum_{m=0}^{k-1}\slashed{\mathcal{L}}_{R_{i_{k}}}...\slashed{\mathcal{L}}_{R_{i_{k-m+1}}}(\leftexp{(R_{i_{k-m}})}{\slashed{\pi}_{1}}
\cdot\slashed{d}(R_{i_{k-m-1}}...R_{i_{1}}h))
\end{align*}
and
\begin{align*}
 \slashed{\mathcal{L}}_{R_{i_{k}}}...\slashed{\mathcal{L}}_{R_{i_{1}}}\check{\slashed{D}}\vartheta
-\check{\slashed{D}}\slashed{\mathcal{L}}_{R_{i_{k}}}...\slashed{\mathcal{L}}_{R_{i_{1}}}\vartheta\\\notag
=-\sum_{m=0}^{k-1}\slashed{\mathcal{L}}_{R_{i_{k}}}...\slashed{\mathcal{L}}_{R_{i_{k-m+1}}}(\leftexp{(R_{i_{k-m}})}{\slashed{\pi}}_{1}
\cdot\slashed{\mathcal{L}}_{R_{i_{k-m-1}}}...\slashed{\mathcal{L}}_{R_{i_{1}}}\vartheta)
\end{align*}
The following assumption, which is another version of $\textbf{H2}$, is used in the proof of the following lemmas.

$\textbf{H2}^{\prime}$ There is a constant $C$ (independent of $s$) such that for any 2-covariant symmetric $S_{t,u}$ tensorfield $\vartheta$ we have, pointwise:
\begin{align*}
 |\slashed{D}\vartheta|^{2}\leq C(1+t)^{-2}\{\sum_{j}|\slashed{\mathcal{L}}_{R_{j}}\vartheta|^{2}+|\vartheta|^{2}\}
\end{align*}

$\textbf{Lemma 10.10}$ Let the hypotheses $\textbf{H0}$, $\textbf{H1}$, $\textbf{H2}^{\prime}$, and (6.177) hold. Let also the bootstrap assumptions 
$\slashed{\textbf{E}}_{[l+1]}$, $\slashed{\textbf{E}}_{[l]}^{Q}$ and $\slashed{\textbf{X}}_{[l]}$ hold, for some positive integer $l$. Then if $\delta_{0}$
is suitably small (depending on $l$) we have:
\begin{align*}
 \max_{j;i_{1}...i_{k}}\|\slashed{\mathcal{L}}_{R_{i_{k}}}...\slashed{\mathcal{L}}_{R_{i_{1}}}\leftexp{(R_{j})}{\slashed{\pi}}_{1}\|
_{L^{\infty}(\Sigma_{t}^{\epsilon_{0}})}\leq C_{l}\delta_{0}(1+t)^{-2}[1+\log(1+t)]
\end{align*}
for all $k=0,...,l-1$.

$\textbf{Lemma 10.11}$ Let the hypotheses $\textbf{H0}$, $\textbf{H1}$, $\textbf{H2}^{\prime}$, and (6.177) hold. Let also the bootstrap assumptions 
$\slashed{\textbf{E}}_{[l_{*}+1]}$, $\slashed{\textbf{E}}_{[l_{*}]}^{Q}$ and $\slashed{\textbf{X}}_{[l_{*}]}$ hold, for some positive integer $l$. Then if $\delta_{0}$
is suitably small (depending on $l$) we have:
\begin{align*}
 \max_{j;i_{1}...i_{k}}\|\slashed{\mathcal{L}}_{R_{i_{k}}}...\slashed{\mathcal{L}}_{R_{i_{1}}}\leftexp{(R_{j})}{\slashed{\pi}}_{1}\|
_{L^{2}(\Sigma_{t}^{\epsilon_{0}})}\leq C_{l}(1+t)^{-1}\{\mathcal{Y}_{0}+(1+t)\mathcal{A}_{[l]}+\mathcal{W}_{[l+1]}\}
\end{align*}
for all $k=0,...,l-1$.

Using these lemmas, it is not difficult to obtain the following $L^{2}$ estimate for the second term in (10.143):
\begin{align}
 \max_{i_{1}...i_{l}}\|\sum_{k=0}^{l-1}(\slashed{\mathcal{L}}_{R_{i_{l}}}+\frac{1}{2}\textrm{tr}\leftexp{(R_{i_{l}})}{\slashed{\pi}})
...(\slashed{\mathcal{L}}_{R_{i_{l-k+1}}}+\frac{1}{2}\textrm{tr}\leftexp{(R_{i_{l-k+1}})}{\slashed{\pi}})\leftexp{(i_{1}...i_{l-k})}{q}_{l-k}\|
_{L^{2}(\Sigma_{t}^{\epsilon_{0}})}\\\notag
\leq C_{l}\delta_{0}(1+t)^{-3}[1+\log(1+t)]\{\mathcal{Y}_{0}+(1+t)\mathcal{A}_{[l]}+\mathcal{W}_{[l+1]}\}
\end{align}
Combining this with (10.145) and (10.142) we obtain:

$\textbf{Proposition 10.3}$ Under the assumptions of Lemma 10.14 we have:
\begin{align*}
 \max_{i_{1}...i_{l}}\|\leftexp{(i_{1}...i_{l})}{i}_{l}\|_{L^{2}(\Sigma_{t}^{\epsilon_{0}})}
\leq C_{l}\delta_{0}(1+t)^{-3}[1+\log(1+t)]\{\mathcal{Y}_{0}+(1+t)\mathcal{A}_{[l]}+\mathcal{W}_{[l+1]}\}
\end{align*}
provided that $\delta_{0}$ is suitably small (depending on $l$).

     We turn to the $S_{t,u}$ 1-form $\leftexp{(i_{1}...i_{l})}{\dot{g}}_{l}$:
\begin{align}
 \leftexp{(i_{1}...i_{l})}{\dot{g}}_{l}=\slashed{\mathcal{L}}_{R_{i_{l}}}...\slashed{\mathcal{L}}_{R_{i_{1}}}g_{0}-\leftexp{(i_{1}...i_{l})}{w}_{l}
+\sum_{k=0}^{l-1}\slashed{\mathcal{L}}_{R_{i_{l}}}...\slashed{\mathcal{L}}_{R_{i_{l-k+1}}}\leftexp{(i_{1}...i_{l-k})}{y}_{l-k}
\end{align}
where $\leftexp{(i_{1}...i_{j})}{y}_{j}$ is given in Proposition 8.3, and $\leftexp{(i_{1}...i_{l})}{w}_{l}$ is
\begin{align}
\leftexp{(i_{1}...i_{l})}{w}_{l}=\sum_{k=0}^{l-1}\slashed{\mathcal{L}}_{R_{i_{l}}}...\slashed{\mathcal{L}}_{R_{i_{l-k+1}}}
\slashed{\mathcal{L}}_{\leftexp{(R_{i_{l-k}})}{Z}}\leftexp{(i_{1}...i_{l-k-1})}{x}_{l-k-1}
-\sum_{k=0}^{l-1}\slashed{\mathcal{L}}_{\leftexp{(R_{i_{l-k}})}{Z}}\leftexp{(i_{1}\overset{>i_{l-k}<}{...}i_{l})}{x}_{l-1}
\end{align}

    We now introduce bootstrap assumptions in regard to the angular derivatives of the function $\mu$. Given a positive integer $l$ let us denote by 
$\slashed{\textbf{M}}_{l}$ the bootstrap assumption that there is a constant $C$ independent of $s$ such that for all $t\in[0,s]$:
\begin{align*}
 \slashed{\textbf{M}}_{l}\quad:\quad \max_{i_{1}...i_{l}}\|R_{i_{l}}...R_{i_{1}}\mu\|_{L^{\infty}(\Sigma_{t}^{\epsilon_{0}})}\leq C\delta_{0}[1+\log(1+t)]
\end{align*}
Let us also denote by $\slashed{\textbf{M}}_{0}$ the bootstrap assumption:
\begin{align*}
 \slashed{\textbf{M}}_{0}\quad:\quad |\mu-1|\leq C\delta_{0}[1+\log(1+t)]
\end{align*}
We then denote by $\slashed{\textbf{M}}_{[l]}$ the bootstrap assumption
\begin{align*}
 \slashed{\textbf{M}}_{[l]}\quad:\quad \slashed{\textbf{M}}_{0}\textrm{and}...\textrm{and}\slashed{\textbf{M}}_{l}
\end{align*}
We also introduce the quantities:
\begin{align*}
 \mathcal{B}_{0}=\|\mu-1\|_{L^{2}(\Sigma_{t}^{\epsilon_{0}})}
\end{align*}
and for $k\slashed{=}0$:
\begin{align*}
 \mathcal{B}_{k}=\max_{i_{1}...i_{k}}\|R_{i_{k}}...R_{i_{1}}\mu\|_{L^{2}(\Sigma_{t}^{\epsilon_{0}})}\quad \textrm{and}:\quad 
\mathcal{B}_{[l]}=\sum_{k=0}^{l}\mathcal{B}_{k}
\end{align*}
Then by the expression (8.27) and (8.29) we obtain:

$\textbf{Lemma 10.12}$ Let hypothesis $\textbf{H0}$ and the estimate (6.177) hold. Let also the bootstrap assumptions $\slashed{\textbf{E}}^{T}_{[l_{*}+1]}$,
$\slashed{\textbf{E}}^{Q}_{[l_{*}+1]}$, $\slashed{\textbf{E}}_{[l_{*}+2]}$, as well as $\slashed{\textbf{X}}_{[l_{*}]}$ and $\slashed{\textbf{M}}_{[l_{*}]}$ hold,
for some positive integer $l$. We then have:
\begin{align*}
 \max_{i_{1}...i_{l}}\|\slashed{d}R_{i_{l}}...R_{i_{1}}\check{g}\|_{L^{2}(\Sigma_{t}^{\epsilon_{0}})}\\\notag
\leq C_{l}\delta_{0}(1+t)^{-3}\{\mathcal{W}^{T}_{[l+1]}+\mathcal{W}^{Q}_{[l+1]}+(1+t)^{-1}[1+\log(1+t)]\mathcal{W}_{[l+2]}\\\notag
+\delta_{0}(1+t)^{-1}[1+\log(1+t)][\mathcal{Y}_{0}+(1+t)\mathcal{A}_{[l]}]+\delta_{0}(1+t)^{-2}\mathcal{B}_{[l+1]}\}
\end{align*}
provided that $\delta_{0}$ is suitably small (depending on $l$). The factor $[1+\log(1+t)]$ comes from the product of $\mu$ and a regular term.

Here, the bootstrap assumption $\slashed{\textbf{E}}^{T}_{[l]}$ is defined as follows:
\begin{align*}
 \slashed{\textbf{E}}^{T}_{[l]}\quad:\quad \slashed{\textbf{E}}^{T}_{0}\textrm{and}...\textrm{and}\slashed{\textbf{E}}^{T}_{l}
\end{align*}
where $\slashed{\textbf{E}}^{T}_{l}$ is the assumption:
\begin{align*}
 \slashed{\textbf{E}}_{l}^{T}\quad:\quad\max_{i_{1}...i_{l}}\max_{\alpha}\|R_{i_{l}}...R_{i_{1}}T\psi_{\alpha}\|_{L^{\infty}(\Sigma_{t}^{\epsilon_{0}})}
\leq C\delta_{0}(1+t)^{-1}
\end{align*}

     We now return to the equation (8.37) for the $S_{t,u}$ 1-form $g_{0}$, to estimate the remaining terms. By the equation 
\begin{align*}
 L\mu=m+\mu e
\end{align*}
and (8.27), we have:
\begin{align}
 \max_{i_{1}...i_{l}}\|\slashed{\mathcal{L}}_{R_{i_{l}}}...\slashed{\mathcal{L}}_{R_{i_{1}}}[\textrm{tr}\chi\slashed{d}(\check{f}+2L\mu)]\|
_{L^{2}(\Sigma_{t}^{\epsilon_{0}})}\\\notag
\leq C_{l}(1+t)^{-3}\{[1+\log(1+t)][\mathcal{W}^{Q}_{[l+1]}+\delta_{0}(1+t)^{-1}(\mathcal{W}_{[l+1]}+\mathcal{Y}_{0}+(1+t)\mathcal{A}_{[l]})]\\\notag
+\delta_{0}(1+t)^{-1}\mathcal{B}_{[l+1]}\}
\end{align}
under the assumptions of Lemma 10.12 except for the assumption $\slashed{\textbf{E}}^{T}_{[l_{*}+1]}$, due to the cancellation between $m$ and 
$\check{f}$. Also due to this cancellation, the coefficient of $\mathcal{A}_{[l]}$ has the decay factor $(1+t)^{-3}[1+\log(1+t)]$.

    Next, we consider the last term in (8.37), which is, by (8.38),
\begin{align*}
 (\slashed{d}\mu)(\mu^{-1}(L\mu)\textrm{tr}\chi-\textrm{tr}\alpha)
\end{align*}
then under the assumptions of Lemma 10.12 except for $\slashed{\textbf{E}}^{T}_{[l_{*}+1]}$, we have:
\begin{align}
 \max_{i_{1}...i_{l}}\|\slashed{\mathcal{L}}_{R_{i_{l}}}...\slashed{\mathcal{L}}_{R_{i_{1}}}
[(\slashed{d}\mu)(\mu^{-1}(L\mu)\textrm{tr}\chi-\textrm{tr}\alpha)]\|_{L^{2}(\Sigma_{t}^{\epsilon_{0}})}\\\notag
\leq C_{l}(1+t)^{-3}\{[1+\log(1+t)][\mathcal{W}^{Q}_{[l+1]}+\delta_{0}(1+t)^{-1}(\mathcal{W}_{[l+2]}+\mathcal{Y}_{0}+(1+t)\mathcal{A}_{[l]})]\\\notag
+\delta_{0}(1+t)^{-1}\mathcal{B}_{[l+1]}\}
\end{align}
The calculation of this can be found in Chapter 12.

Combining the above results, we have:

$\textbf{Lemma 10.13}$ Let the hypotheses $\textbf{H0}$, $\textbf{H1}$, $\textbf{H2}^{\prime}$ and the estimate (6.177) hold. Let also the bootstrap assumptions
$\slashed{\textbf{E}}_{[l_{*}+2]}$, $\slashed{\textbf{E}}_{[l_{*}+1]}^{Q}$ and $\slashed{\textbf{X}}_{[l_{*}]}$, $\slashed{\textbf{M}}_{[l_{*}+1]}$ hold, for some 
positive integer $l$. Then if $\delta_{0}$ is suitably small (depending on $l$) we have:
\begin{align*}
 \max_{i_{1}...i_{l}}\|\slashed{\mathcal{L}}_{R_{i_{l}}}...\slashed{\mathcal{L}}_{R_{i_{1}}}g_{0}\|_{L^{2}(\Sigma_{t}^{\epsilon_{0}})}\leq C_{l}(1+t)^{-3}\cdot\\\notag
\{[1+\log(1+t)](\mathcal{W}^{Q}_{[l+1]}+\delta_{0}(1+t)^{-1}\mathcal{W}_{[l+2]})+\delta_{0}\mathcal{W}_{[l+1]}^{T}\\\notag
+\delta_{0}(1+t)^{-1}[1+\log(1+t)](\mathcal{Y}_{0}+(1+t)\mathcal{A}_{[l]})+\delta_{0}(1+t)^{-1}\mathcal{B}_{[l+1]}\}
\end{align*} \vspace{7mm}

     We turn to the third term in the formula (10.150) for $\leftexp{(i_{1}...i_{l})}{\dot{g}}_{l}$:
\begin{align}
 \sum_{k=0}^{l-1}\slashed{\mathcal{L}}_{R_{i_{l}}}...\slashed{\mathcal{L}}_{R_{i_{l-k+1}}}\leftexp{(i_{1}...i_{l-k})}{y}_{l-k}
\end{align}
where the $S_{t,u}$ 1-form $\leftexp{(i_{1}...i_{j})}{y}_{j}$, $j=1,...,l$ is given in the statement of Proposition 8.3:
\begin{align}
 \leftexp{(i_{1}...i_{j})}{y}_{j}=(R_{i_{j}}\mu)\leftexp{(i_{1}...i_{j-1})}{a}_{j-1}\\\notag
+(\mu R_{i_{j}}(\textrm{tr}\chi-e)-R_{i_{j}}m+\leftexp{(R_{i_{j}})}{Z}\mu-eR_{i_{j}}\mu)\cdot\slashed{d}(R_{i_{j-1}}...R_{i_{1}}\textrm{tr}\chi)\\\notag
+\frac{1}{2}(R_{i_{j}}\textrm{tr}\chi)\slashed{d}\leftexp{(i_{1}...i_{j-1})}{\check{f}}_{j-1}
\end{align}
for $j=1,...,l$.

The $S_{t,u}$ 1-form $\leftexp{(i_{1}...i_{j-1})}{a}_{j-1}$ is given by (8.78):
\begin{align}
 \leftexp{(i_{1}...i_{j-1})}{a}_{j-1}=\slashed{d}(R_{i_{j-1}}...R_{i_{1}}f_{0})+\leftexp{(i_{1}...i_{j-1})}{b}_{j-1}
\end{align}
where
\begin{align*}
 f_{0}=L\textrm{tr}\chi+|\chi|^{2}
\end{align*}
which we have dealt with in deriving (10.153), and:
\begin{align}
\leftexp{(i_{1}...i_{j-1})}{b}_{j-1}=\sum_{m=0}^{j-2}\slashed{d}(R_{i_{j-1}}...R_{i_{j-m}}\leftexp{(R_{i_{j-m-1}})}{Z}R_{i_{j-m-2}}...R_{i_{1}}\textrm{tr}\chi)\\\notag
+\sum_{m=1}^{j-1}\sum_{k_{1}<...<k_{m}=1}^{j-1}(R_{i_{k_{m}}}...R_{i_{k_{1}}}\textrm{tr}\chi)(\slashed{d}R_{i_{j-1}}\overset{>i_{k_{m}}...i_{k_{1}}<}{...}
R_{i_{1}}\textrm{tr}\chi) 
\end{align}
for $j\geq2$; $b_{0}=0$.

To estimate the sum (10.154) in $L^{2}(\Sigma_{t}^{\epsilon_{0}})$, we must estimate:
\begin{align}
 \max_{i_{1}...i_{j+k}}\|\slashed{\mathcal{L}}_{R_{i_{j+k}}}...\slashed{\mathcal{L}}_{R_{i_{j+1}}}\leftexp{(i_{1}...i_{j})}{y}_{j}\|_{L^{2}(\Sigma_{t}^{\epsilon_{0}})}
\quad:\quad j+k=l
\end{align}
To estimate this requires estimating:
\begin{align}
 \max_{i_{1}...i_{j-1+k}}\|\slashed{\mathcal{L}}_{R_{i_{j-1+k}}}...\slashed{\mathcal{L}}_{R_{i_{j}}}\leftexp{(i_{1}...i_{j-1})}{a}_{j-1}\|_{L^{2}(\Sigma_{t}^{\epsilon_{0}})}
\quad:\quad j+k\leq l
\end{align}
as well as:
\begin{align}
 \max_{i_{1}...i_{j-1+k}}\|\slashed{\mathcal{L}}_{R_{i_{j-1+k}}}...\slashed{\mathcal{L}}_{R_{i_{j}}}\leftexp{(i_{1}...i_{j-1})}{a}_{j-1}\|
_{L^{\infty}(\Sigma_{t}^{\epsilon_{0}})}
\quad:\quad j+k\leq l_{*}
\end{align}
Since we have dealt with the first term on the right in (10.156) in deriving (10.153), we just need to estimate:
\begin{align}
 \max_{i_{1}...i_{j-1+k}}\|\slashed{\mathcal{L}}_{R_{i_{j-1+k}}}...\slashed{\mathcal{L}}_{R_{i_{j}}}\leftexp{(i_{1}...i_{j-1})}{b}_{j-1}\|_{L^{2}(\Sigma_{t}^{\epsilon_{0}})}
\quad:\quad j+k\leq l
\end{align}
as well as:
\begin{align}
 \max_{i_{1}...i_{j-1+k}}\|\slashed{\mathcal{L}}_{R_{i_{j-1+k}}}...\slashed{\mathcal{L}}_{R_{i_{j}}}\leftexp{(i_{1}...i_{j-1})}{b}_{j-1}\|
_{L^{\infty}(\Sigma_{t}^{\epsilon_{0}})}
\quad:\quad j+k\leq l_{*}
\end{align}
These quantities 
they enjoy the bounds:
\begin{align}
 \max_{i_{1}...i_{j-1+k}}\|\slashed{\mathcal{L}}_{R_{i_{j-1+k}}}...\slashed{\mathcal{L}}_{R_{i_{j}}}\leftexp{(i_{1}...i_{j-1})}{b}_{j-1}\|
_{L^{\infty}(\Sigma_{t}^{\epsilon_{0}})}\leq C_{l}\delta_{0}^{2}(1+t)^{-5}[1+\log(1+t)]^{2}
\end{align}
for $j+k\leq l_{*}$,

and:
\begin{align}
 \max_{i_{1}...i_{j-1+k}}\|\slashed{\mathcal{L}}_{R_{i_{j-1+k}}}...\slashed{\mathcal{L}}_{R_{i_{j}}}\leftexp{(i_{1}...i_{j-1})}{b}_{j-1}\|
_{L^{2}(\Sigma_{t}^{\epsilon_{0}})}\leq C_{l}\delta_{0}(1+t)^{-4}[1+\log(1+t)]\cdot\\\notag
\{(1+t)\mathcal{A}_{[l]}+\mathcal{Y}_{0}+\mathcal{W}_{[l]}\}
\end{align}
for $j+k\leq l$ under the assumptions of Proposition 10.1 with $l_{*}+1$ in the role of $l$.

The above estimates together with those for $f_{0}$ yield the following lemma.

$\textbf{Lemma 10.14}$ Let the hypotheses $\textbf{H0}$, $\textbf{H1}$, $\textbf{H2}^{\prime}$, and the estimate (6.177) hold. Let also the bootstrap assumptions
$\slashed{\textbf{E}}_{[l_{*}+2]}$, $\slashed{\textbf{E}}^{Q}_{[l_{*}+1]}$ and $\slashed{\textbf{X}}_{[l_{*}]}$ hold, for some positive integer $l$. Then provided
that $\delta_{0}$ is suitably small (depending on $l$), the $S_{t,u}$ 1-form $\leftexp{(i_{1}...i_{j-1})}{a}_{j-1}$ satisfies the following estimates:
\begin{align*}
 \max_{i_{1}...i_{j-1+k}}\|\slashed{\mathcal{L}}_{R_{i_{j-1+k}}}...\slashed{\mathcal{L}}_{R_{i_{j}}}\leftexp{(i_{1}...i_{j-1})}{a}_{j-1}\|
_{L^{\infty}(\Sigma_{t}^{\epsilon_{0}})}\leq C_{l}\delta_{0}(1+t)^{-4}
\end{align*}
for $j+k\leq l_{*}$, and:
\begin{align*}
 \max_{i_{1}...i_{j-1+k}}\|\slashed{\mathcal{L}}_{R_{i_{j-1+k}}}...\slashed{\mathcal{L}}_{R_{i_{j}}}\leftexp{(i_{1}...i_{j-1})}{a}_{j-1}\|
_{L^{2}(\Sigma_{t}^{\epsilon_{0}})}\leq C_{l}(1+t)^{-3}\{\mathcal{W}_{[l+2]}+\delta_{0}(1+t)^{-1}\mathcal{W}^{Q}_{[l+1]}+\mathcal{W}^{Q}_{[l]}\\\notag
+\delta_{0}(1+t)^{-1}[1+\log(1+t)][\mathcal{Y}_{0}+(1+t)\mathcal{A}_{[l]}]\}
\end{align*}
for $j+k\leq l$. \vspace{7mm}

Using Lemma 10.14 and the bootstrap assumption $\slashed{\textbf{M}}_{[l_{*}+1]}$ we readily obtain:
\begin{align}
 \|\slashed{\mathcal{L}}_{R_{i_{j+k}}}...\slashed{\mathcal{L}}_{R_{i_{j+1}}}((R_{i_{j}}\mu)\leftexp{(i_{1}...i_{j-1})}{a}_{j-1}\|_{L^{2}(\Sigma_{t}^{\epsilon_{0}})}\\\notag
\leq C_{l}\delta_{0}(1+t)^{-3}[1+\log(1+t)]\cdot\\\notag
\{\mathcal{W}_{[l+2]}+\delta_{0}(1+t)^{-1}\mathcal{W}_{[l+1]}^{Q}+\mathcal{W}^{Q}_{[l]}\\\notag
+\delta_{0}(1+t)^{-1}[1+\log(1+t)][\mathcal{Y}_{0}+(1+t)\mathcal{A}_{[l]}]\}+C_{l}\delta_{0}(1+t)^{-4}\mathcal{B}_{[l]}
\end{align}
for $j+k=l$.

    We consider next the second term on the right in (10.155). We have the following estimates: 
Here, we also need the assumption $\slashed{\textbf{X}}_{[(l+1)_{*}]}$.
\begin{align*}
\|\slashed{\mathcal{L}}_{R_{i_{j+1}}}...\slashed{\mathcal{L}}_{R_{i_{j+k}}}((\mu R_{i_{j}}(\textrm{tr}\chi-e)-R_{i_{j}}m+\leftexp{(R_{i_{j}})}{Z}\mu-eR_{i_{j}}\mu)\cdot\slashed{d}(R_{i_{j-1}}...R_{i_{1}}\textrm{tr}\chi))\|_{L^{2}(\Sigma_{t}^{\epsilon_{0}})}\\
\leq C_{l}\delta_{0}(1+t)^{-3}\{(1+t)\mathcal{A}_{[l]}+\mathcal{W}^{T}_{[l]}+\mathcal{Y}_{0}
+\mathcal{W}_{[l]}\\
+[1+\log(1+t)]\mathcal{W}_{[l]}^{Q}+\delta_{0}(1+t)^{-1}[1+\log(1+t)]\mathcal{W}_{[l+1]}\}
\end{align*}\vspace{7mm}

     We turn to the last term on the right in (10.155). Recall that:
\begin{align}
 \check{f}=\mu f=-\frac{1}{2}\frac{dH}{dh}\tau_{\underline{L}}
\end{align}
Then by (8.24) and $\slashed{\textbf{M}}_{[l_{*}]}$ we obtain:
\begin{align}
 \max_{i_{1}...i_{k}}\|R_{i_{k}}...R_{i_{1}}\check{f}\|_{L^{\infty}(\Sigma_{t}^{\epsilon_{0}})}\leq C_{l}\delta_{0}(1+t)^{-1}\quad:\quad k\leq l_{*}
\end{align}
and also:
\begin{align}
 \max_{i_{1}...i_{k}}\|R_{i_{k}}...R_{i_{1}}\check{f}\|_{L^{2}(\Sigma_{t}^{\epsilon_{0}})}\\\notag
\leq C_{l}\{\mathcal{W}^{T}_{[l]}+\delta_{0}(1+t)^{-1}(\mathcal{W}_{[l]}+(1+t)^{-1}\mathcal{B}_{[l]})+(1+t)^{-1}[1+\log(1+t)]\mathcal{W}^{Q}_{[l]}\}
\end{align}
for $k\leq l$.

So with the assumption $\slashed{\textbf{X}}_{[(l+1)_{*}]}$, we have:
\begin{align}
 \|\slashed{\mathcal{L}}_{R_{i_{j+k}}}...\slashed{\mathcal{L}}_{R_{i_{j+1}}}((R_{i_{j}}\textrm{tr}\chi)(\slashed{d}R_{i_{j-1}}...R_{i_{1}}\check{f}))\|
_{L^{2}(\Sigma_{t}^{\epsilon_{0}})}\\\notag
\leq C_{l}\delta_{0}(1+t)^{-3}\{\mathcal{Y}_{0}+(1+t)\mathcal{A}_{[l]}+\mathcal{W}_{[l]}\\\notag
+[1+\log(1+t)][\mathcal{W}^{T}_{[l]}+(1+t)^{-1}[1+\log(1+t)]\mathcal{W}^{Q}_{[l]}+\delta_{0}(1+t)^{-2}\mathcal{B}_{[l]}]\}
\end{align}
for $j+k=l$.

    Combining the results we obtain the following lemma:

$\textbf{Lemma 10.15}$ Let the hypotheses $\textbf{H0}$, $\textbf{H1}$, $\textbf{H2}^{\prime}$, and the estimate (6.177) hold. Let also the bootstrap 
assumptions $\slashed{\textbf{E}}_{[l_{*}+2]}$, $\slashed{\textbf{E}}^{Q}_{[l_{*}+1]}$, $\slashed{\textbf{E}}^{T}_{[l_{*}]}$, $\slashed{\textbf{X}}_{[(l+1)_{*}]}$
and $\slashed{\textbf{M}}_{[l_{*}+1]}$ hold, for some positive integer $l$. Then provided that $\delta_{0}$ is suitably small (depending on $l$), we have:
\begin{align*}
\max_{i_{1}...i_{l}}\|\sum_{k=0}^{l-1}\slashed{\mathcal{L}}_{R_{i_{l}}}...\slashed{\mathcal{L}}_{R_{i_{l-k+1}}}\leftexp{(i_{1}...i_{l-k})}{y}_{l-k}\|
_{L^{2}(\Sigma_{t}^{\epsilon_{0}})}\\\notag
\leq C_{l}\delta_{0}(1+t)^{-3}\{[1+\log(1+t)]\times(\mathcal{W}_{[l+2]}+\delta_{0}(1+t)^{-1}\mathcal{W}^{Q}_{[l+1]}\\\notag
+\mathcal{W}^{T}_{[l]}+(1+t)^{-1}[1+\log(1+t)]\mathcal{W}^{Q}_{[l]})+\mathcal{W}_{[l]}\\\notag
+\mathcal{Y}_{0}+(1+t)\mathcal{A}_{[l]}+(1+t)^{-1}\mathcal{B}_{[l]}\}
\end{align*}

    To estimate the $S_{t,u}$ 1-form $\leftexp{(i_{1}...i_{l})}{\dot{g}}_{l}$, given by (10.150), it remains for us to consider the $S_{t,u}$ 
1-form $\leftexp{(i_{1}...i_{l})}{w}_{l}$, defined by (10.151). A treatment similar 
to the above yields:

$\textbf{Lemma 10.16}$ Let hypotheses $\textbf{H0}$ and (6.177) hold. Let also the bootstrap assumptions $\slashed{\textbf{E}}_{[l_{*}+2]}$, 
$\slashed{\textbf{E}}^{Q}_{[l_{*}+1]}$, $\slashed{\textbf{E}}^{T}_{[l_{*}+1]}$, $\slashed{\textbf{X}}_{[(l+1)_{*}]}$ and $\slashed{\textbf{M}}_{[l_{*}+1]}$ hold,
for some positive integer $l$. Then, provided that $\delta_{0}$ is suitably small (depending on $l$), we have:
\begin{align*}
 \max_{i_{1}...i_{l}}\|\leftexp{(i_{1}...i_{l})}{w}_{l}\|_{L^{2}(\Sigma_{t}^{\epsilon_{0}})}\leq C_{l}\delta_{0}(1+t)^{-3}\cdot\\\notag
\{\mathcal{Y}_{0}+(1+t)\mathcal{A}_{[l]}+\mathcal{W}_{[l+1]}+[1+\log(1+t)]\mathcal{W}^{T}_{[l]}\\\notag
+(1+t)^{-1}[1+\log(1+t)]^{2}(\mathcal{W}^{Q}_{[l]}+\delta_{0}(1+t)^{-1}\mathcal{B}_{[l]})\}
\end{align*}
    Combining Lemma 10.13, Lemma 10.15 and Lemma 10.16, we arrive at the following proposition.

$\textbf{Proposition 10.4}$ Let the hypotheses $\textbf{H0}$, $\textbf{H1}$, $\textbf{H2}^{\prime}$ and (6.177) hold. Let also the bootstrap assumptions 
$\slashed{\textbf{E}}_{[l_{*}+2]}$,$\slashed{\textbf{E}}^{Q}_{[l_{*}+1]}$ and $\slashed{\textbf{E}}^{T}_{[l_{*}+1]}$, as well as $\slashed{\textbf{X}}_{[(l+1)_{*}]}$
and $\slashed{\textbf{M}}_{[l_{*}+1]}$ hold, for some positive integer $l$. Then if $\delta_{0}$ is suitably small (depending on $l$) we have:
\begin{align*}
 \max_{i_{1}...i_{l}}\|\leftexp{(i_{1}...i_{l})}{\dot{g}}_{l}\|_{L^{2}(\Sigma_{t}^{\epsilon_{0}})}\leq C_{l}(1+t)^{-3}\cdot\\\notag
\{[1+\log(1+t)][\mathcal{W}^{Q}+\delta_{0}(\mathcal{W}^{T}_{[l+1]}+\mathcal{W}_{[l+2]})]\\\notag
+\delta_{0}(\mathcal{Y}_{0}+(1+t)\mathcal{A}_{[l]}+(1+t)^{-1}\mathcal{B}_{[l+1]})\}
\end{align*}\vspace{7mm}

     Let us now consider the additional terms
\begin{align}
 2\mu\slashed{D}\hat{\chi}\cdot\leftexp{(i_{1}...i_{l})}{\hat{\chi}}_{l}+\mu\slashed{d}\leftexp{(i_{1}...i_{l})}{\dot{h}}_{l}
\end{align}
Using the same methods we establish:

$\textbf{Lemma 10.17}$ Let $\textbf{H0}$ and (6.177) hold. Let also the bootstrap assumptions $\slashed{\textbf{E}}_{[l_{*}+1]}$, 
$\slashed{\textbf{E}}^{Q}_{[l_{*}]}$ and $\slashed{\textbf{X}}_{[(l+1)_{*}]}$ hold for some positive integer $l$. Then we have:
\begin{align*}
 \max_{i_{1}...i_{l}}\|\mu\slashed{d}\leftexp{(i_{1}...i_{l})}{\dot{h}}_{l}\|_{L^{2}(\Sigma_{t}^{\epsilon_{0}})}\\\notag
\leq C_{l}\delta_{0}(1+t)^{-4}[1+\log(1+t)]^{2}\{\mathcal{Y}_{0}+(1+t)\mathcal{A}_{[l]}+\mathcal{W}_{[l+1]}\}
\end{align*}
provided that $\delta_{0}$ is suitably small (depending on $l$).\vspace{7mm}

     Finally the first term in (10.170) can be estimated directly:
\begin{align}
 \max_{i_{1}...i_{l}}\|\mu\slashed{D}\hat{\chi}\cdot\leftexp{(i_{1}...i_{l})}{\hat{\chi}}_{l}\|_{L^{2}(\Sigma_{t}^{\epsilon_{0}})}\\\notag
\leq C_{l}\delta_{0}(1+t)^{-4}[1+\log(1+t)]^{2}\{\mathcal{Y}_{0}+(1+t)\mathcal{A}_{[l]}+\mathcal{W}_{[l]}\}
\end{align}

     Proposition 10.4 together with Lemma 10.17 and (10.170) gives us an appropriate bound for $\|\leftexp{(i_{1}...i_{l})}{\ddot{g}}_{l}\|
_{L^{2}(\Sigma_{t}^{\epsilon_{0}})}$. To get the estimate for $\leftexp{(i_{1}...i_{l})}{Q}_{l}(t)$, it remains for us to bound the quantities
$\|\slashed{d}\leftexp{(i_{1}...i_{l})}{\check{f}}_{l}\|_{L^{2}(\Sigma_{t}^{\epsilon_{0}})}$ and $\|\leftexp{(i_{1}...i_{l})}{z}_{l}\|
_{L^{2}(\Sigma_{t}^{\epsilon_{0}})}$. From (10.168) we have:
\begin{align}
 \max_{i_{1}...i_{l}}\|\slashed{d}\leftexp{(i_{1}...i_{l})}{\check{f}}_{l}\|_{L^{2}(\Sigma_{t}^{\epsilon_{0}})}\\\notag
\leq C_{l}(1+t)^{-1}\{\mathcal{W}_{[l+1]}^{T}+\delta_{0}(1+t)^{-1}[\mathcal{W}_{[l+1]}+(1+t)^{-1}\mathcal{B}_{[l+1]}]
+(1+t)^{-1}[1+\log(1+t)]\mathcal{W}^{Q}_{[l+1]}\}
\end{align}

     The non-negative function $\leftexp{(i_{1}...i_{l})}{z}_{l}$ is given by (8.392). Using the same methods we deduce: 
\begin{align}
 \max_{i_{1}...i_{l}}\|\leftexp{(i_{1}...i_{l})}{z}_{l}\|_{L^{2}(\Sigma_{t}^{\epsilon_{0}})}
\leq C_{l}\delta_{0}(1+t)^{-2}[1+\log(1+t)]\cdot\\\notag
\{\mathcal{W}_{[l+1]}^{T}+(1+t)^{-1}[1+\log(1+t)](\mathcal{W}^{Q}_{[l+1]}+\delta_{0}(1+t)^{-1}\mathcal{W}_{[l+2]})\\\notag
+(1+t)^{-1}[1+\log(1+t)](\mathcal{Y}_{0}+(1+t)\mathcal{A}_{[l]})+\delta_{0}(1+t)^{-2}\mathcal{B}_{[l+1]}\}
\end{align}
We now can write down the estimate for $\leftexp{(i_{1}...i_{l})}{Q}_{l}$:

$\textbf{Proposition 10.5}$ Let the hypotheses $\textbf{H0}$, $\textbf{H1}$, $\textbf{H2}^{\prime}$ and the estimate (6.177) hold. Let also the bootstrap assumptions 
$\slashed{\textbf{E}}_{[l_{*}+2]}$, $\slashed{\textbf{E}}^{Q}_{[l_{*}+1]}$ and $\slashed{\textbf{E}}^{T}_{[l_{*}+1]}$, as well as $\slashed{\textbf{X}}_{[(l+1)_{*}]}$
and $\slashed{\textbf{M}}_{[l_{*}+1]}$ hold, for some positive integer $l$. Then we have:
\begin{align*}
 \max_{i_{1}...i_{l}}\leftexp{(i_{1}...i_{l})}{Q}_{l}\leq C_{l}(1+t)^{-4}\cdot\\\notag
\{[1+\log(1+t)][\mathcal{W}^{Q}_{[l+1]}+\delta_{0}(\mathcal{W}^{T}_{[l+1]}+\mathcal{W}_{[l+2]})]\\\notag
+\delta_{0}(\mathcal{Y}_{0}+(1+t)\mathcal{A}_{[l]}+(1+t)^{-1}\mathcal{B}_{[l+1]})\}
\end{align*}
provided that $\delta_{0}$ is suitably small (depending on $l$).

\subsection{Estimates for $P_{l}$}   
       We now consider the principal term in the defining expression (8.403) for $\leftexp{(i_{1}...i_{l})}{B}_{l}$, namely, the term:
\begin{align*}
 C(1+t)^{-2}\{\leftexp{(i_{1}...i_{l})}{\bar{P}}^{(0)}_{l,a}(t)+(1+t)^{-1/2}\leftexp{(i_{1}...i_{l})}{\bar{P}}^{(1)}_{l,a}(t)\}\bar{\mu}_{m}^{-a}(t)
\end{align*}
Here, the quantities $\leftexp{(i_{1}...i_{l})}{\bar{P}}^{(0)}_{l,a}$, $\leftexp{(i_{1}...i_{l})}{\bar{P}}^{(1)}_{l,a}$ are defined by
(8.348) and (8.349). Recall the definition:
\begin{equation}
 \leftexp{(i_{1}...i_{l})}{P}_{l}(t)=(1+t)^{2}\||\slashed{d}\leftexp{(i_{1}...i_{l})}{\check{f}}_{l}(t)|\|_{L^{2}([0,\epsilon_{0}]\times S^{2})}
\end{equation}
Here we need a more precise estimate for $\|\slashed{d}\leftexp{(i_{1}...i_{l})}{\check{f}}_{l}\|_{L^{2}(\Sigma_{t}^{\epsilon_{0}})}$.

By (8.333) and $\textbf{H0}$ we have:
\begin{equation}
 \leftexp{(i_{1}...i_{l})}{P}_{l}\leq C\sum_{j}\|R_{j}\leftexp{(i_{1}...i_{l})}{\check{f}}_{l}\|_{L^{2}(\Sigma_{t}^{\epsilon_{0}})}
=C\sum_{j}\|R_{j}R_{i_{l}}...R_{i_{1}}\check{f}\|_{L^{2}(\Sigma_{t}^{\epsilon_{0}})}
\end{equation}
Here, we must recall the definition of $\check{f}$ in chapter 8:
\begin{equation}
 \check{f}=\mu f=-\frac{1}{2}\frac{dH}{dh}\tau_{\underline{L}}
\end{equation}
We can use the definition of the acoustical metric and (8.24) to write the above as follows:
\begin{equation}
 \check{f}=-2m-\frac{1}{2}\alpha^{-1}\kappa\frac{dH}{dh}Lh:=-2m-\mu\bar{e}
\end{equation}
where we have used the fact that:
\begin{align*}
 \alpha=\eta,\quad \partial_{t}-\psi_{i}\partial_{i}=B=L+\alpha\kappa^{-1}T
\end{align*}
We can express:
\begin{equation}
 m=m^{\alpha}_{T}(T\psi_{\alpha})
\end{equation}
where:
\begin{equation}
 m^{0}_{T}=\frac{1}{2}\frac{dH}{dh},\quad m^{i}_{T}=-\frac{1}{2}\frac{dH}{dh}\psi_{i}
\end{equation}
and
\begin{equation}
 \bar{e}=\bar{e}^{\alpha}_{L}(L\psi_{\alpha})
\end{equation}
where
\begin{equation}
 \bar{e}^{0}_{L}=\frac{1}{2}\alpha^{-2}\frac{dH}{dh},\quad \bar{e}^{i}_{L}=-\frac{1}{2}\alpha^{-2}\frac{dH}{dh}\psi_{i}
\end{equation}
Then
\begin{align}
 R_{j}R_{i_{l}}...R_{i_{1}}\check{f}=-2m^{\alpha}_{T}(R_{j}R_{i_{l}}...R_{i_{1}}T\psi_{\alpha})-\mu\bar{e}^{\alpha}_{L}
(R_{j}R_{i_{l}}...R_{i_{1}}L\psi_{\alpha})+\leftexp{(i_{1}...i_{l}j)}{n}_{l+1}
\end{align}
Here, $\leftexp{(i_{1}...i_{l}j)}{n}_{l+1}$ is the lower order term: 
\begin{align}
 \leftexp{(i_{1}...i_{l}j)}{n}_{l+1}=-\sum_{|s_{1}|\geq 1}((R)^{s_{1}}\mu)((R)^{s_{2}}\bar{e})\\\notag
-2\sum_{|s_{1}|\geq 1}((R)^{s_{1}}m^{\alpha}_{T})((R)^{s_{2}}T\psi_{\alpha})\\\notag
-\sum_{|s_{1}|\geq 1}\mu((R)^{s_{1}}\bar{e}^{\alpha}_{L})((R)^{s_{2}}L\psi_{\alpha})
\end{align}
We have:
\begin{equation}
 |m^{0}_{T}-\frac{1}{2}\ell|\leq C\delta_{0}(1+t)^{-1},\quad |m^{i}_{T}|\leq C\delta_{0}(1+t)^{-1}
\end{equation}
and also
\begin{equation}
 |\bar{e}^{0}_{L}|\sim|\ell| 
\end{equation}
From (10.182)-(10.184) we have:
\begin{align}
 \|R_{j}R_{i_{l}}...R_{i_{1}}\check{f}\|_{L^{2}(\Sigma_{t}^{\epsilon_{0}})}\leq (|\ell|+C\delta_{0}(1+t)^{-1})
\sum_{\alpha}\|R_{j}R_{i_{l}}...R_{i_{1}}T\psi_{\alpha}\|_{L^{2}(\Sigma_{t}^{\epsilon_{0}})}\\\notag
+C\sum_{\alpha}\|\mu R_{j}R_{i_{l}}...R_{i_{1}}L\psi_{\alpha}\|_{L^{2}(\Sigma_{t}^{\epsilon_{0}})}+\|\leftexp{(i_{1}...i_{l}j)}{n}_{l+1}\|
_{L^{2}(\Sigma_{t}^{\epsilon_{0}})}
\end{align}
where $C$ is a constant which is $independent$ of $l$. Under the assumptions of Proposition 10.5 we readily deduce:
\begin{align}
 \max_{i_{1}...i_{l}j}\|\leftexp{(i_{1}...i_{l}j)}{n}_{l+1}\|_{L^{2}(\Sigma_{t}^{\epsilon_{0}})}\\\notag
\leq C_{l}\delta_{0}(1+t)^{-1}\{\mathcal{W}^{T}_{[l]}+[1+\log(1+t)]\mathcal{W}^{Q}_{[l]}+\mathcal{W}_{[l+1]}+(1+t)^{-1}\mathcal{B}_{[l+1]}\}
\end{align}
To proceed, we need the following lemma, which can be proved in a similar way as Lemma 8.2.

$\textbf{Lemma 10.18}$ Let $Y$ be an arbitrary $S_{t,u}$-tangential vectorfield on the spacetime domain $W^{*}_{\epsilon_{0}}$. We have:
\begin{align*}
 [T,Y]=\leftexp{(Y)}{\Theta}
\end{align*}
where $\leftexp{(Y)}{\Theta}$ is an $S_{t,u}$-tangential vectorfield, associated to $Y$, and defined by the condition that for any vector $V\in TW^{*}_{\epsilon_{0}}$:
\begin{align*}
 g(\leftexp{(Y)}{\Theta},V)=\leftexp{(Y)}{\pi}(T,\Pi V)
\end{align*}
In terms of the $(L,T,X_{1},X_{2})$ frame,
\begin{align*}
 \leftexp{(Y)}{\Theta}=\leftexp{(Y)}{\Theta}^{A}X_{A},\quad \leftexp{(Y)}{\Theta}^{A}=\leftexp{(Y)}{\pi}_{TB}(\slashed{g}^{-1})^{AB}
\end{align*}\vspace{7mm}

    We can express, in rectangular coordinates:
\begin{align}
 [T,R_{i}]=\{T(R_{i})^{j}-R_{i}T^{j}\}\partial_{j}
\end{align}
Since $T=\kappa\hat{T}$, we have:
\begin{align*}
 R_{i}T^{j}=(R_{i}\kappa)\hat{T}^{j}+\kappa R_{i}\hat{T}^{j}
\end{align*}
Using (10.8) we have:
\begin{align}
 R_{i}T^{j}=(R_{i}\kappa)\hat{T}^{j}+\kappa\slashed{q}_{i}\cdot\slashed{d}x^{j}
\end{align}
Projecting to $S_{t,u}$:
\begin{align}
 \Pi(R_{i}\hat{T}^{j}\partial_{j})=\kappa\slashed{q}_{i}
\end{align}

On the other hand, from (10.21), (10.23),
\begin{align*}
 (R_{i})^{j}=\epsilon_{imj}x^{m}-\lambda_{i}\hat{T}^{j}
\end{align*}
hence:
\begin{align}
 T(R_{i})^{j}=\epsilon_{imj}T^{m}-(T\lambda_{i})\hat{T}^{j}-\lambda_{i}T(\hat{T}^{j})
\end{align}
According to (3.192):
\begin{align}
 T(\hat{T}^{j})=q^{j}_{T}
\end{align}
Here $q_{T}$ is an $S_{t,u}$-tangential vectorfield given by:
\begin{align}
 q_{T}=(q_{T})_{b}\cdot \slashed{g}^{-1}
\end{align}
where
\begin{align}
 (q_{T})_{b}=-\slashed{d}\kappa
\end{align}
Substituting (10.193) in (10.192) we obtain:
\begin{align}
 T(R_{i})^{j}\partial_{j}=\kappa v_{i}-(T\lambda_{i})\hat{T}-\lambda_{i}q_{T}
\end{align}
where $v_{i}$ is defined by (10.47). Projecting to $S_{t,u}$:
\begin{align}
 \Pi(T(R_{i})^{j}\partial_{j})=\kappa\Pi v_{i}-\lambda_{i}q_{T}
\end{align}
So we obtain:
\begin{align}
 [T,R_{i}]=\kappa\Pi v_{i}-\lambda_{i}q_{T}-\kappa\slashed{q}_{i}
\end{align}
Defining the vectorfields $v^{\prime}_{i}$ as in (10.50), equation (10.52) holds, hence:
\begin{align}
 \Pi v_{i}=w_{i}-(1-u+t)^{-1}R_{i}
\end{align}
Therefore by (10.43) and (10.14), we obtain:
\begin{align}
 [T,R_{i}]=-\kappa\tilde{\slashed{q}}^{\prime}-\lambda_{i}q_{T}
\end{align}
We have thus proved the following lemma.

$\textbf{Lemma 10.19}$ We have:
\begin{align*}
 [T,R_{i}]=\leftexp{(R_{i})}{\Theta}=-\kappa\tilde{\slashed{q}}^{\prime}-\lambda_{i}q_{T}
\end{align*}\vspace{7mm}
    Now, under the assumptions of Proposition 10.1 with $l_{*}+1$ in the role of $l$ augmented by assumption 
$\slashed{\textbf{M}}_{[l_{*}+1]}$, we have, by (10.194) and (10.195),
\begin{align}
 \max_{i_{1}...i_{k}}\|\slashed{\mathcal{L}}_{R_{i_{k}}}...\slashed{\mathcal{L}}_{R_{i_{1}}}q_{T}\|_{L^{\infty}(\Sigma_{t}^{\epsilon_{0}})}
\leq C_{l}\delta_{0}(1+t)^{-1}[1+\log(1+t)]
\end{align}
for $k\leq l_{*}$.
Also by Corollary 10.2.d, with $l$ replaced by $l+1$,
we have:
\begin{align}
 \max_{i_{1}...i_{k}}\|\slashed{\mathcal{L}}_{R_{i_{k}}}...\slashed{\mathcal{L}}_{R_{i_{1}}}q_{T}\|_{L^{2}(\Sigma_{t}^{\epsilon_{0}})}\\\notag
\leq C_{l}\{(1+t)^{-1}\mathcal{B}_{[l+1]}+\delta_{0}(1+t)^{-1}[1+\log(1+t)][\mathcal{Y}_{0}+(1+t)\mathcal{A}_{[l-1]}+\mathcal{W}_{[l]}]\}
\end{align}
By (10.101) and (10.41)-(10.42), Proposition 10.1, with $l_{*}+1$ in the role of $l$ we have:
\begin{align}
 \max_{i;i_{1}...i_{k}}\|\slashed{\mathcal{L}}_{R_{i_{k}}}...\slashed{\mathcal{L}}_{R_{i_{1}}}\tilde{\slashed{q}}^{\prime}_{i}\|_{L^{\infty}(\Sigma_{t}^{\epsilon_{0}})}
\leq C_{l}\delta_{0}(1+t)^{-1}[1+\log(1+t)]
\end{align}
for $k\leq l_{*}$, and by (10.120) and (10.123) with $l+1$ in the role of $l$:
\begin{align}
 \max_{i;i_{1}...i_{k}}\|\slashed{\mathcal{L}}_{R_{i_{k}}}...\slashed{\mathcal{L}}_{R_{i_{1}}}\tilde{\slashed{q}}^{\prime}_{i}\|_{L^{2}(\Sigma_{t}^{\epsilon_{0}})}
\\\notag
\leq C_{l}\{\mathcal{Y}_{0}+(1+t)\mathcal{A}_{[l]}+\mathcal{W}_{[l+1]}\}
\end{align}
for $k\leq l$.

The above estimates together with Corollary 10.1.a with $l_{*}$ in the role of $l$ and Corollary 10.2.a yield the following lemma.

$\textbf{Lemma 10.20}$ Under the assumptions of Proposition 10.1 with $l_{*}+1$ in the role of $l$ and $\slashed{\textbf{M}}_{[l_{*}+1]}$, we have:
\begin{align*}
 \max_{i;i_{1}...i_{k}}\|\slashed{\mathcal{L}}_{R_{i_{k}}}...\slashed{\mathcal{L}}_{R_{i_{1}}}\leftexp{(R_{i})}{\Theta}\|_{L^{\infty}(\Sigma_{t}^{\epsilon_{0}})}
\leq C_{l}\delta_{0}(1+t)^{-1}[1+\log(1+t)]^{2}
\end{align*}
for $k\leq l_{*}$, and:
\begin{align*}
 \max_{i;i_{1}...i_{k}}\|\slashed{\mathcal{L}}_{R_{i_{k}}}...\slashed{\mathcal{L}}_{R_{i_{1}}}\leftexp{(R_{i})}{\Theta}\|_{L^{2}(\Sigma_{t}^{\epsilon_{0}})}\\\notag
\leq C_{l}[1+\log(1+t)]\{\mathcal{Y}_{0}+(1+t)\mathcal{A}_{[l]}+\delta_{0}(1+t)^{-1}\mathcal{B}_{[l+1]}+\mathcal{W}_{[l+1]}\}
\end{align*}
for $k\leq l$.

Provided that $\delta_{0}$ is suitably small (depending on $l$).\vspace{7mm}

     Consider now the commutator
\begin{align}
 TR_{i_{l}}...R_{i_{1}}-R_{i_{l}}...R_{i_{1}}T=\leftexp{(i_{1}...i_{l})}{A}^{T}_{l}
\end{align}
By Lemma 10.19,
\begin{align}
 \leftexp{(i_{1}...i_{l})}{A}^{T}_{l}=R_{i_{l}}\leftexp{(i_{1}...i_{l-1})}{A}^{T}_{l-1}+\leftexp{(R_{i_{l}})}{\Theta}R_{i_{l-1}}...R_{i_{1}}
\end{align}
Applying Proposition 8.2 to this recursion yields:
\begin{align}
 \leftexp{(i_{1}...i_{l})}{A}^{T}_{l}=\sum_{k=0}^{l-1}R_{i_{l}}...R_{i_{l-k+1}}\leftexp{(R_{i_{l-k}})}{\Theta}R_{i_{l-k-1}}...R_{i_{1}}
\end{align}
So we have:
\begin{align}
 R_{i_{l+1}}R_{i_{l}}...R_{i_{1}}T\psi_{\alpha}=TR_{i_{l+1}}R_{i_{l}}...R_{i_{1}}\psi_{\alpha}\\\notag
-\sum_{k=0}^{l}R_{i_{l+1}}...R_{i_{l-k+2}}\leftexp{(R_{i_{l-k+1}})}{\Theta}R_{i_{l-k}}...R_{i_{1}}\psi_{\alpha}
\end{align}
We shall now use the following lemma, which is proved in a similar way as Lemma 8.5

$\textbf{Lemma 10.21}$ For any $S_{t,u}$-tangential vectorfield $X$ and any function $\phi$ and positive integer $l$ we have:
\begin{align*}
 [X,R_{i_{l}}...R_{i_{1}}]\phi=-\sum_{j=1}^{l}\sum_{k_{1}<...<k_{j}=1}^{l}\leftexp{(i_{k_{1}}...i_{k_{j}})}{Y}R_{i_{l}}
\overset{>i_{k_{j}}...i_{k_{1}}<}{...}R_{i_{1}}\phi
\end{align*}
where $\leftexp{(i_{k_{1}}...i_{k_{j}})}{Y}$ is the $S_{t,u}$-tangential vectorfield:
\begin{align*}
 \leftexp{(i_{k_{1}}...i_{k_{j}})}{Y}=\slashed{\mathcal{L}}_{R_{i_{k_{j}}}}...\slashed{\mathcal{L}}_{R_{i_{k_{1}}}}X
\end{align*}\vspace{7mm}

    Consider then an arbitrary term of the sum on the right in (10.207), corresponding to some $k\in\{0,...,l\}$. Setting
\begin{align}
 X=\leftexp{(R_{i_{l-k+1}})}{\Theta},\quad \phi=R_{i_{l-k}}...R_{i_{1}}\psi_{\alpha}
\end{align}
we write this term in the form, using Lemma 10.21,
\begin{align}
 R_{i_{l+1}}...R_{i_{l-k+2}}X\phi=X\cdot\slashed{d}(R_{i_{l+1}}...R_{i_{l-k+2}}\phi)
+\sum_{|s_{1}|>0}(\slashed{\mathcal{L}}_{R})^{s_{1}}X\cdot\slashed{d}(R)^{s_{2}}\phi
\end{align}
By Lemma 10.20, the first term on the right in (10.209) is bounded in $L^{2}(\Sigma_{t}^{\epsilon_{0}})$ by:
\begin{align}
 C_{l}\delta_{0}(1+t)^{-2}[1+\log(1+t)]^{2}\mathcal{W}_{[l+1]}
\end{align}
A term in the sum in (10.209) is the product of an angular derivative of $X$ of order $|s_{1}|$ with an angular derivative of $\psi_{\alpha}$ of order 
$|s_{2}|+l-k+1$, and we have $|s_{1}|+|s_{2}|=k$. Thus the sum in (10.209) contains angular derivatives of $X$ of order at most $k\leq l$ and angular 
derivatives of $\psi_{\alpha}$ of order at most $l$ (for, $|s_{1}|\geq 1$).

     We have the following two cases.

     Case 1:$|s_{1}|\leq l_{*}$  and:  Case 2:$|s_{1}|\geq l_{*}+1$

     In Case 1 we use the first part of Lemma 10.20 to place the first factor in $L^{\infty}$, obtaining an $L^{2}$ bound by:
\begin{align}
 C_{l}\delta_{0}(1+t)^{-2}[1+\log(1+t)]^{2}\mathcal{W}_{[l]}
\end{align}
     In Case 2 we have $|s_{2}|\leq k-l_{*}-1$, hence:
\begin{align*}
 |s_{2}|+l-k+1\leq l-l_{*}\leq l_{*}+1
\end{align*}
and the bootstrap assumption $\slashed{\textbf{E}}_{[l_{*}+1]}$ allows us to place the second factor in $L^{\infty}$. Using also the second part of 
Lemma 10.20 in placing the first factor in $L^{2}$, we obtain an $L^{2}$ bound by:
\begin{align}
C_{l}(1+t)^{-2}[1+\log(1+t)]\{\mathcal{Y}_{0}+(1+t)\mathcal{A}_{[l]}+\delta_{0}(1+t)^{-1}\mathcal{B}_{[l+1]}+\mathcal{W}_{[l+1]}\}
\end{align}
 Combining the above results (10.210)-(10.212) we conclude that:
\begin{align}
 \|\sum_{k=0}^{l}R_{i_{l+1}}...R_{i_{l-k+2}}\leftexp{(R_{i_{l-k+1}})}{\Theta}R_{i_{l-k}}...R_{i_{1}}\psi_{\alpha}\|_{L^{2}(\Sigma_{t}^{\epsilon_{0}})}\\\notag
\leq C_{l}\delta_{0}(1+t)^{-2}[1+\log(1+t)]\{\mathcal{Y}_{0}+(1+t)\mathcal{A}_{[l]}+\delta_{0}(1+t)^{-1}\mathcal{B}_{[l+1]}\\\notag
+[1+\log(1+t)]\mathcal{W}_{[l+1]}\}
\end{align}
for $k=0,...,l$.

    A commutation relation similar to (10.207) with $\leftexp{(R_{i})}{Z}$ in the role of $\leftexp{(R_{i})}{\Theta}$ holds for $L$ in the role of $T$, that is:
\begin{align}
 R_{i_{l+1}}R_{i_{l}}...R_{i_{1}}L\psi_{\alpha}=LR_{i_{l+1}}R_{i_{l}}...R_{i_{1}}\psi_{\alpha}\\\notag
-\sum_{k=0}^{l}R_{i_{l+1}}...R_{i_{l-k+2}}\leftexp{(R_{i_{l-k+1}})}{Z}R_{i_{l-k}}...R_{i_{1}}\psi_{\alpha}
\end{align}
Using Corollary 10.1.i with $l_{*}+1$ in the role of $l$ and Corollary 10.2.i with $l+1$ in the role of $l$, in place of the estimates of Lemma 10.20, we deduce that:
\begin{align}
 \|\mu\sum_{k=0}^{l}R_{i_{l+1}}...R_{i_{l-k+2}}\leftexp{(R_{i_{l-k+1}})}{Z}R_{i_{l-k}}...R_{i_{1}}\psi_{\alpha}\|_{L^{2}(\Sigma_{t}^{\epsilon_{0}})}\\\notag
\leq C_{l}\delta_{0}(1+t)^{-2}[1+\log(1+t)]\{\mathcal{Y}_{0}+(1+t)\mathcal{A}_{[l]}+[1+\log(1+t)]\mathcal{W}_{[l+1]}\}
\end{align}
for $k=0,...,l$.

     Now we have, from Chapter 5, the following hold for any variation $\psi$:
\begin{align}
 \|T\psi\|_{L^{2}(\Sigma_{t}^{\epsilon_{0}})}\leq C\sqrt{\mathcal{E}_{0}[\psi]}\\\notag
\|\sqrt{\mu}(L\psi+\nu \psi)\|_{L^{2}(\Sigma_{t}^{\epsilon_{0}})}\leq C(1+t)^{-1}\sqrt{\mathcal{E}^{\prime}_{1}[\psi]}
\end{align}
while:
\begin{align}
 \|\mu\nu R_{j}R_{i_{l}}...R_{i_{1}}\psi_{\alpha}\|_{L^{2}(\Sigma_{t}^{\epsilon_{0}})}\leq C(1+t)^{-1}[1+\log(1+t)]\mathcal{W}_{[l+1]}
\end{align}
Thus, up to the estimates for the commutators, the first term on the right in (10.186) is bounded by:
\begin{align}
 (|\ell|+C\delta_{0}(1+t)^{-1})\sqrt{\sum_{\alpha}\mathcal{E}_{0}[R_{j}R_{i_{l}}...R_{i_{1}}\psi_{\alpha}]}
\end{align}
and the second term on the right in (10.186) is bounded by:
\begin{align}
 C(1+t)^{-1}[1+\log(1+t)]^{1/2}\sqrt{\sum_{\alpha}\mathcal{E}^{\prime}_{1}[R_{j}R_{i_{l}}...R_{i_{1}}\psi_{\alpha}]}+C(1+t)^{-1}[1+\log(1+t)]\mathcal{W}_{[l+1]}
\end{align}
 Actually, taking into account the commutator estimates (10.214), (10.216), we have the following proposition.

$\textbf{Proposition 10.6}$ Under the assumptions of Proposition 10.5, we have:
\begin{align*}
 \leftexp{(i_{1}...i_{l})}{P}_{l}\leq \leftexp{(i_{1}...i_{l})}{P}^{(0)}_{l}+\leftexp{(i_{1}...i_{l})}{P}^{(1)}_{l}
\end{align*}
where:
\begin{align*}
 \leftexp{(i_{1}...i_{l})}{P}_{l}^{(0)}=|\ell|\sqrt{\sum_{j,\alpha}\mathcal{E}_{0}[R_{j}R_{i_{l}}...R_{i_{1}}\psi_{\alpha}]}
\end{align*}
and:
\begin{align*}
 \leftexp{(i_{1}...i_{l})}{P}^{(1)}_{l}=C\delta_{0}(1+t)^{-1}\sqrt{\sum_{j,\alpha}\mathcal{E}_{0}[R_{j}R_{i_{l}}...R_{i_{1}}\psi_{\alpha}]}\\
+C(1+t)^{-1}[1+\log(1+t)]^{1/2}\sqrt{\sum_{j,\alpha}\mathcal{E}^{\prime}_{1}[R_{j}R_{i_{l}}...R_{i_{1}}\psi_{\alpha}]}\\
+C(1+t)^{-1}[1+\log(1+t)]\mathcal{W}_{[l+1]}+C_{l}\delta_{0}(1+t)^{-1}\{\mathcal{W}^{T}_{[l]}+[1+\log(1+t)]\mathcal{W}^{Q}_{[l]}\\
+(1+t)^{-1}[[1+\log(1+t)](\mathcal{Y}_{0}+(1+t)\mathcal{A}_{[l]})+\mathcal{B}_{[l+1]}]\}
\end{align*}

\chapter{Control of the Spatial Derivatives of the First Derivatives of the $x^{i}$.\\
Assumptions and Estimates in Regard to $\mu$}

\section{Estimates for $T\hat{T}^{i}$}
\subsection{Basic Lemmas}
The first part of the present chapter is concerned with the derivation of estimates for the spatial derivatives of the $\hat{T}x^{i}=\hat{T}^{i}$, of which at least
one is a $T$-derivative. By spatial derivatives we mean the derivatives with respect to $T$ and the rotation fields $R_{j}$ : $j=1,2,3$. Combining with the estimates of 
Chapter 10 we then obtain estimates for all spatial derivatives of the first derivatives of the $x^{i}$ with respect to all five commutation fields. The derivation of 
the estimates of this chapter is based on bootstrap assumptions in regard to $\mu$, in addition to the bootstrap assumptions of Chapter 10 in regard to $\chi$. The estimates 
are then used to obtain estimates for the spatial derivatives of the deformation tensors of the commutation fields.

     We begin by defining certain norms. In the following $\xi$ is an arbitrary $S_{t,u}$ tensorfield defined on $\Sigma_{t}^{\epsilon_{0}}$. We first define $L^{\infty}$
and $L^{2}$ norms which consider only angular derivatives. Given a non-negative integer $l$ we define:
\begin{align}
 \|\xi\|_{\infty,[l],\Sigma_{t}^{\epsilon_{0}}}=\max_{0\leq n\leq l}\max_{i_{1}...i_{n}}\|\slashed{\mathcal{L}}_{R_{i_{n}}}...\slashed{\mathcal{L}}_{R_{i_{1}}}\xi\|
_{L^{\infty}(\Sigma_{t}^{\epsilon_{0}})}\\
\|\xi\|_{2,[l],\Sigma_{t}^{\epsilon_{0}}}=\sum_{n=0}^{l}\max_{i_{1}...i_{n}}\|\slashed{\mathcal{L}}_{R_{i_{n}}}...\slashed{\mathcal{L}}_{R_{i_{1}}}\xi\|_{L^{2}(\Sigma_{t}
^{\epsilon_{0}})}
\end{align}
 We then define $L^{\infty}$ and $L^{2}$ norms which take into consideration all spatial derivatives.

     Given a pair of non-negative integers $k,l$ with $k\leq l$, we define:
\begin{align}
 \|\xi\|_{\infty,[k,l],\Sigma_{t}^{\epsilon_{0}}}=\max_{0\leq m\leq k}\|(\slashed{\mathcal{L}}_{T})^{m}\xi\|_{\infty,[l-m],\Sigma_{t}^{\epsilon_{0}}}\\
\|\xi\|_{2,[k,l],\Sigma_{t}^{\epsilon_{0}}}=\sum_{m=0}^{k}\|(\slashed{\mathcal{L}}_{T})^{m}\xi\|_{2,[l-m],\Sigma_{t}^{\epsilon_{0}}}
\end{align}
Also, we shall denote:
\begin{align}
 \|\xi\|_{\infty,\{l\},\Sigma_{t}^{\epsilon_{0}}}=\|\xi\|_{\infty,[l,l],\Sigma_{t}^{\epsilon_{0}}},\quad 
\|\xi\|_{2,\{l\},\Sigma_{t}^{\epsilon_{0}}}=\|\xi\|_{2,[l,l],\Sigma_{t}^{\epsilon_{0}}}
\end{align}
According to the above definitions $\|\xi\|_{\infty,[k,l],\Sigma_{t}^{\epsilon_{0}}}$ is the maximum of the quantities:
\begin{align}
 \max_{i_{1}...i_{n}}\|\slashed{\mathcal{L}}_{R_{i_{n}}}...\slashed{\mathcal{L}}_{R_{i_{1}}}(\slashed{\mathcal{L}}_{T})^{m}\xi\|_{L^{\infty}(\Sigma_{t}^{\epsilon_{0}})}
\end{align}
over
\begin{align}
 \{(m,n)\quad:\quad 0\leq n\leq l-m, 0\leq m\leq k\}
\end{align}
Similarly for $\|\xi\|_{2,[k,l],\Sigma_{t}^{\epsilon_{0}}}$.

Also, $\|\xi\|_{\infty,\{l\},\Sigma_{t}^{\epsilon_{0}}}$ is the maximum of the quantities (11.6) over
\begin{align}
 \{(m,n)\quad:\quad m,n\geq 0, m+n\leq l\}
\end{align}
Similarly for $\|\xi\|_{2,\{l\},\Sigma_{t}^{\epsilon_{0}}}$.

      We now introduce some bootstrap assumptions. In the following, $C$ is a constant independent of $s$, and $t\in[0,s]$.
\begin{align*}
 \textbf{E}_{m,n}\quad:\quad \max_{\alpha;i_{1}...i_{n}}\|R_{i_{n}}...R_{i_{1}}(T)^{m}\psi_{\alpha}\|_{L^{\infty}(\Sigma_{t}^{\epsilon_{0}})}\leq C\delta_{0}(1+t)^{-1}
\end{align*}
In particular:
\begin{align*}
 \textbf{E}_{0,0}\quad:\quad\max_{\alpha}\|\psi_{\alpha}\|_{L^{\infty}(\Sigma_{t}^{\epsilon_{0}})}
\leq C\delta_{0}(1+t)^{-1}
\end{align*}
We denote by $\textbf{E}_{\{l\}}$ the conjunction of assumptions $\textbf{E}_{m,n}$ corresponding to (11.8). The constant $C$ depends on $l$ only. Note that 
$\textbf{E}_{0,n}$ coincide with $\slashed{\textbf{E}}_{n}$ and that $\textbf{E}_{\{l\}}$ contain $\slashed{\textbf{E}}_{[l]}$.

      Given non-negative integers $m,n$ not both $0$, we denote by $\mathcal{W}_{m,n}$ the quantity:
\begin{align}
 \mathcal{W}_{m,n}=\max_{\alpha;i_{1}...i_{n}}\|R_{i_{n}}...R_{i_{1}}(T)^{m}\psi_{\alpha}\|_{L^{2}(\Sigma_{t}^{\epsilon_{0}})}
\end{align}
In particular we denote by $\mathcal{W}_{0,0}$ the quantity:
\begin{align}
 \mathcal{W}_{0,0}=\max_{\alpha}\|\psi_{\alpha}\|_{L^{2}(\Sigma_{t}^{\epsilon_{0}})}
\end{align}
We then denote by $\mathcal{W}_{\{l\}}$ the sum of the quantities $\mathcal{W}_{m,n}$ corresponding to (11.8). Note that $\mathcal{W}_{0,n}$ coincide with 
$\mathcal{W}_{n}$, and that $\mathcal{W}_{\{l\}}$ dominate $\mathcal{W}_{[l]}$. Moreover, we have:
\begin{align}
 \max_{\alpha}\|\psi_{\alpha}\|_{2,\{l\},\Sigma_{t}^{\epsilon_{0}}}\leq \mathcal{W}_{\{l\}}
\end{align}

The following lemma, which can be proved by a direct calculation, extends Lemma 10.1 to all spatial derivatives.

$\textbf{Lemma 11.1}$ Let $G$ be a smooth function of the $(\psi_{\alpha}\quad :\quad \alpha=0,1,2,3)$ defined in a neighborhood of $(0,0,0,0)$, and let $G_{0}$ be the constant:
\begin{align*}
 G_{0}=G(0,0,0,0)
\end{align*}
Suppose that the bootstrap assumption $\textbf{E}_{\{l_{*}\}}$ holds for some positive integer $l$. Then, if $\delta_{0}$ is suitably small (depending on $l$), 
we have:
\begin{align*}
 \|G-G_{0}\|_{2,\{l\},\Sigma_{t}^{\epsilon_{0}}}\leq C\mathcal{W}_{\{l\}}
\end{align*}
where $C$ is a constant which is independent of $l$.

     We continue with a product lemma for the norms (11.3) and (11.4), which shall be applied repeatedly in the sequel.

$\textbf{Lemma 11.2}$

    $\textbf{a}$. Let $\xi_{1},...,\xi_{N}$ be arbitrary $S_{t,u}$ tensorfields, defined on $\Sigma_{t}^{\epsilon_{0}}$ and let us denote by
\begin{align*}
 \xi_{1}...\xi_{N}
\end{align*}
an arbitrary tensor product with contractions. Let also $k,l$ be non-negative integers, $k\leq l$. We then have:
\begin{align*}
 \|\xi_{1}...\xi_{N}\|_{2,[k,l],\Sigma_{t}^{\epsilon_{0}}}\leq C_{l}\sum_{i=1}^{N}(\prod_{j\slashed{=}i}\|\xi_{j}\|_{\infty,[k,l_{*}],\Sigma_{t}^{\epsilon_{0}}})
\|\xi_{i}\|_{2,[k,l],\Sigma_{t}^{\epsilon_{0}}}
\end{align*}
$\textbf{b}$. Let $\xi_{1},...,\xi_{N}$ be as above, and $\vartheta$ another arbitrary $S_{t,u}$ tensorfield defined on $\Sigma_{t}^{\epsilon_{0}}$. Let $k,l$ be
positive integers, $k\leq l$. We then have:
\begin{align*}
 \|\xi_{1}...\xi_{N}\cdot\vartheta\|_{2,[k-1,l-1],\Sigma_{t}^{\epsilon_{0}}}\\
\leq C_{l}\|\vartheta\|_{\infty,[k-1,l_{*}-1],\Sigma_{t}^{\epsilon_{0}}}\sum_{i=1}^{N}(\prod_{j\slashed{=}i}\|\xi_{j}\|_{\infty,[k-1,l_{*}],\Sigma_{t}^{\epsilon_{0}}})
\|\xi_{i}\|_{2,[k-1,l-1],\Sigma_{t}^{\epsilon_{0}}}\\
+C_{l}(\prod_{i=1}^{N}\|\xi_{i}\|_{\infty,[k-1,l_{*}],\Sigma_{t}^{\epsilon_{0}}})\|\vartheta\|_{2,[k-1,l-1],\Sigma_{t}^{\epsilon_{0}}}
\end{align*}
$Proof$. It is a direct calculation. $\qed$\vspace{7mm}

      Part $\textbf{a}$ shall be applied to $S_{t,u}$ tensorfields $\xi_{1}$,...,$\xi_{N}$ of order at most 1, and part $\textbf{b}$ shall be applied when the 
$\xi_{1}$,...,$\xi_{N}$ are of order at most 1, but the $S_{t,u}$ tensorfield $\vartheta$ is of order 2.

      As in Chapter 10 a primary role was played by the estimates, given by Propositions 10.1 and 10.2, for the angular derivatives of the rectangular components 
$\hat{T}^{i}$, or equivalently of the functions $y^{i}$, a primary role shall be played in the present chapter by the estimates for the spartial derivatives of the 
$\hat{T}^{i}$ of which at least one is a $T$-derivative. From (10.193) we have:
\begin{align}
 T\hat{T}^{i}=q_{T}\cdot\slashed{d}x^{i}
\end{align}
where
\begin{align}
 q_{T}=(q_{T})_{b}\cdot\slashed{g}^{-1}
\end{align}
and
\begin{align}
 (q_{T})_{b}=-\slashed{d}\kappa
\end{align}
The higher order $T$-derivatives of the $\hat{T}^{i}$ shall be expressed recursively through the following lemma.

$\textbf{Lemma 11.3}$ For any non-negative integer $m$ we have:
\begin{align*}
 (T)^{m+1}\hat{T}^{i}=p_{T,m}\hat{T}^{i}+q_{T,m}\cdot\slashed{d}x^{i}+\sum_{n=0}^{m-1}r^{n}_{T,m}\cdot\slashed{d}(T)^{n}\hat{T}^{i}
\end{align*}
Here the $p_{T,m}$ are the functions and the $q_{T,m}$ are $S_{t,u}$-tangential vectorfields determined by the recursion relations:
\begin{align*}
 p_{T,m}=Tp_{T,m-1}+\slashed{d}\kappa\cdot q_{T,m-1}\\
q_{T,m}=\slashed{\mathcal{L}}_{T}q_{T,m-1}+q_{T}p_{T,m-1}
\end{align*}
and the initial conditions:
\begin{align*}
 p_{T,0}=0,\quad q_{T,0}=q_{T}
\end{align*}
Also, the $r^{n}_{T,m}$ are $S_{t,u}$-tangential vectorfields determined by the recursion relations:
\begin{align*}
 r^{0}_{T,m}=\slashed{\mathcal{L}}_{T}r^{0}_{T,m-1}+\kappa q_{T,m-1}\\
r^{n}_{T,m}=\slashed{\mathcal{L}}_{T}r^{n}_{T,m-1}+r^{n-1}_{T,m-1}\quad:\quad \textrm{for}\quad n\in\{1,...,m-2\}\\
r^{m-1}_{T,m}=r^{m-2}_{T,m-1}
\end{align*}
and the initial condition:
\begin{align*}
 r^{0}_{T,0}=0
\end{align*}
The proof is a direct calculation. $\qed$\vspace{7mm}

      Let us now investigate the recursions of the above lemma. First we have the recursion satisfied by:
\begin{align*}
 \begin{bmatrix}
  p_{T,m}\\
q_{T,m}
 \end{bmatrix}
\end{align*}
Here we considering columns
\begin{align*}
 \begin{bmatrix}
  \phi\\
    X
 \end{bmatrix}
\end{align*}
where $\phi$ is a function and $X$ is a $S_{t,u}$-tangential vectorfield. We define the operator $\textbf{A}$ acting on such columns by:
\begin{align}
 \textbf{A}=\slashed{\mathcal{L}}_{T}+\textbf{B}
\end{align}
where, naturally,
\begin{align*}
 \slashed{\mathcal{L}}_{T}\begin{bmatrix}
                       \phi\\
                        X
                      \end{bmatrix}
=\begin{bmatrix}
  T\phi\\
\slashed{\mathcal{L}}_{T}X
 \end{bmatrix}
\end{align*}
and $\textbf{B}$ is the multiplication operator:
\begin{align}
 \textbf{B}\begin{bmatrix}
            \phi\\
X
           \end{bmatrix}
=\begin{bmatrix}
  0&\slashed{d}\kappa\\
q_{T}&0
 \end{bmatrix}
\cdot\begin{bmatrix}
      \phi\\
X
     \end{bmatrix}
=\begin{bmatrix}
  \slashed{d}\kappa\cdot X\\
q_{T}\phi
 \end{bmatrix}
\end{align}
In terms of the operator $\textbf{A}$ the recursion in question takes the form, simply:
\begin{align}
 \begin{bmatrix}
  p_{T,m}\\
q_{T,m}
 \end{bmatrix}
=\textbf{A}\begin{bmatrix}
            p_{T,m-1}\\
q_{T,m-1}
           \end{bmatrix}
\end{align}
and the solution satisfying the initial condition
\begin{align}
 \begin{bmatrix}
  p_{T,0}\\
q_{T,0}
 \end{bmatrix}
=\begin{bmatrix}
  0\\
q_{T}
 \end{bmatrix}
\end{align}
is:
\begin{align}
 \begin{bmatrix}
  p_{T,m}\\
q_{T,m}
 \end{bmatrix}
=\textbf{A}^{m}\begin{bmatrix}
                0\\
q_{T}
               \end{bmatrix}
\end{align}

     We turn to the 2-dimensional recursion satisfied by the $S_{t,u}$-tangential vectorfields $r^{n}_{T,m}$. The first recursion is
\begin{align}
 r^{0}_{T,m}=\slashed{\mathcal{L}}_{T}r^{0}_{T,m-1}+\kappa q_{T,m-1}
\end{align}
by the initial condition
\begin{align}
 r^{0}_{T,0}=0
\end{align}
we obtain:
\begin{align}
 r^{0}_{T,m}=\sum_{j=0}^{m-1}(\slashed{\mathcal{L}}_{T})^{j}(\kappa q_{T,m-1-j})
\end{align}
Next, we consider the third recursion relation. This is:
\begin{align}
 r^{m-1}_{T,m}=r^{m-2}_{T,m-1}
\end{align}
It follows that:
\begin{align}
 r^{m-1}_{T,m}=r^{0}_{T,1}=\kappa q_{T}
\end{align}
Finally, we consider the second recursion relation, which takes the form (setting $k=m-n$):
\begin{align}
 r^{m-k}_{T,m}=\slashed{\mathcal{L}}_{T}r^{m-k}_{T,m-1}+r^{m-1-k}_{T,m-1}
\end{align}
for $k=2,...,m-1$.

Given any positive integer $k$ let us denote by $\mathcal{N}_{k}$ the set of positive integers:
\begin{align}
 \mathcal{N}_{k}=\{m\geq k\}
\end{align}
On $\mathcal{N}_{k}$ we define the function $x_{k}$ with values in the space $\mathcal{X}$ of $S_{t,u}$-tangential vectorfields by:
\begin{align}
 x_{k}(m)=r^{m-k}_{T,m}
\end{align}
Let $\textbf{L}_{k}$ be the linear map taking $\mathcal{X}$-valued functions defined on $\mathcal{N}_{k-1}$ into $\mathcal{X}$-valued functions defined on 
$\mathcal{N}_{k}$ by:
\begin{align}
 (\textbf{L}_{k}f)(m)=\slashed{\mathcal{L}}_{T}(f(m-1))\quad:\quad \forall m\in\mathcal{N}_{k}
\end{align}
Then (11.25) takes the form:
\begin{align}
 x_{k}(m)-x_{k}(m-1)=(\textbf{L}_{k}x_{k-1})(m)\quad:\quad \forall m\in\mathcal{N}_{k+1}, k\geq 2
\end{align}
Setting:
\begin{align}
 y_{k}=\textbf{L}_{k}x_{k-1}
\end{align}
we can write (11.29) as:
\begin{align}
 x_{k}(m)-x_{k}(m-1)=y_{k}(m)\quad:\quad \forall m\in\mathcal{N}_{k+1}, k\geq 2
\end{align}
Let us sum (11.31) over $m\in\{k+1,...,j\}$. We obtain:
\begin{align*}
 \sum_{m=k+1}^{j}(x_{k}(m)-x_{k}(m-1))=x_{k}(j)-x_{k}(k)
\end{align*}
i.e.
\begin{align}
 x_{k}(m)=x_{k}(k)+\sum_{j=k+1}^{m}y_{k}(j)\quad:\quad \forall m\in\mathcal{N}_{k}, k\geq 2
\end{align}
Now, from (11.27) with $m=k$,
\begin{align}
 x_{k}(k)=c_{k}=r^{0}_{T,k}
\end{align}
is given by (11.22), while from (11.30) and (11.28) the sum on the right is:
\begin{align}
 \sum_{j=k+1}^{m}y_{k}(j)=\sum_{j=k+1}^{m}(\textbf{L}_{k}x_{k-1})(j)=\sum_{j=k+1}^{m}\slashed{\mathcal{L}}_{T}(x_{k-1}(j-1))
\end{align}
Let us define a linear operator $\textbf{M}_{k}$ acting on $\mathcal{X}$-valued functions $f$ defined on $\mathcal{N}_{k}$ by:
\begin{align}
 (\textbf{M}_{k}f)(m)=\sum_{j=k}^{m-1}(\slashed{\mathcal{L}}_{T}f)(j)
\end{align}
In terms of the operator $\textbf{M}_{k}$ and the constants $c_{k}$, (11.32) becomes the following recursion in $k$ for the $\mathcal{X}$-valued
functions $x_{k}$:
\begin{align}
 x_{k}=c_{k}+\textbf{M}_{k}x_{k-1}\quad:\quad k\geq 2
\end{align}
Here we consider $c_{k}$ as constant functions defined on $\mathcal{N}_{k}$. Now, by (11.27) and (11.24) for $k=1$ we have:
\begin{align}
 x_{1}(m)=r^{m-1}_{T,m}=\kappa q_{T}\quad:\quad \forall m\in\mathcal{N}_{1}
\end{align}
hence, by (11.33) for $k=1$,
\begin{align}
 x_{1}=c_{1}
\end{align}
Thus (11.36) holds also for $k=1$ if we set $x_{0}=0$:
\begin{align}
 x_{k}=c_{k}+\textbf{M}_{k}x_{k-1}\quad:\quad k\geq 1\\\notag
x_{0}=0
\end{align}
To this recursion Proposition 8.2 applies, yielding:
\begin{align}
 x_{k}=\sum_{i=0}^{k-1}\textbf{M}_{k}...\textbf{M}_{k-i+1}c_{k-i}
\end{align}
Now by (11.35),
\begin{align}
 (\textbf{M}_{k}c)(m)=(m-k)\slashed{\mathcal{L}}_{T}c\quad:\quad \forall m\in\mathcal{N}_{k}
\end{align}
Let us define, for each non-negative integer $l$, the function $N_{l}$ on the set of non-negative integers, with the value in the same set,
recursively by:
\begin{align}
 N_{0}(n)=1\quad:\textrm{for all non-negative integers}\quad n\\\notag
N_{l}(n)=\sum_{m=1}^{n}N_{l-1}(m)
\end{align}
for all non-negative integers $n$ and all positive integers $l$.

In particular, we have:
\begin{align}
 N_{l}(0)=0\quad:\textrm{for all positive integers} \quad l
\end{align}
In fact, it is readily seen that:
\begin{align}
 N_{l}(n)=\frac{n...(n+l-1)}{l!}\quad:\textrm{for all non-negative integers}\quad n\quad \textrm{and all positive integers} \quad l
\end{align}
We can then show, by induction on $j$ 
, that for all non-negative integers $j$ we have:
\begin{align}
 (\textbf{M}_{k-i+j}...\textbf{M}_{k-i+1}c_{k-i})(m)=N_{j}(m-k+i-j)(\slashed{\mathcal{L}}_{T})^{j}c_{k-i}
\end{align}
Setting $j=i$ in the above we obtain:
\begin{align}
 (\textbf{M}_{k}...\textbf{M}_{k-i+1}c_{k-i})(m)=N_{i}(m-k)(\slashed{\mathcal{L}}_{T})^{i}c_{k-i}\quad:\forall m\in\mathcal{N}_{k}
\end{align}
Substituting in (11.40) then yields:
\begin{align}
 x_{k}(m)=\sum_{i=0}^{k-1}N_{i}(m-k)(\slashed{\mathcal{L}}_{T})^{i}c_{k-i}\quad:\forall m\in\mathcal{N}_{k}, \forall k\geq 1
\end{align}
     We summarize the above results in the following lemma.

$\textbf{Lemma 11.4}$ Let $\textbf{A}$ be the operator
\begin{align*}
 \textbf{A}=\slashed{\mathcal{L}}_{T}+\textbf{B}\quad\textrm{acting on columns}\quad \begin{bmatrix}
                                                                                     \phi\\
X
                                                                                    \end{bmatrix}
\end{align*}
where $\phi$ is a function and $X$ a $S_{t,u}$-tangential vectorfield, $\textbf{B}$ being the multiplication operator:
\begin{align*}
 \textbf{B}=\begin{bmatrix}
             0&\slashed{d}\kappa\\
q_{T}&0
            \end{bmatrix}
\end{align*}
Then functions $p_{T,m}$ and $S_{t,u}$-tangential vectorfields $q_{T,m}$ are given by:
\begin{align*}
 \begin{bmatrix}
  p_{T,m}\\
q_{T,m}
 \end{bmatrix}
=\textbf{A}^{m}\begin{bmatrix}
                p_{T}\\
q_{T}
               \end{bmatrix}
\end{align*}
Moreover, the $S_{t,u}$-tangential vectorfields $r^{n}_{T,m-1}$, $n\in\{0,...,m-1\}$, are given by:
\begin{align*}
 r^{n}_{T,m}=\sum_{i=0}^{m-1-n}N_{i}(n)\sum_{j=i}^{m-1-n}(\slashed{\mathcal{L}}_{T})^{j}(\kappa q_{T,m-1-n-j})
\end{align*}
where $N_{i}(n)$ are the non-negative integers:
\begin{align*}
 N_{i}(n)=1\quad :\textrm{if} \quad i=0\\
=\frac{n...(n+i-1)}{i!}\quad: \textrm{if}\quad i\geq 1
\end{align*}\vspace{7mm}

      The $S_{t,u}$ 1-form $(q_{T})_{b}$ can be directly estimated, however to estimate the $S_{t,u}$-tangential vectorfield $q_{T}$, which enters all the above
formulas, we must estimate the $T$-derivatives of $\slashed{g}^{-1}$. Now,
\begin{align}
 \slashed{\mathcal{L}}_{T}\slashed{g}^{-1}=-\slashed{g}^{-1}\cdot\leftexp{(T)}{\slashed{\pi}}\cdot\slashed{g}^{-1},\quad \leftexp{(T)}{\slashed{\pi}}=
\slashed{\mathcal{L}}_{T}\slashed{g}
\end{align}
 Thus estimating $(\slashed{\mathcal{L}}_{T})^{m+1}\slashed{g}^{-1}$ reduces to estimating $(\slashed{\mathcal{L}}_{T})^{m}\leftexp{(T)}{\slashed{\pi}}$. 
Now from (3.27) we have:
\begin{align}
 \leftexp{(T)}{\slashed{\pi}}=2\kappa\theta=-2\alpha^{-1}\kappa\chi+2\kappa\slashed{k}
\end{align}
Let us define:
\begin{align}
 \chi^{\prime}=\chi-\frac{\slashed{g}}{1-u+t}
\end{align}
and:
\begin{align}
 \leftexp{(T)}{\slashed{\pi}}^{\prime}=-2\alpha^{-1}\kappa\chi^{\prime}+2\kappa\slashed{k}
\end{align}
We then have:
\begin{align}
 \leftexp{(T)}{\slashed{\pi}}=\leftexp{(T)}{\slashed{\pi}}^{\prime}+\lambda(1-u+t)^{-1}\slashed{g}
\end{align}
where $\lambda$ is the function:
\begin{align}
 \lambda=-2(\alpha^{-1}\kappa)
\end{align}

     We shall be able to estimate in a more direct manner the $T$-derivatives of $\leftexp{(T)}{\slashed{\pi}}^{\prime}$. From these estimates for the $T$-derivatives 
of $\leftexp{(T)}{\slashed{\pi}}$, we shall use the auxiliary symmetric 2-covariant $S_{t,u}$ tensorfields $\leftexp{(m;T)}{\slashed{\pi}}$ defined below. Their 
definition requires the polynomials $p_{m}(x)$, which we shall presently define.

     For each non-negative integer $m$ we define:
\begin{align}
 p_{m}(x)=1\quad: \textrm{for}\quad m=0\\\notag
    =(m+x)...(1+x)\quad: \textrm{for}\quad m\geq 1
\end{align}
We can see inductively, that
\begin{align}
 m!\{1+x\sum_{n=0}^{m-1}\frac{p_{n}(x)}{(n+1)!}\}=p_{m}(x)
\end{align}
     We define, for any non-negative integer $m$:
\begin{align}
 \leftexp{(m;T)}{\slashed{\pi}}=(\slashed{\mathcal{L}}_{T})^{m}\leftexp{(T)}{\slashed{\pi}}-\lambda p_{m}(\lambda)(1-u+t)^{-m-1}\slashed{g}
\end{align}
Note that:
\begin{align}
 \leftexp{(0;T)}{\slashed{\pi}}=\leftexp{(T)}{\slashed{\pi}}^{\prime}
\end{align}
Considering the fact, from (11.53), that $\lambda=-2$ in the constant state, we define:
\begin{align}
 \lambda^{\prime}=\lambda+2
\end{align}
$\textbf{Lemma 11.5}$ Let $l$ be a positive integer and $m$ a non-negative integer, $m\leq l$. Suppose that:
\begin{align*}
 \max_{i}\|\leftexp{(R_{i})}{\slashed{\pi}}\|_{\infty,[l-1],\Sigma_{t}^{\epsilon_{0}}}\leq C_{l}\delta_{0}(1+t)^{-1}[1+\log(1+t)]
\end{align*}
and:
\begin{align*}
 \|\lambda^{\prime}\|_{\infty,[m,l],\Sigma_{t}^{\epsilon_{0}}}\leq C_{l}\delta_{0}[1+\log(1+t)]
\end{align*}
Then if $\delta_{0}$ is suitably small (depending on $l$) we have:
\begin{align*}
 \max_{0\leq k\leq m}\|\leftexp{(k;T)}{\slashed{\pi}}\|_{\infty,[l-k],\Sigma_{t}^{\epsilon_{0}}}
\leq C_{l}\|\leftexp{(T)}{\slashed{\pi}}^{\prime}\|_{\infty,[m,l],\Sigma_{t}^{\epsilon_{0}}}+C_{l}\delta_{0}(1+t)^{-1}[1+\log(1+t)]
\end{align*}
$Proof$. We apply $(\slashed{\mathcal{L}}_{T})^{k}$ to (11.52), for $k\in\{0,...,m\}$. A direct calculation gives:
\begin{align}
 \leftexp{(k;T)}{\slashed{\pi}}=(\slashed{\mathcal{L}}_{T})^{k}\leftexp{(T)}{\slashed{\pi}}^{\prime}\\\notag
+\sum_{k_{1}+k_{3}=k, k_{3}\geq 1}\frac{k!}{k_{3}!}(1-u+t)^{-k_{1}-1}\lambda\leftexp{(k_{3}-1;T)}{\slashed{\pi}}\\\notag
+\sum_{k_{1}+k_{2}=k, k_{2}\geq 1}\frac{k!}{k_{2}!}(1-u+t)^{-k_{1}-1}((T)^{k_{2}}\lambda)\slashed{g}\\\notag
+\sum_{k_{1}+k_{2}+k_{3}=k, k_{2},k_{3}\geq 1}\frac{k!}{k_{2}!k_{3}!}(1-u+t)^{-k_{1}-k_{3}-1}((T)^{k_{2}}\lambda)\lambda p_{k_{3}-1}(\lambda)\slashed{g}\\\notag
+\sum_{k_{1}+k_{2}+k_{3}=k, k_{2},k_{3}\geq 1}\frac{k!}{k_{2}!k_{3}!}(1-u+t)^{-k_{1}-1}((T)^{k_{2}}\lambda)\leftexp{(k_{3}-1;T)}{\slashed{\pi}}
\end{align}
We apply $\slashed{\mathcal{L}}_{R_{i_{n}}}...\slashed{\mathcal{L}}_{R_{i_{1}}}$ to (11.59) with $n\in\{0,...,l-k\}$. The first term on the right is bounded by:
\begin{align}
 \|(\slashed{\mathcal{L}}_{T})^{k}\leftexp{(T)}{\slashed{\pi}}^{\prime}\|_{\infty,[l-k],\Sigma_{t}^{\epsilon_{0}}}
\end{align}
We decompose the first sum on the right according to:
\begin{align}
 \lambda=-2+\lambda^{\prime}
\end{align}
Setting:
\begin{align}
 k^{\prime}=k_{3}-1
\end{align}
the contribution of the term $-2$ is:
\begin{align}
 2\sum_{k^{\prime}=0}^{k-1}\frac{k!}{(k^{\prime}+1)!}(1-u+t)^{k^{\prime}-k}\|\leftexp{(k^{\prime};T)}{\slashed{\pi}}\|
_{\infty,[l-k],\Sigma_{t}^{\epsilon_{0}}}
\end{align}
While the term $\lambda^{\prime}$ is bounded by:
\begin{align}
 C_{l}\delta_{0}[1+\log(1+t)]\sum_{k^{\prime}=0}^{k-1}\frac{k!}{(k^{\prime}+1)!}(1-u+t)^{k^{\prime}-k}\|\leftexp{(k^{\prime};T)}{\slashed{\pi}}
\|_{\infty,[l-k],\Sigma_{t}^{\epsilon_{0}}}
\end{align}
The last sum is bounded by:
\begin{align}
 C_{l}\delta_{0}[1+\log(1+t)]\\\notag
\times \sum_{k^{\prime}=0}^{k-2}\sum_{j=0}^{k-2-k^{\prime}}\frac{k!}{(j+1)!(k^{\prime}+1)!}(1-u+t)^{k^{\prime}-k+1+j}\|\leftexp{(k^{\prime};T)}{\slashed{\pi}}\|
_{\infty,[l-k],\Sigma_{t}^{\epsilon_{0}}}
\end{align}
where we have set:
\begin{align}
 k_{2}=1+j
\end{align}
Next, we consider the second sum, which is:
\begin{align}
 \sum_{j=0}^{k-1}\frac{k!}{(j+1)!}(1-u+t)^{j-k}((T)^{j+1}\lambda^{\prime})\slashed{g}
\end{align}
By the first assumption,
\begin{align}
 \|((T)^{j+1}\lambda^{\prime})\slashed{g}\|_{\infty,[l-k],\Sigma_{t}^{\epsilon_{0}}}\leq C\|(T)^{j+1}\lambda^{\prime}\|_{\infty,[l-k],\Sigma_{t}^{\epsilon_{0}}}
\end{align}
for $k\geq 1$, provided that $\delta_{0}$ is suitably small (depending on $l$). Then by the second assumption we know that the second sum is bounded by:
\begin{align}
 C_{l}\delta_{0}(1+t)^{-1}[1+\log(1+t)]
\end{align}
Finally, we have to deal with the third sum. We express:
\begin{align}
 \lambda p_{k^{\prime}}(\lambda)=-2p_{k^{\prime}}(-2)+q_{k^{\prime}+1}(\lambda^{\prime})
\end{align}
where $q_{k^{\prime}+1}(\lambda^{\prime})$ is a polynomial of degree $k^{\prime}+1$ in $\lambda^{\prime}$ with vanishing constant term.
By the second assumption of the lemma, the constant term is bounded by:
\begin{align}
 C_{l}\delta_{0}(1+t)^{-2}[1+\log(1+t)]
\end{align}
while the contribution of the polynomial $q_{k^{\prime}+1}(\lambda^{\prime})$ is bounded by:
\begin{align}
 C_{l}\delta_{0}[1+\log(1+t)]\sum_{k^{\prime}=0}^{k-2}\{\sum_{i=1}^{k^{\prime}+1}(\delta_{0}[1+\log(1+t)])^{i}\}(1+t)^{-2-k^{\prime}}
\end{align}
A simple calculation yields that this is bounded by:
\begin{align}
 C_{l}\delta_{0}^{2}(1+t)^{-2}[1+\log(1+t)]^{2}
\end{align}
provided that $\delta_{0}$ is suitably small.

     Combining the results (11.60), (11.63), (11.64), (11.65), (11.69), (11.71) and (11.73), we conclude that
\begin{align}
 \|\leftexp{(k;T)}{\slashed{\pi}}\|_{\infty,[l-k],\Sigma_{t}^{\epsilon_{0}}}\leq \|(\slashed{\mathcal{L}}_{T})^{k}\leftexp{(T)}{\slashed{\pi}}^{\prime}\|
_{\infty,[l-k],\Sigma_{t}^{\epsilon_{0}}}\\\notag
+C_{l}(1+t)^{-1}[1+\log(1+t)]\sum_{k^{\prime}=0}^{k-1}\|\leftexp{(k^{\prime};T)}{\slashed{\pi}}\|_{\infty,[l-k],\Sigma_{t}^{\epsilon_{0}}}\\\notag
+C_{l}\delta_{0}(1+t)^{-1}[1+\log(1+t)]
\end{align}
Since:
\begin{align}
 \sum_{k'=0}^{k-1}\|\leftexp{(k';T)}{\slashed{\pi}}\|_{\infty,[l-k],\Sigma_{t}^{\epsilon_{0}}}\leq 
\sum_{k'=0}^{k-1}\|\leftexp{(k';T)}{\slashed{\pi}}\|_{\infty,[l-1-k'],\Sigma_{t}^{\epsilon_{0}}}\leq 
\sum_{k'=0}^{k-1}\|\leftexp{(k';T)}{\slashed{\pi}}\|_{\infty,[l-k'],\Sigma_{t}^{\epsilon_{0}}}
\end{align}
(11.74) reduce to:
\begin{align}
 x_{k}\leq c_{k}+a\sum_{k'=0}^{k-1}x_{k'}\quad: \textrm{for}\quad k=0,...,m
\end{align}
where
\begin{align}
 x_{k}=\|\leftexp{(k;T)}{\slashed{\pi}}\|_{\infty,[l-k],\Sigma_{t}^{\epsilon_{0}}}\\\notag
c_{k}=\|(\slashed{\mathcal{L}}_{T})^{k}\leftexp{(T)}{\slashed{\pi}}^{\prime}\|_{\infty,[l-k],\Sigma_{t}^{\epsilon_{0}}}
+C_{l}\delta_{0}(1+t)^{-1}[1+\log(1+t)]\\\notag
a=C_{l}(1+t)^{-1}[1+\log(1+t)]
\end{align}
It follows by direct calculation that
\begin{align}
 x_{k}\leq c_{k}+a\sum_{k'=0}^{k-1}(1+a)^{k-1-k'}c_{k'}\quad: \textrm{for} \quad k=0,...,m
\end{align}
hence:
\begin{align}
 \max_{0\leq k\leq m}x_{k}\leq \max_{0\leq k\leq m}c_{k}\{1+a\sum_{k'=0}^{m-1}(1+a)^{m-1-k'}\}=(1+a)^{m}\max_{0\leq k\leq m}c_{k}
\end{align}
Substituting (11.77) and recalling that $(1+t)^{-1}[1+\log(1+t)]\leq 1$ we obtain:
\begin{align*}
 \max_{0\leq k\leq m}\|\leftexp{(k;T)}{\slashed{\pi}}\|_{\infty,[l-k],\Sigma_{t}^{\epsilon_{0}}}\leq C_{l}
\{\|\leftexp{(T)}{\slashed{\pi}}^{\prime}\|_{\infty,[m,l],\Sigma_{t}^{\epsilon_{0}}}+C_{l}\delta_{0}(1+t)^{-1}[1+\log(1+t)]\}
\end{align*}
$\qed$

$\textbf{Lemma 11.6}$ Let $l$ be a positive integer and $m$ a non-negative integer, $m\leq l$. Suppose that:
\begin{align*}
 \max_{i}\|\leftexp{(R_{i})}{\slashed{\pi}}\|_{\infty,[l_{*}-1],\Sigma_{t}^{\epsilon_{0}}}\leq C_{l}\delta_{0}(1+t)^{-1}[1+\log(1+t)]\\
\|\lambda^{\prime}\|_{\infty, [m,l_{*}], \Sigma_{t}^{\epsilon_{0}}}\leq C_{l}\delta_{0}[1+\log(1+t)]
\end{align*}
and:
\begin{align*}
 \|\leftexp{(T)}{\slashed{\pi}}^{\prime}\|_{\infty,[m-1,l_{*}-1],\Sigma_{t}^{\epsilon_{0}}}\leq C_{l}\delta_{0}(1+t)^{-1}[1+\log(1+t)]
\end{align*}
Then, if $\delta_{0}$ is suitably small (depending on $l$) we have:
\begin{align*}
 \sum_{k=0}^{m}\|\leftexp{(k;T)}{\slashed{\pi}}\|_{2,[l-k],\Sigma_{t}^{\epsilon_{0}}}\leq C_{l}\|\leftexp{(T)}{\slashed{\pi}}^{\prime}\|
_{2,[m,l],\Sigma_{t}^{\epsilon_{0}}}\\
+C_{l}(1+t)^{-1}\{\|\lambda^{\prime}\|_{2,[m,l],\Sigma_{t}^{\epsilon_{0}}}+\delta_{0}[1+\log(1+t)]\max_{i}\|\leftexp{(R_{i})}{\slashed{\pi}}\|_
{2,[l-2],\Sigma_{t}^{\epsilon_{0}}}\}
\end{align*}
$Proof$. We use (11.59). The proof is quite similar to that of Lemma 11.5. As in Chapter 10, when we estimate a product in $L^{2}$, we put 
the lower order factors in $L^{\infty}$, and the highest order factor in $L^{2}$. $\qed$ \vspace{7mm}

     We use Lemmas 11.5 and 11.6 to derive $L^{\infty}$ and $L^{2}$ estimates for the $T$-derivatives of $\leftexp{(T)}{\slashed{\pi}}$ in terms of 
estimates for the $T$-derivatives of $\leftexp{(T)}{\slashed{\pi}}^{\prime}$.

$\textbf{Lemma 11.7}$ Under the assumptions of Lemma 11.5 we have:
\begin{align*}
 \|\leftexp{(T)}{\slashed{\pi}}+2(1-u+t)^{-1}\slashed{g}\|_{\infty,[l],\Sigma_{t}^{\epsilon_{0}}}\leq C_{l}\|\leftexp{(T)}{\slashed{\pi}}^{\prime}\|
_{\infty,[0,l],\Sigma_{t}^{\epsilon_{0}}}\\+C_{l}\delta_{0}(1+t)^{-1}[1+\log(1+t)]\\
\|\slashed{\mathcal{L}}_{T}\leftexp{(T)}{\slashed{\pi}}-2(1-u+t)^{-2}\slashed{g}\|_{\infty,[l-1],\Sigma_{t}^{\epsilon_{0}}}\leq C_{l}
\|\leftexp{(T)}{\slashed{\pi}}^{\prime}\|_{\infty,[1,l],\Sigma_{t}^{\epsilon_{0}}}\\
+C_{l}\delta_{0}(1+t)^{-1}[1+\log(1+t)]
\end{align*}
and:
\begin{align*}
 \max_{2\leq k\leq m}\|(\slashed{\mathcal{L}}_{T})^{k}\leftexp{(T)}{\slashed{\pi}}\|_{\infty,[l-k],\Sigma_{t}^{\epsilon_{0}}}\leq C_{l}
\|\leftexp{(T)}{\slashed{\pi}}^{\prime}\|_{\infty,[m,l],\Sigma_{t}^{\epsilon_{0}}}\\
+C_{l}\delta_{0}(1+t)^{-1}[1+\log(1+t)]
\end{align*}
$Proof$. We consider formula (11.56) with $m$ replaced by $k\in\{0,...,m\}$:
\begin{align}
 (\slashed{\mathcal{L}}_{T})^{k}\leftexp{(T)}{\slashed{\pi}}=\leftexp{(k;T)}{\slashed{\pi}}+\lambda p_{k}(\lambda)(1-u+t)^{-k-1}\slashed{g}
\end{align}
We express $\lambda p_{k}(\lambda)$ as in (11.70) and we bring the contribution of the constant term $-2p_{k}(-2)$ on the left. We obtain:
\begin{align}
 \leftexp{(T)}{\slashed{\pi}}+2(1-u+t)^{-1}\slashed{g}=\leftexp{(0;T)}{\slashed{\pi}}+q_{1}(\lambda^{\prime})(1-u+t)^{-1}\slashed{g}\\\notag
\slashed{\mathcal{L}}_{T}\leftexp{(T)}{\slashed{\pi}}-2(1-u+t)^{-2}\slashed{g}=\leftexp{(1;T)}{\slashed{\pi}}+q_{2}(\lambda^{\prime})(1-u+t)^{-2}\slashed{g}\\\notag
(\slashed{\mathcal{L}}_{T})^{k}\leftexp{(T)}{\slashed{\pi}}=\leftexp{(k;T)}{\slashed{\pi}}+q_{k+1}(\lambda^{\prime})(1-u+t)^{-k-1}\slashed{g}
\end{align}
for $k\geq 2$.

      We first take the $\|\quad\|_{\infty,[l],\Sigma_{t}^{\epsilon_{0}}}$-norm of each side of the first of (11.81). Since:
\begin{align}
 \|q_{1}(\lambda^{\prime})\slashed{g}\|_{\infty,[l],\Sigma_{t}^{\epsilon_{0}}}\leq C\|\lambda^{\prime}\|_{\infty,[l],\Sigma_{t}^{\epsilon_{0}}}\\\notag
+C_{l}\|\lambda^{\prime}\|_{\infty,[l-1],\Sigma_{t}^{\epsilon_{0}}}\max_{j}\|\leftexp{(R_{i})}{\slashed{\pi}}\|_{\infty,[l-1],\Sigma_{t}^{\epsilon_{0}}}
\leq C_{l}\delta_{0}[1+\log(1+t)]
\end{align}
(provided that $\delta_{0}$ is suitably small), the first statement of the lemma results by substituting the estimate for $\|\leftexp{(0;T)}{\slashed{\pi}}\|
_{\infty,[l],\Sigma_{t}^{\epsilon_{0}}}$ of Lemma 11.5 with $m=0$. The next two cases in (11.81) can be dealt in a similar way, then the lemma is proved. $\qed$ 
\vspace{7mm}

$\textbf{Lemma 11.8}$ Let the assumptions of Lemma 11.6 hold. Let us define:
\begin{align*}
 \mathcal{T}_{[0,l]}=\|\leftexp{(T)}{\slashed{\pi}}+2(1-u+t)^{-1}\slashed{g}\|_{2,[l],\Sigma_{t}^{\epsilon_{0}}}\\
\mathcal{T}_{[1,l]}=\|\leftexp{(T)}{\slashed{\pi}}+2(1-u+t)^{-1}\slashed{g}\|_{2,[l],\Sigma_{t}^{\epsilon_{0}}}\\
+\|\slashed{\mathcal{L}}_{T}\leftexp{(T)}{\slashed{\pi}}-2(1-u+t)^{-2}\slashed{g}\|_{2,[l-1],\Sigma_{t}^{\epsilon_{0}}}
\end{align*}
and for $m\geq 2$:
\begin{align*}
 \mathcal{T}_{[m,l]}=\|\leftexp{(T)}{\slashed{\pi}}+2(1-u+t)^{-1}\slashed{g}\|_{2,[l],\Sigma_{t}^{\epsilon_{0}}}\\
+\|\slashed{\mathcal{L}}_{T}\leftexp{(T)}{\slashed{\pi}}-2(1-u+t)^{-2}\slashed{g}\|_{2,[l-1],\Sigma_{t}^{\epsilon_{0}}}\\
+\sum_{k=2}^{m}\|(\slashed{\mathcal{L}}_{T})^{k}\leftexp{(T)}{\slashed{\pi}}\|_{2,[l-k],\Sigma_{t}^{\epsilon_{0}}}
\end{align*}
We then have:
\begin{align*}
 \mathcal{T}_{[m,l]}\leq C_{l}\|\leftexp{(T)}{\slashed{\pi}}^{\prime}\|_{2,[m,l],\Sigma_{t}^{\epsilon_{0}}}\\
+C_{l}(1+t)^{-1}\{\|\lambda^{\prime}\|_{2,[m,l],\Sigma_{t}^{\epsilon_{0}}}+\delta_{0}[1+\log(1+t)]\max_{i}\|\leftexp{(R_{i})}{\slashed{\pi}}\|
_{2,[l-1],\Sigma_{t}^{\epsilon_{0}}}\}
\end{align*}
$Proof$. Again, the proof is similar to that of Lemma 11.7 and Lemma 11.6. $\qed$
\vspace{7mm}    
      In establishing the primary propositions of the present chapter, which concern estimates for the spatial derivatives of the $\hat{T}^{i}$ of which at least 
one is a $T$-derivative, we shall have to estimate, for given $S_{t,u}$ 1-forms $\xi$, the commutators:
\begin{align}
 \leftexp{(i_{1}...i_{n})}{c}_{m,n}[\xi]=\slashed{\mathcal{L}}_{R_{i_{n}}}...\slashed{\mathcal{L}}_{R_{i_{1}}}(\slashed{\mathcal{L}}_{T})^{m}\slashed{D}\xi
-\slashed{D}\slashed{\mathcal{L}}_{R_{i_{n}}}...\slashed{\mathcal{L}}_{R_{i_{1}}}(\slashed{\mathcal{L}}_{T})^{m}\xi
\end{align}
$\textbf{Lemma 11.9}$ The commutators $\leftexp{(i_{1}...i_{n})}{c}_{m,n}[\xi]$ are given by:
\begin{align*}
 \leftexp{(i_{1}...i_{n})}{c}_{m,n}[\xi]=\slashed{\mathcal{L}}_{R_{i_{n}}}...\slashed{\mathcal{L}}_{R_{i_{1}}}c_{m,0}[\xi]\\
-\sum_{j=0}^{n-1}\slashed{\mathcal{L}}_{R_{i_{n}}}...\slashed{\mathcal{L}}_{R_{i_{n-j+1}}}\{\leftexp{(R_{i_{n-j}})}{\slashed{\pi}}_{1}\cdot
\slashed{\mathcal{L}}_{R_{i_{n-j-1}}}...\slashed{\mathcal{L}}_{R_{i_{1}}}(\slashed{\mathcal{L}}_{T})^{m}\xi\}
\end{align*}
where $c_{m,0}[\xi]$ are given by:
\begin{align*}
 c_{m,0}[\xi]=-\sum_{i=0}^{m-1}(\slashed{\mathcal{L}}_{T})^{i}\{\leftexp{(T)}{\slashed{\pi}}_{1}\cdot(\slashed{\mathcal{L}}_{T})^{m-i-1}\xi\}
\end{align*}
Here $\leftexp{(T)}{\slashed{\pi}}_{1}$ is defined by (9.161) and (9.162).

$Proof$. The proof is a direct calculation which relies on Proposition 8.2. $\qed$
\vspace{7mm}      

     We must now obtain estimates for $\leftexp{(T)}{\slashed{\pi}}_{1}$. In terms of the operator $\check{\slashed{D}}$ defined by (10.147), we can express:
\begin{align}
 \leftexp{(T)}{\slashed{\pi}}_{1}=\leftexp{(T)}{(\slashed{\pi}_{1})}_{b}\cdot\slashed{g}^{-1},\quad 
\leftexp{(T)}{(\slashed{\pi}_{1})}_{b}=\check{\slashed{D}}\leftexp{(T)}{\slashed{\pi}}
\end{align}
We shall need the following analogue of Lemma 10.9, which can be also proved in a straightforward manner.

$\textbf{Lemma 11.10}$ Let $\vartheta$ be an arbitrary symmetric 2-covariant $S_{t,u}$ tensorfield. The following commutation formula holds for any non-negative 
integer $m$:
\begin{align*}
 (\slashed{\mathcal{L}}_{T})^{m}\check{\slashed{D}}\vartheta-\check{\slashed{D}}(\slashed{\mathcal{L}}_{T})^{m}\vartheta=-
\sum_{k=0}^{m-1}(\slashed{\mathcal{L}}_{T})^{k}(\leftexp{(T)}{\slashed{\pi}}_{1}\cdot(\slashed{\mathcal{L}}_{T})^{m-k-1}\vartheta)
\end{align*}
We proceed to:

$\textbf{Lemma 11.11}$ Let $l$ be a positive integer and $m$ a non-negative integer $m\leq l-1$. Let also hypothesis $\textbf{H2}^{\prime}$ hold. 
Suppose that:
\begin{align*}
 \|\lambda^{\prime}\|_{\infty,[m,l],\Sigma_{t}^{\epsilon_{0}}}\leq C_{l}\delta_{0}[1+\log(1+t)]\\
\|\leftexp{(T)}{\slashed{\pi}}^{\prime}\|_{\infty,[m,l],\Sigma_{t}^{\epsilon_{0}}}\leq C_{l}\delta_{0}(1+t)^{-1}[1+\log(1+t)]\\
\max_{i}\|\leftexp{(R_{i})}{\slashed{\pi}}\|_{\infty,[l-1],\Sigma_{t}^{\epsilon_{0}}}\leq C_{l}\delta_{0}(1+t)^{-1}[1+\log(1+t)]
\end{align*}
and, if $l\geq 2$,
\begin{align*}
 \max_{i}\|\leftexp{(R_{i})}{\slashed{\pi}}_{1}\|_{\infty,[l-2],\Sigma_{t}^{\epsilon_{0}}}\leq C_{l}\delta_{0}(1+t)^{-2}[1+\log(1+t)]
\end{align*}
Then, if $\delta_{0}$ is suitably small (depending on $l$) we have:
\begin{align*}
 \|\leftexp{(T)}{\slashed{\pi}}_{1}\|_{\infty,[m,l-1],\Sigma_{t}^{\epsilon_{0}}}\leq C_{l}\delta_{0}(1+t)^{-2}[1+\log(1+t)]
\end{align*}
$Proof$. To prove the lemma, we just need to estimate 
\begin{align}
 \|(\slashed{\mathcal{L}}_{T})^{k}\leftexp{(T)}{\slashed{\pi}}_{1}\|_{\infty,[l-1-k],\Sigma_{t}^{\epsilon_{0}}}
\end{align}
To do this, we proceed inductively on $k$. Since under the assumptions of the present lemma, Lemma 11.7 holds, substituting the bound for 
\begin{align*}
 \|\leftexp{(T)}{\slashed{\pi}^{\prime}}\|_{\infty,[m,l],\Sigma_{t}^{\epsilon_{0}}}
\end{align*}
yields:
\begin{align}
 \|\leftexp{(T)}{\slashed{\pi}}+2(1-u+t)^{-1}\slashed{g}\|_{\infty,[l],\Sigma_{t}^{\epsilon_{0}}}\leq C_{l}\delta_{0}(1+t)^{-1}[1+\log(1+t)]\\\notag
\|\slashed{\mathcal{L}}_{T}\leftexp{(T)}{\slashed{\pi}}-2(1-u+t)^{-2}\slashed{g}\|_{\infty,[l-1],\Sigma_{t}^{\epsilon_{0}}}
\leq C_{l}\delta_{0}(1+t)^{-1}[1+\log(1+t)]\\\notag
\max_{2\leq k\leq m}\|(\slashed{\mathcal{L}}_{T})^{k}\leftexp{(T)}{\slashed{\pi}}\|_{\infty,[l-k],\Sigma_{t}^{\epsilon_{0}}}\leq C_{l}\delta_{0}
(1+t)^{-1}[1+\log(1+t)]
\end{align}
We can write
\begin{align}
 \leftexp{(T)}{(\slashed{\pi}_{1})}_{b}=\check{\slashed{D}}(\leftexp{(T)}{\slashed{\pi}}+2(1-u+t)^{-1}\slashed{g})
\end{align}
Then by Lemma 10.9 and a direct calculation, we obtain the conclusion when $k=0$.

     Next, we apply Lemma 11.9 to obtain:
\begin{align}
 (\slashed{\mathcal{L}}_{T})^{k}\leftexp{(T)}{(\slashed{\pi}_{1})}_{b}=\check{\slashed{D}}(\slashed{\mathcal{L}}_{T})^{k}\leftexp{(T)}{\slashed{\pi}}+d^{\prime}_{k}
\end{align}
The estimate for $k=1$ is also straightforward, we just use Lemma 11.9 to obtain an estimate for $d^{\prime}_{1}$.

Then we assume, for $k\in\{2,...,m\}$, that
\begin{align}
\|(\slashed{\mathcal{L}}_{T})^{k'}\leftexp{(T)}{\slashed{\pi}}_{1}\|_{\infty,[l-1-k'],\Sigma_{t}^{\epsilon_{0}}}\leq C_{l}\delta_{0}(1+t)^{-2}[1+\log(1+t)]
\end{align}
holds for all $k'\leq k-1$.

To obtain the estimates in the case $k'=k$, we again use Lemma 11.9, and then use the induction hypothesis to estimate $d^{\prime}_{k}$. Then the lemma follows. $\qed$ \vspace{7mm}

      To proceed with $L^{2}$ estimates for $\leftexp{(T)}{\slashed{\pi}}_{1}$, we shall need the following lemma. This lemma shall also be used in the sequel.

$\textbf{Lemma 11.12}$ Let $k,l$ be positive integers $k\leq l$ and let:
\begin{align*}
 \|\lambda^{\prime}\|_{\infty,[k-1,l_{*}-1],\Sigma_{t}^{\epsilon_{0}}}\leq C_{l}\delta_{0}[1+\log(1+t)]\\
\max_{i}\|\leftexp{(R_{i})}{\slashed{\pi}}\|_{\infty,[l_{*}-1],\Sigma_{t}^{\epsilon_{0}}}\leq C_{l}\delta_{0}(1+t)^{-1}[1+\log(1+t)]\\
\|\leftexp{(T)}{\slashed{\pi}}^{\prime}\|_{\infty,[k-1,l_{*}-1],\Sigma_{t}^{\epsilon_{0}}}\leq C_{l}\delta_{0}(1+t)^{-1}[1+\log(1+t)]
\end{align*}
 Then provided that $\delta_{0}$ is suitably small (depending on $l$), we have, for any $S_{t,u}$ tensorfield $\xi$:
\begin{align*}
 \|\xi\cdot\slashed{g}\|_{2,[k-1,l-1],\Sigma_{t}^{\epsilon_{0}}}, \|\xi\cdot\slashed{g}^{-1}\|_{2,[k-1,l-1],\Sigma_{t}^{\epsilon_{0}}}\\
\leq C_{l}\{\|\xi\|_{2,[k-1,l-1],\Sigma_{t}^{\epsilon_{0}}}+\|\xi\|_{\infty,[k-1,l_{*}-1],\Sigma_{t}^{\epsilon_{0}}}(\max_{i}\|\leftexp{(R_{i})}{\slashed{\pi}}\|
_{2,[l-2],\Sigma_{t}^{\epsilon_{0}}}+\mathcal{T}_{[k-2,l-2]})\}
\end{align*}
where the $\mathcal{T}_{[k,l]}$ are the quantities defined in the statement of Lemma 11.8 (the term $\mathcal{T}_{[k-2,l-2]}$ is present only for $k\geq 2$). 

$Proof$. Under the assumptions of the present lemma, the conclusions of Lemma 11.7 hold with $(k-1,l_{*}-1)$ in the role of $(m,l)$. Then substituting in these the 
bound for 
\begin{align*}
 \|\leftexp{(T)}{\slashed{\pi}}^{\prime}\|_{\infty,[k-1,l_{*}-1],\Sigma_{t}^{\epsilon_{0}}}
\end{align*}
yields:
\begin{align}
 \|\leftexp{(T)}{\slashed{\pi}}+2(1-u+t)^{-1}\slashed{g}\|_{\infty,[l_{*}-1],\Sigma_{t}^{\epsilon_{0}}}\leq C_{l}\delta_{0}(1+t)^{-1}[1+\log(1+t)]\\\notag
\|\slashed{\mathcal{L}}_{T}\leftexp{(T)}{\slashed{\pi}}-2(1-u+t)^{-2}\slashed{g}\|_{\infty,[l_{*}-2],\Sigma_{t}^{\epsilon_{0}}}
\leq C_{l}\delta_{0}(1+t)^{-1}[1+\log(1+t)]\\\notag
\max_{2\leq i\leq k-1}\|(\slashed{\mathcal{L}}_{T})^{i}\leftexp{(T)}{\slashed{\pi}}\|_{\infty,[l_{*}-1-i],\Sigma_{t}^{\epsilon_{0}}}\leq C_{l}\delta_{0}
(1+t)^{-1}[1+\log(1+t)]
\end{align}
We shall estimate:
\begin{align}
 \|(\slashed{\mathcal{L}}_{T})^{j-1}(\xi\cdot\slashed{g})\|_{2,[l-j],\Sigma_{t}^{\epsilon_{0}}}\quad: \textrm{for} j=1,...,k
\end{align}
The case $j=1$ is easy:
\begin{align}
 \|\xi\cdot\slashed{g}\|_{2,[l-1],\Sigma_{t}^{\epsilon_{0}}}\leq C_{l}\{\|\xi\|_{2,[l-1],\Sigma_{t}^{\epsilon_{0}}}+\|\xi\|_{\infty,[l_{*}-1],\Sigma_{t}^{\epsilon_{0}}}
\max_{i}\|\leftexp{(R_{i})}{\slashed{\pi}}\|_{2,[l-2],\Sigma_{t}^{\epsilon_{0}}}\}
\end{align}
For $j\geq 2$ we express:
\begin{align}
 (\slashed{\mathcal{L}}_{T})^{j-1}(\xi\cdot\slashed{g})=\xi_{j-1}\cdot\slashed{g}\\\notag
+(j-1)((\slashed{\mathcal{L}}_{T})^{j-2}\xi)\cdot(\leftexp{(T)}{\slashed{\pi}}+2(1-u+t)^{-1}\slashed{g})\\\notag
+\frac{(j-1)(j-2)}{2}((\slashed{\mathcal{L}}_{T})^{j-3}\xi)\cdot(\slashed{\mathcal{L}}_{T}\leftexp{(T)}{\slashed{\pi}}-2(1-u+t)^{-2}\slashed{g})\\\notag
+\sum_{i=2}^{j-2}\frac{(j-1)!}{(i+1)!(j-2-i)!}((\slashed{\mathcal{L}}_{T})^{j-2-i}\xi)\cdot((\slashed{\mathcal{L}}_{T})^{i}\leftexp{(T)}{\slashed{\pi}})
\end{align}
Here:
\begin{align}
 \xi_{j-1}=(\slashed{\mathcal{L}}_{T})^{j-1}\xi-2(j-1)(1-u+t)^{-1}(\slashed{\mathcal{L}}_{T})^{j-2}\xi\\\notag
+(j-1)(j-2)(1-u+t)^{-2}(\slashed{\mathcal{L}}_{T})^{j-3}\xi
\end{align}
Only the first two terms in (11.93) and (11.94) are present for $j=2$. Also, the sum in (11.93) is present only for $j\geq 4$. What follows is just a 
direct calculation. The estimate for $\|\xi\cdot\slashed{g}^{-1}\|_{2,[k-1,l-1],\Sigma_{t}^{\epsilon_{0}}}$ follows in a similar manner. $\qed$ \vspace{7mm}

     We proceed with $L^{2}$ estimates for $\leftexp{(T)}{\slashed{\pi}}_{1}$.

$\textbf{Lemma 11.13}$ Let $l$ be a positive integer and $m$ a non-negative integer $m\leq l-1$. Let also hypothesis $\textbf{H2}^{\prime}$ hold. Suppose that:
\begin{align*}
 \|\lambda^{\prime}\|_{\infty,[m,l_{*}],\Sigma_{t}^{\epsilon_{0}}}\leq C_{l}\delta_{0}[1+\log(1+t)]\\
\max_{i}\|\leftexp{(R_{i})}{\slashed{\pi}}\|_{\infty,[l_{*}-1],\Sigma_{t}^{\epsilon_{0}}}\leq C_{l}\delta_{0}(1+t)^{-1}[1+\log(1+t)]\\
\|\leftexp{(T)}{\slashed{\pi}}^{\prime}\|_{\infty,[m,l_{*}],\Sigma_{t}^{\epsilon_{0}}}\leq C_{l}\delta_{0}(1+t)^{-1}[1+\log(1+t)]
\end{align*}
and (if $l_{*}\geq 2$):
\begin{align*}
 \max_{i}\|\leftexp{(R_{i})}{\slashed{\pi}}_{1}\|_{\infty,[l_{*}-2],\Sigma_{t}^{\epsilon_{0}}}\leq C_{l}\delta_{0}(1+t)^{-2}[1+\log(1+t)]
\end{align*}
Then if $\delta_{0}$ is suitably small (depending on $l$) we have:
\begin{align*}
 \|\leftexp{(T)}{\slashed{\pi}}_{1}\|_{2,[m,l-1],\Sigma_{t}^{\epsilon_{0}}}\\
\leq C_{l}(1+t)^{-1}\{\mathcal{T}_{[m,l]}+\delta_{0}(1+t)^{-1}[1+\log(1+t)]\cdot\\
(\max_{i}\|\leftexp{(R_{i})}{\slashed{\pi}}\|_{2,[l-2],\Sigma_{t}^{\epsilon_{0}}}+(1+t)\max_{i}\|\leftexp{(R_{i})}{\slashed{\pi}}_{1}\|_{2,[l-2],\Sigma_{t}^{\epsilon_{0}}})\}
\end{align*}
$Proof$. The proof is quite similar to that of Lemma 11.11. The only difference is that here we must use Lemma 11.12 to connect the $L^{2}$ estimates of 
$\leftexp{(T)}{\slashed{\pi}}_{1}$ and $\leftexp{(T)}{(\slashed{\pi}_{1})}_{b}$. $\qed$\vspace{7mm}

     We now apply Lemma 11.11 and 11.13 to obtain estimates for the commutators $\leftexp{(i_{1}...i_{n})}{c}_{m,n}[\xi]$, defined by (11.83).

$\textbf{Lemma 11.14}$ Let $k,l$ be positive integers, $k\leq l$. Let also hypothesis $\textbf{H2}^{\prime}$ hold. Suppose that:
\begin{align*}
 \|\lambda^{\prime}\|_{\infty,[k-1,l],\Sigma_{t}^{\epsilon_{0}}}\leq C_{l}\delta_{0}[1+\log(1+t)]\\
\max_{i}\|\leftexp{(R_{i})}{\slashed{\pi}}\|_{\infty,[l-1],\Sigma_{t}^{\epsilon_{0}}}\leq C_{l}\delta_{0}(1+t)^{-1}[1+\log(1+t)]\\
\|\leftexp{(T)}{\slashed{\pi}}^{\prime}\|_{\infty,[k-1,l],\Sigma_{t}^{\epsilon_{0}}}\leq C_{l}\delta_{0}(1+t)^{-1}[1+\log(1+t)]
\end{align*}
and:
\begin{align*}
 \max_{i}\|\leftexp{(R_{i})}{\slashed{\pi}}_{1}\|_{\infty,[l-1],\Sigma_{t}^{\epsilon_{0}}}\leq C_{l}\delta_{0}(1+t)^{-2}[1+\log(1+t)]
\end{align*}
Then for any $S_{t,u}$ 1-form $\xi$ we have:
\begin{align*}
 \max_{m\leq k,n\leq l-m}\max_{i_{1}...i_{n}}\|\leftexp{(i_{1}...i_{n})}{c}_{m,n}[\xi]\|_{L^{\infty}(\Sigma_{t}^{\epsilon_{0}})}
\leq C_{l}\delta_{0}(1+t)^{-2}[1+\log(1+t)]\|\xi\|_{\infty,[k,l-1],\Sigma_{t}^{\epsilon_{0}}}
\end{align*}
provided that $\delta_{0}$ is suitably small (depending on $l$).

$Proof$. The assumptions of the present lemma include those of Lemma 11.11 with $k-1$ in the role of $m$. Therefore the conclusion of Lemma 11.11 holds with
$m$ replaced by $k-1$, that is, we have:
\begin{align}
 \|\leftexp{(T)}{\slashed{\pi}}_{1}\|_{\infty,[k-1,l-1],\Sigma_{t}^{\epsilon_{0}}}\leq C_{l}\delta_{0}(1+t)^{-2}[1+\log(1+t)]
\end{align}
Then we just use the expressions given by Lemma 11.9. After a direct calculation, the lemma is seen to follow. $\qed$

$\textbf{Lemma 11.15}$ Let $k,l$ be positive integers, $k\leq l$. Let also hypothesis $\textbf{H2}^{\prime}$ hold. Suppose that:
\begin{align*}
 \|\lambda^{\prime}\|_{\infty,[k-1,l_{*}],\Sigma_{t}^{\epsilon_{0}}}\leq C_{l}\delta_{0}[1+\log(1+t)]\\
\max_{i}\|\leftexp{(R_{i})}{\slashed{\pi}}\|_{\infty,[l_{*}-1],\Sigma_{t}^{\epsilon_{0}}}\leq C_{l}\delta_{0}(1+t)^{-1}[1+\log(1+t)]\\
\|\leftexp{(T)}{\slashed{\pi}}^{\prime}\|_{\infty,[k-1,l_{*}],\Sigma_{t}^{\epsilon_{0}}}\leq C_{l}\delta_{0}(1+t)^{-1}[1+\log(1+t)]
\end{align*}
and:
\begin{align*}
 \max_{i}\|\leftexp{(R_{i})}{\slashed{\pi}}_{1}\|_{\infty,[l_{*}-1],\Sigma_{t}^{\epsilon_{0}}}\leq C_{l}\delta_{0}(1+t)^{-2}[1+\log(1+t)]
\end{align*}
Then for any $S_{t,u}$ 1-form $\xi$ we have:
\begin{align*}
 \max_{m\leq k,n\leq l-m}\max_{i_{1}...i_{n}}\|\leftexp{(i_{1}...i_{n})}{c}_{m,n}[\xi]\|_{L^{2}(\Sigma_{t}^{\epsilon_{0}})}\\
\leq C_{l}\delta_{0}(1+t)^{-2}[1+\log(1+t)]\|\xi\|_{2,[k,l-1],\Sigma_{t}^{\epsilon_{0}}}\\
+C_{l}(1+t)^{-1}\|\xi\|_{\infty,[k,l_{*}],\Sigma_{t}^{\epsilon_{0}}}\{\mathcal{T}_{[k-1,l]}+(1+t)\max_{i}\|\leftexp{(R_{i})}{\slashed{\pi}}_{1}\|
_{2,[l-1],\Sigma_{t}^{\epsilon_{0}}}\\
+\delta_{0}(1+t)^{-1}[1+\log(1+t)]\max_{i}\|\leftexp{(R_{i})}{\slashed{\pi}}\|_{2,[l-2],\Sigma_{t}^{\epsilon_{0}}}\}
\end{align*}
provided that $\delta_{0}$ is suitably small (depending on $l$).
     
$Proof$. The assumptions of the present lemma include those of Lemma 11.11 with $(k-1,l_{*})$ in the role of $(m,l)$. Therefore we have:
\begin{align}
 \|\leftexp{(T)}{\slashed{\pi}}_{1}\|_{\infty,[k-1,l_{*}-1],\Sigma_{t}^{\epsilon_{0}}}\leq C_{l}\delta_{0}(1+t)^{-2}[1+\log(1+t)]
\end{align}
Also, the assumptions of the present lemma include those of Lemma 11.13 with $k-1$ in role of $m$. Therefore we have:
\begin{align}
 \|\leftexp{(T)}{\slashed{\pi}}_{1}\|_{2,[k-1,l-1],\Sigma_{t}^{\epsilon_{0}}}\\\notag
\leq C_{l}\delta_{0}(1+t)^{-1}\{\mathcal{T}_{[k-1,l]}+\delta_{0}(1+t)^{-1}[1+\log(1+t)]\cdot\\\notag
(\max_{i}\|\leftexp{(R_{i})}{\slashed{\pi}}\|_{2,[l-2],\Sigma_{t}^{\epsilon_{0}}}+(1+t)\max_{i}\|\leftexp{(R_{i})}{\slashed{\pi}}_{1}\|
_{2,[l-2],\Sigma_{t}^{\epsilon_{0}}})\}
\end{align}
Then, as in establishing Lemma 11.14, we just use the expressions given by Lemma 11.9 and proceed to a direct calculation. The lemma then follows. $\qed$

\subsection{$L^{\infty}$ Estimates for $T\hat{T}^{i}$}
      Given non-negative integers $m,n$, we denote by $\textbf{E}^{Q}_{m,n}$ the bootstrap assumption that there is a constant $C$ independent of $s$ such
that for all $t\in[0,s]$:
\begin{align*}
 \textbf{E}^{Q}_{m,n}\quad:\quad \max_{\alpha;i_{1}...i_{n}}\|R_{i_{n}}...R_{i_{1}}(T)^{m}Q\psi_{\alpha}\|_{L^{\infty}(\Sigma_{t}^{\epsilon_{0}})}
\leq C\delta_{0}(1+t)^{-1}
\end{align*}
We then denote by $\textbf{E}^{Q}_{\{l\}}$ the conjunction of the assumptions $\textbf{E}^{Q}_{m,n}$ corresponding to the triangle (11.8). The constant $C$ then 
depends on $l$ only and shall be denoted by $C_{l}$. Note that the assumptions $\textbf{E}^{Q}_{0,n}$ coincide with the assumptions $\slashed{\textbf{E}}^{Q}_{n}$,
and $\textbf{E}^{Q}_{\{l\}}$ contain $\slashed{\textbf{E}}^{Q}_{[l]}$.

      Given non-negative integers $m,l$, $m\leq l$, we denote by $\textbf{M}_{[m,l]}$ the bootstrap assumption that there is a constant $C_{l}$ independent
of $s$ such that for all $t\in[0,s]$:
\begin{align*}
 \textbf{M}_{[m,l]}\quad:\quad \|\mu-1\|_{\infty,[m,l],\Sigma_{t}^{\epsilon_{0}}}\leq C_{l}\delta_{0}[1+\log(1+t)]
\end{align*}
This is equivalent, modulo $\textbf{E}_{\{l\}}$, to the assumption that there is a constant $C_{l}$ independent of $s$ such that for all $t\in[0,s]$:
\begin{align}
 \|\kappa-1\|_{\infty,[m,l],\Sigma_{t}^{\epsilon_{0}}}\leq C_{l}\delta_{0}[1+\log(1+t)]
\end{align}

$\textbf{Proposition 11.1}$ Let hypotheses $\textbf{H0}$, $\textbf{H1}$, $\textbf{H2}^{\prime}$ and the estimate (6.177) hold. Let also the bootstrap assumptions 
$\textbf{E}_{\{l+1\}}$, $\textbf{E}^{Q}_{\{l\}}$, and $\slashed{\textbf{X}}_{[l]}$, hold for some non-negative integer $l$. Moreover, let the bootstrap
assumption $\textbf{M}_{[m,l+1]}$ hold for some non-negative integer $m\leq l$. Then if $\delta_{0}$ is suitably small (depending on $l$) we have:
\begin{align*}
 \max_{i}\|(T)^{k+1}\hat{T}^{i}\|_{\infty,[l-k],\Sigma_{t}^{\epsilon_{0}}}\leq C_{l}\delta_{0}(1+t)^{-1}[1+\log(1+t)]
\end{align*}
for all $k=0,...,m$.

$Proof$. The assumptions of the present proposition include those of Proposition 10.1 with $l+1$ in the role of $l$. Therefore by the corollaries of Proposition 
10.1, we have, with $l+1$ in the role of $l$:
\begin{align}
\|\psi_{L}-h_{0}\|_{\infty,[l+1],\Sigma_{t}^{\epsilon_{0}}}\leq C_{l}\delta_{0}(1+t)^{-1}\\
 \|\psi_{\hat{T}}\|_{\infty,[l+1],\Sigma_{t}^{\epsilon_{0}}}\leq C_{l}\delta_{0}(1+t)^{-1}\\
\|\slashed{\psi}\|_{\infty,[l+1],\Sigma_{t}^{\epsilon_{0}}}\leq C_{l}\delta_{0}(1+t)^{-1}
\end{align}
and:
\begin{align}
 \|\omega_{L\hat{T}}\|_{\infty,[l],\Sigma_{t}^{\epsilon_{0}}}\leq C_{l}\delta_{0}(1+t)^{-2}\\
\|\slashed{\omega}_{L}\|_{\infty,[l],\Sigma_{t}^{\epsilon_{0}}}\leq C_{l}\delta_{0}(1+t)^{-2}\\
\|\slashed{\omega}\|_{\infty,[l],\Sigma_{t}^{\epsilon_{0}}}\leq C_{l}\delta_{0}(1+t)^{-2}
\end{align}
Moreover, from the definition:
\begin{align*}
 \slashed{\omega}_{\hat{T}}=\hat{T}^{i}\slashed{d}\psi_{i}
\end{align*}
we readily obtain:
\begin{align}
 \|\slashed{\omega}_{\hat{T}}\|_{\infty,[l],\Sigma_{t}^{\epsilon_{0}}}\leq C_{l}\delta_{0}(1+t)^{-2}\\
\|\omega_{T\hat{T}}\|_{\infty,[l],\Sigma_{t}^{\epsilon_{0}}}\leq C_{l}\delta_{0}(1+t)^{-1}
\end{align}
By the assumption $\textbf{M}_{[0,l+1]}$, we have:
\begin{align}
 \|(q_{T})_{b}\|_{\infty,[l],\Sigma_{t}^{\epsilon_{0}}}\leq C_{l}\delta_{0}(1+t)^{-1}[1+\log(1+t)]
\end{align}
(provided that $\delta_{0}$ is suitably small depending on $l$). This together with Corollary 10.1.d implies:
\begin{align}
 \|q_{T}\|_{\infty,[l],\Sigma_{t}^{\epsilon_{0}}}\leq C_{l}\delta_{0}(1+t)^{-1}[1+\log(1+t)]
\end{align}
In view of (11.12), the statement of the proposition for $k=0$ follows.

      We proceed by induction on $k$ for fixed $m$ and $l$. Let $k\in\{1,...,m\}$ and let:
\begin{align}
 \max_{i}\|(T)^{k'+1}\hat{T}^{i}\|_{\infty,[l-k'],\Sigma_{t}^{\epsilon_{0}}}\leq C_{l}\delta_{0}(1+t)^{-1}[1+\log(1+t)]
\end{align}
hold for all $k'=0,...,k-1$.   

By this induction hypothesis and $\textbf{E}_{\{l+1\}}$ we can easily obtain:
\begin{align}
 \|(T)^{k'}\psi_{\hat{T}}\|_{\infty,[l+1-k'],\Sigma_{t}^{\epsilon_{0}}}\leq C_{l}\delta_{0}(1+t)^{-1}
\end{align}
for all $k'\in\{1,...,k\}$.

This together with (11.100) implies:
\begin{align}
 \|\psi_{\hat{T}}\|_{\infty,[k,l+1],\Sigma_{t}^{\epsilon_{0}}}\leq C_{l}\delta_{0}(1+t)^{-1}
\end{align}
Since
\begin{align*}
 \omega_{T\hat{T}}=\hat{T}^{i}(T\psi_{i}),\quad \omega_{L\hat{T}}=\hat{T}^{i}(L\psi_{i}),\quad \slashed{\omega}_{\hat{T}}=\hat{T}^{i}(\slashed{d}\psi_{i})
\end{align*}
these are analogous to $\psi_{\hat{T}}$ with $T\psi_{i}$, $L\psi_{i}$, $\slashed{d}\psi_{i}$, in the role of $\psi_{i}$, respectively.
We thus obtain:
\begin{align}
 \|\omega_{T\hat{T}}\|_{\infty,[k,l],\Sigma_{t}^{\epsilon_{0}}}\leq C_{l}\delta_{0}(1+t)^{-1}\\
\|\omega_{L\hat{T}}\|_{\infty,[k,l],\Sigma_{t}^{\epsilon_{0}}}\leq C_{l}\delta_{0}(1+t)^{-2}\\
\|\slashed{\omega}_{\hat{T}}\|_{\infty,[k,l],\Sigma_{t}^{\epsilon_{0}}}\leq C_{l}\delta_{0}(1+t)^{-2}
\end{align}
Similarly, by the relation:
\begin{align}
 L^{i}=-\alpha\hat{T}^{i}-\psi_{i}
\end{align}
and $\textbf{E}_{\{l+1\}}$, we have:
\begin{align}
 \max_{i}\|(T)^{k'+1}L^{i}\|_{\infty,[l-k'],\Sigma_{t}^{\epsilon_{0}}}\leq C_{l}\delta_{0}(1+t)^{-1}[1+\log(1+t)]
\end{align}
for all $k'=0,...,k-1$.

Then we have:
\begin{align}
 \|(T)^{k'}\psi_{L}\|_{\infty,[l+1-k'],\Sigma_{t}^{\epsilon_{0}}}\leq C_{l}\delta_{0}(1+t)^{-1}
\end{align}
for all $k^{\prime}=0,...,k$，
which, together with (11.99), yields:
\begin{align}
 \|\psi_{L}\|_{\infty,[k,l+1],\Sigma_{t}^{\epsilon_{0}}}\leq C_{l}\delta_{0}(1+t)^{-1}
\end{align}
and then:
\begin{align}
 \|\slashed{\omega}_{L}\|_{\infty,[k,l],\Sigma_{t}^{\epsilon_{0}}}\leq C_{l}\delta_{0}(1+t)^{-2}
\end{align}
Also we have:
\begin{align}
 \|\slashed{\psi}\|_{\infty,[k,l+1],\Sigma_{t}^{\epsilon_{0}}}\leq C_{l}\delta_{0}(1+t)^{-1}\\
\|\slashed{\omega}\|_{\infty,[k,l],\Sigma_{t}^{\epsilon_{0}}}\leq C_{l}\delta_{0}(1+t)^{-2}
\end{align}
Then we have:
\begin{align}
 \|\slashed{k}\|_{\infty,[k,l],\Sigma_{t}^{\epsilon_{0}}}\leq C_{l}\delta_{0}(1+t)^{-2}\\
\|\kappa^{-1}\zeta\|_{\infty,[k,l],\Sigma_{t}^{\epsilon_{0}}}\leq C_{l}\delta_{0}(1+t)^{-2}
\end{align}
Therefore
\begin{align}
 \|\kappa\slashed{k}\|_{\infty,[k,l],\Sigma_{t}^{\epsilon_{0}}}\leq C_{l}\delta_{0}(1+t)^{-2}[1+\log(1+t)]\\
\|\zeta\|_{\infty,[k,l],\Sigma_{t}^{\epsilon_{0}}}\leq C_{l}\delta_{0}(1+t)^{-2}[1+\log(1+t)]
\end{align}
By $\textbf{E}_{\{l+1\}}$ and $\textbf{M}_{[k,l+1]}$ we obtain:
\begin{align}
 \|(q_{T})_{b}\|_{\infty,[k,l],\Sigma_{t}^{\epsilon_{0}}}\leq C_{l}\delta_{0}(1+t)^{-1}[1+\log(1+t)]
\end{align}
To complete the inductive step we need a similar estimate for $q_{T}$. This estimate would follow from Lemma 11.7 if we can show that:
\begin{align}
 \|\leftexp{(T)}{\slashed{\pi}}^{\prime}\|_{\infty,[k-1,l-1],\Sigma_{t}^{\epsilon_{0}}}\leq C_{l}\delta_{0}(1+t)^{-1}[1+\log(1+t)]
\end{align}
because the assumptions of Lemma 11.7, which are those of Lemma 11.5, hold by virtue of Corollary 10.1.d and $\textbf{M}_{[k,l+1]}$.
From the definition (11.51) and the estimate (11.124) we see that (11.127) would in turn follow if we can show that:
\begin{align}
 \|\chi'\|_{\infty,[k-1,l-1],\Sigma_{t}^{\epsilon_{0}}}\leq C_{l}\delta_{0}(1+t)^{-1}
\end{align}
A stronger estimate is in fact provided by the following lemma:

$\textbf{Lemma 11.16}$ Let the assumptions of Proposition 11.1 hold. Let also the inductive hypothesis (11.109) hold, for some $k\in\{1,...,m\}$. Then,
provided that $\delta_{0}$ is suitably small (depending on $l$), the following estimate holds:
\begin{align*}
 \|\chi'\|_{\infty,[k,l],\Sigma_{t}^{\epsilon_{0}}}\leq C_{l}\delta_{0}(1+t)^{-2}[1+\log(1+t)]
\end{align*}
$Proof$. As remarked above, under the present assumptions Proposition 10.1 and all its corollaries hold with $l+1$ in the role of $l$. From the definition 
(11.50) and Corollary 10.1.d we obtain:
\begin{align}
 \|\chi'\|_{\infty,[l],\Sigma_{t}^{\epsilon_{0}}}\leq C_{l}\delta_{0}(1+t)^{-2}[1+\log(1+t)]
\end{align}
We proceed, by induction on $k'$, to show that:
\begin{align}
 \|(\slashed{\mathcal{L}}_{T})^{k'}\chi'\|_{\infty,[l-k'],\Sigma_{t}^{\epsilon_{0}}}\leq C_{l}\delta_{0}(1+t)^{-2}[1+\log(1+t)]
\end{align}
for all $k'\in\{1,...,k\}$，
which, together with (11.129), yields the conclusion of the lemma. In fact we shall show that (11.109) together with the estimate:
\begin{align}
 \|\chi'\|_{\infty,[k-1,l],\Sigma_{t}^{\epsilon_{0}}}\leq C_{l}\delta_{0}(1+t)^{-2}[1+\log(1+t)]
\end{align}
implies
\begin{align}
 \|(\slashed{\mathcal{L}}_{T})^{k}\chi'\|_{\infty,[l-k],\Sigma_{t}^{\epsilon_{0}}}\leq C_{l}\delta_{0}(1+t)^{-2}[1+\log(1+t)]
\end{align}
     From the definition (11.50) we have:
\begin{align}
 \slashed{\mathcal{L}}_{T}\chi'=\slashed{\mathcal{L}}_{T}\chi-\frac{\leftexp{(T)}{\slashed{\pi}}}{1-u+t}-\frac{\slashed{g}}{(1-u+t)^{2}}
\end{align}
Now, $\slashed{\mathcal{L}}_{T}\chi$ is given by the formula (3.125)-(3.126). Substituting in (11.133) this formula, and also substituting for 
$\leftexp{(T)}{\slashed{\pi}}$ in terms of $\chi'$ from (11.51)-(11.52), we obtain a formula for $\slashed{\mathcal{L}}_{T}\chi'$. We apply
$(\slashed{\mathcal{L}}_{T})^{k-1}$ to this formula and then take the $\|\quad\|_{\infty,[l-k],\Sigma_{t}^{\epsilon_{0}}}$ norm. The contribution
of any given term on the right-hand side, will then be bounded by the $\|\quad\|_{\infty,[k-1,l-1],\Sigma_{t}^{\epsilon_{0}}}$ norm of that term. 
Now from (11.51), (11.131) together with (11.124) and $\textbf{M}_{[k,l+1]}$ imply
\begin{align}
 \|\leftexp{(T)}{\slashed{\pi}}^{\prime}\|_{\infty,[k-1,l],\Sigma_{t}^{\epsilon_{0}}}\leq C_{l}\delta_{0}(1+t)^{-2}[1+\log(1+t)]^{2}
\end{align}
By Lemma 11.7 with $k-1$ in the role of $m$ we then obtain:
\begin{align}
 \|\leftexp{(T)}{\slashed{\pi}}+2(1-u+t)^{-1}\slashed{g}\|_{\infty,[l],\Sigma_{t}^{\epsilon_{0}}}\leq C_{l}\delta_{0}(1+t)^{-1}[1+\log(1+t)]\\\notag
\|\slashed{\mathcal{L}}_{T}\leftexp{(T)}{\slashed{\pi}}-2(1-u+t)^{-1}\slashed{g}\|_{\infty,[l-1],\Sigma_{t}^{\epsilon_{0}}}\leq C_{l}\delta_{0}
(1+t)^{-1}[1+\log(1+t)]\\\notag
\max_{2\leq k'\leq k-1}\|(\slashed{\mathcal{L}}_{T})^{k'}\leftexp{(T)}{\slashed{\pi}}\|_{\infty,[l-k'],\Sigma_{t}^{\epsilon_{0}}}
\leq C_{l}\delta_{0}(1+t)^{-1}[1+\log(1+t)]
\end{align}
These estimates, together with Corollary 10.1.d, imply that for any $S_{t,u}$ tensorfield $\xi$ we have:
\begin{align}
 \|\xi\cdot\slashed{g}\|_{\infty,[k-1,l-1],\Sigma_{t}^{\epsilon_{0}}}, \|\xi\cdot\slashed{g}^{-1}\|_{\infty,[k-1,l-1],\Sigma_{t}^{\epsilon_{0}}}
\leq C_{l}\|\xi\|_{\infty,[k-1,l-1],\Sigma_{t}^{\epsilon_{0}}}
\end{align}
provided that $\delta_{0}$ is suitably small (depending on $l$).
This together with the estimates from (11.99) to (11.125) as well as Lemma 11.14 implies (11.132). So the lemma follows. $\qed$

     We now resume the proof of Proposition 11.1. Using Lemma 11.16 we obtain the estimate:
\begin{align}
 \|q_{T}\|_{\infty,[k,l],\Sigma_{t}^{\epsilon_{0}}}\leq C_{l}\delta_{0}(1+t)^{-1}[1+\log(1+t)]
\end{align}
Using this and $\textbf{M}_{[k,l+1]}$ we deduce, from the first statement of Lemma 11.4, the estimates:
\begin{align}
 \|p_{T,j}\|_{\infty,[k-j,l-j],\Sigma_{t}^{\epsilon_{0}}}\leq C_{l}\delta_{0}^{2}(1+t)^{-2}[1+\log(1+t)]^{2}
\end{align}
for all $j\in\{0,...,k\}$
\begin{align}
 \|q_{T,j}\|_{\infty,[k-j,l-j],\Sigma_{t}^{\epsilon_{0}}}\leq C_{l}\delta_{0}(1+t)^{-1}[1+\log(1+t)]
\end{align}
for all $j\in\{0,...,k\}$

We write the second statement of Lemma 11.4 in the form:
\begin{align}
 r^{n}_{T,k}=\sum_{i=0}^{k-1-n}N_{i}(n)\sum_{j'=0}^{k-1-n-i}(\slashed{\mathcal{L}}_{T})^{k-1-n-j'}(\kappa q_{T,j'})
\end{align}
Here $n\in\{0,...,k-1\}$. Taking the $\|\quad\|_{\infty,[l-k+n],\Sigma_{t}^{\epsilon_{0}}}$ norm of this we obtain, using (11.139) and $\textbf{M}_{[k,l+1]}$,
\begin{align}
 \|r^{n}_{T,k}\|_{\infty,[l-k],\Sigma_{t}^{\epsilon_{0}}}\\\notag
\leq \sum_{i=0}^{k-1-n}N_{i}(n)\sum_{j'=0}^{k-1-n-i}\|\kappa q_{T,j'}\|_{\infty,[k-1-n-j',l-1-j'],\Sigma_{t}^{\epsilon_{0}}}\\\notag
\leq C_{l}\max_{j'\in\{0,...,k-1-n\}}\|\kappa q_{T,j'}\|_{\infty,[k-1-n-j',l-1-j'],\Sigma_{t}^{\epsilon_{0}}}\\\notag
\leq C_{l}\delta_{0}(1+t)^{-1}[1+\log(1+t)]^{2}
\end{align}
for all $n\in\{0,...,k-1\}$.

Now according to Lemma 11.3:
\begin{align}
 (T)^{k+1}\hat{T}^{i}=p_{T,k}\hat{T}^{i}+q_{T,k}\cdot\slashed{d}x^{i}+\sum_{n=0}^{k-1}r^{n}_{T,k}\cdot\slashed{d}(T)^{n}\hat{T}^{i}
\end{align}
We take the $\|\quad\|_{\infty,[l-k],\Sigma_{t}^{\epsilon_{0}}}$ norm of this. From the above estimates, we have:
\begin{align}
 \max_{i}\|p_{T,k}\hat{T}\|_{\infty,[l-k],\Sigma_{t}^{\epsilon_{0}}}\leq C_{l}\delta_{0}^{2}(1+t)^{-2}[1+\log(1+t)]^{2}\\
\max_{i}\|q_{T,k}\cdot\slashed{d}x^{i}\|_{\infty,[l-k],\Sigma_{t}^{\epsilon_{0}}}\leq C_{l}\delta_{0}(1+t)^{-1}[1+\log(1+t)]\\
\max_{i}\|\sum_{n=0}^{k-1}r^{n}_{T,k}\cdot\slashed{d}(T)^{n}\hat{T}^{i}\|_{\infty,[l-k],\Sigma_{t}^{\epsilon_{0}}}
\leq C_{l}\delta_{0}(1+t)^{-2}[1+\log(1+t)]^{3}
\end{align}
These imply:
\begin{align}
 \max_{i}\|(T)^{k+1}\hat{T}^{i}\|_{\infty,[l-k],\Sigma_{t}^{\epsilon_{0}}}\leq C_{l}\delta_{0}(1+t)^{-1}[1+\log(1+t)]
\end{align}
This completes the proof of Proposition 11.1. $\qed$

     The foregoing proposition has the following corollaries.

$\textbf{Corollary 11.1.a}$ Under the assumptions of Proposition 11.1 the coefficients of the expression for $(T)^{k+1}\hat{T}^{i}$, 
$k=0,...,m$, of Lemma 11.3 satisfy:
\begin{align*}
 \|p_{T,k}\|_{\infty,[m-k,l-k],\Sigma_{t}^{\epsilon_{0}}}\leq C_{l}\delta_{0}^{2}(1+t)^{-2}[1+\log(1+t)]^{2}
\end{align*}
for all $k\in\{0,...,m\}$
\begin{align*}
 \|q_{T,k}\|_{\infty,[m-k,l-k],\Sigma_{t}^{\epsilon_{0}}}\leq C_{l}\delta_{0}(1+t)^{-1}[1+\log(1+t)]
\end{align*}
for all $k\in\{0,...,m\}$，
and:
\begin{align*}
 \|r^{n}_{T,k}\|_{\infty,[l-k+n],\Sigma_{t}^{\epsilon_{0}}}\leq C_{l}\delta_{0}(1+t)^{-1}[1+\log(1+t)]^{2}
\end{align*}
for all $k\in\{0,...,m\}$, $n\in\{0,...,k-1\}$.

Also, 
\begin{align*}
 \max_{i}\|(T)^{k+1}L^{i}\|_{\infty,[l-k],\Sigma_{t}^{\epsilon_{0}}}\leq C_{l}\delta_{0}(1+t)^{-1}[1+\log(1+t)]
\end{align*}
for all $k=0,...,m$.

$\textbf{Corollary 11.1.b}$ Under the assumptions of Proposition 11.1 the following estimates hold:
\begin{align*}
 \|\psi_{\hat{T}}\|_{\infty,[m+1,l+1],\Sigma_{t}^{\epsilon_{0}}}\leq C_{l}\delta_{0}(1+t)^{-1}\\
\|\psi_{L}\|_{\infty,[m+1,l+1],\Sigma_{t}^{\epsilon_{0}}}\leq C_{l}\delta_{0}(1+t)^{-1}\\
\|\slashed{\psi}\|_{\infty,[m+1,l+1],\Sigma_{t}^{\epsilon_{0}}}\leq C_{l}\delta_{0}(1+t)^{-1}
\end{align*}
and:
\begin{align*}
 \|\omega_{T\hat{T}}\|_{\infty,[m,l],\Sigma_{t}^{\epsilon_{0}}}\leq C_{l}\delta_{0}(1+t)^{-1}\\
\|\omega_{L\hat{T}}\|_{\infty,[m+1,l+1],\Sigma_{t}^{\epsilon_{0}}}\leq C_{l}\delta_{0}(1+t)^{-2}\\
\|\slashed{\omega}_{\hat{T}}\|_{\infty,[m,l],\Sigma_{t}^{\epsilon_{0}}}\leq C_{l}\delta_{0}(1+t)^{-2}\\
\|\slashed{\omega}_{L}\|_{\infty,[m,l],\Sigma_{t}^{\epsilon_{0}}}\leq C_{l}\delta_{0}(1+t)^{-2}\\
\|\slashed{\omega}\|_{\infty,[m,l],\Sigma_{t}^{\epsilon_{0}}}\leq C_{l}\delta_{0}(1+t)^{-2}
\end{align*}
and:
\begin{align*}
 \|\slashed{k}\|_{\infty,[m,l],\Sigma_{t}^{\epsilon_{0}}}\leq C_{l}\delta_{0}(1+t)^{-2}\\
\|\kappa^{-1}\zeta\|_{\infty,[m,l],\Sigma_{t}^{\epsilon_{0}}}\leq C_{l}\delta_{0}(1+t)^{-2}\\
\|\kappa\slashed{k}\|_{\infty,[m,l],\Sigma_{t}^{\epsilon_{0}}}\leq C_{l}\delta_{0}(1+t)^{-2}[1+\log(1+t)]\\
\|\zeta\|_{\infty,[m,l],\Sigma_{t}^{\epsilon_{0}}}\leq C_{l}\delta_{0}(1+t)^{-2}[1+\log(1+t)]
\end{align*}
    The next corollary is obtained from Lemma 11.16 by substituting Proposition 11.1 itself in place of the inductive hypothesis (11.109), and then substituting 
the result in Lemma 11.7 and Lemma 11.11.

$\textbf{Corollary 11.1.c}$ Under the assumptions of Proposition 11.1 we have:
\begin{align*}
 \|\chi'\|_{\infty,[m,l],\Sigma_{t}^{\epsilon_{0}}}\leq C_{l}\delta_{0}(1+t)^{-2}[1+\log(1+t)]
\end{align*}
Moreover, we have:
\begin{align*}
 \|\leftexp{(T)}{\slashed{\pi}}+2(1-u+t)^{-1}\slashed{g}\|_{\infty,[l],\Sigma_{t}^{\epsilon_{0}}}\leq C_{l}\delta_{0}(1+t)^{-1}[1+\log(1+t)]\\
\|\slashed{\mathcal{L}}_{T}\leftexp{(T)}{\slashed{\pi}}-2(1-u+t)^{-2}\slashed{g}\|_{\infty,[l-1],\Sigma_{t}^{\epsilon_{0}}}
\leq C_{l}\delta_{0}(1+t)^{-1}[1+\log(1+t)]\\
\max_{2\leq k\leq m}\|(\slashed{\mathcal{L}}_{T})^{k}\leftexp{(T)}{\slashed{\pi}}\|_{\infty,[l-k],\Sigma_{t}^{\epsilon_{0}}}
\leq C_{l}\delta_{0}(1+t)^{-1}[1+\log(1+t)]
\end{align*}
Also:
\begin{align*}
 \|\Lambda\|_{\infty,[m,l],\Sigma_{t}^{\epsilon_{0}}}\leq C_{l}\delta_{0}(1+t)^{-1}[1+\log(1+t)]
\end{align*}
and:
\begin{align*}
 \|\leftexp{(T)}{\slashed{\pi}}_{1}\|_{\infty,[m,l-1],\Sigma_{t}^{\epsilon_{0}}}\leq C_{l}\delta_{0}(1+t)^{-2}[1+\log(1+t)]
\end{align*}

\subsection{$L^{2}$ Estimates for $T\hat{T}^{i}$}
     We proceed with the $L^{2}$ estimates. Given non-negative integers $m$, $n$, we denote by $\mathcal{W}^{Q}_{m,n}$ the quantity:
\begin{align}
 \mathcal{W}^{Q}_{m,n}=\max_{\alpha;i_{1}...i_{n}}\|R_{i_{n}}...R_{i_{1}}(T)^{m}Q\psi_{\alpha}\|_{L^{2}(\Sigma_{t}^{\epsilon_{0}})}
\end{align}
We then denote by $\mathcal{W}^{Q}_{\{l\}}$ the sum of the quantities $\mathcal{W}^{Q}_{m,n}$ corresponding to the triangle (11.8). Note that the quantities
$\mathcal{W}^{Q}_{0,n}$ coincide with the quantities $\mathcal{W}^{Q}_{n}$, and $\mathcal{W}^{Q}_{\{l\}}$ dominate $\mathcal{W}^{Q}_{[l]}$. 
Moreover, we have:
\begin{align}
 \max_{\alpha}\|L\psi_{\alpha}\|_{2,\{l\},\Sigma_{t}^{\epsilon_{0}}}\leq (1+t)^{-1}\mathcal{W}^{Q}_{\{l\}}
\end{align}
Given non-negative integers $m$, $l$, $m\leq l$ we denote:
\begin{align}
 \mathcal{B}_{[m,l]}=\|\mu-1\|_{2,[m,l],\Sigma_{t}^{\epsilon_{0}}}
\end{align}
Note that $\mathcal{B}_{[0,l]}$ coincide with $\mathcal{B}_{[l]}$. Since $\mu=\alpha\kappa$, it follows from Lemma 11.1 that under the bootstrap
assumptions $\textbf{E}_{\{l_{*}\}}$ and $\textbf{M}_{[m,l_{*}]}$ we have:
\begin{align}
 \|\kappa-1\|_{2,[m,l],\Sigma_{t}^{\epsilon_{0}}}\leq C_{l}\{\mathcal{B}_{[m,l]}+[1+\log(1+t)]\mathcal{W}_{\{l\}}\}
\end{align}

$\textbf{Proposition 11.2}$ Let the hypotheses $\textbf{H0}$, $\textbf{H1}$, $\textbf{H2}^{\prime}$ and the estimate (6.177) hold. Let also the 
bootstrap assumptions $\textbf{E}_{\{l_{*}+1\}}$, $\textbf{E}^{Q}_{\{l_{*}\}}$, and $\slashed{\textbf{X}}_{[l_{*}]}$, hold for some non-negative
integer $l$. Moreover, let the bootstrap assumption $\textbf{M}_{[m,l_{*}+1]}$ hold for some non-negative integer $m\leq l$ (where by definition 
$\textbf{M}_{[m,l_{*}+1]}$ coincides with $\textbf{M}_{\{l_{*}+1\}}=\textbf{M}_{[l_{*}+1,l_{*}+1]}$ for $m\geq l_{*}+1$). Then if $\delta_{0}$
is suitably small  (depending on $l$) we have:
\begin{align*}
 \sum_{k=0}^{m}\max_{i}\|(T)^{k+1}\hat{T}^{i}\|_{2,[l-k],\Sigma_{t}^{\epsilon_{0}}}\\
\leq C_{l}(1+t)^{-1}\{\mathcal{B}_{[m,l+1]}+\delta_{0}[1+\log(1+t)][\mathcal{Y}_{0}+(1+t)\mathcal{A}_{[l-1]}]\\
+[1+\log(1+t)][\mathcal{W}_{\{l+1\}}+\delta_{0}(1+t)^{-2}[1+\log(1+t)]^{2}\mathcal{W}^{Q}_{\{l-1\}}]\}
\end{align*}
$Proof$. The assumptions of the present proposition coincide with those of Proposition 11.1 with $(m,l_{*})$ in the role of $(m,l)$. Therefore under 
the present assumptions Proposition 11.1 and its corollaries hold with $(m,l_{*})$ in the role of $(m,l)$. Moreover, the present assumptions include
those of Proposition 10.1 with $l_{*}+1$ in the role of $l$, as well as those of Proposition 10.2 with $l+1$ in the role of $l$. Therefore under the 
present assumptions Proposition 10.1 and its all corollaries follow with $l_{*}+1$ in the role of $l$, and Proposition 10.2 and its all corollaries 
follow with $l+1$ in the role of $l$. In particular, we have (Corollaries 10.2.b and 10.2.g with $l$ replaced by $l+1$):
\begin{align}
 \|\psi_{L}\|_{2,[l+1],\Sigma_{t}^{\epsilon_{0}}}\leq C_{l}\{\mathcal{W}_{[l+1]}+\delta_{0}(1+t)^{-1}[\mathcal{Y}_{0}+(1+t)\mathcal{A}_{[l]}]\}\\
\|\psi_{\hat{T}}\|_{2,[l+1],\Sigma_{t}^{\epsilon_{0}}}\leq C_{l}\{\mathcal{W}_{[l+1]}+\delta_{0}(1+t)^{-1}[\mathcal{Y}_{0}+(1+t)\mathcal{A}_{[l]}]\}\\
\|\slashed{\psi}\|_{2,[l+1],\Sigma_{t}^{\epsilon_{0}}}\leq C_{l}\{\mathcal{W}_{[l+1]}+\delta_{0}(1+t)^{-1}[\mathcal{Y}_{0}+(1+t)\mathcal{A}_{[l]}]\}
\end{align}
and:
\begin{align}
 \|\omega_{L\hat{T}}\|_{2,[l],\Sigma_{t}^{\epsilon_{0}}}\leq C_{l}(1+t)^{-1}\{\mathcal{W}^{Q}_{[l]}+\delta_{0}(1+t)^{-1}[\mathcal{Y}_{0}
+(1+t)\mathcal{A}_{[l-1]}+\mathcal{W}_{[l]}]\}\\
\|\slashed{\omega}_{L}\|_{2,[l],\Sigma_{t}^{\epsilon_{0}}}\leq C_{l}(1+t)^{-1}\{\mathcal{W}_{[l+1]}+\delta_{0}(1+t)^{-1}
[\mathcal{Y}_{0}+(1+t)\mathcal{A}_{[l-1]}]\}\\
\|\slashed{\omega}\|_{2,[l],\Sigma_{t}^{\epsilon_{0}}}\leq C_{l}(1+t)^{-1}\{\mathcal{W}_{[l+1]}+\delta_{0}(1+t)^{-1}
[\mathcal{Y}_{0}+(1+t)\mathcal{A}_{[l-1]}]\}
\end{align}
Moreover, we readily obtain:
\begin{align}
 \|\slashed{\omega}_{\hat{T}}\|_{2,[l],\Sigma_{t}^{\epsilon_{0}}}\leq C_{l}(1+t)^{-1}\{\mathcal{W}_{[l+1]}+\delta_{0}(1+t)^{-1}
[\mathcal{Y}_{0}+(1+t)\mathcal{A}_{[l-1]}]\}\\
\|\omega_{T\hat{T}}\|_{2,[l],\Sigma_{t}^{\epsilon_{0}}}\leq C_{l}\{\mathcal{W}_{[l+1]}+\delta_{0}(1+t)^{-1}[\mathcal{Y}_{0}+(1+t)
\mathcal{A}_{[l-1]}]\}
\end{align}

     Let us define:
\begin{align}
 \mathcal{U}_{k,l}=\sum_{k'=0}^{k}\max_{i}\|(T)^{k'+1}\hat{T}^{i}\|_{2,[l-k'],\Sigma_{t}^{\epsilon_{0}}}
\end{align}
By the definition (11.14),
\begin{align}
 \|(q_{T})_{b}\|_{2,[l],\Sigma_{t}^{\epsilon_{0}}}\leq C(1+t)^{-1}\|\kappa-1\|_{2,[l+1],\Sigma_{t}^{\epsilon_{0}}}
\end{align}
This together with (10.145) as well as the $L^{\infty}$ estimates (11.126) and Corollary 10.1.d, in turn yields:
\begin{align}
 \|q_{T}\|_{2,[l],\Sigma_{t}^{\epsilon_{0}}}\leq C_{l}(1+t)^{-1}\{\|\kappa-1\|_{2,[l+1],\Sigma_{t}^{\epsilon_{0}}}\\\notag
+\delta_{0}[1+\log(1+t)][\mathcal{Y}_{0}+(1+t)\mathcal{A}_{[l-1]}+\mathcal{W}_{[l]}]\}
\end{align}
From this and Corollary 10.2.a,
\begin{align}
 \mathcal{U}_{0,l}\leq C_{l}(1+t)^{-1}\{\|\kappa-1\|_{2,[l+1],\Sigma_{t}^{\epsilon_{0}}}\\\notag
+\delta_{0}[1+\log(1+t)][\mathcal{Y}_{0}+(1+t)\mathcal{A}_{[l-1]}+\mathcal{W}_{[l]}]\}
\end{align}

     We shall prove the proposition by deriving a recursive, in $k$, inequality, for the quantities $\mathcal{U}_{k,l}$, for fixed $l$. We first establish $L^{2}$ 
estimates for the $T$-derivatives of up to the $k$th order of the $\psi$ and $\omega$ components, in terms of $\mathcal{U}_{k-1,l}$. As we did in proving 
Proposition 11.1, we make use of the following expressions:
\begin{align}
 (T)^{k'}\psi_{\hat{T}}=\hat{T}^{i}(T)^{k'}\psi_{i}+\sum_{j=0}^{k'-1}\frac{k'!}{(j+1)!(k'-1-j)!}((T)^{j+1}\hat{T}^{i})((T)^{k'-1-j}\psi_{i})\\
(T)^{k'}\psi_{L}=(T)^{k'}\psi_{0}+L^{i}(T)^{k'}\psi_{i}+\sum_{j=0}^{k'-1}\frac{k'!}{(j+1)!(k'-1-j)!}((T)^{j+1}L^{i})((T)^{k'-1-j}\psi_{i})\\
(\slashed{\mathcal{L}}_{T})^{k'}\slashed{\psi}=(\slashed{d}x^{i})(T)^{k'}\psi_{i}+\sum_{j=0}^{k'-1}\frac{k'!}{(j+1)!(k'-1-j)!}
(\slashed{d}(T)^{j+1}x^{i})(T)^{k'-1-j}\psi_{i}
\end{align}
We now take the $\|\quad\|_{2,[l+1-k'],\Sigma_{t}^{\epsilon_{0}}}$ norm of the above. 
This together with (11.151)-(11.158), Proposition 10.2 as well as its corollaries and Proposition 11.1 implies:
\begin{align}
 \|\psi_{\hat{T}}\|_{2,[k,l+1],\Sigma_{t}^{\epsilon_{0}}}\leq C_{l}\{\mathcal{W}_{\{l+1\}}+\delta_{0}(1+t)^{-1}[\mathcal{U}_{k-1,l}
+\mathcal{Y}_{0}+(1+t)\mathcal{A}_{[l]}]\}
\end{align}
and:
\begin{align}
 \|\omega_{T\hat{T}}\|_{2,[k,l],\Sigma_{t}^{\epsilon_{0}}}\leq C_{l}\{\mathcal{W}_{\{l+1\}}+\delta_{0}(1+t)^{-1}[\mathcal{U}_{k-1,l-1}
+\mathcal{Y}_{0}+(1+t)\mathcal{A}_{[l-1]}]\}\\
\|\omega_{L\hat{T}}\|_{2,[k,l],\Sigma_{t}^{\epsilon_{0}}}\\\notag
\leq C_{l}(1+t)^{-1}\{\mathcal{W}^{Q}_{\{l\}}+\delta_{0}(1+t)^{-1}[\mathcal{U}_{[k-1,l-1]}
+\mathcal{Y}_{0}+(1+t)\mathcal{A}_{[l-1]}]+\mathcal{W}_{\{l\}}\}\\
\|\slashed{\omega}_{\hat{T}}\|_{2,[k,l],\Sigma_{t}^{\epsilon_{0}}}\leq C_{l}(1+t)^{-1}\{\mathcal{W}_{\{l+1\}}+\delta_{0}(1+t)^{-1}
[\mathcal{U}_{[k-1,l-1]}+\mathcal{Y}_{0}+(1+t)\mathcal{A}_{[l-1]}]\}
\end{align}
Also:
\begin{align}
 \|\psi_{L}\|_{2,[k,l+1],\Sigma_{t}^{\epsilon_{0}}}\leq C_{l}\{\mathcal{W}_{\{l+1\}}+\delta_{0}(1+t)^{-1}
[\mathcal{U}_{k-1,l}+\mathcal{Y}_{0}+(1+t)\mathcal{A}_{[l]}]\}\\
\|\slashed{\omega}_{L}\|_{2,[k,l],\Sigma_{t}^{\epsilon_{0}}}\leq C_{l}(1+t)^{-1}\{\mathcal{W}_{\{l+1\}}+
\delta_{0}(1+t)^{-1}[\mathcal{U}_{k-1,l-1}+\mathcal{Y}_{0}+(1+t)\mathcal{A}_{[l-1]}]\}
\end{align}
and finally:
\begin{align}
 \|\slashed{\psi}\|_{2,[k,l+1],\Sigma_{t}^{\epsilon_{0}}}\leq C_{l}\{\mathcal{W}_{\{l+1\}}+\delta_{0}(1+t)^{-2}\|\kappa-1\|_{2,[k-1,l+1],\Sigma_{t}^{\epsilon_{0}}}\\\notag
+\delta_{0}(1+t)^{-2}[1+\log(1+t)]\mathcal{U}_{k-2,l}\\\notag
+\delta_{0}(1+t)^{-1}[\mathcal{Y}_{0}+(1+t)\mathcal{A}_{[l]}]\}
\end{align}
as well as 
\begin{align}
 \|\slashed{\omega}\|_{2,[k,l],\Sigma_{t}^{\epsilon_{0}}}\leq C_{l}(1+t)^{-1}\{\mathcal{W}_{\{l+1\}}+\delta_{0}(1+t)^{-2}\|\kappa-1\|
_{2,[k-1,l],\Sigma_{t}^{\epsilon_{0}}}\\\notag+\delta_{0}(1+t)^{-2}[1+\log(1+t)]\mathcal{U}_{k-2,l-1}\\\notag
+\delta_{0}(1+t)^{-1}[\mathcal{Y}_{0}+(1+t)\mathcal{A}_{[l-1]}]\}
\end{align}
     (11.173) together with Lemma 11.2.a implies that:
\begin{align}
 \|\slashed{k}\|_{2,[k,l],\Sigma_{t}^{\epsilon_{0}}}\leq C_{l}(1+t)^{-1}\{\mathcal{W}_{\{l+1\}}+\delta_{0}(1+t)^{-2}\|\kappa-1\|
_{2,[k-1,l],\Sigma_{t}^{\epsilon_{0}}}\\\notag+\delta_{0}(1+t)^{-2}[1+\log(1+t)]\mathcal{U}_{k-2,l-1}\\\notag
+\delta_{0}(1+t)^{-1}[\mathcal{Y}_{0}+(1+t)\mathcal{A}_{[l-1]}]\}
\end{align}
Similarly, (11.169) implies:
\begin{align}
 \|\kappa^{-1}\zeta\|_{2,[k,l],\Sigma_{t}^{\epsilon_{0}}}\leq C_{l}(1+t)^{-1}\{\delta_{0}(1+t)^{-1}\mathcal{U}_{k-1,l-1}\\\notag
+\delta_{0}(1+t)^{-1}[\mathcal{Y}_{0}+(1+t)\mathcal{A}_{[l-1]}]+\mathcal{W}_{\{l+1\}}\}
\end{align}
So we have:
\begin{align}
 \|\kappa\slashed{k}\|_{2,[k,l],\Sigma_{t}^{\epsilon_{0}}}\leq C_{l}(1+t)^{-1}\{\delta_{0}(1+t)^{-2}[1+\log(1+t)]^{2}\mathcal{U}_{[k-2,l-1]}\\\notag
+\delta_{0}(1+t)^{-1}\|\kappa-1\|_{2,[k,l],\Sigma_{t}^{\epsilon_{0}}}\\\notag
+\delta_{0}(1+t)^{-1}[1+\log(1+t)][\mathcal{Y}_{0}+(1+t)\mathcal{A}_{[l-1]}]+[1+\log(1+t)]\mathcal{W}_{\{l+1\}}\}\\
\|\zeta\|_{2,[k,l],\Sigma_{t}^{\epsilon_{0}}}\leq C_{l}(1+t)^{-1}\{\delta_{0}(1+t)^{-1}[1+\log(1+t)]\mathcal{U}_{k-1,l-1}\\\notag
+\delta_{0}(1+t)^{-1}\|\kappa-1\|_{2,[k,l],\Sigma_{t}^{\epsilon_{0}}}\\\notag
+\delta_{0}(1+t)^{-1}[1+\log(1+t)][\mathcal{Y}_{0}+(1+t)\mathcal{A}_{[l-1]}]+[1+\log(1+t)]\mathcal{W}_{\{l+1\}}\}
\end{align}
From the definition, we have:
\begin{align}
 \|(q_{T})_{b}\|_{2,[k,l],\Sigma_{t}^{\epsilon_{0}}}\leq C_{l}(1+t)^{-1}\|\kappa-1\|_{2,[k,l+1],\Sigma_{t}^{\epsilon_{0}}}
\end{align}

     To derive a recursive inequality for the quantities $\mathcal{U}_{[k,l]}$, we need to derive a similar estimate for $q_{T}$ in terms of 
$\mathcal{U}_{[k-1,l-1]}$. This requires an appropriate estimate for the quantity $\mathcal{T}_{[k-1,l-1]}$ defined in the statement of Lemma 11.8. Now, 
the assumptions of Lemma 11.8, which are those of Lemma 11.6, all follow under the assumptions of the present proposition, by Corollary 10.1.d, with $l_{*}$ 
in the role of $l$, assumption $\textbf{M}_{[m,l_{*}]}$, and Corollary 11.1.c with $l_{*}$ in the role of $l$. Lemma 11.8 gives a bound for the quantity 
$\mathcal{T}_{[k,l]}$ in terms of $\|\leftexp{(T)}{\slashed{\pi}}^{\prime}\|_{2,[k,l],\Sigma_{t}^{\epsilon_{0}}}$, 
$\|\lambda^{\prime}\|_{2,[k,l],\Sigma_{t}^{\epsilon_{0}}}$ and $\max_{i}\|\leftexp{(R_{i})}{\slashed{\pi}}\|_{2,[l-1],\Sigma_{t}^{\epsilon_{0}}}$. From 
(11.51) we obtain, using the estimate (11.176) and Corollary 11.1.c with $l_{*}$ in the role of $l$,
\begin{align}
 \|\leftexp{(T)}{\slashed{\pi}}^{\prime}\|_{2,[k,l],\Sigma_{t}^{\epsilon_{0}}}\leq C[1+\log(1+t)]\|\chi^{\prime}\|_{2,[k,l],\Sigma_{t}^{\epsilon_{0}}}\\\notag
+C_{l}(1+t)^{-1}[1+\log(1+t)]\{\delta_{0}(1+t)^{-2}[1+\log(1+t)]\mathcal{U}_{k-2,l-1}\\\notag
+\delta_{0}(1+t)^{-1}[\|\kappa-1\|_{2,[k,l],\Sigma_{t}^{\epsilon_{0}}}+\mathcal{Y}_{0}+(1+t)\mathcal{A}_{[l-1]}]+\mathcal{W}_{\{l+1\}}\}
\end{align}
Also, by Lemma 11.1 and Lemma 11.2.a,
\begin{align}
 \|\lambda^{\prime}\|_{2,[k,l],\Sigma_{t}^{\epsilon_{0}}}\leq C_{l}\{\|\kappa-1\|_{2,[k,l],\Sigma_{t}^{\epsilon_{0}}}+[1+\log(1+t)]\mathcal{W}_{\{l\}}\}
\end{align}
Also, by (10.145) we have:
\begin{align}
 \max_{i}\|\leftexp{(R_{i})}{\slashed{\pi}}\|_{2,[l-1],\Sigma_{t}^{\epsilon_{0}}}\leq C_{l}\{\mathcal{Y}_{0}+(1+t)\mathcal{A}_{[l-1]}+\mathcal{W}_{[l]}\}
\end{align}
Substituting in the conclusion of Lemma 11.8 the estimates (11.179)-(11.181) yields the following bound:
\begin{align}
 \mathcal{T}_{[k,l]}\leq C_{l}\{[1+\log(1+t)]\|\chi'\|_{2,[k,l],\Sigma_{t}^{\epsilon_{0}}}\\\notag
+\delta_{0}(1+t)^{-3}[1+\log(1+t)]^{2}\mathcal{U}_{k-2,l-1}\\\notag
+(1+t)^{-1}\|\kappa-1\|_{2,[k,l],\Sigma_{t}^{\epsilon_{0}}}+\delta_{0}(1+t)^{-1}[1+\log(1+t)][\mathcal{Y}_{0}+(1+t)\mathcal{A}_{[l-1]}]\\\notag
+(1+t)^{-1}[1+\log(1+t)]\mathcal{W}_{\{l+1\}}\}
\end{align}
So now we need an estimate for $\|\chi'\|_{2,[k,l],\Sigma_{t}^{\epsilon_{0}}}$, which can be proved similarly as Lemma 11.16:

$\textbf{Lemma 11.17}$ Let the assumptions of Proposition 11.2 hold. Then, provided that $\delta_{0}$ is suitably small (depending on $l$), the following estimate 
holds, for each $k\in\{0,...,m\}$:
\begin{align}
 \|\chi'\|_{2,[k,l],\Sigma_{t}^{\epsilon_{0}}}\leq C_{l}(1+t)^{-1}\{\delta_{0}(1+t)^{-2}[1+\log(1+t)]\mathcal{U}_{[k-2,l-1]}\\\notag
+(1+t)^{-1}\|\kappa-1\|_{2,[k-1,l+1],\Sigma_{t}^{\epsilon_{0}}}+[\mathcal{Y}_{0}+(1+t)\mathcal{A}_{[l]}]\\\notag
+\mathcal{W}_{\{l+1\}}+(1+t)^{-1}[1+\log(1+t)]\mathcal{W}^{Q}_{\{l\}}\}
\end{align}

     We now resume the proof of Proposition 11.2. Substituting the result of Lemma 11.17, yields the following estimate for $\mathcal{T}_{[k,l]}$:
\begin{align}
 \mathcal{T}_{[k,l]}\leq C_{l}(1+t)^{-1}\{\delta_{0}(1+t)^{-2}[1+\log(1+t)]^{2}\mathcal{U}_{k-2,l-1}\\\notag
+\|\kappa-1\|_{2,[k,l+1],\Sigma_{t}^{\epsilon_{0}}}+[1+\log(1+t)][\mathcal{Y}_{0}+(1+t)\mathcal{A}_{[l]}]\\\notag
+[1+\log(1+t)][\mathcal{W}_{\{l+1\}}+(1+t)^{-1}[1+\log(1+t)]\mathcal{W}^{Q}_{\{l\}}]\}
\end{align}
Consider now Lemma 11.12 with $(k+1,l+1)$ in the role of $(k,l)$. Since $(l+1)_{*}\leq l_{*}+1$, the assumptions hold if:
\begin{align}
 \|\lambda^{\prime}\|_{\infty,[k,l_{*}],\Sigma_{t}^{\epsilon_{0}}}\leq C_{l}\delta_{0}[1+\log(1+t)]\\\notag
\max_{i}\|\leftexp{(R_{i})}{\slashed{\pi}}\|_{\infty, [l_{*}],\Sigma_{t}^{\epsilon_{0}}}\leq C_{l}\delta_{0}(1+t)^{-1}[1+\log(1+t)]\\\notag
\|\leftexp{(T)}{\slashed{\pi}}^{\prime}\|_{\infty,[k,l_{*}],\Sigma_{t}^{\epsilon_{0}}}\leq C_{l}\delta_{0}(1+t)^{-1}[1+\log(1+t)]
\end{align}
These will follow from the assumptions of Proposition 11.2. So we have:
\begin{align}
 \|\xi\cdot\slashed{g}\|_{2,[k,l],\Sigma_{t}^{\epsilon_{0}}}, \|\xi\cdot\slashed{g}^{-1}\|_{2,[k,l],\Sigma_{t}^{\epsilon_{0}}}\\\notag
\leq C_{l}\{\|\xi\|_{2,[k,l],\Sigma_{t}^{\epsilon_{0}}}+\|\xi\|_{\infty,[k,l_{*}],\Sigma_{t}^{\epsilon_{0}}}
(\max_{i}\|\leftexp{(R_{i})}{\slashed{\pi}}\|_{2,[l-1],\Sigma_{t}^{\epsilon_{0}}}+\mathcal{T}_{[k-1,l-1]})\}
\end{align}
Substituting from (10.145) and the bound (11.184) this becomes:
\begin{align}
 \|\xi\cdot\slashed{g}\|_{2,[k,l],\Sigma_{t}^{\epsilon_{0}}}, \|\xi\cdot\slashed{g}^{-1}\|_{2,[k,l],\Sigma_{t}^{\epsilon_{0}}}\leq C_{l}\|\xi\|
_{2,[k,l],\Sigma_{t}^{\epsilon_{0}}}\\\notag
+C_{l}\|\xi\|_{\infty,[k,l_{*}],\Sigma_{t}^{\epsilon_{0}}}\cdot\\\notag
\cdot\{\delta_{0}(1+t)^{-3}[1+\log(1+t)]^{2}\mathcal{U}_{k-3,l-2}\\\notag
+(1+t)^{-1}\|\kappa-1\|_{2,[k-1,l],\Sigma_{t}^{\epsilon_{0}}}+[\mathcal{Y}_{0}+(1+t)\mathcal{A}_{[l-1]}]\\\notag
+\mathcal{W}_{\{l\}}+(1+t)^{-2}[1+\log(1+t)]^{2}\mathcal{W}^{Q}_{\{l-1\}}\}
\end{align}
Applying this to $\xi=(q_{T})_{b}$, we obtain, by (11.178) and (11.107),
\begin{align}
\|q_{T}\|_{2,[k,l],\Sigma_{t}^{\epsilon_{0}}}\leq C_{l}(1+t)^{-1}\{\delta_{0}^{2}(1+t)^{-1}\mathcal{U}_{k-1,l-1}+\|\kappa-1\|_{2,[k,l+1],\Sigma_{t}^{\epsilon_{0}}}\\\notag
+\delta_{0}[1+\log(1+t)][\mathcal{Y}_{0}+(1+t)\mathcal{A}_{[l-1]}]\\\notag
+\delta_{0}[1+\log(1+t)][\mathcal{W}_{\{l\}}+(1+t)^{-2}[1+\log(1+t)]^{2}\mathcal{W}^{Q}_{\{l-1\}}]\}
\end{align}

     Consider now the first two recursion relations of Lemma 11.3, replacing $m$ by $j\in\{0,...,k\}$. Taking the $\|\quad\|_{2,[k-j,l-j],\Sigma_{t}^{\epsilon_{0}}}$
norm we obtain:
\begin{align}
 \|p_{T,j}\|_{2,[k-j,l-j],\Sigma_{t}^{\epsilon_{0}}}\leq \|p_{T,j-1}\|_{2,[k+1-j,l+1-j],\Sigma_{t}^{\epsilon_{0}}}\\\notag
+C_{l}\{\|\slashed{d}\kappa\|_{\infty,[k-j,l_{*}],\Sigma_{t}^{\epsilon_{0}}}\|q_{T,j-1}\|_{2,[k-j,l-j],\Sigma_{t}^{\epsilon_{0}}}\\\notag
+\|\slashed{d}\kappa\|_{2,[k-j,l-j],\Sigma_{t}^{\epsilon_{0}}}\|q_{T,j-1}\|_{\infty,[k-j,l_{*}-j],\Sigma_{t}^{\epsilon_{0}}}
\end{align}
\begin{align}
 \|q_{T,j}\|_{2,[k-j,l-j],\Sigma_{t}^{\epsilon_{0}}}\leq \|q_{T,j-1}\|_{2,[k+1-j,l+1-j],\Sigma_{t}^{\epsilon_{0}}}\\\notag
+C_{l}\{\|q_{T}\|_{\infty,[k-j,l_{*}],\Sigma_{t}^{\epsilon_{0}}}\|p_{T,j-1}\|_{2,[k-j,l-j],\Sigma_{t}^{\epsilon_{0}}}\\\notag
+\|q_{T}\|_{2,[k-j,l-j],\Sigma_{t}^{\epsilon_{0}}}\|p_{T,j-1}\|_{\infty,[k-j,l_{*}-j],\Sigma_{t}^{\epsilon_{0}}}\}
\end{align}
Setting, for fixed $(k,l)$ and with $j\in\{0,...,k\}$,
\begin{align}
 x_{j}=\|p_{T,j}\|_{2,[k-j,l-j],\Sigma_{t}^{\epsilon_{0}}}\\\notag
y_{j}=\|q_{T,j}\|_{2,[k-j,l-j],\Sigma_{t}^{\epsilon_{0}}}
\end{align}
and:
\begin{align}
 a=x_{0}=0\\\notag
b=\|q_{T}\|_{2,[k,l],\Sigma_{t}^{\epsilon_{0}}}=y_{0}\\\notag
c=\|\slashed{d}\kappa\|_{2,[k,l],\Sigma_{t}^{\epsilon_{0}}}
\end{align}
(11.189)-(11.190) imply by assumptions $\textbf{M}_{[k,l_{*}+1]}$ and $\textbf{E}_{\{l_{*}+1\}}$ and Corollary 11.1.a with $l_{*}$ in the role of $l$,
\begin{align}
 x_{j}\leq x_{j-1}+\epsilon(y_{j-1}+c)\\\notag
y_{j}\leq y_{j-1}+\epsilon(x_{j-1}+\epsilon b)
\end{align}
where:
\begin{align}
 \epsilon=C_{l}\delta_{0}(1+t)^{-1}[1+\log(1+t)]
\end{align}
The inequalities (11.193) can be written in the form:
\begin{align}
 \begin{bmatrix}
  x_{j}\\
y_{j}
 \end{bmatrix}
\leq (\textbf{I}+\textbf{E})\begin{bmatrix}
                             x_{j-1}\\
y_{j-1}
                            \end{bmatrix}
+\textbf{a}
\end{align}
where
\begin{align}
 \textbf{E}=\begin{bmatrix}
             0&\epsilon\\
\epsilon&0
            \end{bmatrix}
\quad \textbf{a}=\begin{bmatrix}
                  \epsilon c\\
\epsilon^{2} b
                 \end{bmatrix}
\end{align}
It follows that:
\begin{align}
 \begin{bmatrix}
  x_{j}\\
y_{j}
 \end{bmatrix}
\leq (\textbf{I}+\textbf{E})^{j}\begin{bmatrix}
                                 x_{0}\\
y_{0}
                                \end{bmatrix}
+(\sum_{i=0}^{j-1}(\textbf{I}+\textbf{E})^{i})\textbf{a}
\end{align}
Now, in matrix norm:
\begin{align*}
 |\textbf{E}|\leq \epsilon
\end{align*}
and (11.197) implies, for $\epsilon\leq1$,
\begin{align}
 x_{j}\leq x_{0}+C_{j}\epsilon(x_{0}+y_{0})+C_{j}(\epsilon c+\epsilon^{2}b)\\\notag
y_{j}\leq y_{0}+C_{j}\epsilon(x_{0}+y_{0})+C_{j}(\epsilon c+\epsilon^{2}b)
\end{align}
In view of (11.191) and (11.192), substituting for $b$ the bounds (11.188), and noting that:
\begin{align*}
 c\leq C(1+t)^{-1}\|\kappa-1\|_{2,[k,l+1],\Sigma_{t}^{\epsilon_{0}}}
\end{align*}
yields the estimates:
\begin{align}
 \|p_{T,j}\|_{2,[k-j,l-j],\Sigma_{t}^{\epsilon_{0}}}\\\notag
\leq C_{l}\delta_{0}(1+t)^{-2}[1+\log(1+t)]\{\delta_{0}(1+t)^{-1}\mathcal{U}_{[k-1,l-1]}\\\notag
+\|\kappa-1\|_{2,[k,l+1],\Sigma_{t}^{\epsilon_{0}}}
+\delta_{0}[1+\log(1+t)][\mathcal{Y}_{0}+(1+t)\mathcal{A}_{[l-1]}]\\\notag
+\delta_{0}[1+\log(1+t)][\mathcal{W}_{\{l\}}+(1+t)^{-2}[1+\log(1+t)]^{2}\mathcal{W}^{Q}_{\{l-1\}}]\}
\end{align}
for all $j\in\{0,...,k\}$.

\begin{align}
 \|q_{T,j}\|_{2,[k-j,l-j],\Sigma_{t}^{\epsilon_{0}}}\\\notag
\leq C_{l}(1+t)^{-1}\{\delta_{0}(1+t)^{-1}\mathcal{U}_{[k-1,l-1]}\\\notag
+\|\kappa-1\|_{2,[k,l+1],\Sigma_{t}^{\epsilon_{0}}}
+\delta_{0}[1+\log(1+t)][\mathcal{Y}_{0}+(1+t)\mathcal{A}_{[l-1]}]\\\notag
+\delta_{0}[1+\log(1+t)][\mathcal{W}_{\{l\}}+(1+t)^{-2}[1+\log(1+t)]^{2}\mathcal{W}^{Q}_{\{l-1\}}]\}
\end{align}
for all $j\in\{0,...,k\}$

Next, we consider the formula (11.140). Taking the $\|\quad\|_{2,[l-k+n],\Sigma_{t}}$ norm we obtain,
using (11.200),
\begin{align}
 \|r^{n}_{T,k}\|_{2,[l-k+n],\Sigma_{t}^{\epsilon_{0}}}\\\notag
\leq \sum_{i=0}^{k-1-n}N_{i}(n)\sum_{j'=0}^{k-n-1-i}\|\kappa q_{T,j'}\|_{2,[k-1-n-j',l-1-n],\Sigma_{t}^{\epsilon_{0}}}\\\notag
\leq C_{l}\max_{j'\in\{0,...,k-1-n\}}\|\kappa q_{T,j'}\|_{2,[k-1-j',l-1-j'],\Sigma_{t}^{\epsilon_{0}}}\\\notag
\leq C_{l}(1+t)^{-1}[1+\log(1+t)]\{\delta_{0}(1+t)^{-1}\mathcal{U}_{[k-2,l-2]}\\\notag
+\|\kappa-1\|_{2,[k-1,l],\Sigma_{t}^{\epsilon_{0}}}
+\delta_{0}[1+\log(1+t)][\mathcal{Y}_{0}+(1+t)\mathcal{A}_{[l-2]}]\\\notag
+\delta_{0}[1+\log(1+t)][\mathcal{W}_{\{l-1\}}+(1+t)^{-2}[1+\log(1+t)]^{2}\mathcal{W}^{Q}_{\{l-2\}}]\}
\end{align}
for all $n\in\{0,...,k-1\}$.

Finally, we take the $\|\quad\|_{2,[l-k],\Sigma_{t}^{\epsilon_{0}}}$ norm of (11.142). We then obtain on the left:
\begin{align*}
 \|(T)^{k+1}\hat{T}^{i}\|_{2,[l-k],\Sigma_{t}^{\epsilon_{0}}}=\mathcal{U}_{k,l}-\mathcal{U}_{k-1,l}
\end{align*}
By (11.199) and (11.200) with $j=k$ and Proposition 10.2 together with Proposition 10.1 and Corollary 11.1.a with $l_{*}$ in the role of $l$:
\begin{align}
 \max_{i}\|p_{T,k}\hat{T}^{i}\|_{2,[l-k],\Sigma_{t}^{\epsilon_{0}}}\\\notag
\leq C_{l}\delta_{0}(1+t)^{-2}[1+\log(1+t)]\{\delta_{0}(1+t)^{-1}\mathcal{U}_{[k-1,l-1]}\\\notag
+\|\kappa-1\|_{2,[k,l+1],\Sigma_{t}^{\epsilon_{0}}}
+\delta_{0}[1+\log(1+t)][\mathcal{Y}_{0}+(1+t)\mathcal{A}_{[l-1]}]\\\notag
+\delta_{0}[1+\log(1+t)][\mathcal{W}_{\{l\}}+(1+t)^{-2}[1+\log(1+t)]^{2}\mathcal{W}^{Q}_{\{l-1\}}]\}
\end{align}
\begin{align}
 \max_{i}\|q_{T,k}\cdot\slashed{d}x^{i}\|_{2,[l-k],\Sigma_{t}^{\epsilon_{0}}}\\\notag
\leq C_{l}(1+t)^{-1}\{\delta_{0}(1+t)^{-1}\mathcal{U}_{[k-1,l-1]}\\\notag
+\|\kappa-1\|_{2,[k,l+1],\Sigma_{t}^{\epsilon_{0}}}
+\delta_{0}[1+\log(1+t)][\mathcal{Y}_{0}+(1+t)\mathcal{A}_{[l-1]}]\\\notag
+\delta_{0}[1+\log(1+t)][\mathcal{W}_{\{l\}}+(1+t)^{-2}[1+\log(1+t)]^{2}\mathcal{W}^{Q}_{\{l-1\}}]\}
\end{align}
By (11.201) with $n=0$ and Proposition 10.2 ($k\geq 1$) together with Proposition 10.1 with $l_{*}+1$ in the role of $l$ and Corollary 11.1.a with 
$l_{*}$ in the role of $l$:
\begin{align}
 \|r^{0}_{T,k}\cdot\slashed{d}\hat{T}^{i}\|_{2,[l-k],\Sigma_{t}^{\epsilon_{0}}}\\\notag
\leq C_{l}(1+t)^{2}[1+\log(1+t)]\cdot\{\delta_{0}^{2}(1+t)^{-1}\mathcal{U}_{k-2,l-2}+\|\kappa-1\|_{2,[k-1,l],\Sigma_{t}^{\epsilon_{0}}}\\\notag
+\delta_{0}[1+\log(1+t)][\mathcal{Y}_{0}+(1+t)\mathcal{A}_{[l-1]}]\\\notag
+\delta_{0}[1+\log(1+t)][\mathcal{W}_{\{l\}}+(1+t)^{-2}[1+\log(1+t)]^{2}\mathcal{W}^{Q}_{\{l-1\}}]\}
\end{align}
Also by (11.201) together with Proposition 11.1 and Corollary 11.1.a with $l_{*}$ in the role of $l$:
\begin{align}
 \max_{i}\|\sum_{n=1}^{k-1}r^{n}_{T,k}\cdot\slashed{d}(T)^{n}\hat{T}^{i}\|_{2,[l-k],\Sigma_{t}^{\epsilon_{0}}}\\\notag
\leq C_{l}(1+t)^{-1}\{\max_{n\in\{1,...,k-1\}}\max_{i}\|r^{n}_{T,k}\|_{\infty,[l_{*}-k+n],\Sigma_{t}^{\epsilon_{0}}}\mathcal{U}_{k-2,l-1}\\\notag
+(\sum_{n=1}^{k-1}\|r^{n}_{T,k}\|_{2,[l-k],\Sigma_{t}^{\epsilon_{0}}})\max_{n\in\{1,...,k-1\}}\|(T)^{n}\hat{T}^{i}\|
_{\infty,[l_{*}+1-n],\Sigma_{t}^{\epsilon_{0}}}\}\\\notag
\leq C_{l}\delta_{0}(1+t)^{-2}[1+\log(1+t)]^{2}\mathcal{U}_{k-2,l-1}\\\notag
+C_{l}\delta_{0}(1+t)^{-3}[1+\log(1+t)]^{2}\{\|\kappa-1\|_{2,[k-1,l-1],\Sigma_{t}^{\epsilon_{0}}}\\\notag
+\delta_{0}[1+\log(1+t)][\mathcal{Y}_{0}+(1+t)\mathcal{A}_{[l-2]}]\\\notag
+\delta_{0}[1+\log(1+t)][\mathcal{W}_{\{l-1\}}+(1+t)^{-2}[1+\log(1+t)]^{2}\mathcal{W}^{Q}_{\{l-2\}}]\}
\end{align}
Combining (11.204) and (11.205) we obtain:
\begin{align}
 \max_{i}\|\sum_{n=0}^{k-1}r^{n}_{T,k}\cdot\slashed{d}(T)^{n}\hat{T}^{i}\|_{2,[l-k],\Sigma_{t}^{\epsilon_{0}}}\\\notag
\leq C_{l}(1+t)^{-2}[1+\log(1+t)]\cdot\\\notag
\{\delta_{0}[1+\log(1+t)]\mathcal{U}_{k-2,l-1}+\|\kappa-1\|_{2,[k-1,l],\Sigma_{t}^{\epsilon_{0}}}\\\notag
+\delta_{0}[1+\log(1+t)][\mathcal{Y}_{0}+(1+t)\mathcal{A}_{[l-1]}]\\\notag
+\delta_{0}[1+\log(1+t)][\mathcal{W}_{\{l\}}+(1+t)^{-2}[1+\log(1+t)]^{2}\mathcal{W}^{Q}_{\{l-1\}}]\}
\end{align}
In view of (11.202), (11.203) and (11.206), we arrive at the following recursive inequality:
\begin{align}
 \mathcal{U}_{k,l}-\mathcal{U}_{k-1,l}\leq C_{l}\delta_{0}(1+t)^{-2}[1+\log(1+t)]^{2}\mathcal{U}_{k-1,l-1}+b_{k,l}
\end{align}
where:
\begin{align}
 b_{k,l}=C_{l}(1+t)^{-1}\{\|\kappa-1\|_{2,[k,l+1],\Sigma_{t}^{\epsilon_{0}}}\\\notag
+\delta_{0}[1+\log(1+t)][\mathcal{Y}_{0}+(1+t)\mathcal{A}_{[l-1]}\\\notag
+\mathcal{W}_{\{l\}}+(1+t)^{-2}[1+\log(1+t)]^{2}\mathcal{W}^{Q}_{\{l-1\}}]\}
\end{align}
(11.207) implies:
\begin{align}
 \mathcal{U}_{k,l}\leq (1+C_{l}\delta_{0}(1+t)^{-2}[1+\log(1+t)]^{2})\mathcal{U}_{k-1,l}+b_{k,l}
\end{align}
which yields:
\begin{align}
 \mathcal{U}_{k,l}\leq (1+C_{l}\delta_{0}(1+t)^{-2}[1+\log(1+t)]^{2})^{k}\mathcal{U}_{0,l}\\\notag
+\sum_{j=1}^{k}(1+C_{l}\delta_{0}(1+t)^{-2}[1+\log(1+t)]^{2})^{k-j}b_{j,l}
\end{align}
Since $b_{k,l}$ are non-decreasing in $k$, this in turn implies:
\begin{align}
 \mathcal{U}_{k,l}\leq C_{l}\{\mathcal{U}_{0,l}+b_{k,l}\}
\end{align}
Substituting finally in (11.211) the estimate (11.162) and the definition (11.208) we obtain:
\begin{align}
 \mathcal{U}_{k,l}\leq C_{l}(1+t)^{-1}\{\|\kappa-1\|_{2,[k,l+1],\Sigma_{t}^{\epsilon_{0}}}\\\notag
+\delta_{0}[1+\log(1+t)][\mathcal{Y}_{0}+(1+t)\mathcal{A}_{[l-1]}\\\notag
+\mathcal{W}_{\{l\}}+(1+t)^{-2}[1+\log(1+t)]^{2}\mathcal{W}^{Q}_{\{l-1\}}]\}
\end{align}
Since we have:
\begin{align}
 \|\kappa-1\|_{2,[k,l+1],\Sigma_{t}^{\epsilon_{0}}}\leq C_{l}\{\mathcal{B}_{[k,l+1]}+[1+\log(1+t)]\mathcal{W}_{\{l+1\}}\}
\end{align}
The proposition follows. $\qed$\vspace{7mm}

As before, we have many corollaries.

$\textbf{Corollary 11.2.a}$ Under the assumptions of Proposition 11.2 the coefficients of the expression for $(T)^{k+1}\hat{T}^{i}$, $k=0,...,m$, of Lemma 11.3 
satisfy:
\begin{align*}
\|p_{T,k}\|_{2,[m-k,l-k],\Sigma_{t}^{\epsilon_{0}}}\\\notag
\leq C_{l}\delta_{0}(1+t)^{-2}[1+\log(1+t)]\{\mathcal{B}_{[m,l+1]}\\\notag
+\delta_{0}[1+\log(1+t)](\mathcal{Y}_{0}+(1+t)\mathcal{A}_{[l-1]})\\\notag
+[1+\log(1+t)][\mathcal{W}_{\{l+1\}}+\delta_{0}(1+t)^{-2}[1+\log(1+t)]^{2}\mathcal{W}^{Q}_{\{l-1\}}]\}
\end{align*}
for all $k\in\{0,...,m\}$.

\begin{align*}
 \|q_{T,k}\|_{2,[m-k,l-k],\Sigma_{t}^{\epsilon_{0}}}\\\notag
\leq C_{l}(1+t)^{-1}\{\mathcal{B}_{[m,l+1]}\\\notag
+\delta_{0}[1+\log(1+t)](\mathcal{Y}_{0}+(1+t)\mathcal{A}_{[l-1]})\\\notag
+[1+\log(1+t)][\mathcal{W}_{\{l+1\}}+\delta_{0}(1+t)^{-2}[1+\log(1+t)]^{2}\mathcal{W}^{Q}_{\{l-1\}}]\}
\end{align*}
for all $k\in\{0,...,m\}$

and:
\begin{align*}
 \|r^{n}_{T,k}\|_{2,[l-k+n],\Sigma_{t}^{\epsilon_{0}}}\\\notag
\leq C_{l}(1+t)^{-1}[1+\log(1+t)]\cdot\\\notag
\{\mathcal{B}_{[m,l]}+\delta_{0}[1+\log(1+t)](\mathcal{Y}_{0}+(1+t)\mathcal{A}_{[l-2]})\\\notag
+[1+\log(1+t)][\mathcal{W}_{\{l\}}+\delta_{0}(1+t)^{-2}[1+\log(1+t)]^{2}\mathcal{W}^{Q}_{\{l-2\}}]\}
\end{align*}
for all $k\in\{0,...,m\}$, $n\in\{0,...,k-1\}$.

Also,
\begin{align*}
 \sum_{k=0}^{m}\max_{i}\|(T)^{k+1}L^{i}\|_{2,[l-k],\Sigma_{t}^{\epsilon_{0}}}\\\notag
\leq C_{l}(1+t)^{-1}\{\mathcal{B}_{[m,l+1]}+\delta_{0}[1+\log(1+t)][\mathcal{Y}_{0}+(1+t)\mathcal{A}_{[l-1]}]\\\notag
+[1+\log(1+t)][\mathcal{W}_{\{l+1\}}+\delta_{0}(1+t)^{-2}[1+\log(1+t)]^{2}\mathcal{W}^{Q}_{\{l-1\}}]\}
\end{align*}

$\textbf{Corollary 11.2.b}$ Under the assumptions of Proposition 11.2 the following estimates hold:
\begin{align*}
 \|\psi_{\hat{T}}\|_{2,[m+1,l+1],\Sigma_{t}^{\epsilon_{0}}}\leq C_{l}\{\mathcal{W}_{\{l+1\}}+\\\notag
+\delta_{0}(1+t)^{-1}[(1+t)^{-1}\mathcal{B}_{[m,l+1]}+\mathcal{Y}_{0}+(1+t)\mathcal{A}_{[l]}+\delta_{0}(1+t)^{-2}\mathcal{W}^{Q}_{\{l\}}]\}
\end{align*}
\begin{align*}
 \|\psi_{L}\|_{2,[m+1,l+1],\Sigma_{t}^{\epsilon_{0}}}\leq C_{l}\{\mathcal{W}_{\{l+1\}}+\\\notag
+\delta_{0}(1+t)^{-1}[(1+t)^{-1}\mathcal{B}_{[m,l+1]}+\mathcal{Y}_{0}+(1+t)\mathcal{A}_{[l]}+\delta_{0}(1+t)^{-2}\mathcal{W}^{Q}_{\{l\}}]\}
\end{align*}
\begin{align*}
 \|\slashed{\psi}\|_{2,[m+1,l+1],\Sigma_{t}^{\epsilon_{0}}}\leq C_{l}\{\mathcal{W}_{\{l+1\}}+\\\notag
+\delta_{0}(1+t)^{-1}[(1+t)^{-1}\mathcal{B}_{[m,l+1]}+\mathcal{Y}_{0}+(1+t)\mathcal{A}_{[l]}+\delta_{0}(1+t)^{-2}\mathcal{W}^{Q}_{\{l\}}]\}
\end{align*}
and:
\begin{align*}
 \|\omega_{T\hat{T}}\|_{2,[m,l],\Sigma_{T}^{\epsilon_{0}}}\leq C_{l}\{\mathcal{W}_{\{l+1\}}+\\\notag
\delta_{0}(1+t)^{-1}[(1+t)^{-1}\mathcal{B}_{[m-1,l]}+\mathcal{Y}_{0}+(1+t)\mathcal{A}_{[l-1]}+\delta_{0}(1+t)^{-2}\mathcal{W}^{Q}_{\{l-1\}}]\}
\end{align*}
\begin{align*}
 \|\omega_{L\hat{T}}\|_{2,[m,l],\Sigma_{t}^{\epsilon_{0}}}\leq C_{l}(1+t)^{-1}\{\mathcal{W}^{Q}_{\{l\}}+\\\notag
\delta_{0}(1+t)^{-1}[(1+t)^{-1}\mathcal{B}_{[m-1,l]}+\mathcal{Y}_{0}+(1+t)\mathcal{A}_{[l-1]}+\mathcal{W}_{\{l\}}]\}
\end{align*}
\begin{align*}
 \|\slashed{\omega}_{\hat{T}}\|_{2,[m,l],\Sigma_{t}^{\epsilon_{0}}}\leq C_{l}(1+t)^{-1}\{\mathcal{W}_{\{l+1\}}+\\\notag
+\delta_{0}(1+t)^{-1}[(1+t)^{-1}\mathcal{B}_{[m-1,l]}+\mathcal{Y}_{0}+(1+t)\mathcal{A}_{[l-1]}+\delta_{0}(1+t)^{-2}\mathcal{W}^{Q}_{\{l-1\}}]\}
\end{align*}
\begin{align*}
\|\slashed{\omega}_{L}\|_{2,[m,l],\Sigma_{t}^{\epsilon_{0}}}\leq C_{l}(1+t)^{-1}\{\mathcal{W}_{\{l+1\}}+\\\notag
 \delta_{0}(1+t)^{-1}[(1+t)^{-1}\mathcal{B}_{[m-1,l]}+\mathcal{Y}_{0}+(1+t)\mathcal{A}_{[l-1]}+\delta_{0}(1+t)^{-2}\mathcal{W}^{Q}_{\{l-1\}}]\}
\end{align*}
\begin{align*}
 \|\slashed{\omega}\|_{2,[m,l],\Sigma_{t}^{\epsilon_{0}}}\leq C_{l}(1+t)^{-1}\{\mathcal{W}_{\{l+1\}}+\\\notag
\delta_{0}(1+t)^{-1}[(1+t)^{-1}\mathcal{B}_{[m-1,l]}+\mathcal{Y}_{0}+(1+t)\mathcal{A}_{[l-1]}+\delta_{0}(1+t)^{-2}\mathcal{W}^{Q}_{\{l-1\}}]\}
\end{align*}
Also, we have:
\begin{align*}
 \|\slashed{k}\|_{2,[m,l],\Sigma_{t}^{\epsilon_{0}}}\leq C_{l}(1+t)^{-1}\{\mathcal{W}_{\{l+1\}}+\delta_{0}(1+t)^{-2}\mathcal{W}^{Q}_{\{l\}}\\\notag
+\delta_{0}(1+t)^{-1}[(1+t)^{-1}\mathcal{B}_{[m-1,l]}+\mathcal{Y}_{0}+(1+t)\mathcal{A}_{[l-1]}]\}
\end{align*}
\begin{align*}
 \|\kappa^{-1}\zeta\|_{2,[m,l],\Sigma_{t}^{\epsilon_{0}}}\leq C_{l}(1+t)^{-1}\{\mathcal{W}_{\{l+1\}}+\delta_{0}(1+t)^{-2}\mathcal{W}^{Q}_{\{l\}}\\\notag
+\delta_{0}(1+t)^{-1}[(1+t)^{-1}\mathcal{B}_{[m-1,l]}+\mathcal{Y}_{0}+(1+t)\mathcal{A}_{[l-1]}]\}
\end{align*}
\begin{align*}
 \|\kappa\slashed{k}\|_{2,[m,l],\Sigma_{t}^{\epsilon_{0}}}\leq C_{l}(1+t)^{-1}\{[1+\log(1+t)][\mathcal{W}_{\{l+1\}}+\delta_{0}(1+t)^{-2}\mathcal{W}^{Q}_{\{l\}}]\\\notag
+\delta_{0}(1+t)^{-1}[\mathcal{B}_{[m,l]}+[1+\log(1+t)](\mathcal{Y}_{0}+(1+t)\mathcal{A}_{[l-1]})]\}
\end{align*}
\begin{align*}
 \|\zeta\|_{2,[m,l],\Sigma_{t}^{\epsilon_{0}}}\leq C_{l}(1+t)^{-1}\{[1+\log(1+t)][\mathcal{W}_{\{l+1\}}+\delta_{0}(1+t)^{-2}\mathcal{W}^{Q}_{\{l\}}]\\\notag
+\delta_{0}(1+t)^{-1}[\mathcal{B}_{[m,l]}+[1+\log(1+t)](\mathcal{Y}_{0}+(1+t)\mathcal{A}_{[l-1]})]\}
\end{align*}

     The next Corollary is obtained from Lemma 11.17 by substituting Proposition 11.2.

$\textbf{Corollary 11.2.c}$ Under the assumptions of Proposition 11.2 we have:
\begin{align*}
 \|\chi'\|_{2,[m,l],\Sigma_{t}^{\epsilon_{0}}}\leq C_{l}(1+t)^{-1}\{(1+t)^{-1}\mathcal{B}_{[m-1,l+1]}+\mathcal{Y}_{0}+(1+t)\mathcal{A}_{[l]}\\\notag
+\mathcal{W}_{\{l+1\}}+(1+t)^{-1}[1+\log(1+t)]\mathcal{W}^{Q}_{\{l\}}\}
\end{align*}
Moreover by (11.184) we have:
\begin{align*}
 \mathcal{T}_{[m,l]}\leq C_{l}(1+t)^{-1}\{\mathcal{B}_{[m,l+1]}+[1+\log(1+t)][\mathcal{Y}_{0}+(1+t)\mathcal{A}_{[l]}\\\notag
+\mathcal{W}_{\{l+1\}}+(1+t)^{-1}[1+\log(1+t)]\mathcal{W}^{Q}_{\{l\}}]\}
\end{align*}
and obviously,
\begin{align*}
 \|\leftexp{(T)}{\slashed{\pi}}+2(1-u+t)^{-1}\slashed{g}\|_{2,[m,l],\Sigma_{t}^{\epsilon_{0}}}\leq C_{l}\mathcal{T}_{[m,l]}
\end{align*}
Also from the estimate for $\zeta$ and (11.187),
\begin{align*}
 \|\Lambda\|_{2,[m,l],\Sigma_{t}^{\epsilon_{0}}}\leq C_{l}(1+t)^{-1}\{\mathcal{B}_{[m,l+1]}+\delta_{0}[1+\log(1+t)](\mathcal{Y}_{0}+(1+t)\mathcal{A}_{[l-1]})\\\notag
+[1+\log(1+t)][\mathcal{W}_{\{l+1\}}+\delta_{0}(1+t)^{-2}[1+\log(1+t)]\mathcal{W}^{Q}_{\{l\}}]\}
\end{align*}
and by Lemma 11.13 and Lemma 10.11,
\begin{align*}
 \|\leftexp{(T)}{\slashed{\pi}}_{1}\|_{2,[m,l-1],\Sigma_{t}^{\epsilon_{0}}}\leq C_{l}(1+t)^{-2}\{\mathcal{B}_{[m,l+1]}+[1+\log(1+t)]
[\mathcal{Y}_{0}+(1+t)\mathcal{A}_{[l]}\\\notag
+\mathcal{W}_{\{l+1\}}+(1+t)^{-1}[1+\log(1+t)]\mathcal{W}^{Q}_{\{l\}}]\}
\end{align*}\vspace{7mm}

     We now revisit Lemma 11.12. Under the assumptions of Proposition 11.2 the following bounds hold:
\begin{align}
 \|\lambda^{\prime}\|_{\infty,[m,l_{*}],\Sigma_{t}^{\epsilon_{0}}}\leq C_{l}\delta_{0}[1+\log(1+t)]\\\notag
\max_{i}\|\leftexp{(R_{i})}{\slashed{\pi}}\|_{\infty,[l_{*}],\Sigma_{t}^{\epsilon_{0}}}\leq C_{l}\delta_{0}(1+t)^{-1}[1+\log(1+t)]\\\notag
\|\leftexp{(T)}{\slashed{\pi}}^{\prime}\|_{\infty,[m,l_{*}],\Sigma_{t}^{\epsilon_{0}}}\leq C_{l}\delta_{0}(1+t)^{-1}[1+\log(1+t)]
\end{align}
These bounds imply the assumptions of Lemma 11.12 with $(m,l)$ in the role of $(k-1,l-1)$. Then we obtain:
\begin{align}
 \|\xi\cdot\slashed{g}\|_{2,[m,l],\Sigma_{t}^{\epsilon_{0}}}, \|\xi\cdot\slashed{g}^{-1}\|_{2,[m,l],\Sigma_{t}^{\epsilon_{0}}}\\\notag
\leq C_{l}\{\|\xi\|_{2,[m,l],\Sigma_{t}^{\epsilon_{0}}}+\|\xi\|_{\infty,[m,l_{*}-1],\Sigma_{t}^{\epsilon_{0}}}(
\max_{i}\|\leftexp{(R_{i})}{\slashed{\pi}}\|_{2,[l-1],\Sigma_{t}^{\epsilon_{0}}}+\mathcal{T}_{[m-1,l-1]})\}
\end{align}
So by Corollary 10.1.d and Corollary 11.2.c, we obtain:

$\textbf{Corollary 11.2.d}$ Under the assumptions of Proposition 11.2, we have:
\begin{align*}
 \|\xi\cdot\slashed{g}\|_{2,[m,l],\Sigma_{t}^{\epsilon_{0}}}, \|\xi\cdot\slashed{g}^{-1}\|_{2,[m,l],\Sigma_{t}^{\epsilon_{0}}}\\\notag
\leq C_{l}\{\|\xi\|_{2,[m,l],\Sigma_{t}^{\epsilon_{0}}}\\\notag
+\|\xi\|_{\infty,[m,l_{*}-1],\Sigma_{t}^{\epsilon_{0}}}[\mathcal{W}_{\{l\}}+(1+t)^{-2}[1+\log(1+t)]\mathcal{W}^{Q}_{\{l-1\}}\\\notag
+(1+t)^{-1}\mathcal{B}_{[m-1,l]}+\mathcal{Y}_{0}+(1+t)\mathcal{A}_{[l-1]}]\}
\end{align*}

\section{Bounds for Quantities $Q^{\prime}_{m,l}$ and $P^{\prime}_{m,l}$}
The object of this section is to obtain appropriate bounds for the quantities $\leftexp{(i_{1}...i_{l})}{Q}^{\prime}_{m,l}$ and 
$\leftexp{(i_{1}...i_{l})}{P}^{\prime}_{m,l}$, which through $\leftexp{(i_{1}...i_{l})}{B}^{\prime}_{m,l}$, enter the final estimates for the functions
$\leftexp{(i_{1}...i_{l})}{x}^{\prime}_{m,l}$. 
\subsection{Bounds for $Q^{\prime}_{m,l}$}
The quantities $\leftexp{(i_{1}...i_{l})}{Q}^{\prime}_{m,l}(t)$ are defined by (9.266). We first estimate the 
last term:
\begin{align}
 \|\leftexp{(i_{1}...i_{l})}{\dot{g}^{\prime}}_{m,l}(t)\|_{L^{2}([0,\epsilon_{0}]\times S^{2})}
\end{align}
$\leftexp{(i_{1}...i_{l})}{\dot{g}}^{\prime}_{m,l}$ are defined in terms of $\leftexp{(i_{1}...i_{l})}{g}^{\prime}_{m,l}$ by (9.252) and (9.253). 
$\leftexp{(i_{1}...i_{l})}{g}^{\prime}_{m,l}$ are in turn defined in Proposition 9.1 and Proposition 9.2. The first term in the expression for 
$\leftexp{(i_{1}...i_{l})}{g}^{\prime}_{m,l}$ is:
\begin{align*}
 R_{i_{l}}...R_{i_{1}}(T)^{m}\check{g}^{\prime}
\end{align*}
where $\check{g}^{\prime}$ is defined by (9.63). Note that for $m=l=0$ we have:
\begin{align}
 g^{\prime}_{0,0}=\check{g}^{\prime}
\end{align}
while according to (9.252),
\begin{align}
 \dot{g}^{\prime}_{0,0}=g^{\prime}_{0,0}-\xi\cdot x_{0}
\end{align}
where $\xi$ is defined by (9.70). Defining:
\begin{align}
 \check{\dot{g}}^{\prime}:=\dot{g}_{0,0}^{\prime}
\end{align}
Our objective in the following is to estimate:
\begin{align*}
 \|\check{\dot{g}}^{\prime}\|_{2,[m,l],\Sigma_{t}^{\epsilon_{0}}}
\end{align*}
Recalling from Chapter 8 that:
\begin{align*}
 x_{0}=\mu\slashed{d}\textrm{tr}\chi+\slashed{d}\check{f}
\end{align*}
we can write:
\begin{align}
 \check{\dot{g}}^{\prime}=\check{g}^{\prime}-\xi\cdot(\mu\slashed{d}\textrm{tr}\chi+\slashed{d}\check{f})
\end{align}
Substituting the expression (9.63) for $\check{g}^{\prime}$ in (11.220):
\begin{align}
 \check{\dot{g}}^{\prime}=-\xi\cdot\slashed{d}\check{f}-2\mu(\slashed{d}\mu)\cdot(i-\slashed{d}e)\\\notag
+\frac{1}{2}\frac{d\log\Omega}{dh}(Lh)f^{\prime}_{0}-\frac{1}{2}(\mu\textrm{tr}\chi+2m+2\mu e)\mu f^{\prime}_{1}\\\notag
+\dot{n}^{\prime\prime}_{0}+\mu(\dot{n}^{\prime}_{1}+\dot{n}^{\prime}_{2})
\end{align}
 where
\begin{align}
 \dot{n}^{\prime\prime}_{0}=-\frac{1}{2}L(\frac{dH}{dh})(\underline{L}Th)+\frac{1}{2}\frac{dH}{dh}(\ddot{\tau}+\underline{\nu}(LTh)+2\zeta\cdot\slashed{d}Th)
\end{align}
$\ddot{\tau}$ is defined by a formula similar to (9.38) but with
\begin{align}
 \dot{v}^{\prime}=v^{\prime}+\Omega\kappa^{2}(\slashed{d}h)\cdot\slashed{d}\textrm{tr}\chi
\end{align}
in the role of $v^{\prime}$, that is:
\begin{align}
 \ddot{\tau}=\Omega^{-1}(\dot{v}^{\prime}+T\tilde{\tau}+\leftexp{(T)}{\delta}\tilde{\tau})
\end{align}
By (9.35) we know that there is no acoustical principal term in $\dot{n}^{\prime\prime}_{0}$. Also
$\dot{n}_{1}^{\prime}$ is just ${n}_{1}^{\prime}$, while $\dot{n}^{\prime}_{2}$ is defined by:
\begin{align}
\dot{n}_{2}^{\prime}=\mu\dot{n}_{2}-L(\frac{\mu}{2\eta^{2}}(\frac{\rho}{\rho^{\prime}})^{\prime})\slashed{\Delta}h-L(\mu\eta^{-1}\hat{T}^{i})\slashed{\Delta}\psi_{i}
\end{align}
with $\dot{n}_{2}$ defined by:
\begin{align}
 \dot{n}_{2}=n_{2}-(\frac{1}{2\eta^{2}}(\frac{\rho}{\rho^{\prime}})^{\prime}\slashed{d}h+\eta^{-1}\slashed{d}\psi_{\hat{T}}
-\eta^{-2}(L\psi_{i})\slashed{d}x^{i})\cdot\slashed{d}\textrm{tr}\chi
\end{align}
By (9.68), $\dot{n}_{2}$ does not contain principal acoustical terms.

     Now $v^{\prime}$ is defined by (9.34):
\begin{align}
 v^{\prime}=v-\Omega(Th)\slashed{\Delta}\mu
\end{align}
We thus have:
\begin{align}
 \dot{v}^{\prime}=v-\Omega(Th)\slashed{\Delta}\mu+\Omega\kappa^{2}(\slashed{d}h)\cdot\slashed{d}\textrm{tr}\chi
\end{align}
$v$ is given by (9.18)-(9.24):
\begin{align}
 v=v_{1}+v_{2}+v_{3}
\end{align}
Defining:
\begin{align}
 \dot{v}_{1}=\Omega^{-1}v_{1}\\\notag
\dot{v}_{2}=\Omega^{-1}v_{2}-(Th)\slashed{\Delta}\mu+\kappa^{2}(\slashed{d}h)\cdot\slashed{d}\textrm{tr}\chi\\\notag
\dot{v}_{3}=\Omega^{-1}v_{3}
\end{align}
and:
\begin{align}
 \dot{v}=\dot{v}_{1}+\dot{v}_{2}+\dot{v}_{3}
\end{align}
Then according to (9.33), $\dot{v}_{1}$, $\dot{v}_{2}$ and $\dot{v}_{3}$ do not contain principal acoustical terms, and we have:
\begin{align}
 \dot{v}^{\prime}=\Omega\dot{v}
\end{align}
Our first objective is to estimate:
\begin{align}
 \|\dot{v}\|_{2,[m,l],\Sigma_{t}^{\epsilon_{0}}}
\end{align}
Here we just state the results:

$\textbf{Lemma 11.18}$ Under the assumptions of Proposition 11.2 augmented by the assumptions $\textbf{E}_{\{l_{*}+2\}}$, $\textbf{E}^{Q}_{\{l_{*}+1\}}$,
$\textbf{E}^{QQ}_{\{l_{*}\}}$ and $\textbf{M}_{[m+1,l_{*}+1]}$ we have the $L^{2}$ estimate:
\begin{align*}
 \|\dot{v}\|_{2,[m,l],\Sigma_{t}^{\epsilon_{0}}}\\
\leq C_{l}(1+t)^{-2}\{[1+\log(1+t)](\mathcal{W}_{\{l+2\}}+\mathcal{W}^{Q}_{\{l+1\}}+\mathcal{W}^{QQ}_{\{l\}})\\
+\delta_{0}(1+t)^{-1}[\mathcal{B}_{[m+1,l+1]}+[1+\log(1+t)](\mathcal{Y}_{0}+(1+t)\mathcal{A}_{[l]})]\}
\end{align*}
as well as $L^{\infty}$ estimate:
\begin{align*}
 \|\dot{v}\|_{\infty,[m,l_{*}],\Sigma_{t}^{\epsilon_{0}}}\leq C_{l}\delta_{0}(1+t)^{-3}[1+\log(1+t)]
\end{align*}
where we denote by $\textbf{E}^{QQ}_{m,n}$ the bootstrap assumption that there is a constant $C$ independent of $s$ such that for all $t\in[0,s]$:
\begin{align*}
 \textbf{E}^{QQ}_{m,n}\quad:\quad \max_{\alpha;i_{1}...i_{n}}\|R_{i_{n}}...R_{i_{1}}(T)^{m}(Q)^{2}\psi_{\alpha}\|
_{L^{\infty}(\Sigma_{t}^{\epsilon_{0}}}\leq C\delta_{0}(1+t)^{-1}
\end{align*}
We then denote by $\textbf{E}^{QQ}_{\{l\}}$ the conjunction of $\textbf{E}^{QQ}_{m,n}$ corresponding to (11.8). The constant $C$ then depends on $l$ only.
Here we also have introduced the quantities $\mathcal{W}^{QQ}_{m,n}$:
\begin{align}
 \mathcal{W}^{QQ}_{m,n}=\max_{\alpha;i_{1}...i_{n}}\|R_{i_{n}}...R_{i_{1}}(T)^{m}(Q)^{2}\psi_{\alpha}\|_{L^{2}(\Sigma_{t}^{\epsilon_{0}})}
\end{align}
We then denote by $\mathcal{W}^{QQ}_{\{l\}}$ the sum of $\mathcal{W}^{QQ}_{m,n}$ corresponding to (11.8). Moreover, we have:
\begin{align}
 \max_{\alpha}\|(L)^{2}\psi_{\alpha}\|_{2,\{l\},\Sigma_{t}^{\epsilon_{0}}}\leq (1+t)^{-2}(\mathcal{W}^{QQ}_{\{l\}}+\mathcal{W}^{Q}_{\{l\}})
\end{align}

     We now consider $\ddot{\tau}$, defined by (11.224). Setting:
\begin{align}
 \dot{\tau}=\Omega^{-1}\tilde{\tau}
\end{align}
(11.224) becomes, in view of (11.232):
\begin{align}
 \ddot{\tau}=\dot{v}+T\dot{\tau}+(T\log\Omega+\leftexp{(T)}{\delta})\dot{\tau}
\end{align}
By (9.37) we have:
\begin{align}
 \dot{\tau}=-\sum_{i}\{-(L\psi_{i})(\underline{L}\psi_{i})+\mu\slashed{d}\psi_{i}\cdot\slashed{d}\psi_{i}\}
\end{align}
Here we must use Lemma 11.12 with $(k,l)$ replaced by $(m+2,l+2)$ to estimate the second term on the right of (11.237):
\begin{align}
 \|\dot{\tau}\|_{2,[m+1,l+1],\Sigma_{t}^{\epsilon_{0}}}\\\notag
\leq C_{l}\delta_{0}(1+t)^{-2}\{\mathcal{W}_{\{l+2\}}+\mathcal{W}^{Q}_{\{l+1\}}\\\notag
+\delta_{0}(1+t)^{-2}[\mathcal{B}_{[m+1,l+1]}+[1+\log(1+t)](\mathcal{Y}_{0}+(1+t)\mathcal{A}_{[l]})]\}
\end{align}
Also, we have:
\begin{align}
 \|\dot{\tau}\|_{\infty,[m+1,l_{*}+1],\Sigma_{t}^{\epsilon_{0}}}\leq C_{l}\delta_{0}^{2}(1+t)^{-3}
\end{align}
Then combining with Lemma 11.18, we obtain:
\begin{align}
 \|\ddot{\tau}\|_{2,[m,l],\Sigma_{t}^{\epsilon_{0}}}\\\notag
\leq C_{l}(1+t)^{-2}\{[1+\log(1+t)](\mathcal{W}_{\{l+2\}}+\mathcal{W}^{Q}_{\{l+1\}}+\mathcal{W}^{QQ}_{\{l\}})\\\notag
+\delta_{0}(1+t)^{-1}[\mathcal{B}_{[m+1,l+1]}+[1+\log(1+t)](\mathcal{Y}_{0}+(1+t)\mathcal{A}_{[l]})]\}
\end{align}
as well as 
\begin{align}
 \|\ddot{\tau}\|_{2,[m,l_{*}],\Sigma_{t}^{\epsilon_{0}}}\leq C_{l}\delta_{0}(1+t)^{-3}[1+\log(1+t)]
\end{align}
By the above results and a direct calculation, we get an estimate:
\begin{align}
 \|\dot{n}^{\prime\prime}_{0}\|_{2,[m,l],\Sigma_{t}^{\epsilon_{0}}}\\\notag
\leq C_{l}(1+t)^{-2}\{[1+\log(1+t)](\mathcal{W}_{\{l+2\}}+\mathcal{W}^{Q}_{\{l+1\}}+\mathcal{W}^{QQ}_{\{l\}})\\\notag
+\delta_{0}(1+t)^{-1}[\mathcal{B}_{[m+1,l+1]}+[1+\log(1+t)](\mathcal{Y}_{0}+(1+t)\mathcal{A}_{[l]})]\}
\end{align}

     We proceed to estimate $\dot{n}_{1}^{\prime}$ and $\dot{n}^{\prime}_{2}$, given by (9.52) and (11.225). The estimate for $\dot{n}^{\prime}_{1}$ is 
straightforward. We use Lemma 11.12 to obtain:
\begin{align}
 \|\dot{n}_{1}^{\prime}\|_{2,[m,l],\Sigma_{t}^{\epsilon_{0}}}\leq C_{l}\delta_{0}(1+t)^{-3}\{\mathcal{W}_{\{l+2\}}+\delta_{0}(1+t)^{-1}\cdot\\\notag
[(1+t)^{-1}\mathcal{B}_{[m-1,l+1]}+\mathcal{Y}_{0}+(1+t)\mathcal{A}_{[l]}+(1+t)^{-2}[1+\log(1+t)]\mathcal{W}^{Q}_{\{l+1\}}]\}
\end{align}
and also, obviously,
\begin{align}
 \|\dot{n}^{\prime}_{1}\|_{2,[m,l_{*}],\Sigma_{t}^{\epsilon_{0}}}\leq C_{l}\delta_{0}^{2}(1+t)^{-4}
\end{align}

     To estimate $\dot{n}^{\prime}_{2}$ is a bit more complicated, because we have the term $\slashed{\Delta}\hat{T}^{i}$ in the expression (9.55)
for $n_{2}$. We must first derive an estimate for this. Recalling (9.66) and (9.67), we consider:
\begin{align}
 U^{i}=\slashed{\Delta}\hat{T}^{i}+\eta^{-1}\slashed{d}\textrm{tr}
 \chi\cdot\slashed{d}x^{i}+(1-u+t)^{-1}\eta^{-1}\slashed{\Delta}x^{i}
\end{align}
The last term in the above comes from the term $\slashed{q}\cdot\slashed{D}^{2}x^{i}$. We can write:
\begin{align}
 U^{i}=b_{1}\cdot\slashed{d}x^{i}+b_{2}\cdot\slashed{D}^{2}x^{i}
\end{align}
where by (9.66),
\begin{align}
 b_{1}=\slashed{\textrm{div}}\slashed{k}-\eta^{-1}i
 +\eta^{-2}\chi\cdot(\slashed{d}\eta)
\end{align}
and
\begin{align}
 b_{2}=\slashed{k}-\eta^{-1}\chi^{\prime}
\end{align}
To estimate
\begin{align*}
 \|U^{i}\|_{2,[m,l],\Sigma_{t}^{\epsilon_{0}}}
\end{align*}
we need the following two lemmas:

$\textbf{Lemma 11.19}$ Under the assumptions of Proposition 11.2, for an arbitrary $S_{t,u}$ tensorfield $\xi$ we have:
\begin{align*}
 \|\xi\cdot\slashed{d}x^{i}\|_{2,[m,l],\Sigma_{t}^{\epsilon_{0}}}\leq C_{l}\{\|\xi\|_{2,[m,l],\Sigma_{t}^{\epsilon_{0}}}\\\notag
+\|\xi\|_{\infty,[m,l_{*}-1],\Sigma_{t}^{\epsilon_{0}}}[(1+t)^{-1}\mathcal{B}_{[m-1,l]}+\mathcal{Y}_{0}+(1+t)\mathcal{A}_{[l-1]}\\\notag
+\mathcal{W}_{\{l\}}]\}
\end{align*}
$\textbf{Lemma 11.20}$ Under the assumptions of Proposition 11.2, for an arbitrary $S_{t,u}$ tensorfield $\xi$ we have:
\begin{align*}
 \|\xi\cdot\slashed{D}^{2}x^{i}\|_{2,[m,l],\Sigma_{t}^{\epsilon_{0}}}\leq C_{l}(1+t)^{-1}\{\|\xi\|_{2,[m,l],\Sigma_{t}^{\epsilon_{0}}}\\
+\|\xi\|_{\infty,[k,l_{*}-1],\Sigma_{t}^{\epsilon_{0}}}\{(1+t)^{-1}\mathcal{B}_{[m-1,l+1]}+\mathcal{Y}_{0}+(1+t)\mathcal{A}_{[l]}\\
+\mathcal{W}_{\{l+1\}}+(1+t)^{-2}[1+\log(1+t)]^{2}\mathcal{W}^{Q}_{\{l\}}]\}
\end{align*}
By Proposition 10.1/11.1, 10.2/11.2 as well as their corollaries, the proof of these two lemmas is almost a direct calculation.

Using these two lemmas and (11.247)-(11.249), we deduce:

$\textbf{Lemma 11.21}$ Under the assumptions of Proposition 11.2 augmented by the assumptions $\textbf{E}_{\{l_{*}+2\}}$, $\textbf{E}^{Q}_{\{l_{*}+1\}}$,
we have:
\begin{align*}
 \|U^{i}\|_{2,[m,l],\Sigma_{t}^{\epsilon_{0}}}\leq C_{l}(1+t)^{-2}\cdot\\
\{\mathcal{W}_{\{l+2\}}+(1+t)^{-1}[1+\log(1+t)]\mathcal{W}^{Q}_{\{l+1\}}\\
+(1+t)^{-1}\mathcal{B}_{[m-1,l+1]}+\mathcal{Y}_{0}+(1+t)\mathcal{A}_{[l]}\}
\end{align*}
and:
\begin{align*}
 \|U^{i}\|_{\infty,[m,l_{*}],\Sigma_{t}^{\epsilon_{0}}}\leq C_{l}\delta_{0}(1+t)^{-3}[1+\log(1+t)]
\end{align*}\vspace{7mm}

     Next we consider the function $n_{2}^{\prime}$ given by (11.226). In this equation the function $n_{2}$ is given by (9.55)-(9.57). Substituting
in (11.226) we obtain the following expression for $\dot{n}_{2}$:
\begin{align}
 \dot{n}_{2}=\frac{1}{2\eta^{2}}(\frac{\rho}{\rho^{\prime}})^{\prime}[\textrm{tr}\chi\slashed{\Delta}h+2\hat{\chi}\cdot\hat{\slashed{D}^{2}}h
+2i\cdot\slashed{d}h]\\\notag
+\eta^{-1}[\textrm{tr}\chi\hat{T}^{i}\slashed{\Delta}\psi_{i}+2\hat{\chi}\cdot
\hat{T}^{i}\hat{\slashed{D}^{2}}\psi_{i}+2i
\cdot\hat{T}^{i}\slashed{d}\psi_{i}]\\\notag
+\slashed{d}(\frac{1}{\eta^{2}}(\frac{\rho}{\rho^{\prime}})^{\prime})\cdot\slashed{d}Lh+2\slashed{d}(\eta^{-1}\hat{T}^{i})\cdot\slashed{d}(L\psi_{i})\\\notag
+\slashed{\Delta}(\frac{1}{2\eta^{2}}(\frac{\rho}{\rho^{\prime}})^{\prime})Lh+\slashed{\Delta}(\eta^{-1}\hat{T}^{i})(L\psi_{i})\\\notag
+\eta^{-2}(L\psi_{i})\slashed{d}x^{i}\cdot\slashed{d}\textrm{tr}\chi
\end{align}
By Lemma 11.21 and (11.250), we obtain the following estimates for $\dot{n}_{2}$:
\begin{align}
 \|\dot{n}_{2}\|_{2,[m,l],\Sigma_{t}^{\epsilon_{0}}}\leq C_{l}(1+t)^{-3}\{\mathcal{W}_{\{l+2\}}+\mathcal{W}^{Q}_{\{l+1\}}\\\notag
+\delta_{0}(1+t)^{-1}[(1+t)^{-1}\mathcal{B}_{[m-1,l+1]}+\mathcal{Y}_{0}+(1+t)\mathcal{A}_{[l]}]\}
\end{align}
and:
\begin{align}
 \|\dot{n}_{2}\|_{\infty,[m,l_{*}],\Sigma_{t}^{\epsilon_{0}}}\leq C_{l}\delta_{0}(1+t)^{-4}
\end{align}

     We turn to the function $\dot{n}_{2}^{\prime}$ given in terms of $\dot{n}_{2}$ by (11.225). We just use Corollary 11.1.b and Corollary 11.2.b as well
as a direct calculation to obtain:
\begin{align}
 \|\dot{n}_{2}^{\prime}\|_{2,[m,l],\Sigma_{t}^{\epsilon_{0}}}\leq C_{l}(1+t)^{-3}\{[1+\log(1+t)](\mathcal{W}_{\{l+2\}}+\mathcal{W}^{Q}_{\{l+1\}})\\\notag
+\delta_{0}(1+t)^{-1}[\mathcal{B}_{[m,l]}+[1+\log(1+t)](\mathcal{Y}_{0}+(1+t)\mathcal{A}_{[l]})]\}
\end{align}
and:
\begin{align}
 \|\dot{n}_{2}^{\prime}\|_{\infty,[m,l_{*}],\Sigma_{t}^{\epsilon_{0}}}\leq C_{l}(1+t)^{-4}[1+\log(1+t)]
\end{align}

     We now turn to the remaining terms in the expression (11.221) for $\check{\dot{g}}^{\prime}$. We first consider the functions $f^{\prime}_{0}$ and
$f^{\prime}_{1}$, given by (9.48) and (9.59) respectively:
\begin{align}
 f^{\prime}_{0}=\frac{1}{2}\frac{dH}{dh}(\underline{L}Th)\\
f^{\prime}_{1}=\frac{1}{2\eta^{2}}(\frac{\rho}{\rho^{\prime}})^{\prime}(\slashed{\Delta}h)+\eta^{-1}\hat{T}^{i}(\slashed{\Delta}\psi_{i})
\end{align}
Obviously, by the expression:
\begin{align*}
 \underline{L}Th=2(T)^{2}h+\eta^{-1}\kappa TLh+\eta^{-1}\kappa\Lambda\cdot\slashed{d}h
\end{align*}
we have:
\begin{align}
 \|f^{\prime}_{0}\|_{\infty,[m,l_{*}],\Sigma_{t}}\leq C_{l}\delta_{0}(1+t)^{-1}
\end{align}
and:
\begin{align}
 \|f^{\prime}_{0}\|_{2,[m,l],\Sigma_{t}^{\epsilon_{0}}}\leq C_{l}\{\mathcal{W}_{\{l+2\}}+(1+t)^{-1}[1+\log(1+t)]\mathcal{W}^{Q}_{\{l+1\}}\\\notag
+\delta_{0}(1+t)^{-2}[\mathcal{B}_{[m,l+1]}+\delta_{0}(\mathcal{Y}_{0}+(1+t)\mathcal{A}_{[l-1]})]\}
\end{align}
Using Proposition 11.2 we deduce:
\begin{align}
 \|f^{\prime}_{1}\|_{2,[m,l],\Sigma_{t}^{\epsilon_{0}}}\\\notag
\leq C_{l}(1+t)^{-2}\{\mathcal{W}_{\{l+2\}}+\delta_{0}(1+t)^{-3}[1+\log(1+t)]^{2}\mathcal{W}^{Q}_{\{l\}}\\\notag
+\delta_{0}(1+t)^{-1}[(1+t)^{-1}\mathcal{B}_{[m-1,l+1]}+\mathcal{Y}_{0}+(1+t)\mathcal{A}_{[l]}]\}
\end{align}
and also:
\begin{align}
 \|f^{\prime}_{1}\|_{\infty,[m,l_{*}],\Sigma_{t}^{\epsilon_{0}}}\leq C_{l}\delta_{0}(1+t)^{-3}
\end{align}

     With the help of (11.257)-(11.260) we can estimate, in $\|\quad\|_{2,[m,l],\Sigma_{t}^{\epsilon_{0}}}$ norm, the third and fourth terms on the right in (11.221) in 
a straightforward manner. Concerning the second term, we just use the estimates for $e$ and $i$ which we can obtain by using Proposition 11.1 and 11.2:
\begin{align}
 \|e\|_{2,[m+1,l+1],\Sigma_{t}^{\epsilon_{0}}}\leq C_{l}(1+t)^{-1}\{\mathcal{W}^{Q}_{\{l+1\}}\\\notag
+\delta_{0}(1+t)^{-1}[(1+t)^{-1}\mathcal{B}_{[m,l+1]}+\mathcal{Y}_{0}+(1+t)\mathcal{A}_{[l]}+\mathcal{W}_{\{l+1\}}]\}\\
\|e\|_{\infty,[m,l_{*}],\Sigma_{t}^{\epsilon_{0}}}\leq C_{l}\delta_{0}(1+t)^{-2}
\end{align}
and:
\begin{align}
 \|i\|_{2,[m,l],\Sigma_{t}^{\epsilon_{0}}}\\\notag
\leq C_{l}(1+t)^{-2}\{\mathcal{W}_{\{l+2\}}+\delta_{0}(1+t)^{-2}[1+\log(1+t)]\mathcal{W}^{Q}_{\{l+1\}}\\\notag
+\delta_{0}(1+t)^{-1}[(1+t)^{-1}\mathcal{B}_{[m-1,l+1]}+\mathcal{Y}_{0}+(1+t)\mathcal{A}_{[l]}]\}\\
\|i\|_{\infty,[m,l_{*}],\Sigma_{t}^{\epsilon_{0}}}\leq C_{l}\delta_{0}^{2}(1+t)^{-4}
\end{align}
Finally, we have:
\begin{align*}
 -\xi\cdot\slashed{d}\check{f}
\end{align*}
the first term on the right in (11.221). Here $\xi$ is given by (9.70):
\begin{align}
 \xi=-\slashed{d}\mu +(\frac{\mu}{2\eta^{2}}(\frac{\rho}{\rho^{\prime}})^{\prime}-\frac{1}{2}\frac{dH}{dh}\eta^{-1}\kappa)\slashed{d}h+T^{i}\slashed{d}\psi_{i}
-\eta^{-1}\kappa(L\psi_{i})\slashed{d}x^{i}
\end{align}
and $\check{f}$ is the function defined by (8.27):
\begin{align}
 \check{f}=-\frac{1}{2}\frac{dH}{dh}\tau_{\underline{L}}
\end{align}
So we obtain:
\begin{align}
 \|\slashed{d}\check{f}\|_{2,[m,l],\Sigma_{t}^{\epsilon_{0}}}\\\notag
\leq C_{l}(1+t)^{-1}\{\mathcal{W}_{\{l+2\}}+(1+t)^{-1}[1+\log(1+t)]\mathcal{W}^{Q}_{\{l+1\}}+\delta_{0}(1+t)^{-2}\mathcal{B}_{[m,l+1]}\}
\end{align}
and also:
\begin{align}
 \|\slashed{d}\check{f}\|_{\infty,[m,l_{*}],\Sigma_{t}^{\epsilon_{0}}}\leq C_{l}\delta_{0}(1+t)^{-2}
\end{align}
Next, using Corollary 11.1.b with $l_{*}$ in the role of $l$ and Corollary 11.2.b we deduce:
\begin{align}
 \|\xi\|_{\infty,[m,l_{*}],\Sigma_{t}^{\epsilon_{0}}}\leq C_{l}(1+t)^{-1}[1+\log(1+t)]
\end{align}
and:
\begin{align}
 \|\xi\|_{2,[m,l],\Sigma_{t}^{\epsilon_{0}}}\leq C_{l}(1+t)^{-1}\{[1+\log(1+t)](\mathcal{W}_{\{l+1\}}+\delta_{0}^{2}(1+t)^{-2}\mathcal{W}^{Q}_{\{l\}})\\\notag
+\mathcal{B}_{[m,l+1]}+\delta_{0}(1+t)^{-1}(\mathcal{Y}_{0}+(1+t)\mathcal{A}_{[l-1]})\}
\end{align}
To get the final estimate, we need Corollary 11.2.d to see that:
\begin{align}
 \|\xi\cdot\slashed{g}^{-1}\|_{2,[m,l],\Sigma_{t}^{\epsilon_{0}}}\leq C_{l}(1+t)^{-1}\cdot\\\notag
\{[1+\log(1+t)](\mathcal{W}_{\{l+1\}}+\delta_{0}(1+t)^{-2}[1+\log(1+t)]^{2}\mathcal{W}^{Q}_{\{l\}})\\\notag
+\mathcal{B}_{[m,l+1]}+\delta_{0}[1+\log(1+t)](\mathcal{Y}_{0}+(1+t)\mathcal{A}_{[l-1]})\}
\end{align}
and also:
\begin{align}
 \|\xi\cdot\slashed{g}^{-1}\|_{\infty,[m,l_{*}],\Sigma_{t}^{\epsilon_{0}}}\leq C_{l}\delta_{0}(1+t)^{-1}[1+\log(1+t)]
\end{align}
Then the first term on the right of (11.221) can be estimated by applying (11.267), (11.268), (11.271) and (11.272).

     Combining  finally the above results with the earlier results (11.243), (11.244),(11.245), (11.251), (11.252) and also (11.257)-(11.260), yields, through 
(11.221), the following proposition.

$\textbf{Proposition 11.3}$ Let the hypotheses $\textbf{H0}$, $\textbf{H1}$, $\textbf{H2}^{\prime}$ and the estimate (6.177) hold. Let also the bootstrap 
assumptions $\textbf{E}_{\{l_{*}+2\}}$, $\textbf{E}^{Q}_{\{l_{*}+1\}}$, $\textbf{E}^{QQ}_{\{l_{*}\}}$ and $\slashed{\textbf{X}}_{[l_{*}]}$, hold for some 
non-negative integer $l$. Moreover, let the bootstrap assumption $\textbf{M}_{[m+1,l_{*}+1]}$ hold for some non-negative integer $m\leq l$. Then if $\delta_{0}$
is suitably small (depending on $l$) we have:
\begin{align}
 \|\check{\dot{g}}^{\prime}\|_{2,[m,l],\Sigma_{t}^{\epsilon_{0}}}\\\notag
\leq C_{l}(1+t)^{-2}\{[1+\log(1+t)](\mathcal{W}_{\{l+2\}}+\mathcal{W}^{Q}_{\{l+1\}}+\mathcal{W}^{QQ}_{\{l\}})\\\notag
+\delta_{0}(1+t)^{-1}[\mathcal{B}_{[m+1,l+1]}+[1+\log(1+t)](\mathcal{Y}_{0}+(1+t)\mathcal{A}_{[l]})]\}
\end{align}\vspace{7mm}

     We proceed to estimate in $L^{2}(\Sigma_{t}^{\epsilon_{0}})$ the functions $\leftexp{(i_{1}...i_{n})}{\hat{g}}^{\prime}_{m,n}$, defined by (9.252)-(9.253),
for $n=l-m$. Note that $l$ have different meanings here and Chapter 9. We shall now re-express the functions $\leftexp{(i_{1}...i_{n})}{\dot{g}}^{\prime}_{m,n}$ in
terms of the function $\check{\dot{g}}^{\prime}$, which has been estimated by Proposition 11.3. Let us define the functions:
\begin{align}
 \leftexp{(i_{1}...i_{n})}{w}^{\prime}_{m,n}=\sum_{k=0}^{m-1}R_{i_{n}}...R_{i_{1}}(T)^{k}\Lambda x^{\prime}_{m-k-1,0}
-m\Lambda\leftexp{(i_{1}...i_{n})}{x}^{\prime}_{m-1,n}\\
\leftexp{(i_{1}...i_{n})}{w}^{\prime\prime}_{m,n}=\sum_{k=0}^{n-1}R_{i_{n}}...R_{i_{n-k+1}}\leftexp{(R_{i_{n-k}})}{Z}
\leftexp{(i_{1}...i_{n-k-1})}{x}^{\prime}_{m,n-k-1}\\\notag
-\sum_{k=0}^{n-1}\leftexp{(R_{i_{n-k}})}{Z}\leftexp{(i_{1}\overset{>i_{n-k}<}{...}i_{n})}{x}^{\prime}_{m,n-1}
\end{align}
 According to (9.123), (9.131) and (9.134), (9.137), the functions $\leftexp{(i_{1}...i_{n})}{w}^{\prime}_{m,n}$, $\leftexp{(i_{1}...i_{n})}{w}^{\prime\prime}_{m,n}$
do not contain terms of the top order $l+2$. From Proposition 9.1 and 9.2 we then have:
\begin{align}
 \leftexp{(i_{1}...i_{n})}{g}^{\prime}_{m,n}-m\Lambda\leftexp{(i_{1}...i_{n})}{x}^{\prime}_{m-1,n}-
\sum_{k=0}^{n-1}\leftexp{(R_{i_{n-k}})}{Z}\leftexp{(i_{1}\overset{>i_{n-k}<}{...}i_{n})}{x}^{\prime}_{m,n-1}\\\notag
=R_{i_{n}}...R_{i_{1}}(T)^{m}\check{g}^{\prime}+\leftexp{(i_{1}...i_{n})}{w}^{\prime}_{m,n}+\leftexp{(i_{1}...i_{n})}{w}^{\prime\prime}_{m,n}\\\notag
+\sum_{k=0}^{m-1}R_{i_{n}}...R_{i_{1}}(T)^{k}y^{\prime}_{m-k,0}+\sum_{k=0}^{n-1}R_{i_{n}}...R_{i_{n-k+1}}\leftexp{(i_{1}...i_{n-k})}{y}^{\prime}_{m,n-k}
\end{align}
We now substitute for $\check{g}^{\prime}$ in terms of $\check{\dot{g}}^{\prime}$ from (11.220). Consider first the case $m=0$. We have:
\begin{align}
 R_{i_{n}}...R_{i_{1}}\check{g}^{\prime}=R_{i_{n}}...R_{i_{1}}\check{\dot{g}}^{\prime}+\xi\cdot\leftexp{(i_{1}...i_{n})}{x}_{n}+\leftexp{(i_{1}...i_{n})}{u}_{n}
\end{align}
where 
\begin{align}
\leftexp{(i_{1}...i_{n})}{u}_{n}=\xi\cdot\sum_{s_{1}\slashed{=}\emptyset}((R)^{s_{1}}\mu)\slashed{d}(R)^{s_{2}}\textrm{tr}\chi+\sum_{s_{1}\slashed{=}\emptyset}
(\slashed{\mathcal{L}}_{R})^{s_{1}}\xi\cdot(\slashed{\mathcal{L}}_{R})^{s_{2}}x_{0}
\end{align}
Then comparing with (9.252) we obtain:
\begin{align}
 \leftexp{(i_{1}...i_{n})}{\dot{g}}^{\prime}_{0,n}=R_{i_{n}}...R_{i_{1}}\check{\dot{g}}^{\prime}+\leftexp{(i_{1}...i_{n})}{u}_{n}
+\leftexp{(i_{1}...i_{n})}{w}^{\prime\prime}_{0,n}\\\notag
+\sum_{k=0}^{n-1}R_{i_{n}}...R_{i_{n-k+1}}\leftexp{(i_{1}...i_{n-k})}{y}^{\prime}_{0,n-k}
\end{align}
A similar calculation yields:
\begin{align}
 \leftexp{(i_{1}...i_{n})}{\dot{g}}^{\prime}_{m,n}=R_{i_{n}}...R_{i_{1}}(T)^{m}\check{\dot{g}}^{\prime}\\\notag
+\leftexp{(i_{1}...i_{n})}{u}^{\prime}_{m-1,n}+\leftexp{(i_{1}...i_{n})}{w}^{\prime}_{m,n}+\leftexp{(i_{1}...i_{n})}{w}^{\prime\prime}_{m,n}\\\notag
+\sum_{k=0}^{m-1}R_{i_{n}}...R_{i_{1}}(T)^{k}y^{\prime}_{m-k,0}+\sum_{k=0}^{n-1}R_{i_{n}}...R_{i_{n-k+1}}\leftexp{(i_{1}...i_{n-k})}{y}^{\prime}_{m,n-k}
\end{align}
where
\begin{align}
 \leftexp{(i_{1}...i_{n})}{u}^{\prime}_{m-1,n}=R_{i_{n}}...R_{i_{1}}(T)^{m-1}u^{\prime}_{0,0}+\leftexp{(i_{1}...i_{n})}{u}^{\prime\prime}_{m-1,n}
\end{align}
\begin{align}
 \leftexp{(i_{1}...i_{n})}{u}^{\prime\prime}_{m-1,n}=R_{i_{n}}...R_{i_{1}}u^{\prime\prime}_{m-1,0}\\\notag
+\xi\cdot\sum_{s_{1}\slashed{=}\emptyset}\slashed{d}(((R)^{s_{1}}\mu)(R)^{s_{2}}(T)^{m-1}\slashed{\Delta}\mu)\\\notag
+\sum_{s_{1}\slashed{=}\emptyset}(\slashed{\mathcal{L}}_{R})^{s_{1}}\xi\cdot\slashed{d}(R)^{s_{2}}x^{\prime}_{m-1,0}
\end{align}
\begin{align}
 u^{\prime\prime}_{m-1,0}=\sum_{k=1}^{m-1}\frac{(m-1)!}{k!(m-1-k)!}\{\xi\cdot\slashed{d}((T)^{k}\mu)(T)^{m-1-k}\slashed{\Delta}\mu)\\\notag
+((\slashed{\mathcal{L}}_{T})^{k}\xi)\cdot\slashed{d}(T)^{m-1-k}x^{\prime}_{0,0}\}
\end{align}
\begin{align}
u^{\prime}_{0,0}=\xi\cdot\dot{x}_{0}+(\slashed{\mathcal{L}}_{T}\xi)\cdot x_{0} 
\end{align}
where
\begin{align}
 \dot{x}_{0}=\mu\slashed{d}(T\textrm{tr}\chi^{\prime}-\slashed{\Delta}\mu)+(T\mu)\slashed{d}\textrm{tr}\chi^{\prime}-(\slashed{d}\mu)\slashed{\Delta}\mu
+\slashed{d}(T\check{f}+\check{f}^{\prime})
\end{align}
Then by the earlier results in this chapter, we finally obtain:

$\textbf{Proposition 11.4}$ Under the assumptions of Proposition 11.3, augmented by the assumptions $\slashed{X}_{[(l+1)_{*}]}$ and $\textbf{M}_{[m,(l+1)_{*}+1]}$
we have:
\begin{align}
 \max_{i_{1}...i_{l-m}}\|\leftexp{(i_{1}...i_{l-m})}{\dot{g}}^{\prime}_{m,l-m}\|_{L^{2}(\Sigma_{t}^{\epsilon_{0}})}\leq \\\notag
C_{l}(1+t)^{-2}\{[1+\log(1+t)](\mathcal{W}_{\{l+2\}}+\mathcal{W}^{Q}_{\{l+1\}}+\mathcal{W}^{QQ}_{\{l\}})\\\notag
+\delta_{0}[(1+t)^{-1}\mathcal{B}_{[m+1,l+1]}+\mathcal{Y}_{0}+(1+t)\mathcal{A}_{[l]}]\}
\end{align}
In deriving this proposition, we made use of the following lemma:

$\textbf{Lemma 11.22}$ Let $A_{1}$,...,$A_{n}$ and $B$ be linear operators acting on some space $X$. We then have:
\begin{align*}
 [A_{n}...A_{1},B]=\sum_{m=0}^{n-1}A_{n}...A_{n-m+1}[A_{n-m},B]A_{n-m-1}...A_{1}
\end{align*}
$Proof$. We just applying induction on $m$, and then use Proposition 8.2. $\qed$
\vspace{7mm}

     Now we can give an estimate for $\leftexp{(i_{1}...i_{l-m})}{Q}^{\prime}_{m,l-m}$ defined by (9.266). Here we just need to estimate the third and the fourth 
term on the right of (9.266). This is straightforward. We obtain:

$\textbf{Proposition 11.5}$ Under the assumptions of Proposition 11.3, augmented by the assumptions $\slashed{X}_{[(l+1)_{*}]}$ and $\textbf{M}_{[m,(l+1)_{*}+1]}$,
we have:
\begin{align*}
 \max_{(i_{1}...i_{l-m})}{Q}^{\prime}_{m,l-m}(t)\\\notag
\leq C_{l}(1+t)^{-3}\{[1+\log(1+t)](\mathcal{W}_{\{l+2\}}+\mathcal{W}^{Q}_{\{l+1\}}+\mathcal{W}^{QQ}_{\{l\}})\\\notag
+\delta_{0}[(1+t)^{-1}\mathcal{B}_{[m+1,l+1]}+\mathcal{Y}_{0}+(1+t)\mathcal{A}_{[l]}]\}
\end{align*}

\subsection{Bounds for $P^{\prime}_{m,l}$}
     We now consider the principal term in the defining expression (9.269) for the quantity $\leftexp{(i_{1}...i_{n})}{B}^{\prime}_{m,n}$, namely, the term:
\begin{align*}
 C(1+t)^{-1}\{\leftexp{(i_{1}...i_{n})}{\bar{P}'^{(0)}}_{m,n,a}(t)+(1+t)^{-1/2}\leftexp{(i_{1}...i_{n})}{\bar{P}}^{\prime(1)}_{m,n,a}(t)\}\bar{\mu}^{-a}_{m}(t)
\end{align*}
Here, the quantities $\leftexp{(i_{1}...i_{n})}{\bar{P}}^{\prime(0)}_{m,n,a}$, $\leftexp{(i_{1}...i_{n})}{\bar{P}}^{\prime(1)}_{m,n,a}$ are defined by (9.225) and (9.226),
in terms of\\ $\leftexp{(i_{1}...i_{n})}{P}^{\prime(0)}_{m,n}$ and $\leftexp{(i_{1}...i_{n})}{P}^{\prime(1)}_{m,n}$, whose sum bounds the quantity $\leftexp{(i_{1}...i_{n})}
{P}^{\prime}_{m,n}$ itself defined by (9.223):
\begin{equation}
 \leftexp{(i_{1}...i_{n})}{P}^{\prime}_{m,n}(t)=(1+t)\|\leftexp{(i_{1}...i_{n})}{\check{f}}^{\prime}_{m,n}(t)\|_{L^{2}([0,\epsilon_{0}]\times S^{2})}
\end{equation}
We shall derive a bound for $\leftexp{(i_{1}...i_{n})}{P}^{\prime}_{m,n}(t)$, with $n=l-m$. Since
\begin{align*}
 \leftexp{(i_{1}...i_{n})}{\check{f}}^{\prime}_{m,n}=R_{i_{n}}...R_{i_{1}}(T)^{m}\check{f}^{\prime}
\end{align*}
by (8.333) we have:
\begin{align}
 \leftexp{(i_{1}...i_{l-m})}{P}^{\prime}_{m,l-m}\leq C\|\leftexp{(i_{1}...i_{l-m})}{\check{f}'}_{m,l-m}\|_{L^{2}(\Sigma_{t}^{\epsilon_{0}})}
=C\|R_{i_{l-m}}...R_{i_{1}}(T)^{m}\check{f}^{\prime}\|_{L^{2}(\Sigma_{t}^{\epsilon_{0}})}
\end{align}
Let us recall (9.62):
\begin{align*}
 \check{f}^{\prime}=f^{\prime}_{0}+\mu^{2}f^{\prime}_{1}
\end{align*}
where
\begin{align*}
 f^{\prime}_{0}=\frac{1}{2}\frac{dH}{dh}(\underline{L}Th)\\
f^{\prime}_{1}=\frac{1}{2\eta^{2}}(\frac{\rho}{\rho^{\prime}})^{\prime}(\slashed{\Delta}h)+\eta^{-1}\hat{T}^{i}(\slashed{\Delta}\psi_{i})
\end{align*}
 First, we consider the contribution from $f'_{0}$. Since:
\begin{align*}
 Th=T\psi_{0}-\sum_{i}\psi_{i}T\psi_{i}\\
\underline{L}Th=\underline{L}T\psi_{0}-\sum_{i}\underline{L}\psi_{i}T\psi_{i}-\sum_{i}\psi_{i}\underline{L}T\psi_{i}
\end{align*}
defining as in Chapter 10,
\begin{equation}
 m^{0}_{T}=\frac{1}{2}\frac{dH}{dh},\quad m^{i}_{T}=-\frac{1}{2}\frac{dH}{dh}\psi_{i}
\end{equation}
we have:
\begin{equation}
  m=m^{\alpha}_{T}(T\psi_{\alpha})
\end{equation}
We decompose:
\begin{align}
 f^{\prime}_{0}=f^{\prime}_{0,P}+f^{\prime}_{0,N}
\end{align}
where $f^{\prime}_{0,P}$ is the principal part of $f^{\prime}_{0}$:
\begin{equation}
 f^{\prime}_{0,P}=m^{\alpha}_{T}(\eta^{-1}\kappa LT\psi_{\alpha}+2T^{2}\psi_{\alpha})
\end{equation}
and the lower order part:
\begin{align}
 f^{\prime}_{0,N}=-\frac{1}{2}\frac{dH}{dh}\sum_{i}\underline{L}\psi_{i}T\psi_{i}
\end{align}
We have:
\begin{align}
 R_{i_{l-m}}...R_{i_{1}}(T)^{m}f^{\prime}_{0}=2m^{\alpha}_{T}R_{i_{l-m}}...R_{i_{1}}(T)^{m+2}\psi_{\alpha}\\\notag
+\alpha^{-1}\kappa m_{T}^{\alpha}R_{i_{l-m}}...R_{i_{1}}(T)^{m}LT\psi_{\alpha}+ \leftexp{(i_{1}...i_{l-m})}{n}^{\prime}_{m,l-m}
\end{align}
where $\leftexp{(i_{1}...i_{l-m})}{n}^{\prime}_{m,l-m}$, of order $l+1$, is given by:
\begin{align}
\leftexp{(i_{1}...i_{l-m})}{n}^{\prime}_{m,l-m}=2R_{i_{l-m}}...R_{i_{1}}\{\sum_{k=1}^{m}\frac{m!}{k!(m-k)!}((T)^{k}m^{\alpha}_{T})((T)^{m-k+2}\psi_{\alpha})\}\\\notag
+2\sum_{s_{1}\slashed{=}\emptyset}((R)^{s_{1}}m^{\alpha}_{T})((R)^{s_{2}}(T)^{m+2}\psi_{\alpha})\\\notag
+R_{i_{l-m}}...R_{i_{1}}\{\sum_{k=1}^{m}\frac{m!}{k!(m-k)!}((T)^{k}(\alpha^{-1}\kappa m^{\alpha}_{T}))((T)^{m-k}LT\psi_{\alpha})\}\\\notag
+\sum_{s_{1}\slashed{=}\emptyset}((R)^{s_{1}}(\alpha^{-1}\kappa m^{\alpha}_{T}))((R)^{s_{2}}(T)^{m}LT\psi_{\alpha})\\\notag
+R_{i_{l-m}}...R_{i_{1}}(T)^{m}f^{\prime}_{0,N}
\end{align}
From (11.294), taking into account (10.184), we have:
\begin{align}
 \|R_{i_{l-m}}...R_{i_{1}}(T)^{m}f'_{0}\|_{L^{2}(\Sigma_{t}^{\epsilon_{0}})}\leq (|\ell|+C\delta_{0}(1+t)^{-1})
\sum_{\alpha}\|R_{i_{l-m}}...R_{i_{1}}(T)^{m+2}\psi_{\alpha}\|_{L^{2}(\Sigma_{t}^{\epsilon_{0}})}\\\notag
+C\sum_{\alpha}\|\mu R_{i_{l-m}}...R_{i_{1}}(T)^{m}LT\psi_{\alpha}\|_{L^{2}(\Sigma_{t}^{\epsilon_{0}})}+\|\leftexp{(i_{1}...i_{l-m})}{n}^{\prime}_{m,l-m}\|
_{L^{2}(\Sigma_{t}^{\epsilon_{0}})} 
\end{align}
Here, the constant $C$ does not depend on $\ell$. Moreover, using the estimate for $\Lambda$ in Corollary 11.2.c, we readily deduce:
\begin{align}
 \max_{i_{1}...i_{l-m}}\|\leftexp{(i_{1}...i_{l-m})}{n}^{\prime}_{m,l-m}\|_{L^{2}(\Sigma_{t}^{\epsilon_{0}})}\\\notag
\leq C_{l}\delta_{0}(1+t)^{-1}\{\mathcal{W}_{\{l+1\}}+[1+\log(1+t)]\mathcal{W}^{Q}\\\notag
+(1+t)^{-1}[\mathcal{B}_{[m,l]}+\delta_{0}(\mathcal{Y}_{0}+(1+t)\mathcal{A}_{[l-1]})]\}
\end{align}

    Consider next the contribution of the term $\mu^{2}f^{\prime}_{1}$ in $\check{f}'$ to $R_{i_{l-m}}...R_{i_{1}}(T)^{m}\check{f}^{\prime}$. Also, 
we only consider the principal part. Since:
\begin{align*}
 \slashed{d}h=\slashed{d}\psi_{0}-\sum_{i}\psi_{i}\slashed{d}\psi_{i}\\
\slashed{\Delta}h=\slashed{\Delta}\psi_{0}-\sum_{i}\psi_{i}\slashed{\Delta}\psi_{i}-\sum_{i}\slashed{d}\psi_{i}\cdot\slashed{d}\psi_{i}
\end{align*}
we have:
\begin{align}
 \mu^{2}f^{\prime}_{1}=\frac{\mu^{2}}{2\eta^{2}}(\frac{\rho}{\rho^{\prime}})^{\prime}
(\slashed{\Delta}\psi_{0}-\sum_{i}\psi_{i}\slashed{\Delta}\psi_{i}
-\sum_{i}\slashed{d}\psi_{i}\cdot\slashed{d}\psi_{i})\\\notag
+\eta^{-1}\mu^{2}\hat{T}^{i}(\slashed{\Delta}\psi_{i})
\end{align}
if we define:
\begin{align}
 \tilde{m}^{0}=\frac{1}{2\eta^{2}}(\frac{\rho}{\rho^{\prime}})^{\prime},\quad 
\tilde{m}^{i}=-\frac{1}{2\eta^{2}}(\frac{\rho}{\rho^{\prime}})^{\prime}\psi_{i}+\eta^{-1}\hat{T}^{i}
\end{align}

We have the decomposition:
\begin{align}
 f^{\prime}_{1}=f^{\prime}_{1,P}+f^{\prime}_{1,N}
\end{align}
The principal part of $f'_{1}$ as follows:
\begin{align}
 f'_{1,P}=\tilde{m}^{\alpha}\slashed{\Delta}\psi_{\alpha}
\end{align}
and $f^{\prime}_{1,N}$ is the lower order term:
\begin{align}
 f^{\prime}_{1,N}=-\frac{1}{2\eta^{2}}(\frac{\rho}{\rho^{\prime}})^{\prime}\sum_{i}\slashed{d}\psi_{i}\cdot\slashed{d}\psi_{i}
\end{align}
We then obtain:
\begin{align}
 R_{i_{l-m}}...R_{i_{1}}(T)^{m}(\mu^{2}f^{\prime}_{1})=\mu^{2}\tilde{m}^{\alpha}R_{i_{l-m}}...R_{i_{1}}(T)^{m}\slashed{\Delta}\psi_{\alpha}
+\leftexp{(i_{1}...i_{l-m})}{n}^{\prime\prime}_{m,l-m}
\end{align}
where $\leftexp{(i_{1}...i_{l-m})}{n}^{\prime\prime}_{m,l-m}$ is of order $l+1$:
\begin{align}
 \leftexp{(i_{1}...i_{l-m})}{n}^{\prime\prime}_{m,l-m}=R_{i_{l-m}}...R_{i_{1}}\{\sum_{k=1}^{m}\frac{m!}{k!(m-k)!}((T)^{k}(\mu^{2}\tilde{m}^{\alpha}))
((T)^{m-k}\slashed{\Delta}\psi_{\alpha})\}\\\notag
+\sum_{s_{1}\slashed{=}\emptyset}((R)^{s_{1}}(\mu^{2}\tilde{m}^{\alpha}))((R)^{s_{2}}(T)^{m}\slashed{\Delta}\psi_{\alpha})\\\notag
+R_{i_{l-m}}...R_{i_{1}}(T)^{m}(\mu^{2}f^{\prime}_{1,N})
\end{align}
In view of (10.303) we have:
\begin{align}
 \|R_{i_{l-m}}...R_{i_{1}}(T)^{m}(\mu^{2}f^{\prime}_{1})\|_{L^{2}(\Sigma_{t}^{\epsilon_{0}})}\\\notag
\leq C[1+\log(1+t)]\sum_{\alpha}\|\mu R_{i_{l-m}}...R_{i_{1}}(T)^{m}\slashed{\Delta}\psi_{\alpha}\|_{L^{2}(\Sigma_{t}^{\epsilon_{0}})}
+\|\leftexp{(i_{1}...i_{l-m})}{n}^{\prime\prime}_{m,l-m}\|_{L^{2}(\Sigma_{t}^{\epsilon_{0}})}
\end{align}
where $C$ is independent of $\ell$.

The lower order term $\leftexp{(i_{1}...i_{l-m})}{n}^{\prime\prime}_{m,l-m}$ can be estimated straightforwardly, by using Proposition 11.1 and Proposition 11.2,
as well as their corollaries.
\begin{align}
 \max_{i_{1}...i_{l-m}}\|\leftexp{(i_{1}...i_{l-m})}{n}^{\prime\prime}_{m,l-m}\|_{L^{2}(\Sigma_{t}^{\epsilon_{0}})}
\leq C_{l}(1+t)^{-2}[1+\log(1+t)]\cdot\\\notag
\{[1+\log(1+t)](\mathcal{W}_{\{l+1\}}+\delta_{0}(1+t)^{-3}[1+\log(1+t)]^{2}\mathcal{W}^{Q}_{\{l-1\}})\\\notag
+\delta_{0}(1+t)^{-1}[\mathcal{B}_{[m,l]}+[1+\log(1+t)](\mathcal{Y}_{0}+(1+t)\mathcal{A}_{[l-1]})]\}
\end{align}

     We now return to (11.296). Our object is to express
\begin{align*}
 R_{i_{l-m}}...R_{i_{1}}(T)^{m+2}\psi_{\alpha}\quad\textrm{and}\quad R_{i_{l-m}}...R_{i_{1}}(T)^{m}LT\psi_{\alpha}
\end{align*}
as
\begin{align*}
 TR_{i_{l-m}}...R_{i_{1}}(T)^{m+1}\psi_{\alpha}\quad\textrm{and}\quad LR_{i_{l-m}}...R_{i_{1}}(T)^{m+1}\psi_{\alpha}
\end{align*}
respectively, plus lower order terms. Then by (10.217), we can bound the first two terms on the right in (11.296) by:
\begin{align}
 (|\ell|+C\delta_{0}(1+t)^{-1})\sqrt{\sum_{\alpha}\mathcal{E}_{0}[R_{i_{l-m}}...R_{i_{1}}(T)^{m+1}\psi_{\alpha}]}\\\notag
+C(1+t)^{-1}[1+\log(1+t)]^{1/2}\sqrt{\sum_{\alpha}\mathcal{E}^{\prime}_{1}[R_{i_{l-m}}...R_{i_{1}}(T)^{m+1}\psi_{\alpha}]}\\\notag
+C(1+t)^{-1}[1+\log(1+t)]\mathcal{W}_{\{l+1\}}
\end{align}
up to a bound for the lower order terms of a form similar to (11.297).

    In view of (10.204) and (10.206):
\begin{align}
 TR_{i_{n}}...R_{i_{1}}-R_{i_{n}}...R_{i_{1}}T\\\notag
=\sum_{k=0}^{n-1}R_{i_{n}}...R_{i_{n-k+1}}\leftexp{(R_{i_{n-k}})}\Theta R_{i_{n-k-1}}...R_{i_{1}}
\end{align}
we have:
\begin{align}
 R_{i_{l-m}}...R_{i_{1}}(T)^{m+2}\psi_{\alpha}=TR_{i_{l-m}}...R_{i_{1}}(T)^{m+1}\psi_{\alpha}\\\notag
-\sum_{k=0}^{l-m-1}R_{i_{l-m}}...R_{i_{l-m-k+1}}\leftexp{(R_{i_{l-m-k}})}{\Theta}R_{i_{l-m-k-1}}...R_{i_{1}}(T)^{m+1}\psi_{\alpha}
\end{align}
To estimate in $L^{2}(\Sigma_{t}^{\epsilon_{0}})$ the sum on the right, we observe that it is given by:
\begin{align}
 \sum_{k=0}^{l-m-1}\sum_{|s_{1}|+|s_{2}|=k}(\slashed{\mathcal{L}}_{R})^{s_{1}}\leftexp{(R_{i_{l-m-k}})}{\Theta}\cdot\slashed{d}(R)^{s^{\prime}_{2}}
(T)^{m+1}\psi_{\alpha}
\end{align}
Here, in the inner sum, we are considering all ordered partitions $\{s_{1},s_{2}\}$ of the set $\{l-m-k+1,...,l-m\}$ into two  ordered subsets $s_{1}$, $s_{2}$.
Also,
\begin{align*}
 s^{\prime}_{2}=s_{2}\bigcup\{1,...,l-m-k-1\}
\end{align*}
To estimate in $L^{2}(\Sigma_{t}^{\epsilon_{0}})$ a term in the inner sum in (11.310), we use Lemma 10.20. So we obtain that (11.310) is bounded by:
\begin{align}
 C_{l}\delta_{0}(1+t)^{-2}\{[1+\log(1+t)]^{2}\mathcal{W}_{\{l+1\}}+\max_{i}\|\leftexp{(R_{i})}{\Theta}\|_{2,[l-1],\Sigma_{t}^{\epsilon_{0}}}\}\\\notag
\leq C_{l}\delta_{0}(1+t)^{-2}[1+\log(1+t)]\cdot\\\notag
\{[1+\log(1+t)]\mathcal{W}_{\{l+1\}}+\delta_{0}(1+t)^{-1}\mathcal{B}_{[0,l]}+\mathcal{Y}_{0}+(1+t)\mathcal{A}_{[l-1]}\}
\end{align}
Then it follows through (11.309) that:
\begin{align}
 \sum_{\alpha}\|R_{i_{l-m}}...R_{i_{1}}(T)^{m+2}\psi_{\alpha}\|_{L^{2}(\Sigma_{t}^{\epsilon_{0}})}\\\notag
\leq \sum_{\alpha}\|TR_{i_{l-m}}...R_{i_{1}}(T)^{m+1}\psi_{\alpha}\|_{L^{2}(\Sigma_{t}^{\epsilon_{0}})}\\\notag
+C_{l}\delta_{0}(1+t)^{-2}[1+\log(1+t)]\cdot\\\notag
\{[1+\log(1+t)]\mathcal{W}_{\{l+1\}}+\delta_{0}(1+t)^{-1}\mathcal{B}_{[0,l]}+\mathcal{Y}_{0}+(1+t)\mathcal{A}_{[l-1]}\}
\end{align}

     Similarly, we have:
\begin{align}
 LR_{i_{n}}...R_{i_{1}}-R_{i_{n}}...R_{i_{1}}L=\sum_{k=0}^{n-1}R_{i_{n}}...R_{i_{n-k+1}}\leftexp{(R_{i_{n-k}})}{Z}R_{i_{n-k-1}}...R_{i_{1}}
\end{align}
It follows that, taking $n=l-m$,
\begin{align}
 R_{i_{l-m}}...R_{i_{1}}L(T)^{m+1}\psi_{\alpha}=LR_{i_{l-m}}...R_{i_{1}}(T)^{m+1}\psi_{\alpha}\\\notag
-\sum_{k=0}^{l-m-1}R_{i_{l-m}}...R_{i_{l-m-k+1}}\leftexp{(R_{i_{l-m-k}})}{Z}R_{i_{l-m-k-1}}...R_{i_{1}}(T)^{m+1}\psi_{\alpha}
\end{align}
By Lemma 11.22 we have:
\begin{align}
 [(T)^{m},L]=-\sum_{k=0}^{m-1}(T)^{k}\Lambda(T)^{m-k-1}
\end{align}
hence:
\begin{align}
 R_{i_{l-m}}...R_{i_{1}}(T)^{m}LT\psi_{\alpha}=R_{i_{l-m}}...R_{i_{1}}L(T)^{m+1}\psi_{\alpha}\\\notag
-\sum_{k=0}^{m-1}R_{i_{l-m}}...R_{i_{1}}(T)^{k}\Lambda(T)^{m-k}\psi_{\alpha}
\end{align}
Now, $\mu$ times the sum on the right in (11.314) is bounded in $L^{2}(\Sigma_{t}^{\epsilon_{0}})$ by:
\begin{align}
 C[1+\log(1+t)]\sum_{k=0}^{l-m-1}\|\leftexp{(R_{i_{l-m-k}})}{Z}\cdot\slashed{d}R_{i_{l-m-k-1}}...R_{i_{1}}(T)^{m+1}\psi_{\alpha}\|_{2,[k],\Sigma_{t}^{\epsilon_{0}}}
\end{align}
Using Corollary 10.1.i and Corollary 10.2.i we obtain that this is bounded by:
\begin{align}
 C_{l}\delta_{0}(1+t)^{-2}[1+\log(1+t)]\{[1+\log(1+t)]\mathcal{W}_{\{l+1\}}\\\notag
+\mathcal{Y}_{0}+(1+t)\mathcal{A}_{[l-1]}\}
\end{align}
Also, $\mu$ times the sum on the right in (11.316) is bounded in $L^{2}(\Sigma_{t}^{\epsilon_{0}})$ by:
\begin{align}
 C[1+\log(1+t)]\sum_{k=0}^{m-1}\|\Lambda\cdot\slashed{d}(T)^{m-k}\psi_{\alpha}\|_{2,[k,l-m+k],\Sigma_{t}^{\epsilon_{0}}}
\end{align}
Using Corollary 11.1.c and Corollary 11.2.c we obtain that the above is bounded by:
\begin{align}
 C_{l}\delta_{0}(1+t)^{-2}[1+\log(1+t)]\cdot\\\notag
\{[1+\log(1+t)](\mathcal{W}_{\{l+1\}}+\delta_{0}(1+t)^{-3}[1+\log(1+t)]\mathcal{W}^{Q}_{\{l-1\}})\\\notag
+(1+t)^{-1}[\mathcal{B}_{[m-1,l]}+\delta_{0}[1+\log(1+t)](\mathcal{Y}_{0}+(1+t)\mathcal{A}_{[l-2]})]\}
\end{align}
In view of the estimates (11.318) and (11.320), it follows through (11.314), (11.319) that:
\begin{align}
 \sum_{\alpha}\|\mu R_{i_{l-m}}...R_{i_{1}}(T)^{m}LT\psi_{\alpha}\|_{L^{2}(\Sigma_{t}^{\epsilon_{0}})}\\\notag
\leq \sum_{\alpha}\|\mu LR_{i_{l-m}}...R_{i_{1}}(T)^{m+1}\psi_{\alpha}\|_{L^{2}(\Sigma_{t}^{\epsilon_{0}})}\\\notag
+C_{l}\delta_{0}(1+t)^{-2}[1+\log(1+t)]\cdot\\\notag
\{[1+\log(1+t)]\mathcal{W}_{\{l+1\}}+\delta_{0}(1+t)^{-3}[1+\log(1+t)]^{2}\mathcal{W}^{Q}_{\{l-1\}}\\\notag
+(1+t)^{-1}\mathcal{B}_{[m-1,l]}+\mathcal{Y}_{0}+(1+t)\mathcal{A}_{[l-1]}\}
\end{align}
    In view of (11.307), inequalities (11.312), (11.321), together with the estimate (11.297) imply through (11.296) that:
\begin{align}
 \|R_{i_{l-m}}...R_{i_{1}}(T)^{m}f^{\prime}_{0}\|_{L^{2}(\Sigma_{t}^{\epsilon_{0}})}\\\notag
\leq (|\ell|+\delta_{0}(1+t)^{-1})\sqrt{\sum_{\alpha}\mathcal{E}_{0}[R_{i_{l-m}}...R_{i_{1}}(T)^{m+1}\psi_{\alpha}]}\\\notag
+C(1+t)^{-1}[1+\log(1+t)]^{1/2}\sqrt{\sum_{\alpha}\mathcal{E}^{\prime}_{1}[R_{i_{l-m}}...R_{i_{1}}(T)^{m+1}\psi_{\alpha}]}\\\notag
+C(1+t)^{-1}[1+\log(1+t)]\mathcal{W}_{\{l+1\}}+C_{l}\delta_{0}(1+t)^{-1}\{[1+\log(1+t)]\mathcal{W}^{Q}_{\{l\}}\\\notag
+(1+t)^{-1}[\mathcal{B}_{[m,l]}+[1+\log(1+t)](\mathcal{Y}_{0}+(1+t)\mathcal{A}_{[l-1]})]\}
\end{align}

    We finally turn to (11.305). Our object here is to express $R_{i_{l-m}}...R_{i_{1}}(T)^{m}\slashed{\Delta}\psi_{\alpha}$ as \\
$\slashed{\Delta}(R_{i_{l-m}}...R_{i_{1}}(T)^{m}\psi_{\alpha})$ plus lower order terms. Since:
\begin{align*}
 \slashed{\Delta}\psi_{\alpha}=\textrm{tr}(\slashed{g}^{-1})\cdot(\slashed{D}^{2}\psi_{\alpha})
\end{align*}
we have:
\begin{align}
 R_{i_{l-m}}...R_{i_{1}}(T)^{m}\slashed{\Delta}\psi_{\alpha}=R_{i_{l-m}}...R_{i_{1}}(T)^{m}(\slashed{g}^{-1})^{AB}(\slashed{D}^{2}\psi_{\alpha})_{AB}\\\notag
=(\slashed{g}^{-1})^{AB}(\slashed{\mathcal{L}}_{R_{i_{l-m}}}...\slashed{\mathcal{L}}_{R_{i_{1}}}(\slashed{\mathcal{L}}_{T})^{m}(\slashed{D}^{2}\psi_{\alpha})_{AB})
\\\notag
+\slashed{\mathcal{L}}_{R_{i_{l-m}}}...\slashed{\mathcal{L}}_{R_{i_{1}}}\sum_{k=1}^{m}\frac{m!}{k!(m-k)!}((\slashed{\mathcal{L}}_{T})^{k}(\slashed{g}^{-1}))^{AB}
(\slashed{\mathcal{L}}_{T})^{m-k}(\slashed{D}^{2}\psi_{\alpha})_{AB}\}\\\notag
\sum_{|s_{1}|+|s_{2}|=l-m,s_{1}\slashed{=}\emptyset}((\slashed{\mathcal{L}}_{R})^{s_{1}}(\slashed{g}^{-1}))^{AB}(\slashed{\mathcal{L}}_{R})^{s_{2}}
(\slashed{\mathcal{L}}_{T})^{m}(\slashed{D}^{2}\psi_{\alpha})_{AB}\}
\end{align}
Applying formula (11.83) to the $S_{t,u}$ 1-form $\slashed{d}\psi_{\alpha}$ we obtain:
\begin{align}
 \slashed{\mathcal{L}}_{R_{i_{l-m}}}...\slashed{\mathcal{L}}_{R_{i_{1}}}(\slashed{\mathcal{L}}_{T})^{m}\slashed{D}^{2}\psi_{\alpha}\\\notag
=\slashed{D}^{2}(R_{i_{l-m}}...R_{i_{1}}(T)^{m}\psi_{\alpha})+\leftexp{(i_{1}...i_{l-m})}{c}_{m,l-m}[\slashed{d}\psi_{\alpha}]
\end{align}
Hence the first term on the right in (11.323) is given by:
\begin{align}
 (\slashed{g}^{-1})^{AB}(\slashed{\mathcal{L}}_{R_{i_{l-m}}}...\slashed{\mathcal{L}}_{R_{i_{1}}}(\slashed{\mathcal{L}}_{T})^{m}(\slashed{D}^{2}
\psi_{\alpha})_{AB})\\\notag
=\slashed{\Delta}(R_{i_{l-m}}...R_{i_{1}}(T)^{m}\psi_{\alpha})+(\slashed{g}^{-1})^{AB}(\leftexp{(i_{1}...i_{l-m})}{c}_{m,l-m}[\slashed{d}\psi_{\alpha}])_{AB}
\end{align}
Taking into account of the fact that:
\begin{align*}
 \slashed{\mathcal{L}}_{T}(\slashed{g}^{-1})=-\slashed{g}^{-1}\cdot\leftexp{(T)}{\slashed{\pi}}\cdot\slashed{g}^{-1}
\end{align*}
and using Corollary 10.2.d, Corollary 11.1.c, Corollary 11.2.c and Corollary 11.2.d, we deduce that $\mu$ times the first sum on the right in (11.323) is bounded in 
$L^{2}(\Sigma_{t}^{\epsilon_{0}})$ by:
\begin{align}
 C_{l}[1+\log(1+t)]\{(1+t)^{-1}[1+\log(1+t)]\|\slashed{D}^{2}\psi_{\alpha}\|_{2,[m-1,l-1],\Sigma_{t}^{\epsilon_{0}}}\\\notag
+\delta_{0}(1+t)^{-3}(\mathcal{T}_{[m,l]}+(1+t)^{-2}\|\leftexp{(R_{j})}{\slashed{\pi}}\|_{2,[l-1],\Sigma_{t}^{\epsilon_{0}}}\}\\\notag
\leq C_{l}(1+t)^{-3}[1+\log(1+t)]\cdot\\\notag
\{[1+\log(1+t)]\mathcal{W}_{\{l+1\}}+\delta_{0}(1+t)^{-2}[1+\log(1+t)]^{2}\mathcal{W}^{Q}_{\{l-1\}}\\\notag
+\delta_{0}(1+t)^{-1}[\mathcal{B}_{[m-1,l]}+[1+\log(1+t)](\mathcal{Y}_{0}+(1+t)\mathcal{A}_{[l-1]})]\}
\end{align}
Similarly, taking into account of the fact that:
\begin{align*}
 \slashed{\mathcal{L}}_{R_{i}}(\slashed{g}^{-1})=-\slashed{g}^{-1}\cdot\leftexp{(R_{i})}{\slashed{\pi}}\cdot\slashed{g}^{-1}
\end{align*}
and using Corollary 10.1.d, Corollary 10.2.d and Corollary 11.2.d, we deduce that $\mu$ times the second sum on the right in (11.323) is bounded in 
$L^{2}(\Sigma_{t}^{\epsilon_{0}})$ by:
\begin{align}
 C_{l}[1+\log(1+t)]\{\delta_{0}(1+t)^{-1}[1+\log(1+t)]\|\slashed{D}^{2}\psi_{\alpha}\|_{2,[m,l-1],\Sigma_{t}^{\epsilon_{0}}}\\\notag
+\delta_{0}(1+t)^{-3}\max_{i}\|\leftexp{(R_{i})}{\slashed{\pi}}\|_{2,[l-1],\Sigma_{t}^{\epsilon_{0}}}\}\\\notag
\leq C_{l}\delta_{0}(1+t)^{-3}[1+\log(1+t)]\{[1+\log(1+t)]\mathcal{W}_{\{l+1\}}+\delta_{0}(1+t)^{-2}[1+\log(1+t)]^{2}\mathcal{W}^{Q}_{\{l-1\}}\\\notag
+\delta_{0}(1+t)^{-2}[1+\log(1+t)]\mathcal{B}_{[m-1,l]}+\mathcal{Y}_{0}+(1+t)\mathcal{A}_{[l-1]}\}
\end{align}
Finally, to estimate in $L^{2}(\Sigma_{t}^{\epsilon_{0}})$ $\mu$ times the second term on the right in (11.325), we use the expression in Lemma 11.9 as well
as Lemma 10.10, Lemma 10.11, Corollary 11.1.c and Corollary 11.2.c, to obtain:
\begin{align*}
 \max_{k\leq m,n\leq l-k}\max_{i_{1}...i_{n}}\|\leftexp{(i_{1}...i_{n})}{c}_{k,n}[\slashed{d}\psi_{\alpha}]\|_{L^{2}(\Sigma_{t}^{\epsilon_{0}})}\\
\leq C_{l}(1+t)^{-1}\{\delta_{0}(1+t)^{-1}[1+\log(1+t)]\|\slashed{d}\psi_{\alpha}\|_{2,[m,l-1],\Sigma_{t}^{\epsilon_{0}}}\\
+\|\slashed{d}\psi_{\alpha}\|_{\infty,[m,l_{*}],\Sigma_{t}^{\epsilon_{0}}}[\mathcal{W}_{\{l+1\}}+(1+t)^{-2}[1+\log(1+t)]^{2}\mathcal{W}^{Q}_{\{l\}}\\
+(1+t)^{-1}\mathcal{B}_{[m-1,l+1]}+\mathcal{Y}_{0}+(1+t)\mathcal{A}_{[l]}]\}
\end{align*}
which gives a bound by:
\begin{align}
 C_{l}\delta_{0}(1+t)^{-3}[1+\log(1+t)]\cdot\\\notag
\{[1+\log(1+t)](\mathcal{W}_{\{l+1\}}+(1+t)^{-2}[1+\log(1+t)]\mathcal{W}^{Q}_{\{l\}})\\\notag
+(1+t)^{-1}\mathcal{B}_{[m-1,l+1]}+\mathcal{Y}_{0}+(1+t)\mathcal{A}_{[l]}\}
\end{align}
In view of (11.326)-(11.328), it follows through (11.323), (11.325) that:
\begin{align}
 \sum_{\alpha}\|\mu R_{i_{l-m}}...R_{i_{1}}(T)^{m}\slashed{\Delta}\psi_{\alpha}\|_{L^{2}(\Sigma_{t}^{\epsilon_{0}})}\\\notag
\leq \sum_{\alpha}\|\mu\slashed{\Delta}(R_{i_{l-m}}...R_{i_{1}}(T)^{m}\psi_{\alpha})\|_{L^{2}(\Sigma_{t}^{\epsilon_{0}})}\\\notag
+C_{l}(1+t)^{-3}[1+\log(1+t)]\cdot\\\notag
\{[1+\log(1+t)](\mathcal{W}_{\{l+1\}}+\delta_{0}(1+t)^{-2}[1+\log(1+t)]\mathcal{W}^{Q}_{\{l\}})\\\notag
+\delta_{0}[(1+t)^{-1}\mathcal{B}_{[m,l+1]}+\mathcal{Y}_{0}+(1+t)\mathcal{A}_{[l]}]\}
\end{align}
Moreover, by $\textbf{H1}$ the first term on the right in (11.329) is bounded by:
\begin{align}
 \sum_{\alpha}\|\mu\slashed{\Delta}(R_{i_{l-m}}...R_{i_{1}}(T)^{m}\psi_{\alpha})\|_{L^{2}(\Sigma_{t}^{\epsilon_{0}})}\\\notag
\leq C(1+t)^{-1}\sum_{\alpha,j}\|\mu\slashed{d}(R_{j}R_{i_{l-m}}...R_{i_{1}}(T)^{m}\psi_{\alpha})\|_{L^{2}(\Sigma_{t}^{\epsilon_{0}})}\\\notag
\leq C(1+t)^{-2}[1+\log(1+t)]^{1/2}\sqrt{\sum_{\alpha,j}\mathcal{E}^{\prime}_{1}[R_{j}R_{i_{l-m}}...R_{i_{1}}(T)^{m}\psi_{\alpha}]}
\end{align}
The bounds (11.329)-(11.330) together with (11.306) imply through (11.303) that:
\begin{align}
 \|R_{i_{l-m}}...R_{i_{1}}(T)^{m}(\mu^{2}f^{\prime}_{1})\|_{L^{2}(\Sigma_{t}^{\epsilon_{0}})}\\\notag
\leq C(1+t)^{-2}[1+\log(1+t)]^{3/2}\sqrt{\sum_{\alpha,j}\mathcal{E}^{\prime}_{1}[R_{j}R_{i_{l-m}}...R_{i_{1}}(T)^{m}\psi_{\alpha}]}\\\notag
+C_{l}(1+t)^{-2}[1+\log(1+t)]\cdot\\\notag
\{[1+\log(1+t)](\mathcal{W}_{\{l+1\}}+\delta_{0}(1+t)^{-3}[1+\log(1+t)]^{2}\mathcal{W}^{Q}_{\{l\}})\\\notag
+\delta_{0}(1+t)^{-1}[\mathcal{B}_{[m,l+1]}+[1+\log(1+t)](\mathcal{Y}_{0}+(1+t)\mathcal{A}_{[l]})]\}
\end{align}
The bounds (11.322) and (11.331) yield finally the following proposition:

$\textbf{Proposition 11.6}$ Under the assumptions of Proposition 11.5 we have:
\begin{align*}
 \leftexp{(i_{1}...i_{l-m})}{P}^{\prime}_{m,l-m}\leq \leftexp{(i_{1}...i_{l-m})}{P}^{\prime(0)}_{m,l-m}
+\leftexp{(i_{1}...i_{l-m})}{P}^{\prime(1)}_{m,l-m}
\end{align*}
where:
\begin{align*}
 \leftexp{(i_{1}...i_{l-m})}{P}^{\prime(0)}_{m,l-m}=|\ell|\sqrt{\sum_{\alpha}\mathcal{E}_{0}[R_{i_{l-m}}...R_{i_{1}}(T)^{m+1}\psi_{\alpha}]}
\end{align*}
and:
\begin{align*}
 \leftexp{(i_{1}...i_{l-m})}{P}^{\prime(1)}_{m,l-m}=C\delta_{0}(1+t)^{-1}\sqrt{\sum_{\alpha}\mathcal{E}_{0}[R_{i_{l-m}}...R_{i_{1}}(T)^{m+1}\psi_{\alpha}]}\\\notag
+C(1+t)^{-1}[1+\log(1+t)]^{1/2}\sqrt{\sum_{\alpha}\mathcal{E}^{\prime}_{1}[R_{i_{l-m}}...R_{i_{1}}(T)^{m+1}\psi_{\alpha}]}\\\notag
+C(1+t)^{-2}[1+\log(1+t)]^{3/2}\sqrt{\sum_{\alpha,j}\mathcal{E}^{\prime}_{1}[R_{j}R_{i_{l-m}}...R_{i_{1}}(T)^{m}\psi_{\alpha}]}\\\notag
+(C(1+t)^{-1}[1+\log(1+t)]+C_{l}(1+t)^{-2}[1+\log(1+t)]^{2})\mathcal{W}_{\{l+1\}}\\\notag
+C_{l}\delta_{0}(1+t)^{-1}\{[1+\log(1+t)]\mathcal{W}^{Q}_{\{l\}}+(1+t)^{-1}[\mathcal{B}_{[m,l+1]}\\
+[1+\log(1+t)](\mathcal{Y}_{0}+(1+t)\mathcal{A}_{[l]})]\}
\end{align*}

\chapter{ Recovery of the Acoustical Assumptions.\\
Estimates for Up to the Next to the Top Order Angular Derivatives \\
of $\chi$ and Spatial Derivatives of $\mu$}




In the first part of the present chapter we shall show that assumptions $\textbf{E}_{\{1\}}$ and $\textbf{E}^{Q}_{\{0\}}$ on 
$W^{s}_{\epsilon_{0}}$, together with certain primitive assumptions on the initial conditions on $\Sigma_{0}^{\epsilon_{0}}$,
imply pointwise estimates for $\kappa$, the functions $\lambda_{i}, y'_{i}, y_{i}$, as well as sharp upper and lower bounds for $r$.
Moreover, we shall show that the same assumptions imply hypothesis $\textbf{H0}$ on $W^{s}_{\epsilon_{0}}$.

    In the second part of the present chapter we shall show that the assumptions $\textbf{E}_{\{l+2\}}, \textbf{E}^{Q}_{\{l+1\}}$ and
$\textbf{E}_{\{l\}}^{QQ}$, on $W^{s}_{\epsilon_{0}}$, together with appropriate assumptions on the initial conditions on $\Sigma_{0}^{\epsilon_{0}}$,
imply the acoustical assumptions $\slashed{\textbf{X}}_{[l]}$, as well as $\textbf{H1}, \textbf{H2}$ and $\textbf{H2}^{\prime}$. 
Moreover, we shall show that under some appropriate assumptions on the initial conditions, the acoustical assumptions $\textbf{M}_{[m,l+1]}$,
for $m=0,...,l+1$, also follow.

     Finally, in the third part of this chapter we shall derive estimates for the quantities $\mathcal{A}_{[l]}$ and $\mathcal{B}_{[m+1,l+1]}$,
for $m=0,...,l$, under the assumptions $\textbf{E}_{\{l_{*}+2\}}, \textbf{E}^{Q}_{\{l_{*}+1\}}, \textbf{E}^{QQ}_{\{l_{*}\}}$, 
on $W^{s}_{\epsilon_{0}}$, together with appropriate assumptions on the initial conditions on $\Sigma_{0}^{\epsilon_{0}}$.

\section{Estimates for $\lambda_{i}$, $y^{\prime}_{i}$, $y_{i}$ and $r$. Establishing the Hypothesis $\textbf{H0}$}

$\textbf{Proposition 12.1}$ Let assumptions $\textbf{E}_{\{1\}}$ and $\textbf{E}_{\{0\}}^{Q}$ hold on $W^{s}_{\epsilon_{0}}$.
Then provided that $\delta_{0}$ is suitably small, for all $t\in[0,s]$ we have:
\begin{align*}
 \|\kappa-1\|_{L^{\infty}(\Sigma^{\epsilon_{0}}_{t})}\leq C\delta_{0}[1+\log(1+t)]
\end{align*}
where the constant $C$ is independent of $s$.

$Proof$. Recall from Chapter 3, we have:
\begin{equation}
 L\kappa=\eta^{-1}m+\kappa e^{\prime}
\end{equation}
where:
\begin{equation}
 m=\frac{1}{2}\frac{dH}{dh}T(h)
\end{equation}
and
\begin{equation}
 e'=\frac{1}{\eta}\frac{d\eta}{dh}Lh-\alpha^{-1}\hat{T}^{i}(L\psi_{i})-\alpha^{-1}L\alpha
\end{equation}
Assumption $\textbf{E}_{\{0\}}$ implies:
\begin{equation}
 \|h\|_{L^{\infty}(\Sigma_{t}^{\epsilon_{0}})}\leq C\delta_{0}(1+t)^{-1}
\end{equation}
and:
\begin{equation}
 \|\eta^{-1}-1\|_{L^{\infty}(\Sigma_{t}^{\epsilon_{0}})}\leq C\delta_{0}(1+t)^{-1}
\end{equation}
From Chapter 6, we have:
\begin{equation}
 |\hat{T}^{i}|\leq 1, \quad |L^{i}|\leq C
\end{equation}

So $\textbf{E}_{\{0\}}^{Q}$ implies:
\begin{equation}
 \|\hat{T}^{i}(L\psi_{i})\|_{L^{\infty}(\Sigma_{t}^{\epsilon_{0}})}\leq C\delta_{0}(1+t)^{-2}
\end{equation}
Also from $\textbf{E}_{\{1\}}$ and $\textbf{E}_{\{0\}}^{Q}$ we have:
\begin{equation}
 \|Th\|_{L^{\infty}(\Sigma_{t}^{\epsilon_{0}})}\leq C\delta_{0}(1+t)^{-1},\quad
\|Lh\|_{L^{\infty}(\Sigma_{t}^{\epsilon_{0}})}\leq C\delta_{0}(1+t)^{-2}
\end{equation}
It then follows that:
\begin{equation}
 \alpha^{-1}|m|\leq C\delta_{0}(1+t)^{-1}, \quad |e^{\prime}|\leq C\delta_{0}(1+t)^{-2}
\end{equation}
on $W^{s}_{\epsilon_{0}}$.

Writing:
\begin{align}
 L(\kappa-1)=e'(\kappa-1)+\alpha^{-1}m+e'
\end{align}
and noting that for any function $f$:
\begin{align*}
 L|f|\leq |Lf|
\end{align*}
we obtain from (12.9):
\begin{equation}
 L|\kappa-1|\leq C\delta_{0}(1+t)^{-2}|\kappa-1|+C\delta_{0}(1+t)^{-1}
\end{equation}
Integrating this along the integral curves of $L$ and recalling that on $\Sigma_{0}$,
\begin{align*}
 \kappa=1
\end{align*}
the proposition follows. $\qed$

    Recall from Chapter 6:
\begin{equation}
 \lambda_{i}=\bar{g}(\mathring{R}_{i},\hat{T})
\end{equation}

$\textbf{Lemma 12.1}$ Let assumptions $\textbf{E}_{\{1\}}, \textbf{E}_{\{0\}}$ hold on $W^{s}_{\epsilon_{0}}$.
Moreover, let hypothesis $\textbf{H0}$ hold on $W^{s_{0}}_{\epsilon_{0}}$ for some $s_{0}\in(0,s]$. Then, provided that $\delta_{0}$ is suitably small,
we have, for all $t\in[0,s_{0}]$:
\begin{align*}
\max_{i}\|\lambda_{i}\|_{L^{\infty}(\Sigma_{t}^{\epsilon_{0}})}\leq C\delta_{0}[1+\log(1+t)]
\end{align*}
$Proof$. The result follows in the exactly the same way which we used to derive (6.143). Note that to derive an estimate for 
$\kappa^{-1}\zeta$, we must appeal to the hypothesis $\textbf{H0}$. $\qed$

    Since the result of Proposition 12.1 coincides with (6.97), so we conclude that under the assumptions of Proposition 12.1, 
the lower bounds (6.128) and (6.129) holds:
\begin{align}
r\geq1-u+t-C\delta_{0}u[1+\log(1+t)]
\end{align}
and
\begin{align}
r\geq C^{-1}(1+t)
\end{align}
if $\delta_{0}$ is suitably small.

     Recall from Chapter 6 the projection operator $\Sigma$:
\begin{align}
\Sigma_{j}^{i}=\delta_{j}^{i}-r^{-2}x^{i}x^{j}
\end{align}
and the fact that
\begin{align}
|\Sigma\hat{T}|^{2}=r^{-2}\sum_{i}(\lambda_{i})^{2}
\end{align}
Then Lemma 12.1 together with (12.14) implies:
\begin{align}
|\Sigma\hat{T}|\leq C\delta_{0}(1+t)^{-1}[1+\log(1+t)]
\end{align}
in $W^{s_{0}}_{\epsilon_{0}}$.

    Recall also from Chapter 6:
\begin{align}
N=\frac{x^{i}}{r}\frac{\partial}{\partial x^{i}}
\end{align}
and the fact that:
\begin{equation}
 0\leq 1+\bar{g}(N,\hat{T})\leq C\delta_{0}^{2}(1+t)^{-2}[1+\log(1+t)]^{2}
\end{equation}
Moreover, recall from Chapter 6 that
\begin{align*}
 y^{\prime}=\Sigma\hat{T}+(1+\bar{g}(N,\hat{T}))N
\end{align*}

Following the argument in Chapter 6 ((6.169)-(6.177)) we deduce:
\begin{equation}
 |r-1-t+u|\leq C\delta_{0}u[1+\log(1+t)]
\end{equation}
 \begin{equation}
  |\frac{1}{r}-\frac{1}{1-u+t}|\leq C\delta_{0}(1+t)^{-2}[1+\log(1+t)]
 \end{equation}
and:
\begin{equation}
 |y'|\leq C\delta_{0}(1+t)^{-1}[1+\log(1+t)]
\end{equation}
\begin{equation}
 \max_{i}\|y^{i}\|_{L^{\infty}(\Sigma_{t}^{\epsilon_{0}})}\leq C\delta_{0}(1+t)^{-1}[1+\log(1+t)]
\end{equation} 
Note that we have only used the assumptions in Proposition 12.1 and Lemma 12.1 (valid for $t\in[0,s_{0}]$) to get these results.

    Now we begin the argument which leads to our establishing hypothesis $\textbf{H0}$.

$\textbf{Proposition 12.2}$ Consider a surface $S_{t,u}$ for which the following hold:
\begin{align*}
 \inf_{S_{t,u}}r>0  , \quad  \sup_{S_{t,u}}|y'|<1
\end{align*}
Then for every $S_{t,u}$ 1-form $\xi$ we have:
\begin{align*}
 |\xi|^{2}\leq(1-\sup_{S_{t,u}}|y'|^{2})^{-1}(\inf_{S_{t,u}}r)^{-2}\sum_{i=1}^{3}(\xi(R_{i}))^{2}
\end{align*}

$Proof$. Since $R_{i}=\Pi\mathring{R}_{i}$, we have:
\begin{equation}
 (R_{i})^{a}=\Pi^{a}_{b}(\mathring{R}_{i})^{b}=\Pi_{b}^{a}\epsilon_{ikb}x^{k}
\end{equation}
hence:
\begin{align}
 \sum_{i=1}^{3}(R_{i})^{a}(R_{i})^{c}=\sum_{i=1}^{3}\Pi^{a}_{b}\epsilon_{ikb}x^{k}\Pi_{d}^{c}\epsilon_{ild}x^{l}\\\notag
=(\delta_{kl}\delta_{bd}-\delta_{kd}\delta_{bl})x^{k}x^{l}\Pi_{b}^{a}\Pi^{c}_{d}\\\notag
=(r^{2}\delta_{bd}-x^{b}x^{d})\Pi_{b}^{a}\Pi^{c}_{d}
\end{align}
where we have used the fact that:
\begin{equation}
 \sum_{i=1}^{3}\epsilon_{ikb}\epsilon_{ild}=\delta_{kl}\delta_{bd}-\delta_{kd}\delta_{bl}
\end{equation}
It follows that:
\begin{equation}
 \sum_{i=1}^{3}(\xi(R_{i}))^{2}=(\sum_{i=1}^{3}(R_{i})^{a}(R_{i})^{c})\xi_{a}\xi_{c}
=(r^{2}\delta_{bd}-x^{b}x^{d})(\Pi_{b}^{a}\xi_{a})(\Pi_{d}^{c}\xi_{c})
\end{equation}
Since $\xi$ is an $S_{t,u}$ 1-form, we have:
\begin{align*}
 \Pi_{b}^{a}\xi_{a}=\xi_{b}
\end{align*}
We thus obtain:
\begin{equation}
 \sum_{i=1}^{3}(\xi(R_{i}))^{2}=(r^{2}\delta_{bd}-x^{b}x^{d})\xi_{b}\xi_{d}
=r^{2}(|\xi|^{2}-(\xi(N))^{2})
\end{equation}
Recall (6.174):
\begin{equation}
 y^{\prime}=\hat{T}+N
\end{equation}
In view of this and the fact that $\xi(\hat{T})=0$, we have:
\begin{equation}
 \xi(N)=\xi(y^{\prime})=y^{\prime a}\xi_{a}
\end{equation}
hence:
\begin{equation}
 |\xi(N)|=|\xi(y^{\prime})|\leq |y^{\prime}||\xi|
\end{equation}
and:
\begin{equation}
 |\xi|^{2}-(\xi(N))^{2}\geq(1-\sup_{S_{t,u}}|y^{\prime}|^{2})|\xi|^{2}
\end{equation}
Substituting in (12.28) yields:
\begin{equation}
 \sum_{i=1}^{3}(\xi(R_{i}))^{2}\geq (\inf_{S_{t,u}}r)^{2}(1-\sup_{S_{t,u}}|y^{\prime}|^{2})|\xi|^{2}
\end{equation}
The proposition follows. $\qed$

$\textbf{Corollary 12.2.a}$ Suppose that there are positive constants $\epsilon_{1}$ and $\epsilon_{2}<1$ such that:
\begin{align*}
 \inf_{\Sigma^{\epsilon_{0}}_{t}}r\geq\epsilon_{1}(1+t), \quad  \sup_{\Sigma_{t}^{\epsilon_{0}}}|y'|^{2}\leq 1-\epsilon_{2}
\end{align*}
Then $\textbf{H0}$ holds on $\Sigma_{t}^{\epsilon_{0}}$. In fact, for any differentiable function $f$ on $S_{t,u}$, 
we have, pointwise on $\Sigma_{t}^{\epsilon_{0}}$:
\begin{align*}
 |\slashed{d}f|^{2}\leq C(1+t)^{-2}\sum_{i}(R_{i}f)^{2}
\end{align*}
where
\begin{align*}
 C=\frac{1}{\epsilon_{2}\epsilon_{1}^{2}}
\end{align*}
$Proof$. This follows by applying Proposition 12.2 to the $S_{t,u}$ 1-form $\xi=\slashed{d}f$. $\qed$

Since
\begin{align*}
 \sup_{\Sigma_{t}^{\epsilon_{0}}}u=\epsilon_{0}\leq\frac{1}{2}
\end{align*}
fixing any
\begin{align*}
 \epsilon_{1}<\frac{1}{2}
\end{align*}
then (12.13) implies 
\begin{equation}
 \inf_{\Sigma_{t}^{\epsilon_{0}}}r>\epsilon_{1}(1+t)
\end{equation}
for all $t\in[0,s_{0}]$, provided that $\delta_{0}$ is suitably small.
Also, fixing any $\epsilon_{2}<1$, (12.22) implies
\begin{equation}
 \sup_{\Sigma_{t}^{\epsilon_{0}}}|y^{\prime}|^{2}<1-\epsilon_{2}
\end{equation}
for all $t\in[0,s_{0}]$, provided that $\delta_{0}$ is suitably small.

    Here we have used Lemma 12.1, which is based on $\textbf{H0}$ to derive (12.34)-(12.35). We shall show that $\textbf{H0}$ actually holds on $[0,s]$.
Let $\bar{s}_{0}$ be the maximal value of $s_{0}\in[0,s]$ such that $\textbf{H0}$ holds on $W^{s_{0}}_{\epsilon_{0}}$. We show that in fact 
$\bar{s}_{0}=s$. For, suppose that $\bar{s}_{0}<s$. Since (12.17) (12.19) and (12.22) hold on $W_{\epsilon_{0}}^{\bar{s}_{0}}$, hence also (12.34) and (12.35) hold for all 
$t\in[0,\bar{s}_{0}]$, in particular at $t=\bar{s}_{0}$. However, (12.34) and (12.35) are strict inequalities, by continuity they will hold on some small interval 
$[\bar{s}_{0},s_{*})\subset[\bar{s}_{0},s]$. This contradicts the maximality of $\bar{s}_{0}$. So $\textbf{H0}$ holds on $W^{s}_{\epsilon_{0}}$.
From Lemma 12.1 with $s_{0}=s$, we get the following proposition:

$\textbf{Proposition 12.3}$ Let $\textbf{E}_{\{1\}}, \textbf{E}^{Q}_{\{0\}}$, hold on $W^{s}_{\epsilon_{0}}$.
Then if $\delta_{0}$ is suitably small, $\textbf{H0}$ holds on $W^{s}_{\epsilon_{0}}$. Moreover, we have, on $W^{s}_{\epsilon_{0}}$:
\begin{align*}
 1+t-u-C\delta_{0}u[1+\log(1+t)]\leq r\leq 1+t+\min\{0,-u+C\delta_{0}u[1+\log(1+t)]\}\\
|y'|\leq C\delta_{0}(1+t)^{-1}[1+\log(1+t)]
\end{align*}
and, for all $t\in[0,s]$,
\begin{align*}
 \|\kappa-1\|_{L^{\infty}(\Sigma_{t}^{\epsilon_{0}})}\leq C\delta_{0}[1+\log(1+t)]\\
\max_{i}\|\lambda_{i}\|_{L^{\infty}(\Sigma_{t}^{\epsilon_{0}})}\leq C\delta_{0}[1+\log(1+t)]\\
\max_{i}\|y^{i}\|_{L^{\infty}(\Sigma_{t}^{\epsilon_{0}})}\leq C\delta_{0}(1+t)^{-1}[1+\log(1+t)]
\end{align*}

\section{The Coercivity Hypothesis $\textbf{H1}$, $\textbf{H2}$ and $\textbf{H2}^{\prime}$. Estimates for $\chi^{\prime}$}
     We start this section with a proposition analogous to Proposition 12.2.

$\textbf{Proposition 12.4}$ Consider a surface $S_{t,u}$ for which the following hold:
\begin{align*}
 \inf_{S_{t,u}}r>0,  \quad \sup_{S_{t,u}}|y'|<1
\end{align*}
Then for every symmetric 2-covariant $S_{t,u}$ tensorfield $\vartheta$ we have:
\begin{align*}
 |\vartheta|^{2}\leq (1-\sup_{S_{t,u}}|y'|^{2})^{-2}(\inf_{S_{t,u}}r)^{-4}\sum_{i,j=1}^{3}(\vartheta(R_{i},R_{j}))^{2}
\end{align*}
pointwise on $S_{t,u}$

$Proof$. We have:
\begin{align}
 \sum_{i,j}(\vartheta(R_{i},R_{j}))^{2}=(\sum_{i=1}^{3}(R_{i})^{a}(R_{i})^{c})(\sum_{j=1}^{3}(R_{j})^{b}(R_{j})^{d})\vartheta_{ab}\vartheta_{cd}
\end{align}
We substitute (12.25) in the first two factors on the right. Since $x^{b}=rN^{b}=r(y^{\prime b}-\hat{T}^{b})$, and $\hat{T}^{b}\Pi^{a}_{b}=0$, we have:
\begin{align}
 x^{b}\Pi_{b}^{a}=ry^{\prime b}\Pi_{b}^{a}
\end{align}
(12.25) becomes:
\begin{equation}
 \sum_{i=1}^{3}(R_{i})^{a}(R_{i})^{c}=r^{2}(\delta_{bd}-y^{\prime b}y^{\prime d})\Pi^{a}_{b}\Pi_{d}^{c}
\end{equation}
Substituting (12.38) in (12.36) we get:
\begin{align}
 \sum_{i,j=1}^{3}(\vartheta(R_{i},R_{j}))^{2}=r^{4}(\delta_{ac}-y^{\prime a}y^{\prime c})(\delta_{bd}-y^{\prime b}y^{\prime d})\vartheta_{ab}\vartheta_{cd}
=r^{4}(|\vartheta|^{2}-2|i_{y^{\prime}}\vartheta|^{2}+(\vartheta(y^{\prime},y^{\prime}))^{2})
\end{align}
where
\begin{equation}
 (i_{y'}\vartheta)_{b}=y'^{a}\vartheta_{ab}
\end{equation}
Consider the expression:
\begin{equation}
 |\vartheta|^{2}-2|i_{y^{\prime}}\vartheta|^{2}+(\vartheta(y^{\prime},y^{\prime}))^{2}\notag
\end{equation}
At a given point, by a suitable rotation of the rectangular coordinates we can arrange that $y^{\prime 2}=y^{\prime 3}=0$. Then this expression becomes:
\begin{align*}
 (1-(y^{\prime 1})^{2})^{2}(\vartheta_{11})^{2}+(\vartheta_{22})^{2}+(\vartheta_{33})^{2}
+2(1-(y^{\prime 1})^{2})((\vartheta_{12})^{2}+(\vartheta_{13})^{2})+2(\vartheta_{23})^{2}\\
\geq(1-(y^{\prime 1})^{2})^{2}((\vartheta_{11})^{2}+(\vartheta_{22})^{2}+(\vartheta_{33})^{2}+2(\vartheta_{12})^{2}
+2(\vartheta_{13})^{2}+2(\vartheta_{23})^{2})=(1-|y^{\prime}|^{2})^{2}|\vartheta|^{2}
\end{align*}
We conclude that:
\begin{equation}
 |\vartheta|^{2}-2|i_{y^{\prime}}\vartheta|^{2}+(\vartheta(y^{\prime},y^{\prime}))^{2}\geq (1-|y^{\prime}|^{2})^{2}|\vartheta|^{2}
\end{equation}
So (12.39) implies:
\begin{equation}
 \sum_{i,j=1}^{3}(\vartheta(R_{i},R_{j}))^{2}\geq r^{4}(1-|y^{\prime}|^{2})^{2}|\vartheta|^{2}
\end{equation}
The proposition follows. $\qed$

    We now return to the propagation equation for $\chi$ in Chapter 3:
\begin{equation}
 L\chi_{AB}=e\chi_{AB}+\chi_{A}^{C}\chi_{BC}-\alpha'_{AB}
\end{equation}
 by direct calculation (see Chapter 4):
\begin{align}
 \alpha'_{AB}={\alpha'}^{[P]}_{AB}+{\alpha}^{[N]}_{AB}
\end{align}
where
\begin{align}
{\alpha'}^{[P]}_{AB}=-\frac{1}{2}\frac{dH}{dh}\slashed{D}^{2}h(X_{A},X_{B})
\end{align}
and $\alpha^{[N]}_{AB}$ is given by (4.25) and (4.26).

    Taking into account of the fact that
\begin{align}
 \slashed{\mathcal{L}}_{L}\slashed{g}=2\chi\\\notag
 \chi'=\chi-\frac{\slashed{g}}{1-u+t}
\end{align}
we obtain the following propagation equation for $\chi^{\prime}$:
\begin{equation}
 L\chi'_{AB}=e\chi'_{AB}+(\slashed{g}^{-1})^{CD}\chi'_{AD}\chi'_{BC}+\frac{e\slashed{g}_{AB}}{1-u+t}-\alpha'_{AB}
\end{equation}

$\textbf{Lemma 12.2}$ Let assumptions $\textbf{E}_{\{2\}},\textbf{E}^{Q}_{\{1\}},\textbf{E}^{QQ}_{\{0\}}$, hold on $W^{s}_{\epsilon_{0}}$, 
and let the initial data satisfy:
\begin{align*}
 \|\chi'\|_{L^{\infty}(\Sigma_{0}^{\epsilon_{0}})}\leq C\delta_{0}
\end{align*}
Moreover, let $\textbf{H1}$ hold on $W^{s_{0}}_{\epsilon_{0}}$ for some $s_{0}\in(0,s]$. Then, provided $\delta_{0}$ is suitably small,
we have, for all $t\in[0,s_{0}]$:
\begin{align*}
 \|\chi'\|_{L^{\infty}(\Sigma_{t}^{\epsilon_{0}})}\leq C\delta_{0}(1+t)^{-2}[1+\log(1+t)]
\end{align*}
$Proof$. From the assumptions and the expression for $e$, we deduce:
\begin{equation}
 \|e\|_{L^{\infty}(\Sigma_{t}^{\epsilon_{0}})}\leq C\delta_{0}(1+t)^{-2}
\end{equation}
and from the hypotheses $\textbf{H0}$ and $\textbf{H1}$ and the expression for $\alpha'_{AB}$, we deduce:
\begin{equation}
 \|\alpha'\|_{L^{\infty}(\Sigma_{t}^{\epsilon_{0}})}\leq C\delta_{0}(1+t)^{-3}
\end{equation}
Here we have used the fact that the magnitude of the symmetric 2-covariant $S_{t,u}$ tensorfield $\slashed{\omega}$ is bounded by $C\delta_{0}(1+t)^{-2}$
In fact, since 
\begin{align*}
 \slashed{\omega}_{AB}=X_{A}^{i}(X_{B}\psi_{i})
\end{align*}
we have
\begin{align*}
 |\slashed{\omega}|^{2}=(\slashed{g}^{-1})^{AC}(\slashed{g}^{-1})^{BD}X_{A}^{i}(X_{B}\psi_{i})X_{C}^{i}(X_{D}\psi_{i})\\
=(g^{-1})^{ij}(X_{B}\psi_{i})(X_{D}\psi_{j})(\slashed{g}^{-1})^{BD}-\hat{T}^{i}\hat{T}^{j}(X_{B}\psi_{i})(X_{D}\psi_{j})(\slashed{g}^{-1})^{BD}\\
=\sum_{i}|\slashed{d}\psi_{i}|^{2}-|\slashed{\omega}_{\hat{T}}|^{2}
\leq \sum_{i}|\slashed{d}\psi_{i}|^{2}\leq C\delta_{0}^{2}(1+t)^{-4}
\end{align*}
Also, we have used the fact that by $\textbf{H1}$, 
\begin{align*}
 |\slashed{D}^{2}h|\leq C(1+t)^{-1}\max_{i}|\slashed{\mathcal{L}}_{R_{i}}\slashed{d}h|\quad:\textrm{on}\quad W^{s_{0}}_{\epsilon_{0}}
\end{align*}
hence since $\slashed{\mathcal{L}}_{R_{i}}\slashed{d}h=\slashed{d}R_{i}h$, we have, by $\textbf{H0}$,
\begin{align*}
 |\slashed{\mathcal{L}}_{R_{i}}\slashed{d}h|\leq C(1+t)^{-1}\max_{j}|R_{j}R_{i}h|\quad:\textrm{on}\quad W^{s}_{\epsilon_{0}}
\end{align*}
Thus,
\begin{align*}
 |\slashed{D}^{2}h|\leq C(1+t)^{-2}\max_{i,j}|R_{i}R_{j}h|\quad:\textrm{on}\quad W^{s_{0}}_{\epsilon_{0}}
\end{align*}
Therefore, by $\textbf{E}_{\{2\}}$:
\begin{align*}
 \|\slashed{D}^{2}h\|_{L^{\infty}(\Sigma_{t}^{\epsilon_{0}})}\leq C\delta_{0}(1+t)^{-3}\quad:\textrm{for all}\quad t\in[0,s_{0}]
\end{align*}

We can write the equation (12.47) in the form:
\begin{equation}
 \slashed{\mathcal{L}}_{L}\chi'=\chi'\cdot\chi'+e\chi'+\frac{\slashed{g}e}{1-u+t}-\alpha'
\end{equation}
where
\begin{align*}
 (\chi'\cdot\chi')_{AB}=\chi'_{AC}\chi'_{BD}(\slashed{g}^{-1})^{CD}
\end{align*}

     Let now $\vartheta, \vartheta'$, be a pair of symmetric 2-covariant $S_{t,u}$ tensorfields. Their pointwise inner product, with respect to the induced acoustical
metric $\slashed{g}$, is given by:
\begin{equation}
 (\vartheta,\vartheta')=(\slashed{g}^{-1})^{AC}(\slashed{g}^{-1})^{BD}\vartheta_{AB}\vartheta'_{CD}
\end{equation}
Since
\begin{align*}
 \slashed{\mathcal{L}}_{L}(\slashed{g}^{-1})^{AB}=-2\chi^{AB}
\end{align*}
we have:
\begin{equation}
 L(\vartheta,\vartheta')=(\slashed{\mathcal{L}}_{L}\vartheta,\vartheta')+(\vartheta,\slashed{\mathcal{L}}_{L}\vartheta')
-4\textrm{tr}(\vartheta\cdot\chi\cdot\vartheta')
\end{equation}
Taking in particular $\vartheta'=\vartheta$ we obtain:
\begin{equation}
 |\vartheta|L|\vartheta|=\frac{1}{2}L(\vartheta,\vartheta)=(\vartheta,\slashed{\mathcal{L}}_{L}\vartheta)
-2\textrm{tr}(\vartheta\cdot\chi\cdot\vartheta)
\end{equation}
Writing 
\begin{align*}
 \chi=\frac{\slashed{g}}{1-u+t}+\chi'
\end{align*}
we have:
\begin{align*}
 \textrm{tr}(\vartheta\cdot\chi\cdot\vartheta)=\frac{|\vartheta|^{2}}{1-u+t}+\textrm{tr}(\vartheta\cdot\chi'\cdot\vartheta)
\end{align*}
and (12.53) can be written as:
\begin{equation}
 |\vartheta|L|\vartheta|=(\vartheta,\slashed{\mathcal{L}}_{L}\vartheta)-\frac{2|\vartheta|^{2}}{1-u+t}-
2\textrm{tr}(\vartheta\cdot\chi'\cdot\vartheta)
\end{equation}
Now, the following inequality holds:
\begin{equation}
 |\textrm{tr}(\vartheta\cdot\chi'\cdot\vartheta)|\leq |\chi'||\vartheta|^{2}
\end{equation}
To see this, we can work in an orthonormal frame on $S_{t,u}$ relative to $\slashed{g}$ which is a frame of eigenvectors
of $\chi^{\prime}$. Let $\lambda_{1},\lambda_{2}$ be the eigenvalues of $\chi'$,
then at a given point,
\begin{align*}
 \textrm{tr}(\vartheta\cdot\chi'\cdot\vartheta)=\lambda_{1}(\vartheta_{11})^{2}+
(\lambda_{1}+\lambda_{2})(\vartheta_{12})^{2}+\lambda_{2}(\vartheta_{22})^{2}
\end{align*}
Since
\begin{align*}
 |\vartheta|^{2}=(\vartheta_{11})^{2}+2(\vartheta_{12})^{2}+(\vartheta_{22})^{2}
\end{align*}
It follows that:
\begin{align*}
 |\textrm{tr}(\vartheta\cdot\chi'\cdot\vartheta)|\leq\max\{|\lambda_{1}|,|\lambda_{2}|\}|\vartheta|^{2}
\end{align*}
In view of the fact that:
\begin{align*}
 \max\{|\lambda_{1}|,|\lambda_{2}|\}\leq \sqrt{(\lambda_{1})^{2}+(\lambda_{2})^{2}}=|\chi'|
\end{align*}
then (12.55) follows.

    Taking into account (12.55) as well as the fact that:
\begin{align*}
 |(\vartheta,\slashed{\mathcal{L}}_{L}\vartheta)|\leq |\vartheta||\slashed{\mathcal{L}}_{L}\vartheta|
\end{align*}
we deduce from (12.54) the inequality:
\begin{equation}
 L((1-u+t)^{2}|\vartheta|)\leq (1-u+t)^{2}(2|\chi'||\vartheta|+|\slashed{\mathcal{L}}_{L}\vartheta|)
\end{equation}

     We apply this to the case $\vartheta=\chi'$. In view of the fact that:
\begin{align*}
 |\chi'\cdot\chi'|=\sqrt{(\lambda_{1})^{4}+(\lambda_{2})^{4}}\leq (\lambda_{1})^{2}+(\lambda_{2})^{2}=|\chi'|^{2}
\end{align*}
we obtain, substituting for $\slashed{\mathcal{L}}_{L}\chi^{\prime}$ from (12.50),
\begin{equation}
 L((1-u+t)^{2}|\chi'|)\leq (1-u+t)^{2}((3|\chi'|+|e|)|\chi'|+|b|)
\end{equation}
where
\begin{align*}
 b=\frac{\slashed{g}e}{1-u+t}-\alpha'
\end{align*}

    Let $\mathcal{P}(t)$ be the property:
\begin{align*}
 \mathcal{P}(t): \|\chi'\|_{L^{\infty}(\Sigma_{t}^{\epsilon_{0}})}\leq C_{0}\delta_{0}(1+t')^{-2}[1+\log(1+t')]
\end{align*}
for all $t'\in[0,t]$.

Choosing $C_{0}$ suitably large, we have, by the assumption on the initial data,
\begin{equation}
 \|\chi'\|_{L^{\infty}(\Sigma_{0}^{\epsilon_{0}})}<C_{0}\delta_{0}
\end{equation}
It follows by continuity that $\mathcal{P}(t)$ is true for sufficiently small positive $t$. Let 
$t_{0}$ be the least upper bound of the set of values of $t\in[0,s_{0}]$ for which $\mathcal{P}(t)$ holds. Then by continuity
$\mathcal{P}(t_{0})$ is true. Hence, in $W^{t_{0}}_{\epsilon_{0}}$, we have, in view of (12.48):
\begin{equation}
 3|\chi'|+|e|\leq (3C_{0}+C)\delta_{0}(1+t)^{-2}[1+\log(1+t)]
\end{equation}
in $W^{t_{0}}_{\epsilon_{0}}$.

Substituting in (12.57) and taking into account the estimates (12.48) and (12.49) we obtain:
\begin{equation}
 L((1-u+t)^{2}|\chi'|)\leq (3C_{0}+C)\delta_{0}(1+t)^{-2}[1+\log(1+t)]((1-u+t)^{2}|\chi'|)+C\delta_{0}(1+t)^{-1}
\end{equation}
Setting along an integral curve of $L$:
\begin{equation}
 x(t)=(1-u+t)^{2}|\chi'|
\end{equation}
then (12.60) takes the form:
\begin{equation}
 \frac{dx}{dt}\leq fx+g
\end{equation}
where
\begin{equation}
 f(t)=(3C_{0}+C)\delta_{0}(1+t)^{-2}[1+\log(1+t)],\quad g(t)=C\delta_{0}(1+t)^{-1}
\end{equation}
Integrating from $t=0$ yields:
\begin{equation}
 x(t)\leq e^{\int_{0}^{t}f(t')dt'}\{x(0)+\int_{0}^{t}e^{-\int_{0}^{t'}f(t'')dt''}g(t')dt'\}
\end{equation}
Since
\begin{equation}
 x(0)\leq\|\chi'\|_{L^{\infty}(\Sigma_{0}^{\epsilon_{0}})},
\end{equation}
while
\begin{equation}
 x(t)\geq\frac{1}{4}(1+t)^{2}|\chi'|,
\end{equation}
taking into account the facts that:
\begin{equation}
 \int_{0}^{t}f(t')dt'\leq\int_{0}^{\infty}f(t')dt'\leq 2(3C_{0}+C)\delta_{0},\quad
\int_{0}^{t}g(t')dt'=C\delta_{0}\log(1+t)
\end{equation}
(12.64) yields the bound:
\begin{equation}
 \frac{1}{4}(1+t)^{2}|\chi'|\leq e^{2(3C_{0}+C)\delta_{0}}\{\|\chi'\|_{L^{\infty}(\Sigma_{0}^{\epsilon_{0}})}
+C\delta_{0}\log(1+t)\}
\end{equation}
So we get on $\Sigma_{t}^{\epsilon_{0}}$.
\begin{equation}
 \|\chi'\|_{L^{\infty}(\Sigma_{t}^{\epsilon_{0}})}\leq 4e^{2(3C_{0}+C)\delta_{0}}\{\|\chi'\|_{L^{\infty}(\Sigma_{0}^{\epsilon_{0}})}
+C\delta_{0}\log(1+t)\}(1+t)^{-2}
\end{equation}
This holds for all $t\in[0,t_{0}]$. Let us now fix $C_{0}$ large enough so that:
\begin{equation}
 C_{0}\delta_{0}\geq8\|\chi'\|_{L^{\infty}(\Sigma_{0}^{\epsilon_{0}})}, \quad C_{0}>8C
\end{equation}
Then provided $\delta_{0}$ satisfies:
\begin{equation}
 \delta_{0}\leq\frac{\log2}{2(3C_{0}+C)}
\end{equation}
(12.69) implies:
\begin{equation}
 \|\chi'\|_{L^{\infty}(\Sigma_{t}^{\epsilon_{0}})}<C_{0}\delta_{0}(1+t)^{-2}[1+\log(1+t)]
\end{equation}
for all $t\in[0,t_{0}]$.

It follows by continuity that $\mathcal{P}(t)$ is true for some $t>t_{0}$. So we must have $t_{0}=s_{0}$, and the lemma follows. $\qed$

    Next, we investigate $\textbf{H1}$. Consider $\slashed{\mathcal{L}}_{R_{i}}\xi$ for an arbitrary $S_{t,u}$ 1-form $\xi$.
We have, in rectangular coordinates:
\begin{equation}
 (\slashed{\mathcal{L}}_{R_{i}}\xi)_{a}=(\slashed{D}_{R_{i}}\xi)_{a}+\xi_{m}(\slashed{D}R_{i})^{m}_{a},\quad
(\slashed{D}_{R_{i}}\xi)_{a}=(R_{i})^{m}(\slashed{D}\xi)_{ma}
\end{equation}
where we have used that for any $S_{t,u}$ vectorfield $X$:
\begin{align*}
 (\slashed{\mathcal{L}}_{R_{i}}\xi)(X)=(\slashed{D}_{R_{i}}\xi)(X)+\xi(\slashed{D}_{X}R_{i})
\end{align*}
Hence, we have:
\begin{align}
 \sum_{i}|\slashed{\mathcal{L}}_{R_{i}}\xi|^{2}=\sum_{i}|\slashed{D}_{R_{i}}\xi|^{2}\\\notag
+2\sum_{i}(\bar{g}^{-1})^{ab}(\slashed{D}_{R_{i}}\xi)_{a}\xi_{n}(\slashed{D}R_{i})^{n}_{b}\\\notag
+\sum_{i}(\bar{g}^{-1})^{ab}\xi_{m}\xi_{n}(\slashed{D}R_{i})^{m}_{a}(\slashed{D}R_{i})^{n}_{b}
\end{align}
Consider the first term on the right of (12.74). We have:
\begin{align}
 \sum_{i}|\slashed{D}_{R_{i}}\xi|^{2}=\sum_{i}(\bar{g}^{-1})^{ab}(\slashed{D}_{R_{i}}\xi)_{a}(\slashed{D}_{R_{i}}\xi)_{b}
=(\bar{g}^{-1})^{ab}(\sum_{i}(R_{i})^{m}(R_{i})^{n})(\slashed{D}\xi)_{ma}(\slashed{D}\xi)_{nb}
\end{align}
Substituting from (12.38) and noting that:
\begin{align*}
 \Pi_{c}^{m}\Pi^{n}_{d}(\slashed{D}\xi)_{ma}(\slashed{D}\xi)_{nb}=(\slashed{D}\xi)_{ca}(\slashed{D}\xi)_{db}
\end{align*}
we obtain:
\begin{equation}
 \sum_{i}|\slashed{D}_{R_{i}}\xi|^{2}=r^{2}(\bar{g}^{-1})^{ab}(\delta_{cd}-y'^{c}y'^{d})(\slashed{D}\xi)_{ca}(\slashed{D}\xi)_{db}
\end{equation}
Consider the symmetric 2-covariant $S_{t,u}$ tensorfield $\vartheta$:
\begin{equation}
 \vartheta_{cd}=(\bar{g}^{-1})^{ab}(\slashed{D}\xi)_{ca}(\slashed{D}\xi)_{db}
\end{equation}
Let $\mu_{1},\mu_{2},\mu_{3}$, be the eigenvalues of $\vartheta$ relative to $\bar{g}$. Since $\vartheta$ is positive semi-definite, we have:
\begin{equation}
 \mu_{i}\geq 0: i=1,2,3
\end{equation}
It follows that:
\begin{equation}
 \max_{i}\mu_{i}\leq\sum_{i}\mu_{i}=\textrm{tr}\vartheta
\end{equation}
where
\begin{equation}
 \textrm{tr}\vartheta=(\bar{g}^{-1})^{cd}\vartheta_{cd}=\delta_{cd}\vartheta_{cd}
\end{equation}
Then (12.79) implies:
\begin{equation}
 y'^{c}y'^{d}\vartheta_{cd}\leq |y'|^{2}\textrm{tr}\vartheta
\end{equation}
We then have:
\begin{equation}
 (\bar{g}^{-1})^{ab}(\delta_{cd}-y'^{c}y'^{d})(\slashed{D}\xi)_{ca}(\slashed{D}\xi)_{db}
=(\delta_{cd}-y'^{c}y'^{d})\vartheta_{cd}\geq(1-|y'|^{2})\textrm{tr}\vartheta
\end{equation}
Thus, since
\begin{equation}
 \textrm{tr}\vartheta=(\bar{g}^{-1})^{cd}(\bar{g}^{-1})^{ab}(\slashed{D}\xi)_{ca}(\slashed{D}\xi)_{db}=|\slashed{D}\xi|^{2}
\end{equation}
(12.76) implies:
\begin{equation}
 \sum_{i}|\slashed{D}_{R_{i}}\xi|^{2}\geq r^{2}(1-|y'|^{2})|\slashed{D}\xi|^{2}
\end{equation}

     Consider next the third term on the right of (12.74). We should first obtain a suitable expression for $\slashed{D}R_{i}$.
Recall from (6.45) and (6.56), we have:
\begin{equation}
 g(\slashed{D}_{X_{A}}R_{i},X_{B})=X^{l}_{A}\epsilon_{ilm}X^{m}_{B}+\lambda_{i}(\eta^{-1}\chi_{AB}-\slashed{k}_{AB})
\end{equation}
Let us denote:
\begin{equation}
 (\nu_{i})_{AB}=X^{l}_{A}\epsilon_{ilm}X^{m}_{B}
\end{equation}
 The rectangular components of $\nu_{i}$ are given by:
\begin{equation}
 (\nu_{i})_{lm}=\Pi^{l'}_{l}\Pi^{m'}_{m}\epsilon_{il'm'}
\end{equation}
Let us also denote:
\begin{equation}
 \tau_{i}=\lambda_{i}(\eta^{-1}\chi-\slashed{k})
\end{equation}
Then by (12.85) the rectangular components of $\slashed{D}R_{i}$ are given by:
\begin{equation}
 (\slashed{D}R_{i})^{k}_{l}=\Pi^{k}_{m}\Pi^{l'}_{l}\epsilon_{il'm}+(\tau_{i})^{k}_{l}
\end{equation}
Using (12.89) we obtain, in view of (12.26):
\begin{align}
 \sum_{i}(\slashed{D}R_{i})^{m}_{a}(\slashed{D}R_{i})^{n}_{b}=
\Pi_{m'}^{m}\Pi_{n'}^{n}\Pi^{a'}_{a}\Pi^{b'}_{b}(\delta_{a'b'}\delta_{m'n'}-\delta_{a'n'}\delta_{b'm'})\\\notag
+\sum_{i}(\tau_{i})^{m}_{a}\Pi^{n}_{n'}\Pi^{b'}_{b}\epsilon_{ib'n'}+\sum_{i}(\tau_{i})^{n}_{b}\Pi^{m}_{m'}\Pi^{a'}_{a}\epsilon_{ia'm'}\\\notag
+\sum_{i}(\tau_{i})^{m}_{a}(\tau_{i})^{n}_{b}
\end{align}
Substituting this in the third term on the right in (12.74) yields:
\begin{align}
 \sum_{i}(\bar{g}^{-1})^{ab}\xi_{m}\xi_{n}(\slashed{D}R_{i})_{a}^{m}(\slashed{D}R_{i})^{n}_{b}=
\delta_{a'b'}(\bar{g}^{-1})^{ab}\Pi^{a'}_{a}\Pi_{b}^{b'}|\xi|^{2}-|\xi|^{2}\\\notag
+2\sum_{i}(\bar{g}^{-1})^{ab}(\tau_{i}\cdot\xi)_{a}(v_{i})_{b}+\sum_{i}|\tau_{i}\cdot\xi|^{2}
\end{align}
where $v_{i}$ is the $S_{t,u}$ 1-form:
\begin{equation}
 (v_{i})_{b}=\Pi^{b'}_{b}\epsilon_{ib'n'}\xi_{n'}
\end{equation}
Using (12.26) and (12.92) we obtain:
\begin{align}
 \sum_{i}|v_{i}|^{2}=\sum_{i}(\bar{g}^{-1})^{ab}\Pi^{a'}_{a}\Pi^{b'}_{b}\epsilon_{ia'm'}\epsilon_{ib'n'}\xi_{m'}\xi_{n'}\\\notag
=(\bar{g}^{-1})^{ab}\Pi^{a'}_{a}\Pi^{b'}_{b}(\delta_{a'b'}\delta_{m'n'}-\delta_{a'n'}\delta_{b'm'})\xi_{m'}\xi_{n'}\\\notag
=\delta_{a'b'}(\bar{g}^{-1})^{ab}\Pi^{a'}_{a}\Pi^{b'}_{b}|\xi|^{2}-|\xi|^{2}
\end{align}
Since
\begin{equation}
 (\bar{g}^{-1})^{ab}\Pi^{a'}_{a}\Pi^{b'}_{b}=(\bar{g}^{-1})^{a'b'}-\hat{T}^{a'}\hat{T}^{b'}
\end{equation}
we have:
\begin{equation}
 \delta_{a'b'}(\bar{g}^{-1})^{ab}\Pi^{a'}_{a}\Pi^{b'}_{b}=2
\end{equation}
Hence (12.93) reduces to
\begin{align*}
 \sum_{i}|v_{i}|^{2}=|\xi|^{2}
\end{align*}
while
\begin{align*}
 2\sum_{i}(\bar{g}^{-1})^{ab}(\tau_{i}\cdot\xi)_{a}(v_{i})_{b}\geq -2\sum_{i}|\tau_{i}\cdot\xi||v_{i}|
\geq -2|\xi|\sqrt{\sum_{i}|\tau_{i}|^{2}}\sqrt{\sum_{i}|v_{i}|^{2}}
\end{align*}
So from (12.91) we have:
\begin{align}
 \sum_{i}(\bar{g}^{-1})^{ab}\xi_{m}\xi_{n}(\slashed{D}R_{i})^{m}_{a}(\slashed{D}R_{i})^{n}_{b}=
\sum_{i}|v_{i}|^{2}+2\sum_{i}(\tau_{i}\cdot\xi)\cdot v_{i}+\sum_{i}|\tau_{i}\cdot\xi|^{2}\\\notag
\geq\sum_{i}|v_{i}|^{2}-2\sum_{i}|\tau_{i}\cdot\xi||v_{i}|\\\notag
\geq\sum_{i}|v_{i}|^{2}-2\sqrt{\sum_{i}|\tau_{i}\cdot\xi|^{2}}\sqrt{\sum_{i}|v_{i}|^{2}}
\end{align}
i.e.
\begin{equation}
 \sum_{i}(\bar{g}^{-1})^{ab}\xi_{m}\xi_{n}(\slashed{D}R_{i})^{m}_{a}(\slashed{D}R_{i})^{n}_{b}
\geq |\xi|^{2}(1-2\sqrt{\sum_{i}(\tau_{i})^{2}})
\end{equation}

    Finally, we consider the second term on the right of (12.74):
\begin{equation}
 2\sum_{i}(\bar{g}^{-1})^{ab}(\slashed{D}_{R_{i}}\xi)_{a}\xi_{n}(\slashed{D}R_{i})^{n}_{b}=
2\sum_{i}(\bar{g}^{-1})^{ab}(R_{i})^{m}(\slashed{D}\xi)_{ma}\xi_{n}(\slashed{D}R_{i})^{n}_{b}
\end{equation}
From (12.24), (12.26) and (12.89) we have:
\begin{align}
 \sum_{i}(R_{i})^{m}(\slashed{D}R_{i})^{n}_{b}=\sum_{i}\Pi^{m}_{m'}\epsilon_{ikm'}x^{k}
(\Pi^{n}_{n'}\Pi^{b'}_{b}\epsilon_{ib'n'}+(\tau_{i})^{n}_{b})\\\notag
=(\delta_{kb'}\delta_{m'n'}-\delta_{kn'}\delta_{m'b'})\Pi^{m}_{m'}\Pi^{n}_{n'}\Pi^{b'}_{b}x^{k}
+\Pi^{m}_{m'}x^{k}\sum_{i}\epsilon_{ikm'}(\tau_{i})^{n}_{b}\\\notag
=\Pi^{m}_{n'}\Pi^{n}_{n'}(\Pi^{b'}_{b}x^{b'})-\Pi^{m}_{b}(\Pi^{n}_{n'}x^{n'})+\Pi^{m}_{m'}x^{k}\sum_{i}\epsilon_{ikm'}(\tau_{i})^{n}_{b}
\end{align}
Consider the first term on the right of (12.99). We have:
\begin{align*}
 \Pi^{b'}_{b}x^{b'}=r\Pi^{b'}_{b}N^{b'}=r\bar{g}_{b'c}\Pi^{b'}_{b}N^{c}
\end{align*}
By (12.29) and the fact that $\bar{g}_{b'c}\Pi^{b'}_{b}\hat{T}^{c}=0$, we have:
\begin{align*}
 \bar{g}_{b'c}\Pi^{b'}_{b}N^{c}=\bar{g}_{b'c}\Pi^{b'}_{b}y'^{c}
\end{align*}
Hence, we obtain:
\begin{equation}
 \Pi^{b'}_{b}x^{b'}=r\Pi^{b'}_{b}(\bar{g}_{b'c}y'^{c})
\end{equation}
Introducing the $S_{t,u}$-tangential vectorfield:
\begin{equation}
 \slashed{y}'=\Pi\cdot y'
\end{equation}
then (12.100) takes the form:
\begin{equation}
 \Pi^{b'}_{b}x^{b'}=r\slashed{y}'^{b}
\end{equation}
Similarly, we have:
\begin{equation}
 \Pi^{n}_{n'}x^{n'}=r\slashed{y}'^{n}
\end{equation}
In the above, we have used the fact that $\bar{g}_{ab}=\delta_{ab}$.

In view of (12.102), (12.103), (12.99) reduces to:
\begin{align}
 \sum_{i}(R_{i})^{m}(\slashed{D}R_{i})^{n}_{b}=r
\{\Pi^{m}_{n'}\Pi^{n}_{n'}\slashed{y}'^{b}-\Pi^{m}_{b}\slashed{y}'^{n}+\Pi^{m}_{m'}
N^{k}\sum_{i}\epsilon_{ikm'}(\tau_{i})^{n}_{b}\}
\end{align}
Hence, from (12.98), the second term on the right of (12.74) is given by:
\begin{equation}
 2r\{\zeta_{n'}\xi_{n'}-\textrm{tr}(\slashed{D}\xi)(\xi\cdot\slashed{y}')+\sum_{i}\eta_{i}\cdot(\tau_{i}\cdot\xi)\}
\end{equation}
where
\begin{equation}
 \zeta_{n'}=(\slashed{D}\xi)_{n'a}\slashed{y}'^{a}
\end{equation}
and
\begin{equation}
 (\eta_{i})_{a}=N^{k}\epsilon_{ikm'}(\slashed{D}\xi)_{m'a}
\end{equation}
Now we have:
\begin{equation}
|\zeta_{n'}\xi_{n'}|\leq|\zeta||\xi|\leq|\slashed{D}\xi||\xi||\slashed{y}'|  
\end{equation}
Since for any 2-covariant $S_{t,u}$ tensorfield $M$ we have:
\begin{align*}
 \frac{1}{2}(\textrm{tr}M)^{2}\leq |M|^{2}
\end{align*}
we have:
\begin{equation}
 |\textrm{tr}(\slashed{D}\xi)(\xi\cdot\slashed{y}')|\leq \sqrt{2}|\slashed{D}\xi||\xi||\slashed{y}'|
\end{equation}
Finally,
\begin{align}
 |\sum_{i}\eta_{i}\cdot(\tau_{i}\cdot\xi)|\leq\sum_{i}|\eta_{i}||\tau_{i}||\xi|
\leq|\xi|\sqrt{\sum_{i}|\tau_{i}|^{2}}\sqrt{\sum_{i}|\eta_{i}|^{2}}
\end{align}
Using (12.26) we obtain:
\begin{align}
 \sum_{i}|\eta_{i}|^{2}=\sum_{i}(\bar{g}^{-1})^{ab}N^{k}N^{l}\epsilon_{ikm'}\epsilon_{iln'}(\slashed{D}\xi)_{m'a}(\slashed{D}\xi)_{n'b}\\\notag
=(\bar{g}^{-1})^{ab}N^{k}N^{l}(\delta_{kl}\delta_{m'n'}-\delta_{kn'}\delta_{lm'})(\slashed{D}\xi)_{m'a}(\slashed{D}\xi)_{n'b}\\\notag
=(\bar{g}^{-1})^{ab}((\slashed{D}\xi)_{n'a}(\slashed{D}\xi)_{n'b}-N^{m'}N^{n'}(\slashed{D}\xi)_{m'a}(\slashed{D}\xi)_{n'b})\\\notag
\leq (\bar{g}^{-1})^{ab}(\slashed{D}\xi)_{n'a}(\slashed{D}\xi)_{n'b}=|\slashed{D}\xi|^{2}
\end{align}
Combining (12.108), (12.109) and (12.111) we conclude that the second term on the right of (12.74) is bounded by:
\begin{align}
 2r\{(1+\sqrt{2})|\slashed{y}'|+\sqrt{\sum_{i}|\tau_{i}|^{2}}\}|\xi||\slashed{D}\xi|
\end{align}
Combining (12.84), (12.97) and (12.112) results in the following proposition:

$\textbf{Proposition 12.5}$ Consider a surface $S_{t,u}$ for which the following hold:
\begin{align*}
 \sup_{S_{t,u}}|y'|\leq 1
\end{align*}
Then for every $S_{t,u}$ 1-form $\xi$ we have:
\begin{align}
 \sum_{i}|\slashed{\mathcal{L}}_{R_{i}}\xi|^{2}\geq r^{2}(1-|y'|^{2})|\slashed{D}\xi|^{2}+(1-2\sqrt{\sum_{i}|\tau_{i}|^{2}})|\xi|^{2}
-2r\{(1+\sqrt{2})|\slashed{y}'|+\sqrt{\sum_{i}|\tau_{i}|^{2}}\}|\xi||\slashed{D}\xi|
\end{align}
pointwise on $S_{t,u}$.

$\textbf{Corollary 12.5.a}$ Suppose that:
\begin{align*}
 \sup_{\Sigma_{t}^{\epsilon_{0}}}|y'|\leq C\delta_{0},  \quad \sup_{\Sigma_{t}^{\epsilon_{0}}}
\sqrt{\sum_{i}|\tau_{i}|^{2}}\leq C\delta_{0}
\end{align*}
for some fixed positive constant $C$. Then if $\delta_{0}$ is suitably small, for any $S_{t,u}$ 1-form $\xi$, differentiable
on $S_{t,u}$, we have, pointwise on $\Sigma_{t}^{\epsilon_{0}}$:
\begin{align*}
 \sum_{i}|\slashed{\mathcal{L}}_{R_{i}}\xi|^{2}\geq \frac{1}{2}\{r^{2}|\slashed{D}\xi|^{2}+|\xi|^{2}\}
\end{align*}
If also:
\begin{align*}
 \inf_{\Sigma_{t}^{\epsilon_{0}}}r\geq \epsilon_{1}(1+t)
\end{align*}
then $\textbf{H1}$ holds on $\Sigma_{t}^{\epsilon_{0}}$. In fact, we have:
\begin{align*}
 (1+t)^{-2}|\xi|^{2}+|\slashed{D}\xi|^{2}\leq C(1+t)^{-2}\sum_{i}|\slashed{\mathcal{L}}_{R_{i}}\xi|^{2}
\end{align*}
where
\begin{align*}
 C=\frac{2}{\epsilon_{1}^{2}}
\end{align*}
$Proof$. From the assumptions, we can choose $\delta_{0}$ small enough such that
\begin{align*}
 \{(1+\sqrt{2})|\slashed{y}^{\prime}|+\sqrt{\sum_{i}|\tau_{i}|^{2}}\}\leq \frac{1}{4}
\end{align*}
and
\begin{align*}
 (1-|y'|^{2})\geq \frac{3}{4}, \quad (1-2\sqrt{\sum_{i}|\tau_{i}|})^{2}\geq \frac{3}{4}
\end{align*}
then the corollary follows. $\qed$

    Let us now return to Lemma 12.2 and let $\bar{s}_{0}$ be the maximal value of $s_{0}\in[0,s]$ such that $\textbf{H1}$ holds on 
$W^{s_{0}}_{\epsilon_{0}}$. We shall show that in fact $\bar{s}_{0}=s$. Suppose that $\bar{s}_{0}<s$. By Proposition 12.3, all 
the assumptions of Corollary 12.5.a except the one about $\tau_{i}$ hold on $W^{s}_{\epsilon_{0}}$. So if we can prove that the assumption about $\tau_{i}$
holds on $W^{s}_{\epsilon_{0}}$, then $\textbf{H1}$ holds on $W^{s}_{\epsilon_{0}}$. To do this, we use a continuity argument.
From (12.88), $\textbf{H0}$, $\textbf{E}_{\{1\}}$ and Proposition 12.3, we have:
\begin{align}
 \max_{i}\sup_{\Sigma_{t}^{\epsilon_{0}}}|\tau_{i}-\lambda_{i}\eta^{-1}\chi|\leq C\delta_{0}^{2}(1+t)^{-2}[1+\log(1+t)]
\end{align}
for all $t\in[0,s]$.

Thus, if for some $t\in[0,s]$, $\chi'$ satisfies on $\Sigma_{t}^{\epsilon_{0}}$ the estimate:
\begin{equation}
 \|\chi'\|_{L^{\infty}(\Sigma_{t}^{\epsilon_{0}})}\leq C_{0}(1+t)^{-1}
\end{equation}
for some fixed positive constant $C_{0}$, then we have:
\begin{equation}
 \max_{i}\|\tau_{i}\|_{L^{\infty}(\Sigma_{t}^{\epsilon_{0}})}\leq C\delta_{0}(1+t)^{-1}[1+\log(1+t)]
\end{equation}
so the assumptions about $\tau_{i}$ holds as well. On the other hand, taking $s_{0}=\bar{s}_{0}$ in Lemma 12.2, yields the estimate:
\begin{equation}
 \|\chi'\|_{L^{\infty}(\Sigma_{t}^{\epsilon_{0}})}\leq C\delta_{0}(1+t)^{-2}[1+\log(1+t)]
\end{equation}
for all $t\in[0,\bar{s}_{0}]$. If $\delta_{0}$ is suitably small this implies that:
\begin{align*}
 \|\chi'\|_{L^{\infty}(\Sigma_{t}^{\epsilon_{0}})}\leq\frac{1}{2}C_{0}(1+t)^{-1}
\end{align*}
at $t=\bar{s}_{0}$.

Therefore by continuity (12.115) holds in a suitably small interval $[\bar{s}_{0},s_{*}]\subset[\bar{s}_{0},s]$. So all the assumptions
of Corollary 12.5.a hold for $t\in[0,s_{*}]$, hence $\textbf{H1}$ holds on $W^{s_{*}}_{\epsilon_{0}}$, contradicting the maximality of $\bar{s}_{0}$.
We conclude that $\textbf{H1}$ holds on $W^{s}_{\epsilon_{0}}$. So from Lemma 12.2 with $s_{0}=s$, we have:

$\textbf{Proposition 12.6}$ Let assumptions $\textbf{E}_{\{2\}}$, $\textbf{E}^{Q}_{\{1\}}$, $\textbf{E}^{QQ}_{\{0\}}$, hold on $W^{s}_{\epsilon_{0}}$,
and let the initial data satisfy:
\begin{align*}
 \|\chi'\|_{L^{\infty}(\Sigma_{0}^{\epsilon_{0}})}\leq C\delta_{0}
\end{align*}
Then, if $\delta_{0}$ is suitably small, $\textbf{H1}$ holds on $W^{s}_{\epsilon_{0}}$. Moreover, for all $t\in[0,s]$:
\begin{align*}
 \|\chi'\|_{L^{\infty}(\Sigma_{t}^{\epsilon_{0}})}\leq C\delta_{0}(1+t)^{-2}[1+\log(1+t)]
\end{align*}
\vspace{7mm}

Next, we shall investigate $\textbf{H2}$ and $\textbf{H2}^{\prime}$. We consider $\slashed{\mathcal{L}}_{R_{i}}\vartheta$ for any 
symmetric 2-covariant $S_{t,u}$ tensorfield $\vartheta$. We have:
\begin{align}
 (\slashed{\mathcal{L}}_{R_{i}}\vartheta)_{ab}=(\slashed{D}_{R_{i}}\vartheta)_{ab}
+\vartheta_{mb}(\slashed{D}R_{i})^{m}_{a}+\vartheta_{am}(\slashed{D}R_{i})^{m}_{b}\\\notag
(\slashed{D}_{R_{i}}\vartheta)_{ab}=(R_{i})^{m}(\slashed{D}\vartheta)_{mab}
\end{align}
Hence we have:
\begin{align}
 \sum_{i}|\slashed{\mathcal{L}}_{R_{i}}\vartheta|^{2}=\sum_{i}|\slashed{D}_{R_{i}}\vartheta|^{2}\\\notag
+2\sum_{i}(\bar{g}^{-1})^{ac}(\bar{g}^{-1})^{bd}(\vartheta_{mb}(\slashed{D}R_{i})^{m}_{a}+\vartheta_{am}(\slashed{D}R_{i})^{m}_{b})
(\slashed{D}_{R_{i}}\vartheta)_{cd}\\\notag
+\sum_{i}(\bar{g}^{-1})^{ac}(\bar{g}^{-1})^{bd}(\vartheta_{mb}(\slashed{D}R_{i})^{m}_{a}+\vartheta_{am}(\slashed{D}R_{i})^{m}_{b})
(\vartheta_{nd}(\slashed{D}R_{i})^{n}_{c}+\vartheta_{cn}(\slashed{D}R_{i})^{n}_{d})
\end{align}
Consider the first term on the right of (12.119). We have:
\begin{align}
 \sum_{i}|\slashed{D}_{R_{i}}\vartheta|^{2}=\sum_{i}(\bar{g}^{-1})^{ac}(\bar{g}^{-1})^{bd}
(\slashed{D}_{R_{i}}\vartheta)_{ab}(\slashed{D}_{R_{i}}\vartheta)_{cd}\\\notag
=(\bar{g}^{-1})^{ac}(\bar{g}^{-1})^{bd}(\sum_{i}(R_{i})^{m}(R_{i})^{n})(\slashed{D}\vartheta)_{mab}(\slashed{D}\vartheta)_{ncd}
\end{align}
Substituting from (12.38) and noting that:
\begin{align*}
 \Pi^{m}_{k}\Pi^{n}_{l}(\slashed{D}\vartheta)_{mab}(\slashed{D}\vartheta)_{ncd}=
(\slashed{D}\vartheta)_{kab}(\slashed{D}\vartheta)_{lcd}
\end{align*}
we obtain:
\begin{equation}
 \sum_{i}|\slashed{D}_{R_{i}}\vartheta|^{2}=r^{2}(\bar{g}^{-1})^{ac}(\bar{g}^{-1})^{bd}
(\delta_{kl}-y'^{k}y'^{l})(\slashed{D}\vartheta)_{kab}(\slashed{D}\vartheta)_{lcd}
\end{equation}
Consider the symmetric 2-covariant $S_{t,u}$ tensorfield $\varphi$:
\begin{equation}
 \varphi_{kl}=(\bar{g}^{-1})^{ac}(\bar{g}^{-1})^{bd}(\slashed{D}\vartheta)_{kab}(\slashed{D}\vartheta)_{lcd}
\end{equation}
Since $\varphi$ is positive semi definite, we have as in (12.81):
\begin{equation}
 y'^{k}y'^{l}\varphi_{kl}\leq |y'|^{2}\textrm{tr}\varphi
\end{equation}
Assuming that $|y'|\leq 1$ and noting that:
\begin{equation}
 \textrm{tr}\varphi=(\bar{g}^{-1})^{kl}(\bar{g}^{-1})^{ac}(\bar{g}^{-1})^{bd}
(\slashed{D}\vartheta)_{kab}(\slashed{D}\vartheta)_{lcd}=|\slashed{D}\vartheta|^{2}
\end{equation}
(12.121) implies
\begin{equation}
 \sum_{i}|\slashed{D}_{R_{i}}\vartheta|^{2}\geq r^{2}(1-|y'|^{2})|\slashed{D}\vartheta|^{2}
\end{equation}

      Consider next the third term on the right of (12.119). This is:
\begin{align*}
 2A_{1}+2A_{2}
\end{align*}
where:
\begin{equation}
 A_{1}=(\bar{g}^{-1})^{ac}(\bar{g}^{-1})^{bd}\vartheta_{mb}\vartheta_{nd}(\sum_{i}(\slashed{D}R_{i})^{m}_{a}(\slashed{D}R_{i})_{c}^{n})
\end{equation}
\begin{equation}
 A_{2}=(\bar{g}^{-1})^{ac}(\bar{g}^{-1})^{bd}\vartheta_{am}\vartheta_{nd}(\sum_{i}(\slashed{D}R_{i})^{m}_{b}(\slashed{D}R_{i})^{n}_{c})
\end{equation}
Substituting (12.90) in (12.126) we find:
\begin{align}
 A_{1}=(\bar{g}^{-1})^{ac}(\bar{g}^{-1})^{bd}\vartheta_{m'b}\vartheta_{n'd}
(\delta_{a'c'}\Pi^{a'}_{a}\Pi^{c'}_{c}\delta_{m'n'}-\Pi^{n'}_{a}\Pi^{m'}_{c})\\\notag
+2\sum_{i}(\bar{g}^{-1})^{ac}(\bar{g}^{-1})^{bd}(w_{i})_{ab}(\tau_{i}\cdot\vartheta)_{cd}+\sum_{i}|\tau_{i}\cdot\vartheta|^{2}
\end{align}
where $w_{i}$ are the 2-covariant $S_{t,u}$ tensorfields: 
\begin{equation}
 (w_{i})_{ab}=\Pi^{a'}_{a}\epsilon_{ia'm'}\vartheta_{m'b}
\end{equation}
and:
\begin{equation}
 (\tau_{i}\cdot\vartheta)_{ab}=(\tau_{i})^{m}_{a}\vartheta_{mb}
\end{equation}
Substituting (12.96) in the first term on the right of (12.128) this term becomes:
\begin{equation}
 2\delta_{m'n'}(\bar{g}^{-1})^{bd}\vartheta_{m'b}\vartheta_{n'd}-|\vartheta|^{2}=
2|\vartheta|^{2}-|\vartheta|^{2}=|\vartheta|^{2}
\end{equation}
The second term on the right of (12.128) is bounded in absolute value by:
\begin{equation}
 2\sum_{i}|w_{i}||\tau_{i}\cdot\vartheta|\leq 2|\vartheta|\sqrt{\sum_{i}|w_{i}|^{2}}\sqrt{\sum_{i}|\tau_{i}|^{2}}
\end{equation}
Using (12.26) we obtain:
\begin{align}
 \sum_{i}|w_{i}|^{2}=\sum_{i}(\bar{g}^{-1})^{ac}(\bar{g}^{-1})^{bd}\Pi^{a'}_{a}\Pi^{c'}_{c}\epsilon_{ia'm'}\epsilon_{ic'n'}
\vartheta_{m'b}\vartheta_{n'd}\\\notag
=(\bar{g}^{-1})^{ac}(\bar{g}^{-1})^{bd}\Pi^{a'}_{a}\Pi^{c'}_{c}(\delta_{a'c'}\delta_{m'n'}-\delta_{a'n'}\delta_{c'm'})
\vartheta_{m'b}\vartheta_{n'd}\\\notag
=(\bar{g}^{-1})^{bd}\delta_{a'c'}(\bar{g}^{-1})^{ac}\Pi^{a'}_{a}\Pi^{c'}_{c}\delta_{m'n'}\vartheta_{m'b}\vartheta_{n'd}
-|\vartheta|^{2}
\end{align}
hence, by (12.96),
\begin{equation}
 \sum_{i}|w_{i}|^{2}=2\delta_{m'n'}(\bar{g}^{-1})^{bd}\vartheta_{m'b}\vartheta_{n'd}-|\vartheta|^{2}=|\vartheta|^{2}
\end{equation}
So we conclude that:
\begin{align}
 A_{1}\geq (1-2\sqrt{\sum_{i}|\tau_{i}|^{2}})|\vartheta|^{2}+\sum_{i}|\tau_{i}\cdot\vartheta|^{2}
\end{align}
Substituting (12.90) in (12.127) we find:
\begin{align}
 A_{2}=(\bar{g}^{-1})^{ac}(\bar{g}^{-1})^{bd}\vartheta_{am'}\vartheta_{n'd}(\delta_{b'c'}\delta_{m'n'}\Pi^{b'}_{b}\Pi^{c'}_{c}
-\Pi^{n'}_{b}\Pi^{m'}_{c})\\\notag
+2\sum_{i}(\bar{g}^{-1})^{ac}(\bar{g}^{-1})^{bd}(\tilde{\tau_{i}\cdot\vartheta})_{ab}(w_{i})_{cd}
+\sum_{i}(\bar{g}^{-1})^{ac}(\bar{g}^{-1})^{bd}(\tilde{\tau_{i}\cdot\vartheta})_{ab}(\tau_{i}\cdot\vartheta)_{cd}
\end{align}
where:
\begin{align*}
 (\tilde{\tau_{i}\cdot\vartheta})_{ab}=(\tau_{i}\cdot\vartheta)_{ba}
\end{align*}
The first term on the right of (12.136) is given by:
\begin{align*}
 (\bar{g}^{-1})^{ac}(\bar{g}^{-1})^{bd}\vartheta_{am'}\vartheta_{n'd}(\Pi^{b}_{c}\delta_{m'n'}-\Pi^{n'}_{b}\Pi^{m'}_{c})
=|\vartheta|^{2}-(\textrm{tr}\vartheta)^{2}
\end{align*}
Here, we have used the fact that 
\begin{align*}
 \delta_{b^{\prime}c^{\prime}}\Pi^{b^{\prime}}_{b}\Pi^{c^{\prime}}_{c}=\Pi_{c}^{b}
\end{align*}

Since:
\begin{align*}
 |\tilde{\tau_{i}\cdot\vartheta}|=|\tau_{i}\cdot\vartheta|
\end{align*}
the second term on the right of (12.136) is bounded in absolute value in the same way as (12.132). The third term is bounded in absolute value by
\begin{equation}
 \sum_{i}|\tau_{i}\cdot\vartheta|^{2}
\end{equation}
We conclude that
\begin{align}
 A_{2}\geq (1-2\sqrt{\sum_{i}|\tau_{i}|^{2}})|\vartheta|^{2}-(\textrm{tr}\vartheta)^{2}-\sum_{i}|\tau_{i}\cdot\vartheta|^{2}
\end{align}
Combining (12.135) and (12.138) we get:
\begin{equation}
 A_{1}+A_{2}\geq2(1-2\sqrt{\sum_{i}|\tau_{i}|^{2}})|\vartheta|^{2}-(\textrm{tr}\vartheta)^{2}
\end{equation}

     Consider finally the second term on the right of (12.119). This is:
\begin{align*}
 4B
\end{align*}
where:
\begin{equation}
 B=(\bar{g}^{-1})^{ac}(\bar{g}^{-1})^{bd}\vartheta_{mb}(\sum_{i}(\slashed{D}R_{i})^{m}_{a}(R_{i})^{k})(\slashed{D}\vartheta)_{kcd}
\end{equation}
Substituting (12.104), by direct calculation we get:
\begin{align}
 B=r\{(\bar{g}^{-1})^{bd}\vartheta_{m'b}\iota_{m'd}-(\bar{g}^{-1})^{bd}(\vartheta\cdot\slashed{y}')_{b}(\textrm{tr}\slashed{D}\vartheta)_{d}\\\notag
+\sum_{i}(\bar{g}^{-1})^{ac}(\bar{g}^{-1})^{bd}(\tau_{i}\cdot\vartheta)_{ab}(\varpi_{i})_{cd}\}
\end{align}
where
\begin{equation}
 \iota_{m'd}=(\slashed{D}\vartheta)_{m'cd}\slashed{y}'^{c}
\end{equation}
\begin{equation}
 (\varpi_{i})_{cd}=N^{l}\epsilon_{ilk'}(\slashed{D}\vartheta)_{k'cd}
\end{equation}
Also:
\begin{align*}
 (\textrm{tr}\slashed{D}\vartheta)_{d}=(\bar{g}^{-1})^{ac}(\slashed{D}\vartheta)_{acd},\quad 
(\vartheta\cdot\slashed{y}')_{b}=\vartheta_{bm}\slashed{y}'^{m}
\end{align*}
Since
\begin{align*}
 (\bar{g}^{-1})^{m'n'}=\delta_{m'n'}
\end{align*}
we have:
\begin{equation}
 |(\bar{g}^{-1})^{bd}\vartheta_{m'b}\iota_{m'd}|\leq |\vartheta||\iota|
\end{equation}
Since
\begin{equation}
 |\iota|\leq|\slashed{D}\vartheta||\slashed{y}'|
\end{equation}
the first term in parenthesis in (12.141) is bounded in absolute value by:
\begin{equation}
 |\vartheta||\slashed{D}\vartheta||\slashed{y}'|
\end{equation}
Since
\begin{align*}
 \frac{1}{2}|\textrm{tr}\slashed{D}\vartheta|^{2}\leq |\slashed{D}\vartheta|^{2}
\end{align*}
the second term in parenthesis in (12.141) is bounded in absolute value by
\begin{equation}
 \sqrt{2}|\slashed{y}'||\vartheta||\slashed{D}\vartheta|
\end{equation}
Finally the third term in parenthesis in (12.141) is bounded in absolute value by:
\begin{equation}
 \sum_{i}|\vartheta||\tau_{i}||\varpi_{i}|\leq |\vartheta|\sqrt{\sum_{i}|\tau_{i}|^{2}}\sqrt{\sum_{i}|\varpi_{i}|^{2}}
\end{equation}
Moreover, using (12.26) we obtain:
\begin{align}
 \sum_{i}|\varpi_{i}|^{2}=\sum_{i}(\bar{g}^{-1})^{ac}(\bar{g}^{-1})^{bd}N^{k}N^{l}\epsilon_{ikm'}\epsilon_{iln'}
(\slashed{D}\vartheta)_{m'ab}(\slashed{D}\vartheta)_{n'cd}\\\notag
=(\bar{g}^{-1})^{ac}(\bar{g}^{-1})^{bd}N^{k}N^{l}(\delta_{kl}\delta_{m'n'}-\delta_{kn'}\delta_{lm'})(\slashed{D}\vartheta)_{m'ab}(\slashed{D}\vartheta)_{n'cd}\\\notag
=(\bar{g}^{-1})^{ac}(\bar{g}^{-1})^{bd}((\slashed{D}\vartheta)_{m'ab}(\slashed{D}\vartheta)_{m'cd}-N^{m'}N^{n'}(\slashed{D}\vartheta)_{m'ab}(\slashed{D}\vartheta)_{n'cd})\\\notag
\leq (\bar{g}^{-1})^{ac}(\bar{g}^{-1})^{bd}(\slashed{D}\vartheta)_{m'ab}(\slashed{D}\vartheta)_{m'cd}=|\slashed{D}\vartheta|^{2}
\end{align}

     So we conclude that
\begin{equation}
 |B|\leq r\{(1+\sqrt{2})|\slashed{y}'|+\sqrt{\sum_{i}|\tau_{i}|^{2}}\}|\slashed{D}\vartheta||\vartheta|
\end{equation}

      Combining (12.125), (12.139) and (12.150) we get the following proposition:

$\textbf{Proposition 12.7}$  Consider a surface $S_{t,u}$ for which the following hold:
\begin{align*}
 \sup_{S_{t,u}}|y'|\leq 1
\end{align*}
Then for every symmetric 2-covariant $S_{t,u}$ tensorfield $\vartheta$ we have:
\begin{align*}
 \sum_{i}|\slashed{\mathcal{L}}_{R_{i}}\vartheta|^{2}\geq r^{2}(1-|y'|^{2})|\slashed{D}\vartheta|^{2}\\\notag
+4(1-2\sqrt{\sum_{i}|\tau_{i}|^{2}})|\vartheta|^{2}-2(\textrm{tr}\vartheta)^{2}\\\notag
-4r\{(1+\sqrt{2})|\slashed{y}'|+\sqrt{\sum_{i}|\tau_{i}|^{2}}\}|\slashed{D}\vartheta||\vartheta|
\end{align*}
pointwise on $S_{t,u}$.

$\textbf{Corollary 12.7.a}$ Suppose that:
\begin{align*}
 \sup_{\Sigma_{t}^{\epsilon_{0}}}|y'|\leq C\delta_{0} \quad\textrm{and}\quad \sup_{\Sigma_{t}^{\epsilon_{0}}}\sqrt{\sum_{i}|\tau_{i}|^{2}} \leq C\delta_{0}
\end{align*}
for some fixed constant $C$. Then if $\delta_{0}$ is suitably small, for any symmetric 2-covariant $S_{t,u}$ tensorfield $\vartheta$, differentiable on $S_{t,u}$,
we have:
\begin{align*}
 \sum_{i}|\slashed{\mathcal{L}}_{R_{i}}\vartheta|^{2}\geq \frac{1}{2}r^{2}|\slashed{D}\vartheta|^{2}+2|\vartheta|^{2}-2(\textrm{tr}\vartheta)^{2}
\end{align*}
If also:
\begin{align*}
 \inf_{\Sigma_{t}^{\epsilon_{0}}}r\geq\epsilon_{1}(1+t)
\end{align*}
then $\textbf{H2}$ and $\textbf{H2}^{\prime}$ hold on $\Sigma_{t}^{\epsilon_{0}}$. In fact, we have, pointwise on $\Sigma_{t}^{\epsilon_{0}}$:
\begin{align*}
 |\slashed{D}\vartheta|^{2}\leq C(1+t)^{-2}\{\sum_{i}|\slashed{\mathcal{L}}_{R_{i}}\vartheta|^{2}+2|\vartheta|^{2}\}
\end{align*}
 where
\begin{align*}
 C=\frac{2}{\epsilon_{1}^{2}}
\end{align*}
Moreover, if $\vartheta$ is trace-free, we have:
\begin{align*}
 |\slashed{D}\vartheta|^{2}\leq C(1+t)^{-2}\sum_{i}|\slashed{\mathcal{L}}_{R_{i}}\vartheta|^{2}
\end{align*}
$Proof$. Choosing $\delta_{0}$ sufficiently small, the second result is straightforward. Concerning the first result, we just use the fact that:
\begin{align*}
 (\textrm{tr}\vartheta)^{2}\leq 2|\vartheta|^{2}
\end{align*}
 $\qed$

     We end this part with coercivity inequalities for Lie derivatives.
 
$\textbf{Proposition 12.8}$ Let the assumptions of Proposition 12.6 hold and let $X$ be an arbitrary $S_{t,u}$ -tangential vectorfield defined on $W^{s}_{\epsilon_{0}}$ 
and differentiable on $S_{t,u}$. Then, for any 1-form $\xi$ we have, pointwise:
\begin{align*}
 |\slashed{\mathcal{L}}_{X}\xi|\leq C(1+t)^{-1}\{|X|(\max_{i}|\slashed{\mathcal{L}}_{R_{i}}\xi|+|\xi|)+|\xi|\max_{i}|\slashed{\mathcal{L}}_{R_{i}}X|\}
\end{align*}
Also, for any symmetric 2-covariant $S_{t,u}$ tensorfield $\vartheta$ we have, pointwise:
\begin{align*}
 |\slashed{\mathcal{L}}_{X}\vartheta|\leq C(1+t)^{-1}\{|X|(\max_{i}|\slashed{\mathcal{L}}_{R_{i}}\vartheta|+|\vartheta|)
+|\vartheta|\max_{i}|\slashed{\mathcal{L}}_{R_{i}}X|\}
\end{align*}
$Proof$. Consider first the case of $\xi$. By Proposition 12.2, we have:
\begin{equation}
 |\slashed{\mathcal{L}}_{X}\xi|\leq C(1+t)^{-1}\sum_{i}|(\slashed{\mathcal{L}}_{X}\xi)(R_{i})|
\end{equation}
Since
\begin{align}
 (\slashed{\mathcal{L}}_{X}\xi)(R_{i})=X(\xi(R_{i}))-\xi([X,R_{i}])
=X\cdot\slashed{d}(\xi(R_{i}))+\xi(\slashed{\mathcal{L}}_{R_{i}}X)
\end{align}
we have:
\begin{equation}
 |(\slashed{\mathcal{L}}_{X}\xi)(R_{i})|\leq |X||\slashed{d}(\xi(R_{i}))|+|\xi||\slashed{\mathcal{L}}_{R_{i}}X|
\end{equation}
By $\textbf{H0}$:
\begin{equation}
 |\slashed{d}(\xi(R_{i}))|\leq C(1+t)^{-1}\sum_{j}|R_{j}(\xi(R_{i}))|
\end{equation}
and we have:
\begin{equation}
 R_{j}(\xi(R_{i}))=(\slashed{\mathcal{L}}_{R_{j}}\xi)(R_{i})+\xi([R_{i},R_{j}])
\end{equation}
hence:
\begin{equation}
 |R_{j}(\xi(R_{i}))|\leq |\slashed{\mathcal{L}}_{R_{j}}\xi||R_{i}|+|\xi||[R_{i},R_{j}]|
\end{equation}
By Corollary 10.1.e:
\begin{equation}
 |R_{i}|,\quad |[R_{i},R_{j}]|\leq C(1+t)
\end{equation}
It follows that:
\begin{equation}
 |R_{j}(\xi(R_{i}))|\leq C(1+t)(|\slashed{\mathcal{L}}_{R_{j}}\xi|+|\xi|)
\end{equation}
Substituting this in (12.154) we obtain:
\begin{equation}
 |\slashed{d}(\xi(R_{i}))|\leq C(\max_{j}|\slashed{\mathcal{L}}_{R_{j}}\xi|+|\xi|)
\end{equation}
Substituting this in (12.153), the result for $\xi$ follows through (12.151).

     Consider next the case of a symmetric 2-covariant $S_{t,u}$ tensorfield $\vartheta$. By Proposition 12.4:
\begin{equation}
 |\slashed{\mathcal{L}}_{X}\vartheta|\leq C(1+t)^{-2}\sum_{i,j}|(\slashed{\mathcal{L}}_{X}\vartheta)(R_{i},R_{j})|
\end{equation}
Now, we have:
\begin{align}
 (\slashed{\mathcal{L}}_{X}\vartheta)(R_{i},R_{j})=
X(\vartheta(R_{i},R_{j}))-\vartheta([X,R_{i}],R_{j})-\vartheta(R_{i},[X,R_{j}])\\\notag
=X\cdot\slashed{d}(\vartheta(R_{i},R_{j}))+\vartheta(\slashed{\mathcal{L}}_{R_{i}}X,R_{j})+\vartheta(R_{i},\slashed{\mathcal{L}}_{R_{j}}X)
\end{align}
hence:
\begin{align}
 |(\slashed{\mathcal{L}}_{X}\vartheta)(R_{i},R_{j})|\leq |X||\slashed{d}(\vartheta(R_{i},R_{j}))|
+|\vartheta|(|R_{i}||\slashed{\mathcal{L}}_{R_{j}}X|+|R_{j}||\slashed{\mathcal{L}}_{R_{i}}X|)\\\notag
\leq |X||\slashed{d}(\vartheta(R_{i},R_{j}))|+C(1+t)|\vartheta|\max_{k}|\slashed{\mathcal{L}}_{R_{k}}X|
\end{align}
in view of the first of (12.157). By $\textbf{H0}$:
\begin{equation}
 |\slashed{d}(\vartheta(R_{i},R_{j}))|\leq C(1+t)^{-1}\sum_{k}|R_{k}(\vartheta(R_{i},R_{j}))|
\end{equation}
and from
\begin{equation}
 R_{k}(\vartheta(R_{i},R_{j}))=(\slashed{\mathcal{L}}_{R_{k}}\vartheta)(R_{i},R_{j})+
\vartheta([R_{k},R_{i}],R_{j})+\vartheta(R_{i},[R_{k},R_{j}])
\end{equation}
as well as (12.157), we have:
\begin{equation}
 |R_{k}(\vartheta(R_{i},R_{j}))|\leq C(1+t)^{2}(|\slashed{\mathcal{L}}_{R_{k}}\vartheta|+|\vartheta|)
\end{equation}
Substituting this in (12.163) we have:
\begin{equation}
 |\slashed{d}(\vartheta(R_{i},R_{j}))|\leq C(1+t)(\max_{k}|\slashed{\mathcal{L}}_{R_{k}}\vartheta|+|\vartheta|)
\end{equation}
Substituting this in (12.162), the result for $\vartheta$ follows through (12.160). $\qed$

\section{Estimates for Higher Order Derivatives of $\chi^{\prime}$ and $\mu$}
$\textbf{Lemma 12.3}$ Let $Y$ be an arbitrary $S_{t,u}$-tangential vectorfield and $\vartheta$ an arbitrary 2-covariant $S_{t,u}$ 
tensorfield, defined on $W^{s}_{\epsilon_{0}}$. Then we have:
\begin{align*}
 \slashed{\mathcal{L}}_{L}\slashed{\mathcal{L}}_{Y}\vartheta-\slashed{\mathcal{L}}_{Y}\slashed{\mathcal{L}}_{L}\vartheta
=\slashed{\mathcal{L}}_{\leftexp{(Y)}{Z}}\vartheta
\end{align*}
$Proof$. We can restrict ourselves on $C_{u}$. We extend $\vartheta$ to $TC_{u}$ by the condition
\begin{align*}
 \vartheta(L,W)=\vartheta(W,L)=0
\end{align*}
for any $W\in TC_{u}$. Thus, for arbitrary $W\in TC_{u}$ we have:
\begin{align*}
 (\mathcal{L}_{L}\vartheta)(L,W)=L(\vartheta(L,W))-\vartheta([L,L],W)-\vartheta(L,[L,W])=0
\end{align*}
(Note that $[L,W]$ is tangential to $C_{u}$). Similarly, $(\mathcal{L}_{L}\vartheta)(W,L)=0$. It follows that:
\begin{equation}
 \slashed{\mathcal{L}}_{L}\vartheta=\mathcal{L}_{L}\vartheta
\end{equation}
On the other hand, for any $S_{t,u}$-tangential vectorfield $X$, by Lemma 8.2 we have:
\begin{align*}
 (\mathcal{L}_{Y}\vartheta)(L,X)=Y(\vartheta(L,X))-\vartheta([Y,L],X)-\vartheta(L,[Y,X])
=\vartheta([L,Y],X)=\vartheta(\leftexp{(Y)}{Z},X)
\end{align*}
Similarly, $(\mathcal{L}_{Y}\vartheta)(X,L)=\vartheta(X,\leftexp{(Y)}{Z})$. Also,
\begin{align*}
 (\mathcal{L}_{Y}\vartheta)(L,L)=Y(\vartheta(L,L))-\vartheta([Y,L],L)-\vartheta(L,[Y,L])=0
\end{align*}
It follows that, on the manifold $C_{u}$:
\begin{equation}
 \slashed{\mathcal{L}}_{Y}\vartheta=\mathcal{L}_{Y}\vartheta-(\vartheta\cdot\leftexp{(Y)}{Z})\otimes dt-
dt\otimes(\leftexp{(Y)}{Z}\cdot\vartheta)
\end{equation}
where, by definition:
\begin{align}
 (\vartheta\cdot\leftexp{(Y)}{Z})(X)=\vartheta(X,\leftexp{(Y)}{Z}),\\\notag
(\leftexp{(Y)}{Z}\cdot\vartheta)(X)=\vartheta(\leftexp{(Y)}{Z},X) : \quad\textrm{for all}\quad X\in TS_{t,u}
\end{align}

    Consider now the evaluation of 
\begin{align*}
 \slashed{\mathcal{L}}_{L}\slashed{\mathcal{L}}_{Y}\vartheta-\slashed{\mathcal{L}}_{Y}\slashed{\mathcal{L}}_{L}\vartheta
\end{align*}
on the $S_{t,u}$ frame vectorfields $X_{A}$. This evaluation is
\begin{equation}
 (\mathcal{L}_{L}(\slashed{\mathcal{L}}_{Y}\vartheta)-\mathcal{L}_{Y}(\slashed{\mathcal{L}}_{L}\vartheta))(X_{A},X_{B})
\end{equation}
Substituting (12.167) and (12.168) this becomes:
\begin{align}
 (\mathcal{L}_{L}\mathcal{L}_{Y}\vartheta-\mathcal{L}_{L}(\vartheta\cdot\leftexp{(Y)}{Z})\otimes dt-dt\otimes
\mathcal{L}_{L}(\leftexp{(Y)}{Z}\cdot\vartheta)-\mathcal{L}_{Y}\mathcal{L}_{L}\vartheta)(X_{A},X_{B})\\\notag
=(\mathcal{L}_{L}\mathcal{L}_{Y}\vartheta-\mathcal{L}_{Y}\mathcal{L}_{L}\vartheta)(X_{A},X_{B})
\end{align}
where we have used the facts that $\mathcal{L}_{L}dt=d(Lt)=0$ and $dt(X_{A})=dt(X_{B})=0$. Therefore (12.170) reduces to
\begin{align}
 (\mathcal{L}_{L}\mathcal{L}_{Y}\vartheta-\mathcal{L}_{Y}\mathcal{L}_{L}\vartheta)(X_{A},X_{B})=
(\mathcal{L}_{[L,Y]}\vartheta)(X_{A},X_{B})=(\mathcal{L}_{\leftexp{(Y)}{Z}}\vartheta)(X_{A},X_{B})
\end{align}
Thus, the lemma follows. $\qed$ \vspace{7mm}

     Given a positive integer $l$ and a multi-index $(i_{1}...i_{l})$, we define the symmetric 2-covariant $S_{t,u}$ tensorfield:
\begin{equation}
 \leftexp{(i_{1}...i_{l})}{\chi}^{\prime}_{l}=\slashed{\mathcal{L}}_{R_{i_{l}}}...\slashed{\mathcal{L}}_{R_{i_{1}}}\chi'
\end{equation}
We shall derive from (12.47), a propagation equation for $\leftexp{(i_{1}...i_{l})}{\chi'}$.
Using Lemma 11.22, we obtain:
\begin{equation}
 [\slashed{\mathcal{L}}_{R_{i_{l}}}...\slashed{\mathcal{L}}_{R_{i_{1}}}, \slashed{\mathcal{L}}_{L}]\chi'
=\sum_{k=0}^{l-1}\slashed{\mathcal{L}}_{R_{i_{l}}}...\slashed{\mathcal{L}}_{R_{i_{l-k+1}}}[\slashed{\mathcal{L}}_{R_{i_{l-k}}},
\slashed{\mathcal{L}}_{L}]\slashed{\mathcal{L}}_{R_{i_{l-k-1}}}...\slashed{\mathcal{L}}_{R_{i_{1}}}\chi'
\end{equation}
By Lemma 12.3 applied to the case $Y=R_{i_{l-k}}, \vartheta=\slashed{\mathcal{L}}_{R_{i_{l-k-1}}}...\slashed{\mathcal{L}}_{R_{i_{1}}}\chi'$,
we obtain:
\begin{equation}
 [\slashed{\mathcal{L}}_{R_{i_{l-k}}},
\slashed{\mathcal{L}}_{L}]\slashed{\mathcal{L}}_{R_{i_{l-k-1}}}...\slashed{\mathcal{L}}_{R_{i_{1}}}\chi'
=-\slashed{\mathcal{L}}_{\leftexp{(R_{i_{l-k}})}{Z}}\slashed{\mathcal{L}}_{R_{i_{l-k-1}}}...\slashed{\mathcal{L}}_{R_{i_{1}}}\chi'
\end{equation}
Substituting in (12.174) then yields:
\begin{align}
 \slashed{\mathcal{L}}_{L}\leftexp{(i_{1}...i_{l})}{\chi}^{\prime}_{l}=\sum_{k=0}^{l-1}\slashed{\mathcal{L}}_{R_{i_{l}}}...\slashed{\mathcal{L}}_{R_{i_{l-k+1}}}
\slashed{\mathcal{L}}_{\leftexp{(R_{i_{l-k}})}{Z}}\slashed{\mathcal{L}}_{R_{i_{l-k-1}}}...\slashed{\mathcal{L}}_{R_{i_{1}}}\chi'\\\notag
+\slashed{\mathcal{L}}_{R_{i_{l}}}...\slashed{\mathcal{L}}_{R_{i_{1}}}(\slashed{\mathcal{L}}_{L}\chi')
\end{align}
We shall express the last term using (12.47). First, we express:
\begin{equation}
 \slashed{\mathcal{L}}_{R_{i_{l}}}...\slashed{\mathcal{L}}_{R_{i_{1}}}(\chi'\cdot\chi')
=\leftexp{(i_{1}...i_{l})}{\chi}^{\prime}_{l}\cdot\chi'+\chi'\cdot\leftexp{(i_{1}...i_{l})}{\chi}^{\prime}_{l}+\leftexp{(i_{1}...i_{l})}{r}_{l}
\end{equation}
where
\begin{equation}
 \leftexp{(i_{1}...i_{l})}{r}_{l}=\sum_{|s_{1}|+|s_{2}|+|s_{3}|=l;|s_{1}|,|s_{3}|< l}(\slashed{\mathcal{L}}_{R})^{s_{1}}\chi'
\cdot(\slashed{\mathcal{L}}_{R})^{s_{2}}\slashed{g}^{-1}\cdot(\slashed{\mathcal{L}}_{R})^{s_{3}}\chi'
\end{equation}
Next, we have:
\begin{equation}
 \slashed{\mathcal{L}}_{R_{i_{l}}}...\slashed{\mathcal{L}}_{R_{i_{1}}}(e\chi^{\prime})=e\leftexp{(i_{1}...i_{l})}{\chi}^{\prime}_{l}+
\leftexp{(i_{1}...i_{l})}{s}_{l}
\end{equation}
where:
\begin{equation}
 \leftexp{(i_{1}...i_{l})}{s}_{l}=\sum_{|s_{1}|+|s_{2}|=l,|s_{2}|<l}(\slashed{\mathcal{L}}_{R})^{s_{1}}e
(\slashed{\mathcal{L}}_{R})^{s_{2}}\chi'
\end{equation}
We then obtain from (12.47) we get:
\begin{equation}
 \slashed{\mathcal{L}}_{R_{i_{l}}}...\slashed{\mathcal{L}}_{R_{i_{1}}}(\slashed{\mathcal{L}}_{L}\chi')
=e\leftexp{(i_{1}...i_{l})}{\chi}^{\prime}_{l}+\leftexp{(i_{1}...i_{l})}{\chi}^{\prime}_{l}\cdot\chi'+\chi'\cdot\leftexp{(i_{1}...i_{l})}{\chi}^{\prime}_{l}
+\leftexp{(i_{1}...i_{l})}{b}_{l}
\end{equation}
where
\begin{equation}
 \leftexp{(i_{1}...i_{l})}{b}_{l}=\leftexp{(i_{1}...i_{l})}{s}_{l}+\leftexp{(i_{1}...i_{l})}{r}_{l}+
\slashed{\mathcal{L}}_{R_{i_{l}}}...\slashed{\mathcal{L}}_{R_{i_{1}}}b
\end{equation}
with
\begin{equation}
 b=\frac{e\slashed{g}}{1-u+t}-\alpha'
\end{equation}

     Given a non-negative integer $l$ let us denote by $\slashed{\textbf{X}}^{\prime}_{l}$ the statement that there is a constant $C$ independent of $s$
such that for all $t\in[0,s]$:
\begin{align*}
 \slashed{\textbf{X}}^{\prime}_{l}:\quad \max_{i_{1}...i_{l}}
\|\slashed{\mathcal{L}}_{R_{i_{l}}}...\slashed{\mathcal{L}}_{R_{i_{1}}}\chi'\|_{L^{\infty}(\Sigma_{0}^{\epsilon_{0}})}
\leq C\delta_{0}(1+t)^{-2}[1+\log(1+t)]
\end{align*}
We then denote by $\slashed{\textbf{X}}^{\prime}_{[l]}$ the statement:
\begin{align*}
 \slashed{\textbf{X}}^{\prime}_{[l]}\quad :\quad\slashed{\textbf{X}}^{\prime}_{0}\quad\textrm{and}...\textrm{and}\quad\slashed{\textbf{X}}^{\prime}_{l}
\end{align*}

$\textbf{Lemma 12.4}$ Let $\textbf{H0}$ and (12.23) hold. Let also the bootstrap assumptions $\slashed{\textbf{E}}_{[l]},
\slashed{\textbf{E}}^{Q}_{[l-1]}$ and $\slashed{\textbf{X}}^{\prime}_{[l-1]}$ hold, for some positive integer $l$. Then the assumption 
$\slashed{\textbf{X}}_{[l-1]}$ (assumption of Proposition 10.1) also holds.

$Proof$. The lemma is trivial for $l=1$. Let the lemma hold with $l$ replaced by $l-1$, for some $l\geq 2$. That is, let
$\slashed{\textbf{X}'}_{[l-2]}$ together with $\slashed{\textbf{E}}_{[l-1]},\slashed{\textbf{E}}^{Q}_{[l-2]}, \textbf{H0}$ and (12.23), 
imply $\slashed{\textbf{X}}_{[l-2]}$. Then by Corollary 10.1.d with $l$ replaced by $l-1$, we have:
\begin{equation}
 \max_{i}\|\leftexp{(R_{i})}{\slashed{\pi}}\|_{\infty,[l-2],\Sigma_{t}^{\epsilon_{0}}}
\leq C_{l}\delta_{0}(1+t)^{-1}[1+\log(1+t)] : \quad\textrm{for all} \quad t\in[0,s]
\end{equation}
Then for any positive integer $n$ we have:
\begin{align}
 \slashed{\mathcal{L}}_{R_{i_{n}}}...\slashed{\mathcal{L}}_{R_{i_{1}}}\chi=\slashed{\mathcal{L}}_{R_{i_{n}}}...\slashed{\mathcal{L}}_{R_{i_{1}}}\chi'
+\frac{1}{1-u+t}\slashed{\mathcal{L}}_{R_{i_{n}}}...\slashed{\mathcal{L}}_{R_{i_{2}}}\leftexp{(R_{i_{1}})}{\slashed{\pi}}
\end{align}
hence:
\begin{align}
\max_{i_{1}...i_{n}}\| \slashed{\mathcal{L}}_{R_{i_{n}}}...\slashed{\mathcal{L}}_{R_{i_{1}}}\chi\|_{L^{\infty}(\Sigma_{t}^{\epsilon_{0}})}
\leq \| \slashed{\mathcal{L}}_{R_{i_{n}}}...\slashed{\mathcal{L}}_{R_{i_{1}}}\chi'\|_{L^{\infty}(\Sigma_{t}^{\epsilon_{0}})}+
C(1+t)^{-1}\max_{j}\|\leftexp{(R_{j})}{\slashed{\pi}}\|_{\infty,[n-1],\Sigma_{t}^{\epsilon_{0}}}
\end{align}
Taking $n=l-1$, the lemma follows. $\qed$

$\textbf{Proposition 12.9}$ Let assumptions of Proposition 12.6 hold. Let also $\textbf{E}_{\{l+2\}},\textbf{E}_{\{l+1\}}^{Q},
\textbf{E}_{\{l\}}^{QQ}$ hold in $W^{s}_{\epsilon_{0}}$ for some non-negative integer $l$, and let the initial data satisfy:
\begin{align*}
 \|\chi'\|_{\infty,[l],\Sigma_{0}^{\epsilon_{0}}}\leq C_{l}\delta_{0}
\end{align*}
Then there is a constant $C_{l}$ independent of $s$ such that for all $t\in[0,s]$ we have:
\begin{align*}
 \|\chi'\|_{\infty,[l],\Sigma_{t}^{\epsilon_{0}}}\leq C_{l}\delta_{0}(1+t)^{-2}[1+\log(1+t)]
\end{align*}
that is, $\slashed{\textbf{X}}^{\prime}_{[l]}$ holds in $W^{s}_{\epsilon_{0}}$.

$Proof$. For $l=0$ the proposition reduces to Proposition 12.6. We proceed by induction. Let the proposition hold with $l$ replaced by $l-1$, for 
some $l\geq 1$. Then $\slashed{\textbf{X}}^{\prime}_{[l-1]}$ holds on $W^{s}_{\epsilon_{0}}$, hence by Lemma 12.4, $\slashed{\textbf{X}}_{[l-1]}$
holds on $W^{s}_{\epsilon_{0}}$ as well. Then the Proposition 10.1 and all its corollaries hold.

     By Proposition 12.6 and the estimate (12.48) we have:
\begin{align}
 \|e+2\chi'\|_{L^{\infty}(\Sigma_{t}^{\epsilon_{0}})}\leq C\delta_{0}(1+t)^{-2}[1+\log(1+t)]
\quad : \quad\textrm{for all}\quad t\in[0,s]
\end{align}
We shall estimate $\leftexp{(i_{1}...i_{l})}{b}_{l}$ using the corollaries of Proposition 10.1.
From the expression of $e, \kappa^{-1}\zeta$ and the assumptions of the present proposition, we have:
\begin{align}
 \|e\|_{\infty,[l],\Sigma_{t}^{\epsilon_{0}}}, \|\kappa^{-1}\zeta\|_{\infty,[l],\Sigma_{t}^{\epsilon_{0}}},
\|\slashed{k}\|_{\infty,[l],\Sigma_{t}^{\epsilon_{0}}} 
\leq C_{l}\delta_{0}(1+t)^{-2}
\quad:\quad\textrm{for all}\quad t\in[0,s]
\end{align}
Here, we have used Corollary 10.1.a, and Corollary 10.1.b.

Similarly, from the expression of $\alpha'$, we have:
\begin{equation}
 \|\alpha'\|_{\infty,[l],\Sigma_{t}^{\epsilon_{0}}}\leq C\delta_{0}(1+t)^{-3}[1+\log(1+t)]\max_{j}
|\slashed{\mathcal{L}}_{R_{j}}\slashed{\mathcal{L}}_{R_{i_{l}}}...\slashed{\mathcal{L}}_{R_{i_{2}}}\chi^{\prime}|
+C_{l}\delta_{0}(1+t)^{-3}\\
\end{equation}
For all $t\in[0,s]$.

To see how the principal term in $\alpha^{\prime}$, namely $\slashed{D}^{2}h$, is estimated in the $\|\quad\|_{\infty,[l],\Sigma_{t}^{\epsilon_{0}}}$, we
appeal to Lemma 10.9, $\textbf{H0}, \textbf{H1}, \textbf{H2}$ and the fact that:
\begin{align*}
 \leftexp{(R_{i})}{\pi}_{AB}=-2\lambda_{i}(-\eta^{-1}\chi_{AB}+\slashed{k}_{AB})
\end{align*}
(see (3.27) and (6.59))
So we get:
\begin{align}
 |\slashed{\mathcal{L}}_{R_{i_{l}}}...\slashed{\mathcal{L}}_{R_{i_{1}}}b|\leq C\delta_{0}(1+t)^{-3}[1+\log(1+t)]
\max_{j}|\slashed{\mathcal{L}}_{R_{j}}\slashed{\mathcal{L}}_{R_{i_{l}}}...\slashed{\mathcal{L}}_{R_{i_{2}}}\chi^{\prime}|C_{l}\delta_{0}(1+t)^{-3}\\\notag
 :\quad\textrm{pointwise on}\quad
W^{s}_{\epsilon_{0}}
\end{align}
From the inductive hypothesis and Corollary 10.1.d:
\begin{align*}
 \max_{i;i_{1}...i_{k}}\|\slashed{\mathcal{L}}_{R_{i_{k}}}...\slashed{\mathcal{L}}_{R_{i_{1}}}\leftexp{(R_{i})}{\slashed{\pi}}\|_{L^{\infty}(\Sigma_{t}^{\epsilon_{0}})}
\leq C_{l}\delta_{0}(1+t)^{-1}[1+\log(1+t)]\quad:\quad\textrm{for all}\quad k=0,1,...,l-1
\end{align*}
which is needed to estimate $\leftexp{(i_{1}...i_{l})}{r}_{l}$ given by (12.178). We then deduce:
\begin{equation}
 \|\leftexp{(i_{1}...i_{l})}{s_{l}}\|_{L^{\infty}(\Sigma_{t}^{\epsilon_{0}})}\leq C_{l}\delta_{0}^{2}(1+t)^{-4}[1+\log(1+t)]
\quad:\quad\textrm{for all}\quad t\in[0,s]
\end{equation}
and:
\begin{equation}
 \|\leftexp{(i_{1}...i_{l})}{r_{l}}\|_{L^{\infty}(\Sigma_{t}^{\epsilon_{0}})}\leq C_{l}\delta_{0}^{2}(1+t)^{-4}[1+\log(1+t)]^{2}
\quad:\quad\textrm{for all}\quad t\in[0,s]
\end{equation}
It follows that:
\begin{align}
 |\leftexp{(i_{1}...i_{l})}{b}_{l}|\leq C\delta_{0}(1+t)^{-3}[1+\log(1+t)]\max_{j}|\slashed{\mathcal{L}}_{R_{j}}\slashed{\mathcal{L}}_{R_{i_{l}}}
...\slashed{\mathcal{L}}_{R_{i_{2}}}\chi^{\prime}|+C_{l}\delta_{0}(1+t)^{-3}\\\notag
:\quad\textrm{for all}\quad t\in[0,s]
\end{align}
  From (12.181), (12.187) and (12.193) we then conclude that:
\begin{align}
 |\slashed{\mathcal{L}}_{R_{i_{l}}}...\slashed{\mathcal{L}}_{R_{i_{1}}}(\slashed{\mathcal{L}}_{L}\chi')|
\leq C\delta_{0}(1+t)^{-2}[1+\log(1+t)]\max_{j_{1}...j_{l}}|\leftexp{(j_{1}...j_{l})}{\chi}^{\prime}_{l}|\\\notag+C_{l}\delta_{0}(1+t)^{-3}\quad
:\quad\textrm{pointwise on}\quad W^{s}_{\epsilon_{0}}
\end{align}

      What remains to be estimated is the sum on the right of (12.176). Consider a given term in this sum, corresponding to 
$k\in\{0,...,l-1\}$. We use Lemma 8.5, by taking $\slashed{\mathcal{L}}_{R_{i_{l-k-1}}}...\slashed{\mathcal{L}}_{R_{i_{1}}}\chi'$
in the role of $\xi$, and the index set $\{l-k+1,...,l\}$ in the role of the index set $\{1,...,l\}$; to express the term in question
in the form:
\begin{align}
 \slashed{\mathcal{L}}_{R_{i_{l}}}...\slashed{\mathcal{L}}_{R_{i_{l-k+1}}}\slashed{\mathcal{L}}_{\leftexp{(R_{i_{l-k}})}{Z}}
\slashed{\mathcal{L}}_{R_{i_{l-k-1}}}...\slashed{\mathcal{L}}_{R_{i_{1}}}\chi'=\slashed{\mathcal{L}}_{\leftexp{(R_{i_{l-k}})}{Z}}
\leftexp{(i_{1}\overset{>i_{l-k}<}{...}i_{l})}{\chi}^{\prime}_{l-1}\\\notag
+\sum_{p=1}^{k}\sum_{m_{1}<...<m_{p}=l-k+1}^{l}\slashed{\mathcal{L}}_{\leftexp{(i_{m_{1}}...i_{m_{p}};i_{l-k})}{Z}}
\leftexp{(i_{1}\overset{>i_{m_{1}}...i_{m_{p}}i_{l-k}<}{...}i_{l})}{\chi}^{\prime}_{l-1-p}
\end{align}
where
\begin{equation}
 \leftexp{(i_{m_{1}}...i_{m_{p}};i_{l-k})}{Z}=\slashed{\mathcal{L}}_{R_{i_{m_{p}}}}...\slashed{\mathcal{L}}_{R_{i_{m_{1}}}}\leftexp{(R_{i_{l-k}})}{Z}
\end{equation}
(noting that Lemma 8.5 holds for an arbitrary $S_{t,u}$ tensorfields $\xi$)

To estimate the terms on the right of (12.195), we apply the second statement of Proposition 12.8. We then obtain, for the first term:
\begin{align}
 |\slashed{\mathcal{L}}_{\leftexp{(R_{i_{l-k}})}{Z}}\leftexp{(i_{1}\overset{>i_{l-k}<}{...}i_{l})}{\chi}^{\prime}_{l-1}|\leq
C(1+t)^{-1}\{|\leftexp{(R_{i_{l-k}})}{Z}|(\max_{j}|\leftexp{(i_{1}\overset{>i_{l-k}<}{...}i_{l}j)}{\chi}^{\prime}_{l}|
+|\leftexp{(i_{1}\overset{>i_{l-k}<}{...}i_{l})}{\chi}^{\prime}_{l-1}|)\\\notag
+|\leftexp{(i_{1}\overset{>i_{l-k}<}{...}i_{l})}{\chi}^{\prime}_{l-1}|\max_{j}|\leftexp{(j;i_{l-k})}{Z}|\}
\end{align}
In the case of $l\geq 2$, Corollary 10.1.i together with the inductive hypothesis $\slashed{\textbf{X}}^{\prime}_{[l-1]}$ implies that 
the right hand side of (12.197) is bounded pointwise by:
\begin{align}
 C_{l}\delta_{0}(1+t)^{-2}[1+\log(1+t)]\max_{j_{1}...j_{l}}|\leftexp{(j_{1}...j_{l})}{\chi}^{\prime}_{l}|+C_{l}\delta_{0}^{2}(1+t)^{-4}
[1+\log(1+t)]^{2}
\end{align}
In the case of $l=1$, Corollary 10.1.i can only give the estimate for $\leftexp{(R_{i})}{Z}$, so the term 
$\leftexp{(j;i_{l-k})}{Z}$ can not be estimated in this way.
Recalling that $\leftexp{(R_{i})}{Z}=\leftexp{(R_{i})}{\slashed{\pi}}_{L}\cdot\slashed{g}^{-1}$, we use (6.70) to express:
\begin{equation}
 \leftexp{(R_{i})}{\slashed{\pi}}_{L}+\chi'\cdot R_{i}=\epsilon_{ijm}z^{j}\slashed{d}x^{m}+\lambda_{i}(\kappa^{-1}\zeta)
\end{equation}
By Corollary 10.1.g, Corollary 10.1.a for $(R)^{s}x^{j}$ and $\lambda_{i}$, we have the following estimates:
\begin{align*}
 \max_{i}\|\lambda_{i}\|_{\infty,[l],\Sigma_{t}^{\epsilon_{0}}}\leq C_{l}\delta_{0}[1+\log(1+t)],\quad
\max_{i}\|\slashed{d}x^{i}\|_{\infty,[l],\Sigma_{t}^{\epsilon_{0}}}\leq C,\\\notag
\max_{j}\|z^{j}\|_{\infty,[l],\Sigma_{t}^{\epsilon_{0}}}\leq C\delta_{0}(1+t)^{-1}[1+\log(1+t)]\quad: \textrm{for all} \quad t\in[0,s]
\end{align*}
In view also of the second of (12.188), we then obtain:
\begin{equation}
 \|\leftexp{(R_{i})}{\slashed{\pi}}_{L}+\chi'\cdot R_{i}\|_{\infty,[1],\Sigma_{t}^{\epsilon_{0}}}\leq 
C\delta_{0}(1+t)^{-1}[1+\log(1+t)]\quad : \textrm{for all}\quad t\in[0,s]
\end{equation}
Moreover, we have, by Corollary 10.1.e with $l=1$ and $\slashed{\textbf{X}}^{\prime}_{0}$:
\begin{align}
 |\slashed{\mathcal{L}}_{R_{j}}(\chi'\cdot R_{i})|\leq |\slashed{\mathcal{L}}_{R_{j}}\chi'||R_{i}|+|\chi'||\slashed{\mathcal{L}}_{R_{j}}R_{i}|
\leq C(1+t)|\slashed{\mathcal{L}}_{R_{j}}\chi'|+C\delta_{0}(1+t)^{-1}[1+\log(1+t)]
\end{align}
Hence
\begin{equation}
 |\slashed{\mathcal{L}}_{R_{j}}\leftexp{(R_{i})}{\slashed{\pi}}_{L}|\leq C(1+t)|\slashed{\mathcal{L}}_{R_{j}}\chi'|
+C\delta_{0}(1+t)^{-1}[1+\log(1+t)]
\end{equation}
therefore by Corollary 10.1.d, we have:
\begin{equation}
 |\slashed{\mathcal{L}}_{R_{j}}\leftexp{(R_{i})}{Z}|\leq C(1+t)|\slashed{\mathcal{L}}_{R_{j}}\chi'|
+C\delta_{0}(1+t)^{-1}[1+\log(1+t)]
\end{equation}
pointwise on $W^{s}_{\epsilon_{0}}$.
Using this, we can estimate the last term on the right of (12.197). Thus, also in the case $l=1$, an estimate of the form (12.198) holds.

     Consider next a term in the double sum in (12.195), corresponding to an ordered subset $\{m_{1},...,m_{p}\}\\
\subset\{l-k+1,...,l\},\quad p\in\{1,...,k\}$. This double sum is present only for $l\geq 2$. Applying the second statement of Proposition 12.8 we obtain:
\begin{align}
 |\slashed{\mathcal{L}}_{\leftexp{(i_{m_{1}}...i_{m_{p}};i_{l-k})}{Z}}\leftexp{(i_{1}\overset{>i_{m_{1}}...i_{m_{p}}i_{l-k}<}{...}i_{l})}{\chi}^{\prime}_{l-1-p}|
\leq C(1+t)^{-1}\cdot\\\notag
\{|\leftexp{(i_{m_{1}}...i_{m_{p}};i_{l-k})}{Z}|(\max_{j}|\leftexp{(i_{1}\overset{>i_{m_{1}}...i_{m_{p}}i_{l-k}<}{...}i_{l}j)}{\chi}^{\prime}_{l-p}|
+|\leftexp{(i_{1}\overset{>i_{m_{1}}...i_{m_{p}}i_{l-k}<}{...}i_{l})}{\chi}^{\prime}_{l-1-p}|)\\\notag
+|\leftexp{(i_{1}\overset{>i_{m_{1}}...i_{m_{p}}i_{l-k}<}{...}i_{l})}{\chi}^{\prime}_{l-1-p}|\max_{j}|\leftexp{(i_{m_{1}}...i_{m_{p}}j;i_{l-k})}{Z}|\}
\end{align}
Since $p\geq 1$, by the inductive hypothesis $\slashed{\textbf{X}}^{\prime}_{[l-1]}$:
\begin{align}
 \|\leftexp{(i_{1}\overset{>i_{m_{1}}...i_{m_{p}}i_{l-k}<}{...}i_{l})}{\chi}^{\prime}_{l-1-p}\|_{L^{\infty}(\Sigma_{t}^{\epsilon_{0}})}
\leq C_{l}\delta_{0}(1+t)^{-2}[1+\log(1+t)]\\\notag
\max_{j}\|\leftexp{(i_{1}\overset{>i_{m_{1}}...i_{m_{p}}i_{l-k}<}{...}i_{l}j)}{\chi}^{\prime}_{l-p}\|_{L^{\infty}(\Sigma_{t}^{\epsilon_{0}})}\leq 
C_{l}\delta_{0}(1+t)^{-2}[1+\log(1+t)]\quad : \textrm{for all}\quad t\in[0,s]
\end{align}
Also, since $p\leq k\leq l-1$, we have, by Corollary 10.1.i:
\begin{equation}
 \|\leftexp{(i_{m_{1}}...i_{m_{p}};i_{l-k})}{Z}\|_{L^{\infty}(\Sigma_{t}^{\epsilon_{0}})}\leq C_{l}\delta_{0}(1+t)^{-1}[1+\log(1+t)]
\quad :\textrm{for all}\quad t\in[0,s]
\end{equation}
But we can not estimate $\leftexp{(i_{m_{1}}...i_{m_{p}}j;i_{l-k})}{Z}$ when $p=k=l-1$ in this way for the same reason as above.
However, using (12.199) and the estimates below (12.199), as well as the second of (12.188) and results in Chapter 10
we obtain:
\begin{align}
 \|\leftexp{(R_{i})}{\slashed{\pi}}_{L}+\chi'\cdot R_{i}\|_{\infty,[l],\Sigma_{t}^{\epsilon_{0}}}\leq C_{l}\delta_{0}(1+t)^{-1}[1+\log(1+t)]
\quad : \textrm{for all}\quad t\in[0,s]
\end{align}
Moreover, by Corollary 10.1.e and $\slashed{\textbf{X}}^{\prime}_{[l-1]}$,
\begin{align}
 |\slashed{\mathcal{L}}_{R_{j_{l}}}...\slashed{\mathcal{L}}_{R_{j_{1}}}(\chi'\cdot R_{i})|\leq 
|\slashed{\mathcal{L}}_{R_{j_{l}}}...\slashed{\mathcal{L}}_{R_{j_{1}}}\chi'||R_{i}|+\sum_{|s_{1}|+|s_{2}|=l,s_{2}\geq 1}
|(\slashed{\mathcal{L}}_{R})^{s_{1}}\chi'||(\slashed{\mathcal{L}}_{R})^{s_{2}}R_{i}|\\\notag
\leq C(1+t)|\leftexp{(j_{1}...j_{l})}{\chi}^{\prime}_{l}|+C_{l}\delta_{0}(1+t)^{-1}[1+\log(1+t)]
\end{align}
Combining (12.207) and (12.208) yields:
\begin{equation}
 |\slashed{\mathcal{L}}_{R_{j_{1}}}...\slashed{\mathcal{L}}_{R_{j_{1}}}\leftexp{(R_{i})}{\slashed{\pi}}_{L}|\leq C(1+t)
|\leftexp{(j_{1}...j_{l})}{\chi}^{\prime}_{l}|+C_{l}\delta_{0}(1+t)^{-1}[1+\log(1+t)]
\end{equation}
hence, by Corollary 10.1.d, also:
\begin{align}
 |\leftexp{(j_{1}...j_{l};i)}{Z}|\leq C(1+t)|\leftexp{(j_{1}...j_{l})}{\chi}^{\prime}_{l}|+
C_{l}\delta_{0}(1+t)^{-1}[1+\log(1+t)]
\end{align}
Combining the above results, we obtain:

For $p=k=l-1$: 
\begin{align}
 |\slashed{\mathcal{L}}_{\leftexp{(i_{m_{1}}...i_{m_{p}};i_{l-k})}{Z}}\leftexp{(i_{1}\overset{>i_{m_{1}}...i_{m_{p}}i_{l-k}<}{...}i_{l})}{\chi}^{\prime}_{l-1-p}|\\\notag
\leq C_{l}\delta_{0}(1+t)^{-2}[1+\log(1+t)]\max_{j_{1}...j_{l}}|\leftexp{(j_{1}...j_{l})}{\chi}^{\prime}_{l}|+C_{l}\delta_{0}^{2}(1+t)^{-4}[1+\log(1+t)]^{2}
\end{align}
Otherwise:
\begin{align}
  |\slashed{\mathcal{L}}_{\leftexp{(i_{m_{1}}...i_{m_{p}};i_{l-k})}{Z}}\leftexp{(i_{1}\overset{>i_{m_{1}}...i_{m_{p}}i_{l-k}<}{...}i_{l})}{\chi}^{\prime}_{l-k-p}|
\leq C_{l}\delta_{0}^{2}(1+t)^{-4}[1+\log(1+t)]^{2}
\end{align}
Combining these with the estimate (12.198) for the first term on the right in (12.195) we conclude that:
\begin{align}
|\sum_{k=0}^{l-1}\slashed{\mathcal{L}}_{R_{i_{l}}}...\slashed{\mathcal{L}}_{R_{i_{l-k+1}}}\slashed{\mathcal{L}}_{\leftexp{(R_{i_{l-k}})}{Z}}
\slashed{\mathcal{L}}_{R_{i_{l-k-1}}}...\slashed{\mathcal{L}}_{R_{i_{1}}}\chi'|\\\notag
\leq C_{l}\delta_{0}(1+t)^{-2}[1+\log(1+t)]\max_{j_{1}...j_{l}}|\leftexp{(j_{1}...j_{l})}{\chi}^{\prime}_{l}|+C_{l}\delta_{0}^{2}(1+t)^{-4}[1+\log(1+t)]^{2}
\end{align}
This together with (12.194) yields:
\begin{align}
 |\slashed{\mathcal{L}}_{L}\leftexp{(i_{1}...i_{l})}{\chi}^{\prime}_{l}|\leq C_{l}\delta_{0}(1+t)^{-2}[1+\log(1+t)]\max_{j_{1}...j_{l}}|\leftexp{(j_{1}...j_{l})}
{\chi}^{\prime}_{l}|
+C_{l}\delta_{0}(1+t)^{-3}
\end{align}
Applying then (12.56), taking $\vartheta=\leftexp{(i_{1}...i_{l})}{\chi}^{\prime}_{l}$, and using $\slashed{\textbf{X}}^{\prime}_{0}$ we deduce:
\begin{align}
 L((1-u+t)^{2}|\leftexp{(i_{1}...i_{l})}{\chi}^{\prime}_{l}|)\leq C_{l}\delta_{0}(1+t)^{-2}[1+\log(1+t)]
((1-u+t)^{2}\max_{j_{1}...j_{l}}|\leftexp{(j_{1}...j_{l})}{\chi}^{\prime}_{l}|)\\\notag+C_{l}\delta_{0}(1+t)^{-1}
\end{align}
Setting along a given integral curve of $L$,
\begin{equation}
 \leftexp{(i_{1}...i_{l})}{x}_{l}(t)=(1-u+t)^{2}|\leftexp{(i_{1}...i_{l})}{\chi}^{\prime}_{l}|
\end{equation}
(12.215) becomes:
\begin{equation}
 \frac{d\leftexp{(i_{1}...i_{l})}{x}_{l}}{dt}\leq f_{l}\max_{j_{1}...j_{l}}\leftexp{(j_{1}...j_{l})}{x}_{l}+g_{l}
\end{equation}
where:
\begin{equation}
 f_{l}(t)=C_{l}\delta_{0}(1+t)^{-2}[1+\log(1+t)],\quad g_{l}(t)=C_{l}\delta_{0}(1+t)^{-1}
\end{equation}
Integrating from $t=0$ yields:
\begin{equation}
 \leftexp{(i_{1}...i_{l})}{x}_{l}(t)\leq \leftexp{(i_{1}...i_{l})}{x}_{l}(0)
+\int_{0}^{t}(f_{l}(t')\max_{j_{1}...j_{l}}\leftexp{(j_{1}...j_{l})}{x}_{l}(t')+g_{l}(t'))dt'
\end{equation}
Taking the maximum over $i_{1}...i_{l}$, on both sides, and setting:
\begin{equation}
 \bar{x}_{l}(t)=\max_{i_{1}...i_{l}}\leftexp{(i_{1}...i_{l})}{x}_{l}(t)
\end{equation}
then we obtain
\begin{equation}
 \bar{x}_{l}(t)\leq \bar{x}_{l}(0)+\int_{0}^{t}(f_{l}(t')\bar{x}_{l}(t')+g_{l}(t'))dt'
\end{equation}
which implies:
\begin{equation}
 \bar{x}_{l}(t)\leq e^{\int_{0}^{t}f_{l}(t')dt'}\{\bar{x}_{l}(0)+\int_{0}^{t}g_{l}(t')dt'\}
\end{equation}
Since
\begin{equation}
 \bar{x}_{l}(0)\leq \max_{i_{1}...i_{l}}\|\leftexp{(i_{1}...i_{l})}{\chi}^{\prime}_{l}\|_{L^{\infty}(\Sigma_{0}^{\epsilon_{0}})}\leq C_{l}\delta_{0}
\end{equation}
while
\begin{equation}
 \bar{x}_{l}(t)\geq \frac{1}{4}(1+t)^{2}\max_{i_{1}...i_{l}}|\leftexp{(i_{1}...i_{l})}{\chi}^{\prime}_{l}|\quad (u\leq \epsilon_{0}\leq \frac{1}{2})
\end{equation}
taking into account the facts that:
\begin{align*}
 \int_{0}^{t}f_{l}(t')dt'\leq \int_{0}^{\infty}f_{l}(t')dt'=2C_{l}, \qquad
\int_{0}^{t}g_{l}(t')dt'=C_{l}\delta_{0}\log(1+t)
\end{align*}
then taking the supremum on $\Sigma_{t}^{\epsilon_{0}}$, we conclude that:
\begin{equation}
 \max_{i_{1}...i_{l}}\|\leftexp{(i_{1}...i_{l})}{\chi}^{\prime}_{l}\|_{L^{\infty}(\Sigma_{t}^{\epsilon_{0}})}
\leq C_{l}\delta_{0}(1+t)^{-2}[1+\log(1+t)]\quad : \textrm{for all}\quad t\in[0,s]
\end{equation}
This, together with the inductive hypothesis $\slashed{\textbf{X}}^{\prime}_{[l-1]}$ implies $\slashed{\textbf{X}}^{\prime}_{[l]}$, 
completing the inductive step and therefore the proof of the proposition. $\qed$

      From Lemma 12.4 with $l+1$ in the role of $l$ we conclude that under the assumptions of Proposition 12.9, $\slashed{\textbf{X}}_{[l]}$ holds on $W^{s}_{\epsilon_{0}}$.
So Proposition 10.1 and all its Corollaries hold with $l+1$ in the role of $l$.

$\textbf{Proposition 12.10}$ Let the assumptions of Proposition 12.9 hold for some non-negative integer $l$. Let, in addition the initial data
satisfy:
\begin{align*}
 \|\mu-1\|_{\infty,[m,l+1],\Sigma_{0}^{\epsilon_{0}}}\leq C_{l}\delta_{0}
\end{align*}
for some $m\in\{0,...,l+1\}$. Then there is a constant $C_{l}$ independent of $s$ such that for all $t\in[0,s]$ we have:
\begin{align*}
 \|\mu-1\|_{\infty,[m,l+1],\Sigma_{t}^{\epsilon_{0}}}\leq C_{l}\delta_{0}[1+\log(1+t)]
\end{align*}
that is, $\textbf{M}_{[m,l+1]}$ holds on $W^{s}_{\epsilon_{0}}$.

$Proof:$ We consider first the case $m=0$. Since $M_{[0,0]}$ follows directly from Proposition 12.1, we apply induction on $l$.
Assuming then $\textbf{M}_{[0,l]}$ we are to establish $\textbf{M}_{[0,l+1]}$ under the assumptions of the proposition for $m=0$.

     We consider the equation:
\begin{equation}
 L\mu=m+\mu e
\end{equation}
We apply $R_{i_{l+1}}...R_{i_{1}}$ to this equation. By Lemma 11.22, we obtain:
\begin{align}
 [R_{i_{l+1}}...R_{i_{1}}, L]=-\sum_{k=0}^{l}R_{i_{l+1}}...R_{i_{l+2-k}}\leftexp{(R_{i_{l+1-k}})}{Z}R_{i_{l-k}}...R_{i_{1}}
\end{align}
Defining:
\begin{equation}
 \leftexp{(i_{1}...i_{l+1})}{\mu}_{0,l+1}=R_{i_{l+1}}...R_{i_{1}}\mu
\end{equation}
We get:
\begin{align}
 L\leftexp{(i_{1}...i_{l+1})}{\mu}_{0,l+1}=\sum_{k=0}^{l}R_{i_{l+1}}...R_{i_{i+2-k}}(\leftexp{(R_{i_{l+1-k}})}{Z}\cdot\slashed{d}
\leftexp{(i_{1}...i_{l-k})}{\mu}_{0,l-k})\\\notag
+e\leftexp{(i_{1}...i_{l+1})}{\mu}_{0,l+1}+R_{i_{l+1}}...R_{i_{1}}m+\leftexp{(i_{1}...i_{l+1})}{r}^{\prime}_{0,l+1}
\end{align}
where:
\begin{equation}
 \leftexp{(i_{1}...i_{l+1})}{r}^{\prime}_{0,l+1}=\sum_{|s_{1}|+|s_{2}|=l+1,s_{1}\geq 1}
((R)^{s_{1}}e)((R)^{s_{2}}\mu)
\end{equation}
 
     Now by $\textbf{E}_{\{l+2\}}, \textbf{E}^{Q}_{\{l+1\}}$ and Corollary 10.1.a (remember now that Proposition 10.1 is valid for $l+1$):
\begin{align}
 \|Th\|_{\infty,[l+1],\Sigma_{t}^{\epsilon_{0}}}\leq C_{l}\delta_{0}(1+t)^{-1}\\\notag
\|Lh\|_{\infty,[l+1],\Sigma_{t}^{\epsilon_{0}}}\leq C_{l}\delta_{0}(1+t)^{-2}\\\notag
\textrm{also}:\quad \|\omega_{L\hat{T}}\|_{\infty,[l+1],\Sigma_{t}^{\epsilon_{0}}}\leq C_{l}\delta_{0}(1+t)^{-2}
\quad: \textrm{for all}\quad t\in[0,s]
\end{align}
Thus using the expression for $e$ we obtain:
\begin{align}
 \|e\|_{\infty,[l+1],\Sigma_{t}^{\epsilon_{0}}}\leq C_{l}\delta_{0}(1+t)^{-2}\quad :\textrm{for all}\quad t\in[0,s]
\end{align}
Using the first of (12.231) we obtain:
\begin{equation}
 \|m\|_{\infty,[l+1],\Sigma_{t}^{\epsilon_{0}}}\leq C_{l}\delta_{0}(1+t)^{-1}\quad :\textrm{for all}\quad t\in[0,s]
\end{equation}
The estimate (12.232) together with the induction hypothesis $\textbf{M}_{[0,l]}$ imply:
\begin{align}
 \|\leftexp{(i_{1}...i_{l+1})}{r}^{\prime}_{0,l+1}\|_{L^{\infty}(\Sigma_{t}^{\epsilon_{0}})}
\leq C_{l}\delta_{0}(1+t)^{-2}[1+\log(1+t)]\quad :\textrm{for all}\quad t\in[0,s]
\end{align}
Consider finally a term in the sum on the right in (12.229), corresponding to $k\in\{0,...,l\}$. We express the term as:
\begin{align}
 R_{i_{l+1}}...R_{i_{l+2-k}}(\leftexp{(R_{i_{l+1-k}})}{Z}\cdot\slashed{d}\leftexp{(i_{1}...i_{l-k})}{\mu}_{0,l-k})=
\leftexp{(R_{i_{l+1-k}})}{Z}\cdot\slashed{d}\leftexp{(i_{1}\overset{>i_{l+1-k}<}{...}i_{l+1})}{\mu}_{0,l}\\\notag
+\sum_{|s_{1}|+|s_{2}|=k,|s_{1}|\geq 1}((\slashed{\mathcal{L}}_{R})^{s_{1}}\leftexp{(R_{i_{l+1-k}})}{Z})\cdot
\slashed{d}(R)^{s_{2}}\leftexp{(i_{1}...i_{l-k})}{\mu}_{0,l-k}
\end{align}
Consider first the sum on the right of (12.235). By Corollary 10.1.i with $l+1$ in the role of $l$:
\begin{align}
 \max_{i}\|\leftexp{(R_{i})}{Z}\|_{\infty,[l],\Sigma_{t}^{\epsilon_{0}}}\leq C_{l}\delta_{0}(1+t)^{-1}[1+\log(1+t)]
\quad :\textrm{for all}\quad t\in[0,s]
\end{align}
Since $|s_{1}|\leq k\leq l$, this allows us to bound the first factor in the sum. Moreover, by $\textbf{H0}$ and the induction hypothesis
$\textbf{M}_{[0,l]}$, the second factor in the sum can be bounded in $L^{\infty}(\Sigma_{t}^{\epsilon_{0}})$ by:
\begin{equation}
 C_{l}\delta_{0}(1+t)^{-1}[1+\log(1+t)]\notag
\end{equation}
since $1+|s_{2}|+l-k\leq l$. So the sum in (12.235) is bounded in $L^{\infty}(\Sigma_{t}^{\epsilon_{0}})$ by:
\begin{equation}
 C_{l}\delta_{0}(1+t)^{-2}[1+\log(1+t)]^{2}
\end{equation}
On the other hand, by $\textbf{H0}$ and (12.236), the first term on the right of (12.235) is bounded pointwise by:
\begin{equation}
 C\delta_{0}(1+t)^{-2}[1+\log(1+t)]\max_{j}|\leftexp{(i_{1}\overset{>i_{l+1-k}<}{...}i_{l+1}j)}{\mu}_{0,l+1}|
\end{equation}
Combining (12.237) and (12.238), the sum on the right of (12.229) is bounded by:
\begin{align}
 C_{l}\delta_{0}(1+t)^{-2}[1+\log(1+t)]\max_{j_{1}...j_{l}}|\leftexp{(j_{1}...j_{l+1})}{\mu}_{0,l+1}|
+C_{l}\delta_{0}(1+t)^{-2}[1+\log(1+t)]^{2}
\end{align}

    The estimates (12.232), (12.233), (12.234) and (12.239) imply through (12.229):
\begin{align}
 L(|\leftexp{(i_{1}...i_{l+1})}{\mu}_{0,l+1}|)\leq C_{l}\delta_{0}(1+t)^{-2}[1+\log(1+t)]\max_{j_{1}...j_{l+1}}
|\leftexp{(j_{1}...j_{l+1})}{\mu}_{0,l+1}|+C_{l}\delta_{0}(1+t)^{-1}
\end{align}
Defining, along a given integral curve of $L$, 
\begin{equation}
 \leftexp{(i_{1}...i_{l+1})}{x}^{\prime}_{0,l+1}(t)=|\leftexp{(i_{1}...i_{l+1})}{\mu}_{0,l+1}|
\end{equation}
(12.240) takes the form:
\begin{equation}
 \frac{d\leftexp{(i_{1}...i_{l+1})}{x}^{\prime}_{0,l+1}}{dt}\leq f_{l}\max_{j_{1}...j_{l+1}}\leftexp{(j_{1}...j_{l+1})}{x}^{\prime}_{0,l+1}+g_{l}
\end{equation}
where:
\begin{equation}
 f_{l}(t)=C_{l}\delta_{0}(1+t)^{-2}[1+\log(1+t)],\quad g_{l}(t)=C_{l}\delta_{0}(1+t)^{-1}
\end{equation}
This is similar to (12.217)-(12.218). Setting then:
\begin{equation}
 \bar{x}'_{0,l+1}(t)=\max_{i_{1}...i_{l+1}}\leftexp{(i_{1}...i_{l+1})}{x}^{\prime}_{0,l+1}(t)
\end{equation}
and using the initial condition:
\begin{align*}
 \bar{x}'_{0,l+1}(0)\leq C_{l}\delta_{0}
\end{align*}
we deduce:
\begin{equation}
 \max_{i_{1}...i_{l+1}}\|\leftexp{(i_{1}...i_{l+1})}{\mu}_{0,l+1}\|_{L^{\infty}(\Sigma_{t}^{\epsilon_{0}})}
\leq C_{l}\delta_{0}[1+\log(1+t)]\quad :\textrm{for all}\quad t\in[0,s]
\end{equation}
This proves the proposition in the case of $m=0$.

     To prove the proposition for $m\in\{1,...,l+1\}$ we apply induction on $m$. We assume $\textbf{M}_{[m,l+1]}$ for some 
$m\in\{0,...,1\}$, and we will establish $\textbf{M}_{[m+1,l+1]}$ under the assumption:
\begin{equation}
 \|\mu-1\|_{\infty,[m+1,l+1],\Sigma_{0}^{\epsilon_{0}}}\leq C_{l}\delta_{0}
\end{equation}
     We first apply $T^{m+1}$ to $L\mu$ to obtain, using Lemma 11.22,
\begin{equation}
 L(T)^{m+1}\mu=T^{m+1}L\mu+\sum_{k=0}^{m}(T)^{k}\Lambda(T)^{m-k}\mu
\end{equation}
We then apply $R_{i_{n}}...R_{i_{1}}$ to this equation, to obtain, using again Lemma 11.22,
\begin{align}
 LR_{i_{n}}...R_{i_{1}}(T)^{m+1}\mu =R_{i_{n}}...R_{i_{1}}(T)^{m+1}L\mu +\sum_{k=0}^{m}R_{i_{n}}...R_{i_{1}}(T)^{k}\Lambda(T)^{m-k}\mu\\\notag
+\sum_{k=0}^{n-1}R_{i_{n}}...R_{i_{n-k+1}}\leftexp{(R_{i_{n-k}})}{Z}R_{i_{n-k-1}}...R_{i_{1}}(T)^{m+1}\mu
\end{align}
Let us define as in Chapter 9:
\begin{equation}
 \leftexp{(i_{1}...i_{n})}{\mu}_{m+1,n}=R_{i_{n}}...R_{i_{1}}(T)^{m+1}\mu
\end{equation}
Applying $R_{i_{n}}...R_{i_{1}}(T)^{m+1}$ to equation (12.226), we obtain:
\begin{align}
 R_{i_{n}}...R_{i_{1}}(T)^{m+1}L\mu=e\leftexp{(i_{1}...i_{n})}{\mu}_{m+1,n}+R_{i_{n}}...R_{i_{1}}(T)^{m+1}m
+\leftexp{(i_{1}...i_{n})}{r}^{\prime}_{m+1,n}
\end{align}
where
\begin{align}
 \leftexp{(i_{1}...i_{n})}{r}^{\prime}_{m+1,n}=\sum_{|s_{1}|+|s_{2}|=n, s_{1}\geq 1}((R)^{s_{1}}e)((R)^{s_{2}}(T)^{m+1}\mu)\\\notag
+R_{i_{n}}...R_{i_{1}}(\sum_{k=1}^{m+1}\frac{(m+1)!}{k!(m+1-k)!}((T)^{k}e)((T)^{m+1-k}\mu))
\end{align}
Substituting (12.250) in (12.248) we get:
\begin{align}
 L\leftexp{(i_{1}...i_{n})}{\mu}_{m+1,n}=e\leftexp{(i_{1}...i_{n})}{\mu}_{m+1,n}+R_{i_{n}}...R_{i_{1}}(T)^{m+1}m
+\leftexp{(i_{1}...i_{n})}{r}^{\prime}_{m+1,n}\\\notag
+\sum_{k=0}^{m}R_{i_{n}}...R_{i_{1}}(T)^{k}\Lambda(T)^{m-k}\mu\\\notag
+\sum_{k=0}^{n-1}R_{i_{n}}...R_{i_{n-k+1}}\leftexp{(R_{i_{n-k}})}{Z}R_{i_{n-k-1}}...R_{i_{1}}(T)^{m+1}\mu
\end{align}
By assumptions $\textbf{E}_{\{l+2\}}, \textbf{E}_{\{l+1\}}^{Q}$:
\begin{align}
 \|Th\|_{\infty,[m+1,l+1],\Sigma_{t}^{\epsilon_{0}}}\leq C_{l}\delta_{0}(1+t)^{-1}\\\notag
\|Lh\|_{\infty,[m+1,l+1],\Sigma_{t}^{\epsilon_{0}}}\leq C_{l}\delta_{0}(1+t)^{-2}\quad:\textrm{for all}\quad t\in[0,s]
\end{align}
By Corollary 11.1.b, we have:
\begin{equation}
 \|\omega_{L\hat{T}}\|_{\infty,[m+1,l+1],\Sigma_{t}^{\epsilon_{0}}}\leq C_{l}\delta_{0}(1+t)^{-2}\quad:\textrm{for all}\quad t\in[0,s]
\end{equation}
The bounds (12.253) and (12.254) imply:
\begin{equation}
 \|e\|_{\infty,[m+1,l+1],\Sigma_{t}^{\epsilon_{0}}}\leq C_{l}\delta_{0}(1+t)^{-2}\quad:\textrm{for all}\quad t\in[0,s]
\end{equation}
and (12.253) implies:
\begin{equation}
 \|m\|_{\infty,[m+1,l+1],\Sigma_{t}^{\epsilon_{0}}}\leq C_{l}\delta_{0}(1+t)^{-1}\quad:\textrm{for all}\quad t\in[0,s]
\end{equation}
Also, (12.255) together with $\textbf{M}_{[m,l+1]}$ imply:
\begin{equation}
 \|\leftexp{(i_{1}...i_{n})}{r}^{\prime}_{m+1,n}\|_{L^{\infty}(\Sigma_{t}^{\epsilon_{0}})}\leq C_{l}\delta_{0}(1+t)^{-2}[1+\log(1+t)]
\quad:\textrm{for}\quad n\leq l-m
\end{equation}
provided that either $n=0$, so that the first sum on the right in (12.251) vanishes, or $\textbf{M}_{m+1,m+n}$ holds. Moreover, by Corollary 11.1.c:
\begin{equation}
 \|\Lambda\|_{\infty,[m,l],\Sigma_{t}^{\epsilon_{0}}}\leq C_{l}\delta_{0}(1+t)^{-1}[1+\log(1+t)]\quad:\textrm{for all}\quad t\in[0,s]
\end{equation}

    To establish $\textbf{M}_{[m+1,l+1]}$, we establish $\textbf{M}_{[m+1,m+1+n]}$ for $n=0,...,l-m$ by induction on $n$. We begin
with the case $n=0$. In this case, equation (12.252) reduces to:
\begin{equation}
 L\mu_{m+1,0}=e\mu_{m+1,0}+(T)^{m+1}m+r'_{m+1,0}+\sum_{k=0}^{m}(T)^{k}\Lambda(T)^{m-k}\mu
\end{equation}
Using (12.258) with $l=m$, and $\textbf{M}_{[m,m+1]}$, the sum in (12.259) is bounded in $L^{\infty}(\Sigma_{t}^{\epsilon_{0}})$ by:
\begin{align}
 \sum_{k=0}^{m}\|\Lambda\cdot\slashed{d}(T)^{m-k}\mu\|_{\infty,[k,k],\Sigma_{t}^{\epsilon_{0}}}\leq C_{m}
(1+t)^{-1}\|\Lambda\|_{\infty,[m,m],\Sigma_{t}^{\epsilon_{0}}}\|\mu-1\|_{\infty,[m,m+1],\Sigma_{t}^{\epsilon_{0}}}\\\notag
\leq C_{m}\delta_{0}^{2}(1+t)^{-2}[1+\log(1+t)]^{2}
\end{align}
Appealing also to (12.48) and (12.256) with $l=m$, (12.257) with $n=0$, we deduce:
\begin{equation}
 L(|\mu_{m+1,0}|)\leq C\delta_{0}(1+t)^{-2}|\mu_{m+1,0}|+C_{m}\delta_{0}(1+t)^{-1}
\end{equation}
which together with 
\begin{equation}
 \|\mu-1\|_{\infty,[m+1,m+1],\Sigma_{0}^{\epsilon_{0}}}\leq C_{m}\delta_{0}
\end{equation}
yields:
\begin{equation}
 \|\mu_{m+1,0}\|_{L^{\infty}(\Sigma_{t}^{\epsilon_{0}})}\leq C_{m}\delta_{0}[1+\log(1+t)]
\quad:\textrm{for all}\quad t\in[0,s]
\end{equation}
This estimate together with $\textbf{M}_{[m,m+1]}$ yields $\textbf{M}_{[m+1,m+1]}$.

    Next, let $\textbf{M}_{[m+1,m+n]}$ hold for some $n\in\{1,...,l-m\}$. We shall show that $\textbf{M}_{[m+1,m+1+n]}$ holds. 
Consider the first sum on the right in (12.252). Since $n\leq l-m$ this is bounded in $L^{\infty}(\Sigma_{t}^{\epsilon_{0}})$ by:
\begin{align}
 \sum_{k=0}^{m}\|\Lambda\cdot\slashed{d}(T)^{m-k}\mu\|_{\infty,[k,l+k-m],\Sigma_{t}^{\epsilon_{0}}}
\leq C_{l}(1+t)^{-1}\|\Lambda\|_{\infty,[m,l],\Sigma_{t}^{\epsilon_{0}}}\|\mu-1\|_{\infty,[m,l+1],\Sigma_{t}^{\epsilon_{0}}}\\\notag
\leq C_{l}\delta_{0}^{2}(1+t)^{-2}[1+\log(1+t)]^{2}
\end{align}
by (12.258) and $\textbf{M}_{[m,l+1]}$.

    Consider next the second sum on the right of (12.252). Consider a term in this sum corresponding to $k\in\{0,...,n-1\}$. We express the term as:
\begin{align}
 R_{i_{n}}...R_{i_{n-k+1}}(\leftexp{(R_{i_{n-k}})}{Z}\cdot\slashed{d}R_{i_{n-k-1}}...R_{i_{1}}(T)^{m+1}\mu)=
\leftexp{(R_{i_{n-k}})}{Z}\cdot\slashed{d}\leftexp{(i_{1}\overset{>i_{n-k}<}{...}i_{n})}{\mu}_{m+1,n-1}\\\notag
+\sum_{|s_{1}|+|s_{2}|=k,|s_{1}|\geq 1}((\slashed{\mathcal{L}}_{R})^{s_{1}}\leftexp{(R_{i_{n-k}})}{Z})\cdot
\slashed{d}(R)^{s'_{2}}(T)^{m+1}\mu 
\end{align}
where the sum on the right is over all ordered partitions $\{s_{1},s_{2}\}$ of the set $\{n-k+1,...,n\}$ into two ordered subsets $s_{1}, s_{2}$,
with $s_{1}$ non-empty. Also 
\begin{align*}
 s'_{2}=s_{2}\bigcup\{1,...,n-k-1\}
\end{align*}
Since $|s_{1}|\leq k\leq n-1\leq l-m-1$, the first factor in the sum is bounded by (12.236). Also, since $|s_{2}|\leq k-1$, we have:
\begin{equation}
 |s'_{2}|=|s_{2}|+n-k-1\leq n-2\quad\textrm{and}\quad 1+|s'_{2}|+m+1\leq n+m\notag
\end{equation}
 Therefore, using $\textbf{H0}$ and $\textbf{M}_{m+1,m+n}$ the sum is seen to be bounded by:
 \begin{equation}
  C_{l}\delta_{0}^{2}(1+t)^{-2}[1+\log(1+t)]^{2}
 \end{equation}
On the other hand, by (12.236) and $\textbf{H0}$ the first term on the right of (12.265) is bounded by:
\begin{equation}
 C_{l}\delta_{0}(1+t)^{-2}[1+\log(1+t)]\max_{j}|\leftexp{(i_{1}\overset{>i_{n-k}<}{...}i_{n}j)}{\mu}_{m+1,n}|
\end{equation}
Combining (12.266) and (12.267) we conclude that the second sum on the right of (12.252) is bounded by:
\begin{align}
 C_{l}\delta_{0}(1+t)^{-2}[1+\log(1+t)]\max_{j}|\leftexp{(i_{1}\overset{>i_{n-k}<}{...}i_{n}j)}{\mu}_{m+1,n}|
+C_{l}\delta_{0}^{2}(1+t)^{-2}[1+\log(1+t)]^{2}
\end{align}

     The estimates (12.255)-(12.257), (12.264) and (12.268) imply through (12.252):
\begin{align}
 L(|\leftexp{(i_{1}...i_{l})}{\mu}_{m+1,n}|)\leq C_{l}\delta_{0}(1+t)^{-2}[1+\log(1+t)]\max_{j_{1}...j_{n}}
|\leftexp{(j_{1}...j_{n})}{\mu}_{m+1,n}|+C_{l}\delta_{0}(1+t)^{-1}
\end{align}
 Setting along an integral curve of $L$,
\begin{equation}
 \leftexp{(i_{1}...i_{n})}{x}^{\prime}_{m+1,n}(t)=|\leftexp{(i_{1}...i_{n})}{\mu}_{m+1,n}|
\end{equation}
(12.269) takes the form:
\begin{equation}
 \frac{d\leftexp{(i_{1}...i_{n})}{x}^{\prime}_{m+1,n}}{dt}\leq f_{l}\max_{j_{1}...j_{n}}\leftexp{(j_{1}...j_{n})}{x}^{\prime}_{m+1,n}+g_{l}
\end{equation}
where:
\begin{equation}
 f_{l}(t)=C_{l}\delta_{0}(1+t)^{-2}[1+\log(1+t)],\quad g_{l}(t)=C_{l}\delta_{0}(1+t)^{-1}
\end{equation}
This is identical in form to (12.217)-(12.218). Setting then:
\begin{equation}
 \bar{x}'_{m+1,n}(t)=\max_{i_{1}...i_{n}}\leftexp{(i_{1}...i_{n})}{x}^{\prime}_{m+1,n}(t)
\end{equation}
in view of the fact that
\begin{align*}
 \bar{x}'_{m+1,n}(0)\leq C_{l}\delta_{0}
\end{align*}
we conclude that:
\begin{equation}
 \max_{i_{1}...i_{n}}\|\leftexp{(i_{1}...i_{n})}{\mu}_{m+1,n}\|_{L^{\infty}(\Sigma_{t}^{\epsilon_{0}})}\leq C_{l}
\delta_{0}[1+\log(1+t)]\quad:\textrm{for all}\quad t\in[0,s]
\end{equation}
This together with $\textbf{M}_{[m+1,m+n]}$ yields $\textbf{M}_{[m+1,m+1+n]}$, so the proposition 
follows. $\qed$ \vspace{7mm}

 Let us now define, for any non-negative integer $k$, the quantities:
\begin{equation}
 \mathcal{A}'_{k}=\max_{i_{1}...i_{k}}\|\slashed{\mathcal{L}}_{R_{i_{k}}}...\slashed{\mathcal{L}}_{R_{i_{1}}}\chi'\|
_{L^{2}(\Sigma_{t}^{\epsilon_{0}})}
\end{equation}
and, for any non-negative integer $l$, the quantities:
\begin{equation}
 \mathcal{A}'_{[l]}=\sum_{k=0}^{l}\mathcal{A}'_{k}
\end{equation}
Let us recall the definitions:
\begin{align}
 \mathcal{A}_{0}=\|\chi-\frac{\slashed{g}}{1-u+t}\|_{L^{2}(\Sigma_{t}^{\epsilon_{0}})}
\end{align}
and for $k>0$:
\begin{align}
 \mathcal{A}_{k}=\max_{i_{1}...i_{k}}\|\slashed{\mathcal{L}}_{R_{i_{k}}}...\slashed{\mathcal{L}}_{R_{i_{1}}}\chi\|_{L^{2}(\Sigma_{t}^{\epsilon_{0}})}
\end{align}
and:
\begin{equation}
 \mathcal{A}_{[l]}=\sum_{k=0}^{l}\mathcal{A}_{k}
\end{equation}
Let us also recall:
\begin{equation}
 \mathcal{W}_{0}=\max_{\alpha}\{\|\psi_{\alpha}\|_{L^{2}(\Sigma_{t}^{\epsilon_{0}})}\}
\end{equation}
and for $k>0$:
\begin{equation}
 \mathcal{W}_{k}=\max_{\alpha;i_{1}...i_{k}}\|R_{i_{k}}...R_{i_{1}}\psi_{\alpha}\|_{L^{2}(\Sigma_{t}^{\epsilon_{0}})}
\end{equation}
and:
\begin{equation}
 \mathcal{W}_{[l]}=\sum_{k=0}^{l}\mathcal{W}_{k}
\end{equation}
\begin{equation}
 \mathcal{W}^{Q}_{k}=\max_{\alpha;i_{1}...i_{k}}\|R_{i_{k}}...R_{i_{1}}Q\psi_{\alpha}\|_{L^{2}(\Sigma_{t}^{\epsilon_{0}})},\quad
\mathcal{W}^{Q}_{[l]}=\sum_{k=0}^{l}\mathcal{W}^{Q}_{k}
\end{equation}
For all $k=0,...,l$:
\begin{equation}
 \mathcal{Y}_{k}=\max_{j;i_{1}...i_{k}}\|R_{i_{k}}...R_{i_{1}}y^{i}\|_{L^{2}(\Sigma_{t}^{\epsilon_{0}})},\quad
\mathcal{Y}_{[l]}=\sum_{k=0}^{l}\mathcal{Y}_{k}
\end{equation}
Recall also that given a positive integer $n$ we denote:
\begin{align}
 n_{*}=\frac{n}{2}\quad: \textrm{if n is even}\\\notag
=\frac{n-1}{2}\quad: \textrm{if n is odd}
\end{align}
$\textbf{Lemma 12.5}$ Let $\textbf{H0}$ and (12.23) hold. Let $l$ be a positive integer and let the assumptions 
$\slashed{\textbf{E}}_{[l_{*}]}, \slashed{\textbf{E}}_{[l_{*}-1]}^{Q}$ as well as $\slashed{\textbf{X}}'_{[l_{*}-1]}$ hold. Then we have:
\begin{align*}
 \mathcal{A}_{[l]}\leq \mathcal{A}'_{[l]}+C_{l}\mathcal{A}'_{[l-1]}+C_{l}(1+t)^{-1}
[\mathcal{Y}_{0}+\mathcal{W}_{[l]}]
\end{align*}

$Proof$: From (12.185) we obtain, for $k\geq 1$:
\begin{align}
 \mathcal{A}_{k}\leq \mathcal{A}'_{k}+\frac{1}{1-u+t}\max_{i}\max_{j_{1}...j_{k-1}}\|
\slashed{\mathcal{L}}_{R_{j_{k-1}}}...\slashed{\mathcal{L}}_{R_{j_{1}}}\leftexp{(R_{i})}{\slashed{\pi}}\|_{L^{2}(\Sigma_{t}^{\epsilon_{0}})}
\end{align}
Summing over $k\in\{1,...,l\}$ and adding:
\begin{equation}
 \mathcal{A}_{0}=\mathcal{A}'_{0}
\end{equation}
we deduce, for $l\geq 1$:
\begin{equation}
 \mathcal{A}_{[l]}\leq \mathcal{A}'_{[l]}+\frac{3}{1-u+t}\max_{i}\|\leftexp{(R_{i})}{\slashed{\pi}}\|_{2,[l-1],\Sigma_{t}^{\epsilon_{0}}}
\end{equation}
By Lemma 12.4 with $l_{*}$ in the role of $l$, the assumptions of the lemma imply $\slashed{\textbf{X}}_{[l_{*}-1]}$, Therefore Proposition 10.2 holds.
Then Corollary 10.2.d implies:
\begin{align}
 \max_{i}\|\leftexp{(R_{i})}{\slashed{\pi}}\|_{2,[l-1],\Sigma_{t}^{\epsilon_{0}}}\leq 
C_{l}[\mathcal{Y}_{0}+(1+t)\mathcal{A}_{[l-1]}+\mathcal{W}_{[l]}]
\end{align}
Substituting this in (12.288) we get:
\begin{align}
 \mathcal{A}_{[l]}\leq \mathcal{A}'_{[l]}+C_{l}\mathcal{A}_{[l-1]}+C_{l}(1+t)^{-1}[\mathcal{Y}_{0}+\mathcal{W}_{[l]}]
\quad :\textrm{for all} \quad l\geq 1\\\notag
\mathcal{A}_{[0]}=\mathcal{A}'_{[0]}
\end{align}
The lemma follows. $\qed$

$\textbf{Proposition 12.11}$ Let the assumptions of Lemma 12.4 hold. Let $l$ be a non-negative integer and let the assumptions $\textbf{E}_{\{l_{*}+2\}},
\textbf{E}_{\{l_{*}+1\}}^{Q},\textbf{E}_{\{l_{*}\}}^{QQ}$, hold on $W^{s}_{\epsilon_{0}}$. Moreover let the initial data satisfy:
\begin{align*}
 \|\chi'\|_{\infty,[l_{*}],\Sigma_{0}^{\epsilon_{0}}}\leq C_{l}\delta_{0}
\end{align*}
Then $\slashed{\textbf{X}}'_{[l_{*}]}$ holds on $W^{s}_{\epsilon_{0}}$ and there is a constant $C_{l}$ independent of $s$ such that for all $t\in[0,s]$ 
we have:
\begin{align*}
 \mathcal{A}'_{[l]}(t)\leq C_{l}(1+t)^{-1}\{\mathcal{A}'_{[l]}(0)+\delta_{0}\mathcal{Y}_{0}(0)\\
+\int_{0}^{t}(1+t')^{-1}[\mathcal{W}_{[l+2]}(t')+\mathcal{W}^{Q}_{[l+1]}(t')]dt'\}
\end{align*}
$Proof$: The assumptions of this proposition include those of Proposition 12.9 with $l_{*}$ in the role of $l$, therefore $\slashed{\textbf{X}}'_{[l_{*}]}$ 
holds on $W^{s}_{\epsilon_{0}}$. Then from Lemma 12.4, $\slashed{\textbf{X}}_{[l_{*}]}$ holds on $W^{s}_{\epsilon_{0}}$.
Thus, Proposition 10.1  holds with $l_{*}+1$ in the role of $l$, then Proposition 10.2 holds with $l+1$ in the role of $l$.

    We shall estimate $\leftexp{(i_{1}...i_{l})}{b}_{l}$ (see (12.182)) in $L^{2}(\Sigma_{t}^{\epsilon_{0}})$. By assumptions $\textbf{E}^{Q}_{\{l_{*}\}},
\textbf{E}_{\{l_{*}+1\}}$, we have:
\begin{align}
 \|L\psi_{\mu}\|_{\infty,[l_{*}],\Sigma_{t}^{\epsilon_{0}}},\quad \|\slashed{d}\psi_{\mu}\|_{\infty,[l_{*}],\Sigma_{t}^{\epsilon_{0}}}\leq C_{l}\delta_{0}(1+t)^{-2}
\quad:\textrm{for all}\quad t\in[0,s]
\end{align}
as well as
\begin{align}
 \|Lh\|_{\infty,[l_{*}],\Sigma_{t}^{\epsilon_{0}}},\quad \|\slashed{d}h\|_{\infty,[l_{*}],\Sigma_{t}^{\epsilon_{0}}}\leq C_{l}\delta_{0}(1+t)^{-2}
\quad:\textrm{for all}\quad t\in[0,s]
\end{align}
 So by Corollary 10.1.b we have:
\begin{align}
 \|e\|_{\infty,[l_{*}],\Sigma_{t}^{\epsilon_{0}}},\quad \|\kappa^{-1}\zeta\|_{\infty,[l_{*}],\Sigma_{t}^{\epsilon_{0}}}\leq C_{l}\delta_{0}(1+t)^{-2}
\quad:\textrm{for all}\quad t\in[0,s]
\end{align}
Obviously, by the assumptions of the present proposition we have:
\begin{align}
 \|Lh\|_{2,[l],\Sigma_{t}^{\epsilon_{0}}}\leq C_{l}(1+t)^{-1}[\mathcal{W}^{Q}_{[l]}+\delta_{0}(1+t)^{-1}\mathcal{W}_{[l]}]\\\notag
\|\slashed{d}h\|_{2,[l],\Sigma_{t}^{\epsilon_{0}}}\leq C_{l}(1+t)^{-1}\mathcal{W}_{[l+1]}\quad:\textrm{for all}\quad t\in[0,s]
\end{align}
so by Corollary 10.1.g and Corollary 10.2.g with $l$ replaced by $l+1$,  we have:
\begin{align}
 \|e\|_{2,[l],\Sigma_{t}^{\epsilon_{0}}}\leq C_{l}(1+t)^{-1}\{\mathcal{W}^{Q}_{[l]}+\delta_{0}(1+t)^{-1}[\mathcal{Y}_{0}+(1+t)\mathcal{A}_{[l-1]}
+\mathcal{W}_{[l]}]\}\quad:\textrm{for all}\quad t\in[0,s]
\end{align}
Also from the expression 
\begin{align*}
 \kappa^{-1}\zeta=\alpha\epsilon-\slashed{d}\alpha
\end{align*}
and Corollary 10.2.h,
we have:
\begin{align}
 \|\kappa^{-1}\zeta\|_{2,[l],\Sigma_{t}^{\epsilon_{0}}}\leq C_{l}(1+t)^{-1}\{\mathcal{W}_{[l+1]}\\\notag
+\delta_{0}(1+t)^{-1}[\mathcal{Y}_{0}+(1+t)\mathcal{A}_{[l-1]}]\}\quad:\textrm{for all}\quad t\in[0,s]
\end{align}
From (12.45), using the expression for $L^{i}$:
\begin{align*}
 L^{i}=-\eta\hat{T}^{i}+\psi_{i}
\end{align*}
and Corollary 10.2.g, we have:
\begin{align}
 \|\alpha'\|_{2,[l],\Sigma_{t}^{\epsilon_{0}}}\leq C_{l}\delta_{0}(1+t)^{-2}[\mathcal{W}_{[l+2]}+\delta_{0}(1+t)^{-1}\mathcal{W}^{Q}_{[l+1]}+(1+t)^{-1}\mathcal{Y}_{0}
+\mathcal{A}_{[l]}]
\end{align}
where we have used Lemma 10.8 and Lemma 10.11 to estimate the commutator as well as $\textbf{H0}, \textbf{H1}, \textbf{H2}$ to estimate the $S_{t,u}$
-tangential derivatives.
This and (12.295), together with the estimates for $\leftexp{(R_{i})}{\slashed{\pi}}$ in Chapter 10 (10.145):
\begin{align}
 \|\leftexp{(R_{i})}{\slashed{\pi}}\|_{\infty,[l],\Sigma_{t}^{\epsilon_{0}}}\leq C_{l}\delta_{0}(1+t)^{-1}[1+\log(1+t)]\\\notag
\|\leftexp{(R_{i})}{\slashed{\pi}}\|_{2,[l],\Sigma_{t}^{\epsilon_{0}}}\leq C_{l}\{\mathcal{Y}_{0}+(1+t)\mathcal{A}_{[l]}+\mathcal{W}_{[l+1]}\}
\end{align}
yield the following bound for $b$ defined by (12.183):
\begin{align}
 \|b\|_{2,[l],\Sigma_{t}^{\epsilon_{0}}}\leq C_{l}\delta_{0}(1+t)^{-2}[\mathcal{W}_{[l+2]}+(1+t)^{-1}\mathcal{Y}_{0}+\mathcal{A}_{[l]}+\mathcal{W}^{Q}_{[l+1]}]
\end{align}

    Next, we shall estimate $\leftexp{(i_{1}...i_{l})}{s}_{l}$ and $\leftexp{(i_{1}...i_{l})}{r}_{l}$ in $L^{2}(\Sigma_{t}^{\epsilon_{0}})$. 
From the expression (12.178), bounds for $\leftexp{(R_{i})}{\slashed{\pi}}$ (12.298) as well as the fact that $\slashed{\textbf{X}}'_{[l_{*}]}$ holds, 
we get:
\begin{align}
 \|\leftexp{(i_{1}...i_{l})}{r}_{l}\|_{L^{2}(\Sigma_{t}^{\epsilon_{0}})}\leq C_{l}\delta_{0}(1+t)^{-2}[1+\log(1+t)]\mathcal{A}'_{[l-1]}\\\notag
+C_{l}\delta_{0}^{2}(1+t)^{-4}[1+\log(1+t)]^{2}[\mathcal{Y}_{0}+\mathcal{W}_{[l]}]
\end{align}
Also from (12.180) and the estimates for $e$ we get:
\begin{align}
 \|\leftexp{(i_{1}...i_{l})}{s}_{l}\|_{L^{2}(\Sigma_{t}^{\epsilon_{0}})}\leq C_{l}(1+t)^{-2}\mathcal{A}'_{[l-1]}\\\notag
+C_{l}\delta_{0}(1+t)^{-3}[1+\log(1+t)]\{\mathcal{W}^{Q}_{[l]}+\delta_{0}(1+t)^{-1}[\mathcal{Y}_{0}+\mathcal{W}_{[l]}]\}
\end{align}
Combining (12.299)-(12.301) and Lemma 12.5, we get:
\begin{align}
 \|\leftexp{(i_{1}...i_{l})}{b}_{l}\|_{L^{2}(\Sigma_{t}^{\epsilon_{0}})}\leq C_{l}\delta_{0}\{
(1+t)^{-2}[1+\log(1+t)]\mathcal{A}'_{[l]}\\\notag
+(1+t)^{-2}[\mathcal{W}_{[l+2]}+(1+t)^{-1}\mathcal{Y}_{0}+\mathcal{W}^{Q}_{[l+1]}]\}
\end{align}

Next, we shall estimate the $L^{2}(\Sigma_{t}^{\epsilon_{0}})$ norm of the sum on the right of (12.176). A term $k\in\{0,...,l-1\}$ in this sum is expressed by
(12.195). For the first term on the right of (12.195), we have the pointwise estimate (12.197). This implies:
\begin{align}
 \|\slashed{\mathcal{L}}_{\leftexp{(R_{i_{l-k}})}{Z}}\leftexp{(i_{1}\overset{>i_{l-k}<}{...}i_{l})}{\chi}^{\prime}_{l-1}\|_{L^{2}(\Sigma_{t}^{\epsilon_{0}})}
\leq C(1+t)^{-1}\cdot\\\notag
\{\|\leftexp{(R_{i_{l-k}})}{Z}\|_{L^{\infty}(\Sigma_{t}^{\epsilon_{0}})}(3\max_{j}\|\leftexp{(i_{1}\overset{i_{l-k}}{...}i_{l}j)}{\chi}^{\prime}_{l}\|
_{L^{2}(\Sigma_{t}^{\epsilon_{0}})}+\|\leftexp{(i_{1}\overset{i_{l-k}}{...}i_{l})}{\chi}^{\prime}_{l-1}\|_{L^{2}(\Sigma_{t}^{\epsilon_{0}})})\\\notag
+\|\leftexp{(i_{1}\overset{i_{l-k}}{...}i_{l})}{\chi}^{\prime}_{l-1}\|_{L^{2}(\Sigma_{t}^{\epsilon_{0}})}\max_{j}
\|\leftexp{(j;i_{l-k})}{Z}\|_{L^{\infty}(\Sigma_{t}^{\epsilon_{0}})}\}
\end{align}
If $l\geq 2$, so that $l_{*}\geq 1$, we apply Corollary 10.1.i with $l_{*}+1$ in the role of $l$ to get:
\begin{align}
 \|\leftexp{(R_{i_{l-k}})}{Z}\|_{L^{\infty}(\Sigma_{t}^{\epsilon_{0}})},\quad \max_{j}
\|\leftexp{(j;i_{l-k})}{Z}\|_{L^{\infty}(\Sigma_{t}^{\epsilon_{0}})}\leq C\delta_{0}(1+t)^{-1}[1+\log(1+t)]
\end{align}
We then obtain:
\begin{align}
 \|\slashed{\mathcal{L}}_{\leftexp{(R_{i_{l-k}})}{Z}}\leftexp{(i_{1}\overset{i_{l-k}}{...}i_{l})}{\chi'}_{l-1}\|_{L^{2}(\Sigma_{t}^{\epsilon_{0}})}
\leq C\delta_{0}(1+t)^{-2}[1+\log(1+t)]\mathcal{A}'_{[l]}\quad :\quad l\geq 2
\end{align}
On the other hand, if $l=1$, then $l_{*}=0$, we estimate the last term in parenthesis in (12.197) in $L^{2}$ by:
\begin{equation}
 3\|\chi'\|_{L^{\infty}(\Sigma_{t}^{\epsilon_{0}})}\max_{j}\|\leftexp{(j;i_{1})}{Z}\|_{L^{2}(\Sigma_{t}^{\epsilon_{0}})}
\end{equation}
By Corollary 10.2.i with $l+1$ in the role of $l$, we have:
\begin{align}
 \max_{i}\|\leftexp{(R_{i})}{Z}\|_{2,[l],\Sigma_{t}^{\epsilon_{0}}}\leq C_{l}\{\mathcal{Y}_{0}+(1+t)\mathcal{A}'_{[l]}
+\mathcal{W}_{[l+1]}\}
\end{align}
where we have used Lemma 12.5 to express $\mathcal{A}_{[l]}$ in terms of $\mathcal{A}'_{[l]}$.
Then in view of Proposition 12.6 and (12.307) with $l=1$, (12.306) is bounded by
\begin{equation}
 C\delta_{0}(1+t)^{-2}[1+\log(1+t)]\{\mathcal{Y}_{0}+(1+t)\mathcal{A}'_{[1]}
+\mathcal{W}_{[2]}\}
\end{equation}
Thus, in the case of $l=1$, we have:
\begin{align}
 \|\slashed{\mathcal{L}}_{\leftexp{(R_{i})}{Z}}\chi'\|_{L^{2}(\Sigma_{t}^{\epsilon_{0}})}\leq C\delta_{0}
(1+t)^{-2}[1+\log(1+t)]\mathcal{A}'_{[1]}\\\notag
+C\delta_{0}(1+t)^{-3}[1+\log(1+t)][\mathcal{Y}_{0}+\mathcal{W}_{[2]}]
\end{align}

    Next, we should consider a term in the double sum in (12.195), corresponding to an ordered subset 
$\{m_{1},...,m_{p}\}\subset\{l-k+1,...,l\}$, $p\in\{1,...,k\}$. We have the pointwise estimate (12.204).

    We distinguish three cases:
\begin{align*}
 \textrm{Case 1}:\quad p\leq l_{*}-1,\quad\textrm{Case 2}:\quad p=l_{*},\quad \textrm{Case 3}:\quad p\geq l_{*}+1
\end{align*}

    In Case 1, Corollary 10.1.i with $l_{*}+1$ in the role of $l$ implies:
\begin{align}
 \|\leftexp{(i_{m_{1}}...i_{m_{p}};i_{l-k})}{Z}\|_{L^{\infty}(\Sigma_{t}^{\epsilon_{0}})},\quad 
\max_{j}\|\leftexp{(i_{m_{1}}...i_{m_{p}}j;i_{l-k})}{Z}\|_{L^{\infty}(\Sigma_{t}^{\epsilon_{0}})}\\\notag
\leq C_{l}\delta_{0}(1+t)^{-1}[1+\log(1+t)]
\end{align}
then (12.204) implies:
\begin{align}
 \|\slashed{\mathcal{L}}_{\leftexp{(i_{m_{1}}...i_{m_{p}};i_{l-k})}{Z}}\leftexp{(i_{1}\overset{>i_{m_{1}}...i_{m_{p}}i_{l-k}<}{...}i_{l})}{\chi}^{\prime}_{l-1-p}\|
_{L^{2}(\Sigma_{t}^{\epsilon_{0}})}\leq C(1+t)^{-1}\cdot\\\notag
\{\|\leftexp{(i_{m_{1}}...i_{m_{p}};i_{l-k})}{Z}\|_{L^{\infty}(\Sigma_{t}^{\epsilon_{0}})}(3\max_{j}\|\leftexp{(i_{1}
\overset{>i_{m_{1}}...i_{m_{p}}i_{l-k}<}{...}i_{l}j)}{\chi}^{\prime}_{l-p}\|_{L^{2}(\Sigma_{t}^{\epsilon_{0}})}\\\notag
\|\leftexp{(i_{1}\overset{>i_{m_{1}}...i_{m_{p}}i_{l-k}<}{...}i_{l})}{\chi}^{\prime}_{l-1-p}\|_{L^{2}(\Sigma_{t}^{\epsilon_{0}})})\\\notag
+\max_{j}\|\leftexp{(i_{m_{1}}...i_{m_{p}}j;i_{l-k})}{Z}\|_{L^{\infty}(\Sigma_{t}^{\epsilon_{0}})}
\|\leftexp{(i_{1}\overset{>i_{m_{1}}...i_{m_{p}}i_{l-k}<}{...}i_{l})}{\chi}^{\prime}_{l-1-p}\|_{L^{2}(\Sigma_{t}^{\epsilon_{0}})}\}\\\notag
\leq C_{l}\delta_{0}(1+t)^{-2}[1+\log(1+t)]\mathcal{A}^{\prime}_{[l-1]}
\end{align}
where we have used the fact that $p\geq 1$.

     In Case 2, Corollary 10.1.i with $l_{*}+1$ in the role of $l$ implies the first of (12.310), but not the second. We use instead Corollary 10.2.i 
with $l+1$ in the role of $l$, to get:
\begin{align}
 \max_{j}\|\leftexp{(i_{m_{1}}...i_{m_{p}}j;i_{l-k})}{Z}\|_{L^{2}(\Sigma_{t}^{\epsilon_{0}})}\leq C_{l}\{\mathcal{Y}_{0}
+(1+t)\mathcal{A}'_{[l]}+\mathcal{W}_{[l+1]}\}
\end{align}
where we have used the fact that $p\leq k\leq l-1$.

Since in Case 2 we have $l-1-p=l-1-l_{*}\leq l_{*}$, we may use $\slashed{\textbf{X}}'_{[l_{*}]}$ to estimate:
\begin{align}
 \|\leftexp{(i_{1}\overset{>i_{m_{1}}...i_{m_{p}}i_{l-k}<}{...}i_{l})}{\chi}^{\prime}_{l-1-p}\|_{L^{\infty}(\Sigma_{t}^{\epsilon_{0}})}
\leq C_{l}\delta_{0}(1+t)^{-2}[1+\log(1+t)]
\end{align}
In view of (12.312) and (12.313) we get in Case 2:
\begin{align}
 \|\slashed{\mathcal{L}}_{\leftexp{(i_{m_{1}}...i_{m_{p}};i_{l-k})}{Z}}\leftexp{(i_{1}\overset{>i_{m_{1}}...i_{m_{p}}i_{l-k}<}{...}i_{l})}{\chi'_{l-1-p}}\|
_{L^{2}(\Sigma_{t}^{\epsilon_{0}})}\leq C(1+t)^{-1}\cdot\\\notag
\{\|\leftexp{(i_{m_{1}}...i_{m_{p}};i_{l-k})}{Z}\|_{L^{\infty}(\Sigma_{t}^{\epsilon_{0}})}(3\max_{j}\|\leftexp{(i_{1}
\overset{>i_{m_{1}}...i_{m_{p}}i_{l-k}<}{...}i_{l}j)}{\chi}^{\prime}_{l-p}\|_{L^{2}(\Sigma_{t}^{\epsilon_{0}})}\\\notag
\|\leftexp{(i_{1}\overset{>i_{m_{1}}...i_{m_{p}}i_{l-k}<}{...}i_{l})}{\chi}^{\prime}_{l-1-p}\|_{L^{2}(\Sigma_{t}^{\epsilon_{0}})})\\\notag
+3\|\leftexp{(i_{1}\overset{>i_{m_{1}}...i_{m_{p}}i_{l-k}<}{...}i_{l})}{\chi}^{\prime}_{l-1-p}\|_{L^{\infty}(\Sigma_{t}^{\epsilon_{0}})}
\max_{j}\|\leftexp{(i_{m_{1}}...i_{m_{p}}j;i_{l-k})}{Z}\|_{L^{2}(\Sigma_{t}^{\epsilon_{0}})}\}\\\notag
\leq C_{l}\delta_{0}(1+t)^{-3}[1+\log(1+t)][\mathcal{Y}_{0}+(1+t)\mathcal{A}'_{[l]}+\mathcal{W}_{[l+1]}]
\end{align}

     In Case 3 we apply (12.312) and also (12.307) with $l$ replaced by $l-1$ to estimate:
\begin{align}
 \|\leftexp{(i_{m_{1}}...i_{m_{p}};i_{l-k})}{Z}\|_{L^{2}(\Sigma_{t}^{\epsilon_{0}})}\leq C_{l}\{
\mathcal{Y}_{0}+(1+t)\mathcal{A}'_{[l-1]}+\mathcal{W}_{[l]}\}
\end{align}
Since in Case 3 we have $l-p\leq l-l_{*}-1\leq l_{*}$, we may use $\slashed{\textbf{X}}'_{[l_{*}]}$ to estimate:
\begin{align}
 \max_{j}\|\leftexp{(i_{1}
\overset{i_{m_{1}}...i_{m_{p}}i_{l-k}}{...}i_{l}j)}{\chi}^{\prime}_{l-p}\|_{L^{\infty}(\Sigma_{t}^{\epsilon_{0}})}
\leq C_{l}\delta_{0}(1+t)^{-2}[1+\log(1+t)]
\end{align}
as well as (12.313).

So in Case 3 we have:
\begin{align}
 \|\slashed{\mathcal{L}}_{\leftexp{(i_{m_{1}}...i_{m_{p}};i_{l-k})}{Z}}\leftexp{(i_{1}\overset{>i_{m_{1}}...i_{m_{p}}i_{l-k}<}{...}i_{l})}{\chi'_{l-1-p}}\|
_{L^{2}(\Sigma_{t}^{\epsilon_{0}})}\leq C(1+t)^{-1}\cdot\\\notag
\{\|\leftexp{(i_{m_{1}}...i_{m_{p}};i_{l-k})}{Z}\|_{L^{2}(\Sigma_{t}^{\epsilon_{0}})}(3\max_{j}\|\leftexp{(i_{1}
\overset{>i_{m_{1}}...i_{m_{p}}i_{l-k}<}{...}i_{l}j)}{\chi}^{\prime}_{l-p}\|_{L^{\infty}(\Sigma_{t}^{\epsilon_{0}})}\\\notag
\|\leftexp{(i_{1}\overset{>i_{m_{1}}...i_{m_{p}}i_{l-k}<}{...}i_{l})}{\chi}^{\prime}_{l-1-p}\|_{L^{\infty}(\Sigma_{t}^{\epsilon_{0}})})\\\notag
+3\|\leftexp{(i_{1}\overset{>i_{m_{1}}...i_{m_{p}}i_{l-k}<}{...}i_{l})}{\chi}^{\prime}_{l-1-p}\|_{L^{\infty}(\Sigma_{t}^{\epsilon_{0}})}
\max_{j}\|\leftexp{(i_{m_{1}}...i_{m_{p}}j;i_{l-k})}{Z}\|_{L^{2}(\Sigma_{t}^{\epsilon_{0}})}\}\\\notag
\leq C_{l}\delta_{0}(1+t)^{-3}[1+\log(1+t)][\mathcal{Y}_{0}+(1+t)\mathcal{A}'_{[l]}+\mathcal{W}_{[l+1]}]
\end{align}
Combining (12.311), (12.314) and (12.317) we conclude that the double sum in (12.195) is bounded in $L^{2}(\Sigma_{t}^{\epsilon_{0}})$ by:
\begin{equation}
 C_{l}\delta_{0}(1+t)^{-3}[1+\log(1+t)][\mathcal{Y}_{0}+(1+t)\mathcal{A}'_{[l]}+\mathcal{W}_{[l+1]}]
\end{equation}
Combining this result with (12.305) and (12.309) we then conclude that the sum on the right of (12.176) is bounded in $L^{2}(\Sigma_{t}^{\epsilon_{0}})$ by:
\begin{align}
 \|\sum_{k=0}^{l-1}\slashed{\mathcal{L}}_{R_{i_{l}}}...\slashed{\mathcal{L}}_{R_{i_{l-k+1}}}\slashed{\mathcal{L}}_{\leftexp{(R_{i_{l-k}})}{Z}}
\slashed{\mathcal{L}}_{R_{i_{l-k-1}}}...\slashed{\mathcal{L}}_{R_{i_{1}}}\chi'\|_{L^{2}(\Sigma_{t}^{\epsilon_{0}})}\\\notag
\leq C_{l}\delta_{0}(1+t)^{-3}[1+\log(1+t)][\mathcal{Y}_{0}+(1+t)\mathcal{A}'_{[l]}+\mathcal{W}_{[l+1]}]\\\notag
\end{align}

     In view of (12.302) and (12.319), we obtain through (12.176) and (12.181), recalling also (12.187):
\begin{align}
 \|\slashed{\mathcal{L}}_{L}\leftexp{(i_{1}...i_{l})}{\chi}^{\prime}_{l}\|_{L^{2}(\Sigma_{t}^{\epsilon_{0}})}
\leq C_{l}\delta_{0}(1+t)^{-3}[1+\log(1+t)][\mathcal{Y}_{0}+(1+t)\mathcal{A}'_{[l]}]+(1+t)^{-2}\{\mathcal{W}_{[l+2]}+\mathcal{W}^{Q}_{[l+1]}\}
\end{align}

    Let us define the non-negative functions:
\begin{equation}
 \leftexp{(i_{1}...i_{l})}{\phi}_{l}=(1-u+t)^{2}|\leftexp{(i_{1}...i_{l})}{\chi}^{\prime}_{l}|
\end{equation}
\begin{equation}
 \leftexp{(i_{1}...i_{l})}{\rho}_{l}=(1-u+t)^{2}(2|\chi'||\leftexp{(i_{1}...i_{l})}{\chi}^{\prime}_{l}|
+|\slashed{\mathcal{L}}_{L}\leftexp{(i_{1}...i_{l})}{\chi}^{\prime}_{l}|)
\end{equation}
From (12.56) we have:
\begin{equation}
 L\leftexp{(i_{1}...i_{l})}{\phi}_{l}\leq \leftexp{(i_{1}...i_{l})}{\rho}_{l}
\end{equation}
We pull back $\leftexp{(i_{1}...i_{l})}{\phi}_{l}$ and $\leftexp{(i_{1}...i_{l})}{\rho}_{l}$ to $S^{2}$:
\begin{align}
 \leftexp{(i_{1}...i_{l})}{\phi}_{l}(t,u)=\leftexp{(i_{1}...i_{l})}{\phi}_{l}\circ\Phi_{t,u},\quad
\leftexp{(i_{1}...i_{l})}{\rho}_{l}(t,u)=\leftexp{(i_{1}...i_{l})}{\rho}_{l}\circ\Phi_{t,u}
\end{align}
where $\Phi_{t,u}$ is the diffeomorphism of $S^{2}$ onto $S_{t,u}$ defined in Chapter 8. Then (12.324) reads:
\begin{align*}
 \frac{\partial}{\partial t}\leftexp{(i_{1}...i_{l})}{\phi}_{l}(t,u)\leq \leftexp{(i_{1}...i_{l})}{\rho}_{l}(t,u)
\end{align*}
Integrating we obtain:
\begin{equation}
 \leftexp{(i_{1}...i_{l})}{\phi}_{l}(t,u)\leq \leftexp{(i_{1}...i_{l})}{\phi}_{l}(0,u)+\int_{0}^{t}\leftexp{(i_{1}...i_{l})}{\rho}_{l}(t',u)dt'
\end{equation}
Taking $L^{2}$ norm on $[0,\epsilon_{0}]\times S^{2}$ yields:
\begin{align}
 \|\leftexp{(i_{1}...i_{l})}{\phi}_{l}(t)\|_{L^{2}([0,\epsilon_{0}]\times S^{2})}\leq 
\|\leftexp{(i_{1}...i_{l})}{\phi}_{l}(0)\|_{L^{2}([0,\epsilon_{0}]\times S^{2})}+\int_{0}^{t}\|
\leftexp{(i_{1}...i_{l})}{\rho}_{l}(t')\|_{L^{2}([0,\epsilon_{0}]\times S^{2})}dt'
\end{align}
In view of the comparison inequalities (8.333) we then obtain:
\begin{align}
 (1+t)^{-1}\|\leftexp{(i_{1}...i_{l})}{\phi}_{l}\|_{L^{2}(\Sigma_{t}^{\epsilon_{0}})}\\\notag
\leq C\{\|\leftexp{(i_{1}...i_{l})}{\phi}_{l}\|_{L^{2}(\Sigma_{0}^{\epsilon_{0}})}+\int_{0}^{t}(1+t')^{-1}
\|\leftexp{(i_{1}...i_{l})}{\rho}_{l}\|_{L^{2}(\Sigma_{t'}^{\epsilon_{0}})}dt'\}
\end{align}
Replacing $l$ by $k\in\{0,...,l\}$, taking on both sides the maximum over $i_{1},...,i_{k}$ and then summing over $k\in\{0,...,l\}$ 
yields:
\begin{align}
 (1+t)^{-1}\sum_{k=0}^{l}\max_{i_{1}...i_{k}}\|\leftexp{(i_{1}...i_{k})}{\phi}_{k}\|_{L^{2}(\Sigma_{t}^{\epsilon_{0}})}\\\notag
\leq C\{\sum_{k=0}^{l}\max_{i_{1}...i_{k}}\|\leftexp{(i_{1}...i_{k})}{\phi}_{k}\|_{L^{2}(\Sigma_{0}^{\epsilon_{0}})}
+\int_{0}^{t}(1+t')^{-1}\sum_{k=0}^{l}\max_{i_{1}...i_{k}}\|\leftexp{(i_{1}...i_{k})}{\rho}_{k}\|_{L^{2}(\Sigma_{t'}^{\epsilon_{0}})}dt'\}
\end{align}
From (12.321) we have:
\begin{align}
 \sum_{k=0}^{l}\max_{i_{1}...i_{l}}\|\leftexp{(i_{1}...i_{k})}{\phi}_{k}\|_{L^{2}(\Sigma_{t}^{\epsilon_{0}})}\geq C^{-1}(1+t)^{2}\mathcal{A}'_{[l]}(t)\\\notag
\sum_{k=0}^{l}\max_{i_{1}...i_{k}}\|\leftexp{(i_{1}...i_{k})}{\phi}_{k}\|_{L^{2}(\Sigma_{0}^{\epsilon_{0}})}\leq \mathcal{A}'_{[l]}(0)
\end{align}
and from (12.322), (12.320) and $\slashed{\textbf{X}}^{\prime}_{[0]}$ we have:
\begin{align}
 \sum_{k=0}^{l}\max_{i_{1}...i_{k}}\|\leftexp{(i_{1}...i_{k})}{\rho}_{k}\|_{L^{2}(\Sigma_{t}^{\epsilon_{0}})}\leq\\\notag
C\delta_{0}[1+\log(1+t)]\mathcal{A}'_{[l]}(t)+C(1+t)^{2}\|\slashed{\mathcal{L}}_{L}\leftexp{(i_{1}...i_{l})}{\chi}^{\prime}_{l}\|_{L^{2}(\Sigma_{t}^{\epsilon_{0}})}\\\notag
\leq C_{l}\delta_{0}(1+t)^{-1}[1+\log(1+t)]\{\mathcal{Y}_{0}+(1+t)\mathcal{A}'_{[l]}(t)\}+C_{l}\{\mathcal{W}_{[l+2]}(t)+\mathcal{W}^{Q}_{[l+1]}(t)\}
\end{align}
Substituting (12.329) and (12.330) in (12.328) yields:
\begin{align}
 (1+t)\mathcal{A}'_{[l]}(t)\leq C\mathcal{A}'_{[l]}(0)+\int_{0}^{t}C_{l}\delta_{0}(1+t')^{-2}[1+\log(1+t')]\cdot(1+t')\mathcal{A}'_{[l]}(t')dt'\\\notag
+\int_{0}^{t}C_{l}\delta_{0}(1+t')^{-2}[1+\log(1+t')]\mathcal{Y}_{0}(t')dt'+C_{l}\int_{0}^{t}(1+t)^{-1}\{\mathcal{W}_{[l+2]}(t')+\mathcal{W}^{Q}_{[l+1]}(t')\}dt'
\end{align}
This implies:
\begin{align}
 (1+t)\mathcal{A}'_{[l]}(t)\leq C\mathcal{A}'_{[l]}(0)+
\int_{0}^{t}C_{l}\delta_{0}(1+t')^{-2}[1+\log(1+t')]\mathcal{Y}_{0}(t')dt'\\\notag
+C_{l}\int_{0}^{t}(1+t^{\prime})^{-1}\{\mathcal{W}_{[l+2]}(t')+\mathcal{W}^{Q}_{[l+1]}(t')\}dt'
\end{align}

     Finally, we must estimate $\mathcal{Y}_{0}(t)$ in terms of $\mathcal{Y}_{0}(0)$. To do this, we must derive a propagation equation for 
$y^{i}$. According to (6.60):
\begin{equation}
 \hat{T}^{i}=-\frac{x^{i}}{1-u+t}+y^{i}
\end{equation}
hence:
\begin{equation}
 L\hat{T}^{i}=-\frac{L^{i}}{1-u+t}+\frac{x^{i}}{(1-u+t)^{2}}+Ly^{i}
\end{equation}
Also from (6.64) and (6.65),
\begin{equation}
 L^{i}=\frac{x^{i}}{1-u+t}+\omega^{i}-\eta y^{i}
\end{equation}
where
\begin{equation}
 \omega^{i}=\frac{(\eta-1)x^{i}}{1-u+t}-\psi_{i}
\end{equation}
Substituting (12.335) in (12.334) we get:
\begin{equation}
 L\hat{T}^{i}=\frac{(\eta y^{i}-\omega^{i})}{(1-u+t)}+Ly^{i}
\end{equation}
On the other hand, recalling from Chapter 3 that:
\begin{equation}
 L(\hat{T}^{i})=q_{L}\cdot\slashed{d}x^{i}
\end{equation}
we conclude $y^{i}$ satisfies the propagation equation:
\begin{align}
 Ly^{i}+\frac{y^{i}}{1-u+t}=\tilde{\omega}^{i}
\end{align}
where 
\begin{align}
 \tilde{\omega}^{i}=\frac{\omega^{i}-(\eta-1)y^{i}}{1-u+t}+q_{L}\cdot\slashed{d}x^{i}
\end{align}
From (12.23) and (12.336), the first term on the right of (12.340) is bounded in $L^{2}(\Sigma_{t}^{\epsilon_{0}})$ by:
\begin{equation}
 C(1+t)^{-1}\mathcal{W}_{[0]}(t)
\end{equation}
From (3.183) and (12.333), we deduce:
\begin{equation}
 \|q_{L}\|_{L^{2}(\Sigma_{t}^{\epsilon_{0}})}\leq C(1+t)^{-1}\mathcal{W}_{[1]}(t)
\end{equation}
hence we obtain:
\begin{equation}
 \max_{i}\|\tilde{\omega}^{i}\|_{L^{2}(\Sigma_{t}^{\epsilon_{0}})}\leq C(1+t)^{-1}\mathcal{W}_{[1]}(t)
\end{equation}
Integrating (12.339) we get, in acoustical coordinates:
\begin{equation}
 (1-u+t)y^{i}(t,u,\vartheta)=(1-u)y^{i}(0,u,\vartheta)+\int_{0}^{t}(1-u+t')\tilde{\omega}^{i}(t',u,\vartheta)dt'
\end{equation}
Taking $L^{2}$ norm on $[0,\epsilon_{0}]\times S^{2}$ we then obtain:
\begin{align}
 (1+t)\|y^{i}(t)\|_{L^{2}([0,\epsilon_{0}]\times S^{2})}\leq 
C\{\|y^{i}(0)\|_{L^{2}([0,\epsilon_{0}]\times S^{2})}+\int_{0}^{t}(1+t')\|\tilde{\omega}^{i}\|_{L^{2}([0,\epsilon_{0}]\times S^{2})}dt'\}
\end{align}
This is equivalent to
\begin{equation}
 \|y^{i}\|_{L^{2}(\Sigma_{t}^{\epsilon_{0}})}\leq C\{\|y^{i}\|_{L^{2}(\Sigma_{0}^{\epsilon_{0}})}
+\int_{0}^{t}\|\tilde{\omega}^{i}\|_{L^{2}(\Sigma_{t'}^{\epsilon_{0}})}dt'\}
\end{equation}
Taking the maximum over $i=1,2,3$, yields:
\begin{equation}
 \mathcal{Y}_{0}(t)\leq C\{\mathcal{Y}_{0}(0)+\int_{0}^{t}(1+t')^{-1}\mathcal{W}_{[1]}(t')dt'\}
\end{equation}
This implies:
\begin{align}
 \int_{0}^{t}(1+t')^{-2}[1+\log(1+t')]\mathcal{Y}_{0}(t')dt'\leq
C\{\mathcal{Y}_{0}(0)+\int_{0}^{t}(1+t')^{-2}[1+\log(1+t^{\prime})]\mathcal{W}_{[1]}(t')dt'\}
\end{align}
In fact, let
\begin{align*}
 I_{t}(t^{\prime})=-\int_{t^{\prime}}^{t}(1+s)^{-2}[1+\log(1+s)]ds
\end{align*}
Then
\begin{align*}
 \frac{d}{dt^{\prime}}I_{t}(t^{\prime})=(1+t^{\prime})^{-2}[1+\log(1+t^{\prime})]
\end{align*}
Let also
\begin{align*}
 J(t^{\prime})=\int_{0}^{t^{\prime}}(1+s)^{-1}\mathcal{W}_{[1]}(s)ds
\end{align*}
Then:
\begin{align*}
 \int_{0}^{t}(1+t^{\prime})^{-2}[1+\log(1+t^{\prime})]\int_{0}^{t^{\prime}}(1+s)^{-1}\mathcal{W}_{[1]}(s)dsdt^{\prime}\\
=I_{t}(t)J(t)-I_{t}(0)J(0)-\int_{0}^{t}I_{t}(t^{\prime})\frac{dJ}{dt^{\prime}}(t^{\prime})dt^{\prime}
\end{align*}
But
\begin{align*}
 I_{t}(t)=J(0)=0
\end{align*}
So this reduces to:
\begin{align*}
 -\int_{0}^{t}I_{t}(t^{\prime})\frac{dJ}{dt^{\prime}}(t^{\prime})dt^{\prime}
\end{align*}
Now,
\begin{align*}
 -I_{t}(t^{\prime})=\int_{t^{\prime}}^{t}(1+s)^{-2}[1+\log(1+s)]ds
\leq\int_{t^{\prime}}^{\infty}(1+s)^{-2}[1+\log(1+s)]ds\\
=[-\frac{2}{1+s}-\frac{\log(1+s)}{1+s}]^{\infty}_{t^{\prime}}=\frac{2+\log(1+t^{\prime})}{1+t^{\prime}}
\end{align*}
Therefore
\begin{align*}
 -\int_{0}^{t}I_{t}(t^{\prime})\frac{dJ}{dt^{\prime}}(t^{\prime})dt^{\prime}
\leq 2\int_{0}^{t}(1+t^{\prime})^{-2}[1+\log(1+t^{\prime})]\mathcal{W}_{[1]}(t^{\prime})dt^{\prime}
\end{align*}
Substituting (12.348) in the right in (12.332), the integral on the right in (12.348) is absorbed into the last integral on the right in (12.332).
So the proposition follows. $\qed$

$\textbf{Proposition 12.12}$  Let the assumptions of Proposition 12.11 hold for some non-negative integer $l$. 
Let $m\in\{0,...,l+1\}$, and let the initial data satisfy:
\begin{align*}
 \|\mu-1\|_{\infty,[m,l_{*}+1],\Sigma_{0}^{\epsilon_{0}}}\leq C_{l}\delta_{0}
\end{align*}
Then $\textbf{M}_{[m,l_{*}+1]}$ holds on $W^{s}_{\epsilon_{0}}$ and there is a constant $C_{l}$ independent of $s$ such that for all 
$t\in[0,s]$ we have:
\begin{align*}
 (1+t)^{-1}\mathcal{B}_{[m,l+1]}(t)\leq C\mathcal{B}_{[m,l+1]}(0)+C_{l}\{\mathcal{A}'_{[l]}(0)+\delta_{0}\mathcal{Y}_{0}(0)\\\notag
+\int_{0}^{t}(1+t')^{-1}[\mathcal{W}_{\{l+2\}}(t')+\mathcal{W}^{Q}_{\{l+1\}}(t')]dt'\}
\end{align*}

$Proof$: The assumptions of this proposition include those of Proposition 12.10 with $l_{*}$ in the role of $l$, therefore $\textbf{M}_{[m,l_{*}+1]}$
holds on $W^{s}_{\epsilon_{0}}$. Thus, the assumptions of Proposition 11.1 hold with $l_{*}$ in the role of $l$, hence also the assumptions of 
Proposition 11.2 hold. Therefore, the conclusions of Proposition 11.1 and of all its corollaries hold with $l_{*}$ in the role of $l$, and the conclusions
of Proposition 11.2 and of all its corollaries hold as well.

    We consider the propagation equation (12.252), replacing $m+1$ by $m$, so that $m\in\{0,...,l+1\}$. We must estimate in $L^{2}(\Sigma_{t}^{\epsilon_{0}})$
the terms on the right side.   
From the definition of $m$ and $e$, we have:
\begin{align}
 \|m\|_{2,[m,l+1],\Sigma_{t}^{\epsilon_{0}}}\leq C_{l}\mathcal{W}_{\{l+2\}}
\end{align}
\begin{align}
 \|e\|_{2,[m,l+1],\Sigma_{t}^{\epsilon_{0}}}\leq C_{l}(1+t)^{-1}\{\mathcal{W}^{Q}_{\{l+1\}}+\delta_{0}(1+t)^{-1}[(1+t)^{-1}\\\notag
\mathcal{B}_{[m-1,l+1]}+\mathcal{Y}_{0}+(1+t)\mathcal{A}_{[l]}+\mathcal{W}_{\{l+1\}}]\}
\end{align}
Here, we have also used Corollary 11.2.b (the estimate for $\omega_{L\hat{T}}$; note that actually, only $\textbf{M}_{[m,l_{*}]}$ is used in proving
this estimate).

Also from Corollary 11.1.b with $l_{*}$ in the role of $l$ we obtain:
\begin{equation}
 \|e\|_{\infty,[m,l_{*}],\Sigma_{t}^{\epsilon_{0}}}\leq C_{l}\delta_{0}(1+t)^{-2}
\end{equation}
Consider next $\leftexp{(i_{1}...i_{l})}{r}^{\prime}_{m,n}$, given by (12.251) with $m+1$ replaced by $m$. Here $n\leq l+1-m$.
Consider the summand:
\begin{align}
 ((R)^{s_{1}}e)((R)^{s_{2}}(T)^{m}\mu)\quad :\quad |s_{1}|+|s_{2}|=n, |s_{1}|\geq 1\notag
\end{align}
If $|s_{2}|\leq l_{*}+1-m$ we appeal to $\textbf{M}_{[m,l_{*}+1]}$ to place the second factor in $L^{\infty}(\Sigma_{t}^{\epsilon_{0}})$. 
Otherwise, we have $|s_{1}|\leq n-(l_{*}+2-m)\leq l-l_{*}-1\leq l_{*}$ and we appeal to (12.351) to place the first factor in $L^{\infty}(\Sigma_{t}^{\epsilon_{0}})$.
Noting that $|s_{2}|\leq n-1\leq l-m$, we then obtain an $L^{2}(\Sigma_{t}^{\epsilon_{0}})$ for the first sum by:
\begin{equation}
 C_{l}\delta_{0}(1+t)^{-2}\mathcal{B}_{[m,l]}+C_{l}[1+\log(1+t)]\|e\|_{2,[m,l+1],\Sigma_{t}^{\epsilon_{0}}}
\end{equation}
Consider next the summand:
\begin{align*}
 R_{i_{n}}...R_{i_{1}}(((T)^{k}e)((T)^{m-k}\mu))\quad:\quad k\in\{1,...,m\}
\end{align*}
of the second sum. If at most $l_{*}+1+k-m$ of the $n$ angular derivatives fall on $(T)^{m-k}\mu$ we appeal to $\textbf{M}_{[m,l_{*}+1]}$ to place
the corresponding factor in $L^{\infty}(\Sigma_{t}^{\epsilon_{0}})$. Otherwise, at most
\begin{align*}
 n-(l_{*}+2+k-m)\leq l-l_{*}-1-k\leq l_{*}-k
\end{align*}
angular derivatives fall on $(T)^{k}e$ and we appeal to (12.351) to place the corresponding factor in $L^{\infty}(\Sigma_{t}^{\epsilon_{0}})$. 
Noting that $k\geq 1$, the second term on the right of (12.251) is bounded by:
\begin{equation}
 C_{l}\delta_{0}(1+t)^{-2}\mathcal{B}_{[m-1,l]}+C_{l}[1+\log(1+t)]\|e\|_{2,[m,l+1],\Sigma_{t}^{\epsilon_{0}}}
\end{equation}
Combining (12.352) and (12.353), and substituting the bound (12.350), we obtain:
\begin{align}
 \max_{n\leq l+1-m}\max_{i_{1}...i_{n}}\|\leftexp{(i_{1}...i_{n})}{r}^{\prime}_{m,n}\|_{L^{2}(\Sigma_{t}^{\epsilon_{0}})}
\leq C_{l}(1+t)^{-1}\{\delta_{0}(1+t)^{-1}\mathcal{B}_{[m,l+1]}\\\notag
+[1+\log(1+t)][\mathcal{W}^{Q}_{\{l+1\}}+\delta_{0}(1+t)^{-1}(\mathcal{Y}_{0}+(1+t)\mathcal{A}_{[l]}+\mathcal{W}_{\{l+1\}})]\}
\end{align}
Next, we must estimate the two sums on the right of (12.252) with $m+1$ replaced by $m$. The first sum is bounded in $L^{2}(\Sigma_{t}^{\epsilon_{0}})$
by:
\begin{equation}
 \sum_{k=0}^{m-1}\|\Lambda\cdot\slashed{d}(T)^{m-1-k}\mu\|_{2,[k,l+k+1-m],\Sigma_{t}^{\epsilon_{0}}}
\end{equation}
Here, $\Lambda, \slashed{d}(T)^{m-1-k}\mu$, receives at most $k$ $T$-derivatives, thus $\mu$ receives at most $m-1$ $T$-derivatives, and there are 
$l+k+1-m$ spatial derivatives. By Corollary 11.1.c with $l_{*}$ in the role of $l$:
\begin{equation}
 \|\Lambda\|_{\infty,[m,l_{*}],\Sigma_{t}^{\epsilon_{0}}}\leq C_{l}\delta_{0}(1+t)^{-1}[1+\log(1+t)]
\end{equation}
If at most $l_{*}$ spatial derivatives fall on $\Lambda$ we appeal to this estimate to place the corresponding factor in $L^{\infty}(\Sigma_{t}^{\epsilon_{0}})$.
Otherwise, at most
\begin{align*}
 (l+k+1-m)-(l_{*}+1)\leq l_{*}+k+1-m
\end{align*}
spatial derivatives fall on $\slashed{d}(T)^{m-1-k}\mu$, for a total of at most $l_{*}+1$ spatial derivatives falling on $\mu$, and we appeal to
$\textbf{M}_{[m-1,l_{*}+1]}$ to place the corresponding factor in $L^{\infty}(\Sigma_{t}^{\epsilon_{0}})$. In this way, we deduce that (12.355) is 
bounded by:
\begin{align}
 C_{l}\delta_{0}(1+t)^{-2}[1+\log(1+t)]\mathcal{B}_{[m-1,l+1]}\\\notag
+C_{l}\delta_{0}(1+t)^{-1}[1+\log(1+t)]\|\Lambda\|_{2,[m-1,l],\Sigma_{t}^{\epsilon_{0}}}\\\notag
\leq C_{l}\delta_{0}(1+t)^{-2}[1+\log(1+t)]\cdot\{\mathcal{B}_{[m-1,l+1]}+\delta_{0}[1+\log(1+t)](\mathcal{Y}_{0}+(1+t)\mathcal{A}_{[l-1]})\\\notag
+[1+\log(1+t)][\mathcal{W}_{\{l+1\}}+\delta_{0}(1+t)^{-2}[1+\log(1+t)]^{2}\mathcal{W}^{Q}_{\{l\}}]\}
\end{align}
by Corollary 11.2.c with $m$ replaced by $m-1$.

     Finally, we must consider the second sum on the right of (12.252) with $m+1$ replaced by $m$. Consider a term in this sum, corresponding to 
$k\in\{0,...,n-1\}$:
\begin{align}
 R_{i_{n}}...R_{i_{n-k+1}}(\leftexp{(R_{i_{n-k}})}{Z}\cdot\slashed{d}R_{i_{n-k-1}}...R_{i_{1}}(T)^{m}\mu)\\\notag
=\sum_{|s_{1}|+|s_{2}|=k}((\slashed{\mathcal{L}}_{R})^{s_{1}}\leftexp{(R_{i_{n-k}})}{Z})\cdot\slashed{d}(R)^{s'_{2}}(T)^{m}\mu\\\notag
s'_{2}=s_{2}\bigcup\{1,...,n-k-1\}
\end{align}
If $|s_{1}|\leq l_{*}$ we appeal to Corollary 10.1.i with $l_{*}+1$ in the role of $l$ to place the first factor in $L^{\infty}(\Sigma_{t}^{\epsilon_{0}})$.
Otherwise, $|s_{2}|\leq k-l_{*}-1$ and
\begin{align*}
 |s'_{2}|=|s_{2}|+n-k-1\leq n-l_{*}-2\leq l-l_{*}-1-m\leq l_{*}-m
\end{align*}
hence $1+|s'_{2}|+m\leq l_{*}+1$ and we appeal to $\textbf{M}_{[m,l_{*}+1]}$ to place the second factor in $L^{\infty}(\Sigma_{t}^{\epsilon_{0}})$.
Since also $|s_{1}|\leq k\leq n-1\leq l-m$ and
\begin{align*}
 1+|s'_{2}|+m=|s_{2}|+n-k+m\leq n+m\leq l+1
\end{align*}
we obtain in this way that the second sum in (12.252) is bounded by:
\begin{align}
 C_{l}\delta_{0}(1+t)^{-2}[1+\log(1+t)]\mathcal{B}_{[m,l+1]}+
C_{l}\delta_{0}(1+t)^{-1}[1+\log(1+t)]\max_{i}\|\leftexp{(R_{i})}{Z}\|_{2,[l-m],\Sigma_{t}^{\epsilon_{0}}}\\\notag
\leq C_{l}\delta_{0}(1+t)^{-1}[1+\log(1+t)]\cdot\\\notag
\{(1+t)^{-1}\mathcal{B}_{[m,l+1]}+\mathcal{Y}_{0}+(1+t)\mathcal{A}_{[l]}+\mathcal{W}_{\{l+1\}}\}
\end{align}
by Corollary 10.2.i with $l+1$ in the role of $l$.

     In view of (12.347), (12.352), (12.355) and (12.357) and also substituting for $\mathcal{A}_{[l]}$ in terms of $\mathcal{A}'_{[l]}$ from Lemma 12.5
we deduce:
\begin{align}
 \|L\leftexp{(i_{1}...i_{n})}{\mu}_{m,n}\|_{L^{2}(\Sigma_{t}^{\epsilon_{0}})}\leq C_{l}\delta_{0}(1+t)^{-2}[1+\log(1+t)]\mathcal{B}_{[m,l+1]}\\\notag
+C_{l}\{\mathcal{W}_{\{l+2\}}+(1+t)^{-1}[1+\log(1+t)][\mathcal{W}^{Q}_{\{l+1\}}+\delta_{0}(\mathcal{Y}_{0}+(1+t)\mathcal{A}'_{[l]})]\}
\end{align}
for $n\leq l+1-m$.

Defining:
\begin{align}
 \leftexp{(i_{1}...i_{n})}{\phi}^{\prime}_{m,n}=|\leftexp{(i_{1}...i_{n})}{\mu}_{m,n}|\\
\leftexp{(i_{1}...i_{n})}{\rho}^{\prime}_{m,n}=|L\leftexp{(i_{1}...i_{n})}{\mu}_{m,n}|
\end{align}
obviously, we have:
\begin{equation}
 L\leftexp{(i_{1}...i_{n})}{\phi}^{\prime}_{m,n}\leq \leftexp{(i_{1}...i_{n})}{\rho}^{\prime}_{m,n}
\end{equation}
Following an argument similar to that leading to (12.327), we deduce:
\begin{align}
 (1+t)^{-1}\|\leftexp{(i_{1}...i_{n})}{\phi}^{\prime}_{m,n}\|_{L^{2}(\Sigma_{t}^{\epsilon_{0}})}
\leq C\{\|\leftexp{(i_{1}...i_{n})}{\phi}^{\prime}_{m,n}\|_{L^{2}(\Sigma_{0}^{\epsilon_{0}})}+\int_{0}^{t}(1+t')^{-1}
\|\leftexp{(i_{1}...i_{n})}{\rho}^{\prime}_{m,n}\|_{L^{2}(\Sigma_{t'}^{\epsilon_{0}})}dt'\}
\end{align}
Taking on both sides the maximum over $i_{1}...i_{n}$, replacing $m$ by $k\in\{0,...,m\}$ and then summing over $n\in\{0,...,l+1-k\}$
and over $k\in\{0,...,m\}$ yields:
\begin{align}
 (1+t)^{-1}\sum_{k=0}^{m}\sum_{n=0}^{l+1-k}\max_{i_{1}...i_{n}}\|\leftexp{(i_{1}...i_{n})}{\phi}^{\prime}_{k,n}\|_{L^{2}(\Sigma_{t}^{\epsilon_{0}})}\leq
C\{\sum_{k=0}^{m}\sum_{n=0}^{l+1-k}\|\leftexp{(i_{1}...i_{n})}{\phi}^{\prime}_{k,n}\|_{L^{2}(\Sigma_{0}^{\epsilon_{0}})}\\\notag
+\int_{0}^{t}(1+t')^{-1}
\sum_{k=0}^{m}\sum_{n=0}^{l+1-k}\|\leftexp{(i_{1}...i_{n})}{\rho}^{\prime}_{k,n}\|_{L^{2}(\Sigma_{t'}^{\epsilon_{0}})}dt'\}
\end{align}
From (12.361) we have:
\begin{equation}
 \sum_{k=0}^{m}\sum_{n=0}^{l+1-k}\max_{i_{1}...i_{n}}\|\leftexp{(i_{1}...i_{n})}{\phi}^{\prime}_{k,n}\|_{L^{2}(\Sigma_{t}^{\epsilon_{0}})}
=\mathcal{B}_{[m,l+1]}(t)
\end{equation}
From (12.360) and (12.362) we have:
\begin{align}
 \sum_{k=0}^{m}\sum_{n=0}^{l+1-k}\|\leftexp{(i_{1}...i_{n})}{\rho}^{\prime}_{k,n}\|_{L^{2}(\Sigma_{t}^{\epsilon_{0}})}\leq 
C_{l}\delta_{0}(1+t)^{-2}[1+\log(1+t)]\mathcal{B}_{[m,l+1]}\\\notag
+C_{l}\{\mathcal{W}_{\{l+2\}}+(1+t)^{-1}[1+\log(1+t)][\mathcal{W}^{Q}_{\{l+1\}}+\delta_{0}(\mathcal{Y}_{0}+(1+t)\mathcal{A}'_{[l]})]\}
\end{align}

Substitute (12.366) and (12.367) in (12.365) we obtain:
\begin{align}
 (1+t)^{-1}\mathcal{B}_{[m,l+1]}(t)\leq C\mathcal{B}_{[m,l+1]}(0)
+\int_{0}^{t}C_{l}\delta_{0}(1+t')^{-2}[1+\log(1+t')]\cdot(1+t')^{-1}\mathcal{B}_{[m,l+1]}(t')dt'\\\notag
+\int_{0}^{t}C_{l}\delta_{0}(1+t')^{-2}[1+\log(1+t')]\mathcal{Y}_{0}(t')dt'\\\notag
+\int_{0}^{t}C_{l}(1+t')^{-1}[\mathcal{W}_{\{l+2\}}(t')+(1+t')^{-1}[1+\log(1+t')]\mathcal{W}^{Q}_{\{l+1\}}(t')]dt'\\\notag
+\int_{0}^{t}C_{l}\delta_{0}(1+t')^{-1}[1+\log(1+t')]\mathcal{A}'_{[l]}(t')dt'
\end{align}
This implies:
\begin{align}
 (1+t)^{-1}\mathcal{B}_{[m,l+1]}(t)\leq C\mathcal{B}_{[m,l+1]}(0)
+\int_{0}^{t}C_{l}\delta_{0}(1+t')^{-2}[1+\log(1+t')]\mathcal{Y}_{0}(t')dt'\\\notag
+\int_{0}^{t}C_{l}(1+t')^{-1}[\mathcal{W}_{\{l+2\}}(t')+(1+t')^{-1}[1+\log(1+t')]\mathcal{W}^{Q}_{\{l+1\}}(t')]dt'\\\notag
+\int_{0}^{t}C_{l}\delta_{0}(1+t')^{-1}[1+\log(1+t')]\mathcal{A}'_{[l]}(t')dt'
\end{align}
From Proposition 12.11, we have:
\begin{align}
 \int_{0}^{t}(1+t')^{-1}[1+\log(1+t')]\mathcal{A}'_{[l]}(t')dt'\leq \\\notag
C_{l}\{\mathcal{A}'_{[l]}(0)+\delta_{0}\mathcal{Y}_{0}(0)+\int_{0}^{t}(1+t')^{-1}[\mathcal{W}_{[l+1]}(t')+\mathcal{W}^{Q}_{[l]}(t')]dt'\}
\end{align}
Using then also (12.348), the proposition follows. $\qed$

\chapter{Derivation of the Basic Properties of $\mu$}

In this chapter we shall establish the assumptions $\textbf{C1}, \textbf{C2}, \textbf{C3}$ introduced in Chapter 5, on which
Theorem 5.1 relies. To do this, we shall use Proposition 8.6. The assumption of this proposition are those of Proposition 8.5 together with assumption $\textbf{I0}$ in Chapter 6, which coincide with 
the assumption of Proposition 12.1 on the initial data. While the assumptions of
Proposition 8.5 are the assumptions of Lemma 8.10 together with the bootstrap assumptions $\textbf{E}_{\textbf{LT}}\textbf{3}, 
\textbf{E}_{\textbf{LL}}\textbf{3}$. Finally, the assumptions of Lemma 8.10 are basic bootstrap assumptions $\textbf{A1}, \textbf{A2}
,\textbf{A3}$, introduced in Chapter 5, the bootstrap assumptions $\textbf{E1}, \textbf{E2}, \textbf{F2}$, of Chapter 6, together with
the additional bootstrap assumption $\slashed{\textbf{E}}\textbf{3}_{0}$. Now, $\textbf{A1}$ and $\textbf{A2}$ follow from $\textbf{E1}$,
while $\textbf{A3}$ follows from $\textbf{A2}$ and Proposition 12.1. $\textbf{F2}$ coincides with $\slashed{\textbf{X}}^{\prime}_{0}$ which has been 
established by Proposition 12.6. $\textbf{E1}$ is $\textbf{E}_{\{0\}}$ while $\textbf{E2}$ follows from $\textbf{E}_{\{1\}}, \textbf{E}^{Q}_{\{0\}}$,
 and $\textbf{H0}$, which has been established by Corollary 12.2.a. $\textbf{E}_{\textbf{LL}}\textbf{3}$ follows from $\textbf{E}^{QQ}_{\{0\}}$
and $\textbf{E}^{Q}_{\{0\}}$, while $\textbf{E}_{\textbf{LT}}\textbf{3}$ follows from $\textbf{E}^{Q}_{\{1\}}, \textbf{E}_{\{1\}}$, and the estimates for
$\|\Lambda\|_{L^{\infty}(\Sigma_{t}^{\epsilon_{0}})}$ of Corollary 11.1.c with $m=l=0$ as well as $\textbf{H0}$. The assumptions of Corollary 
11.1.c are those of Proposition 11.1 and follow from Proposition 12.10. Finally, $\slashed{\textbf{E}}\textbf{3}_{0}$ follows from $\textbf{E}_{\{2\}}$ 
and $\textbf{H0}$ and $\textbf{H1}$, $\textbf{H1}$ having been established by Corollary 12.5. We conclude from above that the assumptions of Proposition 8.6 
all follow from the assumptions of Proposition 12.10 with $m=l=0$, namely the assumptions of Proposition 12.9 with $l=0$, which coincide with those
of Proposition 12.6, together with:
\begin{equation}
 \|\mu-1\|_{\infty,[0,1],\Sigma_{0}^{\epsilon_{0}}}\leq C\delta_{0}
\end{equation}

$\textbf{Proposition 13.1}$ Let the assumptions of Proposition 12.6 hold and let the initial data satisfy (13.1). 
Then on $W^{s}_{\epsilon_{0}}$ we have:
\begin{align*}
 \mu^{-1}(L\mu)_{+}\leq (1+t)^{-1}[1+\log(1+t)]^{-1}+A(t)
\end{align*}
where
\begin{align*}
 A(t)=C\delta_{0}(1+t)^{-2}[1+\log(1+t)]
\end{align*}
Since
\begin{align*}
 \int_{0}^{s}A(t)dt\leq C\delta_{0}
\end{align*}
$C$ being a constant independent of $s$, we conclude that $\textbf{C1}$ holds on $W^{s}_{\epsilon_{0}}$.

$Proof$: We appeal to Proposition 8.6 to express, in acoustical coordinates:
\begin{align}
 \mu^{-1}(L\mu)_{+}=\frac{1}{\hat{\mu}_{s}}(\frac{\partial\hat{\mu}_{s}}{\partial t})_{+}
=\frac{(\hat{E}_{s}(u,\vartheta)(1+t)^{-1}+\hat{Q}_{1,s}(t,u,\vartheta))_{+}}
{1+\hat{E}_{s}(u,\vartheta)\log(1+t)+\hat{Q}_{0,s}(t,u,\vartheta)}
\end{align}
Consider a given $(u,\vartheta)\in[0,\epsilon_{0}]\times S^{2}$. There are three cases to consider.

Case 1: $\hat{E}_{s}(u,\vartheta)=0$,  Case 2: $\hat{E}_{s}(u,\vartheta)>0$,  Case 3: $\hat{E}_{s}(u,\vartheta)<0$

In Case 1, (13.2) reduces to:
\begin{equation}
 \frac{1}{\hat{\mu}_{s}}(\frac{\partial\hat{\mu}_{s}}{\partial t})_{+}=
\frac{(\hat{Q}_{1,s}(t,u,\vartheta))_{+}}{1+\hat{Q}_{0,s}(t,u,\vartheta)}
\end{equation}
From Proposition 8.6:
\begin{equation}
 |\hat{Q}_{0,s}(t,u,\vartheta)|\leq C\delta_{0}\frac{[1+\log(1+t)]}{(1+t)}
\end{equation}
\begin{equation}
 |\hat{Q}_{1,s}(t,u,\vartheta)|\leq C\delta_{0}\frac{[1+\log(1+t)]}{(1+t)^{2}}
\end{equation}
these imply:
\begin{equation}
 \frac{1}{\hat{\mu}_{s}}(\frac{\partial\hat{\mu}_{s}}{\partial t})_{+}\leq C\delta_{0}\frac{[1+\log(1+t)]}{(1+t)^{2}}
\end{equation}

   In Case 2, we have:
\begin{align*}
 (\hat{E}_{s}(u,\vartheta)(1+t)^{-1}+\hat{Q}_{1,s}(t,u,\vartheta))_{+}\leq 
\hat{E}_{s}(u,\vartheta)(1+t)^{-1}+|\hat{Q}_{1,s}(t,u,\vartheta)|
\end{align*}
hence, substituting in (13.2),
\begin{align}
 \frac{1}{\hat{\mu}_{s}}(\frac{\partial\hat{\mu}_{s}}{\partial t})_{+}\leq 
\frac{\hat{E}_{s}(u,\vartheta)(1+t)^{-1}+|\hat{Q}_{1,s}(t,u,\vartheta)|}{1+\hat{E}_{s}(u,\vartheta)\log(1+t)+\hat{Q}_{0,s}(t,u,\vartheta)}\\\notag
=\frac{\hat{E}_{s}(u,\vartheta)(1+t)^{-1}}{1+\hat{E}_{s}(u,\vartheta)\log(1+t)}\\\notag
+\frac{\hat{E}_{s}(u,\vartheta)}{(1+t)}\{\frac{1}{1+\hat{E}_{s}(u,\vartheta)\log(1+t)+\hat{Q}_{0,s}(t,u,\vartheta)}
-\frac{1}{1+\hat{E}_{s}(u,\vartheta)\log(1+t)}\}\\\notag
+\frac{|\hat{Q}_{1,s}(t,u,\vartheta)|}{1+\hat{E}_{s}(u,\vartheta)\log(1+t)+\hat{Q}_{0,s}(t,u,\vartheta)}
\end{align}
Since, for suitably small $\delta_{0}$,
\begin{equation}
 |\hat{E}_{s}(u,\vartheta)|\leq C\delta_{0}\leq 1
\end{equation}
the first term on the right of (13.7) is bounded by:
\begin{equation}
 \frac{1}{(1+t)[1+\log(1+t)]}
\end{equation}
The factor in parenthesis in the second term on the right in (13.7) is bounded in absolute value by:
\begin{align*}
 \frac{|\hat{Q}_{0,s}(t,u,\vartheta)|}{[1+\hat{E}_{s}(u,\vartheta)\log(1+t)+\hat{Q}_{0,s}(t,u,\vartheta)][1+\hat{E}_{s}(u,\vartheta)\log(1+t)]}\\\notag
\leq \frac{1}{[1+\hat{E}_{s}(u,\vartheta)\log(1+t)]}\frac{|\hat{Q}_{0,s}(t,u,\vartheta)|}{(1-|\hat{Q}_{0,s}(t,u,\vartheta)|)}
\end{align*}
therefore, from (13.4), taking $\delta_{0}$ suitably small, the second term on the right of (13.7) is bounded in absolute value by:
\begin{equation}
 \frac{1}{(1+t)[1+\log(1+t)]}\frac{|\hat{Q}_{0,s}(t,u,\vartheta)|}{(1-|\hat{Q}_{0,s}(t,u,\vartheta)|)}\leq \frac{C\delta_{0}}{(1+t)^{2}}
\end{equation}
Also, by (13.5) and (13.4), for small $\delta_{0}$, the third term on the right in (13.7) is bounded in absolute value by:
\begin{equation}
 C\delta_{0}\frac{[1+\log(1+t)]}{(1+t)^{2}}
\end{equation}
Combining (13.9)-(13.11), we conclude in Case 2:
\begin{equation}
 \frac{1}{\hat{\mu}_{s}}(\frac{\partial\hat{\mu}_{s}}{\partial t})_{+}\leq \frac{1}{(1+t)[1+\log(1+t)]}
+C\delta_{0}\frac{[1+\log(1+t)]}{(1+t)^{2}}
\end{equation}

    Finally, we consider Case 3. In this case we define:
\begin{equation}
 t_{1}(u,\vartheta)=e^{-\frac{1}{2\hat{E}_{s}(u,\vartheta)}}-1
\end{equation}
    There are two subcases to consider.
    
    Subcase 3a: $t\leq t_{1}(u,\vartheta)$,    Subcase 3b: $t>t_{1}(u,\vartheta)$

    In Subcase 3a we have:
\begin{align*}
 \log(1+t)\leq -\frac{1}{2\hat{E}_{s}(u,\vartheta)}
\end{align*}
hence:
\begin{equation}
 1+\hat{E}_{s}(u,\vartheta)\log(1+t)\geq\frac{1}{2}
\end{equation}
From (13.4), for suitably small $\delta_{0}$:
\begin{equation}
 |\hat{Q}_{0,s}(t,u,\vartheta)|\leq\frac{1}{4}\quad :\textrm{for all}\quad t\in[0,s]
\end{equation}
it follows that:
\begin{equation}
 \hat{\mu}_{s}(t,u,\vartheta)=1+\hat{E}_{s}(u,\vartheta)\log(1+t)+\hat{Q}_{0,s}(t,u,\vartheta)\geq \frac{1}{4}
\end{equation}
hence, from (13.2):
\begin{align}
 \frac{1}{\hat{\mu}_{s}}(\frac{\partial\hat{\mu}_{s}}{\partial t})_{+}\leq 4(\frac{\hat{E}_{s}(u,\vartheta)}{(1+t)}+\hat{Q}_{1,s}(t,u,\vartheta))_{+}
\leq 4|\hat{Q}_{1,s}(t,u,\vartheta)|\leq C\delta_{0}\frac{[1+\log(1+t)]}{(1+t)^{2}}
\end{align}
by (13.5).

    In Subcase 3b we have, by (13.5):
\begin{align}
 \frac{\partial\hat{\mu}_{s}}{\partial t}(t,u,\vartheta)=\frac{\hat{E}_{s}(u,\vartheta)}{(1+t)}+\hat{Q}_{1,s}(t,u,\vartheta)
\leq \frac{1}{(1+t)}(\hat{E}_{s}(u,\vartheta)+C\delta_{0}\frac{[1+\log(1+t)]}{(1+t)})
\end{align}
Let us set 
\begin{equation}
 \tau=\log(1+t)\quad\textrm{so that}\quad t=e^{\tau}-1
\end{equation}
The mapping $t\mapsto\tau$ is an orientation preserving diffeomorphsim of $[0,\infty)$ onto itself. Also,
\begin{equation}
 \frac{[1+\log(1+t)]}{(1+t)}=e^{-\tau}(1+\tau):=f(\tau)
\end{equation}
is a decreasing function of $\tau$. Thus $t>t_{1}(u,\vartheta)$ corresponds to $\tau>\tau_{1}(u,\vartheta)$, where:
\begin{equation}
 \tau_{1}(u,\vartheta)=-\frac{1}{2\hat{E}_{s}(u,\vartheta)}
\end{equation}
and to:
\begin{equation}
 f(\tau)<f(\tau_{1}(u,\vartheta))=(1-\frac{1}{2\hat{E}_{s}(u,\vartheta)})e^{\frac{1}{2\hat{E}_{s}(u,\vartheta)}}
\end{equation}
It then follows from (13.18) that in Subcase 3b:
\begin{align}
 \frac{\partial\hat{\mu}_{s}}{\partial t}(t,u,\vartheta)<\frac{1}{(1+t)}\{\hat{E}_{s}(u,\vartheta)+C\delta_{0}
(1-\frac{1}{2\hat{E}_{s}(u,\vartheta)})e^{\frac{1}{2\hat{E}_{s}(u,\vartheta)}}\}\\\notag
=\frac{\hat{E}_{s}(u,\vartheta)}{(1+t)}(1-2C\delta_{0}g(x))
\end{align}
where we set:
\begin{equation}
 x=-2\hat{E}_{s}(u,\vartheta)>0\quad\textrm{and}\quad g(x)=\frac{1}{x}(1+\frac{1}{x})e^{-\frac{1}{x}} 
\end{equation}
Now, $0<x\leq C\delta_{0}$, the function $g(x)$ is bounded on $(0,1]$, and we have:
\begin{equation}
 \lim_{x\rightarrow0_{+}}g(x)=0
\end{equation}
It follows that, for small $\delta_{0}$:
\begin{align*}
 1-2C\delta_{0}g(x)>0
\end{align*}
hence from (13.23):
\begin{equation}
 \frac{\partial\hat{\mu}_{s}}{\partial t}(t,u,\vartheta)<0
\end{equation}
Therefore in Subcase 3b:
\begin{equation}
 \frac{1}{\hat{\mu}_{s}}(\frac{\partial\hat{\mu}_{s}}{\partial t})_{+}=0
\end{equation}
In view of (13.17) and (13.27), we conclude that in Case 3:
\begin{equation}
 \frac{1}{\hat{\mu}_{s}}(\frac{\partial\hat{\mu}_{s}}{\partial t})_{+}\leq C\delta_{0}\frac{[1+\log(1+t)]}{(1+t)^{2}}
\end{equation}
     Combining (13.6), (13.12), (13.28), we conclude that in $W^{s}_{\epsilon_{0}}$:
\begin{equation}
 \frac{1}{\hat{\mu}_{s}}(\frac{\partial\hat{\mu}_{s}}{\partial t})_{+}\leq
\frac{1}{(1+t)[1+\log(1+t)]}+A(t)
\end{equation}
where:
\begin{equation}
 A(t)=C\delta_{0}\frac{[1+\log(1+t)]}{(1+t)^{2}}
\end{equation}
So the proposition follows. $\qed$

$\textbf{Proposition 13.2}$ Let the assumptions of Proposition 12.10 hold for $m=2, l=1$. Then on $W^{s}_{\epsilon_{0}}$ we have:
\begin{align*}
 \mu^{-1}(T\mu)_{+}\leq B_{s}(t)
\end{align*}
Here
\begin{align*}
 B_{s}(t)=C\sqrt{\delta_{0}}\frac{(1+\tau)}{\sqrt{\sigma-\tau}}+C\delta_{0}(1+\tau)
\end{align*}
where $\tau=\log(1+t), \sigma=\log(1+s)$. We have:
\begin{align*}
 \int_{0}^{s}(1+t)^{-2}[1+\log(1+t)]^{4}B_{s}(t)dt\leq C\sqrt{\delta_{0}}
\end{align*}
where $C$ is a constant independent of $s$. Moreover, by $\textbf{C1}$, which we have just proved, $\textbf{C2}$ holds on $W^{s}_{\epsilon_{0}}$.

$Proof$: We shall again use Proposition 8.6. Denote:
\begin{equation}
 \hat{E}_{s,m}=\min_{(u,\vartheta)\in[0,\epsilon_{0}]\times S^{2}}\hat{E}_{s}(u,\vartheta)
\end{equation}
There are two cases, according as to whether $\hat{E}_{s,m}\geq0$, or $\hat{E}_{s,m}<0$. The first case is easier to treat. In the second case, we set:
\begin{equation}
 \hat{E}_{s,m}=-\delta_{1},\quad \delta_{1}>0
\end{equation}
We define
\begin{align}
 \mathcal{V}'_{s-}=\{(u,\vartheta)\in[0,\epsilon_{0}]\times S^{2}\quad:\quad\hat{E}_{s}(u,\vartheta)<-\frac{\delta_{1}}{2}\}
\end{align}
We denote by $\mathcal{U}'_{s,t-}$ the corresponding open subset of $\Sigma_{t}^{\epsilon_{0}}$, namely the set of points on $\Sigma_{t}^{\epsilon_{0}}$ 
with acoustical coordinates $(t,u,\vartheta)$, such that $(u,\vartheta)\in\mathcal{V}'_{s-}$. Now, from Proposition 8.6,
\begin{equation}
 \hat{\mu}_{s}(s,u,\vartheta)=1+\hat{E}_{s}(u,\vartheta)\log(1+s)
\end{equation}
hence, with $\sigma=\log(1+s)$:
\begin{equation}
 0\leq\min_{(u,\vartheta)\in[0,\epsilon_{0}]\times S^{2}}\hat{\mu}_{s}(s,u,\vartheta)=1-\delta_{1}\sigma
\end{equation}
Consider a point $p\in\Sigma_{t}^{\epsilon_{0}}\setminus\mathcal{U}'_{s,t-}$. The acoustical coordinates of $p$ are $(t,u, \vartheta)$ with 
$(u,\vartheta)\in([0,\epsilon_{0}]\times S^{2})\setminus\mathcal{V}'_{s-}$. With $\tau=\log(1+t)$ we then have, by Proposition 8.6 and (13.35),
since $\tau\leq \sigma$:
\begin{align}
 \hat{\mu}_{s}(t,u,\vartheta)=1+\hat{E}_{s}(u,\vartheta)\log(1+t)+\hat{Q}_{0,s}(t,u,\vartheta)\\\notag
\geq1-\frac{1}{2}\delta_{1}\tau-|\hat{Q}_{0,s}(t,u,\vartheta)|\geq \frac{1}{2}-|\hat{Q}_{0,s}(t,u,\vartheta)|\geq \frac{1}{4}
\end{align}
by (13.15). Since
\begin{equation}
 \mu(p)=\mu_{[1],s}(u,\vartheta)\hat{\mu}_{s}(t,u,\vartheta)
\end{equation}
and by Proposition 8.6:
\begin{equation}
 \inf_{(u,\vartheta)\in[0,\epsilon_{0}]\times S^{2}}\mu_{[1],s}(u,\vartheta)\geq\frac{1}{2}
\end{equation}
it follows that:
\begin{equation}
 \mu(p)\geq\frac{1}{8}
\end{equation}
In the case that $\hat{E}_{s,m}\geq 0$, by Proposition 8.6 and (13.4) we have in this case:
\begin{equation}
 \hat{\mu}_{s}(t,u,\vartheta)\geq1+\hat{Q}_{0,s}(t,u,\vartheta)\geq 1-|\hat{Q}_{0,s}(t,u,\vartheta)|\geq \frac{3}{4}
\end{equation}
hence by (13.37) and (13.38):
\begin{equation}
 \mu(p)\geq \frac{3}{8}
\end{equation}
According to Proposition 12.10 with $m=2, l=1, \textbf{M}_{[2,2]}$ holds on $W^{s}_{\epsilon_{0}}$. In particular, we have:
\begin{equation}
 \sup_{\Sigma_{t}^{\epsilon_{0}}}|T\mu|\leq C\delta_{0}[1+\log(1+t)]\quad:\textrm{for all}\quad t\in[0,s]
\end{equation}
In view of (13.39), (13.41) and (13.42), it follows that:
\begin{equation}
 (\mu^{-1}(T\mu)_{+})(p)\leq C\delta_{0}[1+\log(1+t)]
\end{equation}
for any $p\in\Sigma_{t}^{\epsilon_{0}}, t\in[0,s]$, in the case that $\hat{E}_{s,m}\geq 0$, for any $p\in\Sigma_{t}^{\epsilon_{0}}
\setminus\mathcal{U}'_{s,t-}, t\in[0,s]$, in the case that $\hat{E}_{s,m}<0$.

    It remains for us to consider the case where $\hat{E}_{s,m}<0$ and $p\in\mathcal{U}'_{s,t-}$. Let us set:
\begin{equation}
 \tau_{1}=\frac{1}{2\delta_{1}},\quad t_{1}=e^{\tau_{1}}-1
\end{equation}
There are two subcases to consider, according as to whether $t\leq t_{1}$ or $t>t_{1}$. In the first subcase, 
which corresponds to $\tau\leq\tau_{1}$, we have, by Proposition 8.6 and (13.15):
\begin{align}
 \hat{\mu}_{s}(t,u,\vartheta)=1+\hat{E}_{s}(u,\vartheta)\log(1+t)+\hat{Q}_{0,s}(t,u,\vartheta)\\\notag
\geq1-\delta_{1}\tau-|\hat{Q}_{0,s}(t,u,\vartheta)|\geq\frac{1}{2}-|\hat{Q}_{0,s}(t,u,\vartheta)|\geq \frac{1}{4}
\end{align}
 By (13.37), (13.38), the lower bound (13.39) follows, which, together with (13.42) yields:
\begin{equation}
 (\mu^{-1}(T\mu)_{+})(p)\leq C\delta_{0}[1+\log(1+t)]
\end{equation}
for $p\in\mathcal{U}'_{s,t-}$, in the subcase $t\leq t_{1}$.

    Finally, we must consider the subcase $p\in\mathcal{U}'_{s,t-},t>t_{1}$. If $(T\mu)(p)\leq 0$, then $((T\mu))_{+}(p)=0$ and 
the estimate is trivial. Otherwise $(T\mu)(p)>0$, and we consider the integral curve of $T$, on $\Sigma_{t}^{\epsilon_{0}}$, through $p$.
The intersection of $\mathcal{U}'_{s,t-}$ with this integral curve is an open subset of the integral curve, corresponding to an open subset
$\mathcal{J}$ of the parameter interval $(0,\epsilon_{0}]$ (the integral curves of $T$ are parametrized by $u$). Thus $\mathcal{J}$ is a 
countable union of open intervals. The point $0$ is not a boundary point of $\mathcal{J}$. Let the point $p$ correspond to the parameter
value $u_{*}$. Then $u_{*}\in\mathcal{J}$. Consider the interval $\mathcal{I}$, the component of $\mathcal{J}$ to which $u_{*}$ belongs.
Now $\mathcal{I}$ is of the form $\mathcal{I}=(a,b), a>0, b<\epsilon_{0}$, if $\mathcal{I}$ is not the rightmost interval, while either
$\mathcal{I}=(a,\epsilon_{0}]$ or $\mathcal{I}=(a,\epsilon_{0}), a>0$, if $\mathcal{I}$ is the rightmost interval. Consider the function $\mu$
along the given integral curve as a function of $u$. Then $T\mu$ is the function $\frac{d\mu}{du}$, and $T^{2}\mu$ is the function 
$\frac{d^{2}\mu}{du^{2}}$ along the same integral curve. Consider the subinterval $(a,u_{*}]$ of $\mathcal{I}$. Then either (Case 1) 
$\frac{d\mu}{du}>0$ throughout $(a,u_{*}]$, or (Case 2) there is a rightmost value of $u$ in $(a,u_{*}]$ where $\frac{d\mu}{du}=0$, and, denoting
this rightmost value by $\bar{u}$, we have $\frac{d\mu}{du}>0$ on $(\bar{u},u_{*}]$.

     In Case 1 we have:
\begin{equation}
 \mu(u_{*})>\mu(a)
\end{equation}
 Now $a\slashed{\in}\mathcal{J}$, therefore the point $q$ at parameter value $a$ along the integral curve in question belongs to 
$\Sigma_{t}^{\epsilon_{0}}\setminus\mathcal{U}'_{s,t-}$. Consequently the lower bound (13.39) holds at $q$:
\begin{equation}
 \mu(q)\geq\frac{1}{8}
\end{equation}
Since according to (13.47) $\mu(u_{*})>\mu(a)$, this implies:
\begin{equation}
 \mu(p)>\frac{1}{8}
\end{equation}
which together with (13.42) yields:
\begin{equation}
 (\mu^{-1}(T\mu)_{+})(p)\leq C\delta_{0}[1+\log(1+t)]
\end{equation}

     In Case 2 the mean value of $\frac{d^{2}\mu}{du^{2}}$ on $[\bar{u},u_{*}]$ is positive. Denoting then by:
\begin{equation}
 \gamma(t)=\sup_{\Sigma_{t}^{\epsilon_{0}}}(T)^{2}\mu
\end{equation}
we have $\gamma(t)>0$. Consider then some $u\in[\bar{u},u_{*}]$. We have:
\begin{align*}
 \frac{d\mu}{du}(u_{*})-\frac{d\mu}{du}(u)=\int_{u}^{u_{*}}\frac{d^{2}\mu}{du^{2}}(u')du'
\end{align*}
 that is:
\begin{equation}
 (T\mu)(u_{*})-(T\mu)(u)=\int_{u}^{u_{*}}((T)^{2}\mu)(u')du'\leq \gamma(t)(u_{*}-u)
\end{equation}
Therefore
\begin{equation}
 (T\mu)(u)\geq (T\mu)(u_{*})-\gamma(t)(u_{*}-u)\geq \frac{1}{2}(T\mu)(u_{*}) 
\quad:\textrm{for all}\quad u\in[u_{1},u_{*}]
\end{equation}
where:
\begin{equation}
 u_{1}=u_{*}-\frac{(T\mu)(u_{*})}{2\gamma(t)}
\end{equation}

      This shows in particular that $\bar{u}<u_{1}$, for otherwise (13.53) would apply in the case $u=\bar{u}$, yielding 
$(T\mu)(\bar{u})\geq\frac{1}{2}(T\mu)(u_{*})>0$, contradicting the fact that $(T\mu)(\bar{u})=0$. Since $T\mu>0$ throughout 
$(\bar{u}, u_{*}]$, in particular on $(\bar{u},u_{1}]$, we have:
\begin{equation}
 \mu(u_{1})>\mu(\bar{u})
\end{equation}
Moreover, by (13.53), (13.54), we have:
\begin{align}
 \mu(u_{*})-\mu(u_{1})=\int_{u_{1}}^{u_{*}}(T\mu)(u)du\geq 
\frac{1}{2}(T\mu)(u_{*})(u_{*}-u_{1})=\frac{((T\mu)(u_{*}))^{2}}{4\gamma(t)}
\end{align}
Denoting:
\begin{equation}
 \mu_{m}(t)=\min_{\Sigma_{t}^{\epsilon_{0}}}\mu
\end{equation}
we have:
\begin{equation}
 \mu(\bar{u})\geq\mu_{m}(t)
\end{equation}
Combining (13.55), (13.56) and (13.58), we conclude that:
\begin{equation}
 \mu(p)-\mu_{m}(t)\geq\frac{((T\mu)(p))^{2}}{4\gamma(t)}
\end{equation}
We thus obtain:
\begin{align}
 (\mu^{-1}(T\mu)_{+})(p)=\frac{(T\mu)(p)}{\mu(p)}\leq \frac{(T\mu)(p)}{\frac{((T\mu)(p))^{2}}{4\gamma(t)}+\mu_{m}(t)}
=2\sqrt{\gamma(t)}f_{\epsilon}(x)
\end{align}
where
\begin{equation}
 x=\frac{(T\mu)(p)}{2\sqrt{\gamma(t)}}>0,\quad\epsilon=\sqrt{\mu_{m}(t)}>0
\end{equation}
$f_{\epsilon}(x)$ is the function:
\begin{equation}
 f_{\epsilon}(x)=\frac{x}{x^{2}+\epsilon^{2}}
\end{equation}
Now, the function $f_{\epsilon}(x)$ achieves its maximum at $x=\epsilon$ and the maximum value is $\frac{1}{2\epsilon}$.
Thus:
\begin{equation}
 f_{\epsilon}(x)\leq \frac{1}{2\epsilon}
\end{equation}
hence from (13.60):
\begin{equation}
 (\mu^{-1}(T\mu)_{+})(p)\leq\frac{\sqrt{\gamma(t)}}{\sqrt{\mu_{m}(t)}}
\end{equation}

     Consider now $\mu_{m}(t)$. From Proposition 8.6 we have, by (13.37), (13.38):
\begin{align}
 \mu(t,u,\vartheta)\geq\frac{1}{2}\hat{\mu}_{s}(t,u,\vartheta)
=\frac{1}{2}(1+\hat{E}_{s}(u,\vartheta)\log(1+t)+\hat{Q}_{0,s}(t,u,\vartheta))\\\notag
\geq\frac{1}{2}\{1-\delta_{1}\tau-|\hat{Q}_{0,s}|\}
\end{align}
Since $\delta_{1}\sigma\leq 1$, this implies:
\begin{equation}
 \mu_{m}(t)\geq\frac{1}{2}\{\delta_{1}(\sigma-\tau)-|\hat{Q}_{0,s}|\}
\end{equation}
From Proposition 8.6 we have:
\begin{equation}
 |\hat{Q}_{0,s}(t,u,\vartheta)|\leq C\delta_{0}b(t,s)
\end{equation}
where:
\begin{equation}
 b(t,s)=\frac{[1+\log(1+t)]}{(1+t)}\frac{(s-t)}{(1+s)}
\end{equation}
According to (8.288):
\begin{equation}
 0\leq b(t,s)\leq \frac{(1+\tau)}{e^{\tau}}(\sigma-\tau)
\end{equation}
Substituting (13.67), (13.69), in (13.66) we obtain:
\begin{equation}
 \mu_{m}(t)\geq \frac{1}{2}\delta_{1}(1-C\frac{\delta_{0}}{\delta_{1}}\frac{(1+\tau)}{e^{\tau}})(\sigma-\tau)
\end{equation}
Since we are considering the subcase $\tau>\tau_{1}$, we have:
\begin{align*}
 \frac{(1+\tau)}{e^{\tau}}<\frac{(1+\tau_{1})}{e^{\tau_{1}}}=\frac{(1+\frac{1}{2\delta_{1}})}{e^{\frac{1}{2\delta_{1}}}}
\end{align*}
The factor 
\begin{align*}
 \frac{1}{\delta_{1}}(1+\frac{1}{2\delta_{1}})e^{-\frac{1}{2\delta_{1}}}
\end{align*}
is bounded and tends to zero as $\delta_{1}\rightarrow0$. Therefore, if $\delta_{0}$ is suitably small we have:
\begin{align*}
 1-C\frac{\delta_{0}}{\delta_{1}}\frac{(1+\tau)}{e^{\tau}}\geq \frac{1}{2}
\end{align*}
Since also $2\delta_{1}\geq\frac{1}{\tau}$, it follows from (13.70) that:
\begin{equation}
 \mu_{m}(t)\geq \frac{1}{8}\frac{(\sigma-\tau)}{(1+\tau)}
\end{equation}
On the other hand, by $\textbf{M}_{[2,2]}$ we have:
\begin{equation}
 \gamma(t)\leq C\delta_{0}(1+\tau)
\end{equation}
In view of (13.71), (13.72), we conclude from (13.64):
\begin{equation}
 (\mu^{-1}(T\mu)_{+})(p)\leq C\sqrt{\delta_{0}}\frac{(1+\tau)}{\sqrt{\sigma-\tau}}
\end{equation}
This is the bound in Case 2. Combining this with (13.50), we conclude that for any $p\in\mathcal{U}'_{s,t-}, t\in(t_{1},s]$, we have:
\begin{equation}
 (\mu^{-1}(T\mu)_{+})(p)\leq C\sqrt{\delta_{0}}\frac{(1+\tau)}{\sqrt{\sigma-\tau}}+C\delta_{0}(1+\tau)
\end{equation}
Combining this with (13.46) and (13.43), we conclude that, in general:
\begin{equation}
 \mu^{-1}(T\mu)_{+}\leq B_{s}(t)
\end{equation}
where:
\begin{equation}
 B_{s}(t)=C\sqrt{\delta_{0}}\frac{(1+\tau)}{\sqrt{\sigma-\tau}}+C\delta_{0}(1+\tau)
\end{equation}
We have:
\begin{equation}
 \int_{0}^{s}(1+t)^{-2}[1+\log(1+t)]^{4}B_{s}(t)dt=\int_{0}^{\sigma}(1+\tau)^{4}B_{s}(t)e^{-\tau}d\tau\notag
\end{equation}
and:
\begin{align}
 \int_{0}^{\frac{\sigma}{2}}\frac{(1+\tau)^{5}}{\sqrt{\sigma-\tau}}e^{-\tau}d\tau\leq C\int_{0}^{\frac{\sigma}{2}}
\frac{(1+\tau)^{5}}{\sqrt{\tau}}e^{-\tau}d\tau\leq C\int_{0}^{\infty}\frac{(1+\tau)^{5}}{\sqrt{\tau}}e^{-\tau}d\tau=C\quad
\textrm{independent of}\quad s\notag
\end{align}
while:
\begin{align*}
 \int_{\frac{\sigma}{2}}^{\sigma}\frac{(1+\tau)^{5}}{\sqrt{\sigma-\tau}}e^{-\tau}d\tau\leq (1+\sigma)^{5}e^{-\frac{\sigma}{2}}\int_{\frac{\sigma}{2}}^{\sigma}
\frac{d\tau}{\sqrt{\sigma-\tau}}=\sqrt{2\sigma}(1+\sigma)^{5}e^{-\frac{\sigma}{2}}\leq C\quad\textrm{independent of}\quad s
\end{align*}
It follows that:
\begin{align*}
 \int_{0}^{s}(1+t)^{-2}[1+\log(1+t)]^{4}B_{s}(t)dt\leq C\sqrt{\delta_{0}}
\end{align*}
where $C$ is a constant independent of $s$.

     Consider now $\mu^{-1}(L\mu+\underline{L}\mu)_{+}$. Since $\underline{L}=\eta^{-2}\mu L+2T$, we have:
\begin{equation}
 \mu^{-1}(L\mu+\underline{L}\mu)_{+}\leq \eta^{-2}|L\mu|+\mu^{-1}(L\mu)_{+}+2\mu^{-1}(T\mu)_{+}
\end{equation}
By the propagation equation for $\mu$, the first term on the right in (13.77) is bounded by $C\delta_{0}(1+t)^{-1}$.
The second term on the right is bounded by Proposition 13.1, but we treat Case 2 in a different manner, namely, we use 
(13.8) to estimate the first term on the right of (13.7) by:
\begin{align*}
 C\delta_{0}(1+t)^{-1}
\end{align*}
We then obtain in Case 2 ($\hat{E}_{s}(u,\vartheta)>0$) the bound
\begin{equation}
 \frac{1}{\hat{\mu}_{s}}(\frac{\partial\hat{\mu}_{s}}{\partial t})_{+}\leq \frac{C\delta_{0}}{(1+t)}
\end{equation}
so we have:
\begin{equation}
 \mu^{-1}(L\mu)_{+}\leq C\delta_{0}(1+t)^{-1}
\end{equation}
Together with (13.75) this yields, through (13.77):
\begin{equation}
 \mu^{-1}(L\mu+\underline{L}\mu)_{+}\leq B'_{s}(t)\quad:\textrm{on}\quad W^{s}_{\epsilon_{0}}
\end{equation}
where:
\begin{equation}
 B'_{s}(t)=2B_{s}(t)+C\delta_{0}(1+t)^{-1}
\end{equation}
Noting that:
\begin{equation}
 \int_{0}^{s}(1+t)^{-2}[1+\log(1+t)]^{4}B'_{s}(t)dt\leq C\sqrt{\delta_{0}}
\end{equation}
where $C$ is a constant independent of $s$, we conclude that $\textbf{C2}$ holds on $W^{s}_{\epsilon_{0}}$. 
So Proposition 13.2 holds. $\qed$

$\textbf{Proposition 13.3}$ Let the assumptions of Proposition 13.1 hold and let $\mathcal{U}$ be the region:
\begin{align*}
 \mathcal{U}=\{x\in W^{*}_{\epsilon_{0}}\quad:\quad \mu<\frac{1}{4}\}
\end{align*}
Then there is a positive constant $C$ independent of $s$ such that in $\mathcal{U}\bigcap W^{s}_{\epsilon_{0}}$ we have:
\begin{align*}
 L\mu\leq -C^{-1}(1+t)^{-1}[1+\log(1+t)]^{-1}
\end{align*}
that is, property $\textbf{C3}$ holds on $W^{s}_{\epsilon_{0}}$.

$Proof$. Consider a point $p\in\mathcal{U}\bigcap\Sigma_{t}^{\epsilon_{0}}, t\in[0,s]$, and let its acoustical coordinates be $(t,u,\vartheta)$.
First we show that $\hat{E}_{s}(u,\vartheta)<0$. For, otherwise we have, from Proposition 8.6, for suitably small $\delta_{0}$:
\begin{align*}
 \hat{\mu}_{s}(t,u,\vartheta)\geq 1+\hat{Q}_{0,s}(t,u,\vartheta)\geq 1-C\delta_{0}\geq\frac{1}{2}
\end{align*}
hence:
\begin{align*}
 \mu(t,u,\vartheta)=\mu_{[1],s}(u,\vartheta)\hat{\mu}_{s}(t,u,\vartheta)\geq\frac{1}{4}
\end{align*}
contradicts with the fact that $p\in\mathcal{U}$. Next we show that $\log(1+t)=\tau>\tau_{2}$, where:
\begin{equation}
 \tau_{2}=\frac{1}{4\delta},\quad:\quad \delta=-\hat{E}_{s}(u,\vartheta)>0
\end{equation}
For, otherwise we have, from Proposition 8.6, for suitably small $\delta_{0}$:
\begin{align*}
 \hat{\mu}_{s}(t,u,\vartheta)=1-\delta\tau+\hat{Q}_{0,s}(t,u,\vartheta)\geq\frac{3}{4}-C\delta_{0}\geq\frac{1}{2}
\end{align*}
hence again $\mu(t,u,\vartheta)\geq\frac{1}{4}$ in contradiction with the fact that $p\in\mathcal{U}$.

    Now, from Proposition 8.6:
\begin{equation}
 L\mu=\mu_{[1],s}\frac{\partial\hat{\mu}_{s}}{\partial t}
\end{equation}
and:
\begin{equation}
 \frac{\partial\hat{\mu}_{s}}{\partial t}(t,u,\vartheta)=-\frac{\delta}{(1+t)}+\hat{Q}_{1,s}(t,u,\vartheta)
\end{equation}
From (13.5):
\begin{equation}
 \frac{\partial\hat{\mu}_{s}}{\partial t}\leq -\frac{\delta}{(1+t)}\{1-\frac{C\delta_{0}}{\delta}\frac{(1+\tau)}{e^{\tau}}\}
\end{equation}
Now $\tau>\tau_{2}$ implies:
\begin{align*}
 \frac{(1+\tau)}{e^{\tau}}<\frac{(1+\tau_{2})}{e^{\tau_{2}}}=(1+\frac{1}{4\delta})e^{-\frac{1}{4\delta}}
\end{align*}
Since the function
\begin{align*}
 \frac{1}{\delta}(1+\frac{1}{4\delta})e^{-\frac{1}{4\delta}}
\end{align*}
is bounded and tends to zero as $\delta\rightarrow0$, it follows that:
\begin{align*}
 \frac{C\delta_{0}}{\delta}\frac{(1+\tau)}{e^{\tau}}\leq\frac{1}{2}
\end{align*}
provided that $\delta_{0}$ is suitably small. Hence, from (13.86):
\begin{equation}
 \frac{\partial\hat{\mu}_{s}}{\partial t}(t,u,\vartheta)\leq -\frac{1}{2}\frac{\delta}{(1+t)}
\end{equation}
Since $\log(1+t)=\tau>\tau_{2}=\frac{1}{4\delta}$, this implies:
\begin{equation}
 \frac{\partial\hat{\mu}_{s}}{\partial t}(t,u,\vartheta)\leq -\frac{1}{8}\frac{1}{(1+t)[1+\log(1+t)]}
\end{equation}
Then by (13.84) and the fact that, by Proposition 8.6, $\mu_{[1],s}\geq\frac{1}{2}$, it follows that:
\begin{equation}
 L\mu\leq-C^{-1}(1+t)^{-1}[1+\log(1+t)]^{-1}\quad\textrm{with}\quad C=16
\end{equation}
So the proposition follows. $\qed$

     Next, we shall define a function $\omega$, as needed in Chapter 5.

$\textbf{Proposition 13.4}$ Let the assumptions of Proposition 12.6 hold and let the initial data satisfy (13.1). Let us define:
\begin{align*}
 \omega=2(1+t)
\end{align*}
Then the function $\omega$ has the properties $\textbf{D1},\textbf{D2},\textbf{D3},\textbf{D4}$, and $\textbf{D5}$, required in Chapter 5.

$Proof$: $\textbf{D1}$ is obvious. We have:
\begin{equation}
 L\omega=2,\quad T\omega=0,\quad \slashed{d}\omega=0
\end{equation}
By the first of (13.90):
\begin{align}
 L\omega-\frac{\omega}{(1-u+t)}=2[1-\frac{(1+t)}{(1-u+t)}]=\frac{-2u}{1-u+t}
\end{align}
hence:
\begin{equation}
 |L\omega-\frac{\omega}{1-u+t}|\leq\frac{C}{1+t}
\end{equation}
On the other hand, by (5.16):
\begin{equation}
 \nu=\frac{1}{2}(\textrm{tr}\chi+L\log\Omega)=\frac{1}{1-u+t}+\frac{1}{2}(\textrm{tr}\chi^{\prime}+L\log\Omega)
\end{equation}
hence, by Proposition 12.6:
\begin{equation}
 |\nu-\frac{1}{1-u+t}|\leq C\delta_{0}\frac{[1+\log(1+t)]}{(1+t)^{2}}
\end{equation}
Then by (13.92) and (13.94):
\begin{equation}
 |L\omega-\nu\omega|\leq C(1+t)^{-1}[1+\log(1+t)]
\end{equation}
This implies $\textbf{D2}$. By the first two of (13.90):
\begin{equation}
 \underline{L}\omega=2\eta^{-1}\kappa
\end{equation}
hence, by Proposition 12.1:
\begin{equation}
 |\underline{L}\omega|\leq C[1+\log(1+t)]
\end{equation}
This implies $\textbf{D3}$. $\textbf{D4}$ is obvious from the third of (13.90).

     It remains for us to verify $\textbf{D5}$. To compute $\mu\Box_{\tilde{g}}\omega$ we appeal to (3.155) with $f$ replaced by $\omega$.
By the third of (13.90), this reduces to:
\begin{equation}
 -\mu\Omega\Box_{\tilde{g}}\omega=L(\underline{L}\omega)+\nu\underline{L}\omega+\underline{\nu}L\omega
\end{equation}
From (13.96):
\begin{equation}
 L(\underline{L}\omega)=2L(\eta^{-1}\kappa)
\end{equation}
Also, from (13.96) and the first of (13.90):
\begin{align}
 \nu\underline{L}\omega+\underline{\nu}L\omega=2(\nu\eta^{-1}\kappa+\underline{\nu})\\\notag
=(\eta^{-1}\kappa\textrm{tr}\chi+\textrm{tr}\underline{\chi}+\eta^{-1}\kappa L\log\Omega+\underline{L}\log\Omega)\\\notag
=2(\kappa\textrm{tr}\slashed{k}+\eta^{-1}\kappa L\log\Omega+T\log\Omega)
\end{align}
where we have used the definition of $\nu, \underline{\nu}$ and the fact:
\begin{align*}
 \eta^{-1}\kappa\chi+\underline{\chi}=2\kappa\slashed{k}
\end{align*}
We thus obtain:
\begin{equation}
 -\mu\Omega\Box_{\tilde{g}}\omega=2\{L(\eta^{-1}\kappa)+\kappa\textrm{tr}\slashed{k}+\eta^{-1}\kappa L\log\Omega+T\log\Omega\}
\end{equation}
Proposition 12.6 then implies:
\begin{equation}
 \sup_{\Sigma_{t}^{\epsilon_{0}}}(\mu|\Box_{\tilde{g}}\omega|)\leq C\delta_{0}(1+t)^{-1}
\end{equation}
This implies $\textbf{D5}$. So the Proposition follows. $\qed$

      The remaining assumptions on which Theorem 5.1 relies, assumptions $\textbf{B}$, readily follow from assumptions
$\slashed{\textbf{X}}_{[1]}$ and $\textbf{M}_{[1,2]}$, established by Proposition 12.9 with $l=1$ and 12.10 with $m=l=1$,
respectively. Proposition 10.5, 10.6, 11.5, 11.6, require the assumptions $\slashed{\textbf{X}}_{[(l+1)_{*}]}$ 
and $\textbf{M}_{\{(l+1)_{*}+1\}}$. These assumptions are established by Proposition 12.9 and 12.10 on the basis of the 
bootstrap assumptions $\textbf{E}_{\{(l+1)_{*}+2\}}, \textbf{E}^{Q}_{\{(l+1)_{*}+1\}}, \textbf{E}^{QQ}_{\{(l+1)_{*}\}}$.
We now make these assumptions more precise by setting the constant appearing on the right equal to $1$. We define, more generally,
the bootstrap assumptions $\textbf{E}^{Q...Q}_{\{q\}}$, where the length of the string of $Q$'s is $p$, for $p=0,...,q$ as follows. 
First, the bootstrap assumption $\textbf{E}_{0,0}$ holds on $W^{s}_{\epsilon_{0}}$ if for all $t\in[0,s]$:
\begin{align*}
 \textbf{E}_{0,0}\quad:\quad \max_{\alpha}\|\psi_{\alpha}\|_{L^{\infty}(\Sigma_{t}^{\epsilon_{0}})}\leq\delta_{0}(1+t)^{-1}
\end{align*}
For non-negative integers $p,m,n$ not all zero we say that the bootstrap assumption
\begin{align*}
 \textbf{E}^{Q...Q}_{m,n}
\end{align*}
holds on $W^{s}_{\epsilon_{0}}$ if for all $t\in[0,s]$:
\begin{align*}
 \textbf{E}^{Q...Q}_{m,n}\quad:\quad\max_{\alpha;i_{1}...i_{n}}\|R_{i_{n}}...R_{i_{1}}(T)^{m}(Q)^{p}\psi_{\alpha}\|_{
L^{\infty}(\Sigma_{t}^{\epsilon_{0}})}\leq\delta_{0}(1+t)^{-1}
\end{align*}
We then define $\textbf{E}^{Q...Q}_{\{q\}}$ to be the conjunction of the assumptions $\textbf{E}^{Q...Q}_{m,n}$ corresponding 
to the triangle:
\begin{align*}
 \{(m,n)\quad:\quad m,n\geq 0,\quad m+n\leq q\}
\end{align*}
Finally we define the bootstrap assumption $\textbf{E}_{\{\{k\}\}}$ to be the conjunction of the assumptions:
\begin{align*}
 \textbf{E}^{Q...Q}_{\{k-p\}}\quad:\quad p=0,...,k
\end{align*}
From now on by $\textbf{bootstrap assumption}$ we mean the assumption $\textbf{E}_{\{\{(l+1)_{*}+2\}\}}$. We have seen that
all the assumptions which we have used so far follow from the bootstrap assumption together with appropriate assumptions on 
the initial data.

\chapter{The Error Estimates Involving the Top Order Spatial Derivatives of the Acoustical Entities}

\section{The Error Terms Involving the Top Order Spatial \\Derivatives of the Acoustical Entities}
In this chapter we deal with the error estimates involving the top order spatial derivatives of the 
acoustical entities $\chi$ and $\mu$. In view of the fact that $\slashed{\mathcal{L}}_{L}\chi$ and $L\mu$ are 
expressed by the basic propagation equation in terms of $\chi$ and $\mu$ respectively, the top order derivatives of 
$\chi$ and $\mu$ are in fact the top order spatial derivatives of these entities.

     In the following, given a positive integer $n$, we denote by $\leftexp{(\alpha;I_{1}...I_{n-1})}{\psi_{n}}$ the 
$n$th order variation:
\begin{equation}
 \leftexp{(\alpha;I_{1}...I_{n-1})}{\psi_{n}}=Y_{I_{n-1}}...Y_{I_{1}}\leftexp{(\alpha)}{\psi_{1}}
\end{equation}
where:
\begin{align}
\leftexp{(\alpha)}{\psi}_{1}=\psi_{\alpha}
\end{align}
Here the indices $I_{1},...,I_{n-1}$ range over the set $\{1,2,3,4,5\}$ and we have:
\begin{equation}
 Y_{1}=T,\quad Y_{2}=Q,\quad Y_{i+2}=R_{i}\quad: i=1,2,3
\end{equation}
For $n\geq 2$ we may express (14.1) as:
\begin{equation}
 \leftexp{(\alpha;I_{1}...I_{n-1})}{\psi_{n}}=Y_{I_{n-1}}...Y_{I_{1}}\psi_{\alpha}\quad:\textrm{for}\quad n\geq2
\end{equation}
We shall require the multi-indices $(I_{1}...I_{n-1})$ to be one of the following form. First, there should be a string 
of 2's (which may be empty). Then, there should be a string of 1's (which may also be empty). Finally, there should be a 
multi-index from the set $\{3,4,5\}$.
If the string of 2's has length $p$, the string of 1's has length $m$, and the multi-index from the set $\{3,4,5\}$ has 
length $n$, then we have:
\begin{equation}
 (I_{1}...I_{p+m+n})=(2...21...1i_{1}+2...i_{n}+2)
\end{equation}
where $(i_{1}...i_{n})$ is a multi-index from the set $\{1,2,3\}$. Thus, if $p+m+n>0$ we have,
\begin{equation}
 \leftexp{(\alpha;I_{1}...I_{p+m+n})}{\psi_{p+m+n+1}}=R_{i_{n}}...R_{i_{1}}(T)^{m}(Q)^{p}\psi_{\alpha}
\end{equation}
We shall also use in the following, in place of (14.1), the simpler notation:
\begin{equation}
 \leftexp{(I_{1}...I_{n-1})}{\psi_{n}}=Y_{I_{n-1}}...Y_{I_{1}}\psi_{1}
\end{equation}
when we do not care about which of the 1st order variations we are considering.

     We have the inhomogeneous wave equations:
\begin{equation}
 \Box_{\tilde{g}}\leftexp{(I_{1}...I_{n-1})}{\psi}_{n}=\leftexp{(I_{1}...I_{n-1})}{\rho}_{n}
\end{equation}
Also as in Chapter 7, we define:
\begin{equation}
 \leftexp{(I_{1}...I_{n-1})}{\tilde{\rho}}_{n}=\Omega^{2}\mu\leftexp{(I_{1}...I_{n-1})}{\rho}_{n}
\end{equation}
According to (7.30), for $n\geq2$ we have the recursion formula:
\begin{equation}
 \leftexp{(I_{1}...I_{n-1})}{\tilde{\rho}}_{n}=(Y_{I_{n-1}}+\leftexp{(Y_{I_{n-1}})}{\delta})\leftexp{(I_{1}...I_{n-2})}
{\tilde{\rho}}_{n-1}
+\leftexp{(Y_{I_{n-1}};I_{1}...I_{n-2})}{\sigma}_{n-1}
\end{equation}
and we have:
\begin{equation}
 \tilde{\rho}_{1}=0
\end{equation}
In (14.10) the function $\leftexp{(Y_{I_{n-1}})}{\delta}$ is defined by (7.33) and depends only on the commutation 
field $Y_{I_{n-1}}$, while the function $\leftexp{(Y_{I_{n-1}};I_{1}...I_{n-2})}{\sigma}_{n-1}$ is defined by (7.31)
 and (7.32) and depends on $Y_{I_{n-1}}$ as well as $\leftexp{(I_{1}...I_{n-2})}{\psi}_{n-1}$. We shall apply 
Proposition 8.2 to (14.10). Replacing $n$ by $n+1$ we have:
\begin{equation}
  \leftexp{(I_{1}...I_{n})}{\tilde{\rho}}_{n+1}=(Y_{I_{n}}+\leftexp{(Y_{I_{n}})}{\delta})\leftexp{(I_{1}...I_{n-1})}
{\tilde{\rho}}_{n}+\leftexp{(Y_{I_{n}};I_{1}...I_{n-1})}{\sigma}_{n}
\end{equation}
for $n=1,2,...$. Applying then Proposition 8.2 yields:
\begin{equation}
  \leftexp{(I_{1}...I_{n})}{\tilde{\rho}}_{n+1}=\sum_{m=0}^{n-1}(Y_{I_{n}}+\leftexp{(Y_{I_{n}})}{\delta})...
(Y_{I_{n-m+1}}+\leftexp{(Y_{I_{n-m+1}})}{\delta})\leftexp{(Y_{I_{n-m}};I_{1}...I_{n-m-1})}{\sigma_{n-m}}
\end{equation}
or, setting $n=l+1, m=k$:
\begin{equation}
 \leftexp{(I_{1}...I_{l+1})}{\tilde{\rho}}_{l+2}=\sum_{k=0}^{l}(Y_{I_{l+1}}+\leftexp{(Y_{I_{l+1}})}{\delta})...
(Y_{I_{l-k+2}}+\leftexp{(Y_{I_{l-k+2}})}{\delta})\leftexp{(Y_{I_{l-k+1}};I_{1}...I_{l-k})}{\sigma_{l+1-k}}
\end{equation}
    We shall focus our attention on the contribution of top order spatial derivatives of the acoustical entities, i.e.
 the $l+1$st order spatial derivatives of $\chi$ and the $l+2$nd order spatial derivatives of $\mu$. These appear only 
in the $l+1$st spatial derivatives of $\leftexp{(Y)}{\tilde{\pi}}$. From (7.59) we have:
\begin{equation}
 \leftexp{(Y)}{\sigma}_{n-1}= \leftexp{(Y)}{\sigma}_{1,n-1}+ \leftexp{(Y)}{\sigma}_{2,n-1}+ \leftexp{(Y)}{\sigma}_{3,n-1}
\end{equation}
 where $\leftexp{(Y)}{\sigma}_{1,n-1}, \leftexp{(Y)}{\sigma}_{2,n-1}, \leftexp{(Y)}{\sigma}_{3,n-1}$ are given by 
(7.74)-(7.79). Of these only $\leftexp{(Y)}{\sigma}_{2,n-1}$ 
contains the spatial derivatives of $\leftexp{(Y)}{\tilde{\pi}}$.

     We conclude that the derivatives of top order, $l+1$, of the $\leftexp{(Y)}{\tilde{\pi}}$, occur only in the term 
$k=l$ in the sum in (14.14):
\begin{equation}
 (Y_{I_{l+1}}+\leftexp{(Y_{I_{l+1}})}{\delta})...(Y_{I_{2}}+\leftexp{(Y_{I_{2}})}{\delta})\leftexp{(Y_{I_{1}})}{\sigma}
_{1}
\end{equation}
more precisely, in:
\begin{equation}
 Y_{I_{l+1}}...Y_{I_{2}}\leftexp{(Y_{I_{1}})}{\sigma}_{2,1}
\end{equation}
Furthermore, we can replace $\leftexp{(Y_{I_{1}})}{\sigma_{2,1}}$ by its principal spatial derivative part:
\begin{align}
 (1/2)T(\textrm{tr}\leftexp{(Y_{I_{1}})}{\tilde{\slashed{\pi}}})L\psi_{1}-\slashed{\mathcal{L}}_{T}\leftexp{(Y_{I_{1}})}
{\tilde{Z}}\cdot\slashed{d}\psi_{1}\\\notag
-(1/2)\slashed{\textrm{div}}\leftexp{(Y_{I_{1}})}{\tilde{Z}}\underline{L}\psi_{1}-(1/2)\slashed{\textrm{div}}
\leftexp{(Y_{I_{1}})}{\underline{\tilde{Z}}}L\psi_{1}\\\notag
+(1/2)\slashed{d}\leftexp{(Y_{I_{1}})}{\tilde{\pi}}_{L\underline{L}}\cdot(\slashed{d}\psi_{1})
+\slashed{\textrm{div}}(\mu\leftexp{(Y_{I_{1}})}{\hat{\tilde{\slashed{\pi}}}})\cdot(\slashed{d}\psi_{1})
\end{align}
Here we have expressed the derivatives with respect to $\underline{L}$ in terms of the derivatives with respect to $T$. We remark that the only expressions of the form (14.17), which may contain top order spatial derivatives of acoustical entities are those for which $(I_{2}...I_{l+1})$ does not contain any $2'$s, therefore, either $p=0$, in which case $Y_{I_{1}}$ is one of $T,R_{j}\quad:\quad j=1,2,3$, or $p=1$, in which case $Y_{I_{1}}=Q$. In the following, we denote by $[\quad]_{P.A.}$, the principal acoustical part, expressed in terms of $\chi'$ and $\mu$, the difference $\chi-\chi'=(1-u+t)^{-1}\slashed{g}$ being of lower order.

    Consider first the case $p=1, Y_{I_{1}}=Q$. The components of $\leftexp{(Q)}{\tilde{\pi}}$ are given in Chapter 6. 
Noting that:
\begin{equation}
 [\nu]_{P.A.}=\frac{1}{2}\textrm{tr}\chi',\quad [\Lambda]_{P.A.}=-(\slashed{d}\mu)^{\sharp}
\end{equation}
we have:
\begin{align}
 [\leftexp{(Q)}{\tilde{\pi}}_{L\underline{L}}]_{P.A.}=0\\\notag
[\leftexp{(Q)}{\tilde{Z}}]_{P.A.}=0\\\notag
[\leftexp{(Q)}{\underline{\tilde{Z}}}]_{P.A.}=2\Omega(1+t)(\slashed{d}\mu)^{\sharp}\\\notag
[\textrm{tr}\leftexp{(Q)}{\tilde{\slashed{\pi}}}]_{P.A.}=2\Omega(1+t)\textrm{tr}\chi'\\\notag
[\leftexp{(Q)}{\hat{\tilde{\slashed{\pi}}}}]_{P.A.}=2\Omega(1+t)\hat{\chi}'
\end{align}
From (14.18) we then obtain:
\begin{align}
 [\leftexp{(Q)}{\sigma}_{2,1}]_{P.A.}=[\Omega(1+t)((T\textrm{tr}\chi'-\slashed{\Delta}\mu)L\psi_{1}
+2\mu(\slashed{\textrm{div}}\hat{\chi}^{\prime})\cdot(\slashed{d}\psi_{1}))]_{P.A.}
\end{align}
By (8.99):
\begin{equation}
 [2\slashed{\textrm{div}}\hat{\chi}'-\slashed{d}\textrm{tr}\chi']_{P.A.}=0
\end{equation}
Also by (3.125)-(3.126):
\begin{equation}
 [T\textrm{tr}\chi'-\slashed{\Delta}\mu]_{P.A.}=0
\end{equation}
In view of (14.22) and (14.23), (14.21) reduces to:
\begin{equation}
 [\leftexp{(Q)}{\sigma}_{2,1}]_{P.A.}=\Omega(1+t)\mu(\slashed{d}\textrm{tr}\chi')\cdot(\slashed{d}\psi_{1})
\end{equation}

      Consider next the case $p=0, Y_{I_{1}}=R_{j}$. By (6.8) and (6.180):
\begin{equation}
 \leftexp{(R_{j})}{\tilde{\pi}}_{L\underline{L}}=-2\Omega R_{j}\mu-2\mu R_{j}\Omega
\end{equation}
hence:
\begin{equation}
 [\leftexp{(R_{j})}{\tilde{\pi}}_{L\underline{L}}]_{P.A.}=-2\Omega R_{j}\mu
\end{equation}
From (6.67),
\begin{equation}
 [\leftexp{(R_{j})}{\slashed{\pi}}_{L}]_{P.A.}=-R_{j}\cdot\chi'
\end{equation}
Since
\begin{align*}
 \leftexp{(R_{j})}{\tilde{Z}}=\Omega\leftexp{(R_{j})}{Z}
\end{align*}
it follows that
\begin{equation}
 [\leftexp{(R_{j})}{\tilde{Z}}]_{P.A.}=-\Omega R_{j}\cdot\chi^{\prime\sharp}
\end{equation}
We have:
\begin{align*}
 \leftexp{(R_{j})}{\slashed{\pi}}_{\underline{L}}=\eta^{-1}\kappa\leftexp{(R_{j})}{\slashed{\pi}}_{L}
+2\leftexp{(R_{j})}{\slashed{\pi}}_{T}
\end{align*}
From (6.63) and the fact that
\begin{align*}
 \chi=\eta(\slashed{k}-\theta)
\end{align*}
together with (14.28), we deduce:
\begin{equation}
 [\leftexp{(R_{j})}{\tilde{\underline{Z}}}]_{P.A.}=\Omega\eta^{-1}(\kappa R_{j}\cdot\chi^{\prime\sharp}
+2\lambda_{j}(\slashed{d}\mu)^{\sharp})
\end{equation}

Next, from (6.8) we have:
\begin{equation}
 [\leftexp{(R_{j})}{\tilde{\slashed{\pi}}}]_{P.A}=[\Omega\leftexp{(R_{j})}{\slashed{\pi}}]_{P.A.}
\end{equation}
From (6.59) and (3.27):
\begin{equation}
 [\leftexp{(R_{j})}{\slashed{\pi}}]_{P.A.}=2\lambda_{j}\eta^{-1}\chi'
\end{equation}
Hence:
\begin{equation}
 [\leftexp{(R_{j})}{\tilde{\slashed{\pi}}}]_{P.A}=2\Omega\lambda_{j}\eta^{-1}\chi'
\end{equation}
Then
\begin{equation}
 [\textrm{tr}\leftexp{(R_{j})}{\tilde{\slashed{\pi}}}]_{P.A}=2\Omega\lambda_{j}\eta^{-1}\textrm{tr}\chi'
\end{equation}
and:
\begin{equation}
 [\leftexp{(R_{j})}{\hat{\tilde{\slashed{\pi}}}}]_{P.A}=2\Omega\lambda_{j}\eta^{-1}\hat{\chi}'
\end{equation}
    
     Substituting (14.26), (14.28), (14.29), (14.33), (14.34), in (14.18) we obtain:
\begin{align}
 [\leftexp{(R_{j})}{\sigma}_{2,1}]_{P.A}=[\Omega\lambda_{j}\alpha^{-1}(T\textrm{tr}\chi')L\psi_{1}
+\Omega(R_{j}\cdot\slashed{\mathcal{L}}_{T}\chi')\cdot(\slashed{d}\psi_{1})+(1/2)\Omega\slashed{\textrm{div}}
(R_{j}\cdot\chi')
\underline{L}\psi_{1}\\\notag
-(1/2)\Omega\alpha^{-1}(\kappa\slashed{\textrm{div}}(R_{j}\cdot\chi')+2\lambda_{j}(\slashed{\Delta}\mu))L\psi_{1}
-\Omega(\slashed{d}R_{j}\mu)\cdot(\slashed{d}\psi_{1})+2\Omega\lambda_{j}\kappa\slashed{\textrm{div}}\hat{\chi}'
\cdot(\slashed{d}\psi_{1})]_{P.A.}
\end{align}
Using (14.22), (14.23) and the fact by (3.125)-(3.126):
\begin{equation}
 [\slashed{\mathcal{L}}_{T}\chi'-\slashed{D}^{2}\mu]_{P.A.}=0
\end{equation}
(14.35) reduces to:
\begin{align}
 [\leftexp{(R_{j})}{\sigma}_{2,1}]_{P.A}=[\Omega\slashed{\textrm{div}}(R_{j}\cdot\chi')T\psi_{1}+
\Omega\lambda_{j}\kappa\slashed{d}\textrm{tr}\chi'\cdot(\slashed{d}\psi_{1})\\\notag
+\Omega(R_{j}\cdot\slashed{D}^{2}\mu-\slashed{d}R_{j}\mu)\cdot(\slashed{d}\psi_{1})]_{P.A.}
\end{align}
Consider $\slashed{\textrm{div}}(R_{j}\cdot\chi')$ in components in an arbitrary local frame field for $S_{t,u}$:
\begin{align*}
 \slashed{\textrm{div}}(R_{j}\cdot\chi')=\slashed{D}_{A}(R^{B}_{j}\chi'^{A}_{B})=R^{B}_{j}(\slashed{D}_{A}
\chi^{\prime A}_{B})+(\slashed{D}_{A}R^{B}_{j})\chi^{\prime A}_{B}
\end{align*}
and with $R_{jB}=\slashed{g}_{BC}R^{C}_{j}$:
\begin{align*}
 (\slashed{D}_{A}R^{B}_{j})\chi^{\prime A}_{B}=(\slashed{D}_{A}R_{jB})\chi^{\prime AB}
=\frac{1}{2}(\slashed{D}_{A}R_{jB}+\slashed{D}_{B}R_{jA})\chi^{\prime AB}
=\frac{1}{2}\leftexp{(R_{j})}{\slashed{\pi}}_{AB}\chi^{\prime AB}
\end{align*}
We conclude that:
\begin{equation}
 \slashed{\textrm{div}}(R_{j}\cdot\chi')=R_{j}\cdot\slashed{\textrm{div}}\chi'+\frac{1}{2}\textrm{tr}
(\leftexp{(R_{j})}{\slashed{\pi}}\cdot\chi')
\end{equation}
The second term is of lower order. In view of (14.22) it follows that:
\begin{equation}
 [\slashed{\textrm{div}}(R_{j}\cdot\chi')]_{P.A}=[R_{j}\cdot\slashed{\textrm{div}}\chi']_{P.A}=R_{j}\textrm{tr}\chi'
\end{equation}

Consider finally $(R_{j}\cdot\slashed{D}^{2}\mu-\slashed{d}R_{j}\mu)\cdot(\slashed{d}\psi_{1})$.
In terms of the $S_{t,u}$-tangential vectorfield $X=(\slashed{d}\psi_{1})^{\sharp}$ we have:
\begin{equation}
 R_{j}\cdot\slashed{D}^{2}\mu\cdot(\slashed{d}\psi_{1})=\slashed{D}^{2}\mu\cdot(X,R_{j})
\end{equation}
and:
\begin{align}
 \slashed{D}^{2}\mu(X,R_{j})=XR_{j}\mu-\slashed{d}\mu\cdot(\slashed{D}_{X}R_{j})\\\notag
=XR_{j}\mu-\slashed{d}\mu\cdot(\slashed{D}_{R_{j}}X)+\slashed{d}\mu\cdot[R_{j},X]
\end{align}
Therefore:
\begin{align}
 (R_{j}\cdot\slashed{D}^{2}\mu-\slashed{d}R_{j}\mu)\cdot(\slashed{d}\psi_{1})=\slashed{D}^{2}\mu\cdot(X,R_{j})
-XR_{j}\mu\\\notag
=\slashed{d}\mu\cdot\slashed{\mathcal{L}}_{R_{j}}X-\slashed{d}\mu\cdot(\slashed{D}_{R_{j}}X)
\end{align}
Obviously, (14.42) does not contain acoustical terms of order 2. It follows that:
\begin{equation}
 [\Omega(R_{j}\cdot\slashed{D}^{2}\mu-\slashed{d}R_{j}\mu)\cdot(\slashed{d}\psi_{1})]_{P.A.}=0
\end{equation}
In view of (14.39) and (14.43), (14.37) reduces finally to:
\begin{equation}
 [\leftexp{(R_{j})}{\sigma}_{2,1}]_{P.A.}=\Omega(R_{j}\textrm{tr}\chi')T\psi_{1}
+\Omega\kappa\lambda_{j}\slashed{d}\textrm{tr}\chi'\cdot(\slashed{d}\psi_{1})
\end{equation}

     Finally, we consider the case $p=0,Y_{I_{1}}=T$. The components of $\leftexp{(T)}{\tilde{\pi}}$ are given 
in Chapter 6. We have:
\begin{equation}
 \leftexp{(T)}{\tilde{Z}}=\leftexp{(T)}{\tilde{\slashed{\pi}}}_{L}\cdot\slashed{g}^{-1}=\Omega\Lambda,\quad 
\leftexp{(T)}{\tilde{\underline{Z}}}=\leftexp{(T)}{\tilde{\slashed{\pi}}}_{\underline{L}}\cdot\slashed{g}^{-1}
=\Omega\eta^{-1}\kappa\Lambda
\end{equation}
From the fact that $\underline{\chi}=-\eta^{-1}\kappa\chi+2\kappa\slashed{k}$ we have:
\begin{equation}
 [\underline{\nu}]_{P.A}=[\frac{1}{2}\textrm{tr}\underline{\chi}]_{P.A.}=-\frac{1}{2}\eta^{-1}\kappa\textrm{tr}\chi',
\quad[\hat{\underline{\chi}}]_{P.A.}=-\eta^{-1}\kappa\hat{\chi}'
\end{equation}
Using (14.19), we obtain:
\begin{align}
 [\leftexp{(T)}{\tilde{\pi}}_{L\underline{L}}]_{P.A.}=-2\Omega T\mu\\\notag
[\leftexp{(T)}{\tilde{Z}}]_{P.A.}=-\Omega(\slashed{d}\mu)\\\notag
[\leftexp{(T)}{\tilde{\underline{Z}}}]_{P.A.}=-\Omega\eta^{-1}\kappa(\slashed{d}\mu)\\\notag
[\textrm{tr}\leftexp{(T)}{\tilde{\slashed{\pi}}}]_{P.A}=-2\Omega\eta^{-1}\kappa\textrm{tr}\chi'\\\notag
[\leftexp{(T)}{\hat{\tilde{\pi}}}]_{P.A.}=-2\Omega\eta^{-1}\kappa\hat{\chi}'
\end{align}
Substituting above in (14.18) we obtain:
\begin{align}
 [\leftexp{(T)}{\sigma}_{2,1}]_{P.A.}=[-\Omega\alpha^{-1}\kappa(T\textrm{tr}\chi')L\psi_{1}+\Omega
(\slashed{d}T\mu)\cdot(\slashed{d}\psi_{1})\\\notag
+(1/2)\Omega(\slashed{\Delta}\mu)\underline{L}\psi_{1}+(1/2)\Omega\eta^{-1}\kappa(\slashed{\Delta}\mu)L\psi_{1}\\\notag
-\Omega(\slashed{d}T\mu)\cdot(\slashed{d}\psi_{1})-2\Omega\kappa^{2}\slashed{\textrm{div}}
\hat{\chi}'\cdot(\slashed{d}\psi_{1})]_{P.A}
\end{align}
Using (14.22) and (14.23), (14.48) reduces to:
\begin{equation}
 [\leftexp{(T)}{\sigma}_{2,1}]_{P.A.}=\Omega(\slashed{\Delta}\mu)T\psi_{1}
-\Omega\kappa^{2}\slashed{d}\textrm{tr}\chi'\cdot(\slashed{d}\psi_{1})
\end{equation}

      We conclude from the above that the top order variations $\leftexp{(\alpha;I_{1}...I_{l+1})}{\psi}_{l+2}$ to 
which the re-scaled sources 
$\leftexp{(\alpha;I_{1}...I_{l+1})}{\tilde{\rho}}_{l+2}$ containing the top order spatial derivatives of 
the acoustical entities correspond are as follows.

      In the case $Y_{I_{1}}=Q$, the variations:
\begin{equation}
 \leftexp{(\alpha;21...1i_{1}+2...i_{n}+2)}{\psi}_{l+2}=R_{i_{n}}...R_{i_{1}}(T)^{m}Q\psi_{\alpha}\quad:\quad m+n=l
\end{equation}
where there are $m$ 1's in the superscript. By (14.17) and (14.24) the principal acoustical part of the corresponding
 re-scaled sources is:
\begin{align}
 [\leftexp{(\alpha;21...1i_{1}+2...i_{n}+2)}{\tilde{\rho}}_{l+2}]_{P.A.}\\\notag
=[\Omega(1+t)\mu(\slashed{d}R_{i_{n}}...R_{i_{1}}(T)^{m}\textrm{tr}\chi')\cdot(\slashed{d}\psi_{\alpha})]_{P.A.}\\\notag
=\Omega(1+t)\mu(\slashed{d}R_{i_{n}}...R_{i_{1}}\textrm{tr}\chi')\cdot(\slashed{d}\psi_{\alpha})\quad:\quad m=0\\\notag
\textrm{or}\quad =\Omega(1+t)\mu(\slashed{d}R_{i_{n}}...R_{i_{1}}(T)^{m-1}\slashed{\Delta}\mu)\cdot
(\slashed{d}\psi_{\alpha})\quad:\quad m\geq 1
\end{align}
where again there are $m$ 1's in the superscript, and in the case $m\geq 1$ we have used (14.23).

     In the case $Y_{I_{1}}=R_{j}$, setting $i_{1}=j$ we have:
\begin{align*}
 I_{1}...I_{l+1}=i_{1}+2...i_{l+1}+2,
\end{align*}
the variations:
\begin{equation}
 \leftexp{(\alpha;i_{1}+2...i_{l}+2)}{\psi}_{l+2}=R_{i_{l+1}}...R_{i_{1}}\psi_{\alpha}
\end{equation}
By (14.17) and (14.44) the principal acoustical part of the corresponding re-scaled sources is:
\begin{align}
 [\leftexp{(\alpha;i_{1}+2...i_{l+1}+2)}{\tilde{\rho}}_{l+2}]_{P.A.}=\Omega(R_{i_{l+1}}...R_{i_{1}}\textrm{tr}\chi')
T\psi_{\alpha}+\Omega\kappa\lambda_{i_{1}}(\slashed{d}R_{i_{l+1}}...R_{i_{2}}\textrm{tr}\chi')
\cdot(\slashed{d}\psi_{\alpha})
\end{align}

     In the case $Y_{I_{1}}=T$, the variations:
\begin{equation}
 \leftexp{(\alpha;1...1i_{1}+2...i_{n}+2)}{\psi}_{l+2}=R_{i_{n}}...R_{i_{1}}(T)^{m+1}\psi_{\alpha}\quad:\quad m+n=l
\end{equation}
where there are $m+1$ 1's in the superscript. By (14.17) and (14.49) the principal acoustical part of the 
corresponding re-scaled sources is:
\begin{align}
 [\leftexp{(\alpha;i_{1}+2...i_{l+1}+2)}{\tilde{\rho}}_{l+2}]_{P.A.}=\Omega(R_{i_{n}}...R_{i_{1}}(T)^{m}
\slashed{\Delta}\mu)T\psi_{\alpha}
-\Omega\kappa^{2}(\slashed{d}R_{i_{n}}...R_{i_{1}}\textrm{tr}\chi')\cdot(\slashed{d}\psi_{\alpha})
\quad:\quad m=0\\\notag
\textrm{or}\quad =\Omega(R_{i_{n}}...R_{i_{1}}(T)^{m}\slashed{\Delta}\mu)T\psi_{\alpha}-
\Omega\kappa^{2}(\slashed{d}R_{i_{n}}...R_{i_{1}}(T)^{m-1}\slashed{\Delta}\mu)\cdot(\slashed{d}\psi_{\alpha})
\quad:\quad m\geq 1
\end{align}
where again there are $m+1$ 1's in the superscript, and in the case $m\geq1$ we have used (14.23).

      We see from (14.51), (14.53), (14.55), that the most difficult error integrals are which result in terms:
\begin{equation}
 \Omega(R_{i_{l+1}}...R_{i_{1}}\textrm{tr}\chi')T\psi_{\alpha}
\end{equation}
in (14.53), and
\begin{equation}
 \Omega(R_{i_{n}}...R_{i_{1}}(T)^{m}\slashed{\Delta}\mu)T\psi_{\alpha}\quad:\quad m+n=l
\end{equation}
in (14.55). These are the terms proportional to $T\psi_{\alpha}$ in (14.53) and (14.55). The other terms 
in these two expressions are proportional to 
$\slashed{d}\psi_{\alpha}$. Besides containing at least one factor of $\kappa$, which makes the estimates 
much easier in the region of small $\mu$, these terms 
have a decay factor of $(1+t)^{-2}[1+\log(1+t)]^{2}$ relative to the leading terms (14.56) and (14.57). 
The terms in (14.51) are also proportional to 
$\slashed{d}\psi_{\alpha}$, and, besides containing the factor $\mu$, these terms have a decay factor of
 $(1+t)^{-1}[1+\log(1+t)]$ relative to (14.56) ($m=0$)
and (14.57) ($m\geq 1$). So we may confine ourselves to (14.56) and (14.57).

\section{The Borderline Error Integrals}
      We recall from Chapter 7 that the error integrals corresponding to an $n$th order variation $\psi_{n}$ are 
the integrals:
\begin{equation}
 -\int_{W^{t}_{u}}\tilde{\rho}_{n}(K_{0}\psi_{n})dt'du'd\mu_{\slashed{g}}
\end{equation}
associated to $K_{0}$, and
\begin{equation}
 -\int_{W^{t}_{u}}\tilde{\rho}_{n}(K_{1}\psi_{n}+\omega\psi_{n})dt'du'd\mu_{\slashed{g}}
\end{equation}
associated to $K_{1}$.

Since:
\begin{equation}
 K_{0}\psi_{n}=\underline{L}\psi_{n}+(1+\eta^{-1}\kappa)L\psi_{n}
\end{equation}
while
\begin{equation}
 K_{1}\psi_{n}+\omega\psi_{n}=(\omega/\nu)(L\psi_{n}+\nu\psi_{n})
\end{equation}
the coefficient of $L\psi_{n}$ in (14.60) has a decay factor of $(1+t)^{-2}[1+\log(1+t)]$ relative to the coefficient 
of the same term in (14.61). 
Therefore it suffices to bound the error integral (14.59) and the contribution of $\underline{L}\psi_{n}$ in (14.60) 
to the error integral (14.58). 
This contribution is bounded 
in absolute value by:
\begin{equation}
 \int_{W^{t}_{u}}|\tilde{\rho}_{n}||\underline{L}\psi_{n}|dt'dud\mu_{\slashed{g}}
\end{equation}
We are thus to estimate the contributions of (14.56) and (14.57) to (14.62) and (14.59).

      We begin with the contribution of (14.56) to (14.62). The contribution is:
\begin{align}
 \int_{W^{t}_{\epsilon_{0}}}\Omega|R_{i_{l+1}}...R_{i_{1}}\textrm{tr}\chi'||T\psi_{\alpha}|
|\underline{L}R_{i_{l+1}}...R_{i_{1}}\psi_{\alpha}|dt'dud\mu_{\slashed{g}}\\\notag
\leq C\int_{0}^{t}\sup_{\Sigma_{t'}^{\epsilon_{0}}}(\mu^{-1}|T\psi_{\alpha}|)\|\mu R_{i_{l+1}}...R_{i_{1}}\textrm{tr}\chi'\|
_{L^{2}(\Sigma_{t'}^{\epsilon_{0}})}\|\underline{L}R_{i_{l+1}}...R_{i_{1}}\psi_{\alpha}\|_{L^{2}(\Sigma_{t'}^{\epsilon_{0}})}dt'
\end{align}
Now, in view of (8.40),(8.63),(8.335) and (8.346) we have:
\begin{align}
 \|\mu R_{i_{l+1}}...R_{i_{1}}\textrm{tr}\chi'\|_{L^{2}(\Sigma_{t}^{\epsilon_{0}})}\leq C(1+t)
\|\mu\slashed{d}R_{i_{l}}...R_{i_{1}}\textrm{tr}\chi'\|_{L^{2}(\Sigma_{t}^{\epsilon_{0}})}\\\notag
\leq C(1+t)^{2}\{\||\leftexp{(i_{1}...i_{l})}{x}_{l}(t)|\|_{L^{2}([0,\epsilon_{0}]\times S^{2})}
+\||\slashed{d}\leftexp{(i_{1}...i_{l})}{\check{f}}_{l}(t)|\|_{L^{2}([0,\epsilon_{0}]\times S^{2})}\}\\\notag
=C\{(1+t)^{2}\||\leftexp{(i_{1}...i_{l})}{x}_{l}(t)|\|_{L^{2}([0,\epsilon_{0}]\times S^{2})}+\leftexp{(i_{1}...i_{l})}{P}_{l}(t)\}
\end{align}
Substituting (8.413) we obtain:
\begin{align}
 \|\mu R_{i_{l+1}}...R_{i_{1}}\textrm{tr}\chi'\|_{L^{2}(\Sigma_{t}^{\epsilon_{0}})}\leq C(1+t)^{2}
((1+t)^{-2}\leftexp{(i_{1}...i_{l})}{P}_{l}(t)+\leftexp{(i_{1}...i_{l})}{B}_{l}(t))\\\notag
+C_{l}\delta_{0}(1+t)^{-1}[1+\log(1+t)]^{2}\int_{0}^{t}(1+t')[1+\log(1+t')]B_{l}(t')dt'
\end{align}
Here, according to (8.403):
\begin{align}
 \leftexp{(i_{1}...i_{l})}{B}_{l}(t)=C(1+t)^{-2}\{\leftexp{(i_{1}...i_{l})}{\bar{P}}^{(0)}_{l,a}(t)
+(1+t)^{-1/2}\leftexp{(i_{1}...i_{l})}{\bar{P}}^{(1)}_{l,a}(t)\}\bar{\mu}_{m}^{-a}(t)\\\notag
+C(1+t)^{-3}[1+\log(1+t)]^{2}\int_{0}^{t}(1+t')^{3}\leftexp{(i_{1}...i_{l})}{Q}_{l}(t')dt'\\\notag
+C(1+t)^{-3}[1+\log(1+t)]^{2}\|\leftexp{(i_{1}...i_{l})}{x}_{l}(0)\|_{L^{2}(\Sigma_{0}^{\epsilon_{0}})}
\end{align}
Also (see (8.405))
\begin{equation}
 B_{l}(t)=\max_{i_{1}...i_{l}}\leftexp{(i_{1}...i_{l})}{B}_{l}(t)
\end{equation}
According to (8.348) and (8.349):
\begin{align}
 \leftexp{(i_{1}...i_{l})}{\bar{P}}^{(0)}_{l,a}(t)=\sup_{t'\in[0,t]}\{\bar{\mu}^{a}_{m}(t')\leftexp{(i_{1}...i_{l})}
{P}_{l}^{(0)}(t')\}\\
\leftexp{(i_{1}...i_{l})}{\bar{P}}^{(1)}_{l,a}(t)=\sup_{t'\in[0,t]}\{(1+t')^{1/2}\bar{\mu}^{a}_{m}(t')
\leftexp{(i_{1}...i_{l})}{P}_{l}^{(1)}(t')\}
\end{align}
where $\leftexp{(i_{1}...i_{l})}{P}_{l}^{(0)}$, $\leftexp{(i_{1}...i_{l})}{P}_{l}^{(1)}$ are defined in the statement 
of Proposition 10.6 and we have:
\begin{equation}
 \leftexp{(i_{1}...i_{l})}{P}_{l}(t)\leq \leftexp{(i_{1}...i_{l})}{P}_{l}^{(0)}(t)+\leftexp{(i_{1}...i_{l})}
{P}_{l}^{(1)}(t)
\end{equation}
Note that $\leftexp{(i_{1}...i_{l})}{\bar{P}}^{(0)}_{l,a}(t)$ and $\leftexp{(i_{1}...i_{l})}{\bar{P}}^{(1)}_{l,a}(t)$ 
are non-decreasing.

     The leading contribution to the right hand side of (14.63) comes from the first term on the right in (14.65). 
 That is, from  $(1+t)^{-2}\leftexp{(i_{1}...i_{l})}{P}_{l}(t)$ and the first term in the expression (14.66) 
for $\leftexp{(i_{1}...i_{l})}{B}_{l}(t)$. So what we consider here is:
\begin{equation}
 C\{\leftexp{(i_{1}...i_{l})}{\bar{P}}_{l,a}^{(0)}(t)+(1+t)^{-1/2}\leftexp{(i_{1}...i_{l})}
{\bar{P}}^{(1)}_{l,a}(t)\}\bar{\mu}^{-a}_{m}(t)
\end{equation}
 The actual $borderline$ $contribution$ is the contribution from:
\begin{equation}
 \leftexp{(i_{1}...i_{l})}{P}_{l}^{(0)}(t)=|\ell|\sqrt{\sum_{j,\alpha}\mathcal{E}_{0}[
R_{j}R_{i_{l}}...R_{i_{1}}\psi_{\alpha}](t)}
\end{equation}

\section{Assumption $\textbf{J}$}
To estimate the borderline contributions, we shall use the assumption $\textbf{J}$ below. Recall the operators:
\begin{equation}
 S=x^{\alpha}\frac{\partial}{\partial x^{\alpha}}-1,\quad \mathring{R}_{i}=\frac{1}{2}\epsilon_{ijk}
\frac{\partial}{\partial x^{k}}\quad:i=1,2,3
\end{equation}
associated to the background Galilean spacetime.
The assumption is that there is a positive constant $C$ independent of $s$ such that in $W^{s}_{\epsilon_{0}}$:
\begin{equation}
 \textbf{J}:\quad |S\phi|,|TS\phi|\leq C\delta_{0}(1+t)^{-1},\quad |\mathring{R}_{i}\phi|,|T\mathring{R}_{i}\phi|
\leq C\delta_{0}(1+t)^{-1}\quad: i=1,2,3
\end{equation}
where $\phi$ is the wave function. We shall establish this assumption in the sequel, on the basis of the bootstrap 
assumption. We shall presently use assumption $\textbf{J}$ to derive a pointwise estimate for $T\psi_{\alpha}$ in 
terms of $L\mu$.

      Let us consider the $\Sigma_{t}$-tangential vectorfield:
\begin{equation}
 V=\sum_{i=1}^{3}(T\psi_{i})\frac{\partial}{\partial x^{i}}
\end{equation}
Recalling from Chapter 6 the Euclidean outward unit normal $N$ to the Euclidean coordinate spheres:
\begin{equation}
 N=\sum_{i=1}^{3}\frac{x^{i}}{r}\frac{\partial}{\partial x^{i}}
\end{equation}
We decompose $V$ into its components which are tangential and orthogonal, relative to $\bar{g}$, to the Euclidean 
spheres:
\begin{equation}
 V=V^{\Arrowvert}+V^{\bot},\quad V^{\bot}=\bar{g}(V,N)N,\quad |V|^{2}=|V^{\Arrowvert}|^{2}+|V^{\bot}|^{2}
\end{equation}
Consider the vector $N\times V$. We have:
\begin{equation}
 |N\times V|^{2}=|V^{\Arrowvert}|^{2}
\end{equation}
The $i$th rectangular component of $N\times V$ is given by:
\begin{align}
 r(N\times V)^{i}=\sum_{j,k=1}^{3}\epsilon_{ijk}x^{j}T(\partial_{k}\phi)
=\sum_{j,k=1}^{3}\epsilon_{ijk}\{T(x^{j}\partial_{k}\phi)-(Tx^{j})\partial_{k}\phi\}
\end{align}
Thus we have:
\begin{equation}
 r(N\times V)^{i}=T\mathring{R}_{i}\phi-\sum_{j,k=1}^{3}\epsilon_{ijk}T^{j}\psi_{k}
\end{equation}
By $\textbf{J}$, Proposition 12.1 and (12.6):
\begin{equation}
 r|N\times V|\leq C\delta_{0}(1+t)^{-1}[1+\log(1+t)]\quad:\textrm{in}\quad W^{s}_{\epsilon_{0}}
\end{equation}
From (12.14) and (14.78) we then obtain:
\begin{equation}
 \sup_{\Sigma_{t}^{\epsilon_{0}}}|V^{\Arrowvert}|\leq C\delta_{0}(1+t)^{-2}[1+\log(1+t)]
\end{equation}
By (14.75)-(14.77) and (14.82) it follows that:
\begin{align}
 \sup_{\Sigma_{t}^{\epsilon_{0}}}|T\psi_{i}-\frac{x^{i}}{r}\bar{g}(V,N)|\leq C\delta_{0}(1+t)^{-2}[1+\log(1+t)]
\end{align}
Consider next:
\begin{equation}
 S\phi=x^{\alpha}\partial_{\alpha}\phi-\phi=t\psi_{0}+\sum_{i=1}^{3}x^{i}\psi_{i}-\phi
\end{equation}
We have:
\begin{align*}
 TS\phi=tT\psi_{0}+\sum_{i=1}^{3}x^{i}T\psi_{i}+\sum_{i=1}^{3}\psi_{i}Tx^{i}-T\phi
\end{align*}
The last two terms cancel and we obtain:
\begin{equation}
 TS\phi=tT\psi_{0}+\sum_{i=1}^{3}x^{i}T\psi_{i}=tT\psi_{0}+r\bar{g}(V,N)
\end{equation}
By virtue of $\textbf{J}$ we then have:
\begin{equation}
 \sup_{\Sigma_{t}^{\epsilon_{0}}}|\frac{t}{r}T\psi_{0}+\bar{g}(V,N)|\leq C\delta_{0}(1+t)^{-2}
\end{equation}
Combining (14.83) and (14.86) we conclude that:
\begin{equation}
 \sup_{\Sigma_{t}^{\epsilon_{0}}}|T\psi_{i}+\frac{tx^{i}}{r^{2}}T\psi_{0}|\leq C\delta_{0}(1+t)^{-2}[1+\log(1+t)]
\end{equation}
Recalling the proof of Proposition 8.5, we have:
\begin{align*}
 L\mu=m+e\mu
\end{align*}
and, according to (8.209):
\begin{align*}
 m=m_{0}+m_{1},\quad m_{0}=\frac{1}{2}\ell T\psi_{0}
\end{align*}
while according to (8.214), (8.218) and (8.220):
\begin{align*}
 |m_{1}|\leq C\delta_{0}(1+t)^{-2},\quad |e|\leq C\delta_{0}(1+t)^{-2}
\end{align*}
It follows that:
\begin{equation}
 \sup_{\Sigma_{t}^{\epsilon_{0}}}|L\mu-\frac{1}{2}\ell T\psi_{0}|\leq C\delta_{0}(1+t)^{-2}[1+\log(1+t)]
\end{equation}
In view of (14.87), (12.14) and (14.88), we obtain:
\begin{align*}
 \max_{\alpha}(\mu^{-1}|T\psi_{\alpha}|)\leq C\mu^{-1}|T\psi_{0}|+C\mu^{-1}\delta_{0}(1+t)^{-2}[1+\log(1+t)]\\\notag
\leq \frac{C}{|\ell|}\{\mu^{-1}|L\mu|+C\mu^{-1}\delta_{0}(1+t)^{-2}[1+\log(1+t)]\} 
\end{align*}
Taking supremum on $\Sigma_{t}^{\epsilon_{0}}$ yields the desired estimate:
\begin{align}
 \max_{\alpha}\sup_{\Sigma_{t}^{\epsilon_{0}}}(\mu^{-1}|T\psi_{\alpha}|)\leq\frac{C}{|\ell|}
\{\sup_{\Sigma_{t}^{\epsilon_{0}}}(\mu^{-1}|L\mu|)+C\bar{\mu}_{m}^{-1}\delta_{0}(1+t)^{-2}[1+\log(1+t)]\}
\end{align}
We shall use this to estimate the borderline contribution (14.72), through (14.68), (14.71) and (14.65), to the integral 
on the right of (14.63). On the other hand, in estimating all other contributions, through (14.65), to (14.63), 
we simply use:
\begin{equation}
 \max_{\alpha}\sup_{\Sigma_{t}^{\epsilon_{0}}}(\mu^{-1}|T\psi_{\alpha}|)\leq C\bar{\mu}_{m}^{-1}\delta_{0}(1+t)^{-1}
\end{equation}
implied by $\textbf{E}_{\{1\}}$.

\section{The Borderline Estimates Associated to $K_{0}$}
\subsection{Estimates for the Contribution of (14.56)}
    Let us define, for non-negative real numbers $a$ and $p$, the quantities:
\begin{align}
 \leftexp{(i_{1}...i_{l})}{\mathcal{G}}_{0,l+2;a,p}(t)=\sup_{t'\in[0,t]}
\{[1+\log(1+t')]^{-2p}\bar{\mu}_{m}^{2a}(t')\sum_{j,\alpha}\mathcal{E}_{0}[R_{j}R_{i_{l}}...R_{i_{1}}\psi_{\alpha}](t')\}
\end{align}
These quantities are non-decreasing functions of $t$ and we have:
\begin{align}
 \sqrt{\sum_{j,\alpha}\mathcal{E}_{0}[R_{j}R_{i_{l}}...R_{i_{1}}\psi_{\alpha}](t')}
\leq\bar{\mu}_{m}^{-a}(t')[1+\log(1+t')]^{p}\sqrt{\leftexp{(i_{1}...i_{l})}{\mathcal{G}}_{0,l+2;a,p}(t)}\\\notag
\quad:\quad\textrm{for all}\quad t'\in[0,t]
\end{align}
hence, in view of (14.68), the borderline contribution (14.72) to (14.71) is bounded by:
\begin{equation}
 C|\ell|\bar{\mu}^{-a}_{m}(t)[1+\log(1+t)]^{p}\sqrt{\leftexp{(i_{1}...i_{l})}{\mathcal{G}}_{0,l+2;a,p}(t)}
\end{equation}
Also, in regard to the last factor on the right in (14.63) we have:
\begin{align}
 \|\underline{L}R_{i_{l+1}}...R_{i_{1}}\psi_{\alpha}\|_{L^{2}(\Sigma_{t}^{\epsilon_{0}})}
\leq \sqrt{\mathcal{E}_{0}[R_{i_{l+1}}...R_{i_{1}}\psi_{\alpha}](t)}\\\notag
\leq \bar{\mu}^{-a}_{m}(t)[1+\log(1+t)]^{p}\sqrt{\leftexp{(i_{1}...i_{l})}{\mathcal{G}}_{0,l+2;a,p}(t)}
\end{align}

    Substituting (14.89), (14.93) and (14.94) in the integral on the right of (14.63), the factors $|\ell|$ cancel and 
we obtain that the borderline contribution to the integral in question is bounded by:
\begin{align}
 C\int_{0}^{t}\{\sup_{\Sigma_{t'}^{\epsilon_{0}}}(\mu^{-1}|L\mu|)+C\bar{\mu}^{-1}_{m}\delta_{0}(1+t')^{-2}[1+\log(1+t')]\}
\cdot\\\notag
\bar{\mu}^{-2a}_{m}(t')[1+\log(1+t')]^{2p}\leftexp{(i_{1}...i_{l})}{\mathcal{G}}_{0,l+2;a,p}(t')dt'
\end{align}
Now, the partial contribution of the second term in the first factor is actually not borderline. We shall show how to 
estimate contribution of this type afterwards, in connection with the estimate for the contribution of the term
\begin{equation}
 C(1+t)^{-1}[1+\log(1+t)]^{1/2}\sqrt{\sum_{j,\alpha}\mathcal{E}'_{1}[R_{j}R_{i_{l}}...R_{i_{1}}\psi_{\alpha}]}
\end{equation}
in the expression for $\leftexp{(i_{1}...i_{l})}{P}_{l}^{(1)}$ of Proposition 10.6, through (14.69), (14.71) and (14.65), 
to the integral on the right of (14.63). For the present we focus attention on the $borderline$ $integral$:
\begin{align}
 C\int_{0}^{t}\sup_{\Sigma_{t'}^{\epsilon_{0}}}(\mu^{-1}|L\mu|)\bar{\mu}_{m}^{-2a}(t')[1+\log(1+t')]^{2p}
\leftexp{(i_{1}...i_{l})}{\mathcal{G}}_{0,l+2;a,p}(t')dt'
\end{align}
Since
\begin{align*}
 |L\mu|=\max\{-(L\mu)_{-},(L\mu)_{+}\}\leq -(L\mu)_{-}+(L\mu)_{+}
\end{align*}
we have:
\begin{equation}
 \sup_{\Sigma_{t'}^{\epsilon_{0}}}(\mu^{-1}|L\mu|)\leq \sup_{\Sigma_{t'}^{\epsilon_{0}}}(-\mu^{-1}(L\mu)_{-})
+\sup_{\Sigma_{t'}^{\epsilon_{0}}}(\mu^{-1}(L\mu)_{+})
\end{equation}
Substituting (14.98) in (14.97), we first consider the first term on the right of (14.98). According to (8.249):
\begin{equation}
 \sup_{\Sigma_{t}^{\epsilon_{0}}}(-\mu^{-1}(L\mu)_{-})=M(t)
\end{equation}
therefore the contribution in question is:
\begin{align}
 C\int_{0}^{t}M(t')\bar{\mu}^{-2a}_{m}(t')[1+\log(1+t')]^{2p}\leftexp{(i_{1}...i_{l})}{\mathcal{G}}_{0,l+2;a,p}(t')dt'
\\\notag
\leq C[1+\log(1+t)]^{2p}\leftexp{(i_{1}...i_{l})}{\mathcal{G}}_{0,l+2;a,p}(t)I_{2a}(t)
\end{align}
where $I_{a}(t)$ is the integral of Lemma 8.11.
According to Lemma 8.11 we have:
\begin{equation}
 I_{2a}(t)\leq C(2a)^{-1}\bar{\mu}^{-2a}_{m}(t)
\end{equation}
where $C$ is a constant independent of $a$, provided that $a\geq 2$ and $\delta_{0}$ is suitably small depending on $a$. 
We conclude that (14.100) is bounded by:
\begin{equation}
 \frac{C}{2a}\bar{\mu}^{-2a}_{m}(t)[1+\log(1+t)]^{2p}\leftexp{(i_{1}...i_{l})}{\mathcal{G}}_{0,l+2;a,p}(t)
\end{equation}

    Consider next the contribution of the second term on the right in (14.98). Proposition 13.1 implies:
\begin{equation}
 \sup_{\Sigma_{t}^{\epsilon_{0}}}(\mu^{-1}(L\mu)_{+})\leq C(1+t)^{-1}[1+\log(1+t)]^{-1}
\end{equation}
therefore the contribution in question is bounded by:
\begin{equation}
 C\int_{0}^{t}\bar{\mu}^{-2a}_{m}(t')(1+t')^{-1}[1+\log(1+t')]^{2p-1}\leftexp{(i_{1}...i_{l})}{\mathcal{G}}_{0,l+2;a,p}
(t')dt'
\end{equation}
By Corollary 2 of Lemma 8.11:
\begin{equation}
 \bar{\mu}^{-2a}_{m}(t')\leq C\bar{\mu}^{-2a}_{m}(t)
\end{equation}
where $C$ is independent of $a$, also provided that $a\geq 2$ and $\delta_{0}$ is suitably small depending on $a$. Since 
$\leftexp{(i_{1}...i_{l})}{\mathcal{G}}_{0,l+2;a,p}(t)$ is a non-decreasing function of $t$, it follows that (14.104) 
is bounded by:
\begin{equation}
 C\bar{\mu}^{-2a}_{m}(t)\leftexp{(i_{1}...i_{l})}{\mathcal{G}}_{0,l+2;a,p}(t)\int_{0}^{t}
(1+t')^{-1}[1+\log(1+t')]^{2p-1}dt'
\end{equation}
Since
\begin{align*}
 \int_{0}^{t}(1+t')^{-1}[1+\log(1+t')]^{2p-1}dt'=\frac{1}{2p}([1+\log(1+t)]^{2p}-1)
\end{align*}
this is bounded by:
\begin{equation}
 \frac{C}{2p}\bar{\mu}^{-2a}_{m}(t)[1+\log(1+t)]^{2p}\leftexp{(i_{1}...i_{l})}{\mathcal{G}}_{0,l+2;a,p}(t)
\end{equation}
We have thus estimated the borderline contributions to the integral in the right of (14.63). We proceed to 
estimate the remaining contributions

      Let us define:
\begin{equation}
 \mathcal{E}_{0,[n]}(t)=\sum_{m=1}^{n}\mathcal{E}_{0,m}(t)
\end{equation}
where $\mathcal{E}_{0,n}$ represents the sum of the energies associated to the vectorfield $K_{0}$ of all the $n$th order variations.
We then define:
\begin{equation}
 \mathcal{G}_{0,[n];a,p}(t)=\sup_{t'\in[0,t]}\{[1+\log(1+t')]^{-2p}\bar{\mu}^{2a}_{m}(t')\mathcal{E}_{0,[n]}(t')\}
\end{equation}
Let us also define:
\begin{equation}
 \mathcal{E}'_{1,[n]}(t)=\sum_{m=1}^{n}\mathcal{E}'_{1,m}(t)
\end{equation}
where $\mathcal{E}'_{1,n}$ represents the sum of the energies associated to the vectorfield $K_{1}$ of all $n$th order variations. 
We then define, for $q\geq p$:
\begin{equation}
 \mathcal{G}'_{1,[n];a,q}(t)=\sup_{t'\in[0,t]}\{[1+\log(1+t')]^{-2q}\bar{\mu}^{2a}_{m}(t')\mathcal{E}'_{1,[n]}(t')\}
\end{equation}

     We now consider the contribution of the term (14.96), which comes from the second term on the right of definition 
for $\leftexp{(i_{1}...i_{l})}{P}^{(1)}_{l}$, through (14.69) and (14.65), to the integral on the right of (14.63). 
From (14.110) and (14.111), we have:
\begin{align}
 \sqrt{\sum_{j,\alpha}\mathcal{E}'_{1}[R_{j}R_{i_{l}}...R_{i_{1}}\psi_{\alpha}](t')}\leq \sqrt{\mathcal{E}'
_{1,[l+2]}(t')}
\leq \bar{\mu}^{-a}_{m}(t')[1+\log(1+t')]^{q}\sqrt{\mathcal{G}'_{1,[l+2];a,q}(t)}\\\notag
:\quad\textrm{for all}\quad t'\in[0,t]
\end{align}
hence, in view of (14.69), the contribution of (14.96) to (14.71) is bounded by:
\begin{equation}
 C_{q}(1+t)^{-1/2}\bar{\mu}_{m}^{-a}(t)\sqrt{\mathcal{G}'_{1,[l+2];a,q}(t)}
\end{equation}
where $C_{q}$ is a constant depending on $q$:
\begin{align*}
 C_{q}=C\sup_{t\in[0,\infty)}\{(1+t)^{-1/2}[1+\log(1+t)]^{q+1/2}\}
\end{align*}
Also, by (14.94) and the definitions (14.108), (14.109) we have:
\begin{equation}
\|\underline{L}R_{i_{l+1}}...R_{i_{1}}\psi_{\alpha}\|_{L^{2}(\Sigma_{t}^{\epsilon_{0}})}
\leq \bar{\mu}^{-a}_{m}(t)[1+\log(1+t)]^{p}\sqrt{\mathcal{G}_{0,[l+2];a,p}(t)}
 \end{equation}
Hence, by (14.90), the contribution of (14.96) through (14.65) to the integral on the right in (14.63) is bounded by:
\begin{equation}
 C_{q}\delta_{0}\int_{0}^{t}(1+t')^{-3/2}[1+\log(1+t')]^{p}\bar{\mu}_{m}^{-2a-1}(t')\sqrt{\mathcal{G}_{0,[l+2];a,p}(t')
\mathcal{G}'_{1,[l+2];a,q}(t')}dt'
\end{equation}
To estimate this we consider the two cases occurring in the proof of Lemma 8.11:

Case 1: $\hat{E}_{s,m}\geq 0$.

In this case we have the lower bound (8.258):
\begin{equation}
 \bar{\mu}_{m}(t)\geq 1-C\delta_{0}
\end{equation}
which implies:
\begin{equation}
 \bar{\mu}_{m}^{-2a-1}(t)\leq C
\end{equation}
provided that $\delta_{0}a$ is suitably small (see(8.260)). It follows that in Case 1, (14.115) is bounded by:
\begin{equation}
 C_{q}\delta_{0}\int_{0}^{t}(1+t')^{-3/2}[1+\log(1+t')]^{p}\sqrt{\mathcal{G}_{0,[l+2];a,p}(t')
\mathcal{G}'_{1,[l+2];a,q}(t')}dt'
\end{equation}
     In Case 2, we set:
\begin{equation}
 \hat{E}_{s,m}=-\delta_{1},\quad \delta_{1}>0
\end{equation}
Following the proof of Lemma 8.11 with $a$ relaced by $2a$, we then set:
\begin{equation}
 t_{1}=e^{\frac{1}{4a\delta_{1}}}-1
\end{equation}
and we consider the two subcases:

     Subcase 2a: $t'\leq t_{1}$   Subcase 2b: $t'>t_{1}$

In Subcase 2a we have the lower bound (8.273) with $a$ replaced by $2a$:
\begin{equation}
 \bar{\mu}_{m}(t')\geq 1-\frac{1}{a}
\end{equation}
Since $(1-\frac{1}{a})^{-2a-1}$ is bounded  for $a\in[2,\infty)$, if $t\leq t_{1}$ (14.115) is bounded by:
\begin{equation}
 C_{q}\delta_{0}\int_{0}^{t}(1+t')^{-3/2}[1+\log(1+t')]^{p}\sqrt{\mathcal{G}_{0,[l+2];a,p}(t')\mathcal{G}'
_{1,[l+2];a,q}(t')}dt'
\end{equation}
and if $t>t_{1}$ we have:
\begin{align}
 \int_{0}^{t_{1}}(1+t')^{-3/2}[1+\log(1+t')]^{p}\bar{\mu}^{-2a-1}_{m}(t')\sqrt{\mathcal{G}_{0,[l+2];a,p}(t')
\mathcal{G}'_{1,[l+2];a,q}(t')}dt'\\\notag
\leq  C\int_{0}^{t_{1}}(1+t')^{-3/2}[1+\log(1+t')]^{p}\sqrt{\mathcal{G}_{0,[l+2];a,p}(t')
\mathcal{G}'_{1,[l+2];a,q}(t')}dt'
\end{align}

    In Subcase 2b we have the lower bound (8.303) with $a$ replaced by $2a$:
\begin{equation}
 \bar{\mu}_{m}(t')\geq (1-\frac{1}{a})(1-\delta_{1}\tau'),\quad \tau'=\log(1+t')
\end{equation}
In view of the fact that $\mathcal{G}_{0,[n];a,p}(t)$, $\mathcal{G}'_{1,[n];a,q}(t)$ are non-decreasing functions of
 $t$ it follows that:
\begin{align}
 \int_{t_{1}}^{t}(1+t')^{-3/2}[1+\log(1+t')]^{p}\bar{\mu}^{-2a-1}_{m}(t')\sqrt{\mathcal{G}_{0,[l+2];a,p}(t')
\mathcal{G}'_{1,[l+2];a,q}(t')}dt'\\\notag
\leq C\sqrt{\mathcal{G}_{0,[l+2];a,p}(t)\mathcal{G}'_{1,[l+2];a,q}(t)}\cdot\int_{t_{1}}^{t}(1+t')^{-3/2}[1+\log(1+t')]
^{p}(1-\delta_{1}\tau')^{-2a-1}dt'
\end{align}
The integral in (14.125) is:
\begin{align}
 \int_{t_{1}}^{t}(1+t')^{-3/2}[1+\log(1+t')]^{p}(1-\delta_{1}\tau')^{-2a-1}dt'\\\notag
\leq (1+t_{1})^{-1/2}[1+\log(1+t_{1})]^{p}\int_{\tau_{1}}^{\tau}(1-\delta_{1}\tau')^{-2a-1}d\tau'\\\notag
\leq (1+t_{1})^{-1/2}[1+\log(1+t_{1})]^{p}\frac{1}{2a\delta_{1}}(1-\delta_{1}\tau)^{-2a}\\\notag
\leq 2\varphi_{1+p}(4a\delta_{1})(1-\delta_{1}\tau)^{-2a}
\end{align}
Here for any real number $r$ we denote by $\varphi_{r}(x)$ the following function on the positive real line:
\begin{equation}
 \varphi_{r}(x)=e^{-\frac{1}{2x}}(1+\frac{1}{x})^{r}
\end{equation}
The function $\varphi_{r}(x)$ decreases to $0$ exponentially as $x\longrightarrow 0$. Since $\delta_{1}\leq
 C\delta_{0}$, we have:
\begin{equation}
 \varphi_{r}(4a\delta_{1})\leq \varphi_{r}(Ca\delta_{0})
\end{equation}
provided that $Ca\delta_{0}$ is suitably small. In view also of the lower bound (8.314) with $a$ replaced by $2a$:
\begin{equation}
 \bar{\mu}^{-2a}_{m}(t)\geq\frac{1}{C}(1-\delta_{1}\tau)^{-2a}
\end{equation}
we conclude that the right hand side of (14.125) is bounded by:
\begin{equation}
 C\varphi_{1+p}(Ca\delta_{0})\bar{\mu}^{-2a}_{m}(t)\sqrt{\mathcal{G}_{0,[l+2];a,p}(t)\mathcal{G}'_{1,[l+2];a,q}(t)}
\end{equation}
Combining with (14.123) we obtain that if $t>t_{1}$ (14.115) is bounded by:
\begin{align}
 C_{q}\delta_{0}\bar{\mu}^{-2a}_{m}(t)\{\varphi_{1+p}(Ca\delta_{0})\sqrt{\mathcal{G}_{0,[l+2];a,p}(t)\mathcal{G}'
_{1,[l+2];a,q}(t)}\\\notag
+\int_{0}^{t}(1+t')^{-3/2}[1+\log(1+t')]^{p}\sqrt{\mathcal{G}_{0,[l+2];a,p}(t')\mathcal{G}'_{1,[l+2];a,q}(t')}dt'\}
\end{align}
Combining finally with the earlier results (14.122), for the subcase $t\leq t_{1}$, and (14.118), for Case 1, 
we conclude that the contribution of 
(14.96) through (14.65) to the integral on the right in (14.63) is bounded by:
\begin{align}
 C_{q}\delta_{0}\bar{\mu}^{-2a}_{m}(t)\{\varphi_{1+p}(Ca\delta_{0})\sqrt{\mathcal{G}_{0,[l+2];a,p}(t)
\mathcal{G}'_{1,[l+2];a,q}(t)}\\\notag
+\int_{0}^{t}(1+t')^{-3/2}[1+\log(1+t')]^{p}\sqrt{\mathcal{G}_{0,[l+2];a,p}(t')\mathcal{G}'_{1,[l+2];a,q}(t')}dt'\}
\end{align}

Now, we consider the contribution of the second term in the expression (14.66) for $\leftexp{(i_{1}...i_{l})}{B}_{l}(t)$ 
to the first term on the right of (14.65), namely the contribution of:
\begin{equation}
 C(1+t)^{-1}[1+\log(1+t)]^{2}\int_{0}^{t}(1+t')^{3}\leftexp{(i_{1}...i_{l})}{Q}_{l}(t')dt'
\end{equation}
to the estimate for $\|\mu R_{i_{l+1}}...R_{i_{1}}\textrm{tr}\chi'\|_{L^{2}(\Sigma_{t}^{\epsilon_{0}})}$.

      Here, we will use the estimate for $\max_{i_{1}...i_{l}}\leftexp{(i_{1}...i_{l})}{Q}_{l}$ of Proposition 10.5. 
The principal part of the right-hand-side is:
\begin{align}
 C_{l}(1+t)^{-4}[1+\log(1+t)][\mathcal{W}^{Q}_{[l+1]}+\delta_{0}(\mathcal{W}^{T}_{[l+1]}+\mathcal{W}_{[l+2]})]
\end{align}
We shall estimate the contribution of this principal part, through (14.133) to the integral on the right in (14.63).

      We have, for each non-negative integer $n$:
\begin{align}
 \mathcal{W}_{n+1}\leq C(1+t)\sqrt{\sum_{i_{1}...i_{n},\alpha}\|
\slashed{d}R_{i_{n}}...R_{i_{1}}\psi_{\alpha}\|^{2}_{L^{2}(\Sigma_{t}^{\epsilon_{0}})}}
\leq C\bar{\mu}_{m}^{-1/2}\sqrt{\mathcal{E}'_{1,n+1}}
\end{align}
It follows that:
\begin{equation}
 \mathcal{W}_{[l+2]}\leq C_{l}\bar{\mu}^{-1/2}_{m}\sqrt{\mathcal{E}'_{1,[l+2]}}+\mathcal{W}_{0}
\end{equation}
By Lemma 5.1, we have:
\begin{equation}
 \mathcal{W}_{0}\leq\sqrt{\epsilon_{0}}\max_{\alpha}\sup_{u\in[0,\epsilon_{0}]}\|\psi_{\alpha}\|_{L^{2}(S_{t,u})}
\leq C\epsilon_{0}\sqrt{\sum_{\alpha}\mathcal{E}_{0}[\psi_{\alpha}]}=C\epsilon_{0}\sqrt{\mathcal{E}_{0,[1]}}
\end{equation}
Similarly, for each non-negative integer $n$:
\begin{align}
 \mathcal{W}^{T}_{n+1}\leq C(1+t)\sqrt{\sum_{i_{1}...i_{n},\alpha}\|\slashed{d}R_{i_{n}}...R_{i_{1}}T\psi_{\alpha}\|
^{2}_{L^{2}(\Sigma_{t}^{\epsilon_{0}})}}
\leq C\bar{\mu}^{-1/2}_{m}\sqrt{\mathcal{E}'_{1,n+2}}
\end{align}
It follows that:
\begin{equation}
 \mathcal{W}^{T}_{[l+1]}\leq C_{l}\bar{\mu}^{-1/2}_{m}\sqrt{\mathcal{E}'_{1,[l+2]}}+\mathcal{W}^{T}_{0}
\end{equation}
by Lemma 5.1:
\begin{align}
 \mathcal{W}^{T}_{0}\leq\sqrt{\epsilon_{0}}\max_{\alpha}\sup_{u\in[0,\epsilon_{0}]}\|T\psi_{\alpha}\|_{L^{2}(S_{t,u})}
\leq C\epsilon_{0}\sqrt{\sum_{\alpha}\mathcal{E}_{0}[T\psi_{\alpha}]}\leq C\epsilon_{0}\sqrt{\mathcal{E}_{0,[2]}}
\end{align}
Similarly, we have:
\begin{equation}
 \mathcal{W}^{Q}_{[l+1]}\leq C_{l}\bar{\mu}_{m}^{-1/2}\sqrt{\mathcal{E}'_{1,[l+2]}}+\mathcal{W}^{Q}_{0}
\end{equation}
and by Lemma 5.1:
\begin{align}
 \mathcal{W}^{Q}_{0}\leq\sqrt{\epsilon_{0}}\max_{\alpha}\sup_{u\in[0,\epsilon_{0}]}\|Q\psi_{\alpha}\|_{L^{2}(S_{t,u})}
\leq C\epsilon_{0}\sqrt{\sum_{\alpha}\mathcal{E}_{0}[Q\psi_{\alpha}]}\leq C\epsilon_{0}\sqrt{\mathcal{E}_{0,[2]}}
\end{align}

The above inequalities imply that (14.134) is bounded by:
\begin{equation}
 C_{l}(1+t)^{-4}[1+\log(1+t)]\{\bar{\mu}^{-1/2}_{m}\sqrt{\mathcal{E}'_{1,[l+2]}}+C\epsilon_{0}
\sqrt{\mathcal{E}_{0,[2]}}\}
\end{equation}
therefore, the contribution of (14.134) to (14.133) is bounded by:
\begin{align}
 C_{l}(1+t)^{-1}[1+\log(1+t)]^{2}\cdot\\\notag
\int_{0}^{t}(1+t')^{-1}[1+\log(1+t')]\{\bar{\mu}_{m}^{-1/2}(t')\sqrt{\mathcal{E}'_{1,[l+2]}(t')}
+C\epsilon_{0}\sqrt{\mathcal{E}_{0,[2]}(t')}\}dt'
\end{align}
We shall first estimate the contribution of the principal term in the last factor of the above, namely 
the term $\bar{\mu}_{m}^{-1/2}(t')
\sqrt{\mathcal{E}'_{1,[l+2]}(t')}$. From the definition (14.111) we have:
\begin{equation}
 \bar{\mu}_{m}^{-1/2}(t')\sqrt{\mathcal{E}'_{1,[l+2]}(t')}\leq 
\bar{\mu}_{m}^{-a-1/2}(t')[1+\log(1+t')]^{q}\sqrt{\mathcal{G}'_{1,[l+2];a,q}(t)}
\end{equation}
Therefore, the contribution in question is bounded by:
\begin{equation}
 C_{l}(1+t)^{-1}[1+\log(1+t)]^{2}J_{a,q}(t)\sqrt{\mathcal{G}'_{1,[l+2];a,q}(t)}
\end{equation}
where
\begin{equation}
 J_{a,q}(t)=\int_{0}^{t}(1+t')^{-1}\bar{\mu}_{m}^{-a-1/2}(t')[1+\log(1+t')]^{q+1}dt'
\end{equation}
To estimate $J_{a,q}(t)$, we consider again the two cases occurring in the proof of Lemma 8.11, as in the estimate
 of (14.115).

     In Case 1 we have the lower bound (14.116) hence also (14.117), which implies:
\begin{align}
 J_{a,q}(t)\leq C\int_{0}^{t}(1+t')^{-1}[1+\log(1+t')]^{q+1}dt'=C\int_{0}^{\tau}(1+\tau')^{q+1}d\tau'\\\notag
\leq\frac{C}{(q+2)}(1+\tau)^{q+2}=\frac{C}{(q+2)}[1+\log(1+t)]^{q+2}
\end{align}
Similarly, in Subcase 2a), we obtain if $t\leq t_{1}$:
\begin{equation}
 J_{a,q}(t)\leq\frac{C}{(q+2)}[1+\log(1+t)]^{q+2}
\end{equation}
and if $t>t_{1}$:
\begin{equation}
 J_{a,q}(t_{1})\leq\frac{C}{(q+2)}[1+\log(1+t_{1})]^{q+2}<\frac{C}{(q+2)}[1+\log(1+t)]^{q+2}
\end{equation}
In Subcase 2b) we have the lower bound (14.124) which implies:
\begin{align}
 J_{a,q}(t)-J_{a,q}(t_{1})\leq C[1+\log(1+t)]^{q+1}\int_{t_{1}}^{t}(1-\delta_{1}\tau')^{-a-1/2}(1+t')^{-1}dt'\\\notag
=C(1+\tau)^{q+1}\int_{\tau_{1}}^{\tau}(1-\delta_{1}\tau')^{-a-1/2}d\tau'
\leq\frac{C}{\delta_{1}}\frac{(1+\tau)^{q+1}}{(a-\frac{1}{2})}(1-\delta_{1}\tau)^{-a+1/2}\\\notag
\leq\frac{C}{a\delta_{1}}(1+\tau)^{q+1}\bar{\mu}^{-a+1/2}_{m}(t)
\end{align}
where in the last step we have used (8.312) which implies:
\begin{align*}
 \bar{\mu}^{-a+1/2}_{m}(t)\geq\frac{1}{C}(1-\delta_{1}\tau)^{-a+1/2}
\end{align*}
Since $\tau\geq\tau_{1}=1/4a\delta_{1}$, (14.151) implies:
\begin{equation}
 J_{a,q}(t)-J_{a,q}(t_{1})\leq C[1+\log(1+t)]^{q+2}\bar{\mu}_{m}^{-a+1/2}(t)
\end{equation}
In view of (14.149), (14.150), (14.152), we conclude that, in general:
\begin{equation}
 J_{a,q}(t)\leq C[1+\log(1+t)]^{q+2}\bar{\mu}^{-a+1/2}_{m}(t)
\end{equation}
hence (14.146) is bounded by:
\begin{equation}
 C_{l}(1+t)^{-1}[1+\log(1+t)]^{q+4}\bar{\mu}^{-a+1/2}_{m}(t)\sqrt{\mathcal{G}^{\prime}_{1,[l+2];a,q}(t)}
\end{equation}
This bounds the contribution in question through (14.133), to the estimate for 
\begin{align*}
\|\mu R_{i_{l+1}}...R_{i_{1}}
\textrm{tr}\chi'\|_{L^{2}(\Sigma_{t}^{\epsilon_{0}})}
\end{align*}
In view of (14.114) and (14.90), the corresponding contribution to the integral on the right of (14.63) is 
then bounded by:
\begin{align}
 C_{l}\delta_{0}\int_{0}^{t}(1+t')^{-2}[1+\log(1+t')]^{p+q+4}\bar{\mu}_{m}^{-2a-1/2}(t')
\sqrt{\mathcal{G}_{0,[l+2];a,p}(t')\mathcal{G}'_{1,[l+2];a,q}(t')}dt'
\end{align}

This is estimated by following an argument similar to that used to estimate (14.115). We obtain in this way a bound by:
\begin{align}
 C_{l}\delta_{0}\bar{\mu}^{-2a+1/2}_{m}(t)[1+\log(1+t)]^{2p}\{\varphi'_{5+q-p}(Ca\delta_{0})
\sqrt{\mathcal{G}_{0,[l+2];a,p}(t)\mathcal{G}'_{1,[l+2];a,q}(t)}\\\notag
+\int_{0}^{t}(1+t')^{-2}[1+\log(1+t')]^{4+q-p}\sqrt{\mathcal{G}_{0,[l+2];a,p}(t')\mathcal{G}'_{1,[l+2];a,q}(t')}dt'\}
\end{align}
Here, for any real number $r$ we denote by $\varphi'_{r}(x)$ the following function on the positive real line:
\begin{equation}
 \varphi'_{r}(x)=e^{-\frac{1}{x}}(1+\frac{1}{x})^{r}
\end{equation}
Note that $\varphi'_{r}(x)$ decreases to 0 exponentially as $x\longrightarrow 0$, and that, since $\delta_{1}
\leq C\delta_{0}$,
\begin{align*}
 \varphi'_{r}(4a\delta_{1})\leq\varphi'_{r}(Ca\delta_{0})
\end{align*}
provided that $Ca\delta_{0}$ is suitably small.

     Finally, we consider the contribution of the last term in the expression (14.66) for 
$\leftexp{(i_{1}...i_{l})}{B}_{l}(t)$ to the first term on the right in (14.65), namely, the contribution of:
\begin{equation}
 C(1+t)^{-1}[1+\log(1+t)]^{2}\|\leftexp{(i_{1}...i_{l})}{x}_{l}(0)\|_{L^{2}(\Sigma_{0}^{\epsilon_{0}})}
\end{equation}
to the estimate for $\|\mu R_{i_{l+1}}...R_{i_{1}}\textrm{tr}\chi'\|_{L^{2}(\Sigma_{t}^{\epsilon_{0}})}$. 
In view of (14.90) and (14.114), the corresponding contribution to the integral on the right in (14.63) is bounded by:
\begin{align}
 C\delta_{0}\|\leftexp{(i_{1}...i_{l})}{x}_{l}(0)\|_{L^{2}(\Sigma_{0}^{\epsilon_{0}})}\int_{0}^{t}
(1+t')^{-2}[1+\log(1+t')]^{p+2}\bar{\mu}_{m}^{-a-1}(t')\sqrt{\mathcal{G}_{0,[l+2];a,p}(t')}dt'
\end{align}
This is again estimated by following an argument similar to that used to estimate (14.115). We obtain a bound by:
\begin{align}
 C_{p}\delta_{0}\|\leftexp{(i_{1}...i_{l})}{x}_{l}(0)\|_{L^{2}(\Sigma_{0}^{\epsilon_{0}})}\bar{\mu}^{-a}_{m}(t)\{
\varphi'_{3+p}(Ca\delta_{0})\sqrt{\mathcal{G}_{0,[l+2];a,p}(t)}\\\notag
+\int_{0}^{t}(1+t')^{-2}[1+\log(1+t')]^{2+p}\sqrt{\mathcal{G}_{0,[l+2];a,p}(t')}dt'\}
\end{align}
In regard to the contribution of the last term on the right in (14.65), we remark that it is bounded by:
\begin{align}
 C_{l}\delta_{0}(1+t)^{-1}[1+\log(1+t)]^{4}\sup_{t^{\prime}\in[0,t]}
\{(1+t^{\prime})^{2}\bar{\mu}^{a}_{m}(t^{\prime})B_{l}(t^{\prime})\}\bar{\mu}^{-a}_{m}(t)
\end{align}
Thus, relative to the term $(1+t)^{2}\leftexp{(i_{1}...i_{l})}{B}_{l}$ there is an extra factor of 
$C_{l}\delta_{0}(1+t)^{-1}[1+\log(1+t)]^{4}$, consequently the contribution in question is absorbed in the estimates 
already made.

\subsection{Estimates for the Contribution of (14.57)}
     We now consider the contribution of (14.57) to the corresponding integral (14.62). Recalling that this is
 associated to the variation (14.54), the contribution in question is:
\begin{align}
 \int_{W^{t}_{\epsilon_{0}}}\Omega|R_{i_{l-m}}...R_{i_{1}}(T)^{m}\slashed{\Delta}\mu||T\psi_{\alpha}||\underline{L}
R_{i_{l-m}}...R_{i_{1}}(T)^{m+1}\psi_{\alpha}|dt'dud\mu_{\slashed{g}}\\\notag
\leq C\int_{0}^{t}\sup_{\Sigma_{t'}^{\epsilon_{0}}}(\mu^{-1}|T\psi_{\alpha}|)\|\mu 
R_{i_{l-m}}...R_{i_{1}}(T)^{m}\slashed{\Delta}\mu\|_{L^{2}
(\Sigma_{t'}^{\epsilon_{0}})}\cdot\\\notag
\|\underline{L}R_{i_{l-m}}...R_{i_{1}}(T)^{m+1}\psi_{\alpha}\|_{L^{2}(\Sigma_{t'}^{\epsilon_{0}})}dt'
\end{align}
$m=0,...,l$.

Now in view of (9.60), (9.70) and (9.213), we have:
\begin{align}
 \|\mu R_{i_{l-m}}...R_{i_{1}}(T)^{m}\slashed{\Delta}\mu\|_{L^{2}(\Sigma_{t}^{\epsilon_{0}})}\\\notag
\leq C(1+t)\{\||\leftexp{(i_{1}...i_{l-m})}{x}^{\prime}_{m,l-m}(t)|\|_{L^{2}([0,\epsilon_{0}]\times S^{2})}
+\||\leftexp{(i_{1}...i_{l-m})}{\check{f}}^{\prime}_{m,l-m}(t)|\|_{L^{2}([0,\epsilon_{0}]\times S^{2})}\}\\\notag
=C\{(1+t)\||\leftexp{(i_{1}...i_{l-m})}{x}^{\prime}_{m,l-m}(t)|\|_{L^{2}([0,\epsilon_{0}]\times S^{2})}+
\leftexp{(i_{1}...i_{l-m})}{P}^{\prime}_{m,l-m}(t)\}
\end{align}
Substituting (9.268) and (9.270) we then obtain:
\begin{align}
 \|\mu R_{i_{l-m}}...R_{i_{1}}(T)^{m}\slashed{\Delta}\mu\|_{L^{2}(\Sigma_{t}^{\epsilon_{0}})}\\\notag
\leq C(1+t)((1+t)^{-1}\leftexp{(i_{1}...i_{l-m})}{P}^{\prime}_{m,l-m}(t)+\leftexp{(i_{1}...i_{l-m})}
{B}^{\prime}_{m,l-m}(t))\\\notag
+C_{l}\delta_{0}(1+t)^{-1}[1+\log(1+t)]^{2}\{\int_{0}^{t}(1+t')[1+\log(1+t')]B_{l}(t')dt'\\\notag
+\sum_{k=0}^{m}\int_{0}^{t}[1+\log(1+t')]B'_{k,l-k}(t')dt'\}
\end{align}
for $m=0,...,l$.

Here, the quantities $\leftexp{(i_{1}...i_{l-m})}{B}^{\prime}_{m,l-m}(t)$ are defined by (9.259):
\begin{align}
 \leftexp{(i_{1}...i_{l-m})}{B}^{\prime}_{m,l-m}(t)\\\notag
=C(1+t)^{-1}\{\leftexp{(i_{1}...i_{l-m})}{\bar{P}}^{\prime(0)}_{m,l-m,a}(t)+(1+t)^{-1/2}\leftexp{(i_{1}...i_{l-m})}
{\bar{P}}^{\prime(1)}_{m,l-m,a}(t)\}
\bar{\mu}^{-a}_{m}(t)\\\notag
+C(1+t)^{-2}[1+\log(1+t)]^{2}\int_{0}^{t}(1+t')^{2}\leftexp{(i_{1}...i_{l-m})}{Q}^{\prime}_{m,l-m}(t')dt'\\\notag
+C(1+t)^{-2}[1+\log(1+t)]^{2}\|\leftexp{(i_{1}...i_{l-m})}{x}^{\prime}_{m,l-m}(0)\|_{L^{2}(\Sigma_{0}^{\epsilon_{0}})}
\end{align}

Also (see (9.261))
\begin{equation}
 B'_{m,l-m}(t)=\max_{i_{1}...i_{l-m}}\leftexp{(i_{1}...i_{l-m})}{B}^{\prime}_{m,l-m}(t)
\end{equation}
The quantities $\leftexp{(i_{1}...i_{l-m})}{\bar{P}}^{\prime(0)}_{m,l-m,a}(t)$, $\leftexp{(i_{1}...i_{l-m})}
{\bar{P}}^{\prime(1)}_{m,l-m,a}(t)$ are defined by 
(9.215) and (9.216):
\begin{align}
 \leftexp{(i_{1}...i_{l-m})}{\bar{P}}^{\prime(0)}_{m,l-m,a}(t)=\sup_{t'\in[0,t]}\{\bar{\mu}^{a}_{m}(t')
\leftexp{(i_{1}...i_{l-m})}{P}^{\prime(0)}_{m,l-m}(t')\}\\
\leftexp{(i_{1}...i_{l-m})}{\bar{P}}^{\prime(1)}_{m,l-m,a}(t)=\sup_{t'\in[0,t]}\{(1+t')^{1/2}\bar{\mu}^{a}_{m}(t')
\leftexp{(i_{1}...i_{l-m})}{P}^{\prime(1)}_{m,l-m}(t')\}
\end{align}
The quantities $\leftexp{(i_{1}...i_{l-m})}{P}^{\prime(0)}_{m,l-m}(t)$, 
$\leftexp{(i_{1}...i_{l-m})}{P}^{\prime(1)}_{m,l-m}(t)$ are defined in the statement of Proposition 11.6, and we have:
\begin{equation}
 \leftexp{(i_{1}...i_{l-m})}{P}^{\prime}_{m,l-m}(t)\leq \leftexp{(i_{1}...i_{l-m})}{P}^{\prime(0)}_{m,l-m}(t)+
\leftexp{(i_{1}...i_{l-m})}{P}^{\prime(1)}_{m,l-m}(t)
\end{equation}
Note that the quantities $\leftexp{(i_{1}...i_{l-m})}{\bar{P}}^{\prime(0)}_{m,l-m,a}(t)$, 
$\leftexp{(i_{1}...i_{l-m})}{\bar{P}}^{\prime(1)}_{m,l-m,a}(t)$ are non-decreasing.

    The leading contribution to the right hand side of (14.164) comes from the first term on the right in (14.164). 
That is, from $(1+t)^{-1}\leftexp{(i_{1}...i_{l-m})}{P}^{\prime}_{m,l-m}(t)$ and the first term in the expression 
(14.165) for $\leftexp{(i_{1}...i_{l-m})}{B}^{\prime}_{m,l-m}(t)$. So what we consider here is: 
\begin{align}
 C\{\leftexp{(i_{1}...i_{l-m})}{\bar{P}}^{\prime(0)}_{m,l-m,a}(t)+(1+t)^{-1/2}\leftexp{(i_{1}...i_{l-m})}
{\bar{P}}^{\prime(1)}_{m,l-m,a}(t)\}\bar{\mu}_{m}^{-a}(t)
\end{align}
The actual $borderline$ $contribution$ is the contribution from:
\begin{align}
 \leftexp{(i_{1}...i_{l-m})}{P}^{\prime(0)}_{m,l-m}=|\ell|\sqrt{\sum_{\alpha}\mathcal{E}_{0}[R_{i_{l-m}}
...R_{i_{1}}(T)^{m+1}\psi_{\alpha}]}
\end{align}
We shall appeal to (14.89) in estimating the borderline contribution (14.171), through (14.167), (14.170) and (14.165), 
to the integral on the right in (14.162). On the other hand, in estimating all other contributions, through (14.165), 
to the same integral, we simply use (14.90).

     Let us define, for non-negative real numbers $a$ and $p$, the quantities:
\begin{align}
 \leftexp{(i_{1}...i_{l-m})}{\mathcal{G}}_{0,m,l+2;a,p}(t)
=\sup_{t'\in[0,t]}\{[1+\log(1+t')]^{-2p}\bar{\mu}^{2a}_{m}(t')\sum_{\alpha}\mathcal{E}_{0}[R_{i_{l-m}}...R_{i_{1}}
(T)^{m+1}\psi_{\alpha}](t')\}
\end{align}
These quantities are non-decreasing functions of $t$ and we have:
\begin{align}
 \sqrt{\sum_{\alpha}\mathcal{E}_{0}[R_{i_{l-m}}...R_{i_{1}}(T)^{m+1}\psi_{\alpha}](t')}\\\notag
\leq\bar{\mu}^{-a}_{m}(t')[1+\log(1+t')]^{p}\sqrt{\leftexp{(i_{1}...i_{l-m})}{\mathcal{G}}_{0,m,l+2;a,p}(t)}
\end{align}
for all $t'\in[0,t]$.
hence, in view of (14.167), the borderline contribution (14.171) to (14.170) is bounded by:
\begin{align}
 |\ell|\bar{\mu}^{-a}_{m}(t)[1+\log(1+t)]^{p}\sqrt{\leftexp{(i_{1}...i_{l-m})}{\mathcal{G}}_{0,m,l+2;a,p}(t)}
\end{align}
Also, in regard to the last factor in the integral on the right in (14.162) we have:
\begin{align}
 \|\underline{L}R_{i_{l-m}}...R_{i_{1}}(T)^{m+1}\psi_{\alpha}\|_{L^{2}(\Sigma_{t}^{\epsilon_{0}})}\leq 
\sqrt{\mathcal{E}_{0}[R_{i_{l-m}}...R_{i_{1}}(T)^{m+1}\psi_{\alpha}](t)}\\\notag
\leq\bar{\mu}^{-a}_{m}(t)[1+\log(1+t)]^{p}\sqrt{\leftexp{(i_{1}...i_{l-m})}{\mathcal{G}}_{0,m,l+2;a,p}(t)}
\end{align}
Substituting (14.89), (14.174) and (14.175) in the integral on the right in (14.162), the factors $|\ell|$ cancel 
and we obtain 
that the borderline contribution to the integral in question is bounded by:
\begin{align}
 C\int_{0}^{t}\{\sup_{\Sigma_{t'}^{\epsilon_{0}}}(\mu^{-1}|L\mu|)+C\bar{\mu}^{-1}_{m}\delta_{0}(1+t')^{-2}
[1+\log(1+t')]\}\cdot\\\notag
\bar{\mu}^{-2a}_{m}(t')[1+\log(1+t')]^{2p}\leftexp{(i_{1}...i_{l-m})}{\mathcal{G}}_{0,m,l+2;a,p}(t')dt'
\end{align}
Now, the partial contribution of the second term in the first factor is actually not borderline. We shall show how 
to estimate this afterwards, in connection with the estimate for the contribution of the terms
\begin{align}
 C(1+t)^{-1}[1+\log(1+t)]^{1/2}\sqrt{\sum_{\alpha}\mathcal{E}'_{1}[R_{i_{l-m}}...R_{i_{1}}(T)^{m+1}
\psi_{\alpha}]}\\\notag
+C(1+t)^{-2}[1+\log(1+t)]^{3/2}\sqrt{\sum_{\alpha,j}\mathcal{E}'_{1}[R_{j}R_{i_{l-m}}...R_{i_{1}}(T)^{m}\psi_{\alpha}]}
\end{align}
in the expression for $\leftexp{(i_{1}...i_{l-m})}{P}^{\prime(1)}_{m,l-m}$ of Proposition 11.6. For the present we 
focus on the $borderline$ integral:
\begin{align}
 C\int_{0}^{t}\sup_{\Sigma_{t'}^{\epsilon_{0}}}(\mu^{-1}|L\mu|)\bar{\mu}^{-2a}_{m}(t')[1+\log(1+t')]^{2p}
\leftexp{(i_{1}...i_{l-m})}{\mathcal{G}}_{0,m,l+2;a,p}(t')dt'
\end{align}
This is estimated in exactly the same way as the integral (14.97). We obtain that it is bounded by (see (14.102) 
and (14.107)):
\begin{align}
 C(\frac{1}{2a}+\frac{1}{2p})\bar{\mu}^{-2a}_{m}(t)[1+\log(1+t)]^{2p}\leftexp{(i_{1}...i_{l-m})}
{\mathcal{G}}_{0,m,l+2;a,p}(t)
\end{align}

     We proceed to estimate the remaining contributions to the right of (14.162). We first consider the contribution 
of (14.179). From the definitions (14.110) and (14.111) we have:
\begin{align}
 \sqrt{\sum_{\alpha}\mathcal{E}'_{1}[R_{i_{l-m}}...R_{i_{1}}(T)^{m+1}\psi_{\alpha}](t')}
\leq\sqrt{\mathcal{E}'_{1,[l+2]}(t')}
\leq\bar{\mu}^{-a}_{m}(t')[1+\log(1+t')]^{q}\sqrt{\mathcal{G}'_{1,[l+2];a,q}(t)}\\\notag
\sqrt{\sum_{\alpha}\mathcal{E}'_{1}[R_{j}R_{i_{l-m}}...R_{i_{1}}(T)^{m}\psi_{\alpha}](t')}\leq
\sqrt{\mathcal{E}'_{1,[l+2]}(t')}
\leq\bar{\mu}^{-a}_{m}(t')[1+\log(1+t')]^{q}\sqrt{\mathcal{G}'_{1,[l+2];a,q}(t)}
\end{align}
for all $t'\in[0,t]$.
hence, in view of the definition (14.168), the contribution of (14.177) to (14.170) is bounded by:
\begin{align}
 C_{q}(1+t)^{-1/2}\bar{\mu}^{-a}_{m}(t)\sqrt{\mathcal{G}'_{1,[l+2];a,q}(t)}
\end{align}
Also, we have:
\begin{align}
 \|\underline{L}R_{i_{l-m}}...R_{i_{1}}(T)^{m+1}\psi_{\alpha}\|_{L^{2}(\Sigma_{t}^{\epsilon_{0}})}
\leq\bar{\mu}^{-a}_{m}(t)[1+\log(1+t)]^{p}\sqrt{\mathcal{G}_{0,[l+2];a,p}(t)}
\end{align}
Using (14.90) we conclude that the contribution of (14.177) through (14.164) to the integral on the right in (14.162) 
is bounded by:
\begin{align}
 C_{q}\delta_{0}\int_{0}^{t}(1+t')^{-3/2}[1+\log(1+t')]^{p}\bar{\mu}^{-2a-1}_{m}(t')
\sqrt{\mathcal{G}_{0,[l+2];a,p}(t')\mathcal{G}^{\prime}_{1,[l+2];a,q}(t')}dt'
\end{align}
This is identical in the form to (14.115) and it therefore bounded by (14.122).

    Next, we consider the contribution of the second term in the expression (14.165) 
to the first term on the right of (14.164), namely the contribution of 
\begin{align}
 C(1+t)^{-1}[1+\log(1+t)]^{2}\int_{0}^{t}(1+t')^{2}\leftexp{(i_{1}...i_{l-m})}{Q}^{\prime}_{m,l-m}(t')dt'
\end{align}
to the estimate for $\|\mu R_{i_{l-m}}...R_{i_{1}}(T)^{m}\slashed{\Delta}\mu\|_{L^{2}(\Sigma_{t}^{\epsilon_{0}})}$.
Here, we shall use Proposition 11.5, namely the estimate for 
$\max_{i_{1}...i_{l-m}}\leftexp{(i_{1}...i_{l-m})}{Q}^{\prime}_{m,l-m}$. The principal part is:
\begin{align}
 C_{l}(1+t)^{-3}[1+\log(1+t)](\mathcal{W}_{\{l+2\}}+\mathcal{W}^{Q}_{\{l+1\}}+\mathcal{W}^{QQ}_{\{l\}})
\end{align}
We shall estimate the contribution of this principal part, through (14.184) and (14.164), to the integral on the right 
in (14.162).

     For each pair of non-negative integers $m,n$ we have:
\begin{align}
 \mathcal{W}_{m,n+1}\leq C(1+t)\sqrt{\sum_{i_{1}...i_{n},\alpha}\|\slashed{d}R_{i_{n}}...R_{i_{1}}(T)^{m}
\psi_{\alpha}\|^{2}_{L^{2}(\Sigma_{t}^{\epsilon_{0}})}}\leq C\bar{\mu}^{-1/2}_{m}\sqrt{\mathcal{E}'_{1,m+n+1}}
\end{align}
and, for each non-negative integer $m$ we have:
\begin{align}
 \mathcal{W}_{m+1,0}\leq\sqrt{\sum_{\alpha}\|(T)^{m+1}\psi_{\alpha}\|^{2}_{L^{2}(\Sigma_{t}^{\epsilon_{0}})}}
\leq C\sqrt{\mathcal{E}_{0,m+1}}
\end{align}
It follows that:
\begin{align}
 \mathcal{W}_{\{l+2\}}=\sum_{n+m\leq l+2}\mathcal{W}_{m,n}=
\sum_{n+m\leq l+1}\mathcal{W}_{m,n+1}+\sum_{m\leq l+1}\mathcal{W}_{m+1,0}+\mathcal{W}_{0}\\\notag
\leq C_{l}(\bar{\mu}^{-1/2}_{m}\sqrt{\mathcal{E}'_{1,[l+2]}}+\sqrt{\mathcal{E}_{0,[l+2]}})+\mathcal{W}_{0}
\end{align}
Similarly, we have:
\begin{align}
 \mathcal{W}^{Q}_{\{l+1\}}\leq C_{l}(\bar{\mu}^{-1/2}_{m}\sqrt{\mathcal{E}'_{1,[l+2]}}+\sqrt{\mathcal{E}_{0,[l+2]}})
+\mathcal{W}^{Q}_{0}\\
\mathcal{W}^{QQ}_{\{l\}}\leq C_{l}(\bar{\mu}^{-1/2}_{m}\sqrt{\mathcal{E}'_{1,[l+2]}}+\sqrt{\mathcal{E}_{0,[l+2]}})
+\mathcal{W}^{QQ}_{0}
\end{align}

We then conclude that (14.185) is bounded by:
\begin{align}
 C_{l}(1+t)^{-3}[1+\log(1+t)]\{\bar{\mu}^{-1/2}_{m}\sqrt{\mathcal{E}'_{1,[l+2]}}+\sqrt{\mathcal{E}_{0,[l+2]}}\}
\end{align}
(as we may assume that $l\geq 1$). Consequently, the contribution of (14.185) to (14.184) is bounded by:
\begin{align}
 C_{l}(1+t)^{-1}[1+\log(1+t)]^{2}\int_{0}^{t}(1+t')^{-1}[1+\log(1+t')]
\{\bar{\mu}^{-1/2}_{m}(t')\sqrt{\mathcal{E}'_{1,[l+2]}(t')}+\sqrt{\mathcal{E}_{0,[l+2]}(t')}\}dt'
\end{align}
From the definition (14.109):
\begin{align}
 \sqrt{\mathcal{E}_{0,[l+2]}(t')}\leq\bar{\mu}^{-a}_{m}(t')[1+\log(1+t')]^{p}\sqrt{\mathcal{G}_{0,[l+2];a,p}(t)}
\end{align}
This together with (14.145) implies that (14.192) is bounded by, recalling the definition (14.147):
\begin{align}
 C_{l}(1+t)^{-1}[1+\log(1+t)]^{2}\{J_{a,q}(t)\sqrt{\mathcal{G}'_{1,[l+2];a,q}(t)}+J_{a-1/2,p}(t)
\sqrt{\mathcal{G}_{0,[l+2];a,p}(t)}\}
\end{align}
Substituting (14.153) and the same with $(a,q)$ replaced by $(a-\frac{1}{2},p)$ we conclude that the 
contribution of (14.185), through (14.184) to the estimate for 
\begin{align*}
 \|\mu R_{i_{l-m}}...R_{i_{1}}(T)^{m}\slashed{\Delta}\mu\|_{L^{2}(\Sigma_{t}^{\epsilon_{0}})}
\end{align*}
is bounded by:
\begin{align}
 C_{l}(1+t)^{-1}[1+\log(1+t)]^{q+4}\bar{\mu}_{m}^{-a+1/2}(t)\sqrt{\mathcal{G}'_{1,[l+2];a,q}(t)}\\\notag
+C_{l}(1+t)^{-1}[1+\log(1+t)]^{p+4}\bar{\mu}_{m}^{-a+1}(t)\sqrt{\mathcal{G}_{0,[l+2];a,p}(t)}
\end{align}
In view of (14.182) and (14.90), the corresponding contribution to the integral on the right in (14.162) is then 
bounded by:
\begin{align}
 C_{l}\delta_{0}\int_{0}^{t}(1+t')^{-2}[1+\log(1+t')]^{p+q+4}\bar{\mu}^{-2a-1/2}_{m}(t')
\sqrt{\mathcal{G}_{0,[l+2];a,p}(t')\mathcal{G}'_{1,[l+2];a,q}(t')}dt'\\\notag
+C_{l}\delta_{0}\int_{0}^{t}(1+t')^{-2}[1+\log(1+t')]^{2p+4}\bar{\mu}^{-2a}_{m}(t')\mathcal{G}_{0,[l+2];a,p}(t')dt'
\end{align}
The first terms in (14.196) coincides with (14.155) and is estimated by (14.156) while the second is estimated in 
a similar manner by:
\begin{align}
 C_{l}\delta_{0}\bar{\mu}^{-2a+1}_{m}[1+\log(1+t)]^{2p}\{\varphi'_{5}(Ca\delta_{0})\mathcal{G}_{0,[l+2];a,p}(t)\\\notag
+\int_{0}^{t}(1+t')^{-2}[1+\log(1+t')]^{4}\mathcal{G}_{0,[l+2];a,p}(t')dt'\}
\end{align}

     Let us finally consider the contribution of the last term in the expression (14.167) to the first term on 
the right in (14.164), namely the contribution of:
\begin{align}
 C(1+t)^{-1}[1+\log(1+t)]^{2}\|\leftexp{(i_{1}...i_{l-m})}{x}^{\prime}_{m,l-m}(0)\|_{L^{2}(\Sigma_{0}^{\epsilon_{0}})}
\end{align}
In view of (14.182) and the estimate (14.90), the contribution in question is bounded by:
\begin{align}
 C\delta_{0}\|\leftexp{(i_{1}...i_{l-m})}{x}^{\prime}_{m,l-m}(0)\|_{L^{2}(\Sigma_{0}^{\epsilon_{0}})}
\int_{0}^{t}(1+t')^{-2}[1+\log(1+t')]^{p+2}\bar{\mu}^{-a-1}_{m}(t')\sqrt{\mathcal{G}_{0,[l+2];a,p}(t')}dt'
\end{align}
This is similar to (14.159), and is bounded by:
\begin{align}
 C_{p}\delta_{0}\|\leftexp{(i_{1}...i_{l-m})}{x}^{\prime}_{m,l-m}(0)\|_{L^{2}(\Sigma_{0}^{\epsilon_{0}})}
\bar{\mu}_{m}^{-a}(t)\{\varphi'_{3+p}(Ca\delta_{0})\sqrt{\mathcal{G}_{0,[l+2];a,p}(t)}\\\notag
+\int_{0}^{t}(1+t')^{-2}[1+\log(1+t')]^{2+p}\sqrt{\mathcal{G}_{0,[l+2];a,p}(t')}dt'\}
\end{align}
In regard to the contribution of the last term on the right in (14.164), we remark that it is bounded by:
\begin{align}
 C_{l}\delta_{0}(1+t)^{-1}[1+\log(1+t)]^{4}\{\sup_{t'\in[0,t]}\{(1+t')^{2}\bar{\mu}^{a}_{m}(t')B_{l}(t')\}+
\sum_{k=0}^{m}\sup_{t'\in[0,t]}\{(1+t')\bar{\mu}^{a}_{m}(t')B'_{k,l-k}\}\}\bar{\mu}^{-a}_{m}(t)
\end{align}
Thus, relative to the terms $(1+t)^{2}\leftexp{(i_{1}...i_{l})}{B}_{l}$ and $(1+t)\leftexp{(i_{1}...i_{l-m})}
{B}^{\prime}_{m,l-m}$ there is an extra factor of $C_{l}\delta_{0}(1+t)^{-1}[1+\log(1+t)]^{4}$, consequently the 
contribution in question is absorbed in the estimates already made.

\section{The Borderline Estimates Associated to $K_{1}$}
\subsection{Estimates for the Contribution of (14.56)}
     We now consider the contribution of (14.56) to the corresponding integral (14.59). Recalling that this is 
associated to the variation (14.52), the contribution in question is:
\begin{align}
 -\int_{W^{t}_{u}}(\omega/\nu)\Omega(R_{i_{l+1}}...R_{i_{1}}\textrm{tr}\chi')(T\psi_{\alpha})((L+\nu)
R_{i_{l+1}}...R_{i_{1}}\psi_{\alpha})dt'du'd\mu_{\slashed{g}}
\end{align}
Here, using the flux $\mathcal{F}'_{1}[R_{i_{l+1}}...R_{i_{1}}\psi_{\alpha}]$ does not lead to 
an appropriate estimate. 
Instead, we proceed as follows. First, since
\begin{align*}
 \tilde{\slashed{g}}=\Omega\slashed{g},\quad d\mu_{\tilde{\slashed{g}}}=\Omega d\mu_{\slashed{g}}
\end{align*}
the integral (14.202) is:
\begin{align}
 -\int_{W^{t}_{u}}(\omega/\nu)(T\psi_{\alpha})(R_{i_{l+1}}...R_{i_{1}}\textrm{tr}\chi')((L+\nu)
R_{i_{l+1}}...R_{i_{1}}\psi_{\alpha})dt'du'd\mu_{\tilde{\slashed{g}}}
\end{align}
Let us consider, for an arbitrary function $f$, the integral:
\begin{align*}
 \int_{W^{t}_{u}}(Lf+2\nu f)dt'du'd\mu_{\tilde{\slashed{g}}}
\end{align*}
Let:
\begin{align}
 F(t,u)=\int_{S_{t,u}}fd\mu_{\tilde{\slashed{g}}}
\end{align}
Then we have:
\begin{align}
 \frac{\partial F}{\partial t}=\int_{S_{t,u}}(Lf+2\nu f)d\mu_{\tilde{\slashed{g}}}
\end{align}
We have used this formula in Chapter 5. Therefore:
\begin{align}
 \int_{W^{t}_{u}}(Lf+2\nu f)dt'du'd\mu_{\tilde{\slashed{g}}}=\int_{0}^{u}\int_{0}^{t}\frac{\partial F}
{\partial t'}(t',u')dt'du'\\\notag
=\int_{0}^{u}\{F(t,u')-F(0,u')\}du'=\int_{\Sigma_{t}^{u}}fdu'd\mu_{\tilde{\slashed{g}}}
-\int_{\Sigma_{0}^{u}}fdu'd\mu_{\tilde{\slashed{g}}}
\end{align}
Going back to (14.203), we write the integrand in the form:
\begin{align*}
 -(\omega/\nu)(T\psi_{\alpha})(R_{i_{l+1}}...R_{i_{1}}\textrm{tr}\chi')((L+\nu)R_{i_{l+1}}...R_{i_{1}}\psi_{\alpha})\\
=-(L+2\nu)\{(\omega/\nu)(T\psi_{\alpha})(R_{i_{l+1}}...R_{i_{1}}\textrm{tr}\chi')(R_{i_{l+1}}...R_{i_{1}}
\psi_{\alpha})\}\\
+(R_{i_{l+1}}...R_{i_{1}}\psi_{\alpha})(L+\nu)\{(\omega/\nu)(T\psi_{\alpha})(R_{i_{l+1}}...R_{i_{1}}\textrm{tr}\chi')\}
\end{align*}
By (14.206) with the function $(\omega/\nu)(T\psi_{\alpha})(R_{i_{l+1}}...R_{i_{1}}\textrm{tr}\chi')
(R_{i_{l+1}}...R_{i_{1}}\psi_{\alpha})$ in the role 
of the function $f$, we then conclude that (14.203) equals:
\begin{align}
 -\int_{\Sigma_{t}^{u}}(\omega/\nu)(T\psi_{\alpha})(R_{i_{l+1}}...R_{i_{1}}\textrm{tr}\chi')(R_{i_{l+1}}...R_{i_{1}}
\psi_{\alpha})du'd\mu_{\tilde{\slashed{g}}}\\\notag
+\int_{\Sigma_{0}^{u}}(\omega/\nu)(T\psi_{\alpha})(R_{i_{l+1}}...R_{i_{1}}\textrm{tr}\chi')(R_{i_{l+1}}...R_{i_{1}}
\psi_{\alpha})du'd\mu_{\tilde{\slashed{g}}}\\\notag
+\int_{W^{t}_{u}}(R_{i_{l+1}}...R_{i_{1}}\psi_{\alpha})(L+\nu)\{(\omega/\nu)(T\psi_{\alpha})(R_{i_{l+1}}...R_{i_{1}}
\textrm{tr}\chi')\}dt'du'd\mu_{\tilde{\slashed{g}}}
\end{align}
We first consider the hypersurface integral:
\begin{align}
 -\int_{\Sigma_{t}^{u}}(\omega/\nu)(T\psi_{\alpha})(R_{i_{l+1}}...R_{i_{1}}\textrm{tr}\chi')(R_{i_{l+1}}...R_{i_{1}}
\psi_{\alpha})du'd\mu_{\tilde{\slashed{g}}}
\end{align}
(The other hypersurface integral can be expressed in terms of initial data.) Let $f$, $g$ be arbitrary functions 
defined on $S_{t,u}$ and $X$ an arbitrary vectorfield tangent to $S_{t,u}$. We have:
\begin{align*}
 \int_{S_{t,u}}f(Xg)d\mu_{\tilde{\slashed{g}}}=\int_{S_{t,u}}X(fg)d\mu_{\tilde{\slashed{g}}}-
\int_{S_{t,u}}g(Xf)d\mu_{\tilde{\slashed{g}}}
\end{align*}
and
\begin{align*}
 \int_{S_{t,u}}X(fg)d\mu_{\tilde{\slashed{g}}}=\int_{S_{t,u}}\tilde{\slashed{\textrm{div}}}(fgX)
d\mu_{\tilde{\slashed{g}}}-
\int_{S_{t,u}}(\tilde{\slashed{\textrm{div}}}X)fgd\mu_{\tilde{\slashed{g}}}
\end{align*}
Since
\begin{align*}
 \int_{S_{t,u}}\tilde{\slashed{\textrm{div}}}(fgX)d\mu_{\tilde{\slashed{g}}}=0,\quad
 \tilde{\slashed{\textrm{div}}}X=\frac{1}{2}\textrm{tr}\leftexp{(X)}{\tilde{\slashed{\pi}}}
\end{align*}
we obtain:
\begin{align}
 \int_{S_{t,u}}f(Xg)d\mu_{\tilde{\slashed{g}}}=-\int_{S_{t,u}}\{g(Xf)+\frac{1}{2}\textrm{tr}
\leftexp{(X)}{\tilde{\slashed{\pi}}}fg\}d\mu_{\tilde{\slashed{g}}}
\end{align}
Applying (14.209), taking:
\begin{align*}
 X=R_{i_{l+1}},\quad g=R_{i_{l}}...R_{i_{1}}\textrm{tr}\chi',\quad f=(\omega/\nu)(T\psi_{\alpha})
(R_{i_{l+1}}...R_{i_{1}}\psi_{\alpha})
\end{align*}
we conclude that (14.208) equals the sum:
\begin{align}
 H_{0}+H_{1}+H_{2}
\end{align}
with
\begin{align}
 H_{0}=\int_{\Sigma_{t}^{u}}(\omega/\nu)(T\psi_{\alpha})(R_{i_{l}}...R_{i_{1}}\textrm{tr}\chi')(R_{i_{l+1}}
R_{i_{l+1}}...R_{i_{1}}\psi_{\alpha})du'd\mu_{\tilde{\slashed{g}}}\\
H_{1}=\int_{\Sigma_{t}^{u}}(\omega/\nu)(R_{i_{l+1}}T\psi_{\alpha})(R_{i_{l}}...R_{i_{1}}\textrm{tr}\chi')
(R_{i_{l+1}}...R_{i_{1}}\psi_{\alpha})du'd\mu_{\tilde{\slashed{g}}}\\
H_{2}=\int_{\Sigma_{t}^{u}}(T\psi_{\alpha})(R_{i_{l}}...R_{i_{1}}\textrm{tr}\chi')(R_{i_{l+1}}...R_{i_{1}}
\psi_{\alpha})
\{R_{i_{l+1}}(\omega/\nu)+\frac{1}{2}\textrm{tr}\leftexp{(R_{i_{l+1}})}{\tilde{\slashed{\pi}}}(\omega/\nu)\}
du'd\mu_{\tilde{\slashed{g}}}
\end{align}
Now, by Proposition 12.9 and Corollary 10.1.d with $l=1$, recalling that $\omega$ is constant on each $S_{t,u}$,
\begin{align}
 (\nu/\omega)|R_{i_{l+1}}(\omega/\nu)|+\frac{1}{2}|\textrm{tr}\leftexp{(R_{i_{l+1}})}{\tilde{\slashed{\pi}}}|
\leq C\delta_{0}(1+t)^{-1}[1+\log(1+t)]
\end{align}
It follows that $H_{2}$ has an extra decay factor of $\delta_{0}(1+t)^{-1}[1+\log(1+t)]$ relative to $H_{1}$.
So we confine ourselves to $H_{0}$ and $H_{1}$.

     We first consider $H_{0}$. We have:
\begin{align}
 |H_{0}|\leq C\int_{\Sigma_{t}^{\epsilon_{0}}}(1+t)^{3}|T\psi_{\alpha}||R_{i_{l}}...R_{i_{1}}\textrm{tr}\chi'|
|\slashed{d}R_{i_{l+1}}...R_{i_{1}}\psi_{\alpha}|dud\mu_{\slashed{g}}
\end{align}
Although we have a $L^{2}$ bound for $R_{i_{l}}...R_{i_{1}}\textrm{tr}\chi'$ from Proposition 12.11, here we 
need a more precise estimate. Since
\begin{align}
 L\textrm{tr}\chi'=\textrm{tr}(\slashed{\mathcal{L}}_{L}\chi')-2\textrm{tr}(\chi\cdot\chi')\\\notag
=\textrm{tr}(\slashed{\mathcal{L}}_{L}\chi')-\frac{2}{1-u+t}\textrm{tr}\chi'-2|\chi'|^{2}
\end{align}
taking the trace of (12.47), we obtain:
\begin{align}
 L\textrm{tr}\chi'+\frac{2}{1-u+t}\textrm{tr}\chi'=\rho_{0}
\end{align}
where
\begin{align}
 \rho_{0}=e\textrm{tr}\chi-|\chi'|^{2}+\textrm{tr}b
\end{align}
Applying $R_{i_{l}}...R_{i_{1}}$ to this equation and using the Lemma 11.22 then yields the propagation equation:
\begin{align}
 LR_{i_{l}}...R_{i_{1}}\textrm{tr}\chi'+\frac{2}{1-u+t}R_{i_{l}}...R_{i_{1}}\textrm{tr}\chi'=
\leftexp{(i_{1}...i_{l})}{\rho}_{l}
\end{align}
where
\begin{align}
 \leftexp{(i_{1}...i_{l})}{\rho}_{l}=R_{i_{l}}...R_{i_{1}}\rho_{0}+\sum_{k=0}^{l-1}R_{i_{l}}...R_{i_{l-k+1}}
\leftexp{(R_{i_{l-k}})}{Z}R_{i_{l-k-1}}...R_{i_{1}}\textrm{tr}\chi'
\end{align}
Obviously, the principal part of $\rho_{0}$ is contained in $-\textrm{tr}\alpha'$, and is 
$\frac{1}{2}\frac{dH}{dh}\slashed{\Delta}h$. We can write this in a more precise way:
\begin{align*}
 \frac{1}{2}\frac{dH}{dh}\slashed{\Delta}h=\frac{1}{2}\frac{dH}{dh}(\slashed{\Delta}\psi_{0}-\sum_{i}\psi_{i}
\slashed{\Delta}\psi_{i}-\sum_{i}\slashed{d}\psi_{i}\cdot\slashed{d}\psi_{i})
\end{align*}
So we can express:
\begin{align}
 -R_{i_{l}}...R_{i_{1}}(\textrm{tr}\alpha'^{[P]})=
\frac{1}{2}R_{i_{l}}...R_{i_{1}}[\frac{dH}{dh}(\slashed{\Delta}\psi_{0}-\sum_{i}\psi_{i}\slashed{\Delta}\psi_{i}
-\sum_{i}\slashed{d}\psi_{i}\cdot\slashed{d}\psi_{i})]\\\notag=m_{T}^{\alpha}R_{i_{l}}...R_{i_{1}}\slashed{\Delta}
\psi_{\alpha}+\leftexp{(i_{1}...i_{l})}{\tilde{n}}_{l}
\end{align}
Here, $m_{T}^{\alpha}$ are defined in Chapter 10, and $\leftexp{(i_{1}...i_{l})}{\tilde{n}}_{l}$ is a lower
 order term (of order $l+1$):
\begin{align}
 \leftexp{(i_{1}...i_{l})}{\tilde{n}}_{l}=\sum_{|s_{1}|>0}((R)^{s_{1}}m^{\alpha}_{T})((R)^{s_{2}}\slashed{\Delta}
\psi_{\alpha})+R_{i_{l}}...R_{i_{1}}\tilde{n}_{0}
\end{align}
where
\begin{equation}
 \tilde{n}_{0}=-\frac{1}{2}\frac{dH}{dh}\sum_{i}\slashed{d}\psi_{i}\cdot\slashed{d}\psi_{i}
\end{equation}
From the bootstrap assumption, we have:
\begin{equation}
 \|\leftexp{(i_{1}...i_{l})}{\tilde{n}}_{l}\|_{L^{2}(\Sigma_{t}^{\epsilon_{0}})}\leq C_{l}\delta_{0}
(1+t)^{-3}\mathcal{W}_{\{l+1\}}
\end{equation}
According to (10.185):
\begin{equation}
 |m^{0}_{T}-\frac{1}{2}\ell|\leq C\delta_{0}(1+t)^{-1},\quad |m^{i}_{T}|\leq C\delta_{0}(1+t)^{-1}
\end{equation}
Let us define the functions:
\begin{align}
 \leftexp{(i_{1}...i_{l})}{\rho}^{(0)}_{l}=\frac{1}{2}\ell\slashed{\Delta} R_{i_{l}}...R_{i_{1}}\psi_{0}\\
\leftexp{(i_{1}...i_{l})}{\rho}^{(1)}_{l}=(m^{0}_{T}-\frac{1}{2}\ell)\slashed{\Delta}R_{i_{l}}...R_{i_{1}}
\psi_{0}+m^{i}_{T}\slashed{\Delta}R_{i_{l}}...R_{i_{1}}\psi_{i}
\end{align}
We can then express (14.221) as:
\begin{align}
 -R_{i_{l}}...R_{i_{1}}(\textrm{tr}\alpha'^{[P]})=\leftexp{(i_{1}...i_{l})}{\rho}^{(0)}_{l}
+\leftexp{(i_{1}...i_{l})}{\rho}^{(1)}_{l}+\leftexp{(i_{1}...i_{l})}{\tilde{n}}_{l}
\end{align}
Let us define $\leftexp{(i_{1}...i_{l})}{\rho}^{(2)}_{l}$ by:
\begin{align}
 \leftexp{(i_{1}...i_{l})}{\rho}_{l}=\leftexp{(i_{1}...i_{l})}{\rho}^{(0)}_{l}+\leftexp{(i_{1}...i_{l})}
{\rho}^{(1)}_{l}+\leftexp{(i_{1}...i_{l})}{\rho}^{(2)}_{l}
\end{align}
We shall estimate the $L^{2}(\Sigma_{t}^{\epsilon_{0}})$ norm of $\leftexp{(i_{1}...i_{l})}{\rho}^{(2)}_{l}$. 
First, we estimate the contribution from
\begin{align*}
\sum_{k=0}^{l-1}R_{i_{l}}...R_{i_{l-k+1}}\leftexp{(R_{i_{l-k}})}{Z}R_{i_{l-k-1}}...R_{i_{1}}\textrm{tr}\chi'
\end{align*}
 Obviously, this contribution is bounded by:
\begin{align}
(1+t)^{-2}[1+\log(1+t)]\mathcal{A}'_{[l]}(t)
\end{align}
where we have used the bounds for $\leftexp{(R_{i})}{Z}$ in Corollary 10.1.i.

Next, we should estimate the contribution from $R_{i_{l}}...R_{i_{1}}\rho_{0}$. Recall the definition:
\begin{align*}
 \rho_{0}=e\textrm{tr}\chi'+\frac{2e}{1-u+t}-|\chi'|^{2}-\textrm{tr}\alpha'
\end{align*}
By Proposition 12.6, (12.293) and (12.295), the contribution from the first term of above is bounded by:
\begin{align}
 C_{l}(1+t)^{-3}[1+\log(1+t)]\mathcal{W}^{Q}_{\{l\}}+\delta_{0}(1+t)^{-3}[\mathcal{Y}_{0}+(1+t)
\mathcal{A}'_{[l]}+\mathcal{W}_{[l]}]
\end{align}
Also by (12.295), the contribution from the second term is bounded by:
\begin{align}
 C_{l}(1+t)^{-2}\{\mathcal{W}^{Q}_{\{l\}}+\delta_{0}(1+t)^{-1}[\mathcal{Y}_{0}+(1+t)\mathcal{A}'_{[l]}
+\mathcal{W}_{\{l\}}]\}
\end{align}
By Proposition 12.6, the contribution from the third term is bounded by:
\begin{align}
 (1+t)^{-2}[1+\log(1+t)]\{\mathcal{A}'_{[l]}+\delta_{0}(1+t)^{-2}[1+\log(1+t)]\cdot\\\notag
[\mathcal{W}_{[l]}+\mathcal{Y}_{0}]\}
\end{align}
where we have used Corollary 10.2.d.

Finally, by the expression of $\alpha^{[N]}_{AB}$ and Corollary 10.1.g and Corollary 10.2.g, the contribution from 
the fourth term is bounded by:
\begin{align}
 (1+t)^{-3}[\mathcal{W}_{\{l+1\}}+\mathcal{W}^{Q}_{\{l\}}+(1+t)^{-1}(\mathcal{Y}_{0}+(1+t)\mathcal{A}'_{[l]})]
\end{align}
So we obtain the following estimate for $\leftexp{(i_{1}...i_{l})}{\rho}^{(2)}_{l}$:
\begin{align}
 \|\leftexp{(i_{1}...i_{l})}{\rho}^{(2)}_{l}\|_{L^{2}(\Sigma_{t}^{\epsilon_{0}})}\leq C_{l}(1+t)^{-2}
\mathcal{W}^{Q}_{\{l\}}\\\notag
+\delta_{0}(1+t)^{-3}[1+\log(1+t)](\mathcal{W}_{\{l+1\}}+\mathcal{Y}_{0}+(1+t)\mathcal{A}'_{[l]})
\end{align}
On the other hand, by $\textbf{H1}$ we have, pointwise:
\begin{equation}
 |\leftexp{(i_{1}...i_{l})}{\rho}^{(0)}_{l}|\leq C(1+t)^{-1}|\ell|\sqrt{\sum_{j,\alpha}|
\slashed{d}R_{j}R_{i_{l}}...R_{i_{1}}\psi_{\alpha}|^{2}}
\end{equation}
and 
\begin{equation}
 |\leftexp{(i_{1}...i_{l})}{\rho}^{(1)}_{l}|\leq C\delta_{0}(1+t)^{-2}\sqrt{\sum_{j,\alpha}|\slashed{d}
R_{j}R_{i_{l}}...R_{i_{1}}\psi_{\alpha}|^{2}}
\end{equation}
These inequalities imply:
\begin{align}
 \|\leftexp{(i_{1}...i_{l})}{\rho}^{(0)}_{l}\|_{L^{2}(\Sigma_{t}^{\epsilon_{0}})}\leq C|\ell|(1+t)^{-2}
\bar{\mu}^{-1/2}_{m}(t)\sqrt{\sum_{j,\alpha}\mathcal{E}'_{1}[R_{j}R_{i_{l}}...R_{i_{1}}\psi_{\alpha}](t)}
\end{align}
and
\begin{align}
 \|\leftexp{(i_{1}...i_{l})}{\rho}^{(1)}_{l}\|_{L^{2}(\Sigma_{t}^{\epsilon_{0}})}\leq C\delta_{0}(1+t)^{-3}
\bar{\mu}^{-1/2}_{m}(t)\sqrt{\sum_{j,\alpha}\mathcal{E}'_{1}[R_{j}R_{i_{l}}...R_{i_{1}}\psi_{\alpha}](t)}
\end{align}
while:
\begin{align}
 \|\mu^{1/2}\leftexp{(i_{1}...i_{l})}{\rho}^{(0)}_{l}\|_{L^{2}(\Sigma_{t}^{\epsilon_{0}})}
\leq C|\ell|(1+t)^{-2}\sqrt{\sum_{j,\alpha}\mathcal{E}'_{1}[R_{j}R_{i_{l}}...R_{i_{1}}\psi_{\alpha}](t)}
\end{align}
and:
\begin{align}
 \|\mu^{1/2}\leftexp{(i_{1}...i_{l})}{\rho}^{(1)}_{l}\|_{L^{2}(\Sigma_{t}^{\epsilon_{0}})}\leq C\delta_{0}(1+t)^{-3}
\sqrt{\sum_{j,\alpha}\mathcal{E}'_{1}[R_{j}R_{i_{l}}...R_{i_{1}}\psi_{\alpha}](t)}
\end{align}
     We express the propagation equation (14.219) in acoustical coordinates $(t,u,\vartheta)$, so that $L$ becomes 
$\frac{\partial}{\partial t}$, and integrate it along each integral curve of $L$ from $\Sigma_{0}^{\epsilon_{0}}$ to 
obtain:
\begin{align}
 (1-u+t)^{2}(R_{i_{l}}...R_{i_{1}}\textrm{tr}\chi')(t,u,\vartheta)\\\notag
=(1-u)^{2}(R_{i_{l}}...R_{i_{1}}\textrm{tr}\chi')(0,u,\vartheta)+\int_{0}^{t}(1-u+t')^{2}
\leftexp{(i_{1}...i_{l})}{\rho}_{l}(t',u,\vartheta)dt'
\end{align}
It follows that:
\begin{align}
 |R_{i_{l}}...R_{i_{1}}\textrm{tr}\chi'|(t,u,\vartheta)\\\notag
\leq C(1+t)^{-2}\{|R_{i_{l}}...R_{i_{1}}\textrm{tr}\chi'|(0,u,\vartheta)+\sum_{k=0}^{2}\leftexp{(i_{1}...i_{l})}
{A}^{(k)}_{l}(t,u,\vartheta)\}
\end{align}
where:
\begin{align}
 \leftexp{(i_{1}...i_{l})}{A}^{(k)}(t,u,\vartheta)=\int_{0}^{t}(1+t')^{2}|\leftexp{(i_{1}...i_{l})}{\rho}^{(k)}_{l}|
(t',u,\vartheta)dt'\quad :\quad k=0,1,2
\end{align}
By (8.333) we have:
\begin{align}
 \|\leftexp{(i_{1}...i_{l})}{A}^{(k)}_{l}(t)\|_{L^{2}([0,\epsilon_{0}]\times S^{2})}\leq 
\int_{0}^{t}(1+t')^{2}\|\leftexp{(i_{1}...i_{l})}{\rho}^{(k)}_{l}(t')\|_{L^{2}([0,\epsilon_{0}]\times S^{2})}dt'\\\notag
\leq C\int_{0}^{t}(1+t')\|\leftexp{(i_{1}...i_{l})}{\rho}^{(k)}_{l}\|_{L^{2}(\Sigma_{t'}^{\epsilon_{0}})}dt'
\end{align}
We now substitute the pointwise estimate (14.243) in (14.215). The $borderline$ contribution is from 
$\leftexp{(i_{1}...i_{l})}{A}^{(0)}_{l}$. This contribution is the $borderline$ integral:
\begin{align}
 C\int_{\Sigma_{t}^{\epsilon_{0}}}(1+t)|T\psi_{\alpha}||\leftexp{(i_{1}...i_{l})}{A}^{(0)}_{l}|
|\slashed{d}R_{i_{l+1}}...R_{i_{1}}\psi_{\alpha}|dud\mu_{\slashed{g}}
\end{align}
Recalling from Chapter 8 the partition of $[0,\epsilon_{0}]\times S^{2}$ into the subsets $\mathcal{V}_{s-}$, 
$\mathcal{V}_{s+}$, defined by (8.337), (8.338), respectively, we shall consider separately the integrals over 
the corresponding regions $\mathcal{U}^{-}_{s,t}$, $\mathcal{U}^{+}_{s,t}$ of $\Sigma_{t}^{\epsilon_{0}}$:
\begin{align}
 \mathcal{U}^{-}_{s,t}=\{(t,u,\vartheta)\in\Sigma^{\epsilon_{0}}_{t}\quad:\quad (u,\vartheta)\in\mathcal{V}_{s-}\}\\
\mathcal{U}^{+}_{s,t}=\{(t,u,\vartheta)\in\Sigma^{\epsilon_{0}}_{t}\quad:\quad (u,\vartheta)\in\mathcal{V}_{s+}\}
\end{align}
   
    We consider first the integral over $\mathcal{U}^{-}_{s,t}$. We have:
\begin{align}
 \int_{\mathcal{U}^{-}_{s,t}}(1+t)|T\psi_{\alpha}||\leftexp{(i_{1}...i_{l})}{A}^{(0)}_{l}|
|\slashed{d}R_{i_{l+1}}...R_{i_{1}}\psi_{\alpha}|dud\mu_{\slashed{g}}\\\notag
\leq (1+t)\sup_{\mathcal{U}^{-}_{s,t}}|T\psi_{\alpha}|\|\leftexp{(i_{1}...i_{l})}{A}^{(0)}_{l}\|
_{L^{2}(\mathcal{U}^{-}_{s,t})}\|\slashed{d}R_{i_{l+1}}...R_{i_{1}}\psi_{\alpha}\|
_{L^{2}(\Sigma_{t}^{\epsilon_{0}})}\\\notag
\leq C(1+t)\bar{\mu}^{-1/2}_{m}(t)\sup_{\mathcal{U}^{-}_{s,t}}|T\psi_{\alpha}|\|
\leftexp{(i_{1}...i_{l})}{A}^{(0)}_{l}(t)\|_{L^{2}(\mathcal{V}_{s-})}\sqrt{\mathcal{E}'_{1}
[R_{i_{l+1}}...R_{i_{1}}\psi_{\alpha}](t)}
\end{align}
by (8.333). Now by (14.87) and (14.88), we have, pointwise on $\Sigma_{t}^{\epsilon_{0}}$:
\begin{align}
 \max_{\alpha}|T\psi_{\alpha}|\leq \frac{C}{|\ell|}\{|L\mu|+C\delta_{0}(1+t)^{-2}[1+\log(1+t)]\}
\end{align}
In estimating the integral over $\mathcal{U}^{-}_{s,t}$ we can assume that $\mathcal{V}_{s-}$ is non-empty, 
thus, recalling the definitions (8.251) and (8.261),
\begin{align}
 \min_{(u,\vartheta)\in\mathcal{V}_{s-}}\hat{E}_{s}(u,\vartheta)=\hat{E}_{s,m}=-\delta_{1},\quad \delta_{1}>0
\end{align}
Now, from Proposition 8.6 we have:
\begin{align}
 (L\mu)(t,u,\vartheta)=\mu_{[1],s}(u,\vartheta)\{\frac{\hat{E}_{s}(u,\vartheta)}{(1+t)}+\hat{Q}_{1,s}(t,u,\vartheta)\}
\end{align}
It follows that:
\begin{align}
 \sup_{\mathcal{U}^{-}_{s,t}}|L\mu|\leq C\{\frac{\delta_{1}}{(1+t)}+\frac{C\delta_{0}[1+\log(1+t)]}{(1+t)^{2}}\}
\end{align}
hence, substituting in (14.250),
\begin{align}
 \max_{\alpha}\sup_{\mathcal{U}^{-}_{s,t}}|T\psi_{\alpha}|\leq \frac{C}{|\ell|}\{\frac{\delta_{1}}{(1+t)}
+\frac{C\delta_{0}[1+\log(1+t)]}{(1+t)^{2}}\}
\end{align}
On the other hand, from (14.245) for $k=0$ and (14.238),
\begin{align}
 \|\leftexp{(i_{1}...i_{l})}{A}^{(0)}_{l}(t)\|_{L^{2}([0,\epsilon_{0}]\times S^{2})}\leq \\\notag
C|\ell|\int_{0}^{t}(1+t')^{-1}\bar{\mu}^{-1/2}_{m}(t')\sqrt{\sum_{j,\alpha}\mathcal{E}'_{1}
[R_{j}R_{i_{l}}...R_{i_{1}}\psi_{\alpha}](t')}dt'
\end{align}
Substituting (14.254) and (14.255) in (14.249), the factors $|\ell|$ cancel and we obtain that (14.249) is bounded by:
\begin{align}
 C\bar{\mu}^{-1/2}_{m}(t)\{\delta_{1}+\frac{C\delta_{0}[1+\log(1+t)]}{(1+t)}\}\sqrt{\mathcal{E}'_{1}
[R_{i_{l+1}}...R_{i_{1}}\psi_{\alpha}](t)}\\\notag
\cdot\int_{0}^{t}(1+t')^{-1}\bar{\mu}^{-1/2}_{m}(t')\sqrt{\sum_{j,\alpha}\mathcal{E}'_{1}
[R_{j}R_{i_{l}}...R_{i_{1}}\psi_{\alpha}](t')}dt'
\end{align}
We now define
\begin{align}
 \leftexp{(i_{1}...i_{l})}{\mathcal{G}}^{\prime}_{1,l+2;a,q}(t)=\sup_{t'\in[0,t]}
\{[1+\log(1+t)]^{-2q}\bar{\mu}^{2a}_{m}(t')\sum_{j,\alpha}\mathcal{E}'_{1}[R_{j}R_{i_{l}}...R_{i_{1}}
\psi_{\alpha}](t')\}
\end{align}
Then (14.256) is bounded by:
\begin{align}
 C[1+\log(1+t)]^{2q}\bar{\mu}_{m}^{-a-1/2}(t)\{\delta_{1}+\frac{C\delta_{0}[1+\log(1+t)]}{(1+t)}\}\\\notag
\cdot\int_{0}^{t}\bar{\mu}^{-a-1/2}_{m}(t')\frac{dt'}{(1+t')}\cdot\leftexp{(i_{1}...i_{l})}
{\mathcal{G}}^{\prime}_{1,l+2;a,q}(t)
\end{align}
To estimate the integral in (14.258), we follow the argument of Lemma 8.11. Since we are assuming that 
$\mathcal{V}_{s-}$ is non-empty, we only have Case 2 to consider. Setting again, as in (8.267):
\begin{align*}
 t_{1}=e^{\frac{1}{2a\delta_{1}}}-1
\end{align*}
we have two subcases of Case 2 as in the proof of Lemma 8.11. In Subcase 2a the lower bound (8.273) holds, 
hence $\bar{\mu}^{-a-1/2}_{m}(t')\leq C$ and
\begin{align}
 \int_{0}^{t}\bar{\mu}^{-a-1/2}_{m}(t')\frac{dt'}{(1+t')}\leq C\int_{0}^{t}\frac{dt'}{(1+t')}=C\log(1+t)
\end{align}
therefore:
\begin{align}
 \{\delta_{1}+\frac{C\delta_{0}[1+\log(1+t)]}{(1+t)}\}\int_{0}^{t}\bar{\mu}^{-a-1/2}_{m}(t')\frac{dt'}{(1+t')}\\\notag
\leq C\{\delta_{1}\log(1+t)+C\delta_{0}\frac{[1+\log(1+t)]^{2}}{(1+t)}\}\leq \frac{C}{a}
\end{align}
for, in Subcase 2a
\begin{align*}
 \delta_{1}\log(1+t)\leq\delta_{1}\log(1+t_{1})=\frac{1}{2a}
\end{align*}
while
\begin{align*}
 \frac{[1+\log(1+t)]^{2}}{(1+t)}\leq C\quad\textrm{and}\quad C\delta_{0}\leq\frac{1}{a}
\end{align*}
In subcase 2b we first estimate the contribution of 
\begin{align*}
 \int_{0}^{t_{1}}\bar{\mu}^{-a-1/2}_{m}(t')\frac{dt'}{(1+t')}\leq C\log(1+t_{1})=\frac{C}{2a\delta_{1}}
\end{align*}
to the left-hand side of (14.260). This contribution is then bounded by (Note that $t\geq t_{1}$):
\begin{align}
 \{\delta_{1}+C\delta_{0}\frac{[1+\log(1+t)]}{(1+t)}\}\frac{C}{2a\delta_{1}}\leq\{\delta_{1}+C\delta_{0}\frac{[1+\log(1+t_{1})]}{(1+t_{1})}\}
\frac{C}{2a\delta_{1}}\\\notag
\leq\frac{C}{2a}\{1+\frac{C\delta_{0}}{\delta_{1}}(1+\frac{1}{2a\delta_{1}})e^{-\frac{1}{2a\delta_{1}}}\}
\leq\frac{C}{2a}(1+C\varphi'_{2}(2a\delta_{1}))\leq\frac{C}{a}
\end{align}
(see (14.157)). Next, we estimate the contribution of:
\begin{equation}
 \int_{t_{1}}^{t}\bar{\mu}^{-a-1/2}_{m}(t')\frac{dt'}{(1+t')}
\end{equation}
By the lower bound (8.303),
\begin{align}
 \int_{t_{1}}^{t}\bar{\mu}^{-a-1/2}_{m}(t')\frac{dt'}{(1+t')}\leq C\int^{\tau}_{\tau_{1}}
(1-\delta_{1}\tau')^{-a-1/2}d\tau'\\\notag
\leq\frac{C}{\delta_{1}}\frac{(1-\delta_{1}\tau)^{-a+1/2}}{(a-1/2)}\leq\frac{C}{\delta_{1}}
\frac{\bar{\mu}^{-a+1/2}_{m}(t)}{(a-1/2)}
\end{align}
where in the last step, we have used (8.312). It follows that the contribution of (14.262) to 
the left-hand side of (14.260) is bounded by:
\begin{align}
 \{\delta_{1}+C\delta_{0}\frac{[1+\log(1+t)]}{(1+t)}\}\frac{C}{\delta_{1}}\frac{\bar{\mu}^{-a+1/2}_{m}(t)}
{(a-1/2)}\\\notag
\leq\frac{C\bar{\mu}^{-a+1/2}_{m}(t)}{(a-1/2)}\{1+\frac{C\delta_{0}}{\delta_{1}}\frac{[1+\log(1+t_{1})]}
{(1+t_{1})}\}\\\notag
\leq\frac{C\bar{\mu}^{-a+1/2}_{m}(t)}{(a-1/2)}(1+C\varphi'_{2}(2a\delta_{1}))\leq\frac{C\bar{\mu}^{-a+1/2}_{m}(t)}
{(a-1/2)}
\end{align}
Combining the results (14.260), (14.264) of the two subcases we conclude that, in general:
\begin{align}
 \{\delta_{1}+C\delta_{0}\frac{[1+\log(1+t)]}{(1+t)}\}\int_{0}^{t}\bar{\mu}^{-a-1/2}_{m}(t')
\frac{dt'}{(1+t')}\leq C\frac{\bar{\mu}^{-a+1/2}_{m}(t)}{(a-1/2)}
\end{align}
hence (14.258), therefore also (14.256) and the integral on the left in (14.249), is bounded by:
\begin{align}
 \frac{C}{(a-1/2)}\bar{\mu}^{-2a}_{m}(t)[1+\log(1+t)]^{2q}\leftexp{(i_{1}...i_{l})}{\mathcal{G}}^{\prime}_{1,l+2;a,q}(t)
\end{align}

     We now consider the integral over $\mathcal{U}^{+}_{s,t}$:
\begin{align}
 \int_{\mathcal{U}^{+}_{s,t}}(1+t)|T\psi_{\alpha}|\leftexp{(i_{1}...i_{l})}{A}^{(0)}_{l}
|\slashed{d}R_{i_{l+1}}...R_{i_{1}}\psi_{\alpha}|dud\mu_{\slashed{g}}
\end{align}
We estimate this by:
\begin{align}
 (1+t)\sup_{\mathcal{U}^{+}_{s,t}}(\mu^{-1}|T\psi_{\alpha}|)\|\mu^{1/2}\leftexp{(i_{1}...i_{l})}
{A}^{(0)}_{l}\|_{L^{2}(\mathcal{U}^{+}_{s,t})}
\|\mu^{1/2}\slashed{d}R_{i_{l+1}}...R_{i_{1}}\psi_{\alpha}\|_{L^{2}(\Sigma_{t}^{\epsilon_{0}})}\\\notag
\leq C(1+t)\sup_{\mathcal{U}^{+}_{s,t}}(\mu^{-1}|T\psi_{\alpha}|)\|(\mu^{1/2}\leftexp{(i_{1}...i_{l})}
{A}^{(0)}_{l})(t)\|_{L^{2}(\mathcal{V}_{s+})}
\sqrt{\mathcal{E}'_{1}[R_{i_{l+1}}...R_{i_{1}}\psi_{\alpha}](t)}
\end{align}
Since by Proposition 8.6, in $\mathcal{U}^{+}_{s,t}$ we have $\mu\geq C^{-1}$, then the pointwise estimate (14.250) 
implies:
\begin{align}
 \max_{\alpha}\sup_{\mathcal{U}^{+}_{s,t}}(\mu^{-1}|T\psi_{\alpha}|)\leq\frac{C}{|\ell|}
\{\sup_{\mathcal{U}^{+}_{s,t}}(\mu^{-1}|L\mu|)+C\delta_{0}(1+t)^{-2}[1+\log(1+t)]\}
\end{align}
 
Let us define:
\begin{align}
 \hat{E}_{s,M}=\max_{(u,\vartheta)\in[0,\epsilon_{0}]\times S^{2}}\hat{E}_{s}(u,\vartheta)
\end{align}
Setting:
\begin{equation}
 \delta_{2}=\hat{E}_{s,M}
\end{equation}
we can assume in estimating (14.269) that $\delta_{2}>0$. From Proposition 8.6 we have:
\begin{align}
 \mu^{-1}L\mu=\frac{\hat{E}_{s}(u,\vartheta)(1+t)^{-1}+\hat{Q}_{1,s}(t,u,\vartheta)}{1+\hat{E}_{s}
(u,\vartheta)\log(1+t)+\hat{Q}_{0,s}(t,u,\vartheta)}
\end{align}
In view of the fact that $\hat{E}_{s}(u,\vartheta)\geq0$ in $\mathcal{V}_{s+}$, (14.272) implies 
that in $\mathcal{U}^{+}_{s,t}$:
\begin{align}
 \mu^{-1}|L\mu|\leq\frac{\hat{E}_{s}(u,\vartheta)(1+t)^{-1}+|\hat{Q}_{1,s}(t,u,\vartheta)|}
{1+\hat{E}_{s}(u,\vartheta)\log(1+t)-|\hat{Q}_{0,s}(t,u,\vartheta)|}\\\notag
\leq\frac{\hat{E}_{s}(u,\vartheta)(1+t)^{-1}+C\delta_{0}(1+t)^{-2}[1+\log(1+t)]}
{1+\hat{E}_{s}(u,\vartheta)\log(1+t)-C\delta_{0}(1+t)^{-1}[1+\log(1+t)]}\\\notag
\leq\frac{\hat{E}_{s}(u,\vartheta)(1+t)^{-1}}{1+\hat{E}_{s}(u,\vartheta)\log(1+t)
-C\delta_{0}(1+t)^{-1}[1+\log(1+t)]}\\\notag
+C'\delta_{0}\frac{[1+\log(1+t)]}{(1+t)^{2}}
\end{align}
With
\begin{align*}
 \eta=\hat{E}_{s}(u,\vartheta)\log(1+t)\geq0,\quad \epsilon=C\delta_{0}(1+t)^{-1}[1+\log(1+t)]>0
\end{align*}
we write
\begin{align*}
 \frac{1}{1+\eta-\epsilon}=\frac{1}{1+\eta}+\frac{\epsilon}{(1+\eta-\epsilon)(1+\eta)}\leq\frac{1}{1+\eta}+C\epsilon
\end{align*}
then since
\begin{align*}
 \hat{E}_{s}(u,\vartheta)(1+t)^{-1}\epsilon\leq C''\delta_{0}^{2}(1+t)^{-2}[1+\log(1+t)]
\end{align*}
(14.273) implies that in $\mathcal{U}^{+}_{s,t}$:
\begin{align}
 \mu^{-1}|L\mu|\leq\frac{\hat{E}_{s}(u,\vartheta)(1+t)^{-1}}{1+\hat{E}_{s}(u,\vartheta)\log(1+t)}
+C\delta_{0}(1+t)^{-2}[1+\log(1+t)]
\end{align}
The first term on the right in (14.274) is
\begin{align}
 \frac{1}{(1+t)\log(1+t)}\frac{\eta}{1+\eta}
\end{align}
This is an increasing function of $\eta$, and achieves its maximum in $\mathcal{U}^{+}_{s,t}$ where $\eta$ 
achieves its maximum value
\begin{align}
 \eta_{M}=\delta_{2}\log(1+t)
\end{align}
in $\mathcal{V}_{s+}$. The supremum of (14.277) in $\mathcal{U}^{+}_{s,t}$ is therefore:
\begin{align}
 \frac{1}{(1+t)\log(1+t)}\frac{\eta_{M}}{1+\eta_{M}}=\frac{1}{(1+t)}\frac{\delta_{2}}{[1+\delta_{2}\log(1+t)]}
\end{align}
It follows that:
\begin{align}
 \sup_{\mathcal{U}^{+}_{s,t}}\mu^{-1}|L\mu|\leq\frac{1}{(1+t)}\frac{\delta_{2}}{[1+\delta_{2}\log(1+t)]}
+C\delta_{0}\frac{[1+\log(1+t)]}{(1+t)^{2}}
\end{align}
Substituting this in (14.269) and the result in (14.268), we conclude that the integral (14.267) is bounded by:
\begin{align}
 \frac{C}{|\ell|}\{\frac{\delta_{2}}{1+\delta_{2}\log(1+t)}+C\delta_{0}\frac{[1+\log(1+t)]}{(1+t)}\}\cdot\\\notag
\|\mu^{1/2}\leftexp{(i_{1}...i_{l})}{A}^{(0)}_{l}(t)\|_{L^{2}(\mathcal{V}_{s+})}\sqrt{\mathcal{E}'_{1}
[R_{i_{l+1}}...R_{i_{1}}\psi_{\alpha}](t)}
\end{align}
Now from (14.244) for $k=0$ we have:
\begin{align}
 (\mu^{1/2}\leftexp{(i_{1}...i_{l})}{A}^{(0)}_{l})(t,u,\vartheta)=
\int_{0}^{t}(1+t')^{2}(\frac{\mu(t,u,\vartheta)}{\mu(t',u,\vartheta)})^{1/2}|(\mu^{1/2}\leftexp{(i_{1}...i_{l})}
{\rho}^{(0)}_{l})(t',u,\vartheta)|dt'
\end{align}
It follows that:
\begin{align}
 \|(\mu^{1/2}\leftexp{(i_{1}...i_{l})}{A}^{(0)}_{l})(t)\|_{L^{2}(\mathcal{V}_{s+})}\\\notag
\leq\int_{0}^{t}(1+t')^{2}[\sup_{(u,\vartheta)\in\mathcal{V}_{s+}}(\frac{\mu(t,u,\vartheta)}
{\mu(t',u,\vartheta)})]^{1/2}
\|(\mu^{1/2}\leftexp{(i_{1}...i_{l})}{\rho}^{(0)}_{l})(t')\|_{L^{2}(\mathcal{V}_{s+})}dt'\\\notag
\leq C\int_{0}^{t}(1+t')[\sup_{(u,\vartheta)\in\mathcal{V}_{s+}}(\frac{\mu(t,u,\vartheta)}
{\mu(t',u,\vartheta)})]^{1/2}\|\mu^{1/2}
\leftexp{(i_{1}...i_{l})}{\rho}^{(0)}_{l}\|_{L^{2}(\Sigma_{t'}^{\epsilon_{0}})}dt'
\end{align}
where in the last step, we have used (8.333). From Proposition 8.6 we have, for $(u,\vartheta)\in\mathcal{V}_{s+}$,
\begin{align}
 \frac{\mu(t,u,\vartheta)}{\mu(t',u,\vartheta)}=\frac{1+\hat{E}_{s}(u,\vartheta)\log(1+t)+\hat{Q}_{0,s}(t,u,\vartheta)}
{1+\hat{E}_{s}(u,\vartheta)\log(1+t')+\hat{Q}_{0,s}(t',u,\vartheta)}\\\notag
\leq\frac{1+\hat{E}_{s}(u,\vartheta)\log(1+t)+C\delta_{0}(1+t)^{-1}[1+\log(1+t)]}
{1+\hat{E}_{s}(u,\vartheta)\log(1+t')-C\delta_{0}(1+t')^{-1}[1+\log(1+t')]}
\end{align}
Since
\begin{align*}
 C\delta_{0}\frac{[1+\log(1+t')]}{(1+t')}\leq C\delta_{0}\leq\frac{1}{2}
\end{align*}
the denominator in the fraction on the right in (14.282) is
\begin{align*}
 \geq\frac{1}{2}+\hat{E}_{s}(u,\vartheta)\log(1+t')\geq\frac{1}{2}(1+\hat{E}_{s}(u,\vartheta)\log(1+t'))
\end{align*}
Similarly, since
\begin{align*}
 C\delta_{0}\frac{[1+\log(1+t)]}{(1+t)}\leq C\delta_{0}\leq\frac{1}{2}
\end{align*}
the numerator in the fraction on the right in (14.282) is
\begin{align*}
 \frac{3}{2}+\hat{E}_{s}(u,\vartheta)\log(1+t)\leq\frac{3}{2}(1+\hat{E}_{s}(u,\vartheta)\log(1+t))
\end{align*}
Therefore, for $(u,\vartheta)\in\mathcal{V}_{s+}$ it holds that
\begin{align}
 \frac{\mu(t,u,\vartheta)}{\mu(t',u,\vartheta)}\leq3\frac{1+\hat{E}_{s}(u,\vartheta)
\log(1+t)}{1+\hat{E}_{s}(u,\vartheta)\log(1+t')}
\end{align}
Now for $a>b>0, x\geq0$, the function
\begin{align*}
 f(x)=\frac{1+ax}{1+bx}
\end{align*}
is increasing in $x$. Hence, taking $a=\log(1+t)$, $b=\log(1+t')$, $x=\hat{E}_{s}(u,\vartheta)$, 
the ratio on the right in (14.283) achieves its maximum in $\mathcal{V}_{s+}$ where $x$ achieves 
its maximum $\delta_{2}$ in $\mathcal{V}_{s+}$. Therefore:
\begin{equation}
 \sup_{(u,\vartheta)\in\mathcal{V}_{s+}}\frac{\mu(t,u,\vartheta)}{\mu(t',u,\vartheta)}
\leq3\frac{1+\delta_{2}\log(1+t)}{1+\delta_{2}\log(1+t')}
\end{equation}
Substituting this and (14.240) in (14.241) we obtain:
\begin{align}
 \|(\mu^{1/2}\leftexp{(i_{1}...i_{l})}{A}^{(0)}_{l})(t)\|_{L^{2}(\mathcal{V}_{s+})}\\\notag
\leq C|\ell|\int_{0}^{t}\sqrt{\frac{1+\delta_{2}\log(1+t)}{1+\delta_{2}\log(1+t')}}
\sqrt{\sum_{j,\alpha}\mathcal{E}'_{1}[R_{j}R_{i_{l}}...R_{i_{1}}\psi_{\alpha}](t')}
\frac{dt'}{(1+t')}\\\notag
\leq C|\ell|\sqrt{\leftexp{(i_{1}...i_{l})}{\mathcal{G}'_{1,l+2;a,q}}}\int_{0}^{t}
\sqrt{\frac{1+\delta_{2}\log(1+t)}{1+\delta_{2}\log(1+t')}}
[1+\log(1+t')]^{q}\bar{\mu}^{-a}_{m}(t')\frac{dt'}{(1+t')}\\\notag
\leq C|\ell|\bar{\mu}^{-a}_{m}(t)\sqrt{\leftexp{(i_{1}...i_{l})}
{\mathcal{G}'_{1,l+2;a,q}}}\int_{0}^{t}\sqrt{\frac{1+\delta_{2}\log(1+t)}{1+\delta_{2}\log(1+t')}}
[1+\log(1+t')]^{q}\frac{dt'}{(1+t')}
\end{align}
where in the last step we have used the fact that by Corollary 2 of Lemma 8.11, 
$\bar{\mu}^{-a}_{m}(t')\leq C\bar{\mu}^{-a}_{m}(t)$. Substituting (14.285) in (14.279) and noting that
\begin{align*}
 \sqrt{\mathcal{E}'_{1}[R_{i_{l+1}}...R_{i_{1}}\psi_{\alpha}](t)}\leq \bar{\mu}^{-a}_{m}(t)[1+\log(1+t)]^{q}
\sqrt{\leftexp{(i_{1}...i_{l})}{\mathcal{G}}^{\prime}_{1,l+2;a,q}(t)}
\end{align*}
the factors $|\ell|$ cancel and we obtain that (14.279) hence also (14.267) is bounded by:
\begin{align}
 C\bar{\mu}^{-2a}_{m}(t)[1+\log(1+t)]^{q}\leftexp{(i_{1}...i_{l})}{\mathcal{G}}^{\prime}_{1,l+2;a,q}(t)\\\notag
\cdot\{\frac{\delta_{2}}{\sqrt{1+\delta_{2}\log(1+t)}}+C\delta_{0}\frac{[1+\log(1+t)]^{3/2}}{(1+t)}\}I_{q;\delta_{2}}(t)
\end{align}
where:
\begin{align}
 I_{q;\delta_{2}}(t)=\int_{0}^{t}\frac{[1+\log(1+t')]^{q}}{\sqrt{1+\delta_{2}\log(1+t')}}\frac{dt'}{(1+t')}
\end{align}
Setting $x=1+\log(1+t')$, the integral $I_{q;\delta_{2}}(t)$ takes the form:
\begin{align}
 I_{q;\delta_{2}}(t)=\int_{1}^{1+\log(1+t)}\frac{x^{q}}{\sqrt{1+\delta_{2}(x-1)}}dx
\end{align}
Since $\delta_{2}\leq 1$ we have:
\begin{align*}
 1+\delta_{2}(x-1)\geq\delta_{2}+\delta_{2}(x-1)=\delta_{2}x
\end{align*}
hence:
\begin{align}
 I_{q;\delta_{2}}(t)\leq\frac{1}{\sqrt{\delta_{2}}}\int_{0}^{1+\log(1+t)}x^{q-1/2}dx=\frac{1}
{\sqrt{\delta_{2}}}\frac{[1+\log(1+t)]^{q+1/2}}{(q+1/2)}
\end{align}
Also, since the denominator of the integrand in (14.290) is $\geq1$,
\begin{align}
 I_{q;\delta_{2}}(t)\leq\int_{0}^{1+\log(1+t)}x^{q}dx=\frac{[1+\log(1+t)]^{q+1}}{(q+1)}
\end{align}
We use (14.289) to estimate the product:
\begin{align}
 \frac{\delta_{2}}{\sqrt{1+\delta_{2}\log(1+t)}}I_{q;\delta_{2}}(t)\leq\frac{\delta_{2}}
{\sqrt{1+\delta_{2}\log(1+t)}}\frac{[1+\log(1+t)]^{q+1/2}}
{\sqrt{\delta_{2}}(q+1/2)}\\\notag
=\frac{\sqrt{\delta_{2}}\sqrt{1+\log(1+t)}}{\sqrt{1+\delta_{2}\log(1+t)}}\frac{[1+\log(1+t)]^{q}}
{(q+1/2)}\leq\frac{[1+\log(1+t)]^{q}}{(q+1/2)}
\end{align}
where in the last step we have used the fact that $\delta_{2}\leq1$. We use (14.290) to estimate the product:
\begin{align}
 C\delta_{0}\frac{[1+\log(1+t)]^{3/2}}{(1+t)}I_{q;\delta_{2}}(t)\leq C\delta_{0}
\frac{[1+\log(1+t)]^{3/2}}{(1+t)}\frac{[1+\log(1+t)]^{q+1}}{(q+1)}\\\notag
=C\delta_{0}\frac{[1+\log(1+t)]^{5/2}}{(1+t)}\frac{[1+\log(1+t)]^{q}}{(q+1)}\leq
\frac{C'\delta_{0}}{(q+1)}[1+\log(1+t)]^{q}
\end{align}
We conclude that (14.286) hence also (14.267) is bounded by:
\begin{align}
 \frac{C}{(q+1/2)}\bar{\mu}^{-2a}_{m}(t)[1+\log(1+t)]^{2q}\leftexp{(i_{1}...i_{l})}{\mathcal{G}}^{\prime}_{1,l+2;a,q}(t)
\end{align}

     Thus, we have estimated the borderline integral (14.246). We proceed to estimate the remaining contributions 
to the hypersurface integral $H_{0}$. These are the contributions from (14.243). To estimate these contributions, 
we simply use the bound:
\begin{align}
 \max_{\alpha}\sup_{\Sigma_{t}^{\epsilon_{0}}}|T\psi_{\alpha}|\leq C\delta_{0}(1+t)^{-1}
\end{align}
These contributions are then bounded by, for $k=1,2$,
\begin{align}
 C\int_{\Sigma_{t}^{\epsilon_{0}}}(1+t)|T\psi_{\alpha}|\leftexp{(i_{1}...i_{l})}{A}^{(k)}_{l}
|\slashed{d}R_{i_{l+1}}...R_{i_{1}}\psi_{\alpha}|dud\mu_{\slashed{g}}\\\notag
\leq C\delta_{0}\int_{\Sigma_{t}^{\epsilon_{0}}}\leftexp{(i_{1}...i_{l})}{A}^{(k)}_{l}
|\slashed{d}R_{i_{l+1}}...R_{i_{1}}\psi_{\alpha}|dud\mu_{\slashed{g}}\\\notag
\leq C\delta_{0}\|\leftexp{(i_{1}...i_{l})}{A}^{(k)}_{l}\|_{L^{2}(\Sigma_{t}^{\epsilon_{0}})}
\|\slashed{d}R_{i_{l+1}}...R_{i_{1}}\psi_{\alpha}\|_{L^{2}(\Sigma_{t}^{\epsilon_{0}})}\\\notag
\leq C\delta_{0}\bar{\mu}^{-1/2}_{m}(t)\|\leftexp{(i_{1}...i_{l})}{A}^{(k)}_{l}(t)\|
_{L^{2}([0,\epsilon_{0}]\times S^{2})}\sqrt{\mathcal{E}'_{1,[l+2]}(t)}
\end{align}
by (8.333).

     Consider first the contribution of $\leftexp{(i_{1}...i_{l})}{A}^{(1)}_{l}$. From (14.245) and the fact that 
by the estimate (14.239):
\begin{align}
 \|\leftexp{(i_{1}...i_{l})}{\rho}^{(1)}_{l}\|_{L^{2}(\Sigma_{t}^{\epsilon_{0}})}\leq 
C\bar{\mu}^{-1/2}_{m}(t)\delta_{0}(1+t)^{-3}\sqrt{\mathcal{E}'_{1,[l+2]}(t)}
\end{align}
we obtain:
\begin{align}
 \|\leftexp{(i_{1}...i_{l})}{A}^{(1)}_{l}(t)\|_{L^{2}([0,\epsilon_{0}]\times S^{2})}\leq C\delta_{0}\int_{0}^{t}
\bar{\mu}^{-1/2}_{m}(t')(1+t')^{-2}\sqrt{\mathcal{E}'_{1,[l+2]}(t')}dt'\\\notag
\leq C\delta_{0}\sqrt{\mathcal{G}'_{1,[l+2];a,q}(t)}J'_{a,q}(t)
\end{align}
where:
\begin{align}
 J'_{a,q}(t)=\int_{0}^{t}\bar{\mu}^{-a-1/2}_{m}(t')(1+t')^{-2}[1+\log(1+t')]^{q}dt'
\end{align}
To estimate this integral we consider again the two cases occurring in the proof of Lemma 8.11. 
In Case 1 we have the upper bound (14.117), which implies:
\begin{align}
 J'_{a,q}(t)\leq C\int_{0}^{t}(1+t')^{-2}[1+\log(1+t')]^{q}dt'\leq C_{q}
\end{align}
Similarly, in Subcase 2a:
\begin{align}
 J'_{a,q}(t_{1})\leq C_{q}
\end{align}
In Subcase 2b we have the lower bound (14.124), which, with
\begin{align}
 C_{q}(t_{1})=\sup_{t'\geq t_{1}}\{(1+t')^{-1}[1+\log(1+t')]^{q}\}
\end{align}
implies:
\begin{align}
 J'_{a,q}(t)-J'_{a,q}(t_{1})\leq C\cdot C_{q}(t_{1})\int_{t_{1}}^{t}(1+t')^{-1}(1-\delta_{1}\tau')^{-a-1/2}dt'\\\notag
\leq C\cdot\frac{C_{q}(t_{1})}{a\delta_{1}}\bar{\mu}^{-a+1/2}_{m}(t)
\end{align}
In the last step we have used the bound (8.312).

We have:
\begin{align*}
 C_{q}(t_{1})=\sup_{\tau'\geq\tau_{1}}\{e^{-\tau'}(1+\tau')^{q}\}\quad \tau_{1}=\log(1+t_{1})
\end{align*}
Now, the function:
\begin{align*}
 f(x)=e^{-x}x^{q}\quad \textrm{on}\quad [1,\infty)
\end{align*}
is decreasing in $x$ for $x\geq q$. Hence if $\tau_{1}=1/2a\delta_{1}\geq q$, which is the case 
if $aq\delta_{0}$ is suitably small, then:
\begin{align}
 C_{q}(t_{1})=e^{-\tau_{1}}(1+\tau_{1})^{q}=e^{-\frac{1}{2a\delta_{1}}}(1+\frac{1}{2a\delta_{1}})^{q}
\end{align}
Therefore, from (14.302),
\begin{align}
 J'_{a,q}(t)-J'_{a,q}(t_{1})\leq\frac{C}{2a\delta_{1}}(1+\frac{1}{2a\delta_{1}})^{q}
e^{-\frac{1}{2a\delta_{1}}}\bar{\mu}^{-a+1/2}_{m}(t)\\\notag
\leq C\varphi'_{q+1}(2a\delta_{1})\bar{\mu}^{-a+1/2}_{m}(t)\leq C\varphi'_{q+1}(Ca\delta_{0})\bar{\mu}^{-a+1/2}_{m}(t)
\end{align}
We conclude that in general:
\begin{align}
 J'_{a,q}(t)\leq C_{q}+C\varphi'_{q+1}(Ca\delta_{0})\bar{\mu}^{-a+1/2}_{m}(t)
\end{align}
Substituting in (14.297) and the result in (14.295) for $k=1$, then yields:
\begin{align}
 C\int_{\Sigma_{t}^{\epsilon_{0}}}(1+t)|T\psi_{\alpha}|\leftexp{(i_{1}...i_{l})}{A}^{(1)}_{l}|
\slashed{d}R_{i_{l+1}}...R_{i_{1}}\psi_{\alpha}|dud\mu_{\slashed{g}}\\\notag
\leq C\delta_{0}^{2}\mathcal{G}'_{1,[l+2];a,q}(t)\{C_{q}\bar{\mu}^{-a-1/2}_{m}(t)
+C\varphi'_{q+1}(Ca\delta_{0})\bar{\mu}^{-2a}_{m}(t)\}[1+\log(1+t)]^{q}
\end{align}

    Next, we consider the contribution of $\leftexp{(i_{1}...i_{l})}{A}^{(2)}_{l}$. We have (14.245) for 
$k=2$ and for \\
$\|\leftexp{(i_{1}...i_{l})}{\rho}^{(2)}_{l}\|_{L^{2}(\Sigma_{t}^{\epsilon_{0}})}$ we have the estimate (14.235). 
Here we shall only consider the contribution of the leading 
term on the right of (14.235) namely the term $C_{l}(1+t)^{-2}\mathcal{W}^{Q}_{\{l\}}$.
The contribution of this term to (14.245) is:
\begin{align}
 C_{l}\int_{0}^{t}(1+t')^{-1}\mathcal{W}^{Q}_{\{l\}}(t')dt'
\end{align}
By Lemma 5.1 we have:
\begin{align}
 \mathcal{W}^{Q}_{\{l\}}(t')\leq C\epsilon_{0}\sqrt{\mathcal{E}_{0,[l+2]}(t')}\\\notag
\leq C\epsilon_{0}\bar{\mu}^{-a}_{m}(t')[1+\log(1+t')]^{p}\sqrt{\mathcal{G}_{0,[l+2];a,p}}
\end{align}
hence (14.307) is bounded by:
\begin{align}
 C_{l}\epsilon_{0}[1+\log(1+t)]^{p}\sqrt{\mathcal{G}_{0,[l+2];a,p}(t)}\int_{0}^{t}(1+t')^{-1}\bar{\mu}^{-a}_{m}(t')dt'
\end{align}
The last integral coincides with the integral $J_{a-1/2,-1}(t)$ (see (14.147)), therefore  by (14.153) it is bounded by:
\begin{align*}
 C[1+\log(1+t)]\bar{\mu}^{-a+1}_{m}
\end{align*}
Consequently (14.309), hence also (14.307), are bounded by:
\begin{align}
 C_{l}\epsilon_{0}[1+\log(1+t)]^{p+1}\bar{\mu}^{-a+1}_{m}(t)\sqrt{\mathcal{G}_{0,[l+2];a,p}(t)}
\end{align}
We conclude that the contribution to (14.295) for $k=2$ is bounded by:
\begin{align}
 C_{l}\epsilon_{0}\delta_{0}\bar{\mu}^{-2a+1/2}_{m}(t)[1+\log(1+t)]^{p+q+1}\sqrt{\mathcal{G}_{0,[l+2];a,p}(t)
\mathcal{G}'_{1,[l+2];a,q}(t)}
\end{align}

     Finally, we consider the contribution of the hypersurface integral $H_{1}$, given by (14.212). We have:
\begin{align}
 |H_{1}|\leq C\int_{\Sigma_{t}^{\epsilon_{0}}}(1+t)^{2}|R_{i_{l+1}}T\psi_{\alpha}||R_{i_{l}}...R_{i_{1}}\textrm{tr}\chi'|
|R_{i_{l+1}}...R_{i_{1}}\psi_{\alpha}|dud\mu_{\slashed{g}}
\end{align}
By virtue of the bound
\begin{align}
 \max_{\alpha}\sup_{\Sigma_{t}^{\epsilon_{0}}}|R_{i_{l+1}}T\psi_{\alpha}|\leq C\delta_{0}(1+t)^{-1}
\end{align}
we have:
\begin{align}
 |H_{1}|\leq C\delta_{0}\int_{\Sigma_{t}^{\epsilon_{0}}}(1+t)|R_{i_{l}}...R_{i_{1}}
\textrm{tr}\chi'||R_{i_{l+1}}...R_{i_{1}}\psi_{\alpha}|dud\mu_{\slashed{g}}\\\notag
C\delta_{0}(1+t)\|R_{i_{l}}...R_{i_{1}}\textrm{tr}\chi'\|_{L^{2}(\Sigma_{t}^{\epsilon_{0}})}
\|R_{i_{l+1}}...R_{i_{1}}\psi_{\alpha}\|_{L^{2}(\Sigma_{t}^{\epsilon_{0}})}\\\notag
\leq C\delta_{0}(1+t)^{2}\|(R_{i_{l}}...R_{i_{1}}\textrm{tr}\chi')(t)\|_{L^{2}([0,\epsilon_{0}]\times S^{2})}
\mathcal{W}_{\{l+1\}}(t)
\end{align}
Now by (14.243), (14.245):
\begin{align}
 (1+t)^{2}\|(R_{i_{l}}...R_{i_{1}}\textrm{tr}\chi')(t)\|_{L^{2}([0,\epsilon_{0}]\times S^{2})}\\\notag
\leq C\{\|R_{i_{l}}...R_{i_{1}}\textrm{tr}\chi'\|_{L^{2}(\Sigma_{0}^{\epsilon_{0}})}+\sum_{k=0}^{2}
\|\leftexp{(i_{1}...i_{l})}{A}^{(k)}_{l}(t)\|_{L^{2}([0,\epsilon_{0}]\times S^{2})}\}\\\notag
\leq C\{\|R_{i_{l}}...R_{i_{1}}\textrm{tr}\chi'\|_{L^{2}(\Sigma_{0}^{\epsilon_{0}})}+\sum_{k=0}^{2}
\int_{0}^{t}(1+t')\|\leftexp{(i_{1}...i_{l})}{\rho}^{(k)}_{l}\|_{L^{2}(\Sigma_{t'}^{\epsilon_{0}})}dt'\}
\end{align}
Here we need only consider the leading principal contribution, namely that of 
$\|\leftexp{(i_{1}...i_{l})}{\rho}^{(0)}_{l}\|_{L^{2}(\Sigma_{t'}^{\epsilon_{0}})}$.
This contribution to (14.315) is bounded by:
\begin{align}
 C|\ell|\int_{0}^{t}(1+t')^{-1}\bar{\mu}^{-1/2}_{m}(t')\sqrt{\sum_{j,\alpha}\mathcal{E}'_{1}[R_{j}R_{i_{l}}...R_{i_{1}}
\psi_{\alpha}](t')}dt'\\\notag
\leq C|\ell|\sqrt{\mathcal{G}'_{1,[l+2];a,q}(t)}\int_{0}^{t}(1+t')^{-1}\bar{\mu}^{-a-1/2}_{m}(t')[1+\log(1+t')]^{q}dt'
\end{align}
The last integral coincides with $J_{a,q-1}$ (see (14.147)), therefore by (14.153) it is bounded by:
\begin{align*}
 C[1+\log(1+t)]^{q+1}\bar{\mu}^{-a+1/2}_{m}(t)
\end{align*}
Consequently the leading principal contribution to (14.315) is bounded by:
\begin{align}
 C|\ell|\bar{\mu}^{-a+1/2}_{m}(t)[1+\log(1+t)]^{q+1}\sqrt{\mathcal{G}'_{1,[l+2];a,q}(t)}
\end{align}
On the other hand, as in (14.308) we have:
\begin{align}
 \mathcal{W}_{\{l+1\}}\leq C\epsilon_{0}\sqrt{\mathcal{E}_{0,[l+2]}(t)}\\\notag
\leq C\epsilon_{0}\bar{\mu}^{-a}_{m}(t)[1+\log(1+t)]^{p}\sqrt{\mathcal{G}_{0,[l+2];a,p}(t)}
\end{align}
Substituting (14.317) and (14.318) in (14.314), we conclude that the leading principal contribution to 
$H_{1}$ is bounded by:
\begin{align}
 C|\ell|\epsilon_{0}\delta_{0}\bar{\mu}^{-2a+1/2}_{m}(t)[1+\log(1+t)]^{p+q+1}
\sqrt{\mathcal{G}_{0,[l+2];a,p}(t)\mathcal{G}'_{1,[l+2];a,q}(t)}
\end{align}
the same in the form as (14.311).

     We turn to the spacetime integral in (14.207) namely:
\begin{align}
 \int_{W^{t}_{u}}(R_{i_{l+1}}...R_{i_{1}}\psi_{\alpha})(L+\nu)\{(\omega/\nu)(T\psi_{\alpha})
(R_{i_{l+1}}...R_{i_{1}}\textrm{tr}\chi')\}dt'du'd\mu_{\tilde{\slashed{g}}}
\end{align}
Here we shall use the fact that $(L+\nu)T\psi_{\alpha}$ decays faster than $(1+t)^{-2}$. 
To establish this fact we consider the wave equation satisfied by $\psi_{\alpha}$.
In analogy with (8.180) and (8.181) we have for $\alpha=0,1,2,3$:
\begin{equation}
 (L+\nu)\underline{L}\psi_{\alpha}=\rho_{\alpha}
\end{equation}
where:
\begin{align}
 \rho_{\alpha}=\mu\slashed{\Delta}\psi_{\alpha}-\underline{\nu}L\psi_{\alpha}-2\zeta\cdot\slashed{d}\psi_{\alpha}+\mu
\frac{d\log\Omega}{dh}\slashed{d}h\cdot\slashed{d}\psi_{\alpha}
\end{align}
Under the same assumptions as those of Lemma 8.10 an estimate similar to 
(8.182) holds for each $\alpha=0,1,2,3$, that is, we have:
\begin{align}
 \max_{\alpha}|\rho_{\alpha}|\leq C\delta_{0}(1+t)^{-3}[1+\log(1+t)]
\end{align}
Since $2T\psi_{\alpha}=\underline{L}\psi_{\alpha}-\eta^{-1}\kappa L\psi_{\alpha}$, (14.321) implies:
\begin{align}
 (L+\nu)T\psi_{\alpha}=\tau_{\alpha}
\end{align}
where:
\begin{align}
 2\tau_{\alpha}=\rho_{\alpha}-(L+\nu)(\eta^{-1}\kappa L\psi_{\alpha})
\end{align}
Now we have:
\begin{align*}
 (L+\nu)(\eta^{-1}\kappa L\psi_{\alpha})=\eta^{-1}\kappa(L)^{2}\psi_{\alpha}+((L+\nu)(\eta^{-1}\kappa))L\psi_{\alpha}
\end{align*}
Using the bound for $L\mu$ and the assumptions $\textbf{E}^{QQ}_{\{0\}}$, $\textbf{E}^{Q}_{\{0\}}$ we deduce:
\begin{align}
 \max_{\alpha}|(L+\nu)(\eta^{-1}\kappa L\psi_{\alpha})|\leq C\delta_{0}(1+t)^{-3}[1+\log(1+t)]
\end{align}
Combining (14.323) and (14.326) we obtain:
\begin{align}
 \max_{\alpha}|\tau_{\alpha}|\leq C\delta_{0}(1+t)^{-3}[1+\log(1+t)]
\end{align}
This estimate actually relies only on the assumptions of Proposition 12.6. By direct calculation, 
under the assumptions of Proposition 12.9 with $l=1$ together with those of Proposition 12.10 with $m=l=0$, we deduce:
\begin{align}
 \max_{\alpha}\|\tau_{\alpha}\|_{\infty,[1],\Sigma_{t}^{\epsilon_{0}}}\leq C\delta_{0}(1+t)^{-3}[1+\log(1+t)]
\end{align}
    Next, we consider the factor $\omega/\nu$. Setting:
\begin{align}
 \nu'=\frac{1}{2}(\textrm{tr}\chi'+L\log\Omega)
\end{align}
we have:
\begin{align}
 \nu=\frac{1}{1-u+t}+\nu'
\end{align}
Recalling that $\omega=2(1+t)$ we then obtain:
\begin{align}
 (L-2\nu)(\omega/\nu)=\gamma\nu^{-2}\omega
\end{align}
where:
\begin{align}
 \gamma=\frac{-u}{(1+t)(1-u+t)^{2}}+(\frac{1}{1+t}-\frac{4}{1-u+t})\nu'-2\nu'^{2}-L\nu'
\end{align}
Proposition 12.9 with $l=1$ implies:
\begin{align}
 \|\gamma\|_{\infty,[1],\Sigma_{t}^{\epsilon_{0}}}\leq C\delta_{0}(1+t)^{-3}[1+\log(1+t)]
\end{align}
Here we also have used (14.217) and (14.219).

     By (14.324) and (14.331) we have,
\begin{align}
 (L+\nu)\{(\omega/\nu)(T\psi_{\alpha})(R_{i_{l+1}}...R_{i_{1}}\textrm{tr}\chi')\}\\\notag
=(\omega/\nu)\{(T\psi_{\alpha})(L+2\nu)(R_{i_{l+1}}...R_{i_{1}}\textrm{tr}\chi')+\tilde{\tau}
_{\alpha}(R_{i_{l+1}}...R_{i_{1}}\textrm{tr}\chi')\}
\end{align}
where:
\begin{align}
 \tilde{\tau}_{\alpha}=\tau_{\alpha}+\gamma\nu^{-1}T\psi_{\alpha}
\end{align}
By (14.328) and (14.333) and the estimate for $\nu^{\prime}$ resulting from Proposition 12.9 with $l=1$ we have:
\begin{align}
 \max_{\alpha}\|\tilde{\tau}_{\alpha}\|_{\infty,[1],\Sigma_{t}^{\epsilon_{0}}}\leq C\delta_{0}(1+t)^{-3}[1+\log(1+t)]
\end{align}
Substituting (14.334) in (14.320), we are thus to estimate the spacetime integral:
\begin{align}
 \int_{W^{t}_{u}}(R_{i_{l+1}}...R_{i_{1}}\psi_{\alpha})(\omega/\nu)\\\notag
\cdot \{(T\psi_{\alpha})(L+2\nu)(R_{i_{l+1}}...R_{i_{1}}\textrm{tr}\chi')
+\tilde{\tau}_{\alpha}(R_{i_{l+1}}...R_{i_{1}}\textrm{tr}\chi')\}dt'du'd\mu_{\tilde{\slashed{g}}}
\end{align}
Writing:
\begin{align}
 (L+2\nu)R_{i_{l+1}}...R_{i_{1}}\textrm{tr}\chi'=R_{i_{l+1}}(L+2\nu)R_{i_{l}}...R_{i_{1}}\textrm{tr}\chi'\\\notag
+\leftexp{(R_{i_{l+1}})}{Z}R_{i_{l}}...R_{i_{1}}\textrm{tr}\chi'-2(R_{i_{l+1}}\nu)(R_{i_{l}}...R_{i_{1}}\textrm{tr}\chi')
\end{align}
we integrate by parts on each $S_{t,u}$ using (14.209), first taking $X=R_{i_{l+1}}$ to obtain:
\begin{align}
 \int_{S_{t,u}}(\omega/\nu)(R_{i_{l+1}}...R_{i_{1}}\psi_{\alpha})(T\psi_{\alpha})R_{i_{l+1}}((L+2\nu)
R_{i_{l}}...R_{i_{1}}\textrm{tr}\chi')d\mu_{\tilde{\slashed{g}}}\\\notag
=-\int_{S_{t,u}}(\omega/\nu)(T\psi_{\alpha})(R_{i_{l+1}}R_{i_{l+1}}...R_{i_{1}}\psi_{\alpha})(L+2\nu)
R_{i_{l}}...R_{i_{1}}\textrm{tr}\chi'd\mu_{\tilde{\slashed{g}}}\\\notag
-\int_{S_{t,u}}(\omega/\nu)(R_{i_{l+1}}T\psi_{\alpha})(R_{i_{l+1}}...R_{i_{1}}\psi_{\alpha})(L+2\nu)
R_{i_{l}}...R_{i_{1}}\textrm{tr}\chi'd\mu_{\tilde{\slashed{g}}}\\\notag
-\int_{S_{t,u}}(T\psi_{\alpha})\{R_{i_{l+1}}(\omega/\nu)+\frac{1}{2}(\omega/\nu)\textrm{tr}
\leftexp{(R_{i_{l+1}})}{\tilde{\slashed{\pi}}}\}\\\notag
\cdot (R_{i_{l+1}}...R_{i_{1}}\psi_{\alpha})(L+2\nu)R_{i_{l}}...R_{i_{1}}\textrm{tr}\chi'd\mu_{\tilde{\slashed{g}}}
\end{align}
and then taking $X=\leftexp{(R_{i_{l+1}})}{Z}$ to obtain:
\begin{align}
 \int_{S_{t,u}}(\omega/\nu)(T\psi_{\alpha})(R_{i_{l+1}}...R_{i_{1}}\psi_{\alpha})\leftexp{(R_{i_{l+1}})}{Z}
(R_{i_{l}}...R_{i_{1}}\textrm{tr}\chi')d\mu_{\tilde{\slashed{g}}}\\\notag
=-\int_{S_{t,u}}(\leftexp{(R_{i_{l+1}})}{Z}\cdot\slashed{d}R_{i_{l+1}}...R_{i_{1}}\psi_{\alpha})(\omega/\nu)
(T\psi_{\alpha})(R_{i_{l}}...R_{i_{1}}\textrm{tr}\chi')d\mu_{\tilde{\slashed{g}}}\\\notag
-\int_{S_{t,u}}(\omega/\nu)(\leftexp{(R_{i_{l+1}})}{Z}\cdot\slashed{d}T\psi_{\alpha})(R_{i_{l+1}}...R_{i_{1}}
\psi_{\alpha})(R_{i_{l}}...R_{i_{1}}\textrm{tr}\chi')d\mu_{\tilde{\slashed{g}}}\\\notag
-\int_{S_{t,u}}(T\psi_{\alpha})\{\leftexp{(R_{i_{l+1}})}{Z}\cdot\slashed{d}(\omega/\nu)+(\omega/\nu)
\tilde{\slashed{\textrm{div}}}\leftexp{(R_{i_{l+1}})}{Z}\}\\\notag
\cdot (R_{i_{l+1}}...R_{i_{1}}\psi_{\alpha})(R_{i_{l}}...R_{i_{1}}\textrm{tr}\chi')d\mu_{\tilde{\slashed{g}}}
\end{align}
Also, taking again $X=R_{i_{l+1}}$,
\begin{align}
 \int_{S_{t,u}}(R_{i_{l+1}}...R_{i_{1}}\psi_{\alpha})(\omega/\nu)\tilde{\tau}_{\alpha}R_{i_{l+1}}...R_{i_{1}}
\textrm{tr}\chi'd\mu_{\tilde{\slashed{g}}}\\\notag
=-\int_{S_{t,u}}(\omega/\nu)\tilde{\tau}_{\alpha}(R_{i_{l+1}}R_{i_{l+1}}...R_{i_{1}}\psi_{\alpha})
(R_{i_{l}}...R_{i_{1}}\textrm{tr}\chi')d\mu_{\tilde{\slashed{g}}}\\\notag
-\int_{S_{t,u}}(\omega/\nu)(R_{i_{l+1}}\tilde{\tau}_{\alpha})(R_{i_{l+1}}...R_{i_{1}}
\psi_{\alpha})(R_{i_{l}}...R_{i_{1}}\textrm{tr}\chi')d\mu_{\tilde{\slashed{g}}}\\\notag
-\int_{S_{t,u}}\tilde{\tau}_{\alpha}(R_{i_{l+1}}...R_{i_{1}}\psi_{\alpha})\{R_{i_{l+1}}
(\omega/\nu)+\frac{1}{2}(\omega/\nu)\textrm{tr}\leftexp{(R_{i_{l+1}})}
{\tilde{\slashed{\pi}}}\}(R_{i_{l}}...R_{i_{1}}\textrm{tr}\chi')d\mu_{\tilde{\slashed{g}}}
\end{align}
In view of (14.339)-(14.341) the spacetime integral (14.337) becomes:
\begin{equation}
 -V_{0}-V_{1}-V_{2}-V_{3}
\end{equation}
where:
\begin{align}
 V_{0}=V_{0,0}+V_{0,1}\\
V_{0,0}=\int_{W^{t}_{u}}(\omega/\nu)(T\psi_{\alpha})(R_{i_{l+1}}R_{i_{l+1}}...R_{i_{1}}\psi_{\alpha})
(L+2\nu)R_{i_{l}}...R_{i_{1}}\textrm{tr}\chi'dt'du'd\mu_{\tilde{\slashed{g}}}\\
V_{0,1}=\int_{W^{t}_{u}}(\omega/\nu)\{\tilde{\tau}_{\alpha}(R_{i_{l+1}}R_{i_{l+1}}...R_{i_{1}}\psi_{\alpha})
+(T\psi_{\alpha})\leftexp{(R_{i_{l+1}})}{Z}\cdot\slashed{d}R_{i_{l+1}}...R_{i_{1}}\psi_{\alpha}\}
R_{i_{l}}...R_{i_{1}}\textrm{tr}\chi'dt'du'd\mu_{\tilde{\slashed{g}}}\\
V_{1}=V_{1,0}+V_{1,1}\\
V_{1,0}=\int_{W^{t}_{u}}(\omega/\nu)(R_{i_{l+1}}T\psi_{\alpha})(R_{i_{l+1}}...R_{i_{1}}\psi_{\alpha})
(L+2\nu)R_{i_{l}}...R_{i_{1}}\textrm{tr}\chi'dt'du'd\mu_{\tilde{\slashed{g}}}\\
V_{1,1}=\int_{W^{t}_{u}}(\omega/\nu)\{R_{i_{l+1}}\tilde{\tau}_{\alpha}+\leftexp{(R_{i_{l+1}})}{Z}\cdot\slashed{d}
T\psi_{\alpha}\}\cdot(R_{i_{l+1}}...R_{i_{1}}\psi_{\alpha})R_{i_{l}}...R_{i_{1}}\textrm{tr}\chi'
dt'du'd\mu_{\tilde{\slashed{g}}}
\end{align}

\begin{align}
 V_{2}=V_{2,0}+V_{2,1}\\
V_{2,0}=\int_{W^{t}_{u}}(T\psi_{\alpha})\{R_{i_{l+1}}(\omega/\nu)+\frac{1}{2}(\omega/\nu)\textrm{tr}
\leftexp{(R_{i_{l+1}})}{\tilde{\slashed{\pi}}}\}
(R_{i_{l+1}}...R_{i_{1}}\psi_{\alpha})(L+2\nu)R_{i_{l}}...R_{i_{1}}\textrm{tr}\chi'dt'du'd\mu_{\tilde{\slashed{g}}}\\
V_{2,1}=\int_{W^{t}_{u}}(R_{i_{l+1}}...R_{i_{1}}\psi_{\alpha})(R_{i_{l}}...R_{i_{1}}\textrm{tr}\chi')\cdot\\\notag
\{\tilde{\tau}_{\alpha}[R_{i_{l+1}}(\omega/\nu)+\frac{1}{2}(\omega/\nu)\textrm{tr}\leftexp{(R_{i_{l+1}})}
{\tilde{\slashed{\pi}}}]+(T\psi_{\alpha})[\leftexp{(R_{i_{l+1}})}{Z}\cdot\slashed{d}(\omega/\nu)+(\omega/\nu)
\tilde{\slashed{\textrm{div}}}\leftexp{(R_{i_{l+1}})}{Z}]\}dt'du'd\mu_{\tilde{\slashed{g}}}
\end{align}
and:
\begin{align}
 V_{3}=\int_{W^{t}_{u}}2(\omega/\nu)(R_{i_{l+1}}...R_{i_{1}}\psi_{\alpha})(T\psi_{\alpha})
(R_{i_{l+1}}\nu)(R_{i_{l}}...R_{i_{1}}\textrm{tr}\chi')dt'du'd\mu_{\tilde{\slashed{g}}}
\end{align}
The last term is a lower order integral which can be easily estimated.

In the following we shall focus attention on the two leading integrals $V_{0,0}$ and $V_{1,0}$, the other integrals 
containing decay factors compared to these two leading integrals.
In fact, comparing $V_{0,0}$ and $V_{0,1}$ as well as $V_{1,0}$ and $V_{1,1}$ by using (14.336) and the estimate
\begin{align}
 \max_{j}\|\leftexp{(R_{j})}{Z}\|_{\infty,[1],\Sigma_{t}^{\epsilon_{0}}}\leq C\delta_{0}(1+t)^{-1}[1+\log(1+t)]
\end{align}
 of Corollary 10.1.i, we see that $V_{0,1}$ and $V_{1,1}$ contain a decay factor of 
$\delta_{0}(1+t)^{-1}[1+\log(1+t)]$ as compared to $V_{0,0}$ and $V_{1,0}$ respectively.

     We first consider $V_{0,0}$. We have:
\begin{align}
 |V_{0,0}|\leq C\int_{W^{t}_{\epsilon_{0}}}(1+t')^{3}|T\psi_{\alpha}||\slashed{d}R_{i_{l+1}}...R_{i_{1}}\psi_{\alpha}|
|(L+2\nu)R_{i_{l}}...R_{i_{1}}\textrm{tr}\chi'|dt'dud\mu_{\slashed{g}}
\end{align}
 From the propagation equation (14.219) we have:
\begin{align}
 (L+2\nu)R_{i_{l}}...R_{i_{1}}\textrm{tr}\chi'=\leftexp{(i_{1}...i_{l})}{\tilde{\rho}}_{l}
\end{align}
where:
\begin{align}
 \leftexp{(i_{1}...i_{l})}{\tilde{\rho}}_{l}=\leftexp{(i_{1}...i_{l})}{\rho}_{l}+2\nu'R_{i_{l}}...R_{i_{1}}
\textrm{tr}\chi'
\end{align}
Now we can see that $V_{0,1}$ and $V_{1,1}$ enjoy the same bounds as the contribution of the 
second term on the right of above to $V_{0,1}$ and $V_{1,0}$.

We write, as in (14.229),
\begin{align}
 \leftexp{(i_{1}...i_{l})}{\tilde{\rho}}_{l}=\leftexp{(i_{1}...i_{l})}{\rho}^{(0)}_{l}+
\leftexp{(i_{1}...i_{l})}{\rho}^{(1)}_{l}+\leftexp{(i_{1}...i_{l})}{\tilde{\rho}}^{(2)}_{l}
\end{align}
where:
\begin{align}
 \leftexp{(i_{1}...i_{l})}{\tilde{\rho}}^{(2)}_{l}=\leftexp{(i_{1}...i_{l})}{\rho}^{(2)}_{l}+2\nu'R_{i_{l}}...R_{i_{1}}\textrm{tr}\chi'
\end{align}
Since
\begin{align*}
 |\nu'|\leq C\delta_{0}(1+t)^{-2}[1+\log(1+t)]
\end{align*}
$\leftexp{(i_{1}...i_{l})}{\tilde{\rho}}^{(2)}_{l}$ enjoys the same bound as 
$\leftexp{(i_{1}...i_{l})}{\rho}^{(2)}_{l}$.
The $borderline$ contribution to (14.354) is the contribution from $\leftexp{(i_{1}...i_{l})}{\rho}_{l}^{(0)}$. 
This is the $borderline$ integral:
\begin{align}
 C\int_{W^{t}_{u}}(1+t')^{3}|T\psi_{\alpha}||\slashed{d}R_{i_{l+1}}...R_{i_{1}}\psi_{\alpha}|
|\leftexp{(i_{1}...i_{l})}{\rho}^{(0)}_{l}|dt'dud\mu_{\slashed{g}}\\\notag
\leq C\int_{0}^{t}(1+t')^{3}\sup_{\Sigma_{t'}^{\epsilon_{0}}}(\mu^{-1}|T\psi_{\alpha}|)\\\notag
\cdot \|\mu^{1/2}\slashed{d}R_{i_{l+1}}...R_{i_{1}}\psi_{\alpha}\|_{L^{2}(\Sigma_{t'}^{\epsilon_{0}})}
\|\mu^{1/2}\leftexp{(i_{1}...i_{l})}{\rho}^{(0)}_{l}\|
_{L^{2}(\Sigma_{t'}^{\epsilon_{0}})}dt'
\end{align}
Here we substitute (14.89) for $\sup_{\Sigma_{t'}^{\epsilon_{0}}}(\mu^{-1}|T\psi_{\alpha}|)$ and the estimate
(14.240) for $\|\mu^{1/2}\leftexp{(i_{1}...i_{l})}{\rho}^{(0)}_{l}\|_{L^{2}(\Sigma_{t'}^{\epsilon_{0}})}$. 
The factors $|\ell|$ then cancel. Now the partial contribution of the second term 
on the right of (14.89) is actually not borderline. We shall show how to estimate contribution of this type afterwards, 
in connection with the estimate for the contribution of $\leftexp{(i_{1}...i_{l})}{\rho}^{(1)}_{l}$. 
For the present we focus on the $borderline$ integral:
\begin{align}
 C\int_{0}^{t}\sup_{\Sigma_{t'}^{\epsilon_{0}}}(\mu^{-1}|L\mu|)\sqrt{\mathcal{E}'_{1}
[R_{i_{l+1}}...R_{i_{1}}\psi_{\alpha}](t')}\sqrt{\sum_{j,\alpha}\mathcal{E}'_{1}[R_{j}R_{i_{l}}...R_{i_{1}}
\psi_{\alpha}](t')}dt'\\\notag
\leq C\int_{0}^{t}\sup_{\Sigma_{t'}^{\epsilon_{0}}}(\mu^{-1}|L\mu|)\sum_{j,\alpha}
\mathcal{E}'_{1}[R_{j}R_{i_{l}}...R_{i_{1}}\psi_{\alpha}](t')dt'\\\notag
\leq C\int_{0}^{t}\sup_{\Sigma_{t'}^{\epsilon_{0}}}(\mu^{-1}|L\mu|)\bar{\mu}^{-2a}_{m}(t')[1+\log(1+t')]^{2q}
\leftexp{(i_{1}...i_{l})}{\mathcal{G}}^{\prime}_{1,l+2;a,q}(t')dt'
\end{align}
This has the same form as (14.97), with $q$ replacing $p$ and 
$\leftexp{(i_{1}...i_{l})}{\mathcal{G}}^{\prime}_{1,l+2;a,q}$ replacing 
$\leftexp{(i_{1}...i_{l})}{\mathcal{G}}_{0,l+2;a,p}$. Thus, from (14.102) and (14.107) we conclude that 
(14.360) is bounded by:
\begin{align}
 C(\frac{1}{2a}+\frac{1}{2q})\bar{\mu}^{-2a}_{m}(t)[1+\log(1+t)]^{2q}\leftexp{(i_{1}...i_{l})}
{\mathcal{G}}^{\prime}_{1,l+2;a,q}(t)
\end{align}

    We proceed to consider the contribution of $\leftexp{(i_{1}...i_{l})}{\rho}^{(1)}_{l}$ to (14.354). To estimate 
this and all remaining contributions we simply use (14.90). Using (14.241) we obtain that the contribution in question 
is bounded by:
\begin{align}
 C\int_{W^{t}_{\epsilon_{0}}}(1+t')^{3}|T\psi_{\alpha}||\slashed{d}R_{i_{l+1}}...R_{i_{1}}
\psi_{\alpha}||\leftexp{(i_{1}...i_{l})}{\rho}^{(1)}_{l}|dt'dud\mu_{\slashed{g}}\\\notag
\leq C\delta_{0}\int_{0}^{t}(1+t')^{2}\bar{\mu}_{m}^{-1}(t')\|\mu^{1/2}\slashed{d}R_{i_{l+1}}...R_{i_{1}}
\psi_{\alpha}\|_{L^{2}(\Sigma_{t'}^{\epsilon_{0}})}
\|\mu^{1/2}\leftexp{(i_{1}...i_{l})}{\rho}^{(1)}_{l}\|_{L^{2}(\Sigma_{t'}^{\epsilon_{0}})}dt'\\\notag
\leq C\delta_{0}^{2}\int_{0}^{t}(1+t')^{-2}\bar{\mu}^{-1}_{m}(t')\mathcal{E}'_{1,[l+2]}(t')dt'\\\notag
\leq C\delta_{0}^{2}\int_{0}^{t}(1+t')^{-2}[1+\log(1+t')]^{2q}\bar{\mu}^{-1-2a}_{m}(t')\mathcal{G}'_{1,[l+2];a,q}(t')dt'
\end{align}
The last integral is similar as (14.115). We find (14.362) is bounded by:
\begin{align}
 C\delta_{0}^{2}\bar{\mu}_{m}^{-2a}(t)[1+\log(1+t)]^{2q}\\\notag
\cdot \{\varphi_{1}(Ca\delta_{0})\mathcal{G}'_{1,[l+2];a,q}(t)+\int_{0}^{t}(1+t')^{-2}
\mathcal{G}'_{1,[l+2];a,q}(t')dt'\}
\end{align}   

      We consider next the contribution to (14.354) of the leading term $C_{l}(1+t)^{-2}\mathcal{W}^{Q}_{\{l\}}$ in 
the estimate for $\|\leftexp{(i_{1}...i_{l})}{\tilde{\rho}}^{(2)}_{l}\|_{L^{2}(\Sigma_{t}^{\epsilon_{0}})}$. 
This contribution is bounded by:
\begin{align}
 C_{l}\delta_{0}\int_{0}^{t}\bar{\mu}^{-1/2}_{m}(t')\|\mu^{1/2}\slashed{d}R_{i_{l+1}}...R_{i_{1}}
\psi_{\alpha}\|_{L^{2}(\Sigma_{t'}^{\epsilon_{0}})}\mathcal{W}^{Q}_{\{l\}}(t')dt'\\\notag
\leq C_{l}\epsilon_{0}\delta_{0}\int_{0}^{t}(1+t')^{-1}\bar{\mu}^{-2a-1/2}_{m}(t')[1+\log(1+t')]^{p+q}\\\notag
\times\sqrt{\mathcal{G}'_{1,[l+2];a,q}(t')\mathcal{G}_{0,[l+2];a,p}(t')}dt'\\\notag
\leq C_{l}\epsilon_{0}\delta_{0}J_{2a,p+q-1}(t)\sqrt{\mathcal{G}'_{1,[l+2];a,q}(t)\mathcal{G}_{0,[l+2];a,p}(t)}
\end{align}
where $J_{a,q}(t)$ is the integral (14.147). By (14.153) we have:
\begin{align}
 J_{2a,p+q-1}(t)\leq C[1+\log(1+t)]^{p+q+1}\bar{\mu}^{-2a+1/2}_{m}(t)
\end{align}
hence (14.364) is bounded by:
\begin{align}
 C_{l}\epsilon_{0}\delta_{0}\bar{\mu}^{-2a+1/2}_{m}(t)[1+\log(1+t)]^{p+q+1}\sqrt{\mathcal{G}'_{1,[l+2];a,q}(t)
\mathcal{G}_{0,[l+2];a,p}(t)}
\end{align}

     Finally, we consider the integral $V_{1,0}$, given by (14.347). By the bound (14.213) we have:
\begin{align}
 |V_{1,0}|\leq C\delta_{0}\int_{W^{t}_{u}}(1+t')|R_{i_{l+1}}...R_{i_{1}}\psi_{\alpha}|
|\leftexp{(i_{1}...i_{l})}{\tilde{\rho}}_{l}|dt'dud\mu_{\slashed{g}}\\\notag
\leq C\delta_{0}\int_{0}^{t}(1+t')\mathcal{W}_{\{l+1\}}(t')\|\leftexp{(i_{1}...i_{l})}
{\tilde{\rho}}_{l}\|_{L^{2}(\Sigma_{t'}^{\epsilon_{0}})}dt'
\end{align}
Here we need only consider the leading principal contribution, namely that of 
$\|\leftexp{(i_{1}...i_{l})}{\tilde{\rho}}^{(0)}_{l}\|_{L^{2}(\Sigma_{t'}^{\epsilon_{0}})}$.
By (14.240) and (14.318) this contribution is bounded by:
\begin{align}
 C\epsilon_{0}\delta_{0}\int_{0}^{t}(1+t')^{-1}\bar{\mu}^{-2a-1/2}_{m}(t')[1+\log(1+t')]^{p+q}
\sqrt{\mathcal{G}_{0,[l+2];a,p}(t')\mathcal{G}'_{1,[l+2];a,q}(t')}dt'
\end{align}
This is identical to (14.364), which has been bounded by (14.366). This completes the estimates for the spacetime 
integrals $V_{0}, V_{1}, V_{2}$. Thus, the estimate of the contribution of (14.56) to (14.59) is completed.

\subsection{Estimates for the Contribution of (14.57)}
     We now consider the contribution of (14.57) to (14.59). Recalling that this is associated to the variation (14.54),
 the contribution in question is:
\begin{align}
 -\int_{W^{t}_{u}}(\omega/\nu)(T\psi_{\alpha})(R_{i_{l-m}}...R_{i_{1}}(T)^{m}
\slashed{\Delta}\mu)(L+\nu)(R_{i_{l-m}}...R_{i_{1}}(T)^{m+1}\psi_{\alpha})dt'du'd\mu_{\tilde{\slashed{g}}}
\end{align}
We deal with integral in the same way as we dealt with (14.202). That is, we write the integrand in the form:
\begin{align*}
 -(\omega/\nu)(T\psi_{\alpha})(R_{i_{l-m}}...R_{i_{1}}(T)^{m}\slashed{\Delta}\mu)
((L+\nu)R_{i_{l-m}}...R_{i_{1}}(T)^{m+1}\psi_{\alpha})=\\
-(L+2\nu)\{(\omega/\nu)(T\psi_{\alpha})(R_{i_{l-m}}...R_{i_{1}}(T)^{m}\slashed{\Delta}\mu)
(R_{i_{l-m}}...R_{i_{1}}(T)^{m+1}\psi_{\alpha})\}\\
+(R_{i_{l-m}}...R_{i_{1}}(T)^{m+1}\psi_{\alpha})(L+\nu)\{(\omega/\nu)(T\psi_{\alpha})
(R_{i_{l-m}}...R_{i_{1}}(T)^{m}\slashed{\Delta}\mu)\}
\end{align*}
By (14.206) with the function
\begin{align*}
 (\omega/\nu)(T\psi_{\alpha})(R_{i_{l-m}}...R_{i_{1}}(T)^{m}\slashed{\Delta}\mu)
(R_{i_{l-m}}...R_{i_{1}}(T)^{m+1}\psi_{\alpha})
\end{align*}
in the role of the function $f$, we conclude that (14.369) equals:
\begin{align}
 -\int_{\Sigma_{t}^{u}}(\omega/\nu)(T\psi_{\alpha})(R_{i_{l-m}}...R_{i_{1}}(T)^{m}\slashed{\Delta}\mu)
(R_{i_{l-m}}...R_{i_{1}}(T)^{m+1}\psi_{\alpha})du'd\mu_{\tilde{\slashed{g}}}\\\notag
+\int_{\Sigma_{0}^{u}}(\omega/\nu)(T\psi_{\alpha})(R_{i_{l-m}}...R_{i_{1}}(T)^{m}\slashed{\Delta}\mu)
(R_{i_{l-m}}...R_{i_{1}}(T)^{m+1}\psi_{\alpha})du'd\mu_{\tilde{\slashed{g}}}\\\notag
+\int_{W^{t}_{u}}(R_{i_{l-m}}...R_{i_{1}}(T)^{m+1}\psi_{\alpha})(L+\nu)\\\notag
\cdot \{(\omega/\nu)(T\psi_{\alpha})(R_{i_{l-m}}...R_{i_{1}}(T)^{m}\slashed{\Delta}\mu)\}
dt'du'd\mu_{\tilde{\slashed{g}}}	
\end{align}
    We first consider the hypersurface integral:
\begin{align}
 -\int_{\Sigma_{t}^{u}}(\omega/\nu)(T\psi_{\alpha})(R_{i_{l-m}}...R_{i_{1}}(T)^{m}\slashed{\Delta}\mu)(R_{i_{l-m}}...R_{i_{1}}(T)^{m+1}\psi_{\alpha})
du'd\mu_{\tilde{\slashed{g}}}
\end{align}
We define the $S_{t,u}$-tangential vectorfields:
\begin{align}
 \leftexp{(i_{1}...i_{l-m})}{Y}_{m,l-m}=(\slashed{d}R_{i_{l-m}}...R_{i_{1}}(T)^{m}\mu)\cdot(\slashed{g}^{-1})
\end{align}
Then with
\begin{align}
 \leftexp{(i_{1}...i_{l-m})}{r}_{m,l-m}=R_{i_{l-m}}...R_{i_{1}}(T)^{m}\slashed{\Delta}\mu
-\slashed{\Delta}R_{i_{l-m}}...R_{i_{1}}(T)^{m}\mu
\end{align}
in view of the fact that, since $\tilde{\slashed{g}}=\Omega\slashed{g}$, for any $S_{t,u}$-tangential vectorfield $X$:
\begin{align}
 \tilde{\slashed{\textrm{div}}}X=\slashed{\textrm{div}}X+X\cdot\slashed{d}\log\Omega,
\end{align}
we have:
\begin{align}
 R_{i_{l-m}}...R_{i_{1}}(T)^{m}\slashed{\Delta}\mu=\tilde{\slashed{\textrm{div}}}\leftexp{(i_{1}...i_{l-m})}
{Y}_{m,l-m}+\leftexp{(i_{1}...i_{l-m})}{\tilde{r}}_{m,l-m}
\end{align}
where:
\begin{align}
 \leftexp{(i_{1}...i_{l-m})}{\tilde{r}}_{m,l-m}=\leftexp{(i_{1}...i_{l-m})}{r}_{m,l-m}
-\leftexp{(i_{1}...i_{l-m})}{Y}_{m,l-m}\cdot\slashed{d}\log\Omega
\end{align}
In analogy with (11.323)-(11.325) we deduce:
\begin{align}
 \leftexp{(i_{1}...i_{l-m})}{r}_{m,l-m}=\textrm{tr}(\leftexp{(i_{1}...i_{l-m})}{c}_{m,l-m}
[\slashed{d}\mu])\cdot(\slashed{g}^{-1})\\\notag
+\textrm{tr}\{\slashed{\mathcal{L}}_{R_{i_{l-m}}}...\slashed{\mathcal{L}}_{R_{i_{1}}}\sum_{k=1}^{m}\frac{m!}{k!(m-k)!}
((\slashed{\mathcal{L}}_{T})^{k}(\slashed{g}^{-1}))\cdot(\slashed{\mathcal{L}}_{T})^{m-k}\slashed{D}^{2}\mu\}\\\notag
+\textrm{tr}\{\sum_{|s_{1}|+|s_{2}|=l-m, |s_{1}|>0}((\slashed{\mathcal{L}}_{R})^{s_{1}}(\slashed{g}^{-1}))
\cdot(\slashed{\mathcal{L}}_{R})^{s_{2}}(\slashed{\mathcal{L}}_{T})^{m}\slashed{D}^{2}\mu\}
\end{align}
Using Corollaries 10.1.d, 10.2.d, 11.1.c, 11.2.c, 11.2.g, we obtain:
\begin{align}
 \|\leftexp{(i_{1}...i_{l-m})}{r}_{m,l-m}\|_{L^{2}(\Sigma_{t}^{\epsilon_{0}})}\leq 
C_{l}(1+t)^{-2}[1+\log(1+t)]\{(1+t)^{-1}\mathcal{B}_{[m,l+1]}\\\notag
+\delta_{0}[\mathcal{Y}_{0}+(1+t)\mathcal{A}'_{[l]}+\mathcal{W}_{\{l+1\}}
+(1+t)^{-2}[1+\log(1+t)]^{2}\mathcal{W}^{Q}_{\{l\}}]\}
\end{align}
Moreover, since
\begin{align}
 \|\leftexp{(i_{1}...i_{l-m})}{Y}_{m,l-m}\|_{L^{2}(\Sigma_{t}^{\epsilon_{0}})}\leq C(1+t)^{-1}\mathcal{B}_{[m,l+1]}
\end{align}
while
\begin{align}
 \|\slashed{d}\log\Omega\|_{L^{\infty}(\Sigma_{t}^{\epsilon_{0}})}\leq C\delta_{0}(1+t)^{-2}
\end{align}
$\leftexp{(i_{1}...i_{l-m})}{\tilde{r}}_{m,l-m}$ enjoys the same bounds as $\leftexp{(i_{1}...i_{l-m})}{r}_{m,l-m}$.

     Now substituting (14.375) in (14.371) the hypersurface integral becomes:
\begin{align}
 -\int_{\Sigma_{t}^{u}}(\omega/\nu)(T\psi_{\alpha})\tilde{\slashed{\textrm{div}}}\leftexp{(i_{1}...i_{l-m})}{Y}_{m,l-m}
(R_{i_{l-m}}...R_{i_{1}}(T)^{m+1}\psi_{\alpha})du'd\mu_{\tilde{\slashed{g}}}\\\notag
-\int_{\Sigma_{t}^{u}}(\omega/\nu)(T\psi_{\alpha})\leftexp{(i_{1}...i_{l-m})}{\tilde{r}}_{m,l-m}
(R_{i_{l-m}}...R_{i_{1}}(T)^{m+1}\psi_{\alpha})du'd\mu_{\tilde{\slashed{g}}}
\end{align}
If $f$ is an arbitrary function defined on $S_{t,u}$ and $X$ an arbitrary vectorfield tangential to $S_{t,u}$, we have:
\begin{align}
 \int_{S_{t,u}}f\tilde{\slashed{\textrm{div}}}Xd\mu_{\tilde{\slashed{g}}}=-\int_{S_{t,u}}X\cdot\slashed{d}f
d\mu_{\tilde{\slashed{g}}}
\end{align}
Applying (14.382) to the first of (14.381), taking $X=\leftexp{(i_{1}...i_{l-m})}{Y}_{m,l-m}$ and \\
$f=(\omega/\nu)(T\psi_{\alpha})(R_{i_{l-m}}...R_{i_{1}}(T)^{m+1}\psi_{\alpha})$, we obtain that (14.381) equals:
\begin{align}
 H'_{0}+H'_{1}+H'_{2}
\end{align}
where:
\begin{align}
 H'_{0}=\int_{\Sigma_{t}^{u}}(\omega/\nu)(T\psi_{\alpha})\leftexp{(i_{1}...i_{l-m})}{Y}_{m,l-m}
\cdot\slashed{d}(R_{i_{l-m}}...R_{i_{1}}(T)^{m+1}\psi_{\alpha})du'd\mu_{\tilde{\slashed{g}}}\\
H'_{1}=\int_{\Sigma_{t}^{u}}(\omega/\nu)(\slashed{d}T\psi_{\alpha})\cdot\leftexp{(i_{1}...i_{l-m})}{Y}_{m,l-m}
(R_{i_{l-m}}...R_{i_{1}}(T)^{m+1}\psi_{\alpha})du'd\mu_{\tilde{\slashed{g}}}\\
H'_{2}=\int_{\Sigma_{t}^{u}}(T\psi_{\alpha})\{\slashed{d}(\omega/\nu)\cdot\leftexp{(i_{1}...i_{l-m})}{Y}_{m,l-m}
-(\omega/\nu)\leftexp{(i_{1}...i_{l-m})}{\tilde{r}}_{m,l-m}\}(R_{i_{l-m}}...R_{i_{1}}(T)^{m+1}\psi_{\alpha})
du'd\mu_{\tilde{\slashed{g}}}
\end{align}

Now by (14.378) and (14.379) (also (14.214)),
\begin{align}
 \|(\nu/\omega)\leftexp{(i_{1}...i_{l-m})}{Y}_{m,l-m}\cdot\slashed{d}(\omega/\nu)-\leftexp{(i_{1}...i_{l-m})}
{\tilde{r}}_{m,l-m}\|_{L^{2}(\Sigma_{t}^{\epsilon_{0}})}\\\notag
\leq C_{l}(1+t)^{-2}[1+\log(1+t)]\{(1+t)^{-1}\mathcal{B}_{[m,l+1]}+\\\notag
\delta_{0}[\mathcal{Y}_{0}+(1+t)\mathcal{A}'_{[l]}+\mathcal{W}_{\{l+1\}}+(1+t)^{-2}[1+\log(1+t)]^{2}
\mathcal{W}^{Q}_{\{l\}}]\}
\end{align}
Comparing on one hand (14.386) with (14.385) and on the other hand (14.387) with (14.379), we see that 
$H_{2}^{\prime}$ has an extra decay factor of 
$(1+t)^{-1}[1+\log(1+t)]$ relative to $H_{1}^{\prime}$.  
So we confine attention to $H'_{0}$ and $H'_{1}$. 

    We first consider the integral $H'_{0}$. We have, recalling the definition (14.372),
\begin{align}
 |H'_{0}|\leq C\int_{\Sigma_{t}^{\epsilon_{0}}}(1+t)^{2}|T\psi_{\alpha}||\slashed{d}R_{i_{l-m}}...R_{i_{1}}(T)^{m}\mu|
|\slashed{d}R_{i_{l-m}}...R_{i_{1}}(T)^{m+1}\psi_{\alpha}|dud\mu_{\slashed{g}}\\\notag
\leq C\int_{\Sigma_{t}^{\epsilon_{0}}}(1+t)|T\psi_{\alpha}|(\sum_{j}|R_{j}R_{i_{l-m}}...R_{i_{1}}(T)^{m}\mu|)
|\slashed{d}R_{i_{l-m}}...R_{i_{1}}(T)^{m+1}\psi_{\alpha}|dud\mu_{\slashed{g}}
\end{align}
Let us set $j=i_{l-m+1}$.

     Here we must obtain a sharper estimate for $R_{i_{l-m+1}}...R_{i_{1}}(T)^{m}\mu$ than that of Proposition 12.12. 
Writing:
\begin{align}
 \leftexp{(i_{1}...i_{l-m+1})}{\mu}_{m,l-m+1}=R_{i_{l-m+1}}...R_{i_{1}}(T)^{m}\mu\quad (m=0,...,l)
\end{align}
 we have the propagation equations (12.229) for $m=0$, and (12.252), for $m\geq 1$, which we can write in 
combined form as:
\begin{align}
 L\leftexp{(i_{1}...i_{l-m+1})}{\mu}_{m,l-m+1}=\leftexp{(i_{1}...i_{l-m+1})}{\rho}^{\prime}_{m,l-m+1}
\end{align}
where:
\begin{align}
 \leftexp{(i_{1}...i_{l-m+1})}{\rho}^{\prime}_{m,l-m+1}=e\leftexp{(i_{1}...i_{l-m+1})}{\mu}_{m,l-m+1}
+R_{i_{l-m+1}}...R_{i_{1}}(T)^{m}m\\\notag
+\mu R_{i_{l-m+1}}...R_{i_{1}}(T)^{m}e+\leftexp{(i_{1}...i_{l-m+1})}{r}^{\prime\prime}_{m,l-m+1}\\\notag
+\sum_{k=0}^{m-1}R_{i_{l-m+1}}...R_{i_{1}}(T)^{k}\Lambda(T)^{m-1-k}\mu\\\notag
+\sum_{k=0}^{l-m}R_{i_{l-m+1}}...R_{i_{l-m-k+2}}\leftexp{(R_{i_{l-m-k+1}})}{Z}R_{i_{l-m-k}}...R_{i_{1}}(T)^{m}\mu
\end{align}
Here we have defined the functions:
\begin{align}
 \leftexp{(i_{1}...i_{l-m+1})}{r}^{\prime\prime}_{m,l-m+1}=\leftexp{(i_{1}...i_{l-m+1})}
{r}^{\prime}_{m,l-m+1}-(R_{i_{l-m+1}}...R_{i_{1}}(T)^{m}e)\mu
\end{align}
where $\leftexp{(i_{1}...i_{l-m+1})}{r}^{\prime}_{m,l-m+1}$ are given by (12.251). 
The functions $\leftexp{(i_{1}...i_{l-m+1})}{r}^{\prime\prime}_{m,l-m+1}$ are of order 
$l+1$ and contain spatial derivatives of $\mu$ of order at most $l$ of which at most $m$ are $T$-derivatives. 
Explicitly, for $m=0$:
\begin{align}
 \leftexp{(i_{1}...i_{l+1})}{r}^{\prime\prime}_{0,l+1}=\sum_{|s_{1}|+|s_{2}|=l+1,|s_{1}|,|s_{2}|>0}((R)^{s_{1}}e)
((R)^{s_{2}}\mu)
\end{align}
and for $m\geq 1$:
\begin{align}
 \leftexp{(i_{1}...i_{l-m+1})}{r}^{\prime\prime}_{m,l-m+1}=\sum_{|s_{1}|+|s_{2}|=l-m+1,|s_{1}|>0}
((R)^{s_{1}}e)((R)^{s_{2}}(T)^{m}\mu)\\\notag
+\sum_{|s_{1}|+|s_{2}|=l-m+1,|s_{2}|>0}((R)^{s_{1}}(T)^{m}e)((R)^{s_{2}}\mu)\\\notag
+R_{i_{l-m+1}}...R_{i_{1}}(\sum_{k=1}^{m-1}\frac{m!}{k!(m-k)!}((T)^{k}e)((T)^{m-k}\mu))
\end{align}
As in Chapter 10 we write:
\begin{align}
 m=m^{\alpha}_{T}(T\psi_{\alpha})
\end{align}
and we express:
\begin{align}
 R_{i_{l-m+1}}...R_{i_{1}}(T)^{m}m=m^{\alpha}_{T}R_{i_{l-m+1}}...R_{i_{1}}(T)^{m+1}\psi_{\alpha}
+\leftexp{(i_{1}...i_{l-m+1})}{\tilde{n}}^{\prime(0)}_{m,l-m+1}
\end{align}
where:
\begin{align}
 \leftexp{(i_{1}...i_{l-m+1})}{\tilde{n}}^{\prime(0)}_{m,l-m+1}=\sum_{|s_{1}|+|s_{2}|=l-m+1,|s_{1}|>0}
((R)^{s_{1}}m^{\alpha}_{T})((R)^{s_{2}}(T)^{m+1}\psi_{\alpha})\\\notag
+R_{i_{l-m+1}}...R_{i_{1}}(\sum_{k=1}^{m}\frac{m!}{k!(m-k)!}((T)^{k}m^{\alpha}_{T})((T)^{m-k+1}\psi_{\alpha}))
\end{align}
Also, we write:
\begin{align}
 e=e^{\alpha}_{L}(L\psi_{\alpha})
\end{align}
where:
\begin{align}
 e_{L}^{0}=\frac{1}{2\eta^{2}}(\eta^{2})^{\prime}\\\notag
e^{i}_{L}=\eta^{-1}\hat{T}^{i}-\frac{1}{2\eta^{2}}(\eta^{2})^{\prime}\psi_{i}
\end{align}
Obviously,
\begin{align}
 |e^{\alpha}_{L}|\leq C
\end{align}
We then express:
\begin{align}
 R_{i_{l-m+1}}...R_{i_{1}}(T)^{m}e=e^{\alpha}_{L}R_{i_{l-m+1}}...R_{i_{1}}(T)^{m}L\psi_{\alpha}
+\leftexp{(i_{1}...i_{l-m+1})}{\tilde{n}}^{\prime(1)}_{m,l-m+1}
\end{align}
where:
\begin{align}
 \leftexp{(i_{1}...i_{l-m+1})}{\tilde{n}}^{\prime(1)}_{m,l-m+1}=\sum_{|s_{1}|+|s_{2}|=l-m+1,|s_{1}|>0}
((R)^{s_{1}}e^{\alpha}_{L})((R)^{s_{2}}(T)^{m}L\psi_{\alpha})\\\notag
+R_{i_{l-m+1}}...R_{i_{1}}(\sum_{k=1}^{m}\frac{m!}{k!(m-k)!}((T)^{k}e^{\alpha}_{L})((T)^{m-k}L\psi_{\alpha}))
\end{align}
We define:
\begin{align}
 \leftexp{(i_{1}...i_{l-m+1})}{\rho}^{\prime(0)}_{m,l-m+1}=\frac{1}{2}\ell R_{i_{l-m+1}}...R_{i_{1}}(T)^{m+1}\psi_{0}
\end{align}
and:
\begin{align}
 \leftexp{(i_{1}...i_{l-m+1})}{\rho}^{\prime(1)}_{m,l-m+1}=(m^{0}_{T}-\frac{1}{2}\ell)R_{i_{l-m+1}}...R_{i_{1}}
(T)^{m+1}\psi_{0}\\\notag
+m^{i}_{T}R_{i_{l-m+1}}...R_{i_{1}}(T)^{m+1}\psi_{i}+(1+t)^{-1}\mu e^{\alpha}_{L}R_{i_{l-m+1}}...R_{i_{1}}
(T)^{m}Q\psi_{\alpha}
\end{align}
We then have:
\begin{align}
 R_{i_{l-m+1}}...R_{i_{1}}(T)^{m}m+\mu R_{i_{l-m+1}}...R_{i_{1}}(T)^{m}e\\\notag
=\leftexp{(i_{1}...i_{l-m})}{\rho}^{\prime(0)}_{m,l-m+1}+\leftexp{(i_{1}...i_{l-m+1})}{\rho}^{\prime(1)}_{m,l-m+1}
+\leftexp{(i_{1}...i_{l-m+1})}{\tilde{n}}^{\prime}_{m,l-m+1}
\end{align}
where:
\begin{align}
 \leftexp{(i_{1}...i_{l-m+1})}{\tilde{n}}^{\prime}_{m,l-m+1}=\leftexp{(i_{1}...i_{l-m+1})}{\tilde{n}}
^{\prime(0)}_{m,l-m+1}+\mu\leftexp{(i_{1}...i_{l-m+1})}{\tilde{n}}^{\prime(1)}_{m,l-m+1}
\end{align}
and (14.391) takes the form:
\begin{align}
 \leftexp{(i_{1}...i_{l-m+1})}{\rho}^{\prime}_{m,l-m+1}=\leftexp{(i_{1}...i_{l-m+1})}{\rho}^{\prime(0)}_{m,l-m+1}+
\leftexp{(i_{1}...i_{l-m+1})}{\rho}^{\prime(1)}_{m,l-m+1}+\leftexp{(i_{1}...i_{l-m+1})}{\rho}^{\prime(2)}_{m,l-m+1}
\end{align}
where:
\begin{align}
 \leftexp{(i_{1}...i_{l-m+1})}{\rho}^{\prime(2)}_{m,l-m+1}=e\leftexp{(i_{1}...i_{l-m+1})}{\mu}_{m,l-m+1}\\\notag
+\leftexp{(i_{1}...i_{l-m+1})}{\tilde{n}}^{\prime}_{m,l-m+1}+\leftexp{(i_{1}...i_{l-m+1})}{r}^{\prime\prime}_{m,l-m+1}
\\\notag
+\sum_{k=0}^{m-1}R_{i_{l-m+1}}...R_{i_{1}}(T)^{k}\Lambda(T)^{m-1-k}\mu\\\notag
+\sum_{k=0}^{l-m}R_{i_{l-m+1}}...R_{i_{l-m-k+2}}\leftexp{(R_{i_{l-m-k+1}})}{Z}R_{i_{l-m-k}}...R_{i_{1}}(T)^{m}\mu
\end{align}

     We shall now obtain an $L^{2}(\Sigma_{t}^{\epsilon_{0}})$ estimate for 
$\leftexp{(i_{1}...i_{l-m+1})}{\rho}^{\prime(2)}_{m,l-m+1}$. The first term on the right in 
(14.408) is bounded in $L^{2}(\Sigma_{t}^{\epsilon_{0}})$ by:
\begin{align}
 C\delta_{0}(1+t)^{-2}\mathcal{B}_{[m,l+1]}
\end{align}
From (14.397) we deduce:
\begin{align}
 \|\leftexp{(i_{1}...i_{l-m+1})}{\tilde{n}}^{\prime(0)}_{m,l-m+1}\|_{L^{2}(\Sigma_{t}^{\epsilon_{0}})}\leq 
C\delta_{0}(1+t)^{-1}\mathcal{W}_{\{l+1\}}
\end{align}
Using Corollary 11.2.b, we deduce that
$\|\leftexp{(i_{1}...i_{l-m+1})}{\tilde{n}}^{\prime(1)}_{m,l-m+1}\|_{L^{2}(\Sigma_{t}^{\epsilon_{0}})}$ is bounded by:
\begin{align}
 C_{l}\delta_{0}(1+t)^{-1}\{\mathcal{W}_{\{l+1\}}+\mathcal{W}^{Q}_{\{l\}}+\delta_{0}(1+t)^{-1}[(1+t)^{-1}
\mathcal{B}_{[m-1,l+1]}+\mathcal{Y}_{0}+(1+t)\mathcal{A}'_{[l]}]\}
\end{align}
We turn to the third term on the right in (14.408). Using Corollary 11.2.b we deduce through (14.393), (14.394),
\begin{align}
 \|\leftexp{(i_{1}...i_{l-m})}{r}^{\prime\prime}_{m,l-m+1}\|_{L^{2}(\Sigma_{t}^{\epsilon_{0}})}\\\notag
\leq C\delta_{0}\{[1+\log(1+t)](1+t)^{-1}[\mathcal{W}^{Q}_{\{l\}}+\delta_{0}(1+t)^{-1}(\mathcal{Y}_{0}
+(1+t)\mathcal{A}'_{[l-1]}+\mathcal{W}_{\{l\}})]\\\notag+(1+t)^{-2}\mathcal{B}_{[m,l]}\}
\end{align}
Finally, the two sums on the right in (14.408) have already been estimated in $L^{2}(\Sigma_{t}^{\epsilon_{0}})$ 
in (12.357) and (12.359). So we conclude that:
\begin{align}
 \sum_{j}\|\leftexp{(i_{1}...i_{l-m}j)}{\rho}^{\prime(2)}_{m,l-m+1}\|_{L^{2}(\Sigma_{t}^{\epsilon_{0}})}\\\notag
\leq C_{l}\delta_{0}(1+t)^{-1}[1+\log(1+t)]\\\notag
\cdot \{\mathcal{W}_{\{l+1\}}+\mathcal{W}^{Q}_{\{l\}}+(1+t)^{-1}\mathcal{B}_{[m,l+1]}+\mathcal{Y}_{0}+(1+t)
\mathcal{A}'_{[l]}\}
\end{align}
Moreover, from (14.403) and (14.404), similarly as before, we have:
\begin{align}
 \sum_{j}\|\leftexp{(i_{1}...i_{l-m}j)}{\rho}^{\prime(0)}_{m,l-m+1}\|_{L^{2}(\Sigma_{t}^{\epsilon_{0}})}
\leq C|\ell|\bar{\mu}^{-1/2}_{m}\sqrt{\sum_{\alpha}\mathcal{E}'_{1}[R_{i_{l-m}}...R_{i_{1}}(T)^{m+1}\psi_{0}](t)}
\end{align}
\begin{align}
 \sum_{j}\|\leftexp{(i_{1}...i_{l-m}j)}{\rho}^{\prime(1)}_{m,l-m+1}\|_{L^{2}(\Sigma_{t}^{\epsilon_{0}})}\\\notag
\leq C\delta_{0}(1+t)^{-1}\bar{\mu}^{-1/2}_{m}\sqrt{\sum_{\alpha}\mathcal{E}'_{1}[R_{i_{l-m}}...R_{i_{1}}(T)^{m+1}
\psi_{\alpha}](t)}\\\notag
+C\delta_{0}(1+t)^{-1}[1+\log(1+t)]^{1/2}\sqrt{\sum_{\alpha}\mathcal{E}'_{1}
[R_{i_{l-m}}...R_{i_{1}}(T)^{m}Q\psi_{\alpha}](t)}
\end{align}
and:
\begin{align}
 \sum_{j}\|\mu^{1/2}\leftexp{(i_{1}...i_{l-m}j)}{\rho}^{\prime(0)}_{m,l-m+1}\|_{L^{2}(\Sigma_{t}^{\epsilon_{0}})}
\leq C|\ell|\sqrt{\sum_{\alpha}\mathcal{E}'_{1}[R_{i_{l-m}}...R_{i_{1}}(T)^{m+1}\psi_{\alpha}](t)}
\end{align}
\begin{align}
 \sum_{j}\|\mu^{1/2}\leftexp{(i_{1}...i_{l-m}j)}{\rho}^{\prime(1)}_{m,l-m+1}\|_{L^{2}(\Sigma_{t}^{\epsilon_{0}})}\\\notag
\leq C\delta_{0}(1+t)^{-1}\sqrt{\sum_{\alpha}\mathcal{E}'_{1}[R_{i_{l-m}}...R_{i_{1}}(T)^{m+1}\psi_{\alpha}](t)}\\\notag
+C\delta_{0}(1+t)^{-1}[1+\log(1+t)]\sqrt{\sum_{\alpha}\mathcal{E}'_{1}[R_{i_{l-m}}...R_{i_{1}}(T)^{m}Q\psi_{\alpha}](t)}
\end{align}

    We now integrate the propagation equation (14.390) along each integral curve of $L$ from 
$\Sigma_{0}^{\epsilon_{0}}$ to obtain:
\begin{align}
 \sum_{j}|\leftexp{(i_{1}...i_{l-m}j)}{\mu}_{m,l-m+1}(t,u,\vartheta)|\\\notag
\leq C\{\sum_{j}|\leftexp{(i_{1}...i_{l-m}j)}{\mu}_{m,l-m+1}(0,u,\vartheta)|+\sum_{k=0}^{2}
\leftexp{(i_{1}...i_{l-m}j)}{A}^{\prime(k)}_{m,l-m}(t,u,\vartheta)\}
\end{align}
where:
\begin{align}
 \leftexp{(i_{1}...i_{l-m})}{A}^{\prime(k)}_{m,l-m}(t,u,\vartheta)=\int_{0}^{t}\sum_{j}|
\leftexp{(i_{1}...i_{l-m}j)}{\rho}^{\prime(k)}_{m,l-m+1}|(t',u,\vartheta)dt'
\end{align}
for $k=0,1,2$.

By (8.333) we have:
\begin{align}
 \|\leftexp{(i_{1}...i_{l-m})}{A}^{\prime(k)}_{m,l-m}(t)\|_{L^{2}([0,\epsilon_{0}]\times S^{2})}\\\notag
\leq\int_{0}^{t}\sum_{j}\|\leftexp{(i_{1}...i_{l-m})}{\rho}^{\prime(k)}_{m,l-m+1}(t')\|
_{L^{2}([0,\epsilon_{0}]\times S^{2})}dt'\\\notag
\leq\int_{0}^{t}(1+t')^{-1}\sum_{j}\|\leftexp{(i_{1}...i_{l-m}j)}{\rho}^{\prime(k)}_{m,l-m+1}\|
_{L^{2}(\Sigma_{t'}^{\epsilon_{0}})}dt'
\end{align}

     We now substitute (14.418) in (14.388), noting that:
\begin{align*}
 \sum_{j}|R_{j}R_{i_{l-m}}...R_{i_{1}}(T)^{m}\mu|=\sum_{j}|\leftexp{(i_{1}...i_{l-m}j)}{\mu}_{m,l-m+1}|
\end{align*}
The $borderline$ contribution is from $\leftexp{(i_{1}...i_{l-m})}{A}^{\prime(0)}_{m,l-m}$. It is the $borderline$ 
integral:
\begin{align}
 C\int_{\Sigma_{t}^{\epsilon_{0}}}(1+t)|T\psi_{\alpha}|\leftexp{(i_{1}...i_{l-m})}{A}^{\prime(0)}_{m,l-m}
|\slashed{d}R_{i_{l-m}}...R_{i_{1}}(T)^{m+1}\psi_{\alpha}|
dud\mu_{\slashed{g}}
\end{align}
This is similar to (14.246) and it is estimated in exactly the same manner. Defining (see (14.257)):
\begin{align}
 \leftexp{(i_{1}...i_{l-m})}{\mathcal{G}}^{\prime}_{1,m,l+2;a,q}(t)\\\notag
=\sup_{t'\in[0,t]}\{[1+\log(1+t')]^{-2q}\bar{\mu}^{2a}_{m}(t')\sum_{\alpha}\mathcal{E}'_{1}
[R_{i_{l-m}}...R_{i_{1}}(T)^{m+1}\psi_{\alpha}](t')\}
\end{align}
then similarly, we obtain (see (14.266), (14.293)) that the borderline integral (14.421) is bounded by:
\begin{align}
 C(\frac{1}{a-1/2}+\frac{1}{q+1/2})\bar{\mu}^{-2a}_{m}(t)[1+\log(1+t)]^{2q}\leftexp{(i_{1}...i_{l-m})}
{\mathcal{G}}^{\prime}_{1,m,l+2;a,q}(t)
\end{align}

    The remaining contributions to (14.388) are estimated in exactly the same manner as the remaining 
contributions to (14.215). In view of (14.294) these contributions are bounded by:
\begin{align}
 C\int_{\Sigma_{t}^{\epsilon_{0}}}(1+t)|T\psi_{\alpha}|\leftexp{(i_{1}...i_{l-m})}{A}^{\prime(k)}_{m,l-m}
|\slashed{d}R_{i_{l-m}}...R_{i_{1}}(T)^{m+1}\psi_{\alpha}|
dud\mu_{\slashed{g}}\\\notag
\leq C\delta_{0}\int_{\Sigma_{t}^{\epsilon_{0}}}\leftexp{(i_{1}...i_{l-m})}{A}^{\prime(k)}_{m,l-m}
|\slashed{d}R_{i_{l-m}}...R_{i_{1}}(T)^{m+1}\psi_{\alpha}|dud\mu_{\slashed{g}}
\\\notag
\leq C\delta_{0}\|\leftexp{(i_{1}...i_{l-m})}{A}^{\prime(k)}_{m,l-m}\|_{L^{2}(\Sigma_{t}^{\epsilon_{0}})}
\|\slashed{d}R_{i_{l-m}}...R_{i_{1}}(T)^{m+1}\psi_{\alpha}\|_{L^{2}(\Sigma_{t}^{\epsilon_{0}})}\\\notag
\leq C\delta_{0}\bar{\mu}^{-1/2}_{m}(t)\|\leftexp{(i_{1}...i_{l-m})}{A}^{\prime(k)}_{m,l-m}(t)\|
_{L^{2}([0,\epsilon_{0}]\times S^{2})}\sqrt{\mathcal{E}'_{1,[l+2]}(t)}
\end{align}
for $k=1,2$,
which is similar to (14.297). Here we need only consider the contribution of\\
$\leftexp{(i_{1}...i_{l-m})}{A}^{\prime(1)}_{m,l-m}$, because $\leftexp{(i_{1}...i_{l-m})}{A}^{\prime(2)}_{m,l-m}$ 
besides being of lower order does not contain any leading terms, as is evident by comparing (14.413) to (14.237).
By (14.415):
\begin{align}
 \sum_{j}\|\leftexp{(i_{1}...i_{l-m}j)}{\rho}^{\prime(1)}_{m,l-m+1}\|_{L^{2}(\Sigma_{t}^{\epsilon_{0}})}\leq 
C\delta_{0}(1+t)^{-1}[1+\log(1+t)]^{1/2}\bar{\mu}^{-1/2}_{m}(t)\sqrt{\mathcal{E}'_{1,[l+2]}(t)}
\end{align}
we obtain:
\begin{align}
 \|\leftexp{(i_{1}...i_{l-m})}{A}^{\prime(1)}_{m,l-m}(t)\|_{L^{2}([0,\epsilon_{0}]\times S^{2})}\\\notag
\leq C\delta_{0}\int_{0}^{t}\bar{\mu}^{-1/2}_{m}(t')(1+t')^{-2}[1+\log(1+t')]^{1/2}
\sqrt{\mathcal{E}^{\prime}_{1,[l+2]}(t')}dt'\\\notag\leq C\delta_{0}J'_{a,q+1/2}(t)\sqrt{\mathcal{G}'_{1,[l+2];a,q}(t)}
\end{align}
where $J'_{a,q}(t)$ is the integral (14.298). Substituting (14.305) with $q$ replaced by $q+1/2$ yields:
\begin{align}
 C\int_{\Sigma_{t}^{\epsilon_{0}}}(1+t)|T\psi_{\alpha}|\leftexp{(i_{1}...i_{l-m})}{A}^{\prime(1)}_{m,l-m}
|\slashed{d}R_{i_{l-m}}...R_{i_{1}}(T)^{m+1}\psi_{\alpha}|dud\mu_{\slashed{g}}\\\notag
\leq C\delta_{0}^{2}\mathcal{G}'_{1,[l+2];a,q}(t)\{C_{q}\bar{\mu}^{-a-1/2}_{m}(t)+C\varphi^{\prime}_{q+3/2}
(Ca\delta_{0})\bar{\mu}^{-2a}_{m}(t)\}[1+\log(1+t)]^{q}
\end{align}

Finally, we consider the hypersurface integral $H'_{1}$, given by (14.387). We have:
\begin{align}
 |H'_{1}|\leq C\int_{\Sigma_{t}^{\epsilon_{0}}}(1+t)|\slashed{d}T\psi_{\alpha}|\sum_{j}
|R_{j}R_{i_{l-m}}...R_{i_{1}}(T)^{m}\mu||R_{i_{l-m}}...R_{i_{1}}(T)^{m+1}\psi_{\alpha}|dud\mu_{\slashed{g}}
\end{align}
Using 
\begin{align}
 \max_{\alpha}\sup_{\Sigma_{t}^{\epsilon_{0}}}|\slashed{d}T\psi_{\alpha}|\leq C\delta_{0}(1+t)^{-2}
\end{align}
we have:
\begin{align}
 |H'_{1}|\leq C\delta_{0}\int_{\Sigma_{t}^{\epsilon_{0}}}(1+t)^{-1}\sum_{j}|R_{j}R_{i_{l-m}}...R_{i_{1}}(T)^{m}\mu|
|R_{i_{l-m}}...R_{i_{1}}(T)^{m+1}\psi_{\alpha}|dud\mu_{\slashed{g}}\\\notag
\leq C\delta_{0}(1+t)^{-1}\sum_{j}\|R_{j}R_{i_{l-m}}...R_{i_{1}}(T)^{m}\mu\|_{L^{2}(\Sigma_{t}^{\epsilon_{0}})}
\|R_{i_{l-m}}...R_{i_{1}}(T)^{m+1}\psi_{\alpha}\|_{L^{2}(\Sigma_{t}^{\epsilon_{0}})}\\\notag
\leq C\delta_{0}\sum_{j}\|(R_{j}R_{i_{l-m}}...R_{i_{1}}(T)^{m}\mu)(t)\|_{L^{2}([0,\epsilon_{0}]\times S^{2})}
\mathcal{W}_{\{l+1\}}(t)
\end{align}
From (14.418) and (14.420):
\begin{align}
 \sum_{j}\|(R_{j}R_{i_{l-m}}...R_{i_{1}}(T)^{m}\mu)(t)\|_{L^{2}([0,\epsilon_{0}]\times S^{2})}\\\notag
\leq C\{\sum_{j}\|R_{j}R_{i_{l-m}}...R_{i_{1}}(T)^{m}\mu\|_{L^{2}(\Sigma_{0}^{\epsilon_{0}})}\\\notag
+\sum_{k=0}^{2}\|\leftexp{(i_{1}...i_{l-m})}{A}^{\prime(k)}_{m,l-m}(t)\|_{L^{2}([0,\epsilon_{0}]\times S^{2})}\}\\\notag
\leq C\{\sum_{j}\|R_{j}R_{i_{l-m}}...R_{i_{1}}(T)^{m}\mu\|_{L^{2}(\Sigma_{0}^{\epsilon_{0}})}\\\notag
+\sum_{k=0}^{2}\int_{0}^{t}(1+t')^{-1}\sum_{j}\|\leftexp{(i_{1}...i_{l-m}j)}{\rho}^{\prime(k)}_{m,l-m+1}\|
_{L^{2}(\Sigma_{t'}^{\epsilon_{0}})}dt'\}
\end{align}
Here we need only consider the leading principal contribution, namely that of 
\begin{align*}
 \sum_{j}\|\leftexp{(i_{1}...i_{l-m}j)}{\rho'}^{(0)}_{m,l-m+1}\|_{L^{2}(\Sigma_{t'}^{\epsilon_{0}})}
\end{align*}
This contribution to (14.431) is bounded by (see (14.414)):
\begin{align}
 C|\ell|\int_{0}^{t}(1+t')^{-1}\bar{\mu}^{-1/2}_{m}(t')\sqrt{\sum_{\alpha}\mathcal{E}'_{1}
[R_{i_{l-m}}...R_{i_{1}}(T)^{m+1}\psi_{\alpha}](t')}dt'\\\notag
\leq C|\ell|\sqrt{\mathcal{G}'_{1,[l+2];a,q}(t)}\int_{0}^{t}(1+t')^{-1}\bar{\mu}^{-a-1/2}_{m}(t')
[1+\log(1+t')]^{q}dt'\\\notag=C|\ell|J_{a,q-1}(t)\sqrt{\mathcal{G}'_{1,[l+2];a,q}(t)}
\end{align}
Substituting the bound (14.153) with $q$ replaced by $q-1$ and also the bound (14.318) for $\mathcal{W}_{\{l+1\}}$, 
we obtain that the leading contribution to $H'_{1}$ is bounded by:
\begin{align} 
 C\epsilon_{0}\delta_{0}\bar{\mu}^{-2a+1/2}_{m}(t)[1+\log(1+t)]^{p+q+1}\sqrt{\mathcal{G}_{0,[l+2];a,p}(t)
\mathcal{G}'_{1,[l+2];a,q}(t)}
\end{align}
the same in form as (14.319).

    We turn to the spacetime integral in (14.370), namely:
\begin{align}
 \int_{W^{t}_{u}}(R_{i_{l-m}}...R_{i_{1}}(T)^{m+1}\psi_{\alpha})(L+\nu)\{(\omega/\nu)
(T\psi_{\alpha})(R_{i_{l-m}}...R_{i_{1}}(T)^{m}\slashed{\Delta}\mu)\}dt'du'd\mu_{\tilde{\slashed{g}}}
\end{align}
As in (14.334) we have, in regard to the integrand in (14.434),
\begin{align}
 (L+\nu)\{(\omega/\nu)(T\psi_{\alpha})(R_{i_{l-m}}...R_{i_{1}}(T)^{m}\slashed{\Delta}\mu)\}\\\notag
=(\omega/\nu)\{(T\psi_{\alpha})(L+2\nu)(R_{i_{l-m}}...R_{i_{1}}(T)^{m}\slashed{\Delta}\mu)
+\tilde{\tau}_{\alpha}(R_{i_{l-m}}...R_{i_{1}}(T)^{m}\slashed{\Delta}\mu)\}
\end{align}
We substitute on the right-hand side the expression (14.375) for $R_{i_{l-m}}...R_{i_{1}}(T)^{m}\slashed{\Delta}\mu$.
Here we apply the following: Let $X$ be an arbitrary $S_{t,u}$-tangential vectorfield. We then have:
\begin{align}
 (L+2\nu)\tilde{\slashed{\textrm{div}}}X=\tilde{\slashed{\textrm{div}}}(\slashed{\mathcal{L}}_{L}X+2\nu X)
\end{align}
To establish this we work in acoustical coordinates $(t,u,\vartheta^{A} : A=1,2)$. We then have:
\begin{align}
 X=X^{A}\frac{\partial}{\partial\vartheta^{A}}\\\notag
\slashed{\textrm{div}}X=\frac{1}{\sqrt{\det\slashed{g}}}\frac{\partial}{\partial\vartheta^{A}}
(\sqrt{\det\slashed{g}}X^{A})=\frac{\partial X^{A}}{\partial \vartheta^{A}}+X^{A}\frac{\partial}
{\partial \vartheta^{A}}\log\sqrt{\det\slashed{g}}
\end{align}
and:
\begin{align}
 L\slashed{\textrm{div}}X=\frac{\partial}{\partial t}\slashed{\textrm{div}}X\\\notag
=\frac{\partial^{2}X^{A}}{\partial t\partial\vartheta^{A}}+\frac{\partial X^{A}}{\partial t}\frac{\partial}
{\partial\vartheta^{A}}\log\sqrt{\det\slashed{g}}+X^{A}\frac{\partial^{2}}{\partial t\partial\vartheta^{A}}
\log\sqrt{\det\slashed{g}}
\end{align}
Now:
\begin{align}
 \slashed{\mathcal{L}}_{L}X=[L,X]=\frac{\partial X^{A}}{\partial t}\frac{\partial}{\partial\vartheta^{A}}
\end{align}
and:
\begin{align}
 \frac{\partial}{\partial t}\log\sqrt{\det\slashed{g}}=\textrm{tr}\chi
\end{align}
We thus obtain:
\begin{align}
 L\slashed{\textrm{div}}X=\slashed{\textrm{div}}(\slashed{\mathcal{L}}_{L}X)+X\cdot\slashed{d}\textrm{tr}\chi
\end{align}
In view of (14.374) and the definition of $\nu$, this is equivalent to:
\begin{align}
 L\tilde{\slashed{\textrm{div}}}X=\tilde{\slashed{\textrm{div}}}(\slashed{\mathcal{L}}_{L}X)+2X\cdot\slashed{d}\nu
\end{align}
which is equivalent to (14.436).

     We appeal to (14.436) taking $X=\leftexp{(i_{1}...i_{l-m})}{Y}_{m,l-m}$ after substituting the expression (14.375) 
in (14.437). Substituting then the result in (14.434) and integrating by parts on each $S_{t,u}$ we obtain:
\begin{align}
 \int_{W^{t}_{u}}(\omega/\nu)(T\psi_{\alpha})(R_{i_{l-m}}...R_{i_{1}}(T)^{m+1}\psi_{\alpha})
(L+2\nu)\tilde{\slashed{\textrm{div}}}\leftexp{(i_{1}...i_{l-m})}{Y}_{m,l-m}dt'du'd\mu_{\tilde{\slashed{g}}}\\\notag
=-\int_{W^{t}_{u}}(\omega/\nu)(T\psi_{\alpha})((\slashed{\mathcal{L}}_{L}+2\nu)\leftexp{(i_{1}...i_{l-m})}{Y}_{m,l-m})
\cdot\slashed{d}R_{i_{l-m}}...R_{i_{1}}(T)^{m+1}\psi_{\alpha}dt'du'd\mu_{\tilde{\slashed{g}}}\\\notag
-\int_{W^{t}_{u}}(\omega/\nu)(R_{i_{l-m}}...R_{i_{1}}(T)^{m+1}\psi_{\alpha})((\slashed{\mathcal{L}}_{L}+2\nu)
\leftexp{(i_{1}...i_{l-m})}{Y}_{m,l-m})\cdot\slashed{d}T\psi_{\alpha}dt'du'd\mu_{\tilde{\slashed{g}}}\\\notag
-\int_{W^{t}_{u}}(T\psi_{\alpha})(R_{i_{l-m}}...R_{i_{1}}(T)^{m+1}\psi_{\alpha})\\\notag
\cdot ((\slashed{\mathcal{L}}_{L}+2\nu)\leftexp{(i_{1}...i_{l-m})}{Y}_{m,l-m})\cdot\slashed{d}(\omega/\nu)
dt'du'd\mu_{\tilde{\slashed{g}}}
\end{align}
and we are left with:
\begin{align}
 \int_{W^{t}_{u}}(\omega/\nu)(T\psi_{\alpha})(R_{i_{l-m}}...R_{i_{1}}(T)^{m+1}\psi_{\alpha})(L+2\nu)
\leftexp{(i_{1}...i_{l-m})}{\tilde{r}}_{m,l-m}dt'du'd\mu_{\tilde{\slashed{g}}}
\end{align}
We also integrate by parts on each $S_{t,u}$:
\begin{align}
 \int_{W^{t}_{u}}(\omega/\nu)\tilde{\tau}_{\alpha}(R_{i_{l-m}}...R_{i_{1}}(T)^{m+1}\psi_{\alpha})
\tilde{\slashed{\textrm{div}}}\leftexp{(i_{1}...i_{l-m})}{Y}_{m,l-m}dt'du'd\mu_{\tilde{\slashed{g}}}\\\notag
=-\int_{W^{t}_{u}}(\omega/\nu)\tilde{\tau}_{\alpha}\leftexp{(i_{1}...i_{l-m})}{Y}_{m,l-m}\cdot\slashed{d}
(R_{i_{l-m}}...R_{i_{1}}(T)^{m+1}\psi_{\alpha})dt'du'd\mu_{\tilde{\slashed{g}}}\\\notag
-\int_{W^{t}_{u}}(\omega/\nu)(R_{i_{l-m}}...R_{i_{1}}(T)^{m+1}\psi_{\alpha})\leftexp{(i_{1}...i_{l-m})}{Y}_{m,l-m}
\cdot\slashed{d}\tilde{\tau}_{\alpha}dt'du'd\mu_{\tilde{\slashed{g}}}\\\notag
-\int_{W^{t}_{u}}\tilde{\tau}_{\alpha}(R_{i_{l-m}}...R_{i_{1}}(T)^{m+1}\psi_{\alpha})
\leftexp{(i_{1}...i_{l-m})}{Y}_{m,l-m}\cdot\slashed{d}(\omega/\nu)dt'du'd\mu_{\tilde{\slashed{g}}}
\end{align}
and we are left with:
\begin{align}
 \int_{W^{t}_{u}}(\omega/\nu)\tilde{\tau}_{\alpha}(R_{i_{l-m}}...R_{i_{1}}(T)^{m+1}\psi_{\alpha})
\leftexp{(i_{1}...i_{l-m})}{\tilde{r}}_{m,l-m}dt'du'd\mu_{\tilde{\slashed{g}}}
\end{align}
In view of (14.443)-(14.446) the spacetime integral (14.434) becomes:
\begin{align}
 -V'_{0}-V'_{1}-V'_{2}
\end{align}
where:
\begin{align}
 V'_{0}=V'_{0,0}+V'_{0,1}
\end{align}
\begin{align}
 V'_{0,0}=\int_{W^{t}_{u}}(\omega/\nu)(T\psi_{\alpha})((\slashed{\mathcal{L}}_{L}+2\nu)
\leftexp{(i_{1}...i_{l-m})}{Y}_{m,l-m})\cdot\slashed{d}R_{i_{l-m}}...R_{i_{1}}(T)^{m+1}\psi_{\alpha}
dt'du'd\mu_{\tilde{\slashed{g}}}
\end{align}
\begin{align}
 V'_{0,1}=\int_{W^{t}_{u}}(\omega/\nu)\tilde{\tau}_{\alpha}\leftexp{(i_{1}...i_{l-m})}
{Y}_{m,l-m}\cdot\slashed{d}R_{i_{l-m}}...R_{i_{1}}(T)^{m+1}\psi_{\alpha}dt'du'd\mu_{\tilde{\slashed{g}}}
\end{align}
\begin{align}
 V'_{1}=V'_{1,0}+V'_{1,1}
\end{align}
\begin{align}
 V'_{1,0}=\int_{W^{t}_{u}}(\omega/\nu)(R_{i_{l-m}}...R_{i_{1}}(T)^{m+1}\psi_{\alpha})((\slashed{\mathcal{L}}
_{L}+2\nu)\leftexp{(i_{1}...i_{l-m})}{Y}_{m,l-m})\cdot\slashed{d}T\psi_{\alpha}dt'du'd\mu_{\tilde{\slashed{g}}}
\end{align}
\begin{align}
 V'_{1,1}=\int_{W^{t}_{u}}(\omega/\nu)(R_{i_{l-m}}...R_{i_{1}}(T)^{m+1}\psi_{\alpha})
\leftexp{(i_{1}...i_{l-m})}{Y}_{m,l-m}\cdot\slashed{d}\tilde{\tau}_{\alpha}dt'du'd\mu_{\tilde{\slashed{g}}}
\end{align}
and:
\begin{align}
 V'_{2}=V'_{2,0}+V'_{2,1}
\end{align}
\begin{align}
 V'_{2,0}=\int_{W^{t}_{u}}(T\psi_{\alpha})(R_{i_{l-m}}...R_{i_{1}}(T)^{m+1}\psi_{\alpha})\\\notag
\{((\slashed{\mathcal{L}}_{L}+2\nu)\leftexp{(i_{1}...i_{l-m})}{Y}_{m,l-m})\cdot\slashed{d}(\omega/\nu)\\\notag
-(\omega/\nu)(L+2\nu)\leftexp{(i_{1}...i_{l-m})}{\tilde{r}}_{m,l-m}\}dt'du'd\mu_{\tilde{\slashed{g}}}
\end{align}
\begin{align}
 V'_{2,1}=\int_{W^{t}_{u}}(R_{i_{l-m}}...R_{i_{1}}(T)^{m+1}\psi_{\alpha})\tilde{\tau}_{\alpha}\\\notag
(\leftexp{(i_{1}...i_{l-m})}{Y}_{m,l-m}\cdot(\slashed{d}(\omega/\nu))-(\omega/\nu)\leftexp{(i_{1}...i_{l-m})}
{\tilde{r}}_{m,l-m})dt'du'd\mu_{\tilde{\slashed{g}}}
\end{align}
In the following, we estimate $V'_{0,0}$ and $V'_{1,0}$, seeing that the other terms are of lower order.

We have:
\begin{align}
 |V'_{0,0}|\leq C\int_{W^{t}_{\epsilon_{0}}}(1+t')^{2}|T\psi_{\alpha}||\slashed{d}R_{i_{l-m}}...R_{i_{1}}(T)^{m+1}
\psi_{\alpha}||(\slashed{\mathcal{L}}_{L}+2\nu)\leftexp{(i_{1}...i_{l-m})}{Y}_{m,l-m}|dt'dud\mu_{\slashed{g}}
\end{align}
Since according to (14.372):
\begin{align}
 \slashed{d}R_{i_{l-m}}...R_{i_{1}}(T)^{m}\mu=\slashed{g}\cdot\leftexp{(i_{1}...i_{l-m})}{Y}_{m,l-m}
\end{align}
and $\slashed{\mathcal{L}}_{L}\slashed{g}=2\chi$, we have:
\begin{align}
 \slashed{d}LR_{i_{l-m}}...R_{i_{1}}(T)^{m}\mu=\slashed{\mathcal{L}}_{L}\slashed{d}R_{i_{l-m}}...R_{i_{1}}(T)^{m}
\mu\\\notag
=2\chi\cdot\leftexp{(i_{1}...i_{l-m})}{Y}_{m,l-m}+\slashed{g}\cdot\slashed{\mathcal{L}}_{L}\leftexp{(i_{1}...i_{l})}
{Y}_{m,l-m}
\end{align}
Since
\begin{align*}
 \chi=\frac{1}{2}\textrm{tr}\chi\slashed{g}+\hat{\chi}=(\nu-\frac{1}{2}L\log\Omega)\slashed{g}+\hat{\chi}^{\prime}
\end{align*}
defining:
\begin{align}
 \tilde{\chi}^{\prime}=\hat{\chi}^{\prime}-\frac{1}{2}(L\log\Omega)\slashed{g}
\end{align}
we then obtain:
\begin{align}
 \slashed{g}\cdot(\slashed{\mathcal{L}}_{L}+2\nu)\leftexp{(i_{1}...i_{l-m})}{Y}_{m,l-m}=
\slashed{d}LR_{i_{l-m}}...R_{i_{1}}(T)^{m}\mu-2\tilde{\chi}^{\prime}\cdot\leftexp{(i_{1}...i_{l-m})}{Y}_{m,l-m}
\end{align}
Substituting for $LR_{i_{l-m}}...R_{i_{1}}(T)^{m}\mu$ from (14.390) with $l+1$ replaced by $l$ then yields:
\begin{align}
 \slashed{g}\cdot(\slashed{\mathcal{L}}_{L}+2\nu)\leftexp{(i_{1}...i_{l-m})}{Y}_{m,l-m}=\slashed{d}
\leftexp{(i_{1}...i_{l-m})}{\rho}^{\prime}_{m,l-m}-2\tilde{\chi}^{\prime}\cdot\slashed{d}
R_{i_{l-m}}...R_{i_{1}}(T)^{m}\mu
\end{align}
which implies 
\begin{align}
 |(\slashed{\mathcal{L}}_{L}+2\nu)\leftexp{(i_{1}...i_{l-m})}{Y}_{m,l-m}|\leq \sum_{k=0}^{2}|\slashed{d}
\leftexp{(i_{1}...i_{l-m})}{\rho}^{\prime(k)}_{m,l-m}|+2|\tilde{\chi}^{\prime}||\slashed{d}
R_{i_{l-m}}...R_{i_{1}}(T)^{m}\mu|
\end{align}
Now from (14.403) with $l+1$ replaced by $l$:
\begin{align}
 \slashed{d}\leftexp{(i_{1}...i_{l-m})}{\rho}^{\prime(0)}_{m,l-m}=\frac{1}{2}\ell\slashed{d}
R_{i_{l-m}}...R_{i_{1}}(T)^{m+1}\psi_{0}
\end{align}
hence:
\begin{align}
 \|\mu^{1/2}\slashed{d}\leftexp{(i_{1}...i_{l})}{\rho}^{\prime(0)}_{m,l-m}\|_{L^{2}(\Sigma_{t}^{\epsilon_{0}})}
\leq C|\ell|(1+t)^{-1}\sqrt{\sum_{\alpha}\mathcal{E}'_{1}[R_{i_{l-m}}...R_{i_{1}}(T)^{m+1}\psi_{\alpha}](t)}
\end{align}
Also, from (14.404) with $l+1$ replaced by $l$:
\begin{align}
 \slashed{d}\leftexp{(i_{1}...i_{l-m})}{\rho}^{\prime(1)}_{m,l-m}=(m^{0}_{T}-\frac{1}{2}\ell)
\slashed{d}R_{i_{l-m}}...R_{i_{1}}(T)^{m+1}\psi_{0}+m^{i}_{T}\slashed{d}
R_{i_{l-m}}...R_{i_{1}}(T)^{m+1}\psi_{i}\\\notag
+(1+t)^{-1}\mu e^{\alpha}_{L}\slashed{d}R_{i_{l-m}}...R_{i_{1}}(T)^{m}Q\psi_{\alpha}\\\notag
+(\slashed{d}m^{\alpha}_{T})(R_{i_{l-m}}...R_{i_{1}}(T)^{m+1}\psi_{\alpha})\\\notag
+(1+t)^{-1}(\slashed{d}(\mu e^{\alpha}_{L}))(R_{i_{l-m}}...R_{i_{1}}(T)^{m}Q\psi_{\alpha})
\end{align}
hence, in view of the facts that:
\begin{align}
 \sup_{\Sigma_{t}^{\epsilon_{0}}}|\slashed{d}m^{\alpha}_{T}|\leq C\delta_{0}(1+t)^{-2},\quad 
\sup_{\Sigma_{t}^{\epsilon_{0}}}|\slashed{d}e^{\alpha}_{L}|\leq C\delta_{0}(1+t)^{-1}
\end{align}
we obtain:
\begin{align}
 \|\mu^{1/2}\slashed{d}\leftexp{(i_{1}...i_{l-m})}{\rho'}^{(1)}_{m,l-m}\|_{L^{2}(\Sigma_{t}^{\epsilon_{0}})}\\\notag
\leq C\delta_{0}(1+t)^{-2}\sqrt{\sum_{\alpha}\mathcal{E}'_{1}[R_{i_{l-m}}...R_{i_{1}}(T)^{m+1}\psi_{\alpha}](t)}\\\notag
+C\delta_{0}(1+t)^{-2}[1+\log(1+t)]\sqrt{\sum_{\alpha}\mathcal{E}'_{1}[R_{i_{l-m}}...R_{i_{1}}(T)^{m}Q\psi_{\alpha}](t)}
\end{align}
Note that (14.465) and (14.468) are compatible with (14.416) and (14.417) respectively, and from (14.413) we have:
\begin{align}
 \||\slashed{d}\leftexp{(i_{1}...i_{l-m})}{\rho}^{\prime(2)}_{m,l-m}|+2|\tilde{\chi}^{\prime}|
|\slashed{d}R_{i_{l-m}}...R_{i_{1}}(T)^{m}\mu|\|_{L^{2}(\Sigma_{t}^{\epsilon_{0}})}
\\\notag
\leq C_{l}\delta_{0}(1+t)^{-2}[1+\log(1+t)]\{\mathcal{W}_{\{l+1\}}+\mathcal{W}^{Q}_{\{l\}}\\\notag
+(1+t)^{-1}\mathcal{B}_{[m,l+1]}+\mathcal{Y}_{0}+(1+t)\mathcal{A}'_{[l]}\}
\end{align}
 
    The $borderline$ contribution to (14.457) is the contribution from 
$\slashed{d}\leftexp{(i_{1}...i_{l-m})}{\rho}^{\prime(0)}_{m,l-m}$. This is the $borderline$ integral:
\begin{align}
 C\int_{W^{t}_{u}}(1+t')^{2}|T\psi_{\alpha}||\slashed{d}R_{i_{l-m}}...R_{i_{1}}(T)^{m+1}
\psi_{\alpha}||\slashed{d}\leftexp{(i_{1}...i_{l-m})}{\rho}^{\prime(0)}_{m,l-m}|dt'dud\mu_{\slashed{g}}\\\notag
\leq C\int_{0}^{t}(1+t')^{2}\sup_{\Sigma_{t'}^{\epsilon_{0}}}(\mu^{-1}|T\psi_{\alpha}|)
\|\mu^{1/2}\slashed{d}R_{i_{l-m}}...R_{i_{1}}(T)^{m+1}\psi_{\alpha}\|_{L^{2}(\Sigma_{t'}^{\epsilon_{0}})}\\\notag
\cdot\|\mu^{1/2}\slashed{d}\leftexp{(i_{1}...i_{l-m})}{\rho}^{\prime(0)}_{m,l-m}\|_{L^{2}(\Sigma_{t'}^{\epsilon_{0}})}dt'
\end{align}
Here we substitute (14.89) for $\sup_{\Sigma_{t'}^{\epsilon_{0}}}(\mu^{-1}|T\psi_{\alpha}|)$ and \\
(14.465) for $\|\mu^{1/2}\slashed{d}\leftexp{(i_{1}...i_{l-m})}
{\rho'}^{(0)}_{m,l-m}\|_{L^{2}(\Sigma_{t'}^{\epsilon_{0}})}$. The factors $|\ell|$ then cancel. 
Now the contribution of the second term on the right in (14.89) is in fact not borderline. 
The actual $borderline$ contribution is:
\begin{align}
 C\int_{0}^{t}\sup_{\Sigma_{t'}^{\epsilon_{0}}}(\mu^{-1}|L\mu|)\sqrt{\mathcal{E}'_{1}
[R_{i_{l-m}}...R_{i_{1}}(T)^{m+1}\psi_{\alpha}](t')}\\\notag
\cdot\sqrt{\sum_{\alpha}\mathcal{E}'_{1}[R_{i_{l-m}}...R_{i_{1}}(T)^{m+1}\psi_{\alpha}](t')}dt'\\\notag
\leq C(\frac{1}{2a}+\frac{1}{2q})\bar{\mu}^{-2a}_{m}(t)[1+\log(1+t)]^{2q}\leftexp{(i_{1}...i_{l-m})}
{\mathcal{G}}^{\prime}_{1,m,l+2;a,q}(t)
\end{align}
(similar to (14.360) and (14.361)).

    We proceed to consider the contribution of $\slashed{d}\leftexp{(i_{1}...i_{l-m})}{\rho}^{\prime(1)}_{m,l-m}$ to 
(14.457). Here we simply appeal to (14.90) for $\sup_{\Sigma_{t'}^{\epsilon_{0}}}(\mu^{-1}|T\psi_{\alpha}|)$. 
Substituting also (14.468) for\\ $\|\mu^{1/2}\slashed{d}\leftexp{(i_{1}...i_{l-m})}{\rho'}^{(1)}_{m,l-m}\|
_{L^{2}(\Sigma_{t'}^{\epsilon_{0}})}$, we obtain that the contribution in question is bounded by:
\begin{align}
 C\int_{W^{t}_{\epsilon_{0}}}(1+t')^{2}|T\psi_{\alpha}||\slashed{d}R_{i_{l-m}}...R_{i_{1}}(T)^{m+1}\psi_{\alpha}|
|\slashed{d}\leftexp{(i_{1}...i_{l-m})}{\rho'}^{(1)}_{m,l-m}|dt'dud\mu_{\slashed{g}}\\\notag
\leq C\delta_{0}\int_{0}^{t}(1+t')\bar{\mu}^{-1}_{m}(t')\|\mu^{1/2}\slashed{d}R_{i_{l-m}}...R_{i_{1}}(T)^{m+1}
\psi_{\alpha}\|_{L^{2}(\Sigma_{t'}^{\epsilon_{0}})}\\\notag
\cdot\|\mu^{1/2}\slashed{d}\leftexp{(i_{1}...i_{l-m})}{\rho'}^{(1)}_{m,l-m}\|_{L^{2}(\Sigma_{t'}^{\epsilon_{0}})}
dt'\\\notag\leq C\delta_{0}^{2}\bar{\mu}^{-2a}_{m}(t)[1+\log(1+t)]^{2q}\\\notag
\cdot\{\varphi_{2}(Ca\delta_{0})\mathcal{G}'_{1,[l+2];a,q}(t)+\int_{0}^{t}(1+t')^{-2}[1+\log(1+t')]
\mathcal{G}'_{1,[l+2];a,q}(t')dt'\}
\end{align}
(similar to (14.362)).

     Next, we consider $V'_{1,0}$, given by (14.452). By (14.429) we have:
\begin{align}
 |V'_{1,0}|\leq C\delta_{0}\int_{W^{t}_{\epsilon_{0}}}|R_{i_{l-m}}...R_{i_{1}}(T)^{m+1}\psi_{\alpha}|
|(\slashed{\mathcal{L}}_{L}+2\nu)\leftexp{(i_{1}...i_{l-m})}{Y}_{m,l-m}|dt'dud\mu_{\slashed{g}}\\\notag
\leq C\delta_{0}\int_{0}^{t}\mathcal{W}_{\{l+1\}}(t')\|(\slashed{\mathcal{L}}_{L}+2\nu)
\leftexp{(i_{1}...i_{l-m})}{Y}_{m,l-m}\|_{L^{2}(\Sigma_{t'}^{\epsilon_{0}})}dt'dud\mu_{\slashed{g}}
\end{align}
Here we need only consider the contribution of $\slashed{d}\leftexp{(i_{1}...i_{l-m})}{\rho}^{\prime(0)}_{m,l-m}$.
 By (14.464),
\begin{align}
 \|\slashed{d}\leftexp{(i_{1}...i_{l-m})}{\rho}^{\prime(0)}_{m,l-m}\|_{L^{2}(\Sigma_{t}^{\epsilon_{0}})}
\leq C\bar{\mu}^{-1/2}_{m}(1+t)^{-1}\sqrt{\mathcal{E}'_{1,[l+2]}(t)}
\end{align}
hence by (14.318) the contribution in question is bounded by 
\begin{align}
 C\epsilon_{0}\delta_{0}\int_{0}^{t}(1+t')^{-1}\bar{\mu}^{-2a-1/2}_{m}(t')[1+\log(1+t')]^{p+q}
\sqrt{\mathcal{G}_{0,[l+2];a,p}(t')\mathcal{G}'_{1,[l+2];a,q}(t')}dt'
\end{align}
This is similar to (14.364), which has already been bounded by (14.366). This completes the estimates for the 
spacetime integrals $V^{\prime}_{0}, V^{\prime}_{1}, V^{\prime}_{2}$. Thus, the estimate of the contribution 
of (14.57) to (14.59) is completed.
\chapter{The Top Order Energy Estimates}

In this chapter, we consider the energy estimates of top order $l+2$. Of these the most delicate 
are those corresponding to the variations (14.52) and (14.54), to which the $borderline$ error integrals are associated. 

\section{Estimates Associated to $K_{1}$}
For each of the variations, 
we first consider the integral identity associated to $K_{1}$ and then the integral identity associated to $K_{0}$. For each variation $\psi$, of any order,
there is an  additional $borderline$ integral associated to $K_{1}$, contributed by $Q_{1,3}[\psi]$ (see (5.114), (5.190)-(5.194)). This is the integral 
$L(t,\epsilon_{0})[\psi]$ defined by (5.257):
\begin{align}
 L(t,\epsilon_{0})[\psi]=\int_{0}^{t}(1+t')^{-1}[1+\log(1+t')]^{-1}\mathcal{E}'_{1}[\psi](t')dt'
\end{align}
     Consider now the integral identities corresponding to the variations (14.52) and the vectorfield $K_{1}$. In each of these we have the $borderline$ 
hypersurface integral (14.248) bounded by (14.268) and (14.295):
\begin{align}
 C(\frac{1}{a-1/2}+\frac{1}{q+1/2})\bar{\mu}^{-2a}_{m}(t)[1+\log(1+t)]^{2q}\leftexp{(i_{1}...i_{l})}{\mathcal{G}}^{\prime}_{1,l+2;a,q}(t)
\end{align}
We also have the $borderline$ spacetime integral (14.361) bounded by (14.363):
\begin{align}
 C(\frac{1}{2a}+\frac{1}{2q})\bar{\mu}^{-2a}_{m}(t)[1+\log(1+t)]^{2q}\leftexp{(i_{1}...i_{l})}{\mathcal{G}}^{\prime}_{1,l+2;a,q}(t)
\end{align}
We also have the remaining integrals, bounded, in the case of (14.308) by:
\begin{align}
 C_{q}\delta_{0}^{2}[1+\log(1+t)]^{2q}\bar{\mu}^{-2a}_{m}(t)\mathcal{G}'_{1,[l+2];a,q}
\end{align}
and in the case of (14.313), by:
\begin{align}
 C_{l}\epsilon_{0}\delta_{0}\bar{\mu}^{-2a}_{m}(t)[1+\log(1+t)]^{2q}\{\mathcal{G}_{0,[l+2];a,p}(t)+\mathcal{G}'_{1,[l+2];a,q}(t)\}
\end{align}
provided that:
\begin{align}
 q\geq p+1
\end{align}
and (14.321) is the same as (14.313). In the case of (14.365) we have a bound by:
\begin{align}
 C\delta_{0}^{2}\bar{\mu}^{-2a}_{m}(t)[1+\log(1+t)]^{2q}\mathcal{G}'_{1,[l+2];a,q}(t)
\end{align}
and in the case of (14.368) a bound by:
\begin{align}
 C_{l}\epsilon_{0}\delta_{0}\bar{\mu}^{-2a}_{m}(t)[1+\log(1+t)]^{2q}\{\mathcal{G}'_{1,[l+2];a,q}(t)+\mathcal{G}_{0,[l+2];a,p}(t)\}
\end{align}
provided that (15.6) holds. Finally (14.370) is the same as (14.368). Combining we see that all the remaining integrals associated to (14.52) and to $K_{1}$ 
containing the top order spatial derivatives of the acoustical entities are bounded by:
\begin{align}
 C_{q,l}\delta_{0}\bar{\mu}^{-2a}_{m}(t)[1+\log(1+t)]^{2q}\sqrt{\mathcal{G}'_{1,[l+2];a,q}(t)+\mathcal{G}_{0,[l+2];a,p}(t)}\\\notag
\cdot\{\delta_{0}\mathcal{Y}_{0}(0)+\mathcal{A}'_{[l]}(0)+\mathcal{B}_{[l+1]}(0)+\sqrt{\mathcal{G}'_{1,[l+2];a,q}(t)+\mathcal{G}_{0,[l+2];a,p}(t)}\}
\end{align}
The bounds in terms of initial data come from the initial hypersurface integrals and the lower order terms which we omitted.

     Consider now the error integrals associated to the same variations, or in fact any of the top order variations, and to $K_{1}$ or to $K_{0}$, which 
contain the lower order spatial derivatives of the acoustical entities and are contributed by the remaining terms in the sum (14.14). Consider an arbitrary 
term in this sum, corresponding to some $k\in\{0,...,l\}$:
\begin{align}
 (Y_{I_{l+1}}+\leftexp{(Y_{I_{l+1}})}{\delta})...(Y_{I_{l+2-k}}+\leftexp{(Y_{I_{l+2-k}})}{\delta})\leftexp{(Y_{I_{l+1-k}};I_{1}...I_{l-k})}{\sigma}_{l+1-k}
\end{align}
There is a total of $k$ derivatives with respect to the commutation fields acting on $\leftexp{(Y)}{\sigma}_{l+1-k}$. In view of the fact that 
$\leftexp{(Y)}{\sigma}_{l+1-k}$ has the structure described in the paragraph following (14.15), in considering the partial contribution of each term 
in $\leftexp{(Y)}{\sigma}_{1,l+1-k}$, if the factor which is a component of $\leftexp{(Y)}{\tilde{\pi}}$ receives more than $(l+1)_{*}$ derivatives with
 respect to the commutation fields, then the factor which is a 2nd derivative of $\psi_{l+1-k}$ receives at most $k-(l+1)_{*}-1$ derivatives of the the 
commutation fields, thus corresponds to a derivative of the $\psi_{\alpha}$ of order at most:
\begin{align*}
 k-(l+1)_{*}+1+l-k=l_{*}+1
\end{align*}
therefore this factor is bounded in $L^{\infty}(\Sigma_{t}^{\epsilon_{0}})$ by the bootstrap assumption. Also, in considering the partial contribution of each term in 
$\leftexp{(Y)}{\sigma}_{2,l+1-k}$, if the factor which is a 1st derivative of $\leftexp{(Y)}{\tilde{\pi}}$ receives more than $(l+1)_{*}-1$ derivatives with respect to the 
commutation fields, then the factor which is a 1st derivative of $\psi_{l+1-k}$ receives at most $k-(l+1)_{*}$ derivatives with respect to the commutation field, thus 
corresponds to a derivative of the $\psi_{\alpha}$ of order at most:
\begin{align*}
 k-(l+1)_{*}+1+l-k=l_{*}+1
\end{align*}
therefore this factor is again bounded in $L^{\infty}(\Sigma_{t}^{\epsilon_{0}})$ by the bootstrap assumption. Similar considerations apply to $\leftexp{(Y)}{\sigma}
_{3,l+1-k}$. We conclude that for all the terms in the sum in (14.14) of which one factor is a derivative of the $\leftexp{(Y)}{\tilde{\pi}}$ of order more than 
$(l+1)_{*}$, the other factor is then a derivative of the $\psi_{\alpha}$ of order at most $l_{*}+1$ and is thus bounded in $L^{\infty}(\Sigma_{t}^{\epsilon_{0}})$ by 
the bootstrap assumption. Of these terms we have already estimated the contribution of those containing the top order spatial derivatives of the acoustical entities. 
The contributions of the remaining terms are then bounded using Proposition 12.11 and 12.12. We obtain a bound for these contributions to the error integrals associated 
to $K_{1}$ and to any of the variations, up to the top order, by (recall (14.59)):
\begin{align}
 C_{l}\delta_{0}\bar{\mu}^{-a}_{m}(t)[1+\log(1+t)]^{q}\cdot\\\notag
\{\delta_{0}\mathcal{Y}_{0}(0)+\mathcal{A}'_{[l]}(0)+\mathcal{B}_{[l+1]}(0)+\sqrt{\mathcal{G}'_{1,[l+2];a,q}(t)+\mathcal{G}_{0,[l+2];a,p}(t)}\}\\\notag
\cdot\{\int_{0}^{\epsilon_{0}}\mathcal{F}'^{t}_{1,[l+2]}(u')du'\}^{1/2}\\\notag
+ C_{l}\delta_{0}\int_{0}^{\epsilon_{0}}\mathcal{F}'^{t}_{1,[l+2]}(u')du'+C_{l}\delta_{0}\{\int_{0}^{\epsilon_{0}}\mathcal{F}'^{t}_{1,[l+2]}(u')du'\}^{1/2}
\{K_{[l+2]}(t,\epsilon_{0})\}^{1/2}
\end{align}

    The last two terms contributed by the terms in (14.14) containing the top order derivatives of the acoustical entities of which however one is a derivative with 
respect to $L$ (hence are expressible in terms of the top order derivatives of the $\psi_{\alpha}$). 

Here, we have defined:
\begin{align}
 \mathcal{F}'^{t}_{1,[n]}(u)=\sum_{m=1}^{n}\mathcal{F}'^{t}_{1,m}(u)
\end{align}
where $\mathcal{F}^{\prime t}_{1,n}(u)$ represents the sum of the fluxes associated to the vectorfield $K_{1}$ of all the $n$th order variations.

    On the other hand, all the other terms in the sum (14.14) contain derivatives of the $\leftexp{(Y)}{\tilde{\pi}}$ of order at most $(l+1)_{*}$,
thus spatial derivatives of $\chi'$ of order at most $(l+1)_{*}$ and spatial derivatives of $\mu$ of order at most $(l+1)_{*}+1$, which are bounded in 
$L^{\infty}(\Sigma_{t}^{\epsilon_{0}})$ by virtue of Proposition 12.9 and 12.10 and the bootstrap assumption. So we can use Lemma 7.6 to bound these 
contributions by:
\begin{align}
 C_{l}\int_{0}^{\epsilon_{0}}\mathcal{F}'^{t}_{1,[l+2]}(u')du'+C_{l}\{\int_{0}^{\epsilon_{0}}\mathcal{F}'^{t}_{1,[l+2]}(u')du'\}^{1/2}
\{K_{[l+1]}(t,\epsilon_{0})\}^{1/2}\\\notag
+C_{l}\{\int_{0}^{\epsilon_{0}}\mathcal{F}'^{t}_{1,[l+2]}(u')du'\}^{1/2}\cdot\\\notag
\{\int_{0}^{t}(1+t')^{-2}[1+\log(1+t')]^{2}(\mathcal{E}'_{1,[l+2]}(t')+\epsilon_{0}^{2}\mathcal{E}_{0,[l+2]}(t'))dt'\}^{1/2}
\end{align}
 Here we have defined:
\begin{align}
 K_{[n]}(t,u)=\sum_{m=1}^{n}K_{m}(t,u)
\end{align}
where $K_{n}(t,u)$ represents the sum of the integrals
\begin{align}
 K[\psi](t,u)=-\int_{W^{t}_{u}}\frac{\Omega}{2}\omega\nu^{-1}\mu^{-1}(L\mu)_{-}|\slashed{d}\psi|^{2}d\mu_{g}
\end{align}
of all the $n$th order variations.

     Finally, for each variation $\psi$ we have the error integrals:
\begin{align}
 \int_{W^{t}_{\epsilon_{0}}}\sum_{k=1}^{7}Q_{0,k}[\psi]d\mu_{g},\quad \int_{W^{t}_{\epsilon_{0}}}\sum_{k=1}^{8}Q_{1,k}[\psi]d\mu_{g}
\end{align}
from the fundamental energy estimates. According to (5.277), we have:
\begin{align}
 \int_{W^{t}_{\epsilon_{0}}}\sum_{k=1}^{8}Q_{1,k}[\psi]d\mu_{g}\leq -\frac{1}{2}K[\psi](t,\epsilon_{0})+CM[\psi](t,\epsilon_{0})
+\frac{3}{2}L[\psi](t,\epsilon_{0})+\int_{0}^{t}\tilde{A}(t')\mathcal{E}'_{1}[\psi](t')dt'
\end{align}
Here $M[\psi](t,\epsilon_{0})$ is given by (5.274):
\begin{align}
 M[\psi](t,\epsilon_{0})=\bar{\mathcal{E}}_{0}[\psi](t)[1+\log(1+t)]^{4}+\int_{0}^{\epsilon_{0}}\mathcal{F}'^{t}_{1}[\psi](u')du'
\end{align}
As in Chapter 5, we will bring the term $-(1/2)K[\psi](t,\epsilon_{0})$ on the left-hand side of the integral inequality associated to $\psi$ and $K_{1}$.

     Let us define, for any variation $\psi$, in analogy with (14.109) and (14.111),
\begin{align}
 \mathcal{G}_{0;a,p}[\psi](t)=\sup_{t'\in[0,t]}\{[1+\log(1+t')]^{-2p}\bar{\mu}^{2a}_{m}(t')\mathcal{E}_{0}[\psi](t')\}\\
\mathcal{G}'_{1;a,q}[\psi](t)=\sup_{t'\in[0,t]}\{[1+\log(1+t')]^{-2q}\bar{\mu}^{2a}_{m}(t')\mathcal{E}'_{1}[\psi](t')\}
\end{align}
Note that $\mathcal{G}_{0;a,p}[\psi](t)$ and $\mathcal{G}'_{1;a,q}[\psi](t)$ are non-decreasing functions of $t$.

     Recalling from Chapter 8 the definition of $\bar{\mu}_{m}$:
\begin{align}
 \bar{\mu}_{m}(t)=\min\{\mu_{m}(t),1\},\quad \mu_{m}(t)=\min_{\Sigma_{t}^{\epsilon_{0}}}\mu=
\min_{(u,\vartheta)\in[0,\epsilon_{0}]\times S^{2}}\mu(t,u,\vartheta)
\end{align}
we define by replacing $[0,\epsilon_{0}]\times S^{2}$ by $[0,u]\times S^{2}$, or $\Sigma_{t}^{\epsilon_{0}}$ by $\Sigma_{t}^{u}$:
\begin{align}
 \bar{\mu}_{m,u}(t)=\min\{\mu_{m,u}(t),1\},\quad \mu_{m,u}(t)=\min_{\Sigma_{t}^{u}}\mu=\min_{(u',\vartheta)\in[0,u]\times S^{2}}
\mu(t,u',\vartheta)
\end{align}
Note that $\bar{\mu}_{m,u}(t)$ is a non-increasing function of $u$ at each $t$. We then define:
\begin{align}
 \mathcal{H}'^{t}_{1;a,q}[\psi](u)=\sup_{t'\in[0,t]}\{[1+\log(1+t')]^{-2q}\bar{\mu}^{2a}_{m,u}(t')\mathcal{F}'^{t'}_{1}[\psi](u)\}\\
V'_{1;a,q}[\psi](t,u)=\int_{0}^{u}\mathcal{H}'_{1;a,q}[\psi](u')du'
\end{align}
Note that $\mathcal{H}'^{t}_{1;a,q}[\psi](u)$ is a non-decreasing function of $t$ at each $u$ while $V'_{1;a,q}[\psi](t,u)$ is a non-decreasing function
of $u$ at each $t$ as well as a non-decreasing function of $t$ at each $u$. We have, for each $u\in[0,\epsilon_{0}]$:
\begin{align}
 \int_{0}^{u}\mathcal{F}_{1}^{\prime t}[\psi](u')du'\leq \int_{0}^{u}\bar{\mu}^{-2a}_{m,u'}(t)[1+\log(1+t)]^{2q}\mathcal{H}^{\prime t}_{1;a,q}[\psi](u')du'\\\notag
\leq\bar{\mu}^{-2a}_{m,u}(t)[1+\log(1+t)]^{2q}V'_{1;a,q}[\psi](t,u)
\end{align}
We also define:
\begin{align}
 \mathcal{H}'^{t}_{1,[n];a,q}(u)=\sup_{t'\in[0,t]}\{[1+\log(1+t')]^{-2q}\bar{\mu}^{2a}_{m,u}(t')\mathcal{F}'^{t'}_{1,[n]}(u)\}\\
V'_{1,[n];a,q}(t,u)=\int_{0}^{u}\mathcal{H}'^{t}_{1,[n];a,q}(u')du'
\end{align}
Then $\mathcal{H}'^{t}_{1,[n];a,q}(u)$ is a non-decreasing function function of $t$ at each $u$, $V'_{1,[n];a,q}(t,u)$ is a non-decreasing function of $u$ 
at each $t$ as well as a non-decreasing function of $t$ at each $u$, and for each $u\in[0,\epsilon_{0}]$:
\begin{align}
 \int_{0}^{u}\mathcal{F}'^{t}_{1,[n]}(u')du'\leq\int_{0}^{u}\bar{\mu}^{-2a}_{m,u'}(t)[1+\log(1+t)]^{2q}\mathcal{H}'^{t}_{1,[n];a,q}(u')du'\\\notag
\leq\bar{\mu}^{-2a}_{m,u}(t)[1+\log(1+t)]^{2q}V'_{1,[n];a,q}(t,u)
\end{align}
Going back to (15.18), we have, in regard to the first term on the right,
\begin{align}
 \bar{\mathcal{E}}_{0}[\psi](t)=\sup_{t'\in[0,t]}\mathcal{E}_{0}[\psi](t')
\leq\sup_{t'\in[0,t]}\{\bar{\mu}_{m}^{-2a}(t')[1+\log(1+t')]^{2p}\mathcal{G}_{0;a,p}[\psi](t')\}\\\notag
\leq C[1+\log(1+t)]^{2p}\bar{\mu}^{-2a}_{m}(t)\mathcal{G}_{0;a,p}[\psi](t)
\end{align}
(by Corollary 2 of Lemma 8.11). At this point we set:
\begin{align}
 q=p+2
\end{align}
which is consistent with (15.6). Then (15.29) implies that the first term on the right in (15.18) is bounded by:
\begin{align}
 C\bar{\mu}^{-2a}_{m}(t)[1+\log(1+t)]^{2q}\mathcal{G}_{0;a,p}[\psi](t)
\end{align}
The second term on the right in (15.18) being bounded according to (15.25) with $u=\epsilon_{0}$, we conclude that:
\begin{align}
 M[\psi](t,\epsilon_{0})\leq\bar{\mu}^{-2a}_{m}(t)[1+\log(1+t)]^{2q}\{C\mathcal{G}_{0;a,p}[\psi](t)+V'_{1;a,q}[\psi](t,\epsilon_{0})\}\\\notag
\leq\bar{\mu}^{-2a}_{m}(t)[1+\log(1+t)]^{2q}\{C\mathcal{G}_{0,[l+2];a,p}(t)+V'_{1,[l+2];a,q}(t,\epsilon_{0})\}
\end{align}
for all variations of order not exceeding $l+2$.

     In regard to the third term on the right in (15.17), which involves the borderline integral (15.1), we have, by Corollary 2 of Lemma 8.11,
\begin{align}
 L[\psi](t,\epsilon_{0})\leq C\bar{\mu}^{-2a}_{m}(t)\mathcal{G}'_{1;a,q}[\psi](t)\int_{0}^{t}(1+t')^{-1}[1+\log(1+t')]^{2q-1}dt'\\\notag
\leq \frac{C}{2q}\bar{\mu}^{-2a}_{m}(t)[1+\log(1+t)]^{2q}\mathcal{G}'_{1;a,q}[\psi](t)
\end{align}
Finally, in regard to the last term on the right in (15.17) we have:
\begin{align}
 \int_{0}^{t}\tilde{A}(t')\mathcal{E}'_{1}[\psi](t')dt'\leq C\bar{\mu}^{-2a}_{m}[1+\log(1+t)]^{2q}\int_{0}^{t}\tilde{A}(t')\mathcal{G}'_{1;a,q}[\psi](t')dt'\\\notag
\leq C\bar{\mu}^{-2a}_{m}(t)[1+\log(1+t)]^{2q}\int_{0}^{t}\tilde{A}(t')\mathcal{G}'_{1,[l+2];a,q}(t')dt'
\end{align}
for all variations of order not exceeding $l+2$.

      We also recall that in (5.73) (with $u=\epsilon_{0}$) corresponding to the variation $\psi$ and to $K_{1}$ we have the hypersurface integrals bounded 
according to (5.267):
\begin{align}
 |\int_{\Sigma_{t}^{\epsilon_{0}}}(1/2)\Omega(\underline{L}\omega+\underline{\nu}\omega)\psi^{2}-
\int_{\Sigma_{0}^{\epsilon_{0}}}(1/2)\Omega(\underline{L}\omega+\underline{\nu}\omega)\psi^{2}|\\\notag
\leq C\bar{\mathcal{E}}_{0}[\psi](t)[1+\log(1+t)]^{4}\\\notag
\leq C\bar{\mu}^{-2a}_{m}(t)[1+\log(1+t)]^{2q}\mathcal{G}_{0;a,p}[\psi](t)
\end{align}
by (15.29), (15.30).

      We turn to the bound (15.13). Let us define:
\begin{align}
 I_{[n];a,q}(t,u)=\sup_{t'\in[0,t]}\{[1+\log(1+t')]^{-2q}\bar{\mu}^{2a}_{m,u}(t')K_{[n]}(t',u)\}
\end{align}
This is a non-decreasing function of $t$ at each $u$. Then, in view of (15.28) with $u=\epsilon_{0}$, we can bound (15.13) by:
\begin{align}
 C_{l}\bar{\mu}^{-2a}_{m}(t)[1+\log(1+t)]^{2q}\cdot\\\notag
\{\int_{0}^{t}(1+t')^{-2}[1+\log(1+t')]^{2}(\mathcal{G}_{0,[l+2];a,p}(t')+\mathcal{G}'_{1;[l+2];a,q}(t'))dt'\\\notag
+V'_{1,[l+2];a,q}(t,\epsilon_{0})+I_{[l+2];a,q}(t,\epsilon_{0})^{1/2}V'_{1,[l+2];a,q}(t,\epsilon_{0})^{1/2}\}
\end{align}
and the last term in parenthesis is bounded by:
\begin{align}
 \frac{\delta}{2}I_{[l+2];a,q}(t,\epsilon_{0})+\frac{1}{2\delta}V'_{1,[l+2];a,q}(t,\epsilon_{0})
\end{align}
for any positive constant $\delta$ ($\delta$ will be chosen below).

      Finally, (15.11) is bounded by:
\begin{align}
 C_{l}\delta_{0}\bar{\mu}^{-2a}_{m}(t)[1+\log(1+t)]^{2q}\{(\mathcal{Y}_{0}(0)+\mathcal{A}'_{[l]}(0)+\mathcal{B}_{[l+1]}(0))^{2}\\\notag
+\delta(\mathcal{G}'_{1,[l+2];a,q}(t)+\mathcal{G}_{0,[l+2];a,p}(t))+\delta I_{[l+2];a,q}(t,\epsilon_{0})+(1+\frac{1}{\delta})V'_{1,[l+2];a,q}(t,\epsilon_{0})\}
\end{align}

     We now consider the integral identity corresponding to the vectorfield $K_{1}$ (5.73) with $u=\epsilon_{0}$ and to the variations (14.52) 
$R_{j}R_{i_{l}}...R_{i_{1}}\psi_{\alpha}$, $j=1,2,3$, $\alpha=0,1,2,3$, for a given multi-index $(i_{1}...i_{l})$. Summing over $j$ and $\alpha$,
we then obtain from (15.2), (15.3), (15.33), in regard to the borderline integrals, and from (15.9), (15.32), (15.34), (15.35), 
(15.37), (15.39), in regard to the remaining contributions, the following:
\begin{align}
 \{\sum_{j,\alpha}\mathcal{E}'_{1}[R_{j}R_{i_{l}}...R_{i_{1}}\psi_{\alpha}](t)+\sum_{j,\alpha}\mathcal{F}^{\prime t}_{1}[R_{j}R_{i_{1}}...R_{i_{1}}
\psi_{\alpha}](\epsilon_{0})\\\notag
+\frac{1}{2}\sum_{j,\alpha}K[R_{j}R_{i_{l}}...R_{i_{1}}\psi_{\alpha}](t,\epsilon_{0})\}\bar{\mu}^{2a}_{m}(t)[1+\log(1+t)]^{-2q}\\\notag
\leq C(\frac{1}{a}+\frac{1}{q})\leftexp{(i_{1}...i_{l})}{\mathcal{G}}'_{1,l+2;a,q}(t)+C_{q,l}\delta_{0}\mathcal{G}'_{1,[l+2];a,q}(t)\\\notag
+C_{q,l}R_{[l+2];a,q}(t,\epsilon_{0})+C\int_{0}^{t}\bar{A}(t')\mathcal{G}'_{1,[l+2];a,q}(t')dt'
\end{align}
Here,
\begin{align}
 R_{[l+2];a,q}(t,\epsilon_{0})=\mathcal{D}_{[l+2]}+\mathcal{G}_{0,[l+2];a,p}(t)+(1+\frac{1}{\delta})V'_{1,[l+2];a,q}(t,\epsilon_{0})
+\delta I_{[l+2];a,q}(t,\epsilon_{0})
\end{align}
and $\mathcal{D}_{[l+2]}$ is the quantity:
\begin{align}
 \mathcal{D}_{[l+2]}=\mathcal{E}_{0,[l+2]}(0)+\mathcal{E}'_{1,[l+2]}(0)+(\mathcal{P}_{[l+2]})^{2}\\\notag
\mathcal{P}_{[l+2]}=\mathcal{Y}_{0}(0)+\mathcal{A}'_{[l]}(0)+\mathcal{B}_{\{l+1\}}(0)\\\notag
+\sum_{i_{1}...i_{l}}\|\leftexp{(i_{1}...i_{l})}{x}_{l}(0)\|_{L^{2}(\Sigma_{0}^{\epsilon_{0}})}
+\sum_{m=0}^{l}\sum_{i_{1}...i_{l-m}}\|\leftexp{(i_{1}...i_{l-m})}{x}^{\prime}_{m,l-m}(0)\|_{L^{2}(\Sigma_{0}^{\epsilon_{0}})}
\end{align}
Also,
\begin{align}
 \bar{A}(t)=\tilde{A}(t)+C_{l}(1+t)^{-2}[1+\log(1+t)]^{2}
\end{align}
It is crucial that the constant $C$ in front of the first term on the right in (15.40) is independent of $a,q,l$.

     Now all the estimates we have made so far depend only on the bootstrap assumption and on assumption $\textbf{J}$, which refer to $W^{s}_{\epsilon_{0}}$,
on the assumptions of Proposition 12.6, 12.9, 12.10 on the initial data on $\Sigma_{0}^{\epsilon_{0}}$, and on the fact that $\epsilon_{0}\leq 1/2$. These
assumptions hold if $\epsilon_{0}$ is replaced by any $u\in(0,\epsilon_{0}]$; therefore all our estimates hold with $\epsilon_{0}$ replaced by any
$u\in(0,\epsilon_{0}]$. With this replacement, $\bar{\mu}_{m}(t)$ is replaced by $\bar{\mu}_{m,u}(t)$ and, for any variation $\psi$, $\mathcal{E}_{0}[\psi](t)$
and $\mathcal{E}'_{1}[\psi](t)$ are replaced by $\mathcal{E}^{u}_{0}[\psi](t)$ and $\mathcal{E}^{\prime u}_{1}[\psi](t)$ respectively. Thus $\mathcal{G}_{0;a,p}[\psi](t)$,
$\mathcal{G}'_{1;a,q}[\psi](t)$ are replaced by:
\begin{align}
 \mathcal{G}^{u}_{0;a,p}[\psi](t)=\sup_{t'\in[0,t]}\{[1+\log(1+t')]^{-2p}\bar{\mu}^{2a}_{m,u}(t')\mathcal{E}^{u}_{0}[\psi](t')\}\\
\mathcal{G}^{\prime u}_{1;a,q}[\psi](t)=\sup_{t'\in[0,t]}\{[1+\log(1+t')]^{-2q}\bar{\mu}^{2a}_{m,u}(t')\mathcal{E}^{\prime u}_{1}[\psi](t')\}
\end{align}
respectively, and $\leftexp{(i_{1}...i_{l})}{\mathcal{G}}_{0,l+2;a,p}(t)$,   $\leftexp{(i_{1}...i_{l})}{\mathcal{G}}_{0,m,l+2;a,p}(t)$,
$\leftexp{(i_{1}...i_{l})}{\mathcal{G}}^{\prime}_{1,l+2;a,p}(t)$, $\leftexp{(i_{1}...i_{l})}{\mathcal{G}}^{\prime}_{1,m,l+2;a,p}(t)$, are replaced by:
\begin{align}
 \leftexp{(i_{1}...i_{l})}{\mathcal{G}}^{u}_{0,l+2;a,p}(t)\\\notag
=\sup_{t'\in[0,t]}\{[1+\log(1+t')]^{-2p}\bar{\mu}^{2a}_{m,u}(t')\sum_{j,\alpha}\mathcal{E}^{u}_{0}[
R_{j}R_{i_{l}}...R_{i_{1}}\psi_{\alpha}](t')\}
\end{align}
\begin{align}
 \leftexp{(i_{1}...i_{l-m})}{\mathcal{G}}^{u}_{0,m,l+2;a,p}(t)\\\notag
=\sup_{t'\in[0,t]}\{[1+\log(1+t')]^{-2p}\bar{\mu}^{2a}_{m,u}(t')\sum_{\alpha}\mathcal{E}^{u}_{0}[R_{i_{l-m}}...R_{i_{1}}(T)^{m+1}\psi_{\alpha}](t')\}
\end{align}
\begin{align}
 \leftexp{(i_{1}...i_{l})}{\mathcal{G}}^{\prime u}_{1,l+2;a,q}(t)\\\notag
=\sup_{t'\in[0,t]}\{[1+\log(1+t')]^{-2q}\bar{\mu}^{2a}_{m,u}(t')\sum_{j,\alpha}\mathcal{E}^{\prime u}_{1}[
R_{j}R_{i_{l}}...R_{i_{1}}\psi_{\alpha}](t')\}
\end{align}
\begin{align}
 \leftexp{(i_{1}...i_{l-m})}{\mathcal{G}}^{\prime u}_{1,m,l+2;a,q}(t)\\\notag
=\sup_{t'\in[0,t]}\{[1+\log(1+t')]^{-2q}\bar{\mu}^{2a}_{m,u}(t')\sum_{\alpha}\mathcal{E'}^{u}_{1}[R_{i_{l-m}}...R_{i_{1}}(T)^{m+1}\psi_{\alpha}](t')\}
\end{align}
respectively. Also, $\mathcal{G}_{0,[n];a,p}$, $\mathcal{G}'_{1,[n];a,q}$, are replaced by:
\begin{align}
 \mathcal{G}^{u}_{0,[n];a,p}(t)=\sup_{t'\in[0,t]}\{[1+\log(1+t')]^{-2p}\bar{\mu}^{2a}_{m,u}(t')\mathcal{E}^{u}_{0,[n]}(t')\}\\
\mathcal{G}^{\prime u}_{1,[n];a,q}(t)=\sup_{t'\in[0,t]}\{[1+\log(1+t')]^{-2q}\bar{\mu}^{2a}_{m,u}(t')\mathcal{E}^{\prime u}_{1,[n]}(t')\}
\end{align}
respectively. Moreover, $\mathcal{D}_{[l+2]}$ is replaced by $\mathcal{D}^{u}_{[l+2]}$. Thus in place of (15.40) we obtain, for all $u\in[0,\epsilon_{0}]$,
and all $t\in[0,s]$:
\begin{align}
 \{\sum_{j,\alpha}\mathcal{E}^{\prime u}_{1}[R_{j}R_{i_{l}}...R_{i_{1}}\psi_{\alpha}](t)+\sum_{j,\alpha}\mathcal{F}^{\prime t}_{1}[R_{j}R_{i_{l}}...R_{i_{1}}
\psi_{\alpha}](u)\\\notag
+\frac{1}{2}\sum_{j,\alpha}K[R_{j}R_{i_{l}}...R_{i_{1}}\psi_{\alpha}](t,u)\}\bar{\mu}^{2a}_{m,u}(t)[1+\log(1+t)]^{-2q}\\\notag
\leq C(\frac{1}{a}+\frac{1}{q})\leftexp{(i_{1}...i_{l})}{\mathcal{G}}^{\prime u}_{1,l+2;a,q}(t)+C_{q,l}\delta_{0}\mathcal{G}^{\prime u}_{1,[l+2];a,q}(t)\\\notag
+C_{q,l}R_{[l+2];a,q}(t,u)+C\int_{0}^{t}\bar{A}(t')\mathcal{G}^{\prime u}_{1,[l+2];a,q}(t')dt'
\end{align}
and we have:
\begin{align}
 R_{[l+2];a,q}(t,u)=\mathcal{D}^{u}_{[l+2]}+\mathcal{G}^{u}_{0,[l+2];a,p}(t)+(1+\frac{1}{\delta})V'_{1,[l+2];a,q}(t,u)\\\notag
+\delta I_{[l+2];a,q}(t,u)
\end{align}
where:
\begin{align}
 \mathcal{D}^{u}_{[l+2]}=\mathcal{E}^{u}_{0,[l+2]}(0)+\mathcal{E}^{\prime u}_{1,[l+2]}(0)+(\mathcal{P}^{u}_{[l+2]})^{2}\\\notag
\mathcal{P}^{u}_{[l+2]}=\mathcal{Y}^{u}_{0}(0)+\mathcal{A}^{\prime u}_{[l]}(0)+\mathcal{B}^{u}_{\{l+1\}}(0)\\\notag
+\sum_{i_{1}...i_{l}}\|\leftexp{(i_{1}...i_{l})}{x}_{l}(0)\|_{L^{2}(\Sigma_{0}^{u})}+\sum_{m=0}^{l}\sum_{i_{1}...i_{l-m}}
\|\leftexp{(i_{1}...i_{l-m})}{x}^{\prime}_{m,l-m}(0)\|_{L^{2}(\Sigma_{0}^{u})}
\end{align}

     Keeping only the term
\begin{align*}
 \sum_{j,\alpha}\mathcal{E}^{\prime u}_{1}[R_{j}R_{i_{l}}...R_{i_{1}}\psi_{\alpha}](t)
\end{align*}
on the left of (15.52), we have:
\begin{align}
 \bar{\mu}^{2a}_{m,u}(t)[1+\log(1+t)]^{-2q}\sum_{j,\alpha}\mathcal{E}^{\prime u}_{1}[R_{j}R_{i_{l}}...R_{i_{1}}\psi_{\alpha}](t)\\\notag
\leq C(\frac{1}{a}+\frac{1}{q})\leftexp{(i_{1}...i_{l})}{\mathcal{G}}^{\prime u}_{1,l+2;a,q}(t)+C_{q,l}\delta_{0}\mathcal{G}^{\prime u}_{1,[l+2];a,q}(t)\\\notag
+C_{q,l}R_{[l+2];a,q}(t,u)+C\int_{0}^{t}\bar{A}(t')\mathcal{G}^{\prime u}_{1,[l+2];a,q}(t')dt'
\end{align}
The same holds with $t$ replaced by $t'\in[0,t]$, Now the right hand side of (15.55) is a non-decreasing function of $t$ at each $u$. The inequality corresponding
to $t^{\prime}$ thus holds a fortiori if we again replace $t^{\prime}$ by $t$ on the right hand side. Taking the supremum over all $t^{\prime}\in[0,t]$ 
on the left hand side we obtain, in view of the definition (15.48):
\begin{align}
 \leftexp{(i_{1}...i_{l})}{\mathcal{G}}^{\prime u}_{1,[l+2];a,q}(t)\leq C(\frac{1}{a}+\frac{1}{q})\leftexp{(i_{1}...i_{l})}{\mathcal{G}}^{\prime u}_{1,l+2;a,q}(t)
+C_{q,l}\delta_{0}\mathcal{G}^{\prime u}_{1,[l+2];a,q}(t)\\\notag+C_{q,l}R_{[l+2];a,q}(t,u)+C\int_{0}^{t}\bar{A}(t')\mathcal{G}^{\prime u}_{1,[l+2];a,q}(t')dt'
\end{align}
If $a$ and $q$ are chosen suitably large so that:
\begin{align}
 C(\frac{1}{a}+\frac{1}{q})\leq\frac{1}{2}
\end{align}
then (15.56) will imply:
\begin{align}
 \frac{1}{2}\leftexp{(i_{1}...i_{l})}{\mathcal{G}}^{\prime u}_{1,l+2;a,q}(t)\leq C_{q,l}\delta_{0}\mathcal{G}^{\prime u}_{1,[l+2];a,q}(t)\\\notag
+C_{q,l}R_{[l+2];a,q}(t,u)+C\int_{0}^{t}\bar{A}(t')\mathcal{G}^{\prime u}_{1,[l+2];a,q}(t')dt'
\end{align}
The actual choice of $a$ and $q$ shall be specified below.

     Consider now the integral identities corresponding to the variations (14.54) and $K_{1}$. In each of these we have the $borderline$ hypersurface integral
(14.423) bounded by (14.425):
\begin{align}
 C(\frac{1}{a-1/2}+\frac{1}{q+1/2})\bar{\mu}_{m}^{-2a}(t)[1+\log(1+t)]^{2q}\leftexp{(i_{1}...i_{l-m})}{\mathcal{G}}^{\prime}_{1,m,l+2;a,q}(t)
\end{align}
We also have the $borderline$ spacetime integral (14.472) bounded by (14.473):
\begin{align}
 C(\frac{1}{2a}+\frac{1}{2q})\bar{\mu}^{-2a}_{m}(t)[1+\log(1+t)]^{2q}\leftexp{(i_{1}...i_{l-m})}{\mathcal{G}}^{\prime}_{1,m,l+2;a,q}(t)
\end{align}
We also have the remaining integrals, bounded, in the case of (14.429) by (15.4), in the case of (14.435) bounded by (15.5), in the case of (14.474) by (15.7),
and in the case of (14.476) by (15.8). Combining we see that all the remaining integrals associated to the variations (14.54) and to $K_{1}$ containing the top
order spatial derivatives of the acoustical entities are bounded by (15.9). On the other hand, by the discussion following (15.9) and (15.13), the error integrals
associated to (14.54) and to $K_{1}$ contributed by all other terms in (14.14) are bounded by (15.11) and (15.13). Finally, we are left with (15.16), which is bounded
in a similar way. We thus obtain:
\begin{align}
 \{\sum_{\alpha}\mathcal{E}'_{1}[R_{i_{l-m}}...R_{i_{1}}(T)^{m+1}\psi_{\alpha}](t)+\sum_{\alpha}\mathcal{F}_{1}^{\prime t}[R_{i_{l-m}}...R_{i_{1}}(T)^{m+1}\psi_{\alpha}]
(\epsilon_{0})\\\notag
+\frac{1}{2}\sum_{j,\alpha}K[R_{i_{l-m}}...R_{i_{1}}(T)^{m+1}\psi_{\alpha}](t,\epsilon_{0})\}\bar{\mu}^{2a}_{m}(t)[1+\log(1+t)]^{-2q}\\\notag
\leq C(\frac{1}{a}+\frac{1}{q})\leftexp{(i_{1}...i_{l-m})}{\mathcal{G}}^{\prime}_{1,m,l+2;a,q}(t)+C_{q,l}\delta_{0}\mathcal{G}'_{1,[l+2];a,q}(t)\\\notag
+C_{q,l}R_{[l+2];a,q}(t,\epsilon_{0})+C\int_{0}^{t}\bar{A}(t')\mathcal{G}'_{1,[l+2];a,q}(t')dt'
\end{align}
and also:
\begin{align}
 \{\sum_{\alpha}\mathcal{E}^{\prime u}_{1}[R_{i_{l-m}}...R_{i_{1}}(T)^{m+1}\psi_{\alpha}](t)+\sum_{\alpha}\mathcal{F}^{\prime t}_{1}
[R_{i_{l-m}}...R_{i_{1}}(T)^{m+1}\psi_{\alpha}](u)\\\notag
+\frac{1}{2}\sum_{j,\alpha}K[R_{i_{l-m}}...R_{i_{1}}(T)^{m+1}\psi_{\alpha}](t,u)\}\bar{\mu}^{2a}_{m,u}(t)[1+\log(1+t)]^{-2q}\\\notag
\leq C(\frac{1}{a}+\frac{1}{q})\leftexp{(i_{1}...i_{l-m})}{\mathcal{G}}^{\prime u}_{1,m,l+2;a,q}(t)+C_{q,l}\delta_{0}\mathcal{G}^{\prime u}_{1,[l+2];a,q}(t)\\\notag
+C_{q,l}R_{[l+2];a,q}(t,u)+C\int_{0}^{t}\bar{A}(t')\mathcal{G}^{\prime u}_{1,[l+2];a,q}(t')dt'
\end{align}

      Keeping only the term 
\begin{align*}
 \sum_{\alpha}\mathcal{E}^{\prime u}_{1}[R_{i_{l-m}}...R_{i_{1}}(T)^{m+1}\psi_{\alpha}](t)
\end{align*}
on the left of (15.62) we have:
\begin{align*}
 \bar{\mu}^{2a}_{m,u}(t)[1+\log(1+t)]^{-2q}\sum_{\alpha}\mathcal{E}^{\prime u}_{1}[R_{i_{l-m}}...R_{i_{1}}(T)^{m+1}\psi_{\alpha}](t)\\
\leq C(\frac{1}{a}+\frac{1}{q})\leftexp{(i_{1}...i_{l-m})}{\mathcal{G}}^{\prime u}_{1,m,l+2;a,q}(t)+C_{q,l}\delta_{0}\mathcal{G}^{\prime u}_{1,[l+2];a,q}(t)\\
+C_{q,l}R_{[l+2];a,q}(t,u)+C\int_{0}^{t}\bar{A}(t^{\prime})\mathcal{G}^{\prime u}_{1,[l+2];a,q}(t^{\prime})dt^{\prime}
\end{align*}
By the same reasoning as that leading from (15.55) to (15.58) we then deduce:
\begin{align}
 \frac{1}{2}\leftexp{(i_{1}...i_{l-m})}{\mathcal{G}}^{\prime u}_{1,m,l+2;a,q}(t)\\\notag
\leq C_{q,l}\delta_{0}\mathcal{G}^{\prime u}_{1,[l+2];a,q}(t)+C_{q,l}R_{[l+2];a,q}(t,u)+C\int_{0}^{t}\bar{A}(t')\mathcal{G}^{\prime u}_{1,[l+2];a,q}(t')dt'
\end{align}

      Consider finally the integral identities corresponding to any of the other variations $\psi$ of order up to $l+2$ and to $K_{1}$. In each of these we
have only one $borderline$ integral (15.1) bounded by (15.33). Since all the remaining integrals in these identities are bounded by:
\begin{align}
 C_{q,l}\delta_{0}\mathcal{G}^{\prime u}_{1,[l+2];a,q}(t)+C_{q,l}R_{[l+2];a,q}(t,u)+C\int_{0}^{t}\bar{A}(t')\mathcal{G}^{\prime u}_{1,[l+2];a,q}(t')dt'
\end{align}
we obtain, for each such variation $\psi$, the inequality:
\begin{align}
 \{\mathcal{E}^{\prime u}_{1}[\psi](t)+\mathcal{F}^{\prime t}_{1}[\psi](u)+\frac{1}{2}K[\psi](t,u)\}\bar{\mu}^{2a}_{m,u}(t)[1+\log(1+t)]^{-2q}\\\notag
\leq \frac{C}{2q}\mathcal{G}^{\prime u}_{1;a,q}[\psi](t)+C_{q,l}\delta_{0}\mathcal{G}^{\prime u}_{1,[l+2];a,q}(t)\\\notag
+C_{q,l}R_{[l+2];a,q}(t,u)+C\int_{0}^{t}\bar{A}(t')\mathcal{G}^{\prime u}_{1,[l+2];a,q}(t')dt'
\end{align}
     Similarly, keeping only
\begin{align*}
 \mathcal{E}^{\prime u}_{1}[\psi](t)
\end{align*}
on the left and noting that by virtue of the condition (15.57) the coefficient $C/2q$ of the first term on the right is not greater than $1/2$, we have:
\begin{align*}
 \bar{\mu}^{2a}_{m,u}(t)[1+\log(1+t)]^{-2q}\mathcal{E}^{\prime u}_{1}[\psi](t)\\
\leq\frac{1}{2}\mathcal{G}^{\prime u}_{1;a,q}[\psi](t)+C_{q,l}\delta_{0}\mathcal{G}^{\prime u}_{1,[l+2];a,q}(t)\\
+C_{q,l}R_{[l+2];a,q}(t,u)+C\int_{0}^{t}\bar{A}(t^{\prime})\mathcal{G}^{\prime u}_{1,[l+2];a,q}(t^{\prime})dt^{\prime}
\end{align*}
The same holds with $t$ replaced by any $t^{\prime}\in[0,t]$. The right hand side being a non-decreasing function of $t$ at each $u$,
the inequality corresponding to $t^{\prime}$ holds fortiori if we again replace $t^{\prime}$ by $t$ on the right. Taking then the supremum
over all $t^{\prime}\in[0,t]$ on the left we obtain, in view of the definition (15.45),
\begin{align}
 \frac{1}{2}\mathcal{G}^{\prime u}_{1;a,q}[\psi](t)\leq C_{q,l}\delta_{0}\mathcal{G}^{\prime u}_{1,[l+2];a,q}(t)\\\notag
+C_{q,l}R_{[l+2];a,q}(t,u)+C\int_{0}^{t}\bar{A}(t')\mathcal{G'}^{u}_{1,[l+2];a,q}(t')dt'
\end{align}
    Now, according to the definition of $\mathcal{E}^{\prime u}_{1,[l+2]}(t)$ we have:
\begin{align}
 \mathcal{E}^{\prime u}_{1,[l+2]}(t)=\sum_{i_{1}...i_{l}}\{\sum_{j,\alpha}\mathcal{E}^{\prime u}_{1}[R_{j}R_{i_{l}}...R_{i_{1}}\psi_{\alpha}](t)\}\\\notag
+\sum_{m=0}^{l}\sum_{i_{1}...i_{l-m}}\{\sum_{\alpha}\mathcal{E}^{\prime u}_{1}[R_{i_{l-m}}...R_{i_{1}}(T)^{m+1}\psi_{\alpha}](t)\}\\\notag
+\sum_{\psi}\mathcal{E}^{\prime u}_{1}[\psi](t)
\end{align}
where the last sum is over all the other variations $\psi$ of order up to $l+2$. Noting that for non-negative functions
$x_{1}(t),...,x_{N}(t)$ we have:
\begin{align}
 \sup_{t'\in[0,t]}\sum_{n=1}^{N}x_{n}(t')\leq\sum_{n=1}^{N}\sup_{t'\in[0,t]}x_{n}(t')
\end{align}
it follows, in view of the definitions (15.45), (15.48), (15.49), (15.51) that:
\begin{align*}
 \mathcal{G}^{\prime u}_{1,[l+2];a,q}(t)\leq\sum_{i_{1}...i_{l}}\leftexp{(i_{1}...i_{l})}{\mathcal{G}}^{\prime u}_{1,l+2;a,q}(t)\\
+\sum_{m=0}^{l}\sum_{i_{1}...i_{l-m}}\leftexp{(i_{1}...i_{l-m})}{\mathcal{G}}^{\prime u}_{1,m,l+2;a,q}(t)\\
\sum_{\psi}\mathcal{G}^{\prime u}_{1;a,q}[\psi](t)
\end{align*}
Thus, summing inequalities (15.58) over $i_{1}...i_{l}$, summing inequalities (15.63) over $i_{1}...i_{l-m}$ and over $m=0,...,l$, and summing inequalities (15.66)
over all the other variations $\psi$ of order up to $l+2$,and then adding the resulting three inequalities we obtain:
\begin{align}
\frac{1}{2}\mathcal{G}^{\prime u}_{1,[l+2];a,q}(t)\leq C_{q,l}\delta_{0}\mathcal{G}^{\prime u}_{1,[l+2];a,q}(t)\\\notag
+C_{q,l}R_{[l+2];a,q}(t,u)+C_{l}\int_{0}^{t}\bar{A}(t')\mathcal{G}^{\prime u}_{1,[l+2];a,q}(t')dt' 
\end{align}
for new constants $C_{l}$, $C_{q,l}$. Requiring then $\delta_{0}$ to satisfy the smallness condition:
\begin{align}
 C_{q,l}\delta_{0}\leq\frac{1}{4}
\end{align}
(15.69) implies:
\begin{align}
 \mathcal{G}^{\prime u}_{1,[l+2];a,q}(t)\leq C_{q,l}R_{[l+2];a,q}(t,u)+C_{l}\int_{0}^{t}\bar{A}(t')\mathcal{G}^{\prime u}_{1,[l+2];a,q}(t')dt'
\end{align}
By (15.43), (5.258) and Proposition 13.1, we have:
\begin{align}
 \int_{0}^{s}\bar{A}(t)dt\leq C_{l}\quad:\quad\textrm{independent of}\quad s
\end{align}
Since $R_{[l+2];a,q}(t,u)$ is a non-decreasing function of $t$ at each $u$, (15.71) implies:
\begin{align}
 \mathcal{G}^{\prime u}_{1,[l+2];a,q}(t)\leq C_{q,l}R_{[l+2];a,q}(t,u)
\end{align}
hence also:
\begin{align}
 \int_{0}^{t}\bar{A}(t')\mathcal{G}^{\prime u}_{1,[l+2];a,q}(t')dt'\leq C_{q,l}R_{[l+2];a,q}(t,u)
\end{align}
for a new constant $C_{q,l}$.

      Substituting (15.73), (15.74) in the right-hand sides of (15.52), (15.62) and (15.65), omitting in each case the first term in parenthesis on the left, 
we obtain:
\begin{align}
 \{\sum_{j,\alpha}\mathcal{F}^{\prime t}_{1}[R_{j}R_{i_{l}}...R_{i_{1}}\psi_{\alpha}](u)+\frac{1}{2}\sum_{j,\alpha}
K[R_{j}R_{i_{l}}...R_{i_{1}}\psi_{\alpha}](t,u)\}\\\notag
\cdot\bar{\mu}^{2a}_{m,u}(t)[1+\log(1+t)]^{-2q}\\\notag
\leq C_{q,l}R_{[l+2];a,q}(t,u)
\end{align}
\begin{align}
 \{\sum_{\alpha}\mathcal{F}^{\prime t}_{1}[R_{i_{l-m}}...R_{i_{1}}(T)^{m+1}\psi_{\alpha}](u)+\frac{1}{2}\sum_{\alpha}
K[R_{i_{l-m}}...R_{i_{1}}(T)^{m+1}\psi_{\alpha}](t,u)\}\\\notag
\cdot\bar{\mu}^{2a}_{m,u}(t)[1+\log(1+t)]^{-2q}\\\notag
\leq C_{q,l}R_{[l+2];a,q}(t,u)
\end{align}
\begin{align}
 \{\mathcal{F}^{\prime t}_{1}[\psi](u)+\frac{1}{2}K[\psi](t,u)\}\bar{\mu}^{2a}_{m,u}(t)[1+\log(1+t)]^{-2q}\\\notag
\leq C_{q,l}R_{[l+2];a,q}(t,u)
\end{align}
Summing then as in (15.67), yields:
\begin{align}
 \{\mathcal{F}^{\prime t}_{1,[l+2]}(u)+\frac{1}{2}K_{[l+2]}(t,u)\}\bar{\mu}^{2a}_{m,u}(t)[1+\log(1+t)]^{-2q}\\\notag
\leq C_{q,l}R_{[l+2];a,q}(t,u)
\end{align}
for a new constant $C_{q,l}$. Now the definition (15.53) of $R_{[l+2];a,q}(t,u)$ reads:
\begin{align}
 R_{[l+2];a,q}(t,u)=\delta I_{[l+2];a,q}(t,u)+N_{[l+2];a,q}(t,u)
\end{align}
where:
\begin{align}
 N_{[l+2];a,q}(t,u)=\mathcal{D}^{u}_{[l+2]}+\mathcal{G}^{u}_{0,[l+2];a,p}(t)+(1+\frac{1}{\delta})V'_{1,[l+2];a,q}(t,u)
\end{align}

     Keeping only the term $(1/2)K_{[l+2]}(t,u)$ on the left in (15.78) we have:
\begin{align}
 \frac{1}{2}\bar{\mu}^{2a}_{m,u}(t)[1+\log(1+t)]^{-2q}K_{[l+2]}(t,u)\leq C_{q,l}\delta I_{[l+2];a,q}(t,u)+C_{q,l}N_{[l+2];a,q}(t,u)
\end{align}
The same holds with $t$ replaced by any $t^{\prime}\in[0,t]$. The right hand side being a non-decreasing function of $t$ at each $u$, the inequality
holds a fortiori if we again replace $t^{\prime}$ by $t$ on the right. Taking then the supremum over all $t^{\prime}\in[0,t]$ on the left we obtain, 
in view of the definition (15.36), 
\begin{align}
 \frac{1}{2}I_{[l+2];a,q}(t,u)\leq C_{q,l}\delta I_{[l+2];a,q}(t,u)+C_{q,l}N_{[l+2];a,q}(t,u)
\end{align}
We choose $\delta$ according to:
\begin{align}
 C_{q,l}\delta=\frac{1}{4}
\end{align}
Then (15.82) implies:
\begin{align}
 I_{[l+2];a,q}(t,u)\leq C_{q,l}N_{[l+2];a,q}(t,u)
\end{align}

    Substituting the estimate (15.84) in (15.79) and the result in (15.78), and keeping only $\mathcal{F}^{\prime t}_{1,[l+2]}(u)$ on the left of (15.78), 
we obtain:
\begin{align}
 \bar{\mu}^{2a}_{m,u}(t)[1+\log(1+t)]^{-2q}\mathcal{F}^{\prime t}_{1,[l+2]}(u)\leq C_{q,l}(Q_{[l+2];a,p}(t,u)+V'_{1,[l+2];a,q}(t,u))
\end{align}
for a new $C_{q,l}$, where:
\begin{align}
 Q_{[l+2];a,p}(t,u)=\mathcal{D}^{u}_{[l+2]}+\mathcal{G}^{u}_{0,[l+2];a,p}(t)
\end{align}
Seeing that the right hand side of (15.85) is a non-decreasing function of $t$ at each $u$, the previous reasoning applies and we deduce:
\begin{align}
 \mathcal{H}^{\prime t}_{1,[l+2];a,q}(u)\leq C_{q,l}(Q_{[l+2];a,p}(t,u)+V'_{1,[l+2];a,q}(t,u))
\end{align}
In view of the definition of $V'_{1,[l+2];a,q}(t,u)$, this reads:
\begin{align}
 \mathcal{H}^{\prime t}_{1,[l+2];a,q}(u)\leq C_{q,l}Q_{[l+2];a,p}(t,u)+C_{q,l}\int_{0}^{u}\mathcal{H}^{\prime t}_{1,[l+2];a,q}(u')du'
\end{align}
 Defining the functions:
\begin{align}
 \bar{\mathcal{G}}^{u}_{0,[n];a,p}(t)=\sup_{u'\in[0,u]}\mathcal{G}^{u'}_{0,[n];a,p}(t)
\end{align}
which are non-decreasing functions of $u$ at each $t$ as well as non-decreasing functions of $t$ at each $u$, we may replace $Q_{[l+2];a,p}(t,u)$
 on the right in (15.88) by:
\begin{align}
 \bar{Q}_{[l+2];a,p}(t,u)=\mathcal{D}^{u}_{[l+2]}+\bar{\mathcal{G}}^{u}_{0,[l+2];a,p}(t)
\end{align}
 which is also a non-decreasing function of $u$ at each $t$ as well as a non-decreasing function of $t$ at each $u$. Recalling that $[0,\epsilon_{0}]$ is a
bounded interval, (15.88) then implies:
\begin{align}
 \mathcal{H}^{\prime t}_{1,[l+2];a,q}(u)\leq C_{q,l}\bar{Q}_{[l+2];a,p}(t,u)
\end{align}
for a new constant $C_{q,l}$. Hence also:
\begin{align}
 V'_{1,[l+2];a,q}(t,u)\leq C_{q,l}\epsilon_{0}\bar{Q}_{[l+2];a,p}(t,u)
\end{align}
Substituting this in (15.80) and the result in (15.84) we obtain:
\begin{align}
 I_{[l+2];a,q}(t,u)\leq C_{q,l}\bar{Q}_{[l+2];a,p}(t,u)
\end{align}
for a new constant $C_{q,l}$. Substituting finally this in (15.79) and the result in (15.73) yields:
\begin{align}
 \mathcal{G}^{\prime u}_{1,[l+2];a,q}(t)\leq C_{q,l}\bar{Q}_{[l+2];a,p}(t,u)
\end{align}
 for a new constant $C_{q,l}$.

\section{Estimates Associated to $K_{0}$}
      Consider now the integral identities corresponding to the variations (14.52) and to the vectorfield $K_{0}$. In each of these we have the $borderline$ integral
(14.97), bounded by the sum of (14.102) and (14.107):
\begin{align}
 C(\frac{1}{2a}+\frac{1}{2p})\bar{\mu}^{-2a}_{m}(t)[1+\log(1+t)]^{2p}\leftexp{(i_{1}...i_{l})}{\mathcal{G}}_{0,l+2;a,p}(t)
\end{align}
We also have the remaining integrals, bounded in the case of (14.115) by (14.132), hence by:
\begin{align}
 C_{q}\delta_{0}\bar{\mu}^{-2a}_{m}(t)[1+\log(1+t)]^{2p}\{\mathcal{G}_{0,[l+2];a,p}(t)+\mathcal{G}^{\prime}_{1,[l+2];a,q}(t)\}
\end{align}
and in the case of (14.157) by (14.158), hence by (15.96) as well. Also, the integral (14.161), bounded by (14.162), hence by:
\begin{align}
 C_{p}\delta_{0}\bar{\mu}^{-2a}_{m}(t)[1+\log(1+t)]^{2p}\mathcal{P}_{[l+2]}\sqrt{\mathcal{G}_{0,[l+2];a,p}(t)}
\end{align}
Combining we see that all the remaining integrals associated to (14.52) and to $K_{0}$ containing the top order spatial derivatives of the acoustical entities are
bounded by:
\begin{align}
 C_{p,l}\delta_{0}\bar{\mu}^{-2a}_{m}(t)[1+\log(1+t)]^{2p}\sqrt{\mathcal{G}_{0,[l+2];a,p}(t)+\mathcal{G}^{\prime}_{1,[l+2];a,q}(t)}\cdot\\\notag
\{\mathcal{P}_{[l+2]}+\sqrt{\mathcal{G}_{0,[l+2];a,p}(t)+\mathcal{G}^{\prime}_{1,[l+2];a,q}(t)}\}
\end{align}

     By the discussion following (15.9), the contribution of the terms in (14.14) of which one factor is a derivative of $\leftexp{(Y)}{\tilde{\pi}}$ of order more than
$(l+1)_{*}$, but which do not contain top order spatial derivatives of acoustical entities, are bounded using Proposition 12.11 and 12.12 and the bootstrap assumption.
We obtain a bound for these contributions to the error integrals associated to $K_{0}$ and to any of the variations, up to the top order, by:
\begin{align}
 C_{l}\delta_{0}\bar{\mu}^{-a}_{m}(t)[1+\log(1+t)]^{p}\{\mathcal{P}_{[l+2]}+\sqrt{\mathcal{G}^{\prime}_{1,[l+2];a,q}(t)+\mathcal{G}_{0,[l+2];a,p}(t)}\}\\\notag
\cdot \{\int_{0}^{t}(1+t')^{-3/2}\mathcal{E}_{0,[l+2]}(t')dt'+\int_{0}^{\epsilon_{0}}\mathcal{F}^{t}_{0,[l+2]}(u')du'\}^{1/2}\\\notag
+C_{l}\delta_{0}\int_{0}^{t}(1+t')^{-3/2}\mathcal{E}_{0,[l+2]}(t')dt'+C_{l}\delta_{0}\int_{0}^{\epsilon_{0}}\mathcal{F}^{t}_{0,[l+2]}(u')du'\\\notag
+C_{l}\{\bar{\mathcal{F}}^{\prime t}_{1,[l+2]}(u')du'+\bar{K}_{[l+2]}(t,\epsilon_{0})\}^{1/2}\{\int_{0}^{t}(1+t')^{-3/2}\mathcal{E}_{0,[l+2]}(t')dt'
+\int_{0}^{\epsilon_{0}}\mathcal{F}^{t}_{0,[l+2]}(u')du'\}^{1/2}
\end{align}
 The last three terms are contributed by the terms in (14.14) containing the top order derivatives of the acoustical entities of which however one is a derivative with 
respect to $L$ (hence are expressible in terms of the top order derivatives of the $\psi_{\alpha}$). Here:
\begin{align}
 \mathcal{F}^{t}_{0,[n]}(u)=\sum_{m=1}^{n}\mathcal{F}^{t}_{0,m}(u)
\end{align}
where $\mathcal{F}^{t}_{0,n}(u)$ represents the sum of the fluxes associated to the vectorfield $K_{0}$ of all the $n$th order variations. Also:
\begin{align}
 \bar{\mathcal{F}}^{\prime t}_{1,[n]}(u)=\sup_{t'\in[0,t]}\{[1+\log(1+t')]^{-4}\mathcal{F'}^{t'}_{1,[n]}(u)\}\\
\bar{K}_{[n]}(t,u)=\sup_{t'\in[0,t]}\{[1+\log(1+t')]^{-4}K_{[n]}(t',u)\}
\end{align}

    On the other hand, all the other terms in the sum (14.14) contain derivatives of the $\leftexp{(Y)}{\tilde{\pi}}$ of order at most $(l+1)_{*}$, thus spatial 
derivatives of $\chi'$ of order at most $(l+1)_{*}$ and spatial derivatives of $\mu$ of order at most $(l+1)_{*}+1$, which are bounded in $L^{\infty}
(\Sigma_{t}^{\epsilon_{0}})$ by Proposition 12.9 and 12.10 and the bootstrap assumption. The contributions of these terms are bounded according to Lemma 7.6:
\begin{align}
 C_{l}\int_{0}^{t}(1+t')^{-3/2}\mathcal{E}_{0,[l+2]}(t')dt'+C_{l}\int_{0}^{\epsilon_{0}}\mathcal{F}^{t}_{0,[l+2]}(u')du'\\\notag
+C_{l}\{\int_{0}^{\epsilon_{0}}\bar{\mathcal{F}}^{\prime t}_{1,[l+2]}(u')du'+\bar{K}_{[l+2]}(t,\epsilon_{0})\}^{1/2}\{\int_{0}^{t}(1+t')^{-3/2}\mathcal{E}_{0,[l+2]}(t')dt'
+\int_{0}^{\epsilon_{0}}\mathcal{F}^{t}_{0,[l+2]}(u')du'\}^{1/2}
\end{align}
From the definitions (15.26),(15.36) and also the relation (15.30), we have, for all $t^{\prime}\in[0,t]$:
\begin{align}
 [1+\log(1+t')]^{-4}\mathcal{F}^{\prime t'}_{1,[n]}(u)\leq\bar{\mu}^{-2a}_{m,u}(t')[1+\log(1+t')]^{2p}\mathcal{H}^{\prime t}_{1,[n];a,q}(u)\\
[1+\log(1+t')]^{-4}K_{[n]}(t',u)\leq\bar{\mu}^{-2a}_{m,u}(t')[1+\log(1+t')]^{2p}I_{[l+2];a,q}(t,u)
\end{align}
taking the supremum over $t'\in[0,t]$, we obtain, in view of Corollary 2 of Lemma 8.11,
\begin{align}
 \bar{\mathcal{F}}^{\prime t}_{1,[n]}(u)\leq C\bar{\mu}^{-2a}_{m,u}(t)[1+\log(1+t)]^{2p}\mathcal{H'}^{t}_{1,[n];a,q}(u)\\
\bar{K}_{[n]}(t,u)\leq C\bar{\mu}^{-2a}_{m,u}(t)[1+\log(1+t)]^{2p}I_{[n];a,q}(t,u)
\end{align}
In view of the definition (15.27) and the fact that $\bar{\mu}_{m,u}(t)$ is a non-increasing function of $u$ at each $t$ we then have:
\begin{align}
 \int_{0}^{u}\bar{\mathcal{F}}^{\prime t}_{1,[n]}(u')du'\leq C\bar{\mu}^{-2a}_{m,u}(t)[1+\log(1+t)]^{2p}V'_{1,[n];a,q}(t,u)
\end{align}
We also define:
\begin{align}
 \mathcal{H}^{t}_{0,[n];a,p}(u)=\sup_{t'\in[0,t]}\{[1+\log(1+t')]^{-2p}\bar{\mu}^{2a}_{m,u}(t')\mathcal{F}^{t'}_{0,[n]}(u)\}\\
V_{0,[n];a,p}(t,u)=\int_{0}^{u}\mathcal{H}^{t}_{0,[n];a,p}(u')du'
\end{align}
Then $\mathcal{H}^{t}_{0,[n];a,p}(u)$ is a non-decreasing function of $t$ at each $u$, $V_{0,[n];a,p}(t,u)$ is a non-decreasing function of $u$ at each $t$ as well as
a non-decreasing function of $t$ at each $u$, and for each $u\in[0,\epsilon_{0}]$:
\begin{align}
 \int_{0}^{u}\mathcal{F}^{t}_{0,[n]}(u')du'\leq\int_{0}^{u}\bar{\mu}^{-2a}_{m,u'}(t)[1+\log(1+t)]^{2p}\mathcal{H}^{t}_{0,[n];a,p}(u')du'\\\notag
\leq\bar{\mu}^{-2a}_{m,u}(t)[1+\log(1+t)]^{2p}V_{0,[n];a,p}(t,u)
\end{align}
So we can bound (15.99) and (15.103) by:
\begin{align}
 C_{l}\delta_{0}\bar{\mu}^{-2a}_{m}(t)[1+\log(1+t)]^{2p}\cdot\{(\mathcal{P}_{[l+2]})^{2}+\delta'(\mathcal{G'}_{1,[l+2];a,q}(t)+\mathcal{G}_{0,[l+2];a,p}(t))\\\notag
+(1+\frac{1}{\delta'})(\int_{0}^{t}(1+t')^{-3/2}\mathcal{G}_{0,[l+2];a,p}(t')dt'+V_{0,[l+2];a,p}(t,\epsilon_{0}))\\\notag
+\delta'(V'_{1,[l+2];a,q}(t,\epsilon_{0})+I_{[l+2];a,q}(t,\epsilon_{0}))\}
\end{align}
for any positive constant $\delta'$ ($\delta^{\prime}$ will be chosen below).

    Finally, we have to deal with the first of (15.16). According to (5.282) we have:
\begin{align}
 \int_{W^{t}_{\epsilon_{0}}}\sum_{k=1}^{7}Q_{0}[\psi]d\mu_{g}\leq \int_{0}^{t}(1+t')^{-2}[1+\log(1+t')]^{4}B_{s}(t')\bar{\mathcal{E}}^{\prime}_{1}[\psi](t')dt'\\\notag
+C\int_{0}^{t}(1+t')^{-1}[1+\log(1+t')]^{-2}(\bar{\mathcal{E}}^{\prime}_{1}[\psi](t')+\bar{\mathcal{E}}_{0}[\psi](t'))dt'\\\notag
+CV_{0}[\psi](t,\epsilon_{0})\\\notag
+C\{V_{0}[\psi](t,\epsilon_{0})+\int_{0}^{t}(1+t')^{-1}[1+\log(1+t')]^{-2}\bar{\mathcal{E}}_{0}[\psi](t')dt'\}^{1/2}(V'_{1}[\psi](t,\epsilon_{0}))^{1/2}\\\notag
+C\int_{0}^{t}(1+t')^{-2}[1+\log(1+t')]^{2}V'_{1}[\psi](t',\epsilon_{0})dt'\\\notag
+C\bar{K}[\psi](t,\epsilon_{0})^{1/2}(\int_{0}^{t}(1+t')^{-1}[1+\log(1+t')]^{-2}\bar{\mathcal{E}}_{0}[\psi](t')dt')^{1/2}\\\notag
+C\bar{K}[\psi](t,\epsilon_{0})^{1/2}(V_{0}[\psi](t,\epsilon_{0}))^{1/2}
\end{align}
It follows that:
\begin{align}
 \int_{W^{t}_{\epsilon_{0}}}\sum_{k=1}^{7}Q_{0}[\psi]d\mu_{g}\leq\int_{0}^{t}\tilde{B}_{s}(t')\bar{\mathcal{E}}^{\prime}_{1}[\psi](t')dt'\\\notag
+C(1+\frac{1}{\delta'})\{\int_{0}^{t}(1+t')^{-1}[1+\log(1+t')]^{-2}\bar{\mathcal{E}}_{0}[\psi](t')dt'+V_{0}[\psi](t,\epsilon_{0})\}\\\notag
+\delta'V'_{1}[\psi](t,\epsilon_{0})\\\notag
+C\int_{0}^{t}(1+t')^{-2}[1+\log(1+t')]^{2}V'_{1}[\psi](t',\epsilon_{0})dt'\\\notag
+\delta'\bar{K}[\psi](t,\epsilon_{0})
\end{align}
where:
\begin{align}
 \tilde{B}_{s}(t)=(1+t)^{-2}[1+\log(1+t)]^{4}B_{s}(t)+C(1+t)^{-1}[1+\log(1+t)]^{-2}
\end{align}
Since
\begin{align}
 \bar{\mathcal{E}}^{\prime}_{1}[\psi](t)=\sup_{t'\in[0,t]}\{[1+\log(1+t')]^{-4}\mathcal{E}^{\prime}_{1}[\psi](t')\}\\\notag
\leq\sup_{t'\in[0,t]}\{\bar{\mu}^{-2a}_{m}(t')[1+\log(1+t')]^{2p}\mathcal{G}^{\prime}_{1;a,q}[\psi](t')\}\\\notag
\leq C\bar{\mu}^{-2a}_{m}(t)[1+\log(1+t)]^{2p}\mathcal{G}^{\prime}_{1;a,q}[\psi](t)
\end{align}
the first term on the right in (15.114) is bounded by:
\begin{align}
 C\bar{\mu}^{-2a}_{m}(t)[1+\log(1+t)]^{2p}\int_{0}^{t}\tilde{B}_{s}(t')\mathcal{G}^{\prime}_{1;a,q}[\psi](t')dt'
\end{align}
In regard to the second term on the right in (15.114) we have:
\begin{align}
 \int_{0}^{t}(1+t')^{-1}[1+\log(1+t')]^{-2}\bar{\mathcal{E}}_{0}[\psi](t')dt'\\\notag
\leq C\bar{\mu}^{-2a}_{m}(t)[1+\log(1+t)]^{2p}\int_{0}^{t}(1+t')^{-1}[1+\log(1+t')]^{-2}\mathcal{G}_{0;a,p}[\psi](t')dt'
\end{align}
Moreover, defining:
\begin{align}
 \mathcal{H}^{t}_{0;a,p}[\psi](u)=\sup_{t'\in[0,t]}\{[1+\log(1+t')]^{-2p}\bar{\mu}^{2a}_{m,u}(t')\mathcal{F}^{t'}_{0}[\psi](u)\}\\
V_{0;a,p}[\psi](t,u)=\int_{0}^{u}\mathcal{H}^{t}_{0;a,p}[\psi](u')du'
\end{align}
we have, in analogy with (15.111):
\begin{align}
 V_{0}[\psi](t,u)=\int_{0}^{u}\mathcal{F}^{t}_{0}(u')du'\leq\bar{\mu}^{-2a}_{m,u}(t)[1+\log(1+t)]^{2p}V_{0;a,p}[\psi](t,u)
\end{align}
In regard to the third and fourth term on the right of (15.114), recalling from Chapter 5 that:
\begin{align}
 V'_{1}[\psi](t,u)=\int_{0}^{u}\bar{\mathcal{F}}^{\prime t}_{1}[\psi](u')du'
\end{align}
we have, using an inequality similar to (15.106),
\begin{align}
 V'_{1}[\psi](t,u)\leq C\bar{\mu}^{-2a}_{m,u}(t)[1+\log(1+t)]^{2p}V'_{1;a,q}[\psi](t,u)
\end{align}
hence also:
\begin{align}
 \int_{0}^{t}(1+t')^{-2}[1+\log(1+t')]^{2}V'_{1}[\psi](t',\epsilon_{0})dt'\\\notag
\leq C\bar{\mu}^{-2a}_{m}(t)[1+\log(1+t)]^{2p}\int_{0}^{t}(1+t')^{-2}[1+\log(1+t')]^{2}V'_{1;a,q}[\psi](t',\epsilon_{0})dt'
\end{align}
Finally, in regard to the last term, we have, in analogy with (15.107):
\begin{align}
 \bar{K}[\psi](t,\epsilon_{0})\leq C\bar{\mu}^{-2a}_{m}(t)[1+\log(1+t)]^{2p}I_{a,q}[\psi](t,\epsilon_{0})
\end{align}
where:
\begin{align}
 I_{a,q}[\psi](t,u)=\sup_{t'\in[0,t]}\{[1+\log(1+t)]^{-2q}\bar{\mu}^{2a}_{m,u}(t')K[\psi](t',u)\}
\end{align}
Substituting the above in (15.114) we obtain:
\begin{align}
 \int_{W^{t}_{\epsilon_{0}}}\sum_{k=1}^{7}Q_{0}[\psi]d\mu_{g}\leq C\bar{\mu}^{-2a}_{m}(t)[1+\log(1+t)]^{2p}\{
\int_{0}^{t}\tilde{B}_{s}(t')\mathcal{G}^{\prime}_{1;a,q}[\psi](t')dt'+(1+\frac{1}{\delta'})\\\notag
[\int_{0}^{t}(1+t')^{-1}[1+\log(1+t')]^{-2}\mathcal{G}_{0;a,p}[\psi](t')dt'+
V_{0;a,p}[\psi](t,\epsilon_{0})]\\\notag
+\delta'V'_{1;a,q}[\psi](t,\epsilon_{0})+\int_{0}^{t}(1+t')^{-2}[1+\log(1+t')]^{2}V'_{1;a,q}[\psi](t',\epsilon_{0})dt'+\delta'I_{a,q}[\psi](t,\epsilon_{0})\}
\end{align}

    We now consider the energy identity corresponding to $K_{0}$ ((5.74) with $u=\epsilon_{0}$) and to the variations (14.52) $R_{j}R_{i_{l}}...R_{i_{1}}\psi_{\alpha}$,
$j=1,2,3$, $\alpha=0,1,2,3$, for a given multi-index $(i_{1}...i_{l})$.

    Summing over $j$ and $\alpha$ we then obtain from (15.95), (15.98), (15.112) and (15.127):
\begin{align}
 \{\sum_{j,\alpha}\mathcal{E}_{0}[R_{j}R_{i_{l}}...R_{i_{1}}\psi_{\alpha}](t)+\sum_{j,\alpha}\mathcal{F}^{t}_{0}[R_{j}R_{i_{l}}...R_{i_{1}}\psi_{\alpha}]
(\epsilon_{0})\}\bar{\mu}^{2a}_{m}(t)[1+\log(1+t)]^{-2p}\\\notag
\leq C(\frac{1}{a}+\frac{1}{p})\leftexp{(i_{1}...i_{l})}{\mathcal{G}}_{0,l+2;a,p}(t)+C_{p,l}\mathcal{D}_{[l+2]}\\\notag
+C_{p,l}\delta_{0}\{\mathcal{G}_{0,[l+2];a,p}(t)+\mathcal{G}^{\prime}_{1,[l+2];a,q}(t)\}\\\notag
+C\int_{0}^{t}\tilde{B}_{s}(t')\mathcal{G}^{\prime}_{1,[l+2];a,q}(t')dt'\\\notag
+C\int_{0}^{t}(1+t')^{-2}[1+\log(1+t')]^{2}V'_{1,[l+2];a,q}(t',\epsilon_{0})dt'\\\notag
+C_{l}(1+\frac{1}{\delta'})\{\int_{0}^{t}(1+t')^{-1}[1+\log(1+t')]^{-2}\mathcal{G}_{0,[l+2];a,p}(t')dt'+V_{0,[l+2];a,p}(t,\epsilon_{0})\}\\\notag
+C_{l}\delta'\{V'_{1,[l+2];a,q}(t,\epsilon_{0})+I_{[l+2];a,q}(t,\epsilon_{0})\}
\end{align}
We substitute (15.92)-(15.94) to obtain:
\begin{align}
 \{\sum_{j,\alpha}\mathcal{E}_{0}[R_{j}R_{i_{l}}...R_{i_{1}}\psi_{\alpha}](t)+\sum_{j,\alpha}\mathcal{F}^{t}_{0}[R_{j}R_{i_{l}}...R_{i_{1}}\psi_{\alpha}]
(\epsilon_{0})\}\bar{\mu}^{2a}_{m}(t)[1+\log(1+t)]^{-2p}\\\notag
\leq C(\frac{1}{a}+\frac{1}{p})\leftexp{(i_{1}...i_{l})}{\mathcal{G}}_{0,l+2;a,p}(t)+C_{p,l}(\delta_{0}+\delta')\bar{\mathcal{G}}_{0,[l+2];a,p}(t)\\\notag
+C_{p,l}(1+\frac{1}{\delta'})\{\mathcal{D}_{[l+2]}+\int_{0}^{t}\tilde{B}_{s}(t')\bar{\mathcal{G}}_{0,[l+2];a,p}(t')dt'+V_{0,[l+2];a,p}(t,\epsilon_{0})\}
\end{align}
Here we have used the fact in Proposition 13.2:
\begin{align}
 \int_{0}^{s}\tilde{B}_{s}(t)dt\leq C\quad:\textrm{independent of} \quad  s
\end{align}
It is crucial that the constant $C$ in front of the first term on the right of (15.129) is independent of $a,p,l$. Moreover, (15.129) holds with $\epsilon_{0}$ replaced by 
any $u\in(0,\epsilon_{0}]$, that is, we have, for all $u\in(0,\epsilon_{0}]$ and all $t\in[0,s]$:
\begin{align}
 \{\sum_{j,\alpha}\mathcal{E}^{u}_{0}[R_{j}R_{i_{l}}...R_{i_{1}}\psi_{\alpha}](t)+\sum_{j,\alpha}\mathcal{F}^{t}_{0}[R_{j}R_{i_{l}}...R_{i_{1}}\psi_{\alpha}]
(u)\}\bar{\mu}^{2a}_{m,u}(t)[1+\log(1+t)]^{-2p}\\\notag
\leq C(\frac{1}{a}+\frac{1}{p})\leftexp{(i_{1}...i_{l})}{\mathcal{G}}^{u}_{0,l+2;a,p}(t)+C_{p,l}(\delta_{0}+\delta')\bar{\mathcal{G}}^{u}_{0,[l+2];a,p}(t)\\\notag
+C_{p,l}(1+\frac{1}{\delta'})\{\mathcal{D}^{u}_{[l+2]}+\int_{0}^{t}\tilde{B}_{s}(t')\bar{\mathcal{G}}^{u}_{0,[l+2];a,p}(t')dt'+V_{0,[l+2];a,p}(t,u)\}
\end{align}

     Keeping only the term
\begin{align}
 \sum_{j,\alpha}\mathcal{E}^{u}_{0}[R_{j}R_{i_{l}}...R_{i_{1}}\psi_{\alpha}]
\end{align}
on the left of (15.131), we have:
\begin{align}
 \bar{\mu}^{2a}_{m,u}(t)[1+\log(1+t)]^{-2p}\sum_{j,\alpha}\mathcal{E}^{u}_{0}[R_{j}R_{i_{l}}...R_{i_{1}}\psi_{\alpha}](t)\\\notag
\leq C(\frac{1}{a}+\frac{1}{p})\leftexp{(i_{1}...i_{l})}{\mathcal{G}}^{u}_{0,l+2;a,p}(t)+C_{p,l}(\delta_{0}+\delta')\bar{\mathcal{G}}^{u}_{0,[l+2];a,p}(t)\\\notag
+C_{p,l}(1+\frac{1}{\delta'})\{\mathcal{D}^{u}_{[l+2]}+\int_{0}^{t}\tilde{B}_{s}(t')\bar{\mathcal{G}}^{u}_{0,[l+2];a,p}(t')dt'+V_{0,[l+2];a,p}(t,u)\}
\end{align}
The same holds with $t$ replaced by any $t^{\prime}\in[0,t]$. Now the right hand side of (15.132) is a non-decreasing function of $t$ at each $u$. The inequality 
corresponding to $t^{\prime}$ thus holds a fortiori if we again replace $t^{\prime}$ by $t$ on the right hand side. Taking then the supremum over all $t^{\prime}\in[0,t]$
on the left we obtain, in view of the definition (15.46),
\begin{align}
 \leftexp{(i_{1}...i_{l})}{\mathcal{G}}^{u}_{0,l+2;a,p}(t)\\\notag
\leq C(\frac{1}{a}+\frac{1}{p})\leftexp{(i_{1}...i_{l})}{\mathcal{G}}^{u}_{0,l+2;a,p}(t)+C_{p,l}(\delta_{0}+\delta')\bar{\mathcal{G}}^{u}_{0,[l+2];a,p}(t)\\\notag
+C_{p,l}(1+\frac{1}{\delta'})\{\mathcal{D}^{u}_{[l+2]}+\int_{0}^{t}\tilde{B}_{s}(t')\bar{\mathcal{G}}^{u}_{0,[l+2];a,p}(t')dt'+V_{0,[l+2];a,p}(t,u)\}
\end{align}
If we choose $a$ and $p$ suitably large so that:
\begin{align}
 C(\frac{1}{a}+\frac{1}{p})\leq\frac{1}{2}
\end{align}
then (15.133) implies:
\begin{align}
 \frac{1}{2}\leftexp{(i_{1}...i_{l})}{\mathcal{G}}^{u}_{0,l+2;a,p}(t)\leq C_{p,l}(\delta_{0}+\delta')\bar{\mathcal{G}}^{u}_{0,[l+2];a,p}(t)\\\notag
+C_{p,l}(1+\frac{1}{\delta'})\{\mathcal{D}^{u}_{[l+2]}+\int_{0}^{t}\tilde{B}_{s}(t')\bar{\mathcal{G}}^{u}_{0,[l+2];a,p}(t')dt'+V_{0,[l+2];a,p}(t,u)\}
\end{align}

     Consider now the integral identities corresponding to the variations (14.54) and $K_{0}$. In each of these we have the $borderline$ integral (14.180),
bounded by (14.181),
\begin{align}
 C(\frac{1}{2a}+\frac{1}{2p})\bar{\mu}^{-2a}_{m}(t)[1+\log(1+t)]^{2p}\leftexp{(i_{1}...i_{l-m})}{\mathcal{G}}_{0,m,l+2;a,p}(t)
\end{align}
We also have the remaining integrals, bounded in the case of (14.185) by (15.96) and in the case of (14.198) by the sum of (14.158) and (14.199), hence also by 
(15.96). We also have the integral (14.201) bounded by (14.202), hence by (15.97). Combining we see that all the remaining integrals associated to (14.54) and 
$K_{0}$ containing the top order spatial derivatives of the acoustical entities are bounded by (15.98). On the other hand, the error integrals associated to 
(14.54) contributed by all the other terms in (14.14) are bounded by the sum of (15.99) and (15.103), hence by (15.112). Finally, the contribution from the first
of (15.16) is bounded by (15.127).
 
Consequently, the integral identities in question yield:
\begin{align}
 \{\sum_{\alpha}\mathcal{E}_{0}[R_{i_{l-m}}...R_{i_{1}}(T)^{m+1}\psi_{\alpha}](t)+\sum_{\alpha}\mathcal{F}^{t}_{0}[R_{i_{l-m}}...R_{i_{1}}(T)^{m+1}
\psi_{\alpha}](\epsilon_{0})\}\\\notag
\cdot\bar{\mu}^{2a}_{m}(t)[1+\log(1+t)]^{-2p}\\\notag
\leq C(\frac{1}{a}+\frac{1}{p})\leftexp{(i_{1}...i_{l-m})}{\mathcal{G}}_{0,m,[l+2];a,p}(t)+C_{p,l}(\delta_{0}+\delta')\bar{\mathcal{G}}_{0,[l+2];a,p}(t)\\\notag
+C_{p,l}(1+\frac{1}{\delta'})\{\mathcal{D}_{[l+2]}+\int_{0}^{t}\tilde{B}_{s}(t')\bar{\mathcal{G}}_{0,[l+2];a,p}(t')dt'+V_{0,[l+2];a,p}(t,\epsilon_{0})\}
\end{align}
and
\begin{align}
 \{\sum_{\alpha}\mathcal{E}^{u}_{0}[R_{i_{l-m}}...R_{i_{1}}(T)^{m+1}\psi_{\alpha}](t)+\sum_{\alpha}\mathcal{F}^{t}_{\alpha}[R_{i_{l-m}}...R_{i_{1}}(T)^{m+1}
\psi_{\alpha}](u)\}\\\notag
\cdot\bar{\mu}^{2a}_{m,u}(t)[1+\log(1+t)]^{-2p}\\\notag
\leq C(\frac{1}{a}+\frac{1}{p})\leftexp{(i_{1}...i_{l-m})}{\mathcal{G}}^{u}_{0,m,[l+2];a,p}(t)+C_{p,l}(\delta_{0}+\delta')\bar{\mathcal{G}}^{u}_{0,[l+2];a,p}(t)\\\notag
+C_{p,l}(1+\frac{1}{\delta'})\{\mathcal{D}^{u}_{[l+2]}+\int_{0}^{t}\tilde{B}_{s}(t')\bar{\mathcal{G}}^{u}_{0,[l+2];a,p}(t')dt'+V_{0,[l+2];a,p}(t,u)\}
\end{align}
Keeping only the term
\begin{align*}
 \sum_{\alpha}\mathcal{E}^{u}_{0}[R_{i_{l-m}}...R_{i_{1}}(T)^{m+1}\psi_{\alpha}]
\end{align*}
on the left in (15.139), we have:
\begin{align}
 \bar{\mu}^{2a}_{m,u}(t)[1+\log(1+t)]^{-2p}\sum_{\alpha}\mathcal{E}_{0}^{u}[R_{i_{l-m}}...R_{i_{1}}(T)^{m+1}\psi_{\alpha}](t)\\\notag
\leq C(\frac{1}{a}+\frac{1}{p})\leftexp{(i_{1}...i_{l-m})}{\mathcal{G}}^{u}_{0,m,l+2;a,p}(t)+C_{p,l}(\delta_{0}+\delta')\bar{\mathcal{G}}^{u}_{0,[l+2];a,p}(t)\\\notag
C_{p,l}(1+\frac{1}{\delta'})\{\mathcal{D}^{u}_{[l+2]}+\int_{0}^{t}\tilde{B}_{s}(t')\bar{\mathcal{G}}^{u}_{0,[l+2];a,p}(t')dt'+V_{0,[l+2];a,p}(t,u)\}
\end{align}
Seeing that the right hand side of (15.140) is a non-decreasing function of $t$ at each $u$, the previous reasoning applies and we deduce:
\begin{align}
 \leftexp{(i_{1}...i_{l-m})}{\mathcal{G}}^{u}_{0,m,l+2;a,p}(t)\\\notag
\leq C(\frac{1}{a}+\frac{1}{p})\leftexp{(i_{1}...i_{l-m})}{\mathcal{G}}^{u}_{0,m,l+2;a,p}(t)+C_{p,l}(\delta_{0}+\delta')\bar{\mathcal{G}}^{u}_{0,[l+2];a,p}(t)\\\notag
C_{p,l}(1+\frac{1}{\delta'})\{\mathcal{D}^{u}_{[l+2]}+\int_{0}^{t}\tilde{B}_{s}(t')\bar{\mathcal{G}}^{u}_{0,[l+2];a,p}(t')dt'+V_{0,[l+2];a,p}(t,u)\}
\end{align}
Choosing $a$ and $p$ sufficiently large so that (15.135) holds, this implies:
\begin{align}
 \frac{1}{2}\leftexp{(i_{1}...i_{l-m})}{\mathcal{G}}^{u}_{0,m,l+2;a,p}(t)\leq C_{p,l}(\delta_{0}+\delta')\bar{\mathcal{G}}^{u}_{0,[l+2];a,p}(t)\\\notag
+C_{p,l}(1+\frac{1}{\delta'})\{\mathcal{D}^{u}_{[l+2]}+\int_{0}^{t}\tilde{B}_{s}(t')\bar{\mathcal{G}}^{u}_{0,[l+2];a,p}(t')dt'+V_{0,[l+2];a,p}(t,u)\}
\end{align}

     Consider finally the integral identities corresponding to any of the other variations $\psi$ of order up to $l+2$, and to $K_{0}$. These identities 
contain no borderline integrals and all error integrals involved are bounded by:
\begin{align}
 C_{p,l}(\delta_{0}+\delta')\bar{\mathcal{G}}^{u}_{0,[l+2];a,p}(t)\\\notag
+C_{p,l}(1+\frac{1}{\delta'})\{\mathcal{D}^{u}_{[l+2]}+\int_{0}^{t}\tilde{B}_{s}(t')\bar{\mathcal{G}}^{u}_{0,[l+2];a,p}(t')dt'+V_{0,[l+2];a,p}(t,u)\}
\end{align}
 We thus obtain, for each such variation $\psi$ the inequality:
\begin{align}
 \{\mathcal{E}^{u}_{0}[\psi](t)+\mathcal{F}^{t}_{0}[\psi](u)\}\bar{\mu}^{2a}_{m,u}(t)[1+\log(1+t)]^{-2p}\\\notag
\leq C_{p,l}(\delta_{0}+\delta')\bar{\mathcal{G}}^{u}_{0,[l+2];a,p}(t)\\\notag
+C_{p,l}(1+\frac{1}{\delta'})\{\mathcal{D}^{u}_{[l+2]}+\int_{0}^{t}\tilde{B}_{s}(t')\bar{\mathcal{G}}^{u}_{0,[l+2];a,p}(t')dt'+V_{0,[l+2];a,p}(t,u)\}
\end{align}
Keeping only the term 
\begin{align*}
 \mathcal{E}^{u}_{0}[\psi]
\end{align*}
on the left in (15.144) and applying the above reasoning we deduce:
\begin{align}
 \mathcal{G}^{u}_{0;a,p}[\psi](t)\leq C_{p,l}(\delta_{0}+\delta')\bar{\mathcal{G}}^{u}_{0,[l+2];a,p}(t)\\\notag
+C_{p,l}(1+\frac{1}{\delta'})\{\mathcal{D}^{u}_{[l+2]}+\int_{0}^{t}\tilde{B}_{s}(t')\bar{\mathcal{G}}^{u}_{0,[l+2];a,p}(t')dt'+V_{0,[l+2];a,p}(t,u)\}
\end{align}

     Now according to the definition of $\mathcal{E}^{u}_{0,[l+2]}(t)$ we have
\begin{align}
 \mathcal{E}^{u}_{0,[l+2]}(t)=\sum_{i_{1}...i_{l}}\{\sum_{j,\alpha}\mathcal{E}^{u}_{0}[R_{j}R_{i_{l}}...R_{i_{1}}\psi_{\alpha}](t)\}\\\notag
+\sum_{m=0}^{l}\sum_{i_{1}...i_{l-m}}\{\sum_{\alpha}\mathcal{E}^{u}_{0}[R_{i_{l-m}}...R_{i_{1}}(T)^{m+1}\psi_{\alpha}](t)\}\\\notag
+\sum_{\alpha}\mathcal{E}^{u}_{0}[\psi](t)
\end{align}
where the last sum is over all the other variations $\psi$ up to order $l+2$. By (15.68) we then have:
\begin{align}
 \mathcal{G}^{u}_{0,[l+2];a,p}(t)\leq\sum_{i_{1}...i_{l}}\leftexp{(i_{1}...i_{l})}{\mathcal{G}}^{u}_{0,l+2;a,p}(t)
+\sum_{m=0}^{l}\sum_{i_{1}...i_{l-m}}\leftexp{(i_{1}...i_{l-m})}{\mathcal{G}}^{u}_{0,m,l+2;a,p}(t)\\\notag
+\sum_{\psi}\mathcal{G}^{u}_{0;a,p}[\psi](t)
\end{align}
Thus, summing (15.135) over $i_{1}...i_{l}$, summing (15.142) over $i_{1}...i_{l-m}$ and over $m=0,...,l$, and summing (15.145) over all the other variations $\psi$ 
of order up to $l+2$, we obtain:
\begin{align}
 \frac{1}{2}\mathcal{G}^{u}_{0,[l+2];a,p}(t)\leq C_{p,l}(\delta_{0}+\delta')\bar{\mathcal{G}}^{u}_{0,[l+2];a,p}(t)\\\notag
+C_{p,l}(1+\frac{1}{\delta'})\{\mathcal{D}^{u}_{[l+2]}+\int_{0}^{t}\tilde{B}_{s}(t')\bar{\mathcal{G}}^{u}_{0,[l+2];a,p}(t')dt'+V_{0,[l+2];a,p}(t,u)\}
\end{align}
for a new constant $C_{p,l}$. Requiring $\delta_{0}$ and $\delta^{\prime}$ to satisfy the smallness condition and equation respectively:
\begin{align}
 C_{p,l}\delta_{0}\leq \frac{1}{8},\quad C_{p,l}\delta'=\frac{1}{8}
\end{align}
(15.148) implies:
\begin{align}
 \frac{1}{2}\mathcal{G}^{u}_{0,[l+2];a,p}(t)\leq\frac{1}{4}\bar{\mathcal{G}}^{u}_{0,[l+2];a,p}(t)\\\notag
+C_{p,l}\{\mathcal{D}^{u}_{[l+2]}+\int_{0}^{t}\tilde{B}_{s}(t')\bar{\mathcal{G}}^{u}_{0,[l+2];a,p}(t')dt'+V_{0,[l+2];a,p}(t,u)\}
\end{align}
for a new constant $C_{p,l}$. The same holds with $u$ replaced by any $u^{\prime}\in[0,u]$. Now the right hand side of (15.150) is a non-decreasing function of $u$ 
at each $t$. The inequality corresponding to $u^{\prime}$ thus holds a fortiori if we again replace $u^{\prime}$ by $u$ on the right. Taking then the supremum
over all $u^{\prime}\in[0,u]$ on the left we obtain, in view of the definition (15.89),
\begin{align}
 \frac{1}{4}\bar{\mathcal{G}}^{u}_{0,[l+2];a,p}(t)\leq C_{p,l}\{\mathcal{D}^{u}_{[l+2]}+\int_{0}^{t}\tilde{B}_{s}(t')\bar{\mathcal{G}}^{u}_{0,[l+2];a,p}
(t')dt'+V_{0,[l+2];a,p}(t,u)\}
\end{align}
In view of (15.130), since $V_{0,[l+2];a,p}(t,u)$ is a non-decreasing function of $t$ at each $u$, (15.151) implies:
\begin{align}
 \bar{\mathcal{G}}^{u}_{0,[l+2];a,p}(t)\leq C_{p,l}\{\mathcal{D}^{u}_{[l+2]}+V_{0,[l+2];a,p}(t,u)\}
\end{align}
hence also:
\begin{align}
 \int_{0}^{t}\tilde{B}_{s}(t')\bar{\mathcal{G}}^{u}_{0,[l+2];a,p}(t')dt'\leq C_{p,l}\{\mathcal{D}^{u}_{[l+2]}+V_{0,[l+2];a,p}(t,u)\}
\end{align}

     Substituting (15.152) and (15.153) in the right-hand sides of (15.131), (15.138) and (15.144), omitting in each case the first term in parenthesis on the left we
obtain:
\begin{align}
 \bar{\mu}^{2a}_{m,u}(t)[1+\log(1+t)]^{-2p}\sum_{j,\alpha}\mathcal{F}^{t}_{0}[R_{j}R_{i_{l}}...R_{i_{1}}\psi_{\alpha}](u)\\\notag
\leq C_{p,l}\{\mathcal{D}^{u}_{[l+2]}+V_{0,[l+2];a,p}(t,u)\}\\
\bar{\mu}^{2a}_{m,u}(t)[1+\log(1+t)]^{-2p}\sum_{j,\alpha}\mathcal{F}^{t}_{0}[R_{j}R_{i_{l-m}}...R_{i_{1}}(T)^{m+1}\psi_{\alpha}](u)\\\notag
\leq C_{p,l}\{\mathcal{D}^{u}_{[l+2]}+V_{0,[l+2];a,p}(t,u)\}\\
\bar{\mu}^{2a}_{m,u}(t)[1+\log(1+t)]^{-2p}\sum_{j,\alpha}\mathcal{F}^{t}_{0}[\psi](u)\\\notag
\leq C_{p,l}\{\mathcal{D}^{u}_{[l+2]}+V_{0,[l+2];a,p}(t,u)\}
\end{align}
Summing as in (15.146) yields:
\begin{align}
 \bar{\mu}^{2a}_{m,u}(t)[1+\log(1+t)]^{-2p}\mathcal{F}^{t}_{0,[l+2]}(u)\leq C_{p,l}\{\mathcal{D}^{u}_{[l+2]}+V_{0,[l+2];a,p}(t,u)\}
\end{align}
Since the inequality holds with $t$ replaced by any $t^{\prime}\in[0,t]$ and the right hand side is a non-decreasing function of $t$
at each $u$, the inequality corresponding to $t^{\prime}$ holds a fortiori if we again replace $t^{\prime}$ by $t$ on the right. 
Taking then the supremum over all $t^{\prime}\in[0,t]$ on the left we obtain, in view of the definition (15.109), 
\begin{align}
 \mathcal{H}^{t}_{0,[l+2];a,p}(u)\leq C_{p,l}\{\mathcal{D}^{u}_{[l+2]}+V_{0,[l+2];a,p}(t,u)\}
\end{align}
Recalling the definition of $V_{0,[l+2];a,p}(t,u)$, this is:
\begin{align}
 \mathcal{H}^{t}_{0,[l+2];a,p}(u)\leq C_{p,l}\mathcal{D}^{u}_{[l+2]}+C_{p,l}\int_{0}^{u}\mathcal{H}^{t}_{0,[l+2];a,p}(u')du'
\end{align}
Since $[0,\epsilon_{0}]$ is a bounded interval and $\mathcal{D}^{u}_{[l+2]}$ is a non-decreasing function of $u$, (15.159) implies:
\begin{align}
 \mathcal{H}^{t}_{0,[l+2];a,p}\leq C_{p,l}\mathcal{D}^{u}_{[l+2]}
\end{align}
hence also:
\begin{align}
 V_{0,[l+2];a,p}(t,u)\leq C_{p,l}\epsilon_{0}\mathcal{D}^{u}_{[l+2]}
\end{align}
Substituting this in (15.152) we obtain:
\begin{align}
 \bar{\mathcal{G}}^{u}_{0,[l+2];a,p}(t)\leq C_{p,l}\mathcal{D}^{u}_{[l+2]}
\end{align}
Substituting finally (15.161) in (15.90) and the result in (15.91), (15.93), (15.94) yields:
\begin{align}
 \mathcal{G}^{\prime u}_{1,[l+2];a,q}(t),\quad \mathcal{H}^{\prime t}_{1,[l+2];a,q}(u),\quad I_{[l+2];a,q}(t,u)\leq C_{p,l}\mathcal{D}^{u}_{[l+2]}
\end{align}
This completes the top order energy estimates.

      At this point we shall specify the choice of $a$ and $p$ (recall from (15.30) that $q=p+2$). We also take $a$ to be the form:
\begin{align}
 a=[a]+\frac{3}{4}
\end{align}
where $[a]$ denotes the integral part of $a$. We then choose $[a]$ to be the smallest positive integer so that with $C$ the largest of the two constants in 
(15.57) and (15.134) we have:
\begin{align}
 \frac{C}{a}\leq\frac{3}{8}
\end{align}
We then choose $p$ so that with the same constant $C$ we have:
\begin{align}
 \frac{C}{p}=\frac{1}{8}
\end{align}
Then the conditions (15.57) and (15.134) are both satisfied and thereafter $a,p$ and $q$ are fixed. We then fix $l$ so that:
\begin{align}
 l_{*}\geq [a]+4
\end{align}
We shall see the reason for this in the following. So $l$ is thenceforth fixed as well, therefore we shall not denote the dependence of the constants 
on $a,p,q,l$. Actually, the constants shall now depend only on the the function $H$. We shall still require $\delta_{0}$ to be suitably small in 
relation to these constants.

\chapter{The Descent Scheme}
     We now consider again the estimates for the energies $\mathcal{E}_{0,[l+1]}$ and $\mathcal{E}^{\prime}_{1,[l+1]}$, associated to $K_{0}$ and $K_{1}$ respectively.
We shall presently obtain improved estimates for these next to the top order energies, using the estimates for the top order energies just obtained. We shall
show that the next to the top order quantities satisfy estimates similar to (15.160), (15.162) and (15.163), but with $a$ replaced by
\begin{align}
 b=a-1
\end{align}
$p$ replaced by $0$, and, accordingly, $q$ replaced by $2$.

     Let us first consider the error integrals associated to the variations of order $l+1$ and to $K_{1}$, involving the highest spatial derivatives of the acoustical 
entities. The leading contributions are:

1) The contribution of (14.56) with $l+1$ replaced by $l$ to the corresponding integral (14.59), namely the integral (14.202) with $l$ replaced by $l-1$:
\begin{align}
 -\int_{W^{t}_{u}}(\omega/\nu)(R_{i_{l}}...R_{i_{1}}\textrm{tr}\chi')(T\psi_{\alpha})((L+\nu)R_{i_{l}}...R_{i_{1}}\psi_{\alpha})dt'du'd\mu_{\tilde{\slashed{g}}}
\end{align}
 2) The contribution of (14.57) with $l$ replaced by $l-1$ to the corresponding integral (14.59), namely the integral (14.369) with $l$ replaced by $l-1$:
\begin{align}
 -\int_{W^{t}_{u}}(\omega/\nu)(T\psi_{\alpha})(R_{i_{l-1-m}}...R_{i_{1}}(T)^{m}\slashed{\Delta}\mu)((L+\nu)R_{i_{l-1-m}}...R_{i_{1}}(T)^{m+1}\psi_{\alpha})dt'du'
d\mu_{\tilde{\slashed{g}}}
\end{align}
for $m=0,...,l-1$.

     By the bounds (14.294) the integral (16.2) is bounded by:
\begin{align}
 C\delta_{0}\int_{W^{t}_{u}}(1+t')|R_{i_{l}}...R_{i_{1}}\textrm{tr}\chi'||(L+\nu)R_{i_{l}}...R_{i_{1}}\psi_{\alpha}|dt'du'd\mu_{\slashed{g}}\\\notag
\leq C\delta_{0}\{\int_{W^{t}_{u}}|R_{i_{l}}...R_{i_{1}}\textrm{tr}\chi'|^{2}dt'du'd\mu_{\slashed{g}}\}^{1/2}\\\notag
\cdot\{\int_{W^{t}_{u}}(1+t')^{2}|(L+\nu)R_{i_{l}}...R_{i_{1}}\psi_{\alpha}|^{2}dt'du'd\mu_{\slashed{g}}\}^{1/2}
\end{align}
Now the last factor is bounded by:
\begin{align}
 C\{\int_{0}^{u}\mathcal{F}^{\prime t}_{1,[l+1]}(u')du'\}^{1/2}
\end{align}
while the first factor is:
\begin{align}
 \{\int_{0}^{t}\|R_{i_{l}}...R_{i_{1}}\textrm{tr}\chi'\|^{2}_{L^{2}(\Sigma_{t'}^{u})}dt'\}^{1/2}
\end{align}
which is bounded in terms of:
\begin{align}
 \{\int_{0}^{t}(\mathcal{A}^{\prime u}_{[l]}(t'))^{2}dt'\}^{1/2}
\end{align}
Here $\mathcal{A}^{\prime u}_{[l]}(t)$ is defined by replacing the $L^{2}$ norms over $\Sigma_{t}^{\epsilon_{0}}$ in the definition of $\mathcal{A}^{\prime}_{[l]}(t)$ 
by $L^{2}$ norms over $\Sigma_{t}^{u}$. Proposition 12.11 applies with $\epsilon_{0}$ replaced by $u$, giving a bound for $\mathcal{A}^{\prime u}_{[l]}(t)$ by:
\begin{align}
 C(1+t)^{-1}\int_{0}^{t}(1+t')^{-1}[\mathcal{W}^{u}_{[l+2]}(t')+\mathcal{W}^{Qu}_{[l+1]}(t')]dt'
\end{align}
Now we have:
\begin{align}
 \mathcal{W}^{u}_{[l+2]}(t)+\mathcal{W}^{Qu}_{[l+1]}(t)\leq C\bar{\mu}^{-1/2}_{m,u}\sqrt{\mathcal{E}^{\prime u}_{1,[l+2]}(t)}\\\notag
\leq C\bar{\mu}^{-a-1/2}_{m,u}(t)[1+\log(1+t)]^{q}\sqrt{\mathcal{G}^{\prime u}_{1,[l+2];a,q}(t)}
\end{align}
Therefore (16.8) is bounded by:
\begin{align}
 C(1+t)^{-1}\sqrt{\mathcal{G}^{\prime u}_{1,[l+2];a,q}(t)}\int_{0}^{t}(1+t')^{-1}\bar{\mu}^{-a-1/2}_{m,u}(t')[1+\log(1+t')]^{q}dt'\\\notag
\leq C(1+t)^{-1}\sqrt{\mathcal{G}^{\prime u}_{1,[l+2];a,q}(t)}J^{u}_{a,q-1}(t)\\\notag
\leq C(1+t)^{-1}[1+\log(1+t)]^{q+1}\bar{\mu}^{-a+1/2}_{m,u}(t)\sqrt{\mathcal{D}^{u}_{[l+2]}}
\end{align}
Here $J^{u}_{a,q}(t)$ is the integral (14.147) with $\epsilon_{0}$ replaced by $u$ and we have used (14.153) with $\epsilon_{0}$ replaced by $u$ as well as the bounds 
(15.163). It follows that the contribution of (16.8) to (16.7) is bounded by:
\begin{align}
 C\sqrt{\mathcal{D}^{u}_{[l+2]}}\{\int_{0}^{t}(1+t')^{-2}[1+\log(1+t')]^{2q+2}\bar{\mu}^{-2a+1}_{m,u}(t')dt'\}^{1/2}
\end{align}
the last integral being estimated by following an argument similar to that used to estimate (14.115) (with $\epsilon_{0}$ replaced by $u$). By (15.28) with
$a$ replaced by $b$ and $q$ by $2$ we have:
\begin{align}
 \int_{0}^{u}\mathcal{F}^{\prime t}_{1,[l+1]}(u')du'\leq \bar{\mu}^{-2b}_{m,u}(t)[1+\log(1+t)]^{4}V'_{1,[l+1];b,2}(t,u)
\end{align}
we conclude that (16.4) is bounded by:
\begin{align}
 C\delta_{0}\bar{\mu}^{-2b}_{m,u}(t)[1+\log(1+t)]^{2}\sqrt{\mathcal{D}^{u}_{[l+2]}}\sqrt{V'_{1,[l+1];b,2}(t,u)}
\end{align}
where $b$ is defined by (16.1).

    By (14.294), (16.3) is bounded by:
\begin{align}
 C\delta_{0}\int_{W^{t}_{u}}(1+t')|R_{i_{l-1-m}}...R_{i_{1}}(T)^{m}\slashed{\Delta}\mu||(L+\nu)R_{i_{l-1-m}}...R_{i_{1}}(T)^{m+1}\psi_{\alpha}|dt'du'd\mu_{\slashed{g}}
\\\notag
\leq C\delta_{0}\{\int_{W^{t}_{u}}|R_{i_{l-1-m}}...R_{i_{1}}(T)^{m}\slashed{\Delta}\mu|^{2}dt'du'd\mu_{\slashed{g}}\}^{1/2}\\\notag
\cdot\{\int_{W^{t}_{u}}(1+t')^{2}|(L+\nu)R_{i_{l-1-m}}...R_{i_{1}}(T)^{m+1}\psi_{\alpha}|^{2}dt'du'd\mu_{\slashed{g}}\}^{1/2}
\end{align}
Now the last factor is bounded by (16.5) while the first factor is:
\begin{align}
 \{\int_{0}^{t}\|R_{i_{l-1-m}}...R_{i_{1}}(T)^{m}\slashed{\Delta}\mu\|^{2}_{L^{2}(\Sigma_{t'}^{u})}dt'\}^{1/2}
\end{align}
which is bounded by:
\begin{align}
 \{\int_{0}^{t}(1+t')^{-4}(\mathcal{B}^{u}_{[m,l+1]}(t'))^{2}dt'\}^{1/2}
\end{align}
Here the quantity $\mathcal{B}^{u}_{[m,l+1]}(t)$ is defined by replacing the $L^{2}$ norms over $\Sigma_{t}^{\epsilon_{0}}$ in the definition of the quantity 
$\mathcal{B}_{[m,l+1]}(t)$ by $L^{2}$ norms over $\Sigma_{t}^{u}$.
Proposition 12.12 applies with $\epsilon_{0}$ replaced by $u$, giving a bound for $\mathcal{B}^{u}_{[m,l+1]}(t)$ in terms of:
\begin{align}
 C(1+t)\int_{0}^{t}(1+t')^{-1}\{\mathcal{W}^{u}_{\{l+2\}}(t')+\mathcal{W}^{Qu}_{\{l+1\}}(t')\}dt'
\end{align}
Now we have:
\begin{align}
 \mathcal{W}^{u}_{\{l+2\}}(t)+\mathcal{W}^{Qu}_{\{l+1\}}(t)\\\notag
\leq C\bar{\mu}^{-1/2}_{m,u}\{\sqrt{\mathcal{E}^{\prime u}_{1,[l+2]}(t)}+\sqrt{\mathcal{E}^{u}_{0,[l+2]}(t)}\}\\\notag
\leq C\bar{\mu}^{-a-1/2}_{m,u}(t)[1+\log(1+t)]^{q}\{\sqrt{\mathcal{G}^{\prime u}_{1,[l+2];a,q}(t)}+\sqrt{\mathcal{G}^{u}_{0,[l+2];a,p}(t)}\}
\end{align}
Therefore (16.17) is bounded by:
\begin{align}
 C(1+t)\{\sqrt{\mathcal{G}^{\prime u}_{1,[l+2];a,q}(t)}+\sqrt{\mathcal{G}^{u}_{0,[l+2];a,p}(t)}\}\\\notag
\cdot \int_{0}^{t}(1+t')^{-1}\bar{\mu}^{-a-1/2}_{m,u}(t')[1+\log(1+t')]^{q}dt'\\\notag
\leq C(1+t)\{\sqrt{\mathcal{G}^{\prime u}_{1,[l+2];a,q}(t)}+\sqrt{\mathcal{G}^{u}_{0,[l+2];a,p}(t)}\}J^{u}_{a,q-1}\\\notag
\leq C(1+t)[1+\log(1+t)]^{q+1}\bar{\mu}^{-a+1/2}_{m,u}(t)\sqrt{\mathcal{D}^{u}_{[l+2]}}
\end{align}
where we have used (15.162), (15.163). It follows that the contribution of (16.17) to (16.16) is bounded by:
\begin{align}
 C\sqrt{\mathcal{D}^{u}_{[l+2]}}\{\int_{0}^{t}(1+t')^{-2}[1+\log(1+t')]^{2q+2}\bar{\mu}^{-2a+1}_{m,u}(t')dt'\}^{1/2}\leq 
\sqrt{\mathcal{D}^{u}_{[l+2]}}\bar{\mu}^{-a+1}_{m,u}(t)
\end{align}
In view of (16.12), (16.14) is bounded by:
\begin{align}
 C\delta_{0}\bar{\mu}^{-2b}_{m,u}(t)[1+\log(1+t)]^{2}\sqrt{\mathcal{D}^{u}_{[l+2]}}\sqrt{V'_{1,[l+1];b,2}(t,u)}
\end{align}
where $b$ is defined by (16.1).

      We proceed to consider the error integrals associated to the variations of order $l+1$ and to $K_{0}$, involving the highest spatial derivatives of the 
acoustical entities. The leading contributions are:

1) The contribution of (14.56) with $l+1$ replaced by $l$ to the corresponding integral (14.62), namely the integral on the left in (14.63) with $l+1$ replaced by $l$
(and $\epsilon_{0}$ by $u$):
\begin{align}
 \int_{W^{t}_{u}}|R_{i_{l}}...R_{i_{1}}\textrm{tr}\chi'||T\psi_{\alpha}||\underline{L}R_{i_{l}}...R_{i_{1}}\psi_{\alpha}|dt'du'd\mu_{\tilde{\slashed{g}}}
\end{align}
2) The contribution of (14.57) with $l$ replaced by $l-1$ to the corresponding integral (14.62), namely the integral (14.162) with $l$ replaced by $l-1$:
\begin{align}
 \int_{W^{t}_{u}}|R_{i_{l-1-m}}...R_{i_{1}}(T)^{m}\slashed{\Delta}\mu||T\psi_{\alpha}||\underline{L}R_{i_{l-1-m}}...R_{i_{1}}(T)^{m+1}\psi_{\alpha}|
dt'du'd\mu_{\tilde{\slashed{g}}}
\end{align}
By (14.294), (16.22) is bounded by:
\begin{align}
 C\delta_{0}\int_{W^{t}_{u}}(1+t')^{-1}|R_{i_{l}}...R_{i_{1}}\textrm{tr}\chi'||\underline{L}R_{i_{l}}...R_{i_{1}}\psi_{\alpha}|dt'du'd\mu_{\slashed{g}}\\\notag
\leq C\delta_{0}\int_{0}^{t}(1+t')^{-1}\|R_{i_{l}}...R_{i_{1}}\textrm{tr}\chi'\|_{L^{2}(\Sigma_{t'}^{u})}\|\underline{L}R_{i_{l}}...R_{i_{1}}\psi_{\alpha}\|
_{L^{2}(\Sigma_{t'}^{u})}dt'
\end{align}
Substituting (16.10) for the first factor in the integrand on the right and noting the definition (16.1) of $b$ and the fact that:
\begin{align}
 \|\underline{L}R_{i_{l}}...R_{i_{1}}\psi_{\alpha}\|_{L^{2}(\Sigma_{t}^{u})}\leq C\sqrt{\mathcal{E}^{u}_{0,[l+1]}(t)}
\leq C\bar{\mu}^{-b}_{m,u}(t)\sqrt{\mathcal{G}^{u}_{0,[l+1];b,0}(t)}
\end{align}
we obtain that (16.24) is bounded by:
\begin{align}
 C\delta_{0}\sqrt{\mathcal{D}^{u}_{[l+2]}}\sqrt{\mathcal{G}^{u}_{0,[l+1];b,0}(t)}\int_{0}^{t}(1+t')^{-2}[1+\log(1+t')]^{q+1}
\bar{\mu}_{m,u}^{-2b-1/2}(t')dt'\\\notag
\leq C\delta_{0}\bar{\mu}^{-2b+1/2}_{m,u}(t)\sqrt{\mathcal{D}^{u}_{[l+2]}}\sqrt{\mathcal{G}^{u}_{0,[l+1];b,0}(t)}
\end{align}
the last integral being estimated in a similar manner as (14.115).

      By (14.294), (15.23) is bounded by:
\begin{align}
 C\delta_{0}\int_{W^{t}_{u}}(1+t')^{-1}|R_{i_{l-1-m}}...R_{i_{1}}(T)^{m}\slashed{\Delta}\mu||\underline{L}R_{i_{l-1-m}}...R_{i_{1}}(T)^{m+1}\psi_{\alpha}|
dt'du'd\mu_{\slashed{g}}\\\notag
\leq C\delta_{0}\int_{0}^{t}(1+t')^{-1}\|R_{i_{l-1-m}}...R_{i_{1}}(T)^{m}\slashed{\Delta}\mu\|_{L^{2}(\Sigma_{t'}^{u})}\\\notag
\cdot\|\underline{L}R_{i_{l-1-m}}...R_{i_{1}}(T)^{m+1}\psi_{\alpha}\|_{L^{2}(\Sigma_{t'}^{u})}dt'
\end{align}
Substituting (16.19) (multiplied by $C(1+t)^{-2}$) for the first factor in the integrand and noting again the definition (16.1) and the fact that:
\begin{align}
 \|\underline{L}R_{i_{l-1-m}}...R_{i_{1}}(T)^{m+1}\psi_{\alpha}\|_{L^{2}(\Sigma_{t'}^{u})}\leq C\sqrt{\mathcal{E}^{u}_{0,[l+1]}(t)}
\leq C\bar{\mu}^{-b}_{m,u}(t)\sqrt{\mathcal{G}^{u}_{0,[l+1];b,0}(t)}
\end{align}
we obtain that (16.27) is bounded by:
\begin{align}
 C\delta_{0}\sqrt{\mathcal{D}^{u}_{[l+2]}}\sqrt{\mathcal{G}^{u}_{0,[l+1];b,0}(t)}\int_{0}^{t}(1+t')^{-2}[1+\log(1+t')]^{q+1}\bar{\mu}^{-2b-1/2}_{m,u}(t')dt'\\\notag
\leq C\delta_{0}\bar{\mu}^{-2b+1/2}_{m,u}(t)\sqrt{\mathcal{D}^{u}_{[l+2]}}\sqrt{\mathcal{G}^{u}_{0,[l+1];b,0}(t)}
\end{align}
a bound of the same form as (16.26).

     Consider now the error integrals associated to the variations of order $l+1$ and to $K_{1}$ or $K_{0}$, which contain the lower order spatial derivatives of 
the acoustical entities and are contributed by the remaining terms in the sum (14.14) with $l+2$ replaced by $l+1$. By the discussion following (15.9), those  
contributed by the terms in which one factor is a derivative of the $\leftexp{(Y)}{\tilde{\pi}}$ of order more than $(l+1)_{*}$ are bounded, using Proposition 12.11 
and 12.12 and the bootstrap assumption, by (15.11), hence by (15.39), with $([l+2];a,q)$ replaced by $([l+1];b,2)$ (and $\epsilon_{0}$ by $u$), thus by:
\begin{align}
 C_{l}\delta_{0}\bar{\mu}^{-2b}_{m,u}(t)[1+\log(1+t)]^{4}\{\mathcal{D}^{u}_{[l+2]}+\delta(\mathcal{G}^{\prime u}_{1,[l+1];b,2}(t)+\mathcal{G}^{u}_{0,[l+1];b,0}(t))\\\notag
+\delta I_{[l+1];b,2}(t,u)+(1+\frac{1}{\delta})V'_{1,[l+1];b,2}(t,u)\}
\end{align}
 in the case of the error integrals associated to $K_{1}$, by (15.99), hence by (15.112), with $([l+2];a,p)$ replaced by $([l+1];b,0)$ (and $\epsilon_{0}$ by $u$),
thus by:
\begin{align}
 C_{l}\delta_{0}\bar{\mu}^{-2b}_{m,u}(t)\{\mathcal{D}^{u}_{[l+2]}+\delta'(\mathcal{G}^{\prime u}_{1,[l+1];b,2}(t)+\mathcal{G}^{u}_{0,[l+1];b,0}(t))\\\notag
+(1+\frac{1}{\delta'})(\int_{0}^{t}(1+t')^{-3/2}\mathcal{G}^{u}_{0,[l+1];b,0}(t')dt'+V_{0,[l+1];b,0}(t,u))\\\notag
+\delta'(V'_{1,[l+1];b,2}(t,u)+I_{[l+1];b,2}(t,u))\}
\end{align}
in the case of the error integrals associated to $K_{0}$. In (16.30), (16.31), the arbitrary positive constants $\delta$, $\delta'$ shall be chosen appropriately 
below. In fact, (16.30) and (16.31) also bound the error integrals associated to $K_{1}$ and to $K_{0}$ respectively, and to any of the variations of order up to
$l$ in which one factor is a derivative of $\leftexp{(Y)}{\tilde{\pi}}$ of order more than $(l+1)_{*}$.

     On the other hand, the error integrals associated to any of the variations of order up to $l+1$ contributed by the terms in the sum (14.14) containing derivatives 
of the $\leftexp{(Y)}{\tilde{\pi}}$ of order at most $(l+1)_{*}$ are bounded by (15.13) with $[l+2]$ replaced by $[l+1]$ (and $\epsilon_{0}$ by $u$):
\begin{align}
 C\int_{0}^{u}\mathcal{F}^{\prime t}_{1,[l+1]}(u')du'+C\{\int_{0}^{u}\mathcal{F}^{\prime t}_{1,[l+1]}(u')du'\}^{1/2}\{K_{[l+1]}(t,u)\}^{1/2}\\\notag
+C\{\int_{0}^{u}\mathcal{F}^{\prime t}_{1,[l+1]}(u')du'\}^{1/2}\cdot \\\notag
\{\int_{0}^{t}(1+t')^{-2}[1+\log(1+t')]^{2}(\mathcal{E}^{\prime u}_{1,[l+1]}(t')+u^{2}\mathcal{E}^{u}_{0,[l+1]}(t'))dt'\}^{1/2}\\\notag
\leq C\bar{\mu}^{-2b}_{m,u}(t)[1+\log(1+t)]^{4}\{V'_{1,[l+1];b,2}+(V'_{1,[l+1];b,2}(t,u))^{1/2}(I_{[l+1];b,2}(t,u))^{1/2}\\\notag
+(V'_{1,[l+1];b,2}(t,u))^{1/2}\cdot\\\notag
(\int_{0}^{t}(1+t')^{-2}[1+\log(1+t')]^{2}(\mathcal{G}^{\prime u}_{1,[l+1];b,2}(t')+\mathcal{G}^{u}_{0,[l+1];b,0}(t'))dt')^{1/2}\}
\end{align}
in the case of $K_{1}$, and by (15.103) with $[l+2]$ replaced by $[l+1]$ (and $\epsilon_{0}$ by $u$):
\begin{align}
 C\int_{0}^{t}(1+t')^{-3/2}\mathcal{E}^{u}_{0,[l+1]}(t')dt'+C\int_{0}^{u}\mathcal{F}^{t}_{0,[l+1]}(u')du'\\\notag
+C\{\int_{0}^{u}\bar{\mathcal{F}}^{\prime t}_{1,[l+1]}(u')du'\}^{1/2}\{\int_{0}^{t}(1+t')^{-3/2}\mathcal{E}^{u}_{0,[l+1]}(t')dt'+\int_{0}^{u}
\mathcal{F}^{t}_{0,[l+1]}(u')du'\}^{1/2}\\\notag
+C(\bar{K}_{[l+1]}(t,u))^{1/2}\{\int_{0}^{t}(1+t')^{-3/2}\mathcal{E}^{u}_{0,[l+1]}(t')dt'+\int_{0}^{u}
\mathcal{F}^{t}_{0,[l+1]}(u')du'\}^{1/2}\\\notag
\leq C\bar{\mu}^{-2b}_{m,u}(t)\{\int_{0}^{t}(1+t')^{-3/2}\mathcal{G}^{u}_{0,[l+1];b,0}(t')dt'+V_{0,[l+1];b,0}(t,u)\\\notag
+(V'_{1,[l+1];b,2}(t,u)+I_{[l+1];b,2}(t,u))^{1/2}(\int_{0}^{t}(1+t')^{-3/2}\mathcal{G}^{u}_{0,[l+1];b,0}(t')dt'+V_{0,[l+1];b,0}(t,u))^{1/2}\}
\end{align}
in the case of $K_{0}$.

    Finally, for each variation $\psi$ of order up to $l+1$, we have the error integrals (15.16). Consider first the second of these integrals, which is 
associated to $K_{1}$. Summing (15.17) over all such variations we obtain:
\begin{align}
 \sum_{\psi}\int_{W^{t}_{u}}\sum_{k=1}^{8}Q_{1,k}[\psi]d\mu_{g}\leq -\frac{1}{2}K_{[l+1]}(t,u)+CM_{[l+1]}(t,u)+\frac{3}{2}L_{[l+1]}(t,u)\\\notag
+\int_{0}^{t}\tilde{A}(t')\mathcal{E}^{\prime u}_{1,[l+1]}(t')dt'
\end{align}
Here,
\begin{align}
 M_{[l+1]}(t,u)=\sum_{\psi}M[\psi](t,u),\quad L_{[l+1]}(t,u)=\sum_{\psi}L[\psi](t,u)
\end{align}
From (15.18) and (15.29) with $(a,p)$ replaced by $(b,0)$, noting the fact that there is a fixed number of variations under consideration, we obtain:
\begin{align}
 M_{[l+1]}(t,u)\leq C[1+\log(1+t)]^{4}\sup_{t'\in[0,t]}\mathcal{E}^{u}_{0,[l+1]}(t')+\int_{0}^{u}\mathcal{F}^{\prime t}_{1,[l+1]}(u')du'\\\notag
\leq \bar{\mu}^{-2b}_{m,u}(t)[1+\log(1+t)]^{4}\{C\mathcal{G}^{u}_{0;b,0}(t)+V'_{1,[l+1];b,2}(t,u)\}
\end{align}
Also, from (15.1):
\begin{align}
 L_{[l+1]}(t,u)=\int_{0}^{t}(1+t')^{-1}[1+\log(1+t')]^{-1}\mathcal{E}^{\prime u}_{1,[l+1]}(t')dt'
\end{align}

    To estimate this borderline integral appropriately, we must use the following variant of Corollary 2 of Lemma 8.11:

$\textbf{Variant of Corollary 2}$: Let $b$ be a positive constant and $k$ any constant greater than $1$. Then if $\delta_{0}$ is suitably small, depending on an upper 
bound for $b$ and a lower bound greater than $1$ for $k$, we have, for all $t'\in[0,t]$, $t\in[1,s]$, and $u\in(0,\epsilon_{0}]$:
\begin{align*}
 \bar{\mu}^{-b}_{m,u}(t')\leq k\bar{\mu}^{-b}_{m,u}(t)
\end{align*}
     To prove this, we revisit the proof of Corollary 2. In the following, we shall take the constant $a$ in the proof of that corollary to be a suitably large constant,
having nothing to do with the choice of $a$ which we made above in relation to the top order energy estimates. In Case 1 of this proof (8.322) holds, hence:
\begin{align}
 \frac{\bar{\mu}^{-b}_{m,u}(t')}{\bar{\mu}^{-b}_{m,u}(t)}\leq (1-C\delta_{0})^{-b}\leq k
\end{align}
provided that $\delta_{0}$ is suitably small, depending on an upper bound for $b$ and a lower bound greater than $1$ for $k$. In Subcase 2a the lower bound (8.273)
holds (with $\epsilon_{0}$ replaced by $u$) for any given $a$, provided that $a\delta_{0}$ is suitably small, hence:
\begin{align}
 \frac{\bar{\mu}^{-b}_{m,u}(t')}{\bar{\mu}^{-b}_{m,u}(t)}\leq(1-\frac{2}{a})^{-b}\leq k
\end{align}
 if $a$ is suitably large, depending on an upper bound for $b$ and a lower bound greater than $1$ for $k$. In Subcase 2b the lower bound (8.303) as well as the upper
bound (8.312) hold (with $\epsilon_{0}$ replaced by $u$) for any given $a$, provided that $a\delta_{0}$ is suitably small, hence:
\begin{align}
 \frac{\bar{\mu}^{-b}_{m,u}(t')}{\bar{\mu}^{-b}_{m,u}(t)}\leq (\frac{1-\frac{2}{a}}{1+\frac{2}{a}})^{-b}(\frac{1-\delta_{1}\tau}{1-\delta_{1}\tau'})^{b}
\leq(\frac{1-\frac{2}{a}}{1+\frac{2}{a}})^{-b}\leq k
\end{align}
if $a$ is suitably large, depending on an upper bound for $b$ and a lower bound greater than $1$ for $k$. We have thus established the variant of Corollary 2.

     We now set:
\begin{align}
 k=\sqrt{\frac{4}{3}}
\end{align}
in the above variant of Corollary 2, and we proceed to estimate $L_{[l+1]}(t,u)$. We have,
from (16.38) and the definition (14.111) with $(a,q)$ replaced by $(b,2)$,
\begin{align}
 L_{[l+1]}(t,u)\leq\int_{0}^{t}(1+t')^{-1}[1+\log(1+t')]^{3}\bar{\mu}^{-2b}_{m,u}\mathcal{G}^{\prime u}_{1,[l+1];b,2}(t')dt'\\\notag
\leq\frac{4}{3}\bar{\mu}^{-2b}_{m,u}(t)\mathcal{G}^{\prime u}_{1,[l+1];b,2}(t)\int_{0}^{t}(1+t')^{-1}[1+\log(1+t')]^{3}dt'\\\notag
\leq\frac{1}{3}\bar{\mu}^{-2b}_{m,u}(t)[1+\log(1+t)]^{4}\mathcal{G}^{\prime u}_{1,[l+1];b,2}(t)
\end{align}

Consider now the integral identity (5.73) corresponding to $K_{1}$ and to the variations $\psi$ of order up to $l+1$. In each of these identities we also 
have the hypersurface integrals bounded according to (15.35) with $(a,p)$ replaced by $(b,0)$ (and $\epsilon_{0}$ by $u$) by:
\begin{align}
 C\bar{\mu}^{-2b}_{m,u}(t)[1+\log(1+t)]^{4}\mathcal{G}^{u}_{0,[l+1];b,0}(t)
\end{align}
Summing over all such variations we then obtain, from (16.42), and from (16.13), (16.21), (16.30), (16.32), (16.34), (16.36), the following:
\begin{align}
 \{\mathcal{E}^{\prime u}_{1,[l+1]}(t)+\mathcal{F}^{\prime t}_{1,[l+1]}(u)+\frac{1}{2}K_{[l+1]}(t,u)\}\bar{\mu}^{2b}_{m,u}(t)[1+\log(1+t)]^{-4}\\\notag
\leq \frac{1}{2}\mathcal{G}^{\prime u}_{1,[l+1];b,2}(t)+C\delta_{0}\delta\mathcal{G}^{\prime u}_{1,[l+1];b,2}(t)\\\notag
+CR'_{[l+1];b,2}(t,u)+C\int_{0}^{t}\bar{A}(t')\mathcal{G}^{\prime u}_{1,[l+1];b,2}(t')dt'
\end{align}
Here, $\delta$ is any positive constant ($\delta$ shall be chosen below), and:
\begin{align}
 R'_{[l+1];b,2}=\mathcal{D}^{u}_{[l+2]}+\mathcal{G}^{u}_{0,[l+1];b,0}(t)+(1+\frac{1}{\delta})V'_{1,[l+1];b,2}(t,u)+\delta I_{[l+1];b,2}(t,u)
\end{align}

     Keeping only the term $\mathcal{E}^{\prime u}_{1,[l+1]}(t)$ in parenthesis on the left in (16.44) we have:
\begin{align}
 \bar{\mu}^{2b}_{m,u}(t)[1+\log(1+t)]^{-4}\mathcal{E}^{\prime u}_{1,[l+1]}(t)\leq\frac{1}{2}\mathcal{G}^{\prime u}_{1,[l+1];b,2}(t)+C\delta_{0}\delta\mathcal{G}^{\prime u}
_{1,[l+1];b,2}(t)\\\notag
+CR'_{[l+1];b,2}(t,u)+C\int_{0}^{t}\bar{A}(t')\mathcal{G}^{\prime u}_{1,[l+1];b,2}(t')dt'
\end{align}
Since the right-hand side of the above is non-decreasing in $t$, applying the previous reasoning we deduce:
\begin{align}
 \frac{1}{3}\mathcal{G}^{\prime u}_{1,[l+1];b,2}(t)\leq CR'_{[l+1];b,2}(t,u)+C\int_{0}^{t}\bar{A}(t')\mathcal{G}^{\prime u}_{1,[l+1];b,2}(t')dt'
\end{align}
 provided that:
\begin{align}
 C\delta_{0}\delta\leq\frac{1}{6}
\end{align}
From this point we proceed as in the argument leading from (15.71) to (15.73), (15.74), to obtain:
\begin{align}
 \mathcal{G}^{\prime u}_{1,[l+1];b,2}(t)\leq CR'_{[l+1];b,2}(t,u)\\
\int_{0}^{t}\bar{A}(t')\mathcal{G}^{\prime u}_{1,[l+1];b,2}(t')dt'\leq CR'_{[l+1];b,2}(t,u)
\end{align}
Substituting these in the right-hand side of (16.44) and omitting the first term in parenthesis on the left, we obtain:
\begin{align}
 \{\mathcal{F}^{\prime t}_{1,[l+1]}(u)+\frac{1}{2}K_{[l+1]}(t,u)\}\bar{\mu}^{2b}_{m,u}(t)[1+\log(1+t)]^{-4}\leq CR'_{[l+1];b,2}(t,u)
\end{align}
From this point we proceed as in the argument leading from (15.78) to (15.91)-(15.94), choosing $\delta$ according to:
\begin{align}
 C\delta=\frac{1}{4}
\end{align}
We obtain in this way the estimates:
\begin{align}
 \mathcal{G}^{\prime u}_{1,[l+1];b,2}(t),\quad \mathcal{H}^{\prime t}_{1,[l+1];b,2}(u),\quad I_{[l+1];b,2}(t,u)\leq C\bar{Q}^{\prime}_{[l+1];b,2}(t,u)\\
V'_{1,[l+1];b,2}(t,u)\leq C\epsilon_{0}\bar{Q}^{\prime}_{[l+1];b,2}(t,u)
\end{align}
where:
\begin{align}
 \bar{Q}^{\prime}_{[l+1];b,2}(t,u)=\mathcal{D}^{u}_{[l+2]}+\bar{\mathcal{G}}^{u}_{0,[l+1];b,0}(t)
\end{align}

     Consider finally the integral identity corresponding to $K_{0}$ and to the variations $\psi$ of order up to $l+1$. In each of these identities we have
the first of (15.16). Summing (15.127) over all such variations we obtain (there is a fixed number of variations under consideration):
\begin{align}
 \bar{\mu}^{2b}_{m,u}(t)\sum_{\psi}\int_{W^{t}_{u}}\sum_{k=1}^{7}Q_{0}[\psi]d\mu_{g}\leq C\int_{0}^{t}\tilde{B}_{s}(t')\mathcal{G}^{\prime u}_{1,[l+1];b,2}(t')dt'\\\notag
+(1+\frac{1}{\delta'})\{\int_{0}^{t}(1+t')^{-1}[1+\log(1+t')]^{-2}\mathcal{G}^{u}_{0,[l+1];b,0}(t')dt'+V_{0,[l+1];b,0}(t,u)\}\\\notag
+C\delta' V'_{1,[l+1];b,2}(t,u)+C\int_{0}^{t}(1+t')^{-2}[1+\log(1+t')]^{-2}V'_{1,[l+1];b,2}(t',u)dt'\\\notag
+C\delta' I_{[l+1];b,2}(t,u)
\end{align}
Here, as in (16.31), $\delta'$ is an arbitrary positive constant which shall be chosen appropriately below. Summing over all such variations we then obtain, from (16.26), 
(16.29), (16.31), (16.33), (16.56), the following:
\begin{align}
 \{\mathcal{E}^{u}_{0,[l+1]}(t)+\mathcal{F}^{t}_{0,[l+1]}(t)\}\bar{\mu}^{2b}_{m,u}(t)\\\notag
\leq C\delta_{0}\delta'(\mathcal{G}^{u}_{0,[l+1];b,0}(t)+\mathcal{G}^{\prime u}_{1,[l+1];b,2}(t))+C(1+\frac{\delta_{0}}{\delta'})\mathcal{D}^{u}_{[l+2]}\\\notag
+C\int_{0}^{t}\tilde{B}_{s}(t')\mathcal{G}^{\prime u}_{1,[l+1];b,2}(t')dt'\\\notag
+C(1+\frac{1}{\delta'})\int_{0}^{t}(1+t')^{-1}[1+\log(1+t')]^{-2}\mathcal{G}^{u}_{0,[l+1];b,0}(t')dt'\\\notag
+C\delta'V'_{1,[l+1];b,2}(t,u)+C\int_{0}^{t}(1+t')^{-2}[1+\log(1+t')]^{2}V'_{1,[l+1];b,2}(t',u)dt'\\\notag
+C(1+\frac{1}{\delta'})V_{0,[l+1];b,0}(t,u)+C\delta' I_{[l+1];b,2}(t,u)
\end{align}
Substituting on the right (16.53)-(16.55) yields:
\begin{align}
 \{\mathcal{E}^{u}_{0,[l+1]}(t)+\mathcal{F}^{t}_{0,[l+1]}(u)\}\bar{\mu}^{2b}_{m,u}(t)\\\notag
\leq C\delta'\bar{\mathcal{G}}^{u}_{0,[l+1];b,0}(t)+C(1+\frac{\delta_{0}}{\delta'})\mathcal{D}^{u}_{[l+2]}+C(1+\frac{1}{\delta'})
\int_{0}^{t}\tilde{B}_{s}(t')\bar{\mathcal{G}}^{u}_{0,[l+1];b,0}(t')dt'\\\notag
+C(1+\frac{1}{\delta'})V_{0,[l+1];b,0}(t,u)
\end{align}
Keeping only the term $\mathcal{E}^{u}_{0,[l+1]}(t)$ in parenthesis on the left in (16.58), and noting that the right-hand side is a non-decreasing function of $t$
at each $u$, we may replace $t$ by any $t^{\prime}\in[0,t]$ on the left while keeping $t$ on the right, thus, taking the supremum over all $t^{\prime}
\in[0,t]$ we obtain, in view of the definition (15.50),
\begin{align}
 \mathcal{G}^{u}_{0,[l+1];b,0}(t)\leq C\delta'\bar{\mathcal{G}}^{u}_{0,[l+1];b,0}(t)+C(1+\frac{\delta_{0}}{\delta'})\mathcal{D}^{u}_{[l+2]}\\\notag
+C(1+\frac{1}{\delta'})\int_{0}^{t}\tilde{B}_{s}(t')\bar{\mathcal{G}}^{u}_{0,[l+1];b,0}(t')dt'+C(1+\frac{1}{\delta'})V_{0,[l+1];b,0}(t,u)
\end{align}
We now choose $\delta'$ such that:
\begin{align}
 C\delta'=\frac{1}{2}
\end{align}
Then in view of the fact that the right-hand side of (16.59) is non-decreasing in $u$, (16.59) implies:
\begin{align}
 \frac{1}{2}\bar{\mathcal{G}}^{u}_{0,[l+1];b,0}(t)\leq C\mathcal{D}^{u}_{[l+2]}+C\int_{0}^{t}\tilde{B}_{s}(t')\bar{\mathcal{G}}^{u}_{0,[l+1];b,0}(t')dt'
+CV_{0,[l+1];b,0}(t,u)
\end{align}
for new constants $C$.
From this point we proceed as in the argument leading from (15.151) to (15.160)-(15.162) to obtain the estimates:
\begin{align}
 \bar{\mathcal{G}}^{u}_{0,[l+1];b,0}(t),\quad \mathcal{H}^{t}_{0,[l+1];b,0}(u)\leq C\mathcal{D}^{u}_{[l+2]}
\end{align}
Substituting finally (16.62) in (16.55) and the result in (16.53) yields:
\begin{align}
 \mathcal{G}^{\prime u}_{1,[l+1];b,2}(t),\quad \mathcal{H}^{\prime t}_{1,[l+1];b,2}(u),\quad I_{[l+1];b,2}(t,u)\leq C\mathcal{D}^{u}_{[l+2]}
\end{align}
We have thus obtained improved energy estimates of the next to the top order, namely of order $l+1$.

     We proceed in this way, taking at the $n$th step
\begin{align}
 a_{n}=a-n\quad\textrm{and}\quad b_{n}=a_{n}-1=a-(n+1)=a_{n+1}
\end{align}
in the role of $a$ and $b$ respectively, the argument beginning in the paragraph containing (16.1) and concluding with (16.62), (16.63), 
being step $0$. For $n\geq 1$ we have $2$ and $0$ in the role of $q$ and $p$ respectively. The $n$th step is otherwise exactly the same as the $0$th 
step as above, as long as $b_{n}>0$, that is as long as $n\leq [a]-1$. In estimating the integrals $J^{u}_{a_{n},1}$ in (16.10), (16.19), as well
as the integrals
\begin{align}
 \int_{0}^{t}(1+t')^{-2}[1+\log(1+t')]^{6}\bar{\mu}^{-2a_{n}+1}_{m,u}(t')dt'
\end{align}
for $n\geq1$, we follow the argument leading from (14.147) to (14.153) in the case of $J^{u}_{a_{n},1}$, the argument leading from (14.115) to (14.131) in the case 
of (16.65). When referring to Lemma 8.11 in these arguments, we now mean Lemma 8.11 with $a$ set equal to $4$ in that lemma. Then in Case 1, (14.116) holds which 
implies:
\begin{align}
 \bar{\mu}^{-2a_{n}-1}_{m,u}(t)\leq C
\end{align}
for all $t\in[0,s]$, as the $a_{n}$ have a fixed upper bound $a$. Also, the subcases of Case 2 are defined according to:
\begin{align}
 t_{1}=e^{\frac{1}{8\delta_{1}}}-1
\end{align}
and in Subcase 2a we have the lower bound:
\begin{align}
 \bar{\mu}_{m,u}(t')\geq \frac{1}{2}
\end{align}
in place of (14.121), hence again:
\begin{align}
 \bar{\mu}^{-2a_{n}-1}_{m,u}(t')\leq C
\end{align}
while in Subcase 2b we have the lower bound:
\begin{align}
 \bar{\mu}_{m,u}\geq\frac{1}{2}(1-\delta_{1}\tau'),\quad \tau'=\log(1+t')
\end{align}
in place of (14.124), as well as the upper bound (8.312):
\begin{align}
 \bar{\mu}_{m}(t)\leq \frac{7}{4}(1-\delta_{1}\tau)
\end{align}
from the proof of Lemma 8.11 with $a$ set to equal to $4$, which implies:
\begin{align}
 \bar{\mu}^{-2a_{n}}_{m,u}\geq\frac{1}{C}(1-\delta_{1}\tau)^{-2a_{n}}\quad \tau=\log(1+t)
\end{align}
Noting again that the $a_{n}$ have a fixed upper bound. We thus obtain in Subcase 2b, in place of (14.151),
\begin{align}
 J^{u}_{a_{n},1}(t)-J^{u}_{a_{n},1}(t_{1})\leq C(1+\tau)^{2}\int_{\tau_{1}}^{\tau}(1-\delta_{1}\tau')^{-a_{n}-1/2}d\tau'\\\notag
\leq\frac{C}{\delta_{1}}\frac{(1+\tau)^{2}}{(a_{n}-\frac{1}{2})}(1-\delta_{1}\tau)^{-a_{n}+1/2}\\\notag
\leq C\frac{[1+\log(1+t)]^{3}}{(a_{n}-\frac{1}{2})}\bar{\mu}^{-a_{n}+1/2}_{m,u}
\end{align}
noting that the role of $a_{n}$ is different from that of $a$.

Since for any $n\leq[a]$:
\begin{align}
 a_{n}-\frac{1}{2}\geq a_{[a]}-\frac{1}{2}=\frac{3}{4}-\frac{1}{2}=\frac{1}{4}
\end{align}
the denominator in (16.74) is bounded from below by a positive constant even for $n=[a]$, hence combining with the results for Subcase 2a and for Case 1 
we conclude that:
\begin{align}
 J^{u}_{a_{n},1}(t)\leq C[1+\log(1+t)]^{3}\bar{\mu}^{-a_{n}+1/2}_{m,u}
\end{align}
for all $n=1,...,[a]$

Also, for $n\leq[a]-1$, so that $a_{n}\geq 7/4$, we have in Subcase 2b, in regard to the integral (16.65), in place of (14.126),
\begin{align}
 \int_{t_{1}}^{t}(1+t')^{-2}[1+\log(1+t')]^{6}(1-\delta_{1}\tau')^{-2a_{n}+1}dt'\\\notag
\leq C(1+t_{1})^{-1}[1+\log(1+t_{1})]^{6}\int_{\tau_{1}}^{\tau}(1-\delta_{1}\tau')^{-2a_{n}+1}d\tau'\\\notag
\leq C(1+t_{1})^{-1}\frac{[1+\log(1+t_{1})]^{6}}{\delta_{1}}\frac{(1-\delta_{1}\tau)^{-2a_{n}+2}}{2(a_{n}-1)}\\\notag
\leq C(1+t_{1})^{-1}[1+\log(1+t_{1})]^{7}\frac{(1-\delta_{1}\tau)^{-2(a_{n}-1)}}{2(a_{n}-1)}\\\notag
\leq C\bar{\mu}^{-2(a_{n}-1)}_{m,u}
\end{align}
hence combining with the results for Subcase 2a and for Case 1 (where (16.65) is simply bounded by a constant) we conclude that (16.65) is bounded by:
\begin{align}
 C\bar{\mu}^{-2(a_{n}-1)}_{m,u}
\end{align}
for all $n=1,...,[a]-1$.

Obviously, in (16.26) and (16.29) we may replace the integrals on the left by the integral (16.65) and we obtain the bounds at the $n$th step by:
\begin{align}
 C\delta_{0}\bar{\mu}^{-2b_{n}}_{m,u}(t)\sqrt{\mathcal{D}^{u}_{[l+2]}}\sqrt{\mathcal{G}^{u}_{0,[l+1-n];b_{n},0}(t)}
\end{align}

    Following then the argument of step 0, we obtain at the $n$th step, for $n=1,...,[a]-1$, the estimates:
\begin{align}
 \bar{\mathcal{G}}^{u}_{0,[l+1-n];b_{n},0}(t),\quad \mathcal{H}^{t}_{0,[l+1-n];b_{n},0}(u)\leq C\mathcal{D}^{u}_{[l+2]}\\\notag
\mathcal{G}^{\prime u}_{1,[l+1-n];b_{n},2}(t),\quad \mathcal{H}^{\prime t}_{1,[l+1-n];b_{n},2}(u),\quad I_{[l+1-n];b_{n},2}(t,u)\leq C\mathcal{D}^{u}_{[l+2]}
\end{align}
for all $(t,u)\in[0,s]\times [0,\epsilon_{0}]$. Here the constants $C$ may depend on $n$, but $n$ is in any case bounded.

     We now make the last step $n=[a]$. In this case, we have $a_{n}=a_{[a]}=\frac{3}{4}$ and (16.75) reads:
\begin{align}
 J^{u}_{3/4,1}(t)\leq C[1+\log(1+t)]^{3}\bar{\mu}^{-1/4}_{m,u}
\end{align}
while in regard to (16.65), we now have, from (16.76), in Subcase 2b, since $-2a_{[a]}+1=-1/2$:
\begin{align}
 \int_{t_{1}}^{t}(1+t')^{-2}[1+\log(1+t')]^{6}(1-\delta_{1}\tau')^{-1/2}dt'\\\notag
\leq C(1+t_{1})^{-1}[1+\log(1+t_{1})]^{6}\int_{\tau_{1}}^{\tau}(1-\delta_{1}\tau')^{-1/2}d\tau'\\\notag
\leq C(1+t_{1})^{-1}\frac{[1+\log(1+t_{1})]^{6}}{\delta_{1}}\frac{(1-\delta_{1}\tau_{1})^{1/2}}{1/2}\\\notag
\leq C(1+t_{1})^{-1}[1+\log(1+t_{1})]^{7}\frac{(7/8)^{1/2}}{1/2}\leq C
\end{align}
hence combining with the results for Subcase 2a and for Case 1 (where (16.65) is in any case bounded by a constant) we conclude that (16.65) is in the 
case $n=[a]$ simply bounded by a constant. We can thus set:
\begin{align}
 b_{n}=b_{[a]}=0
\end{align}
in the last step, and proceed exactly as in the preceding steps. We thus arrive at the estimates:
\begin{align}
\bar{\mathcal{G}}^{u}_{0,[l+1-[a]];0,0}(t),\quad \mathcal{H}^{t}_{0,[l+1-[a]];0,0}(u)\leq C\mathcal{D}^{u}_{[l+2]}\\\notag
\mathcal{G}^{\prime u}_{1,[l+1-[a]];0,2}(t),\quad \mathcal{H}^{\prime t}_{1,[l+1-[a]];0,2}(u),\quad I_{[l+1-[a]];0,2}(t,u)\leq C\mathcal{D}^{u}_{[l+2]}
\end{align}
for all $(t,u)\in[0,s]\times [0,\epsilon_{0}]$.

    These are the desired estimates. Because from the definitions,
\begin{align}
 \bar{\mathcal{G}}^{u}_{0,[l+1-[a]];0,0}(t)=\sup_{t'\in[0,t]}\{\mathcal{E}^{u}_{0,[l+1-[a]]}(t')\}\\
\mathcal{H}^{t}_{0,[l+1-[a]];0,0}(u)=\sup_{t'\in[0,t]}\{\mathcal{F}^{t'}_{0,[l+1-[a]]}(u)\}\\
\mathcal{G}^{\prime u}_{1,[l+1-[a]];0,2}(t)=\sup_{t'\in[0,t]}\{[1+\log(1+t')]^{-4}\mathcal{E}^{\prime u}_{1,[l+1-[a]]}(t')\}\\
\mathcal{H}^{\prime t}_{1,[l+1-[a]];0,2}(u)=\sup_{t'\in[0,t]}\{[1+\log(1+t)]^{-4}\mathcal{F}^{\prime t'}_{1,[l+1-[a]]}(u)\}\\
I_{[l+1-[a]];0,2}(t,u)=\sup_{t'\in[0,t]}\{[1+\log(1+t)]^{-4}K_{[l+1-[a]]}(t',u)\}
\end{align}
and the weights $\bar{\mu}_{m,u}$ have been eliminated.
\chapter{The Isoperimetric Inequality. Recovery of Assumption $\textbf{J}$. Recovery of the Bootstrap Assumption. Proof of the Main Theorem}

\section{Recovery of $\textbf{J}$-Preliminary}
We now establish assumption $\textbf{J}$ on the basis of the bootstrap assumption. We consider the first order variations:
\begin{align}
 \leftexp{(\alpha)}{\tilde{\psi}}_{1}=
\begin{cases}
  S\phi& \alpha=0\\
  \mathring{R}_{i}\phi&\alpha=i=1,2,3
 \end{cases}
\end{align}
To these Theorem 5.1 applies and we obtain, in particular:
\begin{align}
 \mathcal{E}^{u}_{0}[\leftexp{(\alpha)}{\tilde{\psi}}_{1}](t)\leq C\mathcal{E}^{u}_{0}[\leftexp{(\alpha)}{\tilde{\psi}}_{1}](0)
\end{align}
for all $(t,u)\in[0,s]\times [0,\epsilon_{0}]$.

We also consider the higher order variations, corresponding to the first order variations $\leftexp{(\alpha)}{\tilde{\psi}}_{1}$:
\begin{align}
 \leftexp{(\alpha;I_{1}...I_{n-1})}{\tilde{\psi}}_{n}=Y_{I_{n-1}}...Y_{I_{1}}\leftexp{(\alpha)}{\tilde{\psi}}_{1}
\end{align}
We require the multi-indices $(I_{1}...I_{n-1})$ to be of the form specified in the paragraph following (14.4) and concluding with (14.5).
We consider all such variations of order up to $4$. We denote by $\tilde{\mathcal{E}}^{u}_{0,[4]}(t)$ the sum of energies corresponding to $K_{0}$,
$\Sigma_{t}^{u}$, and to all such variations of order up to $4$. Now the integral identities corresponding to the variations in question and to $K_{0}$
and $K_{1}$ contain derivatives of $\leftexp{(Y)}{\tilde{\pi}}$ of order at most 3, which, in view of (15.167), are bounded in $L^{\infty}(\Sigma_{t}^{\epsilon_{0}})$
by Proposition 12.9 and 12.10 and the bootstrap assumption. The error integrals involved are then bounded exactly as in Chapter 7 and in accordance with the 
remark following Lemma 7.6 can be all absorbed in the error integrals of the fundamental energy estimate. Consequently, the conclusions of Theorem 5.1 hold for the 
variations in question and we obtain, in particular:
\begin{align}
 \tilde{\mathcal{E}}^{u}_{0,[4]}(t)\leq C\tilde{D}^{u}_{[4]}
\end{align}
 for all $(t,u)\in[0,s]\times[0,\epsilon_{0}]$.
Here we denote:
\begin{align}
 \tilde{D}^{u}_{[4]}=\tilde{\mathcal{E}}^{u}_{0,[4]}(0)
\end{align}
Applying Lemma 5.1 we then conclude that:
\begin{align}
 \tilde{\mathcal{S}}_{[4]}(t,u)\leq C\epsilon_{0}\tilde{D}^{u}_{[4]}
\end{align}
for all $(t,u)\in[0,s]\times[0,\epsilon_{0}]$.
Here we denote by $\tilde{\mathcal{S}}_{[4]}$ the integral on $S_{t,u}$ (with respect to $d\mu_{\slashed{g}}$) of the sum of the squares of all the variations 
under consideration. In particular, we have:
\begin{align}
 \int_{S_{t,u}}\{|S\phi|^{2}+\sum_{j_{1}}|R_{j_{1}}S\phi|^{2}+\sum_{j_{1},j_{2}}|R_{j_{2}}R_{j_{1}}S\phi|^{2}\}d\mu_{\slashed{g}}\leq \tilde{\mathcal{S}}_{[4]}(t,u)
\end{align}
\begin{align}
  \int_{S_{t,u}}\{|\mathring{R}_{i}\phi|^{2}+\sum_{j_{1}}|R_{j_{1}}\mathring{R}_{i}\phi|^{2}+\sum_{j_{1},j_{2}}|R_{j_{2}}R_{j_{1}}\mathring{R}_{i}\phi|^{2}\}
d\mu_{\slashed{g}}\leq \tilde{\mathcal{S}}_{[4]}(t,u)\quad: i=1,2,3 
\end{align}
\begin{align}
  \int_{S_{t,u}}\{|TS\phi|^{2}+\sum_{j_{1}}|R_{j_{1}}TS\phi|^{2}+\sum_{j_{1},j_{2}}|R_{j_{2}}R_{j_{1}}TS\phi|^{2}\}d\mu_{\slashed{g}}\leq \tilde{\mathcal{S}}_{[4]}(t,u)
\end{align}
\begin{align}
  \int_{S_{t,u}}\{|T\mathring{R}_{i}\phi|^{2}+\sum_{j_{1}}|R_{j_{1}}T\mathring{R}_{i}\phi|^{2}+\sum_{j_{1},j_{2}}|R_{j_{2}}R_{j_{1}}T\mathring{R}_{i}\phi|^{2}\}
d\mu_{\slashed{g}}\leq \tilde{\mathcal{S}}_{[4]}(t,u)\quad: i=1,2,3 
\end{align}

We are now in a position to apply the following lemma.

$\textbf{Lemma 17.1}$ Let $f$ be a function on $S_{t,u}$ with square-integrable derivatives up to 2nd order. We denote:
\begin{align*}
 \mathcal{S}_{[2]}[f]=\int_{S_{t,u}}\{|f|^{2}+\sum_{j_{1}}|R_{j_{1}}f|^{2}+\sum_{j_{1},j_{2}}|R_{j_{2}}R_{j_{1}}f|^{2}\}d\mu_{\slashed{g}}
\end{align*}
Then $f\in L^{\infty}(S_{t,u})$ and there is a positive numerical constant $C$ such that:
\begin{align*}
 \sup_{S_{t,u}}|f|\leq C(1+t)^{-1}(\mathcal{S}_{[2]}[f])^{1/2}
\end{align*}
$Proof$: By (8.332) we have:
\begin{align}
 C^{-1}(1+t)^{2}\leq A(t,u)\leq C(1+t)^{2}
\end{align}
 where $A(t,u)$ is the area of $S_{t,u}$ (by (8.332)). By Corollary 12.2.a we have:
\begin{align}
 \int_{S_{t,u}}(|f|^{2}+\sum_{i}|R_{i}f|^{2})d\mu_{\slashed{g}}+A(t,u)\int_{S_{t,u}}(|\slashed{d}f|^{2}+\sum_{i}|\slashed{d}R_{i}f|^{2})d\mu_{\slashed{g}}\leq 
C\mathcal{S}_{[2]}[f]
\end{align}
To continue with the proof of the lemma, we need the $isoperimetric$ $Sobolev$ $inequality$.

\section{The Isoperimetric Inequality}
       (See \cite{O}, \cite{Cha}) If $g$ is an arbitrary function which is integrable and with integrable derivative on $S_{t,u}$,
then $g$ is square-integrable on $S_{t,u}$ and, denoting by $\bar{g}$ the mean value of $g$ on $S_{t,u}$
\begin{align}
 \bar{g}=\frac{1}{A(t,u)}\int_{S_{t,u}}gd\mu_{\slashed{g}}
\end{align}
we have:
\begin{align}
 \int_{S_{t,u}}(g-\bar{g})^{2}d\mu_{\slashed{g}}\leq I(t,u)(\int_{S_{t,u}}|\slashed{d}g|d\mu_{\slashed{g}})^{2}
\end{align}
where $I(t,u)$ is the $isoperimetric$ $constant$ of $S_{t,u}$:
\begin{align}
 I=\sup_{U}\frac{\min\{\textrm{Area}(U),\textrm{Area}(U^{c})\}}{(\textrm{Perimeter}(\partial U))^{2}}
\end{align}
where the supremum is over all domains $U$ with $C^{1}$ boundary $\partial U$ in $S_{t,u}$, and $U^{c}=S_{t,u}\setminus U$ denotes the complement of $U$ in $S_{t,u}$.

     Since:
\begin{align}
 \|\bar{g}\|_{L^{2}(S_{t,u})}=|\bar{g}|A^{1/2},\quad |\bar{g}|\leq A^{-1}\|g\|_{L^{1}(S_{t,u})}
\end{align}
It follows that:
\begin{align}
 \|g\|_{L^{2}(S_{t,u})}\leq \sqrt{I'}\|g\|_{W^{1}_{1}(S_{t,u})}
\end{align}
where:
\begin{align}
 \|g\|_{W^{1}_{1}(S_{t,u})}=\|\slashed{d}g\|_{L^{1}(S_{t,u})}+A^{-1/2}\|g\|_{L^{1}(S_{t,u})}
\end{align}
and:
\begin{align}
 I'=\max\{I,1\}
\end{align}
We have:
\begin{align}
 \|f^{2}\|_{W^{1}_{1}(S_{t,u})}=\int_{S_{t,u}}|\slashed{d}(f^{2})|d\mu_{\slashed{g}}+A^{-1/2}\int_{S_{t,u}}|f|^{2}d\mu_{\slashed{g}}
\end{align}
and, since $\slashed{d}(f^{2})=2f\slashed{d}f$,
\begin{align}
 \int_{S_{t,u}}|\slashed{d}(f^{2})|d\mu_{\slashed{g}}\leq 2(\int_{S_{t,u}}|f|^{2}d\mu_{\slashed{g}})^{1/2}(\int_{S_{t,u}}|\slashed{d}f|^{2}d\mu_{\slashed{g}})^{1/2}
\leq CA^{-1/2}\mathcal{S}_{[2]}[f]
\end{align}
by (17.12). Since the second term on the right of (17.20) is similarly bounded, we obtain:
\begin{align}
 \|f^{2}\|_{W^{1}_{1}(S_{t,u})}\leq CA^{-1/2}\mathcal{S}_{[2]}[f]
\end{align}
We have:
\begin{align}
 \|\sum_{i}|R_{i}f|^{2}\|_{W^{1}_{1}(S_{t,u})}=\int_{S_{t,u}}|\slashed{d}(\sum_{i}|R_{i}f|^{2})|^{2}d\mu_{\slashed{g}}+A^{-1/2}\int_{S_{t,u}}\sum_{i}|R_{i}f|^{2}
d\mu_{\slashed{g}}
\end{align}
and, since $\slashed{d}(\sum_{i}|R_{i}f|^{2})=2\sum_{i}(R_{i}f)\slashed{d}(R_{i}f)$,
\begin{align}
 \int_{S_{t,u}}|\slashed{d}(\sum_{i}|R_{i}f|^{2})|d\mu_{\slashed{g}}\leq 2(\int_{S_{t,u}}\sum_{i}|R_{i}f|^{2}d\mu_{\slashed{g}})^{1/2}
(\int_{S_{t,u}}\sum_{i}|\slashed{d}R_{i}f|^{2}d\mu_{\slashed{g}})^{1/2}\leq CA^{-1/2}\mathcal{S}_{[2]}[f]
\end{align}
by (17.12). Since the second term on the right in (17.23) is similarly bounded, we obtain:
\begin{align}
 \|\sum_{i}|R_{i}f|^{2}\|_{W^{1}_{1}(S_{t,u})}\leq CA^{-1/2}\mathcal{S}_{[2]}[f]
\end{align}
We now apply (17.17) first taking $g=f^{2}$ and then taking $g=\sum_{i}|R_{i}f|^{2}$ to obtain, in view of (17.22), (17.25),
\begin{align}
 \{\int_{S_{t,u}}\{|f|^{4}+(\sum_{i}|R_{i}f|^{2})^{2}\}d\mu_{\slashed{g}}\}^{1/2}\leq C\sqrt{I'}A^{-1/2}\mathcal{S}_{[2]}[f]
\end{align}
Then by Corollary 12.2.a, 
\begin{align}
 \|f\|^{2}_{W^{4}_{1}(S_{t,u})}\leq C\sqrt{I'}A^{-3/2}\mathcal{S}_{[2]}[f]
\end{align}
where for an arbitrary function $f$ on $S_{t,u}$ we denote:
\begin{align}
 \|f\|_{W^{4}_{1}(S_{t,u})}=\|\slashed{d}f\|_{L^{4}(S_{t,u})}+A^{-1/2}\|f\|_{L^{4}(S_{t,u})}
\end{align}

    We shall now use an argument found in \cite{GT}. We rescale the metric $\slashed{g}$ on $S_{t,u}$, setting:
\begin{align}
 \hat{\slashed{g}}=A^{-1}\slashed{g}\quad \textrm{so that}\quad d\mu_{\hat{\slashed{g}}}=A^{-1}d\mu_{\slashed{g}}
\end{align}
to a metric $\hat{\slashed{g}}$ of unit area. Taking account of the fact that relative to the new metric we have:
\begin{align}
 |\slashed{d}f|^{2}_{\hat{\slashed{g}}}=(\hat{\slashed{g}}^{-1})^{BC}\partial_{B}f\partial_{C}f=
A(\slashed{g}^{-1})^{BC}\partial_{B}f\partial_{C}f=A|\slashed{d}f|^{2}_{\slashed{g}}
\end{align}
we have:
\begin{align}
 \|f\|_{W^{4}_{1}(S_{t,u},\slashed{g})}=A^{-1/4}\|f\|_{W^{4}_{1}(S_{t,u},\hat{\slashed{g}})}
\end{align}
where
\begin{align}
 \|f\|_{W^{p}_{1}(S_{t,u},\hat{\slashed{g}})}=\|\slashed{d}f\|_{L^{p}(S_{t,u},\hat{\slashed{g}})}+\|f\|_{L^{p}(S_{t,u},\hat{\slashed{g}})}
\end{align}
hence (17.27) reads:
\begin{align}
 \|f\|^{2}_{W^{4}_{1}(S_{t,u},\hat{\slashed{g}})}\leq C\sqrt{I'}A^{-1}\mathcal{S}_{[2]}[f]
\end{align}
Moreover, from (17.17), (17.18) relative to $\hat{\slashed{g}}$,
\begin{align}
 \|g\|_{L^{2}(S_{t,u},\hat{\slashed{g}})}\leq\sqrt{I'}\|g\|_{W^{1}_{1}(S_{t,u},\hat{\slashed{g}})}
\end{align}
We now set:
\begin{align}
 \tilde{f}=\frac{1}{\sqrt{I'}}\frac{|f|}{\|f\|_{W^{4}_{1}(S_{t,u},\hat{\slashed{g}})}}
\end{align}
Then $\tilde{f}\geq 0$ and taking $g=\tilde{f}^{k}$, $k>1$, in (17.34) we obtain:
\begin{align}
 \|\tilde{f}^{k}\|_{L^{2}(S_{t,u},\hat{\slashed{g}})}\leq\sqrt{I'}\|\tilde{f}^{k}\|_{W^{1}_{1}(S_{t,u},\hat{\slashed{g}})}
\end{align}
Since $\slashed{d}(\tilde{f}^{k})=k\tilde{f}^{k-1}\slashed{d}\tilde{f}$ we have,
\begin{align}
 \|\slashed{d}(\tilde{f}^{k})\|_{L^{1}(S_{t,u},\hat{\slashed{g}})}\leq
k\|\tilde{f}^{k-1}\|_{L^{4/3}(S_{t,u},\hat{\slashed{g}})}\|\slashed{d}\tilde{f}\|_{L^{4}(S_{t,u},\hat{\slashed{g}})}
\end{align}
and:
\begin{align}
 \|\tilde{f}^{k}\|_{L^{1}(S_{t,u},\hat{\slashed{g}})}=\|\tilde{f}^{k-1}\tilde{f}\|_{L^{1}(S_{t,u},\hat{\slashed{g}})}
\leq\|\tilde{f}^{k-1}\|_{L^{4/3}(S_{t,u},\hat{\slashed{g}})}\|\tilde{f}\|_{L^{4}(S_{t,u},\hat{\slashed{g}})}
\end{align}
hence:
\begin{align}
 \|\tilde{f}^{k}\|_{W^{1}_{1}(S_{t,u},\hat{\slashed{g}})}\leq k\|\tilde{f}^{k-1}\|_{L^{4/3}(S_{t,u},\hat{\slashed{g}})}\|\tilde{f}\|_{W^{4}_{1}(S_{t,u},\hat{\slashed{g}})}
\end{align}
Now by (17.35) we have:
\begin{align}
 \|\tilde{f}\|_{W^{4}_{1}(S_{t,u},\hat{\slashed{g}})}=\frac{1}{\sqrt{I'}}
\end{align}
Substituting (17.39) and using (17.36) yields:
\begin{align}
 \|\tilde{f}^{k}\|_{L^{2}(S_{t,u},\hat{\slashed{g}})}\leq k\|\tilde{f}^{k-1}\|_{L^{4/3}(S_{t,u},\hat{\slashed{g}})}
\end{align}
which is equivalent to:
\begin{align}
 \|\tilde{f}\|_{L^{2k}(S_{t,u},\hat{\slashed{g}})}\leq k^{1/k}\|\tilde{f}\|^{1-(1/k)}_{L^{(4/3)(k-1)}(S_{t,u},\hat{\slashed{g}})}
\end{align}
This implies:
\begin{align}
 \|\tilde{f}\|_{L^{2k}(S_{t,u},\hat{\slashed{g}})}\leq k^{1/k}\|\tilde{f}\|^{1-(1/k)}_{L^{(4/3)k}(S_{t,u},\hat{\slashed{g}})}
\end{align}
This is because, by virtue of the fact that $S_{t,u}$ has unit area with respect to $\hat{\slashed{g}}$, the norm $\|g\|_{L^{p}(S_{t,u},\hat{\slashed{g}})}$
is a non-decreasing function of the exponent $p$, for any given function $g$. The ratio of the exponent on the left in (17.43) to the 
exponent on the right is $3/2$. We now set:
\begin{align}
 k=(\frac{3}{2})^{n}\quad: n=1,2,3,...
\end{align}
For $n=1$ the exponent on the right in (17.43) is $2$, and taking $g=\tilde{f}$ in (17.34) we obtain, by (17.35):
\begin{align}
 \|\tilde{f}\|_{L^{2}(S_{t,u},\hat{\slashed{g}})}\leq\sqrt{I'}\|\tilde{f}\|_{W^{1}_{1}(S_{t,u},\hat{\slashed{g}})}
=\frac{\|\tilde{f}\|_{W^{1}_{1}(S_{t,u},\hat{\slashed{g}})}}{\|\tilde{f}\|_{W^{4}_{1}(S_{t,u},\hat{\slashed{g}})}}\leq 1
\end{align}
Then by (17.45) and induction on $n$ we get:
\begin{align}
 \|\tilde{f}\|_{L^{2(3/2)^{n}}(S_{t,u},\hat{\slashed{g}})}\leq(\frac{3}{2})^{\sum_{m=1}^{n}m(3/2)^{-m}}
\end{align}
Taking the limit $n\longrightarrow\infty$ we obtain:
\begin{align}
 \sup_{S_{t,u}}\tilde{f}\leq c
\end{align}
where
\begin{align*}
 c=(\frac{3}{2})^{\sum_{m=1}^{\infty}m(3/2)^{-m}}=(\frac{3}{2})^{6}
\end{align*}
In view of (17.35), (17.47) is equivalent to:
\begin{align}
 \sup_{S_{t,u}}|f|\leq c\sqrt{I'}\|f\|_{W^{4}_{1}(S_{t,u},\hat{\slashed{g}})}
\end{align}
Substituting finally (17.33) yields:
\begin{align}
 \sup_{S_{t,u}}|f|\leq cCI'^{3/4}A^{-1/2}(\mathcal{S}_{[2]}[f])^{1/2}
\end{align}

     To complete the proof of the Lemma 17.1, what remains to be done is to obtain an upper bound for $I(S_{t,u})$, the isoperimetric constant of $S_{t,u}$. The
integral curves of $T$ on a given $\Sigma_{t}^{\epsilon_{0}}$ define a diffeomorphism of $S_{t,0}$ onto each $S_{t,u}$, $u\in[0,\epsilon_{0}]$. A domain
$U_{u}\subset S_{t,u}$ with $C^{1}$ boundary $\partial U_{u}$ is mapped by the inverse onto a domain $U_{0}\subset S_{t,0}$ with $C^{1}$ boundary $\partial U_{0}$.
Consider the image $U_{u'}$ of the domain $U_{0}$ on each $S_{t,u'}$, $u'\in[0,u]$, under the diffeomorphism. Then from the definition
\begin{align*}
 \slashed{\mathcal{L}}_{T}\slashed{g}=2\kappa\theta
\end{align*}
we have:
\begin{align}
 \frac{d}{du'}\textrm{Area}(U_{u'})=\int_{U_{u'}}\kappa\textrm{tr}\theta d\mu_{\slashed{g}}
\end{align}
and:
\begin{align}
 \frac{d}{du'}\textrm{Perimeter}(\partial U_{u'})=\int_{\partial U_{u'}}\kappa \theta(V,V)ds
\end{align}
where $V$ is the unit tangent vectorfield and $ds$ is the element of arc length of $\partial U_{u'}$. Here, $u$ is the parameter of the integral curves of $T$.
Since the bootstrap assumption implies:
\begin{align}
 \kappa|\theta|\leq C(1+t)^{-1}[1+\log(1+t)]\leq C
\end{align}
from (17.50) and (17.51) we obtain:
\begin{align}
 |\frac{d}{du'}\textrm{Area}(U_{u'})|\leq C\textrm{Area}(U_{u'})\\
|\frac{d}{du'}\textrm{Perimeter}(\partial U_{u'})|\leq C\textrm{Perimeter}(\partial U_{u'})
\end{align}
Therefore integrating with respect to $u'$ on $[0,u]$ yields:
\begin{align}
 C^{-1}\textrm{Area}(U_{0})\leq \textrm{Area}(U_{u})\leq C\textrm{Area}(U_{0})\\\notag
C^{-1}\textrm{Perimeter}(\partial U_{0})\leq \textrm{Perimeter}(\partial U_{u})\leq C\textrm{Perimeter}(\partial U_{0})
\end{align}
for all $u\in[0,\epsilon_{0}]$.

It follows that:
\begin{align}
 C^{-1}I(S_{t,0})\leq I(S_{t,u})\leq C I(S_{t,0})
\end{align}
for all $u\in[0,\epsilon_{0}]$.

Since $S_{t,0}$ is a round sphere in Euclidean space, we have:
\begin{align*}
 I(S_{t,0})=\frac{1}{2\pi}
\end{align*}
So we obtain an upper bound for $I(S_{t,u})$  by a numerical constant. In view of (17.11), the lemma follows from (17.49). $\qed$

\section{Recovery of $\textbf{J}$-Completion}  
     Noting from (17.7)-(17.10) that:
\begin{align}
 \mathcal{S}_{[2]}[S\phi]\leq\tilde{\mathcal{S}}_{[4]}(t,u)\quad \mathcal{S}_{[2]}[\mathring{R}_{i}\phi]\leq\tilde{\mathcal{S}}_{[4]}(t,u)\quad:\quad
i=1,2,3\\\notag
\mathcal{S}_{[2]}[TS\phi]\leq\tilde{\mathcal{S}}_{[4]}(t,u)\quad\mathcal{S}_{[2]}[T\mathring{R}_{i}\phi]\leq\tilde{\mathcal{S}}_{[4]}(t,u)\quad:\quad
i=1,2,3
\end{align}
then by Lemma 17.1 and (17.6) we have:
\begin{align}
 \sup_{\Sigma_{t}^{\epsilon_{0}}}|S\phi|\leq C(1+t)^{-1}\sqrt{{\tilde{D}}^{u}_{[4]}}\quad
\sup_{\Sigma_{t}^{\epsilon_{0}}}|\mathring{R}_{i}\phi|\leq C(1+t)^{-1}\sqrt{{\tilde{D}}^{u}_{[4]}}\quad:\quad i=1,2,3\\\notag
\sup_{\Sigma_{t}^{\epsilon_{0}}}|TS\phi|\leq C(1+t)^{-1}\sqrt{{\tilde{D}}^{u}_{[4]}}\quad
\sup_{\Sigma_{t}^{\epsilon_{0}}}|T\mathring{R}_{i}\phi|\leq C(1+t)^{-1}\sqrt{{\tilde{D}}^{u}_{[4]}}\quad:\quad i=1,2,3
\end{align}
it follows that if:
\begin{align}
 \sqrt{{\tilde{D}}^{u}_{[4]}}\leq C\delta_{0}
\end{align}
then assumption $\textbf{J}$ holds on $W^{s}_{\epsilon_{0}}$. We thus have established assumption $\textbf{J}$.

\section{Recovery of the Final Bootstrap Assumption}
     We are now ready to recover the bootstrap assumption and proceed to the proof of the main theorem of this monograph. Let us denote by $\mathcal{S}_{[n]}(t,u)$
the integral on $S_{t,u}$ (with respect to $d\mu_{\slashed{g}}$) of the sum of the squares of all the variations (14.4), of order up to $n$. By Lemma 5.1 we have:
\begin{align}
 \mathcal{S}_{[l+1-[a]]}(t,u)\leq C\epsilon_{0}\mathcal{E}^{u}_{0,[l+1-[a]]}(t)\quad:\textrm{for all}\quad (t,u)\in[0,s]\times [0,\epsilon_{0}]
\end{align}
Hence, in view of (16.83) and (16.84):
\begin{align}
 \mathcal{S}_{[l+1-[a]]}(t,u)\leq C\epsilon_{0}\mathcal{D}^{u}_{[l+2]}\quad:\textrm{for all}\quad (t,u)\in[0,s]\times [0,\epsilon_{0}]
\end{align}
Now, for any variation $\psi$ of order up to $l-1-[a]$ we have:
\begin{align}
 \mathcal{S}_{[2]}[\psi]\leq\mathcal{S}_{[l+1-[a]]}
\end{align}
This is because $psi, R_{j_{1}}\psi : j_{1}=1,2,3, R_{j_{2}}R_{j_{1}}\psi : j_{1}, j_{2}=1,2,3$ are themselves variations included in $\mathcal{S}_{[l+1-[a]]}$.
Then by Lemma 17.1 and (17.61), (17.62), we obtain:
\begin{align}
 \sup_{S_{t,u}}|\psi|\leq C(1+t)^{-1}\sqrt{\epsilon_{0}}\sqrt{\mathcal{D}^{u}_{[l+2]}}\quad:\textrm{for all}\quad (t,u)\in[0,s]\times [0,\epsilon_{0}]
\end{align}
That is, we obtain:
\begin{align}
 \max_{\alpha}\max_{i_{1}...i_{n}}\|R_{i_{n}}...R_{i_{1}}(T)^{m}(Q)^{p}\leftexp{(\alpha)}{\psi}_{1}\|_{L^{\infty}(\Sigma_{t}^{\epsilon_{0}})}
\leq C(1+t)^{-1}\sqrt{\epsilon_{0}}\sqrt{\mathcal{D}_{[l+2]}}
\end{align}
for all $p+m+n\leq l-2-[a]$ and all $t\in[0,s]$.

From (15.167), we have:
\begin{align*}
 l-2-[a]\geq (l+1)_{*}+2
\end{align*}
since $l=l_{*}+(l+1)_{*}$. So if 
\begin{align}
 C\sqrt{\epsilon_{0}}\sqrt{\mathcal{D}_{[l+2]}}<\delta_{0}
\end{align}
then we recover the bootstrap assumption $\textbf{E}_{\{\{(l+1)_{*}+2\}\}}$.

\section{Completion of the Proof of the Main Theorem}
     Let now $s_{*}$ be the least upper bound of the set of values of $s$ in the interval $[0,t_{*\epsilon_{0}}]$ such that the bootstrap assumption holds on 
$W^{s}_{\epsilon_{0}}$. We recall from Chapter 2 that $t_{*\epsilon_{0}}$ is defined by:
\begin{align}
 t_{*\epsilon_{0}}=\inf_{u\in[0,\epsilon_{0}]}t_{*}(u)
\end{align}
where $t_{*}(u)$ is the greatest lower bound of the extent of the generators of $C_{u}$, in the parameter $t$, in the domain of the maximal solution.
We note here that $C_{u}$ do not contain cut loci. This follows from the fact that the angle between the outward unit normal $-\hat{T}$ to the surface $S_{t,u}$
with respect to the Euclidean metric on $\Sigma_{t}$ and the outward unit normal $N$ to the Euclidean coordinate spheres does not exceed a fixed constant times
 $\delta_{0}$, as follows from the estimate (12.19). In view of the bound on $\chi$, the second fundamental form of the sections $S_{t,u}$ relative to $C_{u}$, 
the $C_{u}$ do not contain focal points (that is, points along a generator of $C_{u}$ which are conjugate to $S_{0,u}$) either. The absence of cut locus or focal 
points implies that a bicharacteristic generator of $C_{u}$ cannot leave the boundary of the domain of dependence, in the domain of maximal solution, of the
 exterior of the surface $S_{0,u}$ in the initial hyperplane $\Sigma_{0}$. (for the notion of cut locus and of focal or conjugate points in Riemannian 
geometry, see \cite{CE}. For the corresponding notions in Lorentzian geometry see \cite{Pen}.) Thus, unless $t_{*\epsilon_{0}}=\infty$, there is on 
$\Sigma_{t_{*\epsilon_{0}}}^{\epsilon_{0}}$ at least one point which belongs to the boundary of the domain of maximal solution and not to the domain itself. We 
shall presently show that in fact $s_{*}$ coincides with $t_{*\epsilon_{0}}$ under the condition (17.65). For, otherwise, that is if $s_{*}<t_{*\epsilon_{0}}$, 
then by continuity the bootstrap assumption holds on $W^{s_{*}}_{\epsilon_{0}}$ as well, hence by the above, (17.64) holds with $l-2-[a]$ replaced by $(l+1)_{*}+2$ 
and $s$ replaced by $s_{*}$. By (17.65) and continuity however the bootstrap assumption must also hold for some $s>s_{*}$, contradicting the definition of $s_{*}$.
 We conclude that $s_{*}=t_{*\epsilon_{0}}$ and (17.64) holds for all $t\in[0,t_{*\epsilon_{0}})$. We thus have uniform pointwise estimates for the $\psi_{\alpha}$ 
on $W^{*}_{\epsilon_{0}}$ up to order $(l+1)_{*}+2$. Proposition 12.9 and 12.10 then give uniform pointwise estimates for $\chi'$ up to order $(l+1)_{*}$ and for 
$\mu$ up to order $(l+1)_{*}+1$. It follows that the $\psi_{\alpha}$, $\chi'$ and $\mu$ thus also $\kappa$ and the induced acoustical metric $\slashed{g}$, 
extend smoothly in $acoustical$ $coordinates$ to $\Sigma_{t_{*\epsilon_{0}}}^{\epsilon_{0}}$ and $W^{t_{*\epsilon_{0}}}_{\epsilon_{0}}$. Also the rectangular 
components $g_{\mu\nu}$ of the acoustical metric, being functions of $\psi_{\alpha}$, likewise extend. Moreover, since the bootstrap assumption holds 
on $W^{*}_{\epsilon_{0}}$, all the estimates we have derived, in particular the energy estimates (16.83), (16.79), (15.160)-(15.163), and the $L^{2}$ acoustical 
estimates of Proposition 12.11, 12.12, as well as the top order acoustical estimates, hold for all $t\in[0,t_{*\epsilon_{0}})$. The estimates not containing the 
weights $\bar{\mu}_{m,u}$, such as the estimates (16.83), then extend to $t=t_{*\epsilon_{0}}$ as well. Now, we must have:
\begin{align}
 \bar{\mu}_{m}(t_{*\epsilon_{0}})=0
\end{align}
 For otherwise the Jacobian determinant $\triangle$ of the transformation from acoustical to rectangular coordinates has a positive minimum on 
$\Sigma_{t_{*\epsilon_{0}}}^{\epsilon_{0}}$, therefore the inverse transform is regular and the $\psi_{\alpha}$ extend smoothly in $rectangular$ $coordinates$ 
to $\Sigma_{t_{*\epsilon_{0}}}^{\epsilon_{0}}$. However, once the $\psi_{\alpha}$ extend to functions of the rectangular coordinates on 
$\Sigma_{t_{*\epsilon_{0}}}^{\epsilon_{0}}$ which belong to the Sobolev space $H_{3}$, then the standard local existence theorem applies and we obtain an 
extension of the solution to a development containing an extension of all the characteristic hypersurfaces $C_{u}$, $u\in[0,\epsilon_{0}]$, up to a value 
$t_{1}$ of $t$ for some $t_{1}>t_{*}$, in contradiction with the definition of $t_{*\epsilon_{0}}$. We conclude that there is at least one point on 
$\Sigma_{t_{*\epsilon_{0}}}^{\epsilon_{0}}$ where $\mu$ vanishes. By the same argument, at each point $x_{*}\in\Sigma_{t_{*\epsilon_{0}}}^{\epsilon_{0}}$ which
 lies on the boundary of the domain of the maximal solution, $\mu(x_{*})=0$. For, otherwise, $\mu$ has a positive minimum in a suitable neighborhood of $x_{*}$, 
hence the solution is locally extendible at $x_{*}$ in contradiction with the fact that $x_{*}$ lies on the boundary of the domain of the maximal solution.

     Now, from Proposition 8.6, taking $t=s$ we have:
\begin{align}
 \hat{\mu}_{s}(s,u,\vartheta)=1+\hat{E}_{s}(u,\vartheta)\log(1+s)
\end{align}
 and we recall that:
\begin{align}
 \hat{\mu}_{s}(s,u,\vartheta)=\frac{\mu(s,u,\vartheta)}{\mu_{[1],s}(u,\vartheta)},\quad \hat{E}_{s}(u,\vartheta)=
\frac{E_{s}(u,\vartheta)}{\mu_{[1],s}(u,\vartheta)}
\end{align}
and that $\mu_{[1],s}(u,\vartheta)$ is bounded from below by $1/2$. We also recall that according to (8.226):
\begin{align}
 |E_{s}(u,\vartheta)-\frac{1}{4}\ell P_{s}(u,\vartheta)|\leq C\delta_{0}(1+s)^{-1}[1+\log(1+s)]
\end{align}
where, from Lemma 8.10:
\begin{align}
 P_{s}(u,\vartheta)=(1+s)(\underline{L}\psi_{0})(s,u,\vartheta)
\end{align}
Now by (17.64) we have:
\begin{align}
 |P_{s}(u,\vartheta)|\leq C\sqrt{\epsilon_{0}}\sqrt{\mathcal{D}_{[l+2]}}
\end{align}
Hence we obtain:
\begin{align}
 |E_{s}(u,\vartheta)|\leq C|\ell|\sqrt{\epsilon_{0}}\sqrt{\mathcal{D}_{[l+2]}}+C\delta_{0}(1+s)^{-1}[1+\log(1+s)]
\end{align}
Since
\begin{align*}
 \frac{\log(1+s)[1+\log(1+s)]}{(1+s)}\leq C
\end{align*}
this implies:
\begin{align}
 \log(1+s)|E_{s}(u,\vartheta)|\leq C|\ell|\sqrt{\epsilon_{0}}\sqrt{\mathcal{D}_{[l+2]}}\log(1+s)+C\delta_{0}
\end{align}
Substituting in (17.68) then yields:
\begin{align}
 \hat{\mu}_{s}(s,u,\vartheta)\geq 1-C\delta_{0}-C|\ell|\sqrt{\epsilon_{0}}\sqrt{\mathcal{D}_{[l+2]}}\log(1+s)
\end{align}
It follows that:
\begin{align}
 \log(1+t_{*\epsilon_{0}})\geq\frac{1}{C|\ell|\sqrt{\epsilon_{0}}\sqrt{\mathcal{D}_{[l+2]}}}
\end{align}
(for a new constant $C$), because otherwise $\mu$ would be bounded from below by a positive constant on $\Sigma_{t_{*\epsilon_{0}}}^{\epsilon_{0}}$
contradicting the conclusions of the previous paragraph.

It remains for us to analyze the requirements on the initial data. We recall from Chapter 2 that $u$ on $\Sigma_{0}$ is defined by:
\begin{align}
 u=1-r
\end{align}
Thus, in rectangular coordinates we have on $\Sigma_{0}$:
\begin{align}
 \partial_{i}u=-N^{i},\quad N^{i}=\frac{x^{i}}{r}
\end{align}
where $N$ is the Euclidean outward unit normal to the Euclidean coordinate spheres. Since $\kappa$ is defined by:
\begin{align}
 \kappa^{-2}=(\bar{g}^{-1})^{ij}\partial_{i}u\partial_{j}u
\end{align}
and $\bar{g}_{ij}=\delta_{ij}$, on $\Sigma_{0}$ we have simply:
\begin{align}
 \kappa=1
\end{align}
and 
\begin{align}
 \mu=\eta
\end{align}

According to (2.27) we have:
\begin{align}
 \hat{T}^{i}=\kappa\partial_{i}u
\end{align}
then on $\Sigma_{0}$,
\begin{align}
 \hat{T}^{i}=-N^{i}=T^{i}
\end{align}
Then by (2.55), $L$ is given on $\Sigma_{0}$ by:
\begin{align}
 L=\frac{\partial}{\partial x^{0}}+(\eta N^{i}-\psi_{i})\frac{\partial}{\partial x^{i}}
\end{align}
On $\Sigma_{0}$, $Q=L$. Consider next the orthogonal projection to $S_{t,u}$ in $\Sigma_{t}$, relative to $\bar{g}$. It is given in rectangular
coordinates by:
\begin{align}
 \Pi^{a}_{b}=\delta^{a}_{b}-\hat{T}^{a}\hat{T}^{b}
\end{align}
The level surfaces $S_{0,u}$ of $u$ on $\Sigma_{0}$ being Euclidean spheres centered at the origin and the vectorfields $\mathring{R}_{i}$ being tangential 
to these spheres, the commutation fields $R_{i}$ coincide on $\Sigma_{0}$ with $\mathring{R}_{i}$:
\begin{align}
 R_{i}=\Pi\mathring{R}_{i}=\mathring{R}_{i}=\epsilon_{ijk}x^{j}\frac{\partial}{\partial x^{k}}
\end{align}
Moreover, the functions $\lambda_{i}$ defined by (12.12) vanishes on $\Sigma_{0}$:
\begin{align}
 \lambda_{i}=-\kappa N^{j}\epsilon_{ijk}x^{k}=0
\end{align}
in view of (17.83).
    
Since the $S_{0,u}$ are Euclidean spheres of radius $r$ with inward unit normal $\hat{T}$, we have:
\begin{align}
 \theta_{ab}=-\frac{1}{r}(\delta_{ab}-\hat{T}^{a}\hat{T}^{b})
\end{align}
Then from (3.27),
\begin{align}
 \chi_{ab}=-\eta(\theta_{ab}-\slashed{k}_{ab})
\end{align}
where
\begin{align}
 \slashed{k}_{ab}=-\eta^{-1}\Pi^{i}_{a}\Pi^{j}_{b}\partial_{i}\psi_{j}
\end{align}

     Consider now the initial data for the nonlinear wave equation (1.23). The functions $\psi_{\alpha}$ : $\alpha=0,1,2,3$ are given in the exterior of 
the sphere of radius $1-\epsilon_{0}$ with center at the origin and satisfy $\partial_{i}\psi_{j}=\partial_{j}\psi_{i}$. Outside the unit sphere with center 
at the origin, the initial data coincide with those of a constant state $\psi_{\alpha}=0$ : $\alpha=0,1,2,3$. We assume that $\psi_{\alpha}$
: $\alpha=0,1,2,3$ are, in rectangular coordinates, functions belonging to the Sobolev space $H_{l+2}(\Sigma_{0}^{\epsilon_{0}})$ with vanishing traces on 
the unit sphere. We set:
\begin{align}
 D_{[l+2]}=\sum_{\alpha,i}\|\partial_{i}\psi_{\alpha}\|^{2}_{H_{l+1}(\Sigma_{0}^{\epsilon_{0}})}
\end{align}
The nonlinear wave equation in the form:
\begin{align}
 (g^{-1})^{\mu\nu}\partial_{\mu}\psi_{\nu}=0,\quad \partial_{\mu}\psi_{\nu}=\partial_{\nu}\psi_{\mu}
\end{align}
allows us to express $\partial_{0}\psi_{0}$ and $\partial_{0}\psi_{i}$ in terms of $\partial_{i}\psi_{0}$ and $\partial_{i}\psi_{j}$, for, $\partial_{0}\psi_{i}
=\partial_{i}\psi_{0}$ and we have:
\begin{align}
 (g^{-1})^{00}=-\eta^{-2}
\end{align}
We can then express recursively $\partial^{k}_{0}\psi_{\alpha}$ : $\alpha=0,1,2,3$ for $k=1,...,l+2$ in terms of the data $\partial_{i_{k}}...\partial_{i_{1}}\psi_{\alpha}$
: $\alpha=0,1,2,3$; $i_{1},...,i_{n}=1,2,3$. Now the standard Sobolev inequalities yield:
\begin{align}
\max_{\alpha}\max_{i_{1}...i_{n}}\|\partial_{i_{n}}...\partial_{i_{1}}\psi_{\alpha}\|_{L^{\infty}(\Sigma_{0}^{\epsilon_{0}})}\\\notag
\leq C\sqrt{\epsilon_{0}}\sqrt{D_{[l+2]}}\quad :\textrm{for}\quad n=1,...,l-1
\end{align}
It then follows that also:
\begin{align}
\max_{\alpha}\max_{i_{1}...i_{n}}\|\partial_{i_{n}}...\partial_{i_{1}}\partial_{0}^{k}\psi_{j}\|_{L^{\infty}(\Sigma_{0}^{\epsilon_{0}})}\\\notag
\leq C\sqrt{\epsilon_{0}}\sqrt{D_{[l+2]}}\quad :\textrm{for}\quad n,k\geq 0, n+k\leq l-1
\end{align}
Here, the constant $\sqrt{\epsilon_{0}}$ comes from using Lemma 5.1.

Now from (17.83), (17.85) and (17.87) the rectangular components of $T$, $L$ and $R_{i}$ are on $\Sigma_{0}^{\epsilon_{0}}$ smooth 
functions of the rectangular coordinates and the $\psi_{\alpha}$. It then follows that:
\begin{align}
 \mathcal{E}_{0,[l+2]}(0)\quad, \quad \mathcal{E}^{\prime}_{1,[l+2]}(0)\leq CD_{[l+2]}
\end{align}
From (17.94) and (17.95) we have:
\begin{align}
 \|\chi'\|_{\infty,\{l-2\},\Sigma_{0}^{\epsilon_{0}}}\leq C\sqrt{\epsilon_{0}}\sqrt{D_{[l+2]}}
\end{align}

The assumptions of Propositions 12.3, 12.6, 12.9, 12.10, on the initial conditions are then recovered provided that:
\begin{align}
 \sqrt{\epsilon_{0}}\sqrt{D_{[l+2]}}\leq C^{-1}\delta_{0}
\end{align}
On $\Sigma_{0}^{\epsilon_{0}}$, we have $y^{i}=0$, hence
\begin{align}
 \mathcal{Y}_{0}(0)=0
\end{align}
Moreover,
\begin{align}
\mathcal{B}_{\{l+1\}}(0)\leq C\sqrt{D_{[l+2]}}\\
\mathcal{A}_{[l]}(0)\leq C\sqrt{D_{[l+2]}}
\end{align}
and by (8.40), (9.71) we have:
\begin{align}
 \sum_{i_{1}...i_{l}}\|\leftexp{(i_{1}...i_{l})}{x}_{l}(0)\|_{L^{2}(\Sigma_{0}^{\epsilon_{0}})}\leq C\sqrt{D_{[l+2]}}\\
\sum_{m=0}^{l}\sum_{i_{1}...i_{l-m}}\|\leftexp{(i_{1}...i_{l-m})}{x}^{\prime}_{m,l-m}(0)\|_{L^{2}(\Sigma_{0}^{\epsilon_{0}})}
\leq C\sqrt{D_{[l+2]}}
\end{align}
 Combining the above then yields:
\begin{align}
 \mathcal{D}_{[l+2]}\leq CD_{[l+2]}
\end{align}
and we can replace $\epsilon_{0}$ by any $u\in(0,\epsilon_{0}]$ to obtain:
\begin{align}
 \mathcal{D}^{u}_{[l+2]}\leq CD^{u}_{[l+2]}
\end{align}
for all $u\in(0,\epsilon_{0}]$, where
\begin{align}
 D^{u}_{[l+2]}=\sum_{\alpha,i}\|\partial_{i}\psi_{\alpha}\|^{2}_{H_{l+1}(\Sigma^{u}_{0})}
\end{align}

     In conclusion we have proved the following theorem, which is the main theorem of this monograph.

$\textbf{Theorem 17.1}$ Let $(p,s,v^{i}\quad:\quad i=1,2,3)$ be initial data for the equation (1.4), (1.10) and (1.12) on $\Sigma_{0}$ which correspond 
to the initial data of a constant state
\begin{align*}
 p=p_{0},\quad s=s_{0},\quad v^{i}=0:\quad i=1,2,3
\end{align*}
 outside the unit sphere with center at the origin in $\Sigma_{0}$. We can adjust the zero-point of the enthalpy $h$ so that $h_{0}=0$, that is, the enthalpy 
vanishes in the constant state. We can also choose the relation of the unit of time to the unit of length so that $\eta_{0}=1$, that is, the sound speed in the 
constant state is equal to unity. Let also the initial data be irrotational and isentropic outside the sphere of radius 
$1-\epsilon_{0}$, $0<\epsilon_{0}\leq 1/2$ with center at the origin in $\Sigma_{0}$. Then we have initial data $(\phi,\partial_{t}\phi)$ for the nonlinear
wave equation (1.23) outside the sphere  of radius $1-\epsilon_{0}$ with center at the origin in $\Sigma_{0}$, where $\partial_{i}\phi=-v^{i}$ and 
$\partial_{t}\phi=h+\frac{1}{2}\sum_{i=1}^{3}(\partial_{i}\phi)^{2}$. The initial data $(\phi,\partial_{t}\phi)$ vanish
outside the unit sphere with center at the origin in $\Sigma_{0}$. Consider the annular region $\Sigma_{0}^{\epsilon_{0}}$ in
$\Sigma_{0}$ bounded by the two concentric spheres. Then there is a positive integer $[a]$ such that the following hold. Let $l$ be a positive integer such that
\begin{align*}
 l_{*}\geq [a]+4
\end{align*}
and suppose that the functions $\psi_{\alpha}$ : $\alpha=0,1,2,3$ corresponding to the initial data on $\Sigma_{0}^{\epsilon_{0}}$ are functions of the rectangular 
coordinates belonging to the Sobolev space $H_{l+2}(\Sigma_{0}^{\epsilon_{0}})$ with vanishing traces on the unit sphere. Then setting
\begin{align*}
 D_{[l+2]}=\sum_{\alpha,i}\|\partial_{i}\psi_{\alpha}\|^{2}_{H_{l+1}(\Sigma_{0}^{\epsilon_{0}})}
\end{align*}
there is a suitably small positive constant $\bar{\delta}_{0}$ and a positive constant $C$, such that for any $\delta_{0}\in(0,\bar{\delta}_{0}]$, if:
\begin{align*}
 C\sqrt{\epsilon_{0}}\sqrt{D_{[l+2]}}<\delta_{0}
\end{align*}
the following conclusions hold:

(i) Let $u$ be the function $1-r$ on $\Sigma_{0}$ and $S_{0,u}$ the spheres of radius $1-u$ with center at the origin, the level surfaces of $u$ in $\Sigma_{0}$.
We consider, in the domain of the maximal solution corresponding to the given initial data, the family $\{C_{u}\quad :\quad u\in[0,\epsilon_{0}]\}$ of outgoing 
characteristic hypersurfaces corresponding to the family $\{S_{0,u}\quad:\quad u\in[0,\epsilon_{0}]\}$:
\begin{align*}
 C_{u}\bigcap \Sigma_{0}=S_{0,u}\quad:\quad\forall u\in[0,\epsilon_{0}]
\end{align*}
with each bicharacteristic generator of each $C_{u}$ extending in the domain of the maximal solution as long as it remains on the boundary of the domain of dependence
of the exterior of $S_{0,u}$ in $\Sigma_{0}$. Then the bicharacteristic generators of each $C_{u}$ have no future end points except on the boundary of the domain of the 
maximal solution. Let $t_{*}(u)$ be the least upper bound of the extent of the generators of $C_{u}$, in the parameter $t$, in the domain of the maximal solution, and
let:
\begin{align*}
 t_{*\epsilon_{0}}=\inf_{u\in[0,\epsilon_{0}]}t_{*}(u)
\end{align*}
 We define for each $(t,u)\in[0,t_{*\epsilon_{0}})\times[0,\epsilon_{0}]$ the closed surface:
\begin{align*}
 S_{t,u}=C_{u}\bigcap \Sigma_{t}
\end{align*}
Then either $t_{*\epsilon_{0}}=\infty$ or there is on $\Sigma_{t_{*\epsilon_{0}}}^{\epsilon_{0}}$ at least one point which belongs to the boundary of the domain of the
maximal solution and not to the domain itself.

(ii) We have the lower bound:
\begin{align*}
 \log(1+t_{*\epsilon_{0}})\geq\frac{1}{C|\ell|\sqrt{\epsilon_{0}}\sqrt{D_{[l+2]}}}
\end{align*}
where $\ell$ is the constant:
\begin{align*}
 \ell=\frac{dH}{dh}(h_{0})
\end{align*}
In particular, if $\ell=0$ we have $t_{*\epsilon_{0}}=\infty$.

(iii) For all $t\in[0,t_{*\epsilon_{0}})$ we have the $L^{\infty}$ bounds:
\begin{align*}
 \max_{\alpha}\max_{i_{1}...i_{n}}\|R_{i_{n}}...R_{i_{1}}(T)^{m}(Q)^{p}\psi_{\alpha}\|_{L^{\infty}(\Sigma_{t}^{\epsilon_{0}})}
\leq C(1+t)^{-1}\sqrt{\epsilon_{0}}\sqrt{D_{[l+2]}}\\
\quad:\textrm{for all}\quad p+m+n\leq l-2-[a]
\end{align*}
and with $\mu$ and $\chi$ the acoustical entities defined by the family $\{C_{u}\quad:\quad u\in[0,\epsilon_{0}]\}$ and for each $C_{u}$ : $u\in[0,\epsilon_{0}]$ 
by the family of sections $\{S_{t,u}\quad:\quad t\in[0,t_{*\epsilon_{0}})\}$:
\begin{align*}
 \|\mu-1\|_{\infty,\{l-3-[a]\},\Sigma_{t}^{\epsilon_{0}}}\leq C\delta_{0}[1+\log(1+t)]\\
\|\chi'\|_{\infty,\{l-4-[a]\},\Sigma_{t}^{\epsilon_{0}}}\leq C\delta_{0}(1+t)^{-2}[1+\log(1+t)]
\end{align*}
where
\begin{align*}
 \chi'=\chi-\frac{\slashed{g}}{1-u+t}
\end{align*}
and $\slashed{g}$ is the induced acoustical metric on $S_{t,u}$.

(iv) The $\psi_{\alpha}$, $\mu$, $\slashed{g}$ and $\chi$, extend continuously with their first $l-2-[a]$ derivatives in the case of $\psi_{\alpha}$,
 $l-3-[a]$ derivatives in the case of $\mu$, $\slashed{g}$ and $l-4-[a]$ derivatives in the case of $\chi$, in acoustical coordinates to 
$\Sigma_{t_{*\epsilon_{0}}}$. The rectangular components $\hat{T}^{i}$ and $L^{i}$ of the vectorfields $T$ and $L$ likewise extend continuously with their first 
$l-3-[a]$ derivatives in acoustical coordinates to $\Sigma_{t_{*\epsilon_{0}}}$, and so do the rectangular components $g_{\mu\nu}$ of the acoustical spacetime metric.
The function $\mu$ so extended  vanishes at each point in $\Sigma_{t*\epsilon_{0}}^{\epsilon_{0}}$ which lies on the boundary of the domain of the maximal
solution, being positive everywhere else. The Jacobian determinant $\triangle$ of the transformation from acoustical to rectangular coordinates also vanishes at the same 
points, being positive everywhere else. The first derivatives $\hat{T}^{i}\partial_{i}\psi_{\alpha}$ blow up at these points. More precisely what blows up is the 
component $\hat{T}^{i}\hat{T}^{j}\partial_{i}\psi_{j}$. Moreover, there is a positive constant $C$ such that in the subdomain $\mathcal{U}$ of $W^{*}_{\epsilon_{0}}$ 
where $\mu<1/4$, which contains a spacetime neighborhood of each of the points in question, we have:
\begin{align*}
 L\mu\leq-C^{-1}(1+t)^{-1}[1+\log(1+t)]^{-1}
\end{align*}

(v) The following energy estimates hold, for all $u\in[0,\epsilon_{0}]$:
\begin{align*}
 \sup_{t\in[0,t_{*\epsilon_{0}}]}\{\mathcal{E}^{u}_{0,[l+1-[a]]}(t)\}\leq CD^{u}_{[l+2]}\\
\mathcal{F}^{t_{*\epsilon_{0}}}_{0,[l+1-[a]]}(u)\leq CD^{u}_{[l+2]}\\
\sup_{t\in[0,t_{*\epsilon_{0}}]}\{[1+\log(1+t)]^{-4}\mathcal{E}^{\prime u}_{1,[l+1-[a]]}(t)\}\leq CD^{u}_{[l+2]}\\
\sup_{t\in[0,t_{*\epsilon_{0}}]}\{[1+\log(1+t)]^{-4}\mathcal{F}^{\prime t}_{1,[l+1-[a]]}(u)\}\leq CD^{u}_{[l+2]}\\
\sup_{t\in[0,t_{*\epsilon_{0}}]}\{[1+\log(1+t)]^{-4}K_{[l+1-[a]]}(t,u)\}\leq CD^{u}_{[l+2]}
\end{align*}
where
\begin{align*}
 D^{u}_{[l+2]}=\sum_{\alpha,i}\|\partial_{i}\psi_{\alpha}\|^{2}_{H_{l+1}(\Sigma^{u}_{t})}
\end{align*}
and $\Sigma_{0}^{u}$ is the annular region on $\Sigma_{0}$ bounded by $S_{0,u}$ and the unit sphere $S_{0,0}$. Moreover, setting
\begin{align*}
 a=[a]+\frac{3}{4}
\end{align*}
we have, for each $n=0,...,[a]-1$, and for all $u\in[0,\epsilon_{0}]$, the estimates:
\begin{align*}
 \sup_{t\in[0,t_{*\epsilon_{0}})}\{\bar{\mu}^{2(a-n-1)}_{m,u}(t)\mathcal{E}^{u}_{0,[l+1-n]}(t)\}\leq CD^{u}_{[l+2]}\\
\sup_{t\in[0,t_{*\epsilon_{0}})}\{\bar{\mu}^{2(a-n-1)}_{m,u}(t)\mathcal{F}^{t}_{0,[l+1-n]}(u)\}\leq CD^{u}_{[l+2]}\\
\sup_{t\in[0,t_{*\epsilon_{0}})}\{\bar{\mu}^{2(a-n-1)}_{m,u}(t)[1+\log(1+t)]^{-4}\mathcal{E}^{\prime u}_{1,[l+1-n]}(t)\}\leq CD^{u}_{[l+2]}\\
\sup_{t\in[0,t_{*\epsilon_{0}})}\{\bar{\mu}^{2(a-n-1)}_{m,u}(t)[1+\log(1+t)]^{-4}\mathcal{F}^{\prime t}_{1,[l+1-n]}(u)\}\leq CD^{u}_{[l+2]}\\
\sup_{t\in[0,t_{*\epsilon_{0}})}\{\bar{\mu}^{2(a-n-1)}_{m,u}(t)[1+\log(1+t)]^{-4}K_{[l+1-n]}(t,u)\}\leq CD^{u}_{[l+2]}\\
\end{align*}
Furthermore, there is a positive real number $p$ such that with $q=p+2$ we have the top order energy estimates:
\begin{align*}
 \sup_{t\in[0,t_{*\epsilon_{0}})}\{\bar{\mu}^{2a}_{m,u}(t)[1+\log(1+t)]^{-2p}\mathcal{E}^{u}_{0,[l+2]}(t)\}\leq CD^{u}_{[l+2]}\\
\sup_{t\in[0,t_{*\epsilon_{0}})}\{\bar{\mu}^{2a}_{m,u}(t)[1+\log(1+t)]^{-2p}\mathcal{F}^{t}_{0,[l+2]}(u)\}\leq CD^{u}_{[l+2]}\\
\sup_{t\in[0,t_{*\epsilon_{0}})}\{\bar{\mu}^{2a}_{m,u}(t)[1+\log(1+t)]^{-2q}\mathcal{E}^{\prime u}_{1,[l+2]}(t)\}\leq CD^{u}_{[l+2]}\\
\sup_{t\in[0,t_{*\epsilon_{0}})}\{\bar{\mu}^{2a}_{m,u}(t)[1+\log(1+t)]^{-2q}\mathcal{F}^{\prime t}_{1,[l+2]}(u)\}\leq CD^{u}_{[l+2]}\\
\sup_{t\in[0,t_{*\epsilon_{0}})}\{\bar{\mu}^{2a}_{m,u}(t)[1+\log(1+t)]^{-2q}K_{[l+2]}(t,u)\}\leq CD^{u}_{[l+2]}
\end{align*}
(vi) The following $L^{2}$ acoustical estimates hold, for all $u\in[0,\epsilon_{0}]$ and all $t\in[0,t_{*\epsilon_{0}}]$:
\begin{align*}
 \mathcal{A'}^{u}_{[l-[a]-1]}(t)\leq C(1+t)^{-1}[1+\log(1+t)]^{3}\sqrt{D^{u}_{[l+2]}}\\
\mathcal{B}^{u}_{[l-[a]]}(t)\leq C(1+t)[1+\log(1+t)]^{3}\sqrt{D^{u}_{[l+2]}}
\end{align*}
Moreover, for each $n=1,...,[a]$ and for all $u\in[0,\epsilon_{0}]$ and all $t\in[0,t_{*\epsilon_{0}})$:
\begin{align*}
 \mathcal{A}^{\prime u}_{[l-n]}(t)\leq C(1+t)^{-1}[1+\log(1+t)]^{3}\bar{\mu}^{-a+n+1/2}_{m,u}(t)\sqrt{D^{u}_{[l+2]}}\\
\mathcal{B}^{u}_{\{l-n+1\}}(t)\leq C(1+t)[1+\log(1+t)]^{3}\bar{\mu}^{-a+n+1/2}_{m,u}(t)\sqrt{D^{u}_{[l+2]}}
\end{align*}
and for $n=0$ and for all $u\in[0,\epsilon_{0}]$ and all $t\in[0,t_{*\epsilon_{0}})$:
\begin{align*}
 \mathcal{A}^{\prime u}_{[l]}(t)\leq C(1+t)^{-1}[1+\log(1+t)]^{q+1}\bar{\mu}^{-a+1/2}_{m,u}(t)\sqrt{D^{u}_{[l+2]}}\\
\mathcal{B}^{u}_{\{l+1\}}(t)\leq C(1+t)[1+\log(1+t)]^{q+1}\bar{\mu}^{-a+1/2}_{m,u}(t)\sqrt{D^{u}_{[l+2]}}
\end{align*}

Finally, for all $u\in[0,\epsilon_{0}]$ and all $t\in[0,t_{*\epsilon_{0}})$ we have the top order acoustical estimates:
\begin{align*}
 \max_{i_{1}...i_{l+1}}\|\mu R_{i_{l+1}}...R_{i_{1}}\textrm{tr}\chi'\|_{L^{2}(\Sigma_{t}^{u})}\leq C\bar{\mu}^{-2a}_{m,u}(t)[1+\log(1+t)]^{2p}\sqrt{D^{u}_{[l+2]}}\\
\sum_{m=0}^{l}\max_{i_{1}...i_{l-m}}\|\mu R_{i_{l-m}}...R_{i_{1}}(T)^{m}\slashed{\Delta}\mu\|_{L^{2}(\Sigma_{t}^{u})}\leq C\bar{\mu}^{-2a}_{m,u}(t)[1+\log(1+t)]^{2p}
\sqrt{D^{u}_{[l+2]}}
\end{align*}

\chapter{Sufficient Conditions on the Initial Data for the Formation of a Shock in the Evolution}

In the present chapter we shall establish sharp sufficient conditions on the initial data for the formation of a shock in the evolution. We just investigate the problem in
the isentropic irrotational case.

     The set up is as in Theorem 17.1. Following this theorem we must find sharp sufficient conditions on the initial data which will guarantee that the function $\mu$ 
becomes zero somewhere on $\Sigma_{t}^{\epsilon_{0}}$ at some finite $t$. That value of $t$ shall then be $t_{*\epsilon_{0}}$. We shall use Lemma 8.10, Proposition 8.5,
Proposition 8.6. Taking $t=s$ in Proposition 8.6 we obtain:
\begin{align}
 \mu(s,u,\vartheta)=\mu_{[1],s}(u,\vartheta)\hat{\mu}_{s}(s,u,\vartheta)
\end{align}
where $\mu_{[1],s}(u,\vartheta)$ satisfies:
\begin{align}
 \frac{1}{2}\leq\mu_{[1],s}(u,\vartheta)\leq\frac{3}{2}
\end{align}
while $\hat{\mu}_{s}(s,u,\vartheta)$ is given by:
\begin{align}
 \hat{\mu}_{s}(s,u,\vartheta)=1+\hat{E}_{s}(u,\vartheta)\log(1+s)
\end{align}
Here:
\begin{align}
 \hat{E}_{s}(u,\vartheta)=\frac{E_{s}(u,\vartheta)}{\mu_{[1],s}(u,\vartheta)}
\end{align}
and according to (8.226):
\begin{align}
 |E_{s}(u,\vartheta)-\frac{1}{4}\ell P_{s}(u,\vartheta)|\leq C\delta_{0}(1+s)^{-1}[1+\log(1+s)]
\end{align}
where $P_{s}(u,\vartheta)$ is the function defined by Lemma 8.10:
\begin{align}
 P_{s}(u,\vartheta)=(1+s)(\underline{L}\psi_{0})(s,u,\vartheta)
\end{align}
In (18.5) and in all the following the positive constant $\delta_{0}$ shall be as in Theorem 17.1. In view of the above, a sharp sufficient condition for $\mu$ to become 
zero somewhere on $\Sigma_{s}^{\epsilon_{0}}$ at some finite $s$, is an upper bound for $\min_{(u,\vartheta)\in[0,\epsilon_{0}]\times S^{2}}P_{s}(u,\vartheta)$ by a 
negative constant for all sufficiently large $s$ in the case that $\ell>0$, a lower bound for $\max_{(u,\vartheta)\in[0,\epsilon_{0}]\times S^{2}}P_{s}(u,\vartheta)$ by 
a positive constant for all sufficiently large $s$ in the case that $\ell <0$. As was already shown in Theorem 17.1 no shocks can form in the case $\ell=0$ as in this case
we have $t_{*\epsilon_{0}}=\infty$. Thus in the following we assume $\ell\slashed{=}0$.

     The sharp sufficient condition just stated is not a condition on the initial data. Let us revisit the wave equation for $\psi_{0}$ in the form encountered in the proof 
of Lemma 8.10 ((8.180) and (8.181)):
\begin{align}
 L(\underline{L}\psi_{0})+\nu\underline{L}\psi_{0}+\underline{\nu}L\psi_{0}=\rho'_{0}
\end{align}
 where:
\begin{align}
 \rho'_{0}=\mu\slashed{\Delta}\psi_{0}+\mu\slashed{g}^{-1}(\slashed{d}\log\Omega,\slashed{d}\psi_{0})-2\slashed{g}^{-1}(\zeta,\slashed{d}\psi_{0})
\end{align}
Consider the function:
\begin{align}
 \underline{\tau}=(1-u+t)\underline{L}\psi_{0}-\psi_{0}
\end{align}
We have:
\begin{align}
 L\underline{\tau}=(1-u+t)L(\underline{L}\psi_{0})-L\psi_{0}+\underline{L}\psi_{0}
\end{align}
Substituting for $L(\underline{L}\psi_{0})$ from (18.7) yields:
\begin{align}
 L\underline{\tau}=\omega
\end{align}
where:
\begin{align}
 \omega=-[(1-u+t)\nu-1]\underline{L}\psi_{0}-[(1-u+t)\underline{\nu}+1]L\psi_{0}+(1-u+t)\rho'_{0}
\end{align}
From the conclusions (iii) of Theorem 17.1 we have, recalling:
\begin{align*}
 \nu=\frac{1}{2}(\textrm{tr}\chi+L\log\Omega)\\
\underline{\nu}+\alpha^{-1}\kappa\nu=\frac{1}{2}\alpha^{-1}\kappa L\log\Omega+\frac{1}{2}\underline{L}\log\Omega+\kappa\textrm{tr}\slashed{k}
\end{align*}
the bounds:
\begin{align}
 |(1-u+t)\nu-1|\leq C\delta_{0}(1+t)^{-1}[1+\log(1+t)]\\
|(1-u+t)\underline{\nu}+1|\leq C\delta_{0}[1+\log(1+t)]
\end{align}
     Let now $f$ be an arbitrary function defined on $W^{s}_{\epsilon_{0}}$. Let us denote by $\bar{f}(t,u)$ the mean value of $f$ on $S_{t,u}$ with respect to 
$d\mu_{\tilde{\slashed{g}}}$:
\begin{align}
 \bar{f}(t,u)=\frac{1}{\tilde{\textrm{Area}}(t,u)}\int_{S_{t,u}}f d\mu_{\tilde{\slashed{g}}},\quad \tilde{\textrm{Area}}(t,u)=\int_{S_{t,u}}d\mu_{\tilde{\slashed{g}}}
\end{align}
From the facts established in Chapter 5:
\begin{align}
 \frac{\partial}{\partial t}\int_{S_{t,u}}fd\mu_{\tilde{\slashed{g}}}=\int_{S_{t,u}}(Lf+2\nu f)d\mu_{\tilde{\slashed{g}}}\\
\frac{\partial}{\partial t}\tilde{\textrm{Area}}(t,u)=\int_{S_{t,u}}2\nu d\mu_{\tilde{\slashed{g}}}=2\bar{\nu}\tilde{\textrm{Area}}(t,u)
\end{align}
Hence by direct calculation:
\begin{align}
 \frac{\partial\bar{f}}{\partial t}(t,u)=\frac{1}{\tilde{\textrm{Area}}(t,u)}\int_{S_{t,u}}(Lf+2\nu f-2\bar{\nu}\bar{f})d\mu_{\tilde{\slashed{g}}}
\end{align}
Noting that for any pair of functions $f,g$ on $S_{t,u}$ we have:
\begin{align}
 \int_{S_{t,u}}(fg-\bar{f}\bar{g})d\mu_{\tilde{\slashed{g}}}=\int_{S_{t,u}}(f-\bar{f})(g-\bar{g})d\mu_{\tilde{\slashed{g}}}
\end{align}
we can write (18.18) in the form:
\begin{align}
 \frac{\partial \bar{f}}{\partial t}(t,u)=\overline{Lf}+2\overline{(\nu-\bar{\nu})(f-\bar{f})}
\end{align}
     We apply the formula (18.20) to the function $\underline{\tau}$ obtaining, in view of (18.11),
\begin{align}
 \frac{\partial\underline{\bar{\tau}}}{\partial t}(t,u)=\bar{\omega}+2\overline{(\nu-\bar{\nu})(\underline{\tau}-\bar{\underline{\tau}})}
\end{align}
Integrating this with respect to $t$ on $[0,s]$, we obtain:
\begin{align}
 \underline{\bar{\tau}}(s,u)=\underline{\bar{\tau}}(0,u)+\int_{0}^{s}[\bar{\omega}+2\overline{(\nu-\bar{\nu})(\underline{\tau}-\underline{\bar{\tau}})}](t,u)dt
\end{align}
Finally we replace $u$ by $u'\in[0,u]$, multiply by $\tilde{\textrm{Area}}(0,u')$, and integrate with respect to $u'$ on $[0,u]$ to obtain:
\begin{align}
 \int_{0}^{u}\underline{\bar{\tau}}(s,u')\tilde{\textrm{Area}}(0,u')du'\\\notag
=\int_{0}^{u}\int_{S_{0,u'}}\underline{\tau}d\mu_{\tilde{\slashed{g}}}du'+\int_{0}^{u}\tilde{\textrm{Area}}(0,u')
\{\int_{0}^{s}[\bar{\omega}+2\overline{(\nu-\bar{\nu})(\underline{\tau}-\underline{\bar{\tau}})}](t,u')dt\}du'
\end{align}
We shall estimate the second integral on the right in (18.23). Interchanging the order of integration this integral is:
\begin{align}
 I(s,u)=\int_{0}^{s}\{\int_{0}^{u}\tilde{\textrm{Area}}(0,u')[\bar{\omega}+2\overline{(\nu-\bar{\nu})(\underline{\tau}-\bar{\underline{\tau}})}](t,u')du'\}dt
\end{align}
Since $\tilde{\textrm{Area}}(0,u')\leq \tilde{\textrm{Area}}(0,0)=4\pi$, we have:
\begin{align}
 |I(s,u)|\leq C\int_{0}^{s}J(t,u)dt
\end{align}
where $J(t,u)$ is the integral:
\begin{align}
 J(t,u)=\int_{0}^{u}[|\bar{\omega}|+2\overline{|\nu-\bar{\nu}||\underline{\tau}-\bar{\underline{\tau}}|}](t,u')du'
\end{align}

    We consider first the contribution of the term $|\bar{\omega}|$, through (18.26), to the integral on the right in (18.25). Now $\omega$ is given by (18.12), and 
we consider first the third term on the right in (18.12). We have:
\begin{align}
 (1-u+t)\bar{\rho}^{\prime}_{0}(t,u)=\frac{(1-u+t)}{\tilde{\textrm{Area}}(t,u)}\int_{S_{t,u}}\rho'_{0}d\mu_{\tilde{\slashed{g}}}
\end{align}
Now for an arbitrary function $f$ on $S_{t,u}$, we have:
\begin{align}
 \tilde{\slashed{\Delta}}f=\frac{1}{\sqrt{\det\tilde{\slashed{g}}}}\frac{\partial}{\partial\vartheta^{A}}((\tilde{\slashed{g}}^{-1})^{AB}\sqrt{\det\tilde{\slashed{g}}}
\frac{\partial f}{\partial\vartheta^{A}})\\\notag
=\frac{1}{\Omega\sqrt{\det\slashed{g}}}\frac{\partial}{\partial \vartheta^{A}}((\slashed{g}^{-1})^{AB}\sqrt{\det\slashed{g}}\frac{\partial f}{\partial\vartheta^{A}})
=\Omega^{-1}\slashed{\Delta}f
\end{align}
Thus, integrating by parts on $S_{t,u}$ we obtain:
\begin{align}
 \int_{S_{t,u}}\mu\slashed{\Delta}\psi_{0}d\mu_{\tilde{\slashed{g}}}=\int_{S_{t,u}}\mu\Omega\tilde{\slashed{\Delta}}\psi_{0}d\mu_{\tilde{\slashed{g}}}\\\notag
=-\int_{S_{t,u}}\{\mu\tilde{\slashed{g}}^{-1}(\slashed{d}\Omega,\slashed{d}\psi_{0})+\Omega\tilde{\slashed{g}}^{-1}(\slashed{d}\mu,\slashed{d}\psi_{0})\}
d\mu_{\tilde{\slashed{g}}}\\\notag
=-\int_{S_{t,u}}\{\mu\slashed{g}^{-1}(\slashed{d}\log\Omega,\slashed{d}\psi_{0})+\slashed{g}^{-1}(\slashed{d}\mu,\slashed{d}\psi_{0})\}d\mu_{\tilde{\slashed{g}}}
\end{align}
Recalling that
\begin{align*}
 -(\slashed{d}\mu+2\zeta)\cdot\slashed{g}^{-1}=-(\eta+\zeta)\cdot\slashed{g}^{-1}
 =\Lambda
\end{align*}
we then conclude that:
\begin{align}
 \int_{S_{t,u}}\rho'_{0}d\mu_{\tilde{\slashed{g}}}
 =\int_{S_{t,u}}\Lambda\cdot\slashed{d}\psi_{0}d\mu_{\tilde{\slashed{g}}}
\end{align}
Now from conclusions (iii) of Theorem 17.1 we have:
\begin{align}
 |\Lambda|\leq C\delta_{0}(1+t)^{-1}[1+\log(1+t)]
\end{align}
Thus (18.30) implies, through (18.27),
\begin{align}
 (1-u+t)|\overline{\rho}^{\prime}_{0}(t,u)|\leq\frac{(1-u+t)}{(\textrm{Area}(t,u))^{1/2}}\sup_{S_{t,u}}|\Lambda|(\int_{S_{t,u}}|\slashed{d}\psi_{0}|^{2}
d\mu_{\slashed{g}})^{1/2}\\\notag
\leq C\delta_{0}(1+t)^{-1}[1+\log(1+t)](\int_{S_{t,u}}|\slashed{d}\psi_{0}|^{2}d\mu_{\tilde{\slashed{g}}})^{1/2}
\end{align}
This contribution to (18.25) is then bounded by
\begin{align}
 C\delta_{0}\sqrt{u}\int_{0}^{s}(1+t)^{-1}[1+\log(1+t)]\{\int_{\Sigma_{t}^{u}}|\slashed{d}\psi_{0}|^{2}\}^{1/2}dt\\\notag
\leq C\delta_{0}\sqrt{u}\{\int_{0}^{s}\int_{\Sigma_{t}^{u}}|\slashed{d}\psi_{0}|^{2}dt\}^{1/2}
\end{align}
The integral on the right in (18.33) is 
\begin{align}
 \int_{W^{s}_{u}}|\slashed{d}\psi_{0}|^{2}dtdu'd\mu_{\slashed{g}}
\end{align}
This integral has already been estimated in Chapter 7 by (7.210) for any variation. In particular for the first  order variation 
$\psi_{0}$
we have:
\begin{align}
 \int_{W^{s}_{u}}|\slashed{d}\psi_{0}|^{2}dtdu'd\mu_{\slashed{g}}\leq C\{\overline{K}[\psi_{0}](s,u)+\int_{0}^{u}\mathcal{F}^{s}_{0}[\psi_{0}](u')du'\}
\leq C\mathcal{E}^{u}_{0}[\psi_{0}](0)
\end{align}
The last step is by virtue of Theorem 5.1. Then (18.33) is bounded by:
\begin{align}
 C\delta_{0}\sqrt{u}\sqrt{\mathcal{E}^{u}_{0}[\psi_{0}](0)}
\end{align}
Next, we consider the first two terms on the right in (18.12), By (18.13) the contribution of the first term is bounded by:
\begin{align}
 \frac{C\delta_{0}(1+t)^{-1}[1+\log(1+t)]}{\tilde{\textrm{Area}}(t,u)}\int_{S_{t,u}}|\underline{L}\psi_{0}|d\mu_{\tilde{\slashed{g}}}\\\notag
\leq\frac{C\delta_{0}(1+t)^{-1}[1+\log(1+t)]}{(\tilde{\textrm{Area}}(t,u))^{1/2}}(\int_{S_{t,u}}|\underline{L}\psi_{0}|^{2}d\mu_{\tilde{\slashed{g}}})^{1/2}\\\notag
\leq C\delta_{0}(1+t)^{-2}[1+\log(1+t)](\int_{S_{t,u}}|\underline{L}\psi_{0}|^{2}d\mu_{\tilde{\slashed{g}}})^{1/2}
\end{align}
hence the corresponding contribution to $J(t,u)$ is bounded by:
\begin{align}
 C\delta_{0}\sqrt{u}(1+t)^{-2}[1+\log(1+t)]\sqrt{\mathcal{E}^{u}_{0}[\psi_{0}](t)}
\end{align}
and to the integral on the right in (18.25) by:
\begin{align}
 C\delta_{0}\sqrt{u}\{\int_{0}^{s}(1+t)^{-2}[1+\log(1+t)]dt\sqrt{\mathcal{E}^{u}_{0}[\psi_{0}](0)}
\leq C\delta_{0}\sqrt{u}\sqrt{\mathcal{E}^{u}_{0}[\psi_{0}](0)}
\end{align}
by Theorem 5.1. The factor $\sqrt{u}$ comes from Holder Inequality when we obtain (18.38). By (18.14) the contribution of the second term on the right in 
(18.12) to $\bar{\omega}$ 
is bounded by:
\begin{align}
 \frac{C\delta_{0}[1+\log(1+t)]}{\tilde{\textrm{Area}}(t,u)}\int_{S_{t,u}}|L\psi_{0}|d\mu_{\tilde{\slashed{g}}}\\\notag
\leq \frac{C\delta_{0}[1+\log(1+t)]}{(\tilde{\textrm{Area}}(t,u))^{1/2}}(\int_{S_{t,u}}|L\psi_{0}|^{2}d\mu_{\tilde{\slashed{g}}})^{1/2}\\\notag
\leq C\delta_{0}(1+t)^{-1}[1+\log(1+t)](\int_{S_{t,u}}|L\psi_{0}|^{2}d\mu_{\tilde{\slashed{g}}})^{1/2}
\end{align}
The contribution of this to the right of (18.25) is then bounded by:
\begin{align}
 C\delta_{0}\sqrt{u}\int_{0}^{s}(1+t)^{-1}[1+\log(1+t)]\{\int_{\Sigma_{t}^{u}}|L\psi_{0}|^{2}\}^{1/2}dt\\\notag
\leq C\delta_{0}\sqrt{u}\{\int_{0}^{s}\int_{\Sigma_{t}^{u}}|L\psi_{0}|^{2}dt\}^{1/2}\\\notag
\leq C\delta_{0}\sqrt{u}\{\int_{0}^{u}\mathcal{F}^{s}_{0}[\psi_{0}](u')du'\}^{1/2}\\\notag
\leq C\delta_{0}\sqrt{u}\sqrt{\mathcal{E}^{u}_{0}[\psi_{0}](0)}
\end{align}
by Theorem 5.1. Combining (18.36), (18.39), (18.41), we conclude that the contribution of $|\overline{\omega}|$ to the right of (18.25) is bounded by:
\begin{align}
 C\delta_{0}\sqrt{u}\sqrt{\mathcal{E}^{u}_{0}[\psi_{0}](0)}
\end{align}

     Finally, we consider the contribution of the term $2\overline{|\nu-\overline{\nu}||\underline{\tau}-\underline{\overline{\tau}}|}$ to the right of (18.25). We have:
\begin{align}
 \overline{|\nu-\overline{\nu}||\underline{\tau}-\underline{\overline{\tau}}|}=\frac{1}{\tilde{\textrm{Area}}(t,u)}\int_{S_{t,u}}
|\nu-\overline{\nu}||\underline{\tau}-\underline{\overline{\tau}}|d\mu_{\tilde{\slashed{g}}}
\end{align}
By (18.13):
\begin{align}
 |\nu-\bar{\nu}|\leq|\nu-\frac{1}{1-u+t}|+|\bar{\nu}-\frac{1}{1-u+t}|\\\notag
\leq C\delta_{0}(1+t)^{-2}[1+\log(1+t)]
\end{align}
Therefore
\begin{align}
 \overline{|\nu-\overline{\nu}||\underline{\tau}-\underline{\overline{\tau}}|}\leq\frac{C\delta_{0}(1+t)^{-2}[1+\log(1+t)]}{\tilde{\textrm{Area}}(t,u)}
\int_{S_{t,u}}|\underline{\tau}-\underline{\overline{\tau}}|d\mu_{\tilde{\slashed{g}}}\\\notag
\leq\frac{C\delta_{0}(1+t)^{-2}[1+\log(1+t)]}{(\tilde{\textrm{Area}}(t,u))^{1/2}}(\int_{S_{t,u}}|\underline{\tau}-\underline{\overline{\tau}}|^{2}
d\mu_{\tilde{\slashed{g}}})^{1/2}
\end{align}
We have:
\begin{align}
 \int_{S_{t,u}}|\underline{\tau}-\underline{\bar{\tau}}|^{2}d\mu_{\tilde{\slashed{g}}}=
\int_{S_{t,u}}\underline{\tau}^{2}d\mu_{\tilde{\slashed{g}}}-\int_{S_{t,u}}\underline{\bar{\tau}}^{2}d\mu_{\tilde{\slashed{g}}}\\\notag
\leq\int_{S_{t,u}}\underline{\tau}^{2}d\mu_{\tilde{\slashed{g}}}\leq C(1+t)^{2}\int_{S_{t,u}}|\underline{L}\psi_{0}|^{2}d\mu_{\slashed{g}}
+C\int_{S_{t,u}}|\psi_{0}|^{2}d\mu_{\slashed{g}}
\end{align}
by (18.9). The contribution of the first integral on the right in (18.46) to the right-hand side of (18.45) is bounded by:
\begin{align}
 C\delta_{0}(1+t)^{-2}[1+\log(1+t)](\int_{S_{t,u}}|\underline{L}\psi_{0}|^{2}d\mu_{\slashed{g}})^{1/2}
\end{align}
This coincides with (18.37), hence the corresponding contribution is bounded by (18.39). Finally, the contribution of the second term on the right in (18.46)
to (18.45) is bounded by:
\begin{align}
 C\delta_{0}(1+t)^{-3}[1+\log(1+t)](\int_{S_{t,u}}|\psi_{0}|^{2}d\mu_{\slashed{g}})^{1/2}\\\notag
\leq C\delta_{0}(1+t)^{-3}[1+\log(1+t)]\sqrt{u}\sqrt{\mathcal{E}^{u}_{0}[\psi_{0}](t)}
\end{align}
by Lemma 5.1, hence the corresponding contribution to (18.25) is also bounded by (18.39). Combining the above results we conclude that the contribution of $
2\overline{|\nu-\bar{\nu}||\underline{\tau}-\overline{\underline{\tau}}|}$ to the right of (18.25) is also bounded by (18.42).

    Then we have:
\begin{align}
 |I(s,u)|\leq C\delta_{0}\sqrt{u}\sqrt{\mathcal{E}^{u}_{0}[\psi_{0}](0)}
\end{align}
It follows from (18.23) that:
\begin{align}
 |\int_{0}^{u}\overline{\underline{\tau}}(s,u')\tilde{\textrm{Area}}(0,u')du'-\int_{\Sigma_{0}^{u}}\underline{\tau}d\mu_{\tilde{\slashed{g}}}du'|
\leq C\delta_{0}\sqrt{u}\sqrt{\mathcal{E}^{u}_{0}[\psi_{0}](0)}
\end{align}
From (18.6) and (18.9) we have:
\begin{align}
 P_{s}(u,\vartheta)=\underline{\tau}+\frac{1}{(1-u+s)}[u\underline{\tau}+(1+s)\psi_{0}]
\end{align}
Let us denote:
\begin{align}
 Q(u)=\int_{\Sigma_{0}^{u}}\underline{\tau}d\mu_{\tilde{\slashed{g}}}du'
\end{align}
Consider the case $\ell>0$. If the initial data satisfy:
\begin{align}
 Q(u)\leq -2C\delta_{0}\sqrt{u}\sqrt{\mathcal{E}^{u}_{0}[\psi_{0}](0)}<0
\end{align}
for some $u\in(0,\epsilon_{0}]$ with the constant $C$ as in (18.50), by (18.50) we have:
\begin{align}
 \int_{0}^{u}\overline{\underline{\tau}}(s,u')\tilde{\textrm{Area}}(0,u')du'\leq \frac{1}{2}Q(u)<0
\end{align}
It follows that:
\begin{align}
 \min_{u'\in[0,\epsilon_{0}]}\overline{\underline{\tau}}(s,u')\leq\frac{(1/2)Q(u)}{\int_{0}^{u}\tilde{\textrm{Area}}(0,u')du'}<0
\end{align}
hence also:
\begin{align}
 \min_{(u',\vartheta)\in[0,u]\times S^{2}}\underline{\tau}(s,u',\vartheta)\leq\frac{(1/2)Q(u)}{\int_{0}^{u}\tilde{\textrm{Area}}(0,u')du'}
\end{align}
It then follows from (18.51), in view of (iii) of Theorem 17.1, that:
\begin{align}
 \min_{(u',\vartheta)\in[0,u]\times S^{2}}P_{s}(u',\vartheta)\leq\frac{(1/4)Q(u)}{\int_{0}^{u}\tilde{\textrm{Area}}(0,u')du'}<0
\end{align}
for all $s$ large enough so that:
\begin{align}
 \frac{1}{2}+s\geq -C\delta_{0}\frac{\int_{0}^{u}\tilde{\textrm{Area}}(0,u')du'}{(1/8)Q(u)}
\end{align}
Thus by the discussion in the beginning of the chapter $t_{*\epsilon_{0}}$ is finite. In fact, (18.57) implies, in view of (18.5), that:
\begin{align}
 \min_{(u',\vartheta)\in[0,u]\times S^{2}}E_{s}(u',\vartheta)\leq\frac{(1/32)\ell Q(u)}{\int_{0}^{u}\tilde{\textrm{Area}}(0,u')du'}<0
\end{align}
for all $s$ large enough so that:
\begin{align}
 \frac{1+s}{1+\log(1+s)}\geq -C\delta_{0}\frac{\int_{0}^{u}\tilde{\textrm{Area}}(0,u')du'}{(1/32)\ell Q(u)}
\end{align}
where the constant $C$ is the one in (18.5). It follows from (18.3) in view of (18.2) that $s$ is bounded from above according to:
\begin{align}
 \log(1+s)\leq-\frac{\int_{0}^{u}\tilde{\textrm{Area}}(0,u')du'}{(1/32)\ell Q(u)}
\end{align}
Obviously, the upper bound in (18.61) exceeds the lower bounds in (18.58) and (18.60), it follows that $t_{*\epsilon_{0}}$ satisfies the upper bound:
\begin{align}
 \log(1+t_{*\epsilon_{0}})\leq -\frac{\int_{0}^{u}\tilde{\textrm{Area}}(0,u')du'}{(1/32)\ell Q(u)}
\end{align}
In the case $\ell<0$ a similar argument shows that if the initial data satisfy:
\begin{align}
 Q(u)\geq 2C\delta_{0}\sqrt{u}\sqrt{\mathcal{E}^{u}_{0}[\psi_{0}](0)}>0
\end{align}
for some $u\in(0,\epsilon_{0}]$ with the constant $C$ as in (18.50), then $t_{*\epsilon_{0}}$ is finite, in fact satisfies the upper bound:
\begin{align}
 \log(1+t_{*\epsilon_{0}})\leq -\frac{\int_{0}^{u}\tilde{\textrm{Area}}(0,u')du'}{(1/32)\ell Q(u)}
\end{align}
In fact, we have:
\begin{align*}
 \min_{(u',\vartheta)\in[0,u]\times S^{2}}E_{s}(u',\vartheta)\leq
\frac{1}{4}\ell\max_{(u',\vartheta)\in[0,u]\times S^{2}}P_{s}(u',\vartheta)+C\delta_{0}(1+s)^{-1}[1+\log(1+s)]
\end{align*}
Then we get a similar bound as (18.59).

$\textbf{Theorem 18.1 A}$ We obtain the first condition to guarantee the formation of shocks:\vspace{6mm}

The quantity $Q(u)$ satisfies
\begin{align*}
Q(u)\geq 2C\delta_{0}\sqrt{u}\sqrt{\mathcal{E}^{u}_{0}[\psi_{0}](0)}>0
\end{align*}                      
for $\ell<0$; and
\begin{align*}
 Q(u)\leq-2C\delta_{0}\sqrt{u}\sqrt{\mathcal{E}^{u}_{0}[\psi_{0}](0)}<0
\end{align*}
for $\ell>0$.\vspace{6mm}

We can modify this condition so that it is expressed in the Euclidean space.
 
Let us replace $d\mu_{\tilde{\slashed{g}}}$ by $d\mu_{\slashed{g}}$, then since the difference of $\Omega$ from unity, its value in the constant state, 
is bounded on $\Sigma_{0}^{u}$ by $C\delta_{0}$, the change in the integral will be bounded by $C\delta_{0}\sqrt{u}\sqrt{\mathcal{E}^{u}_{0}
[\psi_{0}](0)}$. Moreover we may replace $\underline{L}$ by the one corresponding to the constant state:
\begin{align}
 \mathring{\underline{L}}=\frac{\partial}{\partial x^{0}}-N^{i}\frac{\partial}{\partial x^{i}}
\end{align}
Then let us define:
\begin{align*}
 Q_{1}(u)=\int_{\Sigma_{0}^{u}}\underline{\tau}_{0}d^{3}x
\end{align*}
where
\begin{align*}
 \underline{\tau}_{0}=r\mathring{\underline{L}}\psi_{0}-\psi_{0}=
r(\partial_{0}\psi_{0}-N^{i}\partial_{i}\psi_{0})-\psi_{0}
\end{align*}
Since we are now on the initial hypersurface, the difference of 
\begin{align*}
 \int_{\Sigma_{0}^{u}}\underline{\tau}_{0}d\mu_{\slashed{g}}du'
\end{align*}
from $Q(u)$ will be also bounded by $C\delta_{0}\sqrt{u}\sqrt{\mathcal{E}^{u}_{0}[\psi_{0}](0)}$. 

$\textbf{Theorem 18.1 B}$ We obtain the second
condition to guarantee the formation of shocks:\vspace{6mm}

2) The quantity $Q_{1}(u)$ satisfies 
\begin{align*}
 Q_{1}(u)\geq 2C^{\prime}\delta_{0}\sqrt{u}\sqrt{\mathcal{E}^{u}_{0}[\psi_{0}](0)}>0
\end{align*}
for $\ell<0$; and 
\begin{align*}
 Q_{1}(u)\leq -2C^{\prime}\delta_{0}\sqrt{u}\sqrt{\mathcal{E}^{u}_{0}[\psi_{0}](0)}<0
\end{align*}
for $\ell>0$.\vspace{6mm}

The above two conditions depend only on the variation $\psi_{0}$, but the quantities $Q(u)$ and $Q_{1}(u)$ involve the first derivatives of $\psi_{0}$.
Now we shall obtain a third condition which depends on other variations $\psi_{i} : i=1,2,3$ but involves only the variations themselves (actually on the 
enthalpy $h$ and the normal component of the fluid speed $v_{N}$).
 
Consider the wave equation (1.23) on $\Sigma_{0}$, we can write it in the following form:
\begin{align*}
 \partial_{0}\psi_{0}-\partial_{i}\psi_{i}=2\psi_{i}\partial_{0}\psi_{i}-\psi_{i}\psi_{j}\partial_{i}\psi_{j}+(\eta^{2}-1)\Delta\phi
\end{align*}
where the repeated index means summing. By the $L^{\infty}$ bounds for the first order variations established in the previous chapter,
 we can bound the $L^{1}$ norm on $\Sigma_{0}^{u}$ of $r$ times the right-hand side by: $C\delta_{0}\sqrt{u}\sqrt{\mathcal{E}^{u}_{0,[1]}(0)}$.
Here the $\sqrt{u}$ in the difference comes from the using of Holder inequality.

Therefore defining
\begin{align}
 Q_{0}(u)=\int_{\Sigma_{0}^{u}}\{r(\partial_{i}\psi_{i}-N^{i}\partial_{i}\psi_{0})-\psi_{0}\}d\mu_{\slashed{g}}du'
\end{align}
we have:
\begin{align}
 |Q_{0}(u)-Q(u)|\leq C'\sqrt{\epsilon_{0}}\delta_{0}\sqrt{\mathcal{E}_{0,[1]}(0)}
\end{align}
Hence if in the case $\ell>0$ we have:
\begin{align}
 Q_{0}(u)\leq-2(C'+C)\sqrt{\epsilon_{0}}\delta_{0}\sqrt{\mathcal{E}^{u}_{0,[1]}(0)}<0\quad:\quad \textrm{for some}\quad u\in(0,\epsilon_{0}]
\end{align}
where $C$ is in (18.53), then by direct calculation, (18.53) holds and we have:
\begin{align}
 Q(u)\leq\frac{1}{2}Q_{0}(u)
\end{align}
hence, in place of (18.54):
\begin{align}
 \int_{0}^{u}\overline{\underline{\tau}}(s,u')\tilde{\textrm{Area}}(0,u')du'\leq \frac{1}{4}Q_{0}(u)<0
\end{align}
Similarly, if in the case $\ell>0$ we have:
\begin{align}
 Q_{0}(u)\geq 2(C'+C)\sqrt{\epsilon_{0}}\delta_{0}\sqrt{\mathcal{E}^{u}_{0,[1]}(0)}>0\quad:\quad \textrm{for some}\quad u\in(0,\epsilon_{0}]
\end{align}
where $C$ is in (18.63), then (18.63) holds and we have:
\begin{align}
 Q(u)\geq\frac{1}{2}Q_{0}(u)
\end{align}
hence:
\begin{align}
 \int_{0}^{u}\overline{\underline{\tau}}(s,u')\tilde{\textrm{Area}}(0,u')du'\geq \frac{1}{4}Q_{0}(u)>0
\end{align}
Integrating by parts and using divergence theorem,
\begin{align}
 \int_{\Sigma_{0}^{u}}r(\partial_{i}\psi_{i}-N^{i}\partial_{i}\psi_{0})d^{3}x=
\int_{S_{0,u}}r(\psi_{0}-\psi_{N})d\mu_{\slashed{g}}+\int_{\Sigma_{0}^{u}}(3\psi_{0}-\psi_{N})dx^{3}
\end{align}
To see this, we just write:
\begin{align*}
 \int_{\Sigma_{0}^{u}}r\partial_{i}\psi_{i}d^{3}x=\int_{\Sigma_{0}^{u}}\partial_{i}(r\psi_{i})d^{3}x-\int_{\Sigma_{0}^{u}}\psi_{N}d^{3}x
\end{align*}
and
\begin{align*}
 -\int_{\Sigma_{0}^{u}}rN^{i}\partial_{i}\psi_{0}d^{3}x=
\int_{\Sigma_{0}^{u}}r\partial_{i}N^{i}\psi_{0}d^{3}x+\int_{\Sigma_{0}^{u}}\sum_{i}(N^{i})^{2}\psi_{0}d^{3}x\\
-\int_{\Sigma_{0}^{u}}\partial_{i}(rN^{i}\psi_{0})d^{3}x
\end{align*}

$\textbf{Theorem 18.1 C}$ We obtain the following condition to guarantee the formation of shocks:\vspace{6mm}

3) The quantity 
\begin{align*}
 Q_{0}(u)=\int_{S_{0,u}}r(h+v_{N})+\int_{\Sigma_{0}^{u}}(2h+v_{N})
\end{align*}
satisfies 
\begin{align*}
 Q_{0}(u)\geq 2C^{\prime\prime}\delta_{0}\sqrt{u}\sqrt{\mathcal{E}^{u}_{0,[1]}(0)}
\end{align*}
for $\ell<0$; and
\begin{align*}
 Q_{0}(u)\leq -2C^{\prime\prime}\delta_{0}\sqrt{u}\sqrt{\mathcal{E}^{u}_{0,[1]}(0)}
\end{align*}
for $\ell>0$.\vspace{6mm}

Under one of these three conditions, the lifespan $t_{*\epsilon_{0}}$ is finite, and it satisfies one of the following inequalities:
\begin{align*}
\log(1+t_{*\epsilon_{0}})\leq \frac{Cu}{|\ell||Q(u)|}\\
\log(1+t_{*\epsilon_{0}})\leq \frac{C^{\prime}u}{|\ell||Q_{1}(u)|}\\
\log(1+t^{*\epsilon_{0}})\leq \frac{C^{\prime\prime}u}{|\ell||Q_{0}(u)|}
\end{align*}
respectively.
\chapter{The Structure of the Boundary of the Domain of the Maximal Solution}

\section{Nature of Singular Hypersurface in Acoustical Differential Structure}
\subsection{Preliminary}
In this chapter we shall take $\epsilon_{0}$ to be a variable with range $(0,1/2]$, in agreement with the setup of Theorem 17.1. To emphasize that
$\epsilon_{0}$ is now taken to be a variable, we denote it simply by $\epsilon$. Let us denote as in Chapter 2 by $t_{*}(u)$ the greatest lower bound 
of the parameter $t$ of the generators of $C_{u}$ in the domain of maximal solution, and by $t_{*\epsilon}$ the greatest lower bound of $t_{*}(u)$ for $u\in[0,\epsilon]$.
Obviously $t_{*\epsilon}$ is a non-increasing function of $\epsilon$. Then Theorem 17.1 applies for each $\epsilon\in(0,\epsilon_{0})$, where we now denote
by $\epsilon_{0}$ the maximal real number in the interval $(0,1/2]$ such that the smallness condition:
\begin{align*}
 C\sqrt{\epsilon_{0}}\sqrt{D_{[l+2]}(\epsilon_{0})}\leq\bar{\delta}_{0},\quad 
D_{[l+2]}(\epsilon_{0})=\sum_{\alpha,i}\|\partial_{i}\psi_{\alpha}\|^{2}_{H_{l+1}(\Sigma_{0}^{\epsilon_{0}})}
\end{align*}
holds. For, given any such $\epsilon$ we can find a $\delta_{0}\in(0,\bar{\delta}_{0}]$ such that the smallness condition of Theorem 17.1 holds with 
$\epsilon$ in the role of $\epsilon_{0}$, that is, we have:
\begin{align*}
 C\sqrt{\epsilon}\sqrt{D_{[l+2]}(\epsilon)}<\delta_{0}
\end{align*}
($\sqrt{\epsilon}\sqrt{D_{[l+2]}(\epsilon)}$ is an increasing function in $\epsilon$.)

Now according to Theorem 17.1, for each $\epsilon\in(0,\epsilon_{0})$, the solution, which exists in the classical sense in $W^{t_{*\epsilon}}_{\epsilon}\setminus
\Sigma^{\epsilon}_{t_{*\epsilon}}$ extends smoothly in acoustical coordinates to $\Sigma^{\epsilon}_{t_{*\epsilon}}$, but there is a non-empty set of points 
$K_{\epsilon}\subset\Sigma^{\epsilon}_{t_{*\epsilon}}$ where $\mu$ vanishes, being positive on the complement of $K_{\epsilon}$ in $\Sigma^{\epsilon}_{t_{*\epsilon}}$. Thus the set 
$K_{\epsilon}$ is also the set of minima of the function $\mu$ on $\Sigma_{t_{*\epsilon}}^{\epsilon}$. If $q$ is a point of $K_{\epsilon}$, then either $q$ is an 
interior minimum of $\mu$ on $\Sigma_{t_{*\epsilon}}^{\epsilon}$, in which case we have:
\begin{align}
 \mu(q)=0,\quad (\frac{\partial \mu}{\partial u})(q)=0,\quad (\frac{\partial\mu}{\partial\vartheta^{A}})(q)=0\quad:\quad A=1,2,\quad
(\frac{\partial^{2}\mu}{\partial u^{2}})(q)\geq 0 
\end{align}
at an interior zero;

and we shall assume that the non-degeneracy condition
\begin{align}
 (\frac{\partial^{2}\mu}{\partial u^{2}})(q)>0\quad:\textrm{at an interior zero}
\end{align}
holds, or $q$ is a boundary minimum, necessarily on $S_{t_{*\epsilon},\epsilon}$ (For on the other boundary of $\Sigma_{t_{*\epsilon}}^{\epsilon}$, namely on 
$S_{t_{*\epsilon},0}$, the constant state holds.), in which case we have:
\begin{align}
 \mu(q)=0,\quad (\frac{\partial\mu}{\partial u})(q)\leq 0,\quad (\frac{\partial\mu}{\partial\vartheta^{A}})(q)=0\quad:\quad A=1,2
\end{align}
and we shall assume the non-degeneracy condition
\begin{align}
 (\frac{\partial\mu}{\partial u})(q)<0\quad:\textrm{at a boundary zero}
\end{align}
holds.

      Let now $u_{m}(\epsilon)$ be the minimal value of $u\in(0,\epsilon)$ for which there is an interior zero of $\mu$ on $\Sigma_{t_{*\epsilon}}^{u}$.
Then if we replace $\epsilon$ by $\epsilon'\in(0,\epsilon]$, as we decrease $\epsilon'$, it follows directly from the definitions that $t_{*\epsilon'}$ stays
constant up to the point where $\epsilon'$ reaches the value $u_{m}(\epsilon)$, that is , we have:
\begin{align}
 t_{*\epsilon'}=t_{*\epsilon}\quad:\textrm{for all}\quad \epsilon'\in[u_{m}(\epsilon),\epsilon]
\end{align}
 Thereafter $t_{*\epsilon'}$ increases, $\mu$ being everywhere positive on $\Sigma_{t_{*u_{m}(\epsilon)}}^{\epsilon'}$, and for a while $K_{\epsilon'}$ consists only of 
boundary zeros. This situation either persists for all $\epsilon'\in(0,u_{m}(\epsilon))$ or, as we decrease $\epsilon'$ there is a first $\epsilon_{1}\in(0,u_{m}(\epsilon))$
at which one or more interior zeros of $\mu$ appear on $\Sigma^{\epsilon_{1}}_{t_{*\epsilon_{1}}}$, from which the preceding repeat with $\epsilon_{1}$ in the role of 
$\epsilon$.

      In the following, we shall use the conclusion in (iv) of Theorem 17.1, namely, for every $\epsilon\in(0,\epsilon_{0})$ the following upper bound for $L\mu$ holds in
the region $\mathcal{U}_{\epsilon}\subset W^{t_{*\epsilon}}_{\epsilon}$ where $\mu<1/4$:
\begin{align}
 L\mu\leq -C^{-1}(1+t)^{-1}[1+\log(1+t)]^{-1}\quad:\quad \textrm{in} \quad\mathcal{U}_{\epsilon}
\end{align}
What we shall use is just the fact that for finite $t$, $L\mu$ is bounded from above by a negative constant in the region where $\mu<1/4$.

      Let us denote:
\begin{align}
 V_{\epsilon_{0}}=\bigcup_{\epsilon\in(0,\epsilon_{0})}W^{t_{*\epsilon}}_{\epsilon}
\end{align}
This is the domain covered by Theorem 17.1, applied in the above manner (we are disregarding here the exterior domain where the constant state holds). This domain
is contained in the closure of $W_{\epsilon_{0}}$, the domain of maximal solution, and the set:
\begin{align}
 J_{\epsilon_{0}}=\bigcup_{\epsilon\in(0,\epsilon_{0})}K_{\epsilon}\subset V_{\epsilon_{0}}
\end{align}
of zeros of $\mu$ in $V_{\epsilon_{0}}$ is part of the singular boundary of $W_{\epsilon_{0}}$. 

      Our purpose in the following is to describe the acoustical spacetime structure in a neighborhood of a point of the singular boundary of $W_{\epsilon_{0}}$.
Consider now a value $\epsilon\in(0,\epsilon_{0})$ for which there is an interior zero of $\mu$ on $\Sigma_{t_{*\epsilon}}^{\epsilon}$, that is, there is a 
$u_{0}\in(0,\epsilon)$ and a zero of $\mu$ belonging to $S_{t_{*\epsilon},u_{0}}$. By (19.2) the minimum of $\frac{\partial^{2}\mu}{\partial u^{2}}$ over all zeros 
of $\mu$ on $S_{t_{*\epsilon},u_{0}}$ is positive. It follows that there are $u_{1}\in[0,u_{0})$, $u_{2}\in(u_{0},\epsilon]$ such that in the annular region:
\begin{align}
 \bigcup_{u\in(u_{1},u_{2})}S_{t_{*\epsilon},u}
\end{align}
 in $\Sigma_{t_{*\epsilon}}^{\epsilon}$ there are no zeros of $\mu$ except on $S_{t_{*\epsilon},u_{0}}$ itself. If $u_{0}=u_{m}(\epsilon)$, there are no zeros of $\mu$
on $S_{t_{*\epsilon},u}$ for any $u\in[0,u_{0})$ and we can take $u_{1}=0$.

     According to the discussion above, a boundary zero of $\mu$ on $\Sigma_{t*\epsilon}^{\epsilon}$ belongs to the surface $S_{t_{*\epsilon},\epsilon}$ and by (19.4)
there is a $u_{1}<\epsilon$ such that there are no zeros of $\mu$ on $S_{t_{*\epsilon},u}$ for any $u\in(u_{1},\epsilon)$. Moreover, there is a $u_{0}>\epsilon$ such
that for $\epsilon'\in[\epsilon,u_{0}]$ there is a corresponding boundary zero of $\mu$ on $\Sigma^{\epsilon'}_{t_{*\epsilon'}}$ belonging to the surface 
$S_{t_{*\epsilon'},\epsilon'}$ and no interior zeros, and there is a $u_{2}>u_{0}$ such that for $\epsilon'\in(u_{0},u_{2}]$ 
we have $t_{*\epsilon'}=t_{*u_{0}}$ and on $\Sigma_{t_{*\epsilon'}}$ we have an interior zero of $\mu$ on $S_{t_{*\epsilon'},u_{0}}$ and no zeros of $\mu$ on 
$S_{t_{*\epsilon'},u}$, for any $u\in[0,\epsilon'),u\slashed{=}u_{0}$ (otherwise, the solution can not be extended to $t_{*\epsilon}$, which is larger than $t_{*\epsilon
^{\prime}}$.) 
In this case we consider the annular region:
\begin{align}
 \bigcup_{u\in(u_{1},u_{2})}S_{t_{*u_{2}},u}
\end{align}
in $\Sigma_{t_{*u_{2}}}^{u_{2}}$.

\subsection{Intrinsic View Point}  
     Consider the manifold:
\begin{align}
 \mathcal{M}_{\epsilon_{0}}=[0,\infty)\times[0,\epsilon_{0})\times S^{2}
\end{align}
of all possible $(t,u,\vartheta)$ values with $u\in[0,\epsilon_{0})$. Then the domain $V_{\epsilon_{0}}$ in acoustical spacetime corresponds to the domain
\begin{align}
 \mathcal{V}_{\epsilon_{0}}=\{(t,u,\vartheta)\in\mathcal{M}_{\epsilon_{0}}\quad:\quad t\leq t_{*u}\}
\end{align}
in $\mathcal{M}_{\epsilon_{0}}$  and (19.9) corresponds to
\begin{align}
 \{t_{*\epsilon}\}\times(u_{1},u_{2})\times S^{2}
\end{align}
The mapping
\begin{align}
 (t,u,\vartheta)\in\mathcal{V}_{\epsilon_{0}}\mapsto x(t,u,\vartheta)\in V_{\epsilon_{0}},\quad x=(x^{\alpha}\quad:\alpha=0,1,2,3)
\end{align}
where $x^{\alpha}$ are rectangular coordinates in Galileo spacetime, is a homeomorphism of $\mathcal{V}_{\epsilon_{0}}$ onto $V_{\epsilon_{0}}$, which is 
locally also a diffeomorphism except at the points of the subset $\mathcal{J}_{\epsilon_{0}}$, the image of which is $J_{\epsilon_{0}}$, where $\mu$ vanishes.
On the domain (19.12) we have the acoustical metric $g$, given by (2.41):
\begin{align}
 g=-2\mu dtdu+\alpha^{-2}\mu^{2}du^{2}+\slashed{g}_{AB}(d\vartheta^{A}+\Xi^{A}du)(d\vartheta^{B}+\Xi^{B}du)
\end{align}
Recall that $\alpha(t,u,\vartheta)$ is a positive function near $1$ and that at each $(t,u)$, $\slashed{g}_{AB}(t,u,\vartheta)$ are the components of a positive
definite metric on $S^{2}$. The components of $g$ are smooth functions of the acoustical coordinates in the domain $\mathcal{V}_{\epsilon_{0}}$, including the 
boundary of this domain. Moreover, the properties of $\alpha$ and $\slashed{g}_{AB}$ just mentioned hold on the boundary of $\mathcal{V}_{\epsilon_{0}}$ as well.
However on the subset $\mathcal{J}_{\epsilon_{0}}$ of this boundary $\mu$ vanishes, hence, in view of the fact that
\begin{align}
 \sqrt{-\det g}=\mu\sqrt{\det\slashed{g}}
\end{align}
the metric (19.15) degenerates on $\mathcal{J}_{\epsilon_{0}}$.

     We now consider a smooth extension of the acoustical metric $g$ to $t>t_{*}(u)$. The extension is to satisfy the following two conditions:

   (i) The function $\alpha$ remains near $1$, and $\slashed{g}$ remains a positive definite metric on $S^{2}$.

   (ii) The function $\mu$ extends in such a way that $\frac{\partial\mu}{\partial t}$ is bounded from above by a negative constant where $\mu<1/4$.

Then for each $(u,\vartheta)\in(0,\epsilon_{0})\times S^{2}$ the following alternative holds: Either there is a first $t_{*}(u,\vartheta)\geq t_{*}(u)$ where 
$\mu$ vanishes, or $\mu$ is positive for all $t\geq t_{*}(u)$. Consider the subset $\tilde{\mathcal{D}}\subset(0,\epsilon_{0})\times S^{2}$ consisting of those points where
the first alternative holds. Then by Property (ii) of the extension $\tilde{\mathcal{D}}$ is an open set. In the case of an interior zero of $\mu$ on 
$\Sigma_{t_{*\epsilon}}^{\epsilon}$ belonging to $S_{t_{*\epsilon},u_{0}}$ we have $t_{*\epsilon}=t_{*}(u_{0})$, there is a $\vartheta_{0}\in S^{2}$ such that
$\mu(t_{*}(u_{0}),u_{0},\vartheta_{0})=0$. Therefore $(u_{0},\vartheta_{0})\in\tilde{\mathcal{D}}$. In the case of a boundary zero on $\Sigma^{\epsilon}_{t_{*\epsilon}}$,
for each $\epsilon'\in[\epsilon,u_{0}]$ there is a $\vartheta_{*}(\epsilon')\in S^{2}$ such that $\mu(t_{*}(\epsilon'),\epsilon',\vartheta_{*}(\epsilon'))=0$.
Therefore in this case $\tilde{\mathcal{D}}$ contains the curve
\begin{align}
 \{(\epsilon',\vartheta_{*}(\epsilon'))\quad:\quad \epsilon'\in[\epsilon,u_{0}]\}
\end{align}

     Consider:
\begin{align}
 \tilde{\mathcal{H}}=\{(t_{*}(u,\vartheta),u,\vartheta)\quad:\quad (u,\vartheta)\in\tilde{\mathcal{D}}\}
\end{align}
This is the zero level set of the function $\mu$ over the domain $\tilde{\mathcal{D}}\subset(0,\epsilon_{0})\times S^{2}$. Now $\mu$ is a smooth function which by (ii)
has no critical points on its zero level set. Therefore $\tilde{\mathcal{H}}$ is a smooth graph.

   Since $\tilde{\mathcal{H}}$ is a graph over $\tilde{\mathcal{D}}\subset(0,\epsilon_{0})\times S^{2}$, $(u,\vartheta)$ can be used as coordinates on $\tilde{\mathcal{H}}$,
and in these coordinates $g_{*}$, the metric induced on $\tilde{\mathcal{H}}$, is given by, from (19.15):
\begin{align}
 g_{*}=(\slashed{g}_{*})_{AB}(d\vartheta^{A}+\Xi^{A}_{*}du)(d\vartheta^{B}+\Xi^{B}_{*}du)
\end{align}
where:
\begin{align}
 (\slashed{g}_{*})_{AB}(u,\vartheta)=\slashed{g}_{AB}(t_{*}(u,\vartheta),u,\vartheta)
\end{align}
are the components of a positive definite metric on $S^{2}$, and:
\begin{align}
 \Xi^{A}_{*}(u,\vartheta)=\Xi^{A}(t_{*}(u,\vartheta),u,\vartheta)
\end{align}
Now the metric (19.19) is degenerate:
\begin{align}
 \det g_{*}=0
\end{align}

\subsection{Invariant Curves}
We see that although the hypersurface $\tilde{\mathcal{H}}$ is singular, being a hypersurface where the spacetime metric $g$ degenerates, from the point of view of its 
intrinsic geometry it is just like a regular null hypersurface in a regular spacetime. At each point $q\in\tilde{\mathcal{H}}$ there is a unique line $L_{q}\subset
T_{q}\tilde{\mathcal{H}}$ which we may consider to be the linear span of a non-zero null vector $V(q)$. Thus we have a null vectorfield $V$ on $\tilde{\mathcal{H}}$
and we can express it in terms of the coordinates $(u,\vartheta)$ by:
\begin{align}
 V=V^{u}\frac{\partial}{\partial u}+\slashed{V},\quad \slashed{V}=V^{A}\frac{\partial}{\partial\vartheta^{A}}
\end{align}
At each $q=(t_{*}(u,\vartheta),u,\vartheta)\in\tilde{\mathcal{H}}$, the vector $\slashed{V}(q)$ is tangent to $S_{*u}$ of constant $u$ through $q$, which projects
to the domain $\{u\}\times B_{u}$ where $B_{u}$ is the domain:
\begin{align*}
 B_{u}=\{\vartheta\in S^{2}\quad:\quad (u,\vartheta)\in\tilde{\mathcal{D}}\}\subset S^{2}
\end{align*}
Without loss of generality, we may set $V_{u}=1$. By direct calculation we have:
\begin{align}
 0=g(V,V)=g_{*}(V,V)=\slashed{g}_{*}(\slashed{V},\slashed{V})+2\slashed{g}_{*}(\Xi_{*},\slashed{V})+\slashed{g}_{*}(\Xi_{*},\Xi_{*})\\\notag
=\slashed{g}_{*}(\slashed{V}+\Xi_{*},\slashed{V}+\Xi_{*})
\end{align}
Since $\slashed{g}_{*}$ is positive definite, this holds if and only if:
\begin{align}
 \slashed{V}=-\Xi_{*}
\end{align}
Therefore $V$ is expressed in $(u,\vartheta)$ coordinates on $\tilde{\mathcal{H}}$ by:
\begin{align}
 V=\frac{\partial}{\partial u}-\Xi_{*}
\end{align}
Now, let
\begin{align}
 X(q)=X^{t}(q)\frac{\partial}{\partial t}+X^{u}(q)\frac{\partial}{\partial u}+X^{A}\frac{\partial}{\partial \vartheta^{A}}
\end{align}
be an arbitrary vector in $T_{q}\mathcal{M}_{\epsilon_{0}}, q\in\tilde{\mathcal{H}}$. Then $X(q)\in T_{q}\tilde{\mathcal{H}}$ if and only if:
\begin{align}
 X(q)(t-t_{*}(u,\vartheta))=0
\end{align}
Substituting (19.27) this is equivalent to:
\begin{align}
 X^{t}(q)=\frac{\partial t_{*}}{\partial u}X^{u}(q)+\frac{\partial t_{*}}{\partial \vartheta^{A}}X^{A}(q)
\end{align}
In particular the vectorfield $V$ as a vectorfield in spacetime along $\tilde{\mathcal{H}}$ is expressed by:
\begin{align}
 V=(\frac{\partial t_{*}}{\partial u}-\Xi^{A}_{*}\frac{\partial t_{*}}{\partial\vartheta^{A}})\frac{\partial}{\partial t}
+\frac{\partial}{\partial u}-\Xi^{A}_{*}\frac{\partial}{\partial\vartheta^{A}}
\end{align}
Here, $\frac{\partial}{\partial u}$ and $\frac{\partial}{\partial\vartheta^{A}}$ are the acoustical coordinate vectorfields in spacetime, not in $\tilde{\mathcal{H}}$.
In fact, we have:
\begin{align*}
 \frac{\partial}{\partial u}=\frac{\partial t_{*}}{\partial u}\frac{\partial}{\partial t}+\frac{\partial}{\partial u}\\
\frac{\partial}{\partial\vartheta^{A}}=\frac{\partial t_{*}}{\partial\vartheta^{A}}\frac{\partial}{\partial t}+\frac{\partial}{\partial\vartheta^{A}}
\end{align*}
The vectorfields on the left-hand side are the coordinate vectorfields in $\tilde{\mathcal{H}}$.

Given any vector $X(q)\in T_{q}\mathcal{M}_{\epsilon_{0}}, q\in\tilde{\mathcal{H}}$, there is a unique real number $\lambda$ such that
$X(q)-\lambda L(q)\in T_{q}\tilde{\mathcal{H}}$. In fact it is readily seen from (19.29), recalling that $L=\frac{\partial}{\partial t}$, that:
\begin{align*}
 \lambda=X^{t}(q)-\frac{\partial t_{*}}{\partial u}X^{u}(q)-\frac{\partial t_{*}}{\partial \vartheta^{A}}X^{A}(q)
\end{align*}
This defines a projection operator $\Pi_{*}$ from $T_{q}\mathcal{M}_{\epsilon_{0}}$ to $T_{q}\tilde{\mathcal{H}}$ by: $\Pi_{*}X(q)=X(q)-\lambda L(q)$. 
We call $\Pi_{*}$ $the$ $L$ $projection$ $to$ $\tilde{\mathcal{H}}$. Obviously,
\begin{align}
 V=\Pi_{*} T
\end{align}

    We call the integral curves of $V$ the $invariant$ $curves$. The singular surface $\tilde{\mathcal{H}}$ is ruled by these curves. The invariant curves as 1-dimensional
submanifolds of $\tilde{\mathcal{H}}$ are independent of the choice of acoustical function $u$, being the integral manifolds of the 1-dimensional  distribution $\{L_{q}
\quad:\quad q\in\tilde{\mathcal{H}}\}$ on $\tilde{\mathcal{H}}$. The invariant curves have zero arc length.

    Now we may adapt the coordinates $(u,\vartheta)$ on $\tilde{\mathcal{H}}$ so that the coordinate lines $\vartheta=const$ coincide with the invariant curves. Obviously,
this choice is equivalent to the condition:
\begin{align}
 \Xi_{*}=0
\end{align}
We then call the corresponding acoustical coordinates $canonical$. In these coordinates (19.19) takes the form:
\begin{align}
 g_{*}=(\slashed{g}_{*})_{AB}d\vartheta^{A}d\vartheta^{B}
\end{align}
By (19.32), we have, for $t\leq t_{*}(u,\vartheta)$:
\begin{align*}
 \Xi^{A}(t,u,\vartheta)=-\int_{t}^{t_{*}(u,\vartheta)}(\frac{\partial \Xi^{A}}{\partial t})(t',u,\vartheta)dt'
\end{align*}
Since also:
\begin{align*}
 \mu(t,u,\vartheta)=-\int_{t}^{t_{*}(u,\vartheta)}(\frac{\partial\mu}{\partial t})(t',u,\vartheta)dt'
\end{align*}
it follows that the functions:
\begin{align}
 \hat{\Xi}^{A}=\mu^{-1}\Xi^{A}
\end{align}
which a priori are defined only for $t<t_{*}(u,\vartheta)$, are actually given by:
\begin{align}
 \hat{\Xi}^{A}(t,u,\vartheta)=
\frac{\textrm{mean value on} [t,t_{*}(u,\vartheta)]\textrm{of}\{(\partial\Xi^{A}/\partial t)(,u,\vartheta)\}}
{\textrm{mean value on} [t,t_{*}(u,\vartheta)]\textrm{of}\{(\partial\mu/\partial t)(,u,\vartheta)\}}
\end{align}
hence extend smoothly to $t=t_{*}(u,\vartheta)$, that is, to $\tilde{\mathcal{H}}$.

\subsection{Extrinsic View Point}
     We now consider the character of the singular hypersurface $\tilde{\mathcal{H}}$ from the extrinsic point of view. To do this, we consider the reciprocal 
acoustical metric, a quadratic form in the cotangent space to the spacetime manifold at each regular point, given by:
\begin{align}
 g^{-1}=-(1/2\mu)(L\otimes\underline{L}+\underline{L}\otimes L)+(\slashed{g}^{-1})^{AB}X_{A}\otimes X_{B}
\end{align}
Since
\begin{align*}
 L=\frac{\partial}{\partial t},\quad \underline{L}=\alpha^{-1}\kappa L+2T=
\alpha^{-2}\mu\frac{\partial}{\partial t}+2(\frac{\partial}{\partial u}-\mu\hat{\Xi}^{A}\frac{\partial}{\partial\vartheta^{A}})
\end{align*}
(19.36) takes the form in canonical acoustical coordinates:
\begin{align}
 \mu g^{-1}=-(\frac{\partial}{\partial t}\otimes\frac{\partial}{\partial u}+\frac{\partial}{\partial u}\otimes\frac{\partial}{\partial t})
-\mu\alpha^{-2}\frac{\partial}{\partial t}\otimes\frac{\partial}{\partial t}\\\notag
+\mu(\frac{\partial}{\partial t}\otimes\hat{\Xi}^{A}\frac{\partial}{\partial\vartheta^{A}}+\hat{\Xi}^{A}\frac{\partial}{\partial\vartheta^{A}}\otimes
\frac{\partial}{\partial t})+\mu(\slashed{g}^{-1})^{AB}\frac{\partial}{\partial\vartheta^{A}}\otimes\frac{\partial}{\partial\vartheta^{B}}
\end{align}
We see that although, due to the degeneracy of $g$ on $\tilde{\mathcal{H}}$, $g^{-1}$ blows up at $\tilde{\mathcal{H}}$, $\mu g^{-1}$ in fact extends smoothly to 
$\tilde{\mathcal{H}}$. The character of $\tilde{\mathcal{H}}$ from the extrinsic point of view is then determined by the sign of the invariant $\mu(g^{-1})
^{\alpha\beta}\partial_{\alpha}\mu\partial_{\beta}\mu$ on $\tilde{\mathcal{H}}$. We have:
\begin{align}
 \mu(g^{-1})^{\alpha\beta}\partial_{\alpha}\mu\partial_{\beta}\mu=-2\frac{\partial\mu}{\partial t}\frac{\partial\mu}{\partial u}\\\notag
+\mu\{-\alpha^{-2}(\frac{\partial\mu}{\partial t})^{2}+2\frac{\partial\mu}{\partial t}\frac{\partial\mu}{\partial\vartheta^{A}}\hat{\Xi}^{A}
+(\slashed{g}^{-1})^{AB}\frac{\partial\mu}{\partial\vartheta^{A}}\frac{\partial\mu}{\partial\vartheta^{B}}\}
\end{align}
On $\tilde{\mathcal{H}}$, $\mu$ vanishes and this reduces to, simply:
\begin{align}
 \mu(g^{-1})^{\alpha\beta}\partial_{\alpha}\mu\partial_{\beta}\mu=-2\frac{\partial\mu}{\partial t}\frac{\partial\mu}{\partial u}
\end{align}
In view of condition (ii) $\mu(g^{-1})^{\alpha\beta}\partial_{\alpha}\mu\partial_{\beta}\mu$ is $>0$, $=0$, $<0$, at a point $q\in\tilde{\mathcal{H}}$ according as to
whether $(\partial\mu/\partial u)(q)$ is $>0$, $=0$, $<0$. We conclude that $\tilde{\mathcal{H}}$ is spacelike, null, or timelike, at $q$ according to as whether
$(\partial\mu/\partial u)(q)$ is $<0$, $=0$, $>0$. Now the boundary of the domain of the maximal solution cannot be timelike at any point. Therefore the timelike part 
of $\tilde{\mathcal{H}}$ cannot be part of the boundary of the domain of maximal solution. 
So we shall focus on the spacelike part of $\tilde{\mathcal{H}}$, which we denote by $\mathcal{H}$ and its boundary $\partial\mathcal{H}$, 
where $\tilde{\mathcal{H}}$ is null.

     Consider then the open subset $\mathcal{D}\subset\tilde{\mathcal{D}}$ defined by:
\begin{align}
 \mathcal{D}=\{(u,\vartheta)\in\tilde{\mathcal{D}}\quad:\quad (\partial\mu/\partial u)(t_{*}(u,\vartheta),u,\vartheta)<0\}
\end{align}
 (canonical acoustical coordinates). The boundary of $\mathcal{D}$ in $\tilde{\mathcal{D}}$ is given by:
\begin{align}
 \partial\mathcal{D}=\{(u,\vartheta)\in\tilde{\mathcal{D}}\quad:\quad (\partial\mu/\partial u)(t_{*}(u,\vartheta),u,\vartheta)=0\}
\end{align}
Then $\mathcal{H}$ is the graph (19.18) over $\mathcal{D}$:
\begin{align}
 \mathcal{H}=\{(t_{*}(u,\vartheta),u,\vartheta)\quad:\quad (u,\vartheta)\in\mathcal{D}\}
\end{align}
and $\partial\mathcal{H}$ is the graph (19.18) over $\partial\mathcal{D}$:
\begin{align}
 \partial\mathcal{H}=\{(t_{*}(u,\vartheta),u,\vartheta)\quad:\quad(u,\vartheta)\in\partial\mathcal{D}\}
\end{align}
Assuming that the non-degeneracy condition:
\begin{align}
 (u,\vartheta)\in\partial\mathcal{D}\quad\textrm{implies}\quad (\partial^{2}\mu/\partial u^{2})(t_{*}(u,\vartheta),u,\vartheta)\slashed{=}0
\end{align}
holds, $\partial\mathcal{D}$ splits into the disjoint union $\partial\mathcal{D}=\partial_{-}\mathcal{D}\bigcup\partial_{+}\mathcal{D}$ where:
\begin{align}
 \partial_{-}\mathcal{D}=\{(u,\vartheta)\in\partial\mathcal{D}\quad:\quad(\partial^{2}\mu/\partial u^{2})(t_{*}(u,\vartheta),u,\vartheta)>0\}\\\notag
\partial_{+}\mathcal{D}=\{(u,\vartheta)\in\partial\mathcal{D}\quad:\quad(\partial^{2}\mu/\partial u^{2})(t_{*}(u,\vartheta),u,\vartheta)<0\}
\end{align}
The boundary $\partial\mathcal{H}$ of $\mathcal{H}$ in $\tilde{\mathcal{H}}$ similarly splits into the disjoint union $\partial\mathcal{H}=
\partial_{-}\mathcal{H}\bigcup\partial_{+}\mathcal{H}$  where $\partial_{-}\mathcal{H}$ and $\partial_{+}\mathcal{H}$ are the graphs (19.18) over
$\partial_{-}\mathcal{D}$ and $\partial_{+}\mathcal{D}$ respectively:
\begin{align}
 \partial_{-}\mathcal{H}=\{(t_{*}(u,\vartheta),u,\vartheta)\quad:\quad (u,\vartheta)\in\partial_{-}\mathcal{D}\}\\\notag
\partial_{+}\mathcal{H}=\{(t_{*}(u,\vartheta),u,\vartheta)\quad:\quad (u,\vartheta)\in\partial_{+}\mathcal{D}\}
\end{align}

     If there is an interior zero of $\mu$ on $\Sigma_{t_{*\epsilon}}^{\epsilon}$ corresponding to $(u_{0},\vartheta_{0})$, then $(u_{0}\vartheta_{0})$ 
belongs to $\partial_{-}\mathcal{D}$. If there is a boundary zero of $\mu$ on $\Sigma_{t_{*\epsilon}}^{\epsilon}$, then by (19.4), $\mathcal{D}$ contains 
the curve (19.17) except for the point $(u_{0},\vartheta_{*}(u_{0}))$ which belongs to $\partial_{-}\mathcal{D}$.

     Now on $\tilde{\mathcal{H}}$ we have $\mu=0$, that is, we have:
\begin{align}
 \mu(t_{*}(u,\vartheta),u,\vartheta)=0\quad:\quad\textrm{for all}\quad (u,\vartheta)\in\tilde{\mathcal{D}}
\end{align}
 Differentiating this equation implicitly with respect to $u$ we obtain:
\begin{align}
 (\frac{\partial t_{*}}{\partial u})(u,\vartheta)=-(\frac{\partial\mu/\partial u}{\partial\mu/\partial t})(t_{*}(u,\vartheta),u,\vartheta)
\end{align}
By virtue of property (ii) we have:
\begin{align}
 \mathcal{D}=\{(u,\vartheta)\in\tilde{\mathcal{D}}\quad:\quad (\partial t_{*}/\partial u)(u,\vartheta)<0\}
\end{align}
and:
\begin{align}
 \partial\mathcal{D}=\{(u,\vartheta)\in\tilde{\mathcal{D}}\quad:\quad (\partial t_{*}/\partial u)(u,\vartheta)=0\}
\end{align}
(canonical acoustical coordinates).

     Moreover, differentiating again (19.48) with respect to $u$ and evaluating the result on $\partial\mathcal{D}$, we obtain, in view of (19.41) and
(19.50):
\begin{align}
 (\frac{\partial^{2}t_{*}}{\partial u^{2}})(u,\vartheta)=-(\frac{\partial^{2}\mu/\partial u^{2}}{\partial\mu/\partial t})(t_{*}(u,\vartheta),u,\vartheta)
\quad:\quad\textrm{on}\quad \partial\mathcal{D}
\end{align}
Comparing with (19.45) we conclude that:
\begin{align}
 \partial_{-}\mathcal{D}=\{(u,\vartheta)\in\partial\mathcal{D}\quad:\quad (\partial^{2}t_{*}/\partial u^{2})(u,\vartheta)>0\}\\\notag
\partial_{+}\mathcal{D}=\{(u,\vartheta)\in\partial\mathcal{D}\quad:\quad(\partial^{2}t_{*}/\partial u^{2})(u,\vartheta)<0\}
\end{align}
Thus if we consider a connected component of $\mathcal{H}$ and the corresponding components of $\partial_{-}\mathcal{H}$, $\partial_{+}\mathcal{H}$, then 
the component of $\partial_{-}\mathcal{H}$, which is not empty (In fact, the point on $\partial_{-}\mathcal{H}$ corresponds the interior zero. If there is
no interior zero, this point corresponds to $u=\epsilon_{0}$), is the past boundary of $\mathcal{H}$, the component of $\partial_{+}\mathcal{H}$, which may be empty 
($\mathcal{H}$ could be the asymptote of $C_{0}$), its future boundary, the function $t_{*}$ reaching a minimum, along each invariant curve, at $\partial_{-}\mathcal{H}$, 
a maximum at $\partial_{+}\mathcal{H}$.

     Consider now the function $f=\partial\mu/\partial u$ on $\tilde{\mathcal{H}}$. In canonical acoustical coordinates on $\tilde{\mathcal{H}}$ we have:
\begin{align}
 f(u,\vartheta)=(\frac{\partial\mu}{\partial u})(t_{*}(u,\vartheta),u,\vartheta)
\end{align}
Differentiating with respect to $u$, and evaluating the result on $\partial\mathcal{H}$, we obtain, in view of (19.50):
\begin{align}
 (\frac{\partial f}{\partial u})(u,\vartheta)=(\frac{\partial^{2}\mu}{\partial u^{2}})(t_{*}(u,\vartheta),u,\vartheta)\quad:\quad\textrm{on}\quad \partial\mathcal{H}
\end{align}
Thus $\partial\mathcal{H}$ is the zero level set of a smooth function, namely $f$, on the smooth manifold $\tilde{\mathcal{H}}$, and by (19.54) and the non-degeneracy
condition (19.44) this is a non-critical level set. It follows that $\partial\mathcal{H}$ and its two components $\partial_{-}\mathcal{H}$ and $\partial_{+}\mathcal{H}$
are smooth. Moreover, since $\partial_{-}\mathcal{H}$ and $\partial_{+}\mathcal{H}$ intersect the invariant curves, generators of $\tilde{\mathcal{H}}$ transversally 
they are both spacelike surfaces in acoustical spacetime.

      Finally, we note that $\partial_{+}\mathcal{H}$ cannot be part of the boundary of the domain of the maximal solution. This is seen as follows. For any given 
component of $\partial_{+}\mathcal{H}$, there is a component of $\tilde{\mathcal{H}}\setminus\bar{\mathcal{H}}$, the timelike part of $\tilde{\mathcal{H}}$, whose 
future boundary is the given component of $\partial_{+}{\mathcal{H}}$ (
At $\partial_{+}\mathcal{D}$, $\partial^{2}t_{*}/\partial u^{2}<0$, so as $u$ decreases, $\partial t_{*}/\partial u$ becomes positive beyond $\partial_{+}
\mathcal{D}$, thus $t_{*}$ decreases as $u$ decreases). This, or in fact any, component of $\tilde{\mathcal{H}}\setminus\bar{\mathcal{H}}$ has a past boundary, for, 
$\lim_{u\rightarrow0}t_{*}(u)=\infty$, therefore $t_{*}(u,\vartheta)$ must have intervals of increase along each invariant line as we approach $u=0$. The past boundary 
of the component in question is then also the past boundary of the next outward (that is, in the direction of decreasing $u$ along the invariant curves) component of 
$\mathcal{H}$, a component of $\partial_{-}\mathcal{H}$. But the domain of the maximal solution terminates at the incoming characteristic hypersurface $\underline{C}$
generated by the incoming null normals to this component of $\partial_{-}\mathcal{H}$. So the maximal solution cannot reach $\partial_{+}\mathcal{H}$.

     Summarizing, we get:

$\textbf{Proposition 19.1}$ Consider a smooth extension of the acoustical metric $g$ to $t>t_{*}(u)$, satisfying the following two conditions:

  (i) The function $\alpha$ remains near 1, and the metric $\slashed{g}$ remains a positive definite metric on $S^{2}$.

  (ii) The function $\mu$ is extended in such a way that $\partial\mu/\partial t$ is bounded from above by a negative constant where $\mu<1/4$.

Then there is an open subset $\tilde{\mathcal{D}}\subset(0,\epsilon_{0})\times S^{2}$ and a smooth graph
\begin{align*}
 \tilde{\mathcal{H}}=\{(t_{*}(u,\vartheta),u,\vartheta)\quad:\quad (u,\vartheta)\in\tilde{\mathcal{D}}\}
\end{align*}
where $\mu$ vanishes, being positive below this graph. The singular hypersurface $\tilde{\mathcal{H}}$ with its induced metric $g_{*}$ is smooth and has the intrinsic
geometry of a regular null hypersurface in a regular spacetime. It is ruled by invariant curves of vanishing arc length. The tangent line to an invariant 
curve at a point $q$ is the linear span of the vector $\Pi_{*}T(q)$, the $L$ projection of $T(q)$ to $\tilde{\mathcal{H}}$. Canonical acoustical coordinates are defined 
by taking the $\vartheta=const$, coordinate lines on $\tilde{\mathcal{H}}$ to be the invariant curves. On the other hand, from the point of view of how $\tilde{\mathcal{H}}$
is embedded in the acoustical spacetime, the extrinsic point of view, $\tilde{\mathcal{H}}$ is at a point $q$ spacelike, null, or timelike, according as to whether, in 
canonical acoustical coordinates, $\partial\mu/\partial u$ is $<0$, $=0$, $>0$, at $q$, or equivalently, as to whether $\partial t_{*}/\partial u$ is $<0$, $=0$, $>0$,
at $q$. The timelike part of $\tilde{\mathcal{H}}$ cannot be part of the boundary of the domain of the maximal solution. Moreover, denoting by $\mathcal{H}$ the spacelike 
part of $\tilde{\mathcal{H}}$ and by $\partial\mathcal{H}$ its boundary, under the non-degeneracy condition:
\begin{align*}
 \partial^{2}\mu/\partial u^{2}\slashed{=}0\quad:\quad\textrm{on}\quad \partial\mathcal{H}
\end{align*}
the boundary splits into the disjoint union of $\partial_{-}\mathcal{H}$ and $\partial_{+}\mathcal{H}$, where $\partial_{-}\mathcal{H}$ and $\partial_{+}\mathcal{H}$
correspond to $\partial^{2}\mu/\partial u^{2}>0$ and $<0$ respectively, or equivalently $\partial^{2}t_{*}/\partial u^{2}>0$ and $<0$ respectively. Each of $\partial_{-}
\mathcal{H}$, $\partial_{+}\mathcal{H}$, is a smooth spacelike surface in the acoustical spacetime. For each connected component of $\mathcal{H}$ the corresponding components 
of $\partial_{-}\mathcal{H}$ and $\partial_{+}\mathcal{H}$ are respectively its past and future boundaries, the sets of past and future end points of its invariant curves.
Finally, $\partial_{+}\mathcal{H}$, which, in contrast to $\partial_{-}\mathcal{H}$, may be empty, cannot be part of the boundary of the domain of the maximal solution.

\section{The Trichotomy Theorem for Past Null Geodesics Ending at Singular Boundary}
     We now turn to the investigation of the past null geodesic cone of an arbitrary point $q$ on the singular boundary of the domain of the maximal solution.
According to the above, such a point belongs either to $\mathcal{H}$, or to $\partial_{-}\mathcal{H}$. The domain of the maximal solution or maximal development of the 
initial data on $\Sigma_{0}^{\epsilon_{0}}\bigcup\Sigma_{0}^{E}$, where $\Sigma_{0}^{E}$ denotes the exterior of the unit sphere in $\Sigma_{0}$, where the constant state 
holds, has, like any development of the initial data in question, the property that for each point $q$ in this domain each past directed curve issuing at $q$ which 
is causal with respect to the acoustical metric terminates in the past at a point of $\Sigma^{\epsilon_{0}}_{0}\bigcup\Sigma_{0}^{E}$. In particular, each past directed 
null geodesic issuing  at $q$ terminates in the past at a point of $\Sigma_{0}^{\epsilon_{0}}\bigcup\Sigma_{0}^{E}$. Suppose now that $q$ is a zero of $\mu$ on 
$\Sigma_{t_{*\epsilon}}^{\epsilon}$ for some $\epsilon\in(0,\epsilon_{0})$. Then $q\in W^{t_{*\epsilon}}_{\epsilon}$ and the union of $W^{t_{*\epsilon}}_{\epsilon}$ 
with the domain $E^{t_{*\epsilon}}$  in Galileo spacetime exterior to the cone $C_{0}$ and bounded by the hyperplanes $\Sigma_{0}$ and $\Sigma_{t_{*\epsilon}}$, where the constant
state holds, is a development of the restriction of the initial data to $\Sigma_{0}^{\epsilon}\bigcup\Sigma_{0}^{E}$, therefore each past directed null geodesic issuing at 
$q$ remains in $W^{t_{*\epsilon}}_{\epsilon}\bigcup E^{t_{*\epsilon}}$ and terminates in the past at a point of $\Sigma_{0}^{\epsilon}\bigcup\Sigma_{0}^{E}$.

\subsection{Hamiltonian Flow}
     Now, the null geodesic flow of a Lorentzian manifold $(\mathcal{M},g)$ is from the point view of cotangent bundle the Hamiltonian flow on $T^{*}\mathcal{M}$ generated by the Hamiltonain:
\begin{align}
 H=\frac{1}{2}(g^{-1})^{\mu\nu}(q)p_{\mu}p_{\nu}
\end{align}
on the surface $H=0$. 
Here $(q^{\mu}\quad:\quad \mu=1,...,n)$, $n=\textrm{dim}\mathcal{M}$, are local coordinates on 
$\mathcal{M}$ (called in the Hamiltonian context $canonical$ $coordinates$), and expanding $p\in T^{*}_{q}\mathcal{M}$ in the basis $(dq^{1}(q),...,dq^{n}(q))$:
\begin{align*}
 p=p_{\mu}dq^{\mu}(q)
\end{align*}
 the coefficients $(p_{\mu}\quad:\quad \mu=1,...,n)$ of the expansion are linear coordinates on $T^{*}_{q}\mathcal{M}$ (called in the Hamiltonian context $canonical$ 
$momenta$). We thus have local coordinates \\$(q^{\mu}\quad:\quad \mu=1,...,n;\quad p_{\mu}\quad:\quad \mu=1,...,n)$ on $T^{*}\mathcal{M}$ and the equations defining a 
Hamiltonian flow are the canonical equations:
\begin{align}
 \frac{dq^{\mu}}{d\tau}=\frac{\partial H}{\partial p_{\mu}},\quad \frac{dp_{\mu}}{d\tau}=-\frac{\partial H}{\partial q_{\mu}}
\end{align}
Each level set of the Hamiltonian $H$ is invariant by the Hamiltonian flow. We may thus consider the restriction of the flow to any given level set of $H$, in particular
to the surface $H=0$. Set now:
\begin{align}
 H(q,p)=\Omega(q)\tilde{H}(q,p)
\end{align}
where $\Omega$ is a positive function on $\mathcal{M}$. Then the Hamiltonian flow of $\tilde{H}$ on its zero level set is equivalent, up to reparametrization, to the 
Hamiltonian flow of $H$ on its zero level set. More precisely, let $\tau\mapsto(q(\tau),p(\tau))$ be a solution of the canonical equations (19.56) on the surface $H=0$.
Then defining a new parameter $\tilde{\tau}$ by:
\begin{align}
 \frac{d\tilde{\tau}}{d\tau}=\Omega(q(\tau))
\end{align}
$\tilde{\tau}\mapsto(\tilde{q}(\tilde{\tau}),\tilde{p}(\tilde{\tau}))=(q(\tau),p(\tau))$ is a solution of (19.56) with $H$ replaced by $\tilde{H}$ and $\tau$ by 
$\tilde{\tau}$ on the surface $\tilde{H}=0$. This is so because
\begin{align*}
 \frac{d\tilde{q}^{\mu}}{d\tilde{\tau}}=\frac{dq^{\mu}}{d\tau}\frac{d\tau}{d\tilde{\tau}}=\frac{\partial H}{\partial p_{\mu}}\Omega^{-1}=
\frac{\partial\tilde{H}}{\partial p_{\mu}}\\
\frac{d\tilde{p}_{\mu}}{d\tilde{\tau}}=\frac{dp_{\mu}}{d\tau}\frac{d\tau}{d\tilde{\tau}}=-\frac{\partial H}{\partial q^{\mu}}\Omega^{-1}=
-\frac{\partial(\tilde{H}\Omega)}{\partial q^{\mu}}\Omega^{-1}\\
=-\frac{\partial\tilde{H}}{\partial q^{\mu}}-\tilde{H}\Omega^{-1}\frac{\partial\Omega}{\partial q^{\mu}}=-\frac{\partial\tilde{H}}{\partial q^{\mu}}\quad:
\textrm{on}\quad \tilde{H}=0
\end{align*}
In the case that $H$ is the Hamiltonian (19.55), then
\begin{align*}
 \tilde{H}=\frac{1}{2}(\tilde{g}^{-1})^{\mu\nu}(q)p_{\mu}p_{\nu},\quad \tilde{g}_{\mu\nu}=\Omega g_{\mu\nu}
\end{align*}
and the above corresponds to the fact that null geodesics are invariant, up to reparametrization, under conformal transformations of the metric.

     In the case of our acoustical spacetime $(\mathcal{M}_{\epsilon_{0}},g)$, the reciprocal metric $g^{-1}$ is as we have seen singular at $\tilde{\mathcal{H}}$, 
however $\mu g^{-1}$, given in canonical acoustical coordinates by (19.37), extends smoothly to $\tilde{\mathcal{H}}$. So we take the Hamiltonian to be:
\begin{align}
 H=-p_{u}p_{t}+\mu\{-\frac{1}{2}\alpha^{-2}p_{t}^{2}+p_{t}\hat{\Xi}^{A}\slashed{p}_{A}+\frac{1}{2}(\slashed{g}^{-1})^{AB}\slashed{p}_{A}\slashed{p}_{B}\}
\end{align}
Here $p_{t}$, $p_{u}$, and $\slashed{p}_{A}\quad:\quad A=1,2$ are the momenta conjugate to the coordinates $t$, $u$, $\vartheta^{A}\quad:\quad A=1,2$, respectively. 
Thus the $\slashed{p}_{A}$ are the components of the $angular$ $momentum$ $\slashed{p}$. Given a solution $\tau\mapsto(q(\tau),p(\tau)), q=(t,u,\vartheta), p=(p_{t},
p_{u},\slashed{p})$ (This is an affinely parametrized null geodesic of the conformal acoustical metric $\mu^{-1}g$), then defining $s$ by:
\begin{align}
 \frac{ds}{d\tau}=\mu(q(\tau))
\end{align}
and taking the inverse of the mapping $\tau\mapsto s(\tau)$, then $s\mapsto (q(\tau(s)),p(\tau(s)))$ is an affinely parametrized null geodesic of the acoustical metric $g$,
$s$ being the affine parameter, and conversely.

     The canonical equations (19.56) take for the Hamiltonian (19.59) the following form:
\begin{align}
 \frac{dt}{d\tau}=\frac{\partial H}{\partial p_{t}}=-p_{u}+\mu(-\alpha^{-2}p_{t}+\hat{\Xi}^{A}\slashed{p}_{A})\\\notag
\frac{du}{d\tau}=\frac{\partial H}{\partial p_{u}}=-p_{t}\\\notag
\frac{d\vartheta^{A}}{d\tau}=\frac{\partial H}{\partial \slashed{p}_{A}}=
\mu((\slashed{g}^{-1})^{AB}\slashed{p}_{B}+p_{t}\hat{\Xi}^{A})
\end{align}
\begin{align}
 \frac{dp_{t}}{d\tau}=-\frac{\partial H}{\partial t}=-\frac{\partial\mu}{\partial t}\{-\frac{1}{2}\alpha^{-2}p^{2}_{t}+p_{t}\hat{\Xi}^{A}\slashed{p}_{A}+
\frac{1}{2}(\slashed{g}^{-1})^{AB}\slashed{p}_{A}\slashed{p}_{B}\}\\\notag
-\mu\{\alpha^{-3}\frac{\partial\alpha}{\partial t}p_{t}^{2}+p_{t}\frac{\partial\hat{\Xi}^{A}}{\partial t}\slashed{p}_{A}+\frac{1}{2}
\frac{\partial(\slashed{g}^{-1})^{AB}}{\partial t}\slashed{p}_{A}\slashed{p}_{B}\}\\\notag
\frac{dp_{u}}{d\tau}=-\frac{\partial H}{\partial u}=-\frac{\partial\mu}{\partial u}\{-\frac{1}{2}\alpha^{-2}p^{2}_{t}+p_{t}\hat{\Xi}^{A}\slashed{p}_{A}+
\frac{1}{2}(\slashed{g}^{-1})^{AB}\slashed{p}_{A}\slashed{p}_{B}\}\\\notag
-\mu\{\alpha^{-3}\frac{\partial\alpha}{\partial u}p_{t}^{2}+p_{t}\frac{\partial\hat{\Xi}^{A}}{\partial u}\slashed{p}_{A}+\frac{1}{2}
\frac{\partial(\slashed{g}^{-1})^{AB}}{\partial u}\slashed{p}_{A}\slashed{p}_{B}\}\\\notag
\frac{d\slashed{p}_{A}}{d\tau}=-\frac{\partial H}{\partial \vartheta^{A}}=-\frac{\partial\mu}{\partial \vartheta^{A}}\{-\frac{1}{2}\alpha^{-2}p^{2}_{t}+
p_{t}\hat{\Xi}^{B}\slashed{p}_{B}+\frac{1}{2}(\slashed{g}^{-1})^{BC}\slashed{p}_{B}\slashed{p}_{C}\}\\\notag
-\mu\{\alpha^{-3}\frac{\partial\alpha}{\partial \vartheta^{A}}p_{t}^{2}+p_{t}\frac{\partial\hat{\Xi}^{B}}{\partial \vartheta^{A}}\slashed{p}_{B}+\frac{1}{2}
\frac{\partial(\slashed{g}^{-1})^{BC}}{\partial \vartheta^{A}}\slashed{p}_{B}\slashed{p}_{C}\}\\\notag
\end{align}

\subsection{Asymptotic Behavior}
     Given a point $q_{0}\in\mathcal{H}\bigcup\partial_{-}\mathcal{H}$, we shall study the solutions of the system (19.61)-(19.62), subject to the condition $H=0$, 
which end at $q_{0}=(t_{0},u_{0},\vartheta_{0})$, $t_{0}=t_{*}(u_{0},\vartheta_{0})$, at $\tau=0$. The solution are then studied for $\tau\leq 0$. These are the null 
geodesics of the acoustical metric ending at $q_{0}$, that is, the past null geodesic cone of $q_{0}$, with the orientation of each geodesic reversed so that they 
become future-directed ending at $q_{0}$ instead of past directed issuing at $q_{0}$.

     Now the condition $H=0$ defines at a regular point $p\in\mathcal{M}_{\epsilon_{0}}$, where $\mu(p)>0$, a double cone in $T^{*}_{p}\mathcal{M}_{\epsilon_{0}}$ 
((19.59) is a quadratic homogeneous polynomial in $p_{\mu}$, the coordinates of $T^{*}_{p}\mathcal{M}_{\epsilon_{0}}$) and we are considering the backward part of this 
cone. However at the singular point $q_{0}$, where $\mu(q_{0})=0$, the condition $H=0$ reduces to:
\begin{align}
 p_{u}p_{t}=0
\end{align}
Thus the double cone degenerates to the two hyperplanes $p_{t}=0$ and $p_{u}=0$ and the backward part of the cone degenerates to the following three pieces:

(i) The negative half-hyperplane: $p_{t}=0,\quad p_{u}<0$.

(ii) The negative half-hyperplane: $p_{u}=0,\quad p_{t}<0$.

(iii) The plane $p_{t}=p_{u}=0$.

We thus have a $trichotomy$ of the past null geodesic cone of $q_{0}$. The null geodesics ending at $q_{0}$ such that their momentum covector at $q_{0}$ belongs to (i)
are the $outgoing$ null geodesics. These contain the generator of the characteristic
hypersurface $C_{u_{0}}$ through $q_{0}$, parametrized by $t$, which is the following solution of (19.61)-(19.62) and the condition $H=0$:
\begin{align}
 t=t_{0}+\tau,\quad u=u_{0},\quad \vartheta=\vartheta_{0};\quad p_{t}=0,\quad p_{u}=-1,\quad \slashed{p}=0
\end{align}
The null geodesic ending at $q_{0}$ such that their momentum vector at $q_{0}$ belongs to (ii) are the $incoming$ null geodesics. Finally, the null geodesics 
ending at $q_{0}$ such that their momentum vector at $q_{0}$ belongs to (iii) we simply call the $other$ null geodesics.

     To obtain an intuitive picture of how the null cone in $T^{*}_{p}\mathcal{M}_{\epsilon_{0}}$ degenerates as $p$, a regular point, approaches $q_{0}$, a point on the 
singular boundary, we consider the picture in $4$-dimensional Euclidean space with the regular coordinates $x_{1}$, $x_{2}$, $y$, $z$, taking:
\begin{align}
 x^{2}_{1}+x^{2}_{2}=(\slashed{g}^{-1})^{AB}\slashed{p}_{A}\slashed{p}_{B},\quad y=p_{u},\quad z=\alpha^{-1}p_{t}
\end{align}
 An adequate picture of the degeneration is obtained if we assume for simplicity that $\hat{\Xi}=0$ at a given $\Sigma_{t}$. The equation of the double cone $H=0$ then 
becomes:
\begin{align}
 2yz=\kappa(x_{1}^{2}+x_{2}^{2}-z^{2})
\end{align}
Note that the double cone contains the $y$ axis as well as the line:
\begin{align}
 x_{1}=x_{2}=0,\quad y=-\frac{1}{2}\kappa z
\end{align}
Consider the point $P_{0}$ on the positive $y$ axis at Euclidean distance $1$ from the origin. The coordinates of $P_{0}$ are then $(0,0,1,0)$. Consider also the point 
$P_{1}$ on the line (19.67) at Euclidean distance $1$ from the origin in the positive $z$ direction. The coordinates of $P_{1}$ are:
\begin{align*}
 (0,0,-(\kappa/2)/\sqrt{(\kappa/2)^{2}+1}, 1/\sqrt{(\kappa/2)^{2}+1})
\end{align*}
Let $H_{1}$ be the hyperplane:
\begin{align}
 z=-\lambda(y-1),\quad \lambda=\frac{1}{\sqrt{(\kappa/2)^{2}+1}+(\kappa/2)}
\end{align}
passing through the points $P_{0}$, $P_{1}$, and ruled by planes parallel to the $(x_{1},x_{2})$ plane. Then the intersection of $H_{1}$ with the double cone is the 
spheroid on $H_{1}$ which projects to the following spheroid on the $(x_{1},x_{2},y)$ hyperplane:
\begin{align}
 \frac{x^{2}_{1}+x^{2}_{2}}{a^{2}}+\frac{(y-y_{0})^{2}}{b^{2}}=1;\\\notag a=\sqrt{\frac{\lambda}{\kappa(2-\kappa\lambda)}},\quad 
b=\frac{1}{2-\kappa\lambda},\quad y_{0}=\frac{1-\kappa\lambda}{2-\kappa\lambda}
\end{align}
with semimajor axis $a$, semiminor axis $b$, centered at $y_{0}$ on the $y$ axis. The spheroid on $H_{1}$ is centered at the point $(0,0,y_{0},z_{0})$, $z_{0}=
-\lambda(y_{0}-1)$, and has semimajor axis $a'=a$ in a plane through this center parallel to the $(x_{1},x_{2})$ plane and semiminor axis $b'=\sqrt{1+\lambda^{2}}b$ in 
a line through the center in $H_{1}$ orthogonal to this plane. As the regular point $p$ approaches the singular point $q_{0}$, $\kappa\rightarrow0$, hence 
$\lambda\rightarrow 1$, $y_{0}\rightarrow 1/2$, $z_{0}\rightarrow 1/2$, $b'\rightarrow1/\sqrt{2}$, but $a'\rightarrow\infty$, in fact $\sqrt{2\kappa}a'\rightarrow1$.
The hyperplane $H_{1}$ becomes the hyperplane $z=-(y-1)$, and the spheroid on $H_{1}$, the intersection with the double cone, degenerates to the planes $y=1$, $z=0$, 
and $y=0$, $z=1$, parallel to the $(x_{1},x_{2})$ plane, at which the hyperplane $z=-(y-1)$ intersects the hyperplanes $z=0$ and $y=0$. This shows how the double cone 
degenerates to the two hyperplanes $z=0$ and $y=0$.

     We return to the study of the null geodesics ending at $q_{0}$. The value of $p_{0}=((p_{t})_{0}, (p_{u})_{0},\slashed{p}_{0})$ determines in which of the three 
classes the null geodesic belongs. Note that $p_{0}\slashed{=}0$, so that $(q_{0},p_{0})$ is a regular non-critical point of the Hamiltonian system  
(19.61)-(19.62). Now we can impose in each class a normalization condition by making use of the following remark. Let $\tau\mapsto(q(\tau),p(\tau))$ be a solution of 
the canonical equations (19.56) associated to the Hamiltonian (19.55) and corresponding to the conditions $q(0)=q_{0}$, $p(0)=p_{0}$, at $\tau=0$. Then for any positive 
constant $a$, $\tau\mapsto(q(a\tau),ap(a\tau))$ is the solution corresponding to the conditions $q(0)=q_{0}$, $p(0)=ap_{0}$. Thus the new conditions give a null geodesic 
which is simply a reparametrization of the null geodesic arising from the original conditions. This is so because if $\tilde{q}(\tau)=q(a\tau)$, $\tilde{p}(\tau)=ap(a\tau)$,
we have:
\begin{align*}
 (\frac{d\tilde{q}^{\mu}}{d\tau}(\tau)=a(\frac{dq^{\mu}}{d\tau})(a\tau)=a(\frac{\partial H}{\partial p_{\mu}})(q(a\tau),p(a\tau))\\
=a(g^{-1})^{\mu\nu}(q(a\tau))p_{\nu}(a\tau)=(g^{-1})^{\mu\nu}(\tilde{q}(\tau))\tilde{p}_{\nu}(\tau)=(\frac{\partial H}{\partial p_{\mu}})(\tilde{q}(\tau),\tilde{p}(\tau))\\
(\frac{d\tilde{p}_{\mu}}{d\tau})(\tau)=a^{2}(\frac{dp_{\mu}}{d\tau})(a\tau)=-a^{2}(\frac{\partial H}{\partial q^{\mu}})(q(a\tau),p(a\tau))\\
=-\frac{1}{2}a^{2}(\frac{\partial(g^{-1})^{\lambda\nu}}{\partial q^{\mu}})(q(a\tau))p_{\lambda}(a\tau)p_{\nu}(a\tau)\\
=-\frac{1}{2}(\frac{\partial(g^{-1})^{\lambda\nu}}{\partial q^{\mu}})(\tilde{q}(\tau))\tilde{p}_{\lambda}(\tau)\tilde{p}_{\nu}(\tau)
=-(\frac{\partial H}{\partial q^{\mu}})(\tilde{q}(\tau),\tilde{p}(\tau))
\end{align*}
 Given an $outgoing$ null geodesic ending at $q_{0}$ we can, by suitable choice of the scale factor $a$, set $(p_{0})_{u}=-1$. Thus, the outgoing null geodesics ending
 at $q_{0}$ correspond to the plane:
\begin{align}
 (p_{0})_{t}=0,\quad (p_{0})_{u}=-1\quad:\textrm{in}\quad T^{*}_{q_{0}}\mathcal{M}_{\epsilon_{0}}
\end{align}
Given an $incoming$ null geodesic ending at $q_{0}$ we can, by suitable choice of the scale factor $a$, set $(p_{0})_{t}=-1$. Thus, the incoming null geodesics ending at 
$q_{0}$ correspond to the plane:
\begin{align}
 (p_{0})_{u}=0,\quad (p_{0})_{t}=-1\quad:\textrm{in}\quad T^{*}_{q_{0}}\mathcal{M}_{\epsilon_{0}}
\end{align}
Finally, given one of the $other$ null geodesics ending at $q_{0}$ we can, by suitable choice of the scale factor $a$, set $|(\slashed{p})_{0}|=
\sqrt{(\slashed{g}^{-1})^{AB}(q_{0})(\slashed{p}_{0})_{A}(\slashed{p}_{0})_{B}}=1$. Thus, the $other$ null geodesics ending at $q_{0}$ correspond to the circle:
\begin{align}
 (p_{0})_{t}=0,\quad (p_{0})_{u}=0,\quad |\slashed{p}_{0}|=1\quad:\textrm{in}\quad T^{*}_{q_{0}}\mathcal{M}_{\epsilon_{0}}
\end{align}

    Consider the class of $outgoing$ null geodesics ending at $q_{0}$. For this class equations (19.61) evaluated at $\tau=0$ give:
\begin{align}
 (\frac{dt}{d\tau})(0)=1,\quad(\frac{du}{d\tau})(0)=0,\quad  (\frac{d\vartheta^{A}}{d\tau})(0)=0
\end{align}
Also, equations (19.62) give:
\begin{align}
 (\frac{dp_{t}}{d\tau})(0)=-\frac{1}{2}(\frac{\partial\mu}{\partial t})(q_{0})|\slashed{p}_{0}|^{2}\\\notag
(\frac{dp_{u}}{d\tau})(0)=-\frac{1}{2}(\frac{\partial\mu}{\partial u})(q_{0})|\slashed{p}_{0}|^{2}\\\notag
(\frac{d\slashed{p}_{A}}{d\tau})(0)=-\frac{1}{2}(\frac{\partial\mu}{\partial\vartheta^{A}})(q_{0})|\slashed{p}_{0}|^{2}
\end{align}
Differentiating equations (19.61) with respect to $\tau$ and evaluating the result at $\tau=0$ using (19.74) and the fact that by (19.73):
\begin{align}
 (\frac{d\mu(q(\tau))}{d\tau})_{\tau=0}=(\frac{\partial\mu}{\partial t})(q_{0})
\end{align}
we obtain:
\begin{align}
 (\frac{d^{2}t}{d\tau^{2}})(0)=\frac{1}{2}(\frac{\partial\mu}{\partial u})(q_{0})|\slashed{p}_{0}|^{2}+(\frac{\partial\mu}{\partial t})(q_{0})\hat{\Xi}^{A}(q_{0})
(\slashed{p}_{0})_{A}\\\notag
(\frac{d^{2}u}{d\tau^{2}})(0)=\frac{1}{2}(\frac{\partial\mu}{\partial t})(q_{0})|\slashed{p}_{0}|^{2}\\\notag
(\frac{d^{2}\vartheta^{A}}{d\tau^{2}})(0)=(\frac{\partial\mu}{\partial t})(q_{0})(\slashed{g}^{-1})^{AB}(q_{0})(\slashed{p}_{0})_{B}
\end{align}
In view of (19.73), (19.76) the tangent vector to an outgoing null geodesic ending at $q_{0}$ has the following expansion as $\tau\rightarrow0$:
\begin{align}
 (\frac{dt}{d\tau})(\tau)=1+[\frac{1}{2}(\frac{\partial\mu}{\partial u})(q_{0})|\slashed{p}_{0}|^{2}+(\frac{\partial\mu}{\partial t})(q_{0})\hat{\Xi}^{A}(q_{0})
(\slashed{p}_{0})_{A}]\tau +O(\tau^{2})\\\notag
(\frac{du}{d\tau})(\tau)=\frac{1}{2}(\frac{\partial\mu}{\partial t})(q_{0})|\slashed{p}_{0}|^{2}\tau +O(\tau^{2})\\\notag
(\frac{d\vartheta^{A}}{d\tau})(\tau)=(\frac{\partial\mu}{\partial t})(q_{0})(\slashed{g}^{-1})^{AB}(q_{0})(\slashed{p}_{0})_{B}\tau+O(\tau^{2})
\end{align}
We conclude that in canonical acoustical coordinates (See (19.37)) the tangent vectors to all $outgoing$ null geodesics ending at $q_{0}$ tend as we 
approach $q_{0}$ 
to the tangent vector of the generator of $C_{u_{0}}$ through $q_{0}$.

     Consider the class of $incoming$ null geodesics ending at $q_{0}$. For this class equations (19.61) evaluated at $\tau=0$ give:
\begin{align}
 (\frac{dt}{d\tau})(0)=0,\quad(\frac{du}{d\tau})(0)=1,\quad(\frac{d\vartheta^{A}}{d\tau})(0)=0
\end{align}
Also equations (19.62) give:
\begin{align}
(\frac{dp_{t}}{d\tau})(0)=-\frac{1}{2}(\frac{\partial\mu}{\partial t})(q_{0})[-\alpha^{-2}(q_{0})-2\hat{\Xi}^{A}(q_{0})(\slashed{p}_{0})_{A}+|\slashed{p}_{0}|^{2}]\\\notag
(\frac{dp_{u}}{d\tau})(0)=-\frac{1}{2}(\frac{\partial\mu}{\partial u})(q_{0})[-\alpha^{-2}(q_{0})-2\hat{\Xi}^{A}(q_{0})(\slashed{p}_{0})_{A}+|\slashed{p}_{0}|^{2}]\\\notag
(\frac{d\slashed{p}_{A}}{d\tau})(0)=-\frac{1}{2}(\frac{\partial\mu}{\partial\vartheta^{A}})(q_{0})[-\alpha^{-2}(q_{0})-2\hat{\Xi}^{B}(q_{0})(\slashed{p}_{0})_{B}
+|\slashed{p}_{0}|^{2}]
\end{align}
Differentiating equations (19.61) with respect to $\tau$ and evaluating the result at $\tau=0$ using (19.79) and the fact that by (19.78):
\begin{align}
 (\frac{d\mu(q(\tau))}{d\tau})_{\tau=0}=(\frac{\partial\mu}{\partial u})(q_{0})
\end{align}
we obtain:
\begin{align}
 (\frac{d^{2}t}{d\tau^{2}})(0)=\frac{1}{2}(\frac{\partial\mu}{\partial u})(q_{0})(\alpha^{-2}(q_{0})+|\slashed{p}_{0}|^{2})\\\notag
(\frac{d^{2}\vartheta^{A}}{d\tau^{2}})(0)=(\frac{\partial\mu}{\partial u})(q_{0})[(\slashed{g}^{-1})^{AB}(q_{0})(\slashed{p}_{0})_{B}-\hat{\Xi}^{A}(q_{0})]
\end{align}
The right-hand side of the first equation in (19.81) is negative when $q_{0}\in\mathcal{H}$. However in the case that $q_{0}\in\partial_{-}\mathcal{H}$, we have:
\begin{align}
 (\frac{d\mu(q(\tau))}{d\tau})_{\tau=0}=0,\quad (\frac{d^{2}t}{d\tau^{2}})(0)=0,\quad (\frac{d^{2}\vartheta^{A}}{d\tau^{2}})(0)=0
\end{align}
We must then proceed to the third derivatives. Using (19.78) and (19.82) we find in this case:
\begin{align}
 (\frac{d^{2}\mu(q(\tau))}{d\tau^{2}})_{\tau=0}=(\frac{\partial^{2}\mu}{\partial u^{2}})(q_{0})
\end{align}
and differentiating the second equation of (19.62) with respect to $\tau$ and evaluating the result at $\tau=0$ gives:
\begin{align}
 (\frac{d^{2}p_{u}}{d\tau^{2}})(0)=\frac{1}{2}(\frac{\partial^{2}\mu}{\partial u^{2}})(q_{0})[\alpha^{-2}(q_{0})+2\hat{\Xi}^{A}(q_{0})(\slashed{p}_{0})_{A}-
|\slashed{p}_{0}|^{2}]
\end{align}
Differentiating then the second and third equations of (19.61) with respect to $\tau$ a second time and using (15.83) and (15.84) we then obtain that in the case 
$q_{0}\in\partial_{-}\mathcal{H}$:
\begin{align}
 (\frac{d^{3}t}{d\tau^{3}})(0)=\frac{1}{2}(\frac{\partial^{2}\mu}{\partial u^{2}})(q_{0})(\alpha^{-2}(q_{0})+|\slashed{p}_{0}|^{2})\\\notag
(\frac{d^{3}\vartheta^{A}}{d\tau^{3}})(0)=(\frac{\partial^{2}\mu}{\partial u^{2}})(q_{0})[(\slashed{g}^{-1})^{AB}(q_{0})(\slashed{p}_{0})_{B}-\hat{\Xi}^{A}(q_{0})]
\end{align}
Note that the right-hand side of the first equation of (19.85) is positive. The tangent vector to an incoming null geodesic ending at $q_{0}$ has the following 
expansion as $\tau\rightarrow0$:
\begin{align}
 (\frac{du}{d\tau})(\tau)=1+O(\tau)\\\notag
(\frac{dt}{d\tau})(\tau)=\frac{1}{2}(\frac{\partial\mu}{\partial u})(q_{0})(\alpha^{-2}(q_{0})+|\slashed{p}_{0}|^{2})\tau+O(\tau^{2})\\\notag
(\frac{d\vartheta^{A}}{d\tau})(\tau)=(\frac{\partial\mu}{\partial u})(q_{0})[(\slashed{g}^{-1})^{AB}(q_{0})(\slashed{p}_{0})_{B}-\hat{\Xi}^{A}(q_{0})]\tau+O(\tau^{2})
\end{align}
by (19.78), (19.81), in the case that $q_{0}\in\mathcal{H}$, and:
\begin{align}
 (\frac{du}{d\tau})(\tau)=1+O(\tau)\\\notag
(\frac{dt}{d\tau})(\tau)=\frac{1}{4}(\frac{\partial^{2}\mu}{\partial u^{2}})(q_{0})(\alpha^{-2}(q_{0})+|\slashed{p}_{0}|^{2})\tau^{2}+O(\tau^{3})\\\notag
(\frac{d\vartheta^{A}}{d\tau})(\tau)=\frac{1}{2}(\frac{\partial^{2}\mu}{\partial u^{2}})(q_{0})[(\slashed{g}^{-1})^{AB}(q_{0})(\slashed{p}_{0})_{B}
-\hat{\Xi}^{A}(q_{0})]\tau^{2}+O(\tau^{3})
\end{align}
by (19.78), (19.82), (19.85), in the case that $q_{0}\in\partial_{-}\mathcal{H}$.

     Consider finally the class of the $other$ null geodesics ending at $q_{0}$. For this class equations (19.61) and (19.62) evaluated at $\tau=0$ give:
\begin{align}
 (\frac{dt}{d\tau})(0)=0,\quad(\frac{du}{d\tau})(0)=0,\quad (\frac{d\vartheta^{A}}{d\tau})(0)=0
\end{align}
\begin{align}
 (\frac{dp_{t}}{d\tau})(0)=-\frac{1}{2}(\frac{\partial\mu}{\partial t})(q_{0})\\\notag
(\frac{dp_{u}}{d\tau})(0)=-\frac{1}{2}(\frac{\partial\mu}{\partial u})(q_{0})\\\notag
(\frac{d\slashed{p}_{A}}{d\tau})(0)=-\frac{1}{2}\frac{\partial\mu}{\partial\vartheta^{A}}
\end{align}
Note that the right-hand side of the first equation of (19.89) is positive. Differentiating equation (19.61) with respect to $\tau$ and evaluating the result at $\tau=0$
using (19.89) and the fact that by (19.88):
\begin{align}
 (\frac{d\mu(q(\tau))}{d\tau})_{\tau=0}=0
\end{align}
we obtain:
\begin{align}
 (\frac{d^{2}t}{d\tau^{2}})(0)=\frac{1}{2}(\frac{\partial\mu}{\partial u})(q_{0})\\\notag
(\frac{d^{2}u}{d\tau^{2}})(0)=\frac{1}{2}(\frac{\partial\mu}{\partial t})(q_{0})\\\notag
(\frac{d^{2}\vartheta^{A}}{d\tau^{2}})(0)=0
\end{align}
Note that the right-hand side of the second equation of (19.91) is negative, while the right-hand side of the first equation of (19.91) is negative when $q_{0}\in\mathcal{H}$
and vanishes when $q_{0}\in\partial_{-}\mathcal{H}$. Differentiating the third equation of (19.61) a second time with respect to $\tau$ and evaluating the result at 
$\tau=0$ yields:
\begin{align}
 (\frac{d^{3}\vartheta^{A}}{d\tau^{3}})(0)=(\frac{d^{2}\mu(q(\tau))}{d\tau^{2}})_{\tau=0}(\slashed{g}^{-1})^{AB}(q_{0})(\slashed{p}_{0})_{B}
\end{align}
Now from (19.88) and (19.91) we obtain:
\begin{align}
 (\frac{d^{2}\mu(q(\tau))}{d\tau^{2}})_{\tau=0}=(\frac{d^{2}t}{d\tau^{2}})(0)(\frac{\partial\mu}{\partial t})(q_{0})+(\frac{d^{2}u}{d\tau^{2}})(0)
(\frac{\partial\mu}{\partial u})(q_{0})\\\notag
+(\frac{d^{2}\vartheta^{A}}{d\tau^{2}})(0)(\frac{\partial\mu}{\partial\vartheta^{A}})(q_{0})=(\frac{\partial\mu}{\partial t}\frac{\partial\mu}{\partial u})(q_{0})
\end{align}
Substituting in (19.92) then yields:
\begin{align}
 (\frac{d^{3}\vartheta^{A}}{d\tau^{3}})(0)=(\frac{\partial\mu}{\partial t}\frac{\partial\mu}{\partial u})(q_{0})(\slashed{g}^{-1})^{AB}(q_{0})(\slashed{p}_{0})_{B}
\end{align}
In the case that $q_{0}\in\partial_{-}\mathcal{H}$ the right-hand side of the first equation of (19.91), of (19.93), and of (19.94) all vanish:
\begin{align}
 (\frac{d^{2}\mu(q(\tau))}{d\tau^{2}})_{\tau=0}=0,\quad(\frac{d^{2}t}{d\tau^{2}})(0)=0,\quad(\frac{d^{2}\vartheta^{A}}{d\tau^{2}})(0)=
(\frac{d^{3}\vartheta^{A}}{d\tau^{3}})(0)=0
\end{align}
Differentiating the first equation of (19.61) a second time with respect to $\tau$ and evaluating the result at $\tau=0$ then yields:
\begin{align}
 (\frac{d^{3}t}{d\tau^{3}})(0)=-(\frac{d^{2}p_{u}}{d\tau^{2}})(0)
\end{align}
Also, differentiating the second equation of (19.62) with respect to $\tau$ and evaluating the result at $\tau=0$ yields:
\begin{align}
 (\frac{d^{2}p_{u}}{d\tau^{2}})(0)=-\frac{1}{2}(\frac{d(\partial\mu/\partial u)(q(\tau))}{d\tau})_{\tau=0}
\end{align}
However, from (19.88):
\begin{align}
 (\frac{d(\partial\mu/\partial u)(q(\tau))}{d\tau})_{\tau=0}=(\frac{\partial^{2}\mu}{\partial t\partial u})(q_{0})(\frac{dt}{d\tau})(0)\\\notag
+(\frac{\partial^{2}\mu}{\partial u^{2}})(q_{0})(\frac{du}{d\tau})(0)+(\frac{\partial^{2}\mu}{\partial\vartheta^{A}\partial u})(q_{0})(\frac{d\vartheta^{A}}{d\tau})(0)=0
\end{align}
hence:
\begin{align}
 (\frac{d^{3}t}{d\tau^{3}})(0)=(\frac{d^{2}p_{u}}{d\tau^{2}})(0)=0
\end{align}
We must proceed to the fourth derivatives. Differentiating the first equation of (19.61) a third time with respect to $\tau$ and evaluating the result at $\tau=0$ yields:
\begin{align}
 (\frac{d^{4}t}{d\tau^{4}})(0)=-(\frac{d^{3}p_{u}}{d\tau^{3}})(0)+(\frac{d^{3}\mu(q(\tau))}{d\tau^{3}})_{\tau=0}\hat{\Xi}^{A}(q_{0})(\slashed{p}_{0})_{A}
\end{align}
Also, differentiating the second equation of (19.62) a second time with respect to $\tau$ and evaluating the result at $\tau=0$ yields:
\begin{align}
 (\frac{d^{3}p_{u}}{d\tau^{3}})(0)=-\frac{1}{2}(\frac{d^{2}(\partial\mu/\partial u)(q(\tau))}{d\tau^{2}})_{\tau=0}
\end{align}
Now, in view of (19.88) and (19.91):
\begin{align}
 (\frac{d^{2}(\partial\mu/\partial u)(q(\tau))}{d\tau^{2}})_{\tau=0}=(\frac{\partial^{2}\mu}{\partial t\partial u})(q_{0})(\frac{d^{2}t}{d\tau^{2}})(0)\\\notag
+(\frac{\partial^{2}\mu}{\partial u^{2}})(q_{0})(\frac{d^{2}u}{d\tau^{2}})(0)+(\frac{\partial^{2}\mu}{\partial\vartheta^{A}\partial u})(q_{0})
(\frac{d^{2}\vartheta^{A}}{d\tau^{2}})(0)\\\notag
=(\frac{\partial^{2}\mu}{\partial u^{2}})(q_{0})(\frac{d^{2}u}{d\tau^{2}})(0)=\frac{1}{2}(\frac{\partial\mu}{\partial t}\frac{\partial^{2}\mu}{\partial u^{2}})(q_{0})
\end{align}
On the other hand, by (19.88), (19.91), (19.95), (19.98), (19.99), and the fact that $(\partial\mu/\partial u)(q_{0})=0$,
\begin{align}
 (\frac{d^{3}\mu(q(\tau))}{d\tau^{3}})_{\tau=0}=(\frac{\partial\mu}{\partial t})(q_{0})(\frac{d^{3}t}{d\tau^{3}})(0)+(\frac{\partial\mu}{\partial u})(q_{0})
(\frac{d^{3}u}{d\tau^{3}})(0)\\\notag+(\frac{\partial\mu}{\partial\vartheta^{A}})(q_{0})(\frac{d^{3}\vartheta^{A}}{d\tau^{3}})(0)
+2(\frac{d(\partial\mu/\partial t)(q(\tau))}{d\tau})_{\tau=0}(\frac{d^{2}t}{d\tau^{2}})(0)\\\notag
+2(\frac{d(\partial\mu/\partial u)(q(\tau))}{d\tau})_{\tau=0}(\frac{d^{2}u}{d\tau^{2}})(0)+
2(\frac{d(\partial\mu/\partial\vartheta^{A})(q(\tau))}{d\tau})_{\tau=0}(\frac{d^{2}\vartheta^{A}}{d\tau^{2}})(0)=0
\end{align}
Substituting then (19.102) in (19.101) and the result in (19.100), yields, in view of (19.103),
\begin{align}
 (\frac{d^{4}t}{d\tau^{4}})(0)=\frac{1}{4}(\frac{\partial\mu}{\partial t}\frac{\partial^{2}\mu}{\partial u^{2}})(q_{0})
\end{align}
Differentiating the third equation of (19.61) a third time and evaluating the result at $\tau=0$ we obtain, taking into account (19.103),
\begin{align}
 (\frac{d^{4}\vartheta^{A}}{d\tau^{4}})(0)=(\frac{d^{3}\mu(q(\tau))}{d\tau^{3}})_{\tau=0}(\slashed{g}^{-1})^{AB}(q_{0})(\slashed{p}_{0})_{B}=0
\end{align}
Finally, we differentiate the third equation of (19.61) a fourth time and evaluate the result at $\tau=0$ to obtain:
\begin{align}
 (\frac{d^{5}\vartheta^{A}}{d\tau^{5}})(0)=(\frac{d^{4}\mu(q(\tau))}{d\tau^{4}})_{\tau=0}(\slashed{g}^{-1})^{AB}(q_{0})(\slashed{p}_{0})_{B}
\end{align}
and we have:
\begin{align}
 (\frac{d^{4}\mu(q(\tau))}{d\tau^{4}})_{\tau=0}\\\notag
=(\frac{\partial\mu}{\partial t})(q_{0})(\frac{d^{4}t}{d\tau^{4}})(0)+3(\frac{d^{2}(\partial\mu/\partial u)(q(\tau))}{d\tau^{2}})_{\tau=0}(\frac{d^{2}u}{d\tau^{2}})(0)
\\\notag
=((\frac{\partial\mu}{\partial t})^{2}\frac{\partial^{2}\mu}{\partial u^{2}})(q_{0})
\end{align}
by (19.102), (19.104), and the second of (19.91). Substituting then in (19.106) yields:
\begin{align}
 (\frac{d^{5}\vartheta^{A}}{d\tau^{5}})(0)=((\frac{\partial\mu}{\partial t})^{2}\frac{\partial^{2}\mu}{\partial u^{2}})(q_{0})(\slashed{g}^{-1})^{AB}(q_{0})
(\slashed{p}_{0})_{B}
\end{align}
According to the above results the tangent vector to any of the other null geodesics ending at $q_{0}$ has the following expansion as $\tau\rightarrow0$:
\begin{align}
 (\frac{dt}{d\tau})(\tau)=\frac{1}{2}(\frac{\partial\mu}{\partial u})(q_{0})\tau+O(\tau^{2})\\\notag
(\frac{du}{d\tau})(\tau)=\frac{1}{2}(\frac{\partial\mu}{\partial t})(q_{0})\tau+O(\tau^{2})\\\notag
(\frac{d\vartheta^{A}}{d\tau})(\tau)=\frac{1}{2}(\frac{\partial\mu}{\partial t}\frac{\partial\mu}{\partial u})(q_{0})(\slashed{g}^{-1})^{AB}(\slashed{p}_{0})_{B}\tau^{2}
+O(\tau^{3})
\end{align}
by (19.91), (19.94), in the case $q_{0}\in\mathcal{H}$, and:
\begin{align}
 (\frac{dt}{d\tau})(\tau)=\frac{1}{24}(\frac{\partial\mu}{\partial t}\frac{\partial^{2}\mu}{\partial u^{2}})(q_{0})\tau^{3}+O(\tau^{4})\\\notag
(\frac{du}{d\tau})(\tau)=\frac{1}{2}(\frac{\partial\mu}{\partial t})(q_{0})\tau+O(\tau^{2})\\\notag
(\frac{d\vartheta^{A}}{d\tau})(\tau)=\frac{1}{24}((\frac{\partial\mu}{\partial t})^{2}\frac{\partial^{2}\mu}{\partial u^{2}})(q_{0})
(\slashed{g}^{-1})^{AB}(q_{0})(\slashed{p}_{0})_{B}\tau^{4}+O(\tau^{5})
\end{align}
in the case that $q_{0}\in\partial_{-}\mathcal{H}$.

      To obtain the picture in Galileo spacetime, we now investigate the behavior of the tangent vector to a null geodesic ending at $q_{0}$, in each of the three 
classes, as $t\rightarrow t_{0}$, in rectangular coordinates. The rectangular components $\tilde{L}'^{\mu}$ of the tangent vector $\tilde{L}'$ of a
null geodesic ending at $q_{0}$, parametrized by $\tau$, are given by:
\begin{align}
 \tilde{L}'^{\mu}=\frac{dx^{\mu}}{d\tau}=\frac{\partial x^{\mu}}{\partial t}\frac{dt}{d\tau}+\frac{\partial x^{\mu}}{\partial u}\frac{du}{d\tau}+
\frac{\partial x^{\mu}}{\partial\vartheta^{A}}\frac{d\vartheta^{A}}{d\tau}
\end{align}
If the null geodesic is parametrized by $t$ instead of $\tau$ the tangent vector is:
\begin{align}
 L'=\frac{\tilde{L}'}{dt/d\tau}
\end{align}
and its rectangular components are given by:
\begin{align}
 L'^{\mu}=\frac{\partial x^{\mu}}{\partial t}+\frac{\partial x^{\mu}}{\partial u}\frac{du}{dt}+\frac{\partial x^{\mu}}{\partial\vartheta^{A}}
\frac{d\vartheta^{A}}{dt}
\end{align}
where:
\begin{align}
 \frac{du}{dt}=\frac{du/d\tau}{dt/d\tau},\quad \frac{d\vartheta^{A}/d\tau}{dt/d\tau}
\end{align}
Recalling from Chapter 2 that:
\begin{align}
 \frac{\partial x^{\mu}}{\partial t}=L^{\mu},\quad \frac{\partial x^{0}}{\partial u}=\frac{\partial x^{0}}{\partial\vartheta^{A}}=0,\quad
\frac{\partial x^{i}}{\partial u}=T^{i}+\Xi^{A}X^{i}_{A},\quad\frac{\partial x^{i}}{\partial\vartheta^{A}}=X^{i}_{A}
\end{align}
we then obtain:
\begin{align}
 L'^{0}=1,\quad L'^{i}=L^{i}+\mu(\alpha^{-1}\hat{T}^{i}+\hat{\Xi}^{A}X_{A}^{i})\frac{du}{dt}+X^{i}_{A}\frac{d\vartheta^{A}}{dt}
\end{align}

     Consider the class of the $outgoing$ null geodesics ending at $q_{0}\in\mathcal{H}\bigcup\partial_{-}\mathcal{H}$. Then from (19.77) we have:
\begin{align}
 \lim_{t\rightarrow t_{0}}\frac{du}{dt}=0,\quad \lim_{t\rightarrow t_{0}}\frac{d\vartheta^{A}}{dt}=0
\end{align}
Since $\lim_{t\rightarrow t_{0}}\mu=0$, it follows that for all $outgoing$ null geodesics ending at $q_{0}$ we have:
\begin{align}
 \lim_{t\rightarrow t_{0}}L'^{i}=L^{i}(q_{0})
\end{align}
We conclude that in Galileo spacetime the tangent vectors to all $outgoing$ null geodesics ending at $q_{0}$ tend as we approach $q_{0}$ to the tangent vector $L(q_{0})$
of the generator of $C_{u_{0}}$ through $q_{0}$. The vector $L(q_{0})$ is therefore an $invariant$ $vector$ independent of the choice of acoustical function $u$.
From the point of view of Galileo spacetime this vector is equal, up to a multiplicative constant, to the invariant vector $V(q_{0})$, defined by (19.23), since
the vector $T(q_{0})$ vanishes from the point of view of Galileo spacetime (recalling the definition of $\Pi_{*}$).

     Consider the class of the $incoming$ null geodesics ending at $q_{0}\in\mathcal{H}\bigcup\partial_{-}\mathcal{H}$. If $q_{0}\in\mathcal{H}$, then taking into account 
the fact that by (19.80):
\begin{align}
 \mu(q(\tau))=(\frac{\partial\mu}{\partial u})(q_{0})\tau+O(\tau^{2})
\end{align}
we obtain from (19.86):
\begin{align}
 \lim_{t\rightarrow t_{0}}\mu\frac{du}{dt}=\frac{1}{(1/2)(\alpha^{-2}(q_{0})+|\slashed{p}_{0}|^{2})},\quad
\lim_{t\rightarrow t_{0}}\frac{d\vartheta^{A}}{dt}=\frac{(\slashed{g}^{-1})^{AB}(\slashed{p})_{B}-\hat{\Xi}^{A}(q_{0})}{(1/2)(\alpha^{-2}(q_{0})+|\slashed{p}_{0}|^{2})}
\end{align}
If $q_{0}\in\partial_{-}\mathcal{H}$, then taking into account the fact that by (19.82), (19.83):
\begin{align}
 \mu(q(\tau))=\frac{1}{2}(\frac{\partial^{2}\mu}{\partial u^{2}})(q_{0})\tau^{2}+O(\tau^{3})
\end{align}
we obtain from (19.87):
\begin{align}
 \lim_{t\rightarrow t_{0}}\mu\frac{du}{dt}=\frac{1}{(1/2)(\alpha^{-2}(q_{0})+|\slashed{p}_{0}|^{2})},\quad
\lim_{t\rightarrow t_{0}}\frac{d\vartheta^{A}}{dt}=\frac{(\slashed{g}^{-1})^{AB}(\slashed{p})_{B}-\hat{\Xi}^{A}(q_{0})}
{(1/2)(\alpha^{-2}(q_{0})+|\slashed{p}_{0}|^{2})}
\end{align}
which coincides with (19.120). It follows that in both cases we have:
\begin{align}
 \lim_{t\rightarrow t_{0}}L'^{i}=L^{i}(q_{0})+\frac{\alpha^{-1}(q_{0})\hat{T}^{i}(q_{0})}{(1/2)(\alpha^{-2}(q_{0})+|\slashed{p}_{0}|^{2})}
+\frac{(\slashed{g}^{-1})^{AB}(q_{0})X^{i}_{A}(q_{0})(\slashed{p}_{0})_{B}}{(1/2)(\alpha^{-2}(q_{0})+|\slashed{p}_{0}|^{2})}
\end{align}
Note that this limit depends on $\slashed{p}_{0}$. For $\slashed{p}_{0}=0$ the limit vector is $(L+2\alpha\hat{T})(q_{0})=\alpha^{2}(q_{0})\hat{\underline{L}}(q_{0})$,
where $\hat{\underline{L}}(q_{0})$ is the conjugate to $L(q_{0})$, the null vector in the orthogonal complement of $T_{q_{0}}S_{t_{0},u_{0}}$. As $\slashed{p}_{0}$ 
tends to $\infty$ the limit vector tends to $L(q_{0})$.

     Consider finally the class of the $other$ null geodesics ending at $q_{0}\in\mathcal{H}\bigcup\partial_{-}\mathcal{H}$. If $q_{0}\in\mathcal{H}$, then taking into 
account the fact that by (19.90), (19.93):
\begin{align}
 \mu(q(\tau))=\frac{1}{2}(\frac{\partial\mu}{\partial t}\frac{\partial\mu}{\partial u})(q_{0})\tau^{2}+O(\tau^{3})
\end{align}
we obtain from (19.109):
\begin{align}
 \lim_{t\rightarrow t_{0}}\mu\frac{du}{dt}=0,\quad \lim_{t\rightarrow t_{0}}\frac{d\vartheta^{A}}{dt}=0
\end{align}
If $q_{0}\in\partial_{-}\mathcal{H}$, then taking into account the fact that by (19.90), (19.95), (19.103), (19.107):
\begin{align}
 \mu(q(\tau))=\frac{1}{24}((\frac{\partial\mu}{\partial t})^{2}\frac{\partial^{2}\mu}{\partial u^{2}})(q_{0})\tau^{4}+O(\tau^{5})
\end{align}
we obtain from (19.110):
\begin{align}
 \lim_{t\rightarrow t_{0}}\mu\frac{du}{dt}=0,\quad\lim_{t\rightarrow t_{0}}\frac{d\vartheta^{A}}{dt}=0
\end{align}
as well. It follows that in both cases we have simply:
\begin{align}
 \lim_{t\rightarrow t_{0}}L'^{i}=L^{i}(q_{0})
\end{align}
Thus, from the point of view of Galileo spacetime, the tangent vectors to all the other null geodesics ending at $q_{0}$ tend to the vector $L(q_{0})$ as 
we approach $q_{0}$.

      We now consider, for a fixed point $q_{0}\in\mathcal{H}\bigcup\partial_{-}\mathcal{H}$ and each of the three classes of null geodesics ending at $q_{0}$,
the mapping:
\begin{align}
 (\tau; \slashed{p}_{0})\mapsto(t,u,\vartheta)=F(\tau;\slashed{p}_{0}),\quad F=(F^{t},F^{u},\slashed{F})
\end{align}
where for the $outgoing$ and $incoming$ null geodesics: $\slashed{p}_{0}\in T_{\vartheta_{0}}S^{2}$ which may be identified with $\mathbb{R}^{2}$,
while for the $other$ null geodesics $\slashed{p}_{0}$ takes values in the unit circle in $T_{\vartheta_{0}}S^{2}$, which may be identified with $S^{1}\subset\mathbb{R}^{2}$.
In each of the three cases the mapping (19.129) is smooth (by standard o.d.e. theory), the range of $F^{t}$ being restricted to $t\geq0$.

     Let us consider first the mapping:
\begin{align}
 \tau\mapsto F^{t}(\tau;\slashed{p}_{0})
\end{align}
for fixed $\slashed{p}_{0}$. We shall presently show that this mapping is a diffeomorphism of $[\tau_{0},0]$ onto $[0,t_{0}]$ in the case of the $outgoing$ null geodesic,
a diffeomorphism of $[\tau_{0},0)$ onto $[0,t_{0})$ in the case of $incoming$ null geodesics as well as in the case of the $other$ null geodesics, the difference of the latter
from the former being that $(dt/d\tau)(0)=1$ for the $outgoing$ ((19.77)), while $(dt/d\tau)(0)=0$ for the $incoming$ and the $other$ ((19.86), (19.87), (19.109), (19.110)).
Here $\tau_{0}=\tau_{0}(\slashed{p}_{0})$ ($p_{t}$ and $p_{u}$ are constants). Now the expressions (19.77) and (19.86), (19.87), (19.109), (19.110), show that 
$(dt/d\tau)(\tau)>0$ for $\tau\in[\tau_{1},0)$, $\tau_{1}$ suitably small. Thus either $(dt/d\tau)(\tau)>0$ for all $\tau\in[\tau_{0},0)$, or there is a first $\tau_{*}$ 
where $(dt/d\tau)(\tau_{*})=0$. Now in general, for a Hamiltonian of the form (19.55), the condition $H=0$ is equivalent to the following condition on the tangent vector:
\begin{align}
 \frac{1}{2}g_{\mu\nu}\frac{dq^{\mu}}{d\tau}\frac{dq^{\nu}}{d\tau}=0
\end{align}
In the present case, for our Hamiltonian (19.59), the condition $H=0$ is equivalent to the condition:
\begin{align}
 -\frac{dt}{d\tau}\frac{du}{d\tau}+\frac{1}{2}\alpha^{-2}\mu(\frac{du}{d\tau})^{2}+\frac{1}{2}\slashed{g}_{AB}(\frac{d\vartheta^{A}}{d\tau}+\mu\hat{\Xi}^{A}\frac{du}{d\tau})
(\frac{d\vartheta^{B}}{d\tau}+\mu\hat{\Xi}^{B}\frac{du}{d\tau})=0
\end{align}
Obviously, at $\tau_{*}$ we have:
\begin{align*}
 \frac{dt}{d\tau}(\tau_{*})=\frac{du}{d\tau}(\tau_{*})=\frac{d\vartheta^{A}}{d\tau}(\tau_{*})=0
\end{align*}
By (19.61), we have:
\begin{align*}
 p_{t}(\tau_{*})=p_{u}(\tau_{*})=\slashed{p}_{A}(\tau_{*})=0
\end{align*}
Then with: $(t_{*}, u_{*}, \vartheta_{*})=(t(\tau_{*}), u(\tau_{*}), \vartheta(\tau_{*}))$,
\begin{align*}
 ((t,u,\vartheta); (p_{t},p_{u},\slashed{p}))=((t_{*}, u_{*}, \vartheta_{*}); (0,0,0))
\end{align*}
is a solution of the Hamiltonian system (19.61)-(19.62) coinciding with our solution at $\tau=\tau_{*}$. By the uniqueness theorem for the Hamiltonian system (19.61)-(19.62)
our null geodesic must coincide with the above constant solution for all $\tau$, in particular we must have: $t(0)=t_{*}$ in contradiction with the fact that $t(0)=t_{0}
>t_{*}$. This rules out the second alternative above and we conclude that $(dt/d\tau)(\tau)>0$ for all $\tau\in[\tau_{0},0)$. It follows that the mapping (19.130) is a 
diffeomorphism of $[\tau_{0},0]$ onto $[0,t_{0}]$ in the case of the $outgoing$ null geodesics, a diffeomorphism of $[\tau_{0},0)$ onto $[0,t_{0})$ in the cases of the 
$incoming$ and the $other$ null geodesics. We can thus invert the mapping (19.130) obtaining a continuous mapping of $[0,t_{0}]$ onto $[\tau_{0},0]$ by:
\begin{align}
 t\mapsto (F^{t})^{-1}(t;\slashed{p}_{0})
\end{align}
which is smooth on $[0,t_{0}]$ in the case of the $outgoing$ null geodesics, on $[0,t_{0})$ in the case of the $incoming$ and the $other$ null geodesics.

     Consider the mapping:
\begin{align}
 (t;\slashed{p}_{0})\mapsto(t,u,\vartheta)=G(t;\slashed{p}_{0})=F((F^{t})^{-1}(t;\slashed{p}_{0});\slashed{p}_{0})
\end{align}
and for fixed $t\in[0,t_{0}]$ the mapping:
\begin{align}
 \slashed{p}_{0}\mapsto(t,u,\vartheta)=G_{t}(\slashed{p}_{0}) \quad\textrm{where}\quad G_{t}(\slashed{p}_{0})=G(t;\slashed{p}_{0})
\end{align}
where, as above, $\slashed{p}_{0}\in\mathbb{R}^{2}$ for the $outgoing$ and $incoming$ null geodesics, while $\slashed{p}_{0}\in S^{1}$ for the $other$ null geodesics. 
This is the domain of the mapping (19.135), while its range is the hypersurface $\{t\}\times(-\infty,\epsilon_{0})\times S^{2}$ in the extension of the manifold 
$\mathcal{M}_{\epsilon_{0}}$ to:
\begin{align*}
 \mathcal{M}^{e}_{\epsilon_{0}}=\mathcal{M}_{\epsilon_{0}}\bigcup \mathcal{M}^{E}=[0,\infty)\times (-\infty,\epsilon_{0})\times S^{2}
\end{align*}
where
\begin{align*}
 \mathcal{M}^{E}=[0,\infty)\times (-\infty,0)\times S^{2}
\end{align*}
corresponds to the domain of the surrounding constant state. The hypersurface $\{t\}\times (-\infty,\epsilon_{0})\times S^{2}$ corresponds to the domain in the 
hyperplane $\Sigma_{t}$ in Galileo spacetime where $u>\epsilon_{0}$.

     The mapping (19.135) is smooth for each $t\in[0,t_{0})$, and for $t=t_{0}$, $G_{t_{0}}$ is the constant mapping:
\begin{align}
 G_{t_{0}}(\slashed{p}_{0})=G(t_{0};\slashed{p}_{0})=q_{0}
\end{align}
Also, writing:
\begin{align}
 G=(G^{t},G^{u},\slashed{G}),\quad G_{t}=(G^{t}_{t},G^{u}_{t},\slashed{G}_{t})
\end{align}
we have:
\begin{align}
 G^{t}_{t}=G^{t}=t
\end{align}

     We have:
\begin{align}
 (\frac{\partial\slashed{G}^{A}}{\partial t})(t;\slashed{p}_{0})=(\frac{d\vartheta^{A}}{dt})(t)=\frac{(d\vartheta^{A}/d\tau)(\tau)}{(dt/d\tau)(\tau)}
\end{align}
and:
\begin{align}
 (\frac{\partial^{2}\slashed{G}^{A}}{\partial t^{2}})(t;\slashed{p}_{0})=(\frac{d^{2}\vartheta^{A}}{dt^{2}})(t)\\\notag
=(\frac{(d^{2}\vartheta^{A}/d\tau^{2})(dt/d\tau)-(d^{2}t/d\tau^{2})(d\vartheta^{A}/d\tau)}{(dt/d\tau)^{3}})(\tau)
\end{align}

    In the case of the $outgoing$ null geodesics we obtain, from (19.73):
\begin{align}
 (\frac{\partial\slashed{G}^{A}}{\partial t})(t_{0};\slashed{p}_{0})=0
\end{align}
and, from (19.76):
\begin{align}
 (\frac{\partial^{2}\slashed{G}^{A}}{\partial t^{2}})(t_{0};\slashed{p}_{0})=(\frac{\partial\mu}{\partial t})(q_{0})(\slashed{g}^{-1})^{AB}(q_{0})(\slashed{p}_{0})_{B}
\end{align}
Let us consider the matrix $\partial\slashed{G}^{A}_{t}/\partial(\slashed{p}_{0})_{B}$. By (19.141), (19.142):
\begin{align}
 (\frac{\partial}{\partial t}(\frac{\partial\slashed{G}^{A}}{\partial(\slashed{p}_{0})_{B}}))(t_{0};\slashed{p}_{0})=0,\\\notag
(\frac{\partial^{2}}{\partial t^{2}}(\frac{\partial\slashed{G}^{A}}{\partial(\slashed{p}_{0})_{B}}))(t_{0};\slashed{p}_{0})
=(\frac{\partial\mu}{\partial t})(q_{0})(\slashed{g}^{-1})^{AB}(q_{0})
\end{align}
It follows that:
\begin{align}
 (\frac{\partial\slashed{G}^{A}_{t}}{\partial(\slashed{p}_{0})_{B}})(\slashed{p}_{0})=
\frac{1}{2}(\frac{\partial\mu}{\partial t})(q_{0})(\slashed{g}^{-1})^{AB}(q_{0})((t-t_{0})^{2})+O((t-t_{0})^{3})
\end{align}
and:
\begin{align}
 (\det\partial\slashed{G}_{t}/\partial\slashed{p}_{0})=
\frac{1}{4}(\frac{\partial\mu}{\partial t})^{2}(q_{0})\det\slashed{g}^{-1}(t-t_{0})^{4}+O((t-t_{0})^{5})
\end{align}
Therefore in the case of the $outgoing$ null geodesics the mapping (19.135) is of maximal rank at each $\slashed{p}_{0}\in\mathbb{R}^{2}$ for $t$ suitably close to $t_{0}$.
By the implicit function theorem we conclude that in the case of the $outgoing$ null geodesics the image of $\mathbb{R}^{2}$ by the mapping (19.135) is, for $t$ suitably
close to $t_{0}$, an immersed disk in the hypersurface $\{t\}\times (-\infty,\epsilon_{0})\times S^{2}\subset\mathcal{M}^{e}_{\epsilon_{0}}$, which corresponds to the 
hyperplane $\Sigma_{t}$ in Galileo spacetime. To show that this image is in fact an embedded disk we show that the mapping (19.135) is one-to-one for $t$ suitably close 
to $t_{0}$. For $|\slashed{p}_{0}|\leq B$ and any  given $B>0$, this follows from the fact that by (19.141), (19.142):
\begin{align}
 \slashed{G}^{A}_{t}(\slashed{p}_{0})=\frac{1}{2}(\frac{\partial\mu}{\partial t})(q_{0})(\slashed{g}^{-1})^{AB}(q_{0})(\slashed{p}_{0})_{B}(t-t_{0})^{2}
+O((t-t_{0})^{3})
\end{align}
provided that $t$ is suitably close to $t_{0}$ depending on $B$ (This guarantees the first term on the right of (19.146) is the main part of $\slashed{G}^{A}_{t}
(\slashed{p}_{0})$). To handle the case that $|\slashed{p}_{0}|>B$ (Now, $B$ is fixed), we choose a different renormalization of the momentum covector at $q_{0}$ for 
the class of $outgoing$ geodesics ending at $q_{0}$ than that given by (19.70). We now choose the scale factor $a$ so that :
\begin{align}
 |\slashed{p}_{0}|=1
\end{align}
Then $p_{u}=-1/|\slashed{p}_{0}|$. This corresponds to the diffeomorphism:
\begin{align}
 f: \slashed{p}_{0}\mapsto (-\frac{1}{|\slashed{p}_{0}|},\frac{\slashed{p}_{0}}{|\slashed{p}_{0}|})
\end{align}
which maps the exterior of the disk of radius $B$ in $\mathbb{R}^{2}$ onto $(0,-1/B)\times S^{1}$. Using results analogous to (19.73), (19.76), substituted 
in (19.139), (19.140), we can show that the mapping $G_{t}\circ f^{-1}$ is one-to-one on $(0,-1/B)\times S^{1}$. Actually, what we 
have done here is just putting the ``radial variable'' outside the Taylor expansion, so we can proceed as in the case $|\slashed{p}_{0}|< B$. Still, we have to show that
the image of $(0,-1/B)\times S^{1}$ by $G_{t}\circ f^{-1}$ does not intersect the image of the disk of radius $B$ by $G_{t}$. In fact, although we have made a coordinate 
transformation, but the actual Taylor expansion is the same, so the conclusion mentioned above is obvious. 

    In the case of the $other$ null geodesics, $(F^{t})^{-1}$ is not differentiable at $t=t_{0}$, so we first investigate the mapping $F$. We choose a basis $\omega^{1}$,
$\omega^{2}$ for $T_{\vartheta_{0}}S^{2}$ which is orthonormal relative to $\slashed{g}^{-1}(q_{0})$ and express $\slashed{p}_{0}$, which belongs to the unit circle in
$T_{\vartheta_{0}}S^{2}$ in the form:
\begin{align}
 \slashed{p}_{0}=\omega^{1}\cos\varphi+\omega^{2}\sin\varphi
\end{align}
In the case that $q_{0}\in\mathcal{H}$, we have, from (19.88), (19.91), (19.94):
\begin{align}
 F^{t}(0;\varphi)=t_{0},\quad (\frac{\partial F^{t}}{\partial\tau})(0;\varphi)=0,\quad 
(\frac{\partial^{2}F^{t}}{\partial \tau^{2}})(0,\varphi)=\frac{1}{2}(\frac{\partial\mu}{\partial u})(q_{0})
\end{align}
 and:
\begin{align}
 \slashed{F}^{A}(0;\varphi)=\vartheta^{A}_{0},\quad (\frac{\partial\slashed{F}^{A}}{\partial\tau})(0;\varphi)=
(\frac{\partial^{2}\slashed{F}^{A}}{\partial\tau^{2}})(0;\varphi)=0\\\notag
(\frac{\partial^{3}\slashed{F}^{A}}{\partial\tau^{3}})(0;\varphi)=(\frac{\partial\mu}{\partial t}\frac{\partial\mu}{\partial u})(q_{0})
(\slashed{g}^{-1})^{AB}(q_{0})(\slashed{p}_{0})_{B}(\varphi);\\\notag
(\slashed{p}_{0})_{A}=(\omega^{1})_{A}\cos\varphi+(\omega^{2})_{A}\sin\varphi
\end{align}
It follows that in the case that $q_{0}\in\mathcal{H}$ we have:
\begin{align}
 F^{t}(\tau;\varphi)=t_{0}+\frac{1}{4}(\frac{\partial\mu}{\partial u})(q_{0})\tau^{2}+O(\tau^{3})\\\notag
\slashed{F}^{A}(\tau;\varphi)=\vartheta^{A}_{0}+\frac{1}{6}(\frac{\partial\mu}{\partial u}\frac{\partial\mu}{\partial t})(q_{0})
(\slashed{g}^{-1})^{AB}(q_{0})(\slashed{p}_{0})_{B}(\varphi)\tau^{3}+O(\tau^{4})
\end{align}
hence, since $\slashed{G}^{A}_{t}(\varphi)=\slashed{F}^{A}((F^{t})^{-1}(t;\varphi);\varphi)$,
\begin{align}
 \slashed{G}^{A}_{t}(\varphi)=\vartheta^{A}_{0}-\frac{4}{3}\frac{(\partial\mu/\partial t)(q_{0})}{\sqrt{-(\partial\mu/\partial u)(q_{0})}}
(\slashed{g}^{-1})^{AB}(q_{0})(\slashed{p}_{0})_{B}(\varphi)(t_{0}-t)^{3/2}+O((t_{0}-t)^{2})
\end{align}
(We just use the first of (19.152) to express $\tau$ in terms of $t-t_{0}$.)

Moreover, (19.150) implies:
\begin{align}
 (\frac{\partial F^{t}}{\partial\varphi})(0;\varphi)=(\frac{\partial}{\partial\tau}(\frac{\partial F^{t}}{\partial\varphi}))(0;\varphi)=
(\frac{\partial^{2}}{\partial\tau^{2}}(\frac{\partial F^{t}}{\partial\varphi}))(0;\varphi)=0
\end{align}
hence:
\begin{align}
 (\frac{\partial F^{t}}{\partial\varphi})(\tau;\varphi)=O(\tau^{3})
\end{align}
Also, (19.151) implies:
\begin{align}
 (\frac{\partial\slashed{F}^{A}}{\partial\varphi})(0;\varphi)=(\frac{\partial}{\partial\tau}(\frac{\partial\slashed{F}^{A}}{\partial\varphi}))(0;\varphi)
=(\frac{\partial^{2}}{\partial\tau^{2}}(\frac{\partial\slashed{F}^{A}}{\partial\varphi}))(0;\varphi)=0\\\notag
(\frac{\partial^{3}}{\partial\tau^{3}}(\frac{\partial\slashed{F}^{A}}{\partial\varphi}))(0;\varphi)=(\frac{\partial\mu}{\partial t}\frac{\partial\mu}{\partial u})
(q_{0})(\slashed{g}^{-1})^{AB}(q_{0})(-(\omega^{1})_{B}\sin\varphi+(\omega^{2})_{B}\cos\varphi)
\end{align}
hence:
\begin{align}
 (\frac{\partial\slashed{F}^{A}}{\partial\varphi})(\tau;\varphi)=\frac{1}{6}(\frac{\partial\mu}{\partial t}\frac{\partial\mu}{\partial u})(q_{0})
(\slashed{g}^{-1})^{AB}(q_{0})(-(\omega^{1})_{B}\sin\varphi+(\omega^{2})_{B}\cos\varphi)\tau^{3}+O(\tau^{4})
\end{align}
Differentiating the equation $F^{t}((F^{t})^{-1}(t;\varphi);\varphi)=t$ implicitly with respect to $\varphi$ yields:
\begin{align}
 \frac{\partial(F^{t})^{-1}}{\partial\varphi}=-\frac{\partial F^{t}/\partial\varphi}{\partial F^{t}/\partial\tau}
\end{align}
Now (19.150) implies:
\begin{align}
 (\frac{\partial F^{t}}{\partial\tau})(\tau;\varphi)=\frac{1}{2}(\frac{\partial\mu}{\partial u})(q_{0})\tau+O(\tau^{2})
\end{align}
Substituting (19.155) and (19.159) in (19.158) yields:
\begin{align}
 \frac{\partial(F^{t})^{-1}}{\partial\varphi}=O(\tau^{2})
\end{align}
Differentiating $\slashed{G}^{A}(t;\varphi)=\slashed{F}^{A}((F^{t})^{-1}(t;\varphi);\varphi)$ with respect to $\varphi$ we obtain:
\begin{align}
 \frac{\partial\slashed{G}^{A}_{t}}{\partial\varphi}=\frac{\partial\slashed{F}^{A}}{\partial\tau}\frac{\partial(F^{t})^{-1}}{\partial\varphi}
+\frac{\partial\slashed{F}^{A}}{\partial\varphi}
\end{align}
Now (19.151) implies:
\begin{align}
 \frac{\partial\slashed{F}^{A}}{\partial\tau}=\frac{1}{2}(\frac{\partial\mu}{\partial t}\frac{\partial\mu}{\partial u})
(q_{0})(\slashed{g}^{-1})^{AB}(q_{0})(\slashed{p}_{0})_{B}(\varphi)\tau^{2}+O(\tau^{3})
\end{align}
In view of (19.160) and (19.162) the first term on the right in (19.161) is $O(\tau^{4})$, while the second term is given by (19.157). We conclude that:
\begin{align}
 \frac{\partial\slashed{G}^{A}_{t}}{\partial\varphi}=a(\slashed{g}^{-1})^{AB}(q_{0})(-(\omega^{1})_{B}\sin\varphi+(\omega^{2})_{B}\cos\varphi)\tau^{3}
+O(\tau^{4})
\end{align}
where
\begin{align}
 a=\frac{1}{6}(\frac{\partial\mu}{\partial t}\frac{\partial\mu}{\partial u})(q_{0})
\end{align}
which gives:
\begin{align}
 \slashed{g}_{AB}(q_{0})\frac{\partial\slashed{G}^{A}_{t}}{\partial\varphi}\frac{\partial\slashed{G}^{B}_{t}}{\partial\varphi}\\\notag
=a^{2}(\slashed{g}^{-1})^{AB}(-(\omega^{1})_{A}\sin\varphi+(\omega^{2})_{A}\cos\varphi)(-(\omega^{1})_{B}\sin\varphi+(\omega^{2})_{B}\cos\varphi)\tau^{6}
+O(\tau^{7})\\\notag
=a^{2}\tau^{6}+O(\tau^{7})
\end{align}

     The above hold in the case that $q_{0}\in\mathcal{H}$. In the case $q_{0}\in\partial_{-}\mathcal{H}$ similar calculations using (19.88), (19.95), (19.99),
(19.104), (19.105), (19.108), give:
\begin{align}
 \slashed{G}^{A}_{t}(\varphi)=\vartheta^{A}_{0}-\frac{8}{5}\frac{(-(\partial\mu/\partial t)(q_{0}))^{3/4}}{((1/6)(\partial^{2}\mu/\partial u^{2})(q_{0}))^{1/4}}
(t_{0}-t)^{5/4}+O((t_{0}-t)^{3/2})\\
\frac{\partial\slashed{G}^{A}_{t}}{\partial\varphi}=b(\slashed{g}^{-1})^{AB}(q_{0})(-(\omega^{1})_{B}\sin\varphi+(\omega^{2})_{B}\cos\varphi)\tau^{5}
+O(\tau^{6})
\end{align}
where:
\begin{align}
 b=\frac{1}{120}((\frac{\partial\mu}{\partial t})^{2}\frac{\partial^{2}\mu}{\partial u^{2}})(q_{0})
\end{align}
and:
\begin{align}
 \slashed{g}_{AB}(q_{0})\frac{\partial\slashed{G}^{A}_{t}}{\partial\varphi}\frac{\partial\slashed{G}^{B}_{t}}{\partial\varphi}=b^{2}\tau^{10}+O(\tau^{11})
\end{align}
     The formulas (19.165) and (19.169) imply in the case of the other null geodesics ending at $q_{0}$ that the mapping (19.135) is of maximal rank at each 
$\slashed{p}_{0}\in S^{1}$ for $t$ suitably close to $t_{0}$. By the implicit function theorem we conclude that in the case of the other null geodesics the image of 
$S^{1}$ by the mapping (19.135) is, for $t$ suitably close to $t_{0}$, an immersed circle in the hypersurface $\{t\}\times(-\infty,\epsilon_{0})\times S^{2}\subset
\mathcal{M}^{e}_{\epsilon_{0}}$, which corresponds to the hyperplane $\Sigma_{t}$ in Galileo spacetime. Moreover, we deduce from (19.153) and (19.166) that the mapping
(19.135) is one-to-one for $t$ suitably close to $t_{0}$. It then follows that the image of $S^{1}$ by the mapping (19.135) is actually an embedded circle. This circle is 
the boundary of the embedded disk corresponding to the outgoing null geodesics.

     Similar arguments show that in the case of the $incoming$ null geodesics the image of $\mathbb{R}^{2}$ by the mapping (19.135) is an embedded disk for $t$ suitably 
close to $t_{0}$. The boundary of this disk is the circle corresponding to the $other$ null geodesics. We have thus proved the following theorem.

$\textbf{Theorem 19.1}$  Let $q_{0}=(t_{0},u_{0},\vartheta_{0})$ be a point belonging to $\mathcal{H}\bigcup\partial_{-}\mathcal{H}$, the singular boundary of the domain
of the maximal solution. Then the backward past null cone of $q_{0}$ in $T^{*}_{q_{0}}\mathcal{M}_{\epsilon_{0}}$ degenerates to the three pieces:

(i) The negative half-hyperplane: $p_{t}=0$, $p_{u}<0$.

(ii) The negative half hyperplane: $p_{u}=0$, $p_{t}<0$.

(iii) The plane $p_{t}=p_{u}=0$.

There is a corresponding trichotomy of the past null geodesic cone of $q_{0}$ into three classes of null geodesics ending at $q_{0}$. The $outgoing$ null geodesics,
the $incoming$ null geodesics, and the $other$ null geodesics, corresponding to (i), (ii), (iii), respectively. Let all these null geodesics be parametrized by $t$.
Then, in Galileo spacetime, the following hold:

    The tangent vectors to all $outgoing$ null geodesics ending at $q_{0}$ tend as we approach $q_{0}$ to the tangent vector $L(q_{0})$ of the generator of $C_{u_{0}}$ 
through $q_{0}$. The vector $L(q_{0})$ is thus an invariant null vector associated to the singular point $q_{0}$. 

    The tangent vector to an $incoming$ null geodesic ending at $q_{0}$ with angular momentum $\slashed{p}_{0}$ at $q_{0}$, tends as $t\longrightarrow t_{0}$ to the following
limit:
\begin{align*}
 L^{i}(q_{0})+\frac{\alpha^{-1}(q_{0})\hat{T}^{i}(q_{0})}{(1/2)(\alpha^{-2}(q_{0})+|\slashed{p}_{0}|^{2})}+
\frac{(\slashed{g}^{-1})^{AB}(q_{0})X^{i}_{A}(q_{0})(\slashed{p}_{0})_{B}}{(1/2)(\alpha^{-2}(q_{0})+|\slashed{p}_{0}|^{2})}
\end{align*}
which, for $\slashed{p}_{0}=0$ is equal to $\alpha^{2}(q_{0})\hat{\underline{L}}(q_{0})$, while as $\slashed{p}_{0}\longrightarrow \infty$ tends to $L(q_{0})$.

    The tangent vector to all the $other$ null geodesics ending at $q_{0}$ likewise tend as $t\longrightarrow t_{0}$ to $L(q_{0})$.

    Consider for each of the pieces (i), (ii), (iii), above, and for each $t\in[0,t_{0}]$, the mapping:
\begin{align*}
 \slashed{p}_{0}\mapsto G_{t}(\slashed{p}_{0})
\end{align*}
where $G_{t}(\slashed{p}_{0})$ is the point at parameter value $t$ along the null geodesic ending at parameter value $t_{0}$ at the point $q_{0}$ with momentum covector 
$\slashed{p}_{0}$.
Then for $t\slashed{=}t_{0}$ and suitably close to $t_{0}$ the following hold:

     The image of (i), which corresponds to the points of intersection of the outgoing null geodesics ending at $q_{0}$ with the hyperplane $\Sigma_{t}$, is an 
embedded disk.
     
    The image of (ii), which corresponds to the points of intersection of the incoming null geodesics ending at $q_{0}$ with the hyperplane $\Sigma_{t}$, is also an 
embedded disk.
    
    The image of (iii), which corresponds to the points of intersection of the other null geodesics ending at $q_{0}$ with the hyperplane $\Sigma_{t}$, is an embedded circle, 
the common boundary of the two disks.

\section{Transformation of Coordinates}
    Consider now again the mapping (19.129) for the class of $outgoing$ null geodesics ending at $q_{0}\in\mathcal{H}\bigcup\partial_{-}\mathcal{H}$. We now wish to make 
the dependence  on the point $q_{0}\in\mathcal{H}\bigcup\partial_{-}\mathcal{H}$ explicit, so we denote $F(\tau;\slashed{p}_{0})$ by $F(\tau;q_{0};\slashed{p}_{0})$.
What we wish to investigate presently is the change from one acoustical function $u$ to another $u'$ (arising from different initial data for the eikonal equation), such that the 
null geodesic generators of the level sets $C'_{u'}$ of $u'$, the characteristic hypersurfaces corresponding to $u'$, which have future end points on the singular boundary 
of the domain of the maximal solution, belong to the class of $outgoing$ null geodesics ending at these points. We can view this change as generated by a coordinate 
transformation:
\begin{align}
 u_{0}=u_{0}(u'_{0},\vartheta'_{0}),\quad \vartheta_{0}=\vartheta_{0}(\vartheta'_{0})
\end{align}
on $\mathcal{H}\bigcup\partial_{-}\mathcal{H}$, preserving the invariant curves. We also assume:
\begin{align}
 \frac{\partial u_{0}}{\partial u'_{0}}>0,\quad \det\{\frac{\partial\vartheta_{0}}{\partial\vartheta'_{0}}\}>0
\end{align}
If the coordinates of a point $q_{0}\in \mathcal{H}\bigcup\partial_{-}\mathcal{H}$ in the unprimed system are: $(t_{0},u_{0},\vartheta_{0})$, $t_{0}=t_{*}
(u_{0},\vartheta_{0})$, then the coordinates of the same point in the primed system are: $(t'_{0},u'_{0},\vartheta'_{0})$, $t'_{0}=t'_{*}(u'_{0},\vartheta'_{0})$,
where:
\begin{align}
 t'_{*}(u'_{0},\vartheta'_{0})=t_{*}(u_{0},\vartheta_{0})
\end{align}
hence $t'_{0}=t_{0}$.

     We shall determine $\slashed{p}_{0}$ as a function of $(u'_{0},\vartheta'_{0})$, such that the hypersurfaces, given in parametric form by:
\begin{align}
 t=M^{t}(\tau,u'_{0},\vartheta'_{0})=F^{t}(\tau;(u'_{0},\vartheta'_{0});\slashed{p}_{0}(u'_{0},\vartheta'_{0}))\\\notag
u=M^{u}(\tau,u'_{0},\vartheta'_{0})=F^{u}(\tau;(u'_{0},\vartheta'_{0});\slashed{p}_{0}(u'_{0},\vartheta'_{0}))\\\notag
\vartheta=\slashed{M}(\tau,u'_{0},\vartheta'_{0})=\slashed{F}(\tau;(u'_{0},\vartheta'_{0});\slashed{p}_{0}(u'_{0},\vartheta'_{0}))
\end{align}
for constant values of $u'_{0}$, are null hypersurfaces with respect to the acoustical metric, namely the level sets $C'_{u'_{0}}$ of the new acoustical function $u'$.

     In (19.173), $(t'_{*}(u'_{0},\vartheta'_{0}),u'_{0},\vartheta'_{0})$ are the coordinates of a point $q_{0}\in\mathcal{H}\bigcup\partial_{-}\mathcal{H}$ in the primed
system, hence we have:
\begin{align}
 F^{t}(0;(u'_{0},\vartheta'_{0});\slashed{p}_{0})=t_{*}(u_{0},\vartheta_{0})\\\notag
F^{u}(0;(u'_{0},\vartheta'_{0});\slashed{p}_{0})=u_{0}\\\notag
\slashed{F}(0;(u'_{0},\vartheta'_{0});\slashed{p}_{0})=\vartheta_{0}
\end{align}

The hypersurfaces in question are generated by the outgoing null geodesics corresponding to constant values of $(u'_{0},\vartheta'_{0})$, and the hypersurfaces are null
hypersurfaces if and only if these null geodesics are orthogonal, with respect to the acoustical metric, to the sections $S'_{\tau,u'_{0}}$ corresponding to conatant
values of $\tau$. Let then $X^{\prime}_{A}$ be the tangent vectorfield to the $\vartheta'^{A}_{0}$ coordinate lines on these sections and, as before, $\tilde{L}'$ be the tangent 
vectorfield to the null geodesics, as parametrized by $\tau$. Then the orthogonality condition reads:
\begin{align}
 g(\tilde{L}',X'_{A})=0\quad:\quad A=1,2
\end{align}
and we have:
\begin{align}
 \tilde{L}'=\frac{dF^{t}}{d\tau}\frac{\partial}{\partial t}+\frac{dF^{u}}{d\tau}\frac{\partial}{\partial u}+\frac{d\slashed{F}^{A}}{d\tau}\frac{\partial}
{\partial\vartheta^{A}}
\end{align}
and:
\begin{align}
 X'_{A}=X'^{t}_{A}\frac{\partial}{\partial t}+X'^{u}_{A}\frac{\partial}{\partial u}+X'^{B}_{A}\frac{\partial}{\partial\vartheta^{B}}
\end{align}
where:
\begin{align}
 X'^{t}_{A}=\frac{\partial F^{t}}{\partial\vartheta'^{A}_{0}}+\frac{\partial F^{t}}{\partial(\slashed{p}_{0})_{B}}\frac{\partial(\slashed{p}_{0})_{B}}
{\partial\vartheta'^{A}_{0}}\\\notag
X'^{u}_{A}=\frac{\partial F^{u}}{\partial\vartheta'^{A}_{0}}+\frac{\partial F^{u}}{\partial(\slashed{p}_{0})_{B}}\frac{\partial(\slashed{p}_{0})_{B}}
{\partial\vartheta'^{A}_{0}}\\\notag
X'^{C}_{A}=\frac{\partial\slashed{F}^{C}}{\partial\vartheta'^{A}_{0}}+\frac{\partial\slashed{F}^{C}}{\partial(\slashed{p}_{0})_{B}}\frac
{\partial(\slashed{p}_{0})_{B}}{\partial\vartheta'^{A}_{0}}
\end{align}
In fact, (19.173) defines a change of coordinates in the spacetime domain covered by the hypersurfaces in question, from $(t,u,\vartheta)$ coordinates to $(\tau,
u'_{0},\vartheta'_{0})$ coordinates, and in the new coordinates we have, simply:
\begin{align}
 \tilde{L}'=\frac{\partial}{\partial\tau},\quad X'_{A}=\frac{\partial}{\partial\vartheta'^{A}_{0}}
\end{align}
It follows that:
\begin{align}
 [\tilde{L}',X'_{A}]=0
\end{align}
Now by (19.60) the integral curves of the vectorfield $\mu^{-1}\tilde{L}'$ are affinely parametrized null geodesics, hence:
\begin{align}
 D_{\tilde{L}'}\tilde{L}'=\mu^{-1}(\tilde{L}'\mu)\tilde{L}',\quad g(\tilde{L}',\tilde{L}')=0
\end{align}
Let us then define:
\begin{align}
 \iota_{A}=g(\tilde{L}',X'_{A})
\end{align}
What we shall do is to prove $\iota_{A}=0$.

We have, along the integral curves of $\tilde{L}'$:
\begin{align*}
 \frac{d\iota_{A}}{d\tau}=\tilde{L}'(g(\tilde{L}',X'_{A}))=g(D_{\tilde{L}'}\tilde{L'},X'_{A})+g(\tilde{L}',D_{\tilde{L}'}X'_{A})
\end{align*}
and from (19.180): $D_{\tilde{L}'}X'_{A}=D_{X'_{A}}\tilde{L}'$, hence:
\begin{align*}
 g(\tilde{L}',D_{\tilde{L}'}X'_{A})=\frac{1}{2}X'_{A}(g(\tilde{L}',\tilde{L}'))=0
\end{align*}
while by (19.181):
\begin{align*}
 g(D_{\tilde{L}'}\tilde{L}',X'_{A})=\mu^{-1}\frac{d\mu}{d\tau}\iota_{A}
\end{align*}
It follows that:
\begin{align}
 \frac{d}{d\tau}(\mu^{-1}\iota_{A})=0
\end{align}
In view of (19.15) and (19.176),
\begin{align}
 \iota_{A}=-\mu\frac{dF^{t}}{d\tau}X'^{u}_{A}-\mu\frac{dF^{u}}{d\tau}X'^{t}_{A}+\alpha^{-2}\mu^{2}\frac{dF^{u}}{d\tau}X'^{u}_{A}\\\notag
+\slashed{g}_{BC}(\frac{d\slashed{F}^{B}}{d\tau}+\mu\hat{\Xi}^{B}\frac{dF^{u}}{d\tau})(X'^{C}_{A}+\mu\hat{\Xi}^{C}X'^{u}_{A})
\end{align}
At $\tau=0$ this reduces to:
\begin{align}
 (\iota_{A})_{\tau=0}=(\slashed{g}_{BC}\frac{d\slashed{F}^{B}}{d\tau}X'^{C}_{A})_{\tau=0}=0
\end{align}
by (19.73). Moreover, we have:
\begin{align}
 \lim_{\tau\rightarrow0}(\mu^{-1}\iota_{A})=-(\frac{dF^{t}}{d\tau}X'^{u}_{A}+\frac{dF^{u}}{d\tau}X'^{t}_{A})_{\tau=0}
+(\slashed{g}_{BC}X'^{C}_{A})_{\tau=0}\lim_{\tau\rightarrow0}(\mu^{-1}\frac{d\slashed{F}^{B}}{d\tau})
\end{align}
Now by (19.73):
\begin{align}
 (\frac{dF^{t}}{d\tau})_{\tau=0}=1,\quad (\frac{dF^{u}}{d\tau})_{\tau=0}=0
\end{align}
and by (19.75), (19.77):
\begin{align}
 \lim_{\tau\rightarrow0}(\mu^{-1}\frac{d\slashed{F}^{A}}{d\tau})=(\slashed{g}^{-1})^{AB}(q_{0})(\slashed{p}_{0})_{B}
\end{align}
while by (19.174) and (19.178):
\begin{align}
 (X_{A}^{\prime u})_{\tau=0}=\frac{\partial u_{0}}{\partial\vartheta_{0}^{\prime A}},\quad (X_{A}^{\prime B})_{\tau=0}=\frac{\partial\vartheta^{B}_{0}}
{\partial\vartheta_{0}^{\prime A}}
\end{align}
Note that we view $F^{u}$ and $\slashed{F}$ as the functions of $(u'_{0},\vartheta'_{0})$.

Substituting in (19.186) we obtain:
\begin{align}
 \lim_{\tau\rightarrow0}(\mu^{-1}\iota_{A})=-\frac{\partial u_{0}}{\partial\vartheta_{0}^{\prime A}}+(\slashed{p}_{0})_{B}\frac{\partial\vartheta^{B}_{0}}
{\partial\vartheta_{0}^{\prime A}}
\end{align}
Therefore the condition $\lim_{\tau\rightarrow0}(\mu^{-1}\iota_{A})=0$ is equivalent to the following definition of $\slashed{p}_{0}$ as a function of 
$(u'_{0},\vartheta'_{0})$:
\begin{align}
 (\slashed{p}_{0})_{A}=\frac{\partial u_{0}}{\partial\vartheta'^{B}_{0}}\frac{\partial\vartheta'^{B}_{0}}{\partial\vartheta^{A}_{0}}
\end{align}
Once $\slashed{p}_{0}$ is defined according to (19.191), equation (19.183) implies that $\iota_{A}=0$ for all $\tau$, that is, the orthogonality condition (19.175)
is everywhere satisfied and the hypersurfaces of constant $u'_{0}$, given by (19.173) are outgoing characteristic hypersurfaces.

     In the following we denote $(u'_{0},\vartheta'_{0})$ simply by $(u',\vartheta')$. At fixed $q_{0}$, $\slashed{p}_{0}$ we invert the mapping (19.130), to obtain the
mapping (19.133). Let us denote:
\begin{align}
 f(t',u',\vartheta')=(F^{t})^{-1}(t';(u',\vartheta');\slashed{p}_{0}(u',\vartheta'))
\end{align}
We then substitute in (19.173) to obtain the following transformation from the original acoustical coordinates $(t,u,\vartheta)$ to the new acoustical coordinates 
$(t',u',\vartheta')$:
\begin{align}
 t=t'\\\notag
u=N^{u}(t',u',\vartheta')=F^{u}(f(t',u',\vartheta');(u',\vartheta');\slashed{p}_{0}(u',\vartheta'))\\\notag
\vartheta=\slashed{N}(t',u',\vartheta')=\slashed{F}(f(t',u',\vartheta');(u',\vartheta');\slashed{p}_{0}(u',\vartheta'))
\end{align}
The acoustical metric has in the new acoustical coordinates a form similar to (19.15):
\begin{align}
 g=-2\mu'dt'du'+\alpha'^{-2}\mu'^{2}du'^{2}+\slashed{g}'_{AB}(d\vartheta'^{A}+\Xi'^{A}du')(d\vartheta'^{B}+\Xi'^{B}du')
\end{align}
where we have the new coefficients $\alpha'$, $\mu'$, $\Xi'^{A}$, and $\slashed{g}'_{AB}$. We have, in arbitrary local coordinates:
\begin{align}
 \alpha'^{-2}=-(g^{-1})^{\mu\nu}\partial_{\mu}t'\partial_{\nu}t'=-(g^{-1})^{\mu\nu}\partial_{\mu}t\partial_{\nu}t=\alpha^{-2}
\end{align}
so the functions $\alpha$ and $\alpha'$ coincide. The function $\mu'$ is given by (see Chapter 2):
\begin{align}
 L'=\mu'\hat{L}',\quad\hat{L}'^{\mu}=-(g^{-1})^{\mu\nu}\partial_{\nu}u'
\end{align}
and we have (in arbitrary local coordinates):
\begin{align}
 \frac{1}{\mu'}=-(g^{-1})^{\mu\nu}\partial_{\mu}t'\partial_{\nu}u'=-(g^{-1})^{\mu\nu}\partial_{\mu}t\partial_{\nu}u'
\end{align}
Here $L'$, related to $\tilde{L}'$ by (19.112), has the same integral curves as $\tilde{L}'$, but parametrized by $t$ instead of $\tau$, while the integral curves of 
$\hat{L}'$ are also the same but affinely parametrized. As we have shown above, the integral curves of the vectorfield $\mu^{-1}\tilde{L}'$ are likewise the same and 
affinely parametrized. Now two different affine parameters $s$ and $s'$ along the same geodesic, giving the same orientation, satisfy $s'=as+b$ where $a$ and $b$ are 
constants along the geodesic, $a>0$. It follows that there is a function $a(u',\vartheta')>0$ such that:
\begin{align}
 \mu^{-1}\tilde{L}'=a\hat{L}'
\end{align}
On the other hand, from (19.37) and the second of (19.196) we obtain, working in the original acoustical coordinates:
\begin{align}
 \lim_{\tau\rightarrow0}(\mu\hat{L}')=(\frac{\partial u'}{\partial u})_{\tau=0}\frac{\partial}{\partial t}+(\frac{\partial u'}{\partial t})_{\tau=0}
\frac{\partial}{\partial u}
\end{align}
while by the first of (19.73),
\begin{align}
 \lim_{\tau\rightarrow0}\tilde{L}'=\frac{\partial}{\partial t}
\end{align}
Now, from (19.193) we obtain:
\begin{align}
 \frac{\partial t}{\partial t'}=1,\quad\frac{\partial t}{\partial u'}=0,\quad\frac{\partial t}{\partial\vartheta'^{A}}=0\\\notag
\frac{\partial u}{\partial t'}=\frac{\partial F^{u}}{\partial\tau}\frac{\partial f}{\partial t'}, \quad\frac{\partial u}{\partial u'}=\frac{\partial F^{u}}{\partial\tau}
\frac{\partial f}{\partial u'}+\frac{\partial F^{u}}{\partial u'},\quad\frac{\partial u}{\partial\vartheta'^{A}}=\frac{\partial F^{u}}{\partial\tau}
\frac{\partial f}{\partial\vartheta'^{A}}+\frac{\partial F^{u}}{\partial\vartheta'^{A}}\\\notag
\frac{\partial\vartheta^{A}}{\partial t'}=\frac{\partial\slashed{F}^{A}}{\partial\tau}\frac{\partial f}{\partial t'},\quad \frac{\partial\vartheta^{A}}{\partial u'}
=\frac{\partial\slashed{F}^{A}}{\partial\tau}\frac{\partial f}{\partial u'}+\frac{\partial\slashed{F}^{A}}{\partial u'},\quad\frac{\partial\vartheta^{A}}{\partial\vartheta'^{B}}
=\frac{\partial\slashed{F}^{A}}{\partial\tau}\frac{\partial f}{\partial\vartheta'^{B}}+\frac{\partial\slashed{F}^{A}}{\partial\vartheta'^{B}}\\\notag
\end{align}
In the above, we are considering $F^{t}$, $F^{u}$ and $\slashed{F}$ as functions of $\tau$, $u'$, $\vartheta'$, according to:
\begin{align}
 F^{t}(\tau,u',\vartheta')=F^{t}(\tau;(u',\vartheta');\slashed{p}_{0}(u',\vartheta'))\\\notag
F^{u}(\tau,u',\vartheta')=F^{u}(\tau;(u',\vartheta');\slashed{p}_{0}(u',\vartheta'))\\\notag
\slashed{F}(\tau,u',\vartheta')=\slashed{F}(\tau;(u',\vartheta');\slashed{p}_{0}(u',\vartheta'))\\\notag
\end{align}
 At $\tau=0$, which corresponds to $t'=t'_{*}(u',\vartheta')$, so that:
\begin{align}
 f(t'_{*}(u',\vartheta'), u', \vartheta')=0
\end{align}
we have, from (19.73):
\begin{align}
 \frac{\partial F^{t}}{\partial\tau}=1,\quad\frac{\partial F^{u}}{\partial\tau}=0,\quad\frac{\partial\slashed{F}^{A}}{\partial\tau}=0
\end{align}
hence the second and the third of (19.201) reduce at $\tau=0$ to:
\begin{align}
 \frac{\partial u}{\partial t'}=0,\quad \frac{\partial u}{\partial u'}=\frac{\partial F^{u}}{\partial u'},\quad 
\frac{\partial u}{\partial\vartheta'^{A}}=\frac{\partial F^{u}}{\partial\vartheta'^{A}}\\\notag
\frac{\partial\vartheta^{A}}{\partial t'}=0,\quad \frac{\partial\vartheta^{A}}{\partial u'}=\frac{\partial\slashed{F}^{A}}{\partial u'},\quad 
\frac{\partial\vartheta^{A}}{\partial\vartheta'^{B}}=\frac{\partial\slashed{F}^{A}}{\partial\vartheta'^{B}}\\\notag
\end{align}
Now at $\tau=0$ we have, according to (19.170) and (19.174):
\begin{align}
 F^{t}(0,u',\vartheta')=t'_{*}(u',\vartheta')=t_{*}(u_{0}(u',\vartheta'),\vartheta_{0}(\vartheta'))\\\notag
F^{u}(0,u',\vartheta')=u_{0}(u',\vartheta')\\\notag
\slashed{F}(0,u',\vartheta')=\vartheta_{0}(\vartheta')
\end{align}
Hence, at $\tau=0$:
\begin{align}
 \frac{\partial F^{u}}{\partial u'}=\frac{\partial u_{0}}{\partial u'},\quad \frac{\partial F^{u}}{\partial\vartheta'^{A}}
=\frac{\partial u_{0}}{\partial\vartheta'^{A}}\\\notag
\frac{\partial \slashed{F}^{A}}{\partial u'}=0,\quad\frac{\partial\slashed{F}^{A}}{\partial\vartheta'^{B}}=\frac{\partial\vartheta^{A}_{0}}{\partial\vartheta'^{B}}
\end{align}
Substituting in (19.205) we conclude that the Jacobian matrix of the transformation (19.193) is at $\tau=0$ given by:
\begin{align}
 \frac{\partial(t,u,\vartheta)}{\partial(t',u',\vartheta')}=
\begin{bmatrix}
                                                             1&0&0&0\\
0&\partial u_{0}/\partial u'&\partial u_{0}/\partial\vartheta'^{1}&\partial u_{0}/\partial\vartheta'^{2}\\
0&0&\partial\vartheta^{1}_{0}/\partial\vartheta'^{1}&\partial\vartheta^{1}_{0}/\partial\vartheta'^{2}\\
0&0&\partial\vartheta^{2}_{0}/\partial\vartheta'^{1}&\partial\vartheta^{2}_{0}/\partial\vartheta'^{2}
                                                            \end{bmatrix}
\end{align}
Computing the reciprocal matrix we then find:
\begin{align}
 (\frac{\partial u'}{\partial u})_{\tau=0}=\frac{1}{\partial u_{0}/\partial u'},\quad(\frac{\partial u'}{\partial t})_{\tau=0}=0
\end{align}

Substituting in (19.199) and comparing with (19.198) and (19.200) we conclude that:
\begin{align}
 a(u',\vartheta')=\frac{\partial u_{0}}{\partial u'}>0
\end{align}
by (19.171). And from (19.37), (19.112) and (19.198)
\begin{align}
 \frac{\mu'}{\mu}=\frac{a}{dF^{t}/d\tau}
\end{align}
Then we know that $\mu'/\mu$ is bounded from above and below by positive constants. Also, since $\alpha'=\alpha$ we have:
\begin{align}
 \frac{\kappa'}{\kappa}=\frac{\mu'}{\mu}
\end{align}
     Next we consider the new coefficients $\Xi'^{A}$ and $\slashed{g}'_{AB}$. We have, in view of (19.194),
\begin{align}
 \Xi'_{A}=\slashed{g}'_{AB}\Xi^{\prime B}=g'_{Au'}=g_{\mu\nu}\frac{\partial x^{\mu}}{\partial \vartheta'^{A}}\frac{\partial x^{\nu}}{\partial u'}\\\notag
=\mu^{2}(\alpha^{-2}+|\hat{\Xi}|^{2})\frac{\partial u}{\partial\vartheta'^{A}}\frac{\partial u}{\partial u'}+\mu\hat{\Xi}_{B}
(\frac{\partial u}{\partial\vartheta'^{A}}\frac{\partial\vartheta^{B}}{\partial u'}+\frac{\partial\vartheta^{B}}{\partial\vartheta'^{A}}\frac{\partial u}{\partial u'})
+\slashed{g}_{BC}\frac{\partial\vartheta^{B}}{\partial\vartheta'^{A}}\frac{\partial\vartheta^{C}}{\partial u'}
\end{align}
Here,
\begin{align*}
 |\hat{\Xi}|^{2}=\slashed{g}_{AB}\hat{\Xi}^{A}\hat{\Xi}^{B},\quad\hat{\Xi}_{A}=\slashed{g}_{AB}\hat{\Xi}^{B}
\end{align*}
From (19.208) we have:
\begin{align*}
 (\frac{\partial\vartheta^{A}}{\partial u'})_{\tau=0}=0
\end{align*}
In view of (19.213) this implies that $\Xi'$ vanishes at $\tau=0$, that is, on the singular boundary $\mathcal{H}\bigcup\partial_{-}\mathcal{H}$, a reflection of the fact 
that the new acoustical coordinates $(t',u',\vartheta')$ are also canonical. Since, by (19.211) $\partial\mu'/\partial t'<0$ at $\tau=0$, it follows as in (19.35) that
\begin{align}
 \hat{\Xi}'=\mu^{\prime -1}\Xi'
\end{align}
is regular at $\tau=0$. ($a$ does not depend on $t'$.)

     We have:
\begin{align}
 \slashed{g}'_{AB}=g_{\mu\nu}\frac{\partial x^{\mu}}{\partial \vartheta'^{A}}\frac{\partial x^{\nu}}{\partial\vartheta'^{B}}\\\notag
=\mu^{2}(\alpha^{-2}+|\hat{\Xi}|^{2})\frac{\partial u}{\partial\vartheta'^{A}}\frac{\partial u}{\partial\vartheta'^{B}}
+\mu\hat{\Xi}_{C}(\frac{\partial\vartheta^{C}}{\partial\vartheta'^{B}}\frac{\partial u}{\partial\vartheta'^{A}}+\frac{\partial\vartheta^{C}}{\partial\vartheta'^{A}}
\frac{\partial u}{\partial\vartheta'^{B}})+\slashed{g}_{CD}\frac{\partial\vartheta^{C}}{\partial\vartheta'^{A}}\frac{\partial\vartheta^{D}}{\partial\vartheta'^{B}}
\end{align}
Let:
\begin{align}
 Y'=Y'^{A}\frac{\partial}{\partial\vartheta'^{A}},\quad \slashed{g}'(Y',Y')=\slashed{g}'_{AB}Y'^{A}Y'^{B}
\end{align}
Let also:
\begin{align}
 \tilde{\slashed{g}}_{AB}=\slashed{g}_{CD}\frac{\partial\vartheta^{C}}{\partial\vartheta'^{A}}\frac{\partial\vartheta^{D}}{\partial\vartheta'^{B}},\quad
\tilde{\slashed{g}}(Y',Y')=\tilde{\slashed{g}}_{AB}Y'^{A}Y'^{B}
\end{align}
We then deduce from (19.215) the following lower bound for $\slashed{g}'$ in terms of $\slashed{g}$: For any vector $Y'$ as in (19.216) it holds that 
\begin{align}
 \slashed{g}^{\prime}(Y^{\prime},Y^{\prime})\geq \frac{\tilde{\slashed{g}}(Y^{\prime},Y^{\prime})}{1+\alpha^{2}\slashed{g}(\hat{\Xi},\hat{\Xi})}
\end{align}
In fact, by (19.215),
\begin{align*}
 \slashed{g}^{\prime}(Y^{\prime},Y^{\prime})=
\mu^{2}(\alpha^{-2}+|\hat{\Xi}|^{2})(Y^{\prime}u)^{2}+2\mu(\hat{\Xi}^{\prime}\cdot Y^{\prime})(Y^{\prime}u)+\tilde{\slashed{g}}(Y^{\prime},Y^{\prime})
\end{align*}
where
\begin{align*}
 \hat{\Xi}^{\prime}_{B}=\hat{\Xi}_{C}\frac{\partial\vartheta^{C}}{\partial\vartheta^{\prime}_{B}}
\end{align*}
then by Cauchy-Schwarz,
\begin{align*}
 |\hat{\Xi}^{\prime}\cdot Y^{\prime}|\leq \sqrt{\tilde{\slashed{g}}^{-1}(\hat{\Xi}^{\prime},\hat{\Xi}^{\prime})}\sqrt{\tilde{\slashed{g}}(Y^{\prime},Y^{\prime})}
\end{align*}
Since
\begin{align*}
 \tilde{\slashed{g}}^{-1}(\hat{\Xi}^{\prime},\hat{\Xi}^{\prime})=\slashed{g}^{-1}(\hat{\Xi},\hat{\Xi})=|\hat{\Xi}|^{2}
\end{align*}
we obtain:
\begin{align*}
 |\hat{\Xi}^{\prime}\cdot Y^{\prime}|\leq |\hat{\Xi}|\sqrt{\tilde{\slashed{g}}(Y^{\prime},Y^{\prime})}
\end{align*}
Setting:
\begin{align*}
 a=\mu^{2}(\alpha^{-2}+|\hat{\Xi}|^{2}),\quad b=\mu|\hat{\Xi}|,\quad x=|Y^{\prime}u|,\quad y=\sqrt{\tilde{\slashed{g}}(Y^{\prime},Y^{\prime})}
\end{align*}
then
\begin{align*}
 \slashed{g}^{\prime}(Y^{\prime},Y^{\prime})\geq ax^{2}-2bxy+y^{2}
\end{align*}
Setting:
\begin{align*}
 \rho=\frac{x}{y}
\end{align*}
we have:
\begin{align*}
 \frac{\slashed{g}^{\prime}(Y^{\prime},Y^{\prime})}{\tilde{\slashed{g}}(Y^{\prime},Y^{\prime})}
\geq a\rho^{2}-2b\rho+1\geq 1-\frac{b^{2}}{a}=\frac{\mu^{2}\alpha^{-2}}{\mu^{2}(\alpha^{-2}+|\hat{\Xi}|^{2})}
\end{align*}
Reversing the roles of the coordinates $(t,u,\vartheta)$ and $(t^{\prime},u^{\prime},\vartheta^{\prime})$, we obtain in a similar manner a lower bound for
$\slashed{g}$ in terms of $\slashed{g}^{\prime}$. Let:
\begin{align}
 Y=Y^{A}\frac{\partial}{\partial\vartheta^{A}},\quad \slashed{g}(Y,Y)=\slashed{g}_{AB}Y^{A}Y^{B}
\end{align}
Let also:
\begin{align}
 \tilde{\slashed{g}}^{\prime}_{AB}=\slashed{g}^{\prime}_{CD}
 \frac{\partial\vartheta^{\prime C}}{\partial\vartheta^{A}}
\frac{\partial\vartheta^{\prime D}}{\partial\vartheta^{B}},\quad \tilde{\slashed{g}}^{\prime}(Y,Y)=\tilde{\slashed{g}}^{\prime}_{AB}Y^{A}Y^{B}
\end{align}
Then for any $Y$ as in (19.219) it holds that
\begin{align}
 \slashed{g}(Y,Y)\geq \frac{\tilde{\slashed{g}}^{\prime}(Y,Y)}{1+\alpha^{2}\slashed{g}^{\prime}(\hat{\Xi}^{\prime},\hat{\Xi}^{\prime})}
\end{align}

     We summarize the above results in the following proposition.

$\textbf{Proposition 19.2}$  Let $u'$ be any other acoustical function such that the null geodesic generators of the level surfaces $C'_{u'}$ of $u'$ which have future
end points on the singular boundary of the domain of the maximal solution belong to the class of $outgoing$ null geodesics ending at this points. The change from the original
acoustical function $u$ to the new acoustical function $u'$ can be viewed as generated by a coordinate transformation
\begin{align}
 u_{0}=u_{0}(u'_{0},\vartheta'_{0}),\quad\vartheta_{0}=\vartheta_{0}(\vartheta'_{0})
\end{align}
on $\mathcal{H}\bigcup\partial_{-}\mathcal{H}$ preserving the invariant curves, so that both systems of acoustical coordinates are canonical. We require that the 
transformation $\vartheta_{0}\mapsto\vartheta'_{0}$ be orientation preserving and that on each invariant curve the transformation $u_{0}\mapsto u'_{0}$ be also 
orientation preserving. The characteristic hypersurfaces $C'_{u'_{0}}$ are given in parametric form in terms of the original acoustical coordinates $(t,u,\vartheta)$ by:
\begin{align*}
 t=F^{t}(\tau,u'_{0},\vartheta'_{0})\\
u=F^{u}(\tau,u'_{0},\vartheta'_{0})\\
\vartheta=\slashed{F}(\tau,u'_{0},\vartheta'_{0})
\end{align*}
where:
\begin{align*}
 F^{t}(\tau,u'_{0},\vartheta'_{0})=F^{t}(\tau;q_{0};\slashed{p}_{0}(q_{0}))\\
F^{u}(\tau,u'_{0},\vartheta'_{0})=F^{u}(\tau;q_{0};\slashed{p}_{0}(q_{0}))\\
\slashed{F}(\tau,u'_{0},\vartheta'_{0})=\slashed{F}(\tau;q_{0};\slashed{p}_{0}(q_{0}))
\end{align*}
 Here $q_{0}$ is the point on $\mathcal{H}\bigcup\partial_{-}\mathcal{H}$ with coordinates $(t_{0},u_{0},\vartheta_{0})$ in the original system and 
$(t'_{0},u'_{0},\vartheta'_{0})$ in the new system, where
\begin{align*}
 t_{0}=t_{*}(u_{0},\vartheta_{0})=t'_{*}(u'_{0},\vartheta'_{0})=t'_{0}
\end{align*}
Also,
\begin{align*}
 \slashed{p}_{0}(q_{0})=\frac{\partial u_{0}}{\partial\vartheta_{0}^{\prime A}}d\vartheta^{\prime A}_{0}
\end{align*}
is the angular momentum at $q_{0}$. In the above
\begin{align*}
 \tau\mapsto(F^{t}(\tau;q_{0};\slashed{p}_{0}),F^{u}(\tau;q_{0};\slashed{p}_{0}),\slashed{F}(\tau;q_{0};\slashed{p}_{0}))
\end{align*}
is the $outgoing$ null geodesic, parametrized by $\tau$, ending at $\tau=0$ at $q_{0}\in\mathcal{H}\bigcup\partial_{-}\mathcal{H}$ with angular momentum $\slashed{p}_{0}$.
Denoting by $\tau=f(t',u',\vartheta')$ the solution of 
\begin{align*}
 F^{t}(\tau,u',\vartheta')=t'
\end{align*}
for fixed $(u',\vartheta')$, the transformation between the original system of acoustical coordinates $(t,u,\vartheta)$ and the new system of acoustical coordinates 
$(t',u',\vartheta')$ is given by:
\begin{align*}
 t=t'\\
u=N^{u}(t',u',\vartheta')\\
\vartheta=\slashed{N}(t',u',\vartheta')
\end{align*}
where:
\begin{align*}
 N^{u}(t',u',\vartheta')=F^{u}(f(t',u',\vartheta'),u',\vartheta')\\
\slashed{N}(t',u',\vartheta')=\slashed{F}(f(t',u',\vartheta'),u',\vartheta')
\end{align*}
The acoustical metric has the same form in the new acoustical coordinates:
\begin{align*}
 g=-2\mu dtdu+\alpha^{-2}\mu^{2}du^{2}+\slashed{g}_{AB}(d\vartheta^{A}+\mu\hat{\Xi}^{A}du)(d\vartheta^{B}+\mu\hat{\Xi}^{B}du)\\
g=-2\mu' dt'du'+\alpha^{\prime -2}\mu^{\prime 2}du^{\prime 2}+\slashed{g}^{\prime}_{AB}
(d\vartheta^{\prime A}+\mu\hat{\Xi'}^{A}du')(d\vartheta^{\prime B}+\mu\hat{\Xi'}^{B}du')
\end{align*}
where $\alpha=\alpha'$. The ratio $\mu'/\mu$ is bounded above and below by positive constants, while the metric $\slashed{g}'$ and $\slashed{g}$ are uniformly equivalent in
a neighborhood of the singular boundary. \vspace{7mm}

    We now revisit the mapping (19.135) for a fixed point $q_{0}\in\mathcal{H}\bigcup\partial_{-}\mathcal{H}$ and for the class of $outgoing$ null geodesics ending at $q_{0}$.
Setting $t=0$ we have the mapping of $\mathbb{R}^{2}$ into $\Sigma_{0}$ by:
\begin{align}
 \slashed{p}_{0}\mapsto(0,u,\vartheta)=G_{0}(\slashed{p}_{0})\\\notag
G_{0}(\slashed{p}_{0})=G(0,\slashed{p}_{0})=(0,F^{u}(f(\slashed{p}_{0});\slashed{p}_{0}), \slashed{F}(f(\slashed{p}_{0});\slashed{p}_{0}))\\\notag
\textrm{where}\quad \tau=f(\slashed{p}_{0})\quad\textrm{is the solution of}\quad F^{t}(\tau;\slashed{p}_{0})=0
\end{align}
Let $u^{\prime}$ be a function defined on $\Sigma_{0}$, which has a unique maximum
$p_{0}$ on $\Sigma_{0}$, and is smooth and without critical points on $\Sigma_{0}\setminus\{p_{0}\}$. We may consider $u^{\prime}$ as a function of the acoustical coordinates $(0,u,\vartheta)$ on $\Sigma_{0}$, defined on $(-\infty,1)\times S^{2}\setminus\{(u_{0},\vartheta_{0})\}$ where $(u_{0},\vartheta_{0})$ are the coordinates of the maximum point $p_{0}$. 

Consider then the following function on $\mathbb{R}^{2}$:
\begin{align}
 g(\slashed{p}_{0})=u'(G_{0}(\slashed{p}_{0}))
\end{align}
Consider now a maximum of the function $g(\slashed{p}_{0})$ which is achieved at some point $\slashed{p}_{0,M}\in\mathbb{R}^{2}$. Such a maximum point corresponds to 
a point $p_{M}\in\Sigma_{0}$, where:
\begin{align}
 p_{M}=(0,u_{M},\vartheta_{M}),\quad u_{M}=G^{u}_{0}(\slashed{p}_{0,M}),\quad \vartheta_{M}=\slashed{G}_{0}(\slashed{p}_{0,M})
\end{align}
 and we have:
\begin{align}
 \frac{\partial g}{\partial(\slashed{p}_{0})_{A}}=0\quad:\textrm{at}\quad \slashed{p}_{0,M}
\end{align}
that is:
\begin{align}
 \frac{\partial u^{\prime}}{\partial u}\frac{\partial G^{u}_{0}}{\partial(\slashed{p}_{0})_{A}}
+\frac{\partial u^{\prime}}{\partial\vartheta^{B}}\frac{\partial\slashed{G}^{B}_{0}}{\partial(\slashed{p}_{0})_{A}}=0\quad:\textrm{at}\quad \slashed{p}_{0,M}
\end{align}
Consider the level surface $S^{\prime}_{0,u^{\prime}}$ of $u^{\prime}$ on $\Sigma_{0}$ through the point $p_{M}$. We shall show that $L^{\prime}$ is $g$-orthogonal
to $S_{0,u^{\prime}}$ at $p_{M}$. 

This fact is equivalent to the following. Denoting by $C^{-}(q_{0})$ the past null geodesic cone of the point $q_{0}$ on the singular
boundary, and by $C^{-}_{out}(q_{0})$ the part of $C^{-}(q_{0})$ corresponding to the $outgoing$ null geodesics ending at $q_{0}$, then the statement is that the 
surface $C^{-}_{out}(q_{0})\bigcap\Sigma_{0}$ and $S^{\prime}_{0,u^{\prime}}$ have the same tangent plane at $p_{M}$:
\begin{align}
 T_{p_{M}}(C^{-}_{out}(q_{0})\bigcap\Sigma_{0})=T_{p_{M}}S^{\prime}_{0,u^{\prime}}
\end{align}
To show this, we shall produce a basis for the left-hand side and show that this is also a basis for the right-hand side.

      Consider the vectorfields:
\begin{align}
 Y^{A}=\frac{\partial F^{t}}{\partial(\slashed{p}_{0})_{A}}\frac{\partial}{\partial t}+
\frac{\partial F^{u}}{\partial(\slashed{p}_{0})_{A}}\frac{\partial}{\partial u}+\frac{\partial\slashed{F}^{B}}{\partial(\slashed{p}_{0})_{A}}
\frac{\partial}{\partial\vartheta^{B}}\quad :\quad A=1,2
\end{align}
These are tangential to the sections of $C^{-}_{out}(q_{0})$ which correspond to constant values of $\tau$ (recalling the mapping $F$, 
we view $\tau$ and $\slashed{p}_{0}$ are independent coordinates), hence are $g$-orthogonal to $\tilde{L}^{\prime}$. For each point $p\in C^{-}_{out}(q_{0})$,
$(Y^{1}(p), Y^{2}(p))$ is a basis for the tangent plane at $p$ to the $\tau=const$ section of $C^{-}_{out}(q_{0})$ through $p$, provided that the vectors $Y^{1}(p)$
, $Y^{2}(p)$ are linearly independent.
We project the $Y^{A}$ to the $\Sigma_{t}$ as follows. We find suitable functions $\lambda^{A}$, $A=1,2$, such that the vectorfields:
\begin{align}
 Y^{\prime A}=Y^{A}-\lambda^{A}L^{\prime}
\end{align}
are tangential to the $\Sigma_{t}$. The vectorfields $Y^{\prime A}$, being $g$-orthogonal to $L^{\prime}$ hence tangential to $C^{-}_{out}(q_{0})$, shall then be
tangential to the sections $C^{-}_{out}(q_{0})\bigcap\Sigma_{t}$. As $L^{\prime}$ is given by:
\begin{align}
 L^{\prime}=\frac{\partial}{\partial t}+\frac{\partial F^{u}/\partial\tau}{\partial F^{t}/\partial\tau}\frac{\partial}{\partial u}
+\frac{\partial\slashed{F}^{A}/\partial\tau}{\partial F^{t}/\partial\tau}\frac{\partial}{\partial\vartheta^{A}}
\end{align}
we find:
\begin{align}
 \lambda^{A}=\frac{\partial F^{t}}{\partial(\slashed{p}_{0})_{A}}
\end{align}
and:
\begin{align}
 \frac{\partial F^{t}}{\partial\tau}Y^{\prime A}=(\frac{\partial F^{u}}{\partial(\slashed{p}_{0})_{A}}\frac{\partial F^{t}}{\partial\tau}
-\frac{\partial F^{t}}{\partial(\slashed{p}_{0})_{A}}\frac{\partial F^{u}}{\partial\tau})\frac{\partial}{\partial u}\\\notag
+(\frac{\partial\slashed{F}^{B}}{\partial(\slashed{p}_{0})_{A}}\frac{\partial F^{t}}{\partial\tau}-\frac{\partial F^{t}}{\partial(\slashed{p}_{0})_{A}}
\frac{\partial\slashed{F}^{B}}{\partial\tau})\frac{\partial}{\partial\vartheta^{B}}
\end{align}
Then $(Y^{\prime 1}(p_{M}), Y^{\prime 2}(p_{M}))$ is a basis for $T_{p_{M}}(C^{-}_{out}(q_{0})\bigcap\Sigma_{0})$, provided that $p_{M}$ is not conjugate to $q_{0}$.
Now, we have, from (19.223):
\begin{align}
 \frac{\partial G^{u}_{0}}{\partial(\slashed{p}_{0})_{A}}=\frac{\partial F^{u}}{\partial\tau}\frac{\partial f}{\partial(\slashed{p}_{0})_{A}}
+\frac{\partial F^{u}}{\partial(\slashed{p}_{0})_{A}}\\\notag
\frac{\partial\slashed{G}^{B}_{0}}{\partial(\slashed{p}_{0})_{A}}=\frac{\partial\slashed{F}^{B}}{\partial\tau}\frac{\partial f}{\partial(\slashed{p}_{0})_{A}}
+\frac{\partial\slashed{F}^{B}}{\partial(\slashed{p}_{0})_{A}}
\end{align}
Also, differentiating the equation $F^{t}(f(\slashed{p}_{0});\slashed{p}_{0})=0$ implicitly with respect to $\slashed{p}_{0}$ yields:
\begin{align}
 \frac{\partial f}{\partial (\slashed{p}_{0})_{A}}=-\frac{\partial F^{t}/\partial(\slashed{p}_{0})_{A}}{\partial F^{t}/\partial\tau}
\end{align}
Substituting in (19.234) we then obtain:
\begin{align}
 \frac{\partial G^{u}_{0}}{\partial(\slashed{p}_{0})_{A}}=(-\frac{\partial F^{u}}{\partial\tau}\frac{\partial F^{t}}{\partial(\slashed{p}_{0})_{A}}
+\frac{\partial F^{t}}{\partial\tau}\frac{\partial F^{u}}{\partial (\slashed{p}_{0})_{A}})\frac{1}{\partial F^{t}/\partial\tau}\\\notag
\frac{\partial\slashed{G}^{B}_{0}}{\partial(\slashed{p}_{0})_{A}}=(-\frac{\partial\slashed{F}^{B}}{\partial\tau}\frac{\partial F^{t}}{\partial(\slashed{p}_{0})_{A}}
+\frac{\partial F^{t}}{\partial\tau}\frac{\partial\slashed{F}^{B}}{\partial(\slashed{p}_{0})_{A}})\frac{1}{\partial F^{t}/\partial\tau}
\end{align}
Comparing with (19.233) we see that:
\begin{align}
 Y^{\prime A}=\frac{\partial G^{u}_{0}}{\partial (\slashed{p}_{0})_{A}}\frac{\partial}{\partial u}+
\frac{\partial\slashed{G}^{B}_{0}}{\partial(\slashed{p}_{0})_{A}}\frac{\partial}{\partial\vartheta^{B}}
\end{align}
hence (19.227) reads, simply:
\begin{align}
 Y^{\prime A}(p_{M})u^{\prime}=0\quad: A=1,2
\end{align}
We conclude that $(Y^{\prime 1}(p_{M}), Y^{\prime 2}(p_{M}))$ is a basis for $T_{p_{M}}S^{\prime}_{0,u^{\prime}}$, provided that the vectors $Y^{\prime 1}(p_{M})$,
$Y^{\prime 2}(p_{M})$ are linearly independent. We have thus demonstrated (19.228) under the condition that $p_{M}$ is not conjugate to $q_{0}$. However, (19.226)
holds even in the case that $p_{M}$ is conjugate to $q_{0}$, but in this case we must approach $p_{M}$, which corresponds to $\tau=f(\slashed{p}_{0,M})$ along 
the outgoing geodesic ending at $q_{0}$ with angular momentum $\slashed{p}_{0,M}$, by a sequence of non-conjugate points corresponding to $\tau=\tau_{n} \longrightarrow
f(\slashed{p}_{0,M})$ along the same geodesic.

    The equality (19.228) is equivalent to $L^{\prime}$ being $g$-orthogonal to $S^{\prime}_{0,u^{\prime}}$ at $p_{M}$. Thus, extending $u^{\prime}$ to an 
acoustical function so that its level sets are outgoing characteristic hypersurface $C^{\prime}_{u^{\prime}}$ with 
$C^{\prime}_{u^{\prime}}\bigcap\Sigma_{0}=S^{\prime}_{0,u^{\prime}}$, the generator of the outgoing characteristic hypersurfaces $C^{\prime}_{u^{\prime}}$ 
through $p_{M}$ coincides with the $outgoing$ null geodesic ending at $q_{0}$ with angular momentum 
$\slashed{p}_{0,M}$.

     Let us now consider again our main problem from a more general point of view. Our problem is the initial value problem for the equations of motion of a perfect
fluid, with initial data on the hyperplane $\Sigma_{0}$ which coincide with those of a constant state outside a compact subset $P_{0}\subset \Sigma_{0}$. We may take
$S_{0,0}$ to be any sphere in $\Sigma_{0}$ containing $P_{0}$. The construction in Chapter 2 then gives an acoustical function $u$ based on $S_{0,0}$. If the initial
data is isentropic and irrotational, Theorem 17.1 can be applied after suitable translation and rescaling, giving a development of the initial data in a domain
$V_{\epsilon_{0}}$ as in (19.7) corresponding to the annular region on $\Sigma_{0}$ bounded by $S_{0,0}$ and the concentric sphere of radius $1-\epsilon_{0}$ that 
of $S_{0,0}$. Then by the local uniqueness theorem the union of all the developments each of which corresponds to a sphere $S_{0,0}$ containing $P_{0}$ 
is also a development of the initial data, contained in the maximal development, the union of the domains $V_{\epsilon_{0}}$ being contained in the closure of the domain 
of the maximal solution. Consider now a point $q_{0}$ belonging to the singular boundary of the domain of the maximal solution such that $q_{0}\in V^{\prime}_{\epsilon_{0}}$
where $V^{\prime}_{\epsilon_{0}}$ is the domain associated to the acoustical function $u^{\prime}$, corresponding to the sphere $S^{\prime}_{0,0}$. Thus $q_{0}\in
\Sigma_{t^{\prime}_{*\epsilon}}^{\epsilon}$ for some $\epsilon\in(0,\epsilon_{0})$, where $t^{\prime}_{*\epsilon}=\inf_{u^{\prime}\in[0,\epsilon]}t^{\prime}_{*}(u^{\prime})$
and $t^{\prime}_{*}(u^{\prime})$ is the least upper bound of the extent of the generators of $C^{\prime}_{u^{\prime}}$ in the domain of maximal solution. Then for any given 
acoustical function $u$, corresponding to a given sphere $S_{0,0}$, the above argument (with the roles of $u$ and $u^{\prime}$ reversed) shows that there is an outgoing 
characteristic hypersurface $C_{u}$ associated to $u$ and a generator of $C_{u}$ reaching $q_{0}$, thus of extent $t^{\prime}_{*\epsilon}$, contained in 
$V^{\prime}_{\epsilon_{0}}$, hence in the union of the domains $V_{\epsilon_{0}}$. In view of the properties of the transformation from one acoustical function to another,
expounded by Proposition 19.2, Theorem 17.1 covers the union of the domains $V_{\epsilon_{0}}$.

      However we can extend the domain covered by Theorem 17.1 still further in the following manner. We can replace the hyperplanes $\Sigma_{t}$, by the hyperplanes
$\Sigma_{\alpha;t}$ where for any fixed $\alpha=(\alpha_{1},\alpha_{2},\alpha_{3})\in\mathbb{R}^{3}$ such that
\begin{align*}
 \|\alpha\|=\sqrt{\alpha_{1}^{2}+\alpha_{2}^{2}+\alpha_{3}^{2}}<1
\end{align*}
 and any $t\in\mathbb{R}$, $\Sigma_{\alpha;t}$ is the hyperplane in Galileo spacetime given by:
\begin{align*}
 x^{0}=t+\sum_{i=1}^{3}\alpha_{i}x^{i}
\end{align*}
 The hyperplane $\Sigma_{\alpha;0}$ is then space-like relative to the actual acoustical metric $g$ provided that the departure of the data induced on $\Sigma_{\alpha;0}$ 
from those of the constant state is suitably small 
. We then consider as above a sphere $S_{0,0}$ 
in $\Sigma_{0}$ such that the initial data on $\Sigma_{0}$ coincide with those of the constant state outside $S_{0,0}$, and proceed to construct the corresponding
acoustical function $u$. The surfaces $S_{\alpha;t,0}=C_{0}\bigcap\Sigma_{\alpha;t}$ are then round spheres. The first variations are defined exactly as before, for, these are adapted to the rest frame of the constant 
state. Of the commutation fields $T$, $Q$, $R_{i}$ : $i=1,2,3$, which define the higher order variations, $Q$ is defined exactly as before, $T$ is now defined to be
tangential to the hyperplanes $\Sigma_{\alpha;t}$, $g$-orthogonal to the surfaces $S_{\alpha;t,u}=C_{u}\bigcap\Sigma_{\alpha;t}$ and satisfying $Tu=1$, while the 
$R_{i}$ : $i=1,2,3$ are defined to be the $g$-orthogonal projections to the surfaces $S_{\alpha;t,u}$, of the rotations associated to $g_{0}$ which leave 
$\Sigma_{\alpha;t}$ and $S_{\alpha;t,0}$ invariant. We note that, similarly, $\Sigma_{\alpha;t}$ remain space-like relative to the acoustical metric provided that
the initial departure from the constant state is suitably small. We obtain in this way an extension of Theorem 17.1 and a domain $V_{\alpha;\epsilon_{0}}$ defined 
in a similar manner as the domain $V_{\epsilon_{0}}$, but with the hyperplanes $\Sigma_{\alpha;t}$ in the role of the hyperplanes $\Sigma_{t}$. The domain covered 
by Theorem 17.1 is then maximally extended to be the union of all the domains $V_{\alpha;\epsilon_{0}}$, and the remarks of the previous paragragh apply to this
maximal extension.

\section{How $\mathcal{H}$ Looks Like in Rectangular Coordinates\\
in Galileo Spacetime}
     We now consider a given component of $\mathcal{H}$ and its past boundary, the corresponding component of $\partial_{-}\mathcal{H}$. By Proposition 19.1 these are 
smooth 3-dimensional and 2-dimensional submanifolds of $\mathcal{M}_{\epsilon_{0}}$, respectively. We shall now investigate the image of the component of $\mathcal{H}$
as well as $\partial_{-}\mathcal{H}$ in Galileo spacetime. For clarity we distinguish these images from the sets themselves and we denote by $H$ and $\partial_{-}H$
the image in Galileo spacetime of $\mathcal{H}$ and $\partial_{-}\mathcal{H}$ respectively.

     Now by Proposition 19.1 $\partial_{-}\mathcal{H}$ is a space-like surface in $(\mathcal{M}_{\epsilon_{0}},g)$ intersecting the invariant curves transversally. 
Thus there is a smooth function $u_{*}$ defined on a connected domain $\mathcal{B}\subset S^{2}$, such that the component of $\partial_{-}\mathcal{H}$ is given in 
canonical acoustical coordinates by:
\begin{align}
 \{(t_{*}(u_{*}(\vartheta),\vartheta), u_{*}(\vartheta), \vartheta)\quad:\quad \vartheta\in\mathcal{B}\}
\end{align}
 and the corresponding component of $\mathcal{H}$ is given by:
\begin{align}
 \{(t_{*}(u,\vartheta),u,\vartheta)\quad:\quad u_{*}>u>\underline{u}(\vartheta),\quad\vartheta\in\mathcal{B}\} 
\end{align}
where a given invariant curve, parametrized by $u$, either extends to any $u>0$, in which case $\underline{u}(\vartheta)=0$, or ends at $u=\underline{u}(\vartheta)$,
the point of intersection with the incoming characteristic hypersurface $\underline{C}$ generated by the incoming null normals to the next outward component of 
$\partial_{-}\mathcal{H}$.

     Consider then the component of $\partial_{-}H$, which is the image of (19.239) in Galileo spacetime. This is the image in Galileo spacetime of the domain
$\mathcal{B}\subset S^{2}$ under the mapping:
\begin{align}
 \vartheta\longmapsto x=q(\vartheta),\quad\textrm{where}\quad q^{\mu}(\vartheta)=x^{\mu}(t_{*}(u_{*}(\vartheta),\vartheta),u_{*}(\vartheta),\vartheta)
\quad:\quad \mu=0,1,2,3
\end{align}
This mapping is one-to-one (it is simply the transformation between rectangular coordinates and acoustical coordinates). We shall presently show that the mapping is of maximal rank. Consider the 2-dimensional matrix $m$ with entries:
\begin{align}
 m_{AB}=g_{\mu\nu}X^{\mu}_{A}\frac{\partial q^{\nu}}{\partial\vartheta^{B}}
\end{align}
Now, we have (see (19.115)):
\begin{align}
 \frac{\partial x^{0}}{\partial t}=1,\quad\frac{\partial x^{i}}{\partial t}=L^{i}\\\notag
\frac{\partial x^{0}}{\partial u}=\frac{\partial x^{0}}{\partial\vartheta^{A}}=0,\quad\frac{\partial x^{i}}{\partial u}
=\mu(\alpha^{-1}\hat{T}^{i}+\hat{\Xi}^{A}X^{i}_{A}),\quad\frac{\partial x^{i}}{\partial\vartheta^{A}}=X^{i}_{A}
\end{align}
and from (19.241):
\begin{align}
 \frac{\partial q^{0}}{\partial\vartheta^{B}}=\frac{\partial t_{*}}{\partial u}\frac{\partial u_{*}}{\partial\vartheta^{B}}
+\frac{\partial t_{*}}{\partial\vartheta^{B}}\\\notag
\frac{\partial q^{j}}{\partial\vartheta^{B}}=\frac{\partial x^{j}}{\partial t}(\frac{\partial t_{*}}{\partial u}\frac{\partial u_{*}}{\partial\vartheta^{B}}
+\frac{\partial t_{*}}{\partial\vartheta^{B}})+\frac{\partial x^{j}}{\partial u}\frac{\partial u_{*}}{\partial\vartheta^{B}}+\frac{\partial x^{j}}{\partial\vartheta^{B}}
\end{align}
Now on $\mathcal{H}\bigcup\partial_{-}\mathcal{H}$ we have $\mu=0$, therefore, from (19.243):
\begin{align}
 \frac{\partial x^{i}}{\partial u}=0\quad:\quad\textrm{on}\quad\mathcal{H}\bigcup\partial_{-}\mathcal{H}
\end{align}
Moreover by Proposition 19.1:
\begin{align}
 \frac{\partial t_{*}}{\partial u}=0\quad:\quad\textrm{on}\quad\partial_{-}\mathcal{H}
\end{align}
Thus substituting in (19.244) from (19.245) and taking account of (19.245) and (19.246), we obtain:
\begin{align*}
 \frac{\partial q^{0}}{\partial\vartheta^{B}}=\frac{\partial t_{*}}{\partial\vartheta^{B}},\quad
\frac{\partial q^{j}}{\partial\vartheta^{B}}=L^{j}\frac{\partial t_{*}}{\partial\vartheta^{B}}+X^{j}_{B}
\end{align*}
or, simply:
\begin{align}
 \frac{\partial q^{\nu}}{\partial\vartheta^{B}}=L^{\nu}\frac{\partial t_{*}}{\partial\vartheta^{B}}+X^{\nu}_{B}
\end{align}
Substituting this in (19.242) we obtain simply:
\begin{align}
 m_{AB}=g(X_{A},X_{B})=\slashed{g}_{AB}
\end{align}
Hence $\det m=\det\slashed{g}>0$ and the mapping (19.241) is of maximal rank. It follows that this mapping is a diffeomorphism onto its image and its image, the component
of $\partial_{-}H$, is a smooth embedded 2-dimensional submanifold in Galileo spacetime.

     The component of $H$ can similarly be shown to be a smooth embedded 3-dimensional submanifold in Galileo spacetime except along its past boundary (just using the
definition of linearly independent, the three linearly independent vectors are $L\frac{\partial t_{*}}{\partial u}$ and the $L\frac{\partial t_{*}}{\partial\vartheta^{A}}+X_{A}$), the corresponding component of $\partial_{-}H$. To
analyze the behavior near $\partial_{-}H$, we shall derive an analytic expression for an arbitrary invariant curve of $H$ in a neighborhood of 
$\partial_{-}H$. Now the parametric equations of an invariant curve of $H$, in rectangular coordinates in Galileo spacetime, the parameter being $u$,
are (see (19.240)):
\begin{align}
 x^{\mu}=x^{\mu}(t_{*}(u,\vartheta),u,\vartheta)\quad:\quad u_{*}(\vartheta)>u>\underline{u}(\vartheta),\quad\vartheta=const.
\end{align}
or, since $x^{0}=t$,
\begin{align}
 x^{0}=t_{*}(u,\vartheta)\\\notag
x^{i}=x^{i}(t_{*}(u,\vartheta),u,\vartheta)\quad:\quad u_{*}(\vartheta)>u>\underline{u}(\vartheta),\quad\vartheta=const.
\end{align}
We have:
\begin{align*}
 (\frac{dx^{0}}{du})(u)=(\frac{\partial t_{*}}{\partial u})(u,\vartheta)\\\notag
(\frac{dx^{i}}{du})(u)=(\frac{\partial x^{i}}{\partial t})(t_{*}(u,\vartheta),u,\vartheta)(\frac{\partial t_{*}}{\partial u})(u,\vartheta)
+(\frac{\partial x^{i}}{\partial u})(t_{*}(u,\vartheta),u,\vartheta)
\end{align*}
Substituting from (19.243) and taking into account of (19.245) we obtain:
\begin{align}
 (\frac{dx^{0}}{du})(u)=(\frac{\partial t_{*}}{\partial u})(u,\vartheta)\\\notag
(\frac{dx^{i}}{du})(u)=L^{i}(t_{*}(u,\vartheta),u,\vartheta)(\frac{\partial t_{*}}{\partial u})(u,\vartheta)
\end{align}
or simply:
\begin{align}
 (\frac{dx^{\mu}}{du})(u)=L^{\mu}(t_{*}(u,\vartheta),u,\vartheta)(\frac{\partial t_{*}}{\partial u})(u,\vartheta)
\end{align}
We see again that in Galileo spacetime the vectorfield $L$ is tangential to the invariant curves. We also see that at $u=u_{*}(\vartheta)$, that is
on $\partial_{-}H$, the right-hand sides vanish by (19.246), so the past end point of an invariant curve is a singular point.

     To analyze the behavior at the past end point we consider the second and third derivatives at this point. Differentiating (19.251) with respect to $u$ we 
obtain:
\begin{align}
(\frac{d^{2}x^{0}}{du^{2}})(u)=(\frac{\partial^{2}t_{*}}{\partial u^{2}})(u,\vartheta)\\\notag
(\frac{d^{2}x^{i}}{du^{2}})(u)=(\frac{\partial L^{i}}{\partial t})(t_{*}(u,\vartheta),u,\vartheta)(\frac{\partial t_{*}}{\partial u})^{2}(u,\vartheta)\\\notag
+(\frac{\partial L^{i}}{\partial u})(t_{*}(u,\vartheta),u,\vartheta)(\frac{\partial t_{*}}{\partial u})(u,\vartheta)\\\notag
+L^{i}(t_{*}(u,\vartheta),u,\vartheta)(\frac{\partial^{2}t_{*}}{\partial u^{2}})(u,\vartheta)
\end{align}
Also, differentiating (19.48) with respect to $u$ we obtain:
\begin{align*}
 (\frac{\partial^{2}t_{*}}{\partial u^{2}})(u,\vartheta)=-(\frac{\partial\mu}{\partial t})^{-2}(t_{*}(u,\vartheta),u,\vartheta)\\
\cdot\{\frac{\partial\mu}{\partial t}(\frac{\partial^{2}\mu}{\partial u\partial t}\frac{\partial t_{*}}{\partial u}+\frac{\partial^{2}\mu}{\partial u^{2}})
-\frac{\partial\mu}{\partial u}(\frac{\partial^{2}\mu}{\partial t^{2}}\frac{\partial t_{*}}{\partial u}+\frac{\partial^{2}\mu}{\partial u\partial t})\}
\end{align*}
Substituting for $\partial t_{*}/\partial u$ from (19.48) then yields:
\begin{align}
 (\frac{\partial^{2}t_{*}}{\partial u^{2}})(u,\vartheta)\\\notag
=-(\frac{\partial\mu}{\partial t})^{-3}\{(\frac{\partial\mu}{\partial t})^{2}\frac{\partial^{2}\mu}{\partial u^{2}}
-2\frac{\partial\mu}{\partial t}\frac{\partial\mu}{\partial u}\frac{\partial^{2}\mu}{\partial u\partial t}+
(\frac{\partial\mu}{\partial u})^{2}\frac{\partial^{2}\mu}{\partial t^{2}}\}(t_{*}(u,\vartheta),u,\vartheta)
\end{align}
In particular, at $u=u_{*}(\vartheta)$ we obtain, in accordance with (19.51):
\begin{align}
 (\frac{\partial^{2}t_{*}}{\partial u^{2}})(u_{*}(\vartheta),\vartheta)=a(\vartheta)
\end{align}
where:
\begin{align}
 a(\vartheta)=-(\frac{\partial^{2}\mu/\partial u^{2}}{\partial\mu/\partial t})(q(\vartheta))>0
\end{align}
by the first equation of (19.45), and we denote by
\begin{align}
 q(\vartheta)=(t_{*}(u_{*}(\vartheta),\vartheta),u_{*}(\vartheta),\vartheta)
\end{align}
the end point on $\partial_{-}H$. By (19.246) and (19.255), at $u=u_{*}(\vartheta)$ (19.253) becomes:
\begin{align}
 (\frac{d^{2}x^{0}}{du^{2}})(u_{*}(\vartheta))=a(\vartheta)\\\notag
(\frac{d^{2}x^{i}}{du^{2}})(u_{*}(\vartheta))=L^{i}(q(\vartheta))a(\vartheta)
\end{align}
Moreover, differentiating (19.253) with respect to $u$ and evaluating the result at $u=u_{*}(\vartheta)$ we obtain, in view of (19.246):
\begin{align}
 (\frac{d^{3}x^{0}}{du^{3}})(u_{*}(\vartheta))=(\frac{\partial^{3}t_{*}}{\partial u^{3}})(u_{*}(\vartheta),\vartheta)\\\notag
(\frac{d^{3}x^{i}}{du^{3}})(u_{*}(\vartheta))=2(\frac{\partial L^{i}}{\partial u})(q(\vartheta))(\frac{\partial^{2}t_{*}}{\partial u^{2}})
(u_{*}(\vartheta),\vartheta)+L^{i}(q(\vartheta))(\frac{\partial^{3}t_{*}}{\partial u^{3}})(u_{*}(\vartheta),\vartheta)
\end{align}
We shall presently determine the coefficients:
\begin{align*}
 (\frac{\partial^{3}t_{*}}{\partial u^{3}})(u_{*}(\vartheta),\vartheta)\quad\textrm{and}\quad\frac{\partial L^{i}}{\partial u}(q(\vartheta))
=TL^{i}(q(\vartheta))
\end{align*}
In view of the fact that for any smooth function $f(t,u,\vartheta)$ we have:
\begin{align*}
 \frac{d}{du}f(t_{*}(u,\vartheta),u,\vartheta)=\frac{\partial f}{\partial t}\frac{\partial t_{*}}{\partial u}+\frac{\partial f}{\partial u}
\end{align*}
hence by (19.246):
\begin{align*}
 (\frac{d}{du}f(t_{*}(u,\vartheta),u,\vartheta))_{u=u_{*}(\vartheta)}=\frac{\partial f}{\partial u}(q(\vartheta))
\end{align*}
differentiating (19.254) with respect to $u$ and evaluating the result at $u=u_{*}(\vartheta)$ yields:
\begin{align}
 (\frac{\partial^{3}t_{*}}{\partial u^{3}})(u_{*}(\vartheta),\vartheta)=b(\vartheta)
\end{align}
where:
\begin{align}
 b(\vartheta)=(\frac{\partial\mu}{\partial t})^{-2}(3\frac{\partial^{2}\mu}{\partial u^{2}}\frac{\partial^{2}\mu}{\partial t\partial u}
-\frac{\partial\mu}{\partial t}\frac{\partial^{3}\mu}{\partial u^{3}})(q(\vartheta))
\end{align}
Next, $TL^{i}$ is given by (3.185):
\begin{align}
 TL^{i}=a_{T}\hat{T}^{i}+b^{i}_{T}
\end{align}
where $a_{T}$ is given by (3.188) and $b^{i}_{T}$ are the rectangular components of $S_{t,u}$-tangential vectorfield given by (3.190). 
On $H\bigcup\partial_{-}H$, in particular at $q(\vartheta)$, these formulas reduce to:
\begin{align}
 a_{T}=\alpha^{-1}L\mu<0\\\notag
b^{i}_{T}=(\slashed{g}^{-1}\cdot\slashed{d}\mu)^{i}
\end{align}
Substituting the above as well as (19.255) in (19.259) we obtain:
\begin{align}
 (\frac{d^{3}x^{0}}{du^{3}})(u_{*}(\vartheta))=b(\vartheta)\\\notag
(\frac{d^{3}x^{i}}{du^{3}})(u_{*}(\vartheta))=L^{i}(q(\vartheta))b(\vartheta)+2a(\vartheta)
(a_{T}\hat{T}^{i}+b^{i}_{T})(q(\vartheta))
\end{align}
In view of (19.258) and (19.264) we obtain the following analytic expression for an invariant curve of $H$ in a neighborhood of its origin $q(\vartheta)
\in\partial_{-}H$:
\begin{align*}
 x^{0}(u)=x^{0}(q(\vartheta))+\frac{1}{2}a(\vartheta)(u-u_{*}(\vartheta))^{2}+\frac{1}{6}b(\vartheta)(u-u_{*}(\vartheta))^{3}\\
+O((u-u_{*}(\vartheta))^{4})
\end{align*}
\begin{align*}
 x^{i}(u)=x^{i}(q(\vartheta))+\frac{1}{2}a(\vartheta)L^{i}(q(\vartheta))(u-u_{*}(\vartheta))^{2}\\
+\frac{1}{6}[b(\vartheta)L^{i}(q(\vartheta))+2a(\vartheta)(a_{T}\hat{T}^{i}+b^{i}_{T})(q(\vartheta))](u-u_{*}(\vartheta))^{3}\\
+O((u-u_{*}(\vartheta))^{4})
\end{align*}
We summarize the above results in the following proposition.

$\textbf{Proposition 19.3}$ The singular boundary of the domain of the maximal solution in Galileo spacetime is $H\bigcup\partial_{-}H$. Consider a given 
component of $H$ and its past boundary, the corresponding component of $\partial_{-}H$. The component of $\partial_{-}H$ is a smooth 2-dimensional embedded
submanifold in Galileo spacetime, which is space-like with respect to the acoustical metric. The component of $H$ is a smooth embedded 3-dimensional submanifold 
in Galileo spacetime ruled by invariant curves of vanishing arc length with respect to the acoustical metric, having past end points on the corresponding 
component of $\partial_{-}H$. In a neighborhood of its origin at the point $q(\vartheta)\in\partial_{-}H$ an invariant curve is given, in parametric form, by:
\begin{align*}
 x^{0}(u)=x^{0}(q(\vartheta))+\frac{1}{2}a(\vartheta)(u-u_{*}(\vartheta))^{2}+\frac{1}{6}b(\vartheta)(u-u_{*}(\vartheta))^{3}\\
+O((u-u_{*}(\vartheta))^{4})
\end{align*}
\begin{align*}
 x^{i}(u)=x^{i}(q(\vartheta))+\frac{1}{2}a(\vartheta)L^{i}(q(\vartheta))(u-u_{*}(\vartheta))^{2}\\
+\frac{1}{6}[b(\vartheta)L^{i}(q(\vartheta))+2a(\vartheta)(a_{T}\hat{T}^{i}+b^{i}_{T})(q(\vartheta))](u-u_{*}(\vartheta))^{3}\\
+O((u-u_{*}(\vartheta))^{4})
\end{align*}
where:
\begin{align*}
 a(\vartheta)=-(\frac{\partial^{2}\mu/\partial u^{2}}{\partial\mu/\partial t})(q(\vartheta))>0\\
b(\vartheta)=(\frac{\partial\mu}{\partial t})^{-2}(3\frac{\partial^{2}\mu}{\partial u^{2}}\frac{\partial^{2}\mu}{\partial t\partial u}
-\frac{\partial\mu}{\partial t}\frac{\partial^{3}\mu}{\partial u^{3}})(q(\vartheta))\\
a_{T}(q(\vartheta))=(\alpha^{-1}L\mu)(q(\vartheta))<0\\
b_{T}^{i}(q(\vartheta))=(\slashed{g}^{-1}\cdot\slashed{d}\mu)^{i}(q(\vartheta))
\end{align*}

\bibliographystyle{amsplain}
\bibliography{reference}
\end{document}